\documentclass[]{amsart}

\usepackage[Symbol]{upgreek}

\usepackage{graphicx}

\usepackage{amssymb}
\usepackage{epic}
\usepackage{mathrsfs}
\usepackage{accents}
\usepackage{array}
\usepackage{stmaryrd}
\usepackage{enumitem}
\usepackage{longtable}

\usepackage{hyperref}
\usepackage{etoolbox} % for \patchcmd

\usepackage{tikz}

\input{xy}
\xyoption{matrix}
\xyoption{arrow}

\numberwithin{section}{part}
\numberwithin{equation}{section}
\makeatletter
\renewcommand{\tocsection}[3]
 { \indentlabel{\@ifnotempty{#2}{\parbox{2.5em}{\ignorespaces#1 #2.}\quad}}#3}

\makeindex

\newcommand{\bbold}{\mathbb}
\newcommand{\cal}{\mathcal}

\def\R { {\bbold R} }
\def\Q { {\bbold Q} }
\def\Z { {\bbold Z} }
\def\C { {\bbold C} }
\def\N { {\bbold N} }
\def\T { {\bbold T} }
\def\E {{\mathcal E}}
\let\cedille\c
\def\c {\mathcal{C}}
\def\cc {\operatorname{c}}
\def\g {\operatorname{g}}
\def \I{\operatorname{I}}

\def \Ex{\operatorname{E}}
\def \Dx{\operatorname{D}}
\def \Exp{\operatorname{E}}
\def\ta{\operatorname{tail}}
\def\trig{\operatorname{trig}}
\def\tl{\operatorname{tl}}

\def \sol{\operatorname{sol}}
\def \order{\operatorname{order}}
\def \val{\operatorname{mul}}
\newcommand\dval{\operatorname{dmul}}

\def \exc {{\mathscr E}}
\def \ex{\operatorname{e}}
\def \wr {\operatorname{wr}}

\def \Frac {\operatorname{Frac}}

\def \Univ {{\operatorname{U}}}

\renewcommand\epsilon{\varepsilon}

\def \d{\operatorname{d}}

\def \ev{\operatorname{e}}
\def \bar {\overline}
\def \<{\langle}
\def \>{\rangle}

\def \tilde {\widetilde}

\def \mult{\operatorname{mult}}

\def \hat {\widehat}

\def \supp {\operatorname{supp}}

\def \((  {(\!(}
\def \)) {)\!)}

\def \Hom{\operatorname{Hom}}
\def \Li{\operatorname{Li}}
\def \res{\operatorname{res}}

\def \k {{{\boldsymbol{k}}}}
\def \flatter{\mathrel{\prec\!\!\!\prec}}
\DeclareMathSymbol{\precequ}{\mathrel}{symbols}{"16}
\DeclareMathSymbol{\succequ}{\mathrel}{symbols}{"17}
\def \flattereq{\mathrel{\precequ\!\!\!\precequ}}
\def \steeper{\mathrel{\succ\!\!\!\succ}}

\def \comp{\mathrel{-{\hskip0.06em\!\!\!\!\!\asymp}}}

\def \nasymp{\not\asymp}

\renewcommand{\Re}{\operatorname{Re}}
\renewcommand{\Im}{\operatorname{Im}}

\newcommand{\claim}[2][\!\!]{\medskip\noindent {\bf Claim #1:} {\it #2}\medskip}

\newtheorem{theorem}{Theorem}[section]
\newcounter{tmptheorem}
\newtheorem{lemma}[theorem]{Lemma}
\newtheorem{prop}[theorem]{Proposition}
\newtheorem{cor}[theorem]{Corollary}

\newtheorem{corintro}{Corollary}

\theoremstyle{definition}

\newtheorem{definition}[theorem]{Definition}

\theoremstyle{remark}
\newtheorem*{example}{Example}
\newtheorem*{examples}{Examples}
\newtheorem{exampleNumbered}[theorem]{Example}
\newtheorem{examplesNumbered}[theorem]{Examples}
\newtheorem*{notation}{Notation}

\newtheorem*{remarks}{Remarks}
\newtheorem*{remark}{Remark}
\newtheorem*{question}{Question}
\newtheorem*{conjecture}{Conjecture}
\newtheorem{remarkNumbered}[theorem]{Remark}

\newcommand{\abs}[1]{\lvert#1\rvert}
\newcommand{\dabs}[1]{\lVert#1\rVert}

\def \sgn {\operatorname{sign}}

\def \fM {{\mathfrak M}}
\def \fd {{\mathfrak d}}
\def \fm {{\mathfrak m}}
\def \fp {{\mathfrak p}}
\def \fn {{\mathfrak n}}
\def \fv {{\mathfrak v}}
\def \fw {{\mathfrak w}}

\def \Ric{\operatorname{Ri}}

\def \id{\operatorname{id}}

\let\oldi\i
\let\oldj\j
\renewcommand\i{\relax\ifmmode{\boldsymbol{i}}\else\oldi\fi}
\renewcommand\j{\relax\ifmmode{\boldsymbol{j}}\else\oldj\fi}

\renewcommand\leq{\leqslant}
\renewcommand\geq{\geqslant}
\renewcommand\preceq{\preccurlyeq}
\renewcommand\succeq{\succcurlyeq}
\renewcommand\le{\leq}
\renewcommand\ge{\geq}
\renewcommand\frak{\mathfrak}

\DeclareMathAlphabet{\mathbf}{OML}{cmm}{b}{it}

\DeclareFontFamily{U}{fsy}{}
\DeclareFontShape{U}{fsy}{m}{n}{<->s*[.9]psyr}{}
\DeclareSymbolFont{der@m}{U}{fsy}{m}{n}
\DeclareMathSymbol{\der}{\mathord}{der@m}{182}

\DeclareSymbolFont{der@m}{U}{fsy}{m}{n}
\DeclareMathSymbol{\derdelta}{\mathord}{der@m}{100}

\newcommand\dotprec{\mathrel{\dot\prec}}

\newcommand\wt{\operatorname{wt}}

\newcommand\bsigma{\boldsymbol{\sigma}}
\newcommand\m{\mathfrak m}

\newcommand\dwt{\operatorname{dwt}}

\newcommand\nwt{\operatorname{nwt}}
\newcommand\ndeg{\operatorname{ndeg}}
\newcommand\nval{\operatorname{nmul}}

\DeclareSymbolFont{imag@m}{OT1}{cmr}{m}{ui}
\DeclareMathSymbol{\imag}{\mathord}{imag@m}{105}

\DeclareFontFamily{OMS}{smallo}{}
\DeclareFontShape{OMS}{smallo}{m}{n}{<->s*[.65]cmsy10}{}
\DeclareSymbolFont{smallo@m}{OMS}{smallo}{m}{n}
\DeclareMathSymbol{\smallo}{\mathord}{smallo@m}{79}

\DeclareFontFamily{OMS}{largerdot}{}
\DeclareFontShape{OMS}{largerdot}{m}{n}{<->s*[.8]cmsy10}{}
\DeclareSymbolFont{largerdot@m}{OMS}{largerdot}{m}{n}
\DeclareMathSymbol{\largerdot}{\mathord}{largerdot@m}{15}

\newcommand\Aut{\operatorname{Aut}}

\newcommand\End{\operatorname{End}}

\DeclareMathSymbol{\llambda}{\mathord}{der@m}{108}
\DeclareMathSymbol{\rrho}{\mathord}{der@m}{114}

\def \upg{\upgamma}
\def \Upg{\Upgamma}
\def \upl{\uplambda}
\def \Upl{\Uplambda}
\def \upo{\upomega}
\def \Upo{\Upomega}

\def \Upd{\Updelta}
\def \w{\uptau}

\def\HLO{\Upl\Upo}

\newcommand{\equationqed}[1]{\[\pushQED{\qed}#1 \qedhere\popQED\]\let\qed\relax}
\newcommand{\alignqed}[1]{\begin{align*}\pushQED{\qed} #1 \qedhere\popQED\end{align*}\let\qed\relax}

\makeatletter
\newcommand{\dminus}{\mathbin{\text{\@dminus}}}

\newcommand{\@dminus}{%
  \ooalign{\hidewidth\raise1ex\hbox{\bf.}\hidewidth\cr$\m@th-$\cr}%
}
\makeatother

\def\ddeg{\operatorname{ddeg}}
\def\dwm{\operatorname{dwm}}

\def \Car{\mathcal{C}^r_a}
\def \Caz{\mathcal{C}^0_a}
\def \Cao{\mathcal{C}^1_a}
\def \Cat{\mathcal{C}^2_a}
\def \Cainf{\mathcal{C}^{\infty}_a}
\def \Caom{\mathcal{C}^{\omega}_a}

\def \Gr{\mathcal{C}^r}
\def \Gn{\mathcal{C}^n}
\def \Gz{\mathcal{C}^0}
\def \Go{\mathcal{C}^1}
\def \Gt{\mathcal{C}^2}
\def \Gi{\mathcal{C}^{<\infty}}
\def \Ginf{\mathcal{C}^{\infty}}
\def \Gom{\mathcal{C}^{\omega}}
\def \inv{\operatorname{inv}}
\def \Sol{\operatorname{Sol}}

\def \Caj{\mathcal{C}^j_a}
\def \Cajl{\mathcal{C}^{j-1}_a}

\def \Car{\mathcal{C}^r_a}
\def \Carl{\mathcal{C}^{r-1}_a}
\def \Carm{\mathcal{C}^{r+1}_a}
\def \Carmi{\mathcal{C}^{r-i}_a}
\def \Caz{\mathcal{C}^0_a}
\def \Cazr{\mathcal{C}^r_{a_0}}
\def \Cazrl{\mathcal{C}^{r-1}_{a_0}}

\def \Coo{\mathcal{C}^1_0}
\def \Cao{\mathcal{C}^1_a}
\def \Cat{\mathcal{C}^2_a}
\def \Cainf{\mathcal{C}^{\infty}_a}

\def \Can{\mathcal{C}^n_a}
\def \Caom{\mathcal{C}^{\omega}_a}
\def \Calinf{\mathcal{C}^{<\infty}}
\def \Calr{\mathcal{C}^{r}}
\def \Caln{\mathcal{C}^{n}}

\def\b{\operatorname{b}}
\def \inte{\operatorname{int}}
\def \cf{\operatorname{cf}}
\def \O{\mathcal{O}}

\makeatletter
\renewcommand\part{\@startsection{part}{0}%
  \z@{\linespacing\@plus\linespacing}{.5\linespacing}%
  {\normalfont\bfseries\centering}}
\makeatother

\makeatletter
\renewcommand\theindex{\@restonecoltrue\if@twocolumn\@restonecolfalse\fi
  \columnseprule\z@ \columnsep 35\p@
  \twocolumn[\@xp\part\@xp*\@xp{\bf Index}\bigskip]%
  \let\item\@idxitem
  \parindent\z@  \parskip\z@\@plus.3\p@\relax
  \small}  

\renewenvironment{thebibliography}[1]{%
  \@xp\part\@xp*\@xp{\refname}%
  \normalfont\small\labelsep .5em\relax
  \renewcommand\theenumiv{\arabic{enumiv}}\let\p@enumiv\@empty
  \list{\@biblabel{\theenumiv}}{\settowidth\labelwidth{\@biblabel{#1}}%
    \leftmargin\labelwidth \advance\leftmargin\labelsep
    \usecounter{enumiv}}%
  \sloppy \clubpenalty\@M \widowpenalty\clubpenalty
  \sfcode`\.=\@m
}{%
  \def\@noitemerr{\@latex@warning{Empty `thebibliography' environment}}%
  \endlist
}
 
\def\astr{$\,({}^{*})$}
 
\makeatletter
\newcommand{\smallbullet}{} % for safety
\DeclareRobustCommand\smallbullet{%
  \mathord{\mathpalette\smallbullet@{0.6}}%
}
\newcommand{\smallbullet@}[2]{%
  \vcenter{\hbox{\scalebox{#2}{$\m@th#1\bullet$}}}%
}
\makeatother
  
\dedicatory{To the memory of Michael Boshernitzan \textup{(}1950--2019\textup{)}}
   
\begin{document}

\title{Maximal Hardy fields}
\author[Aschenbrenner]{Matthias Aschenbrenner}
\address{Kurt G\"odel Research Center for Mathematical Logic\\
Universit\"at Wien\\
1090 Wien\\ Austria}
\email{matthias.aschenbrenner@univie.ac.at}

\author[van den Dries]{Lou van den Dries}
\address{Department of Mathematics\\
University of Illinois at Urbana-Cham\-paign\\
Urbana, IL 61801\\
U.S.A.}
\email{vddries@math.uiuc.edu}

\author[van der Hoeven]{Joris van der Hoeven}
\address{CNRS, LIX (UMR 7161)\\ 
Campus de l'\'Ecole Polytechnique\\  91120 Palaiseau \\ France}
\email{vdhoeven@lix.polytechnique.fr}

\date{September, 2023}

\begin{abstract} We show that  all maximal Hardy fields are elementarily equivalent as differential fields, and give various applications of this result and its proof. We also answer some questions on Hardy fields posed by Boshernitzan. 
%\medskip
%\noindent
%{\it This document is not intended for publication in its current form. It is not
%a preprint; please do not cite it as such!}
\end{abstract}

\pagestyle{plain}
 
\maketitle

\bigskip

\tableofcontents

\newpage

\begingroup
\setcounter{tmptheorem}{\value{theorem}}% store current value of theorem counter
\setcounter{theorem}{0} %assign desired value to theorem counter
\renewcommand\thetheorem{\Alph{theorem}}

\part*{Preface}

\medskip

\noindent
A Hardy field is said to be {\it maximal}\/ if it has no proper Hardy field extension.  In these notes we show that all  maximal Hardy fields are elementarily equivalent, as ordered differential fields, to the ordered differential field $\T$ of transseries. This is part of our main result, Theorem~\ref{thm:char d-max}.

\medskip
\noindent
We shall depend heavily on our book [ADH], which contains a model-theoretic analysis of $\T$. Besides developing further the asymptotic differential algebra from that book we require also a good dose of analysis. These notes are divided in Parts~\ref{part:preliminaries}--\ref{part:applications}, preceded by a somewhat lengthy Introduction including a sketch of the proof of our main result. Parts~\ref{part:preliminaries}--\ref{part:dents in H-fields} consist of further asymptotic differential algebra and culminates in various normalization theorems for algebraic differential equations over suitable $H$-fields. Parts~\ref{part:Hardy fields univ exp ext} and~\ref{part:Hardy fields} are more analytic and apply the normalization theorems to Hardy fields. Part~\ref{part:applications} consists of applications. We finish with an index and a list of symbols newly introduced in this work.  (All other notation is standard
or comes from [ADH].)

\medskip
\noindent
The present notes are probably not suitable for publication as a journal article, since we took the liberty of including extensive sections with complete proofs on classical topics such as self-adjoint linear differential operators, almost periodic functions,   uniform distribution modulo~$1$, and Bessel functions. This was partly done for our own education, and partly to put things in a form convenient for our purpose. We also took the opportunity to develop some topics a bit further than needed for the main theorem, and in this way we could also answer in Part~5 some questions about Hardy fields raised by Boshernitzan. We have in mind further use of the material here, for example in \cite{ADHld} and in relation to open problems posed in [ADH]. (Our main theorem solves one of those problems.) 

\medskip
\noindent
Readers only interested in  the proof of our main result can skip Sections~\ref{sec:it log derivative}, \ref{sec:self-adjoint}, \ref{sec:upper lower bds}, as well as several subsections of other sections in Parts~\ref{part:preliminaries}--\ref{part:Hardy fields}. %which are not   required towards that goal.
These (sub)sections are marked by an asterisk\astr.

\medskip
\noindent
The main results in these notes are really about {\em differentially algebraic\/}  Hardy field extensions, especially their construction. We complement this in \cite{fgh} with an account of constructing {\em differentially transcendental\/} Hardy field extensions, leading to the result that all maximal Hardy fields are $\eta_1$ in the sense of Hausdorff: in other words, given any Hardy field $H$ and countable subsets $A < B$ in $H$, there is an element $f$ in a Hardy field extension of $H$ such that $A < f < B$. This can be used to show that all maximal Hardy fields are back-and-forth equivalent, which is considerably stronger than their elementary equivalence. We mention this here because the proof of a key ingredient in \cite{fgh} makes essential use of the main result from the present notes.

\medskip
\noindent
We are still deliberating how to publish this material (these notes and \cite{fgh}), but thought it best to make it available for now on the arxiv.

\newpage

\part*{Introduction}

\medskip

\noindent
Du Bois-Reymond's ``orders of infinity'' \cite{dBR71}--\cite{dBR77} were put on a firm
basis by Har\-dy~\cite{Ha}, leading to the notion of a Hardy field (Bourbaki~\cite{Bou}). A {\it Hardy field}\/ is a field~$H$ of germs at $+\infty$ of differentiable real-valued functions
on intervals~$(a,+\infty)$ such that
for any differentiable function whose germ is in~$H$ the germ of its derivative is also in~$H$. (See Section~\ref{sec:Hardy fields} for more precision.) 
Every Hardy field is naturally a differential field, and
an ordered field with the germ of $f$ being $>0$ iff~$f(t)>0$, eventually. 
{\it Hardy fields are the natural domain of asymptotic analysis, where all rules hold, without qualifying conditions}\/ \cite[p.~297]{Rosenlicht83}.
The basic theory of Hardy fields was mostly developed by Boshernitzan~\cite{Boshernitzan81}--\cite{Boshernitzan87} and Rosenlicht~\cite{Rosenlicht83}--\cite{Rosenlicht95}.

\medskip
\noindent
The germs of Hardy's logarithmico-expo\-nen\-tial functions~\cite{Har12a} furnish the classical example of a Hardy field: these
functions are the real-valued functions that can be built from real constants and the identity function $x$, using addition, multiplication,  division,  taking logarithms, and exponentiating. Examples include  the germs of 
the functions~$(0,+\infty)\to\R$ given by
$x^r$ ($r\in\R$), $\ex^{x^2}$, and $\log\log x$.  Other Hardy fields contain (germs of)  differentially transcendental functions, such as the Riemann $\zeta$-function and Euler's $\Gamma$-function~\cite{Rosenlicht83}, and
even functions ultimately growing faster than each iterate of the exponential function~\cite{Boshernitzan86}.  
One source of Hardy fields is o-minimality: every o-minimal structure on the real field naturally gives rise to 
a Hardy field (of germs of definable functions).  
This yields a wealth of  examples such as those obtained from quasi-analytic Denjoy-Carleman classes~\cite{RSW},
or containing certain transition maps of plane analytic vector fields~\cite{KRS},
and explains the role of Hardy fields in model theory and its applications to real analytic geometry and
 dynamical systems~\cite{AvdD4, BMS, Miller}. 
 %For example, they are  a crucial tool for proving Miller's growth dichotomy.
Hardy fields have also   found applications in computer algebra~\cite{SS98, SS99, Shackell}, ergodic theory (see, e.g., \cite{BMR, BoshernitzanWierdl, FW, KoMue}), and other areas of  mathematics~\cite{BB, CS, CMR,  FishmanSimmons, GvdH}. 
%\marginpar{Joris: any other references that we should add?} 
% (Salvy, Shackell, van der Hoeven, \dots).

\medskip
\noindent
In the remainder of this introduction, $H$ is a Hardy field.
Then   $H(\R)$ (obtained by adjoining the germs of the constant functions) is also a Hardy field,
and for any~${h\in H}$, the germ $\ex^h$ generates a Hardy field~$H(\ex^h)$ over $H$, and so does any differentiable germ with derivative
$h$. Moreover,~$H$ has a unique Hardy field extension
that is algebraic over $H$ and real closed.  (See~\cite{Boshernitzan81, Robinson72, Ros} or 
Section~\ref{sec:Hardy fields} below.)  
%In this paper we prove the ultimate extension result of this kind: 
Our main result is Theorem~\ref{thm:char d-max}, and it yields what appears to be the ultimate fact about differentially algebraic Hardy field extensions:

\begin{theorem}\label{thm:A} Let  $P(Y)$ be a differential polynomial in a single
differential indeterminate $Y$ over $H$,  and let $f< g$ in $H$ be such that $P(f) <0 <  P(g)$. Then there is a $y$ in a Hardy field extension of $H$ such that $f < y < g$ and $P(y)=0$.
\end{theorem}

%\medskip
\noindent
By Zorn, every Hardy field extends to a maximal Hardy field, 
so by the theorem above, maximal Hardy fields have the intermediate value property for differential polynomials. 
(In \cite{fgh} we show there are very many maximal Hardy fields, namely $2^{\frak{c}}$ many, where $\frak{c}$ is the cardinality of the continuum.) By the results mentioned earlier, maximal Hardy fields are also Liouville closed $H$-fields in the 
sense of \cite{AvdD2}; thus they contain  the germs of all logarithmico-exponential functions. Hiding behind the
intermediate value property of Theorem~\ref{thm:A} are two more fundamental properties, {\it $\upo$-freeness}\/ and {\it newtonianity,}\/ which are central in our book [ADH]. (Roughly speaking, $\upo$-freeness controls the solvability of second-order homogeneous differential equations, and newtonianity is a strong version of differential-henselianity.) 
%In our proof of Theorem~\ref{thm:A} 
We show that any Hardy field has an $\upo$-free Hardy field extension (Theorem~\ref{upo}), 
and next the much harder result that any $\upo$-free Hardy field
extends to a newtonian $\upo$-free Hardy field: Theorem~\ref{thm:char d-max},
%thm:extend to H-closed}, 
which is really the main result of this paper. It follows that every maximal Hardy field is, in the terminology of \cite{ADH2}, an 
{\it $H$-closed field}\/ with small derivation.\index{closed!$H$-closed}\index{H-closed@$H$-closed}\index{Hardy field!H-closed@$H$-closed} Now 
the elementary theory $T_H$ of $H$-closed fields with small derivation (denoted by $T^{\text{nl}}_{\text{small}}$ in [ADH]) is 
{\em complete}, by [ADH, 16.6.3]. This means in particular that any two maximal Hardy fields are indistinguishable as to their elementary 
properties: 

\begin{corintro} \label{cor:elem equiv} If $H_1$ and $H_2$ are maximal Hardy fields, then
$H_1$ and $H_2$ are elementarily equivalent as ordered differential fields. 
\end{corintro}

\noindent
To derive Theorem~\ref{thm:A} we use also the key results from the book
%\marginpar{we should also prove this result with the methods of the present notes}
 \cite{JvdH} to the effect that $\T_{\text{g}}$, the ordered differential field of grid-based transseries,\index{transseries!grid-based}
 is $H$-closed with small derivation and the intermediate value property for differential polynomials. 
%That the $H$-field~$\T_{\text{g}}$ is $\upo$-free is also clear from~[ADH, 11.7.15], so
In particular, it is a model of the complete
theory $T_H$. Thus maximal Hardy fields have the intermediate value property for differential polynomials as well,
and this amounts to Theorem~\ref{thm:A}, obtained here as a byproduct of more fundamental results.  (A more detailed account of the differential intermediate value property for $H$-fields is in \cite{ADHip}.) 
We sketch the proof of our main result (Theorem~\ref{thm:char d-max})
%thm:extend to H-closed}) l
later in this introduction, after describing further consequences.

%Somewhat earlier than Hardy, and also inspired by du Bois-Reymond, Hausdorff~\cite{Hau09} had already considered the broader class of  what we called {\it Hausdorff fields}\/ in \cite{ADH2}, that is, subfields of  the ring $\c$ of germs at $+\infty$ of merely {\it continuous}\/ functions~$(a,+\infty)\to\R$ ($a\in\R$).  Similarly to Hardy fields, each Hausdorff field  is an ordered field in a natural way. Hausdorff was particularly interested in the maximal  elements of this class (with respect to inclusion).   By Zorn, every Hausdorff field extends to a maximal one. Hausdorff showed (in modern language) that each maximal Hausdorff field is real closed in the sense of Artin-Schreier~\cite{Artin-Schreier}. Therefore, as a consequence of Tarski's famous theorem~\cite{Tarski}, any two maximal Hausdorff fields are elementarily equivalent, that is, indistinguishable   as to their elementary (first-order) properties. In this work we establish an analogue of this fact for maximal Hardy fields (viewed as {\it differential}\/ fields). Note that just like Hausdorff fields, clearly every Hardy field is contained in  a maximal Hardy field.  Each maximal Hardy field  is also real closed, but  maximal Hardy fields have many more closure properties: for example,  every maximal Hardy field $H$ contains the germs of all constant functions,  and for any~$h\in H$, contains the germ $\ex^h$  as well as any differentiable germ with derivative $h$. (See~\cite{Boshernitzan81, Ros} or Section~\ref{sec:Hardy fields} below.)

\subsection*{Further consequences of  our main result}
 In [ADH] we prove more than completeness of $T_H$:   a certain natural extension by definitions of $T_H$ has quantifier elimination.
This leads to a strengthening of Corollary~\ref{cor:elem equiv} by allowing parameters from a common Hardy subfield of~$H_1$ and~$H_2$.
To fully appreciate this statement requires more knowledge of model theory, as in~[ADH, Appendix~B], which we do not assume for this introduction. However,
we can explain a special case in a direct way,
in terms of  solvability of systems of algebraic differential
equations, inequalities, and asymptotic inequalities.
%We do need to recall the notion of a differential polynomial:   Let $R$ be a differential ring and let $Y,Y',\dots,Y^{(n)},\dots$ be distinct indeterminates.  The derivation on $R$ extends uniquely to a derivation on the polynomial ring $R\{Y\}:=R[Y,Y',\dots]$ such that  $(Y^{(n)})'=Y^{(n+1)}$ for all $n$, whose elements  are called  {\it differential polynomials} in the differential indeterminate $Y$ with  coefficients in $R$.  Inductively we put $R\{Y_1,\dots,Y_n\} := R\{Y_1,\dots,Y_{n-1}\}\{Y_n\}$ for distinct differential interderminates~$Y_1,\dots,Y_n$, $n\geq 1$.
Here we find it convenient to use  the notation
for asymptotic relations introduced by du~Bois-Reymond and Hardy  instead of Bachmann-Landau's  $O$-notation:  
  for germs~$f$,~$g$ in a Hardy field set
\begin{align*}  
f\preceq g	&\quad:\Longleftrightarrow\quad f=O(g) &\hskip-2em& :\Longleftrightarrow\quad \abs{f}\leq c\abs{g}\text{ for some real $c>0$,}\\
f\prec g	&\quad:\Longleftrightarrow\quad f=o(g) &\hskip-2em& :\Longleftrightarrow\quad \abs{f} < c\abs{g}\text{ for all real  $c>0$. }
\end{align*}  
Let now~$Y=(Y_1,\dots,Y_n)$ be a tuple of distinct (differential) indeterminates, and consider    a system of the following form:
\begin{equation}\tag{$\ast$}\label{eq:system}\begin{cases} &
\begin{matrix}
P_1(Y) & \varrho_1 & Q_1(Y) \\
\vdots & \vdots    & \vdots \\
P_k(Y) & \varrho_k & Q_k(Y)
\end{matrix}\end{cases}
\end{equation}
Here each $P_i$, $Q_i$ is a differential polynomial in $Y$  (that is, a
polynomial in the indeterminates $Y_j$ and their formal derivatives $Y_j',Y_j'',\dots$) with coefficients in our Hardy field~$H$, and each
$\varrho_i$ is one of the symbols ${=}$, ${\neq}$, ${\leq}$, ${<}$, ${\preccurlyeq}$,  ${\prec}$. Given a Hardy field $E\supseteq H$,
a {\it solution}\/\index{solution!system} of~\eqref{eq:system} in $E$ is an $n$-tuple $y=(y_1,\dots,y_n)\in E^n$ 
such that for~$i=1,\dots,k$, the relation
$P_i(y) \, \varrho_i\, Q_i(y)$ holds in $E$.
Here is a Hardy field analogue of the ``Tarski Principle''
of real algebraic geometry [ADH, B.12.14]: %\cite{Tarski}:

\begin{corintro}\label{cor:systems, 1}
If the system \eqref{eq:system} has a solution in \emph{some}  Hardy field extension of $H$, then \eqref{eq:system} has a solution in \emph{every} maximal
Hardy field extension of $H$. 
\end{corintro}

\noindent
(The symbols $\neq$, $\leq$, $<$, $\preceq$ in \eqref{eq:system} are for convenience only:
their occurrences  can be eliminated at the cost of increasing $m$, $n$.
But $\prec$ is essential; see [ADH, 16.2.6].) Besides the quantifier elimination alluded to, Corollary~\ref{cor:systems, 1} depends on Lemma~\ref{lem:canonical HLO}, which says that for any Hardy field  $H$ all maximal Hardy field extensions of $H$ induce the same $\HLO$-cut on $H$, as defined in [ADH, 16.3]. 

\medskip
\noindent
In particular, taking for $H$ the smallest Hardy field $\Q$, we see that a system \eqref{eq:system} with a solution in some  Hardy field  
has a solution in \emph{every} maximal
Hardy field, thus recovering a special case of our Corollary~\ref{cor:elem equiv}. Call such a system~\eqref{eq:system} over~$\Q$ {\it consistent.}\/ For example, with $X$, $Y$, $Z$ denoting here single distinct differential indeterminates, the system
$$Y'Z\ \preceq\ Z',\qquad Y \preceq 1, \qquad 1\prec Z$$
is inconsistent, whereas for any $Q\in\Q\{Y\}$ and $n\geq 2$ the system
$$X^n Y'\ =\ Q(Y),\qquad X' = 1,\quad Y\prec 1$$
is consistent.
As a consequence of the completeness of~$T_H$ we   obtain the existence of an algorithm (albeit a very impractical one) for deciding whether a system~\eqref{eq:system} over~$\Q$ is consistent, and this opens up the possibility of automating a substantial part of
asymptotic analysis in Hardy fields. 
%In Part~\ref{part:applications} we will also  give an axiomatization of sorts for the class of the consistent systems over~$\Q$. (See Section~\ref{sec:embeddings into T}.)
We  remark that Singer~\cite{Singer} proved the existence of an algorithm for
deciding whether a given system \eqref{eq:system} over $\Q$ without occurrences of
${\preccurlyeq}$ or ${\prec}$ has a solution in {\it some}\/ ordered differential field (and then it will have a solution 
in the ordered differential field of germs of real meromorphic functions at $0$); but  there are such systems, like
$$X'\ =\ 1,\qquad XY^2\ =\ 1-X,$$
which are solvable in an ordered differential field, but not in a Hardy field.
Also, algorithmically deciding the solvability of  a system~\eqref{eq:system} over~$\Q$ in a {\it given}\/ Hardy field
$H$ may be impossible when $H$ is ``too small'': e.g., if~$H=\R(x)$, by~\cite{Denef}.

\medskip
\noindent
As these results suggest, the aforementioned quantifier elimination for $T_H$ yields a kind of  ``resultant'' for systems~\eqref{eq:system} that allows one to make explicit within~$H$ itself for which choices of coefficients of the 
differential polynomials~$P_i$,~$Q_i$  the system~\eqref{eq:system} has a solution in a Hardy field
extension of $H$. Without going into details,
 we only mention here some attractive  consequences   for systems~\eqref{eq:system} depending on  parameters.
For this, let $X_1,\dots, X_m,Y_1,\dots, Y_n$ be distinct indeterminates and~$X=(X_1,\dots,X_m)$, $Y=(Y_1,\dots,Y_n)$,
and consider a system
\begin{equation}\tag{$\ast \ast$}\label{eq:system param, Y}\begin{cases} &
\begin{matrix}
P_1(X,Y) & \varrho_1 & Q_1(X,Y) \\
\vdots & \vdots    & \vdots \\
P_k(X,Y) & \varrho_k & Q_k(X,Y)
\end{matrix}\end{cases}
\end{equation}
where $P_i$, $Q_i$ are now differential polynomials in $(X,Y)$ over $H$, and the $\varrho_i$ are as before.
Specializing $X$ to $c\in \R^m$ then yields a system
\begin{equation}\tag{$\ast c$}\label{eq:system param}\begin{cases} &
\begin{matrix}
P_1(c,Y) & \varrho_1 & Q_1(c,Y) \\
\vdots & \vdots    & \vdots \\
P_k(c,Y) & \varrho_k & Q_k(c,Y)
\end{matrix}\end{cases}
\end{equation}
where $P_i(c,Y)$, $Q_i(c,Y)$ are  differential polynomials in $Y$ with coefficients in the Hardy field~$H(\R)$.
(We only substitute real constants, so may assume that the~$P_i$,~$Q_i$ are {\it polynomial}\/ in $X$, that is, none of the derivatives~$X_j',X_j'',\dots$ occur in the $P_i, Q_i$.)
Using~[ADH, 16.0.2(ii)] we obtain:

\begin{corintro}\label{cor:parametric systems}
The  set of all $c\in\R^m$ such that 
the system \eqref{eq:system param} has a solution in some Hardy field extension of $H$ is semialgebraic.
\end{corintro}

\noindent
Recall: a subset of $\R^m$ is said to be {\it semialgebraic}\/ if it is a finite union of  sets 
$$\big\{ c\in\R^m:\  p(c)=0,\ q_1(c)>0,\dots,q_l(c)>0 \big\}$$
where $p,q_1,\dots,q_l\in\R[X]$  are ordinary polynomials.\index{semialgebraic}
(The topological and geometric properties of semialgebraic sets have been studied extensively~\cite{BCR}. For example, it is well-known that a semialgebraic set can have only have finitely many connected components, and that each such component is itself semialgebraic.) 

\medskip
\noindent
In connection with  Corollary~\ref{cor:parametric systems} we mention that the asymptotics of Hardy field solutions to algebraic differential equations~${Q(Y)=0}$,
where~$Q$ is a differential polynomial with constant real coefficients, has been investigated by Hardy~\cite{Har12}  and   Fowler~\cite{Fowler} in  cases where $\order Q\leq 2$ (see \cite[Chapter~5]{Bellman}), and later by Shackell~\cite{SS95,Shackell93,Shackell95} in general.
Special case of  our corollary: for any differential polynomial~$P(X,Y)$ with constant real coefficients, the set of parameters~$c\in\R^m$ such that the differential
equation~${P(c,Y)=0}$ has a solution~$y$  in some Hardy field, in addition possibly also satisfying  given asymptotic  side conditions (such as~$y\prec 1$),  is semialgebraic. 
Example: the set of real parameters~$(c_1,\dots,c_{m})\in\R^{m}$ for which the homogeneous linear differential equation
$$y^{(m)} + c_{1}y^{(m-1)}+\cdots+ c_m y\  =\ 0$$
has a nonzero solution $y\prec 1$ in a Hardy field is semialgebraic; in fact, it is the
set of all  $(c_1,\dots,c_{m})\in\R^{m}$ such that the polynomial
$Y^m + c_{1}Y^{m-1}+ \cdots+ c_m \in \R[Y]$ has a negative real zero. (Below we discuss more general linear differential equations over Hardy fields.) Nonlinear example:
for $g_2, g_3\in \R$ the differential equation
$$(Y')^2\ =\ 4Y^3-g_2Y-g_3  $$
%of the Weierstrass $\wp$-function 
has a nonconstant solution in a Hardy field iff
 $g_2^3=27g_3^2$   and~$g_3 \leq 0$.
In both cases, the Hardy field solutions are  germs of logarithmico-exponential functions. But the
class of differentially algebraic germs in Hardy fields is much more extensive; for example, the antiderivatives of $\ex^{x^2}$  are not logarithmico-exponential (Liouville).

\medskip
\noindent
Instead of $c\in \R^m$, substitute $h\in H^m$ for $X$ in   \eqref{eq:system param, Y}, resulting  in a system 
\begin{equation}\tag{$\ast h$}\label{eq:system param, h}\begin{cases} &
\begin{matrix}
P_1(h,Y) & \varrho_1 & Q_1(h,Y) \\
\vdots & \vdots    & \vdots \\
P_k(h,Y) & \varrho_k & Q_k(h,Y)
\end{matrix}\end{cases}
\end{equation}
where $P_i(h,Y)$, $Q_i(h,Y)$ are now differential polynomials in $Y$ with coefficients in~$H$.
It is   well-known  that for any semialgebraic set $S\subseteq \R^{m+1}$  
there is a natural number~$B=B(S)$ such that for every $c\in\R^m$, if the section
$\big\{y\in \R: (c,y)\in S\big\}$
%$S\cap (\{c\}\times\R^n)$
 has~$> B$ elements, then this section has nonempty interior in $\R$. In contrast,
the set of solutions of \eqref{eq:system param, h} for $n=1$ in a maximal $H$ can be simultaneously infinite and discrete in the order topology of $H$: this happens precisely if some nonzero one-variable differential polynomial over $H$ vanishes on  this solution set [ADH, 16.6.11].
(Consider the example of the single algebraic differential equation $Y'=0$, which has solution set $\R$ in each maximal Hardy field.)
Nevertheless, we have the  following
uniform finiteness principle for solutions of~\eqref{eq:system param, h};
its proof is
considerably deeper than Corollary~\ref{cor:parametric systems} and
 also draws on results from~\cite{ADHdim}.

\begin{corintro}\label{cor:uniform finiteness}
There is a natural number $B=B\eqref{eq:system param, Y}$ such that for all $h\in H^m$:
if the system~\eqref{eq:system param, h} has $> B$ solutions in some Hardy field extension of $H$,
then~\eqref{eq:system param, h} has continuum many solutions in every maximal Hardy field extension of~$H$.
\end{corintro}

\noindent
Next we turn to issues of smoothness and analyticity in Corollary~\ref{cor:systems, 1}. 
By definition, a Hardy field is a differential subfield of the  differential ring $\Calinf$ consisting of the germs of functions $(a,+\infty)\to\R$ ($a\in\R$)  which are, for each~$n$, eventually $n$-times
continuously differentiable. Now $\Calinf$ has the differential subring $\Ginf$ whose elements are
the germs that are eventually $\c^\infty$. A {\it $\Ginf$-Hardy field}\/ is a Hardy field~$H\subseteq\Ginf$. (See~\cite{Gokhman} for an example
of a Hardy field~${H\not\subseteq\Ginf}$.)
A $\Ginf$-Hardy field  is said to be {\it $\Ginf$-maximal}\/ if it has no proper $\Ginf$-Hardy field
extension. Now $\Ginf$ in turn has 
the differential subring $\Gom$ whose elements are  the germs that are eventually
real analytic, and so we define likewise  $\Gom$-Hardy fields
($\Gom$-maximal Hardy fields, respectively). Our main theorems go through in the $\Ginf$- and $\Gom$-settings;
combined
with model completeness of~$T_H$ shown in~[ADH, 16.2] this ensures the existence of solutions with appropriate smoothness
in Corollary~\ref{cor:systems, 1}:

\begin{corintro}\label{cor:systems, 2}
If $H\subseteq\Ginf$ and the system \eqref{eq:system} has a solution in some Hardy field extension of $H$, then \eqref{eq:system} has a solution in every $\Ginf$-maximal Hardy field extension of $H$. In particular, if $H$ is   $\Ginf$-maximal  and
\eqref{eq:system} has a solution in a Hardy field extension of $H$, then it has a solution in $H$.
\textup{(}Likewise with $\Gom$ in place of $\Ginf$.\textup{)}
\end{corintro}

\noindent
We already mentioned  $\T_{\text{g}}$ as a quintessential example of an $H$-closed field. Its cousin~$\T$, the ordered differential field
of transseries,   extends $\T_{\text{g}}$ and is also $H$-closed with constant field~$\R$~[ADH, 15.0.2].\index{transseries} The elements of $\T$ are certain generalized  
series (in the sense of Hahn) in an indeterminate $x>\R$ with real coefficients,  involving exponential and logarithmic terms, such as
$$ f\  =\  \ex^{\frac{1}{2}\ex^x}-5\ex^{x^2}+\ex^{x^{-1}+2x^{-2}+\cdots}+\sqrt[3]{2}\log x-x^{-1}+\ex^{-x}+\ex^{-2x}+\cdots+5\ex^{-x^{3/2}}.$$
Mathematically  significant examples are the more simply structured trans\-series
\begin{multline*}
\operatorname{Ai}\ =\  \frac{\ex^{-\xi}}{2\sqrt{\pi}x^{1/4}}\sum_n (-1)^n   c_n \xi^{-n},
\quad
\operatorname{Bi}\ =\ \frac{\ex^{\xi}}{\sqrt{\pi}x^{-1/4}}\sum_n c_n \xi^{-n}, \\
\qquad\text{where $\ c_n\ =\ \frac{(2n+1)(2n+3)\cdots (6n-1)}{(216)^n n!}\ $ and $\ \xi\ =\ \frac{2}{3}x^{3/2}$,}
\end{multline*}
which are $\R$-linearly independent solutions of  the Airy equation $Y''=xY$ {\cite[Chap\-ter~11, (1.07)]{Olver}}.
For   information about~$\T$ see [ADH, Appendix~A] or~\cite{Edgar,JvdH}. We just mention here that
like each $H$-field, $\T$ comes equipped with its
own versions of the asymptotic relations $\preceq$, $\prec$, defined as for $H$ above. 
The asymptotic rules valid in all Hardy fields, such as
$$f\preceq 1\ \Rightarrow\ f'\prec 1,\qquad f\preceq g\prec 1\ \Rightarrow\ f'\preceq g', \qquad
f'=f\neq 0\ \Rightarrow\ f\succ x^n $$
also hold in $\T$. 
Here $x$ denotes, depending on the context, the germ of the identity function on $\R$, as well as the element   $x\in\T$.
(We make this precise in Section~\ref{sec:embeddings into T}, where we also give a finite axiomatization of these rules.)

\medskip
\noindent
Now suppose that we are given an embedding~$\iota\colon H\to\T$   of ordered differential fields.
We may view such an embedding as a {\it formal expansion operator}\/ and its inverse as a {\it summation operator.}\/
(See Section~\ref{sec:embeddings into T} below  for an example of a Hardy field, arising from a fairly rich o-minimal structure, which admits such an embedding.)
From~\eqref{eq:system} we obtain a system
\begin{equation}\tag{$\iota\ast$}\label{eq:iota(system)}\begin{cases} &
\begin{matrix}
\iota(P_1)(Y) & \varrho_1 & \iota(Q_1)(Y) \\
\vdots & \vdots    & \vdots \\
\iota(P_k)(Y) & \varrho_m & \iota(Q_k)(Y)
\end{matrix}\end{cases}
\end{equation}
of algebraic differential equations and (asymptotic) inequalities over $\T$,
where~$\iota(P_i)$, $\iota(Q_i)$   denote the differential polynomials over $\T$ obtained by applying $\iota$ to the coefficients
of $P_i$, $Q_i$, respectively. A {\it solution}\/\index{solution!system} of \eqref{eq:iota(system)} is a tuple $y=(y_1,\dots,y_n)\in\T^n$
such that $\iota(P_i)(y) \, \varrho_i\, \iota(Q_i)(y)$ holds in $\T$, for $i=1,\dots,m$.
Differential-difference equations in $\T$ are sometimes amenable to functional-analytic
techniques like fixed point theorems or small (compact-like) operators~\cite{vdH:noeth}, and
the formal nature of transseries also makes it possible to solve
algebraic differential equations in $\T$ by quasi-algorithmic methods~\cite{vdH:PhD,JvdH}. The simple example of the Euler equation
$$Y'+Y\ =\ x^{-1}$$
is instructive: its solutions in $\Calinf$
are given by the germs of 
$$t\mapsto \ex^{-t}\int_1^t \frac{\ex^{s}}{s}\,ds+c\ex^{-t}\colon (1,+\infty)\to\R\qquad (c\in\R),$$
all contained in a common Hardy field extension of $\R(x)$.
The solutions of this differential equation in $\T$ are 
$$\sum_n n!\, x^{-(n+1)}+c\ex^{-x}\qquad (c\in\R),$$
where the particular solution $\sum_n n!\, x^{-(n+1)}$ is obtained as the unique fixed point of the   operator $f\mapsto x^{-1}-f'$ on
the differential subfield
$\R(\!(x^{-1})\!)$ of $\T$ (cf.~[ADH, 2.2.13]).
(Note:  $\sum_n n!\,t^{-(n+1)}$ diverges for each $t>0$.)
In general, the existence of a solution of
\eqref{eq:iota(system)} in $\T$  entails the existence of a solution of~\eqref{eq:system}  in some  Hardy field extension of~$H$ and vice versa; more precisely:

\begin{corintro}\label{cor:systems, 3}
The system \eqref{eq:iota(system)} has a solution in $\T$ iff \eqref{eq:system} has a solution in some Hardy field extension of $H$.
In this case, we can choose a solution  of \eqref{eq:system} in a Hardy field extension $E$ of $H$ for which $\iota$ extends
to an embedding of ordered differential fields~$E\to\T$.
\end{corintro}

\noindent
In particular, a system \eqref{eq:system} over $\Q$ is consistent if and only if it has a solution in~$\T$.
(The ``if'' direction  already follows from  [ADH, Chapter~16] and~\cite{vdH:hfsol};  the latter constructs a summation operator  on the ordered differential subfield~$\T^{\operatorname{da}}\subseteq\T$   of   differentially algebraic transseries.)
 
\medskip
\noindent
It may seem remarkable that a  result about differential polynomials in one differentiable indeterminate, like
Theorem~\ref{thm:A} (or Theorem~\ref{thm:B} below), yields similar facts about \emph{systems}\/ of algebraic differential equations and asymptotic inequalities in several indeterminates over Hardy fields as in the corollaries above; we owe this to the strength of the model-theoretic methods employed in~[ADH]. But our theorem in combination with  [ADH] already has interesting consequences 
for one-variable differential polynomials over $H$ and over its ``complexification''~${K:=H[\imag]}$ (where~$\imag^2=-1$),
which is a differential subfield of the differential ring $\Calinf[\imag]$. Some of these facts are analogous to familiar properties of
ordinary one-variable polynomials over the real or complex numbers.
First, it follows from Theorem~\ref{thm:A} 
that every differential polynomial in a single differential indeterminate over $H$ of  odd degree has a zero in a Hardy field extension of $H$. (See Corollary~\ref{cor:odd degree}.) For example,  a differential polynomial like
$$ (Y'')^5+\sqrt{2}\ex^x (Y'')^4 Y'''-x^{-1}\log x \, Y^2 Y''+ YY'-\Gamma$$
has a zero in every maximal Hardy field extension of the Hardy field $\R\langle \ex^x,\log x,\Gamma\rangle$.
Passing to $K=H[\imag]$ we have:

\begin{corintro}\label{cor:zeros in complexified Hardy field extensions}
For each  differential polynomial~$P\notin K$ in a single differential indeterminate  with coefficients in $K$   there are~$f$, $g$ in a Hardy field extension of $H$  such that~$P(f+g\imag)=0$. 
\end{corintro}

\noindent
In particular, each   linear differential equation
$$y^{(n)}+a_{1}y^{(n-1)}+\cdots+a_n y\ =\  b\qquad (a_1,\dots,a_{n},b\in K)$$
has a solution $y=f+g\imag$ where $f$, $g$ lie  in some Hardy field extension  of $H$. (Of course, if~$b=0$, then we may take here the trivial solution $y=0$.) Although this special case of Corollary~\ref{cor:zeros in complexified Hardy field extensions} concerns
differential polynomials of degree~$1$, it seems hard to obtain this result without recourse to our more general extension theorems: a solution $y$
of a linear differential equation of order $n$ over $K$ as above may simultaneously be a zero of a non-linear differential polynomial $P$ over $K$ of    order~$< n$, and the structure of the differential field extension   of $K$ generated by $y$ is governed by $P$ (when taken of minimal complexity in the sense of [ADH, 4.3]).

\medskip
\noindent
Turning now   to homogeneous linear differential equations over Hardy fields, we first introduce some notation and terminology. Let~$R[\der]$ be the  ring  of linear differential operators over a differential ring $R$:
this ring  is a free left $R$-module with  basis~$\der^n$~($n\in\N$)
such that $\der^0=1$ and~$\der \cdot f=f\der+f'$ for~$f\in R$, where~${\der:=\der^1}$. (See~[ADH, 5.1] or~\cite[2.1]{vdPS}.)
Any  operator $A\in R[\der]$ gives rise to an additive map~${y\mapsto A(y)\colon R\to R}$, with $\der^n(y)=y^{(n)}$ (the $n$th derivative of $y$ in~$R$) and~$r(y)=ry$ for~$r=r\cdot 1\in {R\subseteq R[\der]}$. 
The elements of~$\der^n+R\der^{n-1}+\cdots+R\subseteq R[\der]$ are said to be  {\it monic}\/ of order~$n$.
It is well-known \cite{CamporesiDiScala,Mammana26,Mammana31} that for $R=\Calinf[\imag]$, each monic~$A\in R[\der]$   factors as a product of monic operators of order~$1$ in $R[\der]$; if~$A\in K[\der]$, then such a factorization already happens over the complexification of some Hardy field extension of $H$:

\begin{corintro}\label{cor:factorization intro}
If $H$ is maximal, then each monic operator in $K[\der]$  is a product of  monic operators of order $1$ in $K[\der]$.
\end{corintro}

\noindent
This follows quite easily from Corollary~\ref{cor:zeros in complexified Hardy field extensions} using the Riccati transform [ADH, 5.8].
{\it In the remainder of this subsection we let $A\in K[\der]$ be monic of order~$n$, and we fix a maximal Hardy field extension $E$ of $H$.}\/
The factorization result in Corollary~\ref{cor:factorization intro} gives rise to a description of a fundamental system of solutions for the  homogeneous linear differential equation~$A(y)=0$ in terms of Hardy field germs. Here, of course, complex exponential terms    naturally appear, but only in a controlled way: the $\C$-linear space consisting of all~$y\in\Calinf[\imag]$ with~$A(y)=0$ has a basis of the form
$$f_1\ex^{\phi_1\imag},\ \dots,\ f_n\ex^{\phi_n\imag}$$
where $f_j\in E[\imag]$ and $\phi_j\in E$ with $\phi_j=0$ or~$\abs{\phi_j}>\R$ for $j=1,\dots,n$.  We can arrange
here that for $i,j=1,\dots,n$ we have $\phi_i=\phi_j$ or $\abs{\phi_i-\phi_j}>\R$.
(Note that for $\phi$ in a Hardy field we have $\phi>\R$ iff $\phi(t)\to+\infty$ as~$t\to+\infty$.)
In this case, the basis elements $f_i\ex^{\phi_i\imag}$  for distinct frequencies $\phi_i$  are pairwise orthogonal
in a sense made precise in Section~\ref{sec:lin diff applications}. 

\begin{example}
If $y\in\Calinf[\imag]$ is {\it holonomic,}\/ that is, $L(y)=0$ for some monic $L\in\C(x)[\der]$, then $y$
is a $\C$-linear combination of germs $f\ex^{\phi\imag}$ where $f\in E[\imag]$, $\phi\in E$, and~$\phi=0$ or~$\abs{\phi}>\R$.
Here, more information about the~$f$,~$\phi$ is available (see, e.g.,~\cite[VIII.7]{FS}, \cite[\S{}19.1]{Wasow}).
Many special functions %such as the hypergeometric functions, 
are holonomic~\cite[B.4]{FS}.
\end{example}

\noindent
By the usual correspondence between linear differential operators and matrix differential equations (see, e.g., [ADH, 5.5]),  our results about zeros of linear differential operators  also yield facts about systems~${y'=Ny}$  of linear differential equations over  Hardy fields. If the matrix $N$ has suitable symmetry, we can even guarantee the existence of a nonzero solution~$y$ which lies in $E[\imag]^n$
 (and thus does not exhibit   oscillating behavior).
A sample result, also shown  in Section~\ref{sec:lin diff applications}: every matrix differential equation  $y'=Ny$,  where $N$ is an~$n\times n$ matrix~over~$K$~(${n\geq 1}$),  has a nonzero solution $y\in E[\imag]^n$ provided $n$ is {\it odd}\/ and~$N$ is {\it skew-symmetric.}\/ 
(The study of such matrix differential equations, for~$n=3$, goes back at least to Darboux~\cite[Livre~I, Chapitre~II]{Darboux}.)
For example, for each $c,d\in\C$ there is a nonzero $y\in E[\imag]^3$ such that
$$y'\ =\ \begin{pmatrix}
0						& -c											& \displaystyle-\frac{d}{\ex^x+\ex^{-x}} \\
c						& 0											& \displaystyle\frac{\ex^x-\ex^{-x}}{\ex^x+\ex^{-x}} \\
\displaystyle\frac{d}{\ex^x+\ex^{-x}}& \displaystyle\frac{\ex^{-x}-\ex^{x}}{\ex^x+\ex^{-x}}	& 0\end{pmatrix} y.$$
(For~$c=0$,~$d=2\sqrt{2}$, this equation is studied in \cite{AHJO}.)
%Examples include the time-dependent Schr\"odinger equation for a quantum system with $n$ states and purely imaginary Hamiltonian~\cite[\S{}3.4]{Teschl}.) 

\medskip
\noindent
We now return to the operator setting and focus on the case where~$A$ is real, that is, $A\in H[\der]$. 
Mammana~\cite{Mammana31}  conjectured that   each monic operator in $\Calinf[\der]$ of odd order has a monic factor of order $1$;
this is false in general~\cite{Sansone} but holds in the Hardy field world,
thanks to our ``real'' version of Corollary~\ref{cor:factorization intro}:

\begin{corintro}\label{cor:factorization intro, real}
Suppose   $A\in H[\der]$. Then $A$ is a product of monic operators in~$E[\der]$, each of order~$1$ or
irreducible of order~$2$. %Moreover, for each $b\in H$ there is some~$y\in H$ with $A(y)=b$.
\end{corintro}

\noindent
As a consequence, the $\R$-linear space of zeros of $A\in H[\der]$ in~$\Calinf$  has a basis
$$g_1\cos\phi_1,\ g_1\sin\phi_1,\ \dots,\ g_r\cos\phi_r,\ g_r\sin\phi_r,\ h_1,\ \dots,\ h_s\qquad (2r+s=n)$$
where $g_j,\phi_j\in E$ with $\phi_j>\R$ for $j=1,\dots,r$ and $h_k\in E$ for $k=1,\dots,s$.
In particular, if $n$ is odd, then $A(y)=0$ for some nonzero $y\in E$. 
%Each matrix differential equation $y'=Ny$ as above, where now   $N$ is an $n\times n$ matrix {\it over~$H$}\/ and $n$ is odd, has a nonzero solution $y$ in $E^n$ (regardless whether $N$ is skew-symmetric or not).

\medskip
\noindent
A   function $y\colon [a,+\infty)\to\R$ ($a\in\R$) is {\it non-oscillating}\/ if $\sgn y(t)$ is eventually constant (and otherwise
$y$ {\em oscillates}).
Similarly we define   (non-) oscillation of germs.  No germ in a Hardy field oscillates.
The following corollary  characterizes  when  $A$ in   Corollary~\ref{cor:factorization intro, real} is a product of monic operators of order $1$ in $E[\der]$:

\begin{corintro}\label{cor:unique factorization}
The $($monic$)$ operator $A\in H[\der]$ is  a product of  monic operators of order $1$ in $E[\der]$  iff no zero of
$A$ in $\Calinf$ oscillates. In this case~$E$ contains a basis~$y_1\prec\cdots\prec y_n$ of the $\R$-linear space of  zeros of $A$ in $\Calinf$,
and $$A\ =\ (\der-a_n)\cdots(\der-a_1)$$ for a unique tuple~$(a_1,\dots,a_n)\in E^n$ 
%of elements of $\Calinf$ 
such that for all sufficiently small~$f>\R$ in $E$ we have $a_j+(f''/f') < a_{j+1}$ for~$j=1,\dots,n-1$.
% moreover, these $a_j$ lie in $E$.
%this tuple~$(a_1,\dots,a_n)$ is also contained in $E$.
\end{corintro}

%\begin{corintro}\label{cor:unique factorization}
%Suppose $A\in H[\der]$. Then  $A$ is   a product of  monic operators of order $1$ in $E[\der]$  precisely if
%$A$ has only non-oscillating zeros in $\Calinf$. In this case~$E$ contains a basis $y_1\prec\cdots\prec y_n$ of the $\R$-linear space of  zeros of $A$ in $\Calinf$,
%and $$A=(\der-a_n)\cdots(\der-a_1)$$ for a unique tuple~$(a_1,\dots,a_n)$ of elements of $\Calinf$ with $a_n\geq\cdots\geq a_1$,
%also contained in $E$.
%\end{corintro}

\noindent
Factorizations of linear differential operators as in Corollary~\ref{cor:unique factorization} are closely connected to the classical topic of {\it disconjugacy.}\/
We recall the definition, which arose from the calculus of variations~\cite{Wintner}.  
Let  $f_1,\dots,f_{n}\colon I\to\R$ be continuous, where $I=[a,+\infty)$, $a\in\R$. The linear differential equation
\begin{equation}\tag{L}\label{eq:disconj intro}
y^{(n)}+f_{1}y^{(n-1)}+\cdots+f_ny\ =\ 0
\end{equation}
on $I$ is said to be {\it disconjugate}\/ if every nonzero solution $y\in \mathcal{C}^n(I)$ of \eqref{eq:disconj intro} has at most~$n-1$ zeros, counted with their multiplicities. (For example, $y^{(n)}=0$, on any such $I$, is disconjugate.) 
The solutions of disconjugate linear differential equations are suitable for approximation and interpolation purposes; see  \cite[Chapter~3]{Coppel} and \cite[Chapter~3, \S{}11]{DVL}.
We also say that \eqref{eq:disconj intro} is
{\it eventually disconjugate}\/ if for some~$b\geq a$  the linear differential equation on $J:=[b,+\infty)$
obtained from  \eqref{eq:disconj intro} by restricting~$f_1,\dots,f_n$ to $J$ is disconjugate.   
If \eqref{eq:disconj intro} is eventually disconjugate, then it has no oscillating solutions in $\mathcal{C}^n(I)$. The converse
of this implication holds when~$n\leq 2$ but fails for each $n>2$~\cite{Gustafson}.
There is an extensive literature, mostly dating back to the 1970s, which develops sufficient conditions for (eventual) disconjugacy
of linear differential equations (see, e.g.~\cite{Coppel,Elias77,Elias97, Gregus, Swanson}), 
often by restricting  the growth of the $f_i$;
for example, \eqref{eq:disconj intro} is disconjugate if~$\int_a^\infty \abs{f_i(t)} (t-a)^{i-1}\,dt<\infty$ for~$i=1,\dots,n$ (cf.~\cite{Willett}).
Corollary~\ref{cor:unique factorization} allows us to contribute another natural criterion for 
eventual disconjugacy:

\begin{corintro}\label{cor:disconjugacy}
If the germs of $f_1,\dots,f_{n}$ lie in a Hardy field and~\eqref{eq:disconj intro} has no oscillating solutions in $\mathcal{C}^n(I)$, then 
 \eqref{eq:disconj intro}  is eventually disconjugate.
\end{corintro}

\noindent
A fundamental property of disconjugate linear differential operators is the existence of a
canonical factorization  
discovered by Trench~\cite{Trench}. (See also Proposition~\ref{prop:Trench} below.)
Corollary~\ref{cor:unique factorization} can also be used to strengthen this factorization  
in the situation of Corollary~\ref{cor:disconjugacy}. See Corollary~\ref{cor:Trench} for the details.

\medskip
\noindent
We finish with discussing the instructive   case of an operator $A\in H[\der]$ of order~$2$. % in Corollary~\ref{cor:factorization intro, real}. 
If such $A$ has a  non-oscillating zero $y\neq 0$ in~$\Calinf$, then by Sturm's Oscillation Theorem all  zeros of $A$ in $\Calinf$ are non-oscillating
and hence contained in every maximal Hardy field,
by Corollary~\ref{cor:char osc} or \cite[Theorem~16.7]{Boshernitzan82}, \cite[Corollary~2]{Ros}. 
For example, the germs of the  $\R$-linearly independent solutions $\operatorname{Ai},\operatorname{Bi}\colon\R\to\R$ 
of the Airy equation~$Y''-xY=0$ given by
\begin{align*}
\operatorname{Ai}(t)\ &=\ \frac{1}{\pi}\int_0^\infty \cos\left(\frac{s^3}{3}+st\right)ds, \\
\operatorname{Bi}(t)\ &=\ \frac{1}{\pi}\int_0^\infty \left[\exp\left(-\frac{s^3}{3}+st\right) +\sin\left(\frac{s^3}{3}+st\right)\right]ds
\end{align*}
lie in each maximal Hardy field, with $\operatorname{Ai}\prec 1 \prec \operatorname{Bi}$.
 In the oscillating case, we have:

\begin{corintro}\label{cor:Boshernitzan intro}
If $A\in H[\der]$ of order $2$ has an oscillating zero in $\Calinf$, then there are~$g,\phi\in E$ with $\phi>\R$  such that
the zeros of $A$ in $\Calinf$ are exactly the germs~${cg\cos(\phi+d)}$ \textup{(}$c,d\in\R$\textup{)}. 
%$$A(y)=0\qquad\Longleftrightarrow\qquad \text{$y=cg\cos(\phi+d)$ for some $c,d\in\R$.}$$
\end{corintro}

\noindent
This corollary was announced by Boshernitzan  as part of~\cite[Theorem~5.4]{Boshernitzan87}, but apparently a proof of this theorem never appeared in print. (See also \cite[Conjecture~4 in \S{}20]{Boshernitzan82}.) 
In Section~\ref{sec:perfect applications} below we state and prove a strengthening of his theorem; this includes a criterion for the uniqueness of the  germs~$g$,~$\phi$.
For every~$\phi>\R$ in $E$ there is at most one~$g\in E$ (up to multiplication by a nonzero constant)   such that
the conclusion of Corollary~\ref{cor:Boshernitzan intro} holds; in general, the 
pair $(g,\phi)$ is not unique (up to multiplication of $g$ by a nonzero constant
and addition of a constant to $\phi$),  but it is if the coefficients of $A$ are differentially algebraic
\textup{(}over~$\R$\textup{)}. As a final example, consider the Bessel equation of order $\nu\in\R$ (see, e.g., \cite[VII]{EMOT}, \cite[I, \S{}9]{Olver}, \cite{Watson}): 
$$x^2Y''+x Y'+ (x^2-\nu^2)Y\ =\ 0.$$
It is well-known that 
   each   solution~$y\in\Calinf$ of this equation satisfies
$$y\ =\ cx^{-1/2} \cos (x+d)+o(x^{-1/2})\quad\text{for some~$c,d\in\R$.}$$
(See \cite[Chapter~6, \S{}18]{Bellman}, \cite[Corollary~XI.8.1]{Hartman}, \cite[Example~13.2]{Wasow}.) %\cite[\S{}27, XVII(g)]{Walter})
%that every solution  $y\in\Calinf$  of this equation is oscillating and satisfies
%$\abs{y}\leq C/\sqrt{x}$ for some constant $C\in\R$, and its amplitudes are strictly decreasing.
We give a similar parametrization  using germs in Hardy fields.
More precisely, we show that there is a unique germ~$\phi=\phi_\nu$ in a Hardy field with   $\phi-x\preceq x^{-1}$  such that
every solution~$y\in\Calinf$ of the Bessel equation has   the form
$$y\ =\ cx^{-1/2}g \cos (\phi+d)\quad\text{for some $c,d\in\R$ and $g:=1/\sqrt{\phi'}$.}$$
This explains the
phenomenon, observed in \cite{HBR, HBRV}, that the Bessel equation admits a non-oscillating phase function. 
Knowing that $\phi$ lives in a Hardy field allows one to reprove a number of classical results about Bessel functions in a short transparent way.
Remarkably, every maximal Hardy field contains
the phase function $\phi$, and $\phi$ has an asymptotic expansion
$$\phi\sim x+\frac{\mu-1}{8}x^{-1}+\frac{\mu^2-26\mu+25}{384}x^{-3}+\frac{\mu^3-115\mu^2+1187\mu-1073}{5120} x^{-5}+\cdots$$
where $\mu=4\nu^2$. Only for special choices of~$\nu$ 
is the germ~$\phi$   contained in the Liouville closure of~$\R(x)$, and hence easily obtainable by the classical extension results for Hardy fields from \cite{Boshernitzan81,Har12a,Ros}:  by results of Liouville~\cite{Liouville41}   this holds  precisely if~$\nu\in \frac{1}{2}+\Z$;  see Section~\ref{sec:Bessel}  for a proof.

\medskip
\noindent
Michael Boshernitzan's pa\-pers~\cite{Boshernitzan81}--\cite{BoshernitzanUniform} on Hardy fields have
been a frequent source of inspiration for us, and we dedicate this work to his memory. 
He paid particular attention to the germs in $\Calinf$, such as $\phi_\nu$ above, that lie in {\it every}\/ maximal Hardy field.
They form  a Liouville closed Hardy field $\Ex$ properly containing Hardy's differential field of logarithmico-exponential functions.
In the course of our work below we prove 
Conjecture~1 from~\cite[\S{}10]{Boshernitzan81} and Conjecture~4 from \cite[\S{}20]{Boshernitzan82} about~$\Ex$. We also prove
Conjecture~17.11 from~\cite[\S{}17]{Boshernitzan82} and answer Question~4 from \cite[\S{}7]{Boshernitzan86}. (See Corollaries~\ref{cor:Bosh Conj 1}, Theorems~\ref{thm:Bosh} and~\ref{cor:17.7 generalized},  and Proposition~\ref{prop:translog}, respectively.)  
Section~\ref{sec:holes perfect} contains some additional  observations which may eventually help to shed further light on the nature of the Hardy field~$\Ex$.

\subsection*{Synopsis of the proof of our main theorem} 
In the rest of the paper we assume familiarity with the terminology and concepts of asymptotic differential algebra from our book~[ADH].
(We review some of this in the last subsection of the introduction below.)
The proof of our main result requires, besides differential-algebraic and valuation-theoretic tools from~[ADH],
also   analytic arguments in an essential way. Some of our analytic machinery is obtained by adapting material from \cite{vdH:hfsol} to a more general setting.
As explained earlier, our main Theorem~\ref{thm:A} is a consequence of the following extension theorem:

\begin{theorem}\label{thm:B}
Every $\upo$-free Hardy field has a newtonian  Hardy field extension.
\end{theorem}

\noindent
The proof of this is long, so it may be useful to outline the strategy behind it. 

\subsubsection*{Holes and slots}
For now, let $K$ be an $H$-asymptotic field with rational asymptotic integration.
In Section~\ref{sec:holes} below we introduce the apparatus of {\it holes}\/ in~$K$ as a means to systematize the study of   solutions of  
algebraic differential equations over~$K$ in immediate asymptotic extensions of~$K$: such a hole in $K$ is a triple~$(P,\fm,\hat f)$
where~$P$ is a differential polynomial in a single differential indeterminate $Y$ with coefficients in $K$, $P\ne 0$,  $0\neq \fm\in K$, and $\hat f\notin K$ lies an immediate asymptotic extension  of~$K$
with~${P(\hat f)=0}$ and~$\hat f\prec\fm$. 
%We shall see shortly that
It is sometimes technically convenient to work with the more flexible concept of a {\it slot}\/  in $K$, where instead of~$P(\hat f)=0$  we only require  
$P$ to vanish at $(K,\hat f)$ in the sense of~[ADH, 11.4]. The {\it complexity}\/ of a slot $(P,\fm,\hat f)$ is the complexity of the differential polynomial~$P$ as  in~[ADH, p.~216]. Now if $K$ is $\upo$-free, then by~Lem\-ma~\ref{lem:no hole of order <=r},
$$\text{$K$ is newtonian}\quad\Longleftrightarrow\quad\text{$K$ has no hole.}$$
This equivalence suggests an opening move for proving Theorem~\ref{thm:B} by induction on complexity as follows:
Let $H\supseteq\R$ be an $\upo$-free Liouville closed Hardy field, and suppose $H$ is not newtonian; it is enough to show that then $H$ has a proper Hardy field extension. By the above equivalence, $H$ has a hole
$(P, \fm,\hat f)$, and we can take here $(P,\fm,\hat f)$ to be of minimal complexity among holes in $H$.
This minimality has consequences that are important for us; for example $r:=\order P\geq 1$, $P$ is a minimal annihilator of $\hat f$ over $H$, and $H$ is $(r-1)$-newtonian as defined in~[ADH, 14.2]. We arrange $\fm=1$  by
replacing $(P,\fm,\hat f)$ with the hole~$(P_{\times\fm},1,\hat f/\fm)$ in $H$.

\subsubsection*{Solving algebraic differential equations over Hardy fields}
For Theorem~\ref{thm:B} it is enough to show that under these conditions $P$ is a minimal annihilator
of some germ~${f\in \Calinf}$ that generates a (necessarily proper) Hardy field extension~$H\langle f\rangle$ of~$H$.  So at a minimum, we need to  find    a  solution in $\Calinf$ to the algebraic differential equation~$P(Y)=0$. For this, it is natural to use fixed point techniques as in~\cite{vdH:hfsol}.
Notation: for~${a\in\R}$, let $\c^n_a$ be the $\R$-linear space of functions~$[a,+\infty)\to\R$ which extend to an $n$-times continuously differentiable
function $U\to\R$ on an open subset~$U\supseteq [a,+\infty)$ of~$\R$. For any $a$ and $n$, each germ in $\Calinf$ has representatives in~$\c^n_a$. 

\subsubsection*{A fixed point theorem}
Let $L:=L_P\in H[\der]$ be the linear part  of $P$. Replacing~$(P,1,\hat f)$ with another minimal hole in~$H$ we arrange ${\order L=r}$.
Representing the coefficients of  $P$ (and thus of $L$) by functions in~$\c^0_a$ we obtain  an $\R$-linear operator $y\mapsto L(y)\colon \c^r_a \to\c^0_a$.
For now we make the  bold assumption that   $L\in H[\der]$ splits over  $H$. Using such a splitting and increasing $a$ if necessary, $r$-fold integration yields an $\R$-linear operator $L^{-1}\colon \c^0_a\to\c^r_a$ which is a {\it right-inverse}\/ of~$L\colon \c^r_a \to\c^0_a$, that is, $L\big(L^{-1}(y)\big)=y$ for all $y\in\c^0_a$.
Consider    the (generally non-linear) operator
$$f\mapsto \Phi(f):=L^{-1}\big(R(f)\big)$$ on $\c^r_a$; here~$P=P_1-R$ where~$P_1$ is the homogeneous part of degree~$1$ of~$P$.
% and we assume that $a$ is chosen so that the coefficients of $R$ are represented by functions in $\c^0_a$. Ideally we could hope 
We try to show that $\Phi$ restricts to a contractive operator
 on a closed ball of an appropriate subspace of~$\c^r_a$ equipped with a suitable complete norm,  
 whose fixed points are then solutions to~$P(Y)=0$; this may also involve increasing $a$ again and replacing the coefficient functions of $P$ by their corresponding restrictions. 
 To obtain such contractivity, we would need to ensure that $R$ is asymptotically small compared to~$P_1$ in a certain sense.
This can indeed be achieved by transforming
$(P,1,\hat f)$ into a certain normal form through successive {\it refinements}\/ and ({\it additive}\/, {\it multiplicative}\/, and {\it compositional}\/) {\it conjugations}\/ of the hole~$(P,1,\hat f)$. This normalization
is done under more general algebraic assumptions in Section~\ref{sec:normalization}.  The analytic arguments leading to fixed points are in Sections~\ref{sec:IHF}--\ref{sec:smoothness}.  Developments below involve the algebraic closure~$K:=H[\imag]$ of $H$ and we work more generally with a decomposition~$P=\tilde{P}_1-R$ where $\tilde{P}_1\in K\{Y\}$ is   homogeneous of degree $1$, not necessarily~$\tilde{P}_1=P_1$, such that
$L_{\tilde{P}_1}\in K[\der]$ splits and $R$ is ``small'' compared to $\tilde{P}_1$.

\subsubsection*{Passing to the complex realm} 
In general we are not so lucky that~$L$ splits over~$H$. The minimality of our hole $(P,1,\hat f)$ does not even ensure that $L$ splits over~$K$. At this point we recall from [ADH, 11.7.23] that $K$ is $\upo$-free because $H$ is.  
We can also draw hope from the fact that every nonzero linear differential operator over $K$ would split over $K$ if
$H$ were newtonian [ADH, 14.5.8].  Although $H$ is not newtonian, it is $(r-1)$-newtonian, and $L$ is only of order $r$, so
we optimistically restart our attempt, and instead of a hole of minimal complexity in $H$, we now let~$(P,\fm,\hat f)$ be a hole of minimal complexity in $K$. Again it follows that $r:=\order P\geq 1$,  $P$ is a minimal annihilator of $\hat f$ over $K$, and $K$ is $(r-1)$-newtonian.
As before we arrange that~$\fm=1$ and the linear part~$L_P\in K[\der]$ of $P$ has order~$r$. We can also arrange~$\hat f\in \hat K = \hat H[\imag]$ where~$\hat H$ is an immediate asymptotic extension of~$H$. So~$\hat f=\hat g+\hat h\imag$ where $\hat g,\hat h\in\hat H$ satisfy~$\hat g,\hat h\prec 1$, and~$\hat g\notin H$ or~$\hat h\notin H$, say $\hat g\notin H$.
 Now minimality of $(P,1,\hat f)$ and algebraic closedness of $K$ give that $K$ is $r$-linearly closed, that is, every nonzero~$A\in K[\der]$ of order~$\leq r$ splits over $K$  (Corollary~\ref{corminholenewt}). Then~$L_P$   splits over~$K$ as desired, and a   version of the above fixed point construction with $\c^r_a[\imag]$ in place of~$\c^r_a$ can be carried out successfully to solve~$P(Y)=0$ in the differential ring extension $\Calinf[\imag]$ of $\Calinf$. 

\subsubsection*{Return to the real world} 
But at this point we face another obstacle: even once we have our hands on a zero~$f\in\Calinf[\imag]$ of $P$, it is not clear why~$g:=\Re f$ should generate a proper Hardy field extension of $H$: Let $Q$ be a minimal annihilator of~$\hat g$ over~$H$; we cannot expect that $Q(g)=0$. 
If $L_Q\in H[\der]$ splits over~$K$, then 
we can try to apply  fixed point arguments like the ones above, with~$(P,1,\hat f)$ replaced by the hole $(Q,1,\hat g)$ in~$H$, 
to  find a zero~$y\in\Calinf$ of~$Q$. (We do need to take care that constructed zero is real.)
Unfortunately we can only ascertain  that~$1\leq s\leq 2r$ for~$s:=\order Q$, and  since we may have~${s>r}$,  we cannot leverage the    minimality of~$(P,1,\hat f)$ anymore to ensure that $L_Q$ splits over~$K$, or to normalize~$(Q,1,\hat g)$ in the same way as 
indicated above for~$(P,1,\hat f)$. 
This situation seems hopeless, but 
now a purely differential-algebraic observation  comes to the rescue: although  the linear part~$L_{Q_{+\hat g}}\in\hat H[\der]$ of the differential polynomial~$Q_{+\hat g}\in\hat H\{Y\}$ also has order~$s$ (which may be $>r$),   {\it if $\hat K$ is $r$-linearly closed, then  $L_{Q_{+\hat g}}$ does split over~$\hat K$}; see~[ADH, 5.1.37]. If moreover~$g\in H$ is sufficiently close to~$\hat g$, then the linear part~$L_{Q_{+g}}\in H[\der]$ of~$Q_{+g}\in H\{Y\}$ is close to an operator in $H[\der]$ that does split
over~$K=H[\imag]$, and so
%``almost'' split over~$K=H[\imag]$,
using~$(Q_{+g},1,\hat g-g)$ instead of~$(Q,1,\hat g)$  may offer a way out of this impasse.

\subsubsection*{Approximating $\hat g$}  Suppose for a moment that $H$ is (valuation) dense in $\hat H$. 
Then by extending $\hat H$ we arrange that $\hat H$ is the completion of $H$, and $\hat K$ of $K$ (as in~[ADH, 4.4]). In this case $\hat K$ inherits from $K$ the property of being $r$-linearly closed, by results in~Section~\ref{sec:complements newton}, and the desired approximation of $\hat g$ by $g\in H$ can be achieved. 
%: although it can always be achieved by
%replacing~$\hat H$ with an immediate newtonian extension  of its $\upo$-free  $H$-subfield   $H\langle\hat g,\hat h\rangle$ (see [ADH, 14.1]),   this
%is prohibitive in view of the intended approximation arguments.
We cannot in general expect $H$ to be dense in $\hat H$. But we are saved by Section~\ref{sec:special elements}  to the effect that $\hat g$ can be made {\it special}\/ over~$H$ in the sense of [ADH, 3.4], that is, some nontrivial convex subgroup $\Delta$ of
the value group of $H$ is cofinal in $v(\hat g-H)$. Then passing to the $\Delta$-specializations of the various
valued differential fields encountered above (see [ADH, 9.4]) we regain density and this allows us to implement the desired approximation. The technical details are involved, and are carried out in the first three sections of Part~\ref{part:dents in H-fields}.  
A minor  obstacle to obtain the necessary specialness of~$\hat g$ is that the hole~$(Q,1,\hat g)$ in~$H$ may not be of minimal complexity. This can be  ameliorated by using a differential polynomial of minimal complexity vanishing at $(H,\hat g)$ instead of~$Q$,   in the process replacing the hole~$(Q,1,\hat g)$ in $H$ by a  slot in~$H$, which we then aim to approximate   by  a {\it strongly split-normal}\/ slot in~$H$; see Definition~\ref{def:strongly repulsive-normal}.
Another caveat:  to carry out our approximation scheme  we require $\deg P>1$. %\marginpar{what about $r=1$?} 
Fortunately,  if $\deg P=1$, then necessarily $r=\order P=1$, and this case
can be dealt with through separate arguments: see Section~\ref{sec:d-alg extensions} where we
finish the proof of Theorem~\ref{thm:B}.

\subsubsection*{Enlarging the Hardy field}
Now suppose we have finally arranged things so that our Fixed Point Theorem applies: it delivers $g\in\Calinf$ such that~$Q(g)=0$ and~$g\prec 1$.
(Notation: for a germ~$\phi\in\Calinf[\imag]$ and~$0\neq \fn\in H$  we write $\phi\prec \fn$ if $\phi(t)/\fn(t)\to 0$ as~$t\to+\infty$.)
However, in order that~$g$ generates a proper Hardy field extension~$H\langle g\rangle$ of $H$
isomorphic to~$H\langle\hat g\rangle$ by an isomorphism over~$H$ sending $g$ to $\hat g$ requires that $g$ and $\hat g$ have similar asymptotic properties with respect
to the elements of~$H$. For example, suppose~$h,\fn\in H$  and~$\hat g-h\prec\fn\preceq 1$; then we must show $g-h\prec\fn$. 
(Of course, we need to show much more about the asymptotic behavior of $g$, and this is expressed using the notion of {\it asymptotic similarity}\/: see
Sections~\ref{sec:asymptotic similarity} and~\ref{sec:d-alg extensions}.)
Now the germ $(g-h)/\fn\in\Calinf$ is a zero of the conjugated differential polynomial~$Q_{+h,\times\fn}\in H\{Y\}$, as is the element~$(\hat g-h)/\fn\prec 1$ of~$\hat H$. The Fixed Point Theorem can also be used to  produce a zero~$y\prec 1$ of $Q_{+h,\times \fn}$ in $\Calinf$. Set $g_1:=y\fn+h$; then~$Q(g)=Q(g_1)=0$ and~$g,g_1\prec 1$. We are thus naturally lead to consider the difference $g-g_1$ between
the solutions~$g,g_1\in\Calinf$ of the differential equation (with asymptotic side condition)
\begin{equation}\label{eq:Q=0}\tag{$\operatorname{E}$}
Q(Y)\ =\ 0,\qquad Y\prec 1.
\end{equation} 
If we manage to show~$g-g_1\prec\fn$, then $g-h=(g-g_1)-y\fn\prec\fn$ as required. Simple estimates coming out of the proof of the Fixed Point Theorem are not good enough for this (cf.~Lemma~\ref{lem:close}). 
We need a generalization of the Fixed Point Theorem
for {\it weighted norms}\/ with (the germ of) the relevant  weight function given by $\fn$, shown in Section~\ref{sec:weights}.
To render this generalized version useful, we also have 
to make the  construction of the right-inverse~$A^{-1}$ of the linear differential operator $A\in H[\der]$, which splits over $K$ and approximates $L_Q$ as postulated by strong split-normality,  and which is central for the  definition of the contractive operator used in the Fixed Point Theorem, in some sense uniform in $\fn$.
This is carried out in Section~\ref{sec:repulsive-normal}, refining our approximation arguments  by improving strong split-normality to {\it strong repulsive-normality}\/ as defined in~\ref{def:strongly repulsive-normal}.

\subsubsection*{Exponential sums}
Just for this discussion, call  $\phi\in\Calinf[\imag]$ {\it small}\/ if~$\phi\prec\fn$ for all~${\fn\in H}$ with
$v\fn\in v(\hat g-H)$. Thus our aim is to show that differences between solutions of \eqref{eq:Q=0} in $\Calinf$ are small in this sense. 
We 
%first observe
show that each such difference gives rise to a zero $z\in\Calinf[\imag]$  of~$A$ with $z\prec 1$
whose smallness would imply  the smallness of the difference under consideration. To ensure that every zero $z\prec 1$ of $A$
is indeed small requires us to have performed beforehand yet another (rather unproblematic) normalization procedure on our  slot, transforming it into {\it ultimate}\/ shape.
(See Section~\ref{sec:ultimate}.) 
Recall the special fundamental systems of solutions to linear differential equations over
maximal Hardy fields explained after Corollary~\ref{cor:factorization intro}: since $A$ splits over~$K$, our zero $z$ of $A$ is a $\C$-linear combination of exponential terms.
As a tool for systematically dealing with such exponential sums over $K$ in a formal
way, we introduce the concept of the {\it universal exponential extension}\/ of a differential field. 
Finally, from conditions like~$z\prec 1$ we need to be able to obtain asymptotic information about the summands of $z$
when expressed as an exponential sum in a certain canonical way.
For this we are able to exploit facts about uniform distribution mod~$1$ for germs in Hardy fields due to Boshernitzan~\cite{BoshernitzanUniform}; see Sections~\ref{sec:almost periodic}--\ref{sec:ueeh}. 

\subsection*{Organization of the manuscript}
Part~\ref{part:preliminaries} has preliminaries on linear differential operators and differential polynomials, on the group of logarithmic derivatives, on special elements, and on differential-henselianity and newtonianity.
In Part~\ref{part:universal exp ext} we define the universal exponential extension of a differential field,
and we  consider the eigenvalues of linear differential operators and their connections to splittings.
Part~\ref{part:normalization} then introduces holes and slots, and proves the Normalization Theorem hinted at earlier in this introduction.
In Part~\ref{part:dents in H-fields} we focus on slots in $H$-fields and their algebraic closures, and implement the approximation arguments for obtaining (strongly) split-normal or repulsive-normal slots.
In Part~\ref{part:Hardy fields univ exp ext} we begin the analytic part of the paper, introducing Hardy fields,
showing that maximal Hardy fields are $\upo$-free, and investigating the universal exponential extensions of Hardy fields.
In the final act (Part~\ref{part:Hardy fields})  we prove our Fixed Point Theorem and give the proof of Theorem~\ref{thm:B}.
We finish with a coda (Part~\ref{part:applications}) consisting of applications, including the proof of Theorem~\ref{thm:A} and the corollaries above.
We refer to the introduction of each part for more details about their respective contents.

\subsection*{Previous work}
Theorem~\ref{thm:A} for $P$ of order $1$ is in~\cite{D}. By \cite{vdH:hfsol} there exists a Hardy field $H\supseteq\R$ isomorphic as an ordered differential field to $\mathbb T_{\text{g}}$, so by \cite{JvdH} this~$H$ has the intermediate value property for all differential polynomials over it.  We announced the $\upo$-freeness of maximal Hardy fields already  in~\cite{ADH2}. 
%~\cite{HIVP}.

\subsection*{Notations and terminology} We freely use the notations and conventions from our book~[ADH], and recall here a few. Throughout, $m$,~$n$ range over the set~$\N=\{0,1,2,\dots\}$.
Given an additively written abelian group $A$ we let $A^{\ne}:=A\setminus\{0\}$.
 Rings (usually, but not always, commutative) are associative with  identity~$1$. For a ring~$R$ we let~$R^\times$ be the multiplicative group of units of $R$ (consisting of the $a\in R$ such that~$ab=ba=1$ for some $b\in R$).

\medskip\noindent
A {\em differential ring\/} is a commutative ring $R$ containing
(an isomorphic copy of) $\Q$ as a subring and equipped with a derivation~$\der\colon R \to R$, in which case $C_R:=\ker\der$   is a subring of $R$, called the ring of constants of $R$, and~$\Q\subseteq C_R$.  A {\em differential field\/} is a differential ring 
$K$ whose underlying ring is a field. In this case  $C_K$ as a subfield of $K$, and if $K$ is understood
from the context we often write $C$ instead of~$C_K$.
An {\em ordered differential field}\/ is an ordered field  equipped with a derivation on its underlying field; such an ordered differential field is
in particular a differential ring. 

Often we are given a differential field $H$ in which $-1$ is not a square, and then~$H[\imag]$ is a differential field extension with
$\imag^2=-1$. Then for $z\in H[\imag]$, $z=a+b\imag$, $a,b \in H$ we set $\Re z:=a$, $\Im z:=b$, and $\bar{z}:=a-b\imag$. Hence
$z\mapsto \bar{z}$ is an automorphism of the differential field $H$.
If in addition there is given a differential field extension $F$ of $H$ in which $-1$ is not a square, we always tacitly arrange $\imag$ to be such that $H[\imag]$ is a differential subfield of the differential field extension $F[\imag]$ of $F$. 

\medskip\noindent
Let $R$ be a differential ring and $a\in R$. When its derivation $\der$ is clear from the context we denote $\der(a),\der^2(a),\dots,\der^n(a),\dots$ by $a', a'',\dots, a^{(n)},\dots$, and if $a\in R^\times$, then~$a^\dagger:=a'/a$ denotes the logarithmic derivative  of $a$, so $(ab)^\dagger=a^\dagger + b^\dagger$ for all~$a,b\in R^\times$.  
We have the differential ring $R\{Y\}=R[Y, Y', Y'',\dots]$ of differential polynomials in a differential indeterminate~$Y$ over $R$. %, and we set $R\{Y\}^{\ne}:= R\{Y\}\setminus \{0\}$. 
Given  
$P=P(Y)\in R\{Y\}$, the smallest~$r\in\N$ such that~$P\in R[Y,Y',\dots, Y^{(r)}]$ is called the order of $P$, denoted by~$r=\order(P)$; if $P$ has order $r$, then $P=\sum_{\i}P_{\i}Y^{\i}$, as in [ADH, 4.2],  with~$\i$ ranging over tuples~$(i_0,\dots,i_r)\in \N^{1+r}$, $Y^{\i}:= Y^{i_0}(Y')^{i_1}\cdots (Y^{(r)})^{i_r}$,     coefficients~$P_{\i}$   in $R$, and $P_{\i}\ne 0$ for only finitely many $\i$. 
For~$P\in R\{Y\}$ and~${a\in R}$ we let~$P_{+a}(Y):=P(a+Y)$ and $P_{\times a}(Y):=P(aY)$ be
the {\it additive conjugate}\/ and the {\it multiplicative conjugate}\/ of $P$ by~$a$, respectively.
%For such $\i$ we set 
%$$|\i|\ :=\ i_0+i_1+ \cdots + i_r, \qquad \|\i\|\ :=\ i_1+2i_2 + \dots + ri_r.$$
%The {\em degree\/} of $P\ne 0$ is 
%$$\deg P:=\max\big\{|\i|:\ P_{\i}\ne 0\big\}\in \N,$$ 
%and the {\em weight} of $P$ is
%$$\wt P:=\max\big\{\|\i\|:\ P_{\i}\ne 0\big\}\in \N.$$
For $\phi\in R^\times$ we also let~$R^{\phi}$ be the {\it compositional conjugate of $R$ by~$\phi$}\/: the differential ring
with the same underlying ring as~$R$ but with derivation~$\phi^{-1}\der$ (usually denoted by $\derdelta$) instead of $\der$. 
We have an $R$-algebra isomorphism~$P\mapsto P^\phi\colon R\{Y\}\to R^\phi\{Y\}$ such that $P^\phi(y)=P(y)$ for all $y\in R$;
see~[ADH, 5.7].

\medskip\noindent
For a field $K$ we have $K^\times=K^{\ne}$, and
a (Krull) valuation on $K$ is a surjective map 
$v\colon K^\times \to \Gamma$ onto an ordered abelian 
group $\Gamma$ (additively written) satisfying the usual laws, and extended to
$v\colon K \to \Gamma_{\infty}:=\Gamma\cup\{\infty\}$ by $v(0)=\infty$,
where the ordering on $\Gamma$ is extended to a total ordering
on $\Gamma_{\infty}$ by $\gamma<\infty$ for all 
$\gamma\in \Gamma$. 
A {\em valued field\/} $K$ is a field (also denoted by $K$) together with a valuation ring~$\mathcal O$ of that field,
and the corresponding valuation $v\colon K^\times \to \Gamma$  on the underlying field is
such that $\mathcal O=\{a\in K:va\geq 0\}$ as explained in [ADH, 3.1].

\medskip
\noindent
Let $K$ be a valued field
with valuation ring $\mathcal O_K$ and valuation $v\colon K^\times \to \Gamma_K$. Then~$\mathcal O_K$ is a local ring  with maximal ideal $\smallo_K=\{a\in K:va>0\}$  and residue field $\res(K)=\mathcal{O}_K/\smallo_K$. 
If $\res(K)$ has characteristic zero, then $K$ is said to be of equicharacteristic zero.
When, as here,  we use the capital $K$ for the valued field under consideration, then
we denote $\Gamma_K$, $\mathcal O_K$, $\smallo_K$, by $\Gamma$, $\mathcal O$, $\smallo$, respectively.
A very handy alternative notation system in connection with the valuation is as follows.  With $a,b$ ranging over $K$, set 
\begin{align*} a\asymp b &\ :\Leftrightarrow\ va =vb, & a\preceq b&\ :\Leftrightarrow\ va\ge vb, & a\prec b &\ :\Leftrightarrow\  va>vb,\\
a\succeq b &\ :\Leftrightarrow\ b \preceq a, &
a\succ b &\ :\Leftrightarrow\ b\prec a, & a\sim b &\ :\Leftrightarrow\ a-b\prec a.
\end{align*}
It is easy to check that if $a\sim b$, then $a, b\ne 0$ and $a\asymp b$, and that
$\sim$ is an equivalence relation on $K^\times$. % let $a^{\sim}$ be the equivalence class of an element $a\in E^\times$ with respect to~$\sim$. 
Given a valued field extension $L$ of $K$, we identify in the usual way~$\res(K)$ with a subfield of $\res(L)$, and $\Gamma$ with an ordered subgroup of~$\Gamma_L$.  We use {\em pc-sequence\/} to abbreviate
{\em pseudocauchy sequence}, and~${a_\rho\leadsto a}$ indicates that~$(a_\rho)$ is a  pc-sequence pseudoconverging to~$a$;
here the $a_{\rho}$ and $a$ lie in a valued field understood from the context, see [ADH, 2.2,~3.2].  

\medskip
\noindent
As in~[ADH],  a {\em valued differential field\/} is a valued field of equicharacteristic zero together with a derivation, generally denoted by $\der$, on the underlying field. (Unlike~\cite{VDF} we do not assume in this definition that $\der$ is continuous with respect to the valuation topology.)  A valued differential field $K$ is said to have {\it small derivation}\/ if $\der\smallo\subseteq\smallo$; then also $\der\mathcal O\subseteq\mathcal O$~ [ADH, 4.4.2], and so $\der$ induces a derivation on $\res(K)$ making the residue morphism~$\mathcal O\to\res(K)$ into a morphism of differential rings. 

\medskip
\noindent
We shall also consider various special classes of valued differential fields introduced in [ADH],
such as the class of {\em asymptotic fields}\/ (and their relatives, {\it $H$-asymptotic fields}\/)
and its subclass of {\it pre-$\d$-valued fields}\/, which in turn contains the class of 
{\it $\d$-valued fields}\/ [ADH, 9.1, 10.1]. (As usual in [ADH], the prefix ``$\d$'' abbreviates ``differential''.) Every asymptotic field $K$ has its associated asymptotic couple $(\Gamma,\psi)$,
where $\psi\colon\Gamma^{\neq}\to\Gamma$ satisfies $\psi(vg)=v(g^\dagger)$ for $g\in K^\times$ with $vg\ne 0$. 
See [ADH, 9.1, 9.2] for more on asymptotic couples, in particular the taxonomy of asymptotic fields introduced
via their asymptotic couples: having a {\it gap,}\/ being {\it grounded,}\/ having  {\it asymptotic integration,}\/ and having
{\it rational asymptotic integration}. 

\medskip
\noindent
An  {\em ordered valued differential field\/} is a valued differential field $K$ equipped with an ordering on $K$ making $K$ an ordered field. An ordered differential field $K$ is called an {\it $H$-field}\/ if
 for all $f\in K$ with if $f\succ 1$ we have $f^\dagger>0$, and
  $\mathcal{O}=C+\smallo$ where 
$\mathcal{O}= \big\{g\in K:\text{$\abs{g} \le c$ for some $c\in C$}\big\}$ and $\smallo$ is the maximal ideal of the convex subring $\mathcal{O}$ of $K$.  Thus $K$ equipped with its valuation ring $\mathcal O$ is an ordered valued
differential field. {\it Pre-$H$-fields}\/ are the ordered valued differential subfields of $H$-fields.
 See [ADH, 10.5] for basic facts about (pre-) $H$-fields. An $H$-field $K$ is said to be {\it Liouville closed}\/
 if $K$ is real closed and for all $f,g\in K$ there exists~$y\in K^\times$ with $y'+fy=g$.
 Every $H$-field extends to  a Liouville closed one; see [ADH, 10.6].
 
 \medskip\noindent
 We alert the reader that in a few places we refer to the Liouville closed $H$-field $\T_{\g}$ of grid-based transseries from
 \cite{JvdH}, which is denoted there by $\T$. Here we  adopt the notation of [ADH] where $\T$ is the larger field of logarithmic-exponential series. 
 
%\marginpar{needs to be extended: recall (pre-) $H$-fields etc.} 
 
\subsection*{Acknowledgements}
Various institutions supported this work during its genesis:
Aschenbrenner and van den Dries thank the Institut Henri Poincar\'e, Paris, for its
hospitality  during
the special semester ``Model Theory, Combinatorics and Valued Fields'' in   2018.
Aschenbrenner   also   thanks the Institut f\"ur Ma\-the\-ma\-ti\-sche Logik und Grundlagenforschung, Universit\"at M\"unster, for     a M\"unster Research Fellowship during the spring of~2019, and he    acknowledges partial support from NSF Research Grant DMS-1700439.
All three authors   received support from the Fields Institute for   Research in Mathematical Sciences, Toronto, during its 
Thematic Program ``Tame Geometry, Transseries and Applications to Analysis and Geometry'', Spring~2022.
 %\marginpar{Joris: do you want to add anything here?}

\endgroup

\newpage  

\part{Preliminaries}\label{part:preliminaries}

\medskip

\noindent
After generalities on linear differential operators and differential polynomials in
Section~\ref{sec:diff ops and diff polys}, we  investigate the group of
logarithmic derivatives in valued differential fields of various kinds (Section~\ref{sec:logder}) and recall how
 iterated logarithmic derivatives can be used to study the asymptotic behavior of differential polynomials over
such valued differential fields for ``large'' arguments  (Section~\ref{sec:it log derivative}). We also
assemble some basic preservation theorems for {\it $\upl$-freeness}\/ and {\it $\upo$-freeness}\/ (Section~\ref{sec:uplupo-freeness})
and continue the study of linear differential operators over $H$-asymptotic fields initiated in [ADH, 5.6, 14.2] (Section~\ref{sec:lindiff}).
In our analysis of the solutions of algebraic differential equations over $H$-asymptotic fields in Part~\ref{part:normalization},
special pc-sequences in the sense of~[ADH, 3.4] play an important role; Section~\ref{sec:special elements} explains why.
A cornerstone of~[ADH] is the concept of {\it newtonianity}\/,  an analogue of henselianity
appropriate for $H$-asymptotic fields with asymptotic integration [ADH, Chapter~14].
Related to this is  {\it differential-henselianity}\/~[ADH, Chapter~7], which makes sense for a broader class of
valued differential fields.
In Sections~\ref{sec:completion d-hens} and~\ref{sec:complements newton} we further explore these notions. Among other things, we study their persistence under taking the completion of a  valued differential field with small derivation (as defined in [ADH, 4.4]).

\section{Linear Differential Operators and Differential Polynomials}\label{sec:diff ops and diff polys}

\noindent
This section gathers miscellaneous facts of a general nature 
about linear differential operators and differential polynomials, sometimes in a valued differential setting. We first discuss splittings and least common left multiples of linear differential operators,  then recall
the complexity and the separant of differential polynomials, and finally deduce some useful estimates for derivatives
of exponential terms.

\subsection*{Splittings}
In this subsection $K$ is a differential field. Let $A\in K[\der]^{\neq}$ be monic of order~$r\geq 1$.
A {\bf splitting of $A$ over $K$} is a tuple $(g_1,\dots, g_r)\in K^r$ such that~$A=(\der-g_1)\cdots(\der-g_r)$. If $(g_1,\dots, g_r)$ is a   splitting of
$A$ over $K$ and $\fn\in K^\times$, then~$(g_1-\fn^\dagger,\dots, g_r-\fn^\dagger)$
is a   splitting of $A_{\ltimes \fn}=\fn^{-1}A\fn$ over $K$. \index{splitting}\index{linear differential operator!splitting}

\medskip
\noindent
Suppose $A=A_1\cdots A_m$ where every $A_i\in K[\der]$ is monic of positive order $r_i$ (so~$r=r_1+\cdots + r_m$). Given any splittings 
$$(g_{11},\dots, g_{1r_1}),\ \dots,\ (g_{m1},\dots,g_{mr_m})$$ of $A_1,\dots, A_m$, respectively, we obtain a splitting 
$$\big(g_{11},\dots, g_{1r_1},\ \dots,\  g_{m1},\dots, g_{mr_m}\big)$$
of $A$ by concatenating the given splittings of $A_1,\dots, A_m$ in the order indicated, and call it a splitting of $A$ {\bf induced} by the factorization $A=A_1\cdots A_m$.  \index{splitting!induced by a factorization}
 For $B\in K[\der]$ of order $r\ge 1$ we have $B=bA$ with $b\in K^{\times}$ and monic
$A\in K[\der]$, and then a {\bf  splitting of $B$ over~$K$\/} is by definition a  splitting of $A$ over $K$.
A splitting of~$B$ over $K$ remains a splitting of $aB$ over $K$, for any $a\in K^\times$. Thus:

\begin{lemma}\label{lem:split and twist} 
If $B\in K[\der]$ has order $r\geq 1$, and $(g_1,\dots, g_r)$ is a splitting of
$B$ over $K$ and $\fn\in K^\times$, then~$(g_1-\fn^\dagger,\dots, g_r-\fn^\dagger)$
is a   splitting of $B_{\ltimes\fn}$ over $K$ and  a splitting of $B\fn$ over $K$.
\end{lemma}

\noindent
From [ADH,~5.1,~5.7] we know that if $A\in K[\der]$ splits over $K$, then for any~$\phi\in K^\times$ the
operator~$A^\phi\in K^\phi[\derdelta]$ splits over $K^\phi$; here is how a splitting of $A$ over $K$ transforms into
a splitting of $A^\phi$ over $K^\phi$: 

\begin{lemma}\label{lem:split and compconj}
Let $\phi\in K^\times$ and $$A\ =\ c(\der-a_1)\cdots(\der-a_r)\quad\text{ with $c\in K^\times$ and $a_1,\dots, a_r\in K$.}$$
Then  in $K^\phi[\derdelta]$ we have
$$A^\phi\ =\ c\phi^r(\derdelta-b_1)\cdots(\derdelta-b_r)\quad\text{ where
$b_j:=\phi^{-1}\big(a_j-(r-j)\phi^\dagger\big)$   \textup{(}$j=1,\dots,r$\textup{)}.}$$
\end{lemma}
\begin{proof}
Induction on $r$. The case $r=0$ being obvious, suppose $r\geq 1$, and set~$B:=(\der-a_2)\cdots(\der-a_r)$.  By inductive hypothesis 
$$B^\phi=\phi^{r-1}(\derdelta-b_2)\cdots(\derdelta-b_r)\quad\text{ where
$b_j:=\phi^{-1}\big(a_j-(r-j)\phi^\dagger\big)$  for $j=2,\dots,r$.}$$
Then
$$A^\phi\ =\ c\phi\,\big(\derdelta-(a_1/\phi)\big)\,B^\phi\ =\ c\phi^r \, \big(\derdelta-(a_1/\phi)\big)_{\ltimes \phi^{r-1}}\, (\derdelta-b_2)\cdots(\derdelta-b_r)$$
with $$\big(\derdelta-(a_1/\phi)\big)_{\ltimes \phi^{r-1}} = \derdelta-(a_1/\phi) + (r-1)\phi^\dagger/\phi$$ by [ADH, p.~243].
\end{proof}

\noindent
A different kind of factorization, see for example~\cite{Polya},
reduces the process of solving the differential equation~$A(y)=0$ to repeated multiplication and integration:

\begin{lemma}\label{lem:Polya fact} 
Let $A\in K[\der]^{\neq}$ be monic of order $r\geq 1$. If $b_1,\dots,b_r\in K^\times$ and
$$A\ =\ b_1\cdots b_{r-1}b_r (\der b_r^{-1})( \der b_{r-1}^{-1})\cdots (\der b_1^{-1}),$$ then $(a_r,\dots,a_1)$, where $a_j:=(b_1\cdots b_j)^\dagger$ for $j=1,\dots,r$,
is a splitting of $A$ over~$K$. Conversely, if $(a_r,\dots,a_1)$ is a splitting of $A$ over $K$ and $b_1,\dots,b_r\in K^\times$ are such that
$b_j^\dagger=a_j-a_{j-1}$ for $j=1,\dots,r$ with $a_0:=0$, then $A$ is as in the display.
%$A=b_1\cdots b_r \der b_r^{-1} \cdots \der b_1^{-1}$.
\end{lemma}

\noindent
This follows easily by induction on $r$.

\subsection*{Real splittings} Let $H$ be a differential field in which
$-1$ is not a square. Then we let $\imag$ denote an element in a differential field extension of $H$ with $\imag^2=-1$, and consider the
differential field $K=H[\imag]$. 
Suppose $A\in H[\der]$ is monic of order $2$ and splits over $K$, so 
$$A\ =\ (\der-f)(\der-g),\qquad f,g\in K.$$ Then 
$$A\ =\ \der^2-(f+g)\der+fg-g',$$ and thus $f\in H$ iff $g\in H$. 
One checks easily that if $g\notin H$, then there are unique~$a,b\in H$ with $b\ne 0$ such that 
$$f\ =\  a-b\imag+b^\dagger, \qquad g\ =\ a+b\imag,$$
and thus $$A\  =\  \der^2-(2a+b^\dagger)\der + a^2+ b^2 -a'+ab^\dagger.$$
Conversely, if $a,b\in H$ and $b\ne 0$, then for $f:=a-b\imag+b^\dagger$ and $g:=a+b\imag$ we have~$(\der -f)(\der-g)\in H[\der]$. 

\medskip
\noindent
Let now $A\in H[\der]$ be monic of order $r\ge 1$.
%If~$r=2$ and~$A$ is irreducible with $A=(\der-f)(\der-g)$, $f,g\in K$, then by [ADH, 5.1.34] there are unique $a,b\in H$ with $b\ne 0$ such that 
%$$f\ =\  a-b\imag+b^\dagger, \qquad g\ =\ a+b\imag. $$

\begin{lemma}\label{hkspl} Suppose $A$ splits over $K$. Then
  $A=A_1\cdots A_m$ for some  $A_1,\dots, A_m$ in $H[\der]$ that are monic and irreducible of order $1$ or $2$ and split over $K$. 
\end{lemma}
\begin{proof} By [ADH, 5.1.35], $A=A_1\cdots A_m$, where every $A_i\in H[\der]$ is monic and irreducible of order $1$ or $2$.
By [ADH, 5.1.22], such $A_i$ split over $K$.  
\end{proof}

\begin{definition}\label{def:real splitting}  \index{splitting!real}\index{linear differential operator!real splitting}
A {\bf real splitting of $A$} (over $K$) is a splitting of~$A$ over $K$ that is induced by a factorization $A=A_1\cdots A_m$ where every $A_i\in H[\der]$ is monic of order $1$ or $2$ and splits over $K$. (Note that we do not require the $A_i$ of order $2$ to be irreducible in $H[\der]$.)
\end{definition}

\noindent
Thus if $A$ splits over $K$, then  $A$ has a real splitting over $K$ by Lemma~\ref{hkspl}. Note that if~$(g_1,\dots, g_r)$ is a real splitting of
$A$ and $\fn\in H^\times$, then~$(g_1-\fn^\dagger,\dots, g_r-\fn^\dagger)$
is a real splitting of $A_{\ltimes \fn}$.  

\medskip
\noindent
It is convenient to extend the above slightly: for $B\in H[\der]$ of order $r\ge 1$ we have~$B=bA$ with $b\in H^{\times}$ and monic
$A\in H[\der]$, and then a {\bf real splitting of~$B$\/}~(over~$K$) is by definition a real splitting of~$A$ (over $K$).  

\medskip
\noindent
In later use, $H$ is a valued differential field with small derivation such that $-1$ is not a square in the differential residue field $\res(H)$. For such $H$, let $\mathcal{O}$ be the valuation ring of $H$. We make $K$ a valued differential field extension of $H$ with small derivation by taking 
$\mathcal{O}_K=\mathcal{O}+\mathcal{O}\imag$ as the valuation ring of $K$. We have the residue map
$a\mapsto \res a\colon \mathcal{O}_K\to \res(K)$, so
 $\res(K)=\res(H)[\imag]$, writing here $\imag$ for~$\res\imag$.
We extend this map to a ring morphism $B\mapsto\res B\colon\mathcal{O}_K[\der]\to\res(K)[\der]$ by sending
$\der\in \mathcal{O}[\der]$ to $\der\in \res(K)[\der]$.  

\begin{lemma}\label{lem:lift real splitting}
Suppose $(g_1,\dots,g_r)\in \res(K)^r$ is a real splitting of a monic
operator $D \in\res(H)[\der]$ of order $r\ge 1$.
Then there are $b_1,\dots,b_r\in\mathcal O_K$ such that 
$$B\ :=\ (\der-b_1)\cdots (\der-b_r)\in \mathcal O[\der],$$
$(b_1,\dots,b_r)$ is a real splitting of $B$, 
and $\res b_j=g_j$ for~$j=1,\dots,r$.
\end{lemma}
\begin{proof}
We can assume  
%$D$ is irreducible over $\res(H)$ and
$r\in\{1,2\}$. The case $r=1$ is obvious, so let $r=2$. Then the case
where $g_1, g_2\in \res(H)$ is again obvious, so let $g_1=\res(a)-\res(b)\imag +(\res b)^\dagger$, $g_2=\res(a)+\res(b)\imag$ where
$a,b\in\mathcal O$, $\res b\neq 0$. Set $b_1:=a-b\imag+b^\dagger$, $b_2:=a+b\imag$. Then $b_1,b_2\in\mathcal O_K$ with~$\res b_1=g_1$,~$\res b_2=g_2$, and $B:=(\der-b_1)(\der-b_2)\in\mathcal O[\der]$ have the desired properties.
%Since $\overline{B}=D$ is irreducible,  $B$ is as well,  by [ADH, 5.6.3].
\end{proof}

\subsection*{Least common left multiples and complex conjugation} 
In this subsection~$H$ is a differential field.
Recall from [ADH,~5.1] the definition of the {\it least common left multiple}\/\index{least common left multiple}\index{linear differential operator!least common left multiple}
$\operatorname{lclm}(A_1,\dots,A_m)$ of operators $A_1,\dots,A_m\in H[\der]^{\neq}$, $m\geq 1$:
this is the monic operator $A\in H[\der]$ such that $H[\der]A_1\cap\cdots \cap H[\der]A_m=H[\der]A$.
%Hence~$\operatorname{lclm}(A_1,\dots,A_m)=\operatorname{lclm}\!\big(\operatorname{lclm}(A_1,\dots,A_{m-1}),A_m\big)$ if $m>1$.
For~${A,B\in H[\der]^{\neq}}$ we have:
$$\max\big\{ \!\order(A),\order(B) \big\} \ \leq\ \order\!\big(\!\operatorname{lclm}(A,B)\big)\ \leq\ \order(A)+\order(B).$$
For the inequality on the right, note that the natural $H[\der]$-module morphism 
$$H[\der]\ \to\ \big(H[\der]/H[\der]A\big)\times \big(H[\der]/H[\der]B\big)$$ has kernel $H[\der]\operatorname{lclm}(A,B)$, and
 for  $D\in H[\der]^{\ne}$,
the $H$-linear space $H[\der]/K[\der]D$ has dimension $\order D$. 

%\begin{lemma} \marginpar{new lemma and corollary}
%If $A,B\in H[\der]^{\neq}$ split over $K$, then so does $\operatorname{lclm}(A,B)$.
%\end{lemma}
%\begin{proof}
%By [ADH, 5.1.25], an operator $D\in H[\der]^{\ne}$ splits over $H$ iff the $H[\der]$-mod\-ule~$M:=H[\der]/DH[\der]$ has a composition series $\{0\}=M_0\subseteq M_1\subseteq\cdots\subseteq M_m=M$ where $\dim_H M_i/M_{i-1}=1$ for $i=1,\dots,m$. The lemma follows from this and the isomorphism above.
%\end{proof}

%\begin{cor}\label{cor:As}
%Let $A\in H[\der]^{\neq}$. There is a unique monic operator $A_{\operatorname{s}}\in H[\der]^{\neq}$ of largest order with $A\in H[\der]A_{\operatorname{s}}$ which splits over $H$. If $A\in H[\der]B$ where $B\in H[\der]$ splits over $H$, then~$A_{\operatorname{s}}\in H[\der]B$.
%\end{cor}
%\begin{proof}
%Take a monic operator $A_{\operatorname{s}}\in H[\der]^{\neq}$ of largest order such that $A\in H[\der]A_{\operatorname{s}}$ and $A_{\operatorname{s}}$ splits over $H$. Let $B\in H[\der]$ split over $H$ with $A\in H[\der]B$, and let $L:=\operatorname{lclm}(A_{\operatorname{s}},B)\in H[\der]$. Then $A\in H[\der]A_{\operatorname{s}}\cap H[\der]B=H[\der]L$ and $L$ splits over $H$, hence~$L=A_{\operatorname{s}}$ and $A_{\operatorname{s}}\in H[\der]B$.
%\end{proof}

\medskip
\noindent
We now assume that $-1$ is not a square in $H$;
then we have a differential field extension $H[\imag]$ where $\imag^2=-1$.
The automorphism $a+b\imag\mapsto \bar{a+b\imag}:= a-b\imag$~(${a,b\in H}$) of the differential field $H[\imag]$ extends uniquely to an automorphism
$A\mapsto \bar{A}$ of the ring $H[\imag][\der]$ with $\bar{\der}=\der$. 
Let $A\in H[\imag][\der]$; then
$\bar{A}=A\Longleftrightarrow A\in H[\der]$. 
Hence if $A\neq 0$ is monic, then~$L:=\operatorname{lclm}(A,\overline{A})\in H[\der]$ and thus
$L= BA = \overline{B}\,\overline{A}$ where $B\in H[\imag][\der]$.

\begin{exampleNumbered}\label{ex:lclm compl conj}
Let $A=\der-a$ where $a\in H[\imag]$. If $a\in H$, then $\operatorname{lclm}(A,\overline{A})=A$, and if~$a\notin H$, then $\operatorname{lclm}(A,\overline{A})=(\der-b)(\der-a)=(\der-\overline{b})(\der-\overline{a})$ where $b\in H[\imag]\setminus H$.
\end{exampleNumbered}

\noindent
Let now $F$ be a differential field extension of $H$ in which $-1$ is not a square; we assume that $\imag$ is an element of a differential ring extension of $F$.  

\begin{lemma}\label{lem:lclm compl conj}
Let $A\in H[\imag][\der]^{\neq}$ be monic, $b\in H[\imag]$, and $f\in F[\imag]$  such that ${A(f)=b}$.
Let $B\in H[\imag][\der]$ be such that $L:=\operatorname{lclm}(A,\overline{A})=BA$. Then $L(f)=B(b)$ and hence
$L\big(\!\Re(f)\big)=\Re\!\big(B(b)\big)$ and  $L\big(\!\Im(f)\big)=\Im\!\big(B(b)\big)$.
\end{lemma}

\noindent
In Sections~\ref{secfhhf} and~\ref{sec:d-alg extensions} we need a slight extension of this lemma:

\begin{remarkNumbered}\label{rem:lclm compl conj}
Let $F$ be a differential ring extension of $H$ in which $-1$ is not a square and let $\imag$ be an element of a commutative ring extension of $F$ such that~$\imag^2=-1$ and the $F$-algebra $F[\imag]=F +  F\imag$ is a free
$F$-module with basis $1$,~$\imag$. For $f=g+h\imag\in F[\imag]$ with $g,h\in F$ we set $\Re(f):= g$ and $\Im(f):= h$. We make $F[\imag]$ into a differential ring extension of~$F$ in the only way possible (which has $\imag'=0$). Then Lemma~\ref{lem:lclm compl conj} goes through. 
\end{remarkNumbered}

%\medskip
%\noindent
%In the following we let $K:=H[\imag]$.  

%\begin{lemma}
%Let $A\in K[\der]^{\neq}$. There is a unique monic operator $A_{\operatorname{r}}\in H[\der]^{\neq}$ of largest order with $A\in K[\der]A_{\operatorname{r}}$. If $A\in K[\der]B$ where $B\in H[\der]$, then~$A_{\operatorname{r}}\in H[\der]B$.
%\end{lemma}
%\begin{proof}
%Take a monic operator $A_{\operatorname{r}}\in H[\der]^{\neq}$ of largest order such that $A\in K[\der]A_{\operatorname{r}}$. Let $B\in H[\der]$ with $A\in K[\der]B$, and let $L:=\operatorname{lclm}(A_{\operatorname{r}},B)\in H[\der]$. Then $A\in K[\der]A_{\operatorname{r}}\cap K[\der]B=K[\der]L$ and hence $L=A_{\operatorname{r}}$ and $A_{\operatorname{r}}\in K[\der]B\cap H[\der]=H[\der]B$, using~[ADH, 5.1.11] for the last equality.
%\end{proof}

%\noindent
%Let $A\in K[\der]^{\neq}$ and  $A_{\operatorname{r}}\in H[\der]^{\neq}$ be as in the previous lemma. We call $A_{\operatorname{r}}$ the {\bf largest real right factor} of $A$.  
%Now also take~$A_{\operatorname{i}}\in K[\der]$ with $A=A_{\operatorname{i}} A_{\operatorname{r}}$. If $A$ is monic, then $A_{\operatorname{i}}=1\Leftrightarrow A_{\operatorname{r}}=A\Leftrightarrow A\in H[\der]$. Hence in general $A_{\operatorname{ri}}=1$ and $A_{\operatorname{rr}}=A_{\operatorname{r}}$, and we also have $A_{\operatorname{ir}}=1$, $A_{\operatorname{ii}}=A_{\operatorname{i}}$.

\subsection*{Complexity and the separant}
We recall some definitions and observations from~[ADH, 4.3].
Let $K$ be a differential field and $P\in K\{Y\}$, $P\notin K$, and set 
$r=\order P$, $s=\deg_{Y^{(r)}} P$, and $t=\deg P$. Then the {\it complexity}\/ of $P$ is the triple
$\cc(P)=(r,s,t)\in\N^3$;  \index{complexity!differential polynomial}\label{p:complexity} \index{differential polynomial!complexity}
%for $P\in K$ we set $\cc(P):=(0,0,0)$.  
we order $\N^3$ lexicographically. Let $a\in K$. Then $\cc(P_{+a})=\cc(P)$, and
$\cc(P_{\times a})=\cc(P)$ if $a\neq 0$.
%Suppose now   $P\notin K$ and set $r:=\order P$.
The differential polynomial~$S_P:=\frac{\partial P}{\partial Y^{(r)}}$ is called the {\it separant}\/ of $P$; \index{separant} \index{differential polynomial!separant}\label{p:separant} thus
 $\cc(S_P)<\cc(P)$ (giving complexity~$(0,0,0)$ to elements of $K$), and $S_{aP}=aS_P$ if $a\neq 0$.  Moreover:

\begin{lemma}\label{lem:separant fms}
We have
$$S_{P_{+a}} = (S_P)_{+a},\quad S_{P_{\times a}}=a\cdot(S_P)_{\times a},\quad \text{and}\quad
S_{P^\phi} = \phi^r (S_P)^\phi\text{ for $\phi\in K^\times$.}$$
\end{lemma}
\begin{proof}
For $S_{P_{+a}}$ and $S_{P_{\times a}}$ this is from [ADH, p.~216]; 
for $S_{P^\phi}$, express $P$ as a polynomial in $Y^{(r)}$ and use $(Y^{(r)})^\phi= \phi^r Y^{(r)}+\text{lower order terms}$. 
\end{proof}

\subsection*{Some transformation formulas}
In this subsection $K$ is a differential ring.
Let~$u\in K^\times$. Then in $K[\der]$ we have
\begin{align*}
(\der-u^\dagger)^0\	&=\ 1,\\
(\der-u^\dagger)^1\	&=\ \der-u'u^{-1},\\
(\der-u^\dagger)^2\	&=\ \der^2 - 2u'u^{-1}\der + \big(2(u')^2 - u''u\big)u^{-2}.
\end{align*}
More generally:

\begin{lemma}\label{lem:def of Qnk}
There are differential polynomials $Q^n_k(X)\in\Q\{X\}$ \textup{(}$0\leq k\leq n$\textup{)}, independent of~$K$ and $u$, such that $Q^n_n=1$ and
$$(\der-u^\dagger)^n\	=\  Q^n_n(u)\der^n+Q^n_{n-1}(u)u^{-1}\der^{n-1}+\cdots+Q^n_0(u)u^{-n}.$$
Setting $Q^n_{-1}:=0$,   we have
$$Q^{n+1}_k(X) \ = \ Q^n_k(X)'X + Q^n_k(X)(k-n-1)X' + Q^n_{k-1}(X)\qquad (0\leq k\leq n).$$
Hence every $Q^n_k$ with $0\le k\le n$ has integer coefficients and is homogeneous of degree $n-k$ and isobaric of weight $n-k$.
\end{lemma}
\begin{proof}
By induction on $n$. The case $n=0$ is obvious. Suppose for a certain $n$  we have $Q^n_k$ for $0\leq k\leq n$ as above.
Then
\begin{align*}
(\der-u^\dagger)^{n+1}\	&=\ (\der-u^\dagger) \sum_{k=0}^n Q^n_k(u)u^{k-n}\der^k \\
						&=\ \sum_{k=0}^n \Big(\big( Q^n_k(u)u^{k-n}\big){}' - u^\dagger Q^n_k(u)u^{k-n} \Big)\der^k  + \sum_{k=0}^n  Q^n_k(u)u^{k-n}\der^{k+1}  \\
						&=\ \sum_{k=0}^n \Big( Q^n_k(u)'u + Q^n_k(u)(k-n-1)u' \Big) u^{k-(n+1)}\der^k + {}\\
						&\ \quad \sum_{k=1}^{n+1} Q^n_{k-1}(u)u^{k-(n+1)}\der^{k},
\end{align*}
and this yields the inductive step.
\end{proof}

\noindent
For $f\in K$ we have
$$(fu^{-1})^{(n)}\ =\ (\der^n)_{\ltimes u^{-1}}(f)u^{-1}\ =\ (\der_{\ltimes u^{-1}})^n(f)u^{-1}\ =\ (\der-u^\dagger)^n(f)u^{-1}$$
and hence:

\begin{cor}\label{cor:Qnk}
Let $f\in K$; then
$$(fu^{-1})^{(n)}\ =\ Q^n_n(u) f^{(n)}u^{-1} + Q^n_{n-1}(u) f^{(n-1)}u^{-2} + \cdots + Q_0^n(u) fu^{-(n+1)}.$$
\end{cor}

\subsection*{Estimates for derivatives of exponential terms} 
In this subsection $K$ is an asymptotic differential field with small derivation, and $\phi\in K$.
We  also fix
$\fm\in K^\times$ with~$\fm\prec 1$.
Recall from [ADH, 4.2]  that for $P\in K\{Y\}^{\neq}$ the {\it multiplicity of~$P$ at~$0$}\/
is $\val(P)=\min\{d\in\N:P_d\neq 0\}$, where $P_d$ denotes the homogeneous part of degree $d$ of $P$.
Here is a useful bound:

\begin{lemma} \label{lem:diff operator at small elt}  
Let $r\in\N$ and $y\in K$ satisfy $y\prec \fm^{r+m}\prec 1$.
Then $P(y) \prec \fm^{m\mu} P$ for all  $P\in K\{Y\}^{\neq}$   of order at most~$r$ with~$\mu=\val(P)\geq 1$.
\end{lemma}
\begin{proof} Note that $0\ne \fm \prec 1$ and $r+m\ge 1$. Hence 
$$y'\ \prec\ (\fm^{r+m})'\ =\ (r+m)\fm^{r+m-1}\fm'\ \prec\ \fm^{r-1+m},$$ so by induction $y^{(i)}\prec \fm^{r-i+m}$ for $i=0,\dots,r$.
Hence $y^{\i} \prec \fm^{(r+m)\abs{\i}-\dabs{\i}} \preceq \fm^{m\abs{\i}}$
for nonzero $\i=(i_0,\dots,i_r)\in\N^{1+r}$, which yields the lemma.
\end{proof}

\begin{cor}\label{cor:kth der of f}
If   $f\in K$ and $f\prec \fm^n$, then  $f^{(k)} \prec \fm^{n-k}$ for $k=0,\dots,n$.
\end{cor}
\begin{proof}
This is a special case of 
 Lemma~\ref{lem:diff operator at small elt}.
 \end{proof}
 
 \begin{cor}\label{cor:kth der of f, preceq}
Let   $f\in K^\times$ and $n\ge 1$ be such that $f\preceq \fm^n$. Then  $f^{(k)} \prec  \fm^{n-k}$ for $k=1,\dots,n$.
\end{cor}
\begin{proof}
Note that   $\fm^n\neq 0$, so $f'\preceq (\fm^n)'=n\fm^{n-1}\fm'\prec\fm^{n-1}$ [ADH, 9.1.3].
Now apply Corollary~\ref{cor:kth der of f} with $f'$, $n-1$ in place of $f$, $n$.
\end{proof}

%\marginpar{material commented out below in this subsection has been checked but is no longer needed}
%\begin{lemma}
%Suppose  $\phi' \preceq \fm^{-1}$; then
 %$\phi^{(\i)}\preceq \phi^{i_0}\fm^{-\dabs{\i}}$  for $\i=(i_0,\dots,i_n)\in\N^{1+n}$.
%\end{lemma}
%\begin{proof}
%It is enough  to show $\phi^{(m)} \preceq \fm^{-m}$ for $m\geq 1$. We proceed by induction on 
%$m\geq 1$.
%The case $m=1$ holds by hypothesis, and for the inductive step note that using
%[ADH, 9.1.3(iii), 9.1.4(i)] and smallness of the derivation of $K$, $\phi^{(m)}\preceq\fm^{-m}\not\asymp 1$  yields  $\phi^{(m+1)}\preceq (\fm^{-m})'=-m\fm'\fm^{-(m+1)} \prec \fm^{-(m+1)}$ as required.
%\end{proof}

%\begin{cor}\label{cor:Enk}
%If $\phi' \preceq \fm^{-1}$, then $E^n_k(\phi) \preceq \fm^{k-n}$ for $k=0,\dots,n$.
%\end{cor}

%\begin{proof}
%Use the facts about $E^n_k$ above and the previous lemma.
%\end{proof}

%\begin{lemma}\label{lem:der of gexphi}
%Suppose $\phi' \preceq \fm^{-1}$, and let   $f\in K$, $f\prec\fm^{m+n}$. Then $$(f\ex^\phi)^{(n)}\ =\ g\ex^\phi \text{ where $g\in K$, $g \prec \fm^m$.}$$
%\end{lemma}
%\begin{proof}
%Set $g=\sum_{k=0}^n E^n_k(\phi)f^{(k)}$; then  $(f\ex^\phi)^{(n)}=g\ex^\phi$.
%By Corollary~\ref{cor:kth der of f} we have~$f^{(k)} \prec \fm^{m+n-k}$ for $k=0,\dots,m+n$.
%For $k=0,\dots,n$ we have  $E^n_k(\phi) \preceq \fm^{k-n}$ 
%by  Corollary~\ref{cor:Enk}, so
%$E^n_k(\phi) f^{(k)} \prec \fm^{k-n} \fm^{m+n-k}=\fm^m$.
%Hence   $g   \prec  \fm^m$.
%\end{proof}

\noindent
In the remainder of this subsection we let $\xi\in K^\times$  and
assume $\xi\succ 1$ and $\zeta:=\xi^\dagger\succeq 1$.   

\begin{lemma}\label{lem:xi zeta} 
{\samepage The elements $\xi,\zeta\in K$ have the following asymptotic properties: 
\begin{enumerate}
\item[\textup{(i)}] $\zeta^n \prec \xi$ for all $n$;
\item[\textup{(ii)}] $\zeta^{(n)} \preceq \zeta^2$ for all $n$.
\end{enumerate}
Thus for each $P\in \mathcal O\{Z\}$ there is an $N\in\N$ with $P(\zeta)\preceq\zeta^N$, and hence $P(\zeta)\prec\xi$.}
\end{lemma}
\begin{proof}
Part (i)   follows from  [ADH, 9.2.10(iv)] for $\gamma=v(\xi)$.
As to (ii), if~$\zeta'\preceq\zeta$, then $\zeta^{(n)}\preceq\zeta$ by [ADH, 4.5.3], and we are done.
Suppose $\zeta'\succ\zeta$ and set~$\gamma:=v(\zeta)$. Then $\gamma,\gamma^\dagger < 0$, so $\gamma^\dagger=o(\gamma)$
by [ADH, 9.2.10(iv)] and hence $v(\zeta^{(n)})=\gamma+n\gamma^\dagger > 2\gamma=v(\zeta^2)$  by
[ADH, 6.4.1(iv)].
\end{proof} 

\noindent
Recall from [ADH, 5.8] that for a homogeneous differential polynomial $P\in K\{Y\}$ of degree $d\in\N$ the {\it Riccati transform}\/ $\Ric(P)\in K\{Z\}$ \label{p:Ric} of $P$
satisfies $$\Ric(P)(z)=P(y)/y^d\quad \text{ for $y\in K^\times$, $z=y^\dagger$.}$$

%\begin{cor}\label{cor:xi 0}
%Let $l\in \Z$ and $\i\in\N^{1+n}$. Then $(\xi^l)^{(\i)} \prec \xi^{l\abs{\i}+1}$.
%\end{cor}
%\begin{proof}
%With $R=\Ric(Y^{(\i)})\in\Q\{Z\}$ we have
%$$(\xi^l)^{(\i)}/\xi^{l\abs{\i}}\ =\ R\big((\xi^l)^\dagger\big)\ =\  R(l\zeta) \prec \xi$$ by Lemma~\ref{lem:xi zeta}.
%\end{proof}

%\noindent
%In the three corollaries below $\xi=\phi'$ and $l\in \Z$.

%\begin{cor}
%$E^n_k(\phi) \sim {n\choose k} \xi^{n-k}$ for $0\leq k\leq n$.
%\end{cor}
%\begin{proof}
%Use \eqref{eq:Enk degree} and Corollary~\ref{cor:xi 0} for $l=1$.
%\end{proof}
\noindent
In the next two corollaries, $l\in \Z$, $\xi=\phi'$, and $\ex^\phi$ denotes a unit of a differential ring extension of $K$ with multiplicative inverse $\ex^{-\phi}$ and such that $(\ex^\phi)'=\phi'\ex^\phi$. 
 
\begin{cor}\label{cor:xi 1}
$(\xi^l\ex^\phi)^{(n)}\ =\ \xi^{l+n}(1+\varepsilon)\ex^\phi$ where $\varepsilon\in K$, $\varepsilon\prec 1$.
\end{cor}
\begin{proof}
By Lemma~\ref{lem:xi zeta}(i) we have $l\zeta+\xi\sim\xi\succ 1$.
Now use $(\xi^l\ex^\phi)^{(n)}/(\xi^l\ex^\phi)=R_n(l\zeta+\xi)$ for $R_n=\Ric(Y^{(n)})$ in combination with~[ADH, 11.1.5].
\end{proof}

\noindent
Applying the corollary above with $\phi$, $\xi$ replaced by~$-\phi$,~$-\xi$, respectively, we obtain:

\begin{cor}\label{cor:xi 2}
$(\xi^l\ex^{-\phi})^{(n)}\ =\ (-1)^n\xi^{l+n}(1+\delta)\ex^{-\phi}$  where $\delta\in K$, $\delta\prec 1$.
\end{cor}

\subsection*{Estimates for Riccati transforms}
In this subsection $K$ is a valued differential field with small derivation. For later use we prove variants of~[ADH, 11.1.5].
%(In loc.~cit.~we assume $K$ is asymptotic, but not here.)

\begin{lemma}\label{Riccatipower}
If $z\in K^{\succ 1}$, then $R_n(z)=z^n(1+ \epsilon)$ with 
$v\epsilon\ge v(z^{-1})+ o(vz)>0$.  
\end{lemma}
\begin{proof} This is clear for $n=0$ and $n=1$. Suppose $z\succ 1$, $n\ge 1$, and
$R_n(z)=z^n(1+ \epsilon)$ with $\epsilon$ as in the lemma. As in the proof of [ADH, 11.1.5],
$$ R_{n+1}(z)\ =\  z^{n+1}\left(1 + \epsilon +    n\frac{z^\dagger}{z}(1+\epsilon) +       
\frac{\epsilon'}{z}\right).$$
Now $v(z^\dagger)\ge o(vz)$: this is obvious if $z^\dagger\preceq 1$, and follows from~$\triangledown(\gamma)=o(\gamma)$ for~$\gamma\ne 0$ if $z^\dagger\succ 1$ [ADH, 6.4.1(iii)]. 
This gives the desired result in view of  $\epsilon'\prec 1$. 
\end{proof}

\begin{lemma}\label{Riccatipower+} Suppose $\der\mathcal O\subseteq\smallo$. If 
$z\in K^{\succeq 1}$, then $R_n(z)=z^n(1+ \epsilon)$ with 
$\epsilon \prec 1$.
\end{lemma}
 
\begin{proof}
The case $z\succ 1$ follows from Lemma~\ref{Riccatipower}. For $z\asymp 1$, proceed as in the proof of that lemma, using $\der\mathcal O\subseteq\smallo$.
\end{proof}

\noindent
By [ADH, 9.1.3 (iv)] the condition $\der\mathcal O\subseteq\smallo$ is satisfied
if $K$ is $\d$-valued, or  asymptotic   with $\Psi\cap \Gamma^{>}\ne \emptyset$.

\begin{lemma}\label{lem:Riccati bd flat}
Suppose $K$ is asymptotic, and $z\in K$ with $0\neq z\preceq z'\prec 1$. Then~$R_n(z)\sim z^{(n-1)}$ for $n\geq 1$.
\end{lemma}
\begin{proof}
Induction on $n$ gives $z\preceq z'\preceq\cdots \preceq z^{(n)}\prec 1$     
 for all $n$.
We now show~$R_n(z)\sim z^{(n-1)}$ for $n\geq 1$, also by induction. The case $n=1$ is clear from $R_1=Z$,
so suppose~${n\geq 1}$ and $R_n(z)\sim z^{(n-1)}$. Then
$$R_{n+1}(z)\ =\ zR_n(z)+R_n(z)'$$
where $R_n(z)'\sim z^{(n)}$ by [ADH, 9.1.4(ii)] and $zR_n(z)\asymp zz^{(n-1)} \prec z^{(n-1)}\preceq z^{(n)}$.
Hence $R_{n+1}(z)\sim z^{(n)}$.
\end{proof}

\subsection*{Valued differential fields with very small derivation\astr} The generalities in this subsection will be used in Section~\ref{sec:embeddings into T}. 
Let $K$ be  a valued differential field with derivation $\der$. Recall that if $K$ has small derivation (that is, $\der\smallo\subseteq \smallo$), then also~${\der\mathcal O\subseteq\mathcal O}$ by [ADH, 4.4.2],  so we have a unique derivation on the residue field~$\k:=\mathcal O/\smallo$ that makes the residue morphism
$\mathcal O\to \k$ into a morphism of differential rings (and we call $\k$ with this induced derivation the differential residue field of $K$).
We say that $\der$ is {\bf very small}\index{very small derivation}\index{valued differential field!very small derivation} if
$\der\mathcal O\subseteq\smallo$. So $K$ has  very small derivation iff~$K$ has small derivation and the induced derivation on $\k$ is trivial.
If  $K$ has small derivation and~$\mathcal O=C+\smallo$, then $K$ has very small derivation.
If $K$ has very small derivation, then so does every valued differential subfield of $K$, and if~$L$ is a valued differential field extension of~$K$ with small derivation and $\k_L=\k$,  
then $L$ has very small derivation. Moreover:

\begin{lemma}\label{lem:very small der alg ext}
Let $L$ be a valued differential field extension of $K$, algebraic over~$K$, and suppose $K$ has very small derivation. Then $L$ also has
very small derivation.
\end{lemma}
\begin{proof}
By [ADH, 6.2.1], $L$ has small derivation. The derivation of $\k$ is trivial and~$\k_L$ is algebraic over $\k$ [ADH, 3.1.9],
so the derivation of $\k_L$ is also trivial.
\end{proof}

\noindent
Next we focus on  pre-$\d$-valued fields with very small derivation. First an easy observation about asymptotic couples:

\begin{lemma}\label{lem:small deriv char}
Let $(\Gamma,\psi)$ be an asymptotic couple;  then  
$$\text{$(\Gamma,\psi)$ has gap $0$}\quad\Longleftrightarrow\quad 
\text{$(\Gamma,\psi)$ has small derivation and $\Psi\subseteq\Gamma^{<}$.}$$
In particular, if $(\Gamma,\psi)$ has small derivation and does not have gap $0$, then each asymptotic couple
extending $(\Gamma,\psi)$ has small derivation.
\end{lemma}

\begin{cor}
Suppose $K$ is pre-$\d$-valued   with   small derivation, and suppose~$0$ is not a gap in $K$. Then $K$ has very small derivation.
\end{cor}
\begin{proof}
The  previous lemma gives $g\in K^\times$ with $g\not\asymp 1$ and $g^\dagger\preceq 1$.
Then for each~$f\in K$ with  $f\preceq 1$ we have $f'\prec g^\dagger\preceq 1$.
\end{proof}

\begin{cor}\label{cor:dv(K) very small der}
Suppose $K$ is pre-$\d$-valued of $H$-type with very small derivation. Then the $\d$-valued hull $\operatorname{dv}(K)$ of $K$ has small derivation.
\end{cor}
\begin{proof}
By Lemma~\ref{lem:small deriv char},
if $0$ is not a gap in $K$, then 
every pre-$\d$-valued field extension of $K$ has small derivation.
If $0$ is a gap in $K$, then no $b\asymp 1$ in~$K$ 
satisfies~$b'\asymp 1$, since $K$ has  very small derivation. Thus $\Gamma_{\operatorname{dv}(K)}=\Gamma$ by 
[ADH,~10.3.2(ii)], so $0$ remains a gap in $\operatorname{dv}(K)$. In both cases, $\operatorname{dv}(K)$ has small derivation.
\end{proof}

\noindent
If $K$ is pre-$\d$-valued and ungrounded, then for each $\phi\in K$ which is active  in $K$,   the pre-$\d$-valued field $K^\phi$ (with derivation $\derdelta=\phi^{-1}\der$)
has very small derivation.

\medskip\noindent
Now a fact about $A\in F[\der]$, where $F$ is any differential field. For the definition of~$A^{(n)}$, see [ADH, p.~243]. Recall that $\operatorname{Ri}(A)\in F\{Z\}$. For $P\in F\{Z\}$  we denote by $P_{[0]}$ the isobaric part of $P$ of weight $0$, as in [ADH, p.~212], so $P\in F[Z]$.

\begin{lemma}\label{lem:RiAn} For $P:= \operatorname{Ri}(A)_{[0]}$ we have
$\operatorname{Ri}(A^{(n)})_{[0]} = P^{(n)}$ for all $n$.
\end{lemma}
\begin{proof}
We  treat the case $n=1$; the general case then follows by induction on~$n$.
By $F$-linearity of $A\mapsto \operatorname{Ri}(A)$ we reduce to the case $A=\der^m$, $m\ge 1$, so $P=Z^m$.  
Then~$A'=m\der^{m-1}$, so $\operatorname{Ri}(A')=mR_{m-1}$ and hence $\operatorname{Ri}(A')_{[0]}=mZ^{m-1}=P'$.
\end{proof}

\noindent
We need this for the next lemma, which in turn is used in proving  Corollary~\ref{cor:Ri unique, 2}. As usual, we extend the residue map $a\mapsto\res a\colon\mathcal O\to\k$
to  the ring morphism $$P\mapsto\res P\ \colon\ \mathcal O[Y]\to\k[Y], \qquad Y\mapsto Y.$$

\begin{lemma}\label{RiPQ} Let $K$ have very small derivation, $A\in \mathcal{O}[\der]$, $R:=\operatorname{Ri}(A)$, so~$R\in \mathcal{O}\{Z\}$, and $P:=R_{[0]}\in \mathcal{O}[Z]$. 
Let $a\in\mathcal O$, so  $Q:=(R_{+a})_{[0]}\in\mathcal O[Z]$.
Then $$(\res P)_{+\res a}\ =\ \res Q.$$
\end{lemma}
\begin{proof}
It suffices to show
$P_{+a}-Q\prec 1$. We have~$R(a)\equiv R_{[0]}(a)\bmod\smallo$, as $K$ has very small derivation.
Applying this to~$\operatorname{Ri}(A^{(n)})$ in place of $R=\operatorname{Ri}(A)$ and using Lemma~\ref{lem:RiAn}
yields~$\operatorname{Ri}(A^{(n)})(a)\equiv P^{(n)}(a)\bmod\smallo$ for all $n$.
Now use $P_{+a}=\sum_n \frac{1}{n!}P^{(n)}(a)\,Z^n$ by Taylor expansion and $R_{+a}=\sum_n \frac{1}{n!}\operatorname{Ri}(A^{(n)})(a)\,R_n$ by [ADH, p.~301],  so
$Q=\sum_n \frac{1}{n!}\operatorname{Ri}(A^{(n)})(a)\,Z^n$. 
\end{proof}

\subsection*{Rosenlicht's proof of a result of Kolchin\astr} Corollary~\ref{cor:Kolchin} below will be used
in Section~\ref{sec:lin diff applications}.
Let $K$ be a differential field  and
$m\geq 1$,  $P\in K\{Y\}$, $\deg P < m$.

\begin{lemma}\label{ympy}
Let $L$ be a differential field extension of $K$ and $t\in L$, $t'\in K + tK$,
and suppose $L$ is algebraic over $K(t)$.
If $y^m=P(y)$ for some $y\in L$, then $z^m=P(z)$ for some $z$ in a differential field extension of $K$
which is algebraic over $K$.
\end{lemma}
\begin{proof}
We may assume $t$ is transcendental over $K$.
View $K(t)$ as a subfield of the Laurent series field $F=K(\!(t^{-1})\!)$.  We equip $F$ with the   valuation ring~$\mathcal O_F=K[[t^{-1}]]$ and the continuous derivation extending that of $K(t)$, cf.~[ADH, p.~226]. Then the valued differential field~$F$ is monotone.
Hence the valued differential subfield~$K(t)$ of $F$ is also monotone. We equip $L$ with a valuation ring $\mathcal O_L$ lying over~$\mathcal O_F\cap K(t)$. Then~$L$ is monotone by [ADH, 6.3.10]. 
We  
identify $K$ with its image under the residue morphism $a\mapsto \res a\colon \mathcal O_L\to\k_L:=\res(L)$. Then $K$ is a differential subfield of the differential residue field $\k_L$ of $L$,
and   $\k_L$ is algebraic over~$K$.
Let now~$y\in L$ with $y^m=P(y)$, and towards a contradiction suppose~$y\succ 1$.
Then $y^\dagger\preceq 1$, thus~$y^{(n)}=y\,R_n(y^\dagger)\preceq y$ for all~$n$, and hence $y^m=P(y)\preceq y^d$ where~$d=\deg P<m$, a contradiction. Thus $y\preceq 1$, and $z:=\res y\in\k_L$ has the required property.
\end{proof}

\noindent
We recall from [ADH, p.~462] that  a {\it Liouville extension}\/ of $K$\index{extension!Liouville} \index{Liouville extension}
is a differential field extension $E$ of $K$ such that $C_E$ is algebraic over $C$ and for each $t\in E$ there are~$t_1,\dots,t_n\in L$ such that~$t\in K(t_1,\dots,t_n)$ and for $i=1,\dots,n$:
 \begin{enumerate}
\item $t_i$ is algebraic over $K(t_1,\dots,t_{i-1})$, or
\item $t_i'\in K(t_1,\dots,t_{i-1})$, or
\item $t_i'\in t_i K(t_1,\dots,t_{i-1})$.
\end{enumerate}
%See [ADH, p.~462] for basic properties of Liouville extensions.

\begin{prop}[{Rosenlicht~\cite[p.~371]{Rosenlicht73}}]
Suppose $y^m=P(y)$ for some $y$ in a Liouville extension of~$K$. Then $z^m=P(z)$ for some $z$ in a differential field extension of $K$ which is algebraic over $K$.
\end{prop}
\begin{proof}
A {\it Liouville sequence}\/ over $K$ is by definition a sequence $(t_1,\dots,t_n)$ of elements of a differential field extension $E$ of $K$ such that  for $i=1,\dots,n$:
\begin{enumerate}
\item $t_i$ is algebraic over $K(t_1,\dots,t_{i-1})$, or
\item $t_i'\in K(t_1,\dots,t_{i-1})$, or
\item $t_i'\in t_i K(t_1,\dots,t_{i-1})$.
\end{enumerate}
Note that then $K_i:=K(t_1,\dots,t_{i})$ is a differential subfield of $E$ for $i=1,\dots,n$.
By induction on $d\in\N$ we now show:  if $(t_1,\dots, t_n)$ is a Liouville sequence over~$K$
with~$\operatorname{trdeg}\!\big(K(t_1,\dots,t_n)|K\big) = d$ and $y^m=P(y)$ for
some $y\in K(t_1,\dots,t_n)$, then   the conclusion of the proposition holds. This is obvious for $d=0$, so let $(t_1,\dots, t_n)$ be a Liouville sequence over~$K$
with $\operatorname{trdeg}\!\big(K(t_1,\dots,t_n)|K\big) = d\ge 1$ and $y^m=P(y)$ for
some $y\in K(t_1,\dots,t_n)$.
Take $i\in\{1,\dots,n\}$ maximal such that $t_i$ is transcendental over~$K_{i-1}=K(t_1,\dots,t_{i-1})$.
Then $t_i'\in  K_{i-1}+t_i K_{i-1}$, and $K(t_1,\dots, t_n)$ is algebraic over $K(t_1,\dots, t_i)$. Applying Lemma~\ref{ympy} 
 to $K_{i-1}$, $t_i$ in the role of $K,t$ yields a~$z$ in an algebraic differential field extension of $K_{i-1}$ with  $z^m=P(z)$. Now apply the inductive hypothesis to the Liouville sequence~$(t_1,\dots,t_{i-1},z)$ over~$K$.
\end{proof}

\begin{cor}[Kolchin]\label{cor:Kolchin}
Let $A\in K[\der]^{\neq}$, and suppose $A(y)=0$ for some $y\neq 0$ in a Liouville extension of $K$.
Then $A(y)=0$ for some $y\neq 0$ in a differential field extension of $K$ with $y^\dagger$  algebraic over $K$.
\end{cor}
\begin{proof}
Let $m=\order A$ and
note that $\operatorname{Ri}(A)=Z^m+Q$ with $\deg Q<m$ [ADH, p.~299]. Now apply  the proposition  above.
\end{proof}

\begin{remark}
Corollary~\ref{cor:Kolchin} goes back to Liouville~\cite{Liouville39} in an analytic setting for $A$ of order $2$ and $K=\C(x)$ with $C=\C$, $x'=1$.
\end{remark}

\subsection*{Results of Srinivasan\astr}  
Corollary~\ref{cor:Srinivasan} and Lemma~\ref{lem:Srinivasan}  below
will be used in Section~\ref{sec:Bessel}.
In this subsection $K$ is a differential field and $a_2,\dots,a_n\in K$, $n\geq 2$, $a_n\neq 0$. We investigate the solutions 
of the differential equation
\begin{equation}\label{eq:Abel equ}
y'\ =\ a_2y^2+a_3y^3+\cdots+a_ny^n
\end{equation}
in Liouville extensions of $K$.
For $n=3$ this is a special case of Abel's differential equation of the first kind,
first studied by Abel~\cite{Abel} (cf.~\cite[\S{}A.4.10]{Kamke}).
In the next three lemmas and in Proposition~\ref{prop:Abel equ} we let $y$ be an
element of a differential field extension $L$ of $K$ satisfying~\eqref{eq:Abel equ}. At various places we consider a field
$E(\!(t)\!)$ of Laurent series over a field $E$; it is to be viewed as a valued field in the usual way. 

\begin{lemma}\label{lem:Abel equ, 1}
Suppose $y$ is   transcendental over $K$.
Then $K\langle y\rangle^\dagger \cap K  =  K^\dagger$. % and if $a_2,a_3\notin\der(K)$, then $\der(K\langle y\rangle)\cap K=\der(K)$.
\end{lemma}
\begin{proof}
We view $K\langle y\rangle=K(y)$ as a differential subfield of $K(\!(y)\!)$ equipped with the unique continuous derivation 
extending that of $K(y)$. Let $f=\sum_{j\ge k} f_jy^j\in K(\!(y)\!)$ with $k\in \Z$, all 
$f_j\in K$, and  $f_k\neq 0$. Then
$$f'=f_k' y^k + (f_{k+1}' + k f_k a_2)y^{k+1}+\big(f_{k+2}'+ k f_k a_3 +  (k+1)f_{k+1}a_2\big)y^{k+2}+\cdots.$$
Hence if $f'=af$ where  $a\in K$,
then~$f_k'=af_k$ and so $a\in K^\dagger$. 
%Next suppose $a_2,a_3\notin\der(K)$, and $f'=b\in K^\times$.  We claim that $f_k'\neq 0$: otherwise  $a_2\notin\der(K)$ yields~$kf_ka_2+f_{k+1}'=b$ and $k=-1$, hence~$k f_k a_3 + f_{k+2}'=0$, contradicting~$a_3\notin\der(K)$. Thus~$f_k'\neq 0$ and so $k=0$ and $b=f_k'\in\der(K)$.
\end{proof}

\begin{lemma}\label{lem:Abel equ, 2}
Suppose $y$   transcendental over $K$ 
and   $R$ is a differential subring of~$K$ with $C\subseteq R=\der(R)$ and~$a_2,\dots,a_n\in R$.
Then   $\der(K\langle y\rangle)\cap K=\der(K)$.
\end{lemma}
\begin{proof}
Let $K(\!(y)\!)$ and $f$ be as in the proof of Lemma~\ref{lem:Abel equ, 1}.
Then $$g:=f'= \sum_{j\geq k}g_jy^j\quad\text{where $g_j=f_j'+\sum_{i=1}^{n-1} (j-i) f_{j-i} a_{i+1}$ and $f_l:=0$ for $l<k$.}$$
Towards a contradiction, suppose $f'=a\in K\setminus \der(K)$. By induction on $j\geq k$ we show that then  $f_j\in R$ and $g_j = 0$.
We have $g_k=f_k'$, and
if $f_k'\neq 0$, then~$k=0$ and~$a=f_k'\in\der(K)$, a contradiction. Therefore~$f_k\in C^\times$ and $g_k=0$.
Suppose~$j\geq k+1$ and $f_k,\dots,f_{j-1}\in R$.
Take $h\in R$ with $h'=\sum_{i=1}^{n-1} (j-i) f_{j-i} a_{i+1}$.
Now~$g_j=(f_j+h)'\neq a$ since $a\notin\der(K)$, hence 
$g_j=0$ and thus $f_j\in -h+C\subseteq R$.
\end{proof}

\begin{lemma}\label{lem:Abel equ, 3}
{\samepage Let  $t\in L^\times$ and suppose $L$ is algebraic over $K(t)$. If
\begin{enumerate}
%\item[\textup{(i)}] $t'\in K\setminus\der(K)$  and $a_2,a_3\notin\der(K)$, or
\item[\textup{(i)}] $t'\in K\setminus\der(K)$  and
there is a differential subring of~$K$ with $C\subseteq R=\der(R)$ and $a_2,\dots,a_n\in R$, or  
\item[\textup{(ii)}] $t^\dagger   \in K \setminus \Q K^\dagger$,
\end{enumerate}
then   $y$ is algebraic over $K$.}
\end{lemma}
\begin{proof}
We may assume that $t$ is transcendental over $K$. Suppose $y$ is transcendental over $K$.
Then $t$ is algebraic over~$K(y)=K\langle y\rangle$. 
If~$t'\in K$, then with $R=K\langle y \rangle$, $r=t'$, and $x=t$ in [ADH, 4.6.10] we obtain $t'\in \der(K\langle y\rangle)$. If
$t^\dagger\in K$, then with~$R=K\langle y \rangle$,  $r=t^\dagger$, and $x=t$ in [ADH, 4.6.11]
we get~$t^\dagger\in \Q K\langle y\rangle^\dagger$. Thus (i) contradicts Lemma~\ref{lem:Abel equ, 2}  and (ii) contradicts Lemma~\ref{lem:Abel equ, 1}.
\end{proof}

\begin{prop}\label{prop:Abel equ}
Suppose $C$ is algebraically closed, $R$ is a differential subring of $K$ with $C\subseteq R=\der(R)$ and $a_2,\dots,a_n\in R$, and $L$ is a Liouville extension of $K$. 
Then $y$ is algebraic over~$K$.
\end{prop}
\begin{proof}
By induction on $d\in\N$ we show: if $(t_1,\dots,t_m)\in L^m$ is a Liouville sequence over~$K$ with
$\operatorname{trdeg}\big(K(t_1,\dots,t_m)|K\big)=d$ and $y\in K(t_1,\dots,t_m)$, then $y$ is algebraic over $K$.
This is clear for $d=0$, so let $(t_1,\dots,t_m)\in L^m$ be a Liouville sequence over $K$ with
$\operatorname{trdeg}\big(K(t_1,\dots,t_m)|K\big)=d\geq 1$ and $y\in K(t_1,\dots,t_m)$. Take $i\in\{1,\dots,m\}$ maximal such that~$t_i$ is transcendental over
$K_{i-1}:=K(t_1,\dots,t_{i-1})$. Then $t_i'\in K_{i-1}\setminus\der(K_{i-1})$ or~$t_i^\dagger\in K_{i-1}\setminus\Q K_{i-1}^\dagger$.
Hence $y$ is algebraic over $K_{i-1}$ by Lemma~\ref{lem:Abel equ, 3} applied to~$K_{i-1}$,~$t_i$, $K(t_1,\dots, t_m)$
in the role of $K$, $t$, $L$.  Now apply
the (tacit)  inductive hypothesis to the Liouville sequence~$(t_1,\dots,t_{i-1},y)$ over~$K$.
\end{proof}

\noindent
In the remainder of this subsection $C$ is algebraically closed, $x\in K$, $x'=1$  (so~$x$
is transcendental over $C$), and  $a_2,\dots,a_n\in C[x]$.
Applying Proposition~\ref{prop:Abel equ}  with~$C(x)$,~$C[x]$ in place of
$K$, $R$, respectively, yields  \cite[Proposition~4.1]{Srinivasan20} with a shorter proof:

\begin{cor}[Srinivasan]\label{cor:Srinivasan}
Any $y$ in any Liouville extension of $C(x)$ satisfying~\eqref{eq:Abel equ} is  algebraic over $C(x)$. 
\end{cor}

\noindent
We now assume $a_2,\dots,a_n\in C$ and  put $P:=a_2Y^2+\cdots+a_nY^n\in C[Y]$.
We equip $C(Y)$ with the derivation that is trivial on $C$ and satisfies $Y'=1$.
Thus the field isomorphism~$C(x)\to C(Y)$ over~$C$ with~$x\mapsto Y$ is an isomorphism between the differential subfield $C(x)$ of~$K$
and~$C(Y)$.
Next a special case of \cite[Proposition~3.1]{Srinivasan17}:

\begin{lemma}[Srinivasan]\label{lem:Srinivasan}
The following are equivalent:
\begin{enumerate}
\item[\textup{(i)}] there is a non-constant $y$ in a differential field extension of $C(x)$ such that~$y$ is algebraic over $C(x)$ and $y'=P(y)$;
\item[\textup{(ii)}] there exists $Q\in C(Y)$ such that $Q'=1/P$. % or $R^\dagger=c/(a_2x^2+\cdots+a_nx^n)$ for some $c\in C^\times$.
\end{enumerate}
\end{lemma}
\begin{proof}
Let $y\notin C$ be algebraic over $C(x)$ with $y'=P(y)$.
 Then~$y$ is transcendental over $C$, hence
$x$ is algebraic over $C(y)$ and so $x\in C(y)$ by~[ADH, 4.6.10] applied to $R=C(y)$. Take~$Q\in C(Y)$ such that $x=Q(y)$. Then $1=x'=Q(y)'= Q'(y)P(y)$
and thus $Q'=1/P$.
This shows (i)~$\Rightarrow$~(ii). Conversely, let
$Q\in C(Y)$ be such that~$Q'=1/P$, and let $Q=A/B$ with $A,B\in C[Y]$, $B\ne 0$. By [ADH, 4.6.14] we have $y$ in a differential field extension of $C(x)$ with constant field $C$ such that~$y'=P(y)$ and $B(y)\ne 0$. 
Then~$Q(y)'=Q'(y)y'=(1/P(y))P(y)=1$, so~$Q(y)=x+c\in C(y)$ with $c\in C$. Then $y\notin C$, hence~$y$ is transcendental over $C$ and $y$ is algebraic over $C(x)$.
\end{proof}

\noindent
Here are two applications of Lemma~\ref{lem:Srinivasan}. In the proofs we
extend the derivation of $C(Y)$ to the  continuous $C$-linear derivation on~$C(\!(Y)\!)$ with   $Y'=1$. 

\begin{cor}\label{cor:Abel equ, 2}
Suppose $n\ge 3$,  $a_2,a_3\neq 0$, and $y$ in a Liouville extension of~$C(x)$ satisfies $y'=P(y)$. Then $y\in C$.
\end{cor}
\begin{proof}
In $C(\!(Y)\!)$ we have $1/P=(1/a_2)Y^{-2}-(a_3/a_2^2)Y^{-1}+\cdots$ and hence
$Q'\neq 1/P$ for all $Q\in C(\!(Y)\!)$, so $y\in C$  by Corollary~\ref{cor:Srinivasan} and Lemma~\ref{lem:Srinivasan}.
\end{proof}

\begin{cor}
Suppose $P$ has a simple zero in $C$ and $y$ in a Liouville extension of $C(x)$ satisfies $y'=P(y)$. Then $y\in C$.
\end{cor}
\begin{proof}
Let $c\in C$ with $P(c)=0$, $P'(c)\neq 0$.  
Then in $C(\!(Y)\!)$ we have $1/P(Y+c)\in aY^{-1}+C[[Y]]$ where $a\in C^\times$, hence
$Q'\neq 1/P$ for all $Q\in C(\!(Y)\!)$.  Thus $y\in C$  by Corollary~\ref{cor:Srinivasan} and Lemma~\ref{lem:Srinivasan}. 
\end{proof}

\section{The Group of Logarithmic Derivatives} \label{sec:logder}

\noindent
Let $K$ be a differential field.
The map $y\mapsto y^\dagger \colon K^\times\to K$ is a morphism 
from the multiplicative group of $K$ to the additive group of $K$, with kernel~$C^\times$.
Its image 
$$(K^\times)^\dagger\ =\ \big\{y^\dagger:\,y\in K^\times\big\}$$
is an additive  subgroup of  $K$, which we call the {\bf group of logarithmic derivatives} \index{group of logarithmic derivatives}\label{p:Kdagger}
of $K$. The morphism $y\mapsto y^\dagger$ induces an isomorphism $K^\times/C^\times\to (K^\times)^\dagger$. To shorten notation, set $0^\dagger:=0$, so $K^\dagger=(K^\times)^\dagger$. 
For $\phi\in K^\times$ we have $\phi(K^\phi)^\dagger=  K^\dagger$.
The group $K^\times$ is divisible iff both $C^\times$ and $K^\dagger$ are divisible. 
If $K$ is algebraically
closed, then $K^\times$ and hence
$K^\dagger$ are divisible, making
$K^\dagger$ a  $\Q$-linear subspace of~$K$.
Likewise, if $K$ is real closed, then
the multiplicative subgroup~$K^{>}$ of $K^\times$ is divisible, so
$K^\dagger=(K^{>})^\dagger$ is a $\Q$-linear subspace of $K$.

\begin{lemma}\label{lem:Kdagger alg closure}
Suppose $K^\dagger$ is divisible, $L$ is a differential field extension of $K$ with~$L^\dagger\cap K=K^\dagger$, and
$M$ is a differential field extension of $L$ and algebraic over~$L$. Then $M^\dagger\cap K=K^\dagger$.
\end{lemma}
\begin{proof}
Let $f\in M^\times$ be such that $f^\dagger\in K$. Then $f^\dagger\in L$, so for $n:=[L(f):L]$, 
$$nf^\dagger\ =\ \operatorname{tr}_{L(f)|L}(f^\dagger)\ =\ \operatorname{N}_{L(f)|L}(f)^\dagger  \in L^\dagger$$ by 
an identity in [ADH, 4.4].
Hence $nf^\dagger\in K^\dagger$, and thus $f^\dagger\in K^\dagger$. 
\end{proof}

\noindent
In particular, if $K^\dagger$ is divisible and $M$ is  a differential field extension of $K$ and algebraic over~$K$, then
$M^\dagger\cap K=K^\dagger$.

\medskip
\noindent
In the next two lemmas~$a,b\in K$; distinguishing  whether or not $a\in K^\dagger$ helps to  describe 
the solutions to the differential equation $y'+ay=b$:

\begin{lemma}\label{0K} Suppose $\der K = K$, and let $L$ be differential field extension of $K$ with~$C_L=C$. 
Suppose also $a\in K^\dagger$. Then for some  $y_0\in K^\times$ and
$y_1\in K$,
$$\{y\in L:\, y'+ay=b\}\ =\ \{y\in K:\, y'+ay=b\}\ =\ Cy_0 + y_1.$$
\end{lemma}
\begin{proof} 
Take $y_0\in K^\times$ with $y_0^\dagger=-a$, so $y_0'+ay_0=0$. Twisting $\der+a\in K[\der]$ by~$y_0$ (see~[ADH, p.~243]) transforms the equation $y'+ay=b$ into $z'=y_0^{-1}b$. This gives $y_1\in K$ with $y_1'+ay_1=b$.
Using $C_L=C$, these  $y_0, y_1$ have the desired properties. 
\end{proof} 

\begin{lemma}\label{1K} Let $L$ be a differential field
extension of $K$ with $L^\dagger\cap K=K^\dagger$. Assume
$a\notin K^\dagger$. Then there is at most one $y\in L$ with $y'+ay=b$.
\end{lemma}
\begin{proof} If $y_1$, $y_2$ are distinct solutions in $L$ of the equation
$y'+ay=b$, then we have~$-a=(y_1-y_2)^\dagger\in L^\dagger\cap K=K^\dagger$, contradicting $a\notin K^\dagger$. 
\end{proof} 

\subsection*{Logarithmic derivatives under algebraic closure}
{\it In this subsection $K$ is a differential field.}\/
We describe for real closed $K$ how $K^\dagger$ changes if we pass from~$K$ to its algebraic closure. More generally, suppose the underlying field of $K$ is euclidean; in particular, $-1$ is not a square in~$K$. 
We equip~$K$ with the unique ordering making~$K$ an ordered field.
For $y=a+b\imag \in K[\imag]$ ($a,b\in K$) we let $\abs{y}\in K^{\geq}$ be such that $\abs{y}^2=a^2+b^2$.
Then $y\mapsto \abs{y}\colon K[\imag]\to K^{\geq}$ is an absolute value on $K[\imag]$, i.e., for all $x,y\in K[\imag]$,
$$\abs{x}\ =\ 0\ \Longleftrightarrow\ x=0, \qquad \abs{xy}\ =\ \abs{x}\abs{y}, \qquad \abs{x+y}\ \leq\ \abs{x}+\abs{y}.$$
For $a\in K$ we have $\abs{a}=\max\{a,-a\}$.
We have the subgroup
\[%\begin{equation}\label{eq:S}
S\ :=\ \big\{ y\in K[\imag]: \abs{y}=1\big\}\ =\ \big\{ a+b\imag : a,b\in K,\ a^2+b^2=1 \big\}
\]%\end{equation}
of the multiplicative group $K[\imag]^\times$.
By an easy computation all elements of $K[\imag]$ are squares in $K[\imag]$; hence $K[\imag]^\dagger$ is $2$-divisible.
The next lemma describes $K[\imag]^\dagger$; it partly generalizes~[ADH, 10.7.8]. 

\begin{lemma}\label{lem:logder}
We have $K[\imag]^\times = K^>\cdot S$ with $K^>\cap S=\{1\}$, and
$$K[\imag]^\dagger\ =\  K^\dagger \oplus S^\dagger\qquad\text{\textup{(}internal direct sum of subgroups of $K[\imag]^\dagger$\textup{)}.}$$
For $a,b\in K$ with $a+b\imag\in S$ we have
$(a+b\imag)^\dagger = \wr(a,b)\imag$. Thus $K[\imag]^\dagger\cap K=K^\dagger$. 
\end{lemma}
\begin{proof}
Let $y=a+b\imag\in K[\imag]^\times$ ($a,b\in K$), and take $r\in K^>$ with $r^2=a^2+b^2$; then~$y=r\cdot (y/r)$ with $y/r\in S$. Thus $K[\imag]^\times = K^>\cdot S$, and clearly $K^>\cap S=\{1\}$.
Hence~$K[\imag]^\dagger = K^\dagger + S^\dagger$.
Suppose $a\in K^\times$, $s\in S$, and $a^\dagger=s^\dagger$; then $a=cs$ with~${c\in C_{K[\imag]}}$, and $C_{K[\imag]}=C[\imag]$ by [ADH,~4.6.20] and hence $\max\{a,-a\}=\abs{a}=\abs{c}\in C$, so $a\in C$ and thus $a^\dagger=s^\dagger=0$; therefore the sum is direct.
Now if $a,b\in K$ and~$\abs{a+b\imag}=1$, then
\begin{align*}
(a+b\imag)^\dagger\	&=\ (a'+b'\imag)(a-b\imag) \\
					&=\ (aa'+bb')+(ab'-a'b)\imag \\ 
					&=\ \textstyle\frac{1}{2}\big( a^2+b^2 \big)' + (ab'-a'b)\imag\ =\ (ab'-a'b)\imag\ =\ \wr(a,b)\imag. \qedhere
\end{align*}				
\end{proof}

\begin{cor}\label{cor:logder abs value}
For $y\in K[\imag]^\times$ we have $\Re(y^\dagger)=\abs{y}^\dagger$, and the group morphism~$y\mapsto  \Re(y^\dagger) \colon K[\imag]^\times\to K$ has kernel $C^> S$.
\end{cor}

\noindent
If $K$ is real closed and $\mathcal{O}$ a convex valuation ring of $K$, 
then $\mathcal{O}[\imag]= \mathcal{O}+ \mathcal{O}\imag$ is the unique valuation ring of $K[\imag]$ that lies over $\mathcal{O}$, and so
$S\subseteq\mathcal O[\imag]^\times$, hence
$y\asymp\abs{y}$ for all $y\in K[\imag]^\times$.
Thus by [ADH, 10.5.2(i)] and Corollary~\ref{cor:logder abs value}:

\begin{cor}\label{cor:10.5.2 variant}
If $K$ is a real closed pre-$H$-field, then for all $y,z\in K[\imag]^\times$,
$$y\prec z \quad\Longrightarrow\quad \Re(y^\dagger) < \Re(z^\dagger).$$
\end{cor}

\noindent
We also have a useful decomposition for $S$:

\begin{cor} \label{cor:decomp of S}
Suppose $K$ is a real closed $H$-field. Then $$S\ =\ S_C\cdot \big(S\cap(1+\smallo_{K[\imag]})\big)$$ where 
$S_C := S\cap C[\imag]^\times$ and $S\cap(1+\smallo_{K[\imag]})$ are subgroups of $\mathcal O[\imag]^\times$.  
\end{cor}  
\begin{proof} The inclusion $\supseteq$ is clear. For the reverse inclusion, let $a,b\in K$, $a^2+b^2=1$ and take the unique $c,d\in C$ with $a-c\prec 1$ and $b-d\prec 1$. Then $c^2+d^2=1$ and~$a+b\imag\sim c+d\imag$, and so $(a+b\imag)/(c+d\imag)\in S\cap (1+\smallo_{K[\imag]})$.
% hence for $y,z\in K[\imag]^\times$ we have $y\sim z$ iff $\abs{y}\sim\abs{z}$. (No!!) 
\end{proof}

\subsection*{Logarithmic derivatives in asymptotic fields}
{\em Let $K$ be an asymptotic field.}\/
If $K$  is   henselian and $\k:=\res K$, then by [ADH, remark before 3.3.33],
$K^\times$ is divisible iff the groups $\k^\times$ and $\Gamma$ are both divisible.
  %and then  $K^\dagger$ is divisible. 
Recall that in~[ADH,~14.2] we defined the 
$\mathcal O$-submodule \label{p:I(K)}
$$\I(K)\ =\ \{y\in K:\, \text{$y\preceq f'$ for some $f\in\mathcal O$}\}$$ 
of $K$. We have  $\der\mathcal O\subseteq \I(K)$, hence $(1+\smallo)^\dagger\subseteq(\mathcal O^\times)^\dagger\subseteq\I(K)$.
One easily verifies:

\begin{lemma}\label{pldv}
Suppose $K$ is pre-$\d$-valued. If $\I(K)\subseteq \der K$, then
$\I(K)=\der\mathcal O$. If~$\I(K)\subseteq K^\dagger$, then 
$\I(K)= (\mathcal O^\times)^\dagger$, with $\I(K)=(1+\smallo)^\dagger$
 if $K$   is $\d$-valued.
\end{lemma}

\noindent
If $K$ is $\d$-valued or $K$ is pre-$\d$-valued without a gap, then
$$\I(K)\ =\ \{y\in K:\text{$y\preceq f'$ for some $f\in\smallo$}\}.$$ 
For $\phi\in K^\times$ we have $\phi \I(K^\phi)=\I(K)$.
If  $K$ has asymptotic integration and $L$ is an asymptotic extension of~$K$, then $\I(K)=\I(L)\cap K$.
The following is [ADH, 14.2.5]:
 
\begin{lemma}\label{lem:ADH 14.2.5}
If $K$ is $H$-asymptotic, has asymptotic integration, and is $1$-linearly newtonian, then it is $\d$-valued and $\der\mathcal O = \I(K) = (1 + \smallo)^\dagger$.
\end{lemma}

\noindent
We now turn our attention to the condition $\I(K)\subseteq K^\dagger$.
If $\I(K)\subseteq K^\dagger$, then also~$\I(K^\phi)\subseteq (K^\phi)^\dagger$ for $\phi\in K^\times$, where 
$$(K^\phi)^\dagger\ :=\ \{\phi^{-1}f'/f:\, f\in K^\times\}\ =\ \phi^{-1}K^\dagger.$$
By [ADH,  Section~9.5 and 10.4.3]: 

\begin{lemma}\label{lem:achieve I(K) subseteq Kdagger} 
Let $K$ be of $H$-type. If $K$ is $\d$-valued, or  pre-$\d$-valued without a gap, then $K$ has an immediate henselian asymp\-to\-tic extension~$L$  with $\I(L)\subseteq L^\dagger$. 
%If $K$ is $\d$-valued of $H$-type with $\Gamma\ne \{0\}$, then $K$ has an immediate henselian $H$-asymp\-to\-tic extension~$L$  such that $\I(L)\subseteq L^\dagger$.   
\end{lemma}

%If $K$ is  $H$-asymptotic with $\Gamma\ne \{0\}$, and $K$ is $\d$-valued or pre-$\d$-valued without a gap, then $K$ has an immediate henselian $H$-asymp\-to\-tic extension~$L$  such that $\I(L)\subseteq L^\dagger$.   \marginpar{generalized}
%% that is closed under small exponential integration.
%\end{lemma}

\begin{cor}\label{cor:no new LDs}
Suppose $K$ has asymptotic integration. Let $L$ be an  asymptotic field extension of $K$ such that $L^\times = K^\times C_L^\times (1+\smallo_L)$.
Then $L^\dagger=K^\dagger+(1+\smallo_L)^\dagger$, and
if $\I(K)\subseteq K^\dagger$, then  $L^\dagger\cap K=K^\dagger$. 
\end{cor}
\begin{proof}
Let $f\in L^\times$, and 
take $b\in K^\times$, $c\in C_L^\times$, $g\in\smallo_L$
with $f=bc(1+g)$; then~$f^\dagger=b^\dagger+(1+g)^\dagger$, showing $L^\dagger=K^\dagger+(1+\smallo_L)^\dagger$.
Next, suppose $\I(K)\subseteq K^\dagger$, let 
$b$, $c$, $f$, $g$ be as before, and assume $a:=f^\dagger\in K$;
then $$a-b^\dagger \in (1+\smallo_L)^\dagger \cap K\ \subseteq\ \I(L)\cap K\ =\  \I(K)\ \subseteq\ K^\dagger$$
and hence  $a\in K^\dagger$. This shows $L^\dagger\cap K=K^\dagger$.
\end{proof}

\noindent
Two cases where the assumption on $L$ in Corollary~\ref{cor:no new LDs} is satisfied: (1)  $L$ is an immediate asymptotic field extension of $K$,  because then  $L^\times = K^\times ({1+\smallo_L})$;  and (2)
$L$ is a $\d$-valued field extension of $K$ with $\Gamma=\Gamma_L$.

\medskip\noindent
If $F$ is a henselian valued field 
of residue characteristic $0$, then clearly the subgroup~$1+\smallo_F$ of $F^\times$ is divisible. 
Hence, if $K$ and $L$ are as in Corollary~\ref{cor:no new LDs} and in addition $K^\dagger$ is divisible and $L$ is henselian, then $L^\dagger$ is divisible.

%\noindent
%Suppose $K$ is pre-$\d$-valued. Then for  $y\in K^\times$, the coset $y^\dagger+\big(K^\dagger\cap\I(K)\big)$ of the subgroup $K^\dagger\cap \I(K)$ of $K^\dagger$ only depends on $vy$, and we have a group isomorphism $$vy\mapsto y^\dagger+\big(K^\dagger\cap\I(K)\big)\colon\Gamma\overset{\cong}\longrightarrow K^\dagger/\big(K^\dagger\cap\I(K)\big).$$

\begin{exampleNumbered}\label{ex:Kdagger} 
Let $C$ be a field of characteristic $0$ and $Q$ be a subgroup of $\Q$ with~$1\in Q$. The Hahn field $C(\!( t^{Q} )\!)=C[[x^Q]]$, with $x=t^{-1}$, is given the natural derivation with $c'=0$ for all~$c\in C$ and~$x'=1$: this derivation is defined by 
$$\bigg(\sum_{q\in Q}c_qx^q\bigg)'\ :=\ \sum_{q\in Q} qc_qx^{q-1}\qquad \text{(all $c_q\in C$)}.$$
Then~$C(\!( t^{Q} )\!)$ has constant field $C$, and is $\d$-valued of $H$-type. Thus $K:=C(\!( t^{Q} )\!)$ satisfies $\I(K)\subseteq K^\dagger$
by Lemma~\ref{lem:achieve I(K) subseteq Kdagger}. Hence by Lemma~\ref{pldv},
$$\I(K)\ =\ (1+\smallo)^\dagger\ =\ \big\{f\in K:\, f\prec x^\dagger = t\big\}\ = \ \smallo\, t.$$ 
It follows easily that
$K^\dagger=Q t\oplus\I(K)$ (internal direct sum of subgroups of $K^\dagger$) and   
thus $(K^t)^\dagger=Q\oplus\smallo\subseteq\mathcal O$. In particular, if $Q=\Z$ (so $K=C(\!( t )\!)$), then~$(K^t)^\dagger=\Z\oplus tC[[t]]$. Moreover, if $L:=\operatorname{P}(C)\subseteq C(\!( t^{\Q} )\!)$ is the differential field of Puiseux series over~$C$, then
$(L^t)^\dagger=\Q\oplus\smallo_L$.
\end{exampleNumbered}

\noindent
In the next three corollaries we continue with the $\d$-valued Hahn field  $K=C(\!(t^Q)\!)$ from the example above. So $CK^\dagger=Ct\oplus \I(K)$ (internal direct sum of $C$-linear subspaces of $K$) where $\I(K)=\smallo t$, hence $CK^\dagger=\mathcal Ot$.
For~$f=\sum_{q\in Q} f_q x^q\in K$ (all~$f_q\in C$) we have the ``residue''~$f_{-1}$ of $f$, and we observe that
$f\mapsto f_{-1}\colon K\to C$ is $C$-linear with kernel~$\der(K)$. Thus:

\begin{cor}\label{cor:derK cap CKdagger}
 $\der(K)\cap C K^\dagger=\I(K)$.
\end{cor}

\noindent
This  yields a fact needed in Section~\ref{sec:Bessel}:

\begin{cor}\label{cor:linindep logs}
Let $F:=C(x)\subseteq K$. Then $\der(F)\cap CF^\dagger=\{0\}$.
\end{cor}
\begin{proof}
We arrange that $C$ is algebraically closed.
Let $f\in \der(F)\cap CF^\dagger$. Then $f\in \I(K)=\smallo t$ by Corollary~\ref{cor:derK cap CKdagger},
so it suffices to show $f\in C[x]$.
For~$c\in C$, let~$v_c\colon F^\times\to\Z$ be the valuation on $F$ with $v(C^\times)=\{0\}$ and~${v(x-c)=1}$. 
 Then~$v_c=v\circ\sigma_c$ where
$\sigma_c$ is  the $C$-linear automorphism of the  field $F$ with~${x\mapsto c+t}$.
Hence it suffices that $\sigma_c(f)\preceq 1$ for all $c\in C$.
For $c\in C$, $g\in F$ we have~$\sigma_c(g)' = -t^2 \sigma_c(g')$, so
$-t^2 \sigma_c(f)\in \der(F)\cap CF^\dagger\subseteq\smallo t$, hence
$\sigma_c(f)\prec x$ and thus~$\sigma_c(f)\preceq 1$.
\end{proof}

\noindent 
For the  next corollary, compare \cite[p.~14]{Hardy16} and 
think of $y_j$ as $\log (x-c_j)$.

\begin{cor}[Linear independence of logarithms]
Let $c_1,\dots,c_n\in C$ be distinct, and let
$y_1,\dots,y_n$ in a common differential field extension of $C(x)$ %with constant field $C$ 
be such that
$y_j'=(x-c_j)^{-1}$ for $j=1,\dots,n$. Then for all $a_1,\dots, a_n\in C$, 
   $$a_1y_1+\cdots+a_ny_n\in C(x)\ \Longrightarrow\ a_1=\cdots=a_n=0.$$
\end{cor}
\begin{proof}
Set $F:=C(x)$ and suppose $a_1,\dots, a_n\in C$ and $f:=a_1y_1+\cdots+a_ny_n\in F$. Then $f'=a_1{(x-c_1)^{-1}}+\cdots+a_n(x-c_n)^{-1}\in \der(F)\cap CF^\dagger$, so by Corollary~\ref{cor:linindep logs},
$$a_1(x-c_1)^{-1}+\cdots+a_n{(x-c_n)^{-1}}\ =\ 0.$$
% by the previous corollary. 
Multiplying both sides of this equality by $\prod_{j=1}^n (x-c_j)$ and substituting $c_j$ for $x$ yields $a_j=0$, for $j=1,\dots,n$.
\end{proof}

\subsection*{The real closed case} {\em In this subsection $H$ is a real closed asymptotic field whose valuation ring $\mathcal{O}$ is convex with respect to the ordering of $H$}. (In later use $H$ is often a Hardy field, which is why we use the letter $H$ here.) The valuation ring of the asymptotic field extension $K=H[\imag]$ of $H$ is then $\mathcal{O}_K=\mathcal{O}+\mathcal{O}\imag$, from which we obtain $\I(K)=\I(H) \oplus \I(H)\imag$. Let 
$$S\ :=\ \big\{y\in K:\, |y|=1\big\}, \qquad
W\ :=\ \big\{\!\wr(a,b):\, a,b\in H,\ a^2+b^2=1\big\},$$ 
so $S$ is a subgroup of $\mathcal{O}_K^\times$ with $S^\dagger = W\imag$ and 
$K^\dagger=H^\dagger\oplus W\imag$ by Lemma~\ref{lem:logder}.
Since~$\der\mathcal O\subseteq \I(H)$,  
we have $W\subseteq \I(H)$, and thus: $W = \I(H) \ \Longleftrightarrow\  \I(H)\imag\subseteq K^\dagger$.

\begin{lemma}\label{lem:W and I(F)} The following are equivalent: \begin{enumerate}
\item[\textup{(i)}] $\I(K)\subseteq K^\dagger$;
\item[\textup{(ii)}] $W=\I(H)\subseteq H^\dagger$.
\end{enumerate}
%Therefore, if $K$ is pre-$\d$-valued and $W=\I(K)$ then $\I(K[\imag])=(\mathcal O_{K[\imag]}^\times)^\dagger$.
\end{lemma}
\begin{proof}
Assume (i). 
Then $\I(H)\,\imag\subseteq\I(K)\subseteq K^\dagger$, so $W=\I(H)$ by the equivalence preceding the lemma. Also
$\I(H)\subseteq  \I(K)$ and 
$K^\dagger\cap H=H^\dagger$ (by Lemma~\ref{lem:logder}), hence $\I(H)\subseteq H^\dagger$, so
(ii) holds. 
For the converse, assume (ii). Then \[ \I(K)\ =\ \I(H)\oplus\I(H)\imag\ \subseteq\ H^\dagger\oplus W\imag\ =\ K^\dagger.\qedhere \]
\end{proof}

\noindent
Applying now Lemma~\ref{lem:ADH 14.2.5} we obtain:

\begin{cor}\label{cor:logder}
If $H$ is $H$-asymptotic and has asymptotic integration, and~$K$ is $1$-linearly newtonian, then $K$ is $\d$-valued and $\I(K)\subseteq K^\dagger$; in particular, $W=\I(H)$.
\end{cor}

\begin{cor}\label{cor:logderset ext} 
Suppose $H$  has asymptotic integration and  $W=\I(H)$. 
Let $F$ be a real closed asymptotic extension of $H$ whose valuation ring is convex.  
Then $$F[\imag]^\dagger\cap K\ =\ (F^\dagger\cap H)\oplus\I(H)\imag.$$  If in addition $H^\dagger=H$, then $F[\imag]^\dagger\cap K=H\oplus\I(H)\imag=K^\dagger$.
\end{cor}
\begin{proof} 
We have 
$$F^\dagger\cap H\subseteq F[\imag]^\dagger\cap K\quad\text{ and }\quad
\I(H)\imag = W\imag \subseteq K^\dagger\cap H\imag \subseteq F[\imag]^\dagger\cap K,$$
so~$(F^\dagger\cap H)\oplus\I(H)\imag \subseteq F[\imag]^\dagger\cap K$. For the reverse inclusion, 
$F[\imag]^\dagger=F^\dagger\oplus W_F\imag$, with 
$$W_F\ :=\ \big\{\!\operatorname{wr}(a,b):\, a,b\in F,\ a^2+b^2=1\big\}\ \subseteq\ \I(F),$$ 
hence 
\begin{align*}
F[\imag]^\dagger\cap K &\, =\, (F^\dagger\cap H)\oplus (W_F\cap H)\imag \\
&\, \subseteq\, (F^\dagger\cap H)\oplus \big(\!\I(F)\cap H\big)\imag  \, =\, (F^\dagger\cap H)\oplus\I(H)\imag,
\end{align*} 
using $\I(F)\cap H=\I(H)$, a consequence of $H$ having asymptotic integration. 
If~$H^\dagger = H$ then clearly $F^\dagger\cap H=H$, hence $F[\imag]^\dagger\cap K=K^\dagger$.
\end{proof} 

\subsection*{Trigonometric closure} In this subsection $H$ is a real closed $H$-field. Let $\O$ be its valuation ring and $\smallo$ the maximal ideal of $\O$.
The algebraic closure $K= H[\imag]$ of $H$ is a $\d$-valued $H$-asymptotic extension with valuation ring $\O_K=\O+\O\imag$.
We have the ``complex conjugation'' automorphism $z=a+b\imag\mapsto \bar{z}=a-b\imag$ ($a,b\in H$) of the valued differential field $K$. For such $z$, $a$, $b$ we have 
 $$ |z|\ =\ \sqrt{z\bar{z}}\ =\ \sqrt{a^2+b^2}\ \in\  H^{\ge}.$$ 
 
 \begin{lemma}\label{lemtri1} Suppose $\theta\in H$ and $\theta'\imag\in K^\dagger$. Then $\theta'\in \der\smallo$, and there is a unique~$y\sim 1$ in $K$ such that $y^\dagger=\theta'\imag$. For this $y$ we have $|y|=1$, so $y^{-1}=\bar{y}$.
\end{lemma} 
\begin{proof} From $\theta'\imag \in K^\dagger$ we get $\theta'\in W\subseteq \I(H)$, so $\theta\preceq 1$, hence
$\theta'\in \der\O=\der\smallo$.  Let~$z\in K^\times$ and $z^\dagger=\theta'\imag$. Then  $\Re z^\dagger=0$, so by Corollaries~\ref{cor:logder abs value} and \ref{cor:decomp of S} we have~$z=cy$ with $c\in C_K^\times$ and $y\in S\cap(1+\smallo_K)$ where $S=\{a\in K:\ |a|=1\}$. Hence~$y\sim 1$, $|y|=1$, and $y^\dagger=\theta'\imag$. If also $y_1\in  K$ and $y_1\sim 1$, $y_1^\dagger=\theta'\imag$, then~$y_1=c_1y$ with~$c_1\in C_K^\times$, so $c_1=1$ in view of $y\sim y_1$. 
\end{proof} 

\noindent
By [ADH, 10.4.3], if $y$ in an $H$-asymptotic extension $L$ of $K$ satisfies
$y\sim 1$ and~${y^\dagger\in \der\smallo_K}$, then the asymptotic field $K(y)\subseteq L$ is an immediate extension of $K$, and so is any algebraic asymptotic extension of $K(y)$.

\medskip\noindent
Call $H$ {\bf trigonometrically closed} if
for all $\theta\prec 1$ in $H$ there is a (necessarily unique)  $y\in K$ such that $y\sim 1$ and $y^\dagger=\theta'\imag$.\index{trigonometric!closed}\index{closed!trigonometrically} (By convention ``trigonometrically closed'' includes ``real closed''.) For such $\theta$ and $y$ we think of
$y$ as $\ex^{\imag\theta}$ and accordingly of the elements $\frac{y+\bar{y}}{2}=\frac{y+y^{-1}}{2}$ and  $\frac{y-\bar{y}}{2\imag}=\frac{y-y^{-1}}{2\imag}$ of $H$ as $\cos \theta$ and $\sin \theta$; this explains the terminology. By Lemma~\ref{lemtri1} the restrictions $\theta\prec 1$ and $y\sim 1$ are harmless. Our aim in this subsection is to construct a canonical trigonometric closure
of $H$.

Note that if  $\I(K)\subseteq K^\dagger$, then $H$ is trigonometrically closed. As a partial converse, if
$\I(H)\subseteq H^\dagger\cap \der H$ and $H$ is trigonometrically closed, then $\I(K)\subseteq K^\dagger$; this is an
easy consequence of $\I(K)=\I(H)+\I(H)\imag$. Thus for Liouville closed $H$ we have:
$$H \text{ is trigonometrically closed}\ \Longleftrightarrow\  \I(K)\subseteq K^\dagger.$$ 
Note also that for trigonometrically closed $H$ there is no $y$ in any $H$-asymptotic extension of $K$ such that
$y\notin K$, $y\sim 1$, and $y^\dagger\in (\der\smallo)\imag$.

If $H$ is Schwarz closed, then $H$ is trigonometrically closed by the next lemma: 

\begin{lemma}\label{lem:sc=>tc}
Suppose $H$ is Liouville closed and $\omega(H)$ is downward closed. Then~$H$ is trigonometrically closed.
\end{lemma}
\begin{proof}
Let $0\ne \theta\prec 1$ in $H$. By Lemma~\ref{lemtri1} it suffices to show that then~${\theta'\imag\in K^\dagger}$.
Note that $h:=\theta'\in\I(H)^{\neq}$; we arrange $h>0$. Now
$$f\ :=\ \omega(-h^\dagger)+4h^2\ =\ \sigma(2h), \qquad 2h\in H^>\cap\I(H),$$ 
hence $2h\in H^>\setminus\Upg(H)$
by [ADH, 11.8.19]. So $f\in\omega(H)^\downarrow=\omega(H)$ by [ADH, 11.8.31],  and thus
$\dim_{C_H} \ker 4\der^2+f\ge 1$ by [ADH, p.~ 258].   
Put $A:=\der^2-h^\dagger \der+h^2\in H[\der]$.
The isomorphism $y\mapsto y\sqrt{h}\colon\ker({4\der^2+f})\to\ker A$ of $C_H$-linear spaces~[ADH, 5.1.13] then yields an element of $\ker^{\neq} A$ that for suggestiveness we denote by $\cos \theta$. 
Put~$\sin\theta:=-(\cos\theta)'/h$.
Then
\begin{align*}
(\sin\theta)'\	&=\  -(\cos\theta)''/h+(\cos\theta)'h^\dagger/h\\
	&=\  \big( {-h^\dagger(\cos\theta)'+h^2\cos\theta } \big)/h+(\cos\theta)'h^\dagger/h\ =\ h\cos\theta
\end{align*}
and thus $y^\dagger=\theta'\imag$ for $y:=\cos\theta+\imag\sin\theta\in K^\times$.
\end{proof}

\noindent
If $H$ is $H$-closed, then $H$ is Schwarz closed by [ADH, 14.2.20], and thus trigonometrically closed.  Using also Lemma~\ref{lem:W and I(F)}  and remarks preceding it  this yields:
% in conjunction with Lemmas~\ref{lem:W and I(F)} and~\ref{lem:sc=>tc} and remarks before these lem we obtain: 

\begin{cor}\label{cor:sc=>tc}
If $H$ is   $H$-closed,  then $\I(K) \subseteq K^\dagger = H\oplus\I(H)\imag$.
\end{cor}

\noindent
Suppose now that $H$ is {\it not}\/ trigonometrically closed; so we have $\theta\prec 1$ in $H$ with~$\theta'\imag\notin K^\dagger$. 
Then [ADH, 10.4.3] provides an immediate asymptotic extension~$K(y)$ of $K$ with $y\sim 1$ and
$y^\dagger=\theta'\imag$. To simplify notation and for suggestiveness we set
$$\cos \theta\ :=\ \frac{y+y^{-1}}{2}, \qquad \sin \theta\ :=\ \frac{y-y^{-1}}{2\imag},$$ so $y=\cos\theta+ \imag\sin \theta$ and
$(\cos\theta)^2+ (\sin\theta)^2=1$.  Moreover $(\cos \theta)'=-\theta'\sin \theta$ and~$(\sin \theta)'=\theta'\cos \theta$. It follows that
$H^+:=H(\cos\theta, \sin\theta)$ is a differential subfield of~$K(y)$ with $K(y)=H^+[\imag]$, and  thus $H^+$, as a valued differential
subfield of $H(y)$, is an asymptotic extension of $H$. 

\begin{lemma}\label{lemtri2}$H^+$ is an immediate extension of $H$.
\end{lemma}
\begin{proof} Since $(y^{-1})^\dagger=-\theta'\imag$, the uniqueness property stated in [ADH, 10.4.3] allows us to extend the complex conjugation automorphism
of $K$ (which is the identity on~$H$ and sends $\imag$ to $-\imag$) to 
 an automorphism $\sigma$ of the valued differential field $K(y)$ such that $\sigma(y)=y^{-1}$.
Then $\sigma(\cos \theta)=\cos\theta$ and $\sigma(\sin \theta)=\sin\theta$, so $H^+= \text{Fix}(\sigma)$.
Let~$\k$ be the residue field of $H$; so $\k[\res\imag]$ is the residue field of $K$ and of its immediate extension $K(y)$.
Now $\sigma(\O_{K(y)})=\O_{K(y)}$, so $\sigma$ induces an automorphism of this residue field $\k[\res\imag]$
which is the identity on $\k$ and sends $\res\imag$ to $-\res \imag$. Hence~$\res\imag$ does not lie in the residue field of $H^+$,   so this residue field is just~$\k$. 
\end{proof} 

\noindent
Equip $H^+$ with the unique field ordering making it an ordered field extension of $H$ in which $\O_{H^+}$ is convex; see 
[ADH, 10.5.8].  Then $H^+$ is an $H$-field, and its real closure is an immediate real closed $H$-field extension of $H$.

\begin{lemma}\label{lemtri3} The $H$-field $H^+$ embeds uniquely over $H$ into any trigonometrically
closed $H$-field extension of $H$. 
\end{lemma} 
\begin{proof} Let $H^*$ be a trigonometrically closed $H$-field extension of $H$. Take the unique $z\sim 1$ in $H^*$ such that $z^\dagger=\theta'\imag$. Then any $H$-field embedding $H^+\to H^*$ over $H$ extends to a valued differential field embedding $H^+[\imag]=K(y)\to H^*[\imag]$ sending $\imag\in K$ to $\imag\in H^*[\imag]$, and this extension
must send $y$ to $z$. Hence there is at most one $H$-field embedding $H^+\to H^*$ over $H$. 
For the existence of such an embedding, the uniqueness properties from [ADH, 10.4.3] yield a valued differential field embedding $K(y)\to H^*[\imag]$ over $H$ sending $\imag\in K$ to $\imag\in H^*[\imag]$ and $y$ to $z$. This embedding maps $H^+$ into $H^*$. The uniqueness property of the ordering
on $H^+$ shows that this embedding restricts to an $H$-field embedding $H^+\to H^*$. 
\end{proof}

\noindent
By iterating the extension step that leads from $H$ to $H^{+}$, alternating it with taking real closures, and taking unions at
limit stages we obtain:

\begin{prop}\label{protrig} $H$ has a trigonometrically closed $H$-field extension $H^{\trig}$
that embeds uniquely over $H$ into any trigonometrically closed $H$-field extension of $H$. 
\end{prop} 

\noindent
This is an easy consequence of Lemma~\ref{lemtri3}. Note that the universal property stated in Proposition~\ref{protrig} determines $H^{\trig}$ up-to-unique-isomorphism of $H$-fields over $H$. We refer to such $H^{\trig}$ as the {\bf trigonometric closure} of $H$.\index{trigonometric!closure}\index{closure!trigonometric}\label{p:Htrig} Note that~$H^{\operatorname{trig}}$ is an immediate extension of $H$, by Lemma~\ref{lemtri2},
and that $H^{\operatorname{trig}}[\imag]$ is a Liouville extension of $K$ and thus of $H$.

\medskip\noindent
A {\em trigonometric  extension\/}\index{extension!trigonometric}\index{trigonometric!extension} of $H$ is a real closed $H$-field extension $E$ of $H$ such that for all $a\in E$ there are
real closed $H$-subfields $H_0\subseteq H_1\subseteq \cdots \subseteq H_n$ of $E$ such that \begin{enumerate}
\item $H_0=H$ and $a\in H_n$;
\item for $j=0,\dots, n-1$ there are $\theta_j\in H_j$ and $y_j\in H_{j+1}[\imag]\subseteq E[\imag]$ such that~$y_j\sim 1$, $\theta_j'\imag=y_j^\dagger$, and $H_{j+1}[\imag]$ is algebraic over $H_j[\imag](y_j)$.
\end{enumerate}
If $E$ is a trigonometric extension of $H$, then $E$ is an immediate extension of $H$ and $E[\imag]$ is an immediate Liouville extension of $K$ and thus of $H$.
The next lemma states some further easy consequences of the definition above: 

\begin{lemma} If $E$ is a trigonometric extension of $H$, then $E$  is a trigonometric extension of any 
real closed $H$-subfield $F\supseteq H$ of $E$.
If $H$ is trigonometrically closed, then $H$ has no proper trigonometric extension. 
\end{lemma} 

\noindent
Induction on $m$ shows that if $E$ is a trigonometric extension of $H$, then for any $a_1,\dots, a_m\in E$ there are real closed $H$-subfields
$H_0\subseteq H_1\subseteq \cdots \subseteq H_n$ of $E$ such that $H_0=H$, $a_1,\dots, a_m\in H_n$ and (2) above holds. This helps in proving:

\begin{cor} A trigonometric extension of a trigonometric extension of $H$ is a trigonometric extension of $H$, and $H^{\trig}$ is a trigonometric extension of $H$. 
\end{cor}

\subsection*{Asymptotic fields of Hardy type} 
Let $(\Gamma,\psi)$ be an asymptotic couple, $\Psi:=\psi(\Gamma^{\neq})$, and let~$\gamma$,~$\delta$ range over $\Gamma$.
Recall that
$[\gamma]$
denotes the archimedean class of~$\gamma$ [ADH, 2.4]. 
Following  \cite[Section~3]{Rosenlicht81} we say that   $(\Gamma,\psi)$ is of {\bf Hardy type} if 
for all $\gamma,\delta\neq 0$ we have
$[\gamma] \leq [\delta]  \Longleftrightarrow  \psi(\gamma)\geq\psi(\delta)$. \index{Hardy field!asymptotic couple}\index{Hardy type} \index{asymptotic couple!of Hardy type}
Note that then $(\Gamma,\psi)$ is of $H$-type, and~$\psi$ induces an order-reversing bijection~$[\Gamma^{\neq}]\to\Psi$.
If $\Gamma$ is archimedean, then $(\Gamma,\psi)$ is of Hardy type. If  $(\Gamma,\psi)$  is of Hardy type, then so is
 $(\Gamma,\psi+\delta)$ for each~$\delta$. 
We also say that an asymptotic field is of Hardy type if its asymptotic couple is. Every asymptotic subfield 
and every compositional conjugate of 
an asymptotic field of Hardy type is also of Hardy type. 
Moreover, every Hardy field is of Hardy type~[ADH, 9.1.11].
Let now~$\Delta$ be a convex  subgroup of $\Gamma$. Note that 
$\Delta$ contains the archimedean class~$[\delta]$ of each~$\delta\in\Delta$. Hence, if
 $\delta\in\Delta^{\neq}$ and $\gamma\notin\Delta$, then  $[\delta] < [\gamma]$ and thus:

\begin{lemma}\label{lem:Hardy type}
If $(\Gamma,\psi)$ is of Hardy type and $\gamma\notin\Delta$, $\delta\in\Delta^{\neq}$,  then
$\psi(\gamma)<\psi(\delta)$.
\end{lemma}

\begin{cor}\label{cor:psi(gamma)<0} 
Suppose $(\Gamma,\psi)$ is of Hardy type with small derivation, $\gamma,\delta\neq 0$, $\psi(\delta) \leq 0$, and $[\gamma']>[\delta]$.  Then $\psi(\gamma)<\psi(\delta)$. 
\end{cor}
\begin{proof}
Let $\Delta$ be the smallest convex subgroup of $\Gamma$ with $\delta\in\Delta$; then $\gamma'\notin\Delta$, and~$\psi(\delta)\in\Delta$
by [ADH, 9.2.10(iv)]. Thus $\gamma\notin\Delta$ by [ADH, 9.2.25].
\end{proof}

\noindent
In \cite[Section~7]{AvdD3} we say that an $H$-field $H$ is {\it closed under powers}\/ if for all $c\in C$ and~$f\in H^\times$ there is a $y\in H^\times$
with $y^\dagger=cf^\dagger$. (Think of $y$ as $f^c$.) \index{closed!under powers} Thus if $H$ is Liouville closed, then $H$ is closed under powers.
{\it In the rest of this subsection we let~$H$ be an $H$-field closed under powers, with asymptotic couple $(\Gamma,\psi)$ and constant field~$C$.} %We let~$f$,~$g$ range over $H$  and~$c$  over~$C$.
We recall some basic facts from \cite[Section~7]{AvdD3}.
First, we can make the value group $\Gamma$ into
an ordered vector space over the constant field $C$:

\begin{lemma}\label{OrderedVectorSpace}
For all $c\in C$ and $\gamma=vf$ with $f\in H^\times$ and each $y\in H^{\times}$
with~$y^\dagger = cf^\dagger$, the element $vy\in \Gamma$ only depends on 
$(c,\gamma)$ \textup{(}not on the choice of $f$ and $y$\textup{)}, and is denoted by $c\cdot \gamma$. The scalar multiplication
$(c,\gamma) \mapsto c \cdot\gamma  \colon C\times \Gamma \to \Gamma$
makes~$\Gamma$ into an ordered vector space over the ordered field $C$.
\end{lemma}

\noindent
Let $G$ be an ordered vector space over the ordered field $C$. From [ADH, 2.4] recall that
the $C$-archimedean class  of $a\in G$ is defined
as  $$[a]_C := \big\{ b\in G:\, 
\textstyle\frac{1}{c}|a|\leq |b|\leq c |a| \text{ for some }c\in C^{>}\big\}.$$
Thus if $C=\Q$, then $[a]_\Q$ is just the archimedean class $[a]$ of $a\in G$.
Moreover, if $C^*$ is an ordered subfield of $C$, then $[a]_{C^*}\subseteq [a]_C$ for each $a\in G$,
with equality if $C^*$ is cofinal in $C$.
Hence  if
$C$ is archimedean, then $[a]=[a]_C$ for all $a\in G$.
Put~$[G]_C:=\bigl\{[a]_C:a\in G\bigr\}$ and linearly order $[G]_C$ by
$$[a]_C<[b]_C \quad :\Longleftrightarrow\quad
\text{$[a]_C\neq [b]_C$ and $|a|<|b|$.}$$ 
%Then $G$ is said to be a {\it Hahn space}\/ if for all   $a,b\in G^{\neq}$ we have $$[a]_C=[b]_C \quad\Rightarrow\quad  [a-c b]_C<[b]_C \text{ for some $c$}.$$
Thus $[G]_C$ has smallest element $[0]_C=\{0\}$. We also set $[G^{\neq}]_C:=[G]_C\setminus\big\{[0]_C\big\}$.
From \cite[Proposition~7.5]{AvdD3} we have: 

\begin{prop}\label{prop:Hahn space property}
%The ordered vector space $\Gamma=v(H^\times)$ over $C$ is a Hahn space,
For all $\gamma,\delta\neq 0$ we have
$$[\gamma]_C \leq [\delta]_C \quad\Longleftrightarrow
\quad \psi(\gamma)\geq\psi(\delta).$$ 
Hence $\psi$ induces an order-reversing bijection $[\Gamma^{\neq}]_C\to\Psi=\psi(\Gamma^{\neq})$.
\end{prop}

\noindent
Proposition~\ref{prop:Hahn space property} yields:

\begin{cor}\label{cor:Hardy type arch C}
$H$ is of Hardy type $\Longleftrightarrow$  $[\gamma]=[\gamma]_C$ for all $\gamma$.
Hence if   $C$ is ar\-chi\-me\-dean, then $H$ is of Hardy type; if $\Gamma\neq\{0\}$, then the converse also holds.
\end{cor}

%\noindent
%Every ordered field arises as the constant field of  an  $H$-field of Hardy type.
%(Take~$C$ to be an ordered field in Example~\ref{ex:Kdagger}.) Conversely, we have:

%\begin{cor}\label{cor:Hardy type}
%Every pre-$H$-field with archimedean ordered residue field is of Har\-dy type.
%\end{cor}
%\begin{proof}
%Use the previous corollary and the fact that every pre-$H$-field $F$ extends to   Liouville closed $H$-field  whose
% residue field is the real closure of $\res(F)$.
%\end{proof}

\section{The Valuation of Differential Polynomials at Infinity\astr} \label{sec:it log derivative} 

\noindent
Our   goal in this work is to   solve certain kinds of algebraic differential equations in Hardy fields. 
In this section we review some general facts about the asymptotic behavior of solutions of algebraic differential equations
in $H$-asymptotic fields.
We will not need these results in order to achieve our main objective,  but they will be used at a few points for
applications and corollaries; see Section~\ref{sec:upper lower bds} and Corollary~\ref{cor:odd degree}.
{\it Throughout this section $K$ is an $H$-asymptotic field, and $f$, $g$ range over $K$.}\/

\subsection*{Iterated logarithmic derivatives}
Let $(\Gamma,\psi)$ be an $H$-asymp\-to\-tic couple. As usual we introduce a new symbol $\infty\notin\Gamma$,
extend the ordering of~$\Gamma$ to an ordering on~$\Gamma_\infty=\Gamma\cup\{\infty\}$ such that $\infty>\Gamma$,
%and declare $\gamma+\infty:=\infty$ for $\gamma\in\Gamma_\infty$ and $0\,\infty:=0$, $k\,\infty:=\infty$ for $k\in\Z^{\neq}$.
and extend $\psi\colon\Gamma^{\neq}\to\Gamma$ to a map $\Gamma_\infty\to\Gamma_\infty$
by setting $\psi(0):=\psi(\infty):=\infty$. (See [ADH, 6.5].) We let~$\gamma$ range over~$\Gamma$,
and we define $\gamma^{\langle n\rangle}\in\Gamma_\infty$ inductively by $\gamma^{\langle 0\rangle}:=\gamma$
and $\gamma^{\langle n+1\rangle}:=
%(\gamma^{\langle n\rangle})^\dagger=
\psi(\gamma^{\langle n\rangle})$.
The following is~\cite[Lem\-ma~5.2]{A}; for the convenience of the reader we include a
proof:

\begin{lemma}\label{lem:itpsi}
Suppose that $0\in (\Gamma^<)'$, $\gamma\neq 0$,  and $n\geq 1$. If $\gamma^{\langle n\rangle}<0$, then~$\gamma^{\langle i\rangle}<0$ for $i=1,\dots,n$ and
$[\gamma] > [\gamma^\dagger] >\cdots > [\gamma^{\langle n-1\rangle}] >  [\gamma^{\langle n\rangle}]$.
\end{lemma}
\begin{proof}
By [ADH, 9.2.9], $(\Gamma,\psi)$ has small derivation, hence the case $n=1$
follows from [ADH, 9.2.10(iv)].
Assume inductively that the lemma holds for a certain value of $n\geq 1$, and suppose $\gamma^{\langle n+1\rangle}<0$. Then $\gamma^{\langle n\rangle}\neq 0$, so we can apply the case~$n=1$ to
$\gamma^{\langle n\rangle}$ instead of $\gamma$ and get $[\gamma^{\langle n\rangle} ]>
 [\gamma^{\langle n+1 \rangle}]$. By the inductive assumption the remaining
inequalities will follow from 
$\gamma^{\langle n\rangle}<0$. 
From $0\in (\Gamma^<)'$ we obtain an element~$1$ of~$\Gamma^>$ with $0=(-1)'=-1+1^\dagger$. Suppose $\gamma^{\langle n\rangle}\ge 0$. 
Then $\gamma^{\langle n\rangle}\in \Psi$, thus~$0<\gamma^{\langle n\rangle}<1+1^\dagger=1+1$ and so
$[\gamma^{\langle n\rangle}]\leq [1]$.
Hence $0> \gamma^{\langle n+1\rangle} \geq 1^\dagger=1$, 
a contradiction.
\end{proof}
 
\noindent
Suppose now that $(\Gamma,\psi)$ is the asymptotic couple of $K$. If  $y\in K^\times$  and
$(vy)^{\langle n \rangle}\ne \infty$, then
the $n$th iterated logarithmic derivative~$y^{\langle n\rangle}$  of $y$ is defined (see  [ADH, 4.2]), and
 $v(y^{\langle n\rangle})= (vy)^{\langle n\rangle}\in\Gamma$. 
 Recall from [ADH, p.~383] that for $f,g\neq 0$,
$$f\flatter g\  :\Leftrightarrow\ f^\dagger\prec g^\dagger,\quad f\flattereq g\ :\Leftrightarrow\ f^\dagger\preceq g^\dagger, \quad f\comp g\ :\Leftrightarrow\ f^\dagger\asymp g^\dagger,$$ hence, assuming also $f, g\nasymp 1$,
$$f\flatter g\  \Rightarrow\ [vf] < [vg], \qquad[vf] \leq [vg] \ \Rightarrow \  f\flattereq g.$$ 
{\it In the rest of this  section we are given $x\succ 1$ in~$K$ with $x'\asymp 1$.}\/
Then $0\in (\Gamma^<)'$, so from  the previous lemma we obtain:

\begin{cor}\label{cor:itpsi}
If $y\in K^\times$,  $y\nasymp 1$, $n\ge 1$, and   $(vy)^{\langle n\rangle}<0$, then~$y^{\langle i\rangle}\succ 1$ for~$i=1,\dots,n$ and $[vy] > \big[v(y^\dagger)\big] > \dots > \big[v( y^{\langle n-1\rangle})\big] > \big[v(y^{\langle n \rangle})\big]$.
\end{cor}

\noindent
Let $\i=(i_0,\dots,i_n)\in\Z^{1+n}$ and $y\in K^\times$ be such that $y^{\langle n\rangle}$ is defined;   we   put 
$$y^{\langle\i\rangle}\ :=\ (y^{\langle 0\rangle})^{i_0}\cdots (y^{\langle n\rangle})^{i_n}\in K.$$ 
If $y^{\langle n\rangle}\ne 0$, then $\i\mapsto y^{\langle\i\rangle}\colon\Z^{1+n}\to K^\times$ is a group morphism.
Suppose now that~$y\in K^\times$, $(vy)^{\langle n\rangle}<0$, and $\i=(i_0,\dots,i_n)\in\Z^{1+n}$, 
$\i\neq 0$, and $m\in\{0,\dots,n\}$ is minimal with $i_m\neq 0$.
Then by Corollary~\ref{cor:itpsi}, $\big[v(y^{\langle\i\rangle})\big]=\big[v(y^{\langle m\rangle})\big]$. Thus if~$y\succ 1$, we have the equivalence
$y^{\langle\i\rangle} \succ 1\  \Leftrightarrow\ i_m\geq 1$.  
%Hence if $f\neq 0$ with $y^{\langle m\rangle}\steeper f$ and $i_m\geq 1$,
%then $y^{\langle\i\rangle} \succ f$.
If $K$ is equipped with an ordering making it a pre-$H$-field and $y\succ 1$, then $y^\dagger >0$, so
$y^{\<i\>}>0$ for   $i=1,\dots,n$, and thus~$\sgn y^{\<\i\>} = \sgn y^{i_0}$.

\subsection*{Iterated exponentials}  
{\it In this subsection we assume   that $\Psi$ is downward closed.}\/
For $f\succ 1$ we have $f'\succ f^\dagger$, so we can and do choose $\Exp(f)\in K^\times$ such that~${\Exp(f)\succ 1}$ and~$\operatorname{E}(f)^\dagger \asymp f'$, hence $f\prec \Exp(f)$ and $f \flatter \Exp(f)$. 
Moreover,  if  $f,g \succ 1$, then $${f \prec g} \quad\Longleftrightarrow\quad \Exp(f) \flatter \Exp(g).$$
For~$f\succ 1$
  define~$\Exp_n(f)\in K^{\succ 1}$ inductively
by 
$$ \Exp_0(f)\ :=\ f,\qquad \Exp_{n+1}(f)\ :=\ \Exp\!\big(\!\Exp_n(f)\big),$$ 
and thus by induction
$$\Exp_n(f)\ \prec\ \Exp_{n+1}(f) \quad\text{ and }\quad\Exp_{n}(f)\ \flatter\ \Exp_{n+1}(f)\qquad \text{ for all $n$.}$$
{\it In the rest of this subsection $f\succeq x$, and $y$ ranges over elements of $H$-asymptotic extensions of $K$.}\/
The proof of the next lemma is like that of~\cite[Lem\-ma~1.3(2)]{AvdD3}.

%Let $K$ be a Liouville closed $H$-field.    As in \cite[Section~1]{AvdD} we choose for each $f\in K$ an element $\operatorname{E}(f)\in K^>$ such that $\operatorname{E}(f)^\dagger=f'$. (So for $g\in K^>$ we have $g^\dagger=f'$ iff $g=c\operatorname{E}(f)$ for some $c\in C^>$.) From \cite{AvdD}   recall the following rules about $\operatorname{E}$, for $f,g\in K$: 
%\begin{itemize}
%\item[(E1)] $\Exp(f+g)=c\Exp(f)\Exp(g)$ and $\Exp(-f)=d\Exp(f)^{-1}$, where  $c,d\in C^{>}$;
%\item[(E2)] $f \preceq 1 \Longleftrightarrow \Exp(f) \asymp 1$;
%\item[(E3)] $f>C \Longleftrightarrow \Exp(f)\succ 1$; $f<C \Longleftrightarrow \Exp(f)\prec 1$;
%\item[(E4)] $1 \prec f  \Longrightarrow     f \flatter \Exp(f) $;
%\item[(E5)] $f>C  \Longrightarrow \Exp(f)>f^n $; $f<C  \Longrightarrow   0<\Exp(f)<|f|^{-n}<C^{>}$ ;
%\item[(E6)] If  $f,g \nasymp 1$, then $f \prec g \Longleftrightarrow \Exp(f) \flatter \Exp(g)$.
%\end{itemize}

\begin{lemma}\label{lem:expn, 1}
If $y\succeq\Exp_{n+1}(f)$, $n\geq 1$, then $y\neq 0$ and $y^\dagger\succeq\Exp_n(f)$.
\end{lemma}
\begin{proof}
If $y\succeq\Exp_2(f)$, then
$y\neq 0$, and using $\Exp_2(f)\succ 1$ we obtain
$$y^\dagger\ \succeq\ \Exp_2(f)^\dagger\  \asymp\ \Exp(f)'\  =\  \Exp(f)\Exp(f)^\dagger \ \asymp\  \Exp(f)f' \ \succeq\  \Exp(f),
$$
Thus the lemma holds for $n=1$. In general,~$\Exp_{n-1}(f) \succeq f \succeq x$, hence  the lemma follows from the case $n=1$ applied
to~$\Exp_{n-1}(f)$ in place of $f$.
\end{proof}

\noindent
An obvious induction on $n$ using Lemma~\ref{lem:expn, 1} shows: if $y\succeq \Exp_n(f)$, then $(vy)^{\langle n \rangle}\le vf <0$. We shall use this fact without further reference.

\begin{lemma}\label{lem:expn, 2}
If   $y \succeq \Exp_{n+1}(f)$, then $y^{\langle n\rangle}$ is defined and $y^{\langle n\rangle}  \succeq\Exp(f)$.
\end{lemma}
\begin{proof}
First note that if $y\neq 0$, $n\ge 1$, and $(y^\dagger)^{\langle n-1\rangle}$ is defined,
then $y^{\langle n\rangle}$ is defined and $y^{\langle n\rangle}=(y^\dagger)^{\langle n-1\rangle}$.
Now use induction on $n$ and Lemma~\ref{lem:expn, 1}.
\end{proof}

\begin{lemma}\label{lem:expn, 3}
If
$y\succeq\Exp_n(f^2)$, then  $y^{\langle n\rangle}$ is defined and $y^{\langle n\rangle}  \succeq f$, with
$y^{\langle n\rangle}  \succ f$ if~$f\succ x$.
\end{lemma} 
\begin{proof}
This is clear if $n=0$, so suppose    $y\succeq\Exp_{n+1}(f^2)$. Then by Lemma~\ref{lem:expn, 2} (applied with $f^2$ in place of $f$)  we have $y^{\langle n\rangle}  \succeq\Exp(f^2)\succ 1$, so
$$y^{\<n+1\>}\ =\ (y^{\<n\>})^\dagger\ \succeq\ \Exp(f^2)^\dagger\ \asymp\ (f^2)'\ =\ 2ff'\ \succeq\ f,$$
with~$y^{\<n+1\>} \succ f$ if $f\succ x$,
 as required. 
\end{proof}

\begin{cor}\label{cor:expn}
Suppose  $y\succeq \operatorname{E}_{n}(f^2)$, and let $\i\in\Z^{1+n}$ be such that $\i>0$ lexicographically. Then 
$y^{\<n\>}$ is defined and $y^{\langle \i\rangle}\succeq f$, with
$y^{\langle \i\rangle}\succ f$ if~$f\succ x$.
\end{cor}
\begin{proof}
By Lemma~\ref{lem:expn, 3}, $y^{\langle n\rangle}$ is defined with $y^{\langle n\rangle}  \succeq f$, and $y^{\langle n\rangle}  \succ f$ if $f\succ x$.
Let~$m\in\{0,\dots,n\}$ be minimal such that $i_m\neq 0$; so $i_m\geq 1$.
If $m=n$ then~$y^{\<\i\>}=(y^{\<n\>})^{i_n}\succeq y^{\<n\>}$, hence $y^{\langle \i\rangle}\succeq f$, with $y^{\langle \i\rangle}\succ  f$ if $f\succ x$.
Suppose $m<n$.
Then~$y\succeq\Exp_{m+1}(f^2)$ and hence $y^{\<m\>}\succeq\Exp(f^2)$ by 
Lemma~\ref{lem:expn, 2}. 
Also,   $f\comp f^2\flatter\Exp(f^2)$, thus~$y^{\langle m\rangle}\steeper f$. 
The remarks following Corollary~\ref{cor:itpsi} now yield $y^{\langle \i\rangle}\succ f$.
\end{proof}

\subsection*{Asymptotic behavior of $P(y)$ for large $y$}
In this subsection~$\i$,~$\j$,~$\k$ range over~$\N^{1+n}$.
Let~$P_{\<\i\>}\in K$ be such that $P_{\<\i\>}=0$ for all but finitely many~$\i$ and~${P_{\<\i\>}\neq 0}$ for some~$\i$, and set
$P:=\sum_{\i} P_{\<\i\>} Y^{\<\i\>}\in K\<Y\>$. 
So if
  $P\in K\{Y\}$, then~$P=\sum_{\i} P_{\<\i\>} Y^{\<\i\>}$ is the logarithmic decomposition of the differential polynomial $P$ as defined in
[ADH, 4.2]. 
If $y$ is an element in a differential field extension~$L$ of $K$ such that~$y^{\<n\>}$ is defined, then
we put~$P(y):=\sum_{\i} P_{\<\i\>} y^{\<\i\>}\in L$ (and for~$P\in K\{Y\}$ this has the usual value).
Let $\j$ be lexicographically maximal such that~$P_{\<\j\>}\neq 0$, and 
choose $\k$ so that
$P_{\<\k\>}$ has   minimal valuation.
If $P_{\<\k\>}/P_{\<\j\>} \succ x$,  set
$f:= P_{\<\k\>}/P_{\<\j\>}$; otherwise set $f:=x^2$. Then $f\succ x$ and $f\succeq P_{\<\i\>}/P_{\<\j\>}$ for all $\i$.
The following is a more precise version of  [ADH, 16.6.10] and
\cite[(8.8)]{JvdH}:

\begin{prop}\label{prop:val at infty}
Suppose $\Psi$ is    downward closed, and $y$ in an $H$-asymptotic extension of $K$  satisfies  $y\succeq\Exp_{n}(f^2)$. Then  $y^{\<n\>}$ is defined and~$P(y)\sim P_{\<\j\>}y^{\langle \j\rangle}$.
\end{prop}

\begin{proof}
Let $\i<\j$.
We have   $f\succ x$, so 
$y^{\langle \j-\i\rangle}\succ f \succeq P_{\<\i\>}/P_{\<\j\>}$ by Corollary~\ref{cor:expn}. Hence 
$P_{\<\j\>}y^{\langle \j\rangle} \succ P_{\<\i\>}y^{\langle \i\rangle}$.
\end{proof}

\noindent
From Corollary~\ref{cor:itpsi}, Lemma~\ref{lem:expn, 3}, and Proposition~\ref{prop:val at infty} we obtain:

\begin{cor}\label{cor:val gp at infty}
Suppose $\Psi$ is   downward closed and $y$ in an $H$-asymptotic extension of $K$ satisfies $y\succ K$. Then
  $y$ is $\d$-transcendental over~$K$, and  for all $n$, $y^{\langle n\rangle}$ is defined, $y^{\langle n\rangle}\succ K$, and~$y^{\langle n+1\rangle}\flatter y^{\langle n\rangle}$. The  $H$-asymptotic extension $K\langle y\rangle$ of $K$ has residue field~$\res K\langle y\rangle=\res K$ and value group~$\Gamma_{K\langle y\rangle} = \Gamma \oplus \bigoplus_n \Z v(y^{\langle n\rangle})$ \textup{(}internal direct sum\textup{)}, and $\Gamma_{K\langle y\rangle} $
contains $\Gamma$ as a convex subgroup.
\end{cor}

\noindent
Suppose now that $K$ is equipped with an ordering making it a pre-$H$-field.  From Proposition~\ref{prop:val at infty} we recover~\cite[Theorem~3.4]{AvdD3} in slightly stronger form: 

\begin{cor}\label{cor:val at infty}
Suppose $y$  lies in
a Liouville closed $H$-field extension of~$K$. If~$y\succeq\Exp_{n}(f^2)$, then $y^{\<n\>}$ is defined and $\sgn P(y)=\sgn P_{\<\j\>}y^{j_0}$. 
In particular, if~$y^{\<n\>}$ is defined and~$P(y)=0$, then~$y\prec \Exp_{n}(f^2)$.
\end{cor}

\begin{example}
Suppose $P\in K\{Y\}$. Using [ADH, 4.2, subsection on logarithmic decomposition] we obtain $j_0=\deg P$, and the logarithmic decomposition
$$P(-Y)\ =\ \sum_{\i} P_{\langle\i\rangle} (-1)^{i_0}Y^{\langle\i\rangle}.$$
If $\deg P$ is odd, and  $y>0$ lies in 
a Liouville closed $H$-field extension of $K$ such that~$y\succeq\Exp_{n}(f^2)$, 
then $$\sgn P(y)\ =\ \sgn P_{\<\j\>}, \qquad
\sgn P(-y)\ =\ -\sgn P_{\langle\j\rangle}\ =\ -\sgn P(y).$$
\end{example}

\section{$\upl$-freeness and $\upo$-freeness}\label{sec:uplupo-freeness} 

\noindent
This section contains  preservation results for  the important properties of {\it $\upl$-freeness}\/ and {\it $\upo$-freeness}\/ from [ADH].  Let  $K$ be an ungrounded $H$-asymptotic field such that~${\Gamma\ne \{0\}}$, and
as in [ADH, 11.5],  fix a logarithmic sequence~$(\ell_\rho)$ for $K$ and define the pc-sequences~$(\upl_\rho)=(-\ell_\rho^{\dagger\dagger})$ and~$(\upo_\rho)=(\omega(\upl_\rho))$ in~$K$, where~$\omega(z):=-2z'-z^2$.  
Recall that~$K$ is  {\it $\upl$-free}\/ iff $(\upl_\rho)$ does not have a pseudolimit in $K$, and~$K$ is   {\it $\upo$-free}\/ iff~$(\upo_\rho)$ does not have a pseudolimit in~$K$.
 If $K$ is $\upo$-free, then~$K$ is $\upl$-free.
We refer to~[ADH, 11.6, 11.7] for this and other basic facts about $\upl$-freeness and $\upo$-freeness used below. (For $\upo$-free Hardy fields, see also Section~\ref{sec:upo-free Hardy fields}.) As in [ADH], $L$ being $\upl$-free or $\upo$-free includes $L$ being an ungrounded $H$-asymptotic field with $\Gamma_L\ne \{0\}$.

\subsection*{Preserving $\upl$-freeness and $\upo$-freeness}
{\it In this subsection $K$ is an ungrounded $H$-asymptotic field with $\Gamma\ne \{0\}$, and~$(\ell_\rho)$, $(\upl_\rho)$, 
$(\upo_\rho)$ are as above.}\/
If $K$ has  a $\upl$-free $H$-asymptotic field extension $L$ such that $\Gamma^<$ is cofinal in $\Gamma_L^<$, then $K$ is
$\upl$-free, and similarly with ``$\upo$-free'' in place of ``$\upl$-free'' [ADH, remarks after 11.6.4, 11.7.19].
The property of $\upo$-freeness is very robust; indeed, by~[ADH, 13.6.1]:

\begin{theorem}\label{thm:ADH 13.6.1}
If $K$ is $\upo$-free and $L$ is a pre-$\d$-valued $\d$-algebraic $H$-asymptotic  extension  of $K$, then $L$  is $\upo$-free and $\Gamma^<$ is cofinal in $\Gamma_L^<$.
\end{theorem}

\noindent
In contrast, $\upl$-freeness is more delicate:  Theorem~\ref{thm:ADH 13.6.1} fails with ``$\upl$-free'' in place of ``$\upo$-free'', as the next example shows.

\begin{exampleNumbered}\label{ex:rat as int and cofinality} 
The $H$-field $K=\R\langle\upo\rangle$ from [ADH, 13.9.1] is $\upl$-free, but its $H$-field extension
$L=\R\langle \upl\rangle$ is not, and this extension is $\d$-algebraic: $2\upl'+\upl^2+\upo=0$.
\end{exampleNumbered}

\noindent
In the rest of this subsection we consider cases where parts of Theorem~\ref{thm:ADH 13.6.1} do hold.
Recall from [ADH, 11.6.8]  that if~$K$ is $\upl$-free, then~$K$ has (rational) asymptotic integration, and $K$
is $\upl$-free iff its algebraic closure is $\upl$-free.
Moreover, $\upl$-freeness is preserved under adjunction of constants: 

\begin{prop}\label{prop:const field ext}
Suppose $K$ is $\upl$-free and $L=K(D)$ is an $H$-asymptotic extension of $K$ with  $D\supseteq C$ a subfield of $C_L$. 
Then $L$ is $\upl$-free with $\Gamma_L=\Gamma$.
\end{prop}

\noindent
We are going  to deduce this from the next three lemmas. 
Recall that $K$ is pre-$\d$-valued, by  [ADH, 10.1.3]. Let $\operatorname{dv}(K)$ be the $\d$-valued hull
of $K$ (see [ADH, 10.3]).

\begin{lemma}\label{lem:Gehret}
Suppose $K$ is $\upl$-free. Then  $L:=\operatorname{dv}(K)$ is $\upl$-free and $\Gamma_L=\Gamma$.
\end{lemma}
\begin{proof}
The first statement is  \cite[Theorem~10.2]{Gehret},
%\marginpar{\cite[10.2]{Gehret} accepted on faith for now} 
and the second statement follows from [ADH, 10.3.2(i)].
\end{proof}

\noindent
If $L=K(D)$ is a differential  field extension of $K$
with $D\supseteq C$ a subfield of $C_L$, then~$D=C_L$, and $K$ and $D$ are linearly disjoint over $C$ [ADH, 4.6.20].
If $K$ is $\d$-valued and $L=K(D)$ is an $H$-asymptotic extension of $K$ with $D\supseteq C$ a subfield of $C_L$,
then $L$ is $\d$-valued  and~$\Gamma_L=\Gamma$  [ADH, 10.5.15].

\begin{lemma}\label{lem:const field ext}
Suppose $K$ is $\d$-valued and $\upl$-free, and $L=K(D)$ is an $H$-asymptotic extension of $K$ with $D\supseteq C$ a subfield of $C_L$. Then $L$ is $\upl$-free.
\end{lemma}
\begin{proof} First, $(\upl_{\rho})$ is of transcendental type over $K$: otherwise, [ADH, 3.2.7] would give an algebraic extension of $K$ that is not $\upl$-free. Next, our logarithmic sequence~$(\ell_\rho)$ for $K$ remains a logarithmic sequence for $L$. 

Zorn and the $\forall\exists$-form of the $\upl$-freeness axiom~[ADH, 1.6.1(ii)] reduce us to the case~$D=C(d)$, $d\notin C$, $d$ transcendental over $K$, so $L=K(d)$. 
%Now $L$ is $\upl$-free iff the $\d$-valued subfield $K^{\operatorname{a}}(d)$ of~$L^{\operatorname{a}}$ is $\upl$-free.
%Hence we may also replace $K$ by $K^{\operatorname{a}}$ to arrange that $K$, and hence also~$C\cong\res K$, are algebraically closed.
 Suppose~$L$ is not $\upl$-free. Then~$\upl_\rho\leadsto\upl\in L$, and such $\upl$ is transcendental over $K$ and gives an immediate extension $K(\upl)$ of $K$ by [ADH, 3.2.6].
%  Then~$K(\upl)$ is an immediate extension of~$K$, so we have some $\upl\in L$ with~$\upl_\rho\leadsto\upl$. Then~$K(\upl)$ is an immediate extension of~$K$ and $\upl$ is transcendental over~$K$~[ADH, 3.2.6, 3.2.7,  11.6.8].
Hence $L$ is algebraic over  $K(\upl)$,   so 
$\res L$ is algebraic over~$\res K(\upl)=\res K\cong C$ and thus $d$ is algebraic over $C$, a contradiction.
\end{proof}

\begin{lemma}\label{conpred}
Suppose $K$ is $\upl$-free and $L$ is an $H$-asymptotic extension of $K$, where~$L=K(d)$ with $d\in C_L$.  
Then $L$ is pre-$\d$-valued.
\end{lemma}
\begin{proof} 
Let $L^{\operatorname{a}}$ be an algebraic closure of the $H$-asymptotic field $L$, and let $K^{\operatorname{a}}$ be the algebraic closure~ of~$K$ inside $L^{\operatorname{a}}$. Then $K^{\operatorname{a}}$  is pre-$\d$-valued by [ADH, 10.1.22]. Replacing $K$, $L$ by $K^{\operatorname{a}}$, $K^{\operatorname{a}}(d)$ we arrange that
$K$ is algebraically closed. We may assume $d\notin C$, so $d$ is transcendental over $K$ by [ADH, 4.1.1, 4.1.2].

Suppose first that $\operatorname{res}(d)\in\operatorname{res}(K)\subseteq \operatorname{res}(L)$, and take~$b\in\mathcal O$ such that
$y:=b-d\prec 1$.
Then   $b'\notin\der\smallo$: otherwise $y'=b'=\delta'$ with $\delta\in\smallo$, so $y=\delta\in K$ and hence~$d\in K$,
a contradiction. 
Also~$vb'\in (\Gamma^>)'$: otherwise  $vb'<(\Gamma^>)'$, by [ADH, 9.2.14], and $vb'$ would be a gap in $K$, contradicting
$\upl$-freeness of $K$. Hence~$L=K(y)$ is pre-$\d$-valued by [ADH, 10.2.4, 10.2.5(iii)]  applied to $s:=b'$.  

If $\operatorname{res}(d)\notin\operatorname{res}(K)$, then $\res(d)$ is transcendental over $\res(K)$ by~[ADH, 3.1.17], hence
$\Gamma_L=\Gamma$ by [ADH, 3.1.11], and so $L$ has asymptotic integration and thus is pre-$\d$-valued by [ADH, 10.1.3].
\end{proof}

\begin{proof}[Proof of Proposition~\ref{prop:const field ext}] By Zorn we reduce to the case $L=K(d)$ with $d\in C_L$.
Then $L$ is pre-$\d$-valued by Lemma~\ref{conpred}.
By Lemma~\ref{lem:Gehret}, the $\d$-valued hull $K_1:=\operatorname{dv}(K)$ of $K$ is $\upl$-free with $\Gamma_{K_1}=\Gamma$, and by the universal property
of $\d$-valued hulls we may arrange that $K_1$ is a $\d$-valued subfield of $L_1:=\operatorname{dv}(L)$ [ADH, 10.3.1].
The proof of~[ADH, 10.3.1] gives $L_1=L(E)$ where $E=C_{L_1}$, and so
$L_1=K_1(E)$. Hence by Lemma~\ref{lem:const field ext} and the remarks preceding it, $L_1$ is $\upl$-free with $\Gamma_{L_1}=\Gamma_{K_1}=\Gamma$.
Thus $L$ is $\upl$-free with~$\Gamma_L=\Gamma$.
\end{proof}

\begin{lemma}\label{lemtrigupl} Let $H$ be a $\upl$-free real closed $H$-field. Then the trigonometric closure $H^{\trig} $ of $H$ is $\upl$-free.
\end{lemma} 
\begin{proof} We show that $H^+$ as in Lemma~\ref{lemtri2} is $\upl$-free. There $H^+[\imag]=K(y)$ where~$K$ is the $H$-asymptotic extension $H[\imag]$ of $H$ and $y\sim 1$, $y^\dagger\notin K^\dagger$,  $y^\dagger\in\imag\der\smallo_H$.  Then $K$ is $\upl$-free, so $K(y)$ is $\upl$-free by \cite[Proposition~7.2]{Gehret},   
hence $H^+$ is $\upl$-free. 
\end{proof}

\noindent
In Example~\ref{ex:rat as int and cofinality} we have a $\upl$-free $K$ and an $H$-asymptotic extension $L$ of $K$ that is not $\upl$-free, with $\operatorname{trdeg}(L|K) = 1$.
The next proposition shows that the second part of the conclusion of Theorem~\ref{thm:ADH 13.6.1} nevertheless
holds for such $K,L$. 
%Recall that if $K$ has rational asymptotic integration, then~$K$ has asymptotic integration and 
%hence is pre-$\d$-valued, by [ADH, 9.2.8, 10.1.3]. 

\begin{prop}\label{prop:rat as int and cofinality}
The following are equivalent:
\begin{enumerate}
\item[\textup{(i)}] $K$ has rational asymptotic integration;
\item[\textup{(ii)}] for every $H$-asymptotic extension $L$ of $K$ with $\operatorname{trdeg}(L|K)\le 1$ we have that
$\Gamma^{<}$  is cofinal in $\Gamma_L^{<}$.
\end{enumerate} 
\end{prop} 
\begin{proof} For  (i)~$\Rightarrow$~(ii), assume (i), and let  $L$ be an $H$-asymptotic extension of $K$ with~$\operatorname{trdeg}(L|K)\le 1$. Towards showing that $\Gamma^{<}$ is cofinal in $\Gamma_L^{<}$ we can arrange that $K$ and $L$ are algebraically closed. Suppose towards a contradiction that $\gamma\in \Gamma_L$ and
$\Gamma^{<} < \gamma < 0$. Then $\Psi < \gamma' < (\Gamma^{>})'$, and so $\Gamma$ is dense in
$\Gamma+\Q\gamma'$ by [ADH, 2.4.16, 2.4.17], in particular, $\gamma\notin \Gamma+\Q\gamma'$. Thus
$\gamma$, $\gamma'$ are $\Q$-linearly independent over $\Gamma$, which contradicts $\operatorname{trdeg}(L|K)\le 1$ by [ADH, 3.1.11]. 

As to (ii)~$\Rightarrow$~(i), we prove the contrapositive, so assume $K$ does not have rational asymptotic integration.
We arrange again that $K$ is algebraically closed. 
Then $K$ has a gap $vs$ with $s\in K^\times$, and so [ADH, 10.2.1 and its proof] gives an $H$-asymptotic extension $K(y)$ of $K$ with $y'=s$ and $0 < vy < \Gamma^{>}$. 
\end{proof}

\noindent
Recall from [ADH, 11.6] that Liouville closed $H$-fields are $\upl$-free. 
To prove the next result we also use Gehret's theorem~\cite[Theorem 12.1(1)]{Gehret}  
%\marginpar{Gehret's result taken on faith for now}
that an $H$-field $H$  has up to isomorphism over $H$ exactly one Liouville closure iff $H$ is grounded or $\upl$-free. Here {\it isomorphism\/} means of course {\it isomorphism of $H$-fields\/}, and likewise with the embeddings referred to in the next result:

\begin{prop}\label{proptl} Let $H$ be a grounded or $\upl$-free  $H$-field. Then $H$ has a trigonometrically closed 
and Liouville closed $H$-field extension $H^{\tl}$ that embeds over $H$ into any trigonometrically closed 
Liouville closed $H$-field extension of $H$.
\end{prop} 
\begin{proof} We build real closed $H$-fields $H_0\subseteq H_1\subseteq H_2\subseteq \cdots$ as follows: $H_0$ is a real closure of $H$, and, recursively,
  $H_{2n+1}$ is a Liouville closure of $H_{2n}$, and~$H_{2n+2}:= H_{2n+1}^{\trig}$ is the trigonometric closure of $H_{2n+1}$. 
Then $H^{*}:= \bigcup_n H_n$ is a trigonometrically closed Liouville closed $H$-field extension of $H$. 
Induction using Lemma~\ref{lemtrigupl} shows that all $H_n$ with $n\ge 1$ are $\upl$-free, and that 
$H_{2n}$ has for all~$n$ up to isomorphism over $H$ a unique Liouville closure. Given any trigonometrically closed Liouville closed $H$-field extension $E$ of $H$ we then use the embedding properties of {\it Liouville closure\/} and {\it trigonometric closure\/} to construct by a similar recursion embeddings $H_n\to E$ that extend to an embedding $H^{*}\to E$ over $H$.  
\end{proof} 

\noindent
For $H$ as in Proposition~\ref{proptl}, the $H^*$ constructed in its proof is minimal: Let~${E\supseteq H}$ be any trigonometrically closed Liouville closed $H$-subfield of $H^*$. Then induction on $n$ yields $H_n\subseteq~E$ for all $n$, so $E=H^*$.
It follows that any $H^{\tl}$ as in Proposition~\ref{proptl} is isomorphic over $H$ to $H^*$, and we refer to such
$H^{\tl}$ as a {\bf  trigonometric-Liouville closure\/} of $H$.\index{closure!trigonometric-Liouville}
\index{trigonometric!Liouville closure}\label{p:Htl} Here are some useful facts about $H^{\tl}$:

\begin{cor}\label{cortrig} Let $H$ be a $\upl$-free $H$-field. Then $C_{H^{\tl}}$ is a real closure of $C_H$, the $H$-asymptotic extension $K^{\tl}:= H^{\tl}[\imag]$ of $H^{\tl}$
is a Liouville extension of $H$ with~$\I(K^{\tl})\subseteq (K^{\tl})^\dagger$, and $\Gamma_H^{<}$ is cofinal in $\Gamma_{H^{\tl}}^{<}$. Moreover,
$$ \text{$H$ is $\upo$-free}\ \Longleftrightarrow\  \text{$H^{\tl}$ is $\upo$-free.}$$ 
\end{cor}
\begin{proof} The construction of $H^*$ in the proof of Proposition~\ref{proptl} gives that $C_{H^{*}}$ is a real closure of $C_H$,
and that the $H$-asymptotic extension $K^{*}:=H^*[\imag]$ of~$H^*$
is a Liouville extension of $H$ with $\I(K^{*})\subseteq (K^{*})^\dagger$.  Induction using Lemma~\ref{lemtrigupl} and Proposition~\ref{prop:rat as int and cofinality} shows that $H_n$ is $\upl$-free  and $\Gamma_H^{<}$ is cofinal in $\Gamma_{H_n}^{<}$, for all $n$, so  $\Gamma_H^<$ is cofinal in $\Gamma_{H^*}^{<}$. 

The final equivalence follows from Theorem~\ref{thm:ADH 13.6.1} and  a remark preceding it.  
\end{proof}

%\begin{cor}\label{cortrig} Let $H$ be a $\upl$-free real closed $H$-field. Then $H$ has a Liouville closed $H$-field extension
%$H^*$  such that $H^*$ has the same constant field as $H$, the $H$-asymptotic extension $K^*:= H^*[\imag]$ of $H^*$
%is a Liouville extension of $H$ for which~$\I(K^*)\subseteq (K^*)^\dagger$, and $\Gamma_H^{<}$ is cofinal in $\Gamma_{H^*}^{<}$. 
%For any such $H^*$ we have: 
%$$H \text{ is $\upo$-free}\ \Longleftrightarrow\ H^* \text{ is $\upo$-free}.$$ 
%\end{cor}
%\begin{proof}  We build an increasing chain $H_0\subseteq H_1\subseteq H_2\subseteq \cdots$ of real closed $H$-fields as follows. We set  $H_0:= H$, and, recursively,
 % $H_{2n+1}:= H_{2n}^{\trig}$ is the trigonometric closure of $H_{2n}$, and $H_{2n+2}$ is a Liouville closure of $H_{2n+1}$. 
%Then $H^*:= \bigcup_n H_n$ is a Liouville closed $H$-field extension with the same constant field as $H$, and $K^*$ 
%is a Liouville extension of $H$ with $\I(K^*)\subseteq (K^*)^\dagger$.  Induction using Lemma~\ref{lemtrigupl} and Proposition~\ref{prop:rat as int and cofinality} shows that all $H_n$ are $\upl$-free,  and $\Gamma_H^{<}$ is cofinal in $\Gamma_{H_n}^{<}$ for all $n$, so  $\Gamma_H^<$ is cofinal in $\Gamma_{H^*}^{<}$. 

%The final equivalence follows from Theorem~\ref{thm:ADH 13.6.1} and  a remark preceding it.  
%\end{proof} 

\noindent
Proposition~\ref{prop:rat as int and cofinality} and
[ADH, remarks after~11.6.4 and   after~11.7.19] yield:

\begin{cor}\label{cor:rat as int and cofinality}
Suppose $K$ has rational asymptotic integration, and let $L$ be an $H$-asymptotic extension of $K$ with $\operatorname{trdeg}(L|K)\leq 1$.
If $L$ is $\upl$-free, then so is $K$, and if $L$ is $\upo$-free, then so is $K$.
\end{cor}

\noindent
We also have a similar characterization of $\upl$-freeness:

\begin{prop}\label{prop:upl-free and as int} 
The following are equivalent:
\begin{enumerate}
\item[\textup{(i)}] $K$ is $\upl$-free;
\item[\textup{(ii)}] every $H$-asymptotic extension $L$ of $K$ with $\operatorname{trdeg}(L|K)\le 1$ has asymptotic integration.
\end{enumerate} 
\end{prop}
\begin{proof}
Assume $K$ is $\upl$-free; let $L$ be an $H$-asymptotic extension of $K$ such that $\operatorname{trdeg}(L|K)\le 1$.
By Proposition~\ref{prop:rat as int and cofinality},  $\Gamma^{<}$ is cofinal in $\Gamma_L^{<}$,  so $L$ is ungrounded.
Towards
a contradiction, suppose $vf$ ($f\in L^\times$) is a gap in $L$. %Then $vf$ remains a gap in an algebraic closure $L^{\operatorname{a}}$ of $L$, and the algebraic closure of $K$ in $L^{\operatorname{a}}$ remains $\upl$-free.
Passing to algebraic closures we arrange that $K$ and $L$ are algebraically closed. 
 Set $\upl:=-f^\dagger$. Then for all active $a$ in $L$
we have $\upl+a^\dagger\prec a$ by~[ADH, 11.5.9] and hence $\upl_\rho\leadsto\upl$ by~[ADH, 11.5.6]. By $\upl$-freeness of $K$
and~[ADH, 3.2.6, 3.2.7], the valued field extension~$K(\upl)\supseteq K$ is  immediate of transcendence degree~$1$, so~$L\supseteq K(\upl)$ is algebraic and~$\Gamma=\Gamma_L$. Hence $vf$ is a gap in $K$, a contradiction. This shows~(i)~$\Rightarrow$~(ii).

To show the contrapositive of (ii)~$\Rightarrow$~(i),  
 suppose  $\upl\in K$ is a pseudolimit of~$(\upl_\rho)$. If the algebraic closure $K^{\operatorname{a}}$ of $K$ does not have asymptotic integration,
 then clearly~(ii) fails. If $K^{\operatorname{a}}$ has asymptotic integration,
 then $-\upl$ creates a gap over~$K$ by~[ADH, 11.5.14] applied to $K^{\operatorname{a}}$ in place of $K$, hence (ii) also fails.
\end{proof}

\noindent
The next two lemmas include converses to Lemmas~\ref{lem:Gehret} and~\ref{lem:const field ext}.

\begin{lemma}\label{notupl1} Let $E$ be a pre-$\d$-valued $H$-asymptotic field. Then:\begin{enumerate}
\item[\textup{(i)}] if $E$ is not $\upl$-free, then $\operatorname{dv}(E)$ is not $\upl$-free; 
\item[\textup{(ii)}] if $E$ is not $\upo$-free, then $\operatorname{dv}(E)$ is not $\upo$-free.
\end{enumerate}
\end{lemma}
\begin{proof} This is clear if $E$ has no rational asymptotic integration, because then $\operatorname{dv}(E)$ has no rational asymptotic integration either, by [ADH, 10.3.2]. Assume $E$ has rational asymptotic integration. Then $\operatorname{dv}(E)$
is an immediate extension of $E$ by~[ADH, 10.3.2], and then (i) and (ii) follow from the characterizations of $\upl$-freeness and 
$\upo$-freeness in terms of nonexistence of certain pseudolimits. 
\end{proof} 

\begin{lemma}\label{notupl2} Let $E$ be a $\d$-valued $H$-asymptotic field and $F$ an $H$-asymptotic extension of $E$ such that
$F=E(C_F)$.
Then: \begin{enumerate}
\item[\textup{(i)}] if $E$ is not $\upl$-free, then $F$ is not $\upl$-free; 
\item[\textup{(ii)}] if $E$ is not $\upo$-free, then $F$ is not $\upo$-free.
\end{enumerate}
\end{lemma}
\begin{proof} By [ADH, 10.5.15] $E$ and $F$ have the same value group. The rest of the proof is like that for the previous lemma, with $F$ instead of $\operatorname{dv}(E)$. 
\end{proof}

\noindent
{\it In the rest of this subsection $K$ is in addition a pre-$H$-field and $L$ a pre-$H$-field extension of $K$.}\/
The following is shown in the proof of~\cite[Lemma~12.5]{Gehret}:

\begin{prop}[Gehret] \label{prop:Gehret} 
%\marginpar{accepted for now on faith}
Suppose $K$  is a $\upl$-free $H$-field and $L$ is a Liou\-ville $H$-field extension of $K$.
Then $L$ is $\upl$-free and $\Gamma^<$ is cofinal in~$\Gamma_L^<$. %\textup{(}In particular, if $L$ is $\upo$-free, then so is $H$.\textup{)}
\end{prop}

\begin{exampleNumbered}\label{ex:Gehret} 
Let $K=\R\langle\upo\rangle$ be the $\upl$-free but non-$\upo$-free $H$-field from [ADH, 13.9.1]. Then 
$K$ has a unique Liouville closure $L$, up to isomorphism over $K$, by~\cite[Theorem~12.1(1)]{Gehret}.
%\marginpar{\bf onfaithfornow} 
By Proposition~\ref{prop:Gehret}, $L$
is not $\upo$-free; \cite{ADH3} has another proof of this fact. By [ADH, 13.9.5] we can take
here $K$ to be a Hardy field, and then $L$ is isomorphic  over $K$ to a Hardy field
extension of $K$ [ADH, 10.6.11].

Applying Corollary~\ref{cortrig} to $H:=\R\langle \upo\rangle$ yields a Liouville closed $H$-field~$H^{\tl}$ that is not
$\upo$-free but does satisfy $\I(K^{\tl})\subseteq (K^{\tl})^\dagger$ for $K^{\tl}:= H^{\tl}[\imag]$. 
\end{exampleNumbered}

\noindent
For a pre-$H$-field $H$ we singled out in [ADH, p.~520]  the following subsets: \label{p:special subsets}
$$\Upg(H)\ :=\ (H^{\succ 1})^\dagger,\quad \Upl(H)\ :=\ -(H^{\succ 1})^{\dagger\dagger}, \quad 
\Upd(H)\ :=\ -(H^{\neq,\prec 1})^{\prime\dagger}.$$
 
\begin{lemma}\label{lem:upl in D(H)}
Suppose $K$ is $\upl$-free, $\upl\in \Upl(L)^\downarrow$, $\upo:=\omega(\upl)\in K$, and suppose~$\omega\big(\Upl(K)\big) <  \upo   < \sigma\big(\Upg(K)\big)$. Then $\upl_\rho\leadsto \upl$, and the pre-$H$-sub\-field~$K\langle \upl\rangle=K(\upl)$ of~$L$   is an immediate extension of $K$ \textup{(}and so $K\langle \upl\rangle$ is not $\upl$-free\textup{)}.
\end{lemma}
\begin{proof} From $\Upl(L)<\Upd(L)$ [ADH, p.~522] and $\Upd(K)\subseteq\Upd(L)$ we
obtain $\upl<\Upd(K)$. 
The restriction of  $\omega$ to $\Upl(L)^\downarrow$ is strictly increasing~[ADH, p.~526] and $\Upl(K)\subseteq\Upl(L)$,  so~$\omega\big(\Upl(K)\big)<\upo=\omega(\upl)$ gives $\Upl(K) < \upl$.
Hence $\upl_\rho\leadsto \upl$ by [ADH, 11.8.16]. Also~$\upo_\rho\leadsto\upo$ by [ADH, 11.8.30]. Thus $K\langle \upl\rangle$ is an immediate   extension of~$K$ by~[ADH, 11.7.13].
\end{proof}

\subsection*{Achieving $\upo$-freeness for pre-$H$-fields}
{\it In the rest of this section $H$ is a pre-$H$-field and  $L$
is a Liouville closed $\d$-algebraic   $H$-field extension of $H$.}\/
Thus if $H$ is $\upo$-free, then so is~$L$, by Theorem~\ref{thm:ADH 13.6.1}. 

The lemmas below give conditions guaranteeing that $L$ is $\upo$-free, while $H$ is not.
 
\begin{lemma}\label{lem:Li(H) upo-free} 
Suppose  $H$ is grounded or has a gap.  Then~$L$ is  $\upo$-free.
%$H$ is $\upo$-free in case $H$ has asymptotic integration. Then~$L$ is   $\upo$-free.
\end{lemma}
\begin{proof}
%Let $\widehat H$ be the  $H$-field hull  of $H$ inside $L$; if $H$ is $\upo$-free, then so is $\widehat H$, and 
%if $\widehat H$ has asymptotic integration, then so does $H$.
%Hence replacing $H$ by $\widehat H$ we first arrange that $H$ is an $H$-field.
%Now if $H$ is $\upo$-free, then the lemma follows from Theorem~\ref{thm:ADH 13.6.1}. The other cases successively reduce to this one:
Suppose~$H$ is grounded. Let $H_{\upo}$ be the $\upo$-free pre-$H$-field extension of $H$ introduced in connection with~[ADH,11.7.17] (where we use the letter $F$ instead of~$H$).  Identifying $H_{\upo}$ with its image in $L$ under an embedding $H_{\upo}\to L$ over $H$ of pre-$H$-fields, we apply Theorem~\ref{thm:ADH 13.6.1} to $K:=H_{\upo}$ to conclude that $L$ is $\upo$-free. 

  Next, suppose $H$ has a gap $\beta=vb$, $b\in H^\times$. Take~$a\in  L$
with~$a'=b$ and~$a\nasymp 1$. Then $\alpha:= va$ satisfies $\alpha'=\beta$, and so the pre-$H$-field
$H(a)\subseteq L$ is grounded, by~[ADH,  9.8.2 and remarks following its proof].
 Now apply the previous case to~$H(a)$ in place of $H$.
\end{proof}

\begin{lemma}\label{lem:D(H) upo-free, 2}
Suppose $H$ has asymptotic integration and divisible value group, and~$s\in H$ creates a gap over $H$. Then $L$ is $\upo$-free.
\end{lemma}
\begin{proof}
Take $f\in L^\times$ with $f^\dagger=s$.
Then by [ADH, remark after 11.5.14], $vf$ is a gap in~$H\langle f\rangle=H(f)$, so~$L$ is $\upo$-free by Lemma~\ref{lem:Li(H) upo-free} applied to $H\langle f\rangle$ in place of~$H$. 
\end{proof}

\begin{lemma}\label{lem:D(H) upo-free, 3}
Suppose $H$ is not $\upl$-free. Then $L$ is $\upo$-free.
\end{lemma}
\begin{proof}
By [ADH, 11.6.8], the real closure   $H^{\operatorname{rc}}$ of $H$ inside $L$
is not $\upl$-free, hence replacing $H$ by $H^{\operatorname{rc}}$  we arrange that $H$ is real closed.
If $H$ does not have asymptotic integration, then we are done 
by  Lemma~\ref{lem:Li(H) upo-free}. So suppose  $H$  has asymptotic integration. Then some
$s\in H$ creates a gap over $H$, by [ADH, 11.6.1], so $L$ is $\upo$-free by Lemma~\ref{lem:D(H) upo-free, 2}.
\end{proof}

\begin{cor}\label{cor:D(H) upo-free, 3}
Suppose $H$ is $\upl$-free and $\upl\in\Upl(L)^\downarrow$ is such that $\upo:=\omega(\upl)\in H$ and
$\omega\big(\Upl(H)\big) <  \upo   < \sigma\big(\Upg(H)\big)$. Then $L$ is $\upo$-free.
\end{cor}
\begin{proof}
By Lemma~\ref{lem:upl in D(H)}, the pre-$H$-subfield~$H\langle \upl\rangle=H(\upl)$ of $L$ is
an immediate non-$\upl$-free extension of $H$.  Now apply  Lemma~\ref{lem:D(H) upo-free, 3} to~$H\langle \upl\rangle$ in place of~$H$.
\end{proof}

\section{Complements on Linear Differential Operators}\label{sec:lindiff}

\noindent
In this section we tie up loose ends from the material on linear differential operators in [ADH,~14.2] and \cite[Section~8]{VDF}.
{\em Throughout $K$ is an ungrounded asymptotic field, $a$,~$b$,~$f$,~$g$,~$h$  range over arbitrary elements of $K$, and $\phi$ over those active in $K$, in particular,  $\phi\ne 0$}.
%We let We say that~$\phi$ is {\it active}\/ if~$\phi$ is active in~$K$.
Recall from [ADH, p.~479] our use of the term ``eventually'': a property $S(\phi)$ of elements 
$\phi$ is said to hold {\em eventually\/} if for some active $\phi_0$ in $K$, $S(\phi)$ holds for all $\phi\preceq \phi_0$. 

\medskip
\noindent
We shall consider linear differential operators $A\in K[\der]^{\ne}$ and  set $r:=\order(A)$. 
In~[ADH, Section~11.1] we introduced the set
$$\exc^{\ev}(A)\ =\ \exc^{\ev}_K(A)\ :=\ \big\{ \gamma\in\Gamma:\,  \nwt_A(\gamma)\geq 1  \big\}\ =\ \bigcap_{\phi}\exc(A^\phi)$$
of {\em eventual exceptional values of $A$}. \label{p:excev}\index{linear differential operator!eventual exceptional values}\index{exceptional values!eventual}For $a \neq 0$ we have $\exc^{\ev}(aA) = \exc^{\ev}(A)$ and $\exc^{\ev}(Aa) = \exc^{\ev}(A) - va$. An easy consequence of the definitions: $\exc^{\ev}(A^f)=\exc^{\ev}(A)$ for~$f\ne 0$. 
 A key fact about $\exc^{\ev}(A)$ is that if~$y\in K^\times$,
$vy\notin \exc^{\ev}(A)$, then $A(y)\asymp A^\phi y$, eventually.
Since $A^\phi y\ne 0$ for $y\in K^\times$, this gives
$v(\ker^{\neq} A)\subseteq \exc^{\ev}(A)$. 

\begin{lemma}\label{lemexc} If $L$ is an ungrounded asymptotic extension of $K$, then ${\exc^{\ev}_{L}(A)\cap\Gamma}\subseteq \exc^{\ev}(A)$, with equality if $\Psi$ is cofinal in $\Psi_L$. 
 \end{lemma}
 \begin{proof} For the inclusion, use that $\dwt(A^\phi)$ decreases as $v\phi$ strictly increases~[ADH, 11.1.12]. Thus its eventual value $\nwt(A)$, evaluated in $K$, cannot strictly increase when evaluated in an ungrounded asymptotic extension of $K$. 
 \end{proof} 

\noindent
{\em In the rest of this section we assume in addition that $K$ is $H$-asymptotic with asymptotic integration.}
Then by [ADH, 14.2.8]:

\begin{prop}\label{kerexc}
If $K$ is $r$-linearly newtonian, then $v(\ker^{\neq} A) = \exc^{\ev}(A)$.
\end{prop}

\begin{remarkNumbered}\label{rem:kerexc}
If  $K$ is $\d$-valued, then $\abs{v(\ker^{\neq} A)} = \dim_C \ker A \leq r$ by [ADH, 5.6.6], using a reduction to the case of ``small derivation'' by compositional conjugation. 
\end{remarkNumbered}

\begin{cor} \label{cor:nonexc sol} Suppose $K$ is $\d$-valued, $\exc^{\ev}(A)=v(\ker^{\neq} A)$, and $0\ne f\in A(K)$. Then $A(y)=f$ for some $y\in K$ with $vy\notin \exc^{\ev}(A)$.
\end{cor}
\begin{proof} Let $y\in K$, $A(y)=f$, with $vy$ maximal. Then $vy\notin \exc^{\ev}(A)$:
otherwise we have $z\in \ker A$ with $z\sim y$, so $A(y-z)=f$ and $v(y-z)>vy$.
\end{proof}

\begin{cor}\label{cor:sum of nwts}
Suppose $K$ is $\upo$-free. Then
$\sum_{\gamma\in\Gamma} \nwt_A(\gamma)=\abs{\exc^{\ev}(A)} \leq r$.
\end{cor}
\begin{proof} The remarks following [ADH,~14.0.1] give an immediate asymptotic extension $L$ of~$K$ that is newtonian. Then $L$ is $\d$-valued 
by Lemma~\ref{lem:ADH 14.2.5},   hence~$\abs{\exc^{\ev}(A)} = \abs{\exc_L^{\ev}(A)}\leq r$ by Proposition~\ref{kerexc} and Remark~\ref{rem:kerexc}.
By~[ADH, 13.7.10]  we have $\nwt_{A}(\gamma)\leq 1$ for all $\gamma\in\Gamma$, thus~$\sum_{\gamma\in\Gamma} \nwt_A(\gamma)=\abs{\exc^{\ev}(A)}$.
\end{proof}

%\begin{cor}\label{cor:sum of nwts}
%If $K$ is $\upl$-free, then
%$\sum_{\gamma\in\Gamma} \nwt_A(\gamma)=\abs{\exc^{\ev}(A)}$.
%If $K$ is $\upo$-free, then $\abs{\exc^{\ev}(A)}\leq r$. \marginpar{generalized}
%\end{cor}
%\begin{proof} 
%If $K$ is $\upl$-free, then
%by~[ADH, 13.7.10]  we have $\nwt_{A}(\gamma)\leq 1$ for all $\gamma\in\Gamma$, so~$\sum_{\gamma\in\Gamma} \nwt_A(\gamma)=\abs{\exc^{\ev}(A)}$.
%Now suppose $K$ is $\upo$-free.
%The remarks following~[ADH, 14.0.1] give an immediate  asymptotic extension $L$ of~$K$ that is newtonian. Then $L$ is $\d$-valued 
%by Lemma~\ref{lem:ADH 14.2.5},   hence~$\abs{\exc^{\ev}(A)} = \abs{\exc_L^{\ev}(A)}\leq r$ by Proposition~\ref{kerexc} and the remark following it.
%\end{proof}

\noindent
In [ADH, Section~11.1] we defined ${v_A^{\ev}\colon\Gamma\to\Gamma}$ by requiring that for all $\gamma\in\Gamma$:
\begin{equation}\label{eq:vAev}
v_{A^\phi}(\gamma)\ =\ v_A^{\ev}(\gamma)+\nwt_A(\gamma)v\phi,\qquad\text{ eventually.}
\end{equation}
We recall from that reference that for $a\neq 0$ and $\gamma\in\Gamma$ we have 
$$v_{aA}^{\ev}(\gamma)\ =\ va+v_A^{\ev}(\gamma),\qquad v_{Aa}^{\ev}(\gamma)=v_A^{\ev}(va+\gamma).$$
As an example from [ADH, p.~481], $v_{\der}^{\ev}(\gamma)=\gamma + \psi(\gamma)$ for $\gamma\in \Gamma\setminus \{0\}$ and $v_{\der}^{\ev}(0)=0$.
By [ADH, 14.2.7 and the remark preceding it] we have:

\begin{lemma}\label{lem:ADH 14.2.7}
The restriction of $v_A^{\ev}$ to a  function $\Gamma\setminus\exc^{\ev}(A)\to\Gamma$ is strictly increasing, and
$v\big(A(y)\big) = v_A^{\ev}(vy)$ for all $y\in K$ with $vy\in \Gamma\setminus\exc^{\ev}(A)$.
Moreover, if~$K$ is $\upo$-free, then $v_A^{\ev}\big(\Gamma\setminus \exc^{\ev}(A)\big)=\Gamma$.  
\end{lemma}

\noindent
The following is [ADH, 14.2.10] without the hypothesis of $\upo$-freeness:

\begin{cor}\label{cor:14.2.10, generalized} 
Suppose $K$ is $r$-linearly newtonian. Then for each $f\neq 0$ there exists $y\in K^\times$ such that $A(y)=f$,
$vy\notin\exc^{\ev}(A)$, and $v^{\ev}_A(vy)=vf$.
\end{cor}
\begin{proof} If $r=0$, then $\exc^{\ev}(A)=\emptyset$ and our claim is obviously valid. Suppose $r\ge 1$. Then
 $K$ is $\d$-valued by Lemma~\ref{lem:ADH 14.2.5}, and $v(\ker^{\neq} A) = \exc^{\ev}(A)$
by Proposition~\ref{kerexc}, Moreover, by [ADH, 14.2.2], $K$ is $r$-linearly surjective, hence $f\in A(K)$.
Now Corollary~\ref{cor:nonexc sol} yields $y\in K^\times$ with $A(y)=f$ and $vy\notin \exc^{\ev}(A)$.
By Lemma~\ref{lem:ADH 14.2.7} we have $v^{\ev}_A(vy)=v\big(A(y)\big)=vf$.
\end{proof}

\noindent
From the proof of [ADH, 14.2.10] we extract the following:

\begin{cor}\label{cor:ADH 14.2.10 extract}
Suppose $K$ is $r$-linearly newtonian with small derivation, and $A\in~\mathcal O[\der]$ with $a_0:=A(1)\asymp 1$,
and $f\asymp^\flat 1$. Then there is $y\in K^\times$ such that~${A(y)=f}$ and $y\sim f/a_0$. For any such $y$ we have $vy\notin\exc^{\ev}(A)$
and $v_A^{\ev}(vy)=vf$.
\end{cor}
\begin{proof} The case $r=0$ is trivial. 
Assume $r\geq 1$, so $K$ is $\d$-valued by Lemma~\ref{lem:ADH 14.2.5}. Hence $f^\dagger\prec 1$, that is, $f'\prec f$,
so $f^{(n)}\prec f$ for all $n\ge 1$ by [ADH, 4.4.2].  Then $Af\preceq~f$ by [ADH, (5.1.3), (5.1.2)], and $A(f)\sim a_0f$, so 
$A_{\ltimes f}\in\mathcal O[\der]$ and~$A_{\ltimes f}(1)\sim~a_0$.  Thus we may replace $A$, $f$ by $A_{\ltimes f}$, $1$ to arrange~${f=1}$.
Now~${a_0\asymp 1}$ gives $\operatorname{dwm}(A)=0$, so $\operatorname{dwt}(A^\phi)=0$ eventually, by [ADH, 11.1.11(ii)],
that is, $\nwt(A)=0$. Also $A^\phi(1)=A(1)=a_0\asymp 1$,  so  $v^{\ev}(A)=0$.  Arguing as in the proof of [ADH, 14.2.10] we obtain
$y\in K^\times$ with $A(y)=1$ and~$y\sim 1/a_0$. It is clear that~$vy=0\notin \exc^{\ev}(A)$ and $v_A^{\ev}(vy)=v^{\ev}(A)=0=vf$
for any such $y$.
\end{proof}

\noindent
In the next few subsections below we consider more closely the case of order $r=1$, and in the last subsection the case of arbitrary order.

\subsection*{First-order operators}
{\em In this subsection $A=\der-g$}. By [ADH, p.~481],   
$$\exc^{\ev}(A)\ =\ \exc^{\ev}_K(A)\ =\ \big\{vy:\, y\in K^\times,\ v(g-y^\dagger)>\Psi\big\}$$ 
has at most one element. 
We also have
$\abs{v(\ker^{\neq}A)} = \dim_C\ker A \leq 1$ in view of~$C^\times\subseteq\mathcal O^\times$. Proposition~\ref{kerexc} holds under a weaker assumption on $K$ for $r=1$: 

\begin{lemma}\label{lem:v(ker)=exc, r=1}
Suppose $\I(K)\subseteq K^\dagger$. Then $v(\ker^{\neq} A)=\exc^{\ev}(A)$.
\end{lemma}
\begin{proof}
It remains to show ``$\supseteq$''.
Suppose $\exc^{\ev}(A)=\{0\}$. Then $g-y^\dagger\in \I(K)$ with~${y\asymp 1}$ in $K$, hence
$g\in \I(K)\subseteq K^\dagger$, so $g=h^\dagger$ with $h\asymp 1$, and thus $0=vh\in v(\ker^{\neq} A)$. The general case reduces to the case $\exc^{\ev}(A)=\{0\}$ by twisting.
\end{proof}

\begin{lemma}\label{lemexc, order 1}
Suppose $L$ is an ungrounded $H$-asymptotic extension of $K$. Then  
$\exc^{\ev}_{L}(A) \cap \Gamma  = \exc^{\ev}(A)$.
\end{lemma}
\begin{proof}
Lemma~\ref{lemexc} gives
$\exc^{\ev}_{L}(A) \cap \Gamma \subseteq \exc^{\ev}(A)$.
Next, let $vy\in \exc^{\ev}(A)$, $y\in K^\times$.
Then $v(g-y^\dagger)>\Psi$ and so $v(g-y^\dagger)\in (\Gamma^>)'$ since $K$ has asymptotic integration. 
Hence $v(g-y^\dagger)>\Psi_L$ and thus $vy\in\exc^{\ev}_{L}(A)$, by [ADH, p.~481].
\end{proof}

\noindent
Recall also from [ADH,~9.7] that for an ordered abelian group $G$ and $U \subseteq G$, a function~$\eta\colon U \to G$ is said to be {\it slowly varying}\/ if $\eta(\alpha)-\eta(\beta) = o(\alpha-\beta)$ for all~$\alpha\neq\beta$ in~$U$;\index{function!slowly varying} then the function
$\gamma\mapsto\gamma+\eta(\gamma)\colon U\to G$ is strictly increasing. The quintessential example of a slowly varying function
is $\psi\colon \Gamma^{\neq}\to\Gamma$   [ADH, 6.5.4(ii)]. 

\begin{prop}\label{prop:slow}
There is a unique slowly varying function $\psi_A\colon \Gamma\setminus\exc^{\ev}(A)\to\Gamma$
such that for all $y\in K^\times$ with $vy\notin \exc^{\ev}(A)$ we have $v\big(A(y)\big)=vy+\psi_A(vy)$.
\end{prop}
\begin{proof}
For $\d$-valued $K$, use \cite[8.4]{VDF}. In general, pass to the $\d$-valued hull $L:=\operatorname{dv}(K)$ of $K$ from 
[ADH,~10.3] and use $\Gamma_{L}=\Gamma$ [ADH, 10.3.2].
\end{proof}

\noindent
If $b\neq 0$, then $\exc^{\ev}(A_{\ltimes b})=\exc^{\ev}(A)-vb$ and $\psi_{A_{\ltimes b}}(\gamma)=\psi_A(\gamma+vb)$ for $\gamma\in\Gamma\setminus \exc^{\ev}(A_{\ltimes b})$.

\begin{example}
We have
$\exc^{\ev}(\der)=\{0\}$ and
$\psi_\der=\psi$.
More generally, if $g=b^\dagger$, $b\neq 0$, then~$A_{\ltimes b}=\der$ and so $\exc^{\ev}(A)=\{vb\}$ and
$\psi_A(\gamma)=\psi(\gamma-vb)$ for $\gamma\in\Gamma\setminus\{vb\}$. 
%Hence if $g\in K^\dagger$, then the function $\psi_A$ is quantifier-free definable  in the  asymptotic couple~$(\Gamma,\psi)$ of~$K$. (Here and below, ``definable'' means ``definable with parameters''.)
\end{example}

\noindent
If $\Gamma$ is divisible, then $\Gamma\setminus v\big(A(K)\big)$ has at most one element by [ADH,  11.6.16]. Also,
~$K$ is $\upl$-free iff $v\big(A(K)\big)=\Gamma_\infty$ for all $A=\der-g$ by [ADH, 11.6.17].

\begin{lemma}\label{lem:slow}
Suppose $K$ is $\upl$-free and $f\neq 0$. Then for some $y\in K^\times$ with ${vy\notin\exc^{\ev}(A)}$ we have $A(y)\asymp f$.
\textup{(}Hence $\gamma\mapsto \gamma+\psi_A(\gamma)\colon \Gamma\setminus\exc^{\ev}(A)\to\Gamma$ is surjective.\textup{)}
\end{lemma}
\begin{proof}{} [ADH, 11.6.17] gives $y\in K^\times$ with $A^\phi y\asymp f$ eventually. Now
$$A^\phi y\ =\ \phi y\derdelta-(g-y^\dagger)y\ \text{ in }K^\phi[\derdelta],\qquad \derdelta:=\phi^{-1}\der.$$ 
Since $v(A^{\phi} y)=vf$ eventually, this forces $g-y^\dagger\succ\phi$ eventually, so $vy\notin \exc^{\ev}(A)$.
\end{proof}

%\noindent
%{\em In the rest of this subsection $K$ has small derivation}. (We keep of course the previous assumptions on $K$.)
%Below we assume $A(y)=f$, with $y,f\in K^\times$. 

\noindent
Call $A$ {\bf steep} if $g\succ^\flat 1$, that is, $g\succ 1$ and $g^\dagger\succeq 1$. If $K$ has small derivation and
$A$ is steep, then $g^\dagger\prec g$ by [ADH, 9.2.10]. \index{steep!linear differential operator}\index{linear differential operator!steep}

\begin{lemma}\label{lem:order1, 2}
Suppose $K$ has small derivation,  $A$ is steep, and $y\in K^\times$ such that~$A(y)=f\ne 0$, $g \succ f^\dagger$, and $vy\notin \exc^{\ev}(A)$. Then $y\sim -f/g$.
\end{lemma}
\begin{proof}
We have $$(f/g)^\dagger-g=f^\dagger-g^\dagger-g\sim-g \succ g^\dagger,$$ hence $v(f/g)\notin \exc^{\ev}(A)$, and
$$A(f/g)\ =\ (f/g)'-(f/g)g\ =\  (f/g)\cdot \big( f^\dagger-g^\dagger-g \big)\ \sim\ (f/g)\cdot (-g)\ =\ -f.$$
Since $A(y)=f\sim A(-f/g)$ and $vy,v(f/g)\in\Gamma\setminus\exc^{\ev}(A)$, this gives $y = u\cdot f/g$ where $u\asymp 1$, by Proposition~\ref{prop:slow}.
Now $u^\dagger \prec 1\prec g$ and $(f/g)^\dagger=f^\dagger-g^\dagger\prec g$, 
hence~$y^\dagger\prec g$ and so $$f=A(y)=y\cdot(y^\dagger-g)\sim - y g.$$
Therefore $y\sim -f/g$.
\end{proof}

\begin{lemma}\label{prlemexc}
Suppose $K$ has small derivation and $y\in K^\times$ is such that $A(y)=f\ne 0$, $g-f^\dagger\succ^\flat 1$ and $vy\notin\exc^{\ev}(A)$.  Then
$y\sim f/(f^\dagger-g)$.
\end{lemma}
\begin{proof} From $g-f^\dagger\succ 1$ we get $vf\notin\exc^{\ev}(A)$. Now $A(y)=f\prec f(f^\dagger-g)=A(f)$, so $y\prec f$ by [ADH, 5.6.8], and
$v(y/f)\notin\exc^{\ev}(A_{\ltimes f})=\exc^{\ev}(A)-vf$. Since $A_{\ltimes f}=\der-(g-f^\dagger)$ is steep,
Lemma~\ref{lem:order1, 2} applies to $A_{\ltimes f}$, $y/f$, $1$ in the role of $A$, $y$, $f$.
\end{proof}

\noindent
Suppose $K$ is $\upl$-free and $f\ne 0$.  Then [ADH, 11.6.1] gives an active~$\phi_0$ in $K$  with ${f^\dagger-g-\phi^\dagger}\succeq \phi_0$ 
for all $\phi\prec\phi_0$. The convex subgroups $\Gamma^\flat_{\phi}$ of $\Gamma$ become arbitrarily small as we let $v\phi$ increase cofinally
in $\Psi^{\downarrow}$, so $\phi \prec^\flat_\phi \phi_0$ eventually, and hence~$f^\dagger-g-\phi^\dagger \succ^\flat_\phi \phi$ eventually, that is,
$\phi^{-1}(f/\phi)^\dagger-g/\phi\succ^\flat_\phi 1$ eventually. So replacing $K$ by $K^\phi$, $A$ by 
$\phi^{-1} A^\phi=\derdelta-(g/\phi)$ in $K^\phi[\derdelta]$,
and $f$ and $g$ by~$f/\phi$ and~$g/\phi$, for suitable $\phi$, we arrange $f^\dagger-g\succ^\flat 1$. 
Thus by Lemma~\ref{prlemexc}:

\begin{cor}\label{cor:prlemexc}
If $K$ is $\upl$-free, $y\in K^\times$, $A(y)=f\ne 0$, and $vy\notin\exc^{\ev}(A)$, then~$y\sim f/\big((f/\phi)^\dagger-g\big)$,  eventually.
\end{cor}

\begin{example}
If $K$ is $\upl$-free and $y\in K$, $y'=f\ne 0$ with $y\nasymp 1$, then $y\sim f/(f/\phi)^\dagger$, eventually.
\end{example}

\subsection*{From $K$ to $K[\imag]$} {\em In this subsection $K$ is a real closed $H$-field}. 
Then $K[\imag]$ with~$\imag^2=-1$ is an $H$-asymptotic extension of $K$, with $\Gamma_{K[\imag]}=\Gamma$.
Consider a linear differential operator $B=\der-(g+h\imag)$ over $K[\imag]$.
Note that $g+h\imag\in K[\imag]^\dagger$ iff $g\in K^\dagger$ and~$h\imag\in K[\imag]^\dagger$, by
Lemma~\ref{lem:logder}.
Under further assumptions on $K$, the next two results give explicit descriptions of $\psi_B$ when $g\in K^\dagger$. 

\begin{prop}\label{prop:psiB} Suppose $K[\imag]$ is $1$-linearly newtonian and $g\in K^\dagger$. Then: \begin{enumerate}
\item[$\rm(i)$] if $h\imag\in K[\imag]^\dagger$, then for some $\beta\in \Gamma$ we have $$\exc^{\ev}(B)\ =\ \{\beta\}, \qquad \psi_B(\gamma)\ =\ \psi(\gamma-\beta)\ \text{ for all }\gamma\in\Gamma\setminus\{\beta\};$$
\item[$\rm(ii)$] if $h\imag\notin K[\imag]^\dagger$ and $g=b^\dagger$, 
$b\ne 0$, then 
$$\exc^{\ev}(B)\ =\ \emptyset, \qquad 
\psi_B(\gamma)\ =\ \min\!\big( \psi(\gamma-vb),vh\big)\  \text{ for all }\gamma\in\Gamma.$$
\end{enumerate}
\end{prop}
\begin{proof} As to (i), apply the example following Proposition~\ref{prop:slow} to $K[\imag]$, $B$, $g+h\imag$ in the roles
of $K$, $A$, $g$. For (ii), assume $h\imag\notin K[\imag]^\dagger$, $g=b^\dagger$, $b\ne 0$. 
Replacing $B$ by~$B_{\ltimes b}$ we arrange $g=0$, $b=1$, $B=\der-h\imag$.
Corollary~\ref{cor:logder} gives $K[\imag]^\dagger=K^\dagger\oplus \I(K)\imag$, so~$h\notin \I(K)$, and thus $vh\in \Psi^{\downarrow}$.
Let $y\in K[\imag]^\times$, and take  $z\in K^\times$ and $s\in\I(K)$ with $y^\dagger=z^\dagger+s\imag$. Then
$vh<vs$, hence
$$v(y^\dagger-h\imag)\ =\ \min\!\big( v(z^\dagger), v(s-h) \big)\ =\ \min\!\big( v(z^\dagger), vs, vh \big)\ =\
 \min\!\big( v(y^\dagger), vh \big),$$
 where the last equality uses $v(y^\dagger)=\min \big(v(z^\dagger), vs\big)$. Thus $v(y^\dagger-h\imag)\in \Psi^{\downarrow}$ and
$$v\big(B(y)\big)-vy\ =\ v(y^\dagger-h\imag)\ =\ \min\!\big(v(y^\dagger),vh\big)\ =\ \min\!\big(\psi(vy),vh\big),$$
which gives the desired result.
\end{proof}

\begin{cor}\label{cor:psiB} Suppose $K$ is $\upo$-free, $g\in K^\dagger$, $g=b^\dagger$, $b\ne 0$. Then either for some $\beta\in \Gamma$ we have $\exc^{\ev}(B) = \{\beta\}$ and $\psi_B(\gamma) = \psi(\gamma-\beta)$ for all $\gamma\in\Gamma\setminus\{\beta\}$, or~$\exc^{\ev}(B) = \emptyset$ and
$\psi_B(\gamma) =\min\!\big( \psi(\gamma-vb),vh\big)$  for all $\gamma\in\Gamma$.
\end{cor}
\begin{proof}By [ADH,~14.0.1 and remarks following it] we have an immediate newtonian extension $L$ of~$K$. Then $L$ is still a real closed $H$-field   [ADH, 10.5.8, 3.5.19], and $L[\imag]$ is newtonian   [ADH, 14.5.7], so Proposition~\ref{prop:psiB} applies to $L$ in place of $K$. 
\end{proof}

\subsection*{Higher-order operators} We begin with the following observation:

\begin{lemma}\label{lem:exce product}
Let $B\in K[\der]^{\neq}$ and $\gamma\in\Gamma$. Then
$\nwt_{AB}(\gamma)\geq \nwt_B(\gamma)$, and
$$ \gamma\notin\exc^{\ev}(B)\ \Longrightarrow\ 
 \nwt_{AB}(\gamma)\ =\ \nwt_A\!\big( v^{\ev}_B(\gamma) \big) \text{ and }\ v^{\ev}_{AB}(\gamma)\ =\ v^{\ev}_A\big(v^{\ev}_B(\gamma)\big).$$
\end{lemma}
\begin{proof}
We have 
$\nwt_{AB}(\gamma)=\dwt_{(AB)^\phi}(\gamma)$ eventually, and 
$(AB)^\phi=A^\phi B^\phi$. Hence by [ADH, Section~5.6] and the definition of $v^{\ev}_B(\gamma)$ in  \eqref{eq:vAev}:
\begin{align*}
\nwt_{AB}(\gamma)\	&=\ \dwt_{A^\phi}\!\big(v_{B^\phi}(\gamma)\big)+\dwt_{B^\phi}(\gamma) \\
					&=\ \dwt_{A^\phi}\!\big(v_B^{\ev}(\gamma)+\nwt_B(\gamma)v\phi\big) + \nwt_B(\gamma), \text{ eventually}, 
\end{align*}
so $\nwt_{AB}(\gamma)\geq \nwt_B(\gamma)$. Now suppose  $\gamma\notin \exc^{\ev}(B)$. Then $\nwt_B(\gamma)=0$, so  
$$\nwt_{AB}(\gamma) = \dwt_{A^\phi}\!\big(v_B^{\ev}(\gamma)\big)   = \nwt_A\!\big(v_B^{\ev}(\gamma)\big), \qquad\text{eventually.}$$ 
Moreover, $v_{(AB)^\phi}= v_{A^\phi B^\phi} =v_{A^\phi}\circ v_{B^\phi}$, hence using   \eqref{eq:vAev}: 
$$v_{(AB)^\phi}(\gamma)\ =\ v_{A^\phi}\big(v_{B^\phi}(\gamma)\big)\ =\ v_{A^\phi}\big(v_{B}^{\ev}(\gamma)\big), \text{ eventually},$$
and thus eventually 
\begin{align*}
v^{\ev}_{AB}(\gamma)\ 	&=\  
     v_{(AB)^\phi}(\gamma)-\nwt_{AB}(\gamma)v\phi \\
						&=\  v_{A^\phi}\big(v_{B}^{\ev}(\gamma)\big) - \nwt_A\!\big(v_B^{\ev}(\gamma)\big)v\phi\ =\ 
v^{\ev}_A\big(v^{\ev}_B(\gamma)\big).\qedhere
\end{align*}
\end{proof} 

\noindent
Lemmas~\ref{lem:ADH 14.2.7} and~\ref{lem:exce product} yield: 

\begin{cor}\label{cor:exce product}
Let $B\in K[\der]^{\neq}$. Then
$$\exc^{\ev}(AB)\  =\  (v^{\ev}_B)^{-1}\big( \exc^{\ev}(A) \big) \cup \exc^{\ev}(B)$$
and hence
$\abs{\exc^{\ev}(AB)} \leq \abs{\exc^{\ev}(A)}+\abs{\exc^{\ev}(B)}$, 
with equality if $v^{\ev}_B\big(\Gamma\setminus\exc^{\ev}(B)\big)=\Gamma$.
\end{cor}

\noindent
As an easy consequence we have a variant of Corollary~\ref{cor:sum of nwts}:

\begin{cor}\label{cor:size of excev}
If $A$ splits over $K$, then $\abs{\exc^{\ev}(A)}\leq r$. 
\end{cor}

\noindent
To study $v^{\ev}_A$ in more detail we introduce the function 
$$ \psi_A\ \colon\ \Gamma\setminus\exc^{\ev}(A)\to\Gamma, \qquad \gamma\mapsto v_A^{\ev}(\gamma)-\gamma.$$
For monic $A$ of order $1$ this agrees with $\psi_A$ as defined in Proposition~\ref{prop:slow}.  
For~$A=a$ ($a\neq 0$) we have $\exc^{\ev}(A)=\emptyset$ and  $\psi_A(\gamma)=va$ for all $\gamma\in\Gamma$.

\begin{lemma}\label{lem:psiAB}
Let  $B\in K[\der]^{\neq}$ and $\gamma\in\Gamma\setminus\exc^{\ev}(AB)$. Then  
$$\psi_{AB}(\gamma)\ =\ \psi_A\big(v_B^{\ev}(\gamma)\big)+\psi_B(\gamma).$$
\end{lemma}
\begin{proof}
We have $\gamma\notin\exc^{\ev}(B)$ and $v_B^{\ev}(\gamma)\notin\exc^{\ev}(A)$ by Corollary~\ref{cor:exce product}, hence
$$\psi_{AB}(\gamma)\ =\ v_A^{\ev}\big(v_B^{\ev}(\gamma)\big)-\gamma\ =\  v_B^{\ev}(\gamma)+\psi_A\big(v_B^{\ev}(\gamma)\big) - \gamma\ =\ \psi_A\big(v_B^{\ev}(\gamma)\big)+\psi_B(\gamma)$$
by Lemma~\ref{lem:exce product}. 
\end{proof}

\noindent
Thus for $a\neq 0$ and $\gamma\in\Gamma$ we have
$$\psi_{aA}(\gamma)=va+\psi_A(\gamma)\text{ if $\gamma\notin\exc^{\ev}(A)$,}\quad \psi_{Aa}(\gamma)=\psi_A(va+\gamma)+va \text{ if 
$\gamma\notin\exc^{\ev}(A)-va$.}$$

\begin{example} Suppose $K$ has small derivation and $x\in K$, $x'\asymp 1$. Then $vx<0$ and~$\exc^{\ev}(\der^2)=\{vx,0\}$, and~$\psi_{\der^2}(\gamma)=\psi\big(\gamma+\psi(\gamma)\big)+\psi(\gamma)$ for $\gamma\in\Gamma\setminus\exc^{\ev}(\der^2)$.
\end{example}

\begin{lemma}\label{lem:v(A(y)) convex subgp}
Suppose $\psi_A$ is slowly varying. Let $\Delta$ be a convex subgroup of $\Gamma$ and let  $y,z\in K^\times$ be such that $vy,vz\notin\exc^{\ev}(A)$.
Then  $$v_\Delta (y) < v_\Delta (z) \ \Longleftrightarrow\ v_\Delta\big(A(y)\big) < v_\Delta\big(A(z)\big).$$ 
\end{lemma}
\begin{proof}
By Lemma~\ref{lem:ADH 14.2.7} we have 
$$v\big(A(y)\big)-v\big(A(z)\big)\ =\ v_A^{\ev}(vy)-v_A^{\ev}(vz)\ =\ vy-vz+\psi_A(vy)-\psi_A(vz)$$
and $\psi_A(vy)-\psi_A(vz)=o(vy-vz)$ if $vy\ne vz$.
\end{proof}

\noindent
Call $A$ {\bf asymptotically surjective} \index{asymptotically surjective} \index{linear differential operator!asymptotically surjective}
if $v_A^{\ev}\big( \Gamma\setminus\exc^{\ev}(A) \big) = \Gamma$  and 
$\psi_A$ is slowly varying.
If $A$ is asymptotically surjective, then so are~$aA$ and $Aa$ for $a\neq 0$, and
if $A$ has order $0$, then $A$ is asymptotically surjective.
If $K$ is $\upl$-free and $A$ has order $1$, then~$A$ is asymptotically surjective, thanks to Proposition~\ref{prop:slow} and Lemma~\ref{lem:slow}. 
%Together with Corollary~\ref{cor:exce product}, this yields:

\noindent
The next lemma has an obvious proof.

\begin{lemma}\label{lem:slowly varying}
Let $G$ be an ordered abelian group and $U, V\subseteq G$. If $\eta_1,\eta_2\colon U\to G$ are slowly varying,
then so is $\eta_1+\eta_2$. If $\eta\colon U\to G$ and $\zeta\colon V\to G$ 
are slowly varying and $\gamma+\zeta(\gamma)\in U$ for all $\gamma\in V$, then the function~$\gamma\mapsto \eta\big(\gamma+\zeta(\gamma)\big)\colon V\to G$ is also slowly varying.
\end{lemma}

\begin{lemma}
If $A$ and $B\in K[\der]^{\neq}$ are asymptotically surjective, then so is $AB$. 
\end{lemma}
\begin{proof}
Let $A$,~$B$ be asymptotically surjective and $\gamma\in \Gamma$. This gives $\alpha\in\Gamma\setminus\exc^{\ev}(A)$
with $v_A^{\ev}(\alpha)=\gamma$ and $\beta\in\Gamma\setminus\exc^{\ev}(B)$
with $v_B^{\ev}(\beta)=\alpha$. Then $\beta\notin\exc^{\ev}(AB)$
by Corollary~\ref{cor:exce product}, and $v_{AB}^{\ev}(\beta)=\gamma$ by Lemma~\ref{lem:exce product}. Moreover, $\psi_{AB}$ is slowly varying
by Lemmas~\ref{lem:psiAB} and \ref{lem:slowly varying}. 
\end{proof}

\noindent
A straightforward induction on $r$ using this lemma yields:

\begin{cor}\label{cor:well-behaved}
If $K$ is $\upl$-free and $A$ splits over $K$, then $A$ is asymptotically surjective. 
\end{cor}

\noindent
We can now add to Lemma~\ref{lem:ADH 14.2.7}:

\begin{cor}\label{cor1524}
Suppose $K$ is $\upo$-free. Then $A$ is asymptotically surjective.
\end{cor}
\begin{proof}
%Note that $K$ has rational asymptotic integration  [ADH, p.~515].
By the second part of Lemma~\ref{lem:ADH 14.2.7} it is enough to show that $\psi_A$ is slowly varying.
For this we may replace $K$ by any $\upo$-free $H$-asymptotic extension $L$ of~$K$ with $\Psi$ cofinal in $\Psi_L$.
Thus we can arrange by~[ADH, 14.5.7, remarks following 14.0.1]  that~$K$ is newtonian, and by passing
to the algebraic closure,  algebraically closed. Then~$A$ splits over $K$ by [ADH, 5.8.9, 14.5.3],
and so $A$ is asymptotically surjective by Corollary~\ref{cor:well-behaved}.
\end{proof}

\section{Special Elements}\label{sec:special elements}

\noindent
Let $K$ be a valued field and let $\hat a$ be an element of an immediate extension of~$K$ with~$\hat a\notin K$.
{\samepage Recall that $$v(\hat a-K)\ =\ \big\{v(\hat a-a):a\in K\big\}$$ is a nonempty downward closed subset of $\Gamma:= v(K^\times)$ without a largest element.}
Call~$\hat a$ {\it special}\/ over~$K$ if \index{special} \index{element!special}
some nontrivial  convex subgroup of~$\Gamma$ is cofinal in~$v({\hat a-K})$~[ADH, p.~167].
In this case $v(\hat a-K)\cap\Gamma^>\neq\emptyset$, and there is a unique such nontrivial convex subgroup $\Delta$ of $\Gamma$,
namely
$$\Delta\ =\ \big\{ \delta\in\Gamma:\, \abs{\delta}\in v(\hat a-K) \big\}.$$
We also call~$\hat a$ {\it almost special}\/ over $K$ if $\hat a/\fm$ is special over $K$ for some~$\fm\in K^\times$. \index{special!almost} \index{element!almost special}
If $\Gamma\neq\{0\}$ is archimedean, then $\hat a$ is special over $K$ iff $v(\hat a-K)=\Gamma$, iff $\hat a$ is the limit of a divergent c-sequence in $K$. (Recall that ``c-sequence'' abbreviates ``cauchy sequence'' [ADH, p.~82].)
In the next lemma $a$ ranges over $K$ and $\fm$, $\fn$ over $K^\times$.

\begin{lemma}\label{lem:special refinement}
Suppose $\hat a\prec \fm$ and $\hat a/\fm$ is special over $K$. Then for all $a$, $\fn$, if~$\hat a-a\prec\fn\preceq\fm$, then $(\hat a-a)/\fn$ is special over $K$.
\end{lemma}
\begin{proof}
Replacing $\hat a$, $a$, $\fm$, $\fn$ by $\hat a/\fm$, $a/\fm$, $1$, $\fn/\fm$, respectively, we arrange $\fm=1$. So let $\hat a$ be special over $K$ with $\hat a \prec 1$. It is enough to show: (1)~$\hat a-a$ is special over $K$, for all $a$; (2)~for all $\fn$, if $\hat a\prec\fn\preceq 1$, then $\hat a/\fn$ is special over $K$. Here~(1) follows from $v(\hat a-a-K)=v(\hat a-K)$. For (2), note that if $\hat a\prec\fn\preceq 1$, 
then $v\fn\in\Delta$ with $\Delta$ as above, and so $v(\hat a/\fn-K)=v(\hat a-K)-v\fn=v(\hat a-K)$.
\end{proof}

%\begin{remark} \marginpar{new remark} Let $\hat a$, $\fm$ be as in the previous lemma, and let $\Delta$ be the convex subgroup of $\Gamma$ which is cofinal in $v(\hat a/\fm-K)$. The proof of the previous lemma shows that then $\Delta$ is also cofinal in $v\big((\hat a-a)/\fn-K\big)$. \end{remark}

\noindent
The remainder of this section is devoted to showing that (almost) special elements   arise naturally in the analysis of
certain immediate $\d$-algebraic extensions of valued differential fields.  
We first treat the case of asymptotic fields with small derivation, and then 
focus on the linearly newtonian $H$-asymptotic case. 

We  recall some notation: for an ordered abelian group $\Gamma$ and  $\alpha\in \Gamma_{\infty}$, 
$\beta\in \Gamma$, $\gamma\in \Gamma^{>}$
we mean by ``$\alpha \ge \beta+o(\gamma)$'' that
$\alpha\ge \beta-(1/n)\gamma$ for all $n\ge 1$, while~``$\alpha < \beta +o(\gamma)$'' is its negation, that is,
$\alpha < \beta-(1/n)\gamma$ for some~$n\ge 1$; see [ADH, p.~312].  Here and later inequalities are
in the sense of the ordered divisible hull $\Q\Gamma$ of the relevant $\Gamma$. 

\subsection*{A source of special elements} We recall that a differential field $F$ is said to be $r$-linearly surjective ($r\in \N$) if $A(F)=F$ for every nonzero $A\in F[\der]$ of order at most $r$.  {\it In this subsection $K$ is an asymptotic field with small derivation, value group $\Gamma=v(K^\times)\ne \{0\}$, and differential residue field $\k$; we also let $r\in \N^{\ge 1}$.}\/ Below we use the notion {\em neatly surjective\/} from~[ADH,~5.6]: $A\in K[\der]^{\ne}$ is neatly surjective iff for all $b\in K^\times$ there exists $a\in K^\times$ such that $A(a)=b$ and $v_A(va)=vb$. We often let $\hat f$ be an element in an immediate asymptotic extension $\hat K$ of $K$, but in the statement of the next lemma we take $\hat f\in K$:

\begin{lemma}\label{neat1} Assume $\k$ is $r$-linearly surjective, $A\in K[\der]^{\ne}$ of order $\le r$ is neatly surjective, $\gamma\in \Q\Gamma$, $\gamma>0$, $\hat f\in K^\times$, and $v\big(A(\hat f)\big)\ge v(A\hat f)+\gamma$. Then $A(f)=0$ and~$v(\hat f-f)\ge v(\hat f)+ \gamma+o(\gamma)$ for some $f\in K$. 
\end{lemma}
\begin{proof} Set $B:= g^{-1}A\hat f$, where we take $g\in K^\times$ such that $vg=v(A\hat f)$. Then~$B\asymp 1$, $B$ is still neatly surjective, and $B(1)=g^{-1}A(\hat f)$, $v\big(B(1)\big)\ge \gamma$. It suffices to find~$y\in K$ such that $B(y)=0$ and $v(y-1)\ge \gamma+o(\gamma)$, because then $f:=\hat fy$ has the desired property.
If $B(1)=0$, then $y=1$ works, so assume $B(1)\ne 0$. By~[ADH, 7.2.7] we have an immediate extension $\hat{K}$ of $K$ that is $r$-differential henselian. Then~$\hat{K}$ is asymptotic by [ADH, 9.4.2 and 9.4.5]. 
Set $R(Z):= \Ric(B)\in K\{Z\}$.
Then the proof of [ADH, 7.5.1] applied to $\hat{K}$ and~$B$ in the roles
of $K$ and $A$ yields~$z\prec 1$ in $\hat{K}$ with $R(z)=0$. Now $R(0)=B(1)$, 
%so~$v\big(R(0)\big)\ge \gamma$, and 
hence by [ADH, 7.2.2] we can take such $z$ with $v(z)  \ge  \beta+o(\beta)$ where $\beta:=v\big(B(1)\big)\ge \gamma$. 
%\ \ge\ \gamma+o(\gamma).$$ 
As in the proof of
[ADH, 7.5.1] we next take $y\in \hat{K}$ with $v(y-1)>0$ and
$y^\dagger=z$ to get~$B(y)=0$, and observe that then
$v(y-1)\ge \beta+o(\beta)$,  by [ADH, 9.2.10(iv)], hence $v(y-1)\ge \gamma+o(\gamma)$. It remains to note that $y\in K$ by [ADH, 7.5.7].  
\end{proof} 

\noindent
By a remark following the proof of [ADH, 7.5.1] the assumption that $\k$ is $r$-linearly surjective in the lemma above can be replaced for $r\ge 2$ by the assumption that $\k$ is $(r-1)$-linearly surjective. 

\medskip
\noindent
Next we establish a version of the above with $\hat f$ in an immediate asymptotic extension of $K$. Recall that an asymptotic extension of $K$ with the same value group as~$K$ has small derivation, by [ADH, 9.4.1].  

\begin{lemma}\label{neat2} Assume $\k$ is $r$-linearly surjective, $A\in K[\der]^{\ne}$ of order $\le r$ is neatly surjective, $\gamma\in \Q\Gamma$, $\gamma > 0$, $\hat{K}$ is an immediate asymptotic extension of $K$, $\hat f\in \hat K^\times$, and~$v\big(A(\hat f)\big)\ge v(A\hat f)+\gamma$. Then for some $f\in K$ we have
$$A(f)\ =\ 0, \qquad v(\hat f-f)\ \ge\ v(\hat f)+\gamma+o(\gamma).$$ 
\end{lemma} 
\begin{proof} By extending $\hat{K}$ we can arrange that 
$\hat{K}$ is $r$-differential henselian, 
so $A$ remains neatly surjective as an element of $\hat{K}[\der]$, by [ADH, 7.1.8]. Then by Lem\-ma~\ref{neat1} with $\hat{K}$ in the role of $K$ we get $f\in \hat{K}$ such that $A(f)=0$ and
$v(\hat f-f)\ge v(\hat f)+\gamma+o(\gamma)$. It remains to note that
$f\in K$ by [ADH, 7.5.7]. 
\end{proof}

\noindent
We actually need an inhomogeneous variant of the above: 

\begin{lemma}\label{neat3} Assume $\k$ is $r$-linearly surjective, $A\in K[\der]^{\ne}$ of order $\le r$ is neatly surjective, $b\in K$, $\gamma\in \Q\Gamma$, $\gamma>0$, $v(A)=o(\gamma)$, $v(b)\ge o(\gamma)$, $\hat{K}$ is an immediate asymptotic extension of $K$, $\hat f\in \hat K$,
$\hat f\preceq 1$, and $v\big(A(\hat f)-b\big)\ge \gamma+o(\gamma)$. 
Then 
$$A(f)\ =\ b,\qquad v(\hat f -f)\ \ge\ (1/2)\gamma +o(\gamma)$$
for some $f\in K$. 
\end{lemma}
\begin{proof} Take $y\in K$ with $A(y)=b$ and $v(y)\ge o(\gamma)$. Then $A(\hat g)=A(\hat f)-b$ for~$\hat g:= \hat f-y$, so $v\big(A(\hat g)\big)\ge \gamma+o(\gamma)$ and $v(\hat g)\ge o(\gamma)$. We distinguish two cases:

\medskip\noindent
(1) $v(\hat g)\ge (1/2)\gamma+o(\gamma)$. Then $v(\hat f - y)\ge (1/2)\gamma+o(\gamma)$, so $f:= y$ works. 

\medskip\noindent
(2) $v(\hat g) < (1/2)\gamma+o(\gamma)$. Then by [ADH, 6.1.3],
$$v(A\hat g)\ <\ (1/2)\gamma+o(\gamma), \qquad v\big(A(\hat g)\big)\ \ge\ \gamma+o(\gamma),$$ so
$v\big(A(\hat g)\big)\ge v(A\hat g) + (1/2)\gamma$. Then 
Lemma~\ref{neat2}
gives an element $g\in K$ such that~$A(g)=0$ and 
$v(\hat g-g)\ge (1/2)\gamma+o(\gamma)$. Hence $f:= y+g$ works. 
\end{proof}

\noindent
Recall from [ADH, 7.2] that $\O$ is said to be {\it $r$-linearly surjective} if for every $A$ in~$K[\der]^{\neq}$ of order $r$ with $v(A)=0$ there exists $y\in\O$ with $A(y)=1$.

\begin{prop}\label{propsp} Assume $\O$ is $r$-linearly surjective,
$P\in K\{Y\}$, $\order(P)\le r$, $\ddeg P=1$, 
and $P(\hat a)=0$, where $\hat a\preceq 1$ lies in an immediate asymptotic extension of $K$ and $\hat a\notin K$. Then $\hat a$ is special over $K$. % some nontrivial convex subgroup of $\Gamma$ is cofinal in $v(\hat a - K)$.
\end{prop}
\begin{proof} The hypothesis on $\O$ yields: $\k$ is $r$-linearly surjective and all $A\in K[\der]^{\ne}$
of order $\le r$ are neatly surjective. Let $0 <\gamma\in v(\hat a -K)$; we claim that $v(\hat a -K)$ has an element $\ge (4/3)\gamma$. 
We arrange $P\asymp 1$. Take $a\in K$ with $v(\hat a -a)=\gamma$. Then~${P_{+a}\asymp 1}$, $\ddeg P_{+a}=1$, so 
$$P_{+a,1}\ \asymp\ 1, \quad P_{+a,>1}\ \prec\ 1, \quad P_{+a}\ =\ P(a) + P_{+a,1} + P_{+a,>1}$$
and
$$0\ =\ P(\hat a)\ =\ P_{+a}(\hat a -a)\ =\ P(a) + P_{+a,1}(\hat a -a) + P_{+a,>1}(\hat a -a),$$
with $$v\big(P_{+a,1}(\hat a -a)+P_{+a,>1}(\hat a -a)\big)\ge \gamma+o(\gamma),$$ and thus $v(P(a))\ge \gamma+o(\gamma)$.
Take $g\in K^\times$ with $vg=\gamma$ and set $Q:= g^{-1}P_{+a,\times g}$, so $Q=Q_0+ Q_1 + Q_{>1}$ with
$$Q_0\ =\ Q(0)\ =\ g^{-1}P(a),\quad Q_1\ =\ g^{-1}(P_{+a,1})_{\times g},\quad Q_{>1}\ =\ g^{-1}(P_{+a,>1})_{\times g}, $$
hence 
$$v(Q_0)\ge o(\gamma), \quad v(Q_1)=o(\gamma), \quad v(Q_{>1})\ge \gamma+o(\gamma).$$ 
We set $\hat f:= g^{-1}(\hat a -a)$, so $Q(\hat f)=0$ and $\hat f\asymp 1$, and $A:=L_Q\in K[\der]$. Then~$Q(\hat f)=0$ gives
$$Q_0+A(\hat f)\ =\ Q_0+Q_1(\hat f)\ =\ -Q_{>1}(\hat f),\ \text{ with }v\big(Q_{>1}(\hat f)\big)\ \ge\ \gamma+o(\gamma),$$ so $v\big(Q_0+A(\hat f)\big)\ge \gamma+o(\gamma)$. Since $v(A)=v(Q_1)=o(\gamma)$, Lemma~\ref{neat3} then gives~$f\in K$ with $v(\hat f - f)\ge (1/3)\gamma$.
In view of $\hat a-a=g\hat f$, this yields
$$v\big(\hat a -(a+gf)\big)\ =\ \gamma+v(\hat f - f)\ \ge\ (4/3)\gamma,$$ 
which proves our claim. It gives the desired result.    
\end{proof}

\subsection*{A source of almost special elements} {\it In this subsection
$K$, $\Gamma$, $\k$, and $r$ are as in the previous subsection, and we assume that $\O$ is
$r$-linearly surjective.}\/ (So $\k$ is $r$-linearly surjective, and $\sup \Psi=0$ by [ADH, 9.4.2].) 
Let 
$\hat a$ be an element in an immediate asymptotic extension of $K$ such that $\hat a\notin K$ and 
$K\<\hat a\>$ has transcendence degree $\le r$
over $K$. We shall use Proposition~\ref{propsp} to show: 
%Towards covering 1.11 in Joris' notes (and more) we show:

\begin{prop} \label{propalsp}  
If $\Gamma$ is divisible, then $\hat a$ is almost special over $K$. % $\alpha+\Delta$ is cofinal in $v(\hat a - K)$ for some $\alpha\in \Gamma$ and some nontrivial convex subgroup $\Delta$ of $\Gamma$. 
\end{prop} 

\noindent
Towards the proof we first note that $\hat a$ has a minimal annihilator $P(Y)$
 over $K$ of order $\le r$. We also fix a divergent pc-sequence $(a_{\rho})$ in $K$ such that $a_{\rho}\leadsto \hat a$. We next show how to improve $\hat a$ and $P$ (without assuming divisibility of $\Gamma$):

\begin{lemma}\label{replace} For some $\hat{b}$ in an immediate asymptotic extension of $K$ we have: \begin{enumerate}
\item[\textup{(i)}] $v(\hat a -K)=v(\hat b -K)$;
\item[\textup{(ii)}] $(a_{\rho})$ has a minimal differential polynomial $Q$ over $K$ of order $\le r$ such that~$Q$
is also a minimal annihilator of $\hat b$ over $K$.
\end{enumerate}
\end{lemma}
\begin{proof} By [ADH, 6.8.1, 6.9.2], $(a_{\rho})$ is of $\d$-algebraic type over $K$ with a minimal differential polynomial $Q(Y)$ over $K$ such that $\order Q \le \order P\le r$. 
By [ADH, 6.9.3, 9.4.5] this gives an element $\hat b$ in an immediate asymptotic extension of $K$ such that $Q$ is a minimal annihilator of $\hat b$ over $K$ and $a_{\rho}\leadsto \hat b$. Then $Q$ and $\hat b$ have the desired properties. 
\end{proof}

\begin{proof}[Proof of Proposition~\ref{propalsp}] Replace $\hat a$ and $P$
by $\hat b$ and $Q$ from Lemma~\ref{replace} (and rename) to arrange
that $P$ is a minimal differential polynomial of $(a_{\rho})$ over $K$.
Now assuming $\Gamma$ is divisible,~\cite[Proposition~3.1]{Nigel} gives $a\in K$ and $g\in K^\times$ such that~$\hat a -a\asymp g$ and $\ddeg P_{+a,\times g} = 1$. 

Set $F:= P_{+a,\times g}$ and
$\hat f:= (\hat a  -a)/g$. Then $\ddeg F=1$, $F(\hat f)=0$, and $\hat f\preceq 1$. Applying Proposition~\ref{propsp} to
$F$ and $\hat f$ in the role of $P$ and $\hat a$ yields a
nontrivial convex subgroup $\Delta$ of $\Gamma$ that is
cofinal in $v(\hat f -K)$. Setting $\alpha:= vg$, it follows that~$\alpha+\Delta$ is cofinal in $v\big((\hat a -a)-K\big)=v(\hat a -K)$.   
\end{proof}

\noindent
We can trade the divisibility assumption in Proposition~\ref{propalsp} against a stronger hypothesis on $K$, the proof  using \cite[3.3]{Nigel} instead of \cite[3.1]{Nigel}:

\begin{cor}
If $K$ is henselian and $\k$ is linearly surjective, then $\hat a$ is almost special over~$K$. 
\end{cor}

%\noindent
%The proof is like that of Proposition~\ref{propalsp}, using \cite[3.3]{Nigel} instead of \cite[3.1]{Nigel}.

\subsection*{The linearly newtonian setting} {\em In this subsection $K$ is an $\upo$-free $r$-linearly newtonian $H$-asymptotic field,
$r\ge 1$.}\/ Thus $K$ is
$\d$-valued by Lemma~\ref{lem:ADH 14.2.5}. We let~$\phi$ range over the elements active in $K$.
We now mimick the material in the previous two subsections. Note that for $A\in K[\der]^{\ne}$ and any element $\hat f$ in an asymptotic extension of $K$ we have $A(\hat f)\preceq A^\phi \hat f$, since
$A(\hat f)=A^\phi(\hat f)$.
% and $K^\phi$ has small derivation.  

\begin{lemma}\label{neneat1} Assume that $A\in K[\der]^{\ne}$ has order $\le r$, $\gamma\in \Q\Gamma$, $\gamma>0$, $\hat f\in K^\times$, and $v\big(A(\hat f)\big)\ge v(A^\phi\hat f)+\gamma$, eventually. Then there exists an $f\in K$ such that~$A(f)=0$ and $v(\hat f-f)\ge v(\hat f) +\gamma+o(\gamma)$. 
\end{lemma}
\begin{proof} Take $\phi$ such that $v\phi\ge \gamma^\dagger$ and $v\big(A(\hat f)\big)\ge v(A^\phi\hat f)+\gamma$. 
Next, take $\beta\in \Gamma$ such that  $\beta\ge \gamma$ and $v\big(A(\hat f)\big)\ge v(A^\phi\hat f)+\beta$.
Then $v\phi\ge \beta^\dagger$, so $\beta> \Gamma_{\phi}^{\flat}$, hence the valuation ring of the flattening $(K^\phi, v_{\phi}^\flat)$ is $r$-linearly surjective, by [ADH, 14.2.1]. 
We now apply
Lemma~\ref{neat1} to $$(K^\phi,v_{\phi}^\flat),\quad
A^\phi,\quad \dot{\beta}:= \beta+\Gamma_{\phi}^{\flat}$$ in the role of $K$, $A$, $\gamma$ to give $f\in K$ with
$A(f)=0$ and $v_{\phi}^{\flat}(\hat f -f)\ge v_{\phi}^{\flat}(\hat f) + \dot{\beta} + o(\dot{\beta})$.
Then also $v(\hat f-f)\ge v(\hat f) +\beta+o(\beta)$, and thus $v(\hat f -f)\ge v(\hat f) + \gamma+o(\gamma)$.  
\end{proof}

\begin{lemma}\label{neneat2} Assume $A\in K[\der]^{\ne}$ has order $\le r$,  $\hat{K}$ is an immediate $\d$-algebraic asymptotic extension of $K$, $\gamma\in \Q\Gamma$, $\gamma>0$, $\hat f\in \hat K^\times$, and  $v\big(A(\hat f)\big)\ge v(A^\phi\hat f)+\gamma$ eventually. Then $A(f)=0$ and $v(\hat f-f)\ge v(\hat f) +\gamma+o(\gamma)$ for some $f\in K$.
\end{lemma} 
\begin{proof} Since $K$ is $\upo$-free, so is $\hat K$ by Theorem~\ref{thm:ADH 13.6.1}. By [ADH, 14.0.1 and subsequent remarks] we can extend $\hat{K}$ to arrange that 
$\hat{K}$ is also newtonian. Then by Lemma~\ref{neneat1} with $\hat{K}$ in the role of $K$ we get $f\in \hat{K}$ with $A(f)=0$ and~$v(\hat f-f)\ge v(\hat f)+\gamma+o(\gamma)$. Now use that
$f\in K$ by [ADH, line before 14.2.10]. 
\end{proof}

\begin{lemma}\label{neneat3} Assume $A\in K[\der]^{\ne}$ has order $\le r$, $b\in K$, $\gamma\in \Q\Gamma$, $\gamma>0$, $\hat{K}$ is an immediate $\d$-algebraic asymptotic extension of $K$, and $\hat f\in \hat K$, $v(\hat f)\ge o(\gamma)$. Assume also that
eventually $v(b)\ge v(A^\phi)+o(\gamma)$ and $v\big(A(\hat f)-b\big)\ge v(A^\phi)+\gamma+o(\gamma)$. 
Then for some~$f\in K$ we have $A(f)=b$ and
$v(\hat f -f)\ge (1/2)\gamma+o(\gamma)$.
\end{lemma}
\begin{proof} We take $y\in K$ with $A(y)=b$ as follows: If $b=0$, then $y:=0$. If~$b\ne 0$, then Corollary~\ref{cor:14.2.10, generalized} yields $y\in K^\times$ such that
$A(y)=b$, $vy\notin \exc^{\ev}(A)$, and
$v_A^{\ev}(vy)=vb$. In any case, $vy\ge o(\gamma)$: when $b\ne 0$, the sentence preceding
[ADH, 14.2.7] gives~$v_{A^{\phi}}(vy)=vb$, eventually, 
to which we apply [ADH, 6.1.3].
 
 Now $A(\hat g)=A(\hat f)-b$ for $\hat g:= \hat f-y$, so
$v(\hat g)\ge o(\gamma)$, and eventually $v\big(A(\hat g)\big)\ge v(A^\phi)+\gamma+o(\gamma)$. We distinguish two cases:

\medskip\noindent
(1) $v(\hat g)\ge (1/2)\gamma+o(\gamma)$. Then $v(\hat f - y)\ge (1/2)\gamma+o(\gamma)$, so $f:= y$ works. 

\medskip\noindent
(2) $v(\hat g) < (1/2)\gamma+o(\gamma)$. Then by [ADH, 6.1.3] we have eventually
$$v(A^\phi\hat g)\ <\ v(A^\phi)+(1/2)\gamma+o(\gamma), \qquad v\big(A(\hat g)\big)\ \ge\ v(A^\phi)+\gamma+o(\gamma),$$ so
$v(A(\hat g))\ge v(A^\phi\hat g) + (1/2)\gamma$, eventually.   
Lemma~\ref{neneat2}
gives an ele\-ment~${g\in K}$ with $A(g)=0$ and 
$v(\hat g-g)\ge (1/2)\gamma+o(\gamma)$. Hence $f:= y+g$ works. 
\end{proof}

\begin{prop}\label{nepropsp}\label{prop:hata special} Suppose that
$P\in K\{Y\}$, $\order P\le r$, $\ndeg P=1$, 
and $P(\hat a)=0$, where $\hat a\preceq 1$ lies in an immediate asymptotic extension of $K$ and $\hat a\notin K$. Then $\hat a$ is special over $K$. %some nontrivial convex subgroup of $\Gamma$ is cofinal in $v(\hat a - K)$.
\end{prop}

\noindent
The proof is like that of Proposition~\ref{propsp}, but there are some differences that call for further details.

\begin{proof} Given $0 <\gamma\in v(\hat a -K)$, we claim that $v(\hat a -K)$ has an element $\ge (4/3)\gamma$. 
Take $a\in K$ with $v(\hat a -a)=\gamma$. Then $\ndeg P_{+a}=1$ by [ADH, 11.2.3(i)], so eventually we have
$$P(a)\ \preceq\ P^\phi_{+a,1}\ \succ\ P^\phi_{+a,>1}, \quad P^\phi_{+a}\ =\ P(a) + P^\phi_{+a,1} + P^\phi_{+a,>1}$$
and
\begin{align*} 
0\  =\ P(\hat a) &\ =\ P^\phi_{+a}(\hat a -a) \\ \ &\ =\ P(a) +  P^\phi_{+a,1}(\hat a -a) + P^\phi_{+a,>1}(\hat a -a),\\
&\phantom{=\ P(a)+} v\big(P^\phi_{+a,1}(\hat a -a)+P^\phi_{+a,>1}(\hat a -a)\big)\ \ge\ v(P^\phi_{+a,1})+\gamma+o(\gamma),
\end{align*}
and thus eventually $v\big(P(a)\big)\ \ge\ v(P^\phi_{+a,1})+\gamma+o(\gamma)$.
Take $g\in K^\times$ with $vg=\gamma$ and set $Q:= g^{-1}P_{+a,\times g}$, so $Q=Q_0+ Q_1 + Q_{>1}$ with
$$Q_0\ =\ Q(0)\ =\ g^{-1}P(a),\quad Q_1\ =\ g^{-1}(P_{+a,1})_{\times g},\quad Q_{>1}\ =\ g^{-1}(P_{+a,>1})_{\times g}. $$
Then $v(Q_0)=v\big(P(a)\big)-\gamma\ge v(P^\phi_{+a,1})+o(\gamma)$,
eventually. By [ADH, 6.1.3], $$v(Q_1^\phi)\ =\ v(P_{+a,1}^\phi)+o(\gamma),\qquad  v(Q_{>1}^\phi)\ \ge\ v(P_{+a,>1}^\phi)+\gamma+o(\gamma)$$
for all $\phi$. Since $P^\phi_{+a, >1}\preceq P^\phi_{+a,1}$, eventually, the last two displayed inequalities give~$v(Q^\phi_{>1})\ge v(Q^\phi_1)+\gamma+o(\gamma)$, eventually. We set $\hat f:= g^{-1}(\hat a -a)$, so $Q(\hat f)=0$ and~$\hat f\asymp 1$. Set $A:=L_Q\in K[\der]$. Then $Q(\hat f)=0$ gives
$$Q_0+A(\hat f)\ =\ Q_0+Q_1(\hat f)\ =\ -Q^\phi_{>1}(\hat f),$$
with $v\big(Q^\phi_{>1}(\hat f)\big)\ge v(Q^\phi_1)+\gamma+o(\gamma)$, eventually, so 
$$v\big(Q_0+A(\hat f)\big)\ \ge\ v(A^\phi)+\gamma+o(\gamma),\quad \text{eventually}.$$ 
Moreover, $v(Q_0)\ge v(A^\phi)+o(\gamma)$, eventually. Lemma~\ref{neneat3} then gives $f\in K$ with~$v(\hat f - f)\ge (1/3)\gamma$.
In view of $\hat a-a=g\hat f$, this yields
$$v\big(\hat a -(a+gf)\big)\ =\ \gamma+v(\hat f - f)\ \ge\ (4/3)\gamma,$$ 
which proves our claim.
\end{proof}

%\noindent
%The proof is the same as for Proposition~\ref{propsp}
%except for $\ndeg$ instead of $\ddeg$. Likewise we get the analogue of Proposition~\ref{propalsp}, replacing `Nigel' and [ADH, 6.9.3] in the proof by [ADH, 14.5.1, 11.4.8, 11.4.13]: 

\noindent
In the rest of this subsection we assume that $\hat a\notin K$ lies in an immediate asymptotic extension of $K$ and $K\<\hat a\>$ has transcendence degree $\le r$ over $K$.

\begin{prop}\label{npropalsp}\label{prop:hata almostspecial}  If $\Gamma$ is divisible, then 
$\hat a$ is almost special over $K$.  
\end{prop} 

\noindent
Towards the proof, we fix  a minimal annihilator $P(Y)$ of $\hat a$
 over $K$, so $\order P\le r$. We also fix a divergent pc-sequence $(a_{\rho})$ in $K$ such that $a_{\rho}\leadsto \hat a$. We next show how to improve $\hat a$ and $P$ if necessary:

\begin{lemma}\label{nreplace} For some $\hat{b}$ in an immediate asymptotic extension of $K$ we have: \begin{enumerate}
\item[\textup{(i)}] $v(\hat a -a)=v(\hat b -a)$ for all $a\in K$;
\item[\textup{(ii)}] $(a_{\rho})$ has a minimal differential polynomial $Q$ over $K$ of order $\le r$ such that~$Q$
is also a minimal annihilator of $\hat b$ over $K$.
\end{enumerate}
\end{lemma}
\begin{proof} By the remarks following the proof of [ADH, 11.4.3] we have $P\in Z(K,\hat a)$.  Take $Q\in Z(K,\hat a)$ of minimal complexity.
Then $\order Q \le \order P\le r$, and $Q$ is a minimal differential polynomial of $(a_{\rho})$ over $K$ by [ADH, 11.4.13]. 
By [ADH, 11.4.8 and its proof] this gives an element $\hat b$ in an immediate asymptotic extension of $K$ such that (i) holds and  $Q$ is a minimal annihilator of $\hat b$ over $K$. Then $Q$ and $\hat b$ have the desired properties. 
\end{proof}

\begin{proof}[Proof of Proposition~\ref{npropalsp}] 
Assume $\Gamma$ is divisible. Replace $\hat a$, $P$
by $\hat b$, $Q$ from Lemma~\ref{nreplace} and rename to arrange
that $P$ is a minimal differential polynomial of $(a_{\rho})$ over $K$.
By [ADH, 14.5.1] we have $a\in K$ and $g\in K^\times$ such that $\hat a -a\asymp g$ and $\ndeg P_{+a,\times g} = 1$.
Set $F:= P_{+a,\times g}$ and
$\hat f:= (\hat a  -a)/g$. Then $\ndeg F=1$, $F(\hat f)=0$, and $\hat f\preceq 1$. Applying Proposition~\ref{nepropsp} to
$F$ and $\hat f$ in the role of~$P$ and~$\hat a$ yields a
nontrivial convex subgroup $\Delta$ of $\Gamma$ that is
cofinal in $v(\hat f -K)$. Setting~$\alpha:= vg$, it follows that $\alpha+\Delta$ is cofinal in $v\big((\hat a -a)-K\big)=v(\hat a -K)$.   
\end{proof}

\begin{cor} 
If $K$ is henselian, then $\hat a$ is almost special over $K$.  
\end{cor}

\noindent
The proof is like that of Proposition~\ref{npropalsp}, using \cite[3.3]{Nigel19} instead of~[ADH, 14.5.1].

\subsection*{The case of order $1$} 
We show here that Proposition~\ref{prop:hata special}  goes through in the case of order $1$ under weaker assumptions:
 {\it in this subsection $K$ is a  $1$-linearly newtonian $H$-asymptotic field with asymptotic integration.}\/ Then $K$ is $\d$-valued with~$\I(K)\subseteq K^\dagger$, by Lemma~\ref{lem:ADH 14.2.5}, and $\upl$-free, by [ADH, 14.2.3].
We let $\phi$ range over elements active in~$K$. In the next two lemmas $A\in K[\der]^{\neq}$ has order~$\leq 1$,   $\gamma\in \Q\Gamma$, $\gamma>0$, and $\hat{K}$ is an immediate  asymptotic extension of $K$.

\begin{lemma}\label{neneat2, r=1} Let $\hat f\in \hat K^\times$ be such that  $v\big(A(\hat f)\big)\ge v(A^\phi\hat f)+\gamma$ eventually. Then there exists~$f\in K$ such that
$A(f)=0$ and~$v(\hat f-f)\ge v(\hat f) +\gamma$.
\end{lemma} 
\begin{proof}
Note that $\order(A)=1$; 
we arrange $A=\der-g$ ($g\in K$).
If~$A(\hat f)=0$, then~$\hat f$ is in $K$ [ADH, line before 14.2.10], and $f:=\hat f$ works. 
Assume~$A(\hat f)\neq 0$.
Then $$v\big(A^\phi(\hat f)\big)=v\big(A(\hat f)\big)\ge v(A^\phi\hat f)+\gamma>v(A^\phi\hat f),\quad\text{ eventually,}$$ 
so
$v(\hat f)\in\exc^{\ev}(A)$, and Lemma~\ref{lem:v(ker)=exc, r=1} yields an $f\in K$ with $f\sim \hat f$ and~$A(f)=0$.
We claim that this $f$ has the desired property. 
Set $b:= A(\hat{f})$. By the remarks preceding Corollary~\ref{cor:prlemexc} we can
replace~$K$,~$\hat K$,~$A$,~$b$ by~$K^\phi$,~$\hat K^\phi$,~$\phi^{-1}A^\phi$,~$\phi^{-1}b$, respectively, for suitable $\phi$,
to arrange that $K$ has small derivation and~$b^\dagger-g\succ^\flat 1$.  Using the hypothesis of the lemma we also arrange $vb\geq v(A\hat f)+\gamma$.
It remains to show that for~$\hat g:=\hat f-f\neq 0$ we have $v(\hat g)\ge v(\hat f)+\gamma$. 
Now $A(\hat g)=b$ with~$v(\hat g)\notin\exc^{\ev}(A)$, hence~$\hat g\sim b/(b^\dagger-g) \prec^\flat b$ by Lemma~\ref{prlemexc}, and thus
$v(\hat g)>vb\geq v(A\hat f)+\gamma$, so it is enough to show $v(A\hat f)\ge v(\hat f)$. 
Now $b=A(\hat f)=\hat f(\hat{f}^\dagger-g)$ and  $A\hat f=\hat f\big(\der+{\hat f}^\dagger -g\big)$. As $vb\ge v(A\hat f) + \gamma> v(A\hat f)$, this yields $v(\hat{f}^\dagger-g) > 0$, so $v(A\hat f)=v(\hat f)$.
% where~$v(g-\hat f^\dagger)>\Psi$; in particular $g-\hat f^\dagger \prec 1$,
%and so~$A\hat f\asymp \hat f$.
\end{proof}

\begin{lemma}\label{neneat3, r=1} Let $b\in K$ and $\hat f\in \hat K$ with $v(\hat f)\ge o(\gamma)$. Assume also that
eventually $v(b)\ge v(A^\phi)+o(\gamma)$ and $v\big(A(\hat f)-b\big)\ge v(A^\phi)+\gamma+o(\gamma)$. 
Then for some~$f\in K$ we have $A(f)=b$ and
$v(\hat f -f)\ge (1/2)\gamma+o(\gamma)$.
\end{lemma}

\noindent
The proof is like that of Lemma~\ref{neneat3}, using Lemma~\ref{neneat2, r=1} instead of Lemma~\ref{neneat2}.
In the same way Lemma~\ref{neneat3} gave Proposition~\ref{prop:hata special}, Lemma~\ref{neneat3, r=1} now yields:

\begin{prop}\label{nepropsp, r=1}\label{prop:hata special, r=1} If
$P\in K\{Y\}$, $\order P \le 1$, $\ndeg P=1$, 
and $P(\hat a)=0$, where~$\hat a\preceq 1$ lies in an immediate asymptotic extension of $K$ and $\hat a\notin K$, then~$\hat a$ is special over $K$. %some nontrivial convex subgroup of $\Gamma$ is cofinal in $v(\hat a - K)$.
\end{prop}

\begin{remark} 
Proposition~\ref{prop:hata almostspecial} does not hold for $r=1$ under present assumptions.
%without assuming $\upo$-freeness.
To see this, let $K$ be a Liouville closed $H$-field which is not $\upo$-free, as in Example~\ref{ex:Gehret} or \cite{ADH3}.
Then $K$ is $1$-linearly newtonian by Corollary~\ref{cor:Liouville closed => 1-lin newt} below.
Consider the pc-sequences $(\upl_\rho)$ and $(\upo_\rho)$  in $K$ as in [ADH, 11.7],
 let $\upo\in K$ with~$\upo_\rho\leadsto\upo$, and~$P=2Y'+Y^2+\upo$. 
Then [ADH, 11.7.13] gives an element $\upl$ in an immediate asymptotic extension of $K$ but not in $K$
with $\upl_\rho\leadsto\upl$ and $P(\upl)=0$. However, $\upl$ is not almost special over $K$ [ADH, 3.4.13, 11.5.2].
\end{remark}

\subsection*{Relating $Z(K,\hat a)$ and $v(\hat a-K)$ for special $\hat a$} 
{\em In this subsection $K$ is a valued differential field with small derivation $\der\ne 0$
such that $\Gamma\ne \{0\}$ and $\Gamma^{>}$ has no least element.}\/ We recall from \cite{VDF} that a valued differential field extension $L$ of $K$ is said to be {\em strict\/}\index{extension!strict} if for all $\phi\in K^\times$,
$$\der \smallo\subseteq \phi \smallo\ \Rightarrow\ \der\smallo_L\subseteq \phi \smallo_L, \qquad \der \O\subseteq \phi \smallo\ \Rightarrow\ \der\O_L\subseteq \phi \smallo_L.$$
(If $K$ is asymptotic, then any immediate asymptotic extension of~$K$ is automatically strict, by \cite[1.11]{VDF}.)
Let $\hat a$ lie in an immediate strict extension of $K$ such that~${\hat a\preceq 1}$, $\hat a\notin K$, and~$\hat a$ is special over $K$.
 We adopt from \cite[Sections~2,~4]{VDF} the definitions of~$\ndeg P$ for $P\in K\{Y\}^{\ne}$ and of the set $Z(K,\hat a)\subseteq K\{Y\}^{\ne}$. 
Also recall that~$\Gamma(\der):=\{v\phi:\, \phi\in K^\times,\, \der\smallo\subseteq\phi\smallo\}$.

\begin{lemma}\label{Zp1} Let $P\in Z(K,\hat a)$ and $P\asymp 1$. Then
$v\big(P(\hat a)\big) >  v(\hat a-K)$. 
\end{lemma}
\begin{proof} 
Take a divergent pc-sequence $(a_{\rho})$ in $\mathcal O$ with $a_{\rho} \leadsto \hat a$, and  
as in [ADH, 11.2] let $\boldsymbol  a:=c_K(a_\rho)$. Then 
$\ndeg_{\boldsymbol a} P\geq 1$ by~\cite[4.7]{VDF}. We arrange $\gamma_{\rho}:= v(\hat a-a_{\rho})$ to be strictly increasing as a function of~$\rho$, with $0 < 2\gamma_{\rho} < \gamma_{s(\rho)}$  for all $\rho$. 
Take~$g_{\rho}\in \smallo$ with $g_{\rho} \asymp \hat a - a_{\rho}$; then $1 \leq d:= \ndeg_{\boldsymbol a} P = \ndeg P_{+a_{\rho}, \times g_{\rho}}$ for all sufficiently large~$\rho$, and we arrange that this holds for all $\rho$. 
We have~$\hat a = a_{\rho} + g_{\rho}y_{\rho}$ with $y_{\rho}\asymp 1$, and 
$$P(\hat a)\ =\ P_{+a_{\rho},\times g_{\rho}}(y_{\rho})\ =\ \sum_i (P_{+a_{\rho},\times g_{\rho}})_i(y_{\rho}).$$
Pick for every $\rho$ an element $\phi_{\rho}\in K^\times$ such that $0\le v(\phi_{\rho})\in \Gamma(\der)$ and
$(P^{\phi_{\rho}}_{+a_{\rho},\times g_{\rho}})_i\ \preceq\ (P^{\phi_{\rho}}_{+a_{\rho},\times g_{\rho}})_{d}$
for all $i$. Then for all $\rho$ and $i$, 
\begin{align*} (P_{+a_{\rho},\times g_{\rho}})_i(y_{\rho})\ & = (P^{\phi_{\rho}}_{+a_{\rho},\times g_{\rho}})_i(y_{\rho})\ \preceq\  (P^{\phi_{\rho}}_{+a_{\rho},\times g_{\rho}})_i\ \preceq\ (P^{\phi_{\rho}}_{+a_{\rho},\times g_{\rho}})_{d}\ 
\text{ with }\\
 v\big((P^{\phi_{\rho}}_{+a_{\rho},\times g_{\rho}})_{d}\big)\ &\ge\ d \gamma_{\rho} + o(\gamma_{\rho})\ \ge\ \gamma_{\rho} + o(\gamma_{\rho}),
 \end{align*}
where for the next to last inequality we use [ADH, 11.1.1, 5.7.1, 5.7.5, 6.1.3].  
Hence~$v\big(P(\hat a)\big) \ge \gamma_{\rho} + o(\gamma_{\rho})$ for all $\rho$, and thus $v\big(P(\hat a)\big)> v(\hat a -K)$.
\end{proof}

\noindent
We also have a converse under extra assumptions:

\begin{lemma} \label{lem:ZKhata} Assume $K$ is asymptotic and $\Psi\subseteq v(\hat a-K)$. Let $P\in K\{Y\}$ be such that $P\asymp 1$ and $v\big(P(\hat a)\big)> v(\hat a -K)$. Then $P\in Z(K,\hat a)$.
\end{lemma}
\begin{proof} 
Let $\Delta$ be the nontrivial convex subgroup of $\Gamma$ that is cofinal in $v(\hat a -K)$. Let $\kappa:=\cf(\Delta)$.
Take a divergent pc-sequence $(a_{\rho})_{\rho< \kappa}$ in $K$ such that $a_{\rho}\leadsto \hat a$.   
We arrange $\gamma_{\rho}:= v(\hat a -a_{\rho})$ is strictly increasing as a function of $\rho$, with $\gamma_{\rho}>0$ for all~$\rho$; thus
$a_{\rho}\preceq 1$ for all $\rho$. 
Consider the $\Delta$-coarsening $\dot v=v_{\Delta}$ of the valuation $v$ of~$K$; it has valuation ring $\dot{\O}$ with
differential residue field $\dot K$. Consider likewise the $\Delta$-coarsening of the valuation of
the immediate extension $L=K\<\hat a\>$ of $K$.  Let~$a^*$ be the image of $\hat a$ in the differential residue field
$\dot{L}$ of $(L,\dot v)$. Note that $\dot{L}$ is an immediate extension of $\dot{K}$. The pc-sequence $(a_{\rho})$ then
yields a sequence $(\dot{a_{\rho}})$  in~$\dot{K}$ with
$v(a^*-\dot{a_{\rho}})=\gamma_{\rho}$ for all $\rho$. Thus $(\dot a_{\rho})$ is a c-sequence in $\dot{K}$ with $\dot a_{\rho} \to a^*$, so
$\dot{P}(\dot a_{\rho})\to \dot{P}(a^*)$ by [ADH, 4.4.5]. From $v\big(P(\hat a)\big)>\Delta$ we obtain~$\dot{P}(a^*)=0$, and so $\dot{P}(\dot a_{\rho})\to 0$.
So far we have not used our assumption that $K$ is asymptotic and~$\Psi\subseteq v(\hat a-K)$. Using this now, we note that
for $\alpha\in \Delta^{>}$ we have 
$0 <\alpha'=\alpha+\alpha^\dagger$, so $\alpha' \in \Delta$, hence
the derivation of $\dot{K}$ is nontrivial. Thus we can apply [ADH, 4.4.10] to $\dot{K}$ and modify the $a_{\rho}$ 
without changing $\gamma_{\rho}=v(a^*-\dot a_{\rho})$ to arrange that in addition~$\dot{P}(\dot a_{\rho})\ne 0$
for all $\rho$.  Since $\kappa=\cf(\Delta)$, we can replace $(a_{\rho})$ by a cofinal subsequence so that $P(a_{\rho})\leadsto 0$, hence $P\in Z(K,\hat a)$ by \cite[4.6]{VDF}.
\end{proof}

\noindent
To elaborate on this, let $\Delta$ be a convex subgroup of $\Gamma$ and $\dot K$ the
valued differential residue field of the  $\Delta$-coarsening $v_\Delta$ of the valuation~$v$ of $K$. We view $\dot K$ in the usual way as a valued differential subfield of the valued differential residue field $\dot{\hat K}$ of the $\Delta$-coarsening of the valuation of $\hat K$ by $\Delta$; see [ADH, pp.~159--160 and~4.4.4]. %\marginpar{moved here from a later section}

\begin{cor}\label{cor:ZKhata} 
Suppose $K$ is asymptotic,  $\Psi\subseteq v(\hat a -K)$, and  $\Delta$  is cofinal in~$v({\hat a-K})$. 
 Let $P\in K\{Y\}$ with $P \asymp 1$. Then
$P\in Z(K,\hat a)$ if and only if~$\dot P(\dot{\hat a})=0$ in $\dot{\hat K}$. Also,
 $P$ is an element of $Z(K,\hat a)$ of minimal complexity if and only if $\dot P$ is a minimal
annihilator of $\dot{\hat a}$ over $\dot K$ and $\dot P$ has the same complexity as~$P$.
\end{cor}
\begin{proof}
The first statement is immediate from  Lemmas~\ref{Zp1} and~\ref{lem:ZKhata}.
For the rest use that 
for $R\in\dot{\mathcal O}\{Y\}$ we have $\cc(\dot R)\leq\cc(R)$ and
that for all $Q\in \dot K\{Y\}$ there is an $R\in\dot{\mathcal O}\{Y\}$
with  $Q=\dot R$
and $Q_\i\neq 0$ iff $R_\i\neq 0$ for all $\i$, so
$\cc(\dot R)=\cc(R)$.
\end{proof}

\section{Differential Henselianity of the Completion}\label{sec:completion d-hens}

\noindent
{\em Let $K$ be a valued differential field with small derivation.}\/ We let $\Gamma:= v(K^\times)$ be the value group of $K$ and 
$\k:=\res(K)$ be the differential residue field of $K$, and we let $r\in\N$.  
The following summarizes [ADH, 7.1.1, 7.2.1]:

\begin{lemma}
The valued differential field $K$ is $r$-$\d$-henselian iff  for each 
$P$ in~$K\{Y\}$ of order~$\leq r$ with $\ddeg P=1$  there is a $y\in\mathcal O$ with $P(y)=0$.
\end{lemma} 

\noindent
Recall that the derivation of $K$ being small, it is continuous [ADH, 4.4.6], and hence extends uniquely to a continuous derivation
on  the completion $K^{\cc}$ of the valued field $K$ [ADH, 4.4.11]. We equip  $K^{\cc}$ with this derivation, which remains small~[ADH, 4.4.12],  so $K^{\cc}$ is an immediate valued differential field extension of~$K$ with small derivation, in particular, $\k=\res(K^{\cc})$. \index{valued differential field!completion} \index{completion}

Below we  characterize in a first-order way when 
$K^{\cc}$ is $r$-$\d$-henselian. We shall use tacitly that for $P\in K\{Y\}$ 
we have $P(g)\preceq P_{\times g}$ for all $g\in K$;
to see this, replace~$P$ by $P_{\times g}$ to reduce to $g=1$, and observe that
$P(1)=\sum_{\dabs{\bsigma}=0} P_{[\bsigma]}\preceq P$.

\begin{lemma}\label{lem:Kc 1}
Let  $P\in K^{\cc}\{Y\}$, $b\in K^{\cc}$ with $b\preceq 1$ and $P(b)=0$, and 
$\gamma\in\Gamma^>$. Then there
is an $a\in\mathcal O$ with $v\big(P(a)\big)>\gamma$.
\end{lemma}
\begin{proof}
To find  an $a$ as claimed we take $f\in K$ satisfying $f\asymp P$ and replace~$P$,~$\gamma$ by $f^{-1}P$, $\gamma-vf$, respectively, to arrange $P\asymp 1$ and thus  $P_{+b}\asymp 1$. We also assume~$b\ne 0$.  
Since $K$ is dense in~$K^{\cc}$ we can take $a\in K$ such that $a\sim b$ (so~$a\in\mathcal O$) and
$v(a-b)>2\gamma$. Then  with $g:=a-b$,   using [ADH, 4.5.1(i) and 6.1.4] we conclude
$$v\big(P(a)\big)=v\big(P_{+b}(g)\big) \geq v\big((P_{+b})_{\times g}\big)
\geq v(P_{+b}) + vg + o(vg) = vg+o(vg) >\gamma$$
as required. 
\end{proof}

\noindent
Recall that if $K$ is asymptotic, then so is $K^{\cc}$ by [ADH, 9.1.6]. 

\begin{lemma}\label{lem:Kc 2}
Suppose $K$ is asymptotic, $\Gamma\ne \{0\}$, and for every $P\in K\{Y\}$ of order at most~$r$ with $\ddeg P=1$ and every $\gamma\in\Gamma^>$ there
is an $a\in\mathcal O$ with~${v\big(P(a)\big)>\gamma}$.
Then~$K^{\cc}$ is $r$-$\d$-henselian.
\end{lemma}
\begin{proof}
The hypothesis applied to $P\in\mathcal O\{Y\}$ of order~$\leq r$ with~${\ddeg P=\deg P=1}$
yields that $\k$ is $r$-linearly surjective.
Let now $P\in K^{\cc}\{Y\}$ be of order~$\leq r$ with~$\ddeg P=1$.   We need to show that there exists $b\in K^{\cc}$ with $b\preceq 1$ and~$P(b)=0$. First we arrange $P\asymp 1$.  
By [ADH, remarks after 9.4.11] we can take $b\preceq 1$ in an  immediate $\d$-henselian asymptotic field extension~$L$ 
of~$K^{\cc}$ with $P(b)=0$.  We prove below that~$b\in K^{\cc}$. We may assume $b\notin K$, so
$v(b-K)$ has no largest element, since~$L\supseteq K$ is immediate.
Note also that   $\ddeg P_{+b}=1$ by [ADH, 6.6.5(i)]; since~${P(b)=0}$ we thus have $\ddeg P_{+b,\times g}=1$
for all $g\preceq 1$ in $L^\times$ by [ADH, 6.6.7].

\claim{Let $\gamma\in\Gamma^>$ and $a\in K$ with $v(b-a)\geq 0$.  There is a $y\in\mathcal O$ such that~$v\big(P(y)\big)>\gamma$ and $v(b-y)\geq v(b-a)$.}

\noindent
To prove this claim, take $g\in K^\times$ with $vg=v(b-a)$. Then by [ADH, 6.6.6] and the observation 
preceding the claim we have~${\ddeg P_{+a,\times g}=\ddeg P_{+b,\times g}=1}$. 
Thanks to density of $K$ in $K^{\cc}$ we may take $Q\in K\{Y\}$ of order~$\leq r$ with $P_{+a,\times g}\sim Q$ and $v(P_{+a,\times g}-Q)>\gamma$. Then 
$\ddeg Q=1$, so by the hypothesis of the lemma  we have~${z\in\mathcal O}$ with $v\big(Q(z)\big)>\gamma$.
Set $y:=a+gz\in\mathcal O$; then~$v\big(P(y)\big)=v\big(P_{+a,\times g}(z)\big)>\gamma$ and $v(b-y)=v(b-a-gz)\geq v(b-a)= vg$ as claimed.

\medskip
\noindent
Let now $\gamma\in\Gamma^>$;  to show that $b\in K^{\cc}$,
it is enough by [ADH, 3.2.15, 3.2.16] to show that then $v(a-b)>\gamma$ for some $a\in K$.
Let   $A:=L_{P_{+b}}\in L[\der]$; then $A\asymp 1$. Since $\abs{\exc_L(A)}\leq r$   by [ADH, 7.5.3], the claim gives an $a\in\mathcal O$ with $v\big(P(a)\big)>2\gamma$ and $0<v(b-a)\notin\exc_L(A)$. Put $g:=a-b$ and $R:=(P_{+b})_{>1}$. Then
$R\prec 1$ and
$$P(a)\ =\ P_{+b}(g)\ =\ A(g) + R(g)$$
where by the definition of $\exc_L(A)$ and [ADH, 6.4.1(iii), 6.4.3] we have in $\Q\Gamma$:
$$v\big(A(g)\big)\ =\ v_A(vg)\ =\ vg+o(vg)\ <\ vR + (3/2)vg\   \leq\ v_R(vg)\ \leq\ v\big(R(g)\big)$$
and so $v\big(P(a)\big) = vg+o(vg) > 2\gamma$. Therefore $v(a-b)=vg>\gamma$ as required.
\end{proof}

\noindent
The last two lemmas yield an analogue of  [ADH, 3.3.7] for $r$-$\d$-hensel\-ianity and a partial generalization of [ADH, 7.2.15]:

\begin{cor}\label{cor:Kc d-henselian}
Suppose $K$ is asymptotic and $\Gamma\ne \{0\}$.
Then the following are equivalent:
\begin{enumerate}
\item[$\mathrm{(i)}$] $K^{\cc}$ is $r$-$\d$-henselian;
\item[$\mathrm{(ii)}$] for every $P\in K\{Y\}$ of order at most $r$ with $\ddeg P=1$ and every $\gamma\in\Gamma^>$ there
exists $a\in\mathcal O$ with $v\big(P(a)\big)>\gamma$.
\end{enumerate}
In particular, if $K$ is $r$-$\d$-henselian, then so is $K^{\cc}$.
\end{cor}

%\subsection*{Results from {\tt asymp}} 
%Let $K$ be $H$-asymptotic with rational asymptotic integration, and   $\hat a$ an element in an immediate asymptotic extension $\hat K$ of $K$ with $\hat a\notin K$. Let $r\ge 1$ and assume $P\in K\{Y\}^{\neq}$ has order $\le r$. Then by  [{\tt asymp}, 0.0.37] and  [{\tt asymp}, 0.0.38], respectively:

%\begin{prop}\label{prop:hata special}
%If $K$ is $\upo$-free and $r$-linearly newtonian, $\ndeg P=1$, $\hat a\preceq 1$,  and $P(\hat a)=0$, then $\hat a$ is   special over $K$. 
%\end{prop}

%\begin{prop}\label{prop:hata almostspecial}
%If $\Gamma$ is divisible, $K$ is $\upo$-free and $r$-linearly newtonian, and $P(\hat a)=0$, then $\hat a$ is  almost special over $K$. 
%\end{prop}

%\noindent
%{\em In the rest of this subsection we also assume that $K$ has small derivation}. The following  is  shown in [{\tt asymp}, 0.0.40 and 0.0.41] under more general assumptions.

%\begin{lemma}\label{lem:ZKhata}
%Suppose  $\hat a \preceq 1$ is special over $K$, $P\asymp 1$, and $\Psi\subseteq v(\hat a -K)$. Then:
%$$P\in Z(K,\hat a)\ \Longleftrightarrow\ v\big(P(\hat a)\big)>v(\hat a-K).$$ 
%\end{lemma}

\section{Complements on Newtonianity}\label{sec:complements newton}

\noindent 
{\em In this section $K$ is an ungrounded $H$-asymptotic field with $\Gamma=v(K^\times)\neq\{0\}$}. Note that then $K^{\cc}$ is also $H$-asymptotic. We let $r$ range over $\N$ and $\phi$ over the active elements of $K$.
Our first aim is a newtonian analogue of Corollary~\ref{cor:Kc d-henselian}:
%, whose ``in particular'' statement partly generalizes 

\begin{prop}\label{prop:Kc newtonian}
For $\d$-valued and $\upo$-free $K$, the  following are equivalent:
\begin{enumerate}
\item[$\mathrm{(i)}$] $K^{\cc}$ is $r$-newtonian;
\item[$\mathrm{(ii)}$]  for every $P\in K\{Y\}$ of order at most $r$ with $\ndeg P=1$ and every $\gamma\in\Gamma^>$ there
is an $a\in\mathcal O$ with $v\big(P(a)\big)>\gamma$.
\end{enumerate}
If $K$ is $\d$-valued, $\upo$-free, and $r$-newtonian, then so is $K^{\cc}$.
\end{prop} 

\noindent
The final statement in this proposition extends [ADH, 14.1.5]. Towards the proof we first
state a variant of [ADH, 13.2.2] which follows easily from [ADH, 11.1.4]: 
%isolate some implications used in the proofs of [ADH, 13.1.9, 13.2.2], which follows easily from [ADH, 11.1.4]: \marginpar{changed this sentence a bit}

\begin{lemma}\label{lem:same ndeg} \label{lem:ndeg of nearby diffpoly} 
Assume $K$ has small derivation and let $P,Q\in K\{Y\}^{\neq}$ and 
$\phi\preceq 1$. Then  $P^\phi \asymp^\flat P$, and so we have the three implications
$$P \preceq^\flat Q\ \Longrightarrow\ P^\phi\preceq^\flat Q^\phi,\quad 
P\prec^\flat Q\ \Longrightarrow\ 
P^\phi \prec^\flat Q^\phi,\quad P \sim^\flat Q\ \Longrightarrow\ P^\phi\sim^{\flat} Q^\phi.$$ 
%and similarly for $\prec^\flat$ and $\sim^\flat$ in place of $\preceq^\flat$. 
The last implication gives:  $P\sim^\flat Q\ \Longrightarrow\ \ndeg P=\ndeg  Q\text{ and }\nval P=\nval Q.$
\end{lemma}

\noindent
For $P^\phi\asymp^\flat P$ and the subsequent three implications 
in the lemma above we can drop the assumption that $K$ is ungrounded.

\begin{lemma}\label{lem:Kc 3}
Suppose $K$ is $\d$-valued, $\upo$-free, and for every $P\in K\{Y\}$ of order at most $r$ with $\ndeg P=1$ and every $\gamma\in\Gamma^>$ there
is an $a\in\mathcal O$ with $v\big(P(a)\big)>\gamma$. Then $K^{\cc}$ is $\d$-valued, $\upo$-free, and $r$-newtonian.
\end{lemma}

\begin{proof} By [ADH, 9.1.6 and 11.7.20], $K^{\cc}$ is $\d$-valued and $\upo$-free. 
Let $P\in K^{\cc}\{Y\}$ be of order~$\leq r$ with $\ndeg P=1$.   We need to show that $P(b)=0$ for some $b\preceq 1$ in $K^{\cc}$. 
To find $b$  we may replace $K$,~$P$ by $K^\phi$,~$P^\phi$; in particular we may assume that
$K$ has small derivation and $\Gamma^{\flat}\ne \Gamma$.
By [ADH, 14.0.1 and the remarks following it] we can take
$b\preceq 1$ in 
an immediate newtonian extension $L$ of~$K^{\cc}$  such that $P(b)=0$. We claim that~$b\in K^{\cc}$.
To show this we may assume $b\notin K$, so  $v(b-K)$ does not have a largest element.
%Replacing $K$, $P$  by $K^\phi$, $P^\phi$ for suitable $\phi$ we can assume  $P=L+R$ where $L\in K^{\cc}\{Y\}$ has degree~$1$ and $R\prec^\flat L\asymp P$. We may also assume that $P\asymp 1$. Denoting the dominant degree of 
By [ADH, 11.2.3(i)] we have $\ndeg P_{+b}=1$ and so $\ndeg P_{+b,\times g}=1$
for all $g\preceq 1$ in~$L^\times$ by [ADH, 11.2.5], in view of $P(b)=0$.

\claim{Let $\gamma\in\Gamma^>$ and $a\in K$ with $v(b-a)\geq 0$. There is a $y\in\mathcal O$ such that~$v\big(P(y)\big)>\gamma$ and $v(b-y)\geq v(b-a)$.}

\noindent
The proof is similar to that of the claim in the proof of Lemma~\ref{lem:Kc 2}:
Take $g\in K^\times$ with $vg=v(b-a)$.  Then $\ndeg P_{+a,\times g}=\ndeg P_{+b,\times g}=1$ by [ADH, 11.2.4] and the observation 
preceding the claim. 
Density of $K$ in $K^{\cc}$ yields $Q\in K\{Y\}$ of order~$\leq r$ with $v(P_{+a,\times g}-Q)>\gamma$ and $P_{+a,\times g}\sim^\flat Q$, the latter using $\Gamma^{\flat}\ne \Gamma$. Then 
$\ndeg Q=\ndeg P_{+a,\times g}=1$ by Lemma~\ref{lem:ndeg of nearby diffpoly}, so the hypothesis of the lemma gives $z\in\mathcal O$ with~$v\big(Q(z)\big)>\gamma$.
Setting $y:=a+gz\in\mathcal O$ we have $v\big(P(y)\big)=v\big(P_{+a,\times g}(z)\big)>\gamma$ and $v(b-y)=v(b-a-gz)\geq vg=v(b-a)$.

\medskip
\noindent
Let   $\gamma\in\Gamma^>$;  to get $b\in K^{\cc}$,
it is enough to show that then $v(a-b)>\gamma$
for some~$a\in K$.
Let   $A:=L_{P_{+b}}\in L[\der]$. Since $\abs{\exc^{\ev}_L(A)}\leq r$   by [ADH, 14.2.9], by the claim we can take
$a\in\mathcal O$ with $v\big(P(a)\big)>2\gamma$ and $0<v(b-a)\notin\exc^{\ev}_L(A)$. Now put $g:=a-b$ and
take $\phi$ with $vg\notin\exc_{L^\phi}(A^\phi)$; note that then $A^\phi=L_{P^\phi_{+b}}$.
Replacing~$K$,~$L$,~$P$ by $K^\phi$,~$L^\phi$,~$P^\phi$ we arrange $vg\notin \exc_L(A)$, and (changing $\phi$ if necessary) $\ddeg P_{+b}=1$. We also arrange
$P_{+b}\asymp 1$, and then  $(P_{+b})_{>1}\prec 1$. 
As in the proof of Lemma~\ref{lem:Kc 2} above we now derive   $v(a-b)=vg>\gamma$.
\end{proof}

\noindent
Combining Lemmas~\ref{lem:Kc 1} and \ref{lem:Kc 3} now yields  Proposition~\ref{prop:Kc newtonian}.  \qed

\medskip
\noindent
To show that newtonianity is preserved under specialization, we assume below that~$\Psi\cap \Gamma^{>}\ne \emptyset$, so $K$ has
small derivation.  Let $\Delta\neq\{0\}$ 
be a convex subgroup of~$\Gamma$ with $\psi(\Delta^{\neq})\subseteq\Delta$. Then $1\in \Delta$ where $1$ denotes the unique positive element of $\Gamma$ fixed by the function 
$\psi$: use that $\psi(\gamma)\ge 1$ for $0 < \gamma < 1$.
(Conversely, any convex subgroup $G$ of $\Gamma$ with $1\in G$ satisfies $\psi(G^{\ne})\subseteq G$.)  Let $\dot v$ be
the coarsening of the valuation $v$ of $K$ by $\Delta$, with valuation ring $\dot{\mathcal O}$,
maximal ideal $\dot\smallo$ of $\dot{\mathcal O}$, and residue field  $\dot K=\dot{\mathcal O}/\dot\smallo$.
The derivation of $K$ is small with respect to $\dot{v}$,
and $\dot K$ with the induced valuation $v\colon \dot K^\times\to\Delta$ and induced derivation as in [ADH, p.~405] is an asymptotic field with asymptotic couple $(\Delta,\psi|\Delta^{\neq})$, and so
is of $H$-type with small derivation. If $K$ is $\d$-valued, then so is
$\dot K$ by [ADH, 10.1.8], and if $K$ is $\upo$-free, then so is $\dot K$ by [ADH, 11.7.24].
The residue map
$a\mapsto\dot a:=a+\dot\smallo\colon \dot{\mathcal O}\to\dot K$ is a differential ring morphism, 
extends to a differential ring morphism $P\mapsto\dot P\colon \dot{\mathcal O}\{Y\}\to\dot K\{Y\}$ of differential rings
sending $Y$ to $Y$, and $\ddeg P=\ddeg\dot P$ for $P\in \dot{\mathcal O}\{Y\}$ with $\dot P \ne 0$. 

We now restrict $\phi$ to range over active elements of $\mathcal O$. Then $v\phi\le 1+1$, so~$v\phi\in \Delta$, and hence $\phi$ is a unit of $\dot{\mathcal O}$. 
It follows that $\dot\phi$ is active in $\dot K$, and every active element of $\dot K$ lying in its valuation ring is of this form. Moreover, the differential subrings~$\dot{\mathcal O}$  of $K$ and $\dot{\mathcal O}^\phi:=(\dot{\mathcal O})^\phi$ of $K^\phi$ have the same underlying ring, and the derivation of~$K^\phi$ is small with respect to $\dot{v}$. Thus the differential residue fields $\dot{K}=\dot{\mathcal O}/\dot{\smallo}$ and~$\dot{K}^\phi:=\dot{\mathcal O}^\phi/\dot{\smallo}$ have the same underlying field and are related as follows:
$$\dot{K}^\phi\ =\ (\dot{K})^{\dot\phi}.$$
For $P\in \dot{\mathcal O}\{Y\}$ we have $P^\phi\in \dot{\mathcal O}^\phi\{Y\}$, and the image of
$P^\phi$ under the residue map~$\dot{\mathcal O}^\phi\{Y\}\to\dot K^\phi\{Y\}$ equals
$\dot P^{\dot\phi}$; hence $\ndeg P=\ndeg \dot P$ for $P\in \dot{\mathcal O}\{Y\}$ satisfying~${\dot{P}\ne 0}$. 
These remarks imply:

\begin{lemma}\label{lem:dotK newtonian}
If $K$ is $r$-newtonian, then $\dot K$ is $r$-newtonian.
\end{lemma}

\noindent
Combining Proposition~\ref{prop:Kc newtonian} and Lemmas~\ref{lem:Kc 3} and~\ref{lem:dotK newtonian}  yields:

\begin{cor}\label{cor:Kc newtonian}
Suppose $K$ is $\d$-valued, $\upo$-free, and $r$-newtonian. Then $\dot K$ and its completion are
$\d$-valued, $\upo$-free, and $r$-newtonian.
\end{cor} 

\noindent
We finish with a newtonian analogue of [ADH, 7.1.7]:

\begin{lemma}\label{rdhrnewt}
Suppose $(K,\dot{\mathcal O})$ is $r$-$\d$-henselian and $\dot K$ is $r$-newtonian. Then $K$ is $r$-newtonian.
\end{lemma}
\begin{proof}
Let $P\in K\{Y\}$ be of order~$\leq r$ and $\ndeg P=1$; we need to show the existence of $b\in\mathcal O$ with $P(b)=0$. Replacing $K$ and $P$ by $K^\phi$ and $P^\phi$ for suitable $\phi$ (and renaming) we arrange $\ddeg P =1$; this also uses [ADH, section 7.3, subsection on compositional conjugation]. 
We can further assume that $P\asymp 1$. Put $Q:=\dot P\in\dot K\{Y\}$, so %The coefficients of $Q$ are in the valuation ring of $\dot K$, and 
$\ndeg Q=1$, and thus $r$-newtonianity of $\dot K$ yields an $a\in\mathcal O$ with $Q(\dot a)=0$. Then $P(a)\dotprec 1$, $P_{+a,1}\sim P_1\asymp 1$,
and $P_{+a,>1}\prec 1$. Since $(K,\dot{\mathcal O})$ is $r$-$\d$-henselian, this gives $y\in\dot\smallo$ with $P_{+a}(y)=0$,
and then $P(b)=0$ for $b:=a+y\in\mathcal O$.
\end{proof}

\noindent
Lemmas~\ref{lem:dotK newtonian},~\ref{rdhrnewt}, and  [ADH, 14.1.2] yield:

\begin{cor}
$K$ is $r$-newtonian iff $(K,\dot{\mathcal O})$ is $r$-$\d$-henselian and $\dot K$ is $r$-new\-to\-ni\-an.
\end{cor}

\subsection*{Invariance of Newton quantities} 
{\it In this subsection $P\in K\{Y\}^{\neq}$.}\/
In~[ADH, 11.1] we associated to $P$ its Newton weight $\nwt P$, Newton degree $\ndeg P$, and Newton multiplicity $\nval P$ at~$0$, all elements of $\N$,
as well as the element $v^{\ev}(P)$ of $\Gamma$;
these quantities do not change when passing to an $H$-asymptotic
extension~$L$ of~$K$ with $\Psi$ cofinal in $\Psi_L$, cf.~[ADH, p.~480], where the assumptions on $K$, $L$ are slightly weaker. 
Thus by Theorem~\ref{thm:ADH 13.6.1}, these quantities do not change for $\upo$-free~$K$ in passing to an  $H$-asymptotic pre-$\d$-valued $\d$-algebraic extension of~$K$.
Below we improve on this in several ways. First,  for $\order P\leqslant1$ we can drop $\Psi$ being cofinal in $\Psi_L$ by a strengthening of [ADH, 11.2.13]:

\begin{lemma}\label{lem:11.2.13 invariant}
Suppose $K$ is $H$-asymptotic with rational asymptotic integration and   $P\in K[Y,Y']^{\ne}$.
Then there are $w\in \N$, $\alpha\in \Gamma^{>}$, 
$A\in K[Y]^{\ne}$, and an active~$\phi_0$ in~$K$ such that for every  asymptotic extension $L$ of $K$ and active~$f\preceq\phi_0$ in $L$,
$$P^{f}\ =\  f^w A(Y)(Y')^w + R_{f},\quad R_{f}\in L^{f}[Y, Y'],\quad v(R_{f})\ \geqslant\ v(P^{f}) + \alpha.$$
For such $w, A$ we have for any ungrounded $H$-asymptotic extension $L$ of $K$,
$$ \nwt_L P\ =\ w, \quad \ndeg_L P\ =\ \deg A+w,\quad
 \nval_L P\ =\ \val A+w,\quad\  v^{\ev}_L(P) = v(A).$$
%\textup{(}In particular, $\nwt_L P=w$, $\ndeg_L P=\deg A+w$, and $\nval_L P=\val A+w$.\textup{)}
\end{lemma}
\begin{proof} Let $w$ be as in the proof of [ADH, 11.2.13]. Using its notations, this proof yields an active 
$\phi_0$ in $K$ such that 
\begin{equation}\label{eq:11.2.13 invariant, 1} w\gamma+v(A_w)\  <\  j\gamma + v(A_j)\end{equation}
for all 
$\gamma\geqslant v(\phi_0)$ in $\Psi^{\downarrow}$ and  $j\in J\setminus \{w\}$. This gives
$\beta\in \Q\Gamma$  such that $\beta>\Psi$ and~\eqref{eq:11.2.13 invariant, 1} remains true for all  $\gamma\in \Gamma$ with $v(\phi_0)\leqslant \gamma < \beta$. Since $(\Q\Gamma, \psi)$ has asymptotic integration, $\beta$ is not a gap in  $(\Q\Gamma,\psi)$, so
$\beta>\beta_0> \Psi$ with $\beta_0\in \Q\Gamma$. This yields an element $\alpha\in (\Q\Gamma)^{>}$ such that
for all $\gamma\in \Q\Gamma$ with $v(\phi_0) \le \gamma \le \beta_0$ we have
\begin{equation}\label{eq:11.2.13 invariant, 2} w\gamma + v(A_w) +\alpha\ \leqslant\ j\gamma + v(A_j)\end{equation}
Since $\Gamma$ has no least positive element, we can decrease $\alpha$ to arrange $\alpha\in \Gamma^{>}$. 
Now~\eqref{eq:11.2.13 invariant, 2} remains true for all elements $\gamma$ 
of any divisible ordered abelian group extending $\Q\Gamma$ 
with $v(\phi_0)\leqslant \gamma \le \beta_0$.
Thus $w$, $\alpha$, $A=A_w$, and $\phi_0$ are as required. 

For any ungrounded $H$-asymptotic extension $L$ of $K$ we obtain for active $f\preceq \phi_0$ in $L$ that $v(P^f)=v(A)+wv(f)$,  so $v_L^{\ev}(P)=v(A)$ in view of  the  identity in  [ADH, 11.1.15]
defining $v_L^{\ev}(P)$. 
\end{proof} 

\noindent
For quasilinear~$P$ we have:

\begin{lemma}\label{lem:13.7.10}
Suppose $K$ is $\upl$-free and $\ndeg P=1$. Then there are active~$\phi_0$ in~$K$ and~$a,b\in K$ with $a\preceq b\ne 0$ 
 such that  either \textup{(i)} or \textup{(ii)} below holds:
\begin{enumerate}
\item[\textup{(i)}]
$P^{f}\, \sim_{\phi_0}^\flat\, a+bY$
for all active~$f\preceq\phi_0$ in all  $H$-asymptotic extensions of $K$;
\item[\textup{(ii)}]
$P^{f}\, \sim_{\phi_0}^\flat\, \frac{f}{\phi_0}b\,Y'$ for all active~$f\preceq\phi_0$ in all  $H$-asymptotic extensions of~$K$.
\end{enumerate}
In particular,  for each ungrounded $H$-asymptotic extension $L$  of~$K$,
$$\nwt_L P=\nwt P\leq 1, \quad \ndeg_L P=1,\quad \nval_L P=\nval P, \quad  v_L^{\ev}(P)=v^{\ev}(P).$$
\end{lemma} 
\begin{proof}
From [ADH, 13.7.10] we obtain an active $\phi_0$ in $K$ and $a,b\in K$ with $a\preceq b$ such that in $K^{\phi_0}\{Y\}$, 
either $P^{\phi_0} \sim^\flat_{\phi_0} a+bY$ or $P^{\phi_0}\ \sim^\flat_{\phi_0}\ b\,Y'$ (so $b\ne 0$). 
In the first case, set $A(Y):=a+bY$, $w:=0$; in the second case, set $A(Y):=bY'$, $w:=1$. Then
$P^{\phi_0}  = A + R$ where~$R\prec_{\phi_0}^\flat b\asymp P^{\phi_0}$ in $K^{\phi_0}\{Y\}$. 

Let $L$ be an $H$-asymptotic extension of $K$.
Then~$R\prec_{\phi_0}^\flat P^{\phi_0}$ remains
true in~$L^{\phi_0}\{Y\}$, and if $f\preceq \phi_0$ is active in $L$, then
$P^f = (P^{\phi_0})^{f/\phi_0}=(f/\phi_0)^w A+R^{f/\phi_0}$ where~$R^{f/\phi_0}\prec_{\phi_0}^\flat P^f$  by Lemma~\ref{lem:same ndeg} and the remark following its proof. As to~$v_L^{\ev}(P)=v^{\ev}(P)$ for ungrounded $L$,  the  identity from  [ADH, 11.1.15]
defining these quantities  shows that both are $vb$ in case  (i), and  $v(b)-v(\phi_0)$ in case (ii).
\end{proof}

\noindent
Lemma~\ref{lem:13.7.10}  has the following consequence, partly generalizing Corollary~\ref{cor:sum of nwts}:
%and to be used in \cite{ADHld}:

\begin{cor} \label{cor:13.7.10}
Suppose $K$ is $\upl$-free, $A\in K[\der]^{\neq}$ and $L$ is an ungrounded $H$-asymptotic extension
of~$K$. Then for $\gamma\in\Gamma$ the quantities $\nwt_A(\gamma)\leq 1$ and $v^{\ev}_A(\gamma)$ do not change when passing from $K$ to $L$; in particular,
$$\exc^{\ev}(A)\ =\ \big\{\gamma\in\Gamma: \nwt_A(\gamma)=1\big\}\  =\  \exc^{\ev}_L(A)\cap\Gamma.$$
\end{cor}

\noindent
This leads to a variant of Corollary~\ref{cor:size of excev}:

\begin{cor}\label{cor:size of excev, strengthened} Suppose $K$ is $\upl$-free. Then $|\exc^{\ev}(A)|\leq \order A$ for all $A\in K[\der]^{\ne}$.
\end{cor}
\begin{proof} By [ADH, 10.1.3], $K$ is pre-$\d$-valued, hence by [ADH, 11.7.18] it has an $\upo$-free $H$-asymptotic extension. It remains to appeal to Corollaries~\ref{cor:sum of nwts} and~\ref{cor:13.7.10}.
\end{proof}

\noindent
For completeness we next state a version of Lemma~\ref{lem:13.7.10} for~$\ndeg P=0$;
the proof using [ADH, 13.7.9] is similar, but simpler, and hence omitted.

\begin{lemma}\label{lem: 13.7.9}
Suppose $K$ is $\upl$-free and $\ndeg P=0$. Then there are an active~$\phi_0$ in $K$ and~$a\in K^\times$ 
 such that   $P^f \sim_{\phi_0}^\flat a$ for all active~$f\preceq\phi_0$ in all $H$-asymptotic extensions of $K$.
\end{lemma}

\noindent
In particular, for $K$, $P$ as in Lemma~\ref{lem: 13.7.9}, no $H$-asymptotic extension of $K$ contains any $y\preceq 1$ such that $P(y)=0$. 

For general $P$ and $\upo$-free $K$ we can still do better than stated earlier:

\begin{lemma}\label{lem:13.6.12} 
Suppose $K$ is $\upo$-free. Then there are $w\in\N$, $A\in K[Y]^{\neq}$, and an active $\phi_0$ in $K$ such that for all active
$f\preceq \phi_0$ in all  $H$-asymptotic extensions of~$K$, 
$$P^{f}\ \sim^\flat_{\phi_0} \  (f/\phi_0)^w A(Y)(Y')^w. $$
For such $w$, $A$, $\phi_0$ we have for any ungrounded $H$-asymptotic extension $L$ of $K$,
\begin{align*} \nwt_L P\ &=\ w, &  \ndeg_L P&\ =\ \deg A+w,\\
 \nval_L P\ &=\ \val A+w, &   v^{\ev}_L(P)&\ =\ v(A)-wv(\phi_0).
 \end{align*} 
\end{lemma}
\begin{proof} By [ADH, 13.6.11] we have active $\phi_0$ in $K$ and $A\in K[Y]^{\ne}$ such that
$$P^{\phi_0}\ =\ A\cdot (Y')^w + R,\quad w:= \nwt P, \quad R \in K^{\phi_0}\{Y\},\ R\prec_{\phi_0}^\flat P^{\phi_0}.$$
(Here $\phi_0$ and $A$ are the $e$ and $aA$ in [ADH, 13.6.11].) The rest of the argument is just like in the second part of the proof of
Lemma~\ref{lem:13.7.10}. 
\end{proof}

\subsection*{Remarks on newton position}
For the next lemma we put ourselves in the setting of~[ADH, 14.3]: $K$ is $\upo$-free, $P\in K\{Y\}^{\neq}$,
and~$a$ ranges over~$K$. 
Recall that~$P$ is said to be in {\it newton position at $a$}\/ if $\nval P_{+a}=1$. \index{differential polynomial!newton position} \index{newton position}

Suppose $P$ is in newton position at~$a$; then $A:=L_{P_{+a}}\in K[\der]^{\ne}$. 
Recall the definition of~$v^{\ev}(P,a)=v^{\ev}_K(P,a)\in \Gamma_\infty$:
if $P(a)=0$, then $v^{\ev}(P,a)=\infty$; if $P(a)\neq 0$, then $v^{\ev}(P,a)=vg$ where $g\in K^\times$
satisfies $P(a) \asymp (P_{+a})^\phi_{1,\times g}$ eventually, that is,
$v_{A^\phi}(vg)=v\big(P(a)\big)$ eventually. In the latter case
$\nwt_A(vg)=0$, that is, $vg\notin\exc^{\ev}(A)$, and 
$v^{\ev}_A(vg)=v\big(P(a)\big)$, since $v_{A^\phi}(vg)=v^{\ev}_A(vg)+ \nwt_A(vg)v\phi$ eventually. 
For any~$f\in K^\times$, $P^f$ is also in newton position at $a$, and  $v^{\ev}(P^f,a)=v^{\ev}(P,a)$.
Note also
that~$P_{+a}$ is in newton position
at $0$ and $v^{\ev}(P_{+a},0)=v^{\ev}(P,a)$.
Moreover, in passing from $K$ to  an $\upo$-free extension, $P$  remains in newton position at $a$ and~$v^{\ev}(P,a)$ does not change, by Lemma~\ref{lem:13.6.12}.

\medskip
\noindent
{\it In the rest of this subsection $P$ is in newton position at $a$, and~$\hat a$ is an element of an $H$-asymptotic extension $\hat K$ of $K$  such that $P(\hat a)=0$.}\/ (We allow~$\hat a\in K$.)
We first generalize part of [ADH, 14.3.1], with a similar proof:

\begin{lemma}\label{lem:14.3.1} 
$v^{\ev}(P,a)>0$ and  $v(\hat a-a) \leq v^{\ev}(P,a)$.
\end{lemma}
\begin{proof}
This is clear if $P(a)=0$. Assume $P(a)\neq 0$.
Replace~$P$,~$\hat a$,~$a$ by~$P_{+a}$,~$\hat a-a$,~$0$, respectively, to arrange $a=0$.
Recall that $K^\phi$ has small derivation.
Set~$\gamma:=v^{\ev}(P,0)\in\Gamma$ and take $g\in K$ with~$vg=\gamma$. 
Now $(P_1^\phi)_{\times g} \asymp P_0$, eventually, and~$\nval P=1$ gives $P(0) \prec P_1^\phi$, eventually, hence~$g\prec 1$.
Moreover, for~$j\ge 2$, $P_1^\phi\succeq P_j^\phi$, eventually, so~$(P_1^\phi)_{\times g}\succ (P_j^\phi)_{\times g}$, eventually, by [ADH, 6.1.3]. 
Thus for~${j\geq 1}$ we have~$(P_{\times g}^\phi)_j = (P_j^\phi)_{\times g} \preceq P(0)$, eventually;
in particular, there is no~$y\prec 1$ in any $H$-asymptotic  extension of $K$ with $P_{\times g}(y)=0$.
Since~$P(\hat a)=0$, this yields~$v(\hat a) \leq \gamma=v^{\ev}(P,0)$.
\end{proof}

\noindent
Here is a situation where $v(\hat a-a) = v^{\ev}(P,a)$:

\begin{lemma}\label{lem:14.3 complement} Suppose $\Psi$ is cofinal in $\Psi_{\hat K}$, $\hat a - a\prec 1$, and
 $v(\hat a-a)\notin\exc_{\hat K}^{\ev}(A)$ where~$A:=L_{P_{+a}}$.  Then $v(\hat a-a)=v^{\ev}(P,a)$. 
\end{lemma}
\begin{proof} Note that $\hat K$ is ungrounded, so $\exc_{\hat K}^{\ev}(A)$ is defined, and $\hat{K}$ is pre-$\d$-valued. 
As in the proof of Lemma~\ref{lem:14.3.1} we arrange $a=0$. As an asymptotic subfield  of~$\hat K$, $K\langle \hat a\rangle$
is pre-$\d$-valued. Hence $K\langle \hat a \rangle$ is $\upo$-free by Theorem~\ref{thm:ADH 13.6.1}.
 The remarks
preceding  Lemma~\ref{lem:14.3.1} then allow us to replace~$K$ by $K\langle \hat a\rangle$  to arrange $\hat a\in K$.
The case $\hat a=0$ is trivial, so assume $0\ne \hat a\prec 1$. Now $\nval P=1$ gives for $j\ge 2$ that $P_1^\phi\succeq P_j^\phi$, eventually, hence
$(P_1^\phi)_{\times \hat a}\succ (P_j^\phi)_{\times \hat a}$, eventually, by~[ADH, 6.1.3]. Moreover, $P_1(\hat a)=A(\hat a)= A^\phi(\hat a)\asymp A^\phi \hat a$, eventually, using $v(\hat a)\notin \exc^{\ev}_{\hat K}(A)$ in the last step, so for $j\ge 2$, eventually
$$P_1(\hat a)\ \asymp\ (P_1^\phi)_{\times \hat a}\ \succ\ (P_j^\phi)_{\times \hat a}\ \succeq P_j^\phi(\hat a)\ =\ P_j(\hat a). $$
Also $P_1(\hat a)\ne 0$, since $A^\phi\hat a \ne 0$. Then $P(\hat a)=0$
gives $P(0)\asymp P_1(\hat a)$. 
%and even $P(0)\sim -P_1(\hat a)$, which we don't need.
Thus $v\big(P(0)\big)=v_{A^\phi}\big(v(\hat a)\big)$, eventually, so $v^{\ev}(P,0)=v(\hat a)$ by the definition of $v^{\ev}(P,0)$. 
\end{proof}

\begin{cor}
Suppose $\hat K$  is ungrounded and equipped with an ordering making it a pre-$H$-field, and assume $\hat a-a \prec 1$ and $v(\hat a-a) \notin \exc^{\ev}_{\hat K}(A)$ where $A:=L_{P_{+a}}$. Then $v(\hat a-a)=v^{\ev}(P,a)$.
\end{cor}
\begin{proof}
In view of Lemma~\ref{lemexc} and using [ADH, 14.5.11] we can extend $\hat K$ to arrange that it is an $\upo$-free newtonian Liouville closed $H$-field. Next, let $H$ be the real closure of the $H$-field hull of $K\langle \hat a\rangle$, all inside~$\hat K$. Then $H$ is $\upo$-free, by Theorem~\ref{thm:ADH 13.6.1}, and hence has a Newton-Liouville closure $L$ inside~$\hat K$~[ADH, 14.5]. Since $L\preccurlyeq \hat K$ by [ADH, 16.2.5], we have~$v({\hat a-a}) \notin \exc^{\ev}_{L}(A)$. 
%Hence we can replace $\hat K$ by
Now $L$ is $\d$-algebraic over $K$ by [ADH, 14.5.9], so $\Psi$ is cofinal in $\Psi_{L}$ by  Theorem~\ref{thm:ADH 13.6.1}. It remains to apply Lemma~\ref{lem:14.3 complement}. 
\end{proof}

%\subsection*{Liouville closed implies $1$-linearly newtonian} 
%The next lemma gives that every Liouville closed $H$-field is $1$-linearly newtonian:
%In this subsection $K$ is an ungrounded $H$-asymptotic field with $\Gamma\neq\{0\}$. 

%\begin{lemma}
%Let $P\in K\{Y\}$ satisfy $\order P=\deg P=1$. Then
%$$\ndeg P = 1 \quad\Longleftrightarrow\quad \ddeg P^\phi=1 \text{ eventually}  \quad\Longleftrightarrow\quad \ddeg P^\phi=1 \text{ for all $\phi$.}$$ 
%\end{lemma}
%\begin{proof}
%The first equivalence holds by definition. Suppose $\ndeg P=1$; to show $\ddeg P^\phi=1$ for all $\phi$, we can assume that $P=Y'+aY+b$ where $a,b\in K$. Then $P^\phi=\phi Y'+aY+b$ for each $\phi$; So if $a\succeq b$ then $\ddeg P^\phi=1$ for each $\phi$. If $a\prec b$, then $\ddeg P^\phi=1$ eventually implies $b\prec \phi$ eventually, and hence $b\prec \phi$ for each $\phi$; so again $\ddeg P^\phi=1$ for each $\phi$.
%\end{proof}

\subsection*{Newton position in the order $1$ case} 
{\it In this subsection $K$ is $\upl$-free, ${P\in K\{Y\}}$ has order~$1$, and $a\in K$.}\/ We basically copy here a definition and two lemmas from~[ADH, 14.3] with the $\upo$-free assumption there replaced by the weaker
$\upl$-freeness, at the cost of restricting $P$ to have order $1$. 

Suppose $\nval P=1$, $P_0\neq 0$. Then  [ADH, 11.6.17] yields $g\in K^\times$ with~$P_0\asymp P^\phi_{1,\times g}$, eventually. Since $P_0\prec P_1^\phi$, eventually, we have $g\prec 1$. 
Moreover, if~${i\geq 2}$, then $P_1^\phi\succeq P_i^\phi$, eventually, hence
$P^\phi_{1,\times g} \succ P^\phi_{i,\times g}$, eventually. Thus $\ndeg P_{\times g}=1$.

Define $P$ to be in {\bf newton position at $a$} if $\nval P_{+a}=1$. \index{differential polynomial!newton position} \index{newton position} Suppose $P$ is in newton position at $a$; set $Q:=P_{+a}$, so $Q(0)=P(a)$.
If $P(a)\neq 0$, then the above yields $g\in K^\times$ with
$P(a)=Q(0)\asymp Q^\phi_{1,\times g}$, eventually; as $vg$ does not depend on
the choice of such $g$, we set $v^{\ev}(P,a):=vg$. If $P(a)=0$, then we set~$v^{\ev}(P,a):=\infty\in\Gamma_\infty$.
 In passing from $K$ to  a $\upl$-free extension, $P$  remains in newton position at $a$ and~$v^{\ev}(P,a)$ does not change, by Lemma~\ref{lem:11.2.13 invariant}. {\em In the rest of this subsection we assume $P$ is in newton position at $a$}.

\begin{lemma}\label{newt4, order 1}
If $P(a)\ne 0$, then  there exists~$b\in K$ with the following properties: \begin{enumerate}
\item[\textup{(i)}] $P$ is in newton position at $b$, $v(a-b)= v^{\ev}(P,a)$, and $P(b)\prec P(a)$; 
\item[\textup{(ii)}] for all $b^*\in K$ with $v(a-b^*)\ge v^{\ev}(P,a)$: $\ P(b^*)\prec P(a)\Leftrightarrow a-b\sim a-b^*$;
\item[\textup{(iii)}] for all $b^*\in K$, if $a-b\sim a-b^*$, then $P$ is in newton position at $b^*$ and~$v^{\ev}(P,b^*)> v^{\ev}(P,a)$.
\end{enumerate}
\end{lemma} 

\noindent
This is shown as in [ADH, 14.3.2]. Next an analogue of [ADH, 14.3.3], 
with the same proof, but using Lemma~\ref{newt4, order 1} in place of [ADH, 14.3.2]:

\begin{lemma}\label{newt.imm, order 1} If there is no $b$ with~${P(b)=0}$ and 
$v(a-b) = v^{\ev}(P,a)$, then there is a divergent pc-sequence~$(a_\rho)_{\rho<\lambda}$ in~$K$, indexed by all ordinals $\rho$ smaller than 
some infinite limit ordinal $\lambda$, such that~$a_0=a$, $v(a_\rho-a_{\rho'})=v^{\ev}(P,a_\rho)$ for all $\rho<\rho'<\lambda$, and $P(a_\rho) \leadsto 0$.
\end{lemma}

\noindent
The next result is proved just like Lemma~\ref{lem:14.3.1}:

\begin{lemma}\label{lem:14.3.1, order 1} If $P(\hat a)=0$ with~$\hat a$ in an $H$-asymptotic extension of $K$, then
$v^{\ev}(P,a)>0$ and $v(\hat a-a) \leq v^{\ev}(P,a)$.
\end{lemma}

%In general $K\langle \hat a\rangle$ is not $\upl$-free (see Example~\ref{ex:rat as int and cofinality}), 

\noindent
Next an analogue of Lemma~\ref{lem:14.3 complement} using Propositions~\ref{prop:rat as int and cofinality} and~\ref{prop:upl-free and as int} in its proof: 

\begin{lemma}\label{lem:14.3 complement, order 1} 
Suppose $\hat a$ in an ungrounded $H$-asymptotic extension  $\hat{K}$ of $K$ satisfies $P(\hat a)=0$, $\hat a - a\prec 1$, and
 $v(\hat a-a)\notin\exc^{\ev}_{\hat K}(A)$, where $A:=L_{P_{+a}}$.  Then~$v({\hat a-a})  =  v^{\ev}(P,a)$. 
\end{lemma}
\begin{proof}
We arrange $a=0$ and assume $\hat a\neq 0$.
Then $L:=K\langle\hat a\rangle$ has asymptotic integration, by
Proposition~\ref{prop:upl-free and as int}, and $v(\hat a)\notin \exc^{\ev}_L(A)$ by
Lemma~\ref{lemexc, order 1} (applied with~$L$,~$\hat K$ in place of $K$, $L$).
Moreover, $\Psi$ is cofinal in $\Psi_L$
by Proposition~\ref{prop:rat as int and cofinality}.  As in the proof of Lemma~\ref{lem:14.3 complement} this leads to
$P_1(\hat a)=A(\hat a)=A^\phi(\hat a)\asymp A^{\phi}\hat a$, eventually, and then as in the rest of that proof we derive
$v^{\ev}(P,0)=v(\hat a)$. 
\end{proof}

\subsection*{Zeros of differential polynomials  of order and degree $1$}
{\it In this subsection~$K$ has asymptotic integration.}\/
%and small derivation.}\/ 
We fix a differential polynomial
$$P(Y)\ =\ a(Y'+gY-u)\in K\{Y\}\qquad (a,g,u\in K,\  a\neq 0),$$
and set $A:=L_P=a(\der+g)\in K[\der]$. Section~\ref{sec:logder} gives  for $y\in K$ the equivalence~$y\in \I(K)\Leftrightarrow vy>\Psi$, so by Section~\ref{sec:lindiff}, $\exc^{\ev}(A)=\emptyset\Leftrightarrow g\notin \I(K)+K^\dagger$,
and $v(\ker_{\hat K}^{\neq} A)\subseteq \exc^{\ev}(A)$ for each immediate $H$-asymptotic field extension $\hat K$ of $K$.
Thus:

\begin{lemma}\label{lem:at most one zero}
If $g\notin \I(K) + K^\dagger$, then each immediate $H$-asymptotic extension of $K$ contains at most one~$y$ such that~$P(y)=0$.
\end{lemma}

\noindent
If  $\der K=K$ and $g\in K^\dagger$, then $P(y)=0$ for some $y\in K$, and if moreover $K$ is $\d$-valued, then any $y$ in any immediate $H$-asymptotic extension of $K$ with $P(y)=0$ lies in $K$. (Lemma~\ref{0K}.)
If $y\prec 1$ in an immediate $H$-asymptotic extension of~$K$ satisfies~$P(y)=0$, then by [ADH, 11.2.3(ii), 11.2.1] 
%\marginpar{using errata update for 11.2.3} 
we have
$$ \nval P\ =\ \nval P_{+y}\ =\ \val P_{+y}\ =\ 1.$$
Lemma~\ref{newt.imm, order 1} 
 yields the following partial converse (a variant of \cite[Lemma~8.5]{VDF}):

\begin{cor}\label{cor:nmul=1}
Suppose $K$ is $\upl$-free and $\nval P=1$. Then there is a $y\prec 1$ in an immediate $H$-asymptotic extension of $K$ with $P(y)=0$.
\end{cor}
\begin{proof}
Replacing $K$ by its henselization and using [ADH, 11.6.7], we arrange that~$K$ is henselian.
Suppose that
$P$ has no zero in $\smallo$. Then $P$ is in newton position at~$0$,
and so   Lemma~\ref{newt.imm, order 1} yields a divergent pc-sequence $(a_{\rho})_{\rho<\lambda}$ in $K$, indexed by all ordinals~$\rho$ smaller than 
some infinite limit ordinal $\lambda$,   with $a_0=0$, $v(a_\rho-a_{\rho'})=v^{\ev}(P,a_\rho)$ for all $\rho<\rho'<\lambda$, and $P(a_{\rho}) \leadsto 0$. Since $\deg P=\order P=1$ and $K$ is henselian, $P$ is a minimal differential polynomial of~$(a_\rho)$ over $K$, and $v(a_\rho)=v^{\ev}(P,0)>0$ for all $\rho>0$.
Hence [ADH, 9.7.6] yields a pseudolimit~$y$ of~$(a_\rho)$ in an immediate asymptotic extension of $K$ with $P(y)=0$
and~$y\prec 1$, as required.
\end{proof}

\noindent
We say that $P$ is {\bf proper} if $u\neq 0$ and $g+u^\dagger\succ^\flat 1$. \index{proper}\index{differential polynomial!proper}\index{proper!differential polynomial} If $P$ is proper, then so is $bP$ for each $b\in K^\times$. For $\fm\in K^\times$ we have
$$P_{\times\fm}\ =\ a\fm\big(Y'+(g+\fm^\dagger)Y-u\fm^{-1}\big),$$
hence if $P$ is proper, then so is $P_{\times\fm}$.
If   $u\neq 0$, then
$P$ is proper iff~$a^{-1}A_{\ltimes u}=\der+(g+u^\dagger)$ is steep, as defined in Section~\ref{sec:lindiff}.
Note that
$$P^\phi\  =\ a\phi\big( Y'+(g/\phi)Y-(u/\phi) \big).$$

\begin{lemma}\label{lem:proper compconj}
Suppose $K$ has small derivation, and $P$ is proper. Then $P^\phi$ is proper $($with respect to  $K^\phi)$ for all $\phi\preceq 1$.
\end{lemma}
\begin{proof} Let $\phi\preceq 1$. Then we have
$\phi\asymp^\flat 1$ and hence $\phi^\dagger\asymp^\flat \phi' \preceq 1\prec^\flat g+u^\dagger$.
Thus
$$g +(u/\phi)^\dagger =(g+u^\dagger)-\phi^\dagger \sim^\flat g+u^\dagger\succ^\flat 1\succeq \phi,$$ hence
$(g/\phi)+\phi^{-1}(u/\phi)^\dagger\succ^\flat  1$ and so~${(g/\phi) + \phi^{-1}(u/\phi)^\dagger} \succ^\flat_\phi 1$.
Therefore $P^\phi$ is proper (with respect to $K^\phi$).
\end{proof}

\begin{lemma}\label{lem:proper evt}
Suppose $K$ is $\upl$-free and $u\neq 0$. Then there is an active $\phi_0$ in~$K$ such that for all $\phi\prec\phi_0$, $P^\phi$ is proper with
$g+(u/\phi)^\dagger\sim g+(u/\phi_0)^\dagger$.
\end{lemma}
\begin{proof}
The argument before Corollary~\ref{cor:prlemexc} yields an active $\phi_0$ in $K$ such that $u^\dagger +g-\phi^\dagger\succeq \phi_0$ for all $\phi\prec \phi_0$. For   such~$\phi$ we have
$\phi^\dagger - \phi_0^\dagger\prec \phi_0$  as noted just before~[ADH, 11.5.3], and so $(u/\phi)^\dagger+g \sim (u/\phi_0)^\dagger+g$. 
The argument before Corollary~\ref{cor:prlemexc} also gives $\phi^{-1}(u/\phi)^\dagger+g/\phi\succ^\flat_\phi 1$ eventually,
and if  $\phi^{-1}(u/\phi)^\dagger+g/\phi\succ^\flat_\phi 1$, then~$P^\phi$ is proper.
\end{proof}

\begin{lemma}\label{lem:proper nmul 1}
We have $\nval P=1$ iff $u\prec g$ or $u\in\I(K)$. Moreover, if $K$ is $\upl$-free, $\nval P=1$, and $u\neq 0$, then $u\prec^\flat_\phi g+(u/\phi)^\dagger$, eventually.
\end{lemma}
\begin{proof}
For the equivalence, note that the identity above for $P^\phi$ yields: $$\nval P=0\ \Longleftrightarrow\  u\succeq g, \text{ and $u/\phi\succeq 1$ eventually}.$$
Suppose $K$ is $\upl$-free, $\nval P=1$, and $u\neq 0$.
If $u\in\I(K)$, then~$u\prec \phi\prec^\flat_\phi g+(u/\phi)^\dagger$, eventually,
by Lemma~\ref{lem:proper evt}.
Suppose  $u\notin \I(K)$. Then $v(u)\in \Psi^{\downarrow}$ and $u\prec g$. Hence by [ADH, 9.2.11] we have~$(u/\phi)^\dagger \prec u\prec g$, eventually, and thus $u\prec g\sim g+(u/\phi)^\dagger$, eventually. Thus~$u\prec^\flat_\phi g+(u/\phi)^\dagger$, eventually.
\end{proof}

\noindent
Assume now $P(y)=0$ with $y$ in
an immediate $H$-asymptotic extension  of $K$; so~$A(y)=u$. Note: if~$vy\in\Gamma\setminus\exc^{\ev}(A)$, then $u\neq 0$. 
From Lemma~\ref{prlemexc} we get:

\begin{lemma}\label{lem:proper asymp}
If $K$ has small derivation, $P$ is proper, and $vy\in\Gamma\setminus\exc^{\ev}(A)$, then~$y\sim u/(g+u^\dagger)$.
\end{lemma}

\noindent
By Lemmas~\ref{lem:proper evt} and \ref{lem:proper asymp}, and using Lemma~\ref{lem:proper nmul 1} for the last part:

\begin{cor}\label{cor:proper}
If $K$ is $\upl$-free  and $vy\in\Gamma\setminus\exc^{\ev}(A)$, then  
$$y\sim  u/\big( g+(u/\phi)^\dagger\big) \quad\text{ eventually.}$$
If in addition $\nval P=1$, then $y\prec 1$. 
\end{cor}

%\noindent
%Together with  the order~$1$ version of Lemma~\ref{lem:14.3 complement}, this yields:
 
%\begin{cor}
%Suppose $K$ is $\upl$-free and $P$ is in newton position at $a$, where~$\hat a-a\prec 1$ and $v(\hat a-a)\notin\exc^{\ev}(L_P)$. Then $u_{+a}:=f^{-1}P(a)\neq 0$, and 
%$$v(\hat a-a) = v^{\ev}(P,a)=v(u_{+a})- v\big( g+(u_{+a}/\phi)^\dagger\big),\quad\text{ eventually.}$$
%\end{cor}

\subsection*{A characterization of $1$-linear newtonianity} 
{\it In this subsection $K$ has asymptotic integration.}\/
We first expand~[ADH, 14.2.4]:

\begin{prop}\label{prop:char 1-linearly newt} 
The following are equivalent:
\begin{enumerate}
\item[\textup{(i)}] $K$ is $1$-linearly newtonian;
\item[\textup{(ii)}] every $P\in K\{Y\}$ with $\nval P=\deg P=1$ and $\order P\leq 1$ has a zero in~$\smallo$;
\item[\textup{(iii)}] $K$ is $\d$-valued, $\upl$-free, and $1$-linearly surjective, with $\I(K)\subseteq K^\dagger$.
\end{enumerate}
\end{prop}
\begin{proof} 
The equivalence of (i) and (ii) is [ADH, 14.2.4], and the implication (i)~$\Rightarrow$~(iii) follows from [ADH, 14.2.2, 14.2.3, 14.2.5].
To show (iii)~$\Rightarrow$~(ii), suppose (iii) holds, and
let $g,u\in K$ and $P=Y'+gY-u$ with $\nval P=1$. We need to find~$y\in \smallo$ such that $P(y)=0$.
Corollary~\ref{cor:nmul=1} gives an element~$y\prec 1$ in an immediate $H$-asymptotic extension $L$ of $K$ with $P(y)=0$.
It suffices to show that then
$y\in K$ (and thus~${y\in \smallo}$). 
If $g\notin K^\dagger$, then  this follows from Lemma~\ref{lem:at most one zero}, using $\I(K)\subseteq K^\dagger$ and $1$-linear surjectivity of $K$; if $g\in K^\dagger$, then this follows from Lemma~\ref{0K} and~$\der K=K$.
\end{proof}

\noindent
By the next corollary, each Liouville closed $H$-field is $1$-linearly newtonian:

\begin{cor} \label{cor:Liouville closed => 1-lin newt}
Suppose $K^\dagger=K$. Then the following are equivalent:
\begin{enumerate}
\item[$\mathrm{(i)}$] $K$ is $1$-linearly newtonian;
\item[$\mathrm{(ii)}$] $K$ is $\d$-valued and $1$-linearly surjective;
\item[$\mathrm{(iii)}$] $K$ is $\d$-valued and $\der K=K$.
\end{enumerate}
\end{cor}
\begin{proof}
Note that $K$ is $\upl$-free by [ADH, remarks following 11.6.2].
Hence the equivalence of (i) and (ii) follows from Proposition~\ref{prop:char 1-linearly newt}. For the equivalence of (ii) with~(iii), see [ADH, example following 5.5.22]. 
%If~(iii) holds, then $K$ is $\d$-valued and $K^\dagger=K$ by [ADH, 14.2.5], and $K$ is $1$-linearly surjective by [ADH, 14.2.2], showing (iii)~$\Rightarrow$~(i). Suppose now that (ii) holds; to show (iii), we only need to verify that  $K$ is $1$-linearly newtonian.  Let $\phi$ range over the active elements of $K$. Replacing $K$ by $K^\phi$ for suitable $\phi$ we arrange that~$K$ has small derivation. Let $P\in K\{Y\}$ satisfy $\order P=\deg P=\ndeg P=1$; we need to find $y\in\mathcal O$ with $P(y)=0$. For this we can assume that $P=Y'+aY-g$ where $a,g\in K$, $g\neq 0$. By $P^\phi=\phi Y'+aY-g$ and $\ndeg P=1$ we obtain    $a\succeq g$ or $a\prec g\prec \phi$ eventually; thus $a\succeq g$ or $a\prec g\prec 1$. Let $A:=\der+a$; then $v_A(0)=v(A)=\min\{0,va\}$. Since $K$ is $1$-linearly surjective and $1$-pv-closed, we can take $y,z\in K$ with $A(y)=g$ and $z\neq 0$, $A(z)=0$. If $y\asymp z$, take $c\in C$ with $y\sim cz$ and replace $y$ by $y-cz$; thus we arrange $y\nasymp z$. We have  $\exc(A^\phi) =  vz+\Gamma^\flat_\phi$ for each $\phi$, and $\exc^{\ev}(A)=\bigcap_\phi\exc(A^\phi)=\{vz\}$.  After replacing~$K$,~$P$ by~$K^\phi$,~$\phi^{-1}P^\phi$ for suitable   $\phi$, we can thus assume that $vy\notin\exc(A)=vz+\Gamma^\flat$. If $a\succeq g$, then $v_A(0)\leq va \leq vg=v_A(vy)$, and if $a\prec g\prec 1$, then $v_A(0)=\min\{0,va\}=0<vg=v_A(vy)$. In both cases we obtain $y\in\mathcal O$, since $v_A\colon\Gamma\to\Gamma$ is strictly increasing.
\end{proof}

\subsection*{Linear newtonianity descends} 
{\em In this subsection $H$ is $\d$-valued with valuation ring 
$\mathcal{O}$ and constant field $C$. Let~$r\in\N^{\ge 1}$}.
If $H$ is $\upo$-free, $\Gamma$ is divisible, and~$H$ has a newtonian algebraic extension $K=H(C_K)$, then
$H$ is also newtonian, by [ADH, 14.5.6]. Here is an analogue of this for $r$-linear newtonianity:

\begin{lemma}\label{lem:descent r-linear newt}
Let $K=H(C_K)$ be an algebraic asymptotic extension of $H$ which is $r$-linearly newtonian. Then $H$ is $r$-linearly newtonian.
\end{lemma}
\begin{proof}
Take a basis $B$ of the $C$-linear space $C_K$ with $1\in B$, and let $b$ range over~$B$.
We have $H(C_K)=H[C_K]$, and $H$ is linearly disjoint from $C_K$ over $C$ [ADH, 4.6.16],
so $B$ is a basis of the $H$-linear space $H[C_K]$.
Let $P\in H\{Y\}$ with~$\deg P=1$ and~$\order(P)\leq r$ be quasilinear; then $P$ as element of~$K\{Y\}$ remains
quasilinear, since~${\Gamma_K=\Gamma}$ by  [ADH, 10.5.15].
Let $y\in\mathcal O_K$ be a zero of $P$. Take~$y_b\in H$ ($b\in B$)  with
$y_b=0$ for all but finitely many $b$ and $y=\sum_b y_b\,b$. Then $y_b\in\mathcal O$ for all~$b$, and
$$0\ =\ P(y)\ =\ P_0 + P_1(y)\ =\ P_0 + \sum_b P_1(y_b)b,$$
so $P(y_1) = P_0 + P_1(y_1) = 0$.
\end{proof}
 
\noindent
Thus if $H[\imag]$ with $\imag^2=-1$ is $r$-linearly newtonian, then $H$ is $r$-linearly newtonian. 

\subsection*{Cases of bounded order} 
{\it In the rest of this section $r\in\N^{\ge 1}$}. \index{H-asymptotic field@$H$-asymptotic field!strongly $r$-newtonian} \index{strongly!$r$-newtonian}
Define~$K$ to be {\bf strongly $r$-newtonian\/} if $K$ is $r$-newtonian and for each divergent pc-sequence~$(a_{\rho})$ in $K$ with minimal differential polynomial $G(Y)$ over $K$ of order~$\le r$ we have~$\ndeg_{\boldsymbol a} G =1$,
where ${\boldsymbol a}:=c_K(a_{\rho})$. Given $P\in K\{Y\}^{\ne}$, a {\bf $K$-external zero of $P$\/} \index{zero!K-external@$K$-external}\index{K-external zero@$K$-external zero}\index{differential polynomial!K-external zero@$K$-external zero}
is an element~$\hat{a}$ of some immediate asymptotic extension $\hat{K}$ of $K$ such that~$P(\hat{a})=0$ and $\hat{a}\notin K$. Now [ADH, 14.1.11] extends as follows with the same proof: % the reason we assume $K$ is of $H$-type here is the use of 11.3.8 in that proof. 

\begin{lemma}\label{14.1.11.r} Suppose $K$ has rational asymptotic integration and $K$ is strongly $r$-newtonian. Then no $P\in K\{Y\}^{\ne}$ of order $\le r$ can have a
$K$-external zero.
\end{lemma} 

\noindent
The following is important in certain inductions on the order.

\begin{lemma}\label{rlcrln} Suppose $K$ has asymptotic integration, is $1$-linearly newtonian, and
$r$-linearly closed. Then $K$ is $r$-linearly newtonian.
\end{lemma}
\begin{proof} Note that $K$ is $\upl$-free and $\d$-valued by Proposition~\ref{prop:char 1-linearly newt}. 
Let~$P\in K\{Y\}$ be such that $\operatorname{nmul} P=\deg P=1$ and $\order P\le r$; by~[ADH, 14.2.6] it suffices to show that then $P$ has a zero in $\smallo$. By [ADH, proof of 13.7.10] we can compositionally conjugate, pass to an elementary extension, and multiply by an element of $K^\times$ to arrange that $K$ has small derivation, $P_0\prec^{\flat} 1$, and $P_1\asymp 1$. 
Let $A:= L_P$. The valuation ring of the flattening 
$(K, v^\flat)$ is $1$-linearly surjective by~[ADH, 14.2.1],
so all operators in $K[\der]$ of order $1$ are neatly
surjective in the sense of $(K, v^\flat)$. Since~$A$ splits over $K$,
we obtain from [ADH, 5.6.10(ii)] that~$A$ is neatly surjective in the sense of $(K,v^\flat)$. As $v^{\flat}(A)=0$ and $v^{\flat}(P_0)>0$, this gives~$y\in K$ with
$v^\flat(y)>0$ such that $P_0+A(y)=0$, that is, $P(y)=0$. 
\end{proof}  

\noindent
Using the terminology of $K$-external zeros,
we can add another item to the list of equivalent statements in Proposition~\ref{prop:char 1-linearly newt}:

\begin{lemma}\label{lem:char 1-linearly newt}
Suppose $K$ has asymptotic integration. Then we have:
\begin{align*} \text{$K$ is $1$-linearly newtonian}\ \Longleftrightarrow\ 
&\text{$K$ is $\upl$-free and no $P\in K\{Y\}$ with $\deg P=1$}\\
&\text{and $\order P=1$ has a $K$-external zero.}
\end{align*} 
\end{lemma}
\begin{proof}
Suppose $K$ is $1$-linearly newtonian. Then by (i)~$\Rightarrow$~(iii) in Proposition~\ref{prop:char 1-linearly newt}, $K$ is
$\upl$-free, $\d$-valued, $1$-linearly surjective, and~$\I(K)\subseteq K^\dagger$.
Let $P\in K\{Y\}$ where~$\deg P=\order P=1$ and  $y$ in an immediate asymptotic extension $L$ of~$K$ with~$P(y)=0$. Then [ADH, 9.1.2] and Corollary~\ref{cor:no new LDs} give $L^\dagger\cap K=K^\dagger$, so $y\in K$ by
Lemmas~\ref{0K} and~\ref{1K}. This gives the direction $\Rightarrow$.  
 The converse follows from Corollary~\ref{cor:nmul=1} and  
(ii)~$\Rightarrow$~(i) in Proposition~\ref{prop:char 1-linearly newt}.
\end{proof}

\noindent
Here is a higher-order version of Lemma~\ref{lem:char 1-linearly newt}:

\begin{lemma}\label{lem:char r-linearly newt}
Suppose $K$ is $\upo$-free. Then
\begin{align*} \text{$K$ is $r$-linearly newtonian}\quad \Longleftrightarrow\quad 
&\text{no $P\in K\{Y\}$ with $\deg P=1$ and $\order P\leq r$}\\
&\text{has a $K$-external zero.}
\end{align*}
\end{lemma}
\begin{proof}
Suppose $K$ is $r$-linearly newtonian. Then $K$ is $\d$-valued by Lemma~\ref{lem:ADH 14.2.5}.
Let  $P\in K\{Y\}$ be of degree~$1$ and order~$\leq r$, and let $y$ be in an immediate asymptotic extension $L$ of $K$
with $P(y)=0$. Then $A(y)=b$ for $A:=L_P\in K[\der]$, $b:=-P(0)\in K$. By~[ADH, 14.2.2] there is also a $z\in K$ with $A(z)=b$, hence~$y-z\in\ker_L A=\ker A$ by   [ADH, remarks after 14.2.9] and so $y\in K$. 
This gives the direction~$\Rightarrow$. For the converse note that every quasilinear  $P\in K\{Y\}$   has a zero~${\hat a\preceq 1}$ in an immediate asymptotic extension of $K$ by [ADH, 14.0.1 and subsequent remarks].
\end{proof}

\noindent
We also have the following $r$-version of [ADH, 14.0.1]:

\begin{prop}\label{14.0.1r} If $K$ is $\upl$-free and no $P\in K\{Y\}^{\ne}$ of order $\le r$ has a $K$-external zero, then $K$ is
$\upo$-free and $r$-newtonian.
\end{prop}

\begin{proof}
The $\upo$-freeness follows as before from [ADH, 11.7.13]. 
%The rest of the proof is as in the proof of $\Leftarrow$ in Lemma~\ref{lem:char 1-linearly newt}. \marginpar{shortened proof (no need to redo proof on  [ADH, p.~653] as indicated earlier)}
The rest of the proof is as in [ADH, p.~653] with $P$ restricted to have order $\le r$.
\end{proof}

\subsection*{Application to solving asymptotic equations} {\it Here $K$ is $\d$-valued, $\upo$-free, with small derivation, and $\fM$ is a monomial group of $K$.}\/ We let $a$,~$b$,~$y$ range over~$K$. 
In addition we fix a $P\in K\{Y\}^{\ne}$ of order $\le r$ and
a $\preceq$-closed set $\E\subseteq K^\times$. (Recall that $r\ge 1$.) This gives the asymptotic equation
\begin{equation}\label{eq:asymp equ}\tag{E}
P(Y)=0,\qquad Y\in \E.
\end{equation} 
This gives the following
$r$-version of [ADH, 13.8.8], with basically the same proof: 

\begin{prop}\label{1388} Suppose $\Gamma$ is divisible, no $Q\in K\{Y\}^{\ne}$ of order $\le r$ has a $K$-external zero, $d:=\ndeg_{\E} P\ge 1$, and there is no $f\in \E\cup\{0\}$ with $\operatorname{mul} P_{+f}=d$. Then \eqref{eq:asymp equ} has an unraveler.
\end{prop}

\noindent
Here is an $r$-version of [ADH, 14.3.4] with the same proof:

\begin{lemma} Suppose $K$ is $r$-newtonian. Let $g\in K^\times$ be an approximate zero of~$P$ with $\ndeg P_{\times g}=1$.
Then there exists $y\sim g$ such that $P(y)=0$. 
\end{lemma}

\noindent
For the next three results we assume the following:

\medskip\noindent
{\em  $C$ is algebraically closed, $\Gamma$ is divisible, and 
no $Q\in K\{Y\}^{\ne}$ of order $\le r$ has a $K$-ex\-ter\-nal zero}.

\medskip\noindent
These three results are $r$-versions of [ADH, 14.3.5, 14.3.6, 14.3.7] with the same proofs, using Propositions~\ref{14.0.1r} and~\ref{1388} instead of [ADH, 14.0.1, 13.8.8]:

\begin{prop} If $\ndeg_{\E}P > \operatorname{mul}(P)=0$, then \eqref{eq:asymp equ} has a solution. 
\end{prop}

\begin{cor}\label{wrdc} $K$ is weakly $r$-differentially closed.
\end{cor} 

\begin{cor} Suppose $g\in K^\times$ is an approximate zero of $P$. Then $P(y)=0$ for some $y\sim g$. 
\end{cor} 

\subsection*{A useful equivalence} {\it Suppose $K$ is $\upo$-free.}\/ 
(No small derivation or monomial group assumed.) Recall that $r\ge 1$. Here is an $r$-version of \cite[3.4]{Nigel19}: 

\begin{cor}\label{14.5.2.r} The following are equivalent: 
\begin{enumerate}
\item[\textup{(i)}] $K$ is $r$-newtonian;
\item[\textup{(ii)}] $K$ is strongly $r$-newtonian;
\item[\textup{(iii)}] no $P\in K\{Y\}^{\ne}$ of order $\le r$ has a $K$-external zero.
\end{enumerate}
\end{cor}
\begin{proof} Since $K$ is $\upo$-free  it has rational asymptotic integration [ADH, p.~515]. Also,
if $K$ is $1$-newtonian, then $K$ is henselian [ADH, p.~645] and $\d$-valued [ADH, 14.2.5].
For~(i)~$\Rightarrow$~(ii), use 
\cite[3.3]{Nigel19}, for (ii)~$\Rightarrow$~(iii), use Lemma~\ref{14.1.11.r}, and for (iii)~$\Rightarrow$~(i), use Proposition~\ref{14.0.1r}.
\end{proof}  

\noindent
Next an $r$-version of [ADH, 14.5.3]:

\begin{cor}\label{14.5.3.r} Suppose $K$ is $r$-newtonian, $\Gamma$ is divisible, and $C$ is algebraically closed. Then $K$ is weakly $r$-differentially closed, so $K$ is $(r+1)$-linearly closed and thus~${(r+1)}$-linearly newtonian.  
\end{cor}
\begin{proof} To show that $K$ is weakly $r$-differentially closed we arrange by compositional conjugation and passing to a suitable elementary extension that $K$ has small derivation and $K$ has a monomial group. Then
$K$ is weakly $r$-differentially closed by 
Corollaries~\ref{wrdc} and~\ref{14.5.2.r}. The rest uses [ADH, 5.8.9] and
Lemma~\ref{rlcrln}. 
\end{proof}

\subsection*{Complementing [ADH, 14.2.12]\astr}
In this subsection $P(Y)\in\mathcal O\{Y\}$ has order at most $r\in\N^{\geq 1}$. 

\begin{lemma}\label{lem:14.2.12 compl}
Let  $y\in K^\times$, $y'\preceq y\prec 1$, and $P(0)=P(y)$.
Then $L_{P}(y) \prec y$.
\end{lemma}
\begin{proof}
Induction on $n$ gives $y^{(n)}\preceq y^{(n-1)}\preceq\cdots\preceq y\prec 1$ for all $n$. Hence
 if~$\i=(i_0,\dots,i_r)\in\N^{1+r}$, $\abs{\i}\geq 2$, then $y^{\i}=y^{i_0}(y')^{i_1}\cdots (y^{(r)})^{i_r} \preceq y^{\abs{\i}}\prec y$.
Now
$$P(y)\ =\ P(0)+ L_{P}(y)+\sum_{\abs{\i}\geq 2} P_{\i}\,y^{\i}, $$
so $L_{P}(y)+\sum_{\abs{\i}\geq 2} P_{\i}\,y^{\i}=0$, and thus $L_{P}(y) \prec y$.
\end{proof}

\noindent
We extend the residue map $a\mapsto\res a\colon\mathcal O\to\k:=\res(K)$
to  the ring morphism
$$p\mapsto\res p\ \colon\ \mathcal O[Y]\to\k[Y], \qquad  Y\mapsto Y.$$
For~${w\in\N}$ we let $P_{[w]}$ be  the isobaric part of $P$ of weight $w$, as in [ADH, 4.2].
Thus~$p:=P_{[0]}\in \mathcal O[Y]$.  

\begin{cor}\label{cor:14.2.12 compl}
Suppose  the derivation of $K$ is very small, and let
 $a\in\mathcal O$, $y\in\smallo$ with $P(a)=P(a+y)$ and $(\res p)'(\res a)\neq 0$. Then $y'\succeq y$.
%If $y\in\smallo$ with~$y' \prec y$ and $P(a+y) = 0$ is $y=0$.
%Then for all $a,b\in\mathcal O$ with~$P(a)=P(b)=0$ and $\bar{a}=\bar{b}=\alpha$ we have $v(a-b)>\Delta$.
\end{cor}
\begin{proof}
Put $R:=\sum_{w\geq 1} P_{[w]}=P-p$. Now $a^{(n)}\prec 1$ for all $n\geq 1$, so
$\big(\frac{\partial R}{\partial Y}\big)(a) \prec 1$.
Towards a contradiction, assume $y'\prec y$. Then  $L_{P_{+a}}(y)\prec y$ by Lemma~\ref{lem:14.2.12 compl} applied to $P_{+a}$ in place of $P$. Induction on $n$ gives $y^{(n)}\prec y^{(n-1)}\prec\cdots\prec  y\prec 1$ for all $n$ and so~$L_{R_{+a}}(y)=\sum_n \big(\frac{\partial R}{\partial Y^{(n)}}\big)(a)y^{(n)}\prec y$.
Together with $L_{P_{+a}}(y)=p'(a)y+L_{R_{+a}}(y)$ and $p'(a)\asymp 1$ this yields the desired contradiction.
\end{proof}

\noindent
In the next corollary we assume that $K$ has   asymptotic integration.
We let $\phi$ range over active elements of $K$, and
we let 
$\smallo^\flat_\phi=\{f\in \smallo: f'\succeq f\phi\}$
be the maximal ideal of the flattened valuation ring  of  $K^\phi$; cf.~[ADH, pp.~406--407]. 

\begin{cor}
Suppose $K$ is $r$-newtonian. Let
$u\in\mathcal O$ and $A\in \k[Y]$ be such that~$A(\res u)=0$, $A'(\res u)\neq 0$, and $D_{P^\phi}\in\k^\times A$, eventually. Then $P$ has a  zero in~$u+\smallo$, and for all zeros
$a,b\in u+\smallo$  of $P$ we have:
$a-b\in\smallo^\flat_\phi$, eventually. 
\end{cor}
\begin{proof}
For the first claim, see [ADH, 14.2.12].
Suppose $D_{P^\phi}\in\k^\times A$ and
take~$\fm\in K^\times$ with~$\fm\asymp P^\phi$, so $Q:=\fm^{-1}P^\phi\in K^\phi\{Y\}$ and $q:=Q_{[0]}$ satisfy~$vQ=0$
and~$\res q\in\k^\times A$.
Note that $K$ is $\d$-valued by [ADH, 14.2.5]; hence $K^\phi$ has very small derivation.
Let $a,b\in u+\smallo$ and $P(a)=P(b)=0$; then $y:=b-a\in\smallo$, and so
Corollary~\ref{cor:14.2.12 compl} applied to $K^\phi$, $Q$ in place of $K$, $P$ yields $y'\succeq y\phi$.
\end{proof}

\subsection*{Newton polynomials of Riccati transforms\astr}
{\it In this subsection we assume that $K$ has small derivation and asymptotic integration.}\/
Let  
$$A\ =\ a_0+a_1\der+\cdots+a_r\der^r\in K[\der]\qquad\text{ where~$a_0,\dots,a_r\in K$, $a_r\neq 0$,}$$
with Riccati transform
$$R\ :=\ \operatorname{Ri}(A)\ =\ a_0R_0(Z)+a_1R_1(Z)+\cdots+a_rR_r(Z)\in K\{Z\},$$
and  set
$$P\ :=\ a_0+a_1Z+\cdots+a_rZ^r\in K[Z].$$
We equip the differential fraction field $K\langle Z\rangle$ of
$K\{Z\}$  with the gaussian extension of the valuation of $K$ and likewise with $K^\phi$ instead of $K$. Then $K^\phi\langle Z\rangle$ is a valued differential field with small derivation by [ADH, 6.3].
(Although $K$ is asymptotic, $K\langle Z\rangle$ is not, by [ADH, 9.4.6].)

\begin{lemma}\label{lem:Rphi}
Eventually, $R^\phi \sim P$.
\end{lemma}
\begin{proof} It is enough to show that $R_n(Z)^\phi\sim Z^n$ eventually. For $n=0,1$ we have~$R_n(Z)=Z^n$.
Now $R_{n+1}(Z)=ZR_n(Z)+\der\big(R_n(Z)\big)$, so by [ADH, 5.7.1], $$R_{n+1}(Z)^\phi\ =\ ZR_n(Z)^\phi + \phi\derdelta\big(R_n(Z)^\phi\big),
\qquad \derdelta\ :=\ \phi^{-1}\der.$$
Assuming $R_n(Z)^\phi\sim Z^n$ eventually, this yields $R_{n+1}(Z)^\phi\sim Z^{n+1}$ eventually .
\end{proof}

\begin{remark}
Suppose $K$ is $\d$-valued and equipped with a monomial group.
In [ADH, 13.0.1] we associate to $Q\in K\{Z\}^{\neq}$  its Newton polynomial $N_Q\in C\{Z\}$ such that~$D_{Q^\phi}=N_Q$, eventually.  Then $N_R=D_P\in C[Z]$ by Lemma~\ref{lem:Rphi}.
\end{remark}

\noindent
Next an application of Lemma~\ref{lem:Rphi}. 
For simplicity,   assume $vA=0$, so~$P\in\mathcal O[Z]$, $vP=0$.
We also let $Q\mapsto\res Q\colon \mathcal O[Z]\to\k[Z]$ be the extension of the residue map~$a\mapsto\res a\colon\mathcal O\to\k$ to a ring morphism with~$Z\mapsto Z$. 

\begin{cor}\label{cor:simple Riccati zero}
Suppose $K$ is $(r-1)$-newtonian, $r\ge 1$. Then for all $\alpha\in\k$  with $\res P(\alpha)=0$ and $(\res P)'(\alpha)\neq 0$ there is $a\in \mathcal O$ with $R(a)=0$
and $\res a=\alpha$. 
\end{cor}
\begin{proof} If $r=1$, use $R=P=a_0+a_1Z$. Assume $r\ge 2$. 
By Lemma~\ref{lem:Rphi} we have~$D_{R^\phi}\in\k^\times \cdot\res P$, eventually, so we can apply [ADH, 14.2.12] to
$R$, $\res P$ in the role of $P$, $A$ there. 
\end{proof}

\noindent
In the rest of this subsection we assume $A\in\mathcal O[\der]$ is monic.
To what extent is the zero $a$ of $R$ in Corollary~\ref{cor:simple Riccati zero} unique?   Corollaries~\ref{cor:Ri unique, 1} and~\ref{cor:Ri unique, 2} below give answers to this question.

\begin{lemma}\label{rarby}
Let $a,b\in\mathcal O$ be such that $R(a)=R(b)=0$ and $y:=b-a\prec 1$. Then
$y'\preceq y$.
\end{lemma}
\begin{proof}
Replace $R$ by $R_{+a}$ to arrange $a=0$, $b=y$, so $a_0=0$. 
Note that $r\geq 1$.
Towards a contradiction, assume~$y\prec y'$. Then  $R_n(y)\sim y^{(n-1)}$ for all $n\geq 1$ by Lemma~\ref{lem:Riccati bd flat},
and $y\prec y' \prec \cdots \prec y^{(r-1)}$,  so $$R(y)\ =\ a_1R_1(y)+\cdots+a_{r-1}R_{r-1}(y)+R_r(y)\ \sim\  y^{(r-1)},$$
hence $R(y)\neq 0$, a contradiction.
\end{proof}

\begin{cor}\label{cor:Ri unique, 1}
Suppose $K$ has very small derivation, 
$\alpha\in\k$ is a simple zero of $\res P$, and $a,b\in\mathcal O$, $R(a)=R(b)=0$, and $\res a=\res b=\alpha$.
Then for $y:=b-a$ we have $y'\asymp y$.
\end{cor}
\begin{proof}
We have~$y'\preceq y$ by Lemma~\ref{rarby}, and $y'\succeq y$
by Corollary~\ref{cor:14.2.12 compl} applied to~$R$ in place of $P$, using
$R_{[0]}=P$.
\end{proof}

\noindent
In the next result we assume 
 $K=H[\imag]$ where $H$ is a real closed differential subfield of $K$   such that the valuation ring $\mathcal O_H:=\mathcal O\cap H$ of $H$ is convex with respect to the ordering of $H$  and  $\mathcal O_H=C_H+\smallo_H$. So $C=C_H[\imag]$ and $\mathcal O=C+\smallo$ (cf.~remarks after Corollary~\ref{cor:logder abs value}).
We identify $C$ with $\k$ via the residue morphism $\mathcal O\to\k$. 

\begin{cor}\label{cor:Ri unique, 2}
Let $\alpha\in C$ be a simple zero of $\res P\in C[Z]$ such that for all
zeros $\beta\in C$ of $\res P$ we have $\Re \alpha \leq \Re \beta$. Then there is at most one $a\in\mathcal O$ with~$R(a)=0$ and~$\res a=\alpha$.
\end{cor}
\begin{proof} Let $a\in \mathcal{O}$, $R(a)=0$, and $\res{a}=\alpha$. 
Towards a contradiction, suppose $b\in \mathcal O$, $b\ne a$,  $R(b)=0$, and $\res{b}=\alpha$. 
By Lemma~\ref{RiPQ} we may replace~$a$,~$R$ by~$0$,~$R_{+a}$ to
arrange $a=0$. Then $0\neq b\prec 1$ and $a_0=R(0)=0$, so $P=QZ$ where~$Q\in \mathcal O[Z]$ and all zeros of $\res Q$ in $C$ have
nonnegative real part.  Moreover
$b'\asymp b$ by Corollary~\ref{cor:Ri unique, 1}. Take $c\in C^\times$ with $b'\sim bc$. 
By Lemma~\ref{lem:Riccati bd flat} and [ADH, 9.1.4(ii)]  we get $R_n(b)\sim b^{(n-1)} \sim bc^{n-1}$ for $n\geq 1$.
Now $R=a_1R_1+\cdots + a_rR_r$ gives~$Q=a_1+\cdots + a_rZ^{r-1}$, so $R(b)\in b\cdot\big(Q(c)+\smallo\big)$, 
hence $(\res Q)(c)=0$ in view of $R(b)=0$, and thus $\Re c \geq 0$.  On the other hand, $b\prec 1$ and Corollary~\ref{cor:10.5.2 variant} give~$\Re(b^\dagger)<0$, so $\Re c<0$, a contradiction.
\end{proof}

\newpage  

\part{The Universal Exponential Extension}\label{part:universal exp ext}

\medskip

\noindent
Let $K$ be an algebraically closed differential field.
In Section~\ref{sec:univ exp ext} below we extend~$K$ in a canonical way to a differential integral  domain $\Univ=\Univ_K$
whose differential fraction field has the same constant field $C$ as $K$,
called the {\it universal exponential extension}\/ of~$K$. (The universal exponential extension of $\mathbb T[\imag]$ appeared in \cite{JvdH} in the guise of ``oscillating transseries''; we explain
the connection at the end of Section~\ref{sec:eigenvalues and splitting}.)
The underlying ring of $\Univ$ is a   group ring of a certain abelian group over~$K$, and we therefore first review some relevant basic facts about such group rings
in Section~\ref{sec:group rings}. The main feature of $\Univ$ is that if $K$ is $1$-linearly surjective, then   each $A\in K[\der]$ of order $r\in\N$ which splits over $K$ has~$r$  many $C$-linearly independent zeros in $\Univ$. This is explained in Section~\ref{sec:eigenvalues and splitting}, after some differential-algebraic preliminaries in Sections~\ref{sec:splitting} and~\ref{sec:self-adjoint}, where we consider a novel kind of {\it spectrum}\/ of a linear differential operator over a differential field.
In Section~\ref{sec:valuniv} we introduce for $H$-asymptotic~$K$ with small derivation and asymptotic integration
the  {\it ultimate exceptional values}\/ of a given linear differential operator $A\in K[\der]^{\neq}$. These  
help to isolate the zeros of~$A$ in~$\Univ$ much like
the exceptional values of $A$ help to locate the zeros of~$A$ in immediate asymptotic extensions of $K$ as in~Section~\ref{sec:lindiff}.
In Section~\ref{sec:ueeh} below we discuss the analytic meaning of~$\Univ$ when $K$ is the algebraic closure of a Liouville closed Hardy field containing $\R$ as a subfield. 

\section{Some Facts about Group Rings}\label{sec:group rings}

\noindent
{\it In this section $G$ is a torsion-free abelian group, written multiplicatively, $K$ is a field, and~$\gamma$,~$\delta$ range over $G$.}\/
For use in Section~\ref{sec:univ exp ext} below we recall some facts about the group ring $K[G]$: a commutative $K$-algebra with $1\ne 0$ that contains $G$ as a subgroup of its multiplicative group~$K[G]^\times$ and
which, as a $K$-linear space, decomposes as \index{group ring}
$$K[G]\ =\ \bigoplus_\gamma K\gamma\qquad\text{(internal direct sum)}.$$
Hence for any $f\in K[G]$ we have a unique family $(f_\gamma)$ of elements of $K$, with $f_\gamma=0$  for all but finitely many $\gamma$,
such that 
\begin{equation}\label{eq:elt of group ring}
f\ =\ \sum_\gamma f_\gamma \gamma.
\end{equation}
We define the support of $f\in K[G]$ as above by
$$\supp(f)\ :=\ \{\gamma:\, f_\gamma\neq 0\}\ \subseteq\ G.$$
{\em In the rest of this section~$f$,~$g$,~$h$ range over $K[G]$}.
For any $K$-algebra $R$, every group morphism $G\to R^\times$ extends uniquely to a $K$-algebra morphism $K[G]\to R$. 

\medskip
\noindent
Clearly $K[G]^\times \supseteq K^\times G$; in fact:

\begin{lemma}\label{lem:only trivial units}
The ring $K[G]$ is an integral domain and $K[G]^\times = K^\times G$.
\end{lemma}
\begin{proof} We take an ordering of~$G$ making $G$ into an ordered abelian group; see~[ADH, 2.4]. Let $f,g\neq 0$ and set 
$$\gamma^-:=\min\supp(f),\ \gamma^+:=\max\supp(f),\quad 
 \delta^-:=\min\supp(g),\ \delta^+:=\max\supp(g);$$
so  $\gamma^- \leq \gamma^+$ and $\delta^-\leq \delta^+$.
We have 
$(fg)_{\gamma^- \delta^-} = f_{\gamma^-}g_{\delta^-}\neq 0$, and likewise with~$\gamma^+$,~$\delta^+$ in place of $\gamma^-$, $\delta^-$.
In particular, 
$fg\neq 0$, showing that~$K[G]$ is an integral domain.
Now suppose $fg=1$. Then $\supp(fg)=\{1\}$, hence $\gamma^- \delta^-=1=\gamma^+\delta^+$, so $\gamma^-=\gamma^+$, and thus $f\in K^\times G$.
\end{proof}

\begin{lemma}\label{lem:Frac U not alg closed} 
Suppose $K$ has characteristic $0$ and $G\neq\{1\}$. Then the fraction field $\Omega$ of $K[G]$ is not algebraically closed.
\end{lemma}
\begin{proof} 
Let $\gamma\in G\setminus\{1\}$ and $n\geq 1$. We claim that there is no $y\in\Omega$ with~$y^2=1-\gamma^n$.
For this, first replace $G$ by its divisible hull to arrange that $G$ is divisible. 
Towards a contradiction, suppose $f,g\in K[G]^{\neq}$ and $f^2=g^2(1-\gamma^n)$.
Take a divisible subgroup~$H$ of  $G$ that is complementary to the smallest divisible subgroup $\gamma^{\Q}$ of $G$ containing $\gamma$, so $G=H\gamma^{\Q}$ and $G\cap \gamma^{\Q}=\{1\}$. Then
$K[G]\subseteq K(H)[\gamma^{\Q}]$ (inside $\Omega$), so we may
replace $K$, $G$ by $K(H)$, $\gamma^{\Q}$ to arrange~$G=\gamma^{\Q}$. 
For suitable~$m\geq 1$ we apply the $K$-algebra automorphism of~$K[G]$ given by $\gamma\mapsto \gamma^m$ to arrange
$f,g\in K[\gamma,\gamma^{-1}]$ (replacing $n$ by $mn$). 
Then replace~$f$,~$g$ by~$\gamma^m f$,~$\gamma^m g$ for suitable $m\ge 1$  
to arrange~$f,g\in K[\gamma]$. Now use that~$1-\gamma$ is a prime divisor of~$1-\gamma^n$ of multiplicity~$1$  in the UFD $K[\gamma]$ to get a contradiction.
\end{proof}

\noindent
The $K$-linear map 
$$f\mapsto \operatorname{tr}(f):=f_1\ \colon\  K[G]\to K$$ is called the
{\bf trace} of $K[G]$. \index{group ring!trace} \index{trace} Thus
$$\operatorname{tr}(fg)\ =\ \sum_\gamma f_\gamma g_{\gamma^{-1}}.$$
We claim that $\operatorname{tr}\circ\sigma=\operatorname{tr}$ for every automorphism $\sigma$ of the $K$-algebra $K[G]$. This invariance comes from an intrinsic description of $\operatorname{tr}(f)$ as follows: given $f$ we have a unique finite set $U\subseteq K[G]^\times = K^\times G$
such that $f=\sum_{u\in U} u$ and $u_1/u_2\notin K^\times$ for all distinct
$u_1, u_2\in U$; if $U\cap K^\times=\{c\}$, then~$\operatorname{tr}(f)=c$;  if $U\cap K^\times=\emptyset$, then~$\operatorname{tr}(f)=0$. 
If $G_0$ is a subgroup of $G$ and $K_0$ is a  subfield  of $K$, then~$K_0[G_0]$ is a subring of $K[G]$, and the trace of $K[G]$ extends the trace of $K_0[G_0]$.

\subsection*{The automorphisms of $K[G]$} For a commutative group $H$, 
written multiplicatively, $\Hom(G,H)$ denotes the set of group morphisms $G\to H$, made into a group by pointwise multiplication. Any $\chi\in\Hom(G,K^\times)$---sometimes called a {\em character}---gives  a $K$-algebra automorphism $f\mapsto f_\chi$ of $K[G]$ defined by
\begin{equation}\label{eq:fchi}
f_\chi\ :=\ \sum_\gamma f_\gamma \chi(\gamma)\gamma.
\end{equation}
This yields a group action of $\Hom(G,K^\times)$ on $K[G]$ by $K$-algebra automorphisms:
$$\Hom(G,K^\times)\times K[G]\to K[G], \qquad (\chi, f)\mapsto f_{\chi}.$$ 
Sending $\chi\in\Hom(G,K^\times)$ to  $f\mapsto f_\chi$ yields an   embedding of the group $\Hom(G,K^\times)$ into the group~$\Aut(K[G]|K)$ 
of automorphisms of the $K$-algebra $K[G]$; its 
 image is the (commutative) subgroup of $\Aut(K[G]|K)$ consisting of the $K$-algebra automorphisms~$\sigma$ of $K[G]$ such that $\sigma(\gamma)/\gamma\in K^\times$ for all $\gamma$. Identify $\Hom(G,K^\times)$ with its image under this embedding. From $K[G]^\times=K^\times G$ we obtain $\sigma(K^\times G)=K^\times G$ for all  $\sigma\in\Aut(K[G]|K)$, and using this one verifies easily that $\Hom(G,K^\times)$ is a normal subgroup of $\Aut(K[G]|K)$. 
 We also have the group embedding 
 $$\Aut(G)\ \to\ \Aut(K[G]|K)$$ assigning to each
$\sigma\in \Aut(G)$ the unique automorphism of the $K$-algebra $K[G]$ extending $\sigma$. Identifying $\Aut(G)$ with its image in $\Aut(K[G]|K)$ via this embedding we have $\Hom(G,K^\times)\cap\Aut(G)=\{\text{id}\}$ and 
$\Hom(G,K^\times)\cdot \Aut(G)=\Aut(K[G],| K)$ inside $\Aut(K[G]|K)$, and thus 
 $\Aut(K[G]|K)=\Hom(G,K^\times)\rtimes\Aut(G)$, an internal semidirect product of subgroups of $\Aut(K[G]|K)$.

\subsection*{The gaussian extension}
{\it In this subsection $v\colon K^\times\to\Gamma$ is a valuation on the field $K$.}\/
We extend $v$ to a map $v_{\g}\colon K[G]^{\neq}\to\Gamma$ by setting
\begin{equation}\label{eq:gaussian ext}
v_{\g}f\ :=\ \min_\gamma vf_\gamma\qquad \text{($f\in K[G]^{\neq}$ as in \eqref{eq:elt of group ring}).}
\end{equation}

\begin{prop} The map $v_{\g}\colon K[G]^{\neq}\to\Gamma$ is a valuation on the domain $K[G]$.
\end{prop}
\begin{proof} 
We can reduce to the case that $G$ is finitely generated,
since $K[G]$ is the directed union of its subrings $K[G_0]$ with $G_0$ a finitely generated subgroup of $G$. 
Then we have a group isomorphism~$G\to \Z^n$ inducing a $K$-algebra isomorphism~$K[G]\to K[X_1, X_1^{-1},\dots, X_n, X_n^{-1}]$ (with distinct indeterminates $X_1,\dots,X_n$) under which~$v_{\g}$ corresponds to the gaussian extension of the valuation of $K$ to $K(X)$ restricted to its subring
$K[X_1, X_1^{-1},\dots, X_n, X_n^{-1}]$; see [ADH, Section~3.1]. 
\end{proof}

\noindent
We call $v_{\g}$ the {\bf gaussian extension} of the valuation of $K$ to~$K[G]$.\index{group ring!gaussian extension of a valuation}\index{gaussian extension}\index{valuation!gaussian extension}\index{extension!gaussian} We
denote by~$\preceq_{\g}$ the dominance relation on $\Omega:=\Frac(K[G])$ associated to the extension of~$v_{\g}$ to a valuation on the field $\Omega$  [ADH, (3.1.1)], with corresponding asymptotic relations~$\asymp_{\g}$ and~$\prec_{\g}$. For the subring $\mathcal O[G]$ of $K[G]$ generated by $G$ over $\mathcal O$ we have
$$\mathcal O[G]\ =\ \{ f:\, f\preceq_{\g} 1\}.$$
The residue morphism~$\mathcal O\to\k:=\mathcal O/\smallo$
extends to a surjective ring mor\-phism $\mathcal O[G]\to\k[G]$ with $\gamma\mapsto\gamma$ for all~$\gamma$ and whose kernel
is the ideal $$\smallo[G]:=\ \{f :\, f\prec_{\g} 1\}$$
of~$\mathcal O[G]$. Hence this ring morphism induces an isomorphism  $\mathcal O[G]/\smallo[G]\cong\k[G]$.
If~$G_0$ is subgroup of $G$ and $K_0$ is a valued subfield of $K$, then the restriction of $v_{\g}$ to a valuation on $K_0[G_0]$ is the gaussian extension of the valuation of $K_0$ to $K_0[G_0]$.\label{p:vg} 

\subsection*{An inner product and two norms}
{\it In the rest of this section $H$ is a real closed subfield of $K$ such that $K=H[\imag]$ where $\imag^2=-1$.}\/ In later use $H$ will be a Hardy field, which is why we use the letter $H$ here. Note that the only nontrivial automorphism of the (algebraically closed) field $K$ over $H$ is {\em complex conjugation}\/: 
$$z=a+b\imag \mapsto \overline{z}:=a-b\imag\qquad (a,b\in H).$$ 
For $f$ as in \eqref{eq:elt of group ring} we set
$$ f^*\ :=\ \sum_\gamma \overline{f_\gamma} \gamma^{-1},$$
so $(f^*)^*=f$, and $f\mapsto f^*$ lies in $\Aut\!\big(K[G]|H\big)$. 
We define the function $$(f,g)\mapsto \langle f,g\rangle\ :\ K[G]\times K[G]\to K$$ by
$$\langle f,g\rangle\ :=\ \operatorname{tr}\!\big(f g^*\big)\ =\ \sum_\gamma f_\gamma  \overline{g_\gamma}.$$
One verifies easily that this  is a ``positive definite hermitian form'' on the $K$-linear space $K[G]$: it is additive on the left and on the right, and for all $f$,~$g$ and all~$\lambda\in K$: $\langle \lambda f, g \rangle=\lambda \langle f,g \rangle$, $\langle g,f\rangle=\overline{\langle f,g\rangle}$, 
$\langle f, f\rangle \in H^{\ge}$, and $\langle f, f\rangle=0\Leftrightarrow f=0$, and thus also $\langle f, \lambda g \rangle=\overline{\lambda}\langle f,g \rangle$. (Hermitian forms are usually defined only on $\mathbb{C}$-linear spaces and are $\mathbb{C}$-valued, which is why we used quote marks, as we do below for {\em norm\/} and {\em orthonormal basis}; see \cite[Chapter~XV, \S{}5]{Lang} for
the more general case.)
Note:
$$\langle f,gh\rangle\ =\ \operatorname{tr}\!\big(f(g h)^*\big)\ =\ \big\langle f g^*,h\big\rangle.$$

\begin{lemma}\label{lem:inner prod intrinsic}
Let $u,w\in K[G]^\times$. If $u\notin K^\times w$, then $\langle u,w\rangle = 0$, and if $u\in K^\times w$, then $\langle u,w\rangle = uw^*$.
\end{lemma}
\begin{proof}
Take $a,b\in K^\times$ and $\gamma$, $\delta$ such that $u=a\gamma$, $w=b\delta$.
If $u\notin K^\times w$, then~$\gamma\neq\delta$, so~$\langle u,w\rangle = 0$.
If $u\in K^\times w$, then $\gamma=\delta$, hence $\langle u,w\rangle = a\overline{b} = uw^*$.
\end{proof}

\noindent
For $z\in K$ we set $\abs{z}:=\sqrt{z\overline{z}}\in H^{\ge}$, and then define $\dabs{\,\cdot\,}\colon K[G]\to H^{\geq}$ by
% denote   
% by  that is, 
%with $\abs{z}=\sqrt{z\overline{z}}\in F^{\ge}$ for $z\in K$ we have
$$\dabs{f}^2\ =\ \langle f,f\rangle\ =\ \sum_\gamma\, \abs{f_\gamma}^2.$$
As in the case $H=\R$ and $K=\mathbb{C}$ one derives the Cauchy-Schwarz Inequality:
$$|\langle f,g\rangle|\ \le\ \dabs{f}\cdot\dabs{g}.$$
Thus $\dabs{\,\cdot\,}$ is a ``norm'' on the $K$-linear space $K[G]$: for all $f,g$ and all $\lambda\in K$,  
$$\dabs{f+g}\le \dabs{f} + \dabs{g}, \quad \dabs{\lambda f}=\abs{\lambda}\cdot\dabs{f}, \quad \dabs{f}=0\Leftrightarrow f=0.$$
Note that $G$ is an ``orthonormal basis'' of $K[G]$ with respect to $\langle\ ,\,\rangle$, and  $f_\gamma= \langle f,\gamma\rangle$. 
We also use the function $\dabs{\,\cdot\,}_1\colon K[G]\to H^{\geq}$ given by
$$\dabs{f}_1\ :=\ \sum_\gamma\, \abs{f_\gamma},$$
which is a ``norm'' on $K[G]$ in the sense of obeying the same laws as we
mentioned for $\dabs{\,\cdot\,}$.  \index{group ring!norms}
The two ``norms'' are in some sense equivalent:
\[%\begin{equation}\label{eq:dabs(f)}
\dabs{f}\ \leq\ \dabs{f}_1\ \leq\ \sqrt{n}\dabs{f}\qquad\text{($n:=\abs{\supp(f)}$).}
\]%\end{equation}
where the first inequality follows from the triangle inequality for $\dabs{\,\cdot\,}$ and the second is of Cauchy-Schwarz type. Moreover:

\begin{lemma}\label{lem:submult}
Let $u\in K[G]^\times$. Then
$\dabs{f u} = \dabs{f}\,\dabs{u}$ and $\dabs{f u}_1 = \dabs{f}_1\,\dabs{u}_1$.
\end{lemma}
\begin{proof} %Since $K[G]^\times=K^\times G$, it suffices to prove this for $u\in G$. 
We have
$$\dabs{f\gamma}\ =\ \langle f\gamma, f\gamma\rangle\ =\ \big\langle f\gamma \gamma^*,f \big\rangle\ =\ 
\big\langle f ,f\big\rangle\ =\
\dabs{f} $$
using $\gamma^*=\gamma^{-1}$. 
Together with $K[G]^\times=K^\times G$ this yields
the first claim; the second claim follows easily from the definition of $\dabs{\,\cdot\,}_1$.
\end{proof}

\begin{cor}\label{cor:submult}
$\dabs{fg}\leq\dabs{f}\,\cdot\dabs{g}_1$ and $\dabs{fg}_1\leq\dabs{f}_1\,\cdot\dabs{g}_1$.
\end{cor}
\begin{proof}
By  the triangle inequality for $\dabs{\,\cdot\,}$ and  the previous lemma,
$$\dabs{fg}\  \leq\ \sum_\gamma\, \dabs{f g_\gamma\gamma}\ 
=\ \sum_\gamma\, \dabs{f}\, \dabs{g_\gamma\gamma}\ =\ \dabs{f} \sum_\gamma\, \abs{g_\gamma}\  =\ \dabs{f}\,\dabs{g}_1.$$
The inequality involving $\dabs{fg}_1$ follows likewise.
\end{proof}

\noindent
In the next lemma we let $\chi\in\Hom(G,K^\times)$; recall from \eqref{eq:fchi} the automorphism~$f\mapsto f_\chi$ of the $K$-algebra $K[G]$. 

\begin{lemma}\label{lem:commuting with f*}
 $(f_\chi)^*=(f^*)_\chi$ iff~$\abs{\chi(\gamma)}=1$ for all $\gamma\in\supp(f)$.
\end{lemma}
\begin{proof}
Let $a\in K$; then $\big((a\gamma)_\chi\big){}^*=\overline{a\chi(\gamma)}\gamma^{-1}$ and
$\big((a\gamma)^*\big){}_\chi = \overline{a}\chi(\gamma)^{-1}\gamma^{-1}$.
\end{proof}

\begin{cor}\label{cor:commuting with f*}
Let $\chi\in\Hom(G,K^\times)$ with $\abs{\chi(\gamma)}=1$ for all $\gamma$. Then~$\langle f_\chi,g_\chi \rangle=\langle f,g\rangle$ for all $f$, $g$, and hence~$\dabs{f_\chi}=\dabs{f}$ for all $f$.
\end{cor}
\begin{proof}
Since $\operatorname{tr}\circ\sigma=\operatorname{tr}$ for every automorphism $\sigma$ of the $K$-algebra $K[G]$, 
$$\langle f_\chi,g_\chi\rangle\ =\ \operatorname{tr}\!\big(f_\chi(g_\chi)^*\big)\ =\ \operatorname{tr}\!\big( (fg^*)_\chi \big)\ =\
\operatorname{tr}(fg^*)\ =\ \langle f,g\rangle,$$
where we use Lemma~\ref{lem:commuting with f*} for the second equality.
\end{proof}

\subsection*{Valuation and norm}
Let $v\colon H^\times\to\Gamma$ be a convex valuation on the ordered field~$H$, 
extended uniquely to a valuation~$v\colon K^\times\to\Gamma$ on the field $K=H[\imag]$, so~$a\asymp\abs{a}$ for~$a\in K$. % \marginpar{noted in {\tt mN} just before 6.3} 
 (See the remarks before Corollary~\ref{cor:10.5.2 variant}.) Let  $v_{\g}\colon K[G]^{\neq}\to\Gamma$ 
be the gaussian extension of $v$, given by \eqref{eq:gaussian ext}.   

\begin{lemma}\label{absval}
$\dabs{f}_1\preceq 1\Leftrightarrow f\preceq_{\g} 1$, and  $\dabs{f}_1\prec 1\Leftrightarrow f\prec_{\g} 1$.
\end{lemma}
\begin{proof} Using that the valuation ring of $H$ is convex we have
$$\dabs{f}_1=\sum_\gamma\, \abs{f_\gamma}\preceq 1\ \Longleftrightarrow\  \text{$\abs{f_\gamma}\preceq 1$ for
all $\gamma$} \ \Longleftrightarrow\ \text{$f_\gamma\preceq 1$ for
all $\gamma$} \ \Longleftrightarrow\ f\preceq_{\g} 1.$$
 Likewise one shows:
$\dabs{f}_1\prec 1\Leftrightarrow f\prec_{\g} 1$.
\end{proof}

\begin{cor}\label{cor:valuation and norm}
$\dabs{f}\asymp\dabs{f}_1\asymp_{\g} f$.
\end{cor}
\begin{proof} This is trivial for $f=0$, so assume $f\neq 0$. Take $a\in H^>$ with $a\asymp_{\g} f$, and replace~$f$ by $f/a$, to arrange $f\asymp_{\g} 1$. Then $\dabs{f}\asymp \dabs{f}_1\asymp_{\g} 1$ by Lemma~\ref{absval}. 
\end{proof}

\section{The Universal Exponential Extension}\label{sec:univ exp ext}

\noindent
As in [ADH, 5.9], given a differential ring $K$, a {\it differential $K$-algebra\/} is a
differential ring $R$ with a morphism $K \to R$ of differential rings. If $R$ is a differential ring extension of a differential ring $K$ we consider $R$ as a differential $K$-algebra via the in\-clu\-sion~${K \to R}$. \index{differential algebra}

\subsection*{Exponential extensions}
{\em In this subsection $R$ is a differential ring and $K$ is a differential subring of $R$}.
Call $a \in R$   {\bf exponential over~$K$} if $a' \in a K$. \index{element!exponential over}\index{exponential!element}
Note that if  
$a\in R$ is exponential over $K$, then~$K[a]$ is a differential subring of $R$.
If~$a\in R$ is exponential over $K$ and $\phi\in K^\times$, then $a$, as element of the differential ring extension~$R^\phi$ of~$K^\phi$, is exponential over~$K^\phi$.
Every   $c\in C_R$ is exponential over~$K$, and every
$u\in K^\times$ is exponential over $K$.  
If $a,b\in R$ are exponential over $K$, then so is $ab$,
and if $a\in R^\times$ is exponential over~$K$, then so is~$a^{-1}$. Hence the 
units of~$R$ that are exponential over~$K$ form a subgroup~$E$ of the  group~$R^\times$ of units of~$R$ with~$E\supseteq C_R^\times\cdot K^\times$;
if~${R=K[E]}$, then
we call~$R$   {\bf exponential over~$K$}.\index{extension!exponential}\index{exponential!extension}\index{differential algebra!exponential extension} An 
{\bf exponential extension of~$K$\/} is a differential ring extension of $K$ that is exponential over $K$. 
If $R=K[E]$ where $E$ is a set  of elements of $R^\times$ which are exponential over $K$, then $R$ is exponential over $K$.
If $R$ is an exponential extension of~$K$ and $\phi\in K^\times$, then $R^\phi$ is an exponential extension of~$K^\phi$. 
The following lemma is extracted from the proof of \cite[Theorem~1]{Rosenlicht75}:

\begin{lemma}[Rosenlicht]\label{lem:Rosenlicht lin indep} 
Suppose $K$ is a field and $R$ is an integral domain with
differential fraction field $F$. Let $I\neq R$ be a differential ideal of $R$, and
let~$u_1,\dots,u_n\in R^\times$ \textup{(}$n\geq 1$\textup{)} be exponential
over $K$ with $u_i\notin u_j C_F^\times K^\times$ for $i\neq j$. Then $\sum_i u_i\notin I$.
\end{lemma}
\begin{proof}
Suppose $u_1,\dots, u_n$ is a counterexample with minimal $n\geq 1$.
% for which there are units $u_1,\dots,u_n$ of $R$, exponential over $K$, with $u_i\notin u_j C_F^\times K^\times$ for $i\neq j$ and $\sum_i u_i\in I$.
Then $n\geq 2$ and $\sum_i u_i'\in I$, so 
$$\sum_i u_i'-u_1^\dagger\sum_i u_i\ =\  \sum_{i>1} (u_i/u_1)^\dagger u_i\in I.$$ 
Hence
~$(u_i/u_1)^\dagger=0$ and thus $u_i/u_1\in  C_F^\times$, for all $i>1$, a contradiction.
\end{proof}

\begin{cor}\label{roscor}
Suppose $K$ is a field and~$F=K(E)$ is a differential field extension of $K$ with $C_F=C$, where $E$ is a subgroup of $F^\times$ whose elements are exponential over $K$. Then $\{y\in F^\times:\ y \text{ is exponential over}~$K$\}=K^\times E$.
%Then the subgroup
%of~$F^\times$ consisting of the elements which are exponential over~$K$ is generated by~$K^\times E$.% \textup{(}thus $F^\dagger\cap K=R^\dagger\cap K=K^\dagger+E^\dagger$\textup{)}.
\end{cor}
\begin{proof}
Let $y\in F^\times$ be exponential over~$K$. Take $K$-linearly independent $u_1,\dots,u_n$ in $E$
and $a_1,\dots,a_n,b_1,\dots,b_n\in K$ with $b_j\ne 0$ for some $j$, such that
$$y\ =\ \Big(\textstyle\sum_i a_i u_i\Big)\Big/\left(\textstyle\sum_j b_j u_j\right).$$
Then $\sum_j b_jyu_j-\sum_i a_iu_i=0$, and so Lemma~\ref{lem:Rosenlicht lin indep} applied with~$R=F$,~$I=\{0\}$ gives $b_jyu_j\in a_iu_iK^\times$
for some $i$, $j$ with~$a_i,b_j\neq 0$, and thus~$y\in K^\times E$.
\end{proof}

\begin{remark}
In the context of Corollary~\ref{roscor}, see~\cite[Theorem~1]{Rosenlicht75} for the structure of the group of elements
of $F^\times$ exponential over~$K$, for finitely generated $E$.
\end{remark}

\begin{lemma}\label{lem:exponential map}
Suppose~$C_R^\times$ is divisible and $E$ is a subgroup of $R^\times$ containing $C_R^\times$.
Then there is a group morphism $e\colon E^\dagger\to E$ such that $e(b)^\dagger=b$
for all $b\in E^\dagger$. 
\end{lemma}
\begin{proof}
We have a short exact sequence of commutative groups
$$1 \to C_R^\times \xrightarrow{\ \iota\ } E \xrightarrow{\ \ell\ } E^\dagger \to 0,$$
where $\iota$ is the natural inclusion and $\ell(a):=a^\dagger$ for~$a\in E$. Since $C_R^\times$ is divisible, this sequence splits, which is what we claimed. 
\end{proof}

\noindent
Let $E$, $e$, $R$ be as in the previous lemma. Then  $e$ is injective, and its image is a complement of $C_R^\times$ in $E$.
Moreover, given also a group morphism $\tilde{e}\colon E^\dagger\to E$ such that~$\tilde{e}(b)^\dagger=b$ for all $b\in E^\dagger$, 
the map $b\mapsto e(b) \tilde{e}(b)^{-1}$ is a  group morphism~$E^\dagger\to C_R^\times$.

\medskip
\noindent
{\it In the rest of this section  $K$ is a differential field with algebraically closed constant field $C$ and divisible
group $K^\dagger$ of logarithmic derivatives.}\/ (These conditions are satisfied if $K$ is an algebraically closed differential field.) In the next subsection we show that   up to isomorphism over $K$ there is a unique exponential  extension $R$ of $K$ satisfying~$C_R=C$ and~${(R^\times)^\dagger=K}$.  By Lemma~\ref{lem:exponential map} we must then have a group embedding~$e\colon K\to R^\times$ such that~$e(b)^\dagger=b$
for all $b\in K$; this motivates the construction below.

\subsection*{The universal exponential extension}
We first describe  a certain exponential extension of $K$.
For this,  
take a {\bf complement\/} $\Lambda$ of $K^\dagger$, that is, a $\Q$-linear subspace of $K$ such that $K=K^\dagger\oplus \Lambda$ (internal direct sum of $\Q$-linear subspaces of~$K$).  Below~$\lambda$ ranges over~$\Lambda$.\index{group of logarithmic derivatives!complement}
Let $\ex(\Lambda)$ be a multiplicatively written abelian group,
isomorphic to the additive subgroup $\Lambda$ of~$K$, 
with isomorphism $\lambda\mapsto \ex(\lambda)\colon \Lambda\to \ex(\Lambda)$. 
Put 
$$\Univ\ :=\ K\big[\!\ex(\Lambda)\big],$$
the group ring of $\ex(\Lambda)$ over $K$, an integral domain. %  (Lemma~\ref{lem:only trivial units}).   
As $K$-linear space, 
$$\Univ\ =\ \bigoplus_\lambda K\ex(\lambda)\qquad (\text{an internal direct sum of $K$-linear subspaces}).$$
For every $f\in \Univ$ we have a unique family $(f_\lambda)$ in $K$ such that
$$f\  =\ \sum_{\lambda} f_{\lambda}  \ex(\lambda),$$
with $f_\lambda=0$ for all but finitely many $\lambda$; we call $(f_\lambda)$ the {\bf spectral decomposition} of $f$ (with respect to $\Lambda$).\index{extension!universal exponential}\index{universal exponential extension}\index{universal exponential extension!spectral decomposition}\index{spectral decomposition!of an element}\index{element!spectral decomposition}
We turn $\Univ$ into a differential ring extension of $K$
by 
$$\ex(\lambda)'\ =\ \lambda\ex(\lambda)\qquad \text{for all $\lambda$.}$$
(Think of $\ex(\lambda)$ as $\exp(\int \lambda)$.) Thus for $f\in \Univ$ with spectral decomposition $(f_\lambda)$,
$$f'\ =\ \sum_{\lambda} \big(f_{\lambda}'+\lambda f_{\lambda}\big) \ex(\lambda),$$
so $f'$ has spectral decomposition $(f_\lambda'+\lambda f_\lambda)$.  
Note that $\Univ$ is exponential over $K$ by Lem\-ma~\ref{lem:only trivial units}:  $\Univ^\times=K^\times  \ex(\Lambda)$, so  $(\Univ^\times)^\dagger=K^\dagger+\Lambda=K$. 

\begin{exampleNumbered}\label{ex:Q}
Let $K=C(\!( t^{\Q} )\!)$ be as in Example~\ref{ex:Kdagger}, so $K^\dagger=(\Q \oplus \smallo) t$. Take a $\Q$-linear subspace $\Lambda_{\operatorname{c}}$ of~$C$ with $C=\Q\oplus \Lambda_{\operatorname{c}}$
(internal direct sum of $\Q$-linear subspaces of $C$), and let
$$K_\succ\ :=\ \big\{f\in K:\ \supp(f) \succ 1 \big\},$$ 
a $C$-linear subspace of $K$.
Then $\Lambda:= (K_\succ \oplus \Lambda_{\operatorname{c}})t$ is a complement to $K^\dagger$,  and
hence~$t^{-1}\Lambda= K_\succ \oplus \Lambda_{\operatorname{c}}$ is a complement to $(K^t)^\dagger$ in $K^t$.
Moreover, if $L:=\operatorname{P}(C)\subseteq K$ is the differential field of Puiseux series over $C$
and $L_\succ:=K_\succ\cap L$, then $L_\succ  \oplus \Lambda_{\operatorname{c}}$ is a complement to $(L^t)^\dagger$. 
\end{exampleNumbered}

\noindent
A subgroup $\Lambda_0$ of $\Lambda$ yields a differential subring  $K\big[\!\ex(\Lambda_0)\big]$ of~$\Univ$ that is exponential over $K$ as well.  These differential subrings have a useful property.
Recall from~[ADH, 4.6] that a differential ring is said to be {\it simple\/} if $\{0\}$ is its only proper differential ideal.

\begin{lemma}\label{lem:U simple} 
Let $\Lambda_0$ be a subgroup of $\Lambda$. Then the differential subring
$K\big[\!\ex(\Lambda_0)\big]$ of $\Univ$ is simple. In particular, the differential ring $\Univ$ is simple. 
\end{lemma}
\begin{proof} 
Let $I\ne R$ be a differential ideal of
$R:=K\big[\!\ex(\Lambda_0)\big]$. Let
$f_1,\dots, f_n\in K^\times$ and let $\lambda_1,\dots, \lambda_n\in \Lambda_0$ be distinct
such that $
f=\sum_{i=1}^n  f_i\ex(\lambda_i)\in I$.
If $n\geq 1$, then  Lemma~\ref{lem:Rosenlicht lin indep} yields $i\neq j$
with $\ex(\lambda_i)/\ex(\lambda_j)=cg$ for some constant $c$ in the differential  fraction field of $\Univ$ and
some $g\in K^\times$, so by taking logarithmic derivatives,  $\lambda_i-\lambda_j\in K^\dagger$ and thus $\lambda_i=\lambda_j$,
a contradiction.  Thus $f=0$.
\end{proof}

\begin{cor}\label{cor:U simple}
Any morphism
$K\big[\!\ex(\Lambda_0)\big]\to R$ of differential $K$-algebras, with~$\Lambda_0$   a subgroup of $\Lambda$ and 
$R$ a differential ring extension of $K$, is injective.
\end{cor}

\noindent
The differential ring $\Univ$ is the directed union of its differential subrings of the form~$\Univ_0=K\big[\!\ex(\Lambda_0)\big]$ where $\Lambda_0$ is a finitely generated subgroup of $\Lambda$. These $\Univ_0$ are simple by Lemma~\ref{lem:U simple} and finitely generated as a $K$-algebra, hence their differential fraction fields have constant field $C$ by [ADH, 4.6.12]. Thus the  differential fraction field of $\Univ$  has constant field $C$. 

\begin{lemma}\label{lem:U minimal}
Suppose $R$ is an exponential extension of $K$ and 
$R_0$ is a differential subring of $R$ with $C_R^\times\subseteq C_{R_0}$ and $K\subseteq (R_0^\times)^\dagger$.
Then $R_0=R$.
\end{lemma}
\begin{proof}
Let $E$ be the group of units of $R$ that are exponential over $K$; so $R=K[E]$.
Given $u\in E$ we have $u^\dagger\in K\subseteq (R_0^\times)^\dagger$, hence we have $u_0\in R_0^\times$ with $u^\dagger=u_0^\dagger$, so~$u=cu_0$  with $c\in C_R^\times\subseteq C_{R_0}$. Thus $E\subseteq R_0$ and so $R_0=R$.
\end{proof}

\begin{cor}
Every endomorphism of the differential $K$-algebra $\Univ$ is an au\-to\-mor\-phism.
\end{cor}
\begin{proof}
Injectivity holds by Corollary~\ref{cor:U simple},  and surjectivity by Lemma~\ref{lem:U minimal}.
\end{proof}

\noindent
Every  exponential  extension of $K$ with constant field $C$ embeds into $\Univ$, and hence is an integral domain. More precisely:   

\begin{lemma}\label{lem:embed into U}  
Let $R$ be an exponential extension of $K$ such that~$C_R^\times$ is divisible, and set  $\Lambda_0:=\Lambda\cap (R^\times)^\dagger$, a subgroup of $\Lambda$.  
Then there exists a mor\-phism $K\big[\!\ex(\Lambda_0)\big]\to R$ of differential $K$-algebras. Any such morphism is injective, and if~$C_R=C$, then any such
morphism is an isomorphism.
\end{lemma}
\begin{proof}
Let $E$ be as in the proof of Lemma~\ref{lem:U minimal}, and let $e_E\colon E^\dagger\to E$ be the map~$e$ from Lemma~\ref{lem:exponential map}.  Since $E^\dagger=K^\dagger+\Lambda_0$  
we have 
\begin{equation}\label{eq:embed into U} 
E\ =\ C_R^\times\,  e_E(E^\dagger)\ =\ C_R^\times \, e_E(K^\dagger)\, e_E(\Lambda_0)\ =\ C_R^\times\, K^\times \, e_E(\Lambda_0).
\end{equation}
The group morphism $\ex(\lambda_0)\mapsto e_E(\lambda_0)\colon \ex(\Lambda_0)\to E$ ($\lambda_0\in \Lambda_0$) extends uniquely to a $K$-algebra morphism
$\iota\colon K\big[\!\ex(\Lambda_0)\big]\to R=K[E]$. One verifies easily that~$\iota$ is a  differential ring morphism. The injectivity claim follows from Corollary~\ref{cor:U simple}. 
If~$C_R=C$, then $E=K^\times e_E(\Lambda_0)$ by \eqref{eq:embed into U}, whence surjectivity. 
\end{proof}

\noindent
Recall that $\Univ$ is an exponential extension of $K$ with $C_{\Univ}=C$
and $(\Univ^\times)^\dagger = K$. By Lemma~\ref{lem:embed into U}, this property characterizes $\Univ$ up to isomorphism:

\begin{cor}\label{corcharexp} If $U$ is an exponential extension of $K$ such that $C_U=C$
and~$K\subseteq (U^\times)^\dagger$, then $U$ is isomorphic to $\Univ$
as a differential $K$-algebra.
\end{cor}

\noindent
Now $\Univ$ is also an exponential extension of $K$ with $C_{\Univ}=C$
and with the property that every exponential extension $R$ of $K$ with $C_R=C$ embeds into $\Univ$ as a differential $K$-algebra. This property
determines $\Univ$ up to isomorphism as well:

\begin{cor}\label{cor:UnivK} Suppose $U$ is an exponential extension of $K$ with $C_U=C$ such that every exponential extension $R$ of $K$ with $C_R=C$ embeds into $U$ as a differential $K$-algebra. Then $U$ is isomorphic to $\Univ$
as a differential $K$-algebra.  
\end{cor}
\begin{proof} Any embedding $\Univ\to U$ of
differential $K$-algebras gives $K\subseteq (U^\times)^\dagger$.  
\end{proof}

\noindent
The results above show to what extent $\Univ$ is independent of the choice of $\Lambda$. We call $\Univ$ the {\bf universal exponential extension of $K$}.  If we need to indicate the dependence of $\Univ$ on $K$ we denote it by $\Univ_K$.
By [ADH, 5.1.40] every~$y\in\Univ=K\{\ex(\Lambda)\}$ satisfies a linear differential equation $A(y)=0$ where $A\in K[\der]^{\neq}$; in the next section we isolate conditions on $K$ which ensure that every $A\in K[\der]^{\neq}$ has a zero~$y\in\Univ^\times=K^\times\ex(\Lambda)$.\index{universal exponential extension}\index{differential algebra!universal exponential extension}\label{p:UK}

\medskip
\noindent
Corollary~\ref{corcharexp} gives for $\phi\in K^\times$ an isomorphism
$\Univ_{K^\phi}\cong(\Univ_K)^\phi$ of differential $K^\phi$-algebras.  Next we investigate how $\Univ_K$ behaves when passing from $K$ to a differential field extension. Therefore,
{\it in the rest of this subsection $L$ is a differential field extension of $K$ with algebraically closed constant field $C_L$, and $L^\dagger$ is divisible.}\/ The next lemma relates the universal exponential extension $\Univ_L$ of $L$ to $\Univ_K$:

\begin{lemma}\label{lem:Univ under d-field ext}
The inclusion $K\to L$ extends to an embedding $\iota\colon\Univ_K\to\Univ_L$ of differential rings. The image of any such embedding~$\iota$ is contained in $K[E]$
where~$E:=\{u\in \Univ_L^\times:u^\dagger\in K\}$, and if $C_L=C$, then $\iota(\Univ_K)=K[E]$. 
\end{lemma}
\begin{proof} The differential subring
 $R:=K[E]$ of~$\Univ_L$ is exponential over $K$ with~$(R^\times)^\dagger=K$, hence Lemma~\ref{lem:embed into U} gives an embedding $\Univ_K\to R$ of differential $K$-algebras.
Let~$\iota\colon \Univ_K\to\Univ_L$ be any embedding of differential $K$-algebras. Then $\iota\big(\!\ex(\Lambda)\big)\subseteq E$, so~$\iota(\Univ_K)\subseteq R$; if $C_L=C$, then $\iota(\Univ_K) = R$
by Lemma~\ref{lem:U minimal}.
\end{proof}

\begin{cor}\label{cor:Univ under d-field ext}
If $L^\dagger\cap K=K^\dagger$ and $\iota\colon\Univ_K\to\Univ_L$ is an embedding of differential $K$-algebras,
then $L^\times\cap \iota(\Univ_K^\times) = K^\times$. 
\end{cor}
\begin{proof} Assume $L^\dagger\cap K=K^\dagger$ and identify $\Univ_K$ with a differential $K$-subalgebra of~$\Univ_L$ via an embedding $\Univ_K\to\Univ_L$ of differential $K$-algebras.
Let $a\in L^\times\cap \Univ_K^\times$; then~$a^\dagger\in L^\dagger\cap K=K^\dagger$, so $a = bc$ where $c\in C_L^\times$, $b\in K^\times$.
Now $c=a/b\in C_L^\times\cap\Univ_K^\times=C^\times$, since $\Univ_K$ has ring of constants~$C$. So $a\in K^\times$ as required. 
\end{proof}

\noindent
Suppose $L^\dagger\cap K=K^\dagger$. Then the subspace~$L^\dagger$ of the $\Q$-linear space $L$  has a complement
$\Lambda_L\supseteq \Lambda$.  We fix such $\Lambda_L$ and
extend~$\ex\colon\Lambda\to \ex(\Lambda)$ to a group isomorphism~$\Lambda_L\to\ex(\Lambda_L)$, also denoted by $\ex$, with~$\ex(\Lambda_L)$ a multiplicatively written commutative group  extending $\ex(\Lambda)$.
   Let $\Univ_L:=L\big[\!\ex(\Lambda_L)\big]$ be the corresponding universal exponential extension of~$L$.  
Then the natural inclusion~$\Univ_K\to\Univ_L$
is an embedding of differential $K$-algebras.
% as in Corollary~\ref{cor:Univ under d-field ext}. 

%\subsection*{The trace} \marginpar{this subsection not used later and commented out but checked} Recall from Section~\ref{sec:group rings} how we defined $\operatorname{tr}(f)\in K$ for $f\in K[G]$ and, what is more important, its intrinsic description. Applying this to $\Univ=K\big[\!\ex(\Lambda)\big]$ yields the trace $\operatorname{tr}_{\Univ}\colon \Univ\to K$ of $\Univ$. Note that $\operatorname{tr}(f)'=\operatorname{tr}(f')$ for $f\in\Univ$. 

%Let $L$ be a differential field extension of $K$ such that $C_L$ is algebraically closed and~$L^\dagger$ is divisible with $L^\dagger\cap K=K^\dagger$, and identify $\Univ_K$ with a differential $K$-subalgebra of $\Univ_L$ via an embedding $\Univ_K \to \Univ_L$ of differential $K$-algebras.  Then the trace of $\Univ_L$ extends the trace of $\Univ_K$, by the intrinsic description alluded to and Corollary~\ref{cor:Univ under d-field ext}.

\subsection*{Automorphisms of $\Univ$} 
These are easy to describe: the beginning of Section~\ref{sec:group rings} gives a group embedding
$$\chi\mapsto \sigma_{\chi}\colon \Hom(\Lambda,K^\times)\to \Aut\!\big(K[\ex(\Lambda)]|K\big)$$ into the group of $K$-algebra automorphisms of $K\big[\!\ex(\Lambda)\big]$, given by
$$ \sigma_{\chi}(f)\ := f_{\chi}\ =\ \sum_\lambda f_\lambda \chi(\lambda)\ex(\lambda)\qquad (\chi\in \Hom(\Lambda,K^\times),\ f\in K[\ex(\Lambda)]).$$
It is easy to check that if $\chi\in \Hom(\Lambda,C^\times)\subseteq \Hom(\Lambda,K^\times)$, then $\sigma_\chi\in \Aut_\der(\Univ|K)$, that is, $\sigma_{\chi}$ is a differential $K$-algebra automorphism of $\Univ$. Moreover:

\begin{lemma}\label{autolem} The map $\chi\mapsto \sigma_\chi \colon \Hom(\Lambda,C^\times)\to\Aut_\der(\Univ|K)$ is a group isomorphism. Its inverse assigns to
any $\sigma\in \Aut_\der(\Univ|K)$ the function
$\chi\colon \Lambda \to C^\times$ given by~$\chi(\lambda):=\sigma\big(\!\ex(\lambda)\big)\ex(- \lambda)$. In particular, $\Aut_\der(\Univ|K)$ is commutative.
\end{lemma}
\begin{proof} Let $\sigma\in \Aut_{\der}(\Univ|K)$ and let $\chi\colon \Lambda\to \Univ^\times$ be given by $\chi(\lambda):=\sigma\big(\!\ex(\lambda)\big)\ex(- \lambda)$. 
Then $\chi(\lambda)^\dagger=0$ for all $\lambda$. It follows easily that $\chi\in \Hom(\Lambda,C^\times)$ and $\sigma_{\chi}=\sigma$.
\end{proof} 

\noindent
The proof of the next result uses that the additive group $\Q$ embeds into~$C^\times$. 
%(Identifying the field of real algebraic numbers with a subfield of $C$,   we have the group embedding $q\mapsto 2^q\colon\Q\to C^\times$.)

\begin{cor}\label{cor:fixed field}
If $f\in\Univ$ and $\sigma(f)=f$ for all $\sigma\in\operatorname{Aut}_\der(\Univ|K)$, then $f\in K$.
\end{cor}
\begin{proof} Suppose $f\in U$ and $\sigma(f)=f$ for all $\sigma\in\operatorname{Aut}_\der(\Univ|K)$.
For $\chi\in\operatorname{Hom}(\Lambda,C^\times)$ we have $f_\chi=f$, that is, $f_\lambda\chi(\lambda)=f_\lambda$
for all $\lambda$, so  $\chi(\lambda)=1$ whenever $f_\lambda\neq 0$. Now use that for $\lambda\ne 0$ there exists
$\chi\in\operatorname{Hom}(\Lambda,C^\times)$ such that $\chi(\lambda)\ne 1$, so $f_{\lambda}=0$. 
\end{proof}

%\noindent
%We also note:

\begin{cor}\label{cor:extend autom}
Every automorphism of the differential field $K$ extends to an automorphism of the differential ring~$\Univ$.
\end{cor}
\begin{proof}
Lemma~\ref{lem:exponential map} yields a group morphism $\mu\colon K\to \Univ^\times$ such that $\mu(a)^\dagger=a$ for all $a\in K$.
Let $\sigma\in\operatorname{Aut}_\der(K)$. Then $\sigma$ extends to an endomorphism, denoted also by $\sigma$,
 of the ring $\Univ$, such that $\sigma\big(\!\ex(\lambda)\big)= \mu\big(\sigma(\lambda)\big)$ for all $\lambda$. 
Then  
$$\sigma\big(\!\ex(\lambda)'\big)\ =\ \sigma\big(\lambda\ex(\lambda)\big)\ =\ \sigma(\lambda)\mu\big(\sigma(\lambda)\big)\ =\ \mu\big(\sigma(\lambda)\big)'\ =\ \sigma\big(\!\ex(\lambda)\big)',$$
hence $\sigma$ is an endomorphism of the differential ring $\Univ$.  
By Lemma~\ref{lem:U simple}, $\sigma$ is injective, and by Lemma~\ref{lem:U minimal}, $\sigma$ is surjective.
\end{proof}

\subsection*{The real case} 
{\it In this subsection $K=H[\imag]$ where~$H$ is a real closed differential subfield of $K$ and $\imag^2=-1$.}\/
Set $S_C:=\big\{c\in C:\, |c|=1\big\}$, a subgroup of $C^\times$. Then by Lemmas~\ref{lem:commuting with f*} and~\ref{autolem}:

\begin{cor}\label{autolem, commuting with f*}
For $\sigma\in\Aut_\der(\Univ|K)$ we have the equivalence
$$\sigma(f^*)=\sigma(f)^*\text{ for all $f\in\Univ$}\quad\Longleftrightarrow\quad \sigma=\sigma_\chi\text{ for some $\chi\in\Hom(\Lambda,S_C)$.}$$
\end{cor}

\noindent
Corollaries~\ref{autolem, commuting with f*} and~\ref{cor:commuting with f*} together give:
 
\begin{cor}\label{cor:inner prod invariant}
Let  $\sigma\in\Aut_\der(\Univ|K)$ satisfy $\sigma(f^*)=\sigma(f)^*$   for all $f\in\Univ$. Then~$\big\langle\sigma(f),\sigma(g)\big\rangle=\langle f,g\rangle$ for all $f,g\in\Univ$, hence $\dabs{\sigma(f)}=\dabs{f}$ for all $f\in\Univ$.
\end{cor}

\noindent
Next we consider the subgroup 
$$S\ :=\ \{a+b\imag:\ a,b\in H,\ a^2+b^2=1\}$$ of $K^\times$, which is divisible,
hence so is the subgroup $S^\dagger$ of $K^\dagger$.
Lemma~\ref{lem:logder} yields~$K^\dagger = H^\dagger \oplus S^\dagger$
(internal direct sum of $\Q$-linear subspaces of $K$) and $S^\dagger\subseteq  H \imag$.
Thus we can (and do) take the complement $\Lambda$ of $K^\dagger$ in $K$ so that
$\Lambda= \Lambda_{\operatorname{r}}+\Lambda_{\operatorname{i}}\imag$ where $\Lambda_{\operatorname{r}}, \Lambda_{\operatorname{i}}$ are subspaces of the $\Q$-linear space $H$ with $\Lambda_{\operatorname{r}}$ a complement of $H^\dagger$ in $H$ and 
$\Lambda_{\operatorname{i}}\imag$
 a complement  of $S^\dagger$ in $H\imag$.
The automorphism $a+b\imag\mapsto \bar{a+b\imag}:= a-b\imag$~(${a,b\in H}$) of the differential field $K$ now satisfies in $\Univ=K[\ex(\Lambda)]$ the identity 
$$\ex(\bar{\lambda + \mu})\ =\ \ex(\bar{\lambda})\ex(\bar{\mu})\qquad(\lambda,\mu\in \Lambda),$$ so it extends to
an automorphism $f\mapsto\overline{f}$ of the ring $\Univ$ as follows: 
for $f\in \Univ$ with spectral decomposition~$(f_\lambda)$, set
$$\overline{f}\ :=\ \sum_\lambda \overline{f_\lambda}\ex(\overline{\lambda})\ =\ \sum_{\lambda}\overline{f_{\overline{\lambda}}}\ex(\lambda),  $$
so $\overline{\ex(\lambda)}=\ex(\overline{\lambda})$, and $\overline{f}$ has spectral decomposition  $(\overline{f_{\overline{\lambda}}})$.
We have $\overline{\overline{f}}=f$ for $f\in\Univ$,
and $f\mapsto\overline{f}$ lies in $\operatorname{Aut}_\der(\Univ|H)$. 
If~$H^\dagger=H$, then $\Lambda_r=\{0\}$ and hence~$\overline{f}=f^*$ for~$f\in\Univ$, where $f^*$ is as defined in Section~\ref{sec:group rings}.
For $f\in\Univ$ we set
$$\Re f\ :=\  \textstyle\frac{1}{2}(f+\overline{f}),\qquad  
\Im f\ :=\ \textstyle\frac{1}{2\imag}(f-\overline{f}).$$
(For $f\in K$ these agree with the usual real and imaginary parts of $f$ as an element of $H[\imag]$.)
Consider the differential $H$-subalgebra 
$$\Univ_{\operatorname{r}}\ :=\ \big\{f\in\Univ: \overline{f}=f\big\}$$ 
of~$\Univ$.
For $f\in\Univ$ with spectral decomposition $(f_\lambda)$ we have~$f\in \Univ_{\operatorname{r}}$ iff $f_{\overline{\lambda}}=\overline{f_\lambda}$ for all~$\lambda$;
in particular $\Univ_{\operatorname{r}}\cap K=H$.
For $f\in\Univ$ we have $f=(\Re f)+(\Im f)\imag$ with~$\Re f,\Im f\in \Univ_{\operatorname{r}}$, hence 
$$\Univ\ =\ \Univ_{\operatorname{r}}\oplus \Univ_{\operatorname{r}}\imag\quad\text{ (internal direct sum of $H$-linear subspaces).}$$
Let $D$ be a subfield of $H$ (not necessarily the constant field of $H$), so $D[\imag]$ is a  subfield  of $K$. Let 
$V$ be a $D[\imag]$-linear subspace of $\Univ$; then $V_{\operatorname{r}}:=V\cap \Univ_{\operatorname{r}}$ is a $D$-linear subspace of $V$. 
If~$\overline{V}=V$ (that is,~$V$~is closed under $f\mapsto \overline{f}$), then~$\Re f,\Im f\in V_{\operatorname{r}}$ for all $f\in V$, hence
 $V=V_{\operatorname{r}}\oplus V_{\operatorname{r}}\imag$ (internal direct sum of $D$-linear subspaces of $V$), so 
any basis of the $D$-linear space $V_{\operatorname{r}}$ is a basis of the $D[\imag]$-linear space $V$. 

\medskip
\noindent
Suppose now that $V=\bigoplus_\lambda V_\lambda$
(internal direct sum of subspaces of $V$) where $V_{\lambda}$ is for each $\lambda$
a $D[\imag]$-linear subspace of $K\ex(\lambda)$. Then~$\overline{V}=V$  iff
$V_{\overline{\lambda}}=\overline{V_\lambda}$ for all $\lambda$. Moreover:

\begin{lemma}\label{lem:real basis}
Assume $H=H^\dagger$, $V_0=\{0\}$,
and $\overline{V}=V$. Let $\mathcal V\subseteq \Univ^\times$ be a basis of the subspace~$\sum_{\Im\lambda>0} V_\lambda$ of $V$.
Then the maps $v\mapsto \Re v,\ v\mapsto \Im v\colon \mathcal{V} \to V_{\operatorname{r}}$ are injective, $\Re \mathcal{V}$ and $\Im \mathcal{V}$ are disjoint, and $\Re \mathcal{V}\cup \Im \mathcal{V}$ is  
%the elements $\Re v$, $\Im v$ \textup{(}$v\in \mathcal V$\textup{)} 
a basis of $V_{\operatorname{r}}$.
\end{lemma}
\begin{proof} Note that $\Lambda=\Lambda_{\operatorname{i}}\imag$. Let $\mu$ range over $\Lambda_{\operatorname{i}}^{>}$ and set $\cal{V}_{\mu}=\cal{V}\cap K^\times \ex(\mu\imag)$, a basis of the $D[\imag]$-linear space $V_{\mu\imag}$. Then $\cal{V}=\bigcup_{\mu}\cal{V}_{\mu}$, a disjoint union. For~$v\in \cal{V}_{\mu}$ we have~$v=a\ex(\mu\imag)$ with $a=a_v\in K^\times$, so 
$$\Re v\ =\ \textstyle\frac{a}{2}\ex(\mu\imag) + \textstyle\frac{\bar{a}}{2}\ex(-\mu\imag), \qquad \Im v\ =\ \textstyle\frac{a}{2i}\ex(\mu\imag) - \textstyle\frac{\bar{a}}{2\imag}\ex(-\mu\imag),$$
from which it is clear that the two maps $\cal{V} \to V_{\operatorname{r}}$ in the statement of the lemma are injective. It is also easy to check that $\Re \mathcal{V}$ and $\Im \mathcal{V}$ are disjoint.  

As $\mathcal{V}$ is a basis of the $D[\imag]$-linear space $\sum_{\mu} V_{\mu \imag}=\sum_{\Im\lambda>0} V_\lambda$, its set of conjugates~$\overline{\mathcal V}$ is a basis of the $D[\imag]$-linear space $\sum_{\mu} \overline{V_{\mu \imag}}=\sum_{\mu}V_{-\mu\imag}=\sum_{\Im\lambda<0} V_\lambda$, and so~$\mathcal{V}\cup \overline{\mathcal{V}}$ (a disjoint union) is a basis of $V$. Thus~$\Re \mathcal{V}\cup \Im \mathcal{V}$ is a basis of $V$ as well. As~$\Re \mathcal{V}\cup \Im \mathcal{V}$ is contained in $V_{\operatorname{r}}$, it is a basis of the $D$-linear space $V_{\operatorname{r}}$. 
%To see that the given elements span $V_{\operatorname{r}}=\Re V$ use that  
%$\mathcal W\cup \imag \mathcal W$ is a basis
%of  $V$ as   $C_H$-linear space, and to see that they are linearly independent use that their $C$-linear span contains
% $\mathcal W$.
\end{proof}

\noindent
If $H=H^\dagger$, then $V:=\sum_{\lambda\neq 0}K\ex(\lambda)$ gives $\overline{V}=V$, so Lemma~\ref{lem:real basis} gives then for~$D:= H$ the basis of the $H$-linear space $V_{\operatorname{r}}$ consisting of the elements 
$$\Re\!\big(\!\ex(\lambda)\big)\ =\ \textstyle\frac{1}{2}\big(\!\ex(\lambda)+\ex(\overline{\lambda})\big),\qquad 
\Im\!\big(\!\ex(\lambda)\big)\ =\ \textstyle\frac{1}{2\imag}\big(\!\ex(\lambda)-\ex(\overline{\lambda})\big)\qquad (\Im\lambda>0).$$

\begin{cor}\label{cor:U_r} Suppose $H=H^\dagger$. Set $\operatorname{c}(\lambda):=\Re\!\big(\!\ex(\lambda)\big)$ and $\operatorname{s}(\lambda):=\Im\!\big(\!\ex(\lambda)\big)$,
for $\Im \lambda>0$. Then for $V:=\sum_{\lambda\neq 0}K\ex(\lambda)$ we have $\Univ_{\operatorname{r}}=H+V_{\operatorname{r}}$, so
$$\Univ_{\operatorname{r}}\ =\ H \oplus \bigoplus_{\Im\lambda>0} \big( H\operatorname{c}(\lambda)\oplus H\operatorname{s}(\lambda) \big) 
\quad \text{\textup{(}internal direct sum of $H$-linear subspaces\textup{)},}$$
and thus $\Univ_{\operatorname{r}} = H\big[ \operatorname{c}(\Lambda_{\operatorname{i}}^>\imag)\cup \operatorname{s}(\Lambda_{\operatorname{i}}^>\imag)\big]$.
\end{cor}

%of $\Univ_{\operatorname{r}}$  are $H$-linearly independent.

%For $f=a+b\imag$ ($a,b\in H$) we have
%$$f\ex(\lambda) = (a+b\imag)\big(\!\operatorname{c}(\lambda)+\operatorname{s}(\lambda)\imag\big) =  \big(a \operatorname{c}(\lambda) - b \operatorname{s}(\lambda)\big) + \big(b \operatorname{c}(\lambda) + a \operatorname{s}(\lambda)\big)\imag,$$
%so
%$$\Re  \!\big(f\ex(\lambda)\big) = a \operatorname{c}(\lambda) - b \operatorname{s}(\lambda),\quad 
%  \Im  \!\big(f\ex(\lambda)\big) = b \operatorname{c}(\lambda) + a \operatorname{s}(\lambda).$$
  
%\begin{lemma}
%The elements $\operatorname{c}(\lambda)$, $\operatorname{s}(\lambda)$  of $\Univ_{\operatorname{r}}$ \textup{(}for varying $\lambda$ with $\Im\lambda>0$\textup{)} are $K$-linearly independent.
%\end{lemma}
%\begin{proof}
%It suffices to show the $K$-linear independence of $\operatorname{c}(\lambda)$, $\operatorname{s}(\lambda)$ for a  given $\lambda$ with $\Im\lambda>0$. Suppose $a,b\in K$ satisfy $a\operatorname{c}(\lambda)+ b \operatorname{s}(\lambda)=0$; then $\left(\frac{a+b\imag}{2} \right) \ex(\lambda) + \left(\frac{a-b\imag}{2} \right) \ex(\overline{\lambda}) = 0$ and hence $a=b=0$.
%\end{proof}

\section{The Spectrum of a Differential Operator}\label{sec:splitting}

\noindent
{\it In this section $K$ is a differential field, $a$, $b$ range over $K$, and $A$, $B$  over~$K[\der]$.}\/
This and the next two sections are mainly differential-algebraic in nature, and deal with splittings of linear differential operators.
In the present section we introduce the concept of {\it eigenvalue}\/ of $A$ and  the {\it spectrum}\/ of $A$ (the collection of its eigenvalues).
In Section~\ref{sec:self-adjoint} we give criteria for $A$ to have eigenvalue $0$, and in Section~\ref{sec:eigenvalues and splitting} we show how the eigenvalues of $A$ relate to the behavior of $A$ over the universal exponential extension of $K$.

\subsection*{Twisting}
Let $L$ be a differential field extension  of $K$ with $L^\dagger\supseteq K$.
Let $u\in L^\times$ be such that $u^\dagger=a\in K$. Then the twist \index{twist}\index{linear differential operator!twist}
$A_{\ltimes u} = u^{-1} A u$
of $A$ by $u$ has the same order as $A$ and coefficients in $K$ [ADH,~5.8.8], and only depends on $a$, not on $u$ or~$L$;
in fact, $\Ric(A_{\ltimes u})=\Ric(A)_{+a}$ [ADH, 5.8.5].
% and $(u,A)\mapsto A_{\ltimes u}$ is a left action of $\Univ_K^\times$ on~$K[\der]$.
%For $c\in C_L^\times$ we have $A_{\ltimes c}=A$.
Hence for  each $a$ we may define
$$A_a\ :=\ A_{\ltimes u}\ =\ u^{-1}Au\in K[\der]$$
where $u\in L^\times$ is arbitrary with $u^\dagger=a$. 
The map
$A\mapsto A_{\ltimes u}$ is an automorphism of the ring $K[\der]$ that is
the identity on $K$ (with inverse $B\mapsto B_{\ltimes u^{-1}}$);
so  $A\mapsto A_a$ is an automorphism of the ring $K[\der]$ that is the identity on $K$  (with inverse $B\mapsto B_{-a}$). 
Note that $\der_a=\der+a$, and that
$$(a,A)\mapsto A_a\ :\ K\times K[\der] \to K[\der]$$
is an action of the additive group of $K$ on the set $K[\der]$, in particular, $A_a=A$ for~$a=0$.  
For  $b\neq 0$ we have $(A_a)_{\ltimes b} = A_{a+b^\dagger}$.

\subsection*{Eigenvalues} 
{\it In the rest of this section $A\neq 0$ and $r:=\order(A)$.}\/
We call 
$$\mult_{a}(A)\ :=\ \dim_C \ker_K A_a\in \{0,\dots,r\}$$
the {\bf multiplicity} of $A$ at $a$. 
If $B\neq 0$, then  
$\mult_a(B) \leq \mult_a(AB)$, as well as
 \index{linear differential operator!multiplicity}\index{multiplicity!linear differential operator}\label{p:multa}
\begin{equation}\label{eq:mult(AB)}
\mult_a(AB)\ \leq\ \mult_a(A)+\mult_a(B),
\end{equation}
with equality if and only if $B_a(K)\supseteq  \ker_K A_a$; see [ADH, remarks before 5.1.12]. 
For $u\in K^\times$ we have an isomorphism 
$$y\mapsto yu\ \colon\ \ker_K A_{\ltimes u} \to \ker_K A$$ 
of $C$-linear spaces, hence 
$$\mult_{a}(A)\ =\  \mult_{b}(A)\qquad\text{whenever $a-b\in K^\dagger$.}$$
Thus we may  define the
{\bf multiplicity} of $A$ at the element $[a]:=a+K^\dagger$ of $K/K^\dagger$ 
as $\mult_{[a]}(A):=\mult_{a}(A)$. \label{p:multalpha}

\medskip
\noindent
{\it In the rest of this section $\alpha$ ranges over $K/K^\dagger$.}\/
We say that $\alpha$ is an {\bf eigenvalue} of~$A$ if~$\mult_{\alpha}(A)\geq 1$.  Thus for $B\ne 0$: if
$\alpha$ is an eigenvalue of~$B$ of multiplicity~$\mu$, then
$\alpha$ is an eigenvalue of $AB$ of multiplicity~$\ge \mu$;
if $\alpha$ is an eigenvalue of~$AB$, then it is an eigenvalue
of $A$ or of $B$; and if
$B_a(K)\supseteq\ker_K(A_a)$, then 
 $\alpha=[a]$ is an eigenvalue of $AB$ if and only if it is an eigenvalue
of~$A$ or of~$B$. 

\begin{exampleNumbered}\label{ex:ev order 1}
Suppose $A=\der-a$. 
Then for each element $u\neq 0$ in a differential field extension  of $K$ with $b:=u^\dagger\in K$ we have $A_b=A_{\ltimes u}=\der-(a-b)$,
so~$\mult_b(A)\geq 1$ iff $a-b\in K^\dagger$. Hence the only eigenvalue of $A$ is $[a]$.
\end{exampleNumbered}

\noindent
The {\bf spectrum} of $A$ is the set~$\Sigma(A)=\Sigma_K(A)$ of its eigenvalues. Thus $\Sigma(A)=\emptyset$ if $r=0$, and
for~$b\neq 0$ we have~$\mult_{a}(A)=\mult_{a}(bA)=\mult_{a}(A_{\ltimes b})$,
so $A$, $bA$, and~$Ab=bA_{\ltimes b}$  all have the same spectrum.\index{linear differential operator!eigenvalue}\index{eigenvalue!linear differential operator}\index{linear differential operator!spectrum}\index{spectrum!linear differential operator}\label{p:SigmaA}
By [ADH, 5.1.21] we have
\begin{equation}\label{eq:spec A}
\Sigma(A) = \big\{ \alpha : A \in K[\der](\der-a) \text{ for some $a$ with $[a]=\alpha$} \big\}.
\end{equation}
Hence for irreducible $A$:   $\ \Sigma(A)\neq \emptyset\ \Leftrightarrow\ r=1$. 
From \eqref{eq:mult(AB)} we obtain:

\begin{lemma}\label{lem:spectrum fact}
Suppose~$B\neq 0$ and set $s:=\order B$. Then $$\mult_\alpha(B)\ \leq\ \mult_\alpha(AB)\ \leq\ 
\mult_\alpha(A)+\mult_\alpha(B),$$
where the second inequality is an equality if  $K$ is $s$-linearly surjective. Hence
$$\Sigma(B)\ \subseteq\ \Sigma(AB)\ \subseteq\ \Sigma(A)\cup\Sigma(B).$$
If $K$ is $s$-linearly surjective, 
then $\Sigma(AB) = \Sigma(A)\cup\Sigma(B)$.
\end{lemma}

\begin{example}
For $n\geq 1$ we have  $\Sigma\big( (\der-a)^n \big)=\big\{[a]\big\}$. (By induction on $n$, using  Example~\ref{ex:ev order 1} and Lemma~\ref{lem:spectrum fact}.)
\end{example}

\noindent
It follows from Lemma~\ref{lem:spectrum fact} that $A$ has at most $r$ eigenvalues.  More precisely:
 
\begin{lemma}\label{lem:size of Sigma(A)}
We have $\sum_\alpha \mult_\alpha(A)\leq r$.
If $\sum_\alpha \mult_\alpha(A) = r$, then $A$ splits over~$K$; the converse holds if
$r=1$ or $K$  is $1$-linearly surjective. 
\end{lemma}
\begin{proof}
By induction on $r$. The case $r=0$ is obvious,
so suppose~$r>0$. We may also assume $\Sigma(A)\neq\emptyset$: otherwise  $\sum_\alpha \mult_\alpha(A)=0$ 
and $A$ does not split over~$K$. Now~\eqref{eq:spec A} gives~$a$,~$B$ with 
$A=B(\der-a)$. By Example~\ref{ex:ev order 1}  we have~$\Sigma(\der-a)=\big\{[a]\big\}$ and $\mult_a(\der-a)=1$.
By the inductive hypothesis applied to $B$ and the second inequality in Lemma~\ref{lem:spectrum fact} we thus get~$\sum_\alpha \mult_\alpha(A)\leq r$.

Suppose that $\sum_\alpha \mult_\alpha(A) = r$.
Then    $\sum_\alpha \mult_\alpha(B) = r-1$ by  
Lemma~\ref{lem:spectrum fact} and
the inductive hypothesis applied to $B$.
Therefore $B$ splits over $K$, again by the inductive hypothesis, and so does $A$.
Finally, if $K$ is $1$-linearly surjective and~$A$ splits over $K$, then we arrange that $B$ splits over $K$,
so $\sum_\alpha \mult_\alpha(B) = r-1$ by the inductive hypothesis,
hence $\sum_\alpha \mult_\alpha(A) = r$ by Lemma~\ref{lem:spectrum fact}.
\end{proof}

\noindent
Section~\ref{sec:eigenvalues and splitting} gives a more explicit proof of Lemma~\ref{lem:size of Sigma(A)},
under additional hypotheses on $K$.
Next, let $L$ be a differential field extension of $K$. Then 
$\mult_a(A)$ does not strictly decrease in passing from $K$ to $L$~[ADH, 4.1.13]. Hence the group morphism
$$a+K^\dagger\mapsto a+L^\dagger\colon K/K^\dagger\to L/L^\dagger$$ restricts to a map
$\Sigma_K(A)\to\Sigma_L(A)$; in particular, if $\Sigma_K(A)\neq\emptyset$, then $\Sigma_L(A)\neq\emptyset$.
If~$L^\dagger \cap K=K^\dagger$, then $\abs{\Sigma_K(A)}\leq\abs{\Sigma_L(A)}$, and
$\sum_\alpha \mult_\alpha(A)$ also 
does not strictly decrease if $K$ is replaced by $L$.

\begin{lemma}\label{lem:split evs}
Let $a_1,\dots, a_r\in K$ and 
$$A\ =\ (\der-a_r)\cdots(\der-a_1),\quad
\sum_{\alpha} \mult_{\alpha}(A)\ =\ r.$$
Then the spectrum of~$A$ is $\big\{[a_1],\dots,[a_r]\big\}$, and for all $\alpha$,
$$\mult_{\alpha}(A)\ =\ 
\big| \big\{ i\in\{1,\dots,r\}:\ \alpha=[a_i] \big\}  \big|.$$
\end{lemma}

\begin{proof}
Let  $i$ range over $\{1,\dots,r\}$.
By Lemma~\ref{lem:spectrum fact} and Example~\ref{ex:ev order 1},
$$\mult_\alpha(A)\ \leq\ \sum_i \mult_\alpha(\der-a_i)\ =\ 
\big| \big\{ i : \alpha=[a_i] \big\}  \big|$$
and hence
$$r\ =\ \sum_{\alpha} \mult_{\alpha}(A)\ \leq\ \sum_\alpha \big| \big\{ i : \alpha=[a_i] \big\}  \big|\ =\  r.$$
Thus for each $\alpha$ we have $\mult_\alpha(A) = 
\big| \big\{ i : \alpha=[a_i] \big\}  \big|$ as required.
\end{proof}

\noindent
Recall from [ADH, 5.1.8] that~$D^*\in K[\der]$ denotes the {\it adjoint}\/ of~$D\in K[\der]$, and that the map $D\mapsto D^*$ is an involution of the ring $K[\der]$ with  $a^*=a$ for all~$a$ and~$\der^*=-\der$.\index{linear differential operator!adjoint}\index{adjoint!linear differential operator}\label{p:A*}
If $A$ splits over $K$, then so does $A^*$. Furthermore, $(A_a)^*=(A^*)_{-a}$ for all $a$. By
Lemmas~\ref{lem:size of Sigma(A)} and~\ref{lem:split evs}:

\begin{cor}\label{cor:spectrum, self-adjoint, 1}
Suppose $K$ is $1$-linearly surjective and $\sum_\alpha \mult_\alpha(A)=r$.  Then $\mult_\alpha(A)=\mult_{-\alpha}(A^*)$
for all $\alpha$. In particular,
the map~$\alpha\mapsto-\alpha$ restricts to a bijection~$\Sigma(A)\to\Sigma(A^*)$.
\end{cor}

%\begin{cor}
%If $K$ is $1$-linearly surjective, $K^\dagger$ is $2$-divisible, and $A$ is self-adjoint and splits over~$K$, then $\mult_0(A)\equiv r\bmod 2$; in particular $\ker_K A\neq\{0\}$ if also $r$ is odd.
%\end{cor}

\noindent
Let $\phi\in K^\times$. Then $(A^\phi)_a=(A_{\phi a})^\phi$ and hence
$$\mult_a(A^\phi)\  =\  \mult_{\phi a}(A),$$
so the group isomorphism 
\begin{equation}\label{eq:Kphidagger}
[a]\mapsto [\phi a]\ \colon\  K^\phi/\phi^{-1}K^\dagger\to K/K^\dagger
\end{equation}
maps $\Sigma(A^\phi)$ onto $\Sigma(A)$.

%In the next lemma we assume that $H$ is a real closed differential subfield of $K$ such that $K=H[\imag]$ where $\imag^2=-1$. The complex conjugation automorphism of the field~$K$ over $H$ then extends to a ring automorphism $B\mapsto\overline{B}$ of $K[\der]$ with~$\overline{\der}=\der$; we have $\overline{B(y)}=\overline{B}(\overline{y})$ for each $y\in K$. Also note that we have the automorphism~$\alpha\mapsto\overline{\alpha}$ of the group $K/K^\dagger$ with $\overline{[a]}=[\overline{a}]$ for each $a$.  We have $\overline{A_a}=\overline{A}_{\overline{a}}$.  The map $y\mapsto\overline{y}$ is a sesquilinear automorphism of the $C$-linear space $K$ \cite[Chapter~XIII, \S{}7]{Lang} which restricts to a sesquilinear isomorphism $\ker_K A_a \to \ker_K \overline{A}_{\overline{a}}$ of $C$-linear subspaces of $K$. In particular, $\mult_a(A)=\mult_{\overline{a}}(\overline{A})$, and hence we obtain: 

%\begin{lemma}\label{lem:spectrum comp conj}\marginpar{new lemma}
%For each $\alpha$ we have  $\mult_\alpha(A)=\mult_{\overline{\alpha}}(\overline{A})$;  so we have a bijection~$\alpha\mapsto\overline{\alpha}\colon\Sigma(A)\to\Sigma(\overline{A})$.
%\end{lemma}

\medskip
\noindent
Note that $K[\der]/K\der]A$ as a $K$-linear space has dimension $r=\order A$. 
Recall from~[ADH, 5.1] that  $A$ and $B\neq 0$  are said to {\it have the same type}\/ if the (left) $K[\der]$-modules
$K[\der]/K[\der]A$ and $K[\der]/K[\der]B$ are isomorphic (and so $\order B=r$).  By [ADH, 5.1.19]:\index{linear differential operator!same type}\index{type}

\begin{lemma}\label{lem:same type}
The operators $A$ and $B\neq 0$ have the same type iff $\order B=r$ and there is $R\in K[\der]$ of order~$<r$ with $1\in K[\der]R+K[\der]A$ and $BR\in K[\der]A$.
\end{lemma}

\noindent
Hence if $A$, $B$ have the same type, then they also have the same type as elements of~$L[\der]$, for any differential field
extension $L$ of $K$.   Since $B\mapsto B_a$ is an automorphism of the ring~$K[\der]$, Lemma~\ref{lem:same type} and~[ADH, 5.1.20] yield:

\begin{lemma}\label{lem:same type and eigenvalues}
If $A$ and $B\neq 0$ have the same type, then so do $A_a$,~$B_a$, for all $a$, and thus
$A$, $B$ have the same eigenvalues, with same multiplicity. \end{lemma}

\noindent
By this lemma the spectrum of $A$ depends only on the type of $A$, that is, on the isomorphism type of the $K[\der]$-module
$K[\der]/K[\der]A$, suggesting one might try to associate a spectrum to each differential module over $K$. 
(Recall from [ADH, 5.5] that  a differential module over $K$ is a $K[\der]$-module of finite dimension as $K$-linear space.) Although
our focus is on differential operators, we carry this out
in the next subsection: it  motivates the terminology of ``eigenvalues'' originating in the case of the
differential field of Puiseux series over $\C$ treated in~\cite{vdPS}.
This point of view will be further developed in the projected second volume of [ADH].

\subsection*{The spectrum of a differential module\astr} 
{\it In this subsection $M$ is a differential module over~$K$ and $r=\dim_K M$.}\/  For each $B$ we  let $\ker_M B$ denote the kernel of the $C$-linear map~$y\mapsto By\colon M\to M$.
For~$M=K$ as horizontal differential module   over $K$ [ADH, 5.5.2],   this agrees with the $C$-linear subspace
$$\ker_K B\  =\ \ker B\ =\  \big\{y\in K:B(y)=0\big\}$$ of $K$. Also, for $B=\der$ we obtain the $C$-linear subspace~$\ker_M\der$ of horizontal elements of $M$.
We define the {\bf spectrum} of~$M$ to be the set\label{p:SigmaM}
$$\Sigma(M) := \big\{ \alpha:  \text{$\ker_M (\der-a) \neq \{0\}$ for some $a$  with $[a]=\alpha$} \big\}.$$
The elements of $\Sigma(M)$ are called {\bf eigenvalues} of $M$.\index{differential module!multiplicity}\index{multiplicity!differential module}\index{differential module!eigenvalue}\index{eigenvalue!differential module}\index{differential module!spectrum}\index{spectrum!differential module}
If $M=\{0\}$, then $\Sigma(M)=\emptyset$. Isomorphic differential modules over $K$ have clearly the same spectrum.

Let $\phi\in K^\times$ and $\derdelta=\phi^{-1}\der$. Then $K[\der]=K^\phi[\derdelta]$ as rings, hence $M$ is also a differential module over $K^\phi$ with $\phi^{-1}\der_M$ instead of $\der_M$ as its derivation; we denote it by $M^\phi$ and call it the {\bf compositional conjugate} of $M$ by $\phi$.\index{differential module!compositional conjugate}
Every cyclic vector of $M$  is also a cyclic vector of~$M^\phi$.
 The  group isomorphism~\eqref{eq:Kphidagger} maps $\Sigma(M^\phi)$ onto $\Sigma(M)$.

In the next lemma we assume $r\geq 1$ and denote the~$r\times r$~identity matrix over~$K$ by $I_r$.
Below $P$ is also an $r\times r$ matrix over $K$. The $C$-linear space of solutions to the matrix differential equation $y'=Py$ over $K$
is the set of all column vectors~$e\in K^r$ such that $e'=Pe$, and is denoted by~$\operatorname{sol}(P)$~[ADH, p.~276]. 
Recall that $a$ is said to be an eigenvalue
of $P$ over $K$ if $Pe=ae$ for some nonzero column vector $e\in K^r$. Recall also from [ADH, p.~277] that we associate
to $P$ the differential module~$M_P$ having the space $K^r$ of column vectors as its underlying $K$-linear space and satisfying $\der e= e'-Pe$ for all $e\in K^r$. Thus by [ADH, 5.4.8]: 
%~$\operatorname{ker}(aI_r-P)\neq\{0\}$.

\begin{lemma}
Let $M=M_{-P}$  be the differential module over $K$ associated to~$-P$. Then
$\ker_M (\der-a)  = \operatorname{sol}(aI_r-P)$, so $\dim_C \ker_M(\der-a)\leq r$ and
$$\Sigma(M) = \big\{ \alpha:  \text{$\operatorname{sol}(aI_r-P) \neq \{0\}$ for some $a$  with $[a]=\alpha$} \big\}.$$
\end{lemma}
%\begin{proof}
%The underlying $K$-linear space of $M$ is $K^r$, and $\der e = e'-Ne$ for $e\in K^r$, so this follows from~[ADH,~5.4.8].
%\end{proof}

\noindent
We define $\mult_a(M):=\dim_C \ker_M(\der-a)$; thus   $\mult_a(M)\in \{0,\dots,r\}$ by the previous lemma. \label{p:multaM}
For $b\neq 0$ we have a $C$-linear isomorphism 
$$y\mapsto by\colon\ker_M(\der-a)\to\ker_M(\der-a-b^\dagger).$$
This observation allows us to define 
the {\bf multiplicity} $\mult_\alpha(M)$ of $M$ at $\alpha$ as the quantity~$\mult_a(M)$ \label{p:multalphaM}
where $a$ with $[a]=\alpha$ is arbitrary. Clearly isomorphic differential modules over $K$ have the
same multiplicity at a given~$\alpha$.\index{differential module!multiplicity}\index{multiplicity!differential module}

\begin{lemma}\label{lem:der-a onto}
The following are equivalent:
\begin{enumerate}
\item[\textup{(i)}] $K$ is $r$-linearly surjective;
\item[\textup{(ii)}] for each differential module $N$ over $K$ with $\dim_K N = r$ and every $a$, 
the $C$-linear map $y\mapsto (\der-a) y\colon N\to N$ is surjective;
\item[\textup{(iii)}] for each differential module $N$ over $K$ with $\dim_K N \leq r$ we have $\der N=N$;
\item[\textup{(iv)}] for $n=1,\dots,r$, each matrix differential equation $y'=Fy+g$ with $F$ an $n\times n$ matrix over $K$
and $g\in K^n$
% is a column vector,  
has a solution in~$K$.
\end{enumerate}
\end{lemma} 
\begin{proof}
For (i)~$\Rightarrow$~(ii),
let $K$ be $r$-linearly surjective. The case $r=0$ being trivial, let $r\geq 1$, so $C\neq K$. Let $N$ be a differential module
over $K$ with $\dim_K N=r$. Towards proving that $y\mapsto (\der-a)y\colon N\to N$ is surjective, we can assume by [ADH, 5.5.3] that~$N = K[\der]/K[\der]A$ with $A$ of order $r$. Let $a$,~$B$ be given, and let $y$ range over $K$. It suffices to find~$R\in K[\der]$ and~$y$ such that $(\der-a)R = yA-B$, that is, $yA- B \in (\der -a)K[\der]$ for some $y$, equivalently,
$yA_a-B_a\in \der K[\der]$ for some~$y$. Taking adjoints this amounts to
finding   $y$ such that~$A_a^*y-B_a^*\in   K[\der]\der$, that is, $A_a^*(y)=B_a^*(1)$.
Such $y$ exists because $K$ is $r$-linearly surjective.

For (ii)~$\Rightarrow$~(iii), use that by [ADH, 5.5.2]
each differential module over $K$ of dimension $\leq r$ is a direct summand of a differential module over $K$ of dimension~$r$.
 For (iii)~$\Rightarrow$~(iv), note that
for  an $n\times n$ matrix $F$ over $K$ ($n\geq 1$), with associated differential module $M_F$ over $K$, and $g,y\in K^n$,
we have $y'=Fy+g$ iff
$\der y = g$ in~$M_F$ [ADH, p.~277]. 
For  (iv)~$\Rightarrow$~(i), use
[ADH, remarks before~5.4.3].
 \end{proof}

\noindent
The previous lemma refines [ADH, 5.4.2], and   leads to a more precise version of~[ADH, 5.4.3] with a similar proof:

\begin{cor}\label{corsol2}
Suppose $K$ is $mn$-linearly surjective and $L$ is a differential field extension of $K$ with $[L:K]=m$. Then $L$ is $n$-linearly surjective.
\end{cor}
\begin{proof}
Let $F$ be an $n\times n$ matrix over $L$, $n\ge 1$, and $g\in L^n$; by (iv)~$\Rightarrow$~(i) in Lemma~\ref{lem:der-a onto}
with $L$ in place of $K$
it is enough to show that the equation $y'+Fy=g$ has a solution in $L$.
For this, take  a basis~$e_1,\dots,e_m$ of the $K$-linear space $L$. As in the proof of [ADH, 5.4.3]
we obtain an $mn\times mn$ matrix $F^\diamond$ over $K$ and a column vector $g^\diamond\in K^{mn}$ such that
any solution of $z'=F^{\diamond}z+g^{\diamond}$ in $K$ yields a solution of~$y'=Fy+g$ in $L$.
Such a solution~$z$ exists by (i)~$\Rightarrow$~(iv)  in  Lemma~\ref{lem:der-a onto}.
\end{proof}

%\begin{lemma}\label{lem:der-a onto}
%Suppose $K$ is $r$-linearly surjective where $r=\dim_K M$.
%Then for every $a$, the $C$-linear map $y\mapsto (\der-a)y\colon M\to M$ is surjective.
%\end{lemma}
%\begin{proof} 
%The case $r=0$ being trivial, let $r\geq 1$, so $C\neq K$. Then by [ADH, 5.5.3] we can assume~$M = K[\der]/K[\der]A$ with $A$ of order $r$. Let $a$,~$B$ be given, and let $y$ range over $K$. It suffices to find~$R\in K[\der]$ and~$y$ such that $(\der-a)R = yA-B$, that is, $yA- B \in (\der -a)K[\der]$ for some $y$, equivalently,
%$yA_a-B_a\in \der K[\der]$ for some~$y$. Taking adjoints this amounts to
%finding   $y$ such that~$A_a^*y-B_a^*\in   K[\der]\der$, that is, $A_a^*(y)=B_a^*(1)$.
%Such $y$ exists because $K$ is $r$-linearly surjective.
%\end{proof}

\noindent
Let 
$0 \to M_1 \xrightarrow{\ \iota\ }  M \xrightarrow{\ \pi\ }  M_2 \to 0$
be a short exact sequence of differential modules over $K$, where for notational simplicity we assume that $M_1$ is a 
submodule of $M$ and $\iota$ is the natural inclusion. By restriction we obtain a sequence
\begin{equation}\label{eq:ker(der-a)}
0 \to \ker_{M_1}(\der-a) \xrightarrow{\ \iota_a\ } \ker_M (\der-a) \xrightarrow{\ \pi_a\ } \ker_{M_2}(\der-a)\to 0,
\end{equation}
of $C$-linear maps, not necessarily exact, but with $\operatorname{im}\iota_a=\ker \pi_a$.
Hence   
\begin{equation}\label{eq:mult subadd}
\mult_a(M)\ \leq\ \mult_a(M_1)+\mult_a(M_2).
\end{equation}
Therefore
$\Sigma(M)\subseteq\Sigma(M_1)\cup\Sigma(M_2)$.
If $(\der-a)M\cap M_1 = (\der-a)M_1$, then~$\pi_a$ is surjective, so the
sequence of $C$-linear maps \eqref{eq:ker(der-a)} is exact, hence the inequality \eqref{eq:mult subadd} is an equality.
Thus the next corollary whose hypothesis holds  if $M=M_1\oplus M_2$ (internal direct sum of submodules of $M$)
and $\iota$ and $\pi$ are the natural morphisms, and  also  
if $K$ is $r_1$-linearly surjective for $r_1:=\dim_K M_1$,
by Lemma~\ref{lem:der-a onto}:

\begin{cor}\label{cormultsum} 
Suppose  $(\der-a)M\cap M_1 = (\der-a)M_1$ for each $a$. Then 
$$\mult_\alpha(M) = \mult_\alpha(M_1)+\mult_\alpha(M_2)\quad\text{for each $\alpha$;}$$ 
in particular,
$\Sigma(M) = \Sigma(M_1)\cup\Sigma(M_2)$.
\end{cor}

\noindent
Let $E:=\End_C(M)$ be the $C$-algebra of endomorphisms of the $C$-linear space~$M$. We have a ring morphism $K[\der]\to E$ which assigns to $B\in K[\der]$   the element $y\mapsto By$ of $E$, and  we view $E$ accordingly as  $K[\der]$-module: $(Bf)(y):= B\cdot f(y)$ for $f\in E$, $y\in M$.  In the next corollary of Lemma~\ref{lem:der-a onto} we let $\der-a$ stand for the image of~$\der-a\in K[\der]$ under the above ring morphism $K[\der]\to E$. 

\begin{cor}
If $K$ is $r$-linearly surjective where $r=\dim_K M$, then
$$\Sigma(M) = \big\{ \alpha:  \text{$(\der-a) \notin E^\times$ for some $a$  with $[a]=\alpha$} \big\}.$$
\end{cor}

\begin{remark}
The description of $\Sigma(M)$ in the previous corollary is reminiscent of the
definition of  the {\it spectrum}\/ of an element $x$ of an arbitrary $K$-algebra $E$ with unit 
as the set of all $a$ such that $(x-a)\notin E^\times$, as given in \cite[\S{}1]{BouS}.
(If $C=K$, then~$K^\dagger=\{0\}$, and identifying $K/K^\dagger$ with $K$ in the natural way, $\Sigma(M)$ is
the spectrum of $\der\in E$ in this sense.)
\end{remark}

\noindent
Let now $N$ be a  differential module over $K$ and $s:=\dim_K N$.
From [ADH, p.~279] recall that   the   $K$-linear space $\Hom_K(M,N)$ of all $K$-linear maps $M \to N$ 
(of dimension $\dim_K \Hom_K(M,N) = rs$) is a differential module  over $K$
with
$$(\der\phi)(f)\  :=\  \der(\phi f)-\phi(\der f)\qquad\text{ for $\phi \in \Hom_K(M,N)$ and $f \in M$.}$$
Given a $K[\der]$-linear map $\theta\colon N\to P$ into a differential module $P$ over $K$, this yields a $K[\der]$-linear map
$\Hom_K(M,\theta)\colon \Hom_K(M,N)\to \Hom_K(M,P)$ which sends any $\phi$ in~$\Hom_K(M,N)$ to 
$\theta\circ \phi\in \Hom_K(M,P)$.  
The horizontal elements of~$\Hom_K(M,N)$ are   the $K[\der]$-module mor\-phisms~$M\to N$;
they are the elements of a finite-di\-men\-sion\-al $C$-linear subspace $\Hom_{K[\der]}(M,N)$ of $\Hom_K(M,N)$:

\begin{lemma}\label{dim of ann}
We have
$$\dim_C \Hom_{K[\der]}(M,N)\ \leq\ \dim_K \Hom_K(M,N),$$
with equality iff $\Hom_K(M,N)$ is horizontal. 
\end{lemma}
\begin{proof}
By [ADH, 5.4.8 and remarks before 5.5.2], the dimension of the $C$-linear space
of horizontal elements of $M$ is at most $\dim_K M$, with equality iff $M$ is horizontal. Now apply this with $\Hom_K(M,N)$ in place of $M$.
%We may assume $M,N\neq \{0\}$. Let $f=(f_1,\dots,f_m)^{\operatorname{t}}$ and $g=(g_1,\dots,g_n)^{\operatorname{t}}$ be ordered bases of the $K$-linear spaces $M$ and $N$, respectively, where $m,n\geq 1$. Let~$A$ be the matrix associated to~$M$ with respect to $f$ and $B$ be the matrix associated to~$N$ with respect to $g$, so  $\der f=-Af$ and $\der g=-Bg$. Then each $\phi\in\Hom_K(M,N)$ is represented by the $m\times n$ matrix $\Phi$ over $K$ such that $\Phi g=\big(\phi(f_1),\dots,\phi(f_m)\big)^{\operatorname{t}}$, and $\phi\in\Hom_{K[\der]}(M,N)$ iff $\der\Phi=B\Phi - \Phi A$. Interpreting $\Phi$ as a column vector in~$K^{mn}$, the latter equation can also be written in the form $\der\Phi = H\Phi$ for some $mn\times mn$ matrix $H$ over $K$ (not depending on $\Phi$). The claim now follows from~[ADH, 5.4.8].
\end{proof}

\noindent
Recall:  $M^*:=\Hom_K(M,K)$ is the {\it dual}\/ of $M$; see [ADH, 5.5].\index{differential module!dual}\label{p:M*}
By Lemma~\ref{dim of ann}, the dimension of the $C$-linear subspace $\Hom_{K[\der]}(M,K)=\ker_{M^*}\der$ of $M^*$ is at most~$\dim_K M$. For the differential module $M=K[\der]/K[\der]A$ we can say more: 

\begin{lemma}\label{lem:kerder}
Suppose $M=K[\der]/K[\der]A$ and $e:=1+K[\der]A\in M$. Then for all~$\phi\in \Hom_{K[\der]}(M,K)$ we have $\phi(e)\in \ker A$, and the map $$\phi\mapsto \phi(e)\colon  \Hom_{K[\der]}(M,K) \der\to \ker A$$ is an isomorphism of $C$-linear spaces.
\end{lemma}
\begin{proof}
The first claim follows from $A(\phi(e))=\phi(Ae)=0$, as $Ae=A+K[\der]A$ is the zero element  of $M$. This yields a $C$-linear map as displayed. To show that it is surjective, let $y\in\ker A$ be given. Then $B\mapsto B(y)\colon K[\der]\to K$ is  $K[\der]$-linear with~$K[\der]A$ contained in its kernel, and thus yields $\phi\in \Hom_{K[\der]}(M,K)$ with~$\phi(e)=y$.
Injectivity is clear since $M=K[\der]e$.
\end{proof}

\noindent
Given $a$, the map $\der-a\colon M\to M$ is a $\der$-compatible derivation on the $K$-linear space~$M$~[ADH, 5.5].
Let $M_a$ be the $K$-linear space $M$ equipped with this $\der$-compatible derivation. Thus
$M_a$ is a differential module over $K$ with 
$$\dim_K M_a=\dim_K M=r, \quad \ker_{M^*}(\der -a)=\Hom_{K[\der]}(M_a,K).$$ 
Moreover, if $e$ is a cyclic vector of $M$ with $Ae=0$, then
$e$ is a cyclic vector of $M_a$ with $A_{a}e=0$. Hence by the previous lemma:

\begin{cor}\label{cor:kerder}
Let $A$, $e$, $M$ be as in Lemma~\ref{lem:kerder}. 
Then for each~$\phi\in\ker_{M^*}({\der-a})$ we have $\phi(e)\in \ker A_{a}$, and the map $$\phi\mapsto \phi(e)\colon  \ker_{M^*}(\der-a)\to \ker A_{a}$$ is an isomorphism of $C$-linear spaces.
In particular, $\mult_{\alpha}(M^*) = \mult_{\alpha}(A)$, so~$\alpha$ is an eigenvalue of $M^*$ iff $\alpha$ is an eigenvalue of $A$.
\end{cor}

\noindent
Recall that every differential module $M$ has finite length, denoted by $\ell(M)$ [ADH, pp.~36--38, 251],
with $\ell(M)\leq\dim_K M=r$.
We say that {\bf $M$ splits} if $\ell(M)=r$.\index{differential module!splitting}\index{splitting!differential module}
By~[ADH, 5.1.25], $M=K[\der]/K[\der]A$ splits iff $A$ splits over $K$.
By additivity of $\ell(-)$ and $\dim_K(-)$ on short exact sequences (see~[ADH, 1.2]) we have:

\begin{lemma}\label{lem:split in short exact sequ}
Let $N$ be a differential submodule of $M$. Then $M$ splits iff both~$N$ and~$M/N$ split.
\end{lemma}

\noindent
Hence if   $N$ is a differential module over $K$, then $M\oplus N$ splits iff   $M$ and $N$ split.
Thus the least common left multiple of $A_1,\dots,A_m\in K[\der]^{\neq}$, $m\geq 1$, splits over $K$ iff~$A_1,\dots,A_m$ split over $K$: use that the differential module $$K[\der]/K[\der]\operatorname{lclm}(A_1,\dots, A_m)$$
over $K$ is isomorphic to the image of the natural (diagonal) $K[\der]$-linear map $$K[\der]\to \big(K[\der]/K[\der]A_1\big)\times \cdots \times (K[\der]/K[\der]A_m\big).$$  
A $K[\der]$-linear map $M\to N$ into a differential module
$N$ over $K$ induces a $K[\der]$-linear map $\phi^*\colon N^*\to M^*$ given by $\phi^*(f)=f\circ \phi$, and if $\phi$ is surjective, 
then $\phi^*$ is injective. This gives a  contravariant functor $(-)^*$ from the category of differential modules over $K$ to itself; 
the morphisms of this category are the $K[\der]$-linear maps between differential modules over $K$.  Using $\dim_K M=\dim_K M^*<\infty$ it follows easily from these facts that if $\phi\colon M\to N$ is an injective $K[\der]$-linear map into a differential module $N$, then $\phi^*\colon N^*\to M^*$ is surjective.

\begin{lemma}\label{lem:dual splits} $\ell(M)=\ell(M^*)$, so
if $M$ splits, then $M^*$ splits as well.
\end{lemma}
\begin{proof}
Induction on $\ell(M)$ using the canonical $K[\der]$-linear isomorphism
$M\cong M^{**}$ and what was said about  the functor $(-)^*$  shows $\ell(M)=\ell(M^*)$.
%Since $M\cong M^{**}$ we then obtain $\ell(M^*)\leq\ell(M)$.
\end{proof}

\noindent
Let $L$ be a differential
field extension   of $K$.
Recall from [ADH, 5.9.2]  that the base change $L\otimes_K M$ of~$M$ to $L$  is  a differential module over $L$ with $\dim_L L\otimes_K M=\dim_K M$. A $K[\der]$-linear map $M\to N$ into a differential module
$N$ over $K$ induces an $L[\der]$-linear map 
$$L\otimes_K \phi\ \colon\ L\otimes_K M \to L\otimes_K N, \qquad \lambda\otimes x \mapsto \lambda\otimes \phi(x),$$
and this yields a covariant functor $L\otimes_K-$ from the category of differential modules over $K$ to the category of
differential modules over $L$. Using the above invariance of dimension one checks easily that this functor transforms short exact sequences in the first category to short exact sequences in the second category.  Hence, by an easy induction on 
$\ell(M)$ we have $\ell(M)\le \ell(L\otimes_K M)$.

We say that~{\bf~$M$~splits over $L$} if $L\otimes_K M$ splits.  The $K[\der]$-linear isomorphism
$$x\mapsto 1\otimes x\ \colon\ M\to K\otimes_K M$$ shows that $M$ splits iff $M$ splits over $K$, and then  $M$ splits over each differential field extension  of $K$. Let $E$ be a differential field extension of $L$. Then we have
an $E[\der]$-linear isomorphism
$$E\otimes_K M\to E\otimes_L (L\otimes_K M),\quad  e\otimes x \mapsto e\otimes(1\otimes x),$$
so if $M$ splits over~$L$, then $M$ also splits over $E$.
If $N$ is a differential submodule of $M$, then $M$ splits over~$L$ iff both $N$ and $M/N$ split over $L$.  

\begin{lemma}
If $M$ splits over $L$, then so does $M^*$.
\end{lemma}
\begin{proof} For $\phi\in M^*$ we have the $L$-linear map 
$$\operatorname{id}_L \otimes\phi\ \colon\ L\otimes M\to L\otimes_K K,\quad s\otimes y\mapsto  s\otimes \phi(y) \quad (\lambda, s\in L, \phi\in M^*, y\in M).$$
We also have the $L[\der]$-linear isomorphism $i_L\colon L\otimes_K K\to L$ given by $i_L(s\otimes 1)=s$ for $s\in L$. It is straightforward to check that this yields an $L$-linear isomorphism
$$L\otimes_K M^* \to (L\otimes_K M)^*,\qquad 1\otimes \phi\mapsto  i_L\circ (\operatorname{id}_L\otimes \phi) \quad(\phi\in M^*),$$
and that this map is even $L[\der]$-linear. Now use Lemma~\ref{lem:dual splits}.
\end{proof}

\noindent
Call $M$ {\bf cyclic}\/ if it has a cyclic vector, equivalently, for some $A$ we have\index{differential module!cyclic}
$$M\ \cong\ K[\der]/K[\der]A, \text{  as $K[\der]$-modules.}$$

\begin{cor}\label{Sigmaatmostr} We have $\sum_\alpha \mult_\alpha(M)\leq r$, hence $\abs{\Sigma(M)}\leq r$.
If moreover~$\sum_\alpha \mult_\alpha(M) = r$, then $M$ splits.
Conversely, if $K$ is $1$-linearly surjective and~$M$ splits, then $\sum_\alpha \mult_\alpha(M) = r$.
\end{cor}
\begin{proof} We prove this for $M^*$ instead of $M$. (Then by various results above and the natural $K[\der]$-linear isomorphism $M\cong M^{**}$ it also follows for $M$.) 
If $M$ is cyclic,  then $\sum_\alpha \mult_\alpha(M^*)\leq r$ by Lemma~\ref{lem:size of Sigma(A)} and Corollary~\ref{cor:kerder}, and the rest follows using also Lemma~\ref{lem:dual splits} and remarks following Corollary~\ref{cor:kerder}.  Thus we are done if $C\ne K$, by [ADH, 5.5.3].

If $C=K$, then a differential module over $K$ is just a finite-dimensional vector space $M$ over $K$ equipped with a $K$-linear map $\der\colon M \to M$, so we can use the well-known internal direct sum decomposition
$\sum_{a} \ker_M(\der-a)=\bigoplus_a\ker_M(\der -a)$.
\end{proof}

\noindent
Likewise we obtain from Corollary~\ref{cor:spectrum, self-adjoint, 1}  and [ADH, 5.5.8]: 

\begin{cor}\label{cor:spectrum, self-adjoint, 2}
Suppose $K$ is $1$-linearly surjective and  $\sum_\alpha \mult_\alpha(M)=r$.
Then the map $\alpha\mapsto-\alpha$ restricts to a
bijection $\Sigma(M)\to\Sigma(M^*)$ with $\mult_\alpha(M)= \mult_{-\alpha}(M^*)$ for each~$\alpha$.
\end{cor}

\noindent We now aim for a variant of Corollary~\ref{cor:spectrum, self-adjoint, 2}:
Corollary~\ref{cor:spectrum, self-adjoint, 3} below.

\begin{lemma}\label{lem:dim ker A^*}
Suppose $\dim_C \ker A=r$. Then $\dim_C \ker A^*=r$.
\end{lemma}
\begin{proof}
The case $r=0$ being trivial, suppose~${r\geq 1}$ and set $M:=K[\der]/K[\der]A$, so~$M\neq\{0\}$. 
By Lemma~\ref{lem:kerder}, $M^*$ is horizontal,
hence   $M$ is also horizontal by [ADH, remark after~5.5.5], and therefore~$\dim_C\ker A^*=\dim_C \ker_{M} \der = r$ by  Lem\-ma~\ref{lem:kerder} and [ADH, 5.5.8].
\end{proof}

\begin{cor}\label{corrd}
Let $d:=\dim_C \ker A$ and suppose $K$ is $(r-d)$-linearly surjective. Then $\dim_C \ker A^*=d$.
\end{cor}
\begin{proof}
It suffices to show $\dim_C\ker A^*\geq d$, since then the reverse inequality follows by interchanging the role of $A$ and $A^*$.
Let $y_1,\dots,y_d$ be a basis of the $C$-linear space $\ker A$. Then
$$L(Y) := \wr(Y,y_1,\dots,y_d)\in K\{Y\}$$
is  homogeneous of degree $1$ and order $d$ with zero set $Z(L)=\ker A$ [ADH, 4.1.13]. So with $B\in K[\der]$ the linear part of $L$
we have $A = D B$  where $D\in K[\der]$ has or\-der~$r-d$, by [ADH, 5.1.15(i)], so $A^* = B^*D^*$ and $D^*(K)=K$.  
Hence
$$\dim_C \ker A^*\  =\  \dim_C \ker B^*+\dim_C\ker D^*\  \geq\ \dim_C \ker B^*\ =\ \dim_C\ker B\ =\ d$$
where we used [ADH, remark before 5.1.12] for the first equality and the previous lemma (applied to $B$ in place of $A$) for the second equality.
\end{proof}

\noindent
Suppose now that $r\ge 1$ and  $K$ is $(r-1)$-linearly surjective. Then Corollary~\ref{corrd} and $A^{**}=A$ give
$\dim_C\ker A=\dim_C\ker A^*$ (even when $\dim_C\ker A=0$). Hence for all $a$ we have
$\dim_C \ker A_a = \dim_C \ker (A^*)_{-a}$.  This leads to:
%In a similar way that Corollary~\ref{cor:spectrum, self-adjoint, 1}
%gave rise to Corollary~\ref{cor:spectrum, self-adjoint, 2}, this implies:

\begin{cor}\label{cor:spectrum, self-adjoint, 3}
If $r\ge 1$ and $K$ is $(r-1)$-linearly surjective, then we have a bijection $\alpha\mapsto-\alpha\colon\Sigma(M)\to\Sigma(M^*)$, and $\mult_\alpha(M)\ =\ \mult_{-\alpha}(M^*)$ for all $\alpha$.
\end{cor}

\subsection*{Complex conjugation\astr}
{\it In this subsection $K=H[\imag]$ where $H$ is a  differential subfield of $K$, $\imag^2=-1$, and $\imag\notin H$.}\/
Then $C=C_H[\imag]$. Recall that $A\in K[\der]^{\ne}$ has order $r$. The
complex conjugation automorphism $z=g+h\imag\mapsto \bar{z}:= g-h\imag$ ($g,h\in H$) of the differential field~$K$ induces an automorphism $\alpha\mapsto\bar\alpha$ of the group $K/K^\dagger$ with~$\bar\alpha=[\bar{a}]$ for $\alpha=[a]$, $a\in K$.
The   automorphism~$z\mapsto \bar{z}$ of~$K$ extends uniquely to an automorphism
$D\mapsto \bar{D}$ of the ring $K[\der]$ with~$\bar{\der}=\der$. 
If $A$ and $B\neq 0$ have the same type, then so do $\bar{A}$ and $\bar{B}$. (Lemma~\ref{lem:same type}.)
 %(For $A=\der$ we get $\bar{\der_a}=\der+\bar{a}=\der_{\bar a}$.)
Now~$\overline{A(f)}=\overline{A}(\overline{f})$ for~$f\in K$,
so~$\dim_C \ker_K A=\dim_C \ker_K\overline{A}$. Moreover, $\bar{A_a}=\bar{A}_{\bar a}$,
hence  
$\mult_\alpha(A)=\mult_{\bar{\alpha}}(\bar A)$ for all~$\alpha$; so
$\alpha$ is an eigenvalue of $A$ iff $\bar{\alpha}$ is an eigenvalue of $\bar A$.
Note that~$\bar{A^*}=\bar{A}{}^*$. We call $\bar{A^*}$ the {\bf conjugate adjoint} of $A$.\index{conjugate!adjoint}\index{linear differential operator!conjugate adjoint}\index{adjoint!conjugate} Corollary~\ref{cor:spectrum, self-adjoint, 1} yields:

\begin{cor}
If $K$ is $1$-linearly surjective and  $\sum_\alpha \mult_\alpha(A)=r$, then 
we have a bijection
$\alpha\mapsto-\bar{\alpha}\colon\Sigma(A)\to\Sigma(\bar{A^*})$, with
$\mult_\alpha(A)=\mult_{-\bar{\alpha}}(\bar{A^*})$ for all $\alpha$. 
\end{cor}

\noindent
Next, let $M$ be a (left) $K[\der]$-module. Then we define $\overline{M}$ as the $K[\der]$-module arising from $M$  by replacing its scalar multiplication $(A,f)\mapsto Af\colon K[\der]\times M\to M$ with
$$(A,f)\mapsto \overline{A}f\ \colon\  K[\der]\times \overline{M}\to\overline{M}.$$
We call $\bar M$ the {\bf complex conjugate} of $M$.\index{conjugate!complex}\index{differential module!complex conjugate} Note that $\bar{\bar M}=M$. 
If $\phi\colon M\to N$ is a morphism of $K[\der]$-modules, then $\phi$ is also a morphism
of $K[\der]$-modules~$\bar{M}\to\bar{N}$, which we denote by $\overline{\phi}$.
Hence we have a covariant functor $\bar{(\,\cdot\,)}$ from the  category of
$K[\der]$-modules to itself.
We have $\dim_K M=\dim_K\overline{M}$, hence if $M$ is a differential module over~$K$, then so is~$\overline{M}$.
Thus  $\bar{(\,\cdot\,)}$ restricts to a functor from the category of differential modules over~$K$ to itself.
If $P$ is an $r\times r$ matrix over $K$ and $M=M_P$ is the differential module associated to $P$ [ADH, p.~277],
then the differential modules~$\overline{M}$ and $M_{\overline{P}}$ over $K$ both have underlying additive group $K^r$,
and the map~$e\mapsto \overline{e}\colon \bar{M}\to M_{\bar P}$ is an isomorphism of differential modules over $K$.

\begin{exampleNumbered}\label{ex:M compl conj}
For $M=K[\der]$ we have an isomorphism $B\mapsto\bar{B}\colon M\to \bar{M}$ of $K[\der]$-modules.
For $N=K[\der]/K[\der]A$ we have an isomorphism  
$$B+K[\der]A\mapsto\bar{B}+K[\der]\bar A\ \colon\ \bar N\to K[\der]/K[\der]\bar A$$
of differential modules over $K$. 
\end{exampleNumbered}

\noindent
{\it Below $M$ is a differential module over~$K$ and $r=\dim_K M$.}\/
Then for each $B$ we have $\ker_M B=\ker_{\bar M}\bar B$. Hence $\mult_\alpha(M)=\mult_{\bar{\alpha}}(\bar M)$ for all $\alpha$, so $\alpha$ is an eigenvalue
of $M$ iff $\bar\alpha$ is an eigenvalue of~$\bar M$. 

Next, let $N$ be a differential module over $K$. A map $\phi\colon M\to N$ is $K$-linear if $\phi\colon\bar{M}\to\bar{N}$ is $K$-linear, so
$\operatorname{Hom}_K(M,N)$ and $\operatorname{Hom}_K(\bar{M},\bar{N})$ have the same underlying additive group.
It is easy to check that for the differential module~$P:=\operatorname{Hom}_K(M,N)$ we have $\bar{P}=\operatorname{Hom}_K(\bar{M},\bar{N})$. 
Thus $\bar{M^*}=\operatorname{Hom}_K(\bar{M}, \bar{K})$. In view of the isomorphism $z\mapsto \bar{z}\colon K\to \bar{K}$ of differential modules over $K$ this yields an isomorphism $\bar{M^*}\cong\bar{M}{}^*$ of differential modules over $K$. We call $\bar{M^*}$ the {\bf conjugate dual}\/ of $M$.\index{differential module!conjugate dual}\index{conjugate!dual} From Corollaries~\ref{cor:spectrum, self-adjoint, 2} and~\ref{cor:spectrum, self-adjoint, 3} we obtain:

\begin{cor}\label{cor:spectrum, conj self-adjoint}
Suppose $K$ is $1$-linearly surjective and  $\sum_\alpha \mult_\alpha(M)=r$,
or~$r\ge 1$ and $K$ is $(r-1)$-linearly surjective.
Then the map~$\alpha\mapsto-\bar\alpha$ restricts to a
bijection~$\Sigma(M)\to\Sigma(\bar{M^*})$ with $\mult_\alpha(M)= \mult_{-\bar{\alpha}}(\bar{M^*})$ for all~$\alpha$.
\end{cor}

\noindent
In the remainder of this section we discuss eigenvalues of differential modules over~$K$ in the presence of a valuation on $K$.
This is only used for the proof of Lemma~\ref{lem:mod osc} in Section~\ref{sec:lin diff applications}. 
In preparation for this, we first study lattices over valued fields.

\subsection*{Lattices\astr}
{\it In this subsection $F$ is a valued  field with valuation ring $R$.}\/
Let $L$ be an $R$-module, with its torsion submodule 
$$L_{\operatorname{tor}}=\big\{y\in L:\text{$ry=0$ for some $r\in R^{\neq}$}\big\}.$$
Call $L$ {\it torsion-free}\/ if $L_{\operatorname{tor}}=\{0\}$, and a
{\it torsion module}\/ if $L_{\operatorname{tor}}=L$. 
For the following basic fact, cf.~\cite[VI, \S{}4, Lemme~1]{BourbakiCA}:

\begin{lemma}\label{lem:Bourb}
Every finitely generated torsion-free $R$-module is free.
\end{lemma}
\begin{proof}
Let  $L$ be a finitely generated torsion-free $R$-module. Let $x_1,\dots,x_m\in L$ be distinct
such that~$\{x_1,\dots,x_m\}$ is a minimal set of generators of $L$ [ADH, p.~44].
Towards a contradiction, suppose $r_1x_1+\cdots+r_mx_m=0$ with $r_1,\dots,r_m$ in $R$ not all zero.
By reordering, arrange ${r_j\in r_1 R}$ for $j=2,\dots,m$.  
Torsion-freeness of~$L$ yields~${x_1+s_2x_2+\cdots+s_mx_m=0}$ where~$s_j:=r_j/r_1$ for~$j=2,\dots,m$. Hence~$\{x_2,\dots,x_m\}$ is also a set of generators of~$L$,
contradicting the minimality of~$\{x_1,\dots,x_m\}$. 
Thus~$x_1,\dots,x_m$ are $R$-linearly independent.
\end{proof}

\noindent
Let now $M$ be a finite-dimensional $F$-linear space and $m:=\dim_F M$.

\begin{lemma}\label{lem:char lattices}
Let $L$ be a finitely generated $R$-submodule of $M$. Then $L$ is free of rank~$\leq m$, and
the following are equivalent:
\begin{enumerate}
\item[\textup{(i)}] $L$ has rank $m$;
\item[\textup{(ii)}] $L$ has a basis which is also a basis of the $F$-linear space $M$;
\item[\textup{(iii)}] $L$ contains a basis of $M$;
\item[\textup{(iv)}] $L$ contains a generating set of $M$;
\item[\textup{(v)}] the $R$-module $M/L$ is a torsion module.
\end{enumerate}
\end{lemma}
\begin{proof}
Freeness of $L$ follows from Lemma~\ref{lem:Bourb}.  Every set of $R$-linearly independent elements of $L$
is $F$-linearly independent, so $\operatorname{rank}(L)\leq m$.
Let~$y_1,\dots,y_n$ be a basis of $L$.  
Assuming~$n<m$   yields $z\in M$ such that 
$y_1,\dots,y_n,z$ are $F$-linearly independent, so $M/L$ is not a torsion module. This shows (v)~$\Rightarrow$~(i), and~(i)~$\Rightarrow$~(ii)~$\Rightarrow$~(iii)~$\Rightarrow$~(iv)~$\Rightarrow$~(v)   are   clear.
\end{proof}

\noindent
A finitely generated $R$-submodule $L$ of $M$ is called an {\bf $R$-lattice} in $M$, or just a {\bf lattice} in $M$ if $R$ is understood from the context, if it satisfies one of the equivalent conditions~(i)--(v) in Lemma~\ref{lem:char lattices}.\index{lattice} 
If $L$ is a lattice in $M$, then every basis of $L$ is a basis  of the $F$-linear space $M$, and every lattice of $M$ is of the form~$\sigma(L)$ for some   automorphism $\sigma$ of $M$. Next some easy consequences of Lemma~\ref{lem:char lattices}:
%From \cite[XVI, Proposition~4.1]{Lang} we obtain:

%\noindent
%In the rest of this subsection we collect further constructions with lattices. First,
%the following is \cite[Lemma~2.2.7]{HHM} (and has a more indirect proof there):

\begin{cor}\label{cor:proj lattice}
If $\pi\colon M\to M'$ is a surjective morphism of $F$-linear spaces and~$L$ a lattice in $M$, then $L':=\pi(L)$ is a lattice in $M'$.
\end{cor}
%\begin{proof}
%The $R$-submodule $L'$ of the torsion-free $R$-module  $M'$ is torsion-free. Moreover, $L'$ is  finitely generated
%since $L$ is. Hence $L'$ is free by Lemma~\ref{lem:Bourb}. The $R$-module $M/L$ is  torsion, hence so is the $R$-module $M'/L'$. Thus $L'$ is a lattice in~$M'$, by Lemma~\ref{lem:char lattices}. 
%\end{proof}

\begin{cor}\label{cor:intersect lattice}
Let $N$ be an $F$-linear subspace of $M$ and $L$ a lattice in $M$. Then $L\cap N$ is a lattice in $N$.
\end{cor}
\begin{proof} Take a basis $x_1,\dots, x_n$ of $N$. Lemma~\ref{lem:char lattices}(v) gives $r_1,\dots, r_n\in R^{\ne}$ with 
$r_1x_1,\dots, r_nx_n\in L$. Then $r_1x_1,\dots, r_nx_n$ is a basis of $N$ contained in $L\cap N$. Now
apply Lemma~\ref{lem:char lattices}(iii) to $N$, $L\cap N$ in place of $M$, $L$.
%Let $\pi\colon M^*\to N^*$ be the surjective $F$-linear map $\phi\mapsto\phi|_N$; by Corollary~\ref{cor:proj lattice}, $\pi(L^*)$ is a lattice in $N^*$,
%hence $L':=\pi(L^*)^*$ is a lattice in $N^{**}$. We have the isomorphism 
%$$\Phi\colon N\to N^{**},\qquad \Phi(y)(\psi)=\psi(y)\text{ for $\psi\in N^*$, $y\in N$,}$$
%and it suffices to show that $\Phi(L\cap N)=L'$.
%Let $y\in L\cap N$ and  $\psi\in\pi(L^*)$, say~$\psi=\phi|_N$ where $\phi\in L^*$;
%then $\Phi(y)(\phi)=\psi(y)=\phi(y)\in R$,  showing
%$\Phi(L\cap N)\subseteq L'$.  To prove the reverse inclusion, let~$y\in N$ be such that $\Phi(y)\in L'$.
%Choose a basis~$y_1,\dots,y_m$ of $L$, and for $i\in\{1,\dots,m\}$ let~$\beta_i\in M^*$ be given by $\beta_i(y)=f_i$
%for~$f_1,\dots,f_m\in F$ and~$y=f_1y_1+\cdots+f_my_m$. Then~$\beta_i\in L^*$ and hence $\beta_i|_N\in\pi(L^*)$, thus $\beta_i(y)=\Phi(y)(\beta_i|_N)\in R$ and
%so $y\in L\cap N$.
\end{proof}

\begin{cor}\label{cor:base change lattice}
If $L$ is a lattice in $M$ and $E$ a valued field extension of~$F$ with valuation ring $S$, then
the $E$-linear space $M_{E}:=E\otimes_F M$ has dimension~$m$, and
the $S$-submodule $L_E$  of $M_{E}$ generated by the image of $L$ under the $F$-linear embedding~$y\mapsto 1\otimes y\colon M\to M_{E}$ is an $S$-lattice  in $M_{E}$.
\end{cor}

\noindent
For $i=1,2$ let $M_i$ be a  $F$-linear space with $m_i:=\dim_K M_i<\infty$ and $L_i$ be a lattice in~$M_i$. Then~$L_1\oplus L_2$ is a lattice in $M_1\oplus M_2$, and
the $R$-submodule of $M_1 \otimes_{F} M_2$ generated by the elements~$y_1\otimes y_2$ ($y_1\in L_1$, $y_2\in L_2$) is a lattice in $M_1 \otimes_{F} M_2$.
%(The first statement is clear, and for the second statement use \cite[XVI, Corollary~2.4]{Lang}.)
The $F$-linear space $\operatorname{Hom}(M_1,M_2)=\operatorname{Hom}_F(M_1,M_2)$ of $F$-linear maps $M_1\to M_2$ has dimension $m_1m_2$, and 
the $R$-module  $\operatorname{Hom}(L_1,L_2)$ of $R$-linear maps $L_1\to L_2$ is free of rank $m_1m_2$. 
Each $R$-linear map $\phi\colon L_1\to L_2$ extends uniquely to an $F$-linear map $\hat\phi\colon M_1\to M_2$, and~$\phi\mapsto\hat\phi$ is an embedding of   $\operatorname{Hom}(L_1,L_2)$ into $\operatorname{Hom}(M_1,M_2)$ viewed as
$R$-module. We identify~$\operatorname{Hom}(L_1,L_2)$ via this embedding with its image in $\operatorname{Hom}(M_1,M_2)$;
then $\operatorname{Hom}(L_1,L_2)$ is a lattice in $\operatorname{Hom}(M_1,M_2)$.
In particular, if~$L$ is a lattice in $M$, then $L^*=\operatorname{Hom}(L,R)$ is a lattice in $M^*=\operatorname{Hom}(M,F)$.

\subsection*{Lattices in differential modules\astr}
{\it In the rest of this section  $K$ is equipped with a valuation ring $\mathcal O$ making $K$ a valued differential field with small derivation. We also let
$M$ be a differential module over~$K$.}\/ Thus~$M$ is a $K[\der]$-module which is finite-dimensional as $K$-linear space.
We have the subring~$\mathcal O[\der]$ of $K[\der]$. 
A {\bf lattice} in~$M$ is an $\mathcal O[\der]$-submodule of $M$ that is also an
$\mathcal O$-lattice in the $K$-linear space~$M$.\index{lattice!differential module}\index{differential module!lattice} An $\mathcal{O}$-lattice~$L$ in  the $K$-linear space~$M$  is a lattice in the differential module~$M$ iff~$\der L\subseteq L$, iff
there is a generating set~$S$ of the $\mathcal{O}$-module $L$ with~$\der S\subseteq L$.
If $a\neq 0$ and $a^\dagger\preceq 1$ and $L$ is a lattice in $M$, then $a L$ is also a lattice in $M$.

{\samepage
\begin{examplesNumbered}\mbox{}\label{ex:lattices}

\begin{enumerate}
\item Suppose $M=M_N$ where $N$ is an $n\times n$ matrix over $\mathcal O$ \textup{(}$n\geq 1$\textup{)}. The underlying $K$-linear space of $M$ is $K^n$, and for each $e\in M$ we have~$\der e=e'-Ne$ [ADH, p.~277], so   $L:=\mathcal O^n$ is a lattice in   $M$.
\item 
Suppose $M\cong K[\der]/K[\der]A$ where $A\in\mathcal O[\der]$ is monic.  Let $e$ be a cyclic vector of $M$ with $Ae=0$.
Then the $K$-linear space $M$ has  basis~$e,\der e,\dots,\der^{r-1}e$ and~$L:=\mathcal O e+\mathcal O\der e+\cdots+\mathcal O \der^{r-1}e$
is a lattice in $M$.
\end{enumerate}
\end{examplesNumbered}}

\noindent
For $i=1,2$ let $M_i$ be a differential module over $K$ and $L_i$ be a lattice in $M_i$. Then~$L_1\oplus L_2$ is a lattice in the differential module
$M_1\oplus   M_2$ over $K$, and  the $\mathcal O$-submodule of $M_1\otimes_{K}M_2$ generated by the elements~$y_1\otimes y_2$ ($y_i\in L_i$, $i=1,2$) is a lattice in 
the differential module $M_1\otimes_K M_2$ over $K$.
Also, $\operatorname{Hom}(L_1,L_2)$ is a lattice in the differential module $\operatorname{Hom}_K(M_1,M_2)$ over $K$.

Let $L$ be a lattice in~$M$. If~$\pi\colon M\to N$ is a surjective morphism of differential modules over $K$, then~$\pi(L)$ is a lattice in $N$. If~$N$ is a differential submodule of $M$, then~$L\cap N$ is a lattice in $N$.
Using the notation from Corollary~\ref{cor:base change lattice} we have:

\begin{lemma}\label{lem:bounded base change}
Let $L$ be a lattice in $M$ and $E$ be a valued differential field extension of $K$ with small derivation.  
Then $L_{E}$ is a lattice in the base change~$M_{E}=E\otimes_K M$ of the differential module $M$ over $K$ to 
a differential module over $E$.
\end{lemma}

\noindent
If $A\in \mathcal{O}[\der]$ is monic (of order $r$ by our convention), then $\mathcal{O}[\der]A$ is a left ideal of the ring
$\mathcal{O}[\der]$, and  the resulting left $\mathcal{O}[\der]$-module $\mathcal{O}[\der]/\mathcal{O}[\der]A$ is free on $e, \der e,\dots, \der^{r-1}e$ for $e:= 1+\mathcal{O}[\der]A$, as is easily verified. Conversely, if $L$ is a left $\mathcal{O}[\der]$-module free on $e,\der e,\dots, \der^{r-1}e$, $e\in L$, then the unique  monic $A\in \mathcal{O}[\der]$ (of order $r$) such that~$A e=0$
yields an isomorphism $\mathcal{O}[\der]/\mathcal{O}[\der] A\to L$ sending $1+\mathcal{O}[\der]A$  to $e$. 
 
Next, let $\k=\mathcal O/\smallo$ be the differential residue field of $K$; cf.~[ADH, 4.4]. Here is a version of the cyclic vector theorem [ADH, 5.5.3] for lattices:

\begin{prop}\label{prop:cyclic vec lattice}
Suppose the derivation on $\k$ is nontrivial and $L$ is a lattice in~$M$ and $\dim_K M=r$. Then $L\cong \mathcal O[\der]/\mathcal O[\der]A$ for some monic $A\in\mathcal O[\der]$.
\end{prop}

\noindent
The case $r=0$ being trivial, we assume for the proof below that $r\ge 1$. We now introduce a tuple
$Y=(Y_0,\dots,Y_{r-1})$ of distinct differential indeterminates over $K$, and let  $i$, $j$, $k$, $l$ range   over $\{0,\dots,r-1\}$.

\begin{lemma}
For all $i$, $j$, let $P_{ij}\in Y_i^{(j)} + \sum_{k<r,\ l<j} K\, Y_k^{(l)}$.
  Then the coefficient of~$Y_0Y_1'\cdots Y_{r-1}^{(r-1)}$ in $\det(P_{ij})\in K\{Y\}$ is~$1$.
\end{lemma}
\begin{proof} For $p=0,\dots, r-1$ we prove by induction on $p$ that the coefficient of $Y_0Y_1'\cdots Y_{p}^{(p)}$ in 
$\det (P_{ij})_{i,j\le p}\in K\{Y\}$ is~$1$ (which for $p=r-1$ gives the desired result). The case $p=0$ is clear, so assume
$p\ge 1$. Then $Y_p^{(p)}$ occurs in the matrix $(P_{ij})_{i,j\le p}$ only in the $(p,p)$-entry $P_{pp}\in Y_p^{(p)}+\sum_{k<r,\ l<p}K\, Y_l^{(l)}$, and so the coefficient of~$Y_0Y_1'\cdots Y_{p}^{(p)}$ in 
$\det(P_{ij})_{i,j\le p}$ is the coefficient of $Y_0Y_1'\cdots Y_{p-1}^{(p-1)}$ in~$\det(P_{ij})_{i,j\le p-1}$, and the latter is $1$ by inductive assumption.  
\end{proof}

%\begin{proof}
%We proceed by  induction on $r$. The case $r=1$ is clear. For the inductive step suppose $r\geq 2$ and note that 
%only one entry of the matrix $(P_{ij})$ involves~$Y_{r-1}^{(r-1)}$, namely~$P_{r-1,r-1}$, and
%the coefficient of $Y_{r-1}^{(r-1)}$ in $P$ viewed as element of $(K\{Y_0,\dots,Y_{r-2}\})\{Y_{r-1}\}$ is 
%$\det(P_{ij})_{i,j<r-1}$.
%\end{proof}

\noindent
{\em Proof of Proposition~\ref{prop:cyclic vec lattice}}. Let~$z_0,\dots,z_{r-1}$ be a basis of the $\mathcal{O}$-module $L$.
Consider the base change~$\hat M:=K\{Y\}\otimes_K  M$ of $M$ to the differential $K$-algebra~$K\{Y\}$, cf.~[ADH, p.~304].
So $\hat M$ is a left $K\{Y\}[\der]$-module and
the $K\{Y\}$-module~$\hat M$ is free on $1\otimes z_0,\dots,1\otimes z_{r-1}$.
 Set
$$\hat e\ :=\ Y_0\otimes z_0 + \cdots + Y_{r-1}\otimes z_{r-1}\in \hat M$$
and  let $P_{ij}\in K\{Y\}$ be such that 
$$\der^j \hat e\  =\  P_{0j}\otimes z_0+\cdots+P_{r-1,j}\otimes z_{r-1}.$$
An easy induction on $j$ using $\der L\subseteq L$ shows that 
$P_{ij}\in Y_i^{(j)} + \sum_{k<r,\ l<j} \mathcal O\, Y_k^{(l)}$.
Put~$P:=\det(P_{ij})\in\mathcal O\{Y\}$; then $v(P)=0$ by the   lemma above, so
[ADH, 4.2.1] applied to $\k$ and the image of $P$
under the natural morphism $\mathcal O\{Y\}\to\k\{Y\}$ in place of $K$ and $P$, respectively, yields an~$a\in\mathcal O^r$ such that $P(a)\asymp 1$.
 We   obtain a   $K[\der]$-module morphism 
 $$\phi\colon \hat M\to M\quad\text{ with }\quad
 \phi(Q\otimes z)=Q(a)z \text{ for~$Q\in K\{Y\}$ and~$z\in M$.}$$ 
Put $R:=\mathcal O\{Y\}$, a differential subring of $K\{Y\}$, and
let $\hat L$ be the $R[\der]$-submodule of $\hat M$ generated by
$1\otimes z_0,\dots,1\otimes z_{r-1}$, so
$$ \hat L\  =\  \{ Q_0\otimes z_0+\cdots + Q_{r-1}\otimes z_{r-1}:\  Q_0,\dots, Q_{r-1}\in\mathcal{O}\{Y\}\}.$$
Then $\der^j\hat e\in\hat L$ for all $j$ and $\phi(\hat L)=L$. With $e:=\phi(\hat e)$ we have
$$\der^j e=\phi(\der^j\hat e)=P_{0j}(a)z_0+\cdots+P_{r-1,j}(a)z_{r-1}$$
and $\det(P_{ij}(a))=P(a)\in\mathcal O^\times$, so $L=\mathcal Oe+\mathcal O\der e+\cdots+\mathcal O\der^{r-1}e$. 
By a remark preceding Proposition~\ref{prop:cyclic vec lattice}, this concludes its proof.  \qed

\begin{remark}
Taking $K=\mathcal O$,  Proposition~\ref{prop:cyclic vec lattice} yields another proof of [ADH, 5.5.3], in the spirit of \cite{Cope};
cf.~\cite{ChKo}. Note also that $e$ as constructed in the proof of Proposition~\ref{prop:cyclic vec lattice} is a cyclic vector of $M$, and so yields the isomorphism $M\cong K[\der]/K[\der]A$ sending $e$ to $1+K[\der]A$, with monic $A\in \mathcal{O}$ (of order $r$ by convention) determined by the requirement $Ae=0$.
\end{remark}

\subsection*{Eigenvalues of bounded operators\astr}
{\it In this subsection $A\in\mathcal O[\der]$ is monic, of order $r$ by earlier convention.}\/
Recall that $[a]:=a+K^\dagger$ for $a\in K$. Put 
$$[\mathcal O]:=(\mathcal O+K^\dagger)/K^\dagger=\big\{[a]:a\in\mathcal O\big\}\qquad\text{ (a divisible subgroup of $K/K^\dagger$).}$$
Thus $\Sigma(A) \subseteq  [\mathcal O]$ by~\eqref{eq:spec A} and [ADH, 5.6.3]. 
More precisely, with $\k$ denoting the differential residue field of $K$, recall  from [ADH, 5.6]  that the residue map
${a\mapsto\res a\colon \mathcal O\to \k}$ extends to a 
 ring morphism $B\mapsto \res B\colon\mathcal O[\der]\to \k[\der]$ with~$\der\mapsto \der$. 
For each~$B\in\mathcal O[\der]$ and $y\in\mathcal O$ we have $B(y)\in\mathcal O$ and $\res\!\big( B(y) \big) = (\res B)(\res y)$. Also, $\res A\in \k[\der]$ is monic with $\order \res A=\order A=r$.
By~[ADH, 5.6.3], if~$B,D\in K[\der]$ are monic and $A=BD$, then $B,D\in \mathcal O[\der]$. 
In particular, all $a_1,\dots,a_r\in K$   such that~$A=(\der-a_1)\cdots(\der-a_r)$ are in $\mathcal O$, and~
$$\res A\ =\ (\der-\res a_1)\cdots(\der-\res a_r).$$
Moreover, using also \eqref{eq:spec A}, we conclude:

\begin{lemma}\label{lem:ev bded}
If $A=B(\der-a)$, $B\in K[\der]$, then $a\in\mathcal O$, $B\in\mathcal O[\der]$, and $\res A=(\res B)\cdot(\der-\res a)$.
Hence for each $\alpha\in\Sigma(A)$ there is an $a\in\mathcal O$ such that~$\alpha=[a]$    and~$\res a+\k^\dagger\in\Sigma(\res A)$.
\end{lemma}

\noindent
Suppose now $\der\mathcal O\subseteq\smallo$. 
Then the derivation of  $\k$   is trivial, and we have a $\k$-algebra isomorphism $P(Y)\mapsto P(\der)\colon\k[Y]\to\k[\der]$.
We let the {\bf characteristic polynomial} of $B$ be  the $\chi_B\in\k[Y]$ satisfying~$\chi_B(\der)=\res B$.\index{linear differential operator!characteristic polynomial}\index{characteristic polynomial}  
Then $B\mapsto\chi_B\colon\mathcal O[\der]\to\k[Y]$ is a ring morphism extending the residue morphism $\mathcal O\to\k$ with $\der\mapsto Y$.
Identifying~$\k/\k^\dagger$ with $\k$ in the natural way,
    the set of zeros of $\chi_A$ in $\k$ is $\Sigma(\res A)$.  
Thus by Lemma~\ref{lem:ev bded}:

\begin{cor}\label{cor:ev bded}
If  $A\in K[\der](\der-a)$, then $a\in\mathcal O$ and $\chi_A(\res a)=0$.
Hence for each~$\alpha\in\Sigma(A)$ there is an $a\in\mathcal O$ such that~$\alpha=[a]$    and~$\chi_A(\res a)=0$.
\end{cor}

\noindent
If $\mathcal O=C+\smallo$, then $\der\mathcal O\subseteq\smallo$ and the residue morphism $\mathcal O\to\k$ restricts to an isomorphism~$C\to\k$,   via which we identify $\k$ with $C$ making $\chi_A$  an element of~$C[Y]$.

\medskip  
\noindent 
{\em In the rest of this subsection  $K=H[\imag]$ where $H$ is a real closed differential subfield of $K$   such that the valuation ring $\mathcal O_H:=\mathcal O\cap H$ of $H$ is convex with respect to the ordering of $H$  and  $\mathcal O_H=C_H+\smallo_H$}. Then $C=C_H[\imag]$. A remark after Corollary~\ref{cor:logder abs value} then yields
~$\mathcal O=\mathcal O_H+\mathcal O_H\imag =C+\smallo$. Using that remark and Lemma~\ref{lem:logder} we have~$K^\dagger \subseteq H^\dagger + \smallo_H\imag\subseteq H^\dagger+\mathcal{O}$, and thus:

\begin{lemma}\label{lem:[O]}
$H^\dagger+\mathcal O\ =\  H^\dagger + C_H+ C_H\imag + \smallo\ =\ \big\{ a\in K: [a]\in[\mathcal O] \big\}$.
\end{lemma}

\begin{lemma}\label{lem:mult upper bd}
Suppose that  $C_H\subseteq H^\dagger$, and
let $\alpha\in[\mathcal O]$.   Then there is a unique~$b\in C_H$
such that $\alpha=[b\imag+\varepsilon]$ for some $\varepsilon\in \smallo$. For this  $b$ we have
\begin{equation}\label{eq:multA vs multP}
\operatorname{mult}_{\alpha}(A)\ \leq\ \sum_{c\in C,\  \Im c=b} \operatorname{mult}_c(\chi_A).
\end{equation}
\end{lemma}
\begin{proof}
Lemma~\ref{lem:[O]} and $C_H\subseteq H^\dagger$ yield the existence of $b\in C_H$ such that $\alpha=[b\imag+\varepsilon]$ for some~$\varepsilon\in \smallo$. Since   $K^\dagger\subseteq H + \smallo_H\imag$, 
there is at most one
 $b\in C_H$.  We prove \eqref{eq:multA vs multP}
by induction on $r$. The cases~$r=0$ and $\operatorname{mult}_{\alpha}(A)=0$ being trivial, suppose $r\geq 1$ and $\alpha\in\Sigma(A)$. 
From   
\eqref{eq:spec A} and~[ADH, 5.6.3] we get  $a\in \mathcal O$ and  monic $B\in\mathcal O[\der]$ with $[a]=\alpha$ and~$A=B(\der-a)$.
Then $\operatorname{mult}_{\alpha}(A)\leq \operatorname{mult}_{\alpha}(B)+1$ by Lemma~\ref{lem:spectrum fact}, and with $c\in C$
such that~$a-c\prec 1$ we have $b=\Im c$ and~$\chi_A(c)=0$.  Now apply the inductive hypothesis to $B$ 
in place of $A$. 
\end{proof}
 
\begin{remark}
The inequality \eqref{eq:multA vs multP} is strict in general: 
$H$ can be an $H$-field with an element $x\in H$ such that $x'=1$. Then $x\succ 1$, $1/x\notin \I(H)$, so
$\varepsilon:=\imag/x\in\smallo\setminus K^\dagger$. Then for
$A\ :=\ (\der-(\imag+\varepsilon))(\der-\imag)$
we have 
%$[i]\in \Sigma(A)$ and 
$\mult_{[\imag]}A=1$
%gives $\Sigma(A)=\big\{[\imag+\varepsilon],[\imag]\big\}$ with $\abs{\Sigma(A)}=2$,
while $\chi_A=(Y-\imag)^2$.
\end{remark}

\begin{cor}\label{cor:Ri unique, 3} Suppose $K$ has asymptotic integration and is $(r-1)$-newtonian, $r\ge 1$. Let $c_1,\dots,c_r\in C$ be the zeros of $\chi_A$,
and suppose $c_1,\dots,c_r$ are distinct and $\Re c_1\geq\cdots\geq\Re c_r$.
Then for each splitting  $(a_1,\dots,a_r)$ of $A$ over $K$ we have $a_1,\dots,a_r\in\mathcal O$. Moreover, 
there is a unique such splitting  of $A$ over $K$ such that~$a_1-c_1,\dots,a_r-c_r\prec 1$.
\end{cor}
\begin{proof}
The first claim is immediate from [ADH, 5.6.3]. We prove the second claim
by induction on $r$. The case $r=1$ being trivial, suppose $r>1$.  
Corollary~\ref{cor:simple Riccati zero} yields an $a_r\in\mathcal O$ with
$\operatorname{Ri}(A)(a_r)=0$ and $a_r-c_r\prec 1$, and then~$A=B(\der-a_r)$ where $B\in\mathcal O[\der]$ is monic,
by [ADH, 5.6.3, 5.8.7]. By inductive hypothesis~$B=(\der-a_1)\cdots(\der-a_{r-1})$ where
$a_j\in\mathcal O$ with $a_j-c_j\prec 1$ for $j=1,\dots,r-1$.
This shows existence. Uniqueness follows in a similar way, using
Corollary~\ref{cor:Ri unique, 2}.  
\end{proof}

\subsection*{Bounded differential modules\astr}
{\it In this subsection  $A$ is monic and $M$ is a differential module over $K$.}\/
We call $M$ {\bf bounded} if there exists a lattice in $M$.\index{bounded!differential module}\index{differential module!bounded}
By  remarks in an earlier subsection the class of bounded differential modules over $K$ is quite robust: if~$M_1$,~$M_2$ are bounded differential modules over $K$, then so are the differential modules~$M_1\oplus M_2$, $M_1\otimes_K M_2$, and~$\operatorname{Hom}_K(M_1,M_2)$ over~$K$, 
and if $M$ is bounded, then so is
every differential submodule of $M$, every
image of~$M$ under a morphism of differential modules over $K$, and every base change of~$M$ to
a valued differential field extension of $K$ with small derivation.

\begin{exampleNumbered}\label{ex:bounded dim 1}
Let $u\in K$ and suppose $M=K$ with $\der a=a'+ua$ for all $a$. Then for~$e\in K^\times$, $\mathcal Oe$  is a lattice in $M$ iff
$e'+ue=\der e   \in \mathcal Oe$. Hence $M$ is bounded iff~$u\in\mathcal O+K^\dagger$.
\end{exampleNumbered}

%\noindent
%Let $a\in K$. Since $B\mapsto B_{a}$ is an automorphism of  the ring $K[\der]$, we also have the subring~$\mathcal O[\der]_{a}:=\big\{B_{a}:B\in\mathcal O[\der]\big\}$ of $K[\der]$. 
%Note that~$\mathcal O[\der]_{a}=\mathcal O[\der]$ if   $a\in\mathcal O$; in particular, if the valued differential field $K$ is monotone, then $\mathcal O[\der]_{a}=\mathcal O[\der]$ whenever $a\in K^\dagger$.
%In Lemmas~\ref{lem:bded cyc, 1}--\ref{lem:bded cyc, 3} below we assume that~$M=K[\der]/K[\der]A$ where $A$ is monic of order~$r$.

%\begin{cor}
%The following are equivalent:
%\begin{enumerate}
%\item[\textup{(i)}] $M\cong K[\der]/K[\der]A$ for some monic $A\in\mathcal O[\der]$;
%\item[\textup{(ii)}] there is a monic $B\in\mathcal O[\der]$ and a cyclic vector $e$ of $M$  with $Be=0$.
%\end{enumerate}
%\end{cor}

%\noindent
%Notation: for~$a\neq 0$ we put $\mathcal O[\der]_{\ltimes a}:=\big\{B_{\ltimes a}:B\in\mathcal O[\der]\big\}$, a subring of $K[\der]$. 
%If~$A\in \mathcal O[\der]_{\ltimes a}$ for some $a\neq 0$, then $K[\der]/K[\der]A$ is bounded by Example~\ref{ex:lattices}(2).
\noindent
If $A_{\ltimes a}\in \mathcal{O}[\der]$ for some $a\ne 0$, then $K[\der]/K[\der]A$  is bounded by Example~\ref{ex:lattices}(2).  

\begin{lemma}\label{lem:bded cyc, 3}
Suppose  $M=K[\der]/K[\der]A$ and  $r=1$. Then 
$$M \text{ is bounded}\ \Longleftrightarrow\  A_{\ltimes a}\in\mathcal O[\der] \text{ for some }a\neq 0.$$
\end{lemma}
\begin{proof}
Let  $A=\der-u$, $u\in K$. Identifying $K$ with $M$ via $a\mapsto a+K[\der]A$ we have~$\der a=a'+ua$ for all $a$ in $K=M$, so
if $M$ is bounded, then Example~\ref{ex:bounded dim 1} gives $a\ne 0$ with $u\in \mathcal O+a^\dagger$,  hence
$A_{\ltimes a}\in \mathcal{O}[\der]$. 
\end{proof}

\begin{lemma}\label{lem:bded cyc, 2}
Suppose  the valuation ring $\mathcal O$ is discrete \textup{(}that is, a DVR\textup{)} and~$M=K[\der]/K[\der]A$  is bounded. Then~$A_{\ltimes a^{-1}}\in   \mathcal O[\der]$
for some~$a\in  \mathcal O^{\neq}$. 
\end{lemma}
\begin{proof}
Let $L$ be a lattice in $M$ and $e:=1+K[\der]A$, a cyclic vector of $M$. Since~$M/L$ is a torsion module we get $a\in \mathcal O^{\neq}$ with $f:=ae\in L$.
Because $\mathcal O$ is noetherian,
the submodule of the finitely generated $\mathcal O$-module $L$ generated by~$f,\der f,\der^2 f,\dots$ is itself finitely generated, and this yields $n$ with~$\der^n f\in \mathcal O f+\mathcal O\der f+\cdots+\mathcal O\der^{n-1}f$~\cite[Chap\-ter~X, \S{}1]{Lang}. We   obtain a monic~$B\in\mathcal O[\der]$ of order~$n$ with~$Bf=0$. Then~$B_{\ltimes a}e=0$, so $B\in K[\der]A_{\ltimes a^{-1}}$, and thus $A_{\ltimes a^{-1}}\in\mathcal O[\der]$ by~[ADH, 5.6.3]. 
\end{proof}

\noindent
 If $K$ is monotone $K$, then $v(B_{\ltimes a})=v(B)$ for all $B$ and $a\ne 0$, by [ADH, 4.5.4]. If~$\mathcal{O}$ is discrete, then $K$ is monotone by [ADH, 6.1.2]. Hence
 %If the valued differential field $K$ is monotone, then~$\mathcal O[\der]_{\ltimes a}=\mathcal O[\der]$ for each $a\neq 0$.
by Lemma~\ref{lem:bded cyc, 2}:

\begin{cor}\label{cor:bded cyc 1}
Suppose $\mathcal{O}$ is discrete. Then: 
$$\text{ the differential module }K[\der]/K[\der]A  \text{ over $K$ is bounded}\ \Longleftrightarrow\ A\in\mathcal O[\der].$$
\end{cor}
\begin{remark}
In the case
$(K,\mathcal O)=\big(\C(\!(t)\!),\C[[t]]\big)$ and~$\der=t\frac{d}{dt}$, bounded differential modules over $K$
are called {\it regular singular}\/ in \cite{vdPS}, and Corollary~\ref{cor:bded cyc 1} in this case is implicit in the proof of \cite[Pro\-po\-sition~3.16]{vdPS}.
\end{remark}

\begin{lemma}\label{lem:bded cyc, 4}
Suppose $M\cong K[\der]/K[\der]A$ where $A\in\mathcal O[\der]$. Then  $\Sigma(M)\subseteq [\mathcal O]$.
\end{lemma}
\begin{proof} The case $r=0$ is trivial, so assume $r\ge 1$. 
%By an earlier remark $M$ is bounded. 
  Corollary~\ref{cor:kerder} and  remarks in the last subsection yield~$\Sigma(M^*)=\Sigma(A)\subseteq [\mathcal O]$. By [ADH, 5.5.8] we have~$M^*\cong K[\der]/K[\der]B$ where $B:=(-1)^rA^*\in\mathcal O[\der]$ is monic. Also $M\cong M^{**}$, hence
$\Sigma(M)=\Sigma(M^{**})\subseteq [\mathcal O]$ by the above applied to $M^*$, $B$ in place of~$M$,~$A$.
\end{proof}

\begin{cor}\label{cor:bded cyc 2}
Suppose $M$ is bounded. Assume also that the derivation of $\k$ is nontrivial, or
the derivation of~$K$ is nontrivial and $\mathcal O$ is discrete. Then~$\Sigma(M)\subseteq[\mathcal O]$.
\end{cor}
\begin{proof} The remark following the proof of
Proposition~\ref{prop:cyclic vec lattice} and Lemma~\ref{lem:bded cyc, 2} yield $A\in\mathcal O[\der]$
with $M\cong K[\der]/K[\der]A$, and so $\Sigma(M)\subseteq[\mathcal O]$ by Lemma~\ref{lem:bded cyc, 4}.
\end{proof}

\begin{cor}\label{cor:bded cyc 3}
Suppose $M$ splits and is bounded. Then  $\Sigma(M)\subseteq [\mathcal O]$.
\end{cor}
\begin{proof}
We proceed by induction on $r=\dim_K M=\ell(M)$. If $r=0$, then~${\Sigma(M)=\emptyset}$. Next suppose $r=1$, and take $e\in M$ with $M=Ke$
and $a\in K$ with~$\der e=ae$. Then~$\Sigma(M)=\big\{[a]\big\}$ by Corollary~\ref{Sigmaatmostr}, and $M\cong K[\der]/K[\der](\der-a)$, so $[a]\in [\mathcal O]$ by Lemma~\ref{lem:bded cyc, 3}.
Now suppose  $r>1$. Take a differential submodule $M_1$ of $M$ with~$\ell(M_1)=r-1$, so~$\ell(M_2)=1$ for $M_2:=M/M_1$.
Then   $M_1$, $M_2$ split and are bounded, by Lemma~\ref{lem:split in short exact sequ} and the remarks before Example~\ref{ex:bounded dim 1}. Hence by inductive hypothesis~$\Sigma(M_i)\subseteq [\mathcal O]$ for $i=1,2$, so $\Sigma(M)\subseteq [\mathcal O]$ by the remark after~\eqref{eq:mult subadd}.
\end{proof}

\begin{cor}\label{cor:bded cyc 4}
Let  $H$ be a Liouville closed, trigonometrically closed $H$-field with small derivation, $K=H[\imag]$,
and suppose $M$ is bounded. Then $\Sigma(M)\subseteq[\mathcal O]$.
\end{cor}
\begin{proof}
We have $K^\dagger=H\oplus\I(H)\imag$ and so $\mathcal O+K^\dagger=H\oplus\mathcal O_H\imag$.
Take an $H$-closed field extension~$H_1$ of $H$ and set $K_1:=H_1[\imag]$.  The base change $M_1:=K_1\otimes_K M$ of $M$ to $K_1$ splits and is bounded, so
$\Sigma(M_1)\subseteq[\mathcal O_1]$ by Corollary~\ref{cor:bded cyc 3}. 
We have~$K_1^\dagger=H_1\oplus\I(H_1)\imag$ by Corollary~\ref{cor:sc=>tc}, so
$\mathcal O_1+K_1^\dagger=H_1\oplus\mathcal O_{H_1}\imag$.
This yields~$K_1^\dagger\cap K=K^\dagger$ and~$(\mathcal O_1+K_1^\dagger)\cap K=\mathcal O+K^\dagger$, so
  identifying~$K/K^\dagger$ with its image under the group embedding $a+K^\dagger\mapsto a+K_1^\dagger\colon K/K^\dagger\to K_1/K_1^\dagger$ ($a\in K$)  we have~$\Sigma(M)\subseteq\Sigma(M_1)$  and~$[\mathcal O] = [\mathcal O_1] \cap (K/K^\dagger)$.
Thus $\Sigma(M)\subseteq [\mathcal O]$.
\end{proof}

\begin{question} Does it follow from $M$ being bounded and $C\ne K$ that $\Sigma(M)\subseteq[\mathcal O]$?
                       %Suppose $C\neq K$ and $M$ is bounded; is then $\Sigma(M)\subseteq[\mathcal O]$?
\end{question}

\section{Self-Adjointness and its Variants\astr}\label{sec:self-adjoint} 

\noindent
{\it In this section $K$ is a differential field. We let $A$, $B$ range over $K[\der]$ with $A\neq 0$, and set $r:=\order A$. We also let $\alpha$ range over $K/K^\dagger$.}\/
The material in this section elaborates on Corollaries~\ref{cor:spectrum, self-adjoint, 1} and \ref{cor:spectrum, self-adjoint, 2} and shows how  symmetries of $A$ force it to have eigenvalue~$0\in K/K^\dagger$, mainly by making some classical results (cf.~\cite[Chapitre~V]{Darboux} 
and~\cite[\S{}23--25]{Schlesinger}) precise and putting them into our present context. It can be skipped on first reading, since it is only needed in Section~\ref{sec:lin diff applications} for applications of  our main theorem  to linear differential equations over complexified Hardy fields.

\subsection*{Operators of the same type}
Suppose $B\neq 0$ and set $s:=\order B$. Consider now the $C$-linear subspace
$$\mathcal E(A,B)\ :=\  \big\{ R\in K[\der]: \order R<r \text{ and } BR\in K[\der]A\big\}$$
of $K[\der]$. 
The next lemma and its corollary elaborate on Lemma~\ref{lem:same type}.

\begin{lemma}\label{lem:E(A,B)}
Let $M:=K[\der]/K[\der]A$ and $N:=K[\der]/K[\der]B$. Then
we have an isomorphism 
$$R\ \mapsto\  \phi_R\ \colon\ \mathcal E(A,B)\to  \Hom_{K[\der]}(N,M)$$
 of $C$-linear spaces where 
$$\phi_R(1+K[\der]B)\ =\ R+K[\der]A\qquad\text{ for $R\in\mathcal E(A,B)$.}$$
\end{lemma}
\begin{proof}
Let $R\in \mathcal E(A,B)$. Then the kernel of the
$K[\der]$-linear map $K[\der]\to K[\der]/K[\der]A$ sending $1$ to $R+K[\der]A$
contains $K[\der]B$, hence induces a $K[\der]$-linear map $$\phi_R\colon N=K[\der]/K[\der]B\to K[\der]/K[\der]A=M$$
as indicated.  It is easy to check that $R\mapsto\phi_R$ is $C$-linear.
If $\phi_R=0$, then $\phi_R({1+K[\der]B})=K[\der]A$ and hence $R\in K[\der]A$, so $R=0$, since $\order R<r$.
Given~$\phi\in \Hom_{K[\der]}(N,M)$ we have $\phi(1+K[\der]B) = R+K[\der]A$ where $R\in K[\der]$ has or\-der~$<r$;
then $BR\in K[\der]A$, so~$R\in\mathcal E(A,B)$, and we have $\phi=\phi_R$.
\end{proof}

\noindent
In particular, $0\leq\dim_C \mathcal E(A,B) \leq rs$ by Lemmas~\ref{dim of ann} and~\ref{lem:E(A,B)}. Moreover:

\begin{cor}\label{cor:tildeE}
Suppose $r=s$. Then the isomorphism $R\mapsto\phi_R$ from the previous lemma maps
the subset 
$$\mathcal E(A,B)^\times\  :=\  \big\{ R\in \mathcal E(A,B):  1\in K[\der]R+K[\der]A\big\}$$
of $\mathcal E(A,B)$ bijectively onto the set of $K[\der]$-linear isomorphisms $N\to M$.
\end{cor}

\noindent
Set $M:=K[\der]/K[\der]A$. We make the $C$-module $\End_{K[\der]}(M):=\Hom_{K[\der]}(M,M)$ into a $C$-algebra with its ring multiplication given by
composition. 
We equip the $C$-module $\mathcal E(A):=\mathcal E(A,A)$ with the ring multiplication making  the map $$R\mapsto\phi_R\ \colon\  \mathcal E(A)\to \End_{K[\der]}(M)$$
an isomorphism of $C$-algebras. The $C$-algebra $\mathcal E(A)$ is called the {\bf eigenring} of~$A$; cf.~\cite{Singer96} or \cite[\S{}2.2]{vdPS}.\index{linear differential operator!eigenring}\index{eigenring!linear differential operator}
Note that if $r\geq 1$, then $C\subseteq\mathcal E(A)$. 
If the $K[\der]$-module $M$ is irreducible, then 
$\End_{K[\der]}(M)$ is a division ring, by Schur's Lemma~\cite[Chap\-ter~XVII, Proposition~1.1]{Lang}. Now $M$ is irreducible iff
$A$ is irreducible [ADH, p. 251], hence:

\begin{cor}\label{cor:tildeE(A)=C}
Suppose $A$ is irreducible. Then $\mathcal E(A)$ is a division algebra over~$C$. 
If  $C$ is algebraically closed, then $\mathcal E(A)=C$.
\end{cor}
\begin{proof} As to the second claim, let $e\in \mathcal E(A)$. The elements of $C$ commute with $e$, so we have a commutative domain $C[e]\subseteq \mathcal E(A)$, hence $e$ is algebraic over $C$ in view of $\dim_C\mathcal E(A)\le r^2$, and thus $e\in C$ if $C$ is algebraically closed. 
\end{proof}

\noindent
We may have $\mathcal E(A)=C$ without $A$ being irreducible \cite[Exercise~2.14]{vdPS}. 
%\marginpar{this exercise skipped}
If~$A$,~$B$ have the same type, then the $C$-algebras $\mathcal E(A)$, $\mathcal E(B)$ are isomorphic. By Lem\-ma~\ref{lem:E(A,B)} and Corollary~\ref{cor:tildeE} we have:

\begin{cor}\label{cor:tildeE(A,B)}
Suppose $\mathcal E(A)=C$ and $A$, $B$ have the same type. Then for  some~$e\in \mathcal E(A,B)^{\ne}$ we have
$\mathcal E(A,B)=Ce$, and  $\mathcal E(A,B)^\times=C^\times e$.
\end{cor}

\subsection*{Self-duality}
Let $M$ be  a differential module over $K$. We say that $M$ is {\bf self-dual} if~$M\cong M^*$.\index{differential module!self-dual}\index{self-dual!differential module}
If $M$ is self-dual, then so is of course every isomorphic $K[\der]$-module, in particular $M^*$.
%If $M$ is self-dual, then so is $M^*$, and $M$ is reflexive, that is,~$M\cong M^{**}$. (The latter is automatic if $M$ is a differential module over $K$.)
Given also a differential module $N$ over $K$,
we say that a $K$-bilinear map ${[\  ,\,  ]}\colon M\times N\to K$   is {\bf $\der$-compatible}\index{bilinear form!$\der$-compatible} if 
$$\der [f,g]\ =\ [\der f,g] + [f,\der g]\qquad\text{for all $f\in M$, $g\in N$.}$$
The non-degenerate $K$-bilinear map
$$(\phi,f)\mapsto \langle \phi,f\rangle:=\phi(f)\ \colon\ M^*\times M\to K$$ 
 is
$\der$-compatible by [ADH, (5.5.1)]. One verifies easily:

\begin{lemma}\label{lem:bilin}
 $M$ is self-dual iff there is a non-degenerate $\der$-compatible
$K$-bilinear form on $M$. In more detail, any isomorphism $\iota\colon M\to M^*$ yields a non-degenerate $\der$-compatible
$K$-bilinear form $(f,g)\mapsto \langle \iota(f),g\rangle\colon M\times M \to K$, and every
non-degenerate $\der$-compatible $K$-bilinear form on~$M$ arises in this way from a unique
isomorphism~$\iota\colon M\to M^*$ \textup{(}of differential modules over $K$\textup{)}. 
\end{lemma}

\noindent
In terms of matrices, let $e_1,\dots, e_n$ be a basis of $M$ and let $e_1^*,\dots, e_n^*$ be the dual basis of $M^*$. Let
$\iota\colon M\to M^*$ be an isomorphism $M\to M^*$ with matrix $P$ with respect to these bases. Then for the 
corresponding $K$-bilinear form ${[\ ,\, ]}$ on $M$ from the above lemma we have $[e_i, e_j]=P_{ji}=(P^{\operatorname{t}})_{ij}$. 

Next a consequence of Corollaries~\ref{cor:spectrum, self-adjoint, 2} and~\ref{cor:spectrum, self-adjoint, 3}. It provides useful information about the spectrum of $M$, which explains our interest in self-duality.    

\begin{cor}\label{cor:self-dual d-module} Let $\dim_K M=r$, and
suppose $\sum_\alpha \mult_\alpha(M)=r$ and $K$ is $1$-linearly surjective, or~${r\geq 1}$ and~$K$ is ${(r-1)}$-linearly surjective.
Assume also that~$M$ is self-dual.
Then $\mult_\alpha(M)=\mult_{-\alpha}(M)$ for all~$\alpha$.
Hence if   additionally~$K^\dagger$ is $2$-divisible and $\sum_\alpha \mult_\alpha(M)$ is odd, then~$0\in\Sigma(M)$.
\end{cor}

\noindent
Suppose now that $M=K[\der]/K[\der]A$ and $r\geq 1$. Then $M^*\cong K[\der]/K[\der]A^*$ by~[ADH, 5.5.8], hence $M$ is self-dual iff $A$, $A^*$ have
the same type. By Lemma~\ref{lem:same type} this is the case iff
there are $R,S\in K[\der]$ of order~$<r$ with $1\in K[\der]R+K[\der]A$ and~$A^*R=SA$. 
When $\mathcal E(A)=C$,  we can replace this with a more symmetric condition:

\begin{lemma}
Suppose $\mathcal E(A)=C$. Then $A$, $A^*$ have the same type iff for some~$R\in K[\der]$ of order~$n<r$ we have
$A^*R=(-1)^{n+r}R^*A$ and $1\in K[\der]R+K[\der]A$. 
\end{lemma}
\begin{proof}
Suppose $A$, $A^*$ have the same type. By Corollary~\ref{cor:tildeE(A,B)} we obtain $R,S\in K[\der]$ of order~$<r$ such that
$A^*R=SA$,  $1\in K[\der]R+K[\der]A$, and $\mathcal E(A,A^*)^\times=C^\times R$.
Now taking adjoints yields  $A^*S^*=R^*A$, so $0\ne S^*\in \mathcal E(A,A^*)$, hence $S^*\in \mathcal E(A,A^*)^\times$ and thus $S^*=cR$ ($c\in C^\times$). Comparing coefficients of the highest order terms on both sides of $cA^*R=R^*A$ gives $c=(-1)^{n+r}$.
\end{proof}

\noindent
We say that $A$ is {\bf self-dual} if $A$, $A^*$ have the same type.\index{self-dual!linear differential operator}\index{linear differential operator!self-dual} Thus if $A$ is self-dual, then so is $A^*$, and
so is every operator of the same type as $A$.
Moreover, by Lemma~\ref{lem:same type and eigenvalues}, if~$A$ is self-dual, then $A$, $A^*$ have the same eigenvalues, with the same multiplicities.
Combining Corollary~\ref{cor:tildeE(A)=C} with the previous lemma yields:

\begin{cor}
Suppose $A$ is irreducible and $C$ is algebraically closed. Then $A$ is self-dual iff for some~$R\in K[\der]$ of order~$n<r$ we have
$A^*R=(-1)^{n+r}R^*A$ and~$1\in K[\der]R+K[\der]A$. 
\end{cor}

\noindent
Here is the operator version of Corollary~\ref{cor:self-dual d-module}:

\begin{cor} \label{cor:self-dual operator}
Suppose  $A$ is self-dual, and set $s:=\sum_\alpha \mult_\alpha(A)$. Also assume~$K$ is $1$-linearly surjective and~$s=r$, or 
$r\geq 1$ and~$K$ is ${(r-1)}$-linearly surjective.
Then~$\mult_\alpha(A)=\mult_{-\alpha}(A)$ for each~$\alpha$.
Hence if in addition~$K^\dagger$ is $2$-divisible and~$s$ is odd, then~$0\in\Sigma(A)$.
\end{cor}

\noindent 
Let $\phi\in K^\times$, $B\ne 0$. If $A$, $B$ have the same type, then so do $A^\phi,B^\phi\in K^\phi[\derdelta]$, by Lemma~\ref{lem:same type}.  Hence by the next lemma, if $A$ is self-dual, then so is $A^\phi$.

\begin{lemma}
$(A^\phi)^*=(A^*)^\phi_{\ltimes\phi}$.
\end{lemma}
\begin{proof}
We have 
$$(\der^\phi)^*\ =\ (\phi\derdelta)^*\ =\ -\derdelta\phi\ 
%=-\phi\derdelta-\phi^\dagger
=\ (-\phi\derdelta)_{\ltimes\phi}\ =\ (-\der)^\phi_{\ltimes\phi}\ =\ (\der^*)^\phi_{\ltimes\phi},$$
so the lemma holds for $A=\der$. It remains to note that $B\mapsto (B^\phi)^*$ and $B\mapsto (B^*)^\phi_{\ltimes\phi}$ are ring morphisms $K[\der]\to K^\phi[\derdelta]^{\operatorname{opp}}$
that are the identity on $K$, where
$K^\phi[\derdelta]^{\operatorname{opp}}$ is the opposite ring of $K^\phi[\derdelta]$; cf.~[ADH, proof of 5.1.8].
\end{proof}

\noindent
If $A^*=(-1)^rA_{\ltimes a}$  ($a\in K^\times$), then~$A$ is self-dual, so there is a  non-degenerate $\der$-compatible $K$-bilinear form on $K[\der]/K[\der]A$. The next proposition gives more information.
We say that a  $K$-bilinear form  ${[\ ,\, ]}$ on  a $K$-linear space~$M$ is {\bf $(-1)^{n}$-symmetric} if~$[f,g]=(-1)^{n}[g,f]$ for all $f,g\in M$.\index{bilinear form!$(-1)^{n}$-symmetric}

\begin{prop}[Bogner~{\cite{Bogner}}]\label{prop:Bogner}
{\samepage Suppose   $r\geq 1$, and let $M=K[\der]/K[\der]A$ and~$e=1+K[\der]A\in M$.
Then the following are equivalent:
\begin{enumerate}
\item[\textup{(i)}] $A^* = (-1)^r A_{\ltimes a}$ for some $a\in K^\times$;
\item[\textup{(ii)}] there is a non-degenerate $\der$-compatible  $K$-bilinear form ${[\ ,\, ]}$ on $M$ such that
$$[e,\der^j e]=0\quad\text{for $j=0,\dots,r-2$.}$$
\end{enumerate}
Any  form  ${[\ ,\, ]}$ on $M$ as in \textup{(ii)} is $(-1)^{r-1}$-symmetric.}
\end{prop}
\begin{proof}
We first arrange that $A$ is monic.
Let $e_0^*,\dots,e_{r-1}^*$ be the basis of $M^*$ dual to the basis $e,\der e,\dots,\der^{r-1}e$ of~$M$, so~$e^*:=e_{r-1}^*$ is
a cyclic vector of $M^*$ with $A^*e^*=0$, by [ADH, 5.5.7]. Below we let~$i$,~$j$ range over~${\{0,\dots,r-1\}}$.

Suppose $A^* = (-1)^r A_{\ltimes a}$, $a\in K^\times$. Then $A^*e^*=0$ gives $Aae^*=0$, so
we have a $K[\der]$-linear isomor\-phism $\varphi\colon M\to M^*$ with $\varphi(e)=ae^*$.
Let~${[\ ,\, ]}$ be the non-degenerate $\der$-compatible $K$-bilinear form on $M$
given by $[f,g] = \langle \varphi(f),g\rangle$.
% for all~$f,g\in M$. 
Then
$$[e,\der^j e]=\langle \varphi(e),\der^j e\rangle=a\langle  e^*,\der^j e\rangle\quad\text{ for all $j$,}$$ 
proving~(ii).
Suppose conversely that    ${[\ ,\, ]}$ is as in~(ii). Then~$a:=[e,\der^{r-1}e]\neq 0$ since~${[\ ,\, ]}$ is non-degenerate.
Let $\varphi\colon M\to M^*$ be the isomorphism with~$\varphi(f)=[f,{-}]$ for all~$f\in M$.
Then $[e,\der^j e] = \langle\varphi(e),\der^j e\rangle$ for all~$j$ and thus $\varphi(e)=a\,e^*$.
Hence
$$Aa\,e^*\  =\  A\varphi(e)\  =\  \varphi(Ae)\  =\  \varphi(0)\  =\  0,$$ and this yields $A^* = (-1)^r A_{\ltimes a}$.

Let now ${[\ ,\, ]}$  be as in \textup{(ii)} and set $a:=[e,\der^{r-1}e]$. Induction on $i$ using $\der$-compatibility of ${[\ ,\, ]}$ shows $[\der^i e, \der^j e]=0$ for $i\leq r-2$, $j\leq r-2-i$. In particular, $[\der^i e,e]=0=[e,\der^i e]$ for $i\leq r-2$. 
Induction on $i$ using the second display in the proof of [ADH, 5.5.7]  gives
\begin{align*}  (\der^*)^i e^*\ &\in\  e_{r-1-i}^* + \sum_{r-i\le j\le r-1} Ke_j^*, \text{ and hence}\\ 
(-1)^{r-1}\der^{r-1}e^*\ &\in\  e_0^*+ K e_1^*+\cdots+K e_{r-1}^*.
\end{align*}
It follows that
$[\der^{r-1}e,e]=\langle \der^{r-1}ae^*,e\rangle=(-1)^{r-1}a$.
This covers the base case $i=0$ of an  induction on $i$ showing
$[\der^i e,g]=(-1)^{r-1} [g,\der^i e]$ for all $g\in M$. 
Suppose this identity holds for a certain $i\leq r-2$. Then  
by $\der$-compatibility
$$[\der^{i+1}e,g] 	= \der[\der^i e,g] - [\der^i e,\der g] 
					= (-1)^{r-1}\big(\der[g,\der^i e] - [\der g,\der^i e]\big) = (-1)^{r-1} [g,\der^{i+1}e]$$
as required.
\end{proof}

\noindent
See~\cite{Bogner} for the geometric significance of operators as in this proposition
when $K$ is the differential field of germs of meromorphic functions at $0$.

\medskip
\noindent
Let $r\geq1$, and $A, a$ be as in (i) of  Proposition~\ref{prop:Bogner}, with
$$A\ =\ \der^r+a_{r-1}\der^{r-1}+\cdots + a_0\quad (a_0,\dots, a_{r-1}\in K).$$ 
Then $a^\dagger=-(2/r) a_{r-1}$. Set~$b:=a^\dagger$. Then the operator
$B:=A_{b/2}$ of order~$r$ satisfies~$B^*=(-1)^r B$, so
the cases~$B^*=B$ and $B^*=-B$ deserve particular attention:
 
\begin{definition}\label{defsask}
A linear differential operator $B$ is said to be (formally) {\bf self-adjoint} if $B^*=B$ and   {\bf skew-adjoint} if $B^*=-B$.\index{self-adjoint!linear differential operator}\index{linear differential operator!self-adjoint}\index{skew-adjoint}\index{linear differential operator!skew-adjoint}
\end{definition}

\noindent
We discuss self-adjoint and skew-adjoint operators in more depth after reviewing a useful   identity  relating
a linear differential operator and its adjoint, which is obtained
by transferring [ADH, (5.5.1)] to the level of operators.

\subsection*{The Lagrange Identity}
Let $M$ be a differential module over $K$ with $\dim_K M=n\ge 1$, and suppose $e$ is a cyclic vector of $M$. Then 
$e_0,\dots,e_{n-1}$ with $e_i:=\der^ie$ is a basis of $M$. Let $e_0^*,\dots,e_{n-1}^*$ be the dual basis of $M^*$.
Then $e^*:= e_{n-1}^*$ is a cyclic vector of $M^*$ [ADH, 5.5.7], so~$e^*,\der e^*,\dots,\der^{n-1} e^*$ is a basis of $M^*$.
Let $L\in K[\der]$ be monic of order $n$ such that~$Le=0$.  Then
$$L=a_0+a_1\der+\cdots+a_n\der^n\qquad\text{
where $a_0,\dots,a_n\in K$  (so $a_n=1$).}$$ 
By [ADH, 5.5.7]  and its proof we have $L^*e^*=0$ and
\begin{equation}\label{eq:5.5.7}
e_{n-i-1}^* = L_ie^* \quad\text{where }L_i:= \sum_{j=0}^i (\der^*)^{i-j}a_{n-j}\in K[\der] \quad (i=0,\dots,n-1).
\end{equation}
%Perhaps more convenient:
%$$e_{i-1}^* = \left( \sum_{j=i}^n (\der^*)^{i-j}a_{j}\right) e^*\qquad\text{for $i=1,\dots,n$.}$$
%Or even:
%$$e_{i}^* = \left( \sum_{j=i+1}^n (\der^*)^{i+1-j}a_{j}\right) e^*\qquad\text{for $i=0,\dots,n-1$.}$$
%By [ADH, (5.1.1)], for $i=0,\dots,n-1$ we have
%$$\sum_{j=0}^i  (\der^*)^{i-j}a_{n-j}=\sum_{k=0}^i b_{ik}\, \der^k\quad    \text{where } b_{ik} := \sum_{j=0}^{i-k} (-1)^{i-j} {i-j \choose k} a_{n-j}^{(i-j-k)}\in K.$$
%(So $b_{ii}=(-1)^i$.) Thus
%$$e_{n-i-1}^* = b_{i0}e^*+b_{i1}\der e^*+\cdots+b_{ii}\der^i e^*.$$ 
Let  $d_0,\dots,d_{n-1}$ be the basis of $M$ dual to the basis $e^*,\der e^*,\dots,\der^{n-1}e^*$ of $M^*$.
%Then for $i,j=0,\dots,n-1$ we have
%\begin{equation}\label{lem:coeff}
%\langle e_{n-i-1}^*, d_j \rangle = b_{ij} \text{ if $j\leq i$,} \qquad \langle e_{n-i-1}^*, d_j \rangle = 0 \text{\ if $j>i$.}
%\end{equation}
Then
\begin{equation}\label{lem:coeff, n-1}
\langle e_{n-1}^*, d_j \rangle=\delta_{0j},\quad \langle e^*_{n-i-1}, d_{n-1} \rangle = (-1)^{n-1}\delta_{i,n-1}\quad \text{ for }1\le i\le n-1 
\end{equation}
(Kronecker deltas)  using \eqref{eq:5.5.7}.
Let $y,z\in K$ and set
$$\phi:=y e_0^*+y'  e_1^* +\cdots + y^{(n-1)} e_{n-1}^*\in M^*, \quad f:=z d_0 + z'd_1 + \cdots + z^{(n-1)}d_{n-1}\in M.$$ 
Then 
$$\der\phi\  =\  L(y)e_{n-1}^*, \qquad \der f\  =\  (-1)^n L^*(z) d_{n-1}.$$
For the first equality, use the first display in the proof of [ADH, 5.5.7].
%that by [ADH, proof of 5.5.8] the matrix of $M^*$ with respect to the basis~$e_0^*,\dots,e_{n-1}^*$ is the companion matrix~$A_L$ of $L$.
The second equality follows 
from the first by reversing the roles of $M$ and~$M^*$ and
noting that~$(-1)^nL^*$ is monic of order~$n$ with~$(-1)^nL^*e^*=0$. Hence
$$\langle \der\phi,f\rangle = L(y)\langle e^*,f \rangle = L(y)z, \qquad \langle \phi,\der f\rangle = (-1)^n L^*(z)\langle \phi,d_{n-1}\rangle = - L^*(z)y$$
by \eqref{lem:coeff, n-1} and so by the identity [ADH, (5.5.1)],
$$\der\langle \phi,f\rangle\ =\  \langle \der\phi,f\rangle+ \langle \phi,\der f\rangle\ =\ L(y)z-L^*(z)y.$$
Now $\<\der^ie^*,f\>=z^{(i)}$ for $i<n$, so 
$\langle Be^*, f\rangle = B(z)$ for all $B$ of order $<n$, hence
$$\langle \phi, f \rangle = \sum_{i=0}^{n-1}  y^{(n-i-1)} \langle L_ie^*,f\rangle =   \sum_{0\leq j\leq i< n} y^{(n-i-1)} (-1)^{i-j}(a_{n-j}z)^{(i-j)}= P_L(y,z)$$
where
\begin{equation}\label{eq:P_L}
P_L(Y,Z)\  :=\  \sum_{0\leq i\leq j< n} Y^{(i)}(-1)^{j-i} (a_{j+1} Z)^{(j-i)} \in K\{Y,Z\},
\end{equation}
a homogeneous differential polynomial of degree $2$. These considerations   show:

\begin{prop}[Lagrange Identity] \label{prop:Lagrange}
The map $$(y,z)\mapsto [y,z]_L:=P_L(y,z)\ \colon\ K\times K\to K$$ is $C$-bilinear, and for $y,z\in K$ we have
\begin{equation}\label{eq:Lagrange}
\der\big([y,z]_L\big)\  =\ L(y)z-L^*(z)y.
\end{equation}
\end{prop}

\noindent
We assumed here that $L$ is monic of order $n\geq 1$. For an arbitrary linear differential operator $L=a_0+a_1\der+\cdots+a_n\der^n\in K[\der]$ ($a_0,\dots,a_n\in K$) we define $P_L$ as in~\eqref{eq:P_L} and set~$[y,z]_L:=P_L(y,z)$ for $y,z\in K$. 
Then Proposition~\ref{prop:Lagrange} continues to hold;
to see \eqref{eq:Lagrange} for $L\neq 0$, reduce to the monic case, using that for all $a\in K$ and~$L\in K[\der]$
 we have 
\begin{equation}\label{eq:PaL} 
P_{aL}(Y,Z)=P_L(Y,aZ),\qquad\text{ hence }\qquad [y,z]_{aL}=[y,az]_L\text{ for $y,z\in K$.}
\end{equation} 
The differential polynomial $P_L$  is called the {\bf concomitant} of~$L$; \index{concomitant}\index{linear differential operator!concomitant}\index{bilinear form!concomitant} it does not change when passing from $K$ to a differential field extension.

\begin{lemma}[Hesse]\label{lem:Hesse}
Let $L,L_1,L_2\in K[\der]$; then 
$$P_{L^*}(Y,Z)=-P_L(Z,Y)\quad\text{and}\quad P_{L_1+L_2}(Y,Z)=P_{L_1}(Y,Z)+P_{L_2}(Y,Z).$$
\end{lemma}
\begin{proof}
By \eqref{eq:Lagrange} we have~$\der\big([y,z]_{L^*}\big) = -\der\big([z,y]_L\big)$ for all $y$,~$z$ in every differential field extension of
$K$, hence $\big(P_{L^*}(Y,Z)+P_{L}(Z,Y)\big){}'=0$ in~$K\{Y,Z\}$ and
then~$P_{L^*}(Y,Z)+P_{L}(Z,Y)=0$, since $P_{L^*}$, $P_L$ are homogeneous of degree $2$.
This shows the first identity. The second identity is clear by inspection of \eqref{eq:P_L}.
\end{proof}

\begin{example}
For $L=0, 1, \der, \der^2$ we have
$$ P_0=P_1=0,\quad P_\der = YZ, \quad P_{\der^2}= Y'Z-YZ'$$
so for $y,z\in K$:  
$$[y,z]_0=[y,z]_1=0,\quad [y,z]_{\der}=yz,\quad  [y,z]_{\der^2} = y'z-yz',$$ which for~$L=a\der^2+b\der+c$ ($a,b,c\in K$), using \eqref{eq:PaL}, gives
$$ P_L\ =\ aY'Z-aYZ' + (b-a')YZ, \quad [y,z]_L\ =\ ay'z-ayz'+(b-a')yz.$$
\end{example}

%\noindent {what follows has been checked but has been commented out, since it seems no longer needed} 
%For later use we note the following immediate consequence of the definition of~$P_L$:

%\begin{lemma}\label{lem:PL}
%Let $[ Y']$ be the differential ideal of $K\{Y,Z\}$ generated by $Y'$. Then
%$$P_L(YZ,Z) \equiv YP_L(Z,Z) \equiv P_L(Z,YZ) \mod [ Y'].$$
%\end{lemma}

\noindent
Below we use that evaluating the differential operator $L\in K[\der]$ at
the element~$Y$ of the differential ring extension $K\{Y\}$  of $K$ results in the differential polynomial~$L(Y)\in K\{Y\}$, which is homogeneous of degree $1$.
With this notation, we have~$P_L(Y,Z)'=L(Y)Z-L^*(Z)Y$. The next result  characterizes the concomitant and adjoint of a differential operator accordingly. 
 
\begin{lemma}[Frobenius]\label{lem:Frobenius}
The pair $(P_L, L^*)$ is the only pair $(P,\tilde L)$ with
$P$ in~$K\{Y,Z\}$ homogeneous of degree $2$
and~$\tilde L\in K[\der]$  such that  $$P(Y,Z)'\  =\  L(Y)Z-\tilde L(Z)Y.$$
\end{lemma}
\begin{proof}
If $(P,\tilde L)$ is such a pair, then $P_1:=P_L-P$, $L_1:= L^*-\tilde L$ gives 
$P_1(Y,Z)'=-L_1(Z)Y$,
and from this one can derive $L_1=0$ and then $P_1=0$.
\end{proof}

\noindent
Let now   $L\in K[\der]^{\neq}$ be  of order $n$, and set~$V:=\ker L$, $W:=\ker L^*$. Then for~$y\in V$,~$z\in W$ we have $[y,z]_L\in C$ by~\eqref{eq:Lagrange}; thus  $[\ ,\,]_L$ restricts to a $C$-bilinear map $V\times W\to C$, also denoted by~$[\ ,\,]_L$.

\begin{cor}\label{cor:nondeg pairing}
Suppose   $\dim_C V=n$. Then the pairing $$[\ ,\,]_L\ \colon\ V\times W\to C$$ is non-degenerate.
\end{cor}
\begin{proof}
By Lemma~\ref{lem:dim ker A^*}  we have $\dim_C W=n$.
Let $y\in V^{\neq}$.
Then $P_L(y,Z)\in K\{Z\}$ is homogeneous of degree~$1$ and order~$n-1$, hence cannot vanish on the $C$-linear subspace~$W$ of~$K$
of dimension~$n$~[ADH, 4.1.14]. Similarly with $z\in W^{\neq}$, $P_L(Y,z)\in K\{Y\}$, $V$ in place of $y$, $P_L(y,Z)$, $W$, respectively.
\end{proof}

%\begin{remark}
%Let $y_1,\dots,y_n$ be a basis of the $C$-linear space $V$. One can show (see \cite[Lemmas~4,~5]{Zettl}) that then~$y_n^*,\dots,y_n^*$ where
%$$y_n^*:=\operatorname{wr}(y_1,\dots,y_{j-1},y_{j+1},\dots,y_n)\big/\operatorname{wr}(y_1,\dots,y_n)\quad\text{for $j=1,\dots,n$}$$
%is a basis of $W$ with $[y_i,y_j^*]_L=\delta_{ij}$ for $i,j=1,\dots,n$.
%\end{remark}

\subsection*{The concomitant of operators  which split over $K$}
In this subsection we assume that $A\in K[\der]$ and $a_1,\dots,a_r\in K$ satisfy
$$A=(\der-a_r)\cdots(\der-a_1),\quad\text{so}\quad A^*=(-1)^r(\der+a_1)\cdots(\der+a_r).$$
For $i=0,\dots,r$ we  define
\begin{equation}\label{eq:AiBi}
A_{i} := (\der-a_{i})\cdots (\der-a_1),\qquad B_{i} := (-1)^{r-i}(\der+a_{i+1})\cdots(\der+a_r).
\end{equation}
Thus $A_{i}$ has order $i$, $B_{i}$ has order $r-i$, and
$$A_0=B_r=1,\quad A_r=A,\quad B_0=A^*.$$
We then have the following formula for $P_A$:

\begin{lemma}\label{lem:P_A split}
$P_A(Y,Z) = \displaystyle \sum_{i=0}^{r-1} A_{i}(Y) B_{i+1}(Z)$.
\end{lemma}

\noindent
Towards proving this, take $b_1,\dots,b_r\neq 0$ in a differential field extension $E$ of $K$ with $b_j^\dagger=a_j-a_{j-1}$ for $j=1,\dots,r$, where $a_0:=0$, and
set $b_{r+1}:=(b_1\cdots b_r)^{-1}$.
Lemma~\ref{lem:Polya fact} then gives
$$A\ =\ b_{r+1}^{-1}(\der b_r^{-1}) \cdots (\der b_2^{-1})( \der b_1^{-1}).$$
For $i=0,\dots,r$, 
$$L_{i}\ :=\  b_{i+1}^{-1}(\der b_{i}^{-1})  \cdots(\der b_2^{-1})( \der b_1^{-1})\in E[\der]$$
has order $i$, with
$L_0=b_1^{-1}$, $L_r=A$. Likewise we introduce for $i=0,\dots,r$, 
$$M_i\ :=\ (-1)^{r-i}b_{i+1}^{-1}(\der b_{i+2}^{-1})  \cdots (\der b_{r}^{-1})(\der b_{r+1}^{-1})\in E[\der]$$
of order $r-i$, so 
$M_0=A^*$, $M_r=b_{r+1}^{-1}$.
Note that 
\begin{equation}\label{eq:LM recursion}
L_{i+1}\ =\ b_{i+2}^{-1}\der L_{i},\quad M_{i}=-b_{i+1}^{-1}\der M_{i+1}\qquad (i=0,\dots,r-1).
\end{equation}
With these notations, we have:
 
\begin{lemma}[Darboux]
$P_A(Y,Z)\ =\ \displaystyle\sum_{i=0}^{r-1} L_{i}(Y) M_{i+1}(Z)$ in $E\{Y,Z\}$. 
\end{lemma}
\begin{proof} The cases $r=0$ and $r=1$ are easy to verify directly. Assume $r\ge 2$.
It suffices  to show that the    differential polynomial 
$P(Y,Z)$ on the right-hand side of the claimed equality satisfies~$P(Y,Z)'=A(Y)Z-YA^*(Z)$. 
From~ \eqref{eq:LM recursion} we obtain
$$A(Y)Z\ =\ L_{r-1}(Y)'M_r(Z),\qquad - YA^*(Z)\ =\  -YM_0(Z)\ =\  L_{0}(Y)M_1(Z)'$$ and
\begin{align*}
\big( L_{i}(Y)M_{i+1}(Z) \big){}'\  &=\  L_{i}(Y)'M_{i+1}(Z) + L_{i}(Y)M_{i+1}(Z)' \\
&=\  L_{i}(Y)M_{i+1}(Z)' -  L_{i+1}(Y)M_{i+2}(Z)'\quad\text{for $i=0,\dots,r-2$.}
\end{align*}
Now use the cancellations in $\sum_{i=0}^{r-1} \big(L_{i}(Y) M_{i+1}(Z)\big){}'$. 
\end{proof}

\noindent
This yields Lemma~\ref{lem:P_A split}: For $i=0,\dots,r-1$, we have
$$L_{i}\ =\  (b_{i+1}b_{i}\cdots b_1)^{-1}(\der-a_{i})\cdots (\der-a_1)\ =\ $$
and
%$$B_i = (-1)^{r-i} (b_{i+1}b_{i+2}\cdots b_{r+1})^{-1}(\der+a_{i+1}) \cdots(\der+a_r)$$
%and thus for $i=1,\dots,r$:
$$M_{i+1}=(-1)^{r-i-1}(b_{i+2}\cdots b_{r+1})^{-1}(\der+a_{i+2})\cdots (\der+a_r)$$
and hence
$L_{i}(Y)M_{i+1}(Z) =  A_{i}(Y)B_{i+1}(Z)$. \qed

\subsection*{Self-adjoint  and skew-adjoint operators}
If $A$ is  self-adjoint, then $r=\order A$  is even.   % \cite[\S{}17.5]{Kamke}, \cite[\S{}10]{Poole}, 
Moreover,  for $B\ne 0$:  $A$ is self-adjoint iff~$B^*AB$ is self-adjoint. The self-adjoint operators form a $C$-linear subspace of~$K[\der]$ containing $K$.

\begin{lemma}[Jacobi]\label{lem:Jacobi, 1}
Let $s\in\N$ and suppose $r=2s$. Then $A$ is self-adjoint iff  there are $b_0,\dots,b_s\in K$ such that
$$A\ =\ \der^s b_s \der^s + \der^{s-1} b_{s-1} \der^{s-1} + \cdots + b_0.$$
\end{lemma}
\begin{proof}
If $A$ has the  displayed shape, then evidently $A$ is self-adjoint. We show the converse by induction on $s$.
The case $s=0$ being trivial, suppose~$s\geq 1$. Say~$A=a_r\der^r+\text{lower order terms}$ ($a_r\in K^\times$). Then
$B=A-\der^s a_r \der^s$ is self-adjoint of order~$<r$, hence the inductive hypothesis applies to $B$.
\end{proof}

\begin{example}
If $r=2$, then $A$ is self-adjoint iff $A=a\der^2+a'\der+b$ ($a,b\in K$). In particular, $\der^2+b$ ($b\in K$) is self-adjoint.
\end{example}

\noindent
If $A$ is self-adjoint, then $[y,z]_A=-[z,y]_A$ for all $y,z\in K$, by Lemma~\ref{lem:Hesse}. Thus~$[y,y]_A=0$ for $y\in K$. This fact is used in the proof of the next lemma:

\begin{lemma}\label{lem:Jacobi, 2}
Suppose $A$ is self-adjoint and splits over $K$, and $r=2s,\ s\in\N$.
Then there are~$a\in K^\times$ and $a_1,\dots,a_s\in K$ such that
$$A\ =\ (\der+a_1)\cdots(\der+a_s)a(\der-a_s)\cdots(\der-a_1).$$
If $A$ is monic, then $a=1$ for any such $a$.
\end{lemma}
\begin{proof}
By induction on $s$. The case $s=0$ being trivial, suppose $s\geq 1$. Let~$z\neq 0$ be a zero of $A$ in a differential field extension
$\Omega$ of $K$ with $a_1:=z^\dagger\in K$.
The differential polynomial $P(Y):=P_A(Y,z)$ is homogeneous of degree $1$ and order~$r-1$
with~$P(z)=[z,z]_A=0$; hence  by [ADH, 5.1.21] we obtain $A_0\in K[\der]$ with $L_P = A_0 (\der-a_1)$.
By~\eqref{eq:Lagrange} we have~$zA=\der L_P=\der A_0(\der-a_1)$ and so $$A=z^{-1}\der A_0(\der-a_1)=(\der+a_1)A_1(\der-a_1)\qquad\text{where~$A_1:=z^{-1}A_0\in \Omega[\der]$.}$$ 
The inductive hypothesis applies to $A_1$:  $A_1\in K[\der]$ by [ADH, 5.1.11],   $A_1$ is  self-adjoint of order $r-2$, and $A_1$ splits over $K$ by [ADH, 5.1.22].
\end{proof}

\noindent
This gives rise to the following corollary:

\begin{cor}[Frobenius, Jacobi] \label{cor:Jacobi}
Suppose $A$ is self-adjoint and~$\dim_C\ker A=r=2s$. Then $A=B^*bB$  where $B=\der b_s^{-1}\cdots \der b_1^{-1}$ with $b,b_1,\dots,b_s\in K^\times$. 
\end{cor}
\begin{proof} From 
$\dim_C\ker A=r$ we obtain that $A$ splits over $K$. 
Hence the previous lemma yields
 $a_1,\dots,a_s\in K$, $a\in K^\times$ such that
$$A=(\der+a_1)\cdots(\der+a_s)a(\der-a_s)\cdots(\der-a_1),$$
and $a_1,\dots,a_s\in K^\dagger$ by Lemma~\ref{lem:split evs}.
Lemma~\ref{lem:Polya fact} yields $b_1,\dots,b_s\in K^\times$ with
$$(\der-a_s)\cdots(\der-a_1)=b_1\cdots b_s \der b_s^{-1}\cdots \der b_1^{-1}.$$
Set $B:=\der b_s^{-1}\cdots \der b_1^{-1}$. Then
$A=B^*bB$ for $b:=(-1)^s(b_1\cdots b_s)^2 a$.
\end{proof}

\noindent
Recall that $A$ is called skew-adjoint if $A^*=-A$ (and then $r=\order A$ is odd).
The skew-adjoint operators form a $C$-linear subspace of $K[\der]$. 
For $B\ne 0$, the operator~$B^*AB$ ($B\neq 0$) is skew-adjoint iff $A$ is skew-adjoint.
%As in [ADH, p.~244] we set 
%$$[L_1,L_2]:=L_1L_2-L_2L_1\quad\text{ for $L_1,L_2\in K[\der]$.}$$ 
%Recall that $\order [L_1,L_2]<\order L_1+\order L_2$ for nonzero $L_1$, $L_2$. We also have~$[L_1,L_2]^*=-[L_1^*,L_2^*]$, hence if $L_1$, $L_2$ are both self-adjoint or both skew-adjoint, then  $[L_1,L_2]$ is skew-adjoint.
We have a characterization of skew-adjoint operators analogous to     Lemma~\ref{lem:Jacobi, 1}:

\begin{lemma}\label{lem:Jacobi, skew-adjoint}
Let $s\in\N$ and suppose $r=2s+1$. Then $A$ is skew-adjoint iff  there are $b_0,\dots,b_s\in K$ such that
$$A=(\der^{s+1} b_s \der^{s} + \der^{s} b_{s} \der^{s+1}) + 
(\der^{s} b_{s-1} \der^{s-1} + \der^{s-1} b_{s-1} \der^{s}) + \cdots +
(\der b_0 + b_0\der).$$
\end{lemma}
\begin{proof}
Suppose $A$ is skew-adjoint.
Say $A=a_r\der^r+\text{lower order terms}$ ($a_r\in K^\times$), and set $b_s:=a_r/2$. Then
$A-(\der^{s+1} b_s \der^{s} + \der^{s} b_{s} \der^{s+1})$ is skew-adjoint of order~$<r$.
Hence the forward direction follows by induction on $s$. The converse is obvious.
\end{proof}

\begin{example}
If $r=1$, then $A$ is skew-adjoint iff $A=a\der+(a'/2)$ ($a\in K^\times$).
\end{example}

\noindent
For monic $A$ of order $3$, $A$ is skew-adjoint iff $A=\der^3+f\der+(f'/2)$ for some~$f\in K$.
In the next lemma we consider this case;  
for a more general version of this lemma, see \cite[Proposition~4.26(1)]{vdPS}.

\begin{lemma} \label{lem:Appell}
Let  $f\in K$, $A=\der^3+f\der+(f'/2)$, $B=4\der^2+f$ and $y,z\in \ker B$. Then $yz\in \ker A$. 
Moreover, if~$y$,~$z$ is a basis of the $C$-linear space $\ker B$, then~$y^2, yz, z^2$ is a basis of $\ker A$.
\end{lemma}
\begin{proof}
We have
$$(yz)'\ =\ y'z+yz',\quad (yz)''\ =\ y''z+2y'z'+yz''\ =\ 2y'z'-(f/2)yz,$$
hence
\begin{align*}
(yz)'''\ 	&=\  2y''z'+2y'z''-(f'/2)yz-(f/2)(yz)' \\
		&=\  -(f/2)(yz'+y'z)-(f'/2)yz-(f/2)(yz)' \\
		&=\  -f(yz)'-(f'/2)yz,
\end{align*}
and so $yz\in \ker A$.
Suppose  $ay^2+byz+cz^2=0$ for some $a,b,c\in C$, not all zero;
we claim that then $y$, $z$ are  $C$-linearly dependent.  We have
  $a\neq 0$ or $c\neq 0$, and so we may assume $a\neq 0$, $z\ne 0$. 
Then  $u:=y/z$ satisfies~$au^2+bu+c=0$, so~$u\in C$~[ADH, 4.1.1], hence~$y\in Cz$. 
\end{proof}

\noindent
If $A$ is skew-adjoint, then $P_A(Y,Z)=P_A(Z,Y)$, so 
$$P_A(Y+Z, Y+Z)\ =\  P_A(Y,Y) + P_A(Z, Z) + 2P_A(Y,Z)$$
and $[y,z]_A=[z,y]_A$ for all $y,z\in K$. 

\begin{lemma}\label{lem:char skew-adjoint}
Suppose $r\geq 1$. Then  the following are equivalent:
\begin{enumerate}
\item[\textup{(i)}] $A$ is skew-adjoint;
\item[\textup{(ii)}] there is a homogeneous differential polynomial $Q\in K\{Y\}$ of degree $2$  such that
$A(Y)Y = Q(Y)'$.
\end{enumerate}
Moreover, the differential polynomial $Q$   in \textup{(ii)} is unique, and 
$Q(Y)=\textstyle\frac{1}{2}P_A(Y,Y)$.
%$$Q(Y)=\textstyle\frac{1}{2}P_A(Y,Y)=a_rYY^{(r-1)}+\text{lower order terms},$$ 
%where $a_r\in K^\times$ is such that $A=a_r\der^r+\text{lower order terms}$.
\end{lemma}
\begin{proof}
For (i)~$\Rightarrow$~(ii) take $Q(Y):=\frac{1}{2}P_A(Y,Y)$.
For the converse let $Q$ be as in~(ii). Let $Z$ be a differential indeterminate over $K$ different from $Y$ and $c$ be a 
constant in a differential field extension $\Omega$, with $c$ transcendental over $C$. Then 
$$A(Y+cZ)(Y+cZ)\  =\ Q(Y+cZ)' \ \text{ in }\Omega\{Y,Z\}.$$
Also, in $\Omega\{Y,Z\}$,
$$A(Y+cZ)(Y+cZ)\	=\  A(Y)Y+c\big(A(Y)Z+A(Z)Y\big)+c^2A(Z)Z.$$
Take   $P,R\in K\{Y,Z\}$  such that
$$Q(Y+cZ)\			=\  Q(Y)+cP(Y,Z)+c^2R(Y,Z).$$ 
Comparing the coefficients of~$c$ now yields
$$A(Y)Z+A(Z)Y\ =\  P(Y,Z)'.$$
Now using Lemma~\ref{lem:Frobenius} gives $A^*=-A$, $P=P_A$, proving (i).
%Now let $Q$ be as in~(ii). Then $\order Q=r-1$;  take $b_1,\dots,b_r\in K$  such that 
%$$Q=Y^{(r-1)}\big(b_1Y^{(r-1)}+b_2Y^{(r-2)}+\cdots+b_r Y\big)+\text{terms of order~$<r-1$.}$$
%Then
%$$Q'=Y^{(r)}\big(2b_1Y^{(r-1)}+b_2Y^{(r-2)}+\cdots+b_{r-1}Y'+b_rY\big)+\text{terms of order~$<r$.}$$
%Comparing the coefficients of $Y^{(r)}Y^{(j)}$ on both sides of the equation $A(Y)Y=Q(Y)'$ yields
%$b_1=b_2=\cdots=b_{r-1}=0$, $b_r=a_r$, so $Q$ has the shape as claimed. 
\end{proof}

\noindent
{\em In the rest of this subsection $A$ is skew-adjoint, $r\geq 3$, and $Q(Y):=\frac{1}{2}P_A(Y,Y)$}.

\begin{lemma}\label{lem:zero of Q, 1}
If $\dim_C\ker A \ge 2$ and $C^\times$ is $2$-divisible,
then $A(z)=Q(z)=0$ for some $z\in K^\times$.
\end{lemma}
\begin{proof} 
 Apply~\cite[Chapter~XV, Theorem~3.1]{Lang} to the symmetric bilinear form 
 $$(y,z)\mapsto [y ,z]_A\ : \ \ker A\times \ker A \to C$$ on the $C$-linear space $\ker A$.
 \end{proof}

\begin{lemma}\label{lem:zero of Q, 2}
Suppose $K^\dagger$ is $2$-divisible, and $z\ne 0$ lies in a differential field extension  of $K$
with $A(z)=0$ and $z^\dagger \in K\setminus K^\dagger$.
Then $Q(z)=0$.
\end{lemma}
\begin{proof} From $z'\in Kz$ it follows by induction that $(Kz)^{(n)}\subseteq Kz$ for all~$n$.
Using~\eqref{eq:P_L} this yields $Q(z) = az^2$ for a certain $a\in K$. 
Also $Q(z)'=A(z)z=0$ and so if $a\neq 0$, then $z^\dagger = - \frac{1}{2}a^\dagger\in K^\dagger$, a contradiction.
\end{proof}

\noindent
Let $z\ne 0$ lie in a differential field extension $\Omega$ of $K$ with
 $A(z)=Q(z)=0$.
%\begin{lemma} \marginpar{checked but superfluous}
%For all $y\in K$ we have $\der\big([y,z]_A\big)=A(y)z$.
%\end{lemma}
%\begin{proof} Substitute $y,z$ for $Y,Z$ in $P_A(Y,Z)'=A(Y)Z+ A(Z)Y$.
%\end{proof}
%By \eqref{eq:qf} we have
%$$[y,z]_A = Q(y+z)-Q(y)-Q(z)=Q(y+z)-Q(y)$$
%and hence
%$$\der\big([y,z]_A\big) = (y+z)A(y+z)-yA(y)=(y+z)A(y)-yA(y)=zA(y)$$
%as required.
%\end{proof}
%\noindent
%Thus $\der\big([yz,z]_A\big)=A(yz)z$ for all $y\in K$.
The differential polynomial $P(Y):=P_A(Yz,z)\in \Omega\{Y\}$ 
is homogeneous of degree~$1$ and order $r-1$. Substitution in the identity
$P_A(Y,Z)'=A(Y)Z+ A(Z)Y$ gives~$P(Y)'=zA(Yz)$. The coefficient of $Y$ in $P$ is~$P(1)=P_A(z,z)=0$, hence 
$$P(Y)\ =\ A_0(Y'), \qquad A_0\in \Omega[\der] \text{ of order }r-2.$$

\begin{lemma}\label{zaz} In $\Omega[\der]$ we have $\der A_0 \der=zAz$, so $A_0$ is skew-adjoint.
\end{lemma}
\begin{proof} From $P(Y)'=zA(Yz)$ and $P(Y)=A_0(Y')$ we obtain $A_0(Y')'=zA(Yz)$. In terms of operators this means
$\der A_0 \der=zAz$. 
\end{proof}

%\begin{lemma} \marginpar{lemma and proof checked but replaced by better version above} 
%$A_0$ is skew-adjoint.
%\end{lemma}
%\begin{proof}
%The differential polynomial $R(Y):=YP(Y)-Q(Yz)\in \Omega\{Y\}$
%=YP_A(Yz,z)-\textstyle\frac{1}{2}P_A(Yz,Yz)\in \Omega\{Y\}$$
%is homogeneous of degree $2$. The above identities for $P(Y)$ and $P(Y)'$, and Lemma~\ref{lem:char skew-adjoint}(ii) give
%$R(Y)'=Y'P(Y)= Y'A_0(Y')\in \Omega\{Y'\}$.  It follows that $R(Y)=S(Y')$ with
%$S(Y)\in \Omega\{Y\}$ homogeneous of degree $2$.   Then $S(Y')'=Y'A_0(Y')$, so $S(Y)'=YA_0(Y)$, and thus $A_0$ is
%skew-adjoint by Lemma~\ref{lem:char skew-adjoint}.
%\end{proof}

\noindent
Next we  use these lemmas  to prove a skew-adjoint version of Lemma~\ref{lem:Jacobi, 2}.

\subsection*{Factorization of skew-adjoint operators}
{\it In this subsection $K$ is $1$-linearly surjective,   $K^\dagger$   and $C^\times$ are $2$-divisible, and $A$ is monic.}\/

\begin{prop}\label{prop:Darboux skew-adjoint} Suppose~$A$ is skew-adjoint and splits over $K$.  
Then there are~$a_1,\dots,a_s\in K$, where $r=2s+1$, such that
$$A\ =\ (\der+a_1)\cdots(\der+a_s)\der(\der-a_s)\cdots(\der-a_1).$$
\end{prop}

\begin{proof}
We proceed by induction on $s$. The case $s=0$ is clear (see the example following Lemma~\ref{lem:Jacobi, skew-adjoint}), so
let $s\geq 1$. With $Q$ as in the previous subsection we claim that  $A(z)=Q(z)=0$ and $z^\dagger\in K$ for some $z\neq 0$ in a differential field extension $\Omega$ of $K$.  
If $\dim_C \ker A=r$, then Lemma~\ref{lem:zero of Q, 1} yields such a $z$ in $\Omega=K$.
Otherwise, Lemma~\ref{lem:size of Sigma(A)} gives $a\in K\setminus K^\dagger$ with
$\mult_a(A)\geq 1$, which in turn gives~$z\neq 0$ in a differential field extension $\Omega$ of $K$ with $A(z)=0$ and $z^\dagger\in  a + K^\dagger$, and thus $Q(z)=0$ by Lemma~\ref{lem:zero of Q, 2}. This proves the claim. 

Let $z$ and $\Omega$ be as in the claim, set $a_1:=z^\dagger$, and let~$A_0\in \Omega[\der]$ be the skew-adjoint differential operator from the previous subsection. Then
$$A\ =\ z^{-1}\der A_0\der z^{-1}\ =\ (\der+a_1)A_1(\der-a_1)\qquad\text{where $A_1:=z^{-1}A_0z^{-1}$.}$$
By Lemma~\ref{zaz}, $A_1$ is skew-adjoint of order $r-2$. By~[ADH, 5.1.11, 5.1.22], $A_1\in K[\der]$ is monic and splits over $K$, so the inductive hypothesis applies to $A_1$.
\end{proof}

\begin{cor}[Darboux]
Suppose $A$ is  skew-adjoint with~$\dim_C \ker A=r=2s+1$. Then~$A=B^*\der B$ for some $B$. More precisely, 
there are $b_1,\dots,b_s\in K^\times$  such that~$A=B^*\der B$ for~$B:=\der b_s^{-1}\cdots \der b_1^{-1}$.
\end{cor}
\begin{proof}  Arguing as in the proof of Corollary~\ref{cor:Jacobi}, using Proposition~\ref{prop:Darboux skew-adjoint} instead of
Lemma~\ref{lem:Jacobi, 2}, gives $A= (-1)^sB^* b\der b B$ with $B=\der b_s^{-1}\cdots \der b_1^{-1}$, $b_1,\dots, b_s\in K^\times$, and $b=b_1\cdots b_s$. But $A$ is monic, so $(-1)^sb^2=1$, hence $b\in C$ and $A=B^*\der B$. 
\end{proof}

\begin{cor}\label{cor:Darboux general}
Suppose $A^*=(-1)^rA_{\ltimes a}$ with $a\in K^\times$, and $A$ splits over~$K$. Then there are $a_1,\dots,a_r\in K$ such that
$$A=(\der-a_r)\cdots(\der-a_1)\quad\text{and}\quad a_j+a_{r+1-j}=a^\dagger\text{ for $j=1,\dots,r$.}$$
\end{cor}
\begin{proof} By a remark preceding Definition~\ref{defsask} we have $B^*=(-1)^r B$ where $B:= A_{b/2}$, $b:=a^\dagger$, so
$B=A_{\ltimes d}$ with $d\in K^\times,\ d^2=a$. Suppose $r=2s$ is even. Then $B$ is self-adjoint and
Lemma~\ref{lem:Jacobi, 2} gives 
$$B\ =\ (\der+b_1)\cdots(\der + b_s)(\der-b_s) \cdots(\der-b_1)\quad\text{ with $b_1,\dots, b_s\in K$.}$$
Hence 
$$A\ =\ B_{\ltimes d^{-1}}\ =\ (\der+b_1-d^\dagger)\cdots(\der + b_s-d^\dagger)(\der-b_s-d^\dagger)\cdots(\der-b_1-d^\dagger).
$$
with the desired result for $a_j=b_j+d^\dagger$  and $a_{r+1-j}=-b_j+d^\dagger$, $j=1,\dots,s$. 
The case of odd $r=2s+1$ is handled in the same way, using Proposition~\ref{prop:Darboux skew-adjoint} instead of 
Lemma~\ref{lem:Jacobi, 2}. 
\end{proof}

\subsection*{Eigenrings of matrix differential equations}
{\it In the rest of this section $N$, $N_1$, $N_2$, $P$, range over  $n\times n$ matrices over~$K$ \textup{(}$n\geq 1$\textup{)}.}\/
Associated to the matrix differential equation ${y'=Ny}$ over $K$ we have the differential module~$M_N$ over $K$ with~$\dim_K M=n$ [ADH, 5.5].
Recall that matrix differential equations~$y'=N_1y$ and $y'=N_2y$ over $K$ are said to be {\it equivalent}\/ if~$M_{N_1}\cong M_{N_2}$. 
Let~$K^{n\times n}$ be the $C$-linear space of  all $n\times n$ matrices over~$K$,
and  consider the subspace
$$\mathcal E(N_1,N_2)\ :=\ \big\{ P:\  P'=N_2P-PN_1 \big\}$$
of $K^{n\times n}$. 
Given a differential ring extension $R$ of $K$, each $P\in \mathcal E(N_1,N_2)$ yields a 
$C_R$-linear map $y\mapsto Py\colon\operatorname{sol}_R(N_1)\to\operatorname{sol}_R(N_2)$.
By Lemma~\ref{dim of ann} and the next lemma we have $\dim_C\mathcal E(N_1,N_2)\leq n^2$:

\begin{lemma}
We have an isomorphism
$$P\mapsto \phi_P\  \colon\  \mathcal E(N_1,N_2) \to \Hom_{K[\der]}(M_{N_1},M_{N_2})$$
of $C$-linear spaces given by
$$\phi_P(y)\ =\  Py\quad\text{for $P\in\mathcal E(N_1,N_2)$ and $y\in M_{N_1}$.}$$
\end{lemma}
\begin{proof}
Let $P\in  \mathcal E(N_1,N_2)$, and define   $\phi_P\in\Hom_K(M_{N_1}, M_{N_2})$
by $\phi_P(y)=Py$. Then for $y\in M_{N_1}$ we have
$$\phi_P(\der y)\ =\ Py' - PN_1y\ =\  Py'+(P'y-N_2Py)\ =\ (Py)'-N_2Py\ =\  \der\phi_P(y),$$
hence $\phi_P\in \Hom_{K[\der]}(M_{N_1},M_{N_2})$. The rest follows easily.
\end{proof}

\noindent
One verifies easily that $\mathcal E(N):=\mathcal E(N,N)$ is a subalgebra of the $C$-algebra of $n\times n$-matrices over $K$ and that this yields an isomorphism
$$P\mapsto\phi_P\ \colon\ \mathcal E(N)\to \End_{K[\der]}(M_N)$$
of $C$-algebras. The $C$-algebra $\mathcal E(N)$ is called the {\bf eigenring} of $y'=Ny$.\index{eigenring!matrix differential equation}\index{matrix differential equation!eigenring}
We have~$1\leq\dim_C \mathcal E(N)\leq n^2$, and $C$ is algebraically closed in $K$. It follows that the minimum polynomial of
any $P\in \mathcal E(N)$ over $K$ (that is, the monic polynomial~$f(T)\in K[T]$ of least degree with $f(P)=0$)
 has degree at most $n^2$ and has its coefficients in $C$. In particular, if~$C$ is algebraically closed,
then the eigenvalues of any $P\in \mathcal E(N)$ are in~$C$. If $y'=N_1y$ and $y'=N_2y$ are equivalent, then their eigenrings are isomorphic as $C$-algebras.

\begin{cor} 
The isomorphism $P\mapsto\phi_P$ from the previous lemma restricts to a bijection
between the subset 
$$\mathcal E(N_1,N_2)^\times\ :=\  \operatorname{GL}_n(K)\cap  \mathcal E(N_1,N_2)$$
of $\mathcal E(N_1,N_2)$ and the set of isomorphisms $M_{N_1}\to M_{N_2}$.  
If $\mathcal E(N_1)=C\cdot 1$ and~$P\in \mathcal E(N_1,N_2)^\times$, then $\mathcal E(N_1,N_2)=C \cdot P$
and $\mathcal E(N_1,N_2)^\times=C^\times\cdot P$.
\end{cor}

\noindent
Hence $y'=N_1y$ and $y'=N_2y$ are equivalent iff $\mathcal E(N_1,N_2)^\times\neq\emptyset$, and
in this case~${y'=N_1y}$ is also called a {\bf gauge transform} of $y'=N_2y$.\index{matrix differential equation!gauge transform}

For $P\in\mathcal E(N_1,N_2)^\times$ and each differential ring extension $R$ of $K$ we have the isomorphism
$$y\mapsto Py\colon \operatorname{sol}_R(N_1)\to\operatorname{sol}_R(N_2)$$ of $C_R$-modules, and any fundamental matrix $F$ for $y'=N_1y$ in $R$ yields a fundamental matrix $PF$ for $y'=N_2y$ in $R$. 

We have a right action of the group $\operatorname{GL}_n(K)$ on $K^{n\times n}$   given by
$$(N,P)\mapsto P^{-1}(N):=P^{-1}NP-P^{-1}P'.$$
For each  $N$ and $P\in\operatorname{GL}_n(K)$,  we have $P\in\mathcal E(P^{-1}(N),N)^\times$, so the matrix differential equation $y'=P^{-1}(N)y$  is a gauge transform of $y'=Ny$.

\medskip
\noindent
Next we relate the eigenrings of linear differential operators introduced above with the
eigenrings of matrix differential equations over $K$. We precede this by some generalities about differential modules:
Let  $M, M_1, M_2$ be (left) $K[\der]$-modules. The {\em dual}~$M^*:=\Hom_K(M,K)$ of $M$ is then a
$K[\der]$-module, and  $\big\langle \phi, f\big\rangle := \phi(f)\in K$ for~$\phi\in M^*$, $f\in M$.  This yields the injective $K[\der]$-linear map
$$\alpha\mapsto\alpha^*\colon \Hom_K(M_2,M_1) \to \Hom_K(M_1^*, M_2^*)\quad\text{where $\alpha^*(\phi)=\phi\circ\alpha$ for $\phi\in M_1^*$,}$$
and 
$$\big\langle \alpha^*(\phi),f\big\rangle\ =\ \big\langle \phi,\alpha(f)\big\rangle\quad\text{for $\alpha\in\Hom_K(M_2,M_1)$, $\phi\in M_1^*$, $f\in M_2$.}$$
If~$M_1$,~$M_2$ are differential modules over $K$, then $\alpha\mapsto\alpha^*$  is an isomorphism.
Note that $\Hom_{K[\der]}(M_2,M_1)$ is a $C$-linear subspace of $H:=\Hom_{K}(M_2,M_1)$, with
$$\Hom_{K[\der]}(M_2,M_1)\ =\ \ker_{H}\der.$$
Hence the $K[\der]$-module morphism   $\alpha\mapsto\alpha^*$  restricts to 
a $C$-linear map $$\Hom_{K[\der]}(M_2,M_1)  \to \Hom_{K[\der]}(M_1^*, M_2^*),$$ which is bijective if $M_1$, $M_2$ are differential modules over $K$.  

Let $N$ be the companion matrix of a monic operator $A\in K[\der]$ of order $n$, and set~$M:= K[\der]/K[\der]A$, a differential module over $K$ of dimension $n$, with cyclic vector~$e:= 1+K[\der]A$, $Ae=0$, and with basis $e_0,\dots, e_{n-1}$, $e_j:=\der^j e$ for~$j=0,\dots,n-1$. Then $M^*$ has matrix $N$ with respect to the dual basis
$e_0^*,\dots, e_{n-1}^*$. Accordingly we identify $M^*$ with $M_N$ via the isomorphism $M^*\to M_N$ sending $e_{j-1}^*$  to the $j$th standard basis vector of $K^n$, for $j=1,\dots,n$.

% {\bf checked, but not needed} Moreover, $e^*:= e_{n-1}^*$ is  a cyclic vector for $M^*$ with  $A^*e^*=0$ by [ADH, 5.5.7], so
%we have the isomorphism $M^*\to K[\der]/K[\der]A^*$ sending $\der^ie^*$ to $\der^i+K[\der]A^*$ for $i=0,\dots,n-1$.]

In the following lemma $N_1$, $N_2$ are the companion matrices of monic operators~$A_1,A_2\in K[\der]$ of order $n$, respectively. Set~$M_1:=K[\der]/K[\der]A_1$, $M_2:=K[\der]/K[\der]A_2$ and identify  $M_1^*$, $M_2^*$ with $M_{N_1}$, $M_{N_2}$, as we just indicated for $M$. 
%natural way~[ADH, 5.5.8].
Let $\Phi$ be the isomorphism of $C$-linear spaces making the diagram 
$$\xymatrix@C+3em{\mathcal E(A_1,A_2) \ar[r]^{\Phi} \ar[d]_{R\mapsto\phi_R}  &\ \mathcal E(N_1, N_2) \ar[d]^{P\mapsto\phi_P} \\
\Hom_{K[\der]}(M_2,M_1) \ar[r]^{\alpha\mapsto\alpha^*}  &\ \Hom_{K[\der]}(M_{N_1},M_{N_2})}
$$
commute. 

\begin{lemma}
Let $R=r_0+r_1\der+\cdots+r_{n-1}\der^{n-1}\in \mathcal E(A_1,A_2)$ \textup{(}$r_0,\dots,r_{n-1}\in K$\textup{)}; then the first row of the $n\times n$ matrix $\Phi(R)$ is $(r_0,r_1,\dots,r_{n-1})$.
%$$\Phi(R) = \begin{pmatrix} r_0 \ r_1 \ \cdots \ r_{n-1} \\  \ \boxed{\parbox{6em}{\vskip1em \hskip2.5em\Huge$\ast$\hskip1em\vskip0.5em}}\  \end{pmatrix} \in K^{n\times n}.$$
\end{lemma}
\begin{proof}
Set $P=\Phi(R)$; so $\phi_P=\phi_R^*$. Let $e:=1+K[\der]A_1\in M_1$, and let $e_0^*,\dots,e_{n-1}^*$ be the basis of $M_{N_1}=M_1^*$
dual to the basis $e,\der e,\dots,\der^{n-1}e$ of $M_1$.
Likewise, let~$f:=1+K[\der]A_2\in M_2$, and let $f_0^*,\dots,f_{n-1}^*$ be the basis of $M_{N_2}=M_2^*$
dual to the basis~$f,\der f,\dots,\der^{n-1}f$ of $M_2$. Then for $j=0,\dots,n-1$ we have
$\phi_P(e_j^*)\in M_2^*$, and
$$\big\langle \phi_P(e_j^*), f\big\rangle\ =\  \big\langle\phi_R^*(e_j^*),f\big\rangle\ =\  
\big\langle e_j^*, \phi_R(f)\big\rangle\ =\  r_j.$$
Hence the matrix $P$ of the $K$-linear map~$\phi_P$ with respect to the standard bases of~$M_{N_1}=K^n$ and $M_{N_2}=K^n$
%$f_0^*,\dots,f_{n-1}^*$ of $M_{N_2}$ 
has first row $(r_0,r_1,\dots,r_{n-1})$.
\end{proof}

%\marginpar{papers by Barkatou, Jacobson}

\subsection*{Self-dual matrix differential equations} Recall that $N^*=-N^{\operatorname{t}}$ by [ADH, 5.5.6] and 
$M_{N^*} \cong (M_N)^*$ by [ADH, p.279]. 
The {\it adjoint equation}\/ of $y'=Ny$ is the matrix differential equation~${y'=N^*y}$ over~$K$.\index{adjoint!matrix differential equation}\index{matrix differential equation!adjoint}
We say that~${y'=Ny}$ is {\bf self-dual} if it is equivalent to its adjoint equation [ADH, p.~277].\index{matrix differential equation!self-dual}\index{self-dual!matrix differential equation} Hence~$y'=Ny$ is self-dual iff the differential module $M_N$ over $K$ is self-dual.
Thus if~$y'=Ny$ is self-dual, then so is any matrix differential equation over~$K$ equivalent to $y'=Ny$,   
as is the adjoint equation~$y'=N^*y$ of~$y'=Ny$. 
By  [ADH, 5.5.8, 5.5.9]  we have: 

\begin{cor}\label{cor:self-dual equivalence}
If $C\neq K$, then  every self-dual matrix differential equation~$y'=N_1y$  over $K$ is equivalent to a matrix differential equation $y'=N_2y$ with $N_2$ the companion matrix of a monic self-dual operator in~$K[\der]$. 
\end{cor}

\noindent
We set $\mult_\alpha(N):=\mult_\alpha(M_N)$ and call
$$\Sigma(N)\ :=\ \Sigma(M_N)\ =\ \big\{ \alpha: \mult_\alpha(N)\geq 1\big\}$$  
the {\bf spectrum} of  $y'=Ny$. The elements of $\Sigma(N)$ are the {\bf eigenvalues} of~$y'=Ny$.\index{matrix differential equation!spectrum}\index{spectrum!matrix differential equation}\index{eigenvalue!matrix differential equation}\index{matrix differential equation!eigenvalue}

\begin{lemma}\label{lem:matrix diff equs vs ops}
Suppose $B\in K[\der]$ is monic of order $n$ and $N$ is the companion matrix of $B$.  Then~$\mult_\alpha(B)=\mult_{\alpha}(N)$
for all $\alpha$. In particular, $\alpha$ is an eigenvalue of $B$ iff~$\alpha$ is an eigenvalue of $y'=Ny$.
\end{lemma}
\begin{proof}
Use Corollary~\ref{cor:kerder} and $M_N\cong M^*$ for $M:=K[\der]/K[\der]B$   [ADH, 5.5.8].
\end{proof}

%\begin{lemma}\label{lem:matrix diff equs vs ops}
%Suppose $N$ is the companion matrix of some monic $B\in K[\der]$. Then~$\mult_\alpha(B)=\mult_{-\alpha}(N)$
%for each $\alpha$; in particular, $\alpha$ is an eigenvalue of $B$ iff~$-\alpha$ is an eigenvalue of $y'=Ny$.
%\end{lemma}
%\begin{proof}
%Use Corollary~\ref{cor:kerder} and $M_N\cong M^*$ for $M:=K[\der]/K[\der]B$   [ADH, 5.5.8].
%\end{proof}

\noindent
From Corollary~\ref{cor:self-dual d-module} we  obtain:

\begin{cor}\label{cor:self-dual d-equ}
Assume that $y'=Ny$ is self-dual. Suppose in addition that $\sum_\alpha \mult_\alpha(N)=n$ and $K$ is $1$-linearly surjective, or 
$K$ is ${(n-1)}$-linearly surjective.
Then $\mult_\alpha(N)=\mult_{-\alpha}(N)$ for all~$\alpha$.
Hence, if also $K^\dagger$ is $2$-divisible and~$\sum_\alpha \mult_\alpha(N)$ is odd, then~$0\in\Sigma(N)$.
\end{cor}

\noindent
Note that
$$\mathcal E(N,N^*)\ =\  \big\{ P :P'=N^*P-PN \big\},\qquad \mathcal E(N,N^*)^\times\ =\ \operatorname{GL}_n(K)\cap \mathcal E(N,N^*)$$
are both closed under matrix transpose. The matrix differential equation~${y'=Ny}$ is self-dual iff $\mathcal E(N,N^*)^\times\neq\emptyset$.  Moreover, there is a $(-1)^n$-symmetric
non-degenerate $\der$-compatible $K$-bilinear form on $M_N$ iff $\mathcal E(N,N^*)^\times$ contains a matrix~$P$
such that~$P^{\operatorname{t}}=(-1)^nP$.
One  calls~$y'=Ny$ {\bf self-adjoint} if~$N^*=N$, that is, $N$ is skew-symmetric (in which case,
${y'=Ny}$ is self-dual) .\index{matrix differential equation!self-adjoint}\index{self-adjoint!matrix differential equation} 

\begin{exampleNumbered}\label{ex:Frenet-Serret}
Suppose $n=3m$ and
$$N=\begin{pmatrix} 
  			& \kappa I	&   	  \\
-\kappa I	&  			& \tau I  \\
 			& -\tau I	&  
\end{pmatrix}$$
where   $I$ denotes the $m\times m$ identity matrix
and $\kappa,\tau\in K$. Then $y'=Ny$ is self-adjoint.  
Let $\pi$ be the permutation of $\{1,\dots,n\}$ given for $i=1,\dots,m$ by
$$\pi(i)\ =\ 3i-2,\quad \pi(m+i)\ =\ 3i-1, \quad \pi(2m+i)\ =\ 3i.$$
Then  $P\in \operatorname{GL}_n(K)$ with $Pe_j=e_{\pi(j)}$ ($j=1,\dots,n$) gives $P'=0\in K^{n\times n}$, so 
$$P^{-1}(N)\  =\ \operatorname{diag}(T,\dots,T)\in K^{n\times n}
\quad\text{where $T:=\begin{pmatrix}0	& \kappa	& 0 \\ -\kappa	& 0	&	\tau \\	0	& -\tau	& 0\end{pmatrix}\in K^{3\times 3}$.}$$
%Earlier description of how $N$ and $D:=\operatorname{diag}(T,\dots,T)$ relate is incorrect.  The permutation involved is not an involution for $m\ge 2$. 
By Corollary~\ref{cormultsum},
$\mult_\alpha(N)=m\mult_\alpha(T)$ for all $\alpha$, so
$\Sigma(N)=\Sigma(T)$. If~$F$ is a fundamental matrix for $y'=Ty$, then $G:=\operatorname{diag}(F,\dots,F)\in \operatorname{GL}_n(K)$ is a fundamental matrix for $y'=P^{-1}(N)y$, that is, $G'=P^{-1}NPG$, so $PG$ is a fundamental matrix for~$y'=Ny$.  Suppose now that $K$ is $1$-linearly surjective, $K^\dagger$ is $2$-divisible, and~$\sum_\alpha \mult_\alpha(T) = 3$.
Then $\sum_\alpha \mult_\alpha(N) = n$ and~$\mult_\alpha(T)=\mult_{-\alpha}(T)$ for all~$\alpha$, so $\Sigma(N)=\Sigma(T)=\{\alpha,-\alpha,0\}$ for some $\alpha$.
\end{exampleNumbered}

\begin{lemma}\label{lem:constant norm}
Suppose $y'=Ny$ is self-adjoint and let $y,z\in  \operatorname{sol}(N)$, where $y=(y_1,\dots,y_n)^{\operatorname{t}}$ and $z=(z_1,\dots,z_n)^{\operatorname{t}}$. Then $y_1z_1+\cdots+y_nz_n\in C$.
\end{lemma}
\begin{proof}
With $\langle\, \cdot\, , \,\cdot\, \rangle$ denoting the usual inner product on $K^n$, we have
$$\langle y,z\rangle'\  =\ \langle y',z\rangle + \langle y,z'\rangle\ =\ 
\langle Ny,z\rangle + \langle y,Nz\rangle\ =\ 
\langle y,N^{\operatorname{t}}z \rangle +  \langle y,Nz\rangle\ =\ 
0$$
since $N^{\operatorname{t}}=-N$.
\end{proof}

\noindent
Thus if  $y'=Ny$ is self-adjoint, then we have a symmetric bilinear form
$$(y,z)\mapsto \langle y,z\rangle = y_1z_1+\cdots+y_nz_n \quad (y=(y_1,\dots,y_n)^{\operatorname{t}},\  z=(z_1,\dots,z_n)^{\operatorname{t}})$$
on the $C$-linear subspace $\operatorname{sol}(N)$ of $K^n$ of dimension~$\leq n$.

\medskip\noindent
 A matrix $F\in K^{n\times n}$ is said to be {\em orthogonal\/} if $FF^{\operatorname{t}}=I_n$, where $I_n$ denotes the identity of the ring $K^{n\times n}$ of $n\times n$-matrices over $K$. This yields the subgroup $\operatorname{O}_n(K)$ of 
 $\operatorname{GL}_n(K)$ consisting of the orthogonal matrices $F\in K^{n\times n}$. 

Suppose $F\in\operatorname{GL}_n(K)$ is a  fundamental matrix  for $y'=Ny$. By [ADH, 5.5.12]  this yields a
a fundamental matrix $(F^{\operatorname{t}})^{-1}\in\operatorname{GL}_n(K)$ for~$y'=N^*y$, so
if $F$ is orthogonal, then $y'=Ny$ is self-adjoint. Conversely:  

\begin{lemma}\label{lem:orthogonal fund matrix}
Suppose $y'=Ny$ is self-adjoint, $\dim_C \operatorname{sol}(N)=n$, and $C^\times$ is $2$-divisible. Then $\operatorname{GL}_n(K)$ contains an orthogonal fundamental matrix for~$y'=Ny$.
\end{lemma}
\begin{proof} Take a fundamental matrix $F\in \operatorname{GL}_n(K)$ for $y'=Ny$. Then
$G:=(F^{\operatorname{t}})^{-1}$ is also a fundamental matrix for $y'=Ny$, so
$F^{\operatorname{t}}F=G^{-1}F\in\operatorname{GL}_n(C)$ by~[ADH, 5.5.11]. 
Now the matrix~$F^{\operatorname{t}}F$ is symmetric,  so  \cite[Chapter~XV, Theorem~3.1]{Lang} gives  $D,U$ in $\operatorname{GL}_n(C)$ with diagonal $D$ such that $F^{\operatorname{t}}F=U^{\operatorname{t}}DU$.
Let~$V:=\sqrt{D}U$ where~$\sqrt{D}$ in~$C^{n\times n}$ is diagonal with~$(\sqrt{D})^2=D$. 
Then $F^{\operatorname{t}}F=V^{\operatorname{t}}V$ and so~$FV^{-1}\in \operatorname{GL}_n(K)$ is an orthogonal fundamental matrix for $y'=Ny$ by [ADH, 5.5.11]. 
\end{proof}

\noindent
The skew-symmetric $n\times n$ matrices over $K$ form a Lie subalgebra 
$$ \mathfrak{so}_n(K)\  =\ \{N :\  N^* = N \}$$
of  $K^{n\times n}$ equipped with the Lie bracket $[N_1,N_2]=N_1N_2-N_2N_1$.
Suppose now~$n=2m$ is even, and set $J:=\left(\begin{smallmatrix}   & I_m \\ -I_m &  \end{smallmatrix}\right)$.
Then $J^{\operatorname{t}}=J^{-1}=-J$, and
$$\mathfrak{sp}_n(K)\ =\ \{N :\  N^*  J = JN \}$$
is also a Lie subalgebra of $K^{n\times n}$. The matrices in  $\mathfrak{sp}_n(K)$ are called {\it hamiltonian}\/;
thus $N$ is hamiltonian iff $JN$ is symmetric.  
We say that the matrix differential equation $y'=Ny$ is {\bf hamiltonian} if~$N\in\mathfrak{sp}_n(K)$; in that case~$J\in\mathcal E(N,N^*)^\times$, so~$y'=Ny$ is self-dual.\index{matrix differential equation!hamiltonian}
A matrix~$F\in K^{n\times n}$ is said to be {\it symplectic}\/ if $F^{\operatorname{t}}JF=J$.
The symplectic matrices form a subgroup~$\operatorname{Sp}_n(K)$ of~$\operatorname{GL}_{n}(K)$. 
If $y'=Ny$ has a fundamental matrix $F\in\operatorname{Sp}_n(K)$, then~$y'=Ny$ is hamiltonian. 
In analogy with Lemma~\ref{lem:orthogonal fund matrix} we have a converse:

\begin{lemma}\label{lem:symplectic fund matrix}
Suppose $y'=Ny$ is hamiltonian and  $\dim_C \operatorname{sol}(N)=n$. Then $\operatorname{GL}_n(K)$ contains a symplectic fundamental matrix for $y'=Ny$.
\end{lemma}
\begin{proof} Take a fundamental matrix $F\in \operatorname{GL}_n(K)$ for $y'=Ny$. Then
$JF$ and $G:=(F^{\operatorname{t}})^{-1}$ are fundamental matrices for $y'=N^*y$, so~$F^{\operatorname{t}}JF=G^{-1}JF\in\operatorname{GL}_n(C)$. Now $F^{\operatorname{t}}JF$ is skew-symmetric, 
hence $F^{\operatorname{t}}JF=U^{\operatorname{t}}JU$ with~$U\in\operatorname{Gl}_n(C)$~\cite[Chap\-ter~XV, Corollary~8.2]{Lang}; then~$FU^{-1}$ is a symplectic fundamental matrix for~${y'=Ny}$.
\end{proof}

\noindent
For a hamiltonian analogue of Lemma~\ref{lem:constant norm}, let  $\langle\, \cdot\, , \,\cdot\, \rangle$ denote the usual inner product on $K^n$ and let
$(y,z)\mapsto   \omega(y,z):=\langle y,Jz\rangle$
be the standard symplectic bilinear form on $K^n$.

\begin{lemma}\label{lem:constant norm, symplectic}
Suppose $y'=Ny$ is hamiltonian and  $y,z\in\operatorname{sol}(N)$. Then $\omega(y,z)\in C$.
\end{lemma}
\begin{proof}
We have
$$\omega(y,z)'\ =\ \langle y',Jz\rangle + \langle y,Jz'\rangle\ =\ \langle Ny,Jz\rangle + \langle y,JNz\rangle\ =\
\langle y,N^{\operatorname{t}}Jz\rangle + \langle y,JNz\rangle\ =\ 0$$
where we used $-N^{\operatorname{t}}J=N^*J=JN$ for the last equality.
\end{proof}

\noindent
Note also that $N$ is hamiltonian iff $N=JN^{\operatorname{t}}J$. It follows that $N$ is hamiltonian iff~$N=\left(\begin{smallmatrix} Q & R \\ P & Q^*  \end{smallmatrix}\right)$ where $P,R\in K^{m\times m}$ are symmetric and $Q\in K^{m\times m}$. 

Hamiltonian matrix differential equations appear naturally in the study of more general (non-linear) Hamiltonian systems %\marginpar{\bf rest of subsection skipped} 
(as  the  variational equations along an integral curve of such a system). See, e.g., \cite{Churchill}.
They also arise from self-adjoint linear differential operators: using Lemma~\ref{lem:Jacobi, 1} one can show that
if $A$ is self-adjoint with companion matrix $M$, then with $n:=r$ there is some hamiltonian $N$ such that
$\mathcal E(M,N)^\times\neq\emptyset$; see \cite[p.~76]{Coppel}.
% To see this let $r=2s$ ($s\in\N$) and
%$$A=(-1)^s\der^s b_s \der^s + (-1)^{s-1}\der^{s-1} b_{s-1} \der^{s-1} + \cdots + b_0 \qquad (b_0,\dots,b_s\in K).$$
%(See Lemma~\ref{lem:Jacobi, 1}.) With $n:=r$, consider the hamiltonian matrix
%$N=\left(\begin{smallmatrix} Q & R \\ P & Q^* \end{smallmatrix}\right)$
%where $$P=\operatorname{diag}(b_0,\dots,b_{s-1}),\ Q=b_s^{-1}\operatorname{diag}(0,\dots,0,1),\  
%R=\left(\begin{smallmatrix} 
%0		& 1 & 0 & \cdots	& 0 \\
%0		& 0	& 1	& \cdots	& 0 \\
%\cdot 	& \cdot & \cdot			& \cdots	& \cdot \\
%0		& 0	& 0	& \cdots	& 1 \\
%0		& 0	& 0	& \cdots	& 0
%  \end{smallmatrix}\right).$$
%With $M$ denoting the companion matrix of $A$ we then have
%$$\begin{pmatrix} I_s & \\ & \end{matrix} \in \mathcal E(M,N)^\times,$$

\subsection*{Anti-self-duality}
We now continue in the setting of the subsection {\it Complex conjugation}\/ in Section~\ref{sec:splitting}.
Thus $K=H[\imag]$ where $H$ is a  differential subfield of $K$, $\imag^2=-1$, and $\imag\notin H$. The isomorphisms below are
 of differential modules over $K$. 
Let $M$ be a differential module over $K$. 
We establish some analogues of results above for the conjugate dual of $M$ instead of its dual.

Call $M$ is {\bf anti-self-dual} if $M\cong \bar{M^*}$.\index{differential module!anti-self-dual}\index{anti-self-dual!differential module}
If $M$ is anti-self-dual, then so is every isomorphic $K[\der]$-module, in particular, $\bar{M^*}$.
Here is an analogue of Corollary~\ref{cor:self-dual d-module} which follows immediately from Corollary~\ref{cor:spectrum, conj self-adjoint}:

\begin{cor}\label{cor:anti-self-dual d-module} Let $\dim_K M=r$ and assume $M$ is anti-self-dual.
Suppose also that $K$ is $1$-linearly surjective and $\sum_\alpha \mult_\alpha(M)=r$, or~${r\geq 1}$ and~$K$ is ${(r-1)}$-linearly surjective.
Then $\mult_\alpha(M)=\mult_{-\bar\alpha}(M)$ for all~$\alpha$.
Hence if   additionally~$K^\dagger$ is $2$-divisible and $\sum_\alpha \mult_\alpha(M)$ is odd, then~$[b\imag]\in\Sigma(M)$ for some~$b\in H$.
\end{cor}

\noindent
Suppose now $M=K[\der]/K[\der]A$ and $r\geq 1$. Then $\overline{M^*}\cong K[\der]/ K[\der]\overline{A^*}$ by
[ADH, 5.5.8] and Example~\ref{ex:M compl conj}. Hence $M$ is anti-self-dual iff $A$, $\bar{A^*}$ have the same type.
We say that $A$ is {\bf anti-self-dual} if $A$, $\bar{A^*}$ have the same type.\index{linear differential operator!anti-self-dual}\index{anti-self-dual!linear differential operator} 
If $A$ is anti-self-dual, then so are $\bar{A}$ and $A^*$,
and so is every operator of the same type as $A$. If $A$ is anti-self-dual, then $A$, $\bar{A^*}$ have the same eigenvalues, with the same multiplicities. The previous corollary yields:

\begin{cor} \label{cor:anti-self-dual operator}
Suppose  $A$ is anti-self-dual, and set $s:=\sum_\alpha \mult_\alpha(A)$. Also assume~$K$ is $1$-linearly surjective and~$s=r$, or 
$r\geq 1$ and~$K$ is ${(r-1)}$-linearly surjective.
Then~$\mult_\alpha(A)=\mult_{-\bar{\alpha}}(A)$ for all~$\alpha$.
Hence if in addition~$K^\dagger$ is $2$-divisible and~$s$ is odd, then~$[b\imag]\in\Sigma(A)$ for some $b\in H$.
\end{cor}

\noindent
Later $H$ is usually a Hardy field with $H^\dagger=H$, so~$\alpha=-\bar{\alpha}$ for all $\alpha$. In this case Corollaries~\ref{cor:anti-self-dual d-module} and~\ref{cor:anti-self-dual operator} are less useful than their cousins Corollaries~\ref{cor:self-dual d-module} and~\ref{cor:self-dual operator}.
Note also that  if $A\in H[\der]$, then $A$ is self-dual iff $A$ is anti-self-dual.

We now consider anti-self-duality for a matrix differential equation $y'=Ny$ over~$K$. Recall that $N$ is an $n\times n$-matrix over $K$ with $n\ge 1$, and that if $M=M_N$ is the differential module over $K$ associated to $N$, then $\bar{M}\cong M_{\bar{N}}$ by 
the remarks preceding Example~\ref{ex:M compl conj}, and $M^*\cong M_{N^*}$ by~[ADH, pp.~279--280].
We say that~$y'=Ny$  is {\bf anti-self-dual} if it is equivalent to the matrix differential equation
$y'=\bar{N^*}y$ over $K$.\index{matrix differential equation!anti-self-dual}\index{anti-self-dual!matrix differential equation} (Note: $\bar{N^*}=-\bar{N}{}^{\operatorname{t}}$.)
 Hence~$y'=Ny$ is anti-self-dual iff   $M_N$   is anti-self-dual.
If~$y'=Ny$ is anti-self-dual, then so is any matrix differential equation over~$K$ equivalent to $y'=Ny$,   
as are the matrix differential equations~$y'=N^*y$ and $y'=\bar{N}y$ over $K$. 
If $N\in H^{n\times n}$, then the matrix differential equation $y'=Ny$ over $K$ is self-dual iff it is anti-self-dual. 

\begin{cor}
Suppose $C\neq K$ and $y'=Ny$ is anti-self-dual. Then~$y'=Ny$ is equivalent to a matrix differential equation $y'=A_Ly$ with $A_L$ the companion matrix of a monic anti-self-dual $L\in K[\der]$. 
\end{cor}
\begin{proof} By [ADH, 5.5.9], $y'=Ny$ is equivalent to $y'=A_Ly$ where $A_L$ is the companion matrix of a monic $L\in K[\der]$.
Then $L$ is anti-self-dual by [ADH, 5.5.8] and Example~\ref{ex:M compl conj}.
\end{proof}

\noindent
We say that $y'=Ny$ is {\bf anti-self-adjoint} if $\bar{N^*}=N$, that is, $N^{\operatorname{t}}=-\bar{N}$.\index{matrix differential equation!anti-self-adjoint}\index{anti-self-adjoint} Then~$y'=Ny$ is in particular anti-self-dual. 
If $N\in H^{n\times n}$, then $y'=Ny$ is anti-self-adjoint iff it is self-adjoint. To state an anti-self-adjoint analogue of Lemma~\ref{lem:constant norm} we use the ``hermitian'' inner product $\langle\, \cdot\, , \,\cdot\, \rangle$ on $K^n$ given by
$\langle y,z\rangle=y_1\bar{z}_1+\cdots + y_n\bar{z}_n$ for~$y=(y_1,\dots,y_n)^{\operatorname{t}}\in K^n$ and $z=(z_1,\dots,z_n)^{\operatorname{t}}\in K^n$.

\begin{lemma}\label{lem:constant norm, anti-self-adj}
If $y'=Ny$ is anti-self-adjoint and $y, z\in \operatorname{sol}(N)$, then $\langle y,z\rangle\in C$.
\end{lemma} 
\begin{proof} Assume $y'=Ny$ is anti-self-adjoint. Then $\bar{N}{}^{\operatorname{t}}=-N$, so
$$\langle y,z\rangle'\ =\ \langle y',z\rangle + \langle y,z'\rangle\ =\ 
\langle Ny,z\rangle + \langle y,Nz\rangle\ =\ 
\langle y,\bar{N}{}^{\operatorname{t}}z \rangle +  \langle y,Nz\rangle\ =\ 
0.\qedhere$$
 \end{proof}

\noindent
Thus if $y'=Ny$ is anti-self-adjoint, then we have
a   hermitian form 
$$(y,z)\mapsto \langle y,z\rangle\ =\ y_1\bar{z_1}+\cdots+y_n\bar{z_n} \quad (y=(y_1,\dots,y_n)^{\operatorname{t}},\  z=(z_1,\dots,z_n)^{\operatorname{t}})$$
on the $C$-linear subspace $\operatorname{sol}(N)$ of $K^n$ of dimension~$\leq n$.

\medskip\noindent
A matrix $U\in K^{n\times n}$ is {\it unitary}\/ if $U{}^{\operatorname{t}}\bar{U}=I_n$, equivalently, $\<Ux, Uy\>=\<x,y\>$ for all~$x,y\in K^n$. 
The unitary matrices form a subgroup~$\operatorname{U}_n(K)$ of $\operatorname{GL}_n(K)$.
Suppose~$F\in \operatorname{GL}_n(K)$ is a fundamental matrix for~$y'=Ny$.
Then $\bar F$ is a fundamental matrix for $y'=\bar{N}y$, and so $(\bar{F}{}^{\operatorname{t}})^{-1}\in
\operatorname{GL}_n(K)$ is a fundamental matrix for $y'=\bar{N^*}y$.
So if $F$ is unitary, then
$y'=Ny$ is anti-self-adjoint. Here is a converse, analogous to Lemma~\ref{lem:orthogonal fund matrix}:

\begin{lemma}\label{lem:unitary fund matrix}
Suppose $y'=Ny$ is anti-self-adjoint, $\dim_C \operatorname{sol}(N)=n$, and  $H$ is real closed.
Then $\operatorname{GL}_n(K)$ contains a unitary fundamental matrix for $y'=Ny$.
\end{lemma}
\begin{proof}
Take a fundamental matrix $F\in\operatorname{GL}_n(K)$ for $y'=Ny$. Then $G:=(\bar{F}{}^{\operatorname{t}})^{-1}$ is also
a fundamental matrix for $y'=Ny$, so $\overline{F}{}^{\operatorname{t}}F=G^{-1}F\in\operatorname{GL}_n(C)$.
Now~$P:=\overline{F}{}^{\operatorname{t}}F$ is hermitian (i.e., $\overline{P}{}^{\operatorname{t}}=P$), so 
    \cite[Chapter~XV, \S{}5,~6]{Lang} gives  a diagonal~$D\in\operatorname{GL}_n(C_H)$ and a $U\in\operatorname{U}_n(C)$ with~$P=\bar{U}{}^{\operatorname{t}}DU$. So for $x\in C^n$ and $y:=U^{-1}x\in C^n$,
    $$\langle Dx,x\rangle\ =\ \langle DUy,Uy\rangle\ =\ \langle Py,y \rangle\ =\ \langle Fy,Fy\rangle,$$
    a sum of squares in $H$. As $C_H$ is also real closed,
    all entries of $D$ are squares in $C_H$. 
Take diagonal $E\in C_H^{n\times n}$  with $E^2=D$. Then $V:=EU\in \operatorname{GL}_n(C)$ and
~$P=\bar{V}{}^{\operatorname{t}}V$, so $FV^{-1}\in\operatorname{GL}_n(K)$ is a unitary fundamental matrix
    for $y'=Ny$.
\end{proof}

\section{Eigenvalues and Splittings}\label{sec:eigenvalues and splitting}

\noindent
{\it In this section $K$ is a differential field such that $C$ is algebraically closed and $K^\dagger$ is divisible.}\/
We let $A$, $B$ range over $K[\der]$, and we assume $A\neq 0$ and set $r:=\order A$.
%In the first subsection we   show that  then $A$ has at most~$r$ eigenvalues.

\subsection*{Spectral decomposition of  differential operators}
Fix a complement $\Lambda$ of the subspace $K^\dagger$ of the $\Q$-linear space $K$, let $\Univ:= K\big[\!\ex(\Lambda)\big]$ be the
universal exponential extension of $K$, let $\Omega$ be the differential fraction field of the differential $K$-algebra~$\Univ$, and let $\lambda$ range over $\Lambda$. Then
$$A_\lambda\ =\ A_{\ltimes\!\ex(\lambda)}\ =\ \ex({-\lambda})A\ex(\lambda)\in K[\der].$$ 
Moreover, for every $a\in K$ there is a unique $\lambda$ with $a-\lambda\in K^\dagger$, so 
$\mult_{[a]}(A)=\mult_{\lambda}(A)$. Call $\lambda$ an {\bf eigenvalue} of $A$ with respect 
to our complement~$\Lambda$ of $K^\dagger$ in~$K$ if~$[\lambda]$ is an eigenvalue of $A$; thus the group isomorphism $\lambda\mapsto [\lambda]\colon \Lambda \to K/K^\dagger$ maps the set of eigenvalues of~$A$ with respect to~$\Lambda$ onto the spectrum of $A$.  
For~$f\in\Univ$ with spectral decomposition~$(f_\lambda)$ we have \index{eigenvalue}
$$A(f)\ =\ \sum_\lambda A_\lambda(f_\lambda)\,\ex(\lambda),$$
so $A(\Univ^\times)\subseteq \Univ^\times\cup \{0\}$.
We call the family $(A_\lambda)$ the {\bf spectral decomposition} of $A$ (with respect to $\Lambda$).
Given a $C$-linear subspace $V$ of $\Univ$, we set $V_\lambda := V\cap K\ex(\lambda)$, a $C$-linear subspace  of $V$; the sum $\sum_\lambda V_\lambda$ is direct. For $V:=\Univ$ we have $\Univ_\lambda=K\ex(\lambda)$, and $\Univ=\bigoplus_\lambda \Univ_\lambda$
with $A(\Univ_\lambda)\subseteq\Univ_\lambda$ for all $\lambda$.
Taking $V:=\ker_{\Univ} A$, we obtain~$V_\lambda=(\ker_{K} A_\lambda)\ex(\lambda)$ 
and hence $\dim_C V_\lambda = \mult_\lambda(A)$,
and $V=\bigoplus_\lambda V_\lambda$. Thus  \index{spectral decomposition!of a differential operator} \index{linear differential operator!spectral decomposition}
\begin{equation}\label{eq:bd on sum mults}
\abs{\Sigma(A)}\ \leq\ \sum_\lambda \mult_\lambda(A)\ =\ \dim_C \ker_{\Univ} A\ \leq\ r.
\end{equation}
Moreover:

\begin{lemma}\label{newlembasis} The $C$-linear space $\ker_{\Univ} A$ has a basis contained in $\Univ^\times=K^\times\ex(\Lambda)$.
\end{lemma}

\begin{example} 
We have a $C$-algebra isomorphism $P(Y)\mapsto P(\der)\colon C[Y]\to C[\der]$.
Suppose~$A\in C[\der]\subseteq K[\der]$,  let $P(Y)\in C[Y]$, $P(\der)=A$, and let
$c_1,\dots,c_n\in C$ be the distinct zeros of $P$, of respective multiplicities $m_1,\dots,m_n\in\N^{\geq 1}$ (so $r=\deg P=m_1+\cdots+m_n$).
Suppose also $C\subseteq \Lambda$, and $x\in K$ satisfies $x'=1$. 
(This holds in Example~\ref{ex:Q}.)
Then the $x^i\ex(c_j)\in \Univ$ ($1\leq j\leq n$, $0\leq i<m_j$)
form a basis of the $C$-linear space $\ker_{\Univ} A$ by [ADH, 5.1.18].
So the eigenvalues of~$A$ with respect to $\Lambda$ are $c_1,\dots,c_n$, with respective multiplicities $m_1,\dots,m_n$.
\end{example}

\begin{cor}\label{cor:sum of evs}
Suppose $\dim_C \ker_{\Univ} A = r\ge 1$ and $A=\der^r+a_{r-1}\der^{r-1}+\cdots+a_0$ where~$a_0,\dots,a_{r-1}\in K$.
Then
$$\sum_{\lambda} \mult_{\lambda}(A)\lambda\  \equiv\  -a_{r-1} \mod K^\dagger.$$
In particular,  $\sum_{\lambda} \mult_{\lambda}(A)\lambda=0$ iff $a_{r-1}\in K^\dagger$.
\end{cor}
\begin{proof}
Take a basis $y_1,\dots,y_r$ of $\ker_{\Univ} A$ with $y_j=f_j\ex(\lambda_j)$, $f_j\in K^\times$, $\lambda_j\in\Lambda$. 
The Wronskian matrix $\operatorname{Wr}(y_1,\dots,y_r)$ of $(y_1,\dots,y_r)$ [ADH, p.~206] equals
$$\operatorname{Wr}(y_1,\dots,y_r)\ =\  M\begin{pmatrix} \ex(\lambda_1) & & \\ & \ddots & \\ & & \ex(\lambda_r) \end{pmatrix}\qquad\text{where $M\in\operatorname{GL}_n(K)$.}$$
Then $w:=\operatorname{wr}(y_1,\dots,y_r)=\det \operatorname{Wr}(y_1,\dots,y_r)\neq 0$ by [ADH, 4.1.13] and
 $$-a_{r-1}\ =\ w^\dagger\ =\ (\det M)^\dagger + \lambda_1+\cdots+\lambda_r$$ where we used [ADH, 4.1.17] for the first equality.
\end{proof}

\noindent
If $A$ splits over $K$, then so does $A_\lambda$.
Moreover, if $A_\lambda(K)=K$, then $A(\Univ_\lambda)=\Univ_\lambda$:
for   $f,g\in K$ with $A_\lambda(f)=g$ we have
$A\big(f\ex(\lambda)\big)=g\ex(\lambda)$.
Thus:

\begin{lemma}\label{lem:A(U)=U}
Suppose $K$ is $r$-linearly surjective, or $K$ is $1$-linearly surjective and~$A$ splits over $K$.  Then $A(\Univ_\lambda)=\Univ_\lambda$ for all $\lambda$ and hence $A(\Univ)=\Univ$.
\end{lemma}

\noindent
%If $K$ is $1$-linearly surjective and  $A$ splits over $K$, then $\sum_\lambda \mult_\lambda(A) = r$
%by the remark after \eqref{eq:mult(AB)} and Example~\ref{ex:ev order 1}. 
In the next subsection we  study the 
connection between splittings of $A$ and
bases of the $C$-linear space $\ker_{\Univ} A$ in more detail.

\subsection*{Constructing splittings and bases}
Recall that $\order A=r\in\N$. Set $\Univ=\Univ_K$, so $\Univ^\times =  K^\times\ex(\Lambda)$. Let 
$y_1,\dots,y_r\in \Univ^\times$. We construct a sequence~$A_0,\dots,A_n$  of monic operators in $K[\der]$ with~$n\leq r$ as follows.
First, set $A_0:=1$. Next, given~$A_0,\dots,A_{i-1}$ in $K[\der]^{\ne}$ ($1\leq i\leq r$), set $f_i:=A_{i-1}(y_i)$; if $f_i\ne 0$, then $f_i\in \Univ^\times$, so $f_i^\dagger\in K$, and the next term in the sequence
is
$$A_i\ :=\ (\der-a_i)A_{i-1}, \qquad a_i\ :=\ f_i^\dagger,$$ 
whereas if $f_i=0$, then $n:=i-1$ and the construction is finished.

\begin{lemma}\label{lem:distinguished splitting}
$\ker_{\Univ} A_i=Cy_1\oplus\cdots\oplus Cy_i$ \textup{(}internal direct sum\textup{)} for $i=0,\dots,n$.
\end{lemma}
\begin{proof}
By induction on $i\le n$. The case $i=0$ being trivial, suppose $1\le i\le n$ and the claim holds
for $i-1$ in place of $i$.  Then $A_{i-1}(y_i)=f_i\neq 0$, hence $y_i\notin \ker_{\Univ} A_{i-1}=Cy_1\oplus\cdots\oplus Cy_{i-1}$,
and $A_i=(\der-f_i^\dagger)A_{i-1}$, so by [ADH, 5.1.14(i)] we have
$\ker_{\Univ} A_i=\ker_{\Univ} A_{i-1}\oplus Cy_i=Cy_1\oplus\cdots\oplus Cy_i$.
\end{proof}

\noindent
We denote the tuple $(a_1,\dots,a_n)\in K^n$ just constructed by $\operatorname{split}(y_1,\dots,y_r)$,
so~$A_n=(\der-a_n)\cdots(\der-a_1)$. 
Suppose $r\geq 1$. Then $n\ge 1$, $a_1=y_1^\dagger$, $A_1=\der-a_1$,
$A_1(y_2),\dots, A_1(y_n)\in \Univ^\times$,  and we have
$$  (a_2,\dots,a_n)\ =\ \operatorname{split}\!\big(A_1(y_2),\dots,A_1(y_n)\big).$$
 By Lemma~\ref{lem:distinguished splitting},  $n\leq r$ is maximal such that $y_1,\dots,y_n$ are $C$-linearly independent. 
In particular, $y_1,\dots,y_r$ are $C$-linearly independent iff $n=r$.
%For the definition of the operator $A_{\operatorname{s}}$ used in the next corollary see Corollary~\ref{cor:As}.

 \begin{cor}\label{corbasissplit} 
If $A(y_i)=0$ for $i=1,\dots,n$, then $A\in K[\der]A_n$.
Thus if~$n=r$ and $A(y_i)=0$ for $i=1,\dots,r$, then  
$A=a(\der-a_r)\cdots(\der-a_1)$ where~$a\in K^\times$.
\end{cor}

\noindent
This follows from [ADH, 5.1.15(i)] and Lemma~\ref{lem:distinguished splitting}.

\medskip\noindent
Suppose that $H$ is a differential subfield of $K$ and $y_1^\dagger,\dots,y_r^\dagger\in H$.  Then we have~$\operatorname{split}(y_1,\dots,y_r)\in~H^n$: use that $y'\in Hy$ with $y\in U$ gives $y^{(m)}\in Hy$ for all $m$, so~$B(y)\in Hy$
for all $B\in H[\der]$, hence for such $B$, if $f:= B(y)\neq 0$, then~$f^\dagger\in H$.

%Note that if $y_1,\dots,y_r$ are contained in a differential subfield $H$ of $K$,   
%then $\operatorname{split}(y_1,\dots,y_r)\in H^n$. 

%\begin{cor}\label{corbasissplit} \marginpar{made more precise}
%Suppose  $A(y_i)=0$ for $i=1,\dots,r$. Then~${A_{\operatorname{s}}\in K[\der] A_n}$. \textup{(}In particular, if $A$ is monic and $n=r$, then $A=(\der-a_r)\cdots(\der-a_1)$.\textup{)}
%\end{cor}
%\begin{proof}
%By Lemma~\ref{lem:distinguished splitting} we have $\order A_n=n=\dim\ker_{\Univ} A_n$, thus $\ker_\Omega A_n=\ker_{\Univ} A_n\subseteq \ker_\Omega A$. So $A\in\Omega[\der]A_n\cap K[\der]=K[\der]A_n$ by [ADH, 5.1.5(i), 5.1.11] and hence ${A_{\operatorname{s}}\in K[\der] A_n}$ by  Corollary~\ref{cor:As}.
%\end{proof}

\begin{cor}\label{corbasiseigenvalues} 
Suppose $\dim_C \ker_{\Univ} A = r$. Then $\ker_{\Univ} A = \ker_{\Omega} A$ and $A$ splits over $K$. If 
$A = (\der-a_r)\cdots(\der-a_1),\ a_1,\dots,a_r\in K$,
then the spectrum of~$A$ is~$\big\{[a_1],\dots,[a_r]\big\}$, and for all $\alpha\in K/K^\dagger$,
$$\mult_{\alpha}(A)\  =\ 
\big| \big\{ i\in\{1,\dots,r\}:\ \alpha=[a_i] \big\}  \big|.$$
\end{cor} 
\begin{proof} 
$A$ splits over $K$ by
 Lemma~\ref{newlembasis} and Corollary~\ref{corbasissplit}. 
The rest follows from Lemma~\ref{lem:split evs} in view of $\sum_\lambda \mult_\lambda(A)=\dim_C \ker_{\Univ} A$.
\end{proof}

\noindent
Conversely, we can associate to a given splitting of $A$ over $K$ a basis of $\ker_{\Univ} A$ consisting of $r$ elements
of~$\Univ^\times$, provided $K$ is $1$-linearly surjective when $r\ge 2$: 

\begin{lemma}\label{lem:basis of kerUA}
Assume $K$ is $1$-linearly surjective in case $r\ge 2$.  Let
$$A\ =\ (\der-a_r)\cdots (\der-a_1)\qquad\text{where $a_i=b_i^\dagger+\lambda_i$, $b_i\in K^\times$, $\lambda_i\in \Lambda$ \textup{(}$i=1,\dots,r$\textup{)}.}$$
Then there are $C$-linearly independent $y_1,\dots,y_r\in \ker_{\Univ} A$   
 with $y_i\in K^\times\ex(\lambda_i)$ for~$i=1,\dots,r$ and $\operatorname{split}(y_1,\dots,y_r)=(a_1,\dots,a_r)$. 
 %In particular, $\ker_{\Univ} A=\ker_{\Omega} A$. 
\end{lemma}
\begin{proof}
By induction on $r$. The case $r=0$ is trivial, and for $r=1$ we can take~$y_1=b_1\ex(\lambda_1)$. Let $r\geq 2$ and suppose inductively that for
$$B\ :=\ (\der-a_r)\cdots (\der-a_{2})$$
we have $C$-linearly independent $z_2,\dots,z_{r}\in\ker_{\Univ} B$ with
$z_i \in K^\times\ex(\lambda_{i})$ for $i=2,\dots,r$ and $\operatorname{split}(z_2,\dots,z_r)=(a_2,\dots,a_r)$.
For $i=2,\dots,r$, Lemma~\ref{lem:A(U)=U} gives~$y_i\in K^\times \ex(\lambda_i)$ with $(\der-a_1)(y_i)=z_i$. Set $y_1:=b_1\ex(\lambda_1)$, so $\ker_{\Univ} (\der-a_1)=~Cy_1$.  Then~$y_1,\dots,y_r\in\ker_{\Univ} A$
are $C$-linearly independent with $y_i\in K^\times\ex(\lambda_i)$ for $i=1,\dots,r$, and one
verifies easily that $\operatorname{split}(y_1,\dots,y_r)=(a_1,\dots,a_r)$.
\end{proof}

\begin{cor}\label{cor:basis of kerUA}
Assume $K$ is $1$-linearly surjective when $r\ge 2$. Then
$$\text{$A$ splits over~$K$}\ \Longleftrightarrow\ \dim_C \ker_{\Univ} A\ =\ r.$$
\end{cor}

\begin{remark}
If $\dim_{C} \ker_{\Univ} A=r$ and $\lambda_1,\dots,\lambda_d$ are the eigenvalues of $A$ with respect to $\Lambda$, then the differential subring $K\big[\ex(\lambda_1),\ex(-\lambda_1),\dots,\ex(\lambda_d),\ex(-\lambda_d)\big]$ of $\Univ$ is the Picard-Vessiot ring for $A$ over $K$; see \cite[Section~1.3]{vdPS}.  
%\marginpar{``Remark'' skipped for now}
If $K$ is linearly closed  and linearly surjective,
then~$\Univ$ is by Corollary~\ref{cor:basis of kerUA} the universal Picard-Vessiot ring of the differential field~$K$ as defined in \cite[Chapter~10]{vdPS}.
Our construction of~$\Univ$ above is modeled on the description of the universal Picard-Vessiot ring 
of the algebraic closure of $C(\!( t)\!)$ given in    \cite[Chapter~3]{vdPS}. 
\end{remark}

\noindent
Recalling our convention that $r=\order A$, here is a complement to Lemma~\ref{newlembasis}: 

\begin{cor}\label{cor:newlembasis}
Let $V$ be a $C$-linear subspace of $\Univ$ with $r=\dim_C V$. Then 
there is at most one monic $A$ with $V=\ker_{\Univ} A$. Moreover,
the following are equivalent:
\begin{enumerate}
\item[\textup{(i)}] $V=\ker_{\Univ} A$ for some monic $A$  that splits over $K$;
\item[\textup{(ii)}] $V=\ker_{\Univ} B$ for some $B\neq 0$;
\item[\textup{(iii)}] $V=\sum_\lambda V_\lambda$;
\item[\textup{(iv)}] $V$ has a basis contained in $\Univ^\times$.
\end{enumerate}
\end{cor}
\begin{proof}
The first claim follows from [ADH, 5.1.15] applied to the differential fraction field of $\Univ$ in place of $K$.
The implication (i)~$\Rightarrow$~(ii) is clear,
(ii)~$\Rightarrow$~(iii)  was  noted
before Lem\-ma~\ref{newlembasis}, and (iii)~$\Rightarrow$~(iv) is obvious.
For   (iv)~$\Rightarrow$~(i),
 let~$y_1,\dots,y_r\in \Univ^\times$ be a basis of~$V$. Then $\operatorname{split}(y_1,\dots,y_r)=(a_1,\dots, a_r)\in K^r$,
 so  $V=\ker_{\Univ} A$ for~$A=(\der-a_r)\cdots(\der-a_1)$ by Lemma~\ref{lem:distinguished splitting}, so~(i) holds.
\end{proof}

\noindent
Let $y_1,\dots,y_r\in\Univ^\times$ and $(a_1,\dots, a_n):=\operatorname{split}(y_1,\dots,y_r)$. We finish this subsection with some remarks about~$(a_1,\dots,a_n)$. 
Let~$A_0,\dots,A_n\in K[\der]$ be as above and
recall that $n\leq r$ is maximal such that~$y_1,\dots,y_n$ are $C$-linearly independent.

\begin{lemma}\label{lem:bases with same split} Assume $n=r$.  
Let $z_1,\dots,z_r\in\Univ^\times$. The following are equivalent:
\begin{enumerate}
\item[\textup{(i)}] $z_1,\dots,z_r$ are $C$-linearly independent and $(a_1,\dots, a_r)=\operatorname{split}(z_1,\dots,z_r)$;
\item[\textup{(ii)}] for $i=1,\dots, r$ there are $c_{ii}, c_{i,i-1},\dots, c_{i1}\in C$  such that
$$z_i\ =\ c_{ii}y_i + c_{i,i-1}y_{i-1}+\cdots+c_{i1}y_1\  \text{and}\  c_{ii}\neq 0.$$
\end{enumerate}
\end{lemma}
\begin{proof} The case $r=0$ is trivial. Let $r=1$. If (i) holds, then $y_1^\dagger=a_1=z_1^\dagger$, hence~$z_1\in C^\times\, y_1$, so (ii) holds. The converse is obvious. 
Let $r\ge 2$, and assume (i) holds. Put $\tilde{y}_i:=A_1(y_i)$ and $\tilde{z}_i:=A_1(z_i)$ for $i=2,\dots,r$. Then
$$\operatorname{split}(\tilde{y}_2,\dots,\tilde{y}_r)\ =\ (a_2,\dots,a_r)\ =\ 
\operatorname{split}(\tilde{z}_2,\dots,\tilde{z}_r) ,$$
so we can assume inductively to have $c_{ij}\in C$ ($2\leq j\le i\leq r$)   with 
$$\tilde{z}_i\ =\ c_{ii}\tilde{y}_i+c_{i,i-1}\tilde{y}_{i-1}+\cdots+c_{i2}\tilde{y}_2\quad\text{and}\quad c_{ii}\neq 0 \qquad (2\le i\le r).$$
Hence for $2\le i\le r$,
$$z_i \in c_{ii}y_i + c_{i,i-1}y_{i-1}+\cdots+c_{i2}y_2+\ker_{\Univ} A_1.$$
Now use $\ker_{\Univ} A_1=Cy_1$ to conclude (ii).  For the converse, let
$c_{ij}\in C$  be as in (ii).  Then clearly $z_1,\dots,z_r$ are $C$-linearly independent. 
Let $(b_1,\dots,b_r):=\operatorname{split}(z_1,\dots,z_r)$ and $B_{r-1}:=(\der-b_{r-1})\cdots (\der-b_1)$. 
Then $a_r=f_r^\dagger$ where $f_r=A_{r-1}(y_r)\neq 0$, and $b_r=g_r^\dagger$ where $g_r:=B_{r-1}(z_r)\neq 0$.
Now inductively we have~$a_j=b_j$ for $j=1,\dots,r-1$, so $A_{r-1}=B_{r-1}$, and $A_{r-1}(y_i)=0$ for $i=1,\dots,r-1$   by Lemma~\ref{lem:distinguished splitting}.
Hence $g_r=c_{rr}f_r$, and thus $a_r=b_r$.
\end{proof}

\begin{lemma}\label{lem:split mult conj}
Let   $z\in\Univ^\times$. Then 
$\operatorname{split}(y_1 z,\dots,y_r z) = (a_1+z^\dagger, \dots, a_n+z^\dagger)$.
\end{lemma}
\begin{proof}
Since for $m\leq r$, the units $y_1z,\dots,y_mz$ of $\Univ$ are $C$-linearly independent iff~$y_1,\dots,y_m$ are
$C$-linearly independent, we see that the tuples $\operatorname{split}(y_1 z,\dots,y_r z)$ and
$\operatorname{split}(y_1,\dots,y_r)$ have the same length $n$.
Let  $(b_1,\dots,b_n):=\operatorname{split}(y_1 z,\dots,y_r z)$; we show $(b_1,\dots,b_n)=(a_1+z^\dagger,\dots,a_n+z^\dagger)$
by induction on $n$. The case $n=0$ is obvious, so suppose $n\geq 1$.
Then $a_1=y_1^\dagger$ and~$b_1=(y_1z)^\dagger=a_1+z^\dagger$ as required.
By remarks following the proof of  Lemma~\ref{lem:distinguished splitting}  we have
$$  (a_2,\dots,a_n)\ =\ \operatorname{split}\!\big(A_1(y_2),\dots,A_1(y_n)\big)\qquad\text{where $A_1:=\der-a_1$.}$$ 
Now $B_1:=\der-b_1=(A_1)_{\ltimes z^{-1}}$, so likewise
$$ (b_2,\dots,b_n)\ =\ \operatorname{split}\!\big(B_1(y_2z),\dots,B_1(y_nz)\big) \ =\ 
\operatorname{split}\!\big( A_1(y_2)z,\dots,A_1(y_n)z\big).$$
Hence $b_2=a_2+z^\dagger,\dots,b_n=a_n+z^\dagger$ by our inductive hypothesis.
\end{proof}

\noindent
For $f\in \der K$ we let $\int f$ denote an element of $K$ such that~$(\int f)'=f$. 

\begin{lemma}\label{lem:Polya fact, 2}
Let $g_1,\dots,g_r\in K^\times$ and
$$A\  =\  g_1\cdots g_r (\der g_r^{-1}) (\der g_{r-1}^{-1})\cdots  (\der g_1^{-1}),$$
and suppose the integrals below can be chosen such that
$$y_1\ =\ g_1,\quad y_2\ =\ g_1\textstyle\int g_2,\quad \dots,\quad y_r\ =\ g_1\int (g_2\int g_3(\cdots(g_{r-1}\int g_r)\cdots)),$$
Then $y_1,\dots, y_r\in K^\times$, $n=r$, and  $a_i=(g_1\cdots g_i)^\dagger$ for $i=1,\dots,r$.
\end{lemma}
\begin{proof}
Let $b_i:=(g_1\cdots g_i)^\dagger$ for $i=1,\dots,r$. By induction on $i=0,\dots,r$ we show~$n\geq i$ and
$(a_1,\dots,a_i)=(b_1,\dots,b_i)$. This is clear for $i=0$, so suppose~${i\in\{1,\dots,r\}}$, $n\geq i-1$, and $(a_1,\dots,a_{i-1})=(b_1,\dots,b_{i-1})$. 
Then 
$$A_{i-1}=(\der-a_{i-1})\cdots (\der-a_1)=(\der-b_{i-1})\cdots(\der-b_1)=g_1\cdots g_{i-1} (\der g_{i-1}^{-1} )  \cdots   (\der g_1^{-1}),$$
using Lemma~\ref{lem:Polya fact} for the last equality. So $A_{i-1}(y_{i})=g_1\cdots g_{i}\neq 0$, and thus $n\geq i$ and
 $a_{i}=A_{i-1}(y_i)^\dagger=b_i$. 
\end{proof}

\subsection*{Splittings and derivatives\astr} 
The material in this subsection is only needed for the proof of Lemma~\ref{lem:param zeros of derivatives}, and not for the proof of our main theorem.
{\it In this subsection $A$ is monic and $a_0:=A(1)\neq 0$.}\/ Let $A^\der$ be the unique element of $K[\der]$ such that  $A^\der \der = \der A - a_0^\dagger A$.
Then $A^\der$ is monic of order $r$, and if $A\in H[\der]$ for some differential subfield $H$ of $K$, then also $A^\der\in H[\der]$.

\begin{examples}
If $\order A=0$ then $A^\der=1$, and if $\order A=1$ then $A^\der=\der+(a_0-a_0^\dagger)$.
Next, suppose $A=\der^2+a_1\der+a_0$ ($a_0,a_1\in K$); then
$$A^\der\ =\ \der^2+(a_1-a_0^\dagger)\der+(a_1'+a_0-a_1a_0^\dagger).$$
\end{examples}

\noindent
If  $A(y)=0$ with $y$ in a differential ring extension of $K$, then~$A^\der(y')=0$. Also:

\begin{lemma}\label{lem:A^der}
Let $R$ be a differential integral domain extending $K$. Suppose the differential fraction field of $R$ has constant field $C$, and  $\dim_C \ker_R A=r$.
Then~$\ker_R A^\der=\{y':y\in\ker_R A\}$ and $\dim_C \ker_R A^\der=r$.
\end{lemma}
\begin{proof}
Let $y_1,\dots,y_r$ be a basis of the $C$-linear space $\ker_R A$, and assume towards a contradiction that $c_1y_1'+\cdots+c_ry_r'=0$ with $c_1,\dots,c_r\in C$ not all zero. 
 Then~$y:=c_1y_1+\cdots+c_ry_r\in\ker^{\neq}_R A$ and~$y'=0$. Hence $a_0y=A(y)=0$ and thus $a_0=0$,
a contradiction.
\end{proof}

\noindent
Let $\Univ=\Univ_K$ and $f_1\ex(\lambda_1),\dots,f_r\ex(\lambda_r)\in\Univ^\times$ be a basis of
the $C$-linear space $\ker_{\Univ} A$, where~$f_j\in K^\times$ and 
$\lambda_j\in\Lambda$ for $j=1,\dots,r$. Then by Lemma~\ref{lem:A^der},
$$(f_1'+\lambda_1f_1)\ex(\lambda_1),\dots,(f_r'+\lambda_rf_r)\ex(\lambda_r)\in\Univ^\times$$ is a basis of
the $C$-linear space $\ker_{\Univ} A^\der$. Hence by Corollary~\ref{corbasiseigenvalues}:

\begin{cor}\label{cor:A^der}
Suppose $\dim_C \ker_{\Univ} A=r$. Then
$\mult_\alpha(A)=\mult_\alpha(A^\der)$ for all~$\alpha\in K/K^\dagger$, so $\Sigma(A)=\Sigma(A^\der)$, and both $A$, $A^\der$ split over $K$.
\end{cor}

\noindent
Suppose now that $K$ is $1$-linearly surjective when~$r\geq 2$, and~$A$ splits over $K$. Then~$A^\der$ splits over $K$ by Corollaries~\ref{cor:basis of kerUA} and~\ref{cor:A^der}.
%if $K$ is $1$-linearly surjective when~$r\geq 2$, then~$A$ splits over $K$ iff $A^\der$ splits over $K$.
%(We don't know if this equivalence continues to hold without the hypothesis on $K$.)

\subsection*{Splitting and adjoints\astr} 
In  this subsection $y_1,\dots,y_r\in \Univ^\times$,
 $$(a_1,\dots,a_r)\ =\ \operatorname{split}(y_1,\dots,y_r),\qquad A\ =\ (\der-a_r)\cdots(\der-a_1).$$ 
 So $y_1,\dots,y_r$ is a basis of the $C$-linear space $V:=\ker_{\Univ} A=\ker_{\Omega} A$, 
 $$A^*\ =\ (-1)^r (\der+a_1)\cdots(\der+a_r),$$ and
$\dim_C W=r$ for the $C$-linear space $W:=\ker_{\Univ} A^*=\ker_{\Omega} A^*$  by Lemma~\ref{lem:dim ker A^*}. 
(Recall here that $\Omega$ denotes the differential fraction field of $\Univ$.)
Proposition~\ref{prop:Lagrange} with $\Omega$ instead of $K$ yields the $C$-bilinear map $[\ ,\,]_A\colon \Omega\times \Omega\to \Omega$, which restricts to a perfect pairing $V\times W \to C$ by
Corollary~\ref{cor:nondeg pairing}. We let $j,k$ range over $\{1,\dots,r\}$ and take~$\lambda_j\in\Lambda$ such that~$y_j^\dagger \equiv \lambda_j\bmod K^\dagger$, so $y_j\in K^\times\ex(\lambda_j)$.

% units of $\Univ$  such that
%$A=(\der-a_r)\cdots(\der-a_1)$ for
%$(a_r,\dots,a_1):=\operatorname{split}(y_r,\dots,y_1)$. Note the deliberate inversion of the order: previously
%we worked with $(a_r,\dots,a_1)=\operatorname{split}(y_1,\dots,y_r)$.
%Then~$y_1,\dots,y_r$ is a basis of the $C$-linear space $V:=\ker_{\Univ} A=\ker_\Omega A$, and $(-a_1,\dots,-a_r)$ is a splitting
%of $A^*$.
%In Corollary~\ref{cor:dual split} below we obtain a basis $y_1^*,\dots,y_r^*$ of the $C$-linear space $W:=\ker_{\Univ} A^*=\ker_\Omega A^*$ which 
%is dual to the basis $y_1,\dots,y_r$ of $V$ 
%under the perfect pairing $[\ ,\,]_A\colon V\times W\to C$ from Corollary~\ref{cor:nondeg pairing}  such that
%$\operatorname{split}(y_1^*,\dots,y_r^*)=(-a_1,\dots,-a_r)$.
%Let $k$, $l$ range over~$\{1,\dots,r\}$, and
%let~$\lambda_k\in\Lambda$ such that~$y_k^\dagger \equiv \lambda_k\bmod K^\dagger$.

\begin{lemma}\label{lem:dual split}
Suppose $z_1,\dots,z_r\in \Univ^\times$ are $C$-linearly independent
such that $$\operatorname{split}(z_r,\dots,z_1)\ =\ (-a_r,\dots,-a_1).$$
%  such that $$\operatorname{split}(z_1,\dots,z_r)\ =\ (-a_1,\dots,-a_r).$$
Then $z_1,\dots, z_r$ is a basis of the $C$-linear space $W$, $[y_j,z_k]_A=0$ if $j<k$, and $[y_k,z_k]_A\neq 0$.
Moreover, $z_k\in K^\times \ex(-\lambda_k)$, and
if~$[y_j,z_k]_A\neq 0$, then $\lambda_j=\lambda_k$.
\end{lemma}
\begin{proof}
Let $i$ range over $\{0,\dots,r\}$. As in \eqref{eq:AiBi}, set
$$A_i\ :=\ (\der-a_i)\cdots(\der-a_1),\qquad B_i\ :=\ (-1)^{r-i}(\der+a_{i+1})\cdots(\der+a_r).$$
 %Define $A_i,B_i\in K[\der]$  be as 
 %in \eqref{eq:AiBi}. 
Then by Lemma~\ref{lem:distinguished splitting} we have
$$\ker_{\Univ} A_i\ =\  Cy_{1}\oplus\cdots\oplus Cy_{i}, \qquad
\ker_{\Univ} B_{i}\ =\  Cz_{r} \oplus\cdots\oplus Cz_{i+1}$$
%$$\ker_{\Univ} A_i\ =\  Cy_{r}\oplus\cdots\oplus Cy_{r-i+1}, \qquad
%\ker_{\Univ} B_{i}\ =\  Cz_{1} \oplus\cdots\oplus Cz_{r-i}$$
and thus
$$A_i(y_{j})\ =\ 0\ \text{ if $i\geq j$},\qquad B_{i}(z_{k})\ =\ 0\ \text{ if $i+1\le k$.}$$
Then Lemma~\ref{lem:P_A split}  yields 
$$[y_{j}, z_{k}]_A\ =\ \sum_{i<r} A_i(y_j)B_{i+1}(z_k)\ =\ \sum_{k-2< i<j} A_i(y_{j})B_{i+1}(z_{k}),$$
so $[y_{j}, z_k]_A=0$ whenever $j<k$. Moreover, 
$$[y_{k},z_k]_A\ =\  A_{k-1}(y_{k}) B_{k}(z_k)\  \neq\  0.$$
%and thus
%$$A_i(y_{k})\ =\ 0\ \text{ if $i\geq r-k+1$},\qquad B_{i}(z_{l})\ =\ 0\ \text{ if $i\leq r-l$.}$$
%This yields 
%$$[y_{k}, z_{l}]_A\ =\ \sum_{i=r-l}^{r-k} A_i(y_{k})B_{i+1}(z_{l});$$
%in particular $[y_{k}, z_l]_A=0$ whenever $k>l$.
%Since  $\operatorname{split}(y_r,\dots,y_1)=(a_r,\dots,a_1)$ and~$\operatorname{split}(z_1,\dots,z_r)=(-a_1,\dots,-a_r)$, 
%we also have
%$$[y_{k},z_k]_A\ =\  A_{r-k}(y_{k}) B_{r-k+1}(z_k) \neq 0.$$
Take $\mu_k\in\Lambda$ with 
$z_k^\dagger\equiv \mu_k\bmod K^\dagger$. Then $y_j\in K\ex(\lambda_j)$ and $z_k\in K\ex(\mu_k)$, so~$[y_j,z_k]_A\in C\cap K\ex(\lambda_j+\mu_k)$ by \eqref{eq:P_L}. Hence, if~$[y_j,z_k]_A\neq 0$, then $\lambda_j+\mu_k=0$.
In particular, $\mu_k=-\lambda_k$ and so $z_k\in K^\times\ex(-\lambda_k)$.
\end{proof}

\begin{cor}\label{cor:dual split}
Assume $K$ is $1$-linearly surjective if $r\geq 2$.
Then there is a basis $y_1^*,\dots,y_r^*$ of the $C$-linear space $W$ such  that $[y_j,y_k^*]_A=\delta_{jk}$ for all $j,k$, and 
$$y_j^*\in K^\times\ex(-\lambda_j) \text{ for all }j, \qquad \operatorname{split}(y_r^*,\dots,y_1^*)\ =\ (-a_r,\dots,-a_1).$$
\end{cor}

\begin{proof}
Lemma~\ref{lem:basis of kerUA} gives $C$-linearly independent 
  $z_1,\dots,z_r\in\Univ^\times$ such that $$\operatorname{split}(z_r,\dots,z_1)\ =\ (-a_r,\dots,-a_1).$$ 
Lemma~\ref{lem:dual split} gives constants  $c_k\in C^\times$ such  that~$[y_j,c_kz_k]_A=\delta_{jk}$ for $j\leq k$.  
We now set $y_r^*:= c_rz_r$, so $[y_j, y_r^*]_A=\delta_{jr}$ for all $j$ and $y_r^*\in K^\times \ex(-\lambda_r)$.  Let $1<k\le r$ and assume inductively that we have
$y_k^*,\dots, y_r^*\in W$ such that for $i=k,\dots,r$ we have $y_i^*\in Cz_i+\cdots + Cz_r$, $[y_j, y_i^*]_A=\delta_{ji}$ for all $j$, and $y_i^*\in K^\times \ex(-\lambda_i)$.
Then for
$$ y_{k-1}^*\ :=\  c_{k-1}z_{k-1}-\sum_{i=k}^r[y_i,c_{k-1}z_{k-1}]_Ay_i^*$$
we have 
$$y_{k-1}^*\in Cz_{k-1}+ Cz_k + \cdots + Cz_r,\qquad  [y_j, y_{k-1}^*]_A=\delta_{j,k-1}\text{ for all $j$.}$$
If~$k\le i\le r$ and $[y_i, c_{k-1}z_{k-1}]_A\ne 0$, then $\lambda_{i}=\lambda_{k-1}$ by the last part of Lemma~\ref{lem:dual split}, so
$y_{k-1}^*\in K^\times \ex(-\lambda_{k-1})$ by the inductive assumption and $z_{k-1}\in K\ex(-\lambda_{k-1})$.

This recursive construction yields a basis $y_1^*,\dots, y_r^*$ of the $C$-linear space $W$ such that $[y_j, y_k^*]=\delta_{jk}$ for all $j$,~$k$, and $y_i^*\in K^\times \ex(-\lambda_i)$ for $i=1,\dots,r$. It now follows from Lemma~\ref{lem:bases with same split} that
$\operatorname{split}(y_r^*,\dots,y_1^*)=(-a_r,\dots,-a_1)$.
\end{proof}

\noindent
Lemma~\ref{lem:dual split}  also yields:

\begin{cor}\label{corcordualsplit} If  $(a_1,\dots,a_r)=(-a_r,\dots,-a_1)$,
then $[y_j,y_{r+1-k}]_A=0$ for all~$j<k$, and $[y_k,y_{r+1-k}]_A\neq 0$ for all $k$.
\end{cor}

%\noindent \marginpar{next statement does not fit the conventions from the beginning of this subsection} 
%If $A$ is self-adjoint or skew-adjoint, then by Lemma~\ref{lem:Jacobi, 2} and Proposition~\ref{prop:Darboux skew-adjoint}, 
%there is a splitting $(a_r,\dots,a_1)$ of $A$ satisfying the hypothesis of Corollary~\ref{corcordualsplit}.
%We finish with an observation which follows from Corollary~\ref{cor:sum of evs} and the remark after Proposition~\ref{prop:Bogner}:

%\begin{cor}
%If $A^*=(-1)^rA_{\ltimes a}$ where $a\in K^\times$, then $\lambda_1+\cdots+\lambda_r=0$.
%\end{cor}

\subsection*{The case of real operators} 
We now continue the subsection {\it The real case}\/ of Section~\ref{sec:univ exp ext}.
Thus $K=H[\imag]$ where $H$ is a real closed differential subfield of~$K$ and~${\imag^2=-1}$, and $\Lambda=\Lambda_{\operatorname{r}}+\Lambda_{\operatorname{i}}\imag$ where $\Lambda_{\operatorname{r}}$, $\Lambda_{\operatorname{i}}$ are subspaces
of the $\Q$-linear space~$H$.
The complex conjugation automorphism $z\mapsto \bar{z}$ of the differential field~$K$ extends uniquely to an automorphism
$B\mapsto \bar{B}$ of the ring $K[\der]$ with~$\bar{\der}=\der$. We have~$\overline{A(f)}=\overline{A}(\overline{f})$ for $f\in\Univ$, from which it follows that $\dim_C \ker_K A=\dim_C \ker_K\overline{A}$,
%If $A\in K[\der]$ has spectral decomposition~$(A_\lambda)$, then $\overline{A}$ has spectral decomposition~$
$(\overline{A})_{\lambda}=\overline{(A_{\overline{\lambda}})}$, 
$\mult_\lambda \overline{A} = \mult_{\overline{\lambda}} A$, and $f\mapsto\overline{f}\colon\Univ\to \Univ$ restricts to a $C_H$-linear bijection~$\ker_{\Univ} A \to \ker_{\Univ} \overline{A}$.

\medskip
\noindent
{\it In the rest of this subsection we assume $H=H^\dagger$ $($so $\Lambda=\Lambda_{\operatorname{i}}\imag)$ and $A\in H[\der]$ $($and by earlier conventions, $A\ne 0$ and $r:=\order A)$.}\/
Then $A=\overline{A}$, hence for all $\lambda$ we have $A_{\overline{\lambda}}=\overline{A_\lambda}$ and $\mult_\lambda  A = \mult_{\overline{\lambda}} A$. Thus with $\mu$ ranging over $\Lambda_{\operatorname{i}}^{>}$:
$$\sum_\lambda \mult_\lambda(A)\ =\ \mult_0(A) + 2\sum_{\mu} \mult_{\mu\imag}(A).$$
%Also, if $\lambda$ is an eigenvalue of $A$ with respect to $\Lambda$, then so is $\overline{\lambda}$.  
Note that  $0$ is an eigenvalue of $A$ iff $\ker_H A\neq\{0\}$.

\medskip
\noindent
Let $V:=\ker_{\Univ} A$, a  subspace of the $C$-linear space $\Univ$ with $\overline{V}=V$ and $\dim_C V\leq r$.
Recall that we have the differential $H$-subalgebra $\Univ_{\operatorname{r}}=\{f\in\Univ:\overline{f}=f\}$ of $\Univ$
and the $C_H$-linear subspace $V_{\operatorname{r}} = \ker_{\Univ_{\operatorname{r}}} A$ of $\Univ_{\operatorname{r}}$.
Now
$V=V_{\operatorname{r}}\oplus V_{\operatorname{r}}\imag$ (internal direct sum of $C_H$-linear subspaces), so 
$\dim_C V = \dim_{C_H} V_{\operatorname{r}}$. 
Combining Lemma~\ref{newlembasis} and the remarks preceding it with Lemma~\ref{lem:real basis} and its proof yields:

{\begin{cor}\label{cor:complex and real basis}
{\samepage The $C$-linear space $V$ has a basis 
$$a_1\ex(\mu_1\imag),\,\overline{a_1}\ex(-\mu_1\imag),\ \dots,\ a_m\ex(\mu_m\imag),\,\overline{a_m}\ex(-\mu_m\imag), \ h_1,\ \dots,\ h_n \qquad (2m+n\leq r),$$ 
where $a_1,\dots,a_m\in K^\times$, $\mu_1,\dots,\mu_m\in \Lambda_{\operatorname{i}}^{>}$, $h_1,\dots,h_n\in H^\times$.} For such a basis,
$$\Re\!\big(a_1\ex(\mu_1\imag)\big),\,\Im\!\big(a_1\ex(\mu_1\imag)\big),\ \dots,\ \Re\!\big(a_m\ex(\mu_m\imag)\big),\,\Im\!\big(a_m\ex(\mu_m\imag)\big), \ h_1,\ \dots,\ h_n$$
is a basis of the $C_H$-linear space $V_{\operatorname{r}}$, and $h_1,\dots,h_n$ is a basis of the $C_H$-linear subspace~$\ker_H A=V\cap H$ of $H$. 
\end{cor}

\noindent
%If $\lambda\in\Lambda_{\operatorname{r}}$, then $A_\lambda\in H[\der]$, and if additionally $A$ splits over $H$, then
%so does $A_\lambda$; moreover, if $A_\lambda(H)=H$, then $A\big(H\ex(\lambda)\big)=H\ex(\lambda)$. In particular, if  $H$ is $r$-linearly surjective, or $H$ is $1$-linearly surjective and $A$ splits over $H$, then for   $\lambda\in\Lambda_{\operatorname{r}}$ we have
% $A\big(H\ex(\lambda)\big)=H\ex(\lambda)$ and hence $A(\Univ_\lambda)=\Univ_\lambda$.
Using $H=H^\dagger$, arguments as in the proof of Lemma~\ref{lem:basis of kerUA} show:
 
\begin{lemma} \label{lem:basis of kerUA, real}
Assume $H$ is $1$-linearly surjective when $r\ge 2$.  Let $a_1,\dots, a_r\in H$ be such that
$A =(\der-a_r)\cdots (\der-a_1)$. 
Then the $C_H$-linear space $\ker_H A$ has a basis~$y_1,\dots, y_r$
such that $\operatorname{split}(y_1,\dots,y_r)=(a_1,\dots, a_r)$. 
\end{lemma}
%\begin{proof}
%By induction on $r$. The case $r=0$ is trivial, and for $r=1$ we can take $y_1=b_1\ex(\lambda_1)$. Let $r\geq 2$ and suppose inductively that for 
%$$B\ :=\ (\der-a_r)\cdots (\der-a_{2})$$
%we have $C$-linearly independent $z_2,\dots,z_{r}\in\ker_{\Univ} B$ with $z_i \in H^\times\ex(\lambda_{i})$ for $i=2,\dots,r$ and $\operatorname{split}(z_2,\dots,z_r)=(a_2,\dots,a_r)$. For $i=2,\dots,r$, the remarks preceding the lemma give $y_i\in H^\times \ex(\lambda_i)$ with $(\der-a_1)(y_i)=z_i$. Set $y_1:=b_1\ex(\lambda_1)$, so $\ker_{\Univ} (\der-a_1)=~Cy_1$.  Then $y_1,\dots,y_r\in\ker_{\Univ} A$ are $C$-linearly independent with $y_i\in H^\times\ex(\lambda_i)$ for $i=1,\dots,r$, and one verifies easily that $\operatorname{split}(y_1,\dots,y_r)=(a_1,\dots,a_r)$.
%\end{proof}

%\begin{cor}\label{cor:basis of kerUA, real}
%If $r\geq 2$, suppose $H$ is $1$-linearly surjective   and  $A$ splits over~$H$; then
%$\dim_{C} \ker_{\operatorname{U}} A\ =\ r$.
%\end{cor} 

%Here is how to determine the quantities $m$, $n$ in Corollary~\ref{cor:complex and real basis} from $A$:

%\begin{lemma} \marginpar{new lemma}
%Assume $H$ is $1$-linearly surjective when $r\ge 2$ and $A$ splits over~$K$. Let $a_i$, $\mu_i$, $h_j$ be as in Corollary~\ref{cor:complex and real basis} and $(b_1,\dots,b_n):=\operatorname{split}(h_1,\dots,h_n)$. Then~$A_{\operatorname{s}}=(\der-b_n)\cdots(\der-b_1)$, so $n=\order A_{\operatorname{s}}$.
%\end{lemma}
%\begin{proof}
%\end{proof}

\noindent
Recall from Lemma~\ref{lem:size of Sigma(A)} that
if $r=1$ or   $K$ is   $1$-linearly surjective,   then
$$\text{$A$ splits over $K$}\quad\Longleftrightarrow\quad \sum_\lambda \mult_\lambda(A)=r.$$
Now $\mult_\lambda(A)=\mult_{\overline{\lambda}}(A)$ for all $\lambda$, so if
$\mult_\lambda(A)=r\ge 1$, then $\lambda=0$.
Also, for~$W:=V\cap K=\ker_K A$  and $W_{\operatorname{r}}:=W\cap\Univ_{\operatorname{r}}$ we have~$W_{\operatorname{r}}=\ker_H A$ and
$$W\ =\ W_{\operatorname{r}}\oplus W_{\operatorname{r}}\imag\quad\text{ (internal direct sum of $C_H$-linear subspaces)},$$ so $\mult_0(A)=\dim_C \ker_K A=\dim_{C_H} \ker_H A$. If $y_1,\dots,y_r$ is a basis of the $C_H$-linear space~$\ker_H A$, then  $\operatorname{split}(y_1,\dots,y_r)\in H^r$ in reversed order is a splitting of~$A$ over $H$ by Corollary~\ref{corbasissplit}.
These remarks and Lem\-ma~\ref{lem:basis of kerUA, real} now yield:

\begin{cor}\label{cor:unique eigenvalue, 1}
If 
$\mult_0(A)=r$, then $A$ splits over $H$. The converse
holds if~$H$ is $1$-linearly surjective or $r = 1$. 
\end{cor}

\begin{cor}\label{cor:unique eigenvalue, 2} 
Suppose $r\ge 1$, and $K$ is $1$-linearly surjective if $r\geq 2$.  Then
$$\text{$A$ splits over~$H$} \quad \Longleftrightarrow\quad  \mult_0(A)=r\quad \Longleftrightarrow\quad \abs{\Sigma(A)}=1.$$
\end{cor}

\noindent
We now focus on the order $2$ case:

\begin{lemma}\label{lem:order 2 eigenvalues}
Suppose  $r=2$ and $A$ splits over~$K$ but not over $H$. 
Then  $$\dim_{C} \ker_{\Univ} A\ =\ 2.$$
If $H$ is $1$-linearly surjective, 
then $A$ has  two distinct eigenvalues. 
\end{lemma}

\begin{proof}
We can assume $A$ is monic, so $A=(\der-f)(\der-g)$ with $f,g\in K$ and~$g=a+b\imag$, $a,b\in H$, $b\ne 0$. Then  $g=d^\dagger+\mu\imag$ with $d\in K^\times$ and $\mu\in \Lambda_{\operatorname{i}}$, and so  $d\ex(\mu\imag)\in \ker_{\operatorname{U}}A$. From $A=\bar{A}$ we obtain $\bar{d}\ex(-\mu\imag)\in \ker_{\operatorname{U}}A$. These two elements of~$\ker_{\Univ} A$ are $C$-linearly independent, since 
$$d\ex(\mu\imag)/\bar{d}\ex(-\mu\imag)\ =\ (d/\bar{d})\ex(2\mu\imag)\notin C:$$
this is clear if $\mu\ne 0$, and if $\mu=0$, then $d^\dagger=g$, so $(d/\bar{d})^\dagger=g-\bar{g}=2b\imag\ne 0$, and hence $d/\bar{d}\notin C$.   Thus $\dim_{C} \ker_{\operatorname{U}} A\ =\ 2$, and $\mu\imag$, $-\mu \imag$ are eigenvalues of~$A$ with respect to $\Lambda$. Now assume $H$ is $1$-linearly surjective. Then we claim that~$\mu\neq 0$. To see this note that [ADH, 5.1.21, 5.2.10] and the assumption that $A$ does not split over $H$ yield $ \dim_{C_H} \ker_H A = \dim_C \ker_K A =0$, hence $g\notin K^\dagger$ and thus~$\mu\imag=g-d^\dagger\neq 0$.
\end{proof}

%\begin{example} \marginpar{commented out material checked but no longer needed}
%Suppose $A=4\der^2+f$ where $f\in H$. Let the functions $\omega\colon H \to H$ and~$\sigma\colon H^\times \to H$ be  as in [ADH, 5.2]. Then we have [ADH, (5.2.1)] :
%$$\text{$A$ splits over $K$ but not over $H$}\quad\Longleftrightarrow\quad f\in \sigma(H^\times)\setminus\omega(H).$$
%Suppose $y\in H^\times$ satisfies $\sigma(y)=f\notin\omega(H)$.
%Then by [ADH, p.~262] we have~$A=4(\der+g)(\der-g)$ where $g=\frac{1}{2}(-y^\dagger+y\imag)$.
%If $H$ is $1$-linearly surjective, then the two distinct eigenvalues of $A$
%are $(y/2)\imag+K^\dagger$ and $-(y/2)\imag+K^\dagger$.
%\end{example}

\noindent
Combining Lemmas~\ref{lem:basis of kerUA, real} and~\ref{lem:order 2 eigenvalues}    yields:  

\begin{cor}\label{spldcr2} If $H$ is $1$-linearly surjective, $A$ has order~$2$, and $A$  splits over~$K$, then~$\dim_{C} \ker_{\operatorname{U}} A\ =\ 2$.
\end{cor} 

\noindent
{\it In the rest of this subsection $H$ is $1$-linearly surjective and $A=4\der^2+f$, $f\in H$.}\/
 Let the functions~$\omega\colon H \to H$ and~$\sigma\colon H^\times \to H$ be  as in [ADH, 5.2]. Then we have~[ADH, remarks before~5.2.1, and (5.2.1)]:
\begin{align*}
\text{$A$ splits over $H$} &\quad\Longleftrightarrow\quad f\in\omega(H), \\
\text{$A$ splits over $K$} &\quad\Longleftrightarrow\quad f\in \sigma(H^\times)\cup\omega(H).
\end{align*}
If $A$ splits over $H$, then $\Sigma(A)=\{0\}$  and~$\operatorname{mult}_0(A)=2$, by Corollary~\ref{cor:unique eigenvalue, 2}.
Suppose~$A$ splits over $K$ but not over~$H$, and let~$y\in H^\times$ satisfy~$\sigma(y)=f\notin\omega(H)$.
Then by~[ADH, p.~262] we have~$A=4(\der+g)(\der-g)$ where~$g=\frac{1}{2}(-y^\dagger+y\imag)$.
Hence the two distinct eigenvalues of~$A$
are~$(y/2)\imag+K^\dagger$ and~$-(y/2)\imag+K^\dagger$.
We
consider also the skew-adjoint differential operator
$$B\ :=\ \der^3+f\der+(f'/2)\in H[\der].$$ 
If  $\dim_C \ker_{\Univ}A=2$, then $\dim_C \ker_{\Univ}B=3$
by Lemma~\ref{lem:Appell}.
% applied tothe differential fraction field of $\Univ$ in place of $K$.  
Likewise, $$ \dim_{C_H} \ker_{H}A=2\ \Longrightarrow\ \dim_{C_H} \ker_{H}B=3.$$

\begin{lemma}
If $A$ splits over $K$, then so does~$B$. Likewise with $H$ instead of~$K$.
\end{lemma}
\begin{proof}
If $A$ splits over $K$, then $\dim_C\ker_{\Univ}A=2$ by  Corollary~\ref{spldcr2} and therefore~$\dim_C \ker_{\Univ}B=3$
by the remark preceding the lemma, so $B$ splits over $K$ by Corollary~\ref{corbasiseigenvalues}.
If $A$ splits over $H$, then $\dim_{C_H}\ker_H A=2$ by Lemma~\ref{lem:basis of kerUA, real} and hence~$\dim_{C_H}\ker_H B=3$, so $B$ splits over $H$ by Corollary~\ref{corbasissplit} and the remark following it.
\end{proof}

\begin{lemma}\label{lem:kerB}
Let $y\in H^\times$ with~$\sigma(y)=f\notin\omega(H)$. Then
$\Sigma(B)=\{\beta,0,-\beta\}$ where $\beta:=y\imag+K^\dagger\neq 0$,  and $\dim_{C_H} \ker_H  B=1$.
\end{lemma}
\begin{proof}
Put $g=\frac{1}{2}(-y^\dagger+y\imag)$, so $A=4(\der+g)(\der-g)$, and
take $d\in K^\times$ and $\mu\in\Lambda_{\operatorname{i}}$ with $g=d^\dagger+\mu\imag$.
Then $d\ex(\mu\imag)$, $\bar{d}\ex(-\mu\imag)$ is a basis of $\ker_{\Univ}A$ and $\mu\neq 0$, by the argument in the proof
of Lemma~\ref{lem:order 2 eigenvalues}. Hence
$$d^2\ex(2\mu\imag),\quad \abs{d}^2,\quad \bar{d}{}^2\ex(-2\mu\imag)$$ is a basis of $\ker_{\Univ} B$ by  Lemma~\ref{lem:Appell}, so
$$\big(d^2\ex(2\mu\imag)\big)^\dagger+K^\dagger=2\mu\imag+K^\dagger,\  \big(\abs{d}^2\big)^\dagger+K^\dagger=[0],\ \big(\bar{d}{}^2\ex(-2\mu\imag)\big)^\dagger+K^\dagger=-2\mu\imag+K^\dagger$$
are eigenvalues of $B$. Since $\mu\imag\notin K^\dagger$, these are distinct eigenvalues, and so there are no other
eigenvalues. Note: $g=\frac{1}{2}(-y^\dagger+y\imag)=d^\dagger+\mu\imag$ gives $y\imag+K^\dagger=2\mu\imag + K^\dagger$.  Finally, 
$\dim_{C_H} \ker_H  B=1$ by Corollary~\ref{cor:complex and real basis}.
\end{proof}

\subsection*{Factoring linear differential operators over $H$-fields\astr} {\it In this subsection~$H$ is a real closed $H$-field with $x\in H$, $x'=1$, $x\succ 1$, and~$K=H[\imag]$ where~$\imag^2=-1$.}\/ In the proof of the next lemma we use [ADH, 10.5.2(i)]: 
\begin{equation}\label{eq:10.5.2}
y,z\in H^\times,\ y\prec z\  \Longrightarrow\  y^\dagger<z^\dagger.
\end{equation}

\begin{lemma}\label{lem:prec vs <}
Let $y$, $z$ be $C_H$-linearly independent elements of $H^\times$
and~$(a,b):=\operatorname{split}(y,z)$. 
If    $xy\succ z$,  then $a>b$, and 
if    $xy\prec z$,  then $a<b$.
\end{lemma}
\begin{proof}
Replacing $(y,z)$ by $(1,z/y)$ we arrange $y=1$, $a=0$, by Lem\-ma~\ref{lem:split mult conj}.  Then $z'\neq 0$ and $b=z'{}^\dagger$.  Now $x\succ z$ implies $1\succ z'$,  and
$x\prec z$ implies $1\prec z'$. It remains to use the remark preceding the lemma.
\end{proof}

\noindent
In the next three lemmas $y_1, \dots, y_r\in H$ ($r\in\N$) are $C_H$-linearly independent
and~$(a_1,\dots,a_r):=\operatorname{split}(y_1,\dots,y_r)\in H^r$. We also assume that $H$ is $\upl$-free.

\begin{lemma}\label{lem:prec vs <,r}
$y_1\succ\cdots\succ y_r \ \Rightarrow\ a_1>\cdots>a_r$.
\end{lemma}
\begin{proof}
The cases  $r=0,1$ are trivial, so suppose that $r\geq 2$ and
$y_1\succ\cdots\succ y_r$.
We have~$(a_1,\dots,a_{r-1})=\operatorname{split}(y_1,\dots,y_{r-1})$. Assume 
 $a_1>\cdots>a_{r-1}$ as inductive hypothesis. It remains to show $a_{r-1}>a_r$. Put $B:=({\der-a_{r-2}})\cdots (\der-a_1)$;
so $B=A_{r-2}$ in the notation introduced before Lemma~\ref{lem:distinguished splitting}, and $(a_{r-1},a_r)=\operatorname{split}\!\big(B(y_{r-1}),B(y_r)\big)$. By Lemma~\ref{lem:distinguished splitting}, $y_1,\dots,y_{r-2}$ is a basis of the $C_H$-linear subspace~$\ker_H B$ of $H$, 
and hence $$v(\ker_H^{\neq} B)\ =\  \big\{ v(y_1),\dots,v(y_{r-2})\big\}\ =\  \exc^{\ev}(B)$$ by Corollary~\ref{cor:size of excev}, so $v(y_{r-1}), v(y_r)\notin \exc^{\ev}(B)$. 
Then Lemma~\ref{lem:ADH 14.2.7} gives $B(y_{r-1})\succ B(y_r)$, so $xB(y_{r-1})\succ B(y_r)$. Now Lem\-ma~\ref{lem:prec vs <} yields $a_{r-1}>a_r$.
\end{proof}

\begin{lemma}\label{lem:prec vs <,r, 2}
$y_1\prec^\flat\cdots \prec^\flat y_r \ \Rightarrow\ a_1<\cdots<a_r$.
\end{lemma}
\begin{proof}
Similar  to the proof of Lemma~\ref{lem:prec vs <,r}, using  in the inductive step that $B$ is asymptotically
surjective by Corollary~\ref{cor:well-behaved}, hence
if $y,z\in H^\times$, $vy,vz\notin\exc^{\ev}(B)$, and $y\prec^\flat z$, then $B(y)\prec^\flat B(z)$ by
Lemma~\ref{lem:v(A(y)) convex subgp}, and so $xB(y)\prec  B(z)$.
\end{proof}

\noindent
Along the lines of the proof of Lemma~\ref{lem:prec vs <,r} we obtain:

\begin{lemma}\label{lem:preceq vs leq,r}
Suppose  $y_i\,\nasymp\, y_j$  for all $i,j$ with $1\leq i<j\leq r$. Then $$a_1\leq \cdots \leq a_r \ \Rightarrow\ y_1\prec\cdots\prec y_r.$$
%$a_1\geq \cdots \geq a_r$ and $y_i\nasymp^\flat y_j$ \textup{(}$1\leq i<j\leq r$\textup{)} $\ \Rightarrow\ y_1\succ^\flat\cdots\succ^\flat y_r$.
\end{lemma}

\noindent
Under present assumptions we can strengthen the conclusion of Lemma~\ref{lem:basis of kerUA, real}:

\begin{lemma} \label{lem:basis of kerUA, real, H-field}
Assume $H$ is Liouville closed. Let $a_1,\dots, a_r\in H$ and set $$A\ :=\ (\der-a_r)\cdots (\der-a_1).$$ 
Then the $C_H$-linear space $\ker_H A$ has a basis~$y_1,\dots, y_r$
such that $\operatorname{split}(y_1,\dots,y_r)=(a_1,\dots, a_r)$ and $y_i\nasymp y_j$  for all $i,j$ with $1\leq i<j\leq r$. 
\end{lemma}
\begin{proof}
The case $r=0$ is clear. Let $r\geq 1$ and assume inductively that
$$B\ :=\ (\der-a_r)\cdots (\der-a_{2})$$
has a basis $z_2,\dots,z_{r}$ of $\ker_{H} B$ such that $\operatorname{split}(z_2,\dots,z_r)=(a_2,\dots,a_r)$
and~$z_i\nasymp z_j$ whenever $2\leq i<j\leq r$. 
Take~$y_1\in H^\times$ with $y_1^\dagger=a_1$, so
 $\ker_{H} (\der-a_1)=Cy_1$ and $\exc_H^{\ev}(\der-a_1)=\{vy_1\}$. 
For~$i=2,\dots,r$, Corollary~\ref{cor:nonexc sol} then gives $y_i\in H^\times$ with~$(\der-a_1)(y_i)=z_i$ and~$y_i\nasymp y_1$.   Then $y_1,\dots,y_r\in\ker_{H} A$ and $y_i\nasymp y_j$ for  all~$i\ne j$, 
%by Lemma~\ref{lem:ADH 14.2.7}, ? 
and~$\operatorname{split}(y_1,\dots,y_r)=(a_1,\dots,a_r)$.
\end{proof}

\noindent
The valuation of $H$ being trivial on $C_H$, the proof of the next lemma is obvious.

\begin{lemma}\label{lem:upper triag}
Let~${g_1\prec\cdots\prec g_n}$ in $H$ and let $h_1,\dots,h_n$ be in the $C_H$-linear subspace spanned by
$g_1,\dots, g_n$.
Then the following are equivalent:
\begin{enumerate}
\item[\textup{(i)}] $h_1\prec\dots\prec h_n$;
\item[\textup{(ii)}] for $i=1,\dots, n$ there are $c_{ii}, c_{i,i-1},\dots, c_{i1}\in C_H$ such that
%there are  $c_{ij}\in C$  \textup{(}$1\leq j\le i\leq n$\textup{)}   with
$$h_i\ =\ c_{ii}g_i + c_{i,i-1}g_{i-1}+\cdots+c_{i1}g_1 \  \text{and} \  c_{ii}\neq 0.$$
\end{enumerate}
\end{lemma}

\noindent
Below $A\in H[\der]^{\neq}$ has order $r\geq 1$. Now the main results of this subsection: 

\begin{lemma}\label{lem:distinguished splitting, unique}
There is at most one splitting $(a_r,\dots,a_1)$ of $A$ over $H$ such that~$a_1\leq \cdots\leq a_r$.
\end{lemma}
\begin{proof}
Let  $(a_r,\dots,a_1)$, $(b_r,\dots,b_1)$ be   splittings of $A$ over $H$ with~$a_1\leq \cdots\leq a_r$ and~$b_1\leq \cdots\leq b_r$.  Towards showing that $a_i=b_i$ for $i=1,\dots,r$ we arrange that~$H$ is Liouville closed.
Then Lemma~\ref{lem:basis of kerUA, real, H-field} yields bases $y_1,\dots,y_r$ and $z_1,\dots,z_r$ of~$\ker_H A$
such that 
$\operatorname{split}(y_1,\dots,y_r)=(a_1,\dots,a_r)$,
$\operatorname{split}(z_1,\dots,z_r)=(b_1,\dots,b_r)$ and $y_i\nasymp y_j$, $z_i\nasymp z_j$ whenever $i\neq j$. 
By Lemma~\ref{lem:preceq vs leq,r} we have $y_1\prec\cdots\prec y_r$ and~$z_1\prec\cdots\prec z_r$, and hence
by Lemmas~\ref{lem:bases with same split} and~\ref{lem:upper triag},
\[ (a_1,\dots,a_r)\ =\ \operatorname{split}(y_1,\dots,y_r)\ =\ \operatorname{split}(z_1,\dots,z_r)\ =\ (b_1,\dots,b_r).  \qedhere\]
\end{proof}

\begin{example}%\label{ex:distinguished splitting, 1}
Let $a,b\in H$ in this example. Then
$$(\der-b)(\der-a)=\der^2-\der a-b\der+ab = \der^2-(a+b)\der+(ab-a'),$$
so for $f,g\in H$,
$$(\der-b)(\der-a)=(\der-g)(\der-f) \quad\Longleftrightarrow\quad
\text{$a+b=f+g$ and $ab-a'=fg-f'$.}$$
Now take $A=\der^2$. Then $1$, $x$ is a basis of $\ker_H A$, and
\begin{align*}
A=(\der-b)(\der-a)&\quad\Longleftrightarrow\quad
\text{$a+b=0$ and $ab-a'=0$}  \\
&\quad\Longleftrightarrow\quad \text{$a=-b=(cx+d)^\dagger$ for some $c,d\in C_H$, not both zero.}
\end{align*}
Hence if $(b,a)$ is  any splitting of $A$ over $H$, 
then $a\geq 0\geq b$, % and for each $d\in\R$, setting $a:=\frac{1}{x+d}$ we obtain the splitting $(-a,a)$ of $A$ over $H$.
and the only splitting~$(b,a)$ of $A$ over $H$ with $a\leq b$ is $(b,a)=(0,0)$.
\end{example}

\noindent
We call $A$   {\bf scrambled} if there are $y,z\in\ker_H^{\neq} A$ with $y\nasymp z$ and $y\asymp^\flat z$, and
{\bf unscrambled} otherwise. Hence if $r=1$, then $A$ is unscrambled, whereas $A=\der^2$ is scrambled.\index{linear differential operator!scrambled}\index{linear differential operator!unscrambled}\index{scrambled}\index{unscrambled}
For $a,b\in H^\times$ we have: $A$ is scrambled $\Leftrightarrow$ $aAb$ is scrambled. Moreover:

\begin{lemma} Assume $H$ has asymptotic integration, and let $B\in H[\der]$ have order $s\geq 1$. If $B$ is scrambled, then so is~$AB$.
If~$A$ is scrambled, $B$ is asymptotically surjective, $\exc^{\ev}(B)=v(\ker_H^{\neq} B)$, and $\ker_H A\subseteq B(H)$, then $AB$ is scrambled.
\end{lemma}
\begin{proof}
The first statement is clear since $\ker_H B\subseteq\ker_H AB$.
Suppose    $B$ is asymptotically surjective, $\exc^{\ev}(B)=v(\ker_H^{\neq} B)$, and $\ker_H A\subseteq B(H)$, and
let $f,g\in\ker_H^{\neq} A$ be such that~$f\nasymp g$ and $f\asymp^\flat g$.
Corollary~\ref{cor:nonexc sol} yields $y,z\in H$ with $B(y)=f$, $B(z)=g$
and $vy,vz\notin\exc^{\ev}(B)$. Then $y,z\in \ker_H^{\neq} AB$, and 
$y \nasymp z$, $y \asymp^\flat z$  by   Lemmas~\ref{lem:ADH 14.2.7} and~\ref{lem:v(A(y)) convex subgp}.
\end{proof}

\begin{prop}\label{prop:distinguished splitting}
Suppose  $H$ is Liouville closed, $A$ splits over~$H$, and $A$ is unscrambled.
Then there is a unique splitting $(a_r,\dots,a_1)$ of $A$ over $H$ such that~$a_1\leq \cdots\leq a_r$.
For this splitting we have $a_1<\cdots<a_r$.
\end{prop}
\begin{proof}
We first arrange that $A$~is monic.
The uniqueness part is immediate from Lem\-ma~\ref{lem:distinguished splitting, unique}.
To obtain a splitting  $(a_r,\dots,a_1)$ of $A$ over $H$ with~${a_1<\cdots<a_r}$, 
Lemma~\ref{lem:basis of kerUA, real, H-field} gives a basis $y_1\prec\cdots\prec y_r$ of the $C_H$-linear space $\ker_H A$.
Now set~$(a_1,\dots,a_r):=\operatorname{split}(y_1,\dots,y_r)$. Then $(a_r,\dots,a_1)$  is a splitting of~$A$ over $H$,
by Corollary~\ref{corbasissplit}. Since $A$ is unscrambled, we have $y_1\prec^\flat\cdots\prec^\flat y_r$, so~$a_1<\cdots<a_r$, by Lemma~\ref{lem:prec vs <,r, 2}. 
\end{proof}

\noindent
In \cite[Exercise~7.28]{JvdH} it is claimed that for the Liouville closed $H$-field $H=\mathbb T_{\g}$ of grid-based transseries, 
%\marginpar{maybe $\T_{\text{g}}$?} 
if $A$ splits over $H$, then $A$ always has a  
splitting~$(a_r,\dots,a_1)$ over $H$ with $a_1\leq \cdots\leq a_r$. The next example shows this to be incorrect for $r=2$:

\begin{exampleNumbered}\label{ex:distinguished splitting, 2}
Let $z\in H\setminus C_H$ and suppose  $A=(\der-g)\der$ where~$g:=z'{}^\dagger$.
Then~$1$,~$z$ is a basis of   $\ker_H A$.
With $a,b\in H$, this fact leads to the equivalence
%The computation in Example~\ref{ex:distinguished splitting, 1} yields 
$$A=(\der-b)(\der-a)\quad\Longleftrightarrow\quad \text{$a=g-b=(cz+d)^\dagger$ for some $c,d\in C_H$, not both zero.}$$
%{Next argument checked, but not included. For the direction $\Rightarrow$, let $A=(\der-b)(\der-a)$. Take $u\in U^\times$ with $u^\dagger=a$. Then $A(u)=0$, so $u\in \ker_U(A)=Cz+C$, hence $u=cz+d$ with $c,d\in C$ not both zero. If $c=0$, then $a=0$. If $c\ne 0$, then we can arrange $c=1$, so $a=z'/(z+d)\in H$, and thus $d\in C\cap H=C_H$.] 
Now take $z=x-x^{-1}$, so $z'=1+x^{-2}$, $z''=-2x^{-3}$, and hence
$$g\ =\ z'{}^\dagger\ =\  \frac{z''}{z'}\ =\ -\frac{2}{x(x^2+1)}<0.$$  
Let $(b,a)$ be a splitting of $A$ over $H$. We claim that then $a\geq 0>b$.  This is clear if $a=0$, so   assume $a\neq 0$.   
For $c,d\in C_H$, $c\neq 0$ we have $(cz+d)^\dagger\sim x^{-1}$.  Thus~$a\sim x^{-1}$ and  $b= g - a\sim -x^{-1}$. Hence $a>0>b$ as claimed.
\end{exampleNumbered}

\noindent
{\em In the rest of this subsection $H$ has asymptotic integration, and $\phi$ ranges over the elements of $H^>$ that are active in $H$}.
Note that if $A$ is unscrambled and~$\phi\preceq 1$, then $A^\phi\in H^\phi[\derdelta]$ is also unscrambled. Moreover:

\begin{lemma} 
$A^\phi$ is unscrambled, eventually.
\end{lemma}
\begin{proof}
By Remark~\ref{rem:kerexc} we have $\abs{v(\ker^{\neq}_H A)}\leq r$, thus we can take $\phi$ so that~${\gamma-\delta}\notin\Gamma^\flat_\phi$ for all $\gamma\neq\delta$ in $v(\ker^{\neq}_H A)$.
Now $A^\phi$ is unscrambled since~$\ker_H A=\ker_{H^\phi} A^\phi$.
\end{proof}

\begin{lemma}
Let $(a_r,\dots,a_1)$ be a splitting of $A$ over $H$ such that $a_1\leq \cdots\leq a_r$,   suppose $\phi\prec 1$, and
set $b_j:=\phi^{-1}\big(a_j-(j-1)\phi^\dagger\big)$ for $j=1,\dots,r$. Then
$(b_r,\dots,b_1)$ is a splitting of $A^\phi$ over $H^\phi$ with $b_1<\cdots < b_r$.
\end{lemma}
\begin{proof}
By Lemma~\ref{lem:split and compconj}, $(b_r,\dots,b_1)$ is a splitting of $A^\phi$ over $H^\phi$. Since $\phi^\dagger<0$, we have $b_1<\cdots < b_r$.
\end{proof}

\begin{cor}\label{cor:distinguished splitting}
Suppose   $H$ is Liouville closed and $A$ splits over~$H$.
Then there is a unique splitting $(a_r,\dots,a_1)$ of $A$ over $H$ such that eventually~$a_j+\phi^\dagger<a_{j+1}$ for $j=1,\dots,r-1$.
\end{cor}
\begin{proof}
Let $(a_r,\dots,a_1)$ be a splitting of $A$ over $H$ and $\phi$ so that $a_j+\phi^\dagger<a_{j+1}$ for~$j=1,\dots,r-1$. Set $b_j:=\phi^{-1}\big(a_j-(j-1)\phi^\dagger\big)$ for $j=1,\dots,r$. Then $(b_r,\dots,b_1)$  is a splitting of $A^\phi$ over $H^\phi$ with $b_1<\cdots < b_r$.
Thus by Lemma~\ref{lem:distinguished splitting, unique} there can be at most one splitting $(a_r,\dots, a_1)$
of $A$ over $H$ such that eventually~$a_j+\phi^\dagger<a_{j+1}$ for $j=1,\dots,r-1$. (Here we also use \eqref{eq:10.5.2}.)
% if $\phi_1\in H^\times$, $\phi_1\prec \phi$, then $\phi_1^\dagger < \phi^\dagger$.) 

For existence, take $\phi$ with unscrambled $A^\phi$ and a splitting $(b_r,\dots,b_1)$ of $A^\phi$ over~$H^\phi$
with $b_1<\cdots<b_r$ as in Proposition~\ref{prop:distinguished splitting}. For $j=1,\dots,r$, take $a_j\in H$ such that
$b_j=\phi^{-1}\big(a_j-(j-1)\phi^\dagger\big)$. Then $(a_r,\dots,a_1)$ is a splitting of $A$ over $H$, by Lemma~\ref{lem:split and compconj}, and $a_j+\phi^\dagger<a_{j+1}$ for $j=1,\dots,r-1$.
\end{proof}

\begin{example}
Suppose  $H$ is Liouville closed  and $A$, $g$, $z=x-x^{-1}$ are as in Example~\ref{ex:distinguished splitting, 2}. Then the unique splitting~$(b,a)$ of $A$ over $H$
such that  eventually~$a+\phi^\dagger < b$  is~$(b,a)=(g,0)$. (To see this use that eventually $\phi^\dagger\sim -x^{-1}$.)
\end{example}

%\noindent
%If $A$ splits over $H$, then $A$ splits over each  differential field extension of $H$. Also recall that under the hypotheses of the previous proposition, $A$ splits over $H$ iff  $0$ is an eigenvalue of $A$ of multiplicity $r$, by Corollary~\ref{cor:unique eigenvalue, 1}. We also note:

%\begin{lemma}
%Suppose $H^\dagger=H$ and $\I(K)\subseteq K^\dagger$. If  $0$ is an eigenvalue of $A$ of multiplicity $r$, then $0$ remains an eigenvalue of $A$ of multiplicity $r$ if $H$ is replaced by a real closed $H$-field extension.
%\end{lemma}
%\begin{proof}
%If $F$ is a real closed $H$-field extension of $H$, then  $F[\imag]^\dagger\cap K = K^\dagger$ by Lem\-ma~\ref{lem:W and I(F)} and Corollary~\ref{cor:logderset ext}. Now use the remarks before Lemma~\ref{lem:split evs}.
%\end{proof}

\noindent
We finish this subsection with a variant of Proposition~\ref{prop:distinguished splitting} for  P\'olya-style splittings.
In [ADH, 11.8]   we defined $\Upg(H):=\{h^\dagger:h\in H^{\succ 1}\}$. If $H$ is Liouville closed,
then $H^>\setminus\I(H)=\Upg(H)$ [ADH, p.~520]. 

\begin{prop}\label{prop:distinguished splitting, Trench}
{\samepage Suppose $H$ is Liouville closed and $A$ is monic  and splits over~$H$. Then there are $g_1,\dots,g_r\in H^>$ such that
$$A\  =\  g_1\cdots g_r (\der g_r^{-1})  \cdots (\der g_2^{-1})( \der g_1^{-1})\quad\text{and}\quad 
g_j\in\Upg(H)\text{ for $j=2,\dots,r$.}$$
Such $g_1,\dots, g_r$ are unique up to multiplication by positive constants.  }
\end{prop}
\begin{proof}
Let $(a_r,\dots,a_1)$ be the  splitting of $A$ over $H$ from Corollary~\ref{cor:distinguished splitting}.
Take~$g_j\in H^>$ such that $g_j^\dagger=a_j-a_{j-1}$ for $j=1,\dots,r$, where $a_0:=0$.
Then $$A = g_1\cdots g_r (\der g_r^{-1}) \cdots (\der g_2^{-1})( \der g_1^{-1})$$ by Lemma~\ref{lem:Polya fact}.
For $j=2,\dots,r$ we have $(g_j/\phi)^\dagger > 0$, eventually,    hence~$g_j/\phi \succeq 1$, eventually, so
$g_j\in H^>\setminus \I(H) = \Upg(H)$. 
Suppose now $h_1,\dots,h_r\in H^>$ are such that
$$A\  =\ h_1\cdots h_r (\der h_r^{-1})  \cdots (\der h_2^{-1})( \der h_1^{-1})\quad\text{and}\quad 
h_j\in\Upg(H)\text{ for $j=2,\dots,r$.}$$
Set $b_j:=(h_1\cdots h_j)^\dagger$ for $j=1,\dots,r$. Then $A = (\der-b_r) \cdots (\der-b_1)$ by Lemma~\ref{lem:Polya fact}. 
Also $\int h_j\succ 1$ for $j=2,\dots,r$, for any choice of the integrals in $H$. Let
$$z_1:=h_1,\quad z_2:=h_1\textstyle\int h_2,\quad \dots,\quad z_r:=h_1\int (h_2\int h_3(\cdots(h_{r-1}\int h_r)\cdots))$$
for some choice of the integrals in $H$. Then $z_1,\dots, z_r\in \ker^{\neq}_H A$, and induction on~$j=1,\dots,r$ using [ADH, 9.1.3(iii)] gives $z_1\prec\cdots\prec z_j$.
Then Lemma~\ref{lem:Polya fact, 2} yields   $\operatorname{split}(z_1,\dots,z_r)=(b_1,\dots,b_r)$.
Applied   to $g_1,\dots, g_r$ instead of $h_1,\dots, h_r$, this gives $y_1\prec \cdots \prec y_r\in \ker^{\neq}_H A$ such that 
$\operatorname{split}(y_1,\dots,y_r)=(a_1,\dots,a_r)$.  Now~$y_1,\dots, y_r$ and $z_1,\dots, z_r$ are both bases
of $\ker_H A$, so by Lemmas~\ref{lem:upper triag} and~\ref{lem:bases with same split}: 
$$(a_1,\dots,a_r)\  =\ \operatorname{split}(y_1,\dots,y_r)\  =\  \operatorname{split}(z_1,\dots,z_r)\  =\  (b_1,\dots,b_r).$$
So $g_j^\dagger=a_j-a_{j-1}=b_j-b_{j-1}=h_j^\dagger$ ($b_0:=0$), and thus $g_j \in C_H^>\, h_j$, $j=1,\dots,r$.
\end{proof}

\subsection*{The case of oscillating transseries\astr} 
We now apply the results in this section to the algebraically closed differential field $K=\T[\imag]$. Note that $\T[\imag]$ has constant field $\C$ and extends the (real closed) differential field $\T$ of transseries. \index{transseries!oscillating}   
%\marginpar{subsection not used later in this file, except for an example at beginning of Section~\ref{sec:valuniv}}
%After recalling a fact that is implicit in [ADH], we next derive a consequence of the above material for the algebraic closure $\T[\imag]$ of $\T$.

\begin{lemma}\label{lem:T[imag] lc and ls}
$\T[\imag]$ is linearly closed and linearly surjective.
\end{lemma}
\begin{proof}
By [ADH,  15.0.2], $\T$ is newtonian, so $\T[\imag]$ is newtonian by [ADH,  14.5.7]. Hence $\T[\imag]$ is linearly closed by [ADH,  5.8.9,  14.5.3], 
and $\T[\imag]$ is linearly surjective  by [ADH,  14.2.2].
\end{proof}

\noindent
Now  applying Corollary~\ref{cor:basis of kerUA} and Lemma~\ref{newlembasis} to $\T[\imag]$ gives:

\begin{cor}
For $K=\T[\imag]$, there are $\mathbb C$-linearly independent units $y_1,\dots,y_r$ of $\Univ_{\T[\imag]}$ with $A(y_1)=\cdots = A(y_r)=0$.
\end{cor}

\noindent
Next we describe another incarnation of $\Univ_{\T[\imag]}$, namely as a ring of ``oscillating'' transseries. Towards this goal we first note that by  
 [ADH, 11.5.1, 11.8.2] we have
\begin{align*}
\I(\T)\	&=\ \big\{y\in\T:\ \text{$y\preceq f'$ for some $f\prec 1$ in $\T$} \big\}\\
		&=\ \big\{y\in \T:\ \text{$y\prec 1/(\ell_0\cdots\ell_n)$ for all $n$} \big\},
\end{align*}
so a complement $\Lambda_{\T}$ of $\I(\T)$ in $\T$ is given by
$$
\Lambda_{\T}\	:=\ \big\{y\in\T:\  \text{$\supp(y) \succ 1/(\ell_0\cdots\ell_{n-1}\ell_n^2)$ for all $n$} \big\}.$$
Since $\T^\dagger=\T$ and $\I\!\big(\T[\imag]\big)\subseteq \T[\imag]^\dagger$ we
have $\T[\imag]^\dagger=\T\oplus \I(\T)\imag$ by Lemmas~\ref{lem:logder} and~\ref{lem:W and I(F)}.
We now take $\Lambda=\Lambda_{\T}\imag$ as our complement $\Lambda$ of $\T[\imag]^\dagger$ in $\T[\imag]$ and 
explain how the universal exponential extension~$\Univ$ of $\T[\imag]$ for this~$\Lambda$ was introduced in \cite[Section~7.7]{JvdH} in a different way.
Let $$\T_{\succ}\ :=\ \{f\in\T:\supp f\succ 1\},$$ and similarly with $\prec$ in place of $\succ$; then $\T_{\prec}=\smallo_{\T}$ and~$\T_{\succ}$ are
$\R$-linear subspaces of $\T$, and~$\T$ decomposes as an internal direct sum 
\begin{equation}\label{eq:T decomp}
\T\ =\ \T_{\succ}\oplus \R \oplus \T_{\prec}
\end{equation}
of $\R$-linear subspaces of $\T$.
Let $\ex^{\imag \T_{\succ}}=\{\ex^{\imag f}:f\in\T_{\succ}\}$ be a multiplicative copy of the additive group~$\T_{\succ}$, with isomorphism $f\mapsto \ex^{\imag f}$. Then we have the group ring
$$\mathbb{O}\ :=\ K\big[\ex^{\imag\T_{\succ}}\big]$$
of $\ex^{\imag\T_{\succ}}$ over $K=\T[\imag]$. 
We make $\mathbb{O}$ into a differential ring extension of $K$ by
$$(\ex^{\imag f})'\ =\  \imag f' \ex^{\imag f}\qquad (f\in\T_{\succ}).$$
Hence $\mathbb{O}$ is an exponential extension of $K$.
The elements of $\mathbb{O}$ are called {\it oscillating transseries.}\/ 
For each $f\in\T$ there is a unique $g\in\T$, to be denoted by $\int f$,
such that~$g'=f$ and $g$ has constant term $g_1=0$.
The injective map $\int\colon \T\to\T$ is $\R$-linear; we use this map
% $\T=(\int P)\oplus\R\oplus \T_{\prec}$.  
to show that $\Univ$ and $\mathbb{O}$ are disguised versions of each other:

\begin{prop}\label{prop:osc transseries}
There is a unique isomorphism $\Univ=K\big[\!\ex(\Lambda)\big] \to \mathbb{O}$ of differential $K$-algebras sending $\ex(h\imag)$ to $\ex^{\imag\int h}$ for all $h\in \Lambda_{\T}$.
\end{prop}

\noindent
This requires the next lemma. We assume familiarity with [ADH, Appendix~A], especially with the ordered group $G^{\operatorname{LE}}$ (a subgroup of $\T^\times$) of
logarithmic-exponential monomials and its subgroup
$G^{\operatorname{E}} = \bigcup_n G_n$ of exponential monomials.

\begin{lemma}\label{succP} If $\fm\in G^{\operatorname{LE}}$ and $\fm\succ 1$, then $\supp \fm'\ \subseteq\ \Lambda_{\T}$. 
\end{lemma}
\begin{proof} We first prove by induction on $n$ a fact about elements of $G^{\operatorname{E}}$: 
$$\text{ if $\fm\in G_n$, $\fm\succ 1$,
then $\supp \fm' \succ 1/x$}.$$ For $r\in\R^>$ we have $(x^r)'=rx^{r-1}\succ 1/x$,
so the claim holds for $n=0$. Suppose the claim holds for a certain $n$. Now
$G_{n+1}=G_n\exp(A_n)$, $G_n$ is a convex subgroup of $G_{n+1}$, and 
$$A_n\ =\ \big\{f\in \R[[G_n]]:\ \supp f\succ G_{n-1}\big\}\qquad\text{ (where $G_{-1}:=\{1\}$).}$$
Let $\fm=\fn\exp(a)\in G_{n+1}$ where $\fn\in G_n$, $a\in A_n$; then
$$\fm\succ 1\quad\Longleftrightarrow\quad\text{$a>0$, or $a=0$, $\fn\succ 1$.}$$
Suppose $\fm\succ 1$. If $a=0$, then $\fm=\fn$, and we are done by inductive hypothesis,
so assume $a>0$. Then $\fm' = (\fn'+\fn a')\exp(a)$ and 
$(\fn'+\fn a')\in \R[[G_n]]$, a differential subfield of $\T$, and
$\exp(a) > \R[[G_n]]$,  hence $\supp \fm' \succ 1\succ 1/x$ as required.

Next, suppose $\fm\in G^{\operatorname{LE}}$ and $\fm\succ 1$. Take $n\ge 1$ such that $\fm{\uparrow}^n\in G^{\operatorname{E}}$. 
We have~$(\fm{\uparrow}^n)'= (\fm'\cdot \ell_0\ell_1\cdots\ell_{n-1}){\uparrow}^n$. For $\fn\in\supp\fm'$ and using $\fm{\uparrow}^n\succ 1$
this gives
$$ (\fn\cdot \ell_0\ell_1\cdots\ell_{n-1} ){\uparrow}^n\ \succ\ 1/x$$
by what we proved for monomials in $G^{\operatorname{E}}$.  Applying ${\downarrow}_n$ this yields
$\fn\succ 1/(\ell_0\ell_1\cdots\ell_n)$,
hence $\fn\in \Lambda_{\T}$ as claimed.
\end{proof}

\begin{proof}[Proof of Proposition~\ref{prop:osc transseries}] Applying $\der$ to the decomposition~\eqref{eq:T decomp} gives $$\T\ =\ \der(\T_{\succ}) \oplus \der(\T_{\prec}).$$ Now $\der(\T_{\succ})\subseteq \Lambda_{\T}$ by Lemma~\ref{succP}, and $\der(\T_{\prec})\subseteq \I(\T)$, and so these two inclusions are equalities. Thus $\int \Lambda_{\T}= \T_{\succ}$, from which the proposition follows.
\end{proof}

\begin{prop} There is a unique group morphism $\exp \colon K=\T[\imag]\to \mathbb O^\times$
that extends the given exponential maps $\exp\colon \T\to \T^\times$ and
$\exp\colon \mathbb C \to \mathbb C^\times$, and such that~$\exp(\imag f)=\ex^{\imag f}$ for all $f\in \T_{\succ}$ and $\exp(\epsilon)=\sum_n \frac{\varepsilon^n}{n!}$ for all $\epsilon\in \smallo$. It is surjective, has kernel $2\pi\imag \Z\subseteq \mathbb C$, and satisfies $\exp(f)'=f'\exp(f)$ for all $f\in K$.   
\end{prop}
\begin{proof} The first statement follows easily from the decompositions
$$ K\ =\ \T \oplus \imag \T\ =\ \T \oplus \imag \T_{\succ} \oplus \imag \R \oplus \imag \smallo_{\T},\qquad \mathbb C\ =\ \R\oplus \imag \R, \qquad \smallo\ =\ \smallo_{\T}\oplus \imag \smallo_{\T}$$
of $K$, $\mathbb C$, and $\smallo=\smallo_K$ as internal direct sums of $\R$-linear subspaces. Next, 
$$\mathbb O^\times\ =\  K^\times \ex^{\imag \T_{\succ}}=\ \T^{>}\cdot S_{\mathbb C}\cdot (1+\smallo)\cdot\ex^{\imag\T_{\succ}}, \qquad S_{\mathbb C}\ :=\ \big\{z\in \mathbb{C}:\ |z|=1\big\},$$
by Lemmas~\ref{lem:only trivial units} and~\ref{lem:logder}, and Corollary~\ref{cor:decomp of S}.
Now $\T^{>}=\exp(\T)$ and $S_{\mathbb C}=\exp(\imag \R)$, so surjectivity follows from
$\exp(\smallo)=1+\smallo$, a consequence of
the well-known bijectivity of the map $\epsilon\mapsto \sum_n\frac{\varepsilon^n}{n!}\colon \smallo\to 1+\smallo$, whose inverse is given by
$$1+\delta\mapsto \log(1+\delta)  :=\sum_{n=1}^\infty \frac{(-1)^{n-1}}{n} \delta^n \qquad(\delta\in \smallo).$$
That the kernel is $2\pi\imag\Z$ follows from the initial decomposition of the additive group of $K$ as $\T \oplus \imag \T_{\succ} \oplus \imag \R \oplus \imag\smallo_{\T}$. The identity $\exp(f)'=f'\exp(f)$ for $f\in K$ follows from it being satisfied for $f\in\T$, $f\in\imag \T_{\succ}$, 
$f\in\mathbb C$, and $f\in\smallo$.  
\end{proof}  

\noindent
To integrate oscillating transseries, note first that the
$\R$-linear operator $\int\colon \T\to \T$ extends uniquely to a $\C$-linear operator $\int\colon \T[\imag] \to \T[\imag]$. This in turn extends  uniquely to a~$\C$-linear operator~$\int\colon \mathbb{O} \to \mathbb{O}$ such that~$(\int \Phi)'=\Phi$ for all $\Phi\in \mathbb{O}$ and
$\int \T[\imag] \ex^{\phi\imag}\subseteq \T[\imag]\ex^{\phi\imag}$ for all $\phi\in \T_{\succ}$: given~$\phi\in \T_{\succ}^{\ne}$ and~${g\in \T[\imag]}$, there is a unique
$f\in \T[\imag[$ such that $(f\ex^{\phi\imag})'=g\ex^{\phi\imag}$: existence holds because~${y'+y\phi'\imag=g}$ has a solution in $\T[\imag]$, the latter being linearly surjective, and uniqueness holds by Lemma~\ref{1K} applied to $K=L=\T[\imag]$, because $\phi'\imag\notin \T[\imag]^\dagger$ in view of remarks preceding Lemma~\ref{lem:W and I(F)}. See also Corollary~\ref{cor:Shackell}. 

\medskip
\noindent
The operator $\int$ is a right-inverse of
the linear differential operator~$\der$ on~$\mathbb{O}$. To extend this to other linear differential operators, make the subgroup~$G^{\mathbb{O}}:=G^{\operatorname{LE}}\ex^{\imag \T_{\succ}}$ of~$\mathbb{O}^\times$ into an ordered group so that the ordered subgroup $G^{\operatorname{LE}}$ of $\T^{>}$ is a convex ordered subgroup of $G^{\mathbb{O}}$ and 
$\ex^{\imag\phi} \succ G^{\operatorname{LE}}$ for $\phi>0$ in $\T_{\succ}$. (Possible in only one way.) Next, extend the natural inclusion 
$\T[\imag]\to \C[[G^{\operatorname{LE}}]]$ to a $\C$-algebra embedding~${\mathbb{O}\to \C[[G^{\mathbb{O}}]]}$ by sending
$\ex^{\imag\phi}\in \mathbb{O}$ to $\ex^{\imag\phi}\in G^{\mathbb{O}}\subseteq \C[[G^{\mathbb{O}}]]$. Identify $\mathbb{O}$
with a subalgebra of~$\C[[G^{\mathbb{O}}]] $ via this embedding, so $\supp f \subseteq G^{\mathbb{O}}$ for $f\in \mathbb{O}$. It makes the Hahn space~$\C[[G^{\mathbb{O}}]]$ over $\C$ an immediate extension of its valued subspace $\mathbb{O}$. The latter is in particular also a Hahn space over $\C$. 

%View $\T[\imag]$ as a valued vector space over $\C$ with valuation~$v\colon \T[\imag]\to \Gamma_\infty$. Let $S:= \T_{\succ}\times\Gamma$ equipped with the lexicographic ordering, and turn $\mathbb{O}$ into the
%valued vector space  $(\mathbb{O},S,v)$ over $\C$ with valuation  $v\colon\mathbb{O}\to S_\infty$ such that~$v(f\ex^{\phi\imag})=(\phi,vf)$ for $f\in\T[\imag]^\times$, $\phi\in\T_{\succ}$.
%Then $(\mathbb{O},S,v)$ is a Hahn space over $\C$. Using [ADH, 2.3.12] we identify this Hahn space with a valued subspace of the Hahn product $H[S,\C]$, making the latter an immediate extension of  $(\mathbb{O},S,v)$.

Let $A\in\T[\imag][\der]^{\neq}$. Then $A(\mathbb{O})=\mathbb{O}$ by Lemmas~\ref{lem:A(U)=U},~\ref{lem:T[imag] lc and ls}, and Proposition~\ref{prop:osc transseries}. 
The proof of [ADH, 2.3.22]  now  gives for each
$g\in \mathbb{O}$ a unique element ${f=:A^{-1}(g)\in\mathbb{O}}$ with $A(f)=g$ and
$\supp(f)\cap \fd\big(\!\ker_{\mathbb O}^{\neq} A\big)=\emptyset$.
This requirement on $\supp A^{-1}(g)$ yields a $\C$-linear   operator~$A^{-1}$ on $\mathbb O$  with~$A\circ A^{-1}=\id_{\mathbb O}$; we call it the {\bf distinguished} right-inverse of the operator $A$ on $\mathbb O$.  
With this definition~$\der^{-1}$ is the operator~$\int$ on $\mathbb{O}$ specified earlier.  

\medskip
\noindent
In the next section we explore various valuations on universal exponential extensions (such as $\mathbb{O}$)
with additional properties. 

\section{Valuations on the Universal Exponential Extension}\label{sec:valuniv}

\noindent
{\it In this section $K$ is a valued differential field with algebraically closed constant field~${C\subseteq\mathcal O}$ and divisible group $K^\dagger$ of
logarithmic derivatives.}\/
Then $\Gamma=v(K^\times)$ is also divisible,
since we have a group isomorphism $$va\mapsto a^\dagger+(\mathcal O^\times)^\dagger\ :\ \Gamma\to K^\dagger/(\mathcal O^\times)^\dagger \qquad (a\in K^\times).$$ 
Let $\Lambda$ be a complement of the $\Q$-linear subspace~$K^\dagger$ of $K$, let $\lambda$ range over $\Lambda$, 
let~$\Univ=K\big[\!\ex(\Lambda)\big]$ be the universal exponential extension of $K$ constructed in Section~\ref{sec:univ exp ext}
and set $\Omega:=\Frac(\Univ)$. Thus $\Omega$ is a differential field with
constant field $C$.

\subsection*{The gaussian extension}
We equip   $\Univ$
with the gaussian extension~$v_{\g}$ of the valuation of~$K$ as defined in Section~\ref{sec:group rings}; so for $f\in\Univ$ with spectral decomposition~$(f_\lambda)$: $$v_{\g}(f)\ =\ \min_\lambda v(f_\lambda),$$
and hence 
$$v_{\g}(f')\ =\ \min_\lambda v(f_\lambda'+\lambda f_\lambda).   $$
The field $\Omega$ with the valuation extending $v_{\g}$ is a valued differential field extension of $K$, but it can happen that 
$K$ has small derivation, whereas  $\Omega$  does not:

\begin{example} Let $K=C(\!( t^{\Q} )\!)$ and $\Lambda$ be as in Example~\ref{ex:Q}, so $t\prec 1\prec x=t^{-1}$ and~$t'=-t^2$. Then $K$ is $\d$-valued of $H$-type with small derivation, but in $\Omega$ with the above valuation, $$t\ex(x)\ \prec\ 1, \qquad \big(t\ex(x)\big){}'\ =\ -t^2\ex(x)+\ex(x)\ \sim\ \ex(x)\ \asymp\ 1.$$
To obtain an example where $K=H[\imag]$ for a Liouville closed $H$-field $H$ and $\imag^2=-1$, take $K:=\T[\imag]$ and $\Lambda:=\Lambda_{\T}\imag$ as at the end of Section~\ref{sec:eigenvalues and splitting}. Now $x\in \Lambda_{\T}$ and in $\Omega$ equipped with the above valuation we have for $t:=x^{-1}$:
$$t\ex(x\imag)\ \prec\ 1,\qquad \big(t\ex(x\imag)\big){}'\ =\ -t^2\ex(x\imag)+\imag\ex(x\imag)\ \sim\ \imag\ex(x\imag)\ \asymp\ 1,$$
so $\big(t\ex(x\imag)\big){}' \not\prec t^\dagger$, hence $\Omega$ is neither asymptotic nor has small derivation.
\end{example}

\noindent
However, we show next that under certain assumptions on $K$ with small derivation, $\Omega$ has also a valuation
which does  make $\Omega$ a valued differential field extension of~$K$ with small derivation.
For this we rely on results from [ADH, 10.4].
Although such a valuation is less canonical than $v_{\g}$, it is useful for  harnessing the finiteness statements about the set $\exc^{\ev}(A)$ of eventual exceptional values of $A\in K[\der]^{\neq}$
from  Section~\ref{sec:lindiff} to obtain similar facts about the set of {\it ultimate exceptional values}\/ of~$A$ introduced  later in this section.

\subsection*{Spectral extensions}
{\it In this subsection $K$ is $\d$-valued of $H$-type with $\Gamma\neq\{0\}$ and with small derivation.}\/
 
\begin{lemma}\label{specex}
The valuation of $K$ extends to a valuation on the field~$\Omega$ that makes~$\Omega$ a $\d$-valued extension of $K$ of $H$-type with small derivation.
\end{lemma}
\begin{proof} Applying [ADH, 10.4.7] to an algebraic closure of $K$ gives a $\d$-valued algebraically closed extension $L$ of $K$ of $H$-type with small derivation and $C_L=C$
such that~$L^\dagger\supseteq K$.
Let~$E:=\{y\in L^\times:\, y^\dagger\in K\}$, so $E$ is a subgroup of~$L^\times$, $E^\dagger=K$, and~$K[E]$ is an exponential extension of $K$ with $C_{K[E]}=C$. Then Corollary~\ref{corcharexp}
gives an embedding~${\Univ\to L}$ of differential $K$-algebras with image $K[E]$, which extends to an
embedding $\Omega\to L$ of differential fields. Using this embedding to transfer the valuation of $L$ to $\Omega$ gives a valuation as required.
\end{proof}

\noindent
A {\bf spectral extension} of the valuation of $K$ to $\Omega$ is a valuation on the field $\Omega$ with the properties stated in Lemma~\ref{specex}.\index{valuation!spectral extension}\index{spectral extension}\index{extension!spectral}
If $K$ is $\upo$-free, then so is $\Omega$ equipped with any spectral extension of the valuation of $K$,
by [ADH, 13.6] (and then $\Omega$ has rational asymptotic integration by [ADH, 11.7]).
We do not know whether this goes through with ``$\upl$-free'' instead of ``$\upo$-free''. Here is something weaker:

\begin{lemma}\label{lem:Omega as int} 
Suppose $K$ is algebraically closed and $\upl$-free. Then some spectral extension of the valuation of $K$ to $\Omega$  makes 
$\Omega$ a $\d$-valued field with divisible value group and asymptotic integration.
\end{lemma}
\begin{proof}
Take $L$, $E$ and an embedding $\Omega\to L$   as in the proof  of Lemma~\ref{specex}. Use this embedding  to identify~$\Omega$ with a differential subfield of~$L$, so $U=K[E]$ and $\Omega=K(E)$, and
equip $\Omega$ with the spectral extension of the valuation of $K$ obtained  by restricting the valuation of $L$ to $\Omega$.
Since $L$ is algebraically closed, $E$ is divisible, and $\Gamma_L=\Gamma+v(E)$ by [ADH, 10.4.7(iv)]. So~$\Gamma_\Omega=\Gamma_L$ is divisible.
Let~$a\in K^\times$, $y\in E$. Then $K(y)$  
has  asymptotic integration by Proposition~\ref{prop:upl-free and as int}, hence $v(ay)\in (\Gamma_{K(y)}^{\neq})'\subseteq (\Gamma_\Omega^{\neq})'$. Thus $\Omega$ has asymptotic integration.
\end{proof} 

\noindent
In the rest of this subsection $\Omega$ is equipped with a spectral extension $v$ (with value group $\Gamma_{\Omega}$) of the valuation of~$K$. The proof of Lemma~\ref{specex} and [ADH, 10.4.7] show that we can choose $v$ so that $\Psi_\Omega\subseteq\Gamma$;
but under suitable hypotheses on $K$, this is automatic: 

\begin{lemma}\label{lem:v(ex(Q))}
Suppose   $K$ has asymptotic integration and $\I(K)\subseteq K^\dagger$. 
Then ${\Psi_{\Omega}\subseteq \Gamma}$, the group morphism
\begin{equation}\label{eq:v(ex(q))}
\lambda\mapsto v\big(\!\ex(\lambda)\big)\ \colon\ \Lambda\to\Gamma_\Omega
\end{equation}
is injective, and $\Gamma_\Omega$ is divisible with $\Gamma_\Omega=\Gamma\oplus v\big(\!\ex(\Lambda)\big)$ \textup{(}internal direct sum of $\Q$-linear subspaces of $\Gamma_\Omega$\textup{)}. Moreover,
$\Psi_\Omega=\Psi^{\downarrow}$ in $\Gamma$. 
%is closed downward in $\Gamma$, $\Psi$ is cofinal in $\Psi_\Omega$, and~$\Omega$ is ungrounded.
\end{lemma}
\begin{proof}
For $a\in K^\times$ we have $\big(a\ex(\lambda)\big){}^\dagger=a^\dagger+\lambda\in K$, and 
if $a\ex(\lambda)\asymp 1$, then 
$$a^\dagger+\lambda\ =\ \big(a\ex(\lambda)\big){}^\dagger\in (\mathcal O^\times_\Omega)^\dagger\cap K\ \subseteq\ \I(\Omega)\cap K\ =\ \I(K),$$ 
so $\lambda\in \Lambda\cap\big(\I(K)+K^\dagger\big)=\Lambda\cap K^\dagger=\{0\}$ and  $a\asymp 1$.  Thus for $a_1,a_2\in K^\times$ and 
distinct $\lambda_1, \lambda_2\in \Lambda$ we have $a_1\ex(\lambda_1)\nasymp a_2\ex(\lambda_2)$, and so
for $f\in\Univ$ with spectral decomposition $(f_\lambda)$ we have
$vf=\min_\lambda v\big(f_\lambda\ex(\lambda)\big)$. Hence
$$\Psi_{\Omega}\ \subseteq 
\big\{v(a^\dagger+\lambda):\ a\in K^\times,\ \lambda\in \Lambda\big\}\ =\ v(K)\ =\ \Gamma_{\infty},$$ the map~\eqref{eq:v(ex(q))} is injective and 
$\Gamma\cap v\big(\!\ex(\Lambda)\big) = \{0\}$, and so 
$\Gamma_\Omega=\Gamma\oplus v\big(\!\ex(\Lambda)\big)$ (internal direct sum of subgroups
of $\Gamma_\Omega$). Since $\Gamma$ and $\Lambda$ are divisible, so is $\Gamma_\Omega$. 
Now $\Psi_\Omega=\Psi^{\downarrow}$  follows from $K=(\Univ^\times)^\dagger\subseteq \Omega^\dagger$
 and $K$ having asymptotic integration.
\end{proof}

%\noindent
%If $K$ is $\upo$-free, then so is $\Omega$ by [ADH, 13.6] (and then $\Omega$ has asymptotic integration by [ADH, 11.7]). (Does it follow  from $K$ being $\upl$-free that~$\Omega$ is $\upl$-free?)

%\medskip
%\noindent
%Let $L$ be a $\d$-valued  field extension of $K$ of $H$-type where $C_L$ is algebraically closed  and $L^\dagger$ is divisible   with $L^\dagger\cap K=K^\dagger$. Let $Q_L$ and $\ex_L$ be as in Remark~\ref{rem:Univ under d-field ext}. Then $\Univ=K\big[\!\ex(Q)\big]$ is a differential subring of $\Univ_L:=L\big[\!\ex_L(Q_L)\big]$,  hence $\Omega$ is a differential subfield of $\Omega_L:=\Frac(\Univ_L)$. Each spectral extension of the valuation of~$L$ to a valuation on $\Omega_L$ restricts to a spectral extension of the valuation of $K$ to~$\Omega$.

%\medskip
\noindent
We can now improve on Lemma~\ref{newlembasis}: 

\begin{cor}\label{cor:spectral valuation basis}
Suppose $K$ has asymptotic integration and $\I(K)\subseteq K^\dagger$, and let~$A\in K[\der]^{\neq}$. Then the $C$-linear space $\ker_{\operatorname{U}} A$ has a basis
$\mathcal B\subseteq \operatorname{U}^\times$ such that $v$ is injective on $\mathcal B$ and
$v(\mathcal B)=v(\ker^{\neq}_{\operatorname{U}} A)$, and thus $\abs{v(\ker^{\neq}_{\operatorname{U}} A)} =\dim_C \ker_{\operatorname{U}} A$.  
\end{cor}
\begin{proof} By [ADH, 5.6.6] we have a basis $\mathcal B_\lambda$ of the $C$-linear space $\ker_K A_\lambda$ such that~$v$ is injective on $\mathcal{B}_{\lambda}$ and $v(\mathcal B_\lambda)=v(\ker^{\neq}_K A_\lambda)$.
%=\mult_\lambda(A)$.
Then $\mathcal B:=\bigcup_\lambda \mathcal B_\lambda\ex(\lambda)$ is a basis of $\ker_{\operatorname{U}} A$. It has the desired properties by Lemma~\ref{lem:v(ex(Q))}.
\end{proof}

\begin{cor}\label{cor:Omega as int}
Suppose $K$ is $\upl$-free and $\I(K)\subseteq K^\dagger$. Then  $\Omega$ has asymptotic integration, and so its $H$-asymptotic couple is closed by Lemma~\ref{lem:v(ex(Q))}. 
\end{cor}
\begin{proof}
By Lemma~\ref{lem:v(ex(Q))}, $\Gamma_\Omega= \Gamma+ v\big(\!\ex(\Lambda)\big)$. Using Proposition~\ref{prop:upl-free and as int} as in the proof of Lemma~\ref{lem:Omega as int}, with $\ex(\Lambda)$ in place of $E$, shows $\Omega$ has asymptotic integration.
\end{proof}

\subsection*{An application\astr} We use spectral extensions to prove an analogue of [ADH, 16.0.3]:

\begin{theorem}\label{noextension} If $K$ is an $\upo$-free newtonian $\d$-valued field, then $K$ has no proper $\d$-algebraic 
$\d$-valued field extension $L$ of $H$-type with $C_L=C$ and $L^\dagger\cap K=K^\dagger$. 
\end{theorem} 

\noindent
We retain of course our  assumption that $C$ is algebraically closed and $K^\dagger$ is divisible. In the same way that [ADH,   16.0.3] follows from [ADH,   16.1.1], Theorem~\ref{noextension} follows from an  analogue of [ADH, 16.1.1]:

\begin{lemma}\label{lem:ADH 16.1.1, Calinfi}
Let $K$ be an $\upo$-free newtonian $\d$-valued field, $L$ a $\d$-valued field extension of $K$ of $H$-type with
$C_L=C$ and $L^\dagger\cap K=K^\dagger$, and let $f\in L\setminus K$. Suppose there is no~$y\in K\langle f\rangle\setminus K$ such that $K\langle y\rangle$ is an immediate extension of~$K$.
Then the $\Q$-linear space~$\Q\Gamma_{K\langle f\rangle}/\Gamma$ is infinite-dimensional.
\end{lemma} 

\noindent
The proof of Lemma~\ref{lem:ADH 16.1.1, Calinfi} is much like that of [ADH,   16.1.1], except where
the latter uses that any $b$ in a Liouville closed $H$-field equals $a^\dagger$ for some nonzero $a$ in that field. This might not work with elements of $K$, and the remedy is to take instead for every $b\in K$ an element $a$ in $\Univ^\times$ with $b=a^\dagger$.
The relevant computation should then take place in the differential fraction field $\Omega_L$ of $\Univ_L$ instead of in $L$ where $\Omega_L$ is equipped with a spectral extension of the valuation of $L$. For all this  to make sense, we first take an active $\phi$ in $K$ and
replace $K$ and $L$ by $K^\phi$ and $L^\phi$, arranging in this way that the derivation of $L$ (and of $K$) is small. Next we
replace~$L$ by its algebraic closure, so that $L^\dagger$ is divisible, while preserving $L^\dagger\cap K=K^\dagger$ by Lemma~\ref{lem:Kdagger alg closure}, and also preserving the other conditions on $L$ in Lemma~\ref{lem:ADH 16.1.1, Calinfi}, as well as the derivation of $L$ being small.  This allows us to
identify  $\Univ$ with a differential subring of $\Univ_L$ as in Lemma~\ref{lem:Univ under d-field ext}, and accordingly $\Omega$ with a differential subfield of $\Omega_L$.  We equip~$\Omega_L$ with a spectral extension of the valuation of $L$ (possible by Lemma~\ref{specex}), and make $\Omega$ a valued subfield of $\Omega_L$. Then the valuation of $\Omega$ is a spectral extension of the valuation of $K$ to $\Omega$, so we have
the following inclusions of $\d$-valued fields:
$$\xymatrix{L \ar[r] & \Omega_L \\
K \ar[u] \ar[r] & \Omega \ar[u]}$$
With these preparations we can now give  the proof of Lemma~\ref{lem:ADH 16.1.1, Calinfi}:

\begin{proof} As we just indicated we arrange that $L$ is algebraically closed with small derivation,  and with an inclusion diagram of $\d$-valued fields involving $\Omega$ and $\Omega_L$, as above. (This will not be used until we arrive at the Claim below.)

By [ADH, 14.0.2], $K$ is asymptotically $\d$-algebraically maximal. Using this and the assumption about $K\langle f\rangle$ it follows
as in the proof of [ADH, 16.1.1] that
there is no divergent pc-sequence in $K$
with a pseudolimit in~$K\langle f\rangle$. Thus every~$y$ in~$K\langle f \rangle \setminus K$ has a 
 {\it a  best approximation in $K$,}\/ that is, an element $b\in K$ such that  $v(y-b)=\max v(y-K)$.
For such $b$ we have $v(y-b)\notin \Gamma$, since $C_L=C$.

Now pick a best approximation $b_0$ in $K$ to $f_0:=f$, and set
$f_1:= (f_0-b_0)^\dagger$. Then~$f_1\in K\langle f\rangle \setminus K$, since $L^\dagger\cap K=K^\dagger$ and $C=C_L$. Thus
$f_1$ has a best approximation $b_1$ in $K$, and continuing this way, we
obtain a sequence $(f_n)$ in~${K\langle f \rangle\setminus K}$ and a sequence
$(b_n)$ in $K$, such that $b_n$ is a best approximation in $K$ to~$f_n$ and~$f_{n+1}=(f_n-b_n)^\dagger$ for all $n$. Thus
$v(f_n-b_n)\in \Gamma_{K\langle f \rangle}\setminus \Gamma$ for all $n$. 

\claim{$v(f_0-b_0),v(f_1-b_1), v(f_2-b_2),\dots$ are $\Q$-linearly independent over $\Gamma$.}

\noindent
To prove this claim, take $a_n\in \Univ^\times$ with $a_n^\dagger=b_n$ for
$n\ge 1$. Then in $\Omega_L$,
$$f_n-b_n\ =\ (f_{n-1}-b_{n-1})^\dagger-a_n^\dagger\ =\ \left(\frac{f_{n-1}-b_{n-1}}{a_n}\right)^\dagger \qquad (n\ge 1).$$
With $\psi:= \psi_{\Omega_L}$ and $\alpha_n=v(a_n)\in \Gamma_{\Omega}\subseteq \Gamma_{\Omega_L}$ for $n\ge 1$, we get 
\begin{align*} v(f_n-b_n)\ &=\ \psi\big(v(f_{n-1}-b_{n-1})-\alpha_n\big),\ 
 \text{ so by an easy induction on $n$,}\\
 v(f_n-b_n)\ &=\ \psi_{\alpha_1,\dots, \alpha_n}\big(v(f_0-b_0)\big),  \qquad (n\ge 1).
 \end{align*} 
Suppose towards a contradiction that $v(f_0-b_0), \dots, v(f_n-b_n)$ are $\Q$-linearly dependent
over $\Gamma$. Then we have $m<n$ and $q_1,\dots, q_{n-m}\in \Q$ such that
$$v(f_m-b_m) + q_1v(f_{m+1}-b_{m+1}) + \cdots + q_{n-m}v(f_n-b_n)\in \Gamma.$$ 
For $\gamma:= v(f_m-b_m)\in \Gamma_L\setminus\Gamma$ this gives
$$\gamma + q_1\psi_{\alpha_{m+1}}(\gamma) + \cdots + q_{n-m}\psi_{\alpha_{m+1},\dots,\alpha_n}(\gamma)\in \Gamma.$$
By Lemma~\ref{lem:ADH 14.2.5} we have $\I(K)\subseteq K^\dagger$,  so the $H$-asymptotic couple of $\Omega$ is closed  with $\Psi_\Omega\subseteq\Gamma$,
by Lemma~\ref{lem:v(ex(Q))} and Corollary~\ref{cor:Omega as int}. Hence $\gamma\in\Gamma_\Omega$ by [ADH, 9.9.2].
 Together with $\Psi_\Omega\subseteq\Gamma$ and~$\alpha_{m+1},\dots,\alpha_n\in\Gamma_\Omega$ this gives 
$$\psi_{\alpha_{m+1}}(\gamma),\dots, \psi_{\alpha_{m+1},\dots,\alpha_n}(\gamma)\in\Gamma$$
and thus $\gamma\in\Gamma$, a contradiction.
\end{proof}

\subsection*{Ultimate exceptional values}
{\it In  this subsection $K$ is $H$-asymptotic with small derivation and asymptotic integration.}\/
Also  $A\in K[\der]^{\neq}$ and $r:=\order(A)$, and $\gamma$ ranges over $\Gamma=v(K^\times)$.
We have $v(\ker^{\neq} A_\lambda)\subseteq \exc^{\ev}(A_\lambda)$, so if
$\lambda$ is an eigenvalue of~$A$ with respect to $\lambda$, then $\exc^{\ev}(A_\lambda)\neq\emptyset$.  
We call the elements of the set
$$\exc^{\operatorname{u}}(A)\ =\ \exc^{\operatorname{u}}_{K}(A)\ :=\ \bigcup_\lambda\, \exc^{\ev}(A_\lambda)\ = \ \big\{ \gamma :\  \text{$\nwt_{A_\lambda}(\gamma)\geq 1$ for some $\lambda$} \big\}$$
the {\bf ultimate exceptional values of $A$} with respect to $\Lambda$. \label{p:excu} The definition of $\exc^{\operatorname{u}}_{K}(A)$ involves our choice of $\Lambda$,  but we are leaving this implicit to avoid complicated notation.
In Section~\ref{sec:ultimate} we shall restrict $K$  and $\Lambda$ so that~$\exc^{\operatorname{u}}(A)$ does not depend any longer on the choice of $\Lambda$. There we shall use the following observation: \index{linear differential operator!ultimate exceptional values}\index{exceptional values!ultimate}\index{ultimate!exceptional values}

\begin{lemma}\label{lem:excu for different Q}
Let $a,b\in K$ be such that $a-b\in (\mathcal O^\times)^\dagger$. Then for all $\gamma$ we have~$\nwt_{A_a}(\gamma)=\nwt_{A_b}(\gamma)$; in particular, $\exc^{\ev}(A_a)=\exc^{\ev}(A_b)$.
\end{lemma} 
\begin{proof}
Use that if $u\in\mathcal O^\times$ and $a-b=u^\dagger$, then $A_a=(A_b)_{\ltimes u}$.
\end{proof}

\begin{cor}\label{cor:excu for different Q}
Let $\Lambda^*$ be a complement of the $\Q$-linear subspace $K^\dagger$ of $K$ and let~$\lambda\mapsto \lambda^*\colon \Lambda \to \Lambda^*$ be the group isomorphism with $\lambda- \lambda^*\in K^\dagger$ for all $\lambda$. If $\lambda-\lambda^*\in (\mathcal O^\times)^\dagger$ for all $\lambda$, then $\nwt_{A_\lambda}(\gamma)=\nwt_{A_{\lambda^*}}(\gamma)$ for all $\gamma$, so $\exc^{\operatorname{u}}(A)= \bigcup_\lambda\, \exc^{\ev}(A_{\lambda^*})$. 
\end{cor}

\begin{remarkNumbered}\label{rem:excu for different Q, 1}
For $a\in K^\times$ we have 
$\exc^{\operatorname{u}}(aA)=\exc^{\operatorname{u}}(A)$
and $\exc^{\operatorname{u}}(Aa)=\exc^{\operatorname{u}}(A)-va$.
Note also that~$\exc^{\ev}(A) = \exc^{\ev}(A_0) \subseteq \exc^{\operatorname{u}}(A)$.
Let $\phi\in K^\times$ be active in $K$, and set~$\lambda^\phi:=\phi^{-1}\lambda$.
Then~$\Lambda^\phi:=\phi^{-1}\Lambda$ is a complement of the $\Q$-linear subspace $(K^\phi)^\dagger=\phi^{-1}K^\dagger$  of~$K^\phi$,
and $(A^\phi)_{\lambda^\phi} = (A_\lambda)^\phi$. Hence $\exc^{\operatorname{u}}_{K}(A)$ agrees with   the set~$\exc^{\operatorname{u}}_{K^\phi}(A^\phi)$ of ultimate exceptional values of $A^\phi$ with respect to $\Lambda^\phi$. 
\end{remarkNumbered}

\begin{remarkNumbered}\label{rem:excu for different Q, 2}
Suppose $L$ is an $H$-asymptotic extension of $K$ with asymptotic integration and algebraically closed constant field $C_L$ such that 
$L^\dagger$ is divisible, and~$\Psi$ is cofinal in~$\Psi_L$ or $K$ is $\upl$-free.
Then~$\exc^{\ev}(A_\lambda)=\exc^{\ev}_L(A_\lambda)\cap\Gamma$, by Lemma~\ref{lemexc} and Corollary~\ref{cor:13.7.10}. Hence if $\Lambda_L\supseteq \Lambda$ is a complement of  the  subspace~$L^\dagger$ of the $\Q$-linear space $L$, and~$\exc^{\operatorname{u}}_{L}(A)$ is  the set of ultimate exceptional values of $A$ (viewed as an element of~$L[\der]$) with respect to~$\Lambda_L$, then $\exc^{\operatorname{u}}(A)\subseteq \exc^{\operatorname{u}}_{L}(A)$. (Note that such a complement~$\Lambda_L$ exists iff~$L^\dagger\cap K=K^\dagger$.)
\end{remarkNumbered}

%\begin{remark} \marginpar{the first sentence has been checked. The meaning of the last sentence is unclear}
%By [ADH, 10.1.3], $K$ is pre-$\d$-valued, hence it  has an $\upo$-free algebraically closed $\d$-valued 
%extension $L$ such that $\I(L)\subseteq L^\dagger$ by [ADH, 11.7.18, 11.7.23, remarks after 14.0.1,  14.5.7, 14.2.5].
%It seems plausible that if $\I(K)\subseteq K^\dagger$, then
%$L$ can be chosen so that   additionally  $L^\dagger\cap K=K^\dagger$;
%we have not pursued this here, but under extra conditions Corollary~\ref{cor:LambdaL} below gives more.
%\end{remark}   

\noindent
In the rest of this subsection  we equip $\Univ$ with the gaussian extension $v_{\g}$ of the valuation of~$K$.
Recall that we have a decomposition $\ker_{\Univ} A = \bigoplus_\lambda (\ker A_\lambda)\ex(\lambda)$ of the $C$-linear space~$\ker_{\Univ} A$
as an internal direct sum of subspaces, 
and hence
\begin{equation}\label{eq:vU(ker A)}
v_{\g}(\ker_{\Univ}^{\neq} A)\ =\ \bigcup_\lambda\, v(\ker^{\neq} A_\lambda)\ \subseteq\ \bigcup_\lambda\, \exc^{\ev}(A_\lambda)\ =\ \exc^{\operatorname{u}}(A).
\end{equation}
Here are some consequences:

\begin{lemma}\label{lem:excu 1}
Suppose $K$ is $r$-linearly newtonian. Then $v_{\g}(\ker_{\Univ}^{\neq} A)=\exc^{\operatorname{u}}(A)$.
\end{lemma}
\begin{proof}
By Proposition~\ref{kerexc} we have $v(\ker^{\neq} A_\lambda) = \exc^{\ev}(A_\lambda)$ for each $\lambda$. There\-fore~$v_{\g}(\ker_{\Univ}^{\neq} A)=\exc^{\operatorname{u}}(A)$ by \eqref{eq:vU(ker A)}.
\end{proof}

\begin{lemma}\label{lem:excu 2}
Suppose $K$ is $\d$-valued. Then  $\abs{v_{\g}(\ker_{\Univ}^{\neq} A)}\leq \dim_C \ker_{\Univ} A\leq r$.
\end{lemma}
{\samepage
\begin{proof}
By [ADH, 5.6.6(i)] applied to $A_\lambda$ in place of $A$ we have
$$\abs{v(\ker^{\neq} A_\lambda)}\ =\ \dim_C \ker A_\lambda\ =\ \mult_\lambda(A)\qquad\text{ for all~$\lambda$}$$ 
and thus by \eqref{eq:vU(ker A)},
$$\abs{v_{\g}(\ker_{\Univ}^{\neq} A)}\ \leq\ 
\sum_\lambda \,\abs{v(\ker^{\neq} A_\lambda)}\ =\ 
\sum_\lambda \,\mult_\lambda(A)\ =\ \dim_C \ker_{\Univ} A\ \leq\ r$$
as claimed.
\end{proof}
}

\begin{lemma}\label{lem:excu, r=1} 
Suppose  $\I(K)\subseteq K^\dagger$ and $r=1$. Then $$v_{\operatorname{g}}(\ker^{\neq}_{\Univ} A)\ =\ \exc^{\operatorname{u}}(A), \qquad \abs{\exc^{\operatorname{u}}(A)}\ =\ 1.$$
\end{lemma}
\begin{proof}
Arrange $A=\der-g$, $g\in K$, and take $f\in K^\times$ and $\lambda$ such that~$g=f^\dagger+\lambda$. Then $u:=f\ex(\lambda)\in\Univ^\times$
satisfies $A(u)=0$, hence $\ker^{\neq}_{\Univ}A=Cu$ and thus~$v_{\text{g}}(\ker^{\neq}_{\Univ} A) = \{ vf \}$.
By Lemma~\ref{lem:v(ker)=exc, r=1} we have $v(\ker^{\neq} A_\lambda)=\exc^{\ev}(A_\lambda)$ for all $\lambda$ and
hence $v_{\operatorname{g}}(\ker^{\neq}_{\Univ} A) = \exc^{\operatorname{u}}(A)$ by \eqref{eq:vU(ker A)}.
\end{proof}

\begin{cor}\label{corevisu} If $\I(K)\subseteq K^\dagger$ and $a\in K^\times$, then $\exc^{\ev}(\der-a^\dagger)=\exc^{\operatorname{u}}(\der-a^\dagger)=\{va\}$. 
\end{cor} 

\noindent
Proposition~\ref{prop:finiteness of excu(A)} below partly extends Lemma~\ref{lem:excu, r=1}.

\subsection*{Spectral extensions and ultimate exceptional values}
{\it In this subsection $K$ is $\d$-valued of $H$-type with small derivation, asymptotic integration, and $\I(K)\subseteq K^\dagger$.}\/ Also $A\in K[\der]^{\ne}$ has order $r$ and $\gamma$ ranges over $\Gamma$. 

Suppose~$\Omega$ is equipped with a spectral extension $v$ of the valuation of $K$. Let ${g\in K^\times}$ with $vg=\gamma$.
The Newton weight of $A_\lambda g\in K[\der]$ does not change in passing from~$K$ to $\Omega$, since
$\Psi$ is cofinal in $\Psi_\Omega$ by Lemma~\ref{lem:v(ex(Q))}; see [ADH, 11.1]. 
Thus
$$\nwt_{A_\lambda}(\gamma)\ =\ \nwt(A_\lambda g)\ =\ \nwt\!\big(Ag\ex(\lambda)\big)\ =\ \nwt_A\!\big(v(g\ex(\lambda)\big)\ =\ \nwt_A\!\big(\gamma+v(\ex(\lambda)\big)\big).$$
In particular, using $\Gamma_{\Omega}=\Gamma\oplus v\big(\!\ex(\Lambda)\big)$,
\begin{equation}\label{eq:excevOmega}
\exc^{\ev}_\Omega(A)\ =\ \bigcup_\lambda \, \exc^{\ev}(A_\lambda)+v\big(\!\ex(\lambda)\big)\qquad\text{(a disjoint union)}.
\end{equation}
%Since $\exc^{\operatorname{u}}(A)=\bigcup_\lambda \exc^{\ev}(A_\lambda)$, this gives 
Thus $\exc^{\operatorname{u}}(A)=\pi\big(\exc^{\ev}_\Omega(A)\big)$ where $\pi\colon \Gamma_{\Omega}\to \Gamma$ is given by
~$\pi\big({\gamma+v\big(\!\ex(\lambda)}\big)\big)=\gamma$. 
%the obvious projection map.  

\begin{lemma}\label{lem:finiteness of excu(A)}
We have $\dim_C \ker_{\Univ} A \leq \sum_\lambda\abs{\exc^{\ev}(A_\lambda)}$, and
$$\dim_C \ker_{\Univ} A\ =\  \sum_\lambda\,\abs{\exc^{\ev}(A_\lambda)}
\ \Longleftrightarrow\  v(\ker^{\neq} A_\lambda)\ =\  \exc^{\ev}(A_\lambda)\text{ for all $\lambda$.}$$
Moreover, if  $\dim_C \ker_{\Univ} A =  \sum_\lambda\,\abs{\exc^{\ev}(A_\lambda)}$, then $v_{\g}(\ker_{\Univ}^{\neq} A) = \exc^{\operatorname{u}}(A)$.
\end{lemma}
\begin{proof}
Clearly, $\dim_C \ker_{\Univ} A \leq \dim_C \ker_\Omega A$.
Equip~$\Omega$ with a spectral extension of the valuation of $K$. 
Then $\dim_C\ker_\Omega A=\abs{v(\ker_\Omega^{\neq}A)}$ and $v(\ker_\Omega^{\neq} A)\subseteq \exc^{\ev}_\Omega(A)$ by~[ADH, 5.6.6(i)] and [ADH, p.~481], respectively, applied to $\Omega$ in the role of~$K$.
Also  $|\exc^{\ev}_{\Omega}(A)|=\sum_\lambda|\exc^{\ev}(A_\lambda)|$ (a sum of cardinals) by the remarks preceding the lemma. This yields the first claim of the lemma. 

Next, note that~$v(\ker^{\neq} A_\lambda)\subseteq\exc^{\ev}(A_\lambda)$ for all $\lambda$.
Hence from~\eqref{eq:excevOmega} and 
$$v(\ker_{\Univ}^{\neq} A)\ =\ \bigcup_\lambda v(\ker^{\neq}A_\lambda)+v\big(\!\ex(\lambda)\big)\qquad\text{(a disjoint union)}$$
we obtain: 
$$v(\ker_{\Univ}^{\neq} A)=\exc^{\ev}_\Omega(A)\quad\Longleftrightarrow\quad v(\ker^{\neq} A_\lambda)\ =\  \exc^{\ev}(A_\lambda) \text{ for all $\lambda$.}$$
Also 
 $\abs{v(\ker_{\Univ}^{\neq} A)}=\dim_C \ker_{\Univ} A$ by~[ADH, 2.3.13], 
and $$v(\ker_{\Univ}^{\neq} A)\ \subseteq\ v(\ker_\Omega^{\neq} A)\ \subseteq\ \exc^{\ev}_\Omega(A), \qquad
|\exc^{\ev}_{\Omega}(A)|\ =\ \sum_\lambda|\exc^{\ev}(A_\lambda)|.$$  
This yields the displayed equivalence. 

Suppose $\dim_C \ker_{\Univ} A= \sum_\lambda\abs{\exc^{\ev}(A_\lambda)}$;
we need to show $v_{\g}(\ker_{\Univ}^{\neq} A)= \exc^{\operatorname{u}}(A)$.
We have~$\pi\big(\exc^{\ev}_\Omega(A)\big)=\exc^{\operatorname{u}}(A)$ for the above projection map $\pi\colon \Gamma_{\Omega} \to \Gamma$, so it is enough to show
$\pi\big(v(\ker_{\Univ}^{\neq} A)\big)=v_{\g}(\ker_{\Univ}^{\ne} A)$.
For that, note that for $\mathcal{B}\subseteq K^\times \ex(\Lambda)$ in Corollary~\ref{cor:spectral valuation basis} we have
$$\pi\big(v(\ker_{\Univ}^{\neq} A)\big)\ =\ \pi(v\mathcal B)\ =\ v_{\g}(\mathcal B)\ =\ v_{\g}(\ker_{\Univ}^{\ne} A),$$
using for the last equality the details in the proof of Corollary~\ref{cor:spectral valuation basis}.
\end{proof}

\begin{prop}\label{prop:finiteness of excu(A)}
Suppose $K$ is $\upo$-free. 
Then  $\nwt_{A_\lambda}(\gamma)=0$ for all but finitely many  pairs $(\gamma,\lambda)$  and  
$$\abs{\exc^{\operatorname{u}}(A)}\ \leq\ 
\sum_\lambda\, \abs{\exc^{\ev}(A_\lambda)}\ =\ 
\sum_{\gamma,\lambda} \nwt_{A_\lambda}(\gamma)\ \leq\ r.$$
If $\dim_C \ker_{\Univ}A =r$, then  $\sum_\lambda\, \abs{\exc^{\ev}(A_\lambda)}=r$ and $v_{\g}(\ker_{\Univ}^{\neq} A) = \exc^{\operatorname{u}}(A)$.
% where $v_{\g}$ is the gaussion extension of the valuation of $K$.
\end{prop}
\begin{proof}
Equip $\Omega$ with a spectral extension $v$ of the valuation of $K$. 
Then $\Omega$ is $\upo$-free, so $\sum_\lambda|\exc^{\ev}(A_\lambda)|=|\exc^{\ev}_{\Omega}(A)|\le r$ by the remarks preceding Lemma~\ref{lem:finiteness of excu(A)}  and Corollary~\ref{cor:sum of nwts} applied to $\Omega$ in place of $K$. These remarks also give 
$\nwt_{A_\lambda}(\gamma)=0$ for all but finitely many pairs $(\gamma,\lambda)$, and so
$$\sum_{\gamma,\lambda} \nwt_{A_\lambda}(\gamma)\ =\ \sum_{\gamma,\lambda}  \nwt_A\!\big(\gamma+v(\ex(\lambda)\big)\  =\  \abs{\exc^{\ev}_\Omega(A)}\ \leq\ r.$$
Corollary~\ref{cor:sum of nwts} applied to $A_\lambda$ in place of $A$ yields $\abs{\exc^{\ev}(A_\lambda)}=\sum_\gamma \nwt_{A_\lambda}(\gamma)$ and so~$\sum_\lambda\, \abs{\exc^{\ev}(A_\lambda)} = \sum_{\gamma,\lambda} \nwt_{A_\lambda}(\gamma)$. This proves the first part (including  the display). The rest follows from this and Lemma~\ref{lem:finiteness of excu(A)}.
 %Suppose $\dim_C\ker_{\Univ}A=r$.
%Then $\dim_C \ker_{\Omega} A = r$, 
%so $r = \abs{v(\ker_\Omega^{\neq} A)}$ by [ADH, 5.6.6(i)] applied to $\Omega$ in place of~$K$. Also $v(\ker_\Omega^{\neq} A)\subseteq \exc^{\ev}_\Omega(A)$ by [ADH, p. 481], so $\abs{v(\ker_\Omega^{\neq} A)}\leq \abs{\exc^{\ev}_\Omega(A)}\le r$, so $\sum_\lambda\, \abs{\exc^{\ev}(A_\lambda)} = r$, as claimed.
%It remains to show $v_{\g}(\ker_{\Univ}^{\neq} A)= \exc^{\operatorname{u}}(A)$.
%We have $\abs{v(\ker_\Omega^{\neq} A)}=\abs{\exc^{\ev}_\Omega(A)}=~r$. Now $v(\ker_\Omega^{\neq} A)\subseteq \exc^{\ev}_\Omega(A)$ as already noted, and $\ker_{\Omega}^{\ne} A= \ker_{\Univ}^{\ne} A$, so $v(\ker_{\Univ}^{\neq} A)= \exc^{\ev}_\Omega(A)$. 
%Let $\pi\colon \Gamma_{\Omega} \to \Gamma$ be the projection map
%given by $\pi\big(\gamma+v\big(\!\ex(\lambda)\big)\big)=\gamma$. Then
%$\pi\big(\exc^{\ev}_\Omega(A)\big)=\exc^{\operatorname{u}}(A)$ by \eqref{eq:excevOmega}. It is therefore enough to show
%$\pi\big(v(\ker_{\Univ}^{\neq} A)\big)=v_{\g}(\ker_{\Univ}^{\ne} A)$.
%For that, use Lemma~\ref{lem:v(ex(Q))} and Corollary~\ref{cor:spectral valuation basis}.  
\end{proof} 

\noindent
In the next lemma (to be used in the proof of Proposition~\ref{prop:stability of excu, real}), as well as in Corollary~\ref{cor:excev cap GammaOmega},  $L$ is a $\d$-valued $H$-asymptotic extension of $K$  with 
algebraically closed constant field and asymptotic integration (so $L$ has small derivation), such that $L^\dagger$ is divisible, $L^\dagger\cap K=K^\dagger$, and $\I(L)\subseteq L^\dagger$ .  We also fix there a complement~$\Lambda_L$ of the $\Q$-linear subspace  $L^\dagger$ of $L$ with
$\Lambda\subseteq\Lambda_L$. Let $\Univ_L=L\big[\!\ex(\Lambda_L)\big]$ be
the corresponding universal exponential extension  
of $L$ containing $\Univ=K\big[\!\ex(\Lambda)\big]$ as a differential subring, as  
described in the remarks following Corollary~\ref{cor:Univ under d-field ext},
with differential fraction field $\Omega_L$.
%Assume now that  $\Omega_L$ is equipped with a spectral extension of the valuation of $L$, and restrict this
%spectral extension to $\Omega$; this gives a spectral extension of the valuation of $K$. 

\begin{lemma}\label{lem:excev cap GammaOmega} Assume $C_L=C$. Let $\Omega_L$ be equipped with a spectral extension of the valuation of $L$, and take $\Omega$ as a valued subfield of $\Omega_L$; so the valuation of $\Omega$ is a spectral extension of the valuation of $K$. Suppose $\Psi$ is cofinal in $\Psi_L$ or $K$ is $\upl$-free.
Then $\exc^{\ev}_{\Omega_L}(A)\cap \Gamma_\Omega\  =\  \exc^{\ev}_\Omega(A)$.
\end{lemma} 
\begin{proof}
Let $\mu$ range over $\Lambda_L$.
We have 
$$\Gamma_{\Omega_L}\ =\ \Gamma_L \oplus v\big(\!\ex(\Lambda_L)\big), \qquad \Gamma_\Omega\ =\ \Gamma\oplus v\big(\!\ex(\Lambda)\big)$$ by Lemma~\ref{lem:v(ex(Q))} and
$$\exc^{\ev}_{\Omega_L}\ =\ \bigcup_{\mu} \exc^{\ev}_L(A_\mu)+v\big(\!\ex(\mu)\big), \qquad 
\exc^{\ev}_\Omega\ =\  \bigcup_\lambda \exc^{\ev}(A_\lambda)+v\big(\!\ex(\lambda)\big)$$
by \eqref{eq:excevOmega}.  Hence
$$\exc^{\ev}_{\Omega_L}(A)\cap \Gamma_\Omega\ =\  \bigcup_{\lambda} \big( \exc^{\ev}_L(A_\lambda)\cap\Gamma \big)+v\big(\!\ex(\lambda)\big)\ =\  \bigcup_{\lambda}  \exc^{\ev}(A_\lambda) +v\big(\!\ex(\lambda)\big)\  =\  \exc^{\ev}_\Omega(A),$$
where we used the injectivity of $\mu\mapsto v\big(\!\ex(\mu)\big)$ 
for the first equality and  Remark~\ref{rem:excu for different Q, 2}  for the second.
\end{proof}

\noindent
Call $A$ {\bf terminal} with respect to $\Lambda$
if $\sum_\lambda \abs{\exc^{\ev}(A_\lambda)}=r$.\index{terminal}\index{linear differential operator!terminal}  We omit  ``with respect to~$\Lambda$'' if it is clear from the context what $\Lambda$ is.
In Section~\ref{sec:ultimate} we shall restrict~$K$,~$\Lambda$ anyway so that this dependence on $\Lambda$ disappears.
Recall also that for a given spectral extension of the valuation of $K$ to $\Omega$ we have $\abs{\exc_{\Omega}(A)}= \sum_\lambda \abs{\exc^{\ev}(A_\lambda)}$. 
If $A$ is terminal and $\phi\in K^\times$ is active in $K$, then $A^\phi\in K^\phi[\derdelta]$ is terminal
with respect to $\Lambda^\phi$ (cf.~remarks after Corollary~\ref{cor:excu for different Q}). If $A$ is terminal and
$a\in K^\times$, then $aA$ is terminal. 
If $r=0$, then $A$ is terminal. 

\begin{lemma}\label{lem:ultimate1}
If $r=1$, then $A$ is terminal.
\end{lemma}
\begin{proof} Assume $r=1$. Then $\dim_C \ker_{\Univ}A=1$, so $\sum_\lambda \abs{\exc^{\ev}(A_\lambda)}\ge 1$
by Lemma~\ref{lem:finiteness of excu(A)}. Equip $\Omega$ with a spectral extension of the valuation of $K$. Then
$\Omega$ is ungrounded by  Lemma~\ref{lem:v(ex(Q))}, and $r=1$ gives $\abs{\exc_{\Omega}(A)}\le 1$. Now use
$\abs{\exc_{\Omega}(A)}= \sum_\lambda \abs{\exc^{\ev}(A_\lambda)}$.
\end{proof}

%By Proposition~\ref{prop:finiteness of excu(A)}, if~${\dim_{C} \ker_{\Univ} A=r}$, then~$A$ is ultimate.
%So if~$A$ splits over $K$ and $K$ is $r$-linearly surjective when~${r\geq 2}$, then~$A$ is ultimate, by Corollary~\ref{cor:basis of kerUA}. 

\begin{lemma}\label{lem:ultimate prod}
Suppose $A$ and $B\in K[\der]^{\neq}$ are terminal, and each operator $B_\lambda$ is asymptotically surjective.   Then $AB$ is terminal.
\end{lemma}
\begin{proof}
Use that $(AB)_\lambda=A_\lambda B_\lambda$, and that
$\abs{\exc^{\ev}(A_\lambda B_\lambda)}=\abs{\exc^{\ev}(A_\lambda)}+\abs{\exc^{\ev}(B_\lambda)}$ by Corollary~\ref{cor:exce product}.
\end{proof}

\noindent
Thus if $A$ is terminal and $a\in K^\times$, then $aA, Aa$, and $A_{\ltimes a}$ are terminal. From 
Lemmas~\ref{lem:ultimate1},~\ref{lem:ultimate prod}, and Corollary~\ref{cor:well-behaved} we conclude:

\begin{cor}\label{cor:ultimate prod, 2}
If $K$ is $\upl$-free and $A$ splits over $K$, then $A$ is terminal.
\end{cor}

\begin{cor}\label{cor:ultimate prod, 1}
Suppose $K$ is $\upo$-free and   $B\in K[\der]^{\neq}$. Then $A$ and $B$ are terminal iff~$AB$ is terminal.
\end{cor}
\begin{proof}
The ``only if'' part   follows from Lemma~\ref{lem:ultimate prod} and Corollary~\ref{cor1524}. 
For the ``if'' part, use Corollary~\ref{cor:exce product} and Proposition~\ref{prop:finiteness of excu(A)}.
\end{proof}

%\noindent
%In~\cite{ADHld} we shall  study other situations where $A$ is ultimate.
%By Proposition~\ref{prop:finiteness of excu(A)} and the identities in the proof of Lemma~\ref{lem:excev cap GammaOmega}:

\begin{cor}\label{cor:excev cap GammaOmega}
Suppose $A$ is terminal, $\Psi$ is cofinal in $\Psi_L$ or $K$ is $\upl$-free, and~$L$ is $\upo$-free. Then, with respect to the complement $\Lambda_L$ of $L^\dagger$ in $L$, we have: 
\begin{enumerate}
\item[\textup{(i)}]  as an element of $L[\der]$, $A$ is terminal;
\item[\textup{(ii)}] $\exc^{\ev}(A_\mu)=\emptyset$ for all $\mu\in\Lambda_L\setminus\Lambda$;
\item[\textup{(iii)}] $\exc^{\ev}(A_\lambda)=\exc^{\ev}_L(A_\lambda)$ for all $\lambda$; and
\item[\textup{(iv)}] $\exc^{\operatorname{u}}(A)=\exc^{\operatorname{u}}_L(A)$.
\end{enumerate}
\end{cor}

{\samepage \begin{proof}
By the remarks after Corollary~\ref{cor:excu for different Q} we have $\exc^{\ev}(A_\lambda)\subseteq\exc^{\ev}_L(A_\lambda)$ for each $\lambda$,
and so with $\mu$ ranging over $\Lambda_L$, by Proposition~\ref{prop:finiteness of excu(A)} applied to $L$ in place of $K$,
we have $r=\sum_\lambda\abs{\exc^{\ev}(A_\lambda)}\leq \sum_\mu \abs{\exc^{\ev}_L(A_\mu)}\leq r$. This yields (i)--(iv).
\end{proof}

\noindent
In~\cite{ADHld} we shall  study other situations where $A$ is terminal.
}

\subsection*{The real case}
{\it In this subsection $H$ is a real closed $H$-field with small derivation, asymptotic integration, and $H^\dagger=H$;
also $K=H[\imag]$, $\imag^2=-1$, for our valued
differential field $K$. We also assume $\I(H)\imag\subseteq K^\dagger$.}\/
Then $K$ is $\d$-valued of $H$-type with small derivation, asymptotic integration, $K^\dagger=H+\I(H)\imag$, 
and $\I(K)\subseteq K^\dagger$. Note that $H$ and thus $K$ is $\upl$-free by [ADH, remark after 11.6.2, and 11.6.8].
Let~$A$ in $K[\der]^{\neq}$ have order $r$ and let  $\gamma$ range over $\Gamma$.

\begin{lemma}\label{lem:LambdaL}
If the real closed $H$-asymptotic extension $F$ of~$H$ has asymptotic integration
and convex valuation ring, then~$L^\dagger\cap K=K^\dagger$ for the  algebraically closed
$H$-asymptotic field extension $L:=F[\imag]$ of $K$. 
\end{lemma}
\begin{proof} Use Corollary~\ref{cor:logderset ext} and earlier remarks in the same subsection. 
\end{proof}

%Together with results from [ADH]  it yields:

\begin{cor}\label{cor:LambdaL}
The $H$-field $H$ has an $H$-closed  extension $F$ with $C_F=C_H$, and for any such~$F$, the algebraically closed
$\d$-valued field extension $L:=F[\imag]$ of $H$-type of $K$ is $\upo$-free with~$C_L=C$, $\I(L)\subseteq L^\dagger$, and
 $L^\dagger\cap K=K^\dagger$.
\end{cor}
 \begin{proof} Use [ADH, 16.4.1, 9.1.2] to extend $H$ to an $\upo$-free $H$-field with the same constant field as $H$, next use~[ADH, 11.7.23] to pass to its real closure, and then use~[ADH, 14.5.9]  to extend further to an $H$-closed $F$, still with the same constant field as $H$. 
 %The existence of $F$ follows from [ADH, 11.7.18, 14.5.11]. 
 For any such~$F$, the $\d$-valued field
  $L:=F[\imag]$ of $H$-type is $\upo$-free by~[ADH,  11.7.23] and newtonian by~[ADH, 14.5.7]. Hence~$\I(L)\subseteq L^\dagger$
  by Lemma~\ref{lem:ADH 14.2.5},
  and~$L^\dagger\cap K=K^\dagger$
  by Lemma~\ref{lem:LambdaL}.
 \end{proof}

\noindent
This leads to a variant of Proposition~\ref{prop:finiteness of excu(A)}:

\begin{prop}\label{prop:finiteness of excu(A), real}
The conclusion of Proposition~\ref{prop:finiteness of excu(A)} holds. In particular:
$$  \dim_C\ker_{\Univ}A=r\ \Longrightarrow\ A \text{ is terminal}.$$
\end{prop}
\begin{proof} 
%Assume $H$ is $\upl$-free. Then $K$ is $\upl$-free by [ADH, 11.6.8].
Corollary~\ref{cor:LambdaL} gives  an $H$-closed extension $F$ of~$H$ with $C_F=C_H$, so~$L:=F[\imag]$ is $\upo$-free, 
$C_L=C$, $\I(L)\subseteq L^\dagger$, and $L^\dagger\cap K=K^\dagger$.   Take a complement~$\Lambda_L\supseteq \Lambda$ of the  subspace~$L^\dagger$ of the $\Q$-linear space $L$.
By Remark~\ref{rem:excu for different Q, 2} we have 
$\exc^{\ev}(A_\lambda)=\exc^{\ev}_L(A_\lambda)\cap\Gamma$. Hence
Proposition~\ref{prop:finiteness of excu(A)} applied to~$K$,~$\Lambda$ replaced by~$L$,~$\Lambda_L$, respectively, and $A$ viewed
as element of $L[\der]$, yields
 $\sum_\lambda\, \abs{\exc^{\ev}(A_\lambda)} \leq r$. 
Corollary~\ref{cor:13.7.10} applied to~$A_\lambda$ in place of~$A$ gives~$\abs{\exc^{\ev}(A_\lambda)} = \sum_\gamma \nwt_{A_\lambda}(\gamma)$. 
This yields  the conclusion of Proposition~\ref{prop:finiteness of excu(A)} as in the proof of that proposition. 
\end{proof}

\noindent
Let now $F$ be a Liouville closed $H$-field extension of $H$ and  suppose
$\I(L)\subseteq L^\dagger$ where~$L:=F[\imag]$.
Lemma~\ref{lem:LambdaL} yields $L^\dagger\cap K=K^\dagger$, so $L$ is as described just before Lemma~\ref{lem:excev cap GammaOmega}, and we have
a complement~$\Lambda_L\supseteq\Lambda$ of the subspace $L^\dagger$ of the $\Q$-linear space $L$. 
Note that if $A$ splits over $K$, then $A$ is terminal by Corollary~\ref{cor:ultimate prod, 2}.  
%Arguing as in the proof of Corollary~\ref{cor:excev cap GammaOmega}, using Proposition~\ref{prop:finiteness of excu(A), real} instead of Proposition~\ref{prop:finiteness of excu(A)}, we obtain:

\begin{cor}\label{cor:excev cap GammaOmega, real} 
Suppose $A$ is terminal. Then, with respect to the complement~$\Lambda_L$ of $L^\dagger$ in $L$, the conclusions
\text{\rm{(i)--(iv)}} of Corollary~\ref{cor:excev cap GammaOmega} hold.
\end{cor}
\begin{proof} By the remarks after Corollary~\ref{cor:excu for different Q} we have $\exc^{\ev}(A_\lambda)\subseteq\exc^{\ev}_L(A_\lambda)$ for all $\lambda$, 
and so with $\mu$ ranging over $\Lambda_L$,  Proposition~\ref{prop:finiteness of excu(A), real} applied to $L$ in place of $K$,
we have $r=\sum_\lambda\abs{\exc^{\ev}(A_\lambda)}\leq \sum_\mu \abs{\exc^{\ev}_L(A_\mu)}\leq r$. This yields (i)--(iv).
\end{proof}

\newpage 

\part{Normalizing Holes and Slots}\label{part:normalization}

\medskip

\noindent
In this introduction $K$ is an $H$-asymptotic field with small derivation and rational asymptotic integration.
In Section~\ref{sec:holes} we introduce  {\it holes}\/ in $K$: A {\it hole in $K$\/} is a triple $(P,\fm,\hat a)$  with $P\in K\{Y\}\setminus K$, $\fm\in K^\times$, and $\hat a\in \hat{K}\setminus K$ for some immediate asymptotic extension
$\hat{K}$ of $K$, such that $\hat a \prec \fm$ and $P(\hat a)=0$.
The main goal of Part~\ref{part:normalization}  is a Normalization Theorem, namely Theorem~\ref{mainthm},  that allows us to transform
under reasonable conditions a hole $(P,\fm,\hat a)$ in $K$ into a ``normal'' hole; this helps to pin down the location of
$\hat a$ relative to $K$. 
%This is accomplished by Theorem}. 
The notion of~$(P,\fm,\hat a)$ being {\em normal\/} involves the linear part of
the differential polynomial~$P_{\times\fm}$, in particular the {\it span} of this linear part. We introduce the span in the preliminary Section~\ref{sec:span}.
In Section~\ref{sec:isolated} we  study {\it isolated}\/ holes~$(P,\fm,\hat a)$ in~$K$, which under reasonable conditions ensure the uniqueness of the 
isomorphism type of $K\<\hat a\>$ as a
valued differential field over $K$; see Proposition~\ref{prop:2.12 isolated}. In Section~\ref{sec:holes of c=(1,1,1)} we   focus on holes~$(P,\fm,\hat a)$ in~$K$ where~$\order P=\deg P=1$. 
For technical reasons we actually work in Part~\ref{part:normalization} also with {\em slots\/} in $K$, which are a bit more general than holes in $K$. 

\medskip\noindent
First some notational conventions. Let $\Gamma$ be an ordered abelian group. For $\gamma, \delta\in \Gamma$ with 
$\gamma\ne 0$ the expression ``$\delta=o(\gamma)$'' means ``$n|\delta|< |\gamma|$ for all $n\ge 1$'' according to~[ADH, 2.4], but here we find it convenient to extend this to $\gamma=0$, in which case~``$\delta=o(\gamma)$" means ``$\delta=0$".
% set $o(\gamma):=0\in \Gamma$ for $\gamma=0\in \Gamma$.  
Suppose $\Gamma=v(E^\times)$ is the value group
of a valued field~$E$ and $\fm\in E^\times$. Then we denote the archimedean class
$[v\fm ]\subseteq \Gamma$ of $v\fm\in \Gamma$ by just $[\fm]$. Suppose $\fm \nasymp 1$. Then we have a
proper convex subgroup\label{p:Delta}
$$\Delta(\fm)\ :=\ \big\{\gamma\in \Gamma:\, \gamma=o(v\fm)\big\}\ =\ \big\{\gamma\in \Gamma:\, [\gamma]<[\fm]\big\},$$ 
of~$\Gamma$. 
If  $\fm\asymp_{\Delta(\fm)}\fn\in E$, then~$0\ne \fn \nasymp 1$ and $\Delta(\fm)=\Delta(\fn)$. 
In particular,  if~$\fm\asymp\fn\in E$, then~$0\ne\fn \nasymp 1$ and $\Delta(\fm)=\Delta(\fn)$. Note that for $f,g\in E$ the meaning of ``$f\preceq_{\Delta(\fm)} g$'' does not change in passing to a valued field extension of~$E$, although~$\Delta(\fm)$ can increase as a subgroup of the value group of the extension.  

\section{The Span of a Linear Differential Operator} \label{sec:span}

\noindent
{\em In this section $K$ is a valued differential field with small derivation and   $\Gamma:= v(K^\times)$.  We let~$a$,~$b$, sometimes subscripted, range over $K$, and $\fm$,~$\fn$ over $K^\times$}. 
Consider a linear differential operator
$$A\ =\ a_0+a_1\der+\cdots+a_r\der^r\in K[\der],\qquad   a_r\neq 0.$$
We shall use below the quantities $\dwm(A)$ and $\dwt(A)$ defined in [ADH,  5.6].  We also introduce a measure $\fv(A)$ for the ``lopsidedness'' of $A$ as follows:\label{p:span}
$$\fv(A)\ :=\ a_r/a_m \in\ K^\times\qquad\text{where $m:=\dwt(A)$.}$$  
So $a_r\asymp \fv(A)A$ and $\fv(A)\preceq 1$, with
$$\fv(A)\asymp 1\quad\Longleftrightarrow\quad\dwt(A)=r\quad\Longleftrightarrow\quad\fv(A)=1.$$
Also note that $\fv(aA)=\fv(A)$ for $a\neq 0$. Moreover,  
$$\fv(A_{\ltimes\fn})A_{\ltimes\fn}\ \asymp\  a_r\ \asymp\ \fv(A)A$$ 
since $A_{\ltimes\fn}=a_r\der^r+\text{lower order terms in $\der$}$. 

\begin{example}
$\fv(a+\der)=1$ if  $a\preceq 1$, and $\fv(a+\der)=1/a$ if $a\succ 1$.
\end{example}

\noindent
We call $\fv(A)$ the {\bf span} of $A$.\index{span}\index{linear differential operator!span}
We are mainly interested in the valuation of $\fv(A)$. 
This is related to the gaussian valuation $v(A)$ of $A$: if~$A$ is monic, then $v\big(\fv(A)\big)=-v(A)$. 
An important property of the span of $A$ is that its valuation is not affected by small additive perturbations of $A$: 

\begin{lemma}\label{lem:fv of perturbed op}
Suppose $B\in K[\der]$, $\order(B)\leq r$ and $B\prec \fv(A) A$. Then: \begin{enumerate}
\item[\textup{(i)}] $A+B\sim A$, $\dwm(A+B)=\dwm(A)$, and $\dwt(A+B)=\dwt(A)$;
\item[\textup{(ii)}] $\order(A+B)=r$ and $\fv(A+B) \sim \fv(A)$. 
\end{enumerate}
\end{lemma}
\begin{proof}
From $B\prec \fv(A) A$ and $\fv(A)\preceq 1$ we obtain $B\prec A$, and thus (i).
Set $m:=\dwt(A)$, let~$i$ range over $\{0,\dots,r\}$, and let $B=b_0+b_1\der+\cdots+b_r\der^r$. 
Then $a_i\preceq a_m$ and $b_i\prec \fv(A)A\asymp a_r \preceq a_m$.
Therefore, if $a_i\asymp a_m$, then $a_i+b_i\sim a_i$, and if $a_i\prec a_m$, then $a_i+b_i\prec a_m$. Hence 
$\fv(A+B)=(a_r+b_r)/(a_m+b_m) \sim a_r/a_m=\fv(A)$.
\end{proof}

\noindent
For $b\neq 0$, the valuation of $\fv(Ab)$ only depends on $vb$; it is enough to check this for~$b\asymp 1$. More generally:

\begin{lemma}\label{fvmult}
Let $B\in K[\der]^{\neq}$ and $b\asymp B$. Then $\fv(AB)\asymp\fv(Ab)\fv(B)$.
\end{lemma}
\begin{proof}
Let $B=b_0+b_1\der+\cdots+b_s\der^s$, $b_s\neq 0$. Then
$$AB\ =\ a_rb_s\der^{r+s}+\text{lower order terms in}\ \der,$$ so by [ADH, 5.6.1(ii)] for $\gamma=0$:
\alignqed{v\big(\fv(AB)\big)\ &=\ v(a_rb_s)-v(AB)\ =\ v(a_rb_s) - v(Ab)\\  &=\ v(a_rb)-v(Ab)+v(b_s)-v(B) \\ \ &=\ v\big(\fv(Ab)\fv(B)\big).} 
\end{proof}

\begin{cor}\label{cor:111} 
Let $B\in K[\der]^{\neq}$. If $\fv(AB)=1$, then $\fv(A)=\fv(B)=1$. The converse holds if $B$ is monic. 
%$$\fv(AB) = 1\quad \Longleftrightarrow\quad \fv(A)=\fv(B)=1.$$
\end{cor}

\noindent
This is clear from from Lemma~\ref{fvmult}, and in turn gives:

\begin{cor}\label{corAfv1}
Suppose $A=a(\der-b_1)\cdots(\der-b_r)$. Then
$$\fv(A) =  1 \quad\Longleftrightarrow\quad b_1,\dots, b_r\preceq 1. $$
\end{cor}

\begin{remark}
Suppose $K=C(\!(t)\!)$ with the $t$-adic valuation and derivation $\der=t\frac{d}{dt}$. In the literature, $A$ is called {\it regular singular}\/ if $\fv(A)=1$, and {\it irregular singular}\/ if~$\fv(A)\prec 1$; see \cite[Definition~3.14]{vdPS}.
\end{remark}

\begin{lemma} 
Let $B\in K[\der]^{\neq}$.  Then $\fv(AB)\preceq\fv(B)$, and if  $B$ is monic, then $\fv(AB)\preceq\fv(A)$.
\end{lemma}
\begin{proof}
Lemma~\ref{fvmult} and $\fv(Ab)\preceq 1$ for $b\neq 0$  yields $\fv(AB)\preceq\fv(B)$.
Suppose~$B$ is monic, so $v(B)\leq 0$. To show $\fv(AB)\preceq\fv(A)$ we  arrange that $A$
is also monic. Then~$AB$ is monic, and $\fv(AB)\preceq\fv(A)$ is equivalent to $v(AB)\leq v(A)$. Now
$$v(AB)\ =\ v_{AB}(0)\  =\  v_A\big(v_B(0)\big)\ =\  v_A\big(v(B)\big)\  \leq\  v_A(0)\  =\  v(A)$$
by [ADH, 4.5.1(iii), 5.6.1(ii)].
\end{proof}

%\begin{lemma}\label{lem:bound on linear factors} 
%Suppose $A=(\der-a)B$ with monic $B\in K[\der]^{\neq}$. Then $a\preceq \fv(A)^{-1}$.
%\end{lemma}
%\begin{proof} Take~$b$ with $b\asymp B$; then $\fv(B)\asymp 1/b$ since $B$ is monic. We have $\fv(A) \asymp \fv\big( (\der-a)b\big)\fv(B)$ by Lemma~\ref{fvmult}, hence $\fv(A)\preceq\fv(B)$. If $a\preceq 1$, then the lemma holds by $\fv(A)\preceq 1$, so assume $a\succ 1$ below. If $a\preceq b^\dagger$, then by [ADH, 9.2.10(iv)]
%$$a\ \preceq\ b^\dagger\ \prec b\ \asymp\ \fv(B)^{-1}\ \preceq\ \fv(A)^{-1}$$
%and we are done. If $a\succ b^\dagger$, then $(\der-a)b = b(\der+b^\dagger-a)$ gives $\fv\big( (\der-a)b\big)\asymp 1/a$, hence $\fv(A) \asymp 1/(ab)$, so again $a \preceq   \fv(A)^{-1}$, as claimed.
%\end{proof}

\begin{cor}\label{cor:bound on linear factors}
If $A=a(\der-b_1)\cdots(\der-b_r)$, then  
$b_1,\dots, b_r\ \preceq\ \fv(A)^{-1}.$
\end{cor}  
%\begin{proof}
%We may assume $a=1$, $r\geq 1$.    Lemma~\ref{lem:bound on linear factors} yields $b_1\preceq\fv(A)^{-1}$, and inductively we assume $b_2,\dots,b_r\preceq\fv(B)^{-1}$ where $B:=(\der-b_2)\cdots(\der-b_r)$. Now use that $\fv(B)^{-1}\preceq\fv(A)^{-1}$ by the beginning of the proof of Lemma~\ref{lem:bound on linear factors}. 
%\end{proof}

\noindent
Let $\Delta$ be a convex subgroup of $\Gamma$, let $\dot{\mathcal O}$ be the valuation ring of the coarsening~$v_\Delta$ of the valuation~$v$ of $K$ by $\Delta$, with maximal ideal $\dot{\smallo}$,
and $\dot K=\dot{\mathcal O}/\dot{\smallo}$ be the
valued differential residue field of   $v_\Delta$.
The residue morphism $\dot{\mathcal O}\to\dot K$ extends to the ring morphism $\dot{\mathcal O}[\der]\to\dot K[\der]$ with $\der\mapsto\der$. If $A\in\dot{\mathcal O}[\der]$ and $\dot A\neq 0$,
then $\dwm(\dot A)=\dwm(A)$ and $\dwt(\dot A)=\dwt(A)$.
We set $\fv:=\fv(A)$.

\begin{lemma}\label{lem:dotfv}
If $A\in\dot{\mathcal O}[\der]$ and $\order(\dot A)=r$, then $\fv(\dot A)=\dot \fv$.
\end{lemma}

\subsection*{Behavior of the span under twisting} 
Recall that $o(\gamma):= 0\in \Gamma$ for $\gamma=0\in \Gamma$. 
With this convention, here is a consequence of [ADH, 6.1.3]: 

\begin{lemma}\label{lem:6.1.3 consequ} 
Let $B\in K[\der]^{\neq}$. Then $v(AB)=v(A)+v(B)+o\big(v(B)\big)$.
\end{lemma}
\begin{proof}
Take $b$ with $b\asymp B$. Then 
$$v(AB)=v_{AB}(0)=v_A\big(v_B(0)\big)=v_A(vb)=v(Ab)$$
by~[ADH, 5.6.1(ii)]. Moreover, $v(Ab)=v(A)+vb+o(vb)$, by [ADH, 6.1.3].
\end{proof}

%\begin{cor}\label{cor:6.1.3 consequ}
%If $A$ and $B\in K[\der]^{\neq}$ are monic, then $v(AB) \leq v(A)$ and $v(B)=O\big(v(AB)\big)$.
%\end{cor}

%\begin{lemma}
%Suppose $A=(\der-b_1)\cdots (\der-b_r)$. Then
%$$v(A) = \sum_{b_j\succ 1} v(b_j)+o\big(v(A)\big).$$
%\end{lemma}
%\begin{proof}
%By induction on $r$. The lemma clearly holds for $r=0$. Assume the claim holds for the above $A$ of order $r$ and consider $AB$ with $B=\der-b_{r+1}$. If $b_{r+1}\preceq 1$, then $v(B)=0$ and so $v(AB)=v(A)$ by Lemma~\ref{lem:6.1.3 consequ}, so the desired result trivially also holds for $AB$ in place of $A$. Suppose $b_{r+1}\succ 1$, thus~$B\asymp b_{r+1}$; hence by Lemma~\ref{lem:6.1.3 consequ} and Corollary~\ref{cor:6.1.3 consequ} we get
%\begin{align*}
%v(AB) 	&= v(A)+v(b_{r+1})+o\big(v(B)\big) \\
%		&= \sum_{b_j\succ 1} v(b_j)+o\big(v(A)\big)+o\big(v(B)\big) 
%		 = \sum_{b_j\succ 1} v(b_j)+o\big(v(AB)\big)
%\end{align*}
%as required.
%\end{proof}

\noindent
We have $\fv(A_{\ltimes\fn})=\fv(A\fn)$, so $v(A_{\ltimes\fn})=v(A)+o(v\fn)$
by Lemma~\ref{lem:6.1.3 consequ}. Moreover:

\begin{lemma}\label{lem:An} 
%Suppose $\fn\nasymp 1$. Then 
$v\big(\fv(A\fn)\big)= v\big(\fv(A)\big) + o(v\fn)$.
\end{lemma}
\begin{proof}
Replacing $A$ by $a_r^{-1} A$ we arrange $A$ is monic, so $A_{\ltimes\fn}$ is monic, and thus
$$v\big(\fv(A\fn)\big)\ =\ v\big(\fv(A_{\ltimes\fn})\big)=-v(A_{\ltimes\fn})=-v(A)+o(v\fn)=v\big(\fv(A)\big)+o(v\fn)$$
by remarks preceding the lemma.
%Replacing $A$ by $a_m^{-1} A$ where $m=\dwt(A)$ we arrange $a_m=1$, so $a_r = \fv:=\fv(A)$. Taking $b_0,\dots,b_r$ with $A\fn=b_0+b_1\der+\cdots+b_r\der^r$ we have $b_r= \fn a_r = \fn\fv$. Let $n=\dwt(A\fn)$; then  $$v(b_n)= v( A\fn)=v(A)+v(\fn)+ o\big(v(\fn)\big)=v(\fn)+o\big(v(\fn)\big)$$ by [ADH, 6.1.3]. Using also $\fv(A\fn) = b_r/b_n$, this gives 
%$$ v\big(\fv(A\fn)\big)\ =\ v(b_r)-v(b_n)\ =\ v(\fn)+v(\fv)-\big(v(\fn)+o(v(\fn))\big)\ =\ v(\fv) + o\big(v(\fn)\big),$$
%as promised.   
\end{proof}

\noindent
Recall: we denote the archimedean class
$[v\fn]\subseteq \Gamma$ by $[\fn]$.   Lemma~\ref{lem:An} yields:

\begin{cor}\label{cor:An}
$\big[\fv(A)\big]<[\fn]\ \Longleftrightarrow\ \big[\fv(A\fn)\big]<[\fn]$.
\end{cor}

\noindent
Under suitable conditions on $K$ we can say more about the valuation of 
$\fv(A_{\ltimes\fn})$: Lemma~\ref{lem:Atwist} below. 
%For this we first prove two variants of~[ADH, 11.1.5]: 

%\begin{lemma}\label{Riccatipower}
%If $z\in K^{\succ 1}$, then $R_n(z)=z^n(1+ \epsilon)$ with 
%$v\epsilon\ge v(z^{-1})+ o(vz)>0$.  
%\end{lemma}
%\begin{proof} This is clear for $n=0$ and $n=1$. Suppose $z\succ 1$, $n\ge 1$, and
%$R_n(z)=z^n(1+ \epsilon)$ with $\epsilon$ as in the lemma. As in the proof of [ADH, 11.1.5],
%$$ R_{n+1}(z)\ =\  z^{n+1}\left(1 + \epsilon +    n\frac{z^\dagger}{z}(1+\epsilon) +       
%\frac{\epsilon'}{z}\right).$$
%Now $v(z^\dagger)\ge o(vz)$: this is obvious if $z^\dagger\preceq 1$, and follows from~$\triangledown(\gamma)=o(\gamma)$ for~$\gamma\ne 0$ if $z^\dagger\succ 1$ [ADH, 6.4.1(iii)]. 
%This gives the desired result in view of  $\epsilon'\prec 1$. 
%\end{proof}

%\begin{lemma}\label{Riccatipower+} Suppose $\der\mathcal O\subseteq\smallo$. If 
%$z\in K^{\succeq 1}$, then $R_n(z)=z^n(1+ \epsilon)$ with 
%$\epsilon \prec 1$.
%\end{lemma}
 
%\begin{proof}
%The case $z\succ 1$ follows from Lemma~\ref{Riccatipower}. For $z\asymp 1$, proceed as in the proof of that lemma, using $\der\mathcal O\subseteq\smallo$.
%\end{proof}

%\noindent
%By [ADH, 9.1.3 (iv)] the condition $\der\mathcal O\subseteq\smallo$ is satisfied
%if $K$ is $\d$-valued, or  asymptotic   with $\Psi\cap \Gamma^{>}\ne \emptyset$.

\begin{lemma}\label{lem:nepsilon, first part}
Let $\fn^\dagger\succeq 1$ and $\fm_0,\dots, \fm_r\in K^\times$ be such that
$$v(\fm_i) + v(A)\ =\ \min_{i\leq j\leq r} v(a_j)+(j-i)v(\fn^\dagger).$$
Then with $m:=\dwt(A)$ we have
$$\fm_0\ \succeq\ \cdots\ \succeq\ \fm_r\quad\text{and}\quad
(\fn^\dagger)^m\ \preceq\ \fm_0\ \preceq\ (\fn^\dagger)^r.$$
\textup{(}In particular, $[\fm_0] \leq [\fn^\dagger]$, with equality if~$m>0$.\textup{)}
\end{lemma}
\begin{proof}
From $v(\fn^\dagger)\leq 0$ we obtain $v(\fm_0)\leq \cdots\leq v(\fm_r)$.
We have $0\leq v(a_j/a_m)$ for~$j=0,\dots,r$ and so
$$rv(\fn^\dagger)\ \leq\ \min_{0\leq j\leq r}  v(a_j/a_m)+jv(\fn^\dagger)\ =\ v(\fm_0)\ \leq\ mv(\fn^\dagger)$$
as required.
\end{proof}

\begin{lemma}\label{lem:Atwist} 
Suppose $\der\mathcal O\subseteq\smallo$. Then
$$\fn^\dagger\preceq 1 \ \Longrightarrow\ v(A_{\ltimes \fn})=v(A),\qquad
\fn^\dagger\succ 1\ \Longrightarrow\ 
\abs{v(A_{\ltimes \fn})-v(A)} \leq -rv(\fn^\dagger).$$
\end{lemma}
\begin{proof}
Let $R:=\Ric A$. Then $v(A_{\ltimes\fn})=v(R_{+\fn^\dagger})$ by [ADH, 5.8.11].
If $\fn^\dagger\preceq 1$, then  $v(R_{+\fn^\dagger})=v(R)$
by [ADH, 4.5.1(i)], hence $v(A_{\ltimes\fn})=v(R)=v(A)$ by [ADH, 5.8.10].
Now suppose $\fn^\dagger\succ 1$. {\it Claim}\/: $v(A_{\ltimes\fn})-v(A) \geq rv(\fn^\dagger)$.
To prove this claim we
replace $A$ by $a^{-1}A$, where $a\asymp A$, to arrange~$A\asymp 1$.
Let~$i$,~$j$ range over~$\{0,\dots,r\}$. We   have 
$R_{+\fn^\dagger} = \sum_i b_i R_i$ where
$$b_i\ =\ \sum_{j\geq i} {j\choose i} a_j R_{j-i}(\fn^\dagger).$$
Take $\fm_i\in K^\times$ as in Lemma~\ref{lem:nepsilon, first part}. 
By Lemma~\ref{Riccatipower+} we have $R_n(\fn^\dagger)\sim (\fn^\dagger)^n$ for all~$n$;
hence $v(b_i)\geq v(\fm_i)$ for all $i$. Thus
$$v(A_{\ltimes\fn})-v(A)\ =\ v(A_{\ltimes\fn})\ =\ v(R_{+\fn^\dagger})\ \geq\ \min_i v(b_i)\ \geq\ v(\fm_0)\ \geq\ rv(\fn^\dagger)$$ by Lemma~\ref{lem:nepsilon, first part}, proving
our claim. 
Applying this claim with $A_{\ltimes \fn}$, $\fn^{-1}$ in place of~$A$,~$\fn$ also yields
$v(A_{\ltimes\fn})-v(A) \leq -rv(\fn^\dagger)$,   thus $\abs{v(A_{\ltimes \fn})-v(A)} \leq -rv(\fn^\dagger)$.
\end{proof}

\begin{remark} Suppose that $\der\mathcal O\subseteq\smallo$ and $\fn^\dagger \succ 1$. Then Lemma~\ref{lem:Atwist} improves on Lem\-ma~\ref{lem:An}, since $v(\fn^\dagger)=o(v\fn)$ by [ADH, 6.4.1(iii)].  
\end{remark}

\begin{lemma}\label{twistprec} Suppose $\der\mathcal O\subseteq\smallo$ and $\fn^\dagger\preceq\fv(A)^{-1}$. 
Let $B\in K[\der]$ and $s\in\N$  be such that $\order(B)\leq s$ and~$B\prec \fv(A)^{s+1}A$. Then $B_{\ltimes\fn} \prec \fv(A_{\ltimes\fn})A_{\ltimes\fn}$. 
\end{lemma}
\begin{proof}
We may assume $B\neq 0$ and $s=\order(B)$. It suffices to show $B_{\ltimes\fn} \prec \fv(A)A$.
If $\fn^\dagger\preceq 1$, then Lemma~\ref{lem:Atwist} applied to $B$ in place of $A$ yields
$B_{\ltimes\fn} \asymp B \prec\fv(A)A$.
Suppose $\fn^\dagger\succ 1$. Then Lemma~\ref{lem:Atwist} gives
$\abs{v(B_{\ltimes \fn})-v(B)} \leq -sv(\fn^\dagger)\leq sv(\fv(A))$ and hence
$B_{\ltimes\fn} \preceq \fv(A)^{-s} B \prec \fv(A) A$.
\end{proof}

\noindent
If   $\der\mathcal{O}\subseteq \smallo$, then we have functions $\dwm_A, \dwt_A\colon \Gamma\to \N$ as defined in~[ADH, 5.6]. 
Combining Lemmas~\ref{lem:fv of perturbed op} and~\ref{twistprec} yields a variant of~[ADH, 6.1.7]:

\begin{cor}\label{cor:Atwist} Suppose $\der\mathcal O\subseteq\smallo$ and $\fn^\dagger\preceq\fv(A)^{-1}$. 
Let $B\in K[\der]$ be such that~$\order(B)\leq r$ and~$B\prec \fv(A)^{r+1}A$. 
Then $\dwm_{A+B}(v\fn)=\dwm_A(v\fn)$ and $\dwt_{A+B}(v\fn)=\dwt_A(v\fn)$. In particular, $$v\fn\in\exc(A+B)\ \Longleftrightarrow\ 
v\fn \in \exc(A).$$
%$(A+B)_{\ltimes\fn}\sim A_{\ltimes\fn}$, $\order(A+B)_{\ltimes\fn}=r$, $\dwt (A+B)_{\ltimes\fn} = \dwt(A_{\ltimes\fn})$, and  $\fv(A+B)_{\ltimes\fn} \sim \fv(A_{\ltimes\fn})$.
\end{cor}

\subsection*{About $A(\fn^q)$ and $A\fn^q$} 
Suppose $\fm^l=\pm \fn^k$ where $k,l\in \Z$, $l\ne0$. Then 
$\fm^\dagger=q\fn^\dagger$ with $q=k/l\in \Q$. In particular, if  $K$ is real closed or algebraically closed, then for any $\fn$ and $q\in \Q$ we have $\fm^\dagger=q\fn^\dagger$ for some $\fm$. 

\medskip
\noindent
{\em Below in this subsection~$K$ is $\d$-valued and $\fn$ is such that for all 
$q\in \Q^{>}$ we are given an element of $K^\times$, denoted by~$\fn^q$ for suggestiveness, with $(\fn^q)^\dagger=q\fn^\dagger$.}

\medskip
\noindent
Let $q\in \Q^{>}$; then $v(\fn^q)=qv(\fn)$: to see this we may arrange that~$K$ is algebraically closed 
by [ADH, 10.1.23], and hence contains an~$\fm$ such that $v\fm=q\,v\fn$ and~$\fm^\dagger=q\fn^\dagger=(\fn^q)^\dagger$, and thus $v(\fn^q)=v\fm=q\,v\fn$.   

\begin{lemma}\label{qlA} Suppose $\fn^\dagger\succeq 1$. %  and for every~$q\in \Q^{>}$ there is given an element of $K^\times$, to be denoted by $\fn^q$ for convenience, with $(\fn^q)^\dagger=q\fn^\dagger$. Let an operator $A=a_0+a_1\der + \cdots + a_r\der^r$ be given with $a_0,\dots, a_r\in K$ and $a_r\ne 0$.
Then for all but finitely many $q\in \Q^{>}$,
$$v\big(A(\fn^q)\big)\ =\ v(\fn^q) + \min_j v(a_j)+jv(\fn^\dagger).$$ 
\end{lemma}
\begin{proof} Let $q\in \Q^{>}$ and take $b_0,\dots,b_r\in K$ with $A\fn^q=b_0+b_1\der+\cdots+b_r\der^r$.
Then 
$$b_0\ =\ A(\fn^q)\ =\ \fn^q \big(a_0 R_0(q\fn^\dagger)+a_1 R_1(q\fn^\dagger)+\cdots+
a_r R_r(q\fn^\dagger)\big).$$
Let $i$,~$j$ range over $\{0,\dots,r\}$.  
By Lemma~\ref{Riccatipower+}, $R_i(q\fn^\dagger) \sim q^i (\fn^\dagger)^i$ for all $i$. Take~$\fm$ (independent of $q$) such that
$v(\fm)=\min_j v(a_j)+jv(\fn^\dagger)$, and let $I$ be the nonempty set of
$i$ with $\fm\asymp a_i(\fn^\dagger)^i$. For $i\in I$ we take $c_i\in C^\times$ such that $a_i(\fn^\dagger)^i \sim c_i\fm$, and set $R:=\sum_{i\in I} c_i Y^i\in C[Y]^{\neq}$.
Therefore, if $R(q)\ne 0$, then 
$$\sum_{i\in I} a_iR_i(q\fn^\dagger)\ \sim\ \fm R(q).$$
Assume $R(q)\ne 0$ in what follows. Then
$$\sum_{i=0}^r a_iR_i(q\fn^\dagger)\ \sim\ \sum_{i\in I} a_iR_i(q\fn^\dagger)\ \sim\ \fm R(q)\ \asymp\ \fm,$$
hence $b_0\asymp \fm\fn^q$, in particular, $b_0\ne 0$.
\end{proof}

\begin{lemma}\label{lem:nepsilon} Assume $\fn^\dagger\succeq 1$ and $[\fv]<[\fn]$ for
$\fv:=\fv(A)$. Then 
$\big[\fv(A\fn^q)\big]<[\fn]$ for all $q\in \Q^>$, and for all but finitely many $q\in \Q^{>}$ 
%with at most $r$ exceptions 
we have 
$\fv(A\fn^q)\preceq \fv$, and thus~$[\fv]\le \big[\fv(A\fn^q)\big]$.
\end{lemma}
\begin{proof} Let $q\in \Q^>$. Then $[\fv]<[\fn]=[\fn^q]$, so $[\fv(A\fn^q)]<[\fn^q]=[\fn]$ by Corollary~\ref{cor:An}.
To show the second part, let $m=\dwt(A)$. 
Replacing $A$ by $a_m^{-1} A$ we arrange $a_m=1$,  so  $a_r = \fv$, $A\asymp 1$.
Take $b_0,\dots,b_r$ with $A\fn^q=b_0+b_1\der+\cdots+b_r\der^r$.
As in the proof of Lemma~\ref{qlA} we obtain an $\fm$ and a polynomial $R(Y)\in C[Y]^{\ne}$ (both independent of $q$) 
%of degree at most~$r$ 
such that $v(\fm)=\min_j v(a_j)+jv(\fn^\dagger)$, and  $b_0\asymp \fm\fn^q$ if~$R(q)\ne 0$. 
Assume $R(q)\ne 0$ in what follows; we show that then $\fv(A\fn^q)\preceq \fv$.  For~$n:= \dwt(A\fn^q)$, 
$$ b_0\fv(A\fn^q)\ \preceq\ b_n\fv(A\fn^q)\ =\ b_r\ =\ \fn^q\fv,$$
hence
$\fv(A\fn^q)\preceq \fv/\fm$.
It remains to note that  
$\fm\succeq a_m(\fn^\dagger)^m=(\fn^\dagger)^m\succeq 1$.  
\end{proof}
 
\begin{lemma}\label{lem:nepsilon, refined} Assume $\fn^\dagger\succeq 1$ and $\fm$ satisfies 
$$v\fm+v(A)\ =\ \min_{0\leq j\leq r}  v(a_j)+jv(\fn^\dagger).$$
Then  $[\fm]\leq[\fn^\dagger]$, with equality if $\dwt(A)>0$, and for all but finitely many $q\in\Q^>$, 
$$A\fn^q\ \asymp\ \fm\, \fn^q\, A, \qquad \fv(A)/\fv(A\fn^q)\ \asymp\ \fm.$$
\end{lemma}
\begin{proof}
Replacing $A$ by $a_m^{-1}A$ where $m=\dwt(A)$ we arrange $a_m=1$, so $a_r=\fv:=\fv(A)$ and~$A\asymp 1$.
Let~$i$,~$j$ range over $\{0,\dots,r\}$. Let $q\in\Q^>$, and
take~$b_i\in K$ such that $A\fn^q = \sum_i b_i \der^i$.
By [ADH, (5.1.3)] we have
$$b_i\ =\   \frac{1}{i!} A^{(i)}(\fn^q)\  =\ \fn^q  \frac{1}{i!}\Ric(A^{(i)})(q\fn^\dagger)\ =\ \fn^q \sum_{j\geq i} {j\choose i} a_j R_{j-i}(q\fn^\dagger).$$
Take $\fm_i\in K^\times$ as in Lemma~\ref{lem:nepsilon, first part}. 
Then  $\fm_0\asymp\fm$ (so $[\fm]\leq[\fn^\dagger]$, with equality if~$m>0$), and $\fm_r\asymp \fv$.
Lemma~\ref{qlA} applied to $A^{(i)}/i!$ instead of $A$ gives that for all but finitely $q\in\Q^>$ we have
$b_i \asymp \fm_i \fn^q$ for all~$i$. Assume that~$q\in \Q^{>}$ has this property.
From  $v(\fm)=v(\fm_0)\leq \cdots\leq v(\fm_r)=v(\fv)$ we obtain
$$v(\fm)+qv(\fn)\ =\ v(b_0)\ \leq\ \cdots\ \leq\ v(b_r)\ =\ v(\fv)+q\,v(\fn).$$
With $n=\dwt(A\fn^q)$ this gives $v(b_0)=\cdots=v(b_n)=v(A\fn^q)$. Thus
$$\fv(A\fn^q)\ =\ b_r/b_n\ \asymp\ b_r/b_0\ \asymp\ 
(\fn^q \fv)/(\fn^q\fm)\ =\ \fv/\fm$$
as claimed.
%We  note that if $\fn^\dagger\succ 1$, then $v(\fm_0)< \cdots< v(\fm_r)$ and hence $v(b_0)< \cdots< v(b_r)$, and so $A\fn^q \sim b_0=A(\fn^q)$.
\end{proof}

\noindent
Let $\fv\in K^\times$ with $\fv \nasymp 1$; so we have the proper convex subgroup of $\Gamma$ given by
$$\Delta(\fv)\ =\ \big\{\gamma\in \Gamma:\, \gamma=o(v\fv)\big\}\ =\ \big\{\gamma\in \Gamma:\, [\gamma]<[\fv]\big\}.$$ 
If $K$ is asymptotic of $H$-type, then we also have the convex subgroup
$$\Delta\ =\ \big\{\gamma\in\Gamma:\, \gamma^\dagger>v(\fv^\dagger) \big\}$$
of $\Gamma$ with $\Delta\subseteq\Delta(\fv)$, and   $\Delta = \Delta(\fv)$ if $K$ is of Hardy type (cf.~Section~\ref{sec:logder}).

\begin{cor}\label{cor:nepsilon} Suppose $\fn^\dagger\succeq 1$ and 
$[\fn^\dagger] < [\fv]$ where $\fv:=\fv(A)$ \textup{(}so $0\ne \fv\prec 1)$.  
Let $A_*\in K[\der]$ and $w\geq r$ be such that
$A_* \prec_{\Delta(\fv)} \fv^w A$. 
Then for all but finitely many~$q\in\Q^>$ we have $\fw:=\fv(A\fn^q)\asymp_{\Delta(\fv)}\fv$  and
$A_*\fn^q \prec_{\Delta(\fw)} \fw^w A\fn^q$.
\end{cor}
\begin{proof} The case $A_*=0$ is trivial, so
assume $A_*\neq 0$.  Take $\fm$ as in Lemma~\ref{lem:nepsilon, refined}, and take $\fm_*$ likewise with $A_*$ in place of $A$. By this lemma, $[\fm],[\fm_*]\leq[\fn^\dagger]< [\fv]$, 
hence $\fm, \fm_*\asymp_{\Delta(\fv)} 1$. 
Moreover, for all but finitely many $q\in\Q^>$ we have $A\fn^q \asymp\fm\fn^q A$,
$A_*\fn^q \asymp\fm_*\fn^q A_*$, and $\fv/\fw\asymp\fm$ where $\fw:=\fv(A\fn^q)$; assume that $q\in \Q^{>}$ has these properties. Then 
$A_* \prec_{\Delta(\fv)} \fv^w A$ yields 
$$A_*\fn^q\ \asymp\ \fm_* \fn^q A_*\ 
\prec_{\Delta(\fv)} \fm \fn^q \fv^w  A\ \asymp\ \fv^w A\fn^q.$$
Now $\fm\asymp_{\Delta(\fv)} 1$ gives $\fv\asymp_{\Delta(\fv)} \fw$, hence 
$A_*\fn^q \prec_{\Delta(\fw)} \fw^w A\fn^q$. 
\end{proof}

\subsection*{The behavior of the span under compositional conjugation}  If $K$ is $H$-asymptotic with asymptotic integration, then $\Psi\cap \Gamma^{>}\ne \emptyset$, but it is convenient not to require ``asymptotic integration'' in some lemmas below. Instead: {\em In this subsection~$K$ is $H$-asymptotic and ungrounded  with $\Psi\cap \Gamma^{>}\ne \emptyset$.}\/ 
We let~$\phi$,~$\fv$  range over~$K^\times$. We say that
$\phi$ is {\em active\/} if $\phi$ is active in $K$. 
%We set $\delta:=v(\fv)$, and if 
%$\fv\prec 1$, we use $\preceq_{\delta}$ to denote the coarsening~
%$\preceq_{\Delta}$ for the convex subgroup $\Delta:=\Delta(\fv)=
%\big\{\gamma\in \Gamma:\, \gamma=o(\delta)\big\}$ of $\Gamma$. We use the %notations~$\prec_\delta$,~$\asymp_{\delta}$,~$\sim_{\delta}$ in the same way. 
Recall from [ADH, pp.~290--292] that $\derdelta$ denotes the derivation $\phi^{-1}\der$ of~$K^\phi$, and that
\begin{equation}\label{eq:Aphi}
A^\phi\ =\ a_r\phi^r\derdelta^r+\text{lower order terms in $\derdelta$.}
\end{equation}

\begin{lemma}\label{lem:v(Aphi)}
Suppose $\fv:=\fv(A)\prec^\flat 1$ and $\phi\preceq 1$ is active. Then 
$$A\ \asymp_{\Delta(\fv)}\ A^\phi,\qquad \fv\ \asymp_{\Delta(\fv)}\ \fv(A^\phi)\ \prec^\flat\ 1, \qquad \fv, \fv(A^\phi)\ \prec_{\phi}^\flat\ 1.
$$
\end{lemma}
\begin{proof} From $\phi^\dagger\prec 1 \preceq \fv^\dagger$ we get $[\phi]<[\fv]$, so $\phi\asymp_{\Delta(\fv)} 1$. Hence $A^\phi \asymp_{\Delta(\fv)} A$
by~[ADH, 11.1.4]. For the rest we can arrange $A\asymp 1$, so $A^\phi\asymp_{\Delta(\fv)} 1$ and $\fv\asymp a_r$. In view of
\eqref{eq:Aphi}
this yields $\fv(A^\phi)\asymp_{\Delta(\fv)} a_r\phi^r\asymp_{\Delta(\fv)} \fv$. 
So $\fv(A^\phi)^\dagger\asymp \fv^\dagger\succeq 1$, which gives $\fv(A^\phi) \prec^{\flat} 1$, and also $\fv, \fv(A^\phi) \prec_{\phi}^\flat 1$.
\end{proof}

%\noindent
%Recall from [ADH, p.~479] our use of the term ``eventually'': a property $S(\phi)$ of active elements  $\phi$ is said to hold {\em eventually\/} if there exists an active $\phi_0$ such that~$S(\phi)$ holds for all active $\phi\preceq \phi_0$. 

\begin{lemma}\label{lem:eventual value of fv} 
If $\nwt(A)=r$, then $\fv(A^\phi)=1$ eventually, and if $\nwt(A)<r$, then $\fv(A^\phi) \prec_\phi^\flat 1$ eventually.
\end{lemma}
\begin{proof}
Clearly, if $\nwt(A)=r$, then $\dwt(A^\phi)=r$ and so $\fv(A^\phi)=1$ eventually.
Suppose $\nwt(A)<r$. To show that $\fv(A^\phi) \prec_\phi^\flat 1$ eventually, we may replace $A$ by~$A^{\phi_0}$ for
suitable active $\phi_0$ and assume that $n:=\nwt(A) = \dwt(A^\phi) = \dwm(A^\phi)$ for all active $\phi\preceq 1$.
Thus $v(A^\phi)=v(A)+nv\phi$ for all active $\phi\preceq 1$ by [ADH, 11.1.11(i)]. Using~\eqref{eq:Aphi} 
we therefore obtain for active $\phi\preceq 1$:
$$\fv(A^\phi)\ \asymp\ a_r\phi^r/a_n\phi^n\ =\ \fv(A)\phi^{r-n}\ \preceq\ \phi^{r-n}\ \preceq\ \phi.$$
Take $x\in K^\times$ with $x\nasymp 1$ and $x'\asymp 1$; then $x\succ 1$, so $x^{-1}\asymp x^\dagger\prec 1$ is active. Hence for active $\phi\preceq x^{-1}$ we have $\phi \prec^\flat_\phi 1$ and thus $\fv(A^\phi)\prec^\flat_\phi 1$. 
\end{proof}

\begin{cor}\label{coruplnwteq}
The following conditions on $K$ are equivalent:
\begin{enumerate}
\item[\textup{(i)}] $K$ is $\upl$-free;
\item[\textup{(ii)}]  $\nwt(B)\leq 1$ for all $B\in K[\der]$ \textup{(}so~$\fv(B^\phi) \prec_\phi^\flat 1$ eventually\textup{)};
\item[\textup{(iii)}] $\nwt(B)\leq 1$ for all $B\in K[\der]$ of order $2$.
\end{enumerate}
\end{cor}
\begin{proof}
The implication (i)~$\Rightarrow$~(ii) follows from [ADH, 13.7.10] and Lemma~\ref{lem:eventual value of fv}, and
(ii)~$\Rightarrow$~(iii) is clear. Suppose  $K$ is not $\upl$-free. Take~$\upl\in K$   such that $\phi^\dagger+\upl\prec\phi$ for all active $\phi$   ([ADH, 11.6.1]);
set~$B:=(\der+\upl)\der=\der^2+\upl\der$. Then for active $\phi$ we have
$B^\phi=\phi^2\big(\derdelta^2+(\phi^\dagger+\upl)\phi^{-1}\derdelta\big)$, so
$\dwt(B^\phi)=2$. Thus~(iii)~$\Rightarrow$~(i). % and $\fv(A^\phi)=1$.
\end{proof}

%\begin{lemma} \marginpar{new lemma}
%Suppose
%$$A=f(\der-g_1)\cdots(\der-g_r)\quad\text{ where $f,g_1,\dots,g_r\in K$.}$$ 
%Then 
%$$A^\phi=\phi^r f(\derdelta-g_1^\phi)\cdots(\derdelta-g_r^\phi)\quad\text{ where $g_i^\phi:=\phi^{-1}\big( g_i-(r-i)\phi^\dagger \big)$ for $i=1,\dots,r$.}$$
%Suppose  $K$ is $\upl$-free. Then for $i=1,\dots,r-1$ we have $g_i^\phi\succ_\phi^\flat 1$ eventually, and $g_r^\phi\succ_\phi^\flat 1$ eventually iff $g_r$ is active.
%\end{lemma}
%\begin{proof}
%The first statement follows from Lemma~\ref{lem:split and compconj}. Suppose $K$ is $\upl$-free and $i\in\{1,\dots,r-1\}$. Note that $g_i^\phi\succ_\phi^\flat 1$ iff $\frac{1}{r-i}g_i-\phi^\dagger \succ^\flat_\phi \phi$. By [ADH, 11.6.1] we can take  an active~$\phi_0$ in $K$ such that $\frac{1}{r-i}g_i-\phi^\dagger \succeq \phi_0$  for all active $\phi\prec\phi_0$. By [ADH, 9.2.11] we have $\phi \prec^\flat_\phi \phi_0$ eventually, and so $\frac{1}{r-i}g_i-\phi^\dagger \succ^\flat_\phi \phi$ eventually. Next, we observe that $g_r^\phi=\phi^{-1}g_r$, so $g_r^\phi\succ 1$ eventually iff $g_r$ is active, and in this case $g_r^\phi\succ_\phi^\flat 1$ eventually.
%\end{proof}

\noindent
Lemma~\ref{lem:v(Aphi)} leads to an  ``eventual'' version of Corollary~\ref{cor:Atwist}:

\begin{lemma}\label{cor:excev stability} Suppose $K$ is $\upl$-free and $B\in K[\der]$ is such that $\order(B)\leq r$ and~$B\prec_{\Delta(\fv)} \fv^{r+1}A$, where $\fv:=\fv(A)\prec^\flat 1$.  Then 
$\exc^{\ev}(A+B) = \exc^{\ev}(A)$.
\end{lemma} 
\begin{proof} By [ADH, 10.1.3, 11.7.18] and Corollary~\ref{cor:13.7.10} we can pass to an extension to arrange that $K$ is
$\upo$-free. Next, by [ADH, 11.7.23, remark following 14.0.1] we extend further to arrange that $K$ is algebraically closed and newtonian, and thus $\d$-valued by Lemma~\ref{lem:ADH 14.2.5}. Then $\exc^{\ev}(A)=v(\ker^{\neq} A)$ by Proposition~\ref{kerexc}, and $A$ splits over $K$ by [ADH, 14.5.3, 5.8.9]. It remains to show that $\exc^{\ev}(A)\subseteq \exc^{\ev}(A+B)$: the reverse inclusion then follows by
interchanging $A$ and $A+B$, using $\fv(A)\sim \fv(A+B)$. 
Let~$\gamma\in\exc^{\ev}(A)$. Take $\fn\in\ker^{\neq} A$ with $v\fn=\gamma$.  Then $A\in K[\der](\der-\fn^\dagger)$
by [ADH, 5.1.21] and so $\fn^\dagger\preceq \fv^{-1}$,
by [ADH, 5.1.22] and  Corollary~\ref{cor:bound on linear factors}. Now
$\exc^{\ev}(A)\subseteq \exc(A)$, so $\gamma=v\fn\in \exc(A+B)$ by Corollary~\ref{cor:Atwist}.  Let $\phi\preceq 1$ be active; it remains to show that then $\gamma\in \exc\big((A+B)^\phi\big)$.  By Lemma~\ref{lem:v(Aphi)}, $A^\phi \asymp_{\Delta(\fv)} A$; 
also   $B^\phi\preceq B$ by~[ADH, 11.1.4].
Lemma~\ref{lem:v(Aphi)} gives $\fv \asymp_{\Delta(\fv)} \fv(A^\phi)$, hence
$B^\phi \prec_{\Delta(\fv)} \fv(A^\phi)^{r+1}A^\phi$. Thus with $K^\phi$, $A^\phi$, $B^\phi$ in the role of $K$, $A$, $B$, the above argument leading to $\gamma\in \exc(A+B)$ gives $\gamma\in \exc(A^\phi+B^\phi)=\exc\big((A+B)^\phi\big)$.
\end{proof} 

\noindent
For $r=1$ we can weaken the hypothesis of  $\upl$-freeness:

\begin{cor}\label{cor:excev stability, r=1}
Suppose $K$ has asymptotic integration,  $r=1$,  and $B\in K[\der]$ of order~$\leq 1$ satisfies  $B\prec_{\Delta(\fv)} \fv^{2}A$, where $\fv:=\fv(A)\prec^\flat 1$.  Then $\exc^{\ev}(A+B) = \exc^{\ev}(A)$. 
\end{cor}
\begin{proof}
Using Lemma~\ref{lem:achieve I(K) subseteq Kdagger} we replace $K$ by an immediate extension to arrange $\I(K)\subseteq K^\dagger$.
Then  $\exc^{\ev}(A)=v(\ker^{\neq} A)$ by Lemma~\ref{lem:v(ker)=exc, r=1}. Now argue as in the proof of Lem\-ma~\ref{cor:excev stability}.
\end{proof}

\noindent
{\em In the next proposition and its corollary $K$ is $\d$-valued  with algebraically closed constant field $C$ and divisible group~$K^\dagger$ of logarithmic derivatives}.
We choose a complement $\Lambda$ of the $\Q$-linear subspace $K^\dagger$ of $K$. Then we have the set $\exc^{\operatorname{u}}(A)$ of ultimate exceptional values of $A$ with respect to $\Lambda$. The following stability result  will be crucial in Section~\ref{sec:ultimate}:

\begin{prop}\label{prop:stability of excu} 
Suppose $K$ is $\upo$-free, $\I(K)\subseteq K^\dagger$, and
$B\in K[\der]$ of order $\le r$ satisfies  $B\prec_{\Delta(\fv)} \fv^{r+1}A$, where $\fv:=\fv(A)\prec^\flat 1$. 
Then $\exc^{\operatorname{u}}(A+B)=\exc^{\operatorname{u}}(A)$.
\end{prop}
\begin{proof}
Let $\Omega$ be the differential fraction field of the universal exponential extension $\Univ=K\big[\!\ex(\Lambda)\big]$ 
of $K$ from Section~\ref{sec:univ exp ext}.
Equip  $\Omega$ with a spectral extension of the valuation of $K$; 
see Section~\ref{sec:valuniv}. Apply Lem\-ma~\ref{cor:excev stability}
to $\Omega$ in place of $K$ to get~$\exc^{\ev}_\Omega(A+B)=\exc^{\ev}_\Omega(A)$. 
Hence $\exc^{\operatorname{u}}(A+B)=\exc^{\operatorname{u}}(A)$ by~\eqref{eq:excevOmega}.
\end{proof}

\noindent
In a similar manner we obtain an analogue of Corollary~\ref{cor:excev stability, r=1}:  

\begin{cor}\label{cor:stability of excu} 
Suppose $K$ has asymptotic integration, $\I(K)\subseteq K^\dagger$, $r=1$, and~$B\in K[\der]$   satisfies~$\order(B)\leq 1$ and~$B\prec_{\Delta(\fv)} \fv^{2}A$, where $\fv:=\fv(A)\prec^\flat 1$. 
Then $\exc^{\operatorname{u}}(A+B)=\exc^{\operatorname{u}}(A)$.
\end{cor}
\begin{proof}
Let $\Omega$ be as in the proof of Proposition~\ref{prop:stability of excu}.
Then $\Omega$ is ungrounded by Lemma~\ref{lem:v(ex(Q))}, hence $\abs{\exc^{\ev}_\Omega(A)}\leq 1$ and $v(\ker^{\neq}_\Omega A)\subseteq\exc^{\ev}_\Omega(A)$ by [ADH, p.~481]. But~$\dim_C \ker_\Omega A=1$, so
$v(\ker_\Omega^{\neq} A)=\exc^{\ev}_\Omega(A)$. The proof of Lemma~\ref{cor:excev stability} with $\Omega$ in place of $K$ now gives $\exc^{\ev}_\Omega(A+B)=\exc^{\ev}_\Omega(A)$, so $\exc^{\operatorname{u}}(A+B)=\exc^{\operatorname{u}}(A)$ by~\eqref{eq:excevOmega}.
\end{proof}

\noindent
In the ``real''  case we have the following variant of Proposition~\ref{prop:stability of excu}: 

\begin{prop}\label{prop:stability of excu, real}  
Suppose $K=H[\imag]$, $\imag^2=-1$, where $H$ is a real closed $H$-field with asymptotic integration such that
$H^\dagger=H$ and $\I(H)\imag\subseteq K^\dagger$.
Let $B\in K[\der]$ of order~$\le r$ be such that~${B\prec_{\Delta(\fv)} \fv^{r+1}A}$ with $\fv:=\fv(A)\prec^\flat 1$. Let $\Lambda$ be a complement of
the subspace~$K^\dagger$ of the $\Q$-linear space $K$. 
Then $\exc^{\operatorname{u}}(A+B)=\exc^{\operatorname{u}}(A)$, where the ultimate exceptional values are with respect to $\Lambda$. 
\end{prop}

\begin{proof}
Take an $H$-closed extension $F$ of $H$ with $C_F=C_H$ as in Corollary~\ref{cor:LambdaL}. Then the algebraically closed $\d$-valued $H$-asymptotic extension $L:=F[\imag]$ of $K$ is $\upo$-free,  $C_L=C$, $\I(L)\subseteq L^\dagger$, and
$L^\dagger\cap K=K^\dagger$. 
Take  a complement $\Lambda_L\supseteq \Lambda$ of the subspace $L^\dagger$ of the $\Q$-linear space $L$.  Let $\Univ_L=L\big[\!\ex(\Lambda_L)\big]$ be
 the universal exponential extension  
of $L$ from Section~\ref{sec:univ exp ext}; it has the  universal exponential extension $\Univ:=K\big[\!\ex(\Lambda)\big]$ of $K$ as a differential subring. 
Let $\Omega$, $\Omega_L$ be the
differential fraction fields of $\Univ$, $\Univ_L$, respectively, and
equip  $\Omega_L$ with a spectral extension of the valuation of $L$; then the restriction of this valuation to
$\Omega$ is  a spectral extension of the valuation of $K$ (see remarks preceding Lemma~\ref{lem:excev cap GammaOmega}).  
Lem\-ma~\ref{cor:excev stability} applied
to $\Omega_L$ in place of $K$ yields~$\exc^{\ev}_{\Omega_L}(A+B)=\exc^{\ev}_{\Omega_L}(A)$,
hence $\exc^{\ev}_\Omega(A+B)=\exc^{\ev}_\Omega(A)$ by
Lemma~\ref{lem:excev cap GammaOmega} and thus $\exc^{\operatorname{u}}(A+B)=\exc^{\operatorname{u}}(A)$.
\end{proof}

%\noindent
%Likewise, using Corollary~\ref{cor:excev stability, r=1} in place of
%Lem\-ma~\ref{cor:excev stability}:

%\begin{cor}\label{cor:stability of excu} \marginpar{new corollary, skipped for now, not used so far}
%Suppose  $\I(K)\subseteq K^\dagger$, $r=1$, and
%$B\in K[\der]$   satisfies~$\order(B)\leq 1$ and~$B\prec_{\Delta(\fv)} \fv^{2}A$, where $\fv:=\fv(A)\prec^\flat 1$. 
%Then $\exc^{\operatorname{u}}(A+B)=\exc^{\operatorname{u}}(A)$.
%\end{cor}

\subsection*{The span of the linear part of a differential polynomial} 
In this subsection~$P\in K\{Y\}^{\neq}$ has order $r$. Recall from [ADH, 5.1] that the {\it linear part}\/ of $P$ is the
differential operator \index{differential polynomial!linear part}\index{linear part!differential polynomial} 
$$L_P\ :=\ \sum_n \frac{\partial P}{\partial Y^{(n)}}(0)\,\der^n \in K[\der]$$
of order~$\leq r$.
We have $L_{P_{\times\fm}}=L_P\fm$ [ADH, p.~242]; hence items~\ref{lem:An}, \ref{cor:An} and \ref{lem:Atwist} above yield information about the span of $L_{P_{\times\fm}}$ (provided $L_P\neq 0$). We now want to similarly investigate the span of the linear part 
$$L_{P_{+a}}\ =\ \sum_n \frac{\partial P}{\partial Y^{(n)}}(a)\,\der^n$$  
of the additive conjugate~$P_{+a}$ of $P$ by some $a\prec 1$.
In the next two lemmas we assume $\order(L_P)=r$ \textup{(}in particular, $L_P\neq 0$\textup{)}, 
$\fv(L_P)\prec 1$, and $a\prec 1$, we set
$$L:=L_P,\quad L^+:=L_{P_{+a}},\quad \fv:=\fv(L),$$ 
%we use the conventions concerning  
%$\prec_\delta$,~$\asymp_{\delta}$,~$\sim_{\delta}$ from the previous subsection, 
and set $L_n:=\frac{\partial P}{\partial Y^{(n)}}(0)$ and $L_n^+:=\frac{\partial P}{\partial Y^{(n)}}(a)$, so $L=\sum_n L_n\der^n$, $L^+=\sum_n L_n^+\der^n$.
Recall from~[ADH, 4.2] the decomposition of $P$ into homogeneous parts: $P=\sum_d P_d$ where $P_d=\sum_{\abs{\i}=d} P_{\i}Y^{\i}$;
we set $P_{>1}:=\sum_{d>1} P_d$.

\begin{lemma}\label{lem:linear part, new}
Suppose $P_{>1} \prec_{\Delta(\fv)} \fv P_1$ and  $n\le r$. Then 
\begin{enumerate}
\item[$\mathrm{(i)}$] $L_r^+ \sim_{\Delta(\fv)} L_r$, and thus $\order(L^+)=\order(L)=r$; 
\item[$\mathrm{(ii)}$] if $L_n\asymp_{\Delta(\fv)} L$, then $L_n^+ \sim_{\Delta(\fv)} L_n$, and so $v(L_n^+)=v(L_n)$;
\item[$\mathrm{(iii)}$] if $L_n\prec_{\Delta(\fv)} L$, then $L_n^+ \prec_{\Delta(\fv)} L$, and so $v(L_n^+)> v(L)$.
\end{enumerate}
In particular, $L^+\sim_{\Delta(\fv)} L$, $\dwt L^+=\dwt L$, 
%and $L_r^*\asymp_\fv \fv  L^*$, so 
and $\fv(L^+)\sim_{\Delta(\fv)}\fv$. 
\end{lemma}
\begin{proof}
Take $Q,R\in K\{Y\}$ with $\deg_{Y^{(n)}} Q\leq 0$ and  $R\in Y^{(n)}K\{Y\}$, such that
$$P\ =\ Q+(L_n+R)Y^{(n)},\qquad\text{so}\qquad  \frac{\partial P}{\partial Y^{(n)}}\  =\ \frac{\partial R}{\partial Y^{(n)}}Y^{(n)}+L_n+R.$$
Now $R\prec_{\Delta(\fv)}  \fv P_1$,
so $\frac{\partial P}{\partial Y^{(n)}}-L_n\prec_{\Delta(\fv)} \fv P_1$. In $K[\der]$ we thus have  
 $$ L_n^+-L_n\ =\ \frac{\partial P}{\partial Y^{(n)}}(a)-L_n\ \prec_{\Delta(\fv)}\   \fv L\ \asymp\ L_r.$$
So $L_n^+-L_n\prec_{\Delta(\fv)} L$ and (taking $r=n$)  
$L_r^+-L_r \prec_{\Delta(\fv)}  L_r$. This yields (i)--(iii). 
\end{proof}

\begin{lemma}\label{lem:linear part, split-normal, new}  
Suppose $P_{>1} \prec_{\Delta(\fv)} \fv^{m+1} P_1$, and let 
$A,B\in K[\der]$ be such that~$L=A+B$, $B\prec_{\Delta(\fv)} \fv^{m+1} L$.  
Then  
$$ L^+\ =\ A+B^+\ \text{ where $B^+\in K[\der]$, $B^+\ \prec_{\Delta(\fv)}\ \fv^{m+1} L^+$.}$$
In particular, $L-L^+\prec_{\Delta(\fv)} \fv^{m+1}L$. 
\end{lemma}

\begin{proof}  
Let $A_n,B_n\in K$ be such that $A=\sum_n A_n\der^n$ and $B=\sum_n B_n\der^n$, so $L_n=A_n+B_n$. Let any $n$ (possibly~$>r$) be given and
take $Q,R\in K\{Y\}$ as in the proof of Lemma~\ref{lem:linear part, new}. Then $R\prec_{\Delta(\fv)} \fv^{m+1} P_1$. 
Since $B\prec_{\Delta(\fv)} \fv^{m+1} L$, this yields
$$\frac{\partial P}{\partial Y^{(n)}} - A_n\ =\ \frac{\partial R}{\partial Y^{(n)}}Y^{(n)} + B_n + R\ \prec_{\Delta(\fv)}\ \fv^{m+1} P_1.$$
We have $L_n^+=\frac{\partial P}{\partial Y^{(n)}}(a)$, so
$$L_n^+-A_n\ =\ \frac{\partial P}{\partial Y^{(n)}}(a) - A_n\ \prec_{\Delta(\fv)}\ \fv^{m+1} L.$$
By Lemma~\ref{lem:linear part, new} we have 
$L^+\sim_{\Delta(\fv)} L$, hence $B^+=L^+-A \prec_{\Delta(\fv)} \fv^{m+1} L^+$. 
\end{proof}

\section{Holes and Slots} \label{sec:holes}

\noindent
{\em Throughout this section $K$ is an $H$-asymptotic field with small derivation and with rational asymptotic integration. We set $\Gamma:= v(K^\times)$}. 
So $K$ is pre-$\d$-valued, $\Gamma\ne \{0\}$ has no least positive element, and $\Psi\cap \Gamma^{>}\ne \emptyset$. We let  $a$, $b$, $f$, $g$ range over $K$, and~$\phi$,~$\fm$,~$\fn$,~$\fv$,~ $\fw$ (possibly decorated) over $K^\times$.   As at the end of the previous section we shorten ``active in $K$'' to ``active''. 

\subsection*{Holes}
A {\bf hole}\/ in $K$ is a triple $(P,\fm,\hat a)$ where 
$P\in K\{Y\}\setminus K$ and
$\hat a$ is an element of~$\hat K\setminus K$, for some
immediate asymptotic extension $\hat K$ of $K$,
such that $\hat a\prec\fm$ and~$P(\hat a)=0$. (The extension $\hat K$ may vary with $\hat a$.) 
The {\bf order}\/, {\bf degree}\/, and {\bf complexity}\/ of a hole
$(P,\fm,\hat a)$ in~$K$ are defined as the order, (total) degree, and complexity, respectively, of the differential
polynomial $P$. 
A hole~$(P,\fm, \hat a)$ in~$K$ is called {\bf minimal}\/ if no hole in $K$ has smaller complexity; then $P$ is a minimal annihilator of $\hat a$ over $K$. \index{hole}\index{hole!minimal}\index{hole!complexity}\index{complexity!hole}\index{minimal hole}\label{p:hole}

\medskip
\noindent
If $(P,\fm,\hat a)$ is a hole in $K$, then $\hat a$ is a $K$-external zero of $P$, in the sense of Section~\ref{sec:complements newton}.
Conversely, every $K$-external zero $\hat a$ of a differential polynomial $P\in K\{Y\}^{\neq}$ 
gives for every $\fm\succ \hat a$ a hole $(P,\fm,\hat a)$ in $K$.
By Proposition~\ref{14.0.1r} and Corollary~\ref{14.5.2.r}:

\begin{lemma}\label{lem:no hole of order <=r} 
Let $r\in\N^{\geq 1}$, and suppose $K$ is $\upl$-free. Then 
$$\text{$K$ is $\upo$-free and $r$-newtonian}\quad\Longleftrightarrow\quad\text{$K$ has no hole of order~$\leq r$.}$$
\end{lemma}

\noindent
Thus for $\upo$-free $K$, being newtonian is equivalent to having no holes.   
Recall that~$K$ being henselian is equivalent to $K$ having no proper immediate algebraic valued field extension, and hence 
to $K$ having no hole of order $0$.

\medskip\noindent
Minimal holes are like the ``minimal counterexamples'' in certain combinatorial settings, and we need to understand such holes in a rather detailed way for later use in inductive arguments. Below we also consider the more general notion of {\em $Z$-minimal hole},
which has an important role to play as well.  We recall that $Z(K,\hat a)$ is the set of all $Q\in K\{Y\}^{\ne}$ that vanish at $(K,\hat a)$ as defined in  [ADH, 11.4]

\begin{lemma}\label{lem:Z(K,hat a)}
Let $(P,\fm,\hat a)$ be a hole in $K$. Then $P\in Z(K,\hat a)$. 
If
$(P,\fm,\hat a)$ is minimal, then $P$ is an element of minimal complexity of $Z(K,\hat a)$.
\end{lemma}
\begin{proof}
Let $a$, $\fv$ with $\hat a-a\prec\fv$. Since $\hat a\notin K$ lies in an immediate extension of~$K$ we can take $\fn$ with $\fn\asymp \hat a-a$.
By [ADH, 11.2.1] we then have $\ndeg_{\prec\fv} P_{+a}\geq\ndeg P_{+a,\times\fn}\geq 1$.
Hence $P\in Z(K,\hat a)$. Suppose $P$ is not of minimal complexity 
in~$Z(K,\hat a)$. Take  $Q\in Z(K,\hat a)$ of minimal
complexity. Then [ADH, 11.4.8] yields a $K$-external zero $\hat b$ of $Q$,
and any $\fn\succ \hat b$ gives a hole $(Q,\fn,\hat b)$ in $K$ of smaller complexity than $(P, \fm, \hat a)$.
\end{proof}

\noindent
In connection with the next result, note that $K$ being $0$-newtonian just means that~$K$ is henselian as a valued field. 

\begin{cor}\label{minholenewt} Suppose $K$ is $\upl$-free and has a minimal hole of order $r\ge 1$. Then~$K$ is $(r-1)$-newtonian,
and $\upo$-free if $r\geq 2$.  
%If in addition $C$ is algebraically closed and
%$L_P\ne 0$, then $L_P\in K[\der]$ splits.  
\end{cor}
\begin{proof}
This is clear for $r=1$ (and doesn't need $\upl$-freeness), and for $r\geq 2$ 
follows from Lemma~\ref{lem:no hole of order <=r}.
\end{proof}

%\begin{proof} Let $Q\in K\{Y\}^{\ne}$ have order $< r$. Using Corollary~\ref{14.5.2.r}  when $r\ge 2$, it suffices to show that then $Q$ has no $K$-external zero. If it had such a zero $\hat{b}$, we would have a hole $(Q, \fn, \hat{b})$ in $K$ of smaller complexity than $(P, \fm, \hat{a})$.    
%\end{proof} 

\begin{cor}\label{corminholenewt} Suppose $K$ is $\upo$-free and has a minimal hole of order $r\ge 2$. Assume also that $C$ is algebraically closed and $\Gamma$ is divisible. Then 
$K$ is $\d$-valued, $r$-linearly closed, and $r$-linearly newtonian.
\end{cor} 
\begin{proof} This follows from Lemma~\ref{lem:ADH 14.2.5}, Corollary~\ref{14.5.3.r}, and Corollary~\ref{minholenewt}.
\end{proof}

\noindent
Here is a   linear version of Lemma~\ref{lem:no hole of order <=r}: 

\begin{lemma}\label{lem:no hole of order <=r, deg 1}
If $K$ is $\upl$-free, then
$$\text{$K$ is $1$-linearly newtonian}\ \Longleftrightarrow\ \text{$K$ has no hole of degree~$1$ and order~$1$.}$$
If $r\in\N^{\geq 1}$ and $K$ is $\upo$-free, then
$$\text{$K$ is $r$-linearly newtonian}\ \Longleftrightarrow\ \text{$K$ has no hole of degree~$1$ and order~$\leq r$.}$$
\end{lemma}
\begin{proof}
The first statement follows from Lemma~\ref{lem:char 1-linearly newt}, and the second statement from Lemma~\ref{lem:char r-linearly newt}.
\end{proof}

\begin{cor}\label{degmorethanone} If $K$ is $\upo$-free and has a minimal hole in $K$ of order~$r$ and degree~$>1$, then $K$ is $r$-linearly newtonian.
\end{cor}
%\begin{proof} Let $Q\in K\{Y\}$ of order $\le r$ be quasilinear with $\deg Q=1$. Then $Q$ has a zero in an immediate asymptotic extension of $K$ by [ADH, 14.0.1 and subsequent remarks]. If such a zero is outside $K$, it would give rise to a hole $(Q,\dots)$ of lower complexity than $(P,\fm,\hat a)$. Thus any such zero must lie in $K$.
%\end{proof}  

\begin{lemma}
Suppose $K$ has a hole $(P,\fm,\hat a)$ of degree $1$, and $L_P\in K[\der]^{\neq}$  splits over $K$. Then $K$ has a hole of complexity $(1,1,1)$.
\end{lemma}
\begin{proof}
Let $(P,\fm,\hat a)$ as in the hypothesis have minimal order. Then~${\order P\ge 1}$, so~$\order P = \order L_P$. 
Take $A,B\in K[\der]$ such that $\order A=1$ and  $L_P=AB$.
If $\order B=0$, then $(P,\fm,\hat a)$ has complexity $(1,1,1)$.
Assume $\order B\ge 1$. Then~$B(\hat a)\notin K$: otherwise, taking $Q\in K\{Y\}$ of degree~$1$ with
$L_Q=B$ and~$Q(0)=-B(\hat a)$ yields a hole $(Q,\fm,\hat a)$ in $K$ where $\deg Q=1$ and $L_Q$ splits over $K$,
and~$(Q,\fm,\hat a)$ has smaller order than $(P,\fm,\hat a)$.
Set $\hat b:=B(\hat a)$ 
and take $R\in K\{Y\}$ of degree $1$ with $L_R=A$ and $R(0)=P(0)$. Then
$$R(\hat b)\ =\ R(0)+L_R(\hat b)\ =\ P(0)+L_P(\hat a)\ =\ 
P(\hat a)\ =\ 0,$$
hence for any $\fn\succ\hat b$, $(R,\fn,\hat b)$ is a hole in $K$ of complexity $(1,1,1)$.
\end{proof}

\begin{cor}\label{cor:minhole deg 1}
Suppose $K$ is $\upo$-free,  $C$ is algebraically closed, and $\Gamma$ is divisible.
Then every minimal hole in $K$ of degree~$1$ has order~$1$.
If in addition $K$ is $1$-linearly newtonian, then  every minimal hole in $K$ has degree~$>1$. 
\end{cor}
\begin{proof}
The first statement follows from Corollary~\ref{corminholenewt} and the preceding lemma.
For the second statement, use the first and Lemma~\ref{lem:no hole of order <=r, deg 1}.  
\end{proof}

\noindent
Let $(P,\fm, \hat a)$ be a hole in $K$. We say $(P,\fm, \hat a)$ is {\bf $Z$-minimal} if
$P$ has minimal complexity in $Z(K,\hat a)$. Thus if
$(P,\fm,\hat a)$ is minimal, then it is $Z$-minimal by Lem\-ma~\ref{lem:Z(K,hat a)}.
If $(P,\fm,\hat a)$ is $Z$-minimal, then by [ADH, remarks following 11.4.3],  the differential polynomial~$P$ is a minimal annihilator of $\hat a$ over $K$.
Note also that~$\ndeg P_{\times\fm} \geq 1$ by~[ADH, 11.2.1]. In more detail:   \index{hole!Z-minimal@$Z$-minimal}\index{Z-minimal@$Z$-minimal!hole}

\begin{lemma}\label{lem:lower bd on ddeg}
Let $(P,\fm,\hat a)$ be a hole in $K$. Then for all $\fn$ with $\hat a \prec \fn \preceq \fm$,
$$1\ \le\ \dval P_{\times \fn}\ \le\ 
\ddeg P_{\times \fn}\ \le\ \ddeg P_{\times\fm}.$$ In particular, $\ddeg_{\prec \fm} P\ge 1$.
\end{lemma}
\begin{proof} Assume $\hat a \prec \fn \preceq \fm$. 
Then $\hat a = \fn \hat b$ with $\hat b\prec 1$; 
put $Q:=P_{\times\fn}\in K\{Y\}^{\neq}$. Then $Q(\hat b)=0$,  hence $D_Q(0)=0$ and so $\dval Q=\dval P_{\times \fn}\ge 1$. The rest follows from
[ADH, 6.6.5(ii), 6.6.7, 6.6.9] and $\Gamma^{>}$ having no least element.
\end{proof}

\noindent
In the next lemma, $(\upl_\rho)$, $(\upo_\rho)$ are pc-sequences in $K$ as in [ADH, 11.5, 11.7].
 
\begin{lemma}\label{lem:upl-free, not upo-free}
Suppose $K$ is $\upl$-free and 
$\upo\in K$ is such that $\upo_\rho\leadsto\upo$
\textup{(}so~$K$ is not $\upo$-free\textup{)}.  Then
we have a  hole $(P,\fm,\upl)$ in $K$ where $P= 2Y'+Y^2+\upo$ and~$\upl_\rho\leadsto\upl$, and
each such hole in $K$ is a $Z$-minimal hole in $K$.
\end{lemma}
\begin{proof}
From [ADH, 11.7.13] we obtain $\upl$ in  an immediate asymptotic extension
of~$K$ such that~$\upl_\rho\leadsto\upl$ and $P(\upl)=0$.
Taking any $\fm$ with $\upl\prec\fm$ then yields a hole $(P,\fm,\upl)$ in~$K$ with $\upl_\rho\leadsto\upl$, and each such hole in $K$ 
is a $Z$-minimal hole in~$K$ by~[ADH,  11.4.13, 11.7.12].
\end{proof}

\begin{cor}\label{cor:upl-free, not upo-free}
If $K$ is $\upl$-free but not $\upo$-free, then each minimal hole in $K$ of positive order has complexity~$(1,1,1)$ or complexity $(1,1,2)$. If $K$ is a Liouville closed $H$-field and not $\upo$-free, then  $(P, \fm, \upl)$ is a minimal hole of complexity $(1,1,2)$,
where  $\upo$, $P$, $\upl$, $\fm$ are as in Lemma~\ref{lem:upl-free, not upo-free}.
\end{cor}

\noindent
Here the second part uses Corollary~\ref{cor:Liouville closed => 1-lin newt} and Lemma~\ref{lem:no hole of order <=r, deg 1}.

\subsection*{Slots}
In some arguments the notion of a hole in $K$ turns out to be too stringent. Therefore we introduce a more
flexible version of it:

\begin{definition}
A {\bf slot} in $K$ is a triple $(P,\fm,\hat a)$ where $P\in K\{Y\}\setminus K$ and~$\hat a$ is an element of~$\hat K\setminus K$, for some immediate asymptotic extension $\hat K$ of $K$, such that $\hat a\prec\fm$ and $P\in Z(K,\hat a)$. The {\bf order}, {\bf degree}, and {\bf complexity} of such a slot  in $K$ are defined to be the order, degree, and complexity of the differential polynomial $P$, respectively. A slot in $K$ of degree $1$ is also called a {\bf linear} slot in $K$. 
A slot $(P,\fm,\hat a)$ in $K$ is {\bf $Z$-minimal} if $P$ is of minimal complexity among elements of~$Z(K,\hat a)$. \index{slot}\index{slot!complexity}\index{complexity!slot}\index{slot!linear}\index{slot!Z-minimal@$Z$-minimal}\index{Z-minimal@$Z$-minimal!slot}\label{p:slot}
\end{definition}

\noindent
Thus  by Lemma~\ref{lem:Z(K,hat a)}, holes in $K$ are slots in $K$, and a hole in $K$ is $Z$-minimal iff  
it is $Z$-minimal as a slot in $K$. 
From  [ADH, 11.4.13] we obtain:

\begin{cor}\label{mindivmin} Let $(P,\fm,\hat a)$ be a $Z$-minimal slot in $K$ and $(a_\rho)$ be a divergent pc-sequence in $K$ such that $a_\rho\leadsto \hat a$. Then $P$ is a minimal differential polynomial of $(a_\rho)$ over $K$.
\end{cor}

\noindent
We say that slots  $(P,\fm,\hat a)$ and $(Q,\fn,\hat b)$ in $K$ are {\bf equivalent}   if
$P=Q$, $\fm=\fn$, and $v(\hat a-a)=v(\hat b-a)$ for all $a$; note that then~$Z(K,\hat a)=Z(K,\hat b)$, so 
$(P,\fm,\hat a)$ is $Z$-minimal iff $(P,\fm,\hat b)$ is $Z$-minimal. Clearly this is an equivalence relation on the class of slots in $K$.
The following lemma often allows us to pass from a $Z$-minimal slot to
a $Z$-minimal hole: \index{slot!equivalence}

\begin{lemma}\label{lem:from cracks to holes}
Let $(P,\fm,\hat a)$ be a $Z$-minimal slot in $K$. Then  $(P,\fm, \hat a)$ is equivalent to a  $Z$-minimal hole  in $K$. 
\end{lemma}
\begin{proof}
By  [ADH, 11.4.8] we obtain $\hat b$ in an immediate asymptotic extension of $K$ with $P(\hat b)=0$
and $v(\hat a-a)=v(\hat b-a)$ for all $a$.
In particular  $\hat b\notin K$, $\hat b\prec\fm$, so~$(P,\fm,\hat b)$ is a hole in $K$ equivalent to $(P, \fm, \hat a)$.
\end{proof}

\noindent
By [ADH, 11.4.8] the extension below containing $\hat b$ is not required to be immediate:

\begin{cor}\label{corisomin}  
If  $(P,\fm, \hat{a})$ is a $Z$-minimal hole in $K$ and   $\hat b$ in an asymptotic extension of $K$ satisfies
$P(\hat b)=0$ and $v(\hat a-a)=v(\hat b-a)$ for all $a$, then there is an isomorphism 
$K\<\hat{a}\>\to K\<\hat{b}\>$ of valued differential fields over $K$ sending~$\hat{a}$ to $\hat{b}$. 
\end{cor} 

\noindent
In particular, equivalent $Z$-minimal holes $(P,\fm, \hat{a})$, $(P,\fm,\hat b)$  in $K$
yield an isomorphism $K\<\hat{a}\>\to K\<\hat{b}\>$ of valued differential fields over $K$ sending~$\hat{a}$ to $\hat{b}$.

\medskip
\noindent
From Lemmas~\ref{lem:no hole of order <=r} and \ref{lem:from cracks to holes} we obtain:

\begin{cor}\label{cor:no dent of order <=r} 
Let $r\in\N^{\geq 1}$, and suppose $K$ is $\upo$-free. Then 
$$\text{$K$ is $r$-newtonian}\quad\Longleftrightarrow\quad\text{$K$ has no  slot of order~$\leq r$.}$$
\end{cor}

\noindent
Let $(P,\fm,\hat a)$ be a slot in $K$. Then $(bP, \fm, \hat a)$ for $b\neq 0$ is a slot in $K$ of the same complexity as $(P,\fm,\hat a)$, and
if $(P,\fm, \hat a)$ is $Z$-minimal, then so is~$(bP, \fm, \hat a)$;
likewise with ``hole in $K$'' in place of ``slot in $K$''. 
For active $\phi$ we have the {\bf compositional conjugate}  $(P^\phi,\fm,\hat a)$ by $\phi$ of $(P,\fm,\hat a)$:\index{slot!compositional conjugate}\index{conjugate!compositional} it is a
slot in $K^\phi$ of the same complexity as $(P, \fm, \hat a)$, it is $Z$-minimal if~$(P, \fm, \hat a)$ is, and it is a hole (minimal hole) in $K^\phi$ if~$(P,\fm,\hat a)$ is a hole (minimal hole, respectively) in $K$.
If the slots~$(P,\fm,\hat a)$, $(Q,\fn,\hat b)$ in $K$ are equivalent, then so are 
$(bP,\fm,\hat a)$, $(bQ,\fn,\hat b)$ for $b\neq 0$, as well as the slots~$(P^\phi,\fm,\hat a)$, $(Q^\phi,\fn,\hat b)$ in $K^\phi$ for active $\phi$.

\medskip
\noindent
The following conventions are in force in the rest of this section: 

\medskip\noindent
{\em We let $r$ range over natural numbers $\ge 1$ and let $(P,\fm,\hat a)$ denote a slot in $K$ of order $r$, so $P\notin K[Y]$ has order $r$. We set $w:=\wt(P)$, so $w\geq r\ge  1$.} 

\medskip
\noindent
Thus $\wt(P_{+a})=\wt(P_{\times \fn})=\wt(P^\phi)=w$.

\subsection*{Refinements and multiplicative conjugates}
For $a$,~$\fn$ such that $\hat a-a\prec\fn\preceq\fm$ 
we obtain a slot $(P_{+a},\fn,{\hat a-a})$ in~$K$
of the same complexity as $(P,\fm,\hat a)$ [ADH, 4.3, 11.4].
Slots of this form are said to {\bf refine}\/ $(P, \fm, \hat a)$ and are called {\bf refinements}\/ of~$(P,\fm,\hat a)$.\index{slot!refinement}\index{refinement} 
A refinement of a refinement of~$(P,\fm,\hat a)$ is itself a refinement of~$(P,\fm,\hat a)$.
If $(P, \fm, \hat a)$ is $Z$-minimal, then so is any refinement of~$(P, \fm, \hat a)$. 
If~$(P,\fm,\hat a)$ is a hole in $K$, then so is each of its refinements, and likewise with ``minimal hole'' in place of ``hole''. 
%Moreover, if~$(P,\fm,\hat b)$ is a hole in $K$ and $v(\hat a-a)=v(\hat b-a)$, then $(P_{+a},\fn,\hat a-a)$ is a refinement of the dent~$(P,\fm,\hat a)$ in $K$ iff~$(P_{+a},\fn,\hat b-a)$ is a refinement of  the hole~$(P,\fm,\hat b)$ in~$K$.
For active~$\phi$, $(P_{+a},\fn,\hat a-a)$ refines $(P,\fm,\hat a)$ iff~$(P^\phi_{+a},\fn,\hat a-a)$ refines $(P^\phi,\fm,\hat a)$.
If $(P,\fm,\hat a)$, $(P,\fm,\hat b)$ are equivalent slots in~$K$ and
$(P_{+a},\fn,{\hat a-a})$ refines $(P,\fm,\hat a)$, then
$(P_{+a},\fn,{\hat b-a})$ refines $(P,\fm,\hat b)$, and the slots
$(P_{+a},\fn,{\hat a-a})$, $(P_{+a},\fn,{\hat b-a})$ in $K$ are equivalent. Conversely, if $(P,\fm,\hat a)$ and $(P,\fm,\hat b)$ are slots in $K$ with equivalent refinements, then $(P,\fm,\hat a)$ and $(P,\fm,\hat b)$ are equivalent.

\begin{lemma}\label{lem:refinements linearly ordered} 
Let $(P_{+a},\fn,\hat a-a)$ be a slot in $K$. Then $(P_{+a},\fn,\hat a-a)$ re\-fines~$(P,\fm,\hat a)$, or
$(P,\fm,\hat a)$ refines $(P_{+a},\fn,\hat a-a)$.
\end{lemma}
\begin{proof}
If $\fn\preceq\fm$, then $\hat a-a\prec\fn\preceq\fm$, so $(P_{+a},\fn,\hat a-a)$ refines~$(P,\fm,\hat a)$,
whereas if~$\fm\prec\fn$, then $(\hat a-a)-(-a)=\hat a\prec\fm\preceq\fn$, so  
$$(P,\fm,\hat a) = \big( (P_{+a})_{+(-a)}, \fm, ({\hat a-a})-(-a) \big)$$ refines $(P_{+a},\fn,\hat a-a)$.
\end{proof}

\begin{lemma}\label{lem:notin Z(K,hata)} 
Let $Q\in K\{Y\}^{\neq}$ be such that $Q\notin Z(K,\hat a)$. Then there is a re\-fine\-ment~$(P_{+a},\fn,\hat a-a)$ of $(P,\fm,\hat a)$ such that $\ndeg  Q_{+a,\times\fn}=0$ and $\hat a-a\prec\fn\prec\hat a$.
\end{lemma}

\begin{proof}
Take~$b$, $\fv$ such that $\hat a-b\prec \fv$ and $\ndeg_{\prec\fv} Q_{+b}=0$. 
We shall find an $a$ such that $\ndeg_{\prec\fv} Q_{+a}=0$, $\hat a-a\preceq\hat a$, and $\hat a-a\prec\fv$:
if~$\hat a-b\preceq\hat a$, we take $a:=b$; if~$\hat a-b\succ\hat a$,
then $-b\sim\hat a-b$ and so $\ndeg_{\prec\fv} Q=\ndeg_{\prec\fv} Q_{+b}= 0$ by [ADH, 11.2.7], hence $a:=0$ works. We next arrange $\hat a-a\prec\hat a$: if $\hat a -a\asymp \hat a$, take~$a_1$ with $\hat a-a_1\prec \hat a$, so $a-a_1\prec\fv$, hence
$\ndeg_{\prec\fv} Q_{+a_1}=\ndeg_{\prec\fv} Q_{+a}=0$, and thus $a$ can be replaced by $a_1$.
Since $\Gamma^{>}$ has no least element, we can choose $\fn$ with $\hat a-a\prec\fn\prec\hat a,\fv$, and then $(P_{+a},\fn,\hat a-a)$ refines $(P,\fm,\hat a)$ as desired.   
\end{proof}

\noindent
If $(P_{+a},\fm,\hat a-a)$ refines~$(P,\fm,\hat a)$,
then 
$D_{P_{+a,\times\fm}}=D_{P_{\times\fm,+(a/\fm)}}=D_{P_{\times\fm}}$ by [ADH, 6.6.5(iii)], and thus 
$$\ddeg P_{+a,\times\fm}\ =\ \ddeg P_{\times\fm}, \qquad
\dval P_{+a,\times\fm}\ =\ \dval P_{\times\fm}.$$ 
In combination with Lemma~\ref{lem:lower bd on ddeg} this has some useful consequences:

\begin{cor}\label{cor:ref 1}
Suppose  $(P,\fm,\hat a)$ is a hole in $K$ and $\ddeg P_{\times\fm}=1$.
Then $\ddeg_{\prec\fm} P=1$, and for all $\fn$ with $\hat a \prec \fn\preceq \fm$, $(P,\fn,{\hat a})$ refines $(P, \fm, \hat a)$ with $\ddeg P_{\times \fn}= \dval P_{\times\fn}=1$.  
\end{cor}

\begin{cor}\label{cor:ref 2}
Suppose $(P_{+a},\fn,{\hat a-a})$ refines
the hole $(P,\fm,\hat a)$ in $K$. Then  
$$\ddeg P_{\times\fm}\ =\ 1\ \Longrightarrow\ \ddeg P_{+a,\times\fn}\ =\ \dval P_{+a,\times\fn}\ =\ 1.$$
\end{cor}
\begin{proof} Use $$1\le \dval P_{+a,\times \fn}\le \ddeg P_{+a,\times \fn}\le \ddeg P_{+a,\times\fm}=\ddeg P_{\times \fm},$$ where the first inequality follows from Lemma~\ref{lem:lower bd on ddeg} applied to
$(P_{+a}, \fn, {\hat a -a})$. 
\end{proof}

\noindent
If $(P_{+a},\fm,\hat a-a)$ refines $(P,\fm,\hat a)$, then in analogy with $\ddeg$ and $\dval$, 
$$\ndeg P_{+a,\times\fm}\ =\ \ndeg P_{\times\fm}, \qquad
\nval P_{+a,\times\fm}\ =\ \nval P_{\times\fm}.$$ 
(Use compositional conjugation by active $\phi$.)
Lemma~\ref{lem:lower bd on ddeg} goes through for slots, provided we use
 $\ndeg$ and $\nval$ instead of 
$\ddeg$ and $\dval$: 

\begin{lemma}\label{lem:lower bd on ndeg}
Suppose $\hat a \prec \fn \preceq \fm$. Then
$$1\ \le\ \nval P_{\times \fn}\ \le\ 
\ndeg P_{\times \fn}\ \le\ \ndeg P_{\times\fm}.$$ %In particular, $\ndeg_{\prec \fm} P\ge 1$.
\end{lemma}
\begin{proof}
By [ADH, 11.2.3(iii), 11.2.5] it is enough to show~$\nval P_{\times\fn}\geq 1$.
Replacing~$(P,\fm,\hat a)$ by its refinement $(P,\fn,\hat a)$ we arrange~$\fm=\fn$.
Now $\Gamma^>$ has no  smallest element, so by definition of~$Z(K,\hat a)$ and [ADH, p.~483] we have  
$$1\ \leq\ \ndeg_{\prec\fm} P\ =\ \max\big\{\!\nval P_{\times\fv} : \fv\prec\fm\big\}.$$
Thus by [ADH, 11.2.5] we can take $\fv$ with 
$\hat a\prec\fv\prec\fm$ with $\nval P_{\times\fv}\geq 1$, and hence $\nval P_{\times\fm}\geq 1$, again  by [ADH, 11.2.5].
\end{proof}

\noindent
Lemma~\ref{lem:lower bd on ndeg} yields results 
analogous to Corollaries~\ref{cor:ref 1} and~\ref{cor:ref 2} above:
% and are proved in a similar way, using [ADH, 11.2.1] 
%in place of Lemma~\ref{lem:lower bd on ddeg} and
%[ADH, 11.2.5, 11.2.6] in place of [ADH, 6.6.7, 6.6.9]:

\begin{cor}\label{cor:ref 1n}
If $\ndeg P_{\times\fm}=1$, then 
%$\ndeg_{\prec\fm} P=1$, and 
for all $\fn$ with $\hat a \prec \fn\preceq \fm$, $(P,\fn,{\hat a})$ re\-fines~$(P, \fm, \hat a)$ and $\ndeg P_{\times \fn}= \nval P_{\times\fn}=1$. 
\end{cor}

\begin{cor}\label{cor:ref 2n}
If $(P_{+a},\fn,{\hat a-a})$ refines  $(P,\fm,\hat a)$, then  
$$\ndeg P_{\times\fm}\ =\ 1\ \Longrightarrow\ \ndeg P_{+a,\times\fn}\ =\ \nval P_{+a,\times\fn}\ =\ 1.$$
\end{cor}

\noindent
Any triple $(P_{\times\fn},\fm/\fn,\hat a/\fn)$ is also
a slot in~$K$, with the same complexity as $(P,\fm, \hat a)$; it is called the {\bf multiplicative conjugate}\/ of $(P,\fm,\hat a)$
by~$\fn$.\index{slot!multiplicative conjugate}\index{conjugate!multiplicative} If $(P,\fm,\hat a)$ is $Z$-mi\-ni\-mal, then so is any multiplicative conjugate.
If $(P,\fm,\hat a)$ is a hole in $K$, then so is any  multiplicative conjugate; likewise with ``minimal hole'' in place
of ``hole''. If two slots in  $K$ are equivalent, then
so are their multiplicative conjugates by $\fn$.

\medskip
\noindent
Refinements and multiplicative conjugates interact in the following way: Suppose $(P_{+a},\fn,\hat a-a)$ refines $(P,\fm,\hat a)$.
Multiplicative conjugation of the slot~$(P_{+a},\fn,{\hat a-a})$ in $K$ by $\fv$
then results in the slot~$(P_{+a,\times\fv},\fn/\fv,(\hat a-a)/\fv)$ in~$K$.
On the other hand, first taking the multiplicative conjugate
$(P_{\times\fv},\fm/\fv,\hat a/\fv)$ of~$(P,\fm,\hat a)$ by $\fv$ and
then refining to
$(P_{\times\fv,+a/\fv},\fn/\fv,\hat a/\fv-a/\fv)$
results in the same slot in $K$, thanks to the identity $P_{+a,\times\fv} = P_{\times\fv,+a/\fv}$.

\subsection*{Quasilinear slots}
Note that $\ndeg P_{\times \fm}\ge 1$ by Lemma~\ref{lem:lower bd on ndeg}. 
We call~$(P,\fm,\hat a)$   {\bf qua\-si\-li\-near}\index{slot!quasilinear} if~$P_{\times\fm}$ is quasilinear, that is, $\ndeg P_{\times\fm}=1$. If $(P,\fm, \hat a)$ is quasilinear, then so is any slot in $K$ equivalent to $(P,\fm,\hat a)$, any multiplicative conjugate of~$(P,\fm,\hat a)$, as well as any refinement of $(P,\fm, \hat a)$,  by Corollary~\ref{cor:ref 2n}. 
If $(P,\fm, \hat a)$ is linear, then it is quasilinear by Lemma~\ref{lem:lower bd on ndeg}. 

\medskip\noindent
Let $(a_\rho)$ be a divergent pc-sequence in $K$ with $a_\rho\leadsto\hat a$ and for each index $\rho$, let~$\fm_\rho\in K^\times$ be such that $\fm_\rho\asymp\hat a-a_\rho$. 
%Then the sequence $(v\fm_{\rho})$ in $\Gamma$ is cofinal in $v(\hat a-K)$. 
Take an index $\rho_0$ such that $\fm_\sigma\prec\fm_\rho\prec \fm$ for all $\sigma>\rho\geq\rho_0$, cf.~[ADH, 2.2]. 
%Increasing $\rho_0$ if necessary we also arrange that $\fm_\rho \prec \fm$ for all $\rho\geq\rho_0$.  

 \begin{lemma}\label{lem:pc vs dent} Let $\sigma\geq\rho\geq \rho_0$. Then 
  \begin{enumerate}
  \item[\textup{(i)}]  $(P_{+a_{\rho+1}},\fm_{\rho},\hat a-a_{\rho+1})$ is a refinement of $(P,\fm,\hat a)$;  
  \item[\textup{(ii)}] if $(P_{+a},\fn,\hat a-a)$ is a refinement of $(P,\fm,\hat a)$, then    $\fm_\rho\preceq\fn$ for all sufficiently large $\rho$, and for such $\rho$, $(P_{+a_{\rho+1}},\fm_{\rho},{\hat a-a_{\rho+1}})$ refines~$(P_{+a},\fn,\hat a-a)$;  
    \item[\textup{(iii)}]     $(P_{+a_{\sigma+1}},\fm_{\sigma},\hat a-a_{\sigma+1})$  refines $(P_{+a_{\rho+1}},\fm_{\rho},\hat a-a_{\rho+1})$.
  \end{enumerate}
 \end{lemma}
\begin{proof}
Part (i) follows from $\hat a-a_{\rho+1}\asymp\fm_{\rho+1}\prec\fm_\rho\preceq \fm$. For (ii) 
let~$(P_{+a},\fn,{\hat a-a})$ be a refinement of $(P,\fm,\hat a)$. Since $\hat a-a\prec\fn$,  we have $\fm_\rho\preceq\fn$ for all sufficiently large $\rho$. For such $\rho$,    with~$b:=a_{\rho+1}-a$ we have
$$(P_{+a_{\rho+1}},\fm_{\rho},\hat a-a_{\rho+1})\  =\ \big( (P_{+a})_{+b}, \fm_\rho, (\hat a-a)-b\big)$$
and 
$$(\hat a-a)-b\ =\ \hat a-a_{\rho+1}\ \asymp\ \fm_{\rho+1}\ \prec\ \fm_\rho\ \preceq\ \fn.$$
Hence  $(P_{+a_{\rho+1}},\fm_{\rho},{\hat a-a_{\rho+1}})$ refines~$(P_{+a},\fn,\hat a-a)$. Part (iii) follows from~(i) and~(ii).
\end{proof}

\noindent
Let $\mathbf a=c_K(a_\rho)$ be the cut defined by $(a_\rho)$ in $K$ and $\ndeg_{\mathbf a} P$ be the Newton degree of $P$ in $\mathbf a$ as introduced in [ADH, 11.2].
Then $\ndeg_{\mathbf a} P$ is the eventual value of~$\ndeg P_{+a_\rho,\times\fm_\rho}$.
Increasing $\rho_0$ we arrange that additionally for all 
$\rho\geq\rho_0$ we have~$\ndeg P_{+a_\rho,\times\fm_\rho}=\ndeg_{\mathbf a} P$. 

\begin{cor}\label{cor:quasilinear refinement}
$(P,\fm,\hat a)$ has a quasilinear refinement iff $\ndeg_{\mathbf a} P=1$.
\end{cor}
\begin{proof}
By Lemma~\ref{lem:lower bd on ndeg} and [ADH, 11.2.8] we have
\begin{equation}\label{eq:quasilinear refinement}
1 \leq \ndeg P_{+a_{\rho+1},\times\fm_{\rho}} = \ndeg P_{+a_\rho,\times\fm_\rho}.
\end{equation}
Thus if $\ndeg_{\mathbf a} P=1$,  then   for $\rho\geq\rho_0$,
 the refinement $(P_{+a_{\rho+1}},\fm_{\rho},\hat a-a_{\rho+1})$ of~$(P,\fm,\hat a)$ is quasilinear. Conversely, if $(P_{+a},\fn,\hat a-a)$ is a quasilinear refinement of~$(P,\fm,\hat a)$, then
Lemma~\ref{lem:pc vs dent}(ii) yields a $\rho\geq\rho_0$ such that  $\fm_\rho\preceq\fn$, and
then~$(P_{+a_{\rho+1}},\fm_{\rho},{\hat a-a_{\rho+1}})$ in $K$ refines~$(P_{+a},\fn,\hat a-a)$ and hence  is also quasilinear, so
$\ndeg_{\mathbf a} P=\ndeg P_{+a_\rho,\times\fm_\rho}=1$ by \eqref{eq:quasilinear refinement}.  
\end{proof}

%\noindent
%From the previous corollary we obtain:

\begin{lemma}\label{lem:quasilinear refinement} Assume $K$ is $\d$-valued and $\upo$-free, and $\Gamma$ is divisible.
Then every $Z$-minimal slot in $K$ of positive order has a quasilinear refinement. 
\end{lemma}
\begin{proof}
Suppose $(P,\fm,\hat a)$ is $Z$-minimal. Take a divergent pc-sequence $(a_\rho)$ in $K$ such that 
$a_\rho\leadsto \hat a$.
Then $P$ is a minimal differential polynomial of $(a_\rho)$ over $K$,  by Corollary~\ref{mindivmin}. Hence 
$\ndeg_{\boldsymbol a} P=1$ by [ADH, 14.5.1], where ${\boldsymbol a}:=c_K(a_\rho)$. 
Now Corollary~\ref{cor:quasilinear refinement}  gives a quasilinear refinement of $(P,\fm,\hat a)$.
%This yields an $a\sim \hat a$ and a $g\asymp \hat a -a$ such that  $\ndeg P_{+a,\times g}=1$. Then [ADH, 13.6.15] gives a quasilinear refinement $(P_{+a}, \fn, \hat a -a)$ of  $(P,\fm,\hat a)$.
\end{proof}
 
\begin{remark}
Suppose $K$ is a real closed $H$-field that is  $\upl$-free but not $\upo$-free.
(For example, the real closure of the $H$-field $\R\langle\upo\rangle$
from [ADH, 13.9.1] satisfies these conditions, by~[ADH, 11.6.8, 11.7.23, 13.9.1].)
Take $(P,\fm,\upl)$  as in Lemma~\ref{lem:upl-free, not upo-free}. Then by Corollary~\ref{cor:quasilinear refinement} and
[ADH, 11.7.9],  $(P,\fm,\upl)$ has no quasilinear refinement. Thus Lemma~\ref{lem:quasilinear refinement}
 fails if 
``$\upo$-free'' is replaced by ``$\upl$-free''.
\end{remark}

\begin{lemma}\label{lem:zero of P}   
% \marginpar{not used for Hardy business; needed for definable closure and of independent interest}
Let $L$ be an $r$-newtonian $H$-asymptotic extension of $K$ such that~$\Gamma^<$ is cofinal in $\Gamma_{L}^<$,  and suppose $(P,\fm,\hat a)$ is quasilinear. Then $P(\hat b)=0$ and $\hat b\prec\fm$ for some $\hat b\in L$.
\end{lemma}
\begin{proof}
Lemma~\ref{lem:lower bd on ndeg} and $\ndeg P_{\times\fm}=1$ gives $\fn\prec \fm$ with $\ndeg_{\times\fn} P=1$.
By [ADH, p.~480], $\ndeg P_{\times\fn}$ does not change in passing from $K$ to $L$. As $L$ is $r$-newtonian this yields $\hat b\preceq\fn$ in $L$   with $P(\hat b)=0$.
\end{proof}

\noindent
In the next two corollaries we assume that  $K$ is  $\d$-valued and $\upo$-free, and that $L$ is a newtonian
$H$-asymptotic extension of $K$.

\begin{cor}\label{cor:find zero of P}
If $(P,\fm,\hat a)$ is quasilinear,  then $P(\hat b)=0$, $\hat b\prec\fm$  for some~$\hat b\in L$.
\end{cor}
\begin{proof} By \cite[Theorem~B]{Nigel19}, $K$ has a newtonization $K^*$ inside $L$. Such  $K^*$ is $\d$-algebraic over $K$ by [ADH, remarks after 14.0.1], 
 so $\Gamma^{<}$ is cofinal in $\Gamma_{K^*}^{<}$ by Theorem~\ref{thm:ADH 13.6.1}. Thus 
we can apply Lemma~\ref{lem:zero of P} to $K^*$ in the role of $L$. 
\end{proof}

\noindent
Here is a variant of Lemma~\ref{lem:from cracks to holes}:

\begin{cor}\label{cor:find zero of P, 2} 
Suppose   $\Gamma$ is divisible and  
 $(P,\fm,\hat a)$ is  $Z$-minimal. Then there exists~$\hat b\in L$ such that $K\<\hat b\>$ is an immediate extension of $K$ and~$(P,\fm,\hat b)$ is a hole in $K$ equivalent to~$(P,\fm,\hat a)$. \textup{(}Thus   if  $(P,\fm,\hat a)$ is also a  hole in $K$,
then there is an embedding~$K\langle \hat a\rangle\to L$ of valued differential fields over $K$.\textup{)}
\end{cor}
\begin{proof}
By Lemma~\ref{lem:quasilinear refinement} we may refine $(P,\fm,\hat a)$ to arrange that~$(P,\fm,\hat a)$ is quasilinear.
Then [ADH, 11.4.8] gives $\hat b$ in an immediate $H$-asymptotic extension of~$K$ with~$P(\hat b)=0$ and $v(\hat a-a)=v(\hat b-a)$ for all $a$. So $(P,\fm,\hat b)$ is a hole in $K$
  equivalent to~$(P,\fm,\hat a)$. 
The immediate $\d$-algebraic extension $K\langle \hat b\rangle$ of $K$ is $\upo$-free by Theorem~\ref{thm:ADH 13.6.1}. Then [ADH, remarks following 14.0.1] gives a newtonian $\d$-algebraic immediate extension $M$ of $K\langle \hat b\rangle$ and thus of $K$. Then $M$ is a newtonization of $K$ by [ADH, 14.5.4] and thus embeds over $K$ into $L$. The rest follows from Corollary~\ref{corisomin}.
%Then there is a $K$-embedding of
%  $L$ into $M$. By [ADH, 14.1.8] its image contains $\hat b$.
\end{proof}

\begin{remark}
Lemma~\ref{lem:quasilinear refinement} and Corollary~\ref{cor:find zero of P, 2}   go through with the hypothesis ``$\Gamma$~is divisible'' replaced by ``$K$ is henselian''.
The proofs are the same, using \cite[3.3]{Nigel19} in place of [ADH, 14.5.1] in the proof of Lemma~\ref{lem:quasilinear refinement},  and \cite[3.5]{Nigel19} in place of~[ADH, 14.5.4] in the proof of Corollary~\ref{cor:find zero of P, 2}.
\end{remark}

\noindent
For $r=1$ we can   weaken  the hypothesis of $\upo$-freeness in Corollary~\ref{cor:find zero of P, 2}:

\begin{cor}\label{cor:find zero of P, 3}
Suppose $K$ is $\upl$-free and $\Gamma$ is divisible, and 
$(P,\fm,\hat a)$ is $Z$-minimal of order~$r=1$ with a quasilinear refinement.
Let
 $L$ be a newtonian $H$-asymptotic extension of $K$.
  Then there exists~$\hat b\in L$ such that $K\langle\hat b\rangle$ is an immediate extension of $K$ and
  $(P,\fm,\hat b)$ is a hole in $K$ equivalent to $(P,\fm,\hat a)$. \textup{(}So if 
  $(P,\fm,\hat a)$ is also a hole in $K$, then we have
  an embedding~$K\langle \hat a\rangle\to L$ of valued differential fields over $K$.\textup{)}
\end{cor}
\begin{proof}
%Refining $(P,\fm,\hat a)$ we arrange that it is quasilinear and isolated and satisfies the conclusion of Corollary~\ref{cor:2.12 isolated, 2, r=1}.
Take a divergent pc-sequence $(a_\rho)$ in $K$ with~$a_\rho\leadsto\hat a$. Then $\ndeg_{\mathbf a}P=1$ for~$\mathbf a:=c_K(a_\rho)$, by
Corollary~\ref{cor:quasilinear refinement}, and~$P$ is a minimal differential polynomial of $(a_\rho)$ over~$K$, by [ADH, 11.4.13].
The equality $\ndeg_{\mathbf a}P=1$ remains valid when passing from $K$, $\mathbf a$ to $L$, $c_L(a_\rho)$, respectively,
by Lemma~\ref{lem:11.2.13 invariant}. Hence [ADH, 14.1.10] yields $\hat b\in L$ such that $P(\hat b)=0$ and~$a_\rho\leadsto\hat b$, so $v(\hat a-a)=v(\hat b-a)$ for all $a$.
Then~$K\langle\hat b\rangle$ is an immediate extension of $K$ by [ADH, 9.7.6], so~$(P,\fm,\hat b)$ is a hole in $K$ equivalent to $(P,\fm,\hat a)$.
For the rest use Corollary~\ref{corisomin}. 
\end{proof}

\subsection*{The linear part of a slot}  
We define the {\bf linear part}\/ of  $(P,\fm,\hat a)$ to be the linear part 
$L_{P_{\times\fm}}\in K[\der]$ of~$P_{\times\fm}$.\index{slot!linear part}\index{linear part!slot} By [ADH, p.~242] and Lemma~\ref{lem:separant fms} we have
$$L_{P_{\times\fm}}\ =\ L_P\, \fm\ =\ \sum_{n=0}^r \frac{\partial P_{\times\fm}}{\partial Y^{(n)}}(0)\,\der^n\ =\ \fm S_{P}(0)\der^r +\text{lower order terms in $\der$}.$$
The slot $(P,\fm,\hat a)$ has the same linear part as each of its multiplicative conjugates.
The linear part of a refinement $(P_{+a},\fn,\hat a-a)$ of $(P,\fm,\hat a)$ is
given by
\begin{align*}
L_{P_{+a,\times\fn}}\	=\ L_{P_{+a}}\fn\ 
							&=\ \sum_{m=0}^r \left(\sum_{n=m}^r {n\choose m} \fn^{(n-m)} \frac{\partial P}{\partial Y^{(n)}}(a)\right)\der^m\\ 
							&=\ \fn\,  S_P(a)\,\der^r +\text{lower order terms in $\der$.}
\end{align*}
(See [ADH, (5.1.1)].) 
By [ADH, 5.7.5] we have $(P^\phi)_d=(P_d)^\phi$ for $d\in \N$; in particular $L_{P^\phi}=(L_P)^\phi$
and so $\order(L_{P^\phi})=\order(L_P)$.
A particularly favorable situation occurs when $L_P$ splits over a given differential field extension $E$ of $K$ (which includes requiring $L_P\ne 0$). Typically, $E$ is an algebraic closure of $K$.  In any case, $L_P$ splits over $E$ iff $L_{P_{\times\fn}}$ splits over $E$, iff $L_{P^\phi}$ splits over $E^\phi$. Thus:

\begin{lemma}\label{lem:deg 1 cracks splitting}  
Suppose $\deg P=1$ and $L_P$  splits over $E$. Then the linear part of any refinement of $(P, \fm, \hat a)$ and any
multiplicative conjugate of $(P, \fm, \hat a)$ also splits over $E$, and any compositional conjugate of $(P, \fm, \hat a)$ by an active $\phi$
 splits over~$E^\phi$.
\end{lemma}
%\begin{proof} Use that $L_{P_{+a}}=L_P$, $L_{P_{\times\fn}}=L_P\fn$, and $L_{P^\phi}=(L_P)^\phi$.
%\end{proof}

\noindent
Let $\i=(i_0,\dots,i_r)$ range over $\N^{1+r}$. As in [ADH, 4.2] we set 
$$P_{(\i)}\ :=\ \frac{P^{(\i)}}{\i !} \qquad\text{where
$P^{(\i)}\ :=\ 
\frac{\partial^{|\i|}P}{\partial^{i_0}Y\cdots \partial^{i_r}Y^{(r)}}$.}$$
If $\abs{\i}=i_0+\cdots+i_r\geq 1$, then $\cc(P_{(\i)})<\cc(P)$. Note that for $\i=(0,\dots,0,1)$ we have $P_{(\i)}=S_P\ne 0$, since $\order P =r$. 
We now aim for Corollary~\ref{cor:order L=r}.

%In Corollary~\ref{cor:order L=r} below we show that every $Z$-minimal hole has a re\-fine\-ment $(P,\fm,\hat a)$ such that 
%for every $\i$ with $|\i|\ge 1$ the valuation of the coefficient $(P_{+a})_{\i}$ of~$P_{+a}$ is independent of $a$, for refinements $(P_{+a},\fn,\hat a-a)$ of $(P,\fm,\hat a)$.

\begin{lemma}\label{lem:ndeg coeff stabilizes, 1}
Suppose that $(P,\fm,\hat a)$ is $Z$-minimal. Then $(P,\fm,\hat a)$ has a refinement $(P_{+a},\fn,{\hat a-a})$
% there exists a refinement $(P_{+a},\fn,{\hat a-a})$ of $(P,\fm,\hat a)$
such that for all $\i$ with $\abs{\i}\geq 1$ and $P_{(\i)}\neq 0$,
$${\ndeg\,(P_{(\i)})_{+a,\times\fn}\ =\ 0}.$$ 
\end{lemma}
\begin{proof}
Let $\i$ range over the (finitely many) elements of $\N^{1+r}$ satisfying $\abs{\i}\geq 1$ and~$P_{(\i)}\neq 0$.
%Let $n$ range over the elements of $\{0,\dots,r\}$ satisfying  $\frac{\partial P}{\partial Y^{(n)}}\neq 0$.
Each $P_{(\i)}$ has smaller complexity than $P$, so~$P_{(\i)}\notin Z(K,\hat a)$.
Then~$Q:=\prod_{\i} P_{(\i)}\notin Z(K,\hat a)$ by [ADH, 11.4.4], so Lemma~\ref{lem:notin Z(K,hata)} gives a refinement~$(P_{+a},\fn,\hat a-a)$ of $(P,\fm,\hat a)$ with $\ndeg Q_{+a,\times\fn}=0$.
Then~$\ndeg\, (P_{(\i)})_{+a,\times\fn}=0$ for all $\i$, by [ADH, remarks before 11.2.6]. 
\end{proof}

\noindent
From [ADH, (4.3.3)] we   recall that $(P_{(\i)})_{+a}=(P_{+a})_{(\i)}$. Also recall that $(P_{+a})_{\i}=P_{(\i)}(a)$ by Tay\-lor expansion. In particular, if $P_{(\i)}=0$, then $(P_{+a})_{\i}=0$.  

\begin{lemma}\label{lem:ndeg coeff stabilizes, 2} 
Suppose $(P_{+a},\fn,\hat a-a)$ refines $(P,\fm,\hat a)$ and~$\i$ is such that~${\abs{\i}\geq 1}$, $P_{(\i)}\neq 0$, 
and $\ndeg\, (P_{(\i)})_{\times\fm}=0$.  Then
$$\ndeg\, (P_{(\i)})_{+a,\times\fn}\ =\ 0, \qquad {(P_{+a})_{\i}\ \sim\ P_{\i}}.$$
\end{lemma}
\begin{proof}
Using [ADH, 11.2.4, 11.2.3(iii), 11.2.5] we get
$$\ndeg\, (P_{(\i)})_{+a,\times\fn}\ =\ 
\ndeg\,(P_{(\i)})_{+\hat a,\times\fn}\ \leq\ 
\ndeg\, (P_{(\i)})_{+\hat a,\times\fm}\ =\ 
\ndeg\, (P_{(\i)})_{\times\fm}\ =\ 0,$$
so $\ndeg\,(P_{(\i)})_{+a,\times\fn}=0$. Thus $P_{(\i)}\notin Z(K,\hat a)$, hence $(P_{+a})_{\i}=P_{(\i)}(a)\sim P_{(\i)}(\hat a)$ by~[ADH, 11.4.3];
applying this to $a=0$, $\fn=\fm$ yields $P_{\i}=P_{(\i)}(0)\sim P_{(\i)}(\hat a)$. 
%Since $\order P=r$, we have $\frac{\partial P}{\partial Y^{(r)}}\neq 0$ and thus $\frac{\partial P}{\partial Y^{(r)}}(a)\neq 0$, that is, $\order L_{P+a}=r$. The rest follows easily.
\end{proof}

%\noindent
%Applying this to $\i=(0,\dots,0,1)$ yields:

%\begin{cor}\label{cor:sep} If $\ndeg S_{P_{\times \fm}}=0$, then for every refinement
%$(P_{+a},\fn, \hat a-a)$ of $(P,\fm,\hat a)$ we have
%$\order L_{P_{a,\times \fn}}=\order L_{P_{+a}}=\order L_P=\order L_{P_{\times \fm}}=r$. 
%\end{cor} 

\noindent
Combining Lemmas~\ref{lem:ndeg coeff stabilizes, 1} and~\ref{lem:ndeg coeff stabilizes, 2} gives:  

\begin{cor}\label{cor:order L=r}
Every $Z$-minimal slot in $K$ of order $r$ has a refinement $(P,\fm,\hat a)$ such that for all
refinements $(P_{+a},\fn,\hat a-a)$ of $(P,\fm,\hat a)$ and all $\i$ with $\abs{\i}\geq 1$ and~$P_{(\i)}\neq 0$ we have
$(P_{+a})_{\i}\sim P_{\i}$ \textup{(}and thus $\order L_{P_{+a}}=\order L_P=r$\textup{)}.% and $\fv(L_{P_{+a}})\sim \fv(L_P)$.
\end{cor}

\noindent
Here the condition ``of order $r$'' may seem irrelevant, but
is forced on us because refinements preserve order and by our convention that $P$ has order $r$.

\subsection*{Special slots} 
The slot $(P,\fm, \hat a)$ in $K$ is said to be {\bf special\/} if $\hat a/\fm$ is special over~$K$ in the sense of [ADH, p.~167]: some nontrivial convex subgroup $\Delta$ of $\Gamma$ is cofinal in~$v\big(\frac{\hat{a}}{\fm}-K\big)$.\index{slot!special}\index{special!slot}  
If $(P,\fm,\hat a)$ is special, then so are $(bP,\fm,\hat a)$ for $b\neq 0$, any multiplicative conjugate of $(P,\fm, \hat a)$, any compositional conjugate of $(P,\fm,\hat a)$, and any slot in $K$ equivalent to $(P,\fm,\hat a)$. Also, by Lemma~\ref{lem:special refinement}:

\begin{lemma}\label{speciallemma} If $(P,\fm, \hat a)$ is special, then so is any refinement.
\end{lemma}
%\begin{proof} Assume $\hat a/\fm$ is $\Delta$-special, where~$\Delta$ is a nontrivial convex subgroup of $\Gamma$.
%Let~$(Q, \fn, \hat b)$ be a refinement of $(P, \fm, \hat a)$. {\em Claim}: $\hat b/\fn$ is $\Delta$-special. 
%We have $\hat b=\hat a -a$ with
%$\hat a -a \prec \fn\preceq \fm$, hence 
%$$0\ \le\ v\left(\frac{\fn}{\fm}\right)\ <\ v\left(\frac{\hat a}{\fm}-\frac{a}{\fm}\right)\in \Delta, \qquad \frac{\hat b}{\fn}\ =\ \left(\frac{\hat a}{\fm} -\frac{a}{\fm}\right)\cdot \frac{\fm}{\fn},$$
%so $v(\fm/\fn)\in \Delta$. 
%Using this to obtain the second equality below, we obtain
%$$v\left(\frac{\hat a}{\fm} - K\right)\ =\ v\left(\left(\frac{\hat a}{\fm} -\frac{a}{\fm}\right) - K\right)\ =\ v\left(\left(\frac{\hat a}{\fm} -\frac{a}{\fm}\right)\cdot\frac{\fm}{\fn} -K\right)\ =\ v\left(\frac{\hat b}{\fn}-K\right),$$
%which proves our claim. 
%Next, suppose $(Q, \fn, \hat b)$ is a dent in $K$ which is equivalent to $(P,\fm,\hat a)$; thus $P=Q$, $\fm=\fn$. {\em Claim}: $\hat b/\fm$ is $\Delta$-special. To see this note that 
%$$v\left(\frac{\hat a}{\fm}-a\right)=v(\hat a-a\fm)-v(\fm)=v(\hat b-a\fm)-v(\fm)=v\left(\frac{\hat b}{\fm}-a\right),$$
%so $v\big(\frac{\hat a}{\fm}-K\big)=v\big(\frac{\hat b}{\fm}-K\big)$.
%\end{proof}

\noindent
Here is our main source of special slots: 

\begin{lemma}\label{lem:special dents} 
Let $K$ be $r$-linearly newtonian, and $\upo$-free if~$r>1$. Suppose $(P,\fm, \hat a)$ is  quasilinear, and  $Z$-minimal or a hole in $K$. Then  $(P,\fm, \hat a)$ is special.
\end{lemma}
\begin{proof} 
Use Lemma~\ref{lem:from cracks to holes} to arrange $(P,\fm,\hat a)$ is a hole in~$K$.
Next arrange $\fm=1$ by replacing  $(P,\fm,\hat a)$ with $(P_{\times\fm},1,\hat a/\fm)$. So $\ndeg P=1$, hence $\hat a$ is special over $K$ by Proposition~\ref{nepropsp} (if $r>1$) and~\ref{nepropsp, r=1} (if $r=1$).
\end{proof}

\noindent
Next an approximation result used in the proof of Corollary~\ref{mfhc} in Part~\ref{part:Hardy fields}:  

\begin{lemma}\label{lem:small P(a)}
Suppose $\fm=1$, $(P,1,\hat a)$ is special and $Z$-minimal, and  $\hat a-a\preceq\fn\prec 1$ for some $a$.
Then $\hat a-b\prec\fn^{r+1}$ for some $b$, and $P(b)\prec\fn P$ for any such~$b$.
\end{lemma}

\begin{proof}
Using Lemma~\ref{lem:from cracks to holes} we arrange $P(\hat a)=0$.
The differential po\-ly\-no\-mial $Q(Y):=\sum_{\abs{\i}\geq 1} P_{(\i)}(\hat a)Y^{\i}\in \hat K\{Y\}$ has order~$\leq r$  and $\val(Q)\geq 1$,
and Taylor expansion   yields,  for all $a$:
$$P(a)\ =\ P(\hat a) + \sum_{\abs{\i}\geq 1} P_{(\i)}(\hat a)(a-\hat a)^{\i}\ =\ Q(a-\hat a).$$
Since $\hat a$ is special over $K$, we have $b$ with $\hat a-b\prec\fn^{r+1}$, and then by Lemma~\ref{lem:diff operator at small elt} we have $Q(b-\hat a)\prec\fn Q\preceq\fn P$.
\end{proof}

\section{The Normalization Theorem}\label{sec:normalization}

\noindent
{\em Throughout this section $K$ is an $H$-asymptotic field with small derivation and with rational asymptotic integration. We set $\Gamma:= v(K^\times)$}. 
The notational conventions introduced in the last section remain in force:
$a$,~$b$,~$f$,~$g$ range over $K$; $\phi$,~$\fm$,~$\fn$,~$\fv$,~$\fw$ over $K^\times$. As at the end of Section~\ref{sec:span} we shall frequently use for $\fv\prec 1$ the coarsening of $v$ by the convex subgroup $\Delta(\fv)=\big\{\gamma\in \Gamma:\, \gamma=o(v\fv)\big\}$ of 
$\Gamma$.  
%$\delta:=v(\fv)$, and if 
%$\fv\prec 1$, then $\prec_\delta$, $\asymp_{\delta}$, $\sim_{\delta}$ 
%denote $\prec_\Delta$, $\asymp_\Delta$, $\sim_\Delta$ for 
%$\Delta=\Delta(\fv)=\big\{\gamma\in \Gamma:\, \gamma=o(\delta)\big\}$.

We fix a slot $(P,\fm,\hat a)$ in $K$ of order $r\geq 1$, and set $w:=\wt(P)$ (so $w\geq r\geq 1$).
In the next subsections we introduce various conditions on~$(P,\fm,\hat a)$. These conditions will be shown to be related as follows:
\[
\xymatrix{ \text{ strictly normal } \ar@{=>}[r] & \text{ normal } \ar@{=>}[r] \ar@{=>}[d] & \text{ steep } \\
& \text{ quasilinear } \ar@{<=}[r] &  \text{ deep }  \ar@{=>}[u]}
\] 
%$$\text{strongly normal}\ \Longrightarrow\ \text{normal}\ \Longrightarrow\ \text{steep}\ \Longleftarrow\ \text{deep}\  \Longrightarrow\ \text{quasilinear},$$
%$$\text{normal}\ \Longrightarrow\ \text{quasilinear}.$$
Thus ``deep + strictly normal'' yields the rest. The main results of this section are Theorem~\ref{mainthm} and its variants \ref{cor:mainthm}, \ref{varmainthm}, and \ref{cor:achieve strong normality, 2}.

\subsection*{Steep and deep slots}
In this subsection, if $\order (L_{P_{\times \fm}})=r$, then we set
$$\fv\ :=\ \fv(L_{P_{\times\fm}}).$$ 
The slot~$(P,\fm,\hat a)$ in $K$ is said to be {\bf steep}\index{steep!slot}\index{slot!steep} if $\order(L_{P_{\times \fm}})=r$ and $\fv\prec^\flat 1$. 
Thus 
$$(P,\fm, \hat a) \text{ is steep }\Longleftrightarrow\ (P_{\times\fn},\fm/\fn,\hat a/\fn) \text{ is steep }\Longleftrightarrow\ (bP,\fm,\hat a) \text{ is steep}$$
for $b\neq 0$.   If $(P,\fm,\hat a)$ is steep, then so is any slot in $K$ equivalent to $(P,\fm,\hat a)$. If $(P, \fm, \hat a)$ is steep, then so is any slot $(P^\phi,\fm,\hat a)$ in~$K^\phi$
for active~$\phi\preceq 1$, 
by  Lemma~\ref{lem:v(Aphi)}, and thus $\nwt(L_{P_{\times \fm}})<r$.  
Below we tacitly use that if~$(P,\fm,\hat a)$ is steep, then
$$\fn\asymp_{\Delta(\fv)}\fv\ \Longrightarrow\ [\fn]=[\fv], \qquad \fn\prec 1,\ [\fn]=[\fv]\ \Longrightarrow\ \fn \prec^\flat 1.$$
Note also that if $(P,\fm,\hat a)$ is steep, then $\fv^\dagger \asymp_{\Delta(\fv)} 1$ by [ADH, 9.2.10(iv)].

\begin{lemma}\label{lem:steep1}
Suppose $(P,\fm,\hat a)$ is steep, $\hat a\prec\fn\preceq \fm$ and $[\fn/\fm]\leq [\fv]$. Then 
$$\order(L_{P_{\times\fn}})\ =\ r, \qquad \fv(L_{P_{\times\fn}})\ \asymp_{\Delta(\fv)}\ \fv,$$ 
so $(P,\fn,\hat a)$ is a steep refinement of 
$(P,\fm,\hat a)$.
\end{lemma}
\begin{proof}
Replace $(P,\fm,\hat a)$ and $\fn$ by $(P_{\times\fm},1,\hat a/\fm)$ and $\fn/\fm$ to arrange
$\fm=1$. Set~$L:= L_P$ and $\tilde L:=L_{P_{\times\fn}}$.
Then $\tilde L=L\fn\asymp_{\Delta(\fv)} \fn L$ by [ADH, 6.1.3]. %using $[v(\fn)]\le [\delta]$. 
Hence 
%\begin{equation}\label{eq:normal 1}
$$\tilde L_r\ =\ \fn L_r\ \asymp\ \fn\fv L\ \asymp_{\Delta(\fv)}\ \fv\tilde L.$$
%\end{equation} 
Since $\fv(\tilde L)\tilde L \asymp \tilde L_r$, this gives
$\fv(\tilde L)\tilde L\asymp_{\Delta(\fv)} \fv \tilde L$, and thus
$\fv(\tilde L)\asymp_{\Delta(\fv)} \fv$. 
\end{proof}

%\begin{lemma}\label{lem:steep2}
%Suppose $(P,\fm,\hat a)$ is steep, $\hat a\prec\fn\preceq \fm$ and $\fn/\fm\steeper\fv$. Then for all sufficiently small $q\in\Q^>$, $(P,\fn^q\fm^{1-q},\hat a)$ is a steep refinement of $(P,\fm,\hat a)$ with $\fv(L_{P_{\times\fn^q\fm^{1-q}}})\steepereq \fv$.
%\end{lemma}
%\begin{proof}
%This follows as in the beginning of the proof of Lemma~\ref{lem:normal for small q}.  
%\end{proof}

\noindent
If $(P,\fm,\hat a)$ is steep and linear, then $L_{P_{+a,\times\fm}}=L_{P_{\times\fm,+(a/\fm)}}=L_{P_{\times\fm}}$, so any refinement $(P_{+a},\fm,\hat a-a)$ of $(P,\fm,\hat a)$  is also steep and linear.

\begin{lemma}\label{lem:achieve steep} 
Suppose  $\order L_{P_{\times\fm}}=r$. Then $(P,\fm,\hat a)$ has a refinement $(P,\fn,\hat a)$ such that 
$\nwt L_{P_{\times\fn}}=0$, and $(P^\phi,\fn,\hat a)$
is steep, eventually.
\end{lemma}
\begin{proof}
Replacing $(P,\fm,\hat a)$ by $(P_{\times\fm},1,\hat a/\fm)$ we  arrange $\fm=1$. 
Take $\fn_1$ with $\hat a\prec \fn_1 \prec 1$. Then $\order\, (P_1)_{\times\fn_1}=\order P_1=\order L_P=r$, and 
thus~${(P_1)_{\times\fn_1}\neq 0}$. So~[ADH, 11.3.6] applied to~$(P_1)_{\times\fn_1}$ in place of $P$
yields an $\fn$ with~$\fn_1\prec\fn\prec 1$ and $\nwt\, (P_1)_{\times\fn}=0$, so $\nwt L_{P_{\times\fn}}=0$.
Hence by Lemma~\ref{lem:eventual value of fv},  $(P^\phi,\fn,\hat a)$
is steep, eventually.
\end{proof}

\noindent
Recall that  the separant $S_P=\partial P/\partial Y^{(r)}$ of $P$ has lower complexity than $P$.
Below we sometimes use the identity $S_{P_{\times\fm}^\phi}=\phi^r (S_{P_{\times\fm}})^\phi$ from Lemma~\ref{lem:separant fms}.

\medskip\noindent
The slot $(P,\fm,\hat a)$ in $K$ is said to be {\bf deep} if it is steep and for all active $\phi\preceq 1$,\index{slot!deep}\index{deep}
{\samepage
\begin{itemize}
\item[(D1)] $\ddeg S_{P^\phi_{\times\fm}}=0$ (hence $\ndeg S_{P_{\times\fm}}=0$), and 
\item[(D2)] $\ddeg P^\phi_{\times\fm}=1$ (hence $\ndeg P_{\times \fm}=1$).
\end{itemize}}\noindent
If $\deg P=1$, then (D1) is automatic, for all active $\phi\preceq 1$. 
If~$(P,\fm,\hat a)$ is deep, then so are $(P_{\times\fn},\m/\fn,\hat a/\fn)$ and $(bP,\fm,\hat a)$ for $b\neq 0$, as well as every slot in $K$ equivalent to  $(P,\fm,\hat a)$ and the slot $(P^\phi,\fm,\hat a)$ in $K^\phi$ for active~${\phi\preceq 1}$. Every deep slot in $K$ is quasilinear, by (D2).
If  $\deg P=1$, then~$(P,\fm,\hat a)$ is quasilinear iff~$(P^\phi,\fm,\hat a)$ is deep for some active $\phi\preceq 1$. 
Moreover, if $(P,\fm,\hat a)$ is a deep hole in $K$, then~$\dval P^\phi_{\times \fm}=1$ for all active $\phi\preceq 1$, by
(D2) and Lemma~\ref{lem:lower bd on ddeg}.

\begin{exampleNumbered}\label{ex:order 1 linear steep} 
Suppose $P=Y'+gY-u$ where $g,u\in K$ and $\fm=1$, $r=1$. Set~$L:=L_P=\der+g$ and $\fv:=\fv(L)$. Then 
$\fv=1$ if $g\preceq 1$, and $\fv=1/g$ if $g\succ 1$. Thus 
$$\text{$(P,1,\hat a)$ is steep} \quad\Longleftrightarrow\quad g\succ^\flat 1  \quad\Longleftrightarrow\quad  
%\text{$(P,1,\hat a)$ is normal} \quad\Longleftrightarrow\quad 
\text{$g\succ 1$ and $g^\dagger \succeq 1$.}$$
Note that $(P,1,\hat a)$ is steep iff $L$ is steep as defined in Section~\ref{sec:lindiff}. 
Also, 
$$\text{$(P,1,\hat a)$ is deep} \quad\Longleftrightarrow\quad \text{$(P,1,\hat a)$ is steep and $g\succeq u$.}$$
%as well as
%$$\text{$(P,1,\hat a)$ is strongly normal} \quad\Longleftrightarrow\quad  \text{$(P,1,\hat a)$ is steep and $f\prec_\delta 1/g$.}$$
Hence if $u=0$, then $(P,1,\hat a)$ is deep iff it is steep. 
\end{exampleNumbered}

\begin{lemma}\label{lem:eventually deep}
For steep $(P,\fm,\hat a)$, the following are equivalent:
\begin{enumerate}
\item[$\mathrm{(i)}$]
 $(P^\phi,\fm,\hat a)$ is deep, eventually; 
\item[$\mathrm{(ii)}$] 
$\ndeg S_{P_{\times\fm}}=0$ and $\ndeg P_{\times\fm}=1$.
%\nval P_{\times\fm}=1$.
\end{enumerate}
\end{lemma}

\noindent
Note that if $\ddeg S_{P_{\times \fm}}=0$ or $\ndeg S_{P_{\times \fm}}=0$, then $S_{P_{\times\fm}}(0)\neq 0$, so $\order L_{P_{\times\fm}}=r$.

\begin{lemma}\label{lem:deep 1}
Suppose $(P_{+a},\fn,\hat a-a)$ refines the hole $(P,\fm,\hat a)$ in $K$. Then: 
%is a hole in $K$ and  refines $(P,\fm,\hat a)$. Then: 
\begin{enumerate}
\item[$\mathrm{(i)}$]  $\ddeg S_{P_{\times \fm}}=0\ \Longrightarrow\ \ddeg S_{P_{+a,\times\fn}}=0$;
\item[$\mathrm{(ii)}$] $\ddeg P_{\times \fm}=1\ \Longrightarrow\ \ddeg P_{+a,\times \fn}=1$;
\item[$\mathrm{(iii)}$] $\ndeg S_{P_{\times \fm}}=0\ \Longrightarrow\ S_P(a)\sim S_P(0)$.
\end{enumerate}
%\item[$\mathrm{(iv)}$] $(P,\fm, \hat a)$ is deep and $\nwt (P_{\times \fm})_1=0\ \Longrightarrow\ \nwt (P_{+a,\times\fn})_1= 0$. (??)
Thus if $(P,\fm, \hat a)$ is deep and $(P_{+a},\fn,\hat a-a)$ is steep, then $(P_{+a},\fn,\hat a-a)$ is deep. 
\end{lemma}
\begin{proof} 
Suppose $\ddeg S_{P_{\times \fm}}=0$. Then $\ddeg S_{P_{+a,\times \fn}}=0$ follows from $$\ddeg S_{P_{+a,\times\fn}}\ =\ \ddeg\,(S_P)_{+a,\times\fn}\ \text{ and }\
\ddeg (S_P)_{\times\fm}\ =\ \ddeg S_{P_{\times\fm}}$$
(consequences of Lemma~\ref{lem:separant fms}), and
$$\ddeg\,(S_P)_{+a,\times\fn}\ =\
\ddeg\,(S_P)_{+\hat a,\times\fn}\ \leq\ 
\ddeg\,(S_P)_{+\hat a,\times\fm}\ =\ \ddeg\,(S_P)_{\times\fm}$$
which holds by [ADH, 6.6.7]. This proves (i). Corollary~\ref{cor:ref 2} yields (ii), and (iii) is contained in Lemma~\ref{lem:ndeg coeff stabilizes, 2}. 
\end{proof}

\noindent
Lemmas~\ref{lem:from cracks to holes} and~\ref{lem:deep 1} give:

\begin{cor}\label{cor:steep refinement}
If  $(P,\fm,\hat a)$ is $Z$-minimal and deep, then each steep refinement of~$(P,\fm,\hat a)$ is deep.
\end{cor}
 
\noindent
Here is another sufficient condition on refinements of deep holes to remain deep:

\begin{lemma}\label{lem:deep 2}
Suppose $(P,\fm,\hat a)$ is a deep hole in $K$, and $(P_{+a},\fn,\hat a-a)$ refines~$(P,\fm,\hat a)$ with $[\fn/\fm]\le[\fv]$. Then $(P_{+a},\fn,\hat a-a)$  is deep with $\fv(L_{P_{+a,\times\fn}})\asymp_{\Delta(\fv)}\fv$.
\end{lemma}
\begin{proof}
From $(P,\fm,\hat a)$  we pass to the hole
$(P_{+a}, \fm, \hat a -a)$ and then to $(P_{+a},\fn,\hat a-a)$.  
We first show that $\order L_{P_{+a,\times \fm}}=r$ and 
$\fv(L_{P_{+a, \times \fm}})\sim \fv$, from which it follows that $(P_{+a},\fm,\hat a-a)$ is steep, hence deep by Lemma~\ref{lem:deep 1}.
By Corollary~\ref{cor:ref 2}, 
$$\ddeg P_{+a,\times \fm}\ =\ \dval P_{+a, \times \fm}\ =\ 1,$$ so
$(P_{+a,\times \fm})_1 \sim P_{+a,\times \fm}$. Also
$$(P_{\times \fm})_1\sim P_{\times \fm}\sim P_{\times \fm,+(a/\fm)}= P_{+a,\times \fm},$$  by [ADH, 4.5.1(i)], and
thus $(P_{+a,\times \fm})_1\sim (P_{\times\fm})_1$.
By Lemmas~\ref{lem:separant fms} and~\ref{lem:deep 1}(iii),
$$S_{P_{+a,\times \fm}}(0)\ =\ \fm S_P(a)\ \sim\ \fm S_P(0)\ =\ S_{P_{\times \fm}}(0),$$
so $S_{P_{+a,\times \fm}}(0)\sim S_{P_{\times \fm}}(0)$. This gives
$\fv(L_{P_{+a,\times\fm}})\sim \fv$ as promised.

Next, Lemma~\ref{lem:steep1} applied to $(P_{+a}, \fm, \hat a -a)$
in the role of $(P,\fm, \hat a)$ gives that~$(P_{+a}, \fn, \hat a -a)$ is steep with $\fv(L_{P_{+a,\times \fn}})\asymp_{\Delta(\fv)} \fv$.
Now  Lemma~\ref{lem:deep 1} applied to $(P_{+a}, \fm, {\hat a -a})$ and  $(P_{+a}, \fn, {\hat a -a})$ 
in the role of  $(P,\fm, \hat a)$ and  $(P_{+a}, \fn, {\hat a -a})$, respectively, gives that  $(P_{+a}, \fn, \hat a -a)$ is deep. 
\end{proof}

\noindent
Lemmas~\ref{lem:from cracks to holes} and~\ref{lem:deep 2} give a version for $Z$-minimal slots:

\begin{cor}\label{cor:deep 2, cracks}
Suppose $(P,\fm,\hat a)$ is $Z$-minimal and deep, and $(P_{+a},\fn,{\hat a-a})$ refines $(P,\fm,\hat a)$ with $[\fn/\fm]\leq[\fv]$,
where $\fv:=\fv(L_{P_{\times\fm}})$.
Then~$(P_{+a},\fn,{\hat a-a})$  is deep with $\fv(L_{P_{+a,\times\fn}})\asymp_{\Delta(\fv)}\fv$.
\end{cor}

\noindent
Next we turn to the task of turning $Z$-minimal slots into deep ones.

\begin{lemma}\label{prop:normalize, q-linear} 
Every quasilinear $Z$-minimal slot in $K$ of order $r$ has a refinement~$(P,\fm,\hat a)$ such that:
\begin{enumerate}
\item[$\mathrm{(i)}$]
$\ndeg\,(P_{(\i)})_{\times\fm}=0$  for all $\i$ with $\abs{\i}\geq 1$ and $P_{(\i)}\neq 0$;
\item[$\mathrm{(ii)}$] $\ndeg P_{\times\fm}=\nval P_{\times\fm}=1$, and
\item[$\mathrm{(iii)}$] $\nwt L_{P_{\times\fm}}=0$. %$\fv(L_{P^\phi_{\times\fm}})\prec^\flat_\phi 1$ eventually.
\end{enumerate}
\end{lemma}
\begin{proof}
By Corollary~\ref{cor:ref 1n}, any quasilinear $(P,\fm,\hat a)$  satisfies (ii). 
Any refinement  of a quasilinear $(P,\fm,\hat a)$ remains quasilinear, by Corollary~\ref{cor:ref 2n}. 
By Lemma~\ref{lem:ndeg coeff stabilizes, 1} and a subsequent remark any quasilinear $Z$-minimal slot in $K$ of order $r$ can be refined to a quasilinear $(P,\fm,\hat a)$ that satisfies (i), and by Lemma~\ref{lem:ndeg coeff stabilizes, 2}, any further refinement of such~$(P,\fm,\hat a)$ continues to satisfy (i). Thus to prove the lemma, assume we are given a
quasilinear $(P, \fm,\hat a)$ satisfying (i); it is enough to show that then $(P,\fm,\hat a)$ has a refinement
$(P,\fn,\hat a)$ satisfying (iii) with $\fn$ instead of $\fm$ (and thus also (i) and (ii) with $\fn$ instead of $\fm$).

 Take $\tilde\fm$ with $\hat a \prec\tilde\fm\prec\fm$. Then~$(P_{\times \tilde\fm})_1\ne 0$ by (ii), so 
[ADH, 11.3.6] applied to~$(P_1)_{\times\tilde\fm}$ in place of $P$ yields an~$\fn$ with $\tilde\fm\prec\fn\prec\fm$ and 
$\nwt\, (P_1)_{\times\fn}=0$. Hence the refinement~$(P,\fn,\hat a)$ of $(P,\fm,\hat a)$ 
satisfies (iii) with~$\fn$ instead of $\fm$.
\end{proof}

\begin{cor}\label{cor:deepening, q-linear}
Every quasilinear $Z$-minimal slot in $K$ of order $r$ has a re\-fine\-ment~$(P,\fm,\hat a)$ such that $\nwt L_{P_{\times\fm}}=0$, and~$(P^\phi,\fm,\hat a)$ is deep, eventually.
\end{cor}
\begin{proof} Given a quasilinear $Z$-minimal slot in $K$ of order $r$, we take a re\-fi\-ne\-ment~$(P,\fm,\hat a)$ as in Lemma~\ref{prop:normalize, q-linear}. Then $\ndeg S_{P_{\times \fm}}=0$ by (i) of that lemma, so $\order L_{P_{\times \fm}}=r$ by the remark that precedes Lemma~\ref{lem:deep 1}.
Then (iii) of Lem\-ma~\ref{prop:normalize, q-linear} and Lem\-ma~\ref{lem:eventual value of fv} give that
$(P^\phi, \fm,\hat a)$ is steep, eventually. Using now~$\ndeg S_{P_{\times \fm}}=0$ and~(ii) of Lemma~\ref{prop:normalize, q-linear} we obtain from Lemma~\ref{lem:eventually deep} that~$(P^\phi, \fm,\hat a)$ is deep, eventually.  
\end{proof}

\noindent
Lemma~\ref{lem:quasilinear refinement} and the previous lemma and its corollary now yield:

\begin{lemma}\label{prop:normalize} Suppose $K$ is $\d$-valued and $\upo$-free, and $\Gamma$ is divisible. Then
every $Z$-minimal slot in $K$ of order $r$ has a refinement $(P,\fm,\hat a)$ satisfying \textup{(i)}--\textup{(iii)} in Lemma~\ref{prop:normalize, q-linear}.
\end{lemma}

\begin{cor}\label{cor:deepening} Suppose $K$ is $\d$-valued and $\upo$-free, and $\Gamma$ is divisible. Then
every $Z$-minimal slot in $K$ of order $r$ has a quasilinear refinement $(P,\fm,\hat a)$ such that $\nwt L_{P_{\times\fm}}=0$, and~$(P^\phi,\fm,\hat a)$ is deep, eventually.
\end{cor}

\subsection*{Approximating $Z$-minimal slots} In this subsection we set, as before, 
$$\fv\ :=\ \fv(L_{P_{\times\fm}}), $$ provided $L_{P_{\times \fm}}$ has order $r$. The next lemma is a key approximation result.  

\begin{lemma}\label{lem:good approx to hata}
Suppose  $(P,\fm,\hat a)$ is $Z$-minimal and steep,  and
$$\ddeg P_{\times\fm}\ =\ \ndeg P_{\times\fm}\ =\ 1,\qquad \ddeg S_{P_{\times_\fm}}\ =\ 0.$$ 
Then there exists an $a$ such that $\hat a-a \prec_{\Delta(\fv)} \fm$.
\end{lemma}
\begin{proof}
We can arrange $\fm=1$ and $P\asymp 1$. Then $\ddeg P=1$
gives $P_1\asymp 1$, so~$S_P(0)\asymp\fv$.
Take $Q,R_1,\dots,R_n\in K\{Y\}$ ($n\geq 1$) of order~$<r$ such that $$ P\ =\ Q+R_1Y^{(r)}+\cdots+R_n(Y^{(r)})^n,\qquad S_P\ =\ R_1+\cdots+nR_n(Y^{(r)})^{n-1}.$$
Then $R_1(0)=S_P(0)\asymp\fv$. As $\ddeg S_P=0$, this gives $S_P \sim R_1(0)$,
hence $$R\ :=\ P-Q\ \sim\ R_1(0)Y^{(r)}\ \asymp\ \fv\  \prec_{\Delta(\fv)}\ 1\ \asymp\ P,$$ so $P\sim_{\Delta(\fv)} Q$. Thus $Q\ne 0$, and $Q\notin Z(K,\hat a)$ because $\order Q < r$.
% Lemma~\ref{lem:Z(K,hat a)}. 
Now Lem\-ma~\ref{lem:notin Z(K,hata)} gives a refinement~$(P_{+a},\fn,\hat a-a)$ of $(P, 1, \hat a)$ such that~${\ndeg Q_{+a,\times \fn}=0}$ and $\fn\prec 1$.
We claim that then $\hat a -a\prec_{\Delta(\fv)} 1$. (Establishing this claim finishes the proof.)
Suppose the claim is false. Then $\hat a-a\asymp_{\Delta(\fv)} 1$, so $\fn \asymp_{\Delta(\fv)} 1$,
hence~$Q_{+a,\times\fn} \asymp_{\Delta(\fv)} Q_{+a}\asymp Q$ by [ADH, 4.5.1]. Likewise,
$R_{+a,\times\fn} \asymp_{\Delta(\fv)} R$.
Using~$P_{+a,\times\fn} =  Q_{+a,\times\fn} +  R_{+a,\times\fn}$ gives
$Q_{+a,\times\fn} \sim_{\Delta(\fv)}  P_{+a,\times\fn}$, so
$Q_{+a,\times\fn} \sim^\flat  P_{+a,\times\fn}$. 
Then $\ndeg  Q_{+a,\times\fn}=\ndeg P_{+a,\times\fn}=1$ by Lemma~\ref{lem:same ndeg} and Corollary~\ref{cor:ref 2n}, a contradiction.
\end{proof}

%\marginpar{commented out a technical lemma here from an old version, which was never used later}

%\begin{lemma}\label{lem:deep minimal} \marginpar{check if really needed}
%Suppose $\Gamma$ is divisible and $(P,\fm,\hat a)$ is $Z$-minimal and deep. Then there exists a deep refinement   $(P_{+a},\fn,\hat a-a)$ of $(P,\fm, \hat a)$ such that $\hat a-a\prec_\delta \fn \prec_\delta \fm$, $\ddeg_\delta P_{+a,\times\fn}=\dval_\delta P_{+a,\times\fn}=1$, and $\order L_{P_{+a,\times \fn}}=r$, $\fv(L_{P_{+a,\times\fn}})\asymp_{\delta} \fv$.
%\end{lemma}
%\begin{proof} 
%Using Lemma~\ref{lem:from cracks to holes} we arrange that $(P,\fm,\hat a)$ is a hole in $K$. Lem\-ma~\ref{lem:good approx to hata} gives $a$ with  $\hat a-a \prec_\delta \fm$.  Then by Lemma~\ref{lem:deep 2} the re\-fine\-ment~$(P_{+a}, \fm, {\hat a-a})$ of $(P,\fm, \hat a)$ is deep; in particular, $\ddeg P_{+a,\times\fm}=1$.  Applying Lemma~\ref{keylem:coars} to~$(P_{+a},\fm,\hat a-a)$  in place of  $(P,\fm,\hat a)$ gives $\fn$ with $\hat a-a \prec_\delta \fn \prec_\delta \fm$,  $[\fn/\fm]=[\delta]$, and $\ddeg_\delta P_{+a,\times\fn}=\dval_\delta P_{+a,\times\fn}=1$. The rest follows by applying Lemma~\ref{lem:deep 2} to~$(P,\fm, \hat a)$ with its refinement $(P_{+a}, \fn, \hat a-a)$ of $(P,\fm, \hat a)$. 
%\end{proof}
\noindent
Lemmas~\ref{lem:lower bd on ddeg} and~\ref{lem:good approx to hata}, and a remark following the definition of {\em deep\/}  give:
 
\begin{cor}\label{cor:good approx to hata, deg 1} 
If $(P,\fm,\hat a)$ is $Z$-minimal, steep, and linear, then there exists an $a$ such that $\hat a-a \prec_{\Delta(\fv)} \fm$.
\end{cor}

\begin{cor}\label{specialvariant}
Suppose $(P,\fm,\hat a)$ is $Z$-minimal, deep, and special. Then for all~$n\ge 1$ there is an~$a$ with $\hat a-a\prec \fv^n\fm$.
\end{cor} 
\begin{proof} We arrange $\fm=1$ in the usual way. Let $\Delta$ be the convex subgroup of $\Gamma$ that is cofinal in $v(\hat a - K)$. Lemma~\ref{lem:good approx to hata} gives an element 
$\gamma\in v(\hat a -K)$
with $\gamma\ge \delta/m$ for some~$m\ge 1$. Hence
$v(\hat a -K)$ contains for every $n\ge 1$ an element $>n\delta$.
\end{proof}

\noindent
Combining Lemma~\ref{lem:special dents} with Corollary~\ref{specialvariant} yields:

\begin{cor}\label{cor:closer to minimal holes} 
If $K$ is   $r$-linearly newtonian, $\upo$-free if $r>1$, and
$(P,\fm,\hat a)$ is $Z$-minimal and deep, then
%(P,\fm,\hat a)$ is special and so
for all $n\ge 1$ there is an~$a$ such that~$\hat a-a\prec \fv^n\fm$.
\end{cor}
%\begin{proof} 
%Use Lemma~\ref{lem:from cracks to holes} to arrange $(P,\fm,\hat a)$ is a hole in~$K$. Next arrange $\fm=1$ as usual, so $\ndeg P=1$, hence $\hat a$ is special over $K$ by Proposition~\ref{nepropsp}.
%\end{proof} 

\subsection*{Normal slots} 
We say that our slot~$(P,\fm,\hat a)$ in $K$, with linear part $L$, is {\bf normal}\/ if  $\order L=r$ and, with $\fv:=\fv(L)$ and $w:=\wt(P)$,\index{slot!normal}\index{normal!slot} 
\begin{itemize}
\item[(N1)] $\fv\prec^\flat 1$;
\item[(N2)] $(P_{\times\fm})_{> 1}\prec_{\Delta(\fv)} \fv^{w+1} (P_{\times\fm})_1$.  
\end{itemize}
Note that then $\fv\prec 1$, $\dwt(L)<r$, $(P, \fm, \hat a)$ is steep, and
\begin{equation}\label{eq:N3}
P_{\times\fm}\sim_{\Delta(\fv)} P(0)+(P_{\times\fm})_1\qquad\text{ (so $\ddeg P_{\times\fm} \le 1$).}
\end{equation}
If $\order L=r$, $\fv:=\fv(L)$, and $L$ is monic, then $(P_{\times\fm})_1\asymp\fv^{-1}$, so that~(N2) is then equivalent to:
$(P_{\times\fm})_{> 1}\prec_{\Delta(\fv)} \fv^{w}$. 
If $\deg P=1$, then $\order L=r$ and~(N2) automatically holds,  
hence $(P,\fm,\hat a)$ is normal iff it is steep.  
Thus by Lemma~\ref{lem:eventual value of fv}:

\begin{lemma}\label{lem:deg1 normal} 
If $\deg P=1$ and $\nwt(L)<r$, then $(P^\phi,\fm,\hat a)$ is normal, eventually.
\end{lemma}

\noindent
If $(P,\fm,\hat a)$ is normal, then so are
$(P_{\times\fn},\fm/\fn,\hat a/\fn)$ and
$(bP,\fm,\hat a)$ for  $b\neq 0$.
In particular, $(P,\fm,\hat a)$ is   normal iff $(P_{\times\fm},1,\hat a/\fm)$ is normal.  
If $(P,\fm,\hat a)$ is normal, then so is any equivalent slot. 
Hence by \eqref{eq:N3} and Lemmas~\ref{lem:lower bd on ddeg} and~\ref{lem:from cracks to holes}:

\begin{lemma}\label{lem:ddeg=dmul=1 normal} 
If $(P,\fm,\hat a)$ is normal, and $(P,\fm,\hat a)$ is $Z$-minimal or is a hole in~$K$, then $\ddeg P_{\times\fm} = \dval P_{\times\fm} = 1$.
\end{lemma}

\begin{example} 
Let $K\supseteq\R(\ex^x)$ be an $H$-subfield of $\T$, $\fm=1$, $r=2$. If $P=D+R$ where 
$$D\ =\ \ex^{-x}Y''-Y,\qquad R\ =\ f+\ex^{-4x}Y^5\quad (f\in K),$$
then $\fv=-\ex^{-x}\prec^\flat 1$, $P_1=D\sim -Y$, $w=2$, and
$P_{>1}=\ex^{-4x}Y^5 \prec_{\Delta(\fv)} \ex^{-3x}P_1$, so~$(P,1,\hat a)$ is normal.
However, if $P=D+S$ with $D$ as above and $S=f+\ex^{-3x}Y^5$ ($f\in K$), then
$P_{>1}=\ex^{-3x}Y^5 \succeq_{\Delta(\fv)} \ex^{-3x}P_1$, so $(P,1,\hat a)$ is not normal.
\end{example}

\begin{lemma}\label{lem:normal pos criterion}
Suppose $\order(L)=r$ and $\fv$ is such that \textup{(N1)} and~\textup{(N2)} hold, and $\fv(L)\asymp_{\Delta(\fv)} \fv$. Then 
$(P,\fm,\hat a)$ is   normal.
\end{lemma}
\begin{proof}
Put $\fw:=\fv(L)$. Then $[\fw]=[\fv]$, and so $\fv\prec^{\flat} 1$ gives $\fw\prec^{\flat} 1$. Also, 
$$(P_{\times\fm})_{> 1}\prec_{\Delta(\fv)} \fv^{w+1} (P_{\times\fm})_1\asymp_{\Delta(\fv)}\fw^{w+1}(P_{\times\fm})_1.$$ 
Hence (N1), (N2) hold with~$\fw$ in place of $\fv$. 
\end{proof}

\begin{lemma}\label{lem:normality comp conj} 
Suppose $(P,\fm,\hat a)$ is normal and $\phi\preceq 1$ is active. Then the slot~$(P^\phi,\fm,\hat a)$  in $K^\phi$ is normal.
\end{lemma}
\begin{proof}
We arrange $\fm=1$ and put $\fv:=\fv(L)$, $\fw:=\fv(L_{P^\phi})$.
Now $L_{P^\phi}=L^\phi$, so~$\fv\asymp_{\Delta(\fv)} \fw$ and $\fv\prec_{\phi}^\flat 1$ by Lemma~\ref{lem:v(Aphi)}. By [ADH, 11.1.1], $[\phi]<[\fv]$, and (N2) we have
$$(P^\phi)_{>1}\ =\ (P_{>1})^\phi\ \asymp_{\Delta(\fv)}\ P_{>1}\ \prec_{\Delta(\fv)}\ \fv^{w+1} P_1\ \asymp_{\Delta(\fv)}\ \fv^{w+1} P^\phi_1,$$
which by Lemma~\ref{lem:normal pos criterion} applied to 
$(P^\phi, 1, \hat a)$ in the role of $(P, \fm, \hat a)$
gives that $(P^\phi, 1, \hat a)$ is normal.  
\end{proof}

\begin{cor}\label{cor:normal=>quasilinear}  %\label{cor:normalimpliesquasilinear} 
Suppose $(P,\fm,\hat a)$ is normal.  Then $(P,\fm,\hat a)$ is quasilinear.
\end{cor}
\begin{proof} Lemma~\ref{lem:lower bd on ndeg} gives $\ndeg P_{\times \fm}\ge 1$. The parenthetical remark after~\eqref{eq:N3} above and
Lemma~\ref{lem:normality comp conj} gives $\ndeg P_{\times \fm}\le 1$.  
\end{proof}

\noindent
Combining Lemmas~\ref{lem:ddeg=dmul=1 normal} and \ref{lem:normality comp conj} yields:

\begin{cor} 
If $(P,\fm,\hat a)$ is normal and  linear, and $(P,\fm,\hat a)$ is $Z$-minimal or a hole in~$K$, then
%~$(P,\fm,\hat a)$ is quasilinear,
%and if $(P,\fm,\hat a)$ is linear, then 
$(P,\fm,\hat a)$ is deep.
\end{cor}

\noindent
There are a few occasions later where we need to change the ``monomial'' $\fm$ in $(P,\fm,\hat a)$ while preserving key properties
of this slot. Here is what we need:  

\begin{lemma}\label{ufm} Let $u\in K$, $u\asymp 1$. Then $(P,u\fm,\hat a)$ refines $(P,\fm,\hat a)$, and
if ${(P_{+a},\fn,\hat a -a)}$ refines $(P,\fm,\hat a)$, then so does~$(P_{+a},u\fn, \hat a -a)$. If $(P,\fm,\hat a)$ is quasilinear, respectively deep, respectively normal, then so is~$(P,u\fm,\hat a)$. 
\end{lemma}
\begin{proof} 
The refinement claims are clearly true, and quasilinearity is preserved since $\ndeg P_{\times u\fm}=\ndeg P_{\times \fm}$ by [ADH, 11.2.3(iii)]. 
``Steep'' is preserved by Lemma~\ref{lem:steep1}, and hence ``deep''
is preserved using Lemma~\ref{lem:separant fms}  and [ADH, 6.6.5(ii)]. Normality is preserved because steepness is, 
$$(P_{\times u\fm})_d\ =\ (P_d)_{\times u\fm}\ \asymp\ (P_d)_{\times \fm}\ =\ (P_{\times \fm})_d\quad\text{ for all~$d\in \N$}$$
by [ADH, 4.3, 4.5.1(ii)], and $\fv(L_{P_{\times u\fm}})\asymp \fv(L_{P_{\times\fm}})$ by Lemma~\ref{fvmult}. 
\end{proof} 

\noindent
Here is a useful invariance property of normal slots:

\begin{lemma}\label{lem:excev normal}
Suppose $(P,\fm,\hat a)$ is normal and $a\prec\fm$.  Then $L_P$ and $L_{P_{+a}}$ have or\-der~$r$. 
If  in addition   $K$ is $\upl$-free or $r=1$, then $\exc^{\ev}(L_{P})=\exc^{\ev}(L_{P_{+a}})$. 
\end{lemma}
\begin{proof} $L_{P_{\times \fm}}=L_P\fm$ (so $L_P$ has order $r$), and 
$L_{P_{+a,\times \fm}}=L_{P_{\times \fm, +a/\fm}}= L_{P_{+a}}\fm$. 
The slot $(P_{\times\fm},1,\hat a/\fm)$ in $K$ is normal and $a/\fm\prec 1$.  Thus we can apply Lemma~\ref{lem:linear part, new}(i)  to  $\hat K$, $P_{\times\fm}$, $a/\fm$ in place of $K$, $P$, $a$ to give $\order L_{P_{+a}} = r$. Next, applying
likewise Lemma~\ref{lem:linear part, split-normal, new}  with $L:=L_{P_{\times \fm}}$, $\fv:=\fv(L_{P_{\times \fm}})$, $m=r$, $B=0$,   gives 
$$L_P\fm - L_{P_{+a}}\fm\ =\ L_{P_{\times\fm}} - L_{P_{\times\fm,+a/\fm}}\ \prec_{\Delta(\fv)}\ \fv^{r+1}L_{P}\fm.$$
%so $L_P-L_{P_{+a}} \prec_{\Delta(\fv)}\ \fv^{r+1}L_{P}$.
Hence, if $K$ is $\upl$-free, then $\exc^{\ev}(L_P\fm)\ =\ \exc^{\ev}(L_{P_{+a}}\fm)$ by Lemma~\ref{cor:excev stability}, so
$$ \exc^{\ev}(L_{P})\ =\ \exc^{\ev}(L_P\fm)+v(\fm)\ =\ \exc^{\ev}(L_{P_{+a}}\fm)+v(\fm)\ =\ \exc^{\ev}(L_{P_{+a}}).$$
%$\exc^{\ev}(L_{P})=\exc^{\ev}(L_{P_{+a}})$, and
If $r=1$ we obtain the same equality from Corollary~\ref{cor:excev stability, r=1}.
\end{proof}

\subsection*{Normality under refinements} 
In this subsection we study how normality behaves under more general refinements. This is not needed to prove the main result of this section, Theorem~\ref{mainthm}, but is included to
obtain useful variants of it.

\begin{prop}\label{normalrefine}
Suppose  $(P,\fm,\hat a)$ is normal. Let a refinement
$(P_{+a},\fm,\hat a-a)$ of $(P,\fm,\hat a)$ be given. Then this refinement is also normal. 
\end{prop}
\begin{proof}
By the remarks following the definition of ``multiplicative conjugate'' in Sec\-tion~\ref{sec:holes} and after replacing the slots 
$(P,\fm,\hat a)$ and $(P_{+a},\fm,\hat a-a)$
in~$K$
by $(P_{\times\fm},1,\hat a/\fm)$ and $\big(P_{\times\fm,+a/\fm},1,(\hat a-a)/\fm\big)$,
respectively, we arrange that $\fm=1$.  Let $\fv:=\fv(L_P)$.
By Lemma~\ref{lem:linear part, new} we have $\order(L_{P_{+a}})=r$, $\fv(L_{P_{+a}})\sim_{\Delta(\fv)} \fv$, and~$(P_{+a})_1\sim_{\Delta(\fv)} P_1$.
Using [ADH, 4.5.1(i)] we have for $d>1$ with~$P_d\ne 0$,
$$
(P_{+a})_{d}\ =\ \big((P_{\geq d})_{+a}\big){}_d\ \preceq\ (P_{\geq d})_{+a}\ \sim\  P_{\geq d}\ \preceq\ P_{>1},$$
and using (N2), this yields
$$ (P_{+a})_{>1}\ \preceq\ P_{>1}\ \prec_{\Delta(\fv)}\ \fv^{w+1}P_1\ \asymp\ \fv^{w+1}(P_{+a})_1.$$ 
Hence (N2) holds with $\fm=1$ and with $P$ replaced by $P_{+a}$. Thus $(P_{+a},1,\hat a-a)$
is normal,  by  Lemma~\ref{lem:normal pos criterion}.
\end{proof}

\begin{prop}\label{easymultnormal}
Suppose   $(P,\fm,\hat a)$ is a normal hole in $K$,   $\hat{a}\prec \fn\preceq \fm$, and~$[\fn/\fm]\le\big[\fv(L_{P_{\times\fm}})\big]$. Then the refinement $(P,\fn,\hat a)$ of $(P,\fm,\hat a)$
is also  normal.
\end{prop}
\begin{proof} As in the proof of Lemma~\ref{lem:steep1}
we arrange $\fm=1$ and  
set $L:= L_P$, $\fv:=\fv(L)$, and $\tilde L:=L_{P_{\times\fn}}$, to obtain $[\fn]\le[\fv]$ and
$\fv(\tilde L)\asymp_{\Delta(\fv)} \fv$. 
Recall from [ADH, 4.3] that~$(P_{\times\fn})_{d}=(P_{d})_{\times\fn}$ for $d\in \N$.
For such $d$ we have by [ADH, 6.1.3],
$$
(P_{d})_{\times\fn}\ \asymp_{\Delta(\fv)}\ \fn^d P_d\ \preceq\ \fn^d P_{\geq d}.$$
In particular, $(P_{\times \fn})_1\asymp_{\Delta(\fv)} \fn P_1$.
By (N2) we also have, for $d>1$:
$$P_{\geq d}\ \preceq\ P_{>1}\ \prec_{\Delta(\fv)}\ \fv^{w+1} P_1.$$
By Lemma~\ref{lem:ddeg=dmul=1 normal} we have   $P \sim P_1$. For $d>1$ we have by [ADH, 6.1.3],
$$\fn^d P\ \asymp\ \fn^d P_1\ \asymp_{\Delta(\fv)}\ \fn^{d-1} (P_1)_{\times\fn}\ \preceq\ 
(P_1)_{\times\fn}\ =\ (P_{\times\fn})_1\ \preceq\ P_{\times\fn}$$
and thus 
$$(P_{\times\fn})_{d}\ =\ (P_d)_{\times\fn}\ \preceq_{\Delta(\fv)}\ \fn^d P_{\geq d}\ \prec_{\Delta(\fv)}\  \fv^{w+1} \fn^d P_1\ \preceq_{\Delta(\fv)}\ \fv^{w+1} (P_{\times\fn})_1.$$
Hence (N2) holds with $\fm$ replaced by $\fn$. Thus $(P,\fn,\hat a)$
is normal, using $\fv(\tilde L)\asymp_{\Delta(\fv)} \fv$ and Lem\-mas~\ref{lem:steep1} and~\ref{lem:normal pos criterion}.
\end{proof}

\noindent
From Lemma~\ref{lem:from cracks to holes} and Proposition~\ref{easymultnormal} we obtain:

\begin{cor}\label{corcorcor} Suppose $(P,\fm,\hat a)$ is normal and $Z$-minimal, $\hat a \prec \fn \preceq \fm$, and~$[\fn/\fm]\le \big[\fv(L_{P_{\times \fm}})\big]$. Then the refinement $(P,\fn, \hat a)$ of $(P,\fm, \hat a)$ is also normal.
\end{cor}

\noindent
{\it In the rest of this subsection $\fm=1$, $\hat a \prec \fn\prec 1$, $\order(L_P)=r$, and $[\fv]<[\fn]$ where~$\fv:=\fv(L_P)$.}\/ So $(P, \fn,\hat a)$ refines $(P,1,\hat a)$, $L_{P_{\times \fn}}=L_P\fn$, and
$\order L_{P_{\times \fn}}=r$. 

\begin{lemma}\label{lem:normal from steep}
Suppose $(P,1,\hat a)$ is  steep, $\fv(L_{P_{\times\fn}})\preceq\fv$, and $P_{>1}\preceq P_1$.
Then $(P,\fn,\hat a)$ is normal.
\end{lemma}
\begin{proof}
Put $\fw:=\fv(L_{P_{\times\fn}})$. Then~$[\fw]<[\fn]$ by Corollary~\ref{cor:An}, and $\fw\preceq \fv \prec^{\flat} 1$ gives~$\fw\prec^\flat 1$.  It remains to show that
$(P_{\times\fn})_{>1} \prec_{\Delta(\fw)} \fw^{w+1} (P_{\times\fn})_1$.
Using~$[\fn]>[\fw]$ it is enough that
$(P_{\times\fn})_{>1} \prec_\Delta \fw^{w+1} (P_{\times\fn})_1$, 
where $\Delta:=\Delta(\fn)$.  Since
$\fw\asymp_\Delta 1$, it is even enough that $(P_{\times\fn})_{>1} \prec_{\Delta}  (P_{\times\fn})_1$, to
be derived below. 
Let~$d>1$. Then by~[ADH, 6.1.3] and $P_d\preceq P_{>1}\preceq P_1$ we have
$$(P_{\times\fn})_d\ =\ (P_d)_{\times\fn}\
\asymp_\Delta\  P_d \fn^d\ \preceq\ P_1\,\fn^d.$$
 In view of
$\fn\prec_{\Delta} 1$ and $d>1$ we have
$$ P_1\,\fn^d\ \prec_\Delta\ P_1\, \fn\ \asymp_\Delta\ (P_1)_{\times\fn}\ =\ (P_{\times\fn})_1,$$
using again [ADH, 6.1.3]. Thus $(P_{\times\fn})_d \prec_\Delta (P_{\times\fn})_1$, as promised.
\end{proof}

%\noindent
%Note that
%$$  \dval P=1\quad \Longrightarrow\quad \text{$P_d\preceq P_1$ for all $d>1$}\quad \Longrightarrow\quad \dval P\leq 1.$$
%Hence by Lemmas~\ref{lem:ddeg=dmul=1 normal} and~\ref{lem:normal from steep}:

\begin{cor}\label{cor:normal for small q, prepZ}
If $(P, 1, \hat a)$ is normal and~${\fv(L_{P_{\times\fn}})\preceq\fv}$, then $(P,\fn,\hat a)$ is normal.
\end{cor}

\noindent
In the next lemma and its corollary $K$ is
$\d$-valued and for every $q\in \Q^{>}$ there is given an element $\fn^q$ of $K^\times$ such that $(\fn^q)^\dagger=q\fn^\dagger$;
the remark before Lemma~\ref{qlA} gives $v(\fn^q)=qv(\fn)$ for $q\in \Q^{>}$. 
Hence for $0 < q\leq 1$ in $\Q$ we have~$\hat a\prec\fn\preceq\fn^q\prec 1$, so $(P,\fn^q,\hat a)$ refines~$(P,1,\hat a)$. 

\begin{lemma}\label{lem:normal for small q}
Suppose $(P,1,\hat a)$ is  steep   and $P_{>1}\preceq P_1$.  Then $(P,\fn^q,\hat a)$ is normal,  for all but finitely many $q\in \Q$ with $0<q \le 1$. 
\end{lemma}
\begin{proof} We have $\fn^\dagger\succeq 1$ by $\fn\prec \fv\prec 1$ and $\fv^\dagger\succeq 1$. Hence 
Lemma~\ref{lem:nepsilon} gives 
$\fv(L_{P_{\times\fn^q}})\preceq \fv$ for all but finitely many $q\in \Q^{>}$. Suppose $\fv(L_{P_{\times\fn^q}})\preceq \fv$, $0<q\le 1$ in $\Q$. Then~$(P,\fn^q,\hat a)$ is normal by Lemma~\ref{lem:normal from steep} applied with $\fn^q$ instead of $\fn$.
\end{proof}

\begin{cor}\label{cor:normal for small q}
If $(P, 1, \hat a)$ is normal, 
%and $Z$-minimal or a hole in $K$, 
then $(P,\fn^q,\hat a)$ is  normal for all but finitely many $q\in \Q$ with~$0<q\le 1$.
\end{cor}

%\begin{cor}
%Suppose  $(P,\fm,\hat a)$ is normal, and we have
%$\hat a \prec \fn\preceq \fm$, with
%$\fn/\fm \steeper \fv:= \fv(L_{P_{\times \fm}})$.
%Then $(P,\fn^q\fm^{1-q},\hat a)$  is a normal refinement
%of $(P,\fm,\hat a)$ for all sufficiently small
%$q\in\Q^>$.
%\end{cor} 

\subsection*{Normalizing} If in this subsection
$\order(L_{P_{\times \fm}})=r$, then $\fv:=\fv(L_{P_{\times\fm}})$. 
Towards proving that normality can always be achieved we first show:

\begin{lemma}\label{lem:normalizing}
Suppose $\Gamma$ is divisible, $(P,\fm,\hat a)$ is a deep hole in $K$, and $\hat a-a \prec \fv^{w+2}\fm$ for some $a$. 
Then $(P,\fm,\hat a)$  has a refinement that is deep and normal.
\end{lemma}
\begin{proof}
Replacing $(P,\fm,\hat a)$ by $(P_{\times \fm}, 1, \hat a/\fm)$ and renaming we arrange $\fm=1$. Take~$a$ such that
$\hat a -a\prec \fv^{w+2}$.   
For $e:=w+\frac{3}{2}$, let $\fv^e$ be an element of~$K^\times$ with~$v(\fv^e)=e\,v(\fv)$. {\em Claim}: the refinement
$(P_{+a},\fv^e,\hat a-a)$ of~$(P,1,\hat a)$  is deep and normal.
By Lemma~\ref{lem:deep 2}, $(P_{+a},\fv^e,\hat a-a)$  is deep, so
we do have $\order(L_{P_{+a,\times\fv^e}})=r$ 
and $\fv(L_{P_{+a,\times\fv^e}})\prec^\flat 1$.
Lemma~\ref{lem:deep 2} also yields $\fv(L_{P_{+a,\times\fv^e}}) \asymp_{\Delta(\fv)} \fv$. Since~$\ddeg P=\dval P=1$,
we can use Corollary~\ref{cor:ref 2} for $\fn=\fv^e$ and for $\fn=1$ to obtain 
$$\ddeg P_{+a,\times\fv^e}\ =\ \dval P_{+a,\times\fv^e}\ =\ \ddeg P_{+a}\ =\ \dval P_{+a}=1$$
and thus 
$(P_{+a,\times\fv^e})_1 \sim P_{+a,\times\fv^e}$; also
$P_1\sim P\sim P_{+a}\sim (P_{+a})_1$,  where $P\sim P_{+a}$ follows from
$a\prec 1$ and [ADH, 4.5.1(i)]. 
Now let $d>1$. Then 
\begin{align*}
(P_{+a,\times\fv^e})_d	&\ \asymp_{\Delta(\fv)}\ (\fv^{e})^d (P_{+a})_d 
						 \ \preceq\ (\fv^{e})^d P_{+a}
						 \ \sim\ (\fv^{e})^d(P_{+a})_1 \\
						&\ \asymp_{\Delta(\fv)}\ (\fv^{e})^{d-1} (P_{+a,\times\fv^e})_1\ \prec_{\Delta(\fv)} \fv^{w+1} (P_{+a,\times\fv^e})_1,
\end{align*}
using [ADH, 6.1.3] for $\asymp_{\Delta(\fv)}$. So $(P_{+a},\fv^e,\hat a-a)$ is normal by Lemma~\ref{lem:normal pos criterion}.
\end{proof}

\noindent
We can now finally show:

\begin{theorem}\label{mainthm} 
Suppose $K$ is $\upo$-free and $r$-linearly newtonian, and $\Gamma$ is divisible. Then every $Z$-minimal slot in $K$
of order $r$ has a refinement $(P,\fm,\hat a)$ such that~$(P^\phi,\fm,\hat a)$ is deep and normal, eventually. 
\end{theorem} 
\begin{proof}
By Lemma~\ref{lem:from cracks to holes} it is enough to show this for $Z$-minimal holes in $K$ of or\-der~$r$.
Given such   hole in $K$, use Corollary~\ref{cor:deepening} to refine it to a hole $(P,\fm,\hat a)$ such that $(P^\phi,\fm,\hat a)$ is deep, eventually.  
%Thus if $\deg P=1$, then $(P^\phi,\fm,\hat a)$ is normal, eventually, by Lemma~\ref{lem:deg1 normal}. In the rest of the proof we assume that $\deg P>1$.  
Replacing $(P,\fm,\hat a)$ by~$(P^\phi,\fm,\hat a)$ for a suitable active $\phi\preceq 1$ we arrange that
$(P,\fm,\hat a)$ itself is deep.
% Using Lemma~\ref{lem:deep minimal} to pass to a refinement of
% $(P,\fm,\hat a)$ we arrange in addition that  $\hat a\prec_\fv \fm$ and $\ddeg_\fv P_{\times\fm}=\dval_\fv P_{\times\fm}=1$. 
Then an appeal to Corollary~\ref{cor:closer to minimal holes} followed by an application of Lemma~\ref{lem:normalizing} yields a deep and normal refinement of  $(P,\fm,\hat a)$. Now apply Lemma~\ref{lem:normality comp conj} to this refinement. 
\end{proof}

\noindent
Next we indicate some variants of Theorem~\ref{mainthm}:

\begin{cor}\label{cor:mainthm} Suppose $K$ is $\d$-valued
and $\upo$-free, and $\Gamma$ is divisible. Then every minimal hole in $K$
of order $r$ has a refinement $(P,\fm,\hat a)$ such that $(P^\phi,\fm,\hat a)$ is deep and normal, eventually.
\end{cor}
\begin{proof} Given a minimal hole in $K$ of order $r$, use Corollary~\ref{cor:deepening} to refine it to a hole~$(P,\fm,\hat a)$ in $K$ such that $\nwt L_{P_{\times \fm}}=0$ and $(P^\phi,\fm,\hat a)$ is deep, eventually.
If~${\deg P=1}$, then $(P^\phi,\fm,\hat a)$ is normal, eventually, by Lemma~\ref{lem:deg1 normal}. 
If~${\deg P>1}$, then~$K$ is $r$-linearly newtonian by
Corollary~\ref{degmorethanone}, so we can use Theorem~\ref{mainthm}.
\end{proof}

\noindent
For $r=1$ we can follow the proof of Theorem~\ref{mainthm},  using Corollary~\ref{cor:deepening, q-linear} in place of Corollary~\ref{cor:deepening}, to obtain:

\begin{cor}  \label{mainthm, r=1} 
If $K$ is  $1$-linearly newtonian and $\Gamma$ is divisible, then every quasilinear $Z$-minimal slot in $K$
of order $1$ has a refinement $(P,\fm,\hat a)$ such that~$(P^\phi,\fm,\hat a)$ is deep and normal, eventually. 
\end{cor}

\noindent
Here is another variant of Theorem~\ref{mainthm}:

\begin{prop}\label{varmainthm} If $K$ is $\d$-valued and $\upo$-free, and
$\Gamma$ is divisible, then every $Z$-minimal special slot
in $K$ of order $r$ has a refinement $(P,\fm,\hat a)$ such that $(P^\phi,\fm,\hat a)$ is deep and normal, eventually. 
\end{prop}

\noindent
To establish this proposition we follow the proof of Theorem~\ref{mainthm}, using Lem\-ma~\ref{speciallemma} to preserve specialness in the initial refining. Corollary~\ref{specialvariant} takes over the role of Corollary~\ref{cor:closer to minimal holes}
in that proof.   

For linear slots in $K$ we can   weaken the hypotheses of Theorem~\ref{mainthm}:

\begin{cor}\label{mainthm deg 1}
Suppose   $\deg P=1$.  Then $(P,\fm,\hat a)$ has a refinement~$(P,\fn,\hat a)$ such that $(P^\phi,\fn,\hat a)$ is deep and normal, eventually.
Moreover, if $K$ is $\upl$-free and~${r>1}$, then $(P^\phi,\fm,\hat a)$ is deep and normal, eventually.
\end{cor}
\begin{proof}
By the remarks before Lemma~\ref{lem:deg1 normal}, $(P,\fm,\hat a)$ is normal iff it is steep. Moreover,
if $(P,\fm,\hat a)$ is normal, then it is quasilinear by Corollary~\ref{cor:normal=>quasilinear}, and hence~$(P^\phi,\fm,\hat a)$ is deep and normal, eventually,
by  the remarks before Example~\ref{ex:order 1 linear steep} and Lemma~\ref{lem:normality comp conj}.
By Lemma~\ref{lem:achieve steep}, $(P,\fm,\hat a)$ has a refinement~$(P,\fn,\hat a)$ such that $(P^\phi,\fn,\hat a)$ is steep, eventually.
This yields the first part. The second part follows from Corollary~\ref{coruplnwteq} and Lemma~\ref{lem:deg1 normal}.
\end{proof}

\begin{cor}\label{mainthm order 1}
Suppose $K$ is $\upl$-free, $\Gamma$ is divisible, and $(P, \fm, \hat a)$ is a quasilinear minimal hole in $K$ of order $r=1$.
Then $(P, \fm, \hat a)$ has a refinement~$(Q,\fn,\hat b)$ such that $(Q^\phi,\fn,\hat b)$ is deep and normal, eventually.
\end{cor}
\begin{proof} The case $\deg P=1$ is part of Corollary~\ref{mainthm deg 1}.
If $\deg P>1$, 
then $K$ is $1$-linearly newtonian by Lemma~\ref{lem:no hole of order <=r, deg 1}, so we can use 
Corollary~\ref{mainthm, r=1}.
\end{proof}

\subsection*{Improving normality} 
%We work under the same assumptions and employing the same notations as in the subsection ``Normal holes'' above.
%Thus $L=L_{P_{\times\fm}}$, $\fv=\fv(L)$, $\delta=v(\fv)$, and $w=\wt(P)$.
{\it In this subsection $L:=L_{P_{\times\fm}}$.}\/ 
Note that if $(P,\fm,\hat a)$ is a normal hole in $K$, then $P_{\times\fm}\sim (P_{\times \fm})_1$ by Lemma~\ref{lem:ddeg=dmul=1 normal}.
We call our slot~$(P,\fm, \hat a)$ in $K$ {\bf strictly normal\/} if it is normal, but
with the condition~(N2) replaced by the stronger condition\index{slot!strictly normal}\index{strictly!normal}\index{normal!strictly} 
\begin{enumerate}
\item[(N2s)] $(P_{\times\fm})_{\ne 1}\prec_{\Delta(\fv)} \fv^{w+1} (P_{\times \fm})_1$.
\end{enumerate}
Thus for normal $(P,\fm, \hat a)$ and $\fv=\fv(L)$ we have:
$$(P,\fm,\hat a) \text{ is strictly normal }\Longleftrightarrow\  P(0) \prec_{\Delta(\fv)} \fv^{w+1} (P_{\times \fm})_1.$$ 
So if $(P,\fm, \hat a)$ is  normal and $P(0)=0$, then $(P,\fm,\hat a)$ is strictly normal.
%If $w=1$ and $L$ is monic of order $1$ then  $(P_{\times\fm})_1 \asymp \fv^{-1}$ and 
%$$\text{$(P,\fm, \hat a)$ is strongly normal}\ \Longleftrightarrow\  \text{$(P,\fm, \hat a)$ is  normal and $P(0) \prec_\delta 1$.}$$
Note that if  $(P,\fm, \hat a)$ is strictly normal, then 
$$P_{\times \fm}\sim_{\Delta(\fv)} (P_{\times \fm})_1\qquad
\text{(and hence $\ddeg P_{\times\fm}=1$).}$$
If $(P,\fm,\hat a)$ is strictly normal, then so are
$(P_{\times\fn},\fm/\fn,\hat a/\fn)$ and
$(bP,\fm,\hat a)$ for  $b\neq 0$.
Thus $(P,\fm,\hat a)$ is strictly normal iff $(P_{\times\fm},1,\hat a/\fm)$ is strictly normal. 
If $(P,\fm,\hat a)$ is strictly normal, then so is every equivalent slot in~$K$.
The proof of Lemma~\ref{ufm} shows that if $(P,\fm,\hat a)$ is strictly normal and~${u\in K}$, $u\asymp 1$, then $(P,u\fm,\hat a)$ is also strictly normal. 
The analogue of
Lem\-ma~\ref{lem:normal pos criterion} goes through, with
$(P_{\times \fm})_{ \ne 1}$ instead of $(P_{\times \fm})_{>1}$
in the proof:

\begin{lemma}\label{lem:strongnormal pos criterion}
Suppose $\order(L)=r$ and $\fv$ are such that \textup{(N1)} and~\textup{(N2s)} hold, and $\fv(L)\asymp_{\Delta(\fv)} \fv$. Then 
$(P,\fm,\hat a)$ is strictly  normal.
\end{lemma}

\noindent
Lemma~\ref{lem:normality comp conj} goes likewise through with ``strictly normal'' instead of ``normal'':

\begin{lemma}\label{lem:normality comp conj, strong}
If $(P,\fm,\hat a)$ is strictly normal and $\phi\preceq 1$ is active, then the slot~$(P^\phi,\fm,\hat a)$  in $K^\phi$ is strictly normal.
\textup{(}Hence if $(P,\fm,\hat a)$ is strictly normal, then $(P,\fm,\hat a)$ is quasilinear, and if
in addition $(P,\fm,\hat a)$ is linear, then it is deep.\textup{)} 
\end{lemma}

\noindent
As to Proposition~\ref{normalrefine}, here is a weak version for strict normality:

\begin{lemma} \label{stronglynormalrefine}
Suppose $(P, \fm, \hat a)$ is a strictly normal hole in $K$ and $\hat a -a\prec_{\Delta(\fv)} \fv^{r+w+1}\fm$ where 
$\fv:=\fv(L)$. Then its refinement
$(P_{+a}, \fm, \hat a -a)$ is also strictly normal.
\end{lemma}

\begin{proof} %[Proof of Lemma~\ref{stronglynormalrefine}]
As in the proof of Proposition~\ref{normalrefine} we 
arrange $\fm=1$.
We can also assume~$P_1\asymp 1$. From $P=P(0)+P_1+P_{>1}$ we get 
$$P(a)\ =\ P(0)+ P_1(a) + P_{>1}(a),$$
where $P(0)\prec_{\Delta(\fv)} \fv^{w+1}$ and $P_{>1}(a) \preceq P_{>1}\prec_{\Delta(\fv)}\fv^{w+1}$ by (N2s) and $a\prec 1$;
 we   show that also $P_1(a)\prec_{\Delta(\fv)}\fv^{w+1}$. To see this note that
$$0\ =\ P(\hat a)\ =\ P(0)+ P_1(\hat a)+ P_{>1}(\hat a),$$ 
where as before $P(0), P_{>1}(\hat a) \prec_{\Delta(\fv)} \fv^{w+1}$, so $P_1(\hat a)\prec_{\Delta(\fv)} \fv^{w+1}$. 
Lemma~\ref{lem:diff operator at small elt} applied to~$(\hat K, \preceq_{\Delta(\fv)}, P_1)$ in place of $(K, \preceq, P)$, with $m=w+1$, $y=a-\hat a$,  yields~$P_1(a-\hat a)\prec_{\Delta(\fv)} \fv^{w+1}$,
hence $$P_1(a)\ =\ P_1(a -\hat a) + P_1(\hat a)\ \prec_{\Delta(\fv)}\ \fv^{w+1}$$ as claimed.  It remains to
use $\fv(L_{P_{+a}})\asymp_{\Delta(\fv)} \fv$ and the normality of $(P_{+a}, 1, \hat a -a)$ obtained from
Proposition~\ref{normalrefine} and its proof. 
\end{proof}

\noindent
We also have a version of Lemma~\ref{stronglynormalrefine} for $Z$-minimal slots, obtained from that lemma via Lemma~\ref{lem:from cracks to holes}: 

\begin{lemma} \label{stronglynormalrefine, cracks}
Suppose $(P, \fm, \hat a)$ is $Z$-minimal and strictly normal. Set $\fv:=\fv(L)$,
and suppose $\hat a -a\prec_{\Delta(\fv)} \fv^{r+w+1}\fm$. Then the refinement
$(P_{+a}, \fm, \hat a -a)$ of~$(P, \fm, \hat a)$ is   strictly normal.
\end{lemma}

\noindent
Next two  versions of Proposition~\ref{easymultnormal}:

\begin{lemma}\label{lem:strongly normal refine, 1}
Suppose $(P,\fm,\hat a)$ is a strictly normal hole in $K$, $\hat a\prec\fn\preceq\fm$, and~$[\fn/\fm]< \big[\fv(L)\big]$. Then the refinement $(P,\fn,\hat a)$ of $(P,\fm,\hat a)$
 is strictly normal.
 \end{lemma}
\begin{proof}
As in the proof of Proposition~\ref{easymultnormal} we arrange $\fm=1$ and,
setting  $\fv:=\fv(L)$, $\tilde L:=L_{P_{\times\fn}}$, show
that $\order(\tilde L)=r$, $\fv(\tilde L)\asymp_{\Delta(\fv)} \fv$, and that~(N2) holds with~$\fm$ replaced by $\fn$. 
Now $[\fn]< [\fv]$ yields $\fn\asymp_{\Delta(\fv)} 1$; together with $(P_{\times\fn})_1\asymp_{\Delta(\fv)} \fn P_1$
this gives 
$P(0) \prec_{\Delta(\fv)} \fv^{w+1} P_1 \asymp_{\Delta(\fv)} \fv^{w+1} (P_{\times\fn})_1$. 
Hence (N2s) holds with $\fm$ replaced by~$\fn$. 
Lemma~\ref{lem:strongnormal pos criterion} now yields that  $(P,\fn,\hat a)$
is strictly normal.
\end{proof}

\begin{lemma}\label{lem:strongly normal refine, 2}  
Suppose $(P,\fm,\hat a)$ is a strictly normal hole in $K$ and $\hat a \prec_{\Delta(\fv)} \fm$ where $\fv:=\fv(L)$.
Assume also that for all $q\in\Q^>$ there is given an element $\fv^q$ of~$K^\times$ with $v(\fv^q)=q\,v(\fv)$.
Then for all sufficiently small $q\in\Q^>$ and~$\fn$ with~$\fn\asymp\fv^q\fm$ we have: $\hat a \prec \fn$
and the refinement $(P,\fn,\hat a)$ of $(P,\fm,\hat a)$ is strictly normal.
\end{lemma}
\begin{proof}
We arrange $\fm=1$ as usual, and
take $q_0\in\Q^>$ such that $\hat a \prec \fv^{q_0}$ and  $P(0) \prec_{\Delta(\fv)} \fv^{w+1+q_0} P_1$.
Let  $q\in\Q$, $0<q\leq q_0$, and suppose $\fn\asymp\fv^q$. Then $(P,\fn,\hat a)$ is a refinement of $(P,1,\hat a)$, and the proof of Proposition~\ref{easymultnormal} gives: $\tilde L:=L_{P_{\times\fn}}$ has order~$r$ with $\fv(\tilde L)\asymp_{\Delta(\fv)} \fv$, $\fn P_1\asymp_{\Delta(\fv)}(P_{\times\fn})_1$, and (N2) holds with $\fm$ replaced by~$\fn$. 
Hence $$P(0) \prec_{\Delta(\fv)} \fv^{w+1+q_0} P_1 \preceq \fv^{w+1} \fn P_1  \asymp_{\Delta(\fv)} \fv^{w+1} (P_{\times\fn})_1.$$ 
Now as in the proof of the  previous lemma we conclude that $(P,\fn,\hat a)$ is strictly normal.
\end{proof}

\begin{remarkNumbered} \label{rem:strongly normal refine, 2}
In Lemmas~\ref{lem:strongly normal refine, 1} and \ref{lem:strongly normal refine, 2} we assumed that $(P,\fm,\hat a)$ is a strictly normal hole in $K$. By Lemma~\ref{lem:from cracks to holes} these lemmas go through if this hypothesis is replaced by
``$(P,\fm,\hat a)$ is a strictly normal $Z$-minimal slot in $K$''.
\end{remarkNumbered}

\noindent
We now turn to refining  a given normal hole to a  strictly normal hole. We only do this under additional hypotheses, tailored so that we may employ Lemma~\ref{lem:nepsilon, refined}. Therefore we assume in the rest of this subsection: {\em  $K$ is $\d$-valued and for all $\fv$ and $q\in\Q^>$   we are given
an element $\fv^q$ of $K^\times$ with $(\fv^q)^\dagger=q\fv^\dagger$.} 
Note that then $v(\fv^q)=q\,v(\fv)$ for such $q$. (In particular, $\Gamma$ is divisible.)  
We also adopt the
convention that if
$\order L=r$, then $\fv:=\fv(L)$.

\begin{lemma}\label{lem:achieve strong normality}
Suppose  $(P,\fm,\hat a)$ is a normal hole in $K$ and   $\hat a-a\preceq  \fv^{w+2}\fm$. Then the refinement
$(P_{+a},\fm,\hat a-a)$  of
 $(P,\fm,\hat a)$  is strictly normal.
\end{lemma}
\begin{proof}
As usual we arrange that $\fm=1$.  By Proposition~\ref{normalrefine}, 
$(P_{+a},1,\hat a-a)$ is  normal; the proof of this proposition 
gives
$\order(L_{P_{+a}})=r$, $\fv(L_{P_{+a}})\sim_{\Delta(\fv)} \fv$, $(P_{+a})_1 \sim_{\Delta(\fv)} P_1$,
and (N2) holds with $\fm=1$ and~$P$ replaced by $P_{+a}$.
It  remains to show that
$P_{+a}(0) \prec_{\Delta(\fv)} \fv^{w+1} (P_{+a})_1$, equivalently, $P(a)\prec_{\Delta(\fv)} \fv^{w+1} P_1$.

Let $\hat L:=L_{P_{+\hat a}}\in\hat K[\der]$ and $R:=P_{>1}\in K\{Y\}$;
note that $P_{(\i)}=R_{(\i)}$ for $\abs{\i}>1$ and $R  \prec_{\Delta(\fv)} \fv^{w+1} P_1$. Hence
Taylor expansion and $P(\hat a)=0$ give
\begin{align*}
P(a)	&\ = \ P(\hat a)+\hat L(a-\hat a) + \sum_{\abs{\i} > 1} P_{(\i)}(\hat a)\cdot(a-\hat a)^\i \\
		&\ =\ \hat L(a-\hat a)+\sum_{\abs{\i} > 1} R_{(\i)}(\hat a)\cdot(a-\hat a)^\i \\
		& \qquad\qquad\qquad\text{where $R_{(\i)}(\hat a)\cdot(a-\hat a)^\i \ \prec_{\Delta(\fv)}\ \fv^{w+1} P_1$ for $\abs{\i} > 1$,}
\end{align*} 
so it is enough to show $\hat L(a-\hat a)\prec_{\Delta(\fv)}\fv^{w+1} P_1$.
Lemma~\ref{lem:linear part, new} applied to $(\hat K, \hat a)$ in place of $(K, a)$ gives
$\order \hat L = r$ and $\hat L\sim_{\Delta(\fv)} L$. Since $\hat K$ is $\d$-valued,
Lemma~\ref{lem:nepsilon, refined}
yields a $q\in\Q$ with $w+1<q\leq w+2$ and a $\fw$ such that $\hat L\fv^q \asymp \fw\, \fv^q\, \hat L$ where~$[\fw]\leq [\fv^\dagger]$ 
and 
hence $\fw \asymp_{\Delta(\fv)} 1$ (see the remark before Lemma~\ref{lem:steep1}).
With~$\fn\asymp a-\hat a$ we have $\fn\preceq\fv^{w+2}\preceq\fv^q\prec_{\Delta(\fv)} \fv^{w+1}$ and therefore
$$\hat L(a-\hat a)\ \preceq\ \hat L\fn\ \preceq\ \hat L\fv^q\ \asymp\ \fw\, \fv^q\, \hat L\ \asymp_{\Delta(\fv)}\ \fv^q \hat L\ \prec_{\Delta(\fv)}\ \fv^{w+1} \hat L.$$
Hence $\hat L(a-\hat a)\prec_{\Delta(\fv)}\fv^{w+1} P_1$ as required.
\end{proof}

\noindent
In particular, if  $(P,\fm,\hat a)$ is a normal hole in $K$ and   $\hat a\preceq  \fv^{w+2}\fm$, then  $(P,\fm,\hat a)$ is strictly  normal. 

\begin{cor}\label{cor:achieve strong normality, 1} 
Suppose $(P,\fm,\hat a)$ is $Z$-minimal, deep, and normal. If $(P,\fm,\hat a)$ is special, then
$(P,\fm,\hat a)$ has a deep and 
strictly normal refinement~$(P_{+a},\fm,{\hat a-a})$ where $\hat a-a\prec_{\Delta(\fv)}\fm$ and $\fv(L_{P_{+a,\times\fm}})\asymp_{\Delta(\fv)}\fv$. \textup{(}Note that
if $K$ is $r$-linearly newtonian, and  $\upo$-free if $r>1$, then $(P,\fm,\hat a)$ is special by Lemma~\ref{lem:special dents}.\textup{)} 
\end{cor}
\begin{proof} By Lemma~\ref{lem:from cracks to holes} we arrange that $(P,\fm,\hat a)$ is a hole in $K$.
If $(P,\fm,\hat a)$ is special, Corollary~\ref{specialvariant}
gives an  $a$ such that $\hat a-a\preceq  \fv^{w+2}\fm$,  and then the refinement
$(P_{+a},\fm,\hat a-a)$  of $(P,\fm,\hat a)$  is strictly normal by Lemma~\ref{lem:achieve strong normality}, and deep 
with $\fv(L_{P_{+a,\times\fm}})\asymp_{\Delta(\fv)}\fv$ by Lemma~\ref{lem:deep 2}.
\end{proof}

\noindent
This leads to a useful variant of the Normalization Theorem~\ref{mainthm}:

\begin{cor}\label{cor:achieve strong normality, 2}
Suppose $K$ is $\upo$-free and $r$-linearly newtonian. Then every $Z$-minimal slot in $K$ of order~$r$ has a refinement
$(P,\fm,\hat a)$ such that $(P^\phi,\fm,\hat a)$ is deep and strictly normal, eventually.
\end{cor}
\begin{proof}
Let a $Z$-minimal slot in $K$ of order~$r$ be given. Use Theorem~\ref{mainthm} to refine it to a slot $(P,\fm,\hat a)$ in $K$
with an active $\phi_0$ such that the slot 
$(P^{\phi_0},\fm,\hat a)$ in~$K^{\phi_0}$ is deep and normal.
% Replacing~$(P,\fm,\hat a)$ by 
%$(P^\phi,\fm,\hat a)$  for suitable active $\phi\preceq 1$ we arrange that~$(P,\fm,\hat a)$ is deep and normal. 
Corollary~\ref{cor:achieve strong normality, 1} gives a deep and  strictly normal refinement~$(P^{{\phi}_0}_{+a},\fm,\hat a-a)$ of $(P^{\phi_0},\fm,\hat a)$.  By 
 Lemma~\ref{lem:normality comp conj, strong} the slot 
 $(P^\phi_{+a},\fm, {\hat a -a})$ in~$K^\phi$ is deep and strictly normal, for all active $\phi\preceq \phi_0$ (in $K$). Thus $(P_{+a},\fm, \hat a -a)$ refines the original $Z$-minimal slot in $K$ and has the desired property. 
\end{proof}

\noindent
Corollaries~\ref{degmorethanone} and~\ref{cor:achieve strong normality, 2} have the following consequence:

\begin{cor}\label{cor:achieve strong normality, 3} Suppose $K$ is $\upo$-free. Then
every minimal hole in $K$ of order $r$ and degree~$>1$ has a refinement  $(P,\fm,\hat a)$ such that $(P^\phi,\fm,\hat a)$ is deep and strictly normal, eventually. 
\end{cor}

\noindent
%Similarly to Corollary~\ref{cor:achieve strong normality, 2}, with the role of  Theorem~\ref{mainthm}  in its proof taken over by
%Proposition~\ref{varmainthm}, we obtain from 
Corollary~\ref{cor:achieve strong normality, 1} also gives the following variant of Corollary~\ref{cor:achieve strong normality, 2}, where the role of  Theorem~\ref{mainthm}  in its proof is taken over by
Proposition~\ref{varmainthm}:

\begin{cor}\label{cor:achieve strong normality, 4} Suppose $K$ is $\upo$-free. Then every $Z$-minimal special  slot  in~$K$ of order~$r$ has a refinement~$(P,\fm,\hat a)$ such that $(P^\phi,\fm,\hat a)$
is deep and strictly normal, eventually.
\end{cor}

\section{Isolated Slots}\label{sec:isolated} 

\noindent
In this short section we study the concept of isolation, which plays well together with normality. 
{\em Throughout this section $K$ is an $H$-asymptotic field with small derivation and with rational asymptotic integration.}\/ 
We let $a$, $b$ range over $K$ and~$\phi$, $\fm$, $\fn$, $\fw$ over $K^\times$. We also let
$(P,\fm,\hat a)$ be a slot in $K$ of order $r\geq 1$.
Recall that~$v(\hat a-K)$ is a cut in $\Gamma$ without largest element.
Note that $v\big( (\hat a-a) - K \big) = v(\hat a-K)$
and~$v( \hat a \fn - K) = v(\hat a-K)+v\fn$. 

\begin{definition}\label{def:isolated}
We say that $(P,\fm,\hat a)$  is {\bf isolated}\index{slot!isolated} if for all $a\prec\fm$,
$$\order(L_{P_{+a}})=r\ \text{ and }\ \exc^{\ev}(L_{P_{+a}}) \cap v(\hat a-K)\ <\  v(\hat a-a);$$
equivalently, for all $a\prec \fm$:  $\order(L_{P_{+a}})=r$ and whenever
$\fw \preceq \hat a-a$ is such that~$v(\fw) \in \exc^{\ev}(L_{P_{+a}})$, then 
$\fw\prec \hat a-b$ for all $b$.
\end{definition}

\noindent
In particular,  if $(P,\fm,\hat a)$ is isolated, then $v(\hat a)\notin \exc^{\ev}(L_P)$. If $(P,\fm,\hat a)$ is isolated, then so is every equivalent slot in~$K$,
as well as $(bP,\fm,\hat a)$ for $b\neq 0$ and the slot~$(P^\phi,\fm,\hat a)$ in $K^\phi$ for active $\phi$ in $K$. Moreover:

\begin{lemma}\label{lem:isolated refinement}
If $(P,\fm,\hat a)$ is isolated, then so is any refinement~$(P_{+a},\fn,\hat a-a)$  of~it. 
\end{lemma}
\begin{proof} For the case $\fn=\fm$, use $v\big( (\hat a-a) - K \big) = v(\hat a-K)$. The case $a=0$ is clear. The general case reduces to these two special cases.
\end{proof}

\begin{lemma}\label{lem:isolated}
Suppose  $(P,\fm,\hat a)$ is isolated. Then the multiplicative con\-ju\-gate $(P_{\times\fn},\fm/\fn,\hat a/\fn)$ of 
$(P,\fm,\hat a)$ by $\fn$ is isolated.
\end{lemma}
\begin{proof} Let $a\prec\fm/\fn$. Then $a\fn\prec \fm$, so
$\order(L_{P_{\times \fn, +a}})=\order(L_{P_{+a\fn,\times \fn}})=\order(L_{P_{+a\fn}})=r$. Suppose $\fw\preceq (\hat a/\fn)-a$ and
$ v(\fw) \in \exc^{\ev}\big(L_{P_{\times\fn,+a}}\big)$.
Now $L_{P_{\times\fn,+a}} = L_{P_{+a\fn,\times\fn}} = L_{P_{+a\fn}}\fn$  and  thus $\fw\fn \preceq \hat a-a\fn$, 
$v(\fw\fn)\in \exc^{\ev}\big(P_{+a\fn}\big)$.
But~$(P,\fm,\hat a)$ is isolated, so $v(\fw\fn)> v(\hat a-K)$ and hence~$v(\fw)> v\big((\hat a/\fn)-K\big)$.
Thus $(P_{\times\fn},\fm/\fn,\hat a/\fn)$  is isolated.
\end{proof}

\begin{lemma}\label{lem:isolated normal}
Suppose   $K$ is $\upl$-free or   $r=1$, and $(P,\fm,\hat a)$ is normal. Then  
$$ \text{$(P,\fm,\hat a)$  is isolated} \quad\Longleftrightarrow\quad
\exc^{\ev}(L_P) \cap v(\hat a-K)\ \leq\  v\fm.$$
\end{lemma}
\begin{proof} Use Lemma~\ref{lem:excev normal}; for the direction $\Rightarrow$, use also that  $\hat a-a\prec\fm$ iff $a\prec\fm$.
\end{proof}

\begin{lemma}\label{lem:isolated deg 1} 
Suppose  $\deg P=1$. Then
$$ \text{$(P,\fm,\hat a)$  is isolated} \quad\Longleftrightarrow\quad
\exc^{\ev}(L_P) \cap v(\hat a-K)\  \leq\  v\fm.$$
\end{lemma}
\begin{proof}
Use that $\order L_P=r$ and $L_{P_{+a}}=L_P$ for all $a$.
\end{proof}

\begin{prop}\label{prop:achieve isolated}
Suppose~$K$ is $\upl$-free or $r=1$, and $(P,\fm,\hat a)$ is normal.  
Then $(P,\fm,\hat a)$ has an isolated refinement.
\end{prop}
\begin{proof}
Suppose $(P,\fm,\hat a)$  is not already isolated.
Then Lemma~\ref{lem:isolated normal} gives $\gamma$~with
$$\gamma\in\exc^{\ev}(L_P)\cap v(\hat a-K),\quad 
\gamma>v\fm.$$
We have $|\exc^{\operatorname{\ev}}(L_P)|\le r$,  by [ADH, p.~481] if $r=1$,  and Corollary~\ref{cor:size of excev, strengthened} and $\upl$-freeness of $K$ if $r>1$.
Hence we can take $\gamma:= \max\exc^{\ev}(L_P)\cap v(\hat a-K)$, and then~$\gamma > v\fm$.
Take~$a$ and~$\fn$ with~$v(\hat a-a)>\gamma=v(\fn)$;
then $(P_{+a},\fn,\hat a-a)$ is a refinement of~$(P,\fm,\hat a)$ and~$a\prec\fm$. 
Let $b\prec \fn$; then 
$a+b\prec\fm$, so by~Lemma~\ref{lem:excev normal},
$$\order(L_{(P_{+a})_{+b}})\ =\ r, \qquad
\exc^{\ev}(L_{(P_{+a})_{+b}})\ =\ 
\exc^{\ev}(L_P).$$
Also $v\big((\hat a-a)-b\big)>\gamma$, hence
$$\exc^{\ev}\big(L_{(P_{+a})_{+b}}\big)  \cap v\big((\hat a-a)-K\big)\  =\ 
\exc^{\ev}(L_P)\cap v(\hat a-K)\ \le\ \gamma\ <\ v\big((\hat a-a)-b\big).$$
Thus $(P_{+a},\fn,\hat a-a)$  is isolated. 
\end{proof}

\begin{remarkNumbered}\label{rem:achieve isolated}
Proposition~\ref{prop:achieve isolated}  goes through if instead of assuming that  $(P,\fm,\hat a)$ is normal, we assume that
 $(P,\fm,\hat a)$ is linear. (Same argument,  
using Lem\-ma~\ref{lem:isolated deg 1} in place of Lemma~\ref{lem:isolated normal} and 
$L_{(P_{+a})_{+b}}=L_P$ in place of Lemma~\ref{lem:excev normal}.) 
\end{remarkNumbered}

\begin{cor}\label{cor:isolated r=1}
Suppose $r=1$, and $(P,\fm,\hat a)$ is normal or linear.
If~$\exc^{\ev}(L_P)=\emptyset$, then $(P,\fm,\hat a)$ is isolated.
If $\exc^{\ev}(L_P)\neq\emptyset$, so $\exc^{\ev}(L_P)=\{v\mathfrak g\}$ where $\mathfrak g\in K^\times$, then~$(P,\fm,\hat a)$ is isolated iff $\fm\preceq\mathfrak g$ or $\hat a-K\succ\mathfrak g$.
\end{cor}

\noindent
This follows immediately from Lemmas~\ref{lem:isolated normal} and~\ref{lem:isolated deg 1}.
The results in the rest of this subsection are the {\em raison d'\^etre}\/ of isolated holes:

\begin{prop}\label{prop:2.12 isolated}  
Suppose $K$ is $\upo$-free and $(P,\fm,\hat a)$ is an isolated hole in~$K$ 
which is normal or linear.
Let $\hat b$ in an immediate
asymptotic extension of $K$ satisfy~${P(\hat b)=0}$ and $\hat b\prec\fm$. Then~$v(\hat a-a)=v(\hat b-a)$ for all~$a$, so~$\hat b\notin K$. \end{prop}
\begin{proof}
Replacing $(P,\fm,\hat a)$,~$\hat b$ by~$(P_{\times\fm},1,\hat a/\fm)$,~$\hat b/\fm$, we arrange $\fm=1$. Let~$a$ be given; we show $v(\hat a-a)=v(\hat b -a)$. This is clear if $a\succeq 1$, so assume~${a\prec 1}$. 
Corollary~\ref{cor:normal=>quasilinear} (if $(P,\fm,\hat a)$ is normal)
and Lemma~\ref{lem:lower bd on ndeg} (if $(P,\fm,\hat a)$ is linear) give~${\ndeg P=1}$. Thus $P$ is in newton position at $a$ by Corollary~\ref{cor:ref 2n}. Moreover~$v(\hat a-a) \notin \exc^{\ev}(L_{P_{+a}})$,
hence $v(\hat a-a)=v^{\ev}(P,a)$ by Lemma~\ref{lem:14.3 complement}. 
Likewise, if~$v(\hat b-a)\notin \exc^{\ev}(L_{P_{+a}})$, then 
$v(\hat b-a)=v^{\ev}(P,a)$ by Lemma~\ref{lem:14.3 complement}, so~$v({\hat a-a})=v(\hat b-a)$.

Thus to finish the proof it is enough to show that $\exc^{\ev}(L_{P_{+a}})\cap v(\hat b-K)\leq 0$. 
Now~$\abs{\exc^{\ev}(L_{P_{+a}})} \leq r$ by   Corollary~\ref{cor:sum of nwts}, so
we have $b\prec 1$ such that $$ \exc^{\ev}(L_{P_{+a}}) \cap v(\hat b-K)\ <\ v(\hat b-b),$$
in particular, $v(\hat b -b)\notin \exc^{\ev}(L_{P_{+a}})$. 
If $(P,\fm,\hat a)$ is normal, then
Lemma~\ref{lem:excev normal} gives
$$\exc^{\ev}(L_{P_{+a}})\ =\ \exc^{\ev}(L_P)\ =\ \exc^{\ev}(L_{P_{+b}}),$$
so by the above with $b$ instead of $a$ we have $v(\hat a-b)=v(\hat b-b)$. 
If $(P,\fm,\hat a)$ is linear, then~$L_{P_{+a}}=L_P=L_{P_{+b}}$, and we obtain likewise $v(\hat a-b)=v(\hat b-b)$. 
Hence
$$\exc^{\ev}(L_{P_{+a}})\cap v(\hat b-K)\ \subseteq\ 
\exc^{\ev}(L_{P_{+a}})\cap \Gamma^{< v(\hat a-b)}\ \subseteq\
\exc^{\ev}(L_P)\cap v(\hat a-K)\ \leq\ 0.$$
using Lemmas~\ref{lem:isolated normal} and~\ref{lem:isolated deg 1} for the last step.
%and $$\exc^{\ev}(L_P)\cap v(\hat a-K)\ \subseteq\ \exc^{\ev}(L_P)\cap v(\hat a-K)\ <\ v(\hat a)\ ??<\ 0$$
%since $(P,\fm,\hat a)$ is isolated.
\end{proof}

\noindent
Combining Proposition~\ref{prop:2.12 isolated} with Corollary~\ref{corisomin} yields:

\begin{cor}\label{cor:2.12 isolated}
Let $K$, $(P,\fm,\hat a)$, $\hat b$ be as in Proposition~\ref{prop:2.12 isolated}, and assume also that $(P,\fm,\hat a)$ is $Z$-minimal. Then there is an isomorphism $K\langle\hat a\rangle\to K\langle \hat b\rangle$ of valued differential fields
over $K$ sending~$\hat a$ to $\hat b$.
\end{cor}

\noindent
Using the Normalization Theorem, we now obtain:

\begin{cor}\label{cor:2.12 isolated, 2}
Suppose $K$ is $\upo$-free and $\Gamma$ is divisible. Then every minimal hole in $K$ of order $r$ has an
isolated refinement 
$(P,\fm,\hat a)$ such that for any $\hat b$ in an immediate
asymptotic extension of $K$ with~${P(\hat b)=0}$ and $\hat b\prec\fm$ there is an isomorphism~$K\langle\hat a\rangle\to K\langle \hat b\rangle$ of valued differential fields
over $K$ sending~$\hat a$ to $\hat b$.
\end{cor}
\begin{proof}
Given a minimal linear hole in $K$ of order $r$, use Remark~\ref{rem:achieve isolated}  to refine it to an isolated minimal linear hole 
$(P,\fm,\hat a)$ in $K$ of order $r$, and use Corollary~\ref{cor:2.12 isolated}.
Suppose we are given a minimal non-linear hole in $K$ of order $r$. Then $K$ is $r$-linearly newtonian by Corollary~\ref{degmorethanone}. Then Theorem~\ref{mainthm} yields a refinement~$(Q, \fw, \hat{d})$ of it and an active $\theta$
in $K$ such that  the minimal hole $(Q^\theta, \fw, \hat{d})$ in $K^\theta$ is normal. Proposition~\ref{prop:achieve isolated} gives an
isolated refinement $(Q^\theta_{+d}, \fv , \hat{d}-d)$ of $(Q^\theta, \fw, \hat{d})$. 
Suitably refining $(Q^\theta_{+d}, \fv , \hat{d}-d)$ further followed by compositionally conjugating with a suitable active element of $K^\theta$ yields by Theorem~\ref{mainthm} and Lemma~\ref{lem:isolated refinement} a
refinement~$(P, \fm, \hat a)$ of $(Q, \fw, \hat{d})$ (and thus of the originally given hole) and an active $\phi$ in $K$ such that
 $(P^\phi, \fm, \hat a)$ is both normal and isolated. Then $(P, \fm, \hat a)$ is isolated, and we can apply Corollary~\ref{cor:2.12 isolated} to $K^\phi$ and $(P^\phi, \fm, \hat a)$ in the role of $K$ and $(P,\fm,\hat a)$. 
%Lemma~\ref{lem:isolated refinement}, and Proposition~\ref{prop:achieve isolated} we obtain 
% a normal isolated refinement $(P,\fm,\hat a)$ of this hole, and then Corollary~\ref{cor:2.12 isolated} applies.
 \end{proof}

\noindent
For $r=1$ we can replace ``$\upo$-free'' in  Proposition~\ref{prop:2.12 isolated} and Corollary~\ref{cor:2.12 isolated} by the weaker ``$\upl$-free'' (same proofs,  using Lemma~\ref{lem:14.3 complement, order 1} instead of Lem\-ma~\ref{lem:14.3 complement}):

\begin{prop}\label{prop:2.12 isolated, r=1}
Suppose $K$ is $\upl$-free, $(P,\fm,\hat a)$ is an isolated   hole in~$K$ of order~$r=1$, and suppose $(P,\fm,\hat a)$ is normal or  linear. Let $\hat b$ in an immediate
asymptotic extension of $K$ satisfy~${P(\hat b)=0}$ and $\hat b\prec\fm$. Then~$v(\hat a-a)=v(\hat b-a)$ for all~$a$. \textup{(}Hence if~$(P,\fm,\hat a)$ is $Z$-minimal, then there is an isomorphism $K\langle\hat a\rangle\to K\langle \hat b\rangle$ of valued differential fields
over $K$ sending~$\hat a$ to $\hat b$.\textup{)}
\end{prop}

\noindent
This leads to an analogue of Corollary~\ref{cor:2.12 isolated, 2}: 

\begin{cor}\label{cor:2.12 isolated, 2, r=1}
Suppose $K$ is $\upl$-free and $\Gamma$ is divisible.
 Then every quasilinear minimal  hole in $K$ of order~$r=1$ has an
isolated refinement 
$(P,\fm,\hat a)$ such that for any $\hat b$ in an immediate
asymptotic extension of $K$ with~${P(\hat b)=0}$ and $\hat b\prec\fm$ there is an isomorphism~$K\langle\hat a\rangle\to K\langle \hat b\rangle$ of valued differential fields
over $K$ sending~$\hat a$ to $\hat b$.
\end{cor}
\begin{proof}
Suppose we are given a quasilinear minimal  hole in $K$ of order~$r=1$.
Then Corollary~\ref{mainthm order 1} yields a refinement $(Q,\fw,\hat d)$ of it and an active $\theta$ in $K$
such that the quasilinear minimal hole $(Q^\theta,\fw,\hat d)$ in $K^\theta$ of order $1$ is normal. 
Proposition~\ref{prop:achieve isolated} gives an isolated refinement $(Q^{\theta}_{+d}, \fv, \hat d -d)$ of $(Q^\theta,\fw,\hat d)$,
and then Corollary~\ref{mainthm order 1} yields a refinement~$(P,\fm,\hat a)$ of $(Q,\fw,\hat d)$ and an active $\phi$ in $K$ such
that $(P^\phi,\fm,\hat a)$ is  normal and isolated.  Now apply Proposition~\ref{prop:2.12 isolated, r=1}
with $K^\phi$ and $(P^\phi, \fm, \hat a)$ in the role of $K$ and $(P,\fm,\hat a)$. 
\end{proof}

\noindent
Next a variant of Lemma~\ref{lem:no hole of order <=r}  for $r=1$ without assuming  $\upo$-freeness:

\begin{cor}\label{cor:1-newt => no qlin Zmin dent order 1}
Suppose $K$ is  $1$-newtonian and $\Gamma$ is divisible. Then $K$ has no  quasilinear $Z$-minimal slot  of order $1$.
\end{cor}
\begin{proof}
By Proposition~\ref{prop:char 1-linearly newt}, $K$ is $\upl$-free.
 Towards a contradiction, let $(P,\fm,\hat a)$ be a quasilinear $Z$-minimal slot in~$K$ of order~$1$.
 By Lemma~\ref{lem:from cracks to holes} we arrange that~$(P,\fm,\hat a)$  is a hole in $H$. Using Corollary~\ref{mainthm, r=1}, 
  Lemma~\ref{lem:isolated refinement} and the remark before it,
  and Proposition~\ref{prop:achieve isolated}, we can refine further so that
  $(P^\phi,\fm,\hat a)$ is normal and isolated for some active $\phi$ in $K$.
   Then there is no $y\in K$ with~$P(y)=0$ and~$y\prec \fm$, by
  Proposition~\ref{prop:2.12 isolated, r=1}, contradicting Lemma~\ref{lem:zero of P} for~$L=K$.
\end{proof}

\noindent
Finally, for isolated linear holes, without additional hypotheses:

\begin{lemma}\label{lem:isolated, d=1}
Suppose $(P,\fm,\hat a)$ is an  isolated linear hole in $K$, and~${\hat a-a}\prec\fm$.  Then
$P(a)\neq 0$, and $\gamma=v(\hat a-a)$ is the unique element of $\Gamma\setminus\exc^{\ev}(L_P)$ such that~$v^{\ev}_{L_P}(\gamma)=v\big(P(a)\big)$.  
\end{lemma}
\begin{proof}
By Lemma~\ref{lem:isolated deg 1}, $\gamma:=v(\hat a-a)\in\Gamma\setminus\exc^{\ev}(L_P)$.
Since $\deg P =1$,    $$L_P(\hat a-a)\ =\ L_P(\hat a)-L_P(a)\ =\ -P(0)-L_P(a)\ =\ -P(a),$$ 
so $P(a)\neq 0$. By Lemma~\ref{lem:ADH 14.2.7},  
$v^{\ev}_{L_P}(\gamma)=v\big(L_P(\hat a-a)\big)=v\big(P(a)\big)$.
\end{proof}

\noindent
In \cite{ADHld} we shall prove a version of Proposition~\ref{prop:2.12 isolated} without the hypothesis that~$\hat b$ lies in an
immediate extension of $K$. 
% \marginpar{skipped for now}% (but with extra requirements on $K$).
In Section~\ref{sec:ultimate} below we consider, in a  more restricted setting, a variant of isolated slots, with 
ultimate exceptional values taking over the role played by exceptional values in Definition~\ref{def:isolated}.

\section{Holes of Order and Degree One}\label{sec:holes of c=(1,1,1)} 

\noindent
In this section $K$ is a $\d$-valued field of $H$-type with small derivation and rational asymptotic integration. (Later on we will impose 
additional restrictions on $K$.)   We also let $\hat K$ be an immediate asymptotic extension of $K$.
We focus here on slots of complexity $(1,1,1)$ in~$K$. 
As a byproduct we obtain in Corollary~\ref{cor:achieve strong normality, general} a partial generalization 
of Corollary~\ref{cor:achieve strong normality, 3} to minimal holes in $K$ of arbitrary   degree.
First we establish in the next subsection a useful formal identity.
We  let~$j$,~$k$ range over~$\N$ (in addition to $m$, $n$, as usual).  

\subsection*{An integration identity}
Let $R$ be a differential ring, and let $f$,~$g$,~$h$,  range over~$R$.   We use 
$\int f= g + \int h$ as a suggestive way to express that $f=g'+h$, and likewise, $\int f = g-\int h$ means that $f=g'-h$. For example, 
$$\int f'g\ =\ fg-\int fg'\qquad \text{(integration by parts)}.$$
Let $\ex,\xi\in R^\times$ satisfy
$\ex^\dagger=\xi$. We wish to expand~$\int \ex$ by iterated integration by parts. Now for $g=\ex$ we have $g'\in R^\times$ with $\frac{g}{g'}=\frac{1}{\xi}$, so in view of~$\ex=g'\frac{\ex}{g'}$:
$$\int \ex\ =\ \int g'\frac{\ex}{g'}\ =\ \frac{\ex}{\xi} - \int g\bigg(\frac{\ex}{g'}\bigg)',$$
and $$\left(\frac{\ex}{g'}\right)'= \left(\frac{1}{\xi}\right)'=\frac{-\xi'}{\xi^2}=\frac{-\xi^\dagger}{\xi},$$ and thus  
$$\int \ex\ =\ \frac{\ex}{\xi} + \int \frac{\xi^\dagger}{\xi}\ex.$$
More generally, using the above identities for $g=\ex$,
\begin{align*}  \int f\ex\ &=\ \int g'f\frac{\ex}{g'}\ =\ \frac{f}{\xi}\ex - \int g\bigg(f\frac{\ex}{g'}\bigg)'\\
&=\ \frac{f}{\xi} \ex -\int g\left(f'\frac{\ex}{g'}+f\bigg(\frac{\ex}{g'}\bigg)'\,\right)\ =\ \frac{f}{\xi} \ex -\int \left(\frac{f'}{\xi}\ex + fg\Big(\frac{-\xi^\dagger}{\xi}\Big)\right)\\
&=\ \frac{f}{\xi} \ex -\int \left(\frac{f'}{\xi}\ex +\frac{-f\xi^\dagger}{\xi}\ex \right)\ =\ 
\frac{f}{\xi} \ex +\int\left(\frac{\xi^\dagger f -f'}{\xi}\right)\ex.
\end{align*}
Replacing $f$ by $f/\xi^{k}$ gives the following variant of this identity:
$$\int \frac{f}{\xi^{k}}\ex\ =\ 
\frac{f}{\xi^{k+1}}\ex + 
\int \frac{(k+1)\xi^\dagger f -f'}{\xi^{k+1}}\ex.$$
Induction on $m$ using the last identity yields:

\begin{lemma}\label{fex} Set  $\zeta:= \xi^\dagger$. Then
$$\int f\ex \ =\ \sum_{j=0}^m P_{j}(\zeta,f)\frac{\ex}{\xi^{j+1}} + \int P_{m+1}(\zeta,f)\frac{\ex}{\xi^{m+1}},$$
where the $P_j\in \Q\{Z,V\}=\Q\{Z\}\{V\}$ are independent of $R$, $\ex$, $\xi$: 
$$P_0\ :=\ V, \qquad P_{j+1}\ :=\ (j+1)ZP_j-P_j'.$$
%\begin{align*} 
%P_{0}\ &:=\ V,\qquad\qquad \quad Q_{0}\ :=\ ZV-V',\\
% P_{m+1}\ &:=\ Q_{m}, \qquad Q_{m+1}\ :=\ (m+2)ZQ_{m}-Q_{m}'.
% \end{align*}
Thus $P_j=P_{j0}V+P_{j1}V'+\cdots+P_{jj}V^{(j)}$ 
with all  $P_{jk}\in\Q\{Z\}$ 
%for $k=0,\dots,j-1$,
and $P_{jj}=(-1)^j$.  
% is of order at most $j-i-1$ for $j=1,\dots,i-1$ and . 
\end{lemma}

\noindent
For example,
$$P_0=V, \quad P_1 = ZV-V', \quad P_2 = (2Z^2-Z')V-3ZV'+V''.$$

\subsection*{An asymptotic expansion}
{\it In this subsection $\xi\in K$ and $\xi\succ^\flat 1$}\/; equivalently,~$\xi\in K$ satisfies $\xi\succ 1$ and~$\zeta:=\xi^\dagger\succeq 1$. 
{\it We also assume that $\xi\notin \I(K)+K^\dagger$.}\/  
Since $\hat K$ is $\d$-valued of $H$-type with asymptotic integration, it has by [ADH, 10.2.7] an immediate
asymptotic extension
$\hat{K}(\phi)$ with $\phi'=\xi$.  Then the algebraic closure of $\hat{K}(\phi)$ is still $\d$-valued of $H$-type, by [ADH, 9.5], and so [ADH, 10.4.1] yields a
 $\d$-valued $H$-asymptotic extension~$L$ of this algebraic closure with an element $\ex\ne 0$ such that~$\ex^\dagger=\xi$. 
All we need about $L$ below is that it is a $\d$-valued $H$-asymptotic extension of $\hat{K}$ with elements $\phi$ and $\ex$
such that $\phi'=\xi$ and $\ex\ne 0,\ \ex^\dagger=\xi$. Note that then~$L$  has small derivation, and $\xi\succ^\flat 1$ in $L$. (The element $\phi$ will only play an auxiliary role later in this subsection.) 

\begin{lemma}\label{venotin}
$v(\ex)\notin\Gamma$.  
\end{lemma}
\begin{proof}
Suppose otherwise. Take $a\in K^\times$ with $a\ex\asymp 1$. Then 
$a^\dagger+\xi=(a\ex)^\dagger\in \I(L)\cap K=\I(K)$ and thus
$\xi\in \I(K)+K^\dagger$, a contradiction.
\end{proof}

\noindent
By Lemma~\ref{venotin} there is for each $g\in L$ at most one $\hat f\in \hat K$ with $\big(\hat f\frac{\ex}{\xi}\big)'=g$.
Let~$f\in K^\times$ be given with $f \preceq 1$, and suppose  $\hat f\in \hat K$ satisfies
$\big(\hat f\frac{\ex}{\xi}\big)'=f\ex$. 
Our aim is to show that with   $P_j$   as in Lemma~\ref{fex},
the series $\sum_{j=0}^\infty P_{j}(\zeta,f)\frac{1}{\xi^{j}}$ is a kind of asymptotic expansion of $\hat f$. The partial sums $$f_m:=\sum_{j=0}^m P_{j}(\zeta,f)\frac{1}{\xi^{j}}$$ of this series lie in $K$, with $f_0=f$
and $f_n  - f_m   \prec \xi^{-m}$ for $m<n$,
by Lemma~\ref{lem:xi zeta}.
More precisely,  we show:

\begin{prop}\label{uaintfex} 
We have $\hat f-f_m \prec \xi^{-m}$ for all $m$.
\textup{(}Thus:  $f \asymp 1\,  \Rightarrow\, \hat f \sim f$.\textup{)} 
\end{prop} 

\noindent
Towards the proof, note that by
Lemma~\ref{fex} with~$R=L$,
\begin{align} \hat f\,\frac{\ex}{\xi}\ &=\ \sum_{j=0}^m P_{j}(\zeta,f)\frac{\ex}{\xi^{j+1}} + I_m,\quad I_m\in L, \text{ and thus}\notag 
\\
\hat f\ &=\ f_m + \frac{\xi}{\ex}I_m \label{eq:F= ...}
\end{align}
where $I_m\in \hat K \ex$ satisfies $I_m'= P_{m+1}(\zeta,f)\frac{\ex}{\xi^{m+1}}$, a condition that determines $I_m$ uniquely up to an additive constant from $C_L$.
The proof of Proposition~\ref{uaintfex} now rests on the following lemmas:

\begin{lemma}\label{aexpphiy} 
In $L$ we have
$(\ex \xi^l)^{(k)}\ \sim\ \ex \xi^{l+k}$, for all $l\in\Z$ and all $k$.
\end{lemma}
 
\noindent
This is Corollary~\ref{cor:xi 1} with our $L$ in the role of $K$ there, and taking $\ex^\phi$ there as our $\ex\in L$; note that here our $\phi\in L$ with $\phi'=\xi$ is needed. 

\begin{lemma}\label{exi1} Suppose $\ex\succ \xi^{m+1}$. Then $\frac{\xi}{\ex}I_m\prec  \xi^{-m}$.
\end{lemma}
\begin{proof} This amounts to $I_m \prec \frac{\ex}{\xi^{m+1}}$. Suppose $I_m\succeq \frac{\ex}{\xi^{m+1}}\succ 1$. Then we have $I_m'\succeq \big(\frac{\ex}{\xi^{m+1}}\big)'\sim \frac{\ex}{\xi^m}$
by Lemma~\ref{aexpphiy}, so $P_{m+1}(\zeta,f)\frac{\ex}{\xi^{m+1}}\succeq \frac{\ex}{\xi^m}$, and thus $P_{m+1}(\zeta,f)\succeq \xi$, contradicting Lemma~\ref{lem:xi zeta}.  
\end{proof}

\begin{lemma}\label{exi2} Suppose $\ex\preceq \xi^m$. Then $I_m\prec 1$ and $\frac{\xi}{\ex}I_m\prec \xi^{-m}$.
\end{lemma}
\begin{proof} Lemma~\ref{lem:xi zeta} gives 
$$P_{m+1}(\zeta,f) \frac{\ex}{\xi^{m+1}}\preceq \zeta^N\frac{\ex}{\xi^{m+1}}\preceq \frac{\zeta^N}{\xi}\quad\text{ for some $N\in \N$,}$$ so $v(I_m')>\Psi_L$, and thus
$I_m\preceq 1$.  If   $I_m\asymp 1$, then~$v(\frac{\xi}{\ex}I_m)=v(\xi)-v(\ex)\notin \Gamma$, contradicting~$\frac{\xi}{\ex}I_m = \hat f-f_m\in \hat K$, by \eqref{eq:F= ...}. Thus~$I_m\prec 1$.
 Now assume towards a contradiction that $\frac{\xi}{\ex}I_m\succeq \xi^{-m}$. Then~$\frac{\ex}{\xi^{m+1}}\preceq I_m\prec 1$, so
 $I_m'\succeq \big(\frac{\ex}{\xi^{m+1}}\big)'\sim \frac{\ex}{\xi^m}$
by Lemma~\ref{aexpphiy}, and this yields a contradiction as in the proof of Lemma~\ref{exi1}.   
\end{proof}

\begin{proof}[Proof of Proposition~\ref{uaintfex}] Let $m$ be given. If $\ex\succ \xi^{m+1}$, then 
$\hat f-f_m = \frac{\xi}{\ex}I_m\prec \xi^{-m}$ by  Lemma~\ref{exi1}. Suppose
$\ex\preceq \xi^{m+1}$. Then  Lemma~\ref{exi2} (with $m+1$ instead of $m$) gives $\hat f-f_{m+1}\prec  \xi^{-(m+1)}$, hence  $\hat f-f_m=(\hat f-f_{m+1})+(f_{m+1}-f_m)\prec  \xi^{-m}$.   
\end{proof}

\subsection*{Application to linear differential equations of order $1$}
Proposition~\ref{uaintfex} yields   information about the asymptotics of solutions (in~$
\hat K$) of certain linear differential equations of order $1$ over $K$:

\begin{cor}\label{cor:uaintfex, 1}
Let $f,\xi\in K$, $f \preceq 1$, $\xi\succ^\flat 1$, $\xi\notin\I(K)+K^\dagger$, and suppose~${y\in \hat K}$ satisfies $y'+\xi y=f$. Then there is for every $m$ an element $y_m\in K$ with~${y-y_m\prec\xi^{-m}}$. Also, 
$f\asymp 1\ \Rightarrow\ y\sim f\xi^{-1}$.
\end{cor}
\begin{proof}
Take $L$ and $\ex\in L$  as at the beginning of the previous subsection, and
set~$\hat f:=y\xi\in \hat K$. Then for $A:=\der+\xi$ we have $A(\hat f/\xi)=f$, so
$$ \left( \hat f\,\frac{\ex}{\xi}\right)'\ =\ (\hat f/\xi)'\ex + (\hat f/\xi)\xi\ex\ =\  A(\hat f/\xi) \ex\  =\  f\ex,$$
hence $\hat f$ is as in the previous subsection. Now apply  Proposition~\ref{uaintfex}.
\end{proof}

\begin{cor}
Let $g\in K$, $u\in K^\times$ be such that $g\notin \I(K)+K^\dagger$ and $\xi:=g+u^\dagger \succ^\flat 1$. Suppose~$z\in\hat K$ satisfies $z'+gz=u$.
Then $z\sim u/\xi$, and for every $m$ there is a~$z_m\in K$ such that $z-z_m\prec u\xi^{-m}$. 
\end{cor}
\begin{proof} Set $A:= \der+g$.
Then $A_{\ltimes u}=\der+\xi$, so $A(z)=u$ yields for $y:= z/u$ that~$y'+\xi y=1$. Now observe that $\xi\notin \I(K)+K^\dagger$ and use the previous corollary.
\end{proof}

\subsection*{Slots of order and degree $1$}
In the rest of this section we use the material above to analyze slots of order and degree $1$ in $K$. 
{\it Below $K$ is henselian and~$(P,\fm,\hat f)$ is a slot in~$K$ with $\order P=\deg P=1$ and $\hat f\in \hat K\setminus K$. We  let $f$ range over $K$, $\fn$ over $K^\times$, and $\phi$ over active elements of $K$.}\/
Thus
\begin{align*}
P\				&=\  a(Y'+gY-u)\quad\text{where $a\in K^\times,\ g,u\in K$,}  \\
P_{\times\fn}\	&=\   a\fn\big(Y'+(g+\fn^\dagger)Y-\fn^{-1}u\big).
\end{align*} 
Since $K$ is henselian,   $(P,\fm,\hat f)$  is $Z$-minimal and thus equivalent to a hole in $K$, by Lem\-ma~\ref{lem:from cracks to holes}.
Also, $\nval P_{\times\fm}=\ndeg P_{\times\fm}=1$ by Lemma~\ref{lem:lower bd on ndeg}.
We have~$L_P=a(\der+g)$, so 
$$ g\in K^\dagger\ \Longleftrightarrow\ \ker L_P\neq\{0\},\quad\quad\qquad
g\in \I(K)+K^\dagger\ \Longleftrightarrow\ \exc^{\ev}(L_P)\neq\emptyset, $$
using for the second equivalence the remark on $\exc^{\ev}(A)$ preceding Lemma~\ref{lem:v(ker)=exc, r=1}. 
If~$(P,\fm,\hat f)$ is isolated, then $P(f)\neq 0$ for~$\hat f-f\prec\fm$ by Lemmas~\ref{lem:from cracks to holes} and~\ref{lem:isolated, d=1}, so, taking $f=0$, we have $u\neq 0$.

\begin{lemma}\label{lem:excev empty}
Suppose $\der K=K$ and $\I(K)\subseteq K^\dagger$.
Then  $\exc^{\ev}(L_P)=\emptyset$, so 
$(P,\fm,\hat f)$ is isolated by Lemma~\ref{lem:isolated deg 1}.
\end{lemma}
\begin{proof}
Passing to an equivalent hole in $K$, arrange that  $(P,\fm,\hat f)$  is a hole in~$K$. Since~$\der K=K$ and $\hat f\in \hat K\setminus K$,   the remark following Lemma~\ref{lem:at most one zero}
yields $g\notin K^\dagger=\I(K)+K^\dagger$, therefore~$\exc^{\ev}(L_P)=\emptyset$.
\end{proof}

%\noindent
%Corollary~\ref{cor:proper}  implies:

%\begin{cor}\label{cor:excev empty}
%If $K$ is $\upl$-free and $(P,\fm,\hat f)$ is isolated, then
%$$\hat f\sim  u/\big( g+(u/\phi)^\dagger\big) \prec^\flat_\phi 1,\quad\text{ eventually.}$$
%\end{cor}

\noindent
Set $\fv:=\fv(L_{P_{\times\fm}})$; thus $\fv=1$ if $g+\fm^\dagger\preceq 1$ and $\fv=1/(g+\fm^\dagger)$ otherwise. Hence from Example~\ref{ex:order 1 linear steep} and the remarks before Lemma~\ref{lem:deg1 normal} we obtain:
\begin{align*}
\text{$(P,\fm,\hat f)$ is normal}	&\quad\Longleftrightarrow\quad\text{$(P,\fm,\hat f)$ is steep}
\quad\Longleftrightarrow\quad \fv\prec^\flat 1, \\
\text{$(P,\fm,\hat f)$ is deep}		&\quad\Longleftrightarrow\quad \fv\prec^\flat 1 \text{ and }  u \preceq \fm/\fv.
\end{align*}
We have $P(0)=-au$, and if $\fv\prec 1$, then $(P_{\times\fm})_1 \sim (a\fm/\fv) Y$. Thus
$$\text{$(P,\fm,\hat f)$ is strictly normal}\quad\Longleftrightarrow\quad
\fv\prec^\flat 1 \text{ and } u \prec_{\Delta(\fv)} \fm\fv.$$
We say that $(P,\fm,\hat f)$ is {\bf balanced}\index{slot!balanced} if  $(P,\fm,\hat f)$ is steep and $P(0) \preceq S_{P_{\times\fm}}(0)$, equivalently, $(P,\fm,\hat f)$ is steep and $u \preceq \fm$.
Thus 
$$\text{$(P,\fm, \hat f)$ is strictly normal}\ \Longrightarrow\ \text{$(P,\fm, \hat f)$ is balanced}\ \Longrightarrow\ 
\text{$(P,\fm, \hat f)$ is deep},$$ 
and with $b\in K^\times$,
$$(P,\fm, \hat f) \text{ is balanced }\Longleftrightarrow (P_{\times\fn},\fm/\fn,\hat f/\fn) \text{ is balanced }\Longleftrightarrow (bP,\fm,\hat f) \text{ is balanced}.$$
If $(P,\fm,\hat f)$ is balanced, then so is any slot in $K$ equivalent to~$(P,\fm,\hat f)$.
Moreover, if~$(P,\fm,\hat f)$ is a hole in $K$, then $P(0)=-L_P(\hat f)$, so
$(P,\fm,\hat f)$ is balanced iff it is steep and $L_P(\hat f)\preceq S_{P_{\times\fm}}(0)$.
By Corollary~\ref{cor:good approx to hata, deg 1}, if $(P,\fm,\hat f)$ is steep, then  $\hat f-f\prec_{\Delta(\fv)}\fm$ for some $f$.
For balanced $(P,\fm,\hat f)$ we have a variant of this fact:

\begin{lemma}\label{lem:balanced good approx}
Suppose $(P,\fm,\hat f)$ is balanced and $g\notin \I(K)+K^\dagger$. Then there is for all~$n$ an $f$ such that $\hat f-f\prec\fv^n\fm$.
\end{lemma}
\begin{proof}
Replacing $(P,\fm,\hat f)$ by an equivalent hole in $K$, we arrange that $(P,\fm,\hat f)$ is a hole in $K$, and replacing
$(P,\fm,\hat f)$  by $(P_{\times\fm},1,\hat f/\fm)$, that~$\fm=1$. Then~${\hat f}'+g\hat f=u$ with $g=1/\fv\succ^\flat 1$, $g\notin \I(K)+K^\dagger$, and $u\preceq 1$. Hence the lemma follows from Corollary~\ref{cor:uaintfex, 1}.
\end{proof}

\noindent 
In the next corollary we assume that the subgroup $K^\dagger$ of $K$ is divisible.  (Since $K$ is henselian and $\d$-valued, this holds if the groups~$C^\times$ and $\Gamma$ are divisible.)

\begin{cor}\label{cor:balanced -> strongly normal}
Suppose $(P,\fm,\hat f)$ is balanced and $g\notin \I(K)+K^\dagger$. Then $(P,\fm,\hat f)$ has a strictly normal refinement $(P_{+f},\fm,\hat f-f)$.
\end{cor}
\begin{proof}
First arrange that $(P,\fm,\hat f)$ is a hole in $K$.
The previous lemma yields an~$f$ such that $\hat f-f\preceq\fv^3\fm$.
Then $(P_{+f},\fm,\hat f-f)$ is a strictly normal refinement of~$(P,\fm,\hat f)$, by Lemma~\ref{lem:achieve strong normality} (where the latter uses divisibility of $K^\dagger$).
\end{proof}

%The next lemma and its proof have been checked, but are commented out, since it doesn't seem to be used

\begin{lemma}\label{lem:balanced refinement}
Suppose $(P,\fm,\hat f)$ is  balanced with  $v\hat f\notin \exc^{\ev}(L_P)$ and~$\hat f - f \preceq \hat f$. Then the refinement $(P_{+f},\fm,\hat f-f)$ of $(P,\fm,\hat f)$
is balanced.
\end{lemma}
%\begin{proof}
%We first arrange~$\fm=1$, by replacing $(P,\fm,\hat f)$, $f$ by $(P_{\times\fm},1,\hat f/\fm)$, $f/\fm$, respectively. 
%By the remark after Lemma~\ref{lem:steep1},  $(P_{+f},1,\hat f-f)$ is steep. 
%Take~$\phi$ with $v\hat f\notin\exc\big((L_P)^\phi\big)$, and set $\hat g:=\hat f-f$, so $0\neq\hat g\preceq\hat f$.
%From [ADH, 5.7.5] recall~$L_{P^\phi}=(L_P)^\phi$ and hence
%$L_{P^\phi}(\hat f)=L_P(\hat f)$ and $L_P(\hat g)=L_{P^\phi}(\hat g)$.
%Thus
%$$L_{P_{+f}}(\hat g)=L_P(\hat g) \preceq L_{P^\phi}\hat g \preceq L_{P^\phi}\hat f \asymp L_{P^\phi}(\hat f)   \preceq S_P(0)=S_{P_{+f}}(0),$$
%using [ADH, 4.5.1(iii)] to get the second $\preceq$ and $v\hat f \notin \exc(L_{P^\phi})$  to get   $\asymp$.
%Therefore~$(P_{+f},\fm,\hat g)$ is balanced.
%\end{proof}

\begin{proof} By Lemma~\ref{lem:from cracks to holes}  we arrange $(P,\fm,\hat f)$ is a hole. 
Replacing $(P,\fm,\hat f)$ and $f$ by $(P_{\times\fm},1,\hat f/\fm)$ and $f/\fm$ we arrange next that~$\fm=1$. 
By the remark preceding Lemma~\ref{lem:achieve steep},  $(P_{+f},1,\hat f-f)$ is steep.  
Take~$\phi$ such that $v\hat f\notin\exc\big((L_P)^\phi\big)$, and set~$\hat g:=\hat f-f$, so $0\neq\hat g\preceq\hat f$.
Recall from [ADH, 5.7.5] that~$L_{P^\phi}=(L_P)^\phi$ and hence
$L_{P^\phi}(\hat f)=L_P(\hat f)$ and $L_P(\hat g)=L_{P^\phi}(\hat g)$.
Thus
$$L_{P_{+f}}(\hat g)\ =\ L_P(\hat g)\  \preceq\  L_{P^\phi}\hat g\  \preceq\ L_{P^\phi}\hat f\ \asymp\ L_{P^\phi}(\hat f)\  =L_P(\hat f)\    \preceq\ S_P(0)\ =\ S_{P_{+f}}(0),$$
using [ADH, 4.5.1(iii)] to get the second $\preceq$ and $v\hat f \notin \exc(L_{P^\phi})$  to get   $\asymp$; the last $\preceq$ uses $(P,1,\hat f)$ being a hole.
Therefore~$(P_{+f},1,\hat g)$ is balanced.
\end{proof} 

%Moreover,
%$$L_{P_{+f}}(\hat f-f)\ =\ L_P(\hat f-f)\ \preceq\  L_P(\hat f)\ \preceq\ S_P(0)\ =\ S_{P_{+f}}(0),$$
%using Proposition~\ref{prop:slow} and $\exc^{\ev}(L_P) < v\hat f \leq v(\hat f-f)$ for the first $\preceq$ and $(P,1,\hat f)$ being a hole for the second. Therefore~$(P_{+f},1,\hat f-f)$ is balanced.
%\end{proof}

% Version for $$v\hat f\notin \exc^{\ev}(L_P)< v\hat f$ (not checked):  Set $\hat g:=\hat f-f$. Then $0\neq\hat g\preceq\hat f$, so
%$$L_{P_{+f}}(\hat g)\ =\ L_P(\hat g)\  \preceq\  L_P\hat g\  \preceq\ L_P\hat f\ \asymp(???) \ L_P(\hat f)\ \preceq\ %S_P(0)=S_{P_{+f}}(0),$$
%using $v\hat f \notin \exc^{\ev}(L_P)$  for the second $\preceq$ . 

\noindent
Combining Lemmas~\ref{lem:isolated refinement} and~\ref{lem:balanced refinement} yields:

\begin{cor}\label{cor:balanced refinement}
If $(P,\fm,\hat f)$ is balanced and isolated, and $\hat f -f \preceq \hat f$, then the refinement~$(P_{+f},\fm,\hat f-f)$ of $(P,\fm,\hat f)$ is also balanced and isolated. 
\end{cor}
 
\noindent
We call  $(P,\fm,\hat f)$ {\bf proper}\index{slot!proper}\index{proper!slot} if the differential polynomial $P$ is proper as defined in Section~\ref{sec:complements newton} (that is, $u\neq 0$ and $g+u^\dagger\succ^\flat 1$). If~$(P,\fm,\hat f)$ is proper, then so are~$(bP,\fm,\hat f)$ for $b\neq 0$ and $(P_{\times\fn},\fm/\fn,\hat f/\fn)$, as well as each refinement $(P,\fn,\hat f)$ of~$(P,\fm,\hat f)$ and each slot in $K$ equivalent to $(P,\fm,\hat f)$. By Lemma~\ref{lem:proper compconj},
if $(P,\fm,\hat f)$ is proper, then so is~$(P^\phi,\fm,\hat f)$ for $\phi\preceq 1$.
%[{\bf skipped}If $u\neq 0$, $u^\dagger\preceq 1$,   then~$(P,\fm,\hat f)$ is proper iff it is steep.]

\begin{lemma}\label{lem:proper => balanced}
Suppose $(P,\fm,\hat f)$ is proper and $\fm\asymp u$; then $(P,\fm,\hat f)$ is balanced.
\end{lemma}
\begin{proof}
Replacing $(P,\fm,\hat f)$ by $(P_{\times\fm},1,\hat f/\fm)$, we arrange $\fm=1$.
Then $u\asymp 1$ and thus $(P,1,\hat f)$ is balanced.
\end{proof}

\begin{prop}\label{prop:balanced refinement}
Suppose $(P,\fm,\hat f)$ is proper and $v\hat f\notin\exc^{\ev}(L_P)$. Then~$(P,\fm,\hat f)$ has a balanced re\-fine\-ment. %~$(P_{+f},\fn,\hat f-f)$.  
\end{prop}

\begin{proof}
We arrange $\fm=1$ as usual.
By Lemmas~\ref{lem:proper asymp} and~\ref{lem:from cracks to holes} we have 
$$\hat f\  \sim\  u/(g+u^\dagger) \prec^\flat u.$$ Hence if $u\preceq 1$, then
$(P,u,\hat f)$ refines $(P,1,\hat f)$, and so $(P,u,\hat f)$ is balanced by Lem\-ma~\ref{lem:proper => balanced}.
Assume now that $u\succ 1$. Then $1\prec u\prec g$ by Lemma~\ref{lem:proper nmul 1} and~$\nval P=1$, and hence $u^\dagger\preceq g^\dagger \prec  g$. So~$g\sim g+u^\dagger \succ^\flat 1$, hence~$(P,1,\hat f)$ is steep, and~$\hat f\sim u/g$. Set $f:=u/g\prec 1$; then~$(P_{+f},1,\hat f-f)$ is a steep refinement of~$(P,1,\hat f)$. Moreover
$$P_{+f}(0) =P(f) = af'\prec a=S_{P_{+f}}(0),$$
hence  $(P_{+f},1,\hat f-f)$ is balanced.
\end{proof}

\begin{cor}\label{nonamecor}
Suppose $K$ is $\upl$-free.  Then there exists $\phi\preceq 1$ and  a re\-fine\-ment~$(P_{+f},\fn,{\hat f-f})$ of $(P,\fm,\hat f)$ such that $(P^\phi_{+f},\fn,\hat f-f)$ is balanced.
\end{cor}
\begin{proof}  
Using Remark~\ref{rem:achieve isolated} we can replace~$(P,\fm,\hat f)$ by a refinement to arrange that~$(P,\fm,\hat f)$ is isolated. Then $u\neq 0$
by the remark before Lemma~\ref{lem:excev empty}, so by Lemma~\ref{lem:proper evt} ,  $P^\phi$ is proper, eventually.
Now apply  Proposition~\ref{prop:balanced refinement} to a proper (and isolated) $(P^\phi,\fm,\hat f)$ with $\phi\preceq 1$.
\end{proof}

\noindent
%The last corollary together with Lemmas~\ref{lem:normality comp conj, strong}, \ref{lem:excev empty}, and %Corollary~\ref{cor:balanced -> strongly normal} yields:
 
\begin{cor}\label{cor:strongly normal d=1}
Suppose $K$ is $\upl$-free, $\der K=K$,~$\I(K)\subseteq K^\dagger$, and $K^\dagger$ is divisible. Then $(P,\fm,\hat f)$ has a
refinement~$(P_{+f},\fn,\hat f-f)$ such that  $(P_{+f}^\phi,\fn,\hat f -f)$ is strictly normal for some $\phi\preceq 1$.
\end{cor}
\begin{proof}  Corollary~\ref{nonamecor} yields a refinement $(P_{+f_1}, \fn_1, \hat f-f_1)$ of $(P, \fm,\hat f)$ and a~$\phi\preceq 1$ such that
$(P^\phi_{+f_1},\fn_1,\hat f-f_1)$ is balanced. By  Lemma~\ref{lem:excev empty} with $K^\phi$ in the role of~$K$ and
$(P^\phi_{+f_1}, \fn_1, \hat f-f_1)$ in the role of $(P,\fm, \hat f)$ we can apply Corollary~\ref{cor:balanced -> strongly normal} to~$(P^\phi_{+f_1}, \fn_1, \hat f-f_1)$ to give a strictly normal refinement $(P^{\phi}_{f_1+f_2}, \fn, \hat f-f_1-f_2)$ of it.
Thus for $f:=f_1+f_2$ the refinement $(P_{+f},\fn, \hat f-f)$ of $(P,\fm,\hat f)$ has the property that $(P_{+f}^\phi,\fn,\hat f -f)$ is strictly normal. 
\end{proof}

\noindent
Combining this corollary with Corollaries~\ref{cor:minhole deg 1},~\ref{cor:achieve strong normality, 3}, and Lemma~\ref{lem:normality comp conj, strong} yields:

\begin{cor}\label{cor:achieve strong normality, general}
If $K$ is $\upo$-free and algebraically closed with $\der K=K$ and ${\I(K)\subseteq K^\dagger}$, then every minimal hole in $K$ of order $\ge 1$ has
a re\-fine\-ment~$(Q,\fn,\hat g)$ such that~$(Q^\phi,\fn,\hat g)$ is deep and strictly normal, eventually.  
\end{cor} 

\begin{remark}
Suppose $K$ is $\upl$-free,  with $\der K=K$, $\I(K)\subseteq K^\dagger$, and $K^\dagger$ is divisible.
By Corollary~\ref{mainthm deg 1} every linear slot in $K$ of order $r\geq 1$ has a refinement 
$(Q,\fn,\hat g)$ such that $(Q^\phi,\fn,\hat g)$ is deep and normal, eventually.
We don't know whether every linear minimal hole in $K$ of order~$r\geq 1$ has a refinement
$(Q,\fn,\hat g)$ such that $(Q^\phi,\fn,\hat g)$ is deep and strictly normal, eventually.
(For $r=1$ this holds by Corollary~\ref{cor:strongly normal d=1}.)
\end{remark}

\newpage 

\part{Slots in $H$-Fields}\label{part:dents in H-fields}

\medskip

\noindent
Here we specialize to the case that $K$ is the algebraic closure of a Liouville closed $H$-field $H$ with small derivation.
After the preliminary Sections~\ref{sec:aux} and~\ref{sec:approx linear diff ops} we come in Sections~\ref{sec:split-normal holes}--\ref{sec:repulsive-normal} to 
the technical heart of Part 4.
  Section~\ref{sec:split-normal holes} shows that every minimal hole in $K$ gives rise to 
a {\it split-normal}\/ slot~$(Q,\fn,\hat b)$ in~$H$:
a normal slot  in~$H$
whose linear part~$L_{Q_{\times\fn}}\in H[\der]$ ``asymptotically'' splits over $K$; see
Definition~\ref{SN} for the precise definition, and Theorem~\ref{thm:split-normal} for the detailed statement of the main result of this section.
In the intended setting where~$H$ is a Hardy field, this asymptotic splitting will allow us to define in Part~\ref{part:Hardy fields}
a contractive operator on a space of real-valued functions; this operator then has a fixed point whose germ $y$ satisfies $Q(y)=0$, $y\prec\fn$. 
A main difficulty in that part will lie in guaranteeing that  such germs $y$ have similar asymptotic properties as~$\hat b$.  
Sections~\ref{sec:ultimate} and~\ref{sec:repulsive-normal} prepare the ground for dealing with this:
In Section~\ref{sec:ultimate} we strengthen the concept of isolated slot  to {\it ultimate}\/ slot (in $H$, or in $K$).
This relies on the ultimate exceptional values of linear differential operators over $K$ introduced in Part~\ref{part:universal exp ext}. 
In Section~\ref{sec:repulsive-normal} we single out among
split-normal slots in $H$ those that are {\it repulsive-normal}, 
culminating in the proof of Theorem~\ref{thm:repulsive-normal}: an analogue of Theorem~\ref{thm:split-normal} producing
repulsive-normal   ultimate slots in~$H$ from minimal holes in $K$.

\section{Some Valuation-Theoretic Lemmas}\label{sec:aux}

\noindent
%In the next section we prove Lemma~2.8 from Joris' notes.
The present section contains preliminaries for the next section on approximating  splittings of linear differential operators; these facts
in turn will be used in  Section~\ref{sec:split-normal holes} on split-normality. 
We shall often deal with real closed fields with extra structure, denoted usually by $H$, since the results in this section about such $H$ will later be applied to $H$-fields and Hardy fields. 
We begin by summarizing some purely valuation-theoretic facts.

\subsection*{Completion and specialization of real closed valued fields}
Let $H$ be a real closed valued field whose valuation ring $\mathcal O$ is convex in $H$ (with respect to the unique ordering on $H$ making~$H$ an ordered field). Using [ADH, 3.5.15] we equip the algebraic closure $K=H[\imag]$ \textup{(}$\imag^2=-1$\textup{)} of~$H$ with its unique valuation ring lying over~$\mathcal O$, which is $\mathcal O + \mathcal O\imag$. We set $\Gamma:= v(H^\times)$, so $\Gamma_K=\Gamma$. 

\begin{lemma}\label{lem:Kc real closed}
The completion $H^{\operatorname{c}}$ of the valued field $H$ is real closed, its valuation ring is convex in $H^{\operatorname{c}}$, and there is a unique valued field embedding 
$H^{\operatorname{c}}\to K^{\operatorname{c}}$ over~$H$. Identifying $H^{\operatorname{c}}$ with its image under this embedding we have $H^{\operatorname{c}}[\imag]=K^{\operatorname{c}}$.
\end{lemma}
\begin{proof} For the first two claims, see [ADH, 3.5.20].  By [ADH, 3.2.20] we have a unique valued field embedding 
$H^{\operatorname{c}}\to K^{\operatorname{c}}$ over~$H$, and viewing $H^{\operatorname{c}}$ as a valued subfield of~$K^{\operatorname{c}}$ via this embedding we have $K^{\operatorname{c}}=H^{\operatorname{c}}K=H^{\operatorname{c}}[\imag]$ by [ADH, 3.2.29].
\end{proof}

\noindent
We identify $H^{\operatorname{c}}$ with its image in $K^{\operatorname{c}}$ as in the previous lemma.
Fix a  convex subgroup $\Delta$ of $\Gamma$. Let $\dot{\mathcal O}$ be the valuation ring of the coarsening
of $H$ by $\Delta$, with maximal ideal $\dot{\smallo}$. Then by [ADH, 3.5.11 and subsequent remarks] $\dot{\mathcal O}$ and $\dot{\smallo}$ are convex in $H$, the specialization $\dot H=\dot{\mathcal O}/\dot{\smallo}$ of~$H$ by $\Delta$ is naturally an ordered and valued field, and the valuation ring of $\dot{H}$ is convex in $\dot{H}$. Moreover,
$\dot{H}$ is even real closed by [ADH, 3.5.16]. 
 Likewise, the coarsening of $K$ by $\Delta$ has valuation ring 
$\dot{\mathcal O}_K$ with maximal ideal $\dot{\smallo}_K$
and valued residue field $\dot K$. Thus $\dot{\mathcal O}_K$ lies over $\dot{\mathcal O}$ by [ADH,~3.4, subsection {\it Coarsening and valued field extensions}\/], so
$(K, \dot{\mathcal O}_K)$ is a valued field extension of
$(H, \dot{\mathcal{O}})$. In addition:  

\begin{lemma}\label{lem:dotK real closed}
$\dot K$ is a valued field extension of $\dot H$ and an algebraic closure of $\dot H$.
\end{lemma}
\begin{proof} The second part follows by general valuation theory from $K$ being an algebraic closure of $H$. In fact, with the image of $\imag\in \mathcal{O}_K\subseteq \dot{\mathcal O}_K$
in $\dot{K}$ denoted by the same symbol, we have
$\dot{K}=\dot{H}[\imag]$.
\end{proof}

\noindent
Next, let $\hat H$ be an immediate valued field extension of $H$. We equip $\hat H$ with the unique field ordering making it an ordered field extension of $H$ in which $\mathcal O_{\hat H}$ is convex; see [ADH, 3.5.12].
Choose $\imag$ in a field extension of $\hat H$ with $\imag^2=-1$.  Equip~$\hat H[\imag]$ with the unique valuation ring of $\hat H[\imag]$ that lies over $\mathcal O_{\hat H}$, namely~${\mathcal O_{\hat H} + \mathcal{O}_{\hat H}\imag}$  [ADH, 3.5.15]. Let  $\hat a =\hat b+\hat c\,\imag\in \hat H[\imag]\setminus H[\imag]$ with  $\hat b,\hat c\in \hat H$, and let~$b$,~$c$ range over $H$.
Then
$$v\big(\hat a-(b+c\imag)\big)\ =\ \min\!\big\{ v(\hat b-b),v(\hat c-c) \big\}$$
and thus $v\big( \hat a - H[\imag] \big) \subseteq v(\hat b-H)$ and $v\big( \hat a - H[\imag] \big) \subseteq v(\hat c-H)$.

\begin{lemma}\label{lem:same width}
We have $v(\hat b-H)\subseteq v(\hat c-H)$ or $v(\hat c-H)\subseteq v(\hat b-H)$. Moreover, the following are equivalent:
\begin{enumerate}
\item[$\mathrm{(i)}$] $v(\hat b-H)\subseteq v(\hat c-H)$;
\item[$\mathrm{(ii)}$] for all $b$ there is a $c$ with $v\big(\hat a-(b+c\imag)\big)=v(\hat b-b)$;
\item[$\mathrm{(iii)}$]  $v\big( \hat a - H[\imag] \big) = v(\hat b-H)$.
\end{enumerate}
%\item if $v(\hat c-K)\subseteq v(\hat b-K)$ then for each $c$ there is some $b$ with $v\big(\hat a-(b+c\imag)\big)=v(\hat c-c)$, and thus $v\big( \hat a - K[\imag] \big) = v(\hat c-K)$.
%\end{enumerate}
\end{lemma}
\begin{proof} For the first assertion, use that $v(\hat b-H), v(\hat c-H)\subseteq \Gamma_\infty$ are downward closed.  Suppose $v(\hat b-H)\subseteq v(\hat c-H)$, and let $b$ be given. 
If $\hat c\in H$, then for $c:=\hat c$ we have $v\big(\hat a-(b+c\imag)\big)=v(\hat b-b)$. Suppose $\hat c\notin H$. Then $v(\hat c-H)\subseteq\Gamma$ does not have a largest element and
$v(\hat b-b)\in v(\hat c-H)$, so we have~$c$ with $v(\hat b-b)<v(\hat c-c)$; thus 
$$v\big(\hat a-(b+c\imag)\big)=\min\!\big\{v(\hat b-b),v(\hat c-c)\big\}=v(\hat b-b).$$ 
This shows~(i)~$\Rightarrow$~(ii). Moreover,
(ii)~$\Rightarrow$~(iii) follows from $v\big( \hat a - H[\imag] \big) \subseteq v(\hat b-H)$, and (iii)~$\Rightarrow$~(i) from~$v\big( \hat a - H[\imag] \big) \subseteq v(\hat c-H)$. 
%the remark preceding the lemma.
\end{proof} 

\noindent
So if $v(\hat b-H)\subseteq v(\hat c-H)$, then:  $\hat a$ is special over $H[\imag]\ \Longleftrightarrow\ \hat b$ is special over $H$. 
 
\medskip
\noindent
To apply Lemma~\ref{lem:same width} to $H$-fields we assume in the next lemma more generally that~$H$ is equipped with a derivation making it a $\d$-valued field and that $\hat H$ is equipped with a derivation~$\der$ making it an asymptotic field extension of $H$; then $\hat H$ is also $\d$-valued
with the same constant field as $H$ [ADH, 9.1.2].
 
\begin{lemma}\label{lem:same width, der}
Suppose $H$ is closed under integration. Then we have:
$$v(\hat b-H)\subseteq v(\hat c-H)\ \Longrightarrow\ v(\der\hat b-H)\subseteq v(\der\hat c-H).$$
\end{lemma}
\begin{proof} Assume $v(\hat b-H)\subseteq v(\hat c-H)$.
Let $b\in H$, and take $g\in H$ with $g'=b$;    adding a suitable constant to $g$ we   arrange $\hat b-g\nasymp 1$.
Next, take $h\in H$ with $\hat b-g\asymp \hat c-h$. Then 
$$\der\hat b-b\ =\ \der(\hat b-g)\ \asymp\ \der(\hat c-h)\ =\ \der\hat c -h',$$
so $v(\der\hat b - b)\in v(\der\hat c -H)$.
\end{proof}

\subsection*{Embedding into the completion}
In this subsection $K$ is an asymptotic field, $\Gamma:= v(K^\times)\neq\{0\}$, and 
$L$ is an asymptotic field extension of $K$ such that $\Gamma$ is cofinal in $\Gamma_L$.

\begin{lemma}\label{cdercon}
Let $a\in L$ and let $(a_\rho)$ be a c-sequence in $K$ with $a_\rho\to a$ in $L$. Then for each $n$, $(a_\rho^{(n)})$ is a c-sequence in $K$ with $a_\rho^{(n)}\to a^{(n)}$ in $L$. 
\end{lemma}
\begin{proof}
By induction on $n$ it suffices to treat the case $n=1$. Let
$\gamma\in\Gamma_L$; we need to show the existence of an index $\sigma$ such that $v(a'-a_\rho')>\gamma$ for all $\rho>\sigma$.
By [ADH, 9.2.6]   we have $f\in L^\times$ with $f\prec 1$ and $v(f') \geq \gamma$. 
Take $\sigma$ such that $v(a-a_\rho) > vf$ for all $\rho>\sigma$.
Then $v(a'-a_\rho') > v(f') \geq \gamma$ for $\rho>\sigma$. \end{proof}

\noindent
Let  $K^{\operatorname{c}}$ be the completion of the valued differential field $K$; then 
$K^{\operatorname{c}}$ is   asymptotic   by [ADH, 9.1.6].  Lemma~\ref{cdercon} and [ADH, 3.2.13 and 3.2.15] give:

\begin{cor}\label{klemb}
Let $(a_i)_{i\in I}$ be a family of elements of $L$ such that  $a_i$  is the limit in~$L$
of a c-sequence in $K$, for each $i\in I$. 
Then there is a unique embedding~$K\big\<(a_i)_{i\in I}\big\>\to K^{\operatorname{c}}$ of valued differential fields over $K$.
\end{cor}

\noindent
Next suppose that $H$ is a real closed asymptotic field whose
valuation ring $\mathcal{O}$ is convex in $H$ with $\mathcal{O}\ne H$, the asymptotic extension  $\hat H$ of $H$
is immediate, and~$\imag$ is an element of an asymptotic
extension of~$\hat H$ with $\imag^2=-1$.  Then $\imag\notin \hat H$, and we identify~$H^{\operatorname{c}}$ with a valued subfield of $H[\imag]^{\operatorname{c}}$
as in Lemma~\ref{lem:Kc real closed}, so that
$H^{\operatorname{c}}[\imag]=H[\imag]^{\operatorname{c}}$  as in that lemma.  Using also  Lemma~\ref{cdercon}  we see that
$H^{\operatorname{c}}$ is actually a valued {\em differential\/} subfield of the asymptotic field $H[\imag]^{\operatorname{c}}$, and so $H^{\operatorname{c}}[\imag]=H[\imag]^{\operatorname{c}}$ also as {\em asymptotic\/} fields.   Thus by Corollary~\ref{klemb} applied to $K:=H$ and
$L:=\hat{H}$:

\begin{cor}\label{cor:embed into Kc}
Let $a\in \hat{H}[\imag]$ be the limit in $\hat{H}[\imag]$ of a c-sequence in~$H[\imag]$. Then~$\Re a$,~$\Im a$ are limits in $\hat{H}$ of c-sequences in $H$, hence
there is a unique embedding
$H[\imag]\big\<\!\Re a,\Im a\big\>\to H^{\operatorname{c}}[\imag]$ of valued differential fields over $H[\imag]$.
\end{cor}

\section{Approximating Linear Differential Operators}\label{sec:approx linear diff ops}

\noindent
{\em In this section $K$ is a valued differential field with small derivation,}\/  $\Gamma:= v(K^\times)$.
For later use we prove here Corollaries~\ref{cor:approx LP+f} and~\ref{cor:approx LP+f, real, general} and consider {\em strong splitting}\/.\index{splitting!strong}\index{linear differential operator!strong splitting} %and   {\em Hardy type}\/ asymptotic fields. 
Much of this section rests on the following basic estimate for linear differential operators which split over $K$:

\begin{lemma}\label{lem:approxB}
Let   $b_1,\dots,b_r\in K$ and $n$ be given.
Then there exists $\gamma_0\in\Gamma^{\geq}$ such that  for all $b_1^{\smallbullet},\dots,b_r^{\smallbullet}\in K$ and $\gamma\in\Gamma$ with $\gamma>\gamma_0$ and $v(b_i-b_i^{\smallbullet})\geq (n+r)\gamma$ for~$i=1,\dots,r$, we have $v(B-B^{\smallbullet})  \geq vB+n\gamma$, where
$$B\ :=\ (\der-b_1)\cdots(\der-b_r)\in K[\der],\quad B^{\smallbullet}\ :=\ (\der-b_1^{\smallbullet})\cdots(\der-b_r^{\smallbullet})\in K[\der].$$ 
\end{lemma}
\begin{proof}
By induction on $r\in\N$. The case $r=0$ is clear (any $\gamma_0\in\Gamma^{\geq}$ works).  Suppose the lemma holds for a certain $r$.
Let $b_1,\dots,b_{r+1}\in K$  and $n$ be given. Set~$\beta_i:=vb_i$ ($i=1,\dots,r+1$). Take~$\gamma_0$ as in the lemma applied to
$b_1,\dots,b_r$ and~$n+1$ in place of $n$, and
let $\gamma_1:=\gamma_0$ if $b_{r+1}=0$, $\gamma_1:=\max\big\{\gamma_0,\abs{\beta_{r+1}}\big\}$ otherwise. 
Let $b_1^{\smallbullet},\dots,b_{r+1}^{\smallbullet}\in K$ and $\gamma\in\Gamma$ with
$\gamma>\gamma_1$ and $v(b_i-b_i^{\smallbullet})\geq (n+r+1)\gamma$ for $i=1,\dots,r+1$.
Set 
$$B\ :=\ (\der-b_1)\cdots(\der-b_r),\qquad B^{\smallbullet}\ :=\ (\der-b_1^{\smallbullet})\cdots(\der-b_r^{\smallbullet}),\qquad E\ :=\ B-B^{\smallbullet}.$$
Then
$$B(\der-b_{r+1})\  =\ B^{\smallbullet}(\der-b^{\smallbullet}_{r+1})+ 
 B^{\smallbullet}(b^{\smallbullet}_{r+1}-b_{r+1})+E(\der-b_{r+1}).$$
Inductively we have  $vE \geq vB+(n+1)\gamma$. Suppose $E\ne 0$ and $0\neq b_{r+1}\nasymp 1$. Then by
[ADH, 6.1.5],  
\begin{align*}
v_E(\beta_{r+1})-v_B(\beta_{r+1})\ &=\ vE-vB+o(\beta_{r+1}) \\
												&\geq\ (n+1)\gamma + o(\beta_{r+1})\\
												&\geq\ n\gamma + |\beta_{r+1}|+o(\beta_{r+1})\ >\  n\gamma.
\end{align*}
Hence, using $E(\der-b_{r+1})=E\der - Eb_{r+1}$ and $v(E\der)=v(E)\ne v_E(\beta_{r+1})$,
\begin{align*}
v\big(E(\der-b_{r+1})\big)\ =\ \min\!\big\{vE,v_E(\beta_{r+1})\big\}\ &>\  \min\!\big\{ vB,v_B(\beta_{r+1})\big\}+n\gamma \\
																&=\ v\big(B(\der-b_{r+1})\big)+n\gamma,
\end{align*}
where for the last equality we use $vB\ne v_B(\beta_{r+1})$. Also,
$$v\big( B^{\smallbullet}(b^{\smallbullet}_{r+1}-b_{r+1}) \big) = v_{B^{\smallbullet}}\big(v(b^{\smallbullet}_{r+1}-b_{r+1})\big) \geq v_{B^{\smallbullet}}\big((n+r+1)\gamma\big)=v_B\big((n+r+1)\gamma\big)$$
where we use [ADH, 6.1.7] for the last equality. Moreover,    by [ADH, 6.1.4],
$$v_B\big((n+r+1)\gamma\big)-n\gamma\ \geq\ vB+(r+1)\gamma+o(\gamma)\ >\ vB\ 	\geq\   v\big(B(\der-b_{r+1})\big).$$
This yields the desired result for $E\ne 0$, $0\ne b_{r+1}\nasymp 1$. The cases $E\ne 0$, $b_{r+1}=0$ and $E=0$, $0\ne b_{r+1}\nasymp 1$ are simpler versions of the above, and so is the case~$E\ne 0$, $b_{r+1}\asymp 1$ using
[ADH, 5.6.1(i)]. The remaining cases, $E=0$, $b_{r+1}=0$ and $E=0$, $b_{r+1}\asymp 1$, are even simpler to handle. 
\end{proof}

\begin{cor}\label{cor:approx A}
Let   $a,b_1,\dots,b_r\in K$, $a\neq 0$.
Then there exists $\gamma_0\in\Gamma^{\geq}$ such that  for all $a^{\smallbullet},b_1^{\smallbullet},\dots,b_r^{\smallbullet}\in K$ and $\gamma\in\Gamma$ with $\gamma>\gamma_0$, $v(a-a^{\smallbullet})\geq va+\gamma$, and~$v(b_i-b_i^{\smallbullet})\geq (r+1)\gamma$ for $i=1,\dots,r$, we have
$v(A-A^{\smallbullet}) \geq vA+\gamma$, where
$$A\ :=\ a(\der-b_1)\cdots(\der-b_r)\in K[\der],\quad A^{\smallbullet}\ :=\ a^{\smallbullet}(\der-b_1^{\smallbullet})\cdots(\der-b_r^{\smallbullet})\in K[\der].$$ 
\end{cor}
\begin{proof}
Take $\gamma_0$ as in the previous lemma applied to $b_1,\dots,b_r$ and $n=1$, and let~$B=(\der-b_1)\cdots(\der-b_r)$, $A=aB$. Let $a^{\smallbullet},b_1^{\smallbullet},\dots,b_r^{\smallbullet}\in K$ and  $\gamma\in\Gamma$ be such that~$\gamma>\gamma_0$, $v(a-a^{\smallbullet})\geq va+\gamma$, and $v(b_i-b_i^{\smallbullet})\geq (r+1)\gamma$ for $i=1,\dots,r$. Set~$B^{\smallbullet}:=(\der-b_1^{\smallbullet})\cdots(\der-b_r^{\smallbullet})$, $A^{\smallbullet}:=a^{\smallbullet} B^{\smallbullet}$.
Then 
$$E\ :=\ A-A^{\smallbullet}\ =\ a(B-B^{\smallbullet})+(a-a^{\smallbullet})B^{\smallbullet}.$$
Lemma~\ref{lem:approxB} gives $vB^{\smallbullet}=vB$, and so {\samepage
$$v\big(a(B-B^{\smallbullet})\big)\ \geq\ va+vB+\gamma\ =\ vA+\gamma, \quad
v\big( (a-a^{\smallbullet})B^{\smallbullet} \big)\ =\ v(a-a^{\smallbullet})+vB\ \geq\ vA+\gamma,$$
so $vE\geq vA+\gamma$.}
\end{proof}

\noindent
{\em In the rest of this subsection we assume $P\in K\{Y\}\setminus K$,
set $r:=\order P$, and let~$\i$,~$\j$ range over $\N^{1+r}$}. 

\begin{lemma} For $\delta:= v\big(P-P(0)\big)$ and all $h\in \smallo$ we have $v\big(P_{+h} - P\big) \geq \delta+\textstyle\frac{1}{2}vh$.
\end{lemma}
\begin{proof} Note that $\delta\in\Gamma$ and $v(P_{\j})\ge \delta$ for all $\j$ with $\abs{\j}\geq 1$.
Let $h\in \smallo^{\ne}$ and~$\i$ be given; we claim that  $v\big( (P_{+h})_{\i} - P_{\i}\big) \geq \delta+\textstyle\frac{1}{2}vh$.
By [ADH, (4.3.1)] we have
$$(P_{+h})_{\i}\ 	=\  P_{\i}+   Q(h)\quad\text{where $Q(Y):=
\sum_{\abs{\j}\geq 1}  {\i+\j \choose \i} P_{\i+\j}\,Y^{\j}\in K\{Y\}$.}$$
From $Q(0)=0$ and [ADH, 6.1.4] we obtain
$$v(Q_{\times h})\ \geq\  v(Q)  + vh+o(vh)\ \geq\ \delta+\textstyle\frac{1}{2}vh.$$
Together with $v\big(Q(h)\big)\ge v(Q_{\times h})$ this yields the lemma.
\end{proof}

\begin{cor}\label{cor:approx P+f}
Let $f\in K$. Then there exists $\delta\in\Gamma$ such that for all $f^{\smallbullet}\in K$ with $f-f^{\smallbullet}\prec 1$  we have
$v\big(P_{+f^{\smallbullet}} - P_{+f}\big) \geq \delta+\textstyle\frac{1}{2}v(f^{\smallbullet}-f)$.
\end{cor}
\begin{proof}
Take  $\delta$   as
in the preceding lemma with  $P_{+f}$  in place of $P$ and $h=f^{\smallbullet}-f$.
\end{proof}

\begin{cor}\label{corcorappr}
Let $a,b_1,\dots,b_r,f\in K$ be such that
$$A\ :=\ L_{P_{+f}}\ =\ a(\der-b_1)\cdots(\der-b_r),\qquad a\ \neq\ 0.$$
Then there exists $\gamma_1\in\Gamma^{\geq}$ such that for all $a^{\smallbullet}, b_1^{\smallbullet},\dots,b_r^{\smallbullet},f^{\smallbullet}\in K$ and $\gamma\in\Gamma$, if 
$$\gamma>\gamma_1,\ v(a-a^{\smallbullet})\geq va +\gamma,\ v(b_i-b_i^{\smallbullet})\geq (r+1)\gamma\ (i=1,\dots,r), \text{ and }
v(f-f^{\smallbullet})\geq 4\gamma,$$ 
then
\begin{enumerate}
\item[$\mathrm{(i)}$] $v\big( P_{+f^{\smallbullet}}  - P_{+f}\big) \geq  vA+\gamma$; and
\item[$\mathrm{(ii)}$] $L_{P_{+f^{\smallbullet}}} = a^{\smallbullet}(\der-b_1^{\smallbullet})\cdots(\der-b_r^{\smallbullet}) + E$  where $vE\geq vA+\gamma$.
\end{enumerate}
\end{cor}
\begin{proof}
Take $\gamma_0$ as in Corollary~\ref{cor:approx A} applied to $a,b_1,\dots,b_r$, and
take  $\delta$   as
in  Corollary~\ref{cor:approx P+f}.
Then $\gamma_1:=\max\{\gamma_0,vA-\delta\}$  has the required property.
\end{proof}

\noindent
In the next result $L$ is a valued differential field extension of $K$ with small derivation such that~$\Gamma$ is cofinal in $\Gamma_L$.
Then the natural inclusion~$K\to L$ extends uniquely to an embedding~$K^{\operatorname{c}}\to L^{\operatorname{c}}$ 
of valued fields by [ADH, 3.2.20]. It is easy to check that this is even an embedding of valued {\em differential\/} fields; we identify~$K^{\operatorname{c}}$ with a valued differential subfield of $L^{\operatorname{c}}$ via this embedding.   

\begin{cor}\label{cor:approx LP+f}
Let $a,b_1,\dots,b_r\in L^{\operatorname{c}}$ and $f\in K^{\operatorname{c}}$ be such that in $L^{\operatorname{c}}[\der]$, 
$$A\ :=\ L_{P_{+f}}\ =\ a(\der-b_1)\cdots(\der-b_r),\qquad a,f\neq 0,\quad \fv:=\fv(A)\prec 1,$$
and let $w\in\N$. Then there  are $a^{\smallbullet},b_1^{\smallbullet},\dots,b_r^{\smallbullet}\in L$ and $f^{\smallbullet}\in K$ such that 
$$a^{\smallbullet}\ \sim\ a, \qquad f^{\smallbullet}\ \sim\ f, \qquad A^{\smallbullet}:=L_{P_{+f^{\smallbullet}}}\ \sim\ A,\qquad \order A^{\smallbullet}\ =\ r, \qquad
 \fv(A^{\smallbullet})\ \sim\ \fv,$$
and such that for $\Delta:=\big\{\alpha\in \Gamma_L:\, \alpha=o\big(v(\fv)\big)\big\}$ we have in $L[\der]$,
$$A^{\smallbullet}\ =\ a^{\smallbullet}(\der-b_1^{\smallbullet})\cdots(\der-b_r^{\smallbullet}) + E,\qquad  E\prec_{\Delta} \fv^{w+1} A.$$
\end{cor}
\begin{proof} 
Let $\gamma_1\in \Gamma_L^{\ge}$ be as in Corollary~\ref{corcorappr} applied to $L^{\operatorname{c}}$ in place of $K$, and take~$\gamma_2\in\Gamma$ such that $\gamma_2\geq \max\{\gamma_1,\frac{1}{4}vf\}+vA$ and
$\gamma_2\geq v\big((P_{+f})_{\i} \big)$ for all $\i$ with~$(P_{+f})_{\i}\neq 0$.
Let $\gamma\in\Gamma$ and $\gamma>\gamma_2$. Then 
$\gamma-vA>\gamma_1$. By the
density of~$K$,~$L$ in~$K^{\operatorname{c}}$,~$L^{\operatorname{c}}$, respectively, we can take
%in combination with the previous corollary then yields
$a^{\smallbullet},b_1^{\smallbullet},\dots,b_r^{\smallbullet}\in L$ and $f^{\smallbullet}\in K$ such that
$$v(a-a^{\smallbullet})\ \ge\ va+(\gamma-vA),\qquad v(b_i-b_i^{\smallbullet})\geq (r+1)(\gamma-vA)\ \text{ for $i=1,\dots,r$,}$$ and 
$v(f-f^{\smallbullet})\geq 4(\gamma-vA) > vf$. Then $a^{\smallbullet}\sim a$, $f^{\smallbullet}\sim f$, and by Corollary~\ref{corcorappr},
$$v\big( P_{+f^{\smallbullet}}  - P_{+f}\big)\ \geq\ \gamma, \quad
A^{\smallbullet}:= L_{P_{+f^{\smallbullet}}}\ =\ a^{\smallbullet}(\der-b_1^{\smallbullet})\cdots(\der-b_r^{\smallbullet}) + E, \quad vE\geq \gamma.$$
Hence 
 $(P_{+f^{\smallbullet}})_{\i} \sim (P_{+f})_{\i}$ if $(P_{+f})_{\i}\neq 0$, and
 $v\big((P_{+f^{\smallbullet}})_{\i}\big)> \gamma_2\ge vA$ if $(P_{+f})_{\i}=0$, so~${A^{\smallbullet}\sim A}$, $\order A^{\smallbullet}=r$, and $\fv(A^{\smallbullet})\sim \fv$.
Choosing $\gamma$ so that also $\gamma > v(\fv^{w+1} A)+\Delta$  we achieve in addition that $E\prec_{\Delta} \fv^{w+1} A$. 
\end{proof}

\subsection*{Keeping it real} 
{\em In this subsection $H$ is a real closed $H$-asymptotic field with small derivation whose valuation ring is convex, with 
$\Gamma:= v(H^\times)\ne \{0\}$, and~$K$ is the asymptotic extension $H[\imag]$ of $H$ with $\imag^2=-1$.}\/ Then $H^{\operatorname{c}}$ is real closed and~$H^{\operatorname{c}}[\imag]=K^{\operatorname{c}}$ as valued field extension of $H$ according to Lemma~\ref{lem:Kc real closed}, and as asymptotic field extension of $H$ by the discussion after Corollary~\ref{klemb}.
Using the real splittings from Definition~\ref{def:real splitting} we show here that we can ``preserve the reality of $A$'' in Corollary~\ref{cor:approx LP+f}.

\begin{lemma}\label{lem:approx real} Let $A\in H^{\operatorname{c}}[\der]$ be of order $r\ge 1$ and let $(g_1,\dots, g_r)\in H^{\operatorname{c}}[\imag]^r$
be a real splitting of $A$ over $H^{\operatorname{c}}[\imag]$. Then for every
$\gamma\in \Gamma$ there are $g_1^{\smallbullet},\dots,g_r^{\smallbullet}$
in~$H[\imag]$
such that $v(g_i-g_i^{\smallbullet}) > \gamma$ for $i=1,\dots,r$, 
$$A^{\smallbullet}\ :=\ (\der-g_1^{\smallbullet})\cdots(\der-g_r^{\smallbullet})\in H[\der],$$ 
and~$(g_1^{\smallbullet},\dots, g_r^{\smallbullet})$ is a real splitting of $A^{\smallbullet}$ over $H[\imag]$. 
\end{lemma}  
\begin{proof} We can reduce to the case where $r=1$ or $r=2$.
If $r=1$, then the lemma holds trivially, so suppose $r=2$. Then again the lemma holds trivially if~$g_1, g_2\in H^{\operatorname{c}}$, so we can assume instead that 
$$g_1\ =\ a-b\imag+b^\dagger, \quad g_2\ =\ a+b\imag, \qquad a\in H^{\operatorname{c}},\ b\in (H^{\operatorname{c}})^\times.$$
Let $\gamma\in\Gamma$ be given. The density of $H$ in $H^{\operatorname{c}}$ gives $a^{\smallbullet}\in H$ with $v(a-a^{\smallbullet})\ge\gamma$.
Next, choose $\gamma^{\smallbullet}\in\Gamma$ such that $\gamma^{\smallbullet}\ge \max\{\gamma,vb\}$ and $\alpha' >\gamma$ for all nonzero $\alpha>\gamma^{\smallbullet}-vb$ in~$\Gamma$,
and take $b^{\smallbullet}\in H$ with $v(b-b^{\smallbullet})>\gamma^{\smallbullet}$. Then $v(b-b^{\smallbullet})>\gamma$ and $b\sim b^{\smallbullet}$.
In fact, $b=b^{\smallbullet}(1+\varepsilon)$ where $v\varepsilon+vb=v(b-b^{\smallbullet})> \gamma^{\smallbullet}$ and so
$v\big( (b/b^{\smallbullet})^\dagger \big) = v(\varepsilon')>\gamma$. Set~$g_1^{\smallbullet}:= a^{\smallbullet}-b^{\smallbullet}\imag+{b^{\smallbullet}}^\dagger$ and $g_2^{\smallbullet}:= a^{\smallbullet}+b^{\smallbullet}\imag$. Then
\begin{align*} v(g_1 - g_1^{\smallbullet})\ &=\ 
v\big( a-a^{\smallbullet}+(b/b^{\smallbullet})^\dagger + (b^{\smallbullet}-b)\imag \big)\ >\ \gamma,\quad v(g_2-g_2^{\smallbullet})\ >\ \gamma,\\
( \der- g_1^{\smallbullet}) \cdot( \der-g_2^{\smallbullet})\ &=\  \der^2-\big(2a^{\smallbullet}+ b^{\smallbullet}{}^\dagger\big)\der + \big((-a^{\smallbullet})'+a^{\smallbullet}{}^2 + a^{\smallbullet} b^{\smallbullet}{}^\dagger + b^{\smallbullet}{}^2\big)\ \in\ H[\der].
\end{align*} 
Hence $(g_1^{\smallbullet}, g_2^{\smallbullet})$ is a real splitting of 
$A^{\smallbullet}:=(\der-g_1^{\smallbullet})(\der-g_2^{\smallbullet})\in H[\der]$.
\end{proof}

\noindent 
In the next two corollaries $a\in (H^{\operatorname{c}})^{\times}$ and
$b_1,\dots,b_r\in K^{\operatorname{c}}$ are such that 
$$A\ :=\ a(\der - b_1)\cdots(\der-b_r)\in H^{\operatorname{c}}[\der],$$
$(b_1,\dots,b_r)$ is a real splitting of $A$ over 
$K^{\operatorname{c}}$, and $\fv:=\fv(A)\prec 1$.  We set $\Delta:=\Delta(\fv)$.

\begin{cor}\label{cor:approx LP+f, real}  
Suppose  $A = L_{P_{+f}}$ with $P\in H\{Y\}$
of order $r\ge 1$ and $f$ in~$(H^{\operatorname{c}})^{\times}$. Let $\gamma\in \Gamma$ and $w\in \N$. Then there is $f^{\smallbullet}\in H^\times$ such that $v(f^{\smallbullet}- f) \ge \gamma$, 
\begin{equation}\label{eq:approx LP+f, real}\
f^{\smallbullet}\ \sim\ f, \quad A^{\smallbullet}\ :=\ L_{P_{+f^{\smallbullet}}}\ \sim\ A,\quad \order A^{\smallbullet}\ =\ r, \quad
 \fv(A^{\smallbullet})\ \sim\ \fv,
\end{equation}
and  we have $a^{\smallbullet}\in H^\times$, $b_1^{\smallbullet},\dots, b_r^{\smallbullet}\in K$, and $B^{\smallbullet},E^{\smallbullet}\in H[\der]$ with 
$A^{\smallbullet}= B^{\smallbullet} + E^{\smallbullet}$,  $E^{\smallbullet}\prec_{\Delta} \fv^{w+1} A$, such that
$$B^{\smallbullet}\ =\ a^{\smallbullet}(\der-b_1^{\smallbullet})\cdots(\der-b_r^{\smallbullet}),\qquad
v(a-a^{\smallbullet}),\
v(b_1-b_1^{\smallbullet}),\ \dots\ ,v(b_r-b_r^{\smallbullet})\ \geq\ \gamma,$$
and $(b_1^{\smallbullet},\dots, b_r^{\smallbullet})$ is a real splitting of $B^{\smallbullet}$ over $K$. 
\end{cor}
\begin{proof} We apply Corollary~\ref{cor:approx LP+f} with $H$,~$K$ in the role of $K$,~$L$, and take~$\gamma_1$,~$\gamma_2$ as in the proof of that corollary. We can assume $\gamma>\gamma_2$, so that $\gamma-vA>0$.
The density of $H$ in $H^{\operatorname{c}}$ gives $a^{\smallbullet}\in H$ such that $v(a-a^{\smallbullet})\geq \max\!\big\{va+(\gamma-vA), \gamma\big\}$ (so $a^{\smallbullet}\sim a$),
and Lemma~\ref{lem:approx real} gives $b_1^{\smallbullet},\dots,b_r^{\smallbullet}\in K$ such that
$v(b_i-b_i^{\smallbullet})\geq \max\!\big\{(r+1)(\gamma-vA), \gamma\big\}$ for $i=1,\dots,r$, and 
$(b_1^{\smallbullet},\dots, b_r^{\smallbullet})$ is a real splitting of $$B^{\smallbullet}\ :=\ a^{\smallbullet}(\der-b_1^{\smallbullet})\cdots(\der-b_r^{\smallbullet})\in H[\der]$$ over $K$.
Take $f^{\smallbullet}\in H$ with $v(f-f^{\smallbullet})\geq \max\!\big\{4(\gamma-vA),\gamma\big\}$.
Then \eqref{eq:approx LP+f, real} follows from the proof of Corollary~\ref{cor:approx LP+f}. We can increase $\gamma$ so that $\gamma> v(\fv^{w+1}A)+\Delta$, and then we have
$A^{\smallbullet}-B^{\smallbullet} \prec_\Delta \fv^{w+1} A$.
\end{proof}

\noindent
This result persists after multiplicative conjugation:

\begin{cor}\label{cor:approx LP+f, real, general}
Suppose  $A= L_{P_{+f,\times\fm}}$ with $P\in H\{Y\}$ of order $r\ge 1$, and $f$ in~$(H^{\operatorname{c}})^{\times}$,  $\fm\in H^\times$. Let $\gamma\in\Gamma$, $w\in\N$. Then there is $f^{\smallbullet}\in H^\times$ such that 
$$v(f^{\smallbullet}-f) \geq \gamma,\quad f^{\smallbullet}\ \sim\ f, \quad A^{\smallbullet}\ :=\ L_{P_{+f^{\smallbullet},\times\fm}}\ \sim\ A,\quad 
\order A^{\smallbullet}\ =\ r, \quad
 \fv(A^{\smallbullet})\ \sim\ \fv,$$
and we have $a^{\smallbullet}\in H^\times$, $b_1^{\smallbullet},\dots, b_r^{\smallbullet}\in K$, and $B^{\smallbullet},E^{\smallbullet}\in H[\der]$ with the properties stated in the previous corollary.
\end{cor}
\begin{proof}
Put $Q:=P_{\times\fm}\in H\{Y\}$, $g:=f/\fm\in H^{\operatorname{c}}$; then $Q_{+g}=P_{+f,\times\fm}$. Applying the previous corollary to $Q$, $g$ in place of $P$, $f$ yields
$g^{\smallbullet}\in H^\times$, $a^{\smallbullet}\in H^\times$, and~${b_1^{\smallbullet},\dots,b_r^{\smallbullet}\in K}$  such that $v(g^{\smallbullet}-g) \geq\ \gamma-v\fm$,
$$g^{\smallbullet}\ \sim\ g, \qquad A^{\smallbullet}\ :=\ L_{Q_{+g^{\smallbullet}}}\ \sim\ A,\qquad 
\order A^{\smallbullet}\ =\ r, \qquad
 \fv(A^{\smallbullet})\ \sim\ \fv$$
and
$A^{\smallbullet} = B^{\smallbullet} + E^{\smallbullet}$, with $B^{\smallbullet},E^{\smallbullet}\in H[\der]$, $E^{\smallbullet}\prec_{\Delta} \fv^{w+1} A$,
and   
$$B^{\smallbullet}\ =\ a^{\smallbullet}(\der-b_1^{\smallbullet})\cdots(\der-b_r^{\smallbullet}),\qquad
v(a-a^{\smallbullet}),\
v(b_1-b_1^{\smallbullet}),\ \dots\ ,v(b_r-b_r^{\smallbullet})\ \geq\ \gamma,$$
and $(b_1^{\smallbullet},\dots, b_r^{\smallbullet})$ is a real splitting of $B^{\smallbullet}$ over $K$. 
Therefore  $f^{\smallbullet}:=g^{\smallbullet}\fm\in H^\times$ and~$a^{\smallbullet},b_1^{\smallbullet},\dots,b_r^{\smallbullet}$ have the required properties.
\end{proof}

\subsection*{Strong splitting} 
{\em In this subsection $H$ is a real closed $H$-field with small derivation and asymptotic integration}. Thus $K:=H[\imag]$ is a $\d$-valued extension of $H$.
Let~$A\in K[\der]^{\neq}$ have order $r\ge 1$ and set $\fv:=\fv(A)$, and let  $f$, $g$, $h$ (possibly subscripted) range over $K$. 
Recall from Section~\ref{sec:diff ops and diff polys} that  a splitting of $A$ over $K$ is an $r$-tuple $(g_1,\dots,g_r)$ such that
$$A\ =\ f(\der-g_1)\cdots(\der-g_r)\quad\text{where $f\neq 0$.}$$
We call such a splitting $(g_1,\dots,g_r)$ of $A$ over $K$ {\bf strong} if 
$\Re g_j\succeq \fv^\dagger$ for $j=1,\dots,r$, and
we say that $A$ {\bf  splits strongly over~$K$} if there is a strong splitting of $A$ over~$K$.
This notion is mainly of interest for $\fv\prec 1$, since otherwise $\fv=1$, and then any splitting of 
$A$  over $K$ is a strong splitting of $A$ over $K$. 

\begin{lemma}\label{lem:Ah splits strongly} 
Let $(g_1,\dots,g_r)$ be a strong splitting of $A$ over $K$. If~${h\neq 0}$, then
$(g_1,\dots,g_r)$ is a strong splitting of $hA$ over $K$. If  $h\asymp 1$, then~${(g_1-h^\dagger,\dots,g_r-h^\dagger)}$ is a strong splitting of $Ah$  over $K$.
\end{lemma}
\begin{proof} Suppose $h\asymp 1$. Now use Lemma~\ref{lem:split and twist}, 
 and the fact that if $\fv \prec 1$, then $\Re h^\dagger\preceq h^\dagger\prec \fv^\dagger$. If $\fv=1$, then use that $\fv(Ah)=1$ by Corollary~\ref{cor:111}. 
\end{proof}

\begin{lemma}\label{lem:order 1 splits strongly}
Suppose $g\asymp\Re g$. Then $A=\der-g$ splits strongly over $K$.
\end{lemma}
\begin{proof}
Assuming $\fv\prec 1$ gives $\fv'\prec 1$, so  $\fv^\dagger \prec 1/\fv\asymp g\asymp\Re g$.
\end{proof}

\noindent
In particular, every $A\in H[\der]^{\ne}$ of order~$1$ splits strongly over $K$.

\begin{lemma}\label{lem:split strongly compconj}
Suppose $(g_1,\dots,g_r)$ is a strong splitting of $A$  over $K$ and $\fv\prec^\flat 1$. Let $\phi\preceq 1$ be active in $H$ and set $h_j:=\phi^{-1}\big(g_j-(r-j)\phi^\dagger\big)$   for $j=1,\dots,r$. Then~$(h_1,\dots,h_r)$ is a strong splitting of $A^\phi$ over $K^\phi=H^\phi[\imag]$.
\end{lemma}
\begin{proof}
By Lemma~\ref{lem:split and compconj}, $(h_1,\dots,h_r)$ is a splitting of $A^\phi$ over $K^\phi$. 
We have $\phi^\dagger\prec 1\preceq\fv^\dagger$, 
so $\Re h_j\sim \phi^{-1} \Re g_j\succeq \phi^{-1}\fv^\dagger$ for $j=1,\dots,r$.
Set $\fw:=\fv(A^\phi)$ and~$\derdelta:=\phi^{-1}\der$. 
Lem\-ma~\ref{lem:v(Aphi)} gives
$\fv^\dagger\asymp \fw^\dagger$, so 
$\phi^{-1} \fv^\dagger  \asymp \derdelta(\fw)/\fw$. 
\end{proof}

\noindent
{\em In the next two results we assume that for all $q\in \Q^{>}$ and
$\fn\in H^\times$  there is given an element
$\fn^q\in H^\times$ such that $(\fn^q)^\dagger=q\fn^\dagger$ \textup{(}and thus
$v(\fn^q)=q\,v(\fn)$\textup{)}}. 

\begin{lemma}\label{lem:split strongly multconj} 
Suppose $(g_1,\dots,g_r)$ is a splitting of $A$  over $K$, $\fv\prec 1$, $\fn\in H^\times$, and $[\fv]\le[\fn]$. Then for all $q\in \Q^{>}$ with at most~$r$ exceptions, $(g_1-q\fn^\dagger,\dots,g_r-q\fn^\dagger)$ is a strong splitting of $A\fn^q$   over $K$.  
\end{lemma}
\begin{proof} 
Let $q\in \Q^{>}$. Then
$(g_1-q\fn^\dagger,\dots,g_r-q\fn^\dagger)$ is a splitting of $A\fn^q$ over $K$, by Lemma~\ref{lem:split and twist}. Moreover, $\big[\fv(A\fn^q)\big]\le [\fn]$,  by Lemma~\ref{lem:An}, so
$\fv(A\fn^q)^\dagger \preceq \fn^\dagger$.
Thus if $\Re g_j\not\sim q\fn^\dagger$ for $j=1,\dots,r$, then
$(g_1-q\fn^\dagger,\dots,g_r-q\fn^\dagger)$ is a strong splitting of~$A\fn^q$   over $K$. 
\end{proof}

\begin{cor}\label{cor:split strongly multconj}
Let $(P,\fm,\hat a)$ be a steep slot in $K$ of order $r\ge 1$ whose linear part~$L:=L_{P_{\times\fm}}$ splits over $K$ and such that $\hat a\prec_{\Delta} \fm$ for $\Delta:=\Delta\big(\fv(L)\big)$.  Then for all sufficiently small $q\in \Q^{>}$, any $\fn\asymp|\fv(L)|^q\fm$ in $K^\times$ gives a  steep refinement~$\big(P,\fn,\hat a\big)$ of $(P,\fm,\hat a)$ whose linear part $L_{P_{\times\fn}}$ splits strongly over $K$.
\end{cor}
\begin{proof} Note that $|f|\asymp f$ for all $f$.  Lemma~\ref{lem:steep1} gives $q_0\in \Q^{>}$ such that
for all~$q\in \Q^{>}$ with $q\le q_0$ and any $\fn\asymp|\fv(L)|^q\fm$,  
$(P,\fn, \hat a)$ is a steep refinement of~$(P,\fm, \hat a)$. Now apply Lemma~\ref{lem:split strongly multconj} with $L$, $\fv(L)$, $|\fv(L)|$ in the respective roles
of~$A$,~$\fv$,~$\fn$, and use Lemma~\ref{lem:Ah splits strongly} and the fact that for $\fn\asymp |\fv(L)|^q\fm$ we have~$L_{P_{\times \fn}}=L\cdot\fn/\fm=L|\fv(L)|^qh$ with~$h\asymp 1$. 
 \end{proof} 

\noindent
We finish this section with a useful fact on slots in $K$. 
Given such a slot $(P,\fm,\hat a)$, the element $\hat a$ lies in an immediate asymptotic extension of~$K$ that might not be of the form 
$\hat H[\imag]$ with $\hat H$ an immediate $H$-field extension of $H$. 
By the next lemma we can nevertheless often reduce to this situation, and more: 

\begin{lemma}\label{lem:hole in hat K}
Suppose $H$ is $\upo$-free. Then every $Z$-minimal slot in $K$ of positive order  is equivalent to a hole
$(P,\fm,\hat b)$  in $K$ with $\hat b\in\hat K=\hat H[\imag]$ for some immediate $\upo$-free newtonian $H$-field extension $\hat H$ of $H$.
%If $H$ is $\upl$-free, then the same conclusion holds for every $Z$-minimal dent in~$K$ of order $1$ with a quasilinear refinement.
\end{lemma}
\begin{proof}
Let $(P,\fm,\hat a)$ be a $Z$-minimal slot in $K$ of order $\ge 1$.
Take an immediate $\upo$-free newtonian $H$-field extension $\hat H$ of $H$; such $\hat H$ exists by remarks following~[ADH, 14.0.1]. Then $\hat K=\hat H[\imag]$ is also newtonian by~[ADH, 14.5.7]. 
Now apply  Corollary~\ref{cor:find zero of P, 2} with $L:=\hat K$  to obtain~$\hat b\in\hat K$ such that~$(P,\fm,\hat b)$ is a hole in $K$ equivalent to  $(P,\fm,\hat a)$.
\end{proof}

\section{Split-Normal Slots}\label{sec:split-normal holes}

\noindent
{\em In this section $H$ is a real closed $H$-field with small derivation and  asymptotic integration. We let $\mathcal{O}:= \mathcal{O}_H$ be its valuation ring and  $C:= C_H$ its constant field.  
We fix an immediate asymptotic extension $\hat H$ of~$H$ with valuation ring $\hat{\mathcal{O}}$ and an element~$\imag$ of an asymptotic extension of $\hat H$ with $\imag^2=-1$}. 
Then $\hat H$ is also an $H$-field by [ADH, 10.5.8], $\imag\notin\hat H$ and $K:=H[\imag]$ is an algebraic closure of $H$. With~$\hat K:= \hat H[\imag]$ we have the inclusion diagram
$$\xymatrix{ {\hat H} \ar@{-}[r]  &  \hat K ={\hat H}[\imag]  \\ 
 H \ar@{-}[u] \ar@{-}[r] &    K = H[\imag] \ar@{-}[u] }$$
By [ADH, 3.5.15, 10.5.7],
$K$ and $\hat K$ are $\d$-valued with valuation rings $\mathcal{O}+\mathcal{O}\imag$ and~$\hat{\mathcal{O}}+\hat{\mathcal{O}}\imag$ and
with the same constant field~$C[\imag]$, and $\hat K$ is an immediate extension of~$K$. Thus $H$, $K$, $\hat H$, $\hat K$ have the same
$H$-asymptotic couple $(\Gamma, \psi)$.

\begin{lemma} 
Let  $\hat a\in\hat H\setminus H$. Then $Z(H,\hat a) = Z\big(K,\hat a\big)\cap H\{Y\}$.
\end{lemma}
\begin{proof}
The inclusion ``$\supseteq$'' is obvious since  the Newton degree of a differential polynomial $Q\in H\{Y\}^{\neq}$
does not change when $H$ is replaced by its algebraic closure; see~[ADH, 11.1]. 
Conversely, let $P\in Z(H,\hat a)$. Then for all $\fv\in H^\times$ and~$a\in H$ such that $a-\hat a\prec\fv$ we have $\ndeg_{\prec\fv} H_{+a}\geq 1$.
Let $\fv\in H^\times$ and $z\in K$ be such that $z-\hat a\prec\fv$.
Take $a,b\in H$ such that $z=a+b\imag$. Then $a-\hat a,b\imag\prec\fv$ and hence~$\ndeg_{\prec\fv} P_{+z}=\ndeg_{\prec\fv} P_{+a}\geq 1$, using [ADH, 11.2.7].
Thus $P\in Z\big(K,\hat a\big)$.
\end{proof}

\begin{cor} \label{cor:holes in K vs holes in K[i]} 
Let $(P,\fm,\hat a)$ be a slot in $H$ with $\hat a\in\hat H$. Then 
$(P,\fm,\hat a)$ is also a slot in $K$, and if
   $(P,\fm,\hat a)$  is $Z$-minimal as a slot in $K$, then 
$(P,\fm,\hat a)$  is $Z$-minimal as a slot in $H$.
Moreover, $(P,\fm,\hat a)$ is a hole in $H$ iff  $(P,\fm,\hat a)$ is a hole in $K$, and
if $(P,\fm,\hat a)$ is a minimal hole in $K$, then $(P,\fm,\hat a)$ is a minimal hole in $H$.
\end{cor}
\begin{proof}
The first three claims are obvious from $\hat K$ being an immediate extension of $K$ and the previous lemma.
Suppose $(P,\fm,\hat a)$ is minimal as a hole  in $K$.
Let~$(Q,\fn,\tilde{b})$ be a hole in~$H$;
thus $\tilde{b}\in \widetilde{H}$ where $\widetilde{H}$ is an immediate asymptotic extension of $H$.
By the first part of the corollary applied to $(Q,\fn,\tilde{b})$ and $\widetilde{H}$ in place of
$(P,\fm,\hat a)$ and~$\hat H$,  respectively, $(Q,\fn,\tilde{b})$ is also a hole  in $K$. Hence $\cc(P) \leq \cc(Q)$, proving the last claim. 
\end{proof}

\noindent
In the next subsection we define the notion of a {\it split-normal}\/ slot in $H$. Later in this section we employ the results of Sections~\ref{sec:normalization}--\ref{sec:approx linear diff ops} to show, under suitable hypotheses on $H$, that minimal holes in $K$ of  order~$\geq 1$ give rise to a split-normal $Z$-minimal slots in $H$. (Theorem~\ref{thm:split-normal}.)
We then investigate which kinds of refinements preserve split-normality, and also consider a strengthening of split-normality.

\subsection*{Defining split-normality} {\em In this subsection  $b$ ranges over $H$ and $\fm, \fn$ over $H^\times$. Also, $(P, \fm, \hat a)$ is a slot in $H$ of order $r\ge 1$ with $\hat a\in \hat H\setminus H$ and linear part~$L:=L_{P_{\times \fm}}$. Set $w:=\wt(P)$, so $w\ge r$; if $\order L=r$, we set $\fv:=\fv(L)$}.
%  $\delta:=v(\fv)$}.

\begin{definition}\label{SN}  
We say that $(P,\fm,\hat a)$ is {\bf split-normal} if $\order L=r$, and \index{slot!split-normal}\index{split-normal}
\begin{itemize}
\item[(SN1)] $\fv\prec^\flat 1$; 
\item[(SN2)]  $(P_{\times\fm})_{\geq 1}=Q+R$ where $Q, R\in H\{Y\}$, $Q$ is homogeneous of degree~$1$ and order~$r$,  $L_Q$ splits over $K$, and $R\prec_{\Delta(\fv)} \fv^{w+1} (P_{\times\fm})_1$. 
\end{itemize}
\end{definition}

\noindent
Note that in (SN2) we do not require that $Q=(P_{\times\fm})_1$. 

\begin{lemma}\label{splnormalnormal} Suppose $(P,\fm,\hat a)$ is split-normal.  Then $(P,\fm, \hat a)$ is normal, and with~$Q$,~$R$ as in \textup{(SN2)} we have $(P_{\times \fm})_1-Q\prec_{\Delta(\fv)} \fv^{w+1}(P_{\times \fm})_1$, so $(P_{\times \fm})_1\sim Q$. 
\end{lemma}
\begin{proof} We have $(P_{\times \fm})_1=Q+R_1$ and $R_1\preceq R \prec_{\Delta(\fv)} \fv^{w+1} (P_{\times\fm})_1$, and thus
$(P_{\times \fm})_1-Q\prec_{\Delta(\fv)} \fv^{w+1}(P_{\times \fm})_1$. Now $(P,\fm,\hat a)$ is normal because $(P_{\times \fm})_{>1}= R_{>1} \prec_{\Delta(\fv)} \fv^{w+1}(P_{\times \fm})_1$.
\end{proof}

\noindent
If $(P,\fm,\hat a)$ is normal and $(P_{\times\fm})_{1}=Q+R$  where $Q, R\in H\{Y\}$, $Q$ is homogeneous of degree~$1$ and order $r$,  $L_Q$ splits over $K$, and $R\prec_{\Delta(\fv)} \fv^{w+1} (P_{\times\fm})_1$, then~$(P,\fm, \hat a)$ is split-normal. 
Thus  if $(P,\fm,\hat a)$ is normal and $L$ splits over $K$, then $(P,\fm,\hat a)$ is split-normal; in particular, if
$(P,\fm,\hat a)$ is normal  of order $r=1$, then it is split-normal. 
If~$(P,\fm,\hat a)$ is  split-normal, then so are $(bP,\fm,\hat a)$ for $b\neq 0$ and~$(P_{\times\fn},\fm/\fn,\hat a/\fn)$. 
Note also that if  $(P,\fm,\hat a)$ is split-normal, 
then with $Q$ as in~(SN2) we have $\fv(L)\sim\fv(L_Q)$, by Lemma~\ref{lem:fv of perturbed op}. If $(P,\fm, \hat a)$ is split-normal and $H$ is $\upl$-free, then $\exc^{\ev}(L)=\exc^{\ev}(L_Q)$ with $Q$ as in (SN2),  by
 Lemmas~\ref{splnormalnormal} and~\ref{cor:excev stability}.

\begin{lemma}\label{lem:split-normal comp conj}
Suppose $(P,\fm,\hat a)$ is split-normal and  $\phi\preceq 1$ is active in $H$ and~$\phi>0$ \textup{(}so $H^\phi$ is still an $H$-field\textup{)}. Then the slot $(P^\phi,\fm,\hat a)$ in $H^\phi$ is split-normal.
\end{lemma}
\begin{proof}
We first arrange $\fm=1$. Note that $L_{P^\phi}=L^\phi$ has order $r$. Put $\fw:=\fv(L_{P^\phi})$, and take $Q$, $R$ as in (SN2). Then 
$\fv\asymp_{\Delta(\fv)}\fw\prec^\flat_\phi 1$ by   Lemma~\ref{lem:v(Aphi)}. Moreover, $L_{Q^\phi}=L_Q^\phi$ splits over $K^\phi$; see [ADH, p.~291] or Lemma~\ref{lem:split and compconj}. By
[ADH, 11.1.4],  
$$R^\phi\ \asymp_{\Delta(\fv)} R\ \prec_{\Delta(\fv)}\ \fv^{w+1} P_1\ \asymp_{\Delta(\fv)}\ \fw^{w+1} P_1^\phi,$$
so $(P^\phi,\fm,\hat a)$  is split-normal. 
\end{proof}

\noindent
Since we need to preserve $H$ being an $H$-field when compositionally conjugating, we say:  {\em $(P^\phi,\fm, \hat a)$ is eventually 
split-normal\/} if there exists an active $\phi_0$ in $H$ such that $(P^\phi,\fm, \hat a)$ is split-normal for all active $\phi\preceq \phi_0$ in $H$ with $\phi>0$. 
We use this terminology in a similar way with ``split-normal'' replaced by other properties of slots of  order~$r \geq 1$ in real closed $H$-fields with small derivation and asymptotic integration, such as ``deep'' and ``deep and split-normal''. 

%In the same way we define the notion  {\em $(P^\phi,\fm, \hat a)$ is eventually deep\/} and the notion {\em $(P^\phi,\fm, \hat a)$ is  eventually deep and split-normal}. 

\subsection*{Achieving split-normality} {\em Assume $H$ is $\upo$-free and $(P,\fm,\hat a)$ is a minimal hole in $K=H[\imag]$ of order $r\ge 1$, with  $\fm\in H^\times$ and $\hat a\in\hat K\setminus K$}. Note that then~$K$ is $\upo$-free by [ADH, 11.7.23],  $K$ is
$(r-1)$-newtonian by Corollary~\ref{minholenewt}, and~$K$ is $r$-linearly closed  by Corollary~\ref{corminholenewt}.
In particular, the linear part of $(P,\fm,\hat a)$ is~$0$ or splits over~$K$. 
If  $\deg P=1$, then $r=1$ by
Corollary~\ref{cor:minhole deg 1}. 
If $\deg P >1$, then
$K$ and~$H$ are $r$-linearly newtonian by Corollary~\ref{degmorethanone} and Lemma~\ref{lem:descent r-linear newt}.  
In particular, if~$H$ is $1$-linearly newtonian, then $H$ is $r$-linearly newtonian.  
%We set~$K:=H[\imag]$ and $\hat K:= {\hat H}[\imag]$, so we have the inclusion diagram
%$$\xymatrix{ {\hat H} \ar@{-}[r]  &  \hat K ={\hat H}[\imag]  \\ 
% H \ar@{-}[u] \ar@{-}[r] &    K = H[\imag] \ar@{-}[u] }$$
%of $\d$-valued fields with common $H$-asymptotic couple $(\Gamma,\psi)$. 
{\em In this subsection we let $a$ range over $K$, $b$, $c$ over $H$, and $\fn$ over $H^\times$}.

\begin{lemma}\label{lem:hole in K with same c as P} 
Let~$(Q,\fn,\hat b)$ be a hole in $H$ with $\cc(Q)\le \cc(P)$ and
% of the same complexity as~$(P,\fm,\hat a)$, where 
$\hat b\in \hat H$. Then $\cc(Q)=\cc(P)$, $(Q,\fn,\hat b)$ is minimal and remains a minimal hole in $K$. The linear part of~$(Q,\fn,\hat b)$
is $0$ or splits over~$K$, and $(Q,\fn,\hat b)$ has a refinement~$(Q_{+b},\fp,\hat b-b)$ \textup{(}in $H$\textup{)} such that $(Q^\phi_{+b},\fp,\hat b-b)$ is eventually deep and split-normal. 
\end{lemma} 
\begin{proof}
By Corollary~\ref{cor:holes in K vs holes in K[i]}, 
$(Q,\fn,\hat b)$ is a  hole in $K$, and this hole in $K$ is minimal with $\cc(Q)=\cc(P)$, since $(P,\fm,\hat a)$ is minimal. By Corollary~\ref{cor:holes in K vs holes in K[i]} again, $(Q,\fn,\hat b)$  as a hole in $H$ is also minimal.  Since $K$ is $r$-linearly closed, the  linear part of $(Q,\fn,\hat b)$ is~$0$ or  splits over $K$. 
Corollary~\ref{cor:mainthm} gives a refinement $(Q_{+b},\fp,\hat b-b)$ of the minimal hole $(Q,\fn,\hat b)$ in $H$
such that $(Q_{+b}^\phi,\fp,\hat b-b)$ is deep and normal, eventually. Thus the linear part of $(Q_{+b},\fp,\hat b-b)$ is not $0$, and
as $\cc(Q_{+b})=\cc(P)$, this linear part splits over $K$. Hence for active $\phi$ in $H$ the linear part of
$(Q_{+b}^\phi,\fp,\hat b-b)$ splits over $K^\phi=H^\phi[\imag]$. Thus $(Q^\phi_{+b},\fp,\hat b-b)$ is eventually split-normal.
\end{proof}

\noindent
Now $\hat a=\hat b + \hat c\, \imag$ with $\hat b, \hat c\in \hat H$, and $\hat b, \hat c \prec \fm$. Moreover, $\hat b\notin H$ or $\hat c\notin H$.  
Since $\hat a$  is differentially algebraic over~$H$, so is its conjugate  $\hat b - \hat c\, \imag$, and therefore  its 
real and imaginary parts~$\hat b$ and~$\hat c$ are differentially algebraic over~$H$; thus $Z(\hat b,H)\ne \emptyset$ for $\hat b\notin H$, and
$Z(\hat c,H)\ne \emptyset$ for $\hat c\notin H$. More precisely:

\begin{lemma}\label{kb} We have $\operatorname{trdeg}\!\big(H\<\hat b\>|H\big)\le 2r$. If $\hat b\notin H$, then $Z(H,\hat b)\cap H[Y]=\emptyset$, so
$1\le \order Q \le 2r$ for all $Q\in Z(H,\hat b)$ of minimal complexity. 
These statements also hold for $\hat c$ instead of $\hat b$. 
\end{lemma} 
\begin{proof} The first statement follows from $\hat b \in  H\<\hat b +\hat c\, \imag, \hat b - \hat c\, \imag\>$. Suppose $\hat b\notin H$. If~$Q\in Z(H,\hat b)$ has minimal complexity, then [ADH, 11.4.8] yields an element $f$ in a proper immediate asymptotic extension of $H$ with  $Q(f)=0$, so $Q\notin H[Y]$.
\end{proof}

\begin{lemma}\label{lem:r=deg P=1}
Suppose $\deg P=1$ and $\hat b\notin H$.  
Let $Q\in Z(H,\hat b)$ be of minimal complexity; then either $\order{Q}=1$, or $\order{Q}=2$, $\deg{Q}=1$.
Let $\hat{Q}\in H\{Y\}$ be a minimal annihilator of
$\hat b$ over $H$;
then either $\order{\hat{Q}}=1$, or $\order{\hat{Q}}=2$, $\deg{\hat{Q}}=1$, and $L_{\hat{Q}}\in H[\der]$ splits over $K$.
\end{lemma}
\begin{proof}
Recall that $r=1$ by Corollary~\ref{cor:minhole deg 1}. Example~\ref{ex:lclm compl conj} and
Lemma~\ref{lem:lclm compl conj}  
give a $\tilde Q\in H\{Y\}$ of degree $1$ and order~$1$ or $2$ such that $\tilde Q(\hat b)=0$ and~$L_{\tilde Q}$
splits over $K$. Then $\cc(\tilde Q)=(1,1,1)$ or $\cc(\tilde Q)=(2,1,1)$, which proves the claim about~$Q$, using also Lemma~\ref{kb}. 
Also, $\tilde{Q},\hat{Q}\in Z(H,\hat b)$, hence $\cc(Q)\leq \cc(\hat Q)\leq\cc(\tilde Q)$. If~$\cc(\hat{Q})=\cc(\tilde Q)$, then $\hat{Q}=a\tilde Q$ for some~$a\in H^\times$. The claim about $\hat Q$ now follows easily. 
\end{proof}

\noindent
By Corollary~\ref{cor:mainthm} and Lemma~\ref{ufm}, our minimal hole $(P,\fm,\hat a)$ in $K$ has a refinement $(P_{+a},\fn,\hat a-a)$ such that eventually $(P_{+a}^\phi,\fn,\hat a-a)$ is deep and normal. 
Moreover, as $K$ is $r$-linearly closed, the
linear part of $(P_{+a}^\phi,\fn,\hat a-a)$ (for active~$\phi$ in~$K$)  splits over $K^\phi=H^\phi[\imag]$. Our main goal in this subsection is to prove  analogues
of these facts for  suitable $Z$-minimal slots $(Q,\fm,\hat b)$ or $(R,\fm,\hat c)$ in~$H$:

\begin{theorem}\label{thm:split-normal} 
If $H$ is $1$-linearly newtonian, then one of the following holds:
\begin{enumerate}
\item[$\mathrm{(i)}$] $\hat b\notin H$ and there exists a $Z$-minimal slot $(Q,\fm,\hat b)$ in $H$ with a re\-fine\-ment~${(Q_{+b},\fn,\hat b-b)}$ such that $(Q^\phi_{+b},\fn,\hat b-b)$ is eventually deep and split-normal; 
\item[$\mathrm{(ii)}$] $\hat c\notin H$ and there exists a $Z$-minimal slot $(R,\fm,\hat c)$ in $H$ with a refinement~${(R_{+c},\fn,\hat c-c)}$ such that $(R^\phi_{+c},\fn,\hat c-c)$ is eventually deep and split-normal. 
\end{enumerate}
\end{theorem}

\noindent
Lemmas~\ref{lem:hat c in K}, \ref{lem:hat b in K} and Corollaries~\ref{cor:evsplitnormal, 1}--\ref{cor:r=deg P=1, 2}
%~\ref{cor:evsplitnormal, 1}, \ref{cor:evsplitnormal, 2}, \ref{cor:r=deg P=1, 1},~\ref{cor:r=deg P=1, 2} 
below are more precise (only Corollary~\ref{cor:r=deg P=1, 1} has $H$ being $1$-linearly newtonian as a hypothesis) and together give Theorem~\ref{thm:split-normal}.  
We first deal with the case where $\hat b$ or $\hat c$ is in $H$:

\begin{lemma}\label{lem:hat c in K}
Suppose $\hat c\in H$. Then some hole $(Q,\fm,\hat b)$ in $H$ has the same complexity
as $(P,\fm,\hat a)$. Any such hole $(Q,\fm,\hat b)$ in $H$ is minimal and has a refinement~$(Q_{+b},\fn,\hat b-b)$ such that $(Q^\phi_{+b},\fn,\hat b-b)$ is eventually deep and split-normal. 
\end{lemma}
\begin{proof}
Let $A,B\in H\{Y\}$ be such that $P_{+\hat c\,\imag}(Y)=A(Y)+B(Y)\,\imag$.
Then $A(\hat b)=B(\hat b)=0$. If $A\ne 0$, then $\cc(A)\le \cc(P)$ gives that
$Q:=A$ has the desired property by Lemma~\ref{lem:hole in K with same c as P}. If $B\ne 0$, then likewise $Q:= B$ has the desired property.
% $\cc(P)=\max\big\{\!\cc(A),\cc(B)\big\}$ and . 
%Set $Q:=A$ if $\cc(P)=\cc(A)$ and $Q:=B$ otherwise;
%then $(Q,\fm,\hat b)$ is a hole in $K$  with the same complexity
%as $(P,\fm,\hat a)$. 
The rest also follows from that lemma. \end{proof}

\noindent
Thus if $\hat{c}\in H$, we obtain a strong version of (i) in Theorem~\ref{thm:split-normal}. Likewise, the next lemma gives a strong version of (ii) in Theorem~\ref{thm:split-normal} if $\hat{b}\in H$.

\begin{lemma}\label{lem:hat b in K} 
Suppose $\hat b\in H$. Then there is a hole $(R,\fm,\hat c)$ in $H$ with the same complexity
as $(P,\fm,\hat a)$. Every such hole in $H$ is minimal and   has a refinement~$(R_{+c},\fn,\hat c-c)$ such that $(R^\phi_{+c},\fn,\hat c-c)$ is eventually deep and split-normal. 
\end{lemma}

\noindent
This follows by applying Lemma~\ref{lem:hat c in K}  with $(P,\fm,\hat a)$ replaced by  the minimal hole~$\big(P_{\times\imag},\fm,-\imag\hat a\big)$ in $K$, which has the same complexity as $(P,\fm,\hat a)$.

\medskip
\noindent
{\em We assume in the rest of this subsection that $\hat b, \hat c\notin H$ and that 
$Q\in Z(H,\hat b)$ has minimal complexity}.
Hence~$(Q,\fm,\hat b)$ is a $Z$-minimal slot in $H$, and so is every refinement of $(Q,\fm,\hat b)$.
If~$(P_{+a},\fn,\hat a-a)$ is a refinement of~$(P,\fm,\hat a)$ and $b=\Re a$, then
$(Q_{+b},\fn,\hat b-b)$ is a refinement of $(Q,\fm,\hat b)$.
Conversely,  if~$(Q_{+b},\fn,\hat b-b)$ is a refinement of~$(Q,\fm,\hat b)$ and $v\big({\hat b-H}\big)\subseteq v({\hat c-H})$, then
Lemma~\ref{lem:same width} yields a refinement $(P_{+a},\fn,\hat a-a)$ of~$(P,\fm,\hat a)$ with $\Re a=b$. Recall from that lemma
that~$v({\hat b-H})\subseteq v({\hat c-H})$ is equivalent to~$v({\hat a-K})=v({\hat b-H})$; in this case, 
 $(P,\fm,\hat a)$ is special iff $(Q,\fm,\hat b)$ is special. 
 Recall also that if $(Q,\fm,\hat b)$ is deep, then so is each of its refinements 
$(Q_{+b},\fm,\hat b-b)$, by Corollary~\ref{cor:deep 2, cracks}.

\medskip
\noindent 
Here is a key technical fact underlying Theorem~\ref{thm:split-normal}:

\begin{prop}\label{evsplitnormal}
Suppose  the hole $(P,\fm,\hat a)$ in $K$ is special, the slot~$(Q,\fm,\hat b)$ in~$H$ is normal,
and $v\big({\hat b-H}\big)\subseteq v({\hat c-H})$. Then some refinement $(Q_{+b},\fm,{\hat b-b})$ of $(Q,\fm,\hat b)$ has the property that $(Q^\phi_{+b},\fm,{\hat b-b})$ is eventually split-normal.
\end{prop}

\begin{proof}
Replacing $(P,\fm,\hat a)$, $(Q,\fm,\hat b)$ by $(P_{\times\fm},1,\hat a/\fm)$, $(Q_{\times\fm},1,\hat b/
\fm)$, respectively, we reduce to the case $\fm=1$; then $\hat a, \hat b\prec 1$.
Since $\hat a$ is special over~$K=H[\imag]$,
$$\Delta\ :=\ \big\{ \delta\in\Gamma:\  \abs{\delta}\in v(\hat a-K) \big\}$$
is a convex subgroup of~$\Gamma$ which is cofinal in $v(\hat a-K)$ and hence
in $v(\hat b-H)$, so $\hat b$ is special over $H$.  Compositionally conjugate $H$, $\hat H$, $K$, $\hat K$ 
by a suitable active~$\phi\preceq 1$ in $H^{>}$, and replace $P$, $Q$ by $P^\phi$, $Q^\phi$, 
to arrange  $\Gamma^\flat\subseteq\Delta$; in particular, $\Psi\subseteq v(\hat b-H)$ and  $\psi(\Delta^{\neq})\subseteq\Delta$.
Multiplying $P$, $Q$ by suitable elements of $H^\times$ we also arrange that~$P,Q\asymp 1$. 
By Lemma~\ref{lem:split-normal comp conj} it suffices to show that then~$(Q,1,\hat b)$ has a split-normal refinement
$(Q_{+b},1,\hat b-b)$, and this is what we shall do. 

Note that $H$,~$\hat H$,~$K$,~$\hat K$ have small derivation, so the specializations~$\dot H$, $\dot{\hat H}$, $\dot K$, $\dot{\hat K}$ of~$H$, $\hat H$, $K$, $\hat K$, respectively, by~$\Delta$, are valued differential fields with small derivation. These specializations are asymptotic with 
asymptotic couple~$(\Delta,\psi|\Delta^{\neq})$, and of $H$-type with asymptotic integration, by [ADH, 9.4.12]; in addition they are $\d$-valued, by [ADH, 10.1.8]. 
The natural inclusions $\dot{\mathcal O}\to\dot{\mathcal O}_K$,  
$\dot{\mathcal O}\to\dot{\mathcal O}_{\hat H}$, $\dot{\mathcal O}_{\hat H} \to \dot{\mathcal O}_{\hat K}$, and
$\dot{\mathcal O}_K\to\dot{\mathcal O}_{\hat K}$ 
induce valued differential field embeddings
$\dot H\to\dot K$, $\dot H\to\dot{\hat H}$, $\dot{\hat H}\to\dot{\hat K}$ and $\dot K\to\dot{\hat K}$, which we make into inclusions by
the usual identifications; see [ADH, pp.~405--406].
 By Lemma~\ref{lem:dotK real closed} and the remarks preceding it, $\dot H$ is real closed with convex valuation ring and $\dot K$ is an algebraic closure of $\dot H$. Moreover,~$\dot{\hat H}$ is an immediate
extension of $\dot H$ and $\dot{\hat K}$ is an immediate
extension of $\dot K$.
Denoting the image of $\imag$ under the residue morphism $\dot{{\mathcal O}}_{\hat K}\to \dot{\hat K}$ by the same symbol, we then have $\dot K=\dot H[\imag]$, $\dot{\hat K}=\dot{\hat H}[\imag]$, and $\imag\notin\dot{\hat H}$. This gives the following inclusion diagram:
$$\xymatrix{ \dot{\hat H} \ar@{-}[r]  & \dot{\hat K} =\dot{\hat H}[\imag]  \\ 
\dot H \ar@{-}[u] \ar@{-}[r] &   \dot K =\dot H[\imag] \ar@{-}[u] }$$
Now $\hat a\in\mathcal O_{\hat K}\subseteq\dot{\mathcal O}_{\hat K}$ and $\hat b,\hat c\in\mathcal O_{\hat H}\subseteq\dot{\mathcal O}_{\hat H}$, and 
$\dot{\hat a}=\dot{\hat b} + \dot{\hat c}\,\imag$, $\Re \dot{\hat a}=\dot{\hat b}$, $\Im\dot{\hat a}=\dot{\hat c}$. 
For all~$a\in\dot{\mathcal O}_K$ we have $v(\dot{\hat a}-\dot a)=v(\hat a-a)\in\Delta$, hence
$\dot{\hat a}\notin \dot K$; likewise~${v(\hat b -b)}\in \Delta$
for all $b\in \dot{\mathcal O}$, so $\dot{\hat b}\notin \dot H$. Moreover, for all $\delta\in\Delta$ there is an $a\in\dot{\mathcal O}_K$ with~${v(\dot{\hat a}-\dot a)}=\delta$;
hence~$\dot{\hat a}$ is the limit of a c-sequence in~$\dot K$. This leads us to consider the completions~$\dot{H}^{\operatorname{c}}$ and $\dot{K}^{\operatorname{c}}$ of $\dot{H}$ and $\dot{K}$.
By [ADH, 4.4.11] and Lemma~\ref{lem:Kc real closed}, these  yield an inclusion diagram of valued differential field extensions:
$$\xymatrix{ \dot{H}^{\operatorname{c}} \ar@{-}[r]  & \dot{K}^{\operatorname{c}}=\dot{H}^{\operatorname{c}}[\imag]  \\ 
\dot H \ar@{-}[u] \ar@{-}[r] &   \dot K =\dot H[\imag] \ar@{-}[u] }$$
where $\dot H^{\operatorname{c}}$ is real closed with algebraic closure $\dot K^{\operatorname{c}}=
\dot H^{\operatorname{c}}[\imag]$. These completions are $\d$-valued by [ADH, 9.1.6]. By Corollary~\ref{cor:Kc newtonian},
$\dot{K}$ and $\dot K^{\operatorname{c}}$ are $\upo$-free and ${(r-1)}$-newtonian;  thus~$\dot K^{\operatorname{c}}$ is $r$-linearly closed by Corollary~\ref{14.5.3.r}. We identify the valued differential subfield   $\dot K\big\<\!\Re\dot{\hat a},\Im\dot{\hat a}\big\>$ of $\dot{\hat K}$ with its image under
the embedding into~$\dot K^{\operatorname{c}}$ over $\dot K$ from
Corollary~\ref{cor:embed into Kc}; then $\dot{\hat a}\in \dot{K}^{\operatorname{c}}$ and
$\dot{\hat b}=\Re\dot{\hat a}\in \dot H^{\operatorname{c}}$.  This leads to the next inclusion diagram:
$$\xymatrix{ \dot{H}^{\operatorname{c}} \ar@{-}[r]  & \dot{K}^{\operatorname{c}}   \\ 
\dot H\<\dot{\hat b}\>\ar@{-}[u] & \dot K\<\dot{\hat a}\>\ar@{-}[u] \\
\dot H \ar@{-}[u] \ar@{->}[r] &   \dot K   \ar@{-}[u] }$$
By Corollary~\ref{cor:ZKhata}, 
$\dot P\in\dot K\{Y\}$ is a minimal annihilator of $\dot{\hat a}$ over $\dot K$ and has the same complexity as $P$. Likewise, 
$\dot Q\in\dot H\{Y\}$ is a minimal annihilator of $\dot{\hat b}$ over
$\dot H$ and has the same complexity as $Q$. 
Let $s:=\order Q=\order\dot Q$, so $1\le s\le 2r$ by Lemma~\ref{kb}, and    
the linear part~$A\in \dot{H}^{\operatorname{c}}[\der]$ of $\dot Q_{+\dot{\hat b}}$ has order $s$ as well. By [ADH, 5.1.37] applied to $\dot H^{\operatorname{c}}$, $\dot H$, $\dot P$, $\dot Q$, $\dot{\hat a}$ in the role of $K$, $F$, $P$, $S$, $f$, respectively, $A$ splits over~$\dot K^{\operatorname{c}}=\dot H^{\operatorname{c}}[\imag]$, so Lemma~\ref{hkspl} gives a real splitting $(g_1,\dots, g_s)$ of $A$ over $\dot{K}^{\operatorname{c}}$: 
$$ A\ =\ f(\der-g_1)\cdots(\der-g_s),\qquad
f,g_1,\dots,g_s\in \dot K^{\operatorname{c}},\ f\neq 0.$$
%From $\dot{\hat b}\prec 1$ and [ADH, 4.5.1(i)] we get $\dot Q\sim \dot Q_{+\dot{\hat b}}$ and thus $\fv(A)\prec 1$ in $\dot L^{\operatorname{c}}$.
The slot $(Q,1,\hat b)$ in $H$ is normal, so $\fv(L_{Q_{+\hat b}}) \sim \fv(L_Q)\prec^\flat 1$ by Lemma~\ref{lem:linear part, new}, hence 
$\fv(A)\prec^\flat 1$ in $\dot K^{\operatorname{c}}$ by Lemma~\ref{lem:dotfv}. 
Then 
Corollary~\ref{cor:approx LP+f, real}  gives $a,b\in\dot{\mathcal O}$  and~$b_1,\dots,b_s\in\dot{\mathcal O}_K$ with $\dot a, \dot b\ne 0$ in
$\dot H$ such that for the linear part $\tilde{A}\in \dot{H}[\der]$ of~$\dot Q_{+\dot b}$,
$$\dot b\ \sim\ \dot{\hat b},\qquad \tilde{A}\ \sim\ A,\qquad \order \tilde{A}\ =\ s, \qquad
\fw\ :=\ \fv(\tilde{A})\ \sim\ \fv(A),$$
and such that for $w:=\wt(Q)$ and with $\Delta(\fw)\subseteq \Delta$: 
$$ \tilde{A} =  \tilde{B} +  \tilde{E},\ \tilde{B} = \dot a(\der-\dot b_1)\cdots(\der-\dot b_s)\in \dot{H}[\der],\quad 
\tilde{E}\in \dot{H}[\der],\quad  \tilde{E} \prec_{\Delta(\fw)}  \fw^{w+1} \tilde{A}, $$ 
%   \end{align*}
and $(\dot{b}_1,\dots, \dot{b}_s)$ is a real splitting of $\tilde{B}$ over $\dot{K}$.  Lemma~\ref{lem:lift real splitting} shows that we can change $b_1,\dots, b_s$ if necessary, without changing $\dot{b}_1,\dots, \dot{b}_s$, to arrange that $B:=a(\der-b_1)\cdots (\der-b_s)$ lies in $\dot{\mathcal{O}}[\der]\subseteq H[\der]$ and
$(b_1,\dots, b_s)$ is a real splitting of $B$ over~$K$.  
%\bigskip\noindent
%and such that in $\dot K[\der]$, with $w:=\wt(Q)$ and $\delta:=v(\fw)$: 
%$$ A^*\ =\ \dot a(\der-\dot b_1)\cdots(\der-\dot b_s) + \dot E\qquad\text{  where $\dot{E}\in\dot{H}[\der]$, $\dot E\prec_{\delta} \fw^{w+1} A^*$},$$
%in particular, $(\der-\dot b_1)\cdots(\der-\dot b_s)\in \dot{H}[\der]$.  
Now $\hat b-b\prec \hat b\prec 1$, so
$(Q_{+b},1,\hat b-b)$ is a refinement of the normal slot~$(Q,1,\hat b)$. Hence
$(Q_{+b},1,\hat b-b)$ is normal by Proposition~\ref{normalrefine}, so
$\fv:=\fv(L_{Q_{+ b}})\prec^\flat 1$.
By Lem\-ma~\ref{lem:dotfv} we have $\dot\fv=\fw$, so $\Delta(\fv)=\Delta(\fw)\subseteq \Delta$. Hence in $H[\der]$:
$$L_{Q_{+b}}\ =\ B +   E, \quad E\in \dot{\mathcal O}[\der],\  E\prec_{\Delta(\fv)} \fv^{w+1} L_{Q_{+ b}}.$$
Thus $(Q_{+b},1,\hat b-b)$ is split-normal.
\end{proof}

\noindent
Recall from the beginning of this subsection that if~${\deg P>1}$, 
then $K=H[\imag]$ is $r$-linearly new\-tonian;  this allows
us to remove the assumptions that~$(P,\fm,\hat a)$ is special and~$(Q,\fm,\hat b)$ is normal in Proposition~\ref{evsplitnormal}, by reducing to that case:

\begin{cor}  \label{cor:evsplitnormal, 1}
Suppose  $\deg P>1$ and $v(\hat b -H)\subseteq v(\hat c - H)$.
Then~$(Q,\fm,\hat b)$  has a special refinement~$(Q_{+b},\fn,\hat b-b)$ such that $(Q^\phi_{+b},\fn,{\hat b-b})$ is eventually deep and split-normal. 
\end{cor}
\begin{proof}
By Lemmas~\ref{lem:quasilinear refinement} and~\ref{ufm}, the hole $(P,\fm,\hat a)$ in $K$ has a quasilinear refinement~${(P_{+a},\fn,\hat a-a)}$. (The use of Lemma~\ref{ufm} is because we require~${\fn\in H^\times}$.)
Let $b=\Re a$. Then, using Lemma~\ref{lem:same width} for the second equality,
$$v\big((\hat a-a)-K\big)\ =\ v(\hat a - K)\ =\ v(\hat b -H)\ =\ {v\big((\hat b-b)-H\big)},$$
and~${(Q_{+b},\fn,\hat b-b)}$ is a $Z$-minimal refinement of $(Q,\fm,\hat b)$.
We replace~$(P,\fm,\hat a)$ and~$(Q,\fm,\hat b)$ by
$(P_{+a},\fn,\hat a-a)$ and $(Q_{+b},\fn,\hat b-b)$, respectively, to arrange that the hole~$(P,\fm,\hat a)$ in $K$ is quasilinear.
Then
by Proposition~\ref{prop:hata special} and $K$ being $r$-linearly newtonian,  $(P,\fm,\hat a)$ is special. 
Hence~$(Q,\fm,\hat b)$ is also special, so Proposition~\ref{varmainthm} gives a refinement~$(Q_{+b},\fn,\hat b-b)$ of $(Q,\fm,\hat b)$ and an active~${\phi_0\in H^{>}}$ such that $(Q^{\phi_0}_{+b},\fn,{\hat b-b})$ is deep and normal. 
Refinements of~$(P,\fm,\hat a)$ remain quasilinear, by Corollary~\ref{cor:ref 2n}. Since $v(\hat b-H)\subseteq v(\hat c -H)$ we have
a refinement~$(P_{+a},\fn,\hat a-a)$ of~$(P,\fm,\hat a)$ with~$\Re a=b$.
Then by Lemma~\ref{speciallemma} the minimal hole $(P^{\phi_0}_{+a},\fn,{\hat a-a})$ in~$H^{\phi_0}[\imag]$ is special. Now apply Proposition~\ref{evsplitnormal} with~$H^{\phi_0}$, $(P^{\phi_0}_{+a},\fn,\hat a-a)$, $(Q^{\phi_0}_{+b},\fn,{\hat b-b})$
in place of $H$, $(P,\fm,\hat a)$,   $(Q,\fm,\hat b)$, respectively: it gives~${b_0\in H}$ and a 
refinement 
$$\big( (Q^{\phi_0}_{+b})_{+b_0},\fn, (\hat b-b)-b_0\big)\ =\ 
\big(Q^{\phi_0}_{+(b+b_0)},\fn, {\hat b - (b+b_0)}\big)$$ of 
$(Q^{\phi_0}_{+b},\fn,{\hat b-b})$, and thus a refinement~$\big(Q_{+(b+b_0)},\fn, {\hat b - (b+b_0)}\big)$ of $(Q_{+b},\fn,{\hat b-b})$, such that  $\big(Q^\phi_{+(b+b_0)}, \fn,{ \hat b - (b+b_0)}\big)$  is eventually split-normal. 
By the remark before Proposition~\ref{evsplitnormal}, $\big(Q^{\phi}_{+(b+b_0)},\fn, {\hat b - (b+b_0)}\big)$ is also eventually deep. 
\end{proof}

\noindent
%In this corollary we assumed $v\big(\hat a-K[\imag]\big)=v(\hat b-K)$. Recall from Lemma~\ref{lem:same width} that this is equivalent to $v(\hat b-K)\subseteq v(\hat c-K)$, and that 
Recall that $v(\hat b-H)\subseteq v(\hat c-H)$ or  $v(\hat c-H)\subseteq v(\hat b-H)$. The following corollary concerns the second case:

\begin{cor} \label{cor:evsplitnormal, 2}
If  $\deg P>1$, $v(\hat c-H)\subseteq v(\hat b-H)$, and $R\in Z(H,\hat c)$ has minimal complexity,  then the $Z$-minimal slot~$(R,\fm,\hat c)$ in $H$ has a special refinement~$(R_{+c},\fn,\hat c-c)$   such that $(R^\phi_{+c},\fn,{\hat c-c})$ is eventually deep and split-normal.
\end{cor}
\begin{proof}
Apply Corollary~\ref{cor:evsplitnormal, 1} to the minimal hole $(P_{\times \imag},\fm,-\imag\hat a)$ in~$H[\imag]$. 
\end{proof}

\noindent
In the next two corollaries we handle the case $\deg P=1$.
Recall from Lemma~\ref{lem:r=deg P=1} that then
$\order{Q}=1$ or $\order{Q}=2$, $\deg{Q}=1$. 
Theorem~\ref{mainthm} gives:

\begin{cor}\label{cor:r=deg P=1, 1}
Suppose $H$ is $1$-linearly newtonian and $\order{Q}=1$. 
Then the slot~$(Q,\fm,\hat b)$ in $H$  has a refinement~$(Q_{+b},\fn,\hat b-b)$   such that $(Q^\phi_{+b},\fn,{\hat b-b})$ is eventually deep and split-normal. 
\end{cor}

\begin{cor} \label{cor:r=deg P=1, 2}
Suppose $\deg P=1$ and $\order Q=2$, $\deg Q=1$.  
Let $\hat{Q}\in H\{Y\}$ be a minimal annihilator of $\hat b$ over~$H$.
Then $\big(\hat{Q},\fm,\hat b\big)$ is a $Z$-minimal hole in~$H$ and has a refinement~$\big({\hat{Q}_{+b},\fn,\hat b-b}\big)$  such that $\big(\hat{Q}^\phi_{+b},\fn,{\hat b-b}\big)$ is eventually deep and split-normal. 
\end{cor}
\begin{proof}
By the proof of Lemma~\ref{lem:r=deg P=1} we have $\cc(Q)=\cc(\hat{Q})$ 
(hence  $\big(\hat{Q},\fm,\hat b\big)$  is a $Z$-minimal hole in $H$)
and $L_{\hat{Q}}$ splits over $H[\imag]$. Corollary~\ref{cor:deepening} gives a refine\-ment~$\big({\hat{Q}_{+b},\fn,\hat b-b}\big)$ of $\big(\hat{Q},\fm,\hat b\big)$ whose linear part has Newton weight~$0$ and such that the slot $\big({\hat{Q}^\phi_{+b},\fn,\hat b-b}\big)$ in $H^\phi$ is deep, eventually.  Moreover, by Lemmas~\ref{lem:deg1 normal} and~\ref{lem:deg 1 cracks splitting}, $\big(\hat{Q}^\phi_{+b},\fn,{\hat b-b}\big)$ is normal and its linear part splits over $H^{\phi}[\imag]$, eventually. Thus $\big(\hat{Q}^\phi_{+b},\fn,{\hat b-b}\big)$ is eventually deep and split-normal. 
\end{proof}

\noindent
This concludes the proof of Theorem~\ref{thm:split-normal}.

\subsection*{Split-normality and refinements}
We now study the behavior of split-normality under refinements. In this subsection $a$ ranges over $H$ and $\fm$, $\fn$, $\fv$ range over $H^\times$.
Let $(P,\fm,\hat a)$ be a slot in $H$ of order $r\geq 1$ with $\hat a\in\hat H\setminus H$, and $L:=L_{P_{\times\fm}}$, $w:=\wt(P)$.
Here is the split-normal analogue of Lemma~\ref{lem:normal pos criterion}:

\begin{lemma}\label{lem:split-normal criterion}
Suppose $\order(L)=r$ and $\fv$ is such that \textup{(SN1)} and~\textup{(SN2)} hold, and $\fv(L)\asymp_{\Delta(\fv)} \fv$. Then $(P,\fm,\hat a)$ is split-normal.
\end{lemma}
\begin{proof} Same as that of~\ref{lem:normal pos criterion}, but with $R$ as in (SN2) instead of $(P_{\times \fm})_{>1}$.
\end{proof}

\noindent
Now split-normal analogues of Propositions~\ref{normalrefine} and~\ref{easymultnormal}:

\begin{lemma}\label{splitnormalrefine}
Suppose $(P,\fm,\hat a)$ is split-normal. Let a refinement $(P_{+a},\fm,{\hat a-a})$ of $(P,\fm,\hat a)$ be given. Then $(P_{+a},\fm,\hat a-a)$ is also split-normal.
\end{lemma}
\begin{proof}
As in the proof of Proposition~\ref{normalrefine} we arrange $\fm=1$ and show for $\fv:= \fv(L_P)$, using Lemmas~\ref{lem:linear part, new} and~\ref{splnormalnormal}, that $\order(L_{P_{+a}})=r$ and $$(P_{+a})_1\sim_{\Delta(\fv)} P_1, \quad \fv(L_{P_{+a}})\sim_{\Delta(\fv)}\fv, \quad (P_{+a})_{>1} \prec_{\Delta(\fv)} \fv^{w+1} (P_{+a})_1.$$
Now take $Q$, $R$ as in (SN2) for $\fm=1$. Then $P_1=Q+R_1$, and so by Lemma~\ref{lem:linear part, split-normal, new} for $A=L_Q$ we
obtain $(P_{+a})_1 - Q \prec_{\Delta(\fv)} \fv^{w+1} (P_{+a})_1$, and thus
$(P_{+a})_{\geq 1}-Q \prec_{\Delta(\fv)}\fv^{w+1} (P_{+a})_1$.
Hence (SN2) holds with $\fm=1$ and $P_{+a}$ instead of~$P$. Thus the slot~$(P_{+a},\fm,\hat a-a)$ in $H$ is split-normal by Lemma~\ref{lem:split-normal criterion}.
\end{proof}

\begin{lemma}\label{easymultsplitnormal}
Suppose $(P,\fm,\hat a)$ is split-normal,  $\hat{a}\prec \fn\preceq \fm$, and $[\fn/\fm]\le[\fv]$, $\fv:=\fv(L)$.
% with \big[\fv(L_{P_{\times\fm}})\big]$. 
Then the refinement $(P,\fn,\hat a)$ of $(P,\fm,\hat a)$
is split-normal: if $\fm$, $P$, $Q$, $\fv$ are as in \textup{(SN2)}, then \textup{(SN2)} holds with $\fn$, $Q_{\times \fn/\fm}$, $R_{\times \fn/\fm}$, $\fv(L_{P_{\times \fn}})$ in place of~$\fm$,~$Q$,~$R$,~$\fv$.
\end{lemma} 
%\textup{(SN2)} holds with $\fn, Q_{\times \fn/\fm}, R_{\times \fn/\fm}$ in %place of $\fm, Q, R$. 
%\end{lemma}
\begin{proof} Set $\tilde{L}:= L_{P_{\times\fn}}$. Lemma~\ref{lem:steep1} gives $\order(\tilde{L})=r$ and $\fv(\tilde{L})\asymp_{\Delta(\fv)}\fv$. Thus~$(P_{\times\fn})_{>1} \prec_{\Delta(\fv)} \fv^{w+1} (P_{\times\fn})_1$ by Proposition~\ref{easymultnormal}. Now arrange $\fm=1$ in the usual way, and
take~$Q$,~$R$ as in (SN2) for $\fm=1$.  Then $$(P_{\times\fn})_1\ =\ (P_1)_{\times\fn}\ =\ Q_{\times\fn}+(R_1)_{\times\fn}, \qquad (P_{\times \fn})_{>1}\ =\ (R_{\times \fn})_{>1}\ =\ (R_{>1})_{\times \fn}$$ by~[ADH,~4.3], where
$Q_{\times\fn}$ is homogeneous of degree $1$ and order $r$, and $L_{Q_{\times\fn}}=L_Q\fn$ splits over~$K$.
Using [ADH, 4.3, 6.1.3] and $[\fn]\leq[\fv]$ we obtain
$$(R_1)_{\times\fn}\ \asymp_{\Delta(\fv)}\ \fn R_1\ \preceq\ \fn R\ \prec_{\Delta(\fv)}\ \fn\fv^{w+1}  P_1\ \asymp_{\Delta(\fv)} \fv^{w+1} (P_1)_{\times\fn}\ =\ \fv^{w+1} (P_{\times\fn})_1.$$
Hence (SN2) holds for $\fn, Q_{\times \fn}, R_{\times \fn}, \fv(\tilde{L})$ in place of $\fm, Q, R,\fv$. 
%Using Lemma~\ref{lem:split-normal criterion} it follows 
%that~$(P,\fn,\hat a)$ is split-normal.  
\end{proof}

\noindent 
Recall our standing assumption in this section that $H$ is a real closed $H$-field. Thus~$H$ is $\d$-valued, and
for all $\fn$ and $q\in \Q^{>}$ we have $\fn^q\in H^\times$ such that $(\fn^q)^\dagger=q\fn^\dagger$. {\em In the 
rest of this section we fix such an $\fn^q$ for all $\fn$ and $q\in\Q^{>}$.}\/ 
Now we upgrade Corollary~\ref{cor:normal for small q}  with ``split-normal'' instead of ``normal'':

\begin{lemma}\label{splnq}
Suppose $\fm=1$, $(P,1,\hat a)$ is split-normal, $\hat a \prec \fn\prec 1$, and for~$\fv:=\fv(L_P)$ we have $[\fn^\dagger]<[\fv]<[\fn]$. 
Then $(P,\fn^q,\hat a)$ is a split-normal refinement
of~$(P,1,\hat a)$ for all but finitely many $q\in \Q$ with $0 < q < 1$. 
\end{lemma}
\begin{proof}
Corollary~\ref{cor:normal for small q} gives that $(P,\fn^q,\hat a)$ is a normal refinement of $(P,1,\hat a)$ for all but finitely many $q\in \Q$ with $0<q<1$. 
Take~$Q$,~$R$  as in (SN2) for $\fm=1$. Then $L = L_Q + L_R$ where $L_Q$ splits over~$H[\imag]$ and $L_R \prec_{\Delta(\fv)} \fv^{w+1} L$,
for $\fv:=\fv(L)$.
%By Lemma~\ref{lem:fv of perturbed op} we have $\fv\sim \fv(L_Q)$.
Applying Corollary~\ref{cor:nepsilon} to $A:=L$, $A_*:=L_R$ we obtain: 
$L_R\fn^q \prec_{\Delta(\fw)} \fw^{w+1} L\fn^q$, $\fw:=\fv(L\fn^q)$, for all but finitely many $q\in \Q^>$. 

Let $q\in \Q$ be such that $0<q<1$, $(P,\fn^q,\hat a)$ is a normal refinement of $(P,1,\hat a)$, and $L_R\fn^q \prec_{\Delta(\fw)} \fw^{w+1} L\fn^q$, with $\fw$ as above. 
Then $(P_{\times\fn^q})_1 = Q_{\times\fn^q}+(R_1)_{\times\fn^q}$ where~$Q_{\times\fn^q}$ is homogeneous of degree~$1$ and order $r$,  $L_{Q_{\times\fn^q}}=L_Q\fn^q$ splits over~$H[\imag]$,
and $(R_1)_{\times\fn^q} \prec_{\Delta(\fw)} \fw^{w+1} (P_{\times\fn^q})_1$ for $\fw:=\fv(L_{P_{\times\fn^q}})$. Since  $(P,\fn^q, \hat a)$ is normal, we also have  
$(P_{\times\fn^q})_{>1}\prec_{\Delta(\fw)} \fw^{w+1} (P_{\times\fn^q})_1$. Thus $(P,\fn^q,\hat a)$ is split-normal. 
\end{proof}

\begin{remark} We do not know if in this last lemma we can drop the assumption $[\fn^\dagger]<[\fv]$. 
\end{remark}

\subsection*{Strengthening split-normality} 
{\em In this subsection $a, b$ range over $H$ and $\fm,\fn$ over $H^\times$,
and $(P,\fm,\hat a)$ is a slot in $H$ of order~$r\geq 1$ and weight $w:=\wt(P)$, so~${w\geq 1}$, and $L:=L_{P_{\times\fm}}$. If
$\order L=r$, we set $\fv:= \fv(L)$}. 
%we work under the same assumptions as in the subsection ``Defining split-normality'' above. In particular,~

With an eye towards later use in connection with fixed point theorems over Hardy fields we strengthen here the concept of
split-normality; in the next subsection we  show how to improve Theorem~\ref{thm:split-normal} accordingly.
See the last subsection of Section~\ref{sec:approx linear diff ops} for the notion of strong splitting.

\begin{definition}
{\samepage Call 
$(P,\fm,\hat a)$ {\bf almost strongly split-normal} \index{slot!almost strongly split-normal}\index{almost strongly!split-normal}\index{split-normal!almost strongly} if $\order L=r$, 
$\fv\prec^\flat 1$, and the following strengthening of (SN2) holds:
\begin{list}{*}{\addtolength\itemindent{-3.5em}\addtolength\leftmargin{0.5em}}
\item[(SN2as)]  $(P_{\times\fm})_{\geq 1}=Q+R$ where $Q, R\in H\{Y\}$, $Q$ is homogeneous of degree~$1$ and order~$r$,  $L_Q$ splits strongly over $K$, and $R\prec_{\Delta(\fv)} \fv^{w+1} (P_{\times\fm})_1$. 
\end{list}}\noindent
We say that $(P,\fm,\hat a)$ is {\bf strongly split-normal}
 if $\order L =r$, $\fv\prec^\flat 1$, and the following condition is satisfied:\index{slot!strongly split-normal}\index{strongly!split-normal}\index{split-normal!strongly}
\begin{enumerate}
\item[(SN2s)]  
$P_{\times\fm}=Q+R$ where $Q,R\in H\{Y\}$, $Q$ is homogeneous of degree~$1$ and order~$r$, 
$L_Q$   splits strongly over $K$, and $R\prec_{\Delta(\fv)} \fv^{w+1} (P_{\times\fm})_1$.
\end{enumerate}
\end{definition}

\noindent
To facilitate use of (SN2s) we observe:

\begin{lemma}\label{SN2suse} Suppose $(P,\fm,\hat a)$ is strongly split-normal and $P_{\times\fm}=Q+R$ as in {\rm(SN2s)}. Then $Q\sim (P_{\times \fm})_1$, 
$\fv_Q:= \fv(L_Q)\sim \fv$, so $R\prec_{\Delta(\fv)} \fv_Q^{w+1}Q$.
\end{lemma}
\begin{proof} We have $(P_{\times \fm})_1=Q+R_1$, so 
$Q=(P_{\times \fm})_1-R_1$ with $R_1\prec_{\Delta(\fv)} \fv^{w+1}(P_{\times \fm})_1$. Now apply Lemma~\ref{lem:fv of perturbed op} to $A:=L$ and $B:=-L_{R_1}$. 
\end{proof}

\noindent
If $(P,\fm,\hat a)$  is almost strongly split-normal, then
$(P,\fm,\hat a)$  is split-normal and hence normal by Lemma~\ref{splnormalnormal}. 
If $(P,\fm,\hat a)$ is normal and $L$ splits strongly over~$K$, then~$(P,\fm,\hat a)$ is  almost strongly split-normal;
in particular, if $(P,\fm,\hat a)$ is  normal of order $r=1$, then~$(P,\fm,\hat a)$ is almost strongly split-normal, by Lemma~\ref{lem:order 1 splits strongly}.
Moreover:

\begin{lemma}\label{lem:char strong split-norm}
The following are equivalent:
\begin{enumerate}
\item[\textup{(i)}] $(P,\fm,\hat a)$ is strongly split-normal;
\item[\textup{(ii)}] $(P,\fm,\hat a)$ is almost strongly split-normal and strictly normal;
\item[\textup{(iii)}] $(P,\fm,\hat a)$ is almost strongly split-normal and $P(0)\prec_{\Delta(\fv)} \fv^{w+1} (P_1)_{\times\fm}$.
\end{enumerate}
\end{lemma}
\begin{proof}
Suppose $(P,\fm,\hat a)$ is strongly split-normal, and
let $Q$, $R$ be as in (SN2s). Then~$(P_{\times\fm})_{\geq 1}=Q+R_{\geq 1}$, $L_Q$   splits strongly over $K$, and $R_{\geq 1} \prec_{\Delta(\fv)} \fv^{w+1} (P_{\times\fm})_1$.
Hence  $(P,\fm,\hat a)$  is almost strongly split-normal, and thus normal.
Also~$P(0)=R(0)\prec_{\Delta(\fv)} \fv^{w+1} (P_{\times\fm})_1$, so $(P,\fm,\hat a)$ is strictly normal.
This shows~(i)~$\Rightarrow$~(ii), and~(ii)~$\Rightarrow$~(iii) is clear.
For (iii)~$\Rightarrow$~(i) suppose~$(P,\fm,\hat a)$ is almost strongly split-normal and 
$P(0)\prec_{\Delta(\fv)} \fv^{w+1} (P_1)_{\times\fm}$. Take~$Q$,~$R$ as in~(SN2as). Then $P_{\times\fm} = Q+\tilde R$ where 
$\tilde R:=P(0)+R\prec_{\Delta(\fv)}\fv^{w+1} (P_1)_{\times\fm}$.
Thus $(P,\fm,\hat a)$ is strongly split-normal.
\end{proof}

\begin{cor}\label{cor:strongly splitting => strongly split-normal}
If $L$ splits strongly over $K$, then
$$\text{$(P,\fm,\hat a)$ is strongly split-normal }\ \Longleftrightarrow\   \text{$(P,\fm,\hat a)$ is strictly normal.}$$
\end{cor}

\noindent
The following diagram summarizes some implications between these variants of normality, for slots $(P,\fm,\hat a)$ in $H$
of order $r\ge 1$:  
$$\xymatrix@L=6pt{	\text{strongly split-normal\ } \ar@{=>}[r] \ar@{=>}[d]&  \text{\ almost strongly split-normal} \ar@{=>}[r] &  \text{\ split-normal} \ar@{=>}[d] \\ 
\text{strictly normal\ } \ar@{=>}[rr] & & \text{\ normal}   }$$
If $(P,\fm,\hat a)$ is almost strongly split-normal, then  so are
$(bP,\fm,\hat a)$ for $b\neq 0$ and $(P_{\times\fn},\fm/\fn,\hat a/\fn)$, and likewise with ``strongly''
in place of ``almost strongly''.  

\medskip
\noindent
Here is a version of Lemma~\ref{splitnormalrefine} for (almost) strong split-normality:

\begin{lemma}\label{stronglysplitnormalrefine}
Suppose $(P_{+a},\fm,\hat a-a)$ refines $(P,\fm,\hat a)$.
If $(P,\fm,\hat a)$ is almost strongly split-normal, then so is $(P_{+a},\fm,\hat a-a)$.
If  $(P,\fm,\hat a)$ is   strongly split-normal, $Z$-minimal, and 
$\hat a - a \prec_{\Delta(\fv)} \fv^{r+w+1}\fm$, then $(P_{+a},\fm,\hat a-a)$ is strongly split-normal. 
\end{lemma}
\begin{proof}
The first part follows from Lemma~\ref{splitnormalrefine} and its proof.
In combination with Lemmas~\ref{stronglynormalrefine, cracks} and \ref{lem:char strong split-norm}, this also yields the second part.
\end{proof} 
 
\begin{lemma}\label{stronglysplitnormalrefine, q}
Suppose that $(P,\fm,\hat a)$ is split-normal and  $\hat{a}\prec_{\Delta(\fv)} \fm$.
Then for all sufficiently small  $q\in\Q^>$, any $\fn\asymp\fv^q\fm$ yields  an almost strongly split-normal refinement $(P,\fn,\hat a)$ of~$(P,\fm,\hat a)$.
\end{lemma}
\begin{proof}
We arrange $\fm=1$, so $\hat{a}\prec_{\Delta(\fv)} 1$. Take $Q$, $R$ as in~(SN2) with $\fm=1$, and take~$q_0\in\Q^>$ such that $\hat a\prec \fv^{q_0}\prec 1$. By Lemma~\ref{lem:split strongly multconj} we
can decrease~$q_0$ so that for all $q\in\Q$ with $0<q\leq q_0$ and any $\fn\asymp \fv^q$, $L_{Q_{\times\fn}}=L_Q\fn$ splits strongly over~$K$. 
Suppose $q\in \Q$, $0<q\leq q_0$, and $\fn\asymp\fv^q$.
Then  $(P,\fn,\hat a)$ is an almost strongly split-normal refinement of $(P, 1,\hat a)$, by Lemma~\ref{easymultsplitnormal}. 
%Indeed, let  
%$\fw:=\fv(L_{P_{\times\fv^q}})$, $\varepsilon:=v(\fw)$he proof
%of this lemma shows that  $\fw\prec^\flat 1$ and
%$(P_{\times\fn})_{\geq 1} = Q_{\times\fn} + R_{\times\fn}$ where~$Q_{\times\fv^q}$ is homogeneous of degree~$1$  and order~$r$, and
%$R_{\times\fv^q} \prec_\varepsilon \fw^{w+1} (P_{\times\fv^q})_1$.
%Now $L_{Q_{\times\fv^q}}$ splits strongly over $K[\imag]$.
\end{proof}

\begin{cor}\label{cor:deep and almost strongly split-normal}
Suppose that $(P,\fm,\hat a)$ is $Z$-minimal, deep, and split-normal. Then $(P,\fm,\hat a)$ has a refinement which is deep and
almost strongly split-normal.
\end{cor}
\begin{proof}  
Lemma~\ref{lem:good approx to hata} gives $a$ such that $\hat a - a \prec_{\Delta(\fv)} \fm$. By Corollary~\ref{cor:deep 2, cracks}, the re\-fine\-ment~$(P_{+a},\fm,\hat a-a)$ of
$(P,\fm,\hat a)$ is   deep  with $\fv(L_{P_{+a,\times \fm}})\asymp_{\Delta(\fv)} \fv$, and by Lemma~\ref{splitnormalrefine} it is also
split-normal. 
Now apply Lemma~\ref{stronglysplitnormalrefine, q} to $(P_{+a},\fm,\hat a-a)$ in place of~$(P,\fm,\hat a)$  and again use Corollary~\ref{cor:deep 2, cracks}.
\end{proof}

\noindent
We now turn to the behavior of these properties under compositional conjugation.

\begin{lemma}\label{lem:strongly split-normal compconj}
Let $\phi$ be active in $H$ with $0<\phi\preceq 1$.
If $(P,\fm,\hat a)$ is almost strongly split-normal, then so is 
the slot $(P^\phi,\fm,\hat a)$ in $H^\phi$. Likewise with ``strongly''
in place of ``almost strongly''.
\end{lemma}
\begin{proof}
We arrange $\fm=1$, 
assume $(P,\fm,\hat a)$ is almost strongly split-normal, and take~$Q$,~$R$ as in~(SN2as).
% where $\fm=1$ and $L_Q$ splits strongly over $K[\imag]$.
The proof of Lemma~\ref{lem:split-normal comp conj} shows that with $\fw:=\fv(L_{P^\phi})$
we have $\fw\prec^\flat_\phi 1$ and  $(P^\phi)_{\geq 1} = Q^\phi + R^\phi$ where
$Q^\phi\in H^\phi\{Y\}$ is homogeneous of degree~$1$ and order~$r$, $L_{Q^\phi}$ splits over~$H^\phi[\imag]$, and
$R^\phi \prec_{\Delta(\fw)} \fw^{w+1} (P^\phi)_1$. By Lemma~\ref{lem:split strongly compconj}, $L_{Q^\phi}=L_Q^\phi$ even  splits strongly over $H[\imag]$.
Hence~$(P^\phi,\fm,\hat a)$ is almost strongly split-normal.
The rest follows from Lemma~\ref{lem:char strong split-norm} and 
the fact that  if~$(P,\fm,\hat a)$ is strictly normal, then so is  $(P^\phi,\fm,\hat a)$.
\end{proof}

%\noindent
%We say that {\em $(P^\phi,\fm, \hat a)$ is eventually  strongly split-normal\/} if for some  active~$\phi_0$ in~$H$, the dent $(P^\phi,\fm, \hat a)$ in $H^\phi$ is strongly split-normal for all active $\phi\preceq \phi_0$ in~$H$ with~$\phi>0$; likewise with ``almost strongly'' in place of ``strongly''.

%\begin{lemma}\label{lem:order 1 strongly split-normal}
%Suppose $H$ is $1$-linearly newtonian and $\upo$-free. Then every $Z$-minimal dent in $H$ of order~$r=1$ has a refinement  $(P,\fm,\hat a)$  such that $(P^\phi,\fm,\hat a)$ is eventually deep and strongly split-normal.
%\end{lemma}
%\begin{proof} Let a $Z$-minimal dent in $H$ of order $1$ be given. Using first Theorem~\ref{mainthm} and then Corollary~\ref{cor:achieve strong normality, 1} gives a refinement $(P,\fm,\hat a)$ of this dent and an active $\phi_0$ in $H$ with $0 < \phi_0\preceq 1$ such that the dent $(P^{\phi_0}, \fm, \hat a)$ in $H^{\phi_0}$ is deep and strongly normal. Then by Lemma~\ref{lem:char strong split-norm} and the sentence preceding it, $(P^{\phi_0},\fm,\hat a)$ is strongly split-normal. Thus $(P,\fm,\hat a)$ has the desired property by Lemma~\ref{lem:strongly split-normal compconj}.
%\end{proof} 

\noindent
If $H$ is $\upo$-free and $r$-linearly newtonian, then by Corollary~\ref{cor:achieve strong normality, 2}, every $Z$-minimal slot in~$H$ of order~$r$
has a refinement $(P,\fm,\hat a)$  such that the slot~$(P^\phi,\fm,\hat a)$ in $H^\phi$ is eventually deep and strictly normal.
Corollary~\ref{cor:5.30real} of the next lemma 
is a variant of this fact for strong split-normality.

\begin{lemma}\label{5.30real} 
Assume $H$ is $\upo$-free and $r$-linearly newtonian, and every
$A\in H[\der]$ of order~$r$ splits over $K$.  
Suppose $(P,\fm,\hat a)$ is $Z$-minimal. Then there is a refinement $(P_{+a},\fn,\hat a-a)$  of $(P,\fm,\hat a)$ and an active $\phi$ in $H$
with $0<\phi\preceq 1$   such that  $(P^\phi_{+a},\fn,\hat a-a)$  is  deep and strictly normal, and its   linear part splits strongly over $K^\phi$ \textup{(}so $(P^\phi_{+a},\fn,\hat a-a)$  is strongly split-normal by Corollary~\ref{cor:strongly splitting => strongly split-normal}\textup{)}.
\end{lemma}
\begin{proof} 
For any active $\phi$ in $H$ with $0<\phi\preceq 1$ we may replace~$H$,~$(P,\fm,\hat a)$  by~$H^\phi$,~$(P^\phi,\fm,\hat a)$, respectively. We may
 also replace~$(P,\fm,\hat a)$ by any of its refinements. 
Now Theorem~\ref{mainthm} gives a refinement~$(P_{+a},\fn,\hat a-a)$ of $(P,\fm,\hat a)$ and an active~$\phi$ in $H$ such that~$0<\phi\preceq 1$ and $(P^\phi_{+a},\fn,\hat a-a)$ is deep and   normal.
Replacing~$H$,~$(P,\fm,\hat a)$ by~$H^\phi$,~$(P^\phi_{+a},\fn,{\hat a-a})$, respectively,  
we thus arrange that~$(P,\fm,\hat a)$ itself is deep and   normal. We show that then the lemma holds with~$\phi=1$. 
For this we first replace~$(P,\fm,\hat a)$  by a suitable refinement~$(P_{+a},\fm,{\hat a-a})$ to arrange by
Corollary~\ref{cor:achieve strong normality, 1}  that~$(P,\fm,\hat a)$ is strictly normal and $\hat a \prec_{\Delta(\fv)} \fm$. 
Now $L$ splits over~$K$, so by 
Corollary~\ref{cor:split strongly multconj},
% and~\ref{cor:holes in K vs holes in K[i]}, 
for sufficiently small~${q\in\Q^>}$, any $\fn\asymp|\fv|^q\fm$ gives a  refinement~$(P,\fn,\hat a)$  of~$(P,\fm,\hat a)$ whose linear part~$L_{P_{\times\fn}}$ has order $r$ and splits strongly over $K$.
For each such $\fn$,  $(P,\fn,\hat a)$ is deep  by Corollary~\ref{cor:deep 2, cracks}, and for some such $\fn$,
 $(P,\fn,\hat a)$ is also strictly normal, by Remark~\ref{rem:strongly normal refine, 2}.
\end{proof}

\noindent
The previous lemma in combination with Lemma~\ref{lem:strongly split-normal compconj}  yields:
 
\begin{cor}\label{cor:5.30real}
With the same assumptions on $H$, $K$ as in Lemma~\ref{5.30real},  every $Z$-minimal slot in $H$ of order $r$ has a refinement $(P,\fm,\hat a)$ such that
$(P^\phi,\fm,\hat a)$ is eventually deep and strongly split-normal. 
\end{cor}

\noindent
For $r=1$ the splitting assumption is automatically satisfied (and this is the case most relevant later). 
We do not know whether ``every
$A\in H[\der]^{\ne}$ of order~$\le r$ splits over $K$'' is strictly weaker than ``$K$ is $r$-linearly closed''.

\subsection*{Achieving strong split-normality}
We make the same assumptions as in the subsection {\it Achieving split-normality}\/: {\em $H$ is $\upo$-free and $(P,\fm,\hat a)$ is a minimal hole in $K=H[\imag]$ of order $r\ge 1$, with  $\fm\in H^\times$ and $\hat a\in\hat K\setminus K$.}\/ Recall that $K$ is also
$\upo$-free [ADH, 11.7.23]. 
We have  $$\hat a\ =\ \hat b + \hat c\, \imag,\qquad \hat b, \hat c\in \hat H.$$
We let  $a$ range over $K$, $b$, $c$  over $H$, and $\fn$ over~$H^\times$. 
In connection with the next two lemmas we note that given 
an active $\phi$ in~$H$ with~${0<\phi\preceq 1}$, 
if 
$(P,\fm,\hat a)$ is normal (strictly normal, respectively), then
so is~$(P^\phi,\fm,\hat a)$, by Lemma~\ref{lem:normality comp conj} (Lemma~\ref{lem:normality comp conj, strong}, respectively);
moreover, if  the linear part of $(P,\fm,\hat a)$ splits strongly over $K$, then the linear part of $(P^\phi,\fm,\hat a)$ splits strongly over $K^\phi=H^\phi[\imag]$,
by Lem\-ma~\ref{lem:split strongly compconj}. Here is a ``complex'' version of Lemma~\ref{5.30real}, with a similar proof:

\begin{lemma}\label{lem:achieve strong splitting}  
For some refinement $(P_{+a},\fn,\hat a-a)$  of $(P,\fm,\hat a)$ and active $\phi$ in~$H$ with $0<\phi\preceq 1$,
the hole $(P^\phi_{+a},\fn,\hat a-a)$ in $K^\phi$ is deep and  normal, its  linear part splits 
strongly over $K^\phi$, and it is moreover strictly normal if $\deg P>1$. 
\end{lemma}
\begin{proof} For any active $\phi$ in $H$ with $0<\phi\preceq 1$ 
we may replace 
 $H$ and $(P,\fm,\hat a)$  by~$H^\phi$ and the minimal hole $(P^\phi,\fm,\hat a)$ in $K^\phi$. We may
 also replace~$(P,\fm,\hat a)$   by any of its refinements $(P_{+a},\fn,\hat a-a)$. %(Recall that  $\fn\in K^\times$.)
As noted before Theorem~\ref{thm:split-normal}, Corollary~\ref{cor:mainthm} and Lemma~\ref{ufm} give a refinement 
$(P_{+a},\fn,\hat a-a)$ of $(P,\fm,\hat a)$ and an active~${\phi}$ in $H$
with $0<\phi\preceq 1$ such that   $(P^\phi_{+a},\fn,\hat a-a)$ is deep and normal.
Replacing $H$, $(P,\fm,\hat a)$ by~$H^\phi$,~$(P^\phi_{+a},\fn,{\hat a-a})$, respectively,  
we thus arrange that~$(P,\fm,\hat a)$ itself is deep and normal. We show that then the lemma holds with $\phi=1$. 

Set $L:=L_{P_{\times\fm}}$ and $\fv:=\fv(L)$.
Lemma~\ref{lem:good approx to hata} gives $a$ with~${\hat a - a \prec_{\Delta(\fv)} \fm}$. 
If~${\deg P>1}$, then~$K$ is $r$-linearly newtonian and we use Corollary~\ref{cor:closer to minimal holes} to take~$a$ such that  
even $\hat a - a \preceq \fv^{w+2}\fm$.
Replacing~$(P,\fm,\hat a)$  by $(P_{+a},\fm,\hat a-a)$, we thus arrange by Lemma~\ref{lem:deep 2} and 
Proposition~\ref{normalrefine} that $\hat a \prec_{\Delta(\fv)} \fm$, and also by Lemma~\ref{lem:achieve strong normality} that~$(P,\fm,\hat a)$  is  strictly normal if $\deg P>1$. 
Now $L$ splits over~$K$, since~$K$ is $r$-linearly closed by Corollary~\ref{corminholenewt}.
Then by Corollary~\ref{cor:split strongly multconj}, for sufficiently small~${q\in\Q^>}$, any $\fn\asymp|\fv|^q\fm$ gives a refinement $(P,\fn,\hat a)$   of~$(P,\fm,\hat a)$ whose linear part~$L_{P_{\times\fn}}$ splits strongly over $K$. 
For such $\fn$,  $(P,\fn,\hat a)$ is deep by Lemma~\ref{lem:deep 2} and normal by Proposition~\ref{easymultnormal}. 
If~$(P,\fm,\hat a)$ is strictly normal, then for some such~$\fn$, 
 $(P,\fn,\hat a)$ is also strictly normal, thanks to Lemma~\ref{lem:strongly normal refine, 2}.  % and the remark preceding it.
\end{proof}

\noindent
The following version of Lemma~\ref{lem:achieve strong splitting} also encompasses linear $(P,\fm,\hat a)$: 

\begin{lemma}\label{lem:achieve strong splitting, d=1}
Suppose $\der K=K$ and $\I(K)\subseteq K^\dagger$. Then there is a refinement~$(P_{+a},\fn,\hat a-a)$  of $(P,\fm,\hat a)$ and an active $\phi$ in~$H$ with~$0<\phi\preceq 1$ such that
the hole~$(P^\phi_{+a},\fn,\hat a-a)$ in $K^\phi$ is deep and strictly normal, and its  linear part splits 
strongly over $K^\phi$. 
\end{lemma}
\begin{proof}  
Thanks to Lemma~\ref{lem:achieve strong splitting} we need only consider the case $\deg P=1$. Then we have~$r=1$ by Corollary~\ref{cor:minhole deg 1}. (See now the remark following this proof.) 
As in the proof of Lem\-ma~\ref{lem:achieve strong splitting} we may replace 
 $H$ and $(P,\fm,\hat a)$ for any active $\phi\preceq 1$ in $H^{>}$ by $H^\phi$ and~$(P^\phi, \fm,\hat a)$, and
 also $(P,\fm,\hat a)$   by any of its refinements $(P_{+a},\fn,\hat a-a)$. Recall here that  $\fn\in H^\times$.
 Hence using a remark preceding Lemma~\ref{lem:strongnormal pos criterion} and Corollary~\ref{cor:strongly normal d=1}   
 we   arrange that $(P,\fm,\hat a)$ is strictly normal, and thus balanced and deep. We  show that then the lemma holds with $\phi=1$. 

 Set $L:=L_{P_{\times\fm}}$, $\fv:=\fv(L)$.
 Lemmas~\ref{lem:excev empty} and~\ref{lem:balanced good approx} yield an~$a$ with~$\hat a-a \preceq \fv^4\fm$.
Replacing~$(P,\fm,\hat a)$  by $(P_{+a},\fm,\hat a-a)$ arranges that $\hat a \prec_{\Delta(\fv)} \fm$, by Lemmas~\ref{lem:deep 2} and~\ref{stronglynormalrefine}. As in the proof of Lemma~\ref{lem:achieve strong splitting}, for sufficiently small~${q\in\Q^>}$, any~$\fn\asymp|\fv|^q\fm$ now gives a strictly normal and deep refinement $(P,\fn,\hat a)$   of~$(P,\fm,\hat a)$ whose linear part splits strongly over $K$.
\end{proof}

\begin{remark}
Suppose we replace our standing assumption that $H$ is $\upo$-free and $(P, \fm, \hat a)$ is a minimal hole in $K$ by the assumption that $H$ is $\upl$-free and $(P, \fm, \hat a)$ is a slot in $K$ of order and degree $1$ (so $K$ is $\upl$-free by [ADH, 11.6.8] and $(P, \fm,\hat a)$ is $Z$-minimal). 
Then Lemma~\ref{lem:achieve strong splitting, d=1} goes through with ``hole" replaced by ``slot''. Its proof also goes through with the references to Lemmas~\ref{lem:deep 2} and~\ref{stronglynormalrefine} replaced by references to Corollary~\ref{cor:deep 2, cracks} and Lemma~\ref{stronglynormalrefine, cracks}. The end of that proof refers to the end of the proof of  Lemma~\ref{lem:achieve strong splitting}, and there one should replace Proposition~\ref{easymultnormal}  by Corollary~\ref{corcorcor}, and  Lemma~\ref{lem:strongly normal refine, 2} by  Remark~\ref{rem:strongly normal refine, 2}. 
\end{remark}

\noindent
In the remainder of this subsection we prove the following variant of Theorem~\ref{thm:split-normal}:

\begin{theorem}\label{thm:strongly split-normal} If $H$ is $1$-linearly newtonian, then one of the following holds: 
\begin{enumerate}
\item[\textup{(i)}] $\hat b\notin H$ and there exists a $Z$-minimal slot $(Q,\fm,\hat b)$ in $H$ with a refinement~${(Q_{+b},\fn,\hat b-b)}$ such that $(Q^\phi_{+b},\fn,\hat b-b)$ is eventually deep and almost strongly split-normal; 
\item[\textup{(ii)}] $\hat c\notin H$ and there exists a $Z$-minimal slot $(R,\fm,\hat c)$ in $H$ with a refinement~${(R_{+c},\fn,\hat c-c)}$ such that $(R^\phi_{+c},\fn,\hat c-c)$ is eventually deep and almost strongly split-normal. 
\end{enumerate}
Moreover, if $H$ is $1$-linearly newtonian and either  $\deg P>1$, or $\hat b\notin H$ and   $Z(H,\hat b)$ contains an element of order~$1$, or $\hat c\notin H$ and   $Z(H,\hat c)$ contains an element of order~$1$,
then \textup{(i)} holds with ``almost'' omitted, or  \textup{(ii)} holds with ``almost'' omitted. 
\end{theorem}

\noindent
Towards the proof of this theorem we first show:

\begin{lemma}\label{lem:refine to almost strongly split-normal, Q}
Suppose $\hat b\notin H$ and $(Q,\fm,\hat b)$ is  a $Z$-minimal slot in $H$ with a refinement~${(Q_{+b},\fn,\hat b-b)}$  such that $(Q^\phi_{+b},\fn,\hat b-b)$ is eventually deep and split-normal.
Then $(Q,\fm,\hat b)$ has a refinement~${(Q_{+b},\fn,\hat b-b)}$  such that $(Q^\phi_{+b},\fn,\hat b-b)$ is eventually deep and almost strongly split-normal.
\end{lemma}
\begin{proof}
Let ${(Q_{+b},\fn,\hat b-b)}$ be a refinement of $(Q,\fm,\hat b)$  and
let $\phi_0$ be active in $H$ such that $0<\phi_0\preceq 1$
and $(Q^{\phi_0}_{+b},\fn,\hat b-b)$ is   deep and split-normal.
Then  Corollary~\ref{cor:deep and almost strongly split-normal}
yields a refinement $\big((Q^{\phi_0}_{+b})_{+b_0},\fn_0,(\hat b-b)-b_0\big)$ of $(Q^{\phi_0}_{+b},\fn,\hat b-b)$ which is deep and
almost strongly split-normal. Hence 
$$\big((Q_{+b})_{+b_0},\fn_0,(\hat b-b)-b_0\big)\ =\ \big( Q_{+(b+b_0)},\fn_0,\hat b - (b+b_0) \big)$$
is  a refinement of $(Q,\fm,\hat b)$, and $\big( Q^\phi_{+(b+b_0)},\fn_0,\hat b - (b+b_0) \big)$ is eventually deep and almost strongly split-normal by Lemma~\ref{lem:strongly split-normal compconj}.
\end{proof}

\noindent
Likewise:

\begin{lemma}\label{lem:refine to almost strongly split-normal, R}
Suppose $\hat c\notin H$, and $(R,\fm,\hat c)$ is  a $Z$-minimal slot in $H$ with   a refinement~${(R_{+c},\fn,\hat c-c)}$  such that $(R^\phi_{+c},\fn,\hat c-c)$ is eventually deep and split-normal.
Then $(R,\fm,\hat c)$  has a refinement~${(R_{+c},\fn,\hat c-c)}$  such that $(R^\phi_{+c},\fn,\hat c-c)$ is eventually deep and almost strongly split-normal.
\end{lemma}

\noindent
Theorem~\ref{thm:split-normal} and the two lemmas above give the first part of Theorem~\ref{thm:strongly split-normal}.  
We break up the proof of the ``moreover'' part  into several cases, along the lines of the proof of Theorem~\ref{thm:split-normal}. We begin with the case where $\hat b\in H$ or~$\hat c\in H$.

\begin{lemma}\label{lem:refine to almost strongly split-normal, order(Q)=r}
Suppose $H$ is $1$-linearly newtonian, $\hat b\notin H$, $(Q,\fm,\hat b)$ is a $Z$-minimal slot in $H$ of order $r$, and some 
refinement~${(Q_{+b},\fn,\hat b-b)}$ of $(Q,\fm, \hat b)$ is such that~$(Q^\phi_{+b},\fn,\hat b-b)$ is eventually deep and split-normal.
Then $(Q,\fm,\hat b)$ has a refinement~${(Q_{+b},\fn,\hat b-b)}$  with $(Q^\phi_{+b},\fn,\hat b-b)$ eventually deep and strongly split-normal.
\end{lemma}
\begin{proof}
Lemma~\ref{lem:refine to almost strongly split-normal, Q} gives a refinement 
$(Q_{+b},\fn,\hat b-b)$ of $(Q,\fm,\hat b)$ with $(Q^\phi_{+b},\fn,\hat b-b)$ eventually deep and almost strongly split-normal.
We  upgrade this to ``strongly split-normal'' as follows: Take active $\phi_0$ in $H$ with $0<\phi_0\preceq 1$
such that the slot~$(Q^{\phi_0}_{+b},\fn,\hat b-b)$ in $H^{\phi_0}$ is  deep and almost strongly split-normal.  
Now $H$ is $1$-linearly newtonian, hence $r$-linearly newtonian. Therefore
Corollary~\ref{cor:achieve strong normality, 1} yields a deep and strictly normal refinement
$\big( (Q^{\phi_0}_{+b})_{+b_0},\fn, (\hat b - b)-b_0 \big)$
of $\big( Q^{\phi_0}_{+b},\fn,\hat b - b \big)$. 
By
Lemma~\ref{stronglysplitnormalrefine}, this refinement is still almost  strongly split-normal,   thus strongly split-normal by
Lemma~\ref{lem:char strong split-norm}. 
Then by  Lemma~\ref{lem:strongly split-normal compconj},  
$\big( Q_{+(b+b_0)},\fn,{\hat b - (b+b_0)} \big)$
is  a refinement of $(Q,\fm,\hat b)$ such that $\big( Q^\phi_{+(b+b_0)},\fn,\hat b - (b+b_0) \big)$ is eventually deep and strongly split-normal.
\end{proof}

\noindent
Lemmas~\ref{lem:hat c in K} and~\ref{lem:refine to almost strongly split-normal, order(Q)=r} give the following:

\begin{cor}\label{cor:hat c in K, strong}
Suppose $H$ is $1$-linearly newtonian and $\hat c\in H$. Then there is a hole $(Q,\fm,\hat b)$ in $H$ of the same complexity
as $(P,\fm,\hat a)$. Every such hole $(Q,\fm,\hat b)$ in $H$ is minimal and has a refinement~$(Q_{+b},\fn,\hat b-b)$ such that $(Q^\phi_{+b},\fn,\hat b-b)$ is eventually deep and strongly split-normal.
\end{cor}

\noindent
Just as Lemma~\ref{lem:hat c in K} gave rise to Lemma~\ref{lem:hat b in K},  Corollary~\ref{cor:hat c in K, strong} leads to:

\begin{cor}\label{cor:hat b in K, strong}
Suppose $H$ is $1$-linearly newtonian and $\hat b\in H$. Then there is a hole~$(R,\fm,\hat c)$ in $H$ of the same complexity
as $(P,\fm,\hat a)$. Every such hole in $H$ is minimal and   has a refinement~$(R_{+c},\fn,\hat c-c)$ such that $(R^\phi_{+c},\fn,\hat c-c)$ is eventually deep and strongly split-normal. 
\end{cor}

\noindent
In the following two lemmas we assume that $\hat b,\hat c\notin H$. Let $Q\in Z(H,\hat b)$ be of minimal complexity, so $(Q,\fm,\hat b)$ is a $Z$-minimal slot in $H$, as is each of its refinements. 
The next lemma strengthens Corollary~\ref{cor:evsplitnormal, 1}:

\begin{lemma}\label{lem:evstronlysplitnormal, 1}
Suppose $\deg P>1$ and $v(\hat b-H)\subseteq v(\hat c-H)$. Then
$(Q,\fm,\hat b)$ has a refinement $(Q_{+b},\fn,\hat b-b)$ such that $(Q^\phi_{+b},\fn,\hat b-b)$ is eventually deep and strongly split-normal.
\end{lemma}
\begin{proof}
Corollary~\ref{cor:evsplitnormal, 1} and Lemma~\ref{lem:refine to almost strongly split-normal, Q} give a refinement 
$(Q_{+b},\fn,\hat b-b)$ of $(Q,\fm,\hat b)$ and an active $\phi_0$ in $H$ with $0<\phi_0\preceq 1$
such that the slot $(Q^{\phi_0}_{+b},\fn,{\hat b-b})$ in $H^{\phi_0}$ is  deep and almost strongly split-normal.  From $\deg P >1$ we obtain that~$H$ is $r$-linearly newtonian. Now argue as in the proof of Lemma~\ref{lem:refine to almost strongly split-normal, order(Q)=r}.
%using that $(Q^{\phi_0}_{+b},\fn,\hat b-b)$ is special when applying Corollary~\ref{cor:achieve strong normality, cracks}.
\end{proof}

\noindent
Similarly we obtain a strengthening of Corollary~\ref{cor:evsplitnormal, 2}, using that corollary and Lemma~\ref{lem:refine to almost strongly split-normal, R} in place of Corollary~\ref{cor:evsplitnormal, 1} and Lemma~\ref{lem:refine to almost strongly split-normal, Q} in the proof:

\begin{lemma}\label{lem:evstronlysplitnormal, 2}
If $\deg P >1$, $v(\hat c-H)\subseteq v(\hat b-H)$, and $R\in Z(H,\hat c)$ has minimal complexity, then the $Z$-minimal slot~$(R,\fm,\hat c)$ in $H$ has a refinement~$(R_{+c},\fn,\hat c-c)$   such that $(R^\phi_{+c},\fn,{\hat c-c})$ is eventually deep and strongly split-normal. 
\end{lemma}

\noindent
We now prove the ``moreover'' part of Theorem~\ref{thm:strongly split-normal}.
Thus, suppose $H$ is $1$-linearly newtonian.
If~$\hat b\in H$,   then $\hat c\notin H$ and
Corollary~\ref{cor:hat b in K, strong} yields a strong version of~(ii)
with ``almost'' omitted.  Likewise, if $\hat c\in H$, then $\hat b\notin H$ and Corollary~\ref{cor:hat c in K, strong} yields a strong version of~(i),
with ``almost'' omitted.
In the rest of the proof we assume~${\hat b,\hat c\notin H}$.
By Lemma~\ref{lem:same width} we have $v({\hat b-H})\subseteq v({\hat c-H})$ or  $v(\hat c-H)\subseteq v({\hat b-H})$,
and thus Lemmas~\ref{lem:evstronlysplitnormal, 1} and \ref{lem:evstronlysplitnormal, 2} take care of the case~${\deg P >1}$. 
If~$Z(H,\hat b)$ contains an element of order~$1$, and $Q\in Z(H,\hat b)$ has minimal complexity, then $\order Q =1$ by Lemma~\ref{kb}, so Corollary~\ref{cor:5.30real} and the remark following it
yield (i) with ``almost'' omitted. Likewise, if $Z(H,\hat c)$ contains an element of order~$1$, then~(ii) holds with ``almost'' omitted.
\qed

\section{Ultimate Slots and Firm Slots}\label{sec:ultimate}

\noindent 
{\em In this section $H$ is a Liouville closed $H$-field
with small derivation, $\hat H$ is an immediate asymptotic extension 
of~$H$, and $\imag$ be an element of an asymptotic extension of $\hat H$
with $\imag^2=-1$.} 
Then~$\hat H$ is an $H$-field, $\imag\notin\hat H$,
$K:=H[\imag]$  is an algebraic closure of $H$, and $\hat K:=\hat H[\imag]$ is an immediate $\d$-valued extension of $K$. (See the beginning of Section~\ref{sec:split-normal holes}.) Let $C$ be the constant field of $H$, let $\mathcal{O}$ denote the valuation ring of $H$ and $\Gamma$ its value group. Accordingly, the constant field of $K$ is $C_K=C[\imag]$ and the valuation ring of~$K$ is~$\mathcal{O}_K=\mathcal{O}+\mathcal{O}\imag$.  
Let $\fm$, $\fn$, $\fw$ range over $H^\times$ and $\phi$ over the  elements of $H^>$ which are active in~$H$ (and hence in $K$).

\medskip
\noindent
In Section~\ref{sec:logder} we introduced
$$W\ :=\ \big\{\!\wr(a,b):\ a,b\in H,\ a^2+b^2=1\big\}.$$
Note that $W$ is a subspace of the 
$\Q$-linear space $H$, because $W\imag=S^\dagger$ where 
$$S\ :=\ \{a+b\imag:\ a,b\in H,\ a^2+b^2=1\}$$ is a divisible subgroup of $K^\times$. We have $K^\dagger=H +W\imag$ by Lemma~\ref{lem:logder}. 
Thus there exists a complement $\Lambda$ of the subspace~$K^\dagger$ of $K$ such that $\Lambda\subseteq H\imag$, and in this section we fix such $\Lambda$ and let~$\lambda$ range over $\Lambda$. 
 Let $\Univ=K\big[\! \ex(\Lambda) \big]$ be the universal exponential extension of~$K$ defined in Section~\ref{sec:univ exp ext}.

\medskip
\noindent
For $A\in K[\der]^{\ne}$ we have its set~$\exc^{\operatorname{u}}(A)\subseteq \Gamma$ of ultimate exceptional values, 
which a-priori might depend on our choice of~$\Lambda$. We now make good on a promise from Section~\ref{sec:valuniv} by showing under the mild assumption $\I(K)\subseteq K^\dagger$ and with our restriction $\Lambda\subseteq H\imag$ there is no such dependence:

\begin{cor}\label{prop:excu independent of Q} Suppose $\I(K)\subseteq K^\dagger$. Then for $A\in K[\der]^{\neq}$, the status of $A$ being terminal does not depend on the choice of $\Lambda$, and  the set $\exc^{\operatorname{u}}(A)$ of ultimate exceptional values of $A$  also does not depend on this choice.
%  of $\Lambda$, and $A$ being ultimate doesn't depend on this either. 
\end{cor} 
\begin{proof}
Let $\Lambda^*\subseteq H\imag$ also be a complement of $K^\dagger$. Let   $\lambda\mapsto \lambda^*$ be the $\Q$-linear bijection $\Lambda\to \Lambda^*$ with $\lambda- \lambda^*\in W\imag$ for all $\lambda$. Then by Lemmas~\ref{pldv} and~\ref{lem:W and I(F)},
$$\lambda- \lambda^*\in \I(H)\imag\ \subseteq\ \I(K)\ \subseteq\  (\mathcal O_K^\times)^\dagger$$  for all $\lambda$. Now use  Lemma~\ref{lem:excu for different Q} and Corollary~\ref{cor:excu for different Q}.
\end{proof}

\begin{cor}\label{cor:excu independent of Q} 
Suppose $\I(K)\subseteq K^\dagger$. Let $A=\der-g\in K[\der]$ where $g\in K$ and let~$\mathfrak g\in H^\times$ be such that $\mathfrak g^\dagger=\Re g$. Then
$$\exc^{\operatorname{u}}(A)\ =\ v_{\g}(\ker_{\Univ}^{\neq} A)\ =\ \{v\mathfrak g\}.$$
In particular, if $\Re g\in \I(H)$, then~$\exc^{\operatorname{u}}(A)=\{0\}$.
\end{cor}
\begin{proof}
Let $f\in K^\times$ and $\lambda$ be such that $g=f^\dagger+\lambda$. Then 
$$\exc^{\operatorname{u}}(A)\ =\ v_{\g}(\ker_{\Univ}^{\neq} A)\ =\ \{vf\}$$ by  Lemma~\ref{lem:excu, r=1} and its proof.
Recall that $K^\dagger = H+\I(H)\imag$ by Lemma~\ref{lem:W and I(F)} and remarks preceding it, so $g\in K^\dagger$ iff $\Im g\in\I(H)$.
Consider first the case $g\notin K^\dagger$. Then by Corollary~\ref{prop:excu independent of Q} we can change
$\Lambda$ if necessary to arrange $\lambda:=(\Im g)\imag\in\Lambda$ so that we can take $f:=\mathfrak g$ in the above.
Now suppose~$g\in K^\dagger$. Then $g=(\mathfrak g h)^\dagger$ where $h\in K^\times$, $h^\dagger=(\Im g)\imag$.  Then 
we can take $f:=\mathfrak g h$, $\lambda:=0$, and we have $h\asymp 1$ since~$h^\dagger\in \I(H)\imag\subseteq\I(K)$.
\end{proof}

\begin{cor}\label{cor:LambdaL, purely imag}
Suppose  $\I(K)\subseteq K^\dagger$, and let $F$ be a Liouville closed $H$-field extension of $H$, and $L:=F[\imag]$. 
Then the subspace~$L^\dagger$ of the  $\Q$-linear space $L$ has a
complement~$\Lambda_L$ with $\Lambda\subseteq\Lambda_L\subseteq F\imag$. For any such $\Lambda_L$ and 
  $A\in K[\der]^{\neq}$ we have~~$\exc^{\ev}(A_\lambda)=\exc^{\ev}_L(A_\lambda)\cap\Gamma$ for all $\lambda$, and thus  
  $\exc^{\operatorname{u}}(A) \subseteq \exc^{\operatorname{u}}_L(A)$, where $\exc^{\operatorname{u}}_L(A)$ is the
set of  ultimate exceptional values of $A\in L[\der]^{\ne}$  with respect to $\Lambda_L$.
\end{cor}
\begin{proof}
By the remarks at the beginning of this subsection applied to $F$, $L$ in place of~$H$,~$K$ we have 
$L^\dagger=F+W_F\imag$ where
$W_F$ is a subspace of the $\Q$-linear space~$F$. Also  $K^\dagger=H + \I(H)\imag$ by Lemma~\ref{lem:W and I(F)}, and $L^\dagger\cap K= K^\dagger$ by Lemma~\ref{lem:LambdaL}.  This yields a
complement~$\Lambda_L$ of $L^\dagger$ in $L$ with
$\Lambda\subseteq\Lambda_L\subseteq F\imag$. 
 Since~$H$ is Liouville closed and hence $\upl$-free by [ADH, 11.6.2], its algebraic closure~$K$ is $\upl$-free by [ADH, 11.6.8].   Now the rest follows from  remarks preceding Lemma~\ref{lem:excu 1}. 
 %lem:LambdaL}.
\end{proof}

\noindent
Given $A\in K[\der]^{\neq}$, let $\exc^{\operatorname{u}}(A^\phi)$ be the set of ultimate exceptional values
of the linear differential operator $A^\phi\in K^\phi[\derdelta]$, $\derdelta=\phi^{-1}\der$, with respect to  $\Lambda^\phi=\phi^{-1}\Lambda$. 
We summarize some properties of ultimate exceptional values used later in this section:

\begin{lemma}\label{lem:excu properties}
Let $A\in K[\der]^{\neq}$ have order~$r$. Then for all $b\in K^\times$ and all $\phi$,
$$\exc^{\operatorname{u}}(bA)\ =\ \exc^{\operatorname{u}}(A),\quad  \exc^{\operatorname{u}}(Ab)\ =\ \exc^{\operatorname{u}}(A)-vb,
\quad \exc^{\operatorname{u}}(A^\phi)\ =\ \exc^{\operatorname{u}}(A).$$
Moreover, if  $\I(K)\subseteq K^\dagger$, then: 
\begin{enumerate}
\item[\rm{(i)}] $\abs{\exc^{\operatorname{u}}(A)} \leq r$;
\item[\rm{(ii)}] $\dim_{C[\imag]} \ker_{\Univ} A=r\ \Longrightarrow\ 
\exc^{\operatorname{u}}(A)=v_{\g}(\ker_{\Univ}^{\neq} A)$; 
\item[\rm{(iii)}] under the assumption that $\fv:=\fv(A)\prec^\flat 1$ and  $B\prec_{\Delta(\fv)} \fv^{r+1}A$ where~$B\in K[\der]$ has order~$\leq r$,
we have $\exc^{\operatorname{u}}(A+B)=\exc^{\operatorname{u}}(A)$;
\item[\rm(iv)] for $r=1$ we have $\abs{\exc^{\operatorname{u}}(A)} = 1$ and $ \exc^{\operatorname{u}}(A)=v_{\g}(\ker_{\Univ}^{\neq} A)$.
\end{enumerate}
\end{lemma}

\begin{proof}
For the displayed equalities, see Remark~\ref{rem:excu for different Q, 1}. Now assume~$\I(K)\subseteq K^\dagger$. Then~$K^\dagger=H+\I(H)\imag$, 
%by Lemma~\ref{lem:W and I(F)} and remarks preceding it,  
so (i) and (ii) follow from  Proposition~\ref{prop:finiteness of excu(A), real}
and (iii) from Proposition~\ref{prop:stability of excu, real}. Corollary~\ref{cor:excu independent of Q} yields (iv). 
\end{proof}

%\noindent
%For $A$ of order $1$ we have the following variant thanks to
%Corollaries~\ref{cor:stability of excu} and~\ref{cor:excu independent of Q}:

%\begin{cor}\label{cor:excu properties} 
%Suppose $A\in K[\der]^{\neq}$ has order~$1$ and $\I(K)\subseteq K^\dagger$. Then:
%\begin{enumerate}
%\item[\rm{(i)}] $\abs{\exc^{\operatorname{u}}(A)} = 1$;
%\item[\rm{(ii)}] 
%$v_{\g}(\ker_{\Univ}^{\neq} A) = \exc^{\operatorname{u}}(A)$;  
%\item[\rm{(iii)}] if $\fv:=\fv(A)\prec^\flat 1$ and  $B\prec_{\Delta(\fv)} \fv^{2}A$ where $B\in K[\der]$ has order~$\leq 1$,
%then $\exc^{\operatorname{u}}(A+B)=\exc^{\operatorname{u}}(A)$.
%\end{enumerate}
%\end{cor}

\noindent
Recall from Lemma~\ref{lem:ADH 14.2.5} that if $K$ is $1$-linearly newtonian, then $\I(K)\subseteq K^\dagger$. 

Suppose $\I(K)\subseteq K^\dagger$. Then
$K^\dagger=H+\I(H)\imag$, so our $\Lambda$ has the form~$\Lambda_H\imag$ with~$\Lambda_H$ a complement of $\I(H)$ in $H$. Conversely, any complement~$\Lambda_H$ of~$\I(H)$ in $H$ yields a complement $\Lambda=\Lambda_H\imag$ of $K^\dagger$ in $K$ with $\Lambda\subseteq H\imag$. Now $\I(H)$ is a $C$-linear subspace of $H$, so 
$\I(H)$ has a complement $\Lambda_H$ in $H$ that is a $C$-linear subspace of $H$, and then $\Lambda:=\Lambda_H\imag$ is also a $C$-linear subspace of $K$. 

\begin{lemma}\label{lladd}
Suppose $\I(K)\subseteq K^\dagger$ and $g\in K$, $g-\lambda\in K^\dagger$. Then 
$$\Im g \in \I(H)\ \Longleftrightarrow\  
\lambda=0, \qquad \Im g \notin \I(H)\ \Longrightarrow\ \lambda \sim (\Im g)\imag.$$
\end{lemma}
\begin{proof}
Recall that $\Lambda=\Lambda_H\imag$ where $\Lambda_H$ is a complement of $\I(H)$ in $H$, so $\lambda=\lambda_H\imag$ where $\lambda_H\in\Lambda_H$. Also, $K^\dagger=H\oplus\I(H)\imag$, hence $\Im(g)-\lambda_H\in \I(H)$; this proves the displayed equivalence. Suppose $\Im g \notin \I(H)$; since $\I(H)$  is an $\mathcal O_H$-submodule of $H$ and $\lambda_H\notin\I(H)$, we then have $\Im(g)-\lambda_H\prec\lambda_H$,
so $\lambda =\lambda_H\imag\sim \Im(g)\imag$.
\end{proof}

\begin{cor}\label{cor:evs v-small}
 Suppose $\I(K)\subseteq K^\dagger$,  $A\in K[\der]^{\neq}$ has order $r$, $\dim_{C[\imag]}\ker_{\Univ}A=r$, and
 $\lambda$ is an eigenvalue of $A$ with respect to $\Lambda$. Then $\lambda\preceq\fv(A)^{-1}$.
\end{cor}
\begin{proof}
Take $f\neq 0$ and $g_1,\dots,g_r$ in $K$ with   $A=f(\der-g_1)\cdots(\der-g_r)$.
By Corollary~\ref{cor:bound on linear factors} we have $g_1,\dots,g_r\preceq\fv(A)^{-1}$, and
so Corollary~\ref{corbasiseigenvalues} gives $j\in\{1,\dots,r\}$ with $g_j-\lambda\in K^\dagger$. Now use Lemma~\ref{lladd}.
\end{proof}

\subsection*{Ultimate slots in $H$}
{\it In this subsection  $a$, $b$ range over $H$.}\/
Also, $(P,\fm,\hat a)$ is a slot in $H$ of order $r\geq 1$, where $\hat a\in \hat H\setminus H$. Recall that~$L_{P_{\times\fm}}=L_P\fm$, so
if $(P,\fm,\hat a)$ is normal, then $L_P$ has order $r$. 

\begin{cor} \label{4.7} 
Suppose $\I(K)\subseteq K^\dagger$  and the slot $(P,\fm,\hat a)$ is split-normal with linear part $L:=L_{P_{\times\fm}}$.
Then with~$Q$ and~$R$ as in~\textup{(SN2)} 
%\textup{(}but with~$H$ in the role of~$K$ there\textup{)} 
we have $\exc^{\operatorname{u}}(L)=\exc^{\operatorname{u}}(L_Q)$.
\end{cor}

\noindent
This follows from Lemmas~\ref{splnormalnormal} and~\ref{lem:excu properties}(iii).  In a similar vein we have an analogue of Lemma~\ref{lem:excev normal}:

\begin{lemma}\label{lem:excu normal}
Suppose $(P,\fm,\hat a)$ is normal and $a\prec\fm$.  Then $L_P$ and $L_{P_{+a}}$ have or\-der~$r$, and if $\I(K)\subseteq K^\dagger$, then $\exc^{\operatorname{u}}(L_{P})=\exc^{\operatorname{u}}(L_{P_{+a}})$. 
\end{lemma}
\begin{proof} 
% \marginpar{shortened proof}
%This is shown in a similar way as Lemma~\ref{lem:excev normal}, using
%Lemma~\ref{lem:excu properties}(iii) and Corollary~\ref{cor:excu properties} in place of
%Lemma~\ref{cor:excev stability} and Corollary~\ref{cor:excev stability, r=1}, respectively.
We have $L_{P_{\times \fm}}=L_P\fm$ and  $L_{P_{+a,\times \fm}}=L_{P_{\times \fm, +a/\fm}}= L_{P_{+a}}\fm$.  The slot 
$(P_{\times\fm},1,\hat a/\fm)$ in $H$ is normal and $a/\fm\prec 1$.   Lemma~\ref{lem:linear part, split-normal, new} applied to  $\hat H$, $P_{\times\fm}$, $\hat a/\fm$ in place of $K$, $P$, $a$, respectively,  gives: $L_P$ and $L_{P_{+a}}$ have order $r$, and 
$$L_P\fm - L_{P_{+a}}\fm\ =\ L_{P_{\times\fm}} - L_{P_{\times\fm,+a/\fm}}\ \prec_{\Delta(\fv)}\ \fv^{r+1}L_{P}\fm$$ where $\fv:=\fv(L_{P}\fm)\prec^\flat 1$ by (N1). Suppose now that  $\I(K)\subseteq K^\dagger$. Then
$$\exc^{\operatorname{u}}(L_P)\ =\ \exc^{\operatorname{u}}(L_P\fm)+v(\fm)\ =\ \exc^{\operatorname{u}}(L_{P_{+a}}\fm) + v(\fm)\ =\  \exc^{\operatorname{u}}(L_{P_{+a}})$$
by Lemma~\ref{lem:excu properties}(iii).
\end{proof}

\noindent
The notion introduced below is modeled on that of ``isolated slot'' (Definition~\ref{def:isolated}): 

\begin{definition}\label{def:ultimate}
Call $(P,\fm,\hat a)$   {\bf ultimate} if for all $a\prec\fm$,\index{ultimate!slot}\index{slot!ultimate}
$$\order(L_{P_{+a}})=r\ \text{ and }\ \exc^{\operatorname{u}}(L_{P_{+a}}) \cap v(\hat a-H)\ <\  v(\hat a-a);$$
equivalently, for all $a\prec \fm$:  $\order(L_{P_{+a}})=r$ and whenever
$\fw \preceq \hat a-a$ is such that~$v(\fw) \in \exc^{\operatorname{u}}(L_{P_{+a}})$, then 
$\fw\prec \hat a-b$ for all $b$. (Thus if $(P,\fm,\hat a)$ is ultimate, then it is isolated.) 
\end{definition}

\noindent
If $(P,\fm,\hat a)$ is ultimate, then so is every equivalent slot in $H$ and $(bP,\fm,\hat a)$ for~$b\neq 0$, as well as the slot~$(P^\phi,\fm,\hat a)$ in $H^\phi$ (by Lemma~\ref{lem:excu properties}). The proofs of the next two lemma are like those of their ``isolated'' versions, Lemmas~\ref{lem:isolated refinement} and~\ref{lem:isolated}:

\begin{lemma}\label{lem:ultimate refinement}
If $(P,\fm,\hat a)$ is ultimate, then so is any of its refinements. 
\end{lemma}
%\begin{proof} For the case $\fn=\fm$, use $v\big( (\hat a-a) - H \big) = v(\hat a-H)$. The case $a=0$ is clear. The general case reduces to these two special cases.
%\end{proof}

\begin{lemma}\label{lem:ultmult}
If  $(P,\fm,\hat a)$ is ultimate, then so is any of its multiplicative con\-ju\-gates.
\end{lemma}
%\begin{proof} Let $a\prec\fm/\fn$. Then $a\fn\prec \fm$, so $\order(L_{P_{\times \fn, +a}})=\order(L_{P_{+a\fn,\times \fn}})=\order(L_{P_{+a\fn}})=r$. Suppose $\fw\preceq (\hat a/\fn)-a$ and $v(\fw) \in \exc^{\operatorname{u}}\big(L_{P_{\times\fn,+a}}\big)$. Now $L_{P_{\times\fn,+a}} = L_{P_{+a\fn,\times\fn}} = L_{P_{+a\fn}}\fn$  and  thus $\fw\fn \preceq \hat a-a\fn$,  $v(\fw\fn)\in \exc^{\operatorname{u}}\big(P_{+a\fn}\big)$. But~$(P,\fm,\hat a)$ is ultimate, so $v(\fw\fn)> v(\hat a-H)$ and hence~$v(\fw)> v\big((\hat a/\fn)-H\big)$. Thus $(P_{\times\fn},\fm/\fn,\hat a/\fn)$  is ultimate.
%\end{proof}

\noindent
The ultimate condition is most useful in combination with other properties: 

\begin{lemma}\label{lem:ultimate normal}
If $\I(K)\subseteq K^\dagger$ and $(P,\fm,\hat a)$ is normal, then 
$$ \text{$(P,\fm,\hat a)$  is ultimate} \quad\Longleftrightarrow\quad
\exc^{\operatorname{u}}(L_P) \cap v(\hat a-H) \leq v\fm.$$
\end{lemma}
\begin{proof} Use Lemma~\ref{lem:excu normal} and the equivalence  $\hat a-a\prec\fm\Leftrightarrow a\prec\fm$.
\end{proof}

\noindent
The ``ultimate'' version of Lemma~\ref{lem:isolated deg 1} has the same proof: 

\begin{lemma}\label{lem:ultimate deg 1} 
If  $\deg P=1$, then
$$ \text{$(P,\fm,\hat a)$  is ultimate} \quad\Longleftrightarrow\quad
\exc^{\operatorname{u}}(L_P) \cap v(\hat a-H) \leq v\fm.$$
\end{lemma}
%\begin{proof}
%Use that $\order L_P=r$ and $L_{P_{+a}}=L_P$ for each $a$.
%\end{proof}

\noindent
The next proposition is the ``ultimate'' version of Proposition~\ref{prop:achieve isolated}:

\begin{prop}\label{prop:achieve ultimate}
Suppose $\I(K)\subseteq K^\dagger$, 
%$K$ is $\upo$-free or $r=1$, 
and~$(P,\fm,\hat a)$  is normal. 
Then   $(P,\fm,\hat a)$  has an ultimate refinement.
\end{prop}
\begin{proof}
Suppose $(P,\fm,\hat a)$ is not already ultimate.
Then Lem\-ma~\ref{lem:ultimate normal} gives $\gamma$ with
$$\gamma\in\exc^{\operatorname{u}}(L_P)\cap v(\hat a-H),\quad 
\gamma>v\fm.$$
%If $K$ is $\upo$-free, then $|\exc^{\operatorname{u}}(L_P)|\le r$ by Lemma~\ref{lem:excu properties}(i),  
%and if~$r=1$, then $|\exc^{\operatorname{u}}(L_P)| = 1$ by Lemma~\ref{lem:excu, r=1}. 
Lemma~\ref{lem:excu properties}(i) gives $|\exc^{\operatorname{u}}(L_P)|\le r$, 
so we can take $$\gamma\ :=\  \max\exc^{\operatorname{u}}(L_P)\cap v(\hat a-H),$$ and then $\gamma > v\fm$.
Take~$a$ and $\fn$ with $v(\hat a-a)>\gamma=v(\fn)$;
then $(P_{+a},\fn,\hat a-a)$ is a refinement of~$(P,\fm,\hat a)$ and $a\prec\fm$. 
Let $b\prec \fn$; then 
$a+b\prec\fm$, so by~Lemma~\ref{lem:excu normal},
$$\order(L_{(P_{+a})_{+b}})\ =\ r, \qquad
\exc^{\operatorname{u}}(L_{(P_{+a})_{+b}})\ =\ 
\exc^{\operatorname{u}}(L_P).$$
Also $v\big((\hat a-a)-b\big)>\gamma$, hence
$$\exc^{\operatorname{u}}\big(L_{(P_{+a})_{+b}}\big)  \cap v\big((\hat a-a)-H\big)\  =\ 
\exc^{\operatorname{u}}(L_P)\cap v(\hat a-H)\ \le\ \gamma\ <\ v\big((\hat a-a)-b\big).$$
Thus $(P_{+a},\fn,\hat a-a)$  is ultimate.
\end{proof}

\begin{remarkNumbered}\label{rem:achieve ultimate}
Proposition~\ref{prop:achieve ultimate} goes through if instead of assuming that $(P,\fm,\hat a)$  is normal, we assume
that $(P,\fm,\hat a)$  is linear. (Same argument, using Lem\-ma~\ref{lem:ultimate deg 1} in place of Lemma~\ref{lem:ultimate normal}.)
\end{remarkNumbered}

\noindent
Finally, here is a consequence of Corollaries~\ref{corevisu},~\ref{cor:excu independent of Q}, and   Lem\-ma~\ref{lem:ultimate normal}, where we recall that $\order(L_{P_{\times \fm}})=\order(L_P\fm)=\order(L_P)$: 

\begin{cor}\label{cor:ultimate order=1}
Suppose  $\I(K)\subseteq K^\dagger$ and $(P,\fm,\hat a)$ is normal of order~$r=1$. 
Then~$L_{P}=f(\der-g)$ with $f\in H^\times$, $g\in H$, and for $\mathfrak g\in H^\times$ with $\mathfrak g^\dagger=g$ we have:
$$\text{$(P,\fm,\hat a)$  is ultimate} \quad\Longleftrightarrow\quad 
\text{$(P,\fm,\hat a)$  is isolated} \quad\Longleftrightarrow\quad. 
 \text{$\mathfrak g \succeq \fm$ or $\mathfrak g \prec \hat a-H$.}$$ 
\textup{(}In particular, if~$g\in\I(H)$ and~$\fm\preceq 1$,
then~$(P,\fm,\hat a)$  is ultimate.\textup{)}
\end{cor}

\subsection*{Ultimate slots in $K$}
{\it In this subsection, $a$, $b$ range over $K=H[\imag]$.}\/
Also $(P,\fm,\hat a)$ is a slot in $K$ of order $r\geq 1$, where $\hat a\in \hat K\setminus K$.
Lemma~\ref{lem:excu normal} goes through in this setting, with $H$ in the proof replaced by $K$:

\begin{lemma}\label{lem:excu normal, K} 
Suppose $(P,\fm,\hat a)$ is normal, and $a\prec\fm$.  Then $L_P$ and $L_{P_{+a}}$ have order $r$, and if $\I(K)\subseteq K^\dagger$, then $\exc^{\operatorname{u}}(L_{P})=\exc^{\operatorname{u}}(L_{P_{+a}})$.
\end{lemma}

\noindent
We adapt Definition~\ref{def:ultimate} to slots in $K$:
call $(P,\fm,\hat a)$    {\bf ultimate} if for all $a\prec\fm$ we have $\order(L_{P_{+a}})=r$ and
$\exc^{\operatorname{u}}(L_{P_{+a}}) \cap v(\hat a-K) < v(\hat a-a)$.\index{ultimate!slot}\index{slot!ultimate} If $(P,\fm,\hat a)$ is ultimate, then it is isolated. 
Moreover, if $(P,\fm,\hat a)$ is ultimate, then so is~$(bP,\fm,\hat a)$ for $b\neq 0$  as well as
the slot~$(P^\phi,\fm,\hat a)$ in $K^\phi$. 
Lemmas~\ref{lem:ultimate refinement} and~\ref{lem:ultmult} go through in the present context,  and so do Lemmas~\ref{lem:ultimate normal} and~\ref{lem:ultimate deg 1} with $H$ replaced by $K$.   The analogue of Proposition~\ref{prop:achieve ultimate} 
follows likewise:

\begin{prop}\label{prop:achieve ultimate, K} 
If  $\I(K)\subseteq K^\dagger$ and~$(P,\fm,\hat a)$  is normal, then~$(P,\fm,\hat a)$  has an ultimate refinement.
\end{prop}

\begin{remarkNumbered}\label{rem:achieve ultimate, K} 
Proposition~\ref{prop:achieve ultimate, K} also holds if instead of assuming that $(P,\fm,\hat a)$  is normal, we assume
that $(P,\fm,\hat a)$  is linear. 
\end{remarkNumbered}

\noindent
Corollary~\ref{cor:excu independent of Q} and the $K$-versions of Lem\-mas~\ref{lem:ultimate normal} and~\ref{lem:ultimate deg 1} yield:

\begin{cor}\label{cor:ultimate order=1, K}
Suppose $\I(K)\subseteq K^\dagger$, $r=1$, and $(P,\fm,\hat a)$ is normal or linear. Then
$L_P=f(\der-g)$ with $f\in K^\times, g\in K$, and for $\mathfrak g\in H^\times$ with
$\mathfrak g^\dagger=\Re g$ we have:
$$\text{$(P,\fm,\hat a)$  is ultimate} \quad\Longleftrightarrow\quad  \text{$\mathfrak g\succeq \fm$ or 
$\mathfrak g \prec \hat a-K $.}$$ 
\textup{(}In particular, if~$\Re g\in\I(H)$ and~$\fm\preceq 1$, 
then~$(P,\fm,\hat a)$  is ultimate.\textup{)} %Also, $(P,\fm,\hat a)$ is flabby iff it is ultimate and $\fm\succ\mathfrak g$.
\end{cor}

\subsection*{Using the norm to characterize being ultimate}  We use here the  ``norm''~$\dabs{\,\cdot\,}$ on $\Univ$
and the gaussian extension $v_{\g}$ of the valuation of $K$ from Section~\ref{sec:group rings}. 

\begin{lemma}\label{uudag} For $u\in \Univ^\times$ we have $\dabs{u}^\dagger=\Re u^\dagger$.
\end{lemma}
\begin{proof} For $u=f\ex(\lambda)$, $f\in K^\times$ we have $\dabs{u}=|f|$ and
$u^\dagger=f^\dagger + \lambda$, so 
$$\dabs{u}^\dagger\ =\ |f|^\dagger\ =\ \Re f^\dagger\ =\ \Re u^\dagger,$$
using Corollary~\ref{cor:logder abs value} for the second equality. 
\end{proof}

\noindent
Using  Corollary~\ref{cor:valuation and norm}, Lemma~\ref{uudag}, and [ADH, 10.5.2(i)] we obtain:

\begin{lemma}\label{Wlem}
Let  $\mathfrak W\subseteq H^\times$ be $\preceq$-closed.
% such that $v(\mathfrak W)$ is closed upward in $\Gamma$ without a smallest element. 
Then for all $u\in\Univ^\times$,
$$\dabs{u} \in \mathfrak W \quad \Longleftrightarrow\quad
v_{\g}u \in v(\mathfrak W) \quad \Longleftrightarrow\quad 
\text{$\Re u^\dagger < \fn^\dagger$ for all $\fn\notin \mathfrak W$.}$$
\end{lemma} 

\noindent
Let $(P,\fm,\hat a)$ be a slot in $H$ of order $r\ge 1$. 
Applying Lemma~\ref{Wlem} to the set~$\mathfrak W = \{ \fw:\ \fw\prec\hat a-H\}$---
so $v(\mathfrak W)=\Gamma\setminus v(\hat a-H)$---we obtain a  reformulation of the condition ``$(P,\fm,\hat a)$ is ultimate'' 
in terms of the ``norm'' $\dabs{\,\cdot\,}$ on $\Univ$:

{\samepage
\begin{cor}
The following are equivalent \textup{(}with $a$ ranging over $H$\textup{)}:
\begin{enumerate}
\item[\rm{(i)}] $(P,\fm,\hat a)$  is ultimate;
\item[\rm{(ii)}] for  all $a\prec\fm$: $\order(L_{P_{+a}})=r$ and whenever~$u\in\Univ^\times$, $v_{\g}u \in \exc^{\operatorname{u}}(L_{P_{+a}})$, and
$\dabs{u}\preceq \hat a-a$, then  $\dabs{u}\prec \hat a-H$; 
\item[\rm{(iii)}] for  all $a\prec\fm$: $\order(L_{P_{+a}})=r$ and whenever~$u\in\Univ^\times$, $v_{\g}u \in \exc^{\operatorname{u}}(L_{P_{+a}})$, and 
$\dabs{u}\preceq \hat a-a$, then 
$\Re u^\dagger < \fn^\dagger$ for all $\fn$ with 
 $v(\fn)\in v(\hat a-H)$.
\end{enumerate}
\end{cor}}

\subsection*{Firm slots and flabby slots in $H$\astr} 
Let $(P,\fm,\hat a)$ be a slot in~$H$ of order $r\geq 1$, where~$\hat a\in \hat H\setminus H$.  We let $a$, $b$ range over $H$.

\begin{definition}
We  call~$(P,\fm,\hat a)$  {\bf firm} if for all $a\prec\fm$,
$$\order(L_{P_{+a}})=r \quad\text{and}\quad \exc^{\operatorname{u}}(L_{P_{+a}})\subseteq v(\hat a-H).$$
We call $(P,\fm,\hat a)$ {\bf flabby} if  it is not firm, that is, if there is an $a\prec\fm$ such that $\order(L_{P_{+a}})<r$, or $\order(L_{P_{+a}})=r$ and  $\gamma>v(\hat a-H)$ for some~$\gamma\in\exc^{\operatorname{u}}(L_{P_{+a}})$.\index{slot!firm}\index{slot!flabby}  
\end{definition}

\noindent
If  $(P,\fm,\hat a)$ is firm, then so are~$(bP,\fm,\hat a)$ for~$b\neq 0$ and any slot~$(P,\fm,\hat b)$ in $H$ that is
equivalent to~$(P,\fm,\hat a)$. 
%By Lemma~\ref{lem:excu properties}, 
For any $\phi$, the slot~$(P,\fm,\hat a)$ in $H$ is firm iff the slot~$(P^\phi,\fm,\hat a)$ in $H^\phi$ is firm.

\begin{lemma}\label{lem:firm refinement}
If $(P,\fm,\hat a)$ is firm, then so is any of its refinements.
If  $(P,\fm,\hat a)$ is flabby, then so is any refinement $(P_{+a},\fm,\hat a-a)$ of it.
\end{lemma}

\noindent
The proof is like that of Lemma~\ref{lem:ultimate refinement}.

\begin{lemma}\label{lem:firm mult conj}
Suppose $(P,\fm,\hat a)$ is firm. Then~$(P_{\times\fn},\fm/\fn,\hat a/\fn)$ is firm.
\end{lemma}
\begin{proof}
Let $a\prec\fm/\fn$, so $a\fn\prec\fm$ with  $L_{P_{\times\fn,+a}}=L_{P_{+a\fn}}\fn$. Since $(P,\fm,\hat a)$ is firm, this yields 
$\order(L_{P_{\times\fn,+a}})=\order(L_{P_{+a\fn}})=r$ and 
$$\exc^{\operatorname{u}}(L_{P_{\times\fn,+a}})=
\exc^{\operatorname{u}}(L_{P_{+a\fn}})-v\fn 
\subseteq v(\hat a-H)-v\fn = v\big( (\hat a/\fn)-H \big),$$
using Lemma~\ref{lem:excu properties} for the first equality.
\end{proof}

\noindent
The proofs of the next two lemmas are clear, using Lemma~\ref{lem:excu normal}  for the first one: 

\begin{lemma}\label{lem:firm normal}
If $\I(K)\subseteq K^\dagger$ and $(P,\fm,\hat a)$ is normal,  then
$$ \text{$(P,\fm,\hat a)$  is firm} \quad\Longleftrightarrow\quad
\exc^{\operatorname{u}}(L_P) \subseteq v(\hat a-H).$$
\end{lemma}

\begin{lemma}\label{lem:firm deg 1}
If  $\deg P=1$, then
$$ \text{$(P,\fm,\hat a)$  is firm} \quad\Longleftrightarrow\quad
\exc^{\operatorname{u}}(L_P)  \subseteq v(\hat a-H).$$
\end{lemma}

\begin{remarkNumbered}\label{rem:firm deg 1}
If the hypothesis of Lemma~\ref{lem:firm normal} or Lemma~\ref{lem:firm deg 1} holds, then
$$ \text{$(P,\fm,\hat a)$  is firm and ultimate} \quad\Longleftrightarrow\quad
\exc^{\operatorname{u}}(L_P) \leq v\fm,$$
as a consequence of Lemmas~\ref{lem:ultimate normal} and \ref{lem:ultimate deg 1}. 
\end{remarkNumbered}

\noindent
Lemma~\ref{lem:excu normal} yields:

\begin{cor} \label{cor:flabby refinement}
If the hypothesis of Lemma~\ref{lem:firm normal} or Lemma~\ref{lem:firm deg 1} holds and~$(P,\fm,\hat a)$ is flabby, then so is each refinement  of  $(P,\fm,\hat a)$.
\end{cor}

\subsection*{Firm slots and flabby slots in $K$\astr} 
Let now $(P,\fm,\hat a)$ be a slot in~$K$ of order~${r\geq 1}$ with~$\hat a\in \hat K\setminus K$, and let $a$, $b$ range over $K$.
We  define~$(P,\fm,\hat a)$  to be {\bf firm} if for all~$a\prec\fm$ we have
$\order(L_{P_{+a}})=r$ and $\exc^{\operatorname{u}}(L_{P_{+a}})\subseteq v(\hat a-H)$,
and we say that~$(P,\fm,\hat a)$ is {\bf flabby} if it is not firm.\index{slot!firm}\index{slot!flabby}
The results in the subsection above about a slot $(P,\fm,\hat a)$ in $H$  go through for the slot $(P,\fm,\hat a)$ in $K$, replacing~$H$,~$\hat H$ by~$K$,~$\hat K$ throughout.  
%From Lemma~\ref{lem:excu normal} and and~\ref{cor:ultimate order=1, K}, we obtain: 

\begin{cor}\label{cor:flabby, r=1}
Suppose  $\I(K)\subseteq K^\dagger$, $r=1$, and $(P,\fm,\hat a)$ is normal or linear.
Then $L_P=f(\der-g)$ with $f\in K^\times,\ g\in K$. For $\mathfrak g\in H^\times$ with
$\mathfrak g^\dagger=\Re g$ we have:\begin{enumerate}
\item[\rm{(i)}] $(P,\fm,\hat a)$  is flabby $\quad\Longleftrightarrow\quad \mathfrak g \prec \hat a-K\quad \Longrightarrow\quad (P,\fm,\hat a)$  is  ultimate;
\item[\rm{(ii)}]  $(P,\fm,\hat a)$  is firm and ultimate $\quad\Longleftrightarrow\quad \mathfrak g\succeq \fm$;
\item[\rm{(iii)}] $\mathfrak g\succeq 1\  \Longleftrightarrow\  \Re g\in\I(H) \text{  or } \Re g>0$. 
\end{enumerate} 
\end{cor}
\begin{proof} The equivalence in (i) follows from  Corollary~\ref{cor:excu independent of Q} and the $K$-versions of Lemmas~\ref{lem:firm normal}  and~\ref{lem:firm deg 1}. Corollary~\ref{cor:ultimate order=1, K} yields the last part of (i). 
For (ii), use Corollary~\ref{cor:excu independent of Q} and the $K$-version of the equivalence in Remark~\ref{rem:firm deg 1}. As to (iii), this is an elementary fact about the relation between $\mathfrak{g}\in H^\times$ and
$\mathfrak{g}^\dagger$. 
\end{proof} 

\noindent
For the significance of firm slots in the Hardy field setting, see Section~\ref{sec:holes perfect} below.

\subsection*{Counterexamples\astr} 
Suppose  $\I(K)\subseteq K^\dagger$ and $H$ is not  $\upo$-free. 
%(Lemma~\ref{lem:achieve I(K) subseteq Kdagger, Hfield}, Proposition~\ref{prop:Gehret}, and Example~\ref{ex:Gehret} 
%show how to obtain such $H$. 
(In Example~\ref{ex:counterex} we provide an $H$ with these properties.)
Let $(\upl_{\rho})$ and~$(\upo_{\rho})$ be as in Lemma~\ref{lem:upl-free, not upo-free} with
$H$ in the role of $K$ there. That lemma yields a minimal hole $(P,\fm, \upl)$ in $H$ with $P= 2Y'+Y^2+\upo$ ($\upo\in H$). This is a good source of counterexamples: 

\begin{lemma}\label{lem:counterex}
The minimal hole $(P,\fm,\upl)$ in $H$ is ultimate,  and none of its refinements is quasilinear,   normal, or firm.
\end{lemma}
\begin{proof} 
Let $a\in H$.  Then
$P_{+a} = 2Y'+2aY+Y^2 + P(a)$
and thus~$L_{P_{+a}}=2(\der+a)$, so for $b\in H^\times$ with $b^\dagger=-a$ we have
$\exc^{\operatorname{u}}(L_{P_{+a}})=\{vb\}$, by Corollary~\ref{cor:excu independent of Q}. 
Thus~$(P,\fm,\upl)$ is ultimate iff $\upl-a\prec b$ for all $a\prec\fm$ in $H$ and $b\in H^\times$ with~$b^\dagger=-a$
and $vb\in v(\upl-H)$; the latter holds by~[ADH, 11.5.6] since $v(\upl-H)=\Psi$.
Hence~$(P,\fm,\upl)$ is ultimate. No refinement of $(P,\fm,\upl)$ is    quasilinear  by Corollary~\ref{cor:quasilinear refinement} and
[ADH, 11.7.9], and so by Corollary~\ref{cor:normal=>quasilinear}, no refinement of $(P,\fm,\upl)$ is normal.

It remains to show that no refinement of $(P,\fm,\upl)$ is firm.   
Let $(\ell_\rho)$, $(\upg_\rho)$,  be the sequences from [ADH, 11.5] that give rise to $\upl_{\rho}=-\upg_{\rho}^\dagger$
 with $H$ in place of $K$.
If~$(P,\fm,\upl)$ has a firm refinement, then it has a firm refinement
$(P_{+\upl_{\rho}},\upg_{\rho},\upl-\upl_\rho)$,
by Lemmas~\ref{lem:pc vs dent} and~\ref{lem:firm refinement}, so it suffices
that $(P_{+\upl_{\rho}},\upg_{\rho},\upl-\upl_\rho)$ is flabby for all~$\rho$. 
For $a\in H$ we have $L_{P_{+(\upl_{\rho}+a)}}=2(\der+\upl_{\rho}+a)$, so
$\exc^{\operatorname{u}}(L_{P_{+(\upl_{\rho}+a)}})=\{vb\}$ with $b\in H^\times$, $b^\dagger=-(\upl_{\rho}+a)$.
Also $v\big((\upl-\upl_{\rho})-H\big)=v(\upl-H)=\Psi$. 
Hence $(P_{+\upl_{\rho}},\upg_{\rho},\upl-\upl_\rho)$ is flabby if  there is $a\prec \upg_{\rho}$ in $H$ and $b\in H^\times$, not active in $H$, such that $b^\dagger=-(\upl_\rho+a)$. 
We  take  $a:=2\upg_{\rho+1}$, $b:=\upg_\rho/\ell_{\rho+1}^2$.  Then
$b^\dagger=\upg_\rho^\dagger-2\ell_{\rho+1}^\dagger=-(\upl_\rho+a)$ as required.
Also, $b$ is not active in $H$.
To see this let $\sigma>\rho+1$.
Then $\upg_{\rho}/\upg_{\rho+1}, \upg_{\rho+1}/\upg_{\sigma} \succ 1$ and
$$ (\upg_{\rho}/\upg_{\rho+1})^\dagger\  =\  \upl_{\rho+1}-\upl_{\rho}\ \sim\ \upg_{\rho+1}\ \succ\ \upg_{\rho+2}\ \sim\ \upl_{\sigma}-\upl_{\rho+1}\ =\ 
 (\upg_{\rho+1}/\upg_{\sigma})^\dagger $$
by [ADH, 11.5.2], hence $\upg_{\rho}/\upg_{\rho+1} \succ \upg_{\rho+1}/\upg_{\sigma}$. Also
$\ell_{\rho+1}\asymp\upg_\rho/\upg_{\rho+1}$ by [ADH, proof of 11.5.2] and thus
$b=\upg_\rho/\ell_{\rho+1}^2 \asymp \upg_{\rho+1}^2/\upg_\rho\prec\upg_\sigma$.
\end{proof}

\section{Repulsive-Normal Slots}\label{sec:repulsive-normal}

\noindent
{\it In this section   $H$ is a real closed $H$-field with
small derivation and asymptotic integration, with $\Gamma:= v(H^\times)$.
Also $K:= H[\imag]$ with $\imag^2=-1$ is an algebraic closure of $H$.}\/ We study here the concept of a repulsive-normal slot in $H$,  which strengthens that of a  split-normal slot in $H$.  
Despite their name, 
repulsive-normal slots will turn out to have attractive analytic properties in the realm of Hardy fields.

\subsection*{Attraction and repulsion}
In this subsection~$a$, $b$ range over $H$,  $\fm$,~$\fn$  over $H^\times$, $f$,~$g$,~$h$ (possibly with subscripts)   over $K$, and $\gamma$, $\delta$ over $\Gamma$.\index{element!repulsive}\index{repulsive!element}\index{element!attractive} 
We say that~$f$ is {\bf attractive} if $\Re f\succeq 1$ and $\Re f < 0$, and
{\bf repulsive} if $\Re f\succeq 1$ and~$\Re f > 0$.
%If~$H\supseteq \R$ is Hardy field, this agrees with the terminology introduced at the end of the subsection ``Twisted integration'' of  [{\tt ahardy}, Section~1]. 
If~$\Re f \sim \Re g$, then $f$ is attractive iff $g$ is attractive, and likewise with
``repulsive'' in place of ``attractive''.
Moreover, if $a>0$, $a\succeq 1$, and $f$ is attractive (repulsive), then~$af$ is attractive (repulsive, respectively). 

\begin{definition}
Let $\gamma>0$; we say $f$ is  {\bf $\gamma$-repulsive}\index{element!gamma-repulsive@$\gamma$-repulsive}\index{S-repulsive@$S$-repulsive!element} if $v(\Re f) < \gamma^\dagger$ or~${\Re f > 0}$.
Given $S\subseteq\Gamma$, we say $f$ is {\bf $S$-repulsive} if $f$ is
$\gamma$-repulsive for all $\gamma\in S\cap\Gamma^>$, equivalently,
$\Re f >0$, or $v(\Re f) < \gamma^\dagger$ for all $\gamma\in S\cap \Gamma^{>}$.
\end{definition}

\noindent
Note the following implications for $\gamma>0$:
\begin{align*} \text{$f$ is $\gamma$-repulsive}\quad &\Longrightarrow\quad   \Re f\neq 0,\\
 \text{$f$ is $\gamma$-repulsive, $\Re g\sim\Re f$}\quad &\Longrightarrow\quad  \text{$g$ is $\gamma$-repulsive}.
\end{align*} 
The following is easy to show:

\begin{lemma}\label{lem:repulsive 1}
Suppose $\gamma>0$ and $\Re f \succeq 1$. Then    $f$ is $\gamma$-repulsive iff $v(\Re f) < \gamma^\dagger$ or~$f$ is repulsive.
Hence, if $f$ is repulsive, then  $f$ is $\Gamma$-repulsive; the converse of this implication holds if $\Psi$ is not
bounded from below in $\Gamma$.
\end{lemma}

\noindent  
Let $\gamma, \delta>0$. If $f$ is $\gamma$-repulsive and $a>0$, $a\succeq 1$, then $af$ is $\gamma$-repulsive.
If $f$ is $\gamma$-repulsive and $\delta$-repulsive, then $f$ is $(\gamma+\delta)$-repulsive. If $f$ is $\gamma$-repulsive
and~$\gamma>\delta$, then  $f$ is $(\gamma-\delta)$-repulsive. Moreover:

\begin{lemma}\label{lem:repulsive 2}
Suppose $\gamma\geq\delta=v\fn>0$. Set $g:=f-\fn^\dagger$. Then:
$$\text{$f$ is $\gamma$-repulsive}\quad\Longleftrightarrow\quad\text{$f$ is $\delta$-repulsive and $g$ is $\gamma$-repulsive.}$$
\end{lemma} 
\begin{proof} Note that $\gamma\geq\delta>0$ gives $\gamma^\dagger \leq \delta^\dagger$.
Suppose $f$ is $\gamma$-repulsive;   by our remark, $f$ is $\delta$-repulsive. Now if $v(\Re f) < \gamma^\dagger$, then $\Re g\sim\Re f$,   whereas if~${\Re f>0}$, then~$\Re(g)=\Re(f)-\fn^\dagger > \Re(f)>0$; in both cases, $g$ is $\gamma$-repulsive.
Conversely, suppose $f$ is $\delta$-repulsive and $g$ is $\gamma$-repulsive.
If $\Re f>0$, then clearly~$f$ is $\gamma$-repulsive. Otherwise, $v(\Re f)<\delta^\dagger$, hence $\Re g\sim \Re f$, so $f$ is also $\gamma$-repulsive.
\end{proof}

\noindent
In a similar way we deduce a useful characterization of repulsiveness:

\begin{lemma}\label{lem:repulsive 3}
Suppose $\gamma=v\fm>0$. Set $g:=f-\fm^\dagger$. Then:
$$\text{$f$ is repulsive}\quad\Longleftrightarrow\quad\text{$\Re f\succeq 1$, $f$ is $\gamma$-repulsive, and
 $g$ is repulsive.}$$
\end{lemma}
\begin{proof}
Suppose $f$ is repulsive; then by Lemma~\ref{lem:repulsive 1}, $f$ is $\gamma$-repulsive. Moreover,
$\Re g=\Re(f)-\fm^\dagger > \Re f>0$, hence $\Re g\succeq 1$ and
$\Re g>0$, that is, $g$ is repulsive.
Conversely, suppose $\Re f\succeq 1$, $f$ is $\gamma$-repulsive, and $g$ is repulsive.
If $v(\Re f) < \gamma^\dagger$, then $\Re f\sim \Re g$;
otherwise $\Re f>0$. In both cases,   $f$ is repulsive.
\end{proof}

\begin{cor}\label{cor:repulsive}
Suppose $f$ is $\gamma$-repulsive where $\gamma=v\fm>0$, and $\Re f\succeq 1$. Then~$f$ is repulsive iff $f-\fm^\dagger$ is repulsive, and $f$ is attractive iff $f-\fm^\dagger$ is attractive.
\end{cor}
\begin{proof}
The first equivalence is immediate from Lemma~\ref{lem:repulsive 3};
this equivalence yields 
\begin{align*}
\text{$f$ is attractive} &\ \Longleftrightarrow\  \text{$f$ is not repulsive} 
\ \Longleftrightarrow\  \text{$f-\fm^\dagger$ is not repulsive} \\
&\ \Longleftrightarrow \ \text{$\Re(f)-\fm^\dagger\prec 1$ or $f-\fm^\dagger$ is attractive.}
\end{align*}
Thus if $f-\fm^\dagger $ is attractive, so is $f$. Now assume towards a contradiction that~$f$ is attractive and $f-\fm^\dagger$ is not. 
Then $\Re f <0$ and $\Re(f) -\fm^\dagger \prec 1$ by the above equivalence, so
 $\Re f\sim\fm^\dagger$ thanks to $\Re f\succeq 1$. But $f$ is $\gamma$-repulsive, that is, 
 $\Re f\succ \fm^\dagger$ or $\Re f >0$, a contradiction.
\end{proof}

\begin{lemma}\label{lem:make repulsive}
Suppose  $\gamma=v\fm>0$ and $v(\Re g)\geq\gamma^\dagger$. Then for all sufficiently large~$c\in C^{>}$ we have 
$\Re(g)-c\fm^\dagger > 0$ \textup{(}and hence $g-c\fm^\dagger$ is  $\Gamma$-repulsive\textup{)}.
\end{lemma}
\begin{proof}
If $v(\Re g)>\gamma^\dagger$, then $\Re(g)-c\fm^\dagger\sim-c\fm^\dagger>0$ for all $c\in C^>$.
Suppose~$v(\Re g)=\gamma^\dagger$. Take $c_0\in C^\times$ with $\Re g\sim c_0\fm^\dagger$;
then  $\Re(g)-c\fm^\dagger >0$ for~$c>c_0$.
\end{proof}

\noindent
{\em In the rest of this subsection we assume that $S\subseteq \Gamma$}. If $f$ is $S$-repulsive, then so is~$af$ for $a>0$, $a\succeq 1$.
If $S>0$, $\delta>0$, and $f$ is $S$-repulsive and $\delta$-repulsive, then $f$ is $(S+\delta)$-repulsive.

\begin{lemma}\label{lem:repulsive hata, 1}
Suppose $f$ is $S$-repulsive and $0<\delta=v\fn \in S$. Then
\begin{enumerate}
\item[\textup{(i)}] $f$ is $(S-\delta)$-repulsive;
\item[\textup{(ii)}]   $g:=f-\fn^\dagger$ is $S$-repulsive.
\end{enumerate}
\end{lemma}
\begin{proof}
Let $\gamma\in (S-\delta)$, $\gamma>0$.
Then $\gamma+\delta\in S$,  so
$f$ is $(\gamma+\delta)$-repulsive, hence $\gamma$-repulsive.
This shows (i).
For (ii), suppose   $\gamma\in S$, $\gamma>0$; we  need to show that  $g$ is $\gamma$-repulsive.
If $\gamma\geq\delta$, then $g$ is $\gamma$-repulsive
by Lemma~\ref{lem:repulsive 2}. Taking $\gamma=\delta$ we see that $g$ is $\delta$-repulsive, hence if $\gamma<\delta$, then
$g$ is also $\gamma$-repulsive.
\end{proof}

\noindent
Let $A\in K[\der]^{\neq}$ have order $r\ge 1$. An {\bf $S$-repulsive splitting}\index{linear differential operator!S-repulsive splitting@$S$-repulsive splitting}\index{splitting!S-repulsive@$S$-repulsive}\index{S-repulsive@$S$-repulsive!splitting} of $A$ over $K$ is a splitting $(g_1,\dots,g_r)$ of $A$ over $K$ where  $g_1,\dots,g_r$ are $S$-repulsive.
An $S$-repulsive splitting of $A$ over $K$ remains an
$S$-repulsive splitting of~$hA$ over $K$ for $h\neq 0$.
We say that  {\bf $A$ splits $S$-repulsively} over $K$ if
there is an $S$-repulsive splitting of $A$ over $K$.
%if $\hat a$ is understood from the context, we also just say that ``$A$ splits repulsively over~$K[\imag]$''.
From Lemmas~\ref{lem:split and twist} and~\ref{lem:repulsive hata, 1} we obtain:

\begin{lemma}\label{lem:repulsive hata, 2}
Suppose $(g_1,\dots,g_r)$ is an $S$-repulsive splitting of $A$ over $K$ and ${0<\delta=v\fn \in S}$.
Then $(g_1,\dots,g_r)$ is  an $(S-\delta)$-repulsive splitting of $A$ over~$K$, and~$(h_1,\dots,h_r):={(g_1-\fn^\dagger,\dots,g_r-\fn^\dagger)}$ is an $S$-repulsive splitting of
$A\fn$ over $K$.  
\textup{(}Hence~${(h_1,\dots,h_r)}$  is also an  $(S-\delta)$-re\-pul\-sive splitting of $A\fn$  over $K$.\textup{)}
\end{lemma}
%\begin{proof}
%Take $f$ with $A=f(\der-g_1)\cdots(\der-g_r)$.   Then $$A \fn\ =\ \fn A_{\ltimes \fn}\ =\ f \fn (\der-h_1)\cdots (\der-h_r).$$
%Now use Lemma~\ref{lem:repulsive hata, 1}. 
%\end{proof}

\noindent
Note that if $\phi$ is active in $H$ with $0<\phi\preceq 1$, and $f$ is $\gamma$-repulsive (in~$K$), then~$\phi^{-1}f$ is $\gamma$-repulsive in $K^\phi=H^\phi[\imag]$.

\begin{lemma}\label{lem:repulsive splitting compconj} 
Suppose $(g_1,\dots,g_r)$ is an $S$-repulsive splitting of $A$ over $K$  and ${S\cap \Gamma^{>}}\not\subseteq\Gamma^\flat$.
Let~$\phi$ be active in $H$ with $0<\phi\prec 1$, and set $h_j:=g_j-(r-j)\phi^\dagger$ for~$j=1,\dots,r$. 
Then $(\phi^{-1}h_1,\dots,\phi^{-1}h_r)$ is an $S$-repulsive splitting of $A^\phi$ over~$K^\phi$.
\end{lemma}
\begin{proof}
By Lem\-ma~\ref{lem:split and compconj}, $(\phi^{-1}h_1,\dots,\phi^{-1}h_r)$ is   splitting of $A^\phi$ over~$K^\phi$. 
Let $j\in\{1,\dots,r\}$. If $\Re g_j>0$, then $\phi^\dagger < 0$ yields  $\Re h_j\geq \Re g_j>0$.
Otherwise, 
$v(\Re g_j)<\gamma^\dagger$ whenever $0<\gamma\in S$; in particular, $\Re g_j\succ 1\succ\phi^\dagger$,
so $\Re h_j\sim  \Re g_j$. In both cases $h_j$ is $S$-repulsive, so
$\phi^{-1}h_j$ is $S$-repulsive in $K^\phi$.
\end{proof}

\begin{prop}\label{prop:repulsive splitting}
Suppose $S\cap \Gamma^{>}\neq\emptyset$, $nS\subseteq S$ for all $n\geq 1$,   the ordered constant field~$C$ of~$H$ is ar\-chi\-me\-dean, and 
$(g_1,\dots,g_r)$ is a splitting of $A$ over $K$.  Then there exists~$\gamma\in S\cap \Gamma^{>}$
such that for any $\fm$ with $\gamma=v\fm$: $({g_1-n\fm^\dagger},\dots,{g_r-n\fm^\dagger})$ is an $S$-repulsive splitting of
$A\fm^n$ over $K$, for all big enough $n$.
\end{prop}
\begin{proof}
%Take $f\neq 0$ and $g_1,\dots,g_r$  such that $A=f(\der-g_1)\cdots(\der-g_r)$, and 
Let~$J$ be the set of $j\in\{1,\dots,r\}$ such that
$g_j$ is not $S$-repulsive. If  $\gamma >0$ and $g$ is not $\gamma$-repulsive, then~$g$ is not
$\delta$-repulsive, for all $\delta\geq\gamma$. 
Hence we can take $\gamma\in S\cap \Gamma^{>}$ such that $g_j$ is not $\gamma$-repulsive, for all $j\in J$.  
Suppose $\gamma=v\fm$.
Lemma~\ref{lem:make repulsive} yields $m\geq 1$ such that for all $n\geq m$, setting $\fn:=\fm^n$,
$g_j-\fn^\dagger$ is $\Gamma$-repulsive for all $j\in J$. For such $\fn$ we have $v\fn\in S$, so
by Lemma~\ref{lem:repulsive hata, 1}(ii), $g_j-\fn^\dagger$ is also $S$-repulsive for $j\notin J$. %and thus~$A\fn$ splits $S$-repulsively over $K[\imag]$.
\end{proof}

%\noindent
%If $A$ splits $S$-repulsively over $K[\imag]$, then for any $\widetilde{S}\subseteq S$, $A$ also splits $\widetilde{S}$-repulsively over $K[\imag]$.
%By a {\bf repulsive splitting} of $A$ over $K[\imag]$ we mean a $\Gamma$-repulsive splitting of $A$ over $K[\imag]$, and
%if $A$ splits $\Gamma$-repulsively over $K[\imag]$, then we say that {\bf $A$ splits repulsively} over~$K[\imag]$.

\begin{cor}\label{cor:repulsive splitting}
If  $C$ is archimedean and~$(g_1,\dots,g_r)$ is a splitting of $A$ over $K$,
 then there exists $\gamma>0$ such that for all $\fm$ with $\gamma=v\fm$: ${(g_1-n\fm^\dagger,\dots,g_r-n\fm^\dagger)}$ is a $\Gamma$-repulsive splitting of
$A\fm^n$ over $K$, for all big enough $n$. If $\Gamma\neq\Gamma^\flat$ then we can choose such $\gamma>\Gamma^\flat$.
\end{cor}
\begin{proof} Taking $S=\Gamma$ this
follows from Proposition~\ref{prop:repulsive splitting} and its proof.
%, for the first part applied with
%$S=\Gamma^>$, and for the second part with $S=\Gamma^>\setminus\Gamma^\flat$.
\end{proof}

\noindent
In logical jargon, the condition that $C$ is archimedean is not {\em first-order}. But it is satisfied when $H$ is a Hardy field, the case
where the results of this section will be applied. For other possible uses we indicate here a first-order variant of Proposition~\ref{prop:repulsive splitting} with essentially the same proof: 

\begin{cor} 
Suppose  $(g_1,\dots,g_r)$ is a splitting of $A$ over $K$. Then there exists $\fm\prec 1$ such that for all sufficiently large $c\in C^{>}$ and all~$\fn$, if $\fn^\dagger=c\fm^\dagger$, then~$({g_1-\fn^\dagger},\dots,{g_r-\fn^\dagger})$ is a $\Gamma$-repulsive splitting of
$A\fn$ over $K$.
\end{cor}

\noindent
In connection with this corollary we recall from~\cite[p.~105]{AvdD3}
%Section~\ref{sec:approx linear diff ops}  
that $H$ is said to be {\it closed under powers}\/ if for all $c\in C$ and~$\fm$ there is an $\fn$
with $c\fm^\dagger=\fn^\dagger$.  

\medskip
\noindent
{\it In the rest of this section $\hat H$ is an immediate asymptotic extension of~$H$ and $\imag$ with~$\imag^2=-1$ lies in an asymptotic extension of $\hat H$. Also $K:= H[\imag]$ and $\hat K:=\hat H[\imag]$.}

\medskip
\noindent
Let $\hat a\in\hat H\setminus H$, so 
$v({\hat a-H})$
is a downward closed subset of $\Gamma$. 
We say that~$f$ is~{\bf $\hat a$-repulsive} if $f$ is~$v({\hat a-H})$-repulsive;\index{element!a-repulsive@$\hat a$-repulsive} that is, $\Re f>0$, or $\Re f \succ \fm^\dagger$ for all~$a$,~$\fm$
with $\fm\asymp\hat a-a\prec 1$. 
(Of course, this notion is only interesting if $v({\hat a-H})\cap\Gamma^>\neq \emptyset$, since otherwise every $f$ is $\hat a$-repulsive.)
Various earlier results give:

{\samepage
\begin{lemma}\label{lem:hata-repulsive}
Suppose $f$ is $\hat a$-repulsive.  Then
\begin{enumerate}
\item[\textup{(i)}] $b>0,\ b\succeq 1\ \Longrightarrow\ bf$ is $\hat a$-repulsive; 
\item[\textup{(ii)}] $f$ is   $(\hat a-a)$-repulsive;  
\item[\textup{(iii)}] $\fm \asymp 1\ \Longrightarrow\ f$ is $\hat a\fm$-repulsive; 
\item[\textup{(iv)}] $\fn\asymp \hat a-a\prec 1\ \Longrightarrow\ f$ is $\hat a/\fn$ repulsive and $f-\fn^\dagger$ is $\hat a$-repulsive.
\end{enumerate}
\end{lemma}
}

\noindent
For (iv), use Lemma~\ref{lem:repulsive hata, 1}. 
An {\bf $\hat a$-repulsive splitting} of $A$ over $K$ is a
$v({\hat a-H})$-re\-pul\-sive splitting $(g_1,\dots,g_r)$ of $A$ over $K$:\index{linear differential operator!a-repulsive splitting@$\hat a$-repulsive splitting}\index{splitting!a-repulsive@$\hat a$-repulsive} 
$$A\ =\ f(\der-g_1)\cdots(\der-g_r)\qquad \text{where $f\neq 0$ and $g_1,\dots,g_r$  are $\hat a$-repulsive.}$$
We say
that {$A$ \bf splits $\hat a$-repulsively} over $K$ if it splits $v({\hat a-H})$-repulsively over~$K$.
Thus if $A$ splits $\hat a$-repulsively over $K$, then so does $hA$ ($h\neq 0$), and 
$A$ splits $(\hat a-a)$-repulsively over $K$, and splits $\hat a\fm$-repulsively over $K$ for $\fm\asymp 1$.
Moreover, from Lemma~\ref{lem:repulsive hata, 2} we obtain:

\begin{cor}\label{cor:repulsive hata} 
Suppose $(g_1,\dots,g_r)$ is an $\hat a$-repulsive splitting of $A$ over $K$ and~$\fn\asymp \hat a-a\prec 1$. 
Then $(g_1,\dots,g_r)$ is  an $\hat a/\fn$-repulsive splitting of $A$ over~$K$ and~${(g_1-\fn^\dagger,\dots,g_r-\fn^\dagger)}$ is an $\hat a$-repulsive splitting of
$A\fn$ over $K$.  
%\textup{(}Hence~${(h_1,\dots,h_r)}$  is also an  $(S-\delta)$-re\-pul\-sive splitting of $A\fn$  over $K$.\textup{)}
%Suppose $A$ splits $\hat a$-repulsively over $K$ and $\fn\asymp \hat a-a\prec 1$. Then~$A$ splits $\hat a/\fn$-repulsively over $K$, and $A\fn$ splits $\hat a$-repulsively over $K$.
\end{cor}

\noindent
Proposition~\ref{prop:repulsive splitting} yields:

\begin{cor}\label{cor:repulsive hata, 1} 
If $\hat a\preceq 1$ is special over $H$, $C$ is archimedean,  and $A$ splits over~$K$, then $A\fn$ splits $\hat a$-repulsively over $K$ for some~$a$ and $\fn\asymp\hat a-a\prec 1$.
\end{cor}

\noindent
%An {\bf $\hat a$-repulsive splitting} of $A$ over $K$ is a  $v({\hat a-H})$-repulsive splitting of $A$ over $K$.
Recall that in Section~\ref{sec:approx linear diff ops} we defined a splitting $(g_1,\dots,g_r)$ of $A$ over $K$ to be {\it strong}\/ if~$\Re g_j \succeq\fv(A)^\dagger$ for $j=1,\dots,r$. 
%We say that  {$A$ \bf splits strongly $\hat a$-repulsively} over $K[\imag]$ if there exists a strong $\hat a$-repulsive splitting of
%$A$ over $K[\imag]$.
 
\begin{lemma}\label{lem:hata-repulsive splitting compconj}
Suppose $\hat a-a\prec^\flat 1$ for some $a$. Let $(g_1,\dots,g_r)$ be an $\hat a$-repulsive splitting of $A$ over $K$, let $\phi$ be active in $H$ with $0<\phi\prec 1$, 
and set $$h_j\ :=\ \phi^{-1}\big( g_j - (r-j) \phi^\dagger \big) \qquad (j=1,\dots,r).$$ Then
$(h_1,\dots,h_r)$ is an $\hat a$-repulsive splitting of~$A^\phi$ over $K^\phi=H^\phi[\imag]$.
If $\fv(A)\prec^{\flat} 1$ and $(g_1,\dots,g_r)$ is strong, then $(h_1,\dots,h_r)$ is strong.
\end{lemma} 

\noindent
This follows from Lemmas~\ref{lem:split strongly compconj} and~\ref{lem:repulsive splitting compconj}. 

\begin{lemma}\label{lem:achieve strong repulsive splitting} 
Suppose $\fv:=\fv(A)\prec 1$ and $\hat a \prec_{\Delta(\fv)} 1$.
Let~$(g_1,\dots,g_r)$ be an $\hat a$-repulsive splitting of $A$ over $K$. Then for all sufficiently small $q\in\Q^>$ and any~$\fn\asymp\abs{\fv}^q$, $(g_1-\fn^\dagger,\dots,g_r-\fn^\dagger)$ is a strong $\hat a/\fn$-repulsive splitting of $A\fn$ over $K$.
\end{lemma}
\begin{proof}
Take $q_0\in\Q^>$ with $\hat a\prec |\fv|^{q_0}\prec 1$. Then for any $q\in\Q$ with $0<q\leq q_0$ and any $\fn\asymp\abs{\fv}^q$,
$(g_1-\fn^\dagger,\dots,g_r-\fn^\dagger)$ is an $\hat a/\fn$-repulsive splitting of $A\fn$ over~$K$, by Corollary~\ref{cor:repulsive hata}.
Using Lemmas~\ref{lem:split strongly multconj} and~\ref{lem:Ah splits strongly} (in that order) we can decrease~$q_0$ so that for all~$q\in\Q$ with~$0<q\leq q_0$ and $\fn\asymp\abs{\fv}^q$,
$(g_1-\fn^\dagger,\dots,g_r-\fn^\dagger)$ is also a strong splitting of~$A\fn$ over $K$.
\end{proof}

\noindent
{\it In the rest of this subsection we assume that $H$ is Liouville closed with~$\I(K)\subseteq K^\dagger$.}\/
 We choose a complement~$\Lambda\subseteq H\imag$ of $K^\dagger$ in $K$ as in Section~\ref{sec:ultimate} and set $\Univ := K\big[\!\ex(\Lambda)\big]$.
 We then have the set $\exc^{\operatorname{u}}(A)\subseteq\Gamma$ of ultimate exceptional values of $A$ 
 (which doesn't depend on $\Lambda$ by Corollary~\ref{prop:excu independent of Q}).
Recall from Corollary~\ref{cor:Hardy type arch C} that $H$ is of Hardy type iff $C$ is archimedean.
{\it  We now assume $r=1$ and~$\hat a\prec 1$ is special over $H$, and let 
$\Delta$ be the nontrivial convex subgroup of $\Gamma$ that is cofinal in $v(\hat a - H)$}. 

\begin{lemma}\label{lem:aAH}
Suppose $C$ is archimedean and 
$\exc^{\operatorname{u}}(A) \cap v(\hat a-H) < 0$. Then $A$ splits $\hat a$-repulsively over $K$.
\end{lemma}
\begin{proof}
We may arrange $A=\der-f$. 
%Let $\Delta$ be the convex subgroup of $\Gamma$ which is cofinal in $v(\hat a-H)$.
Take $u\in\Univ^\times$  with $u^\dagger=f$, and set $b:=\dabs{u}\in H^>$. 
Then $\exc^{\operatorname{u}}(A)=\{vb\}$ by Lemma~\ref{lem:excu, r=1} and its proof, hence 
$$\exc^{\operatorname{u}}(A) \cap v(\hat a-H) < 0 \quad\Longleftrightarrow\quad \text{$b\succ 1$ or $vb>\Delta$,}$$ 
and $\Re f =  b^\dagger$ by Lemma~\ref{uudag}.
If $b\succ 1$, then $\Re f>0$, and if $vb>\Delta$, then for all $\delta\in\Delta^{\neq}$ we have $\psi(vb)<\psi(\delta)$ by Lemma~\ref{lem:Hardy type}, so $\Re f \succ \fm^\dagger$ for all $a$, $\fm$ with~$\hat a-a\asymp\fm\prec 1$.
In both cases $A$ splits $\hat a$-repulsively over $K$. 
\end{proof}

\begin{lemma}
Suppose $A\in H [\der]$ and $\fv(A)\prec 1$.
Then $0\notin\exc^{\operatorname{u}}(A)$, and if
$A$ splits $\hat a$-repulsively over $K$, then $\exc^{\operatorname{u}}(A) \cap v(\hat a-H) < 0$.
%$$\text{$A$ splits $\hat a$-repulsively over $K$} \ \Longleftrightarrow\ \exc^{\operatorname{u}}(A) \cap v(\hat a-H) \subseteq \Gamma^<.$$
\end{lemma}
\begin{proof}
We again arrange $A=\der-f$ and
take $u$, $b$ as in the proof of Lemma~\ref{lem:aAH}. Then
$f\in H$ and $b^\dagger=f=-1/\fv(A)\succ 1$, so $b\nasymp 1$, and thus $0\notin \{vb\}=\exc^{\operatorname{u}}(A)$. Now suppose $A$ splits $\hat a$-repulsively over $K$, that is,
$f>0$ or $f\succ\fm^\dagger$ for all $a$, $\fm$ with $\hat a-a\asymp\fm\prec 1$.
In the first case $f=b^\dagger$ and $b\nasymp 1$ yield $b\succ 1$.
In the second case  $\psi(vb)=vf<\psi(\delta)$ for all $\delta\in\Delta^{\neq}$, hence $vb\notin\Delta$.
\end{proof}

\noindent
Combining Lemma~\ref{lem:order 1 splits strongly} with the previous two lemmas yields:

\begin{cor}\label{cor:split rep ultimate}
Suppose $A\in H [\der]$ and $\fv(A)\prec 1$, and $H$ is of Hardy type. Then~$A$ splits strongly over $K$, and we have the equivalence
$$\text{$A$ splits $\hat a$-repulsively over $K$} \ \Longleftrightarrow\ \exc^{\operatorname{u}}(A) \cap v(\hat a-H)\ \leq\ 0.$$
\end{cor}

\subsection*{Defining repulsive-normality}
{\it In this subsection
 $(P, \fm, \hat a)$ is a slot in $H$ of or\-der~$r\ge 1$ with $\hat a\in \hat H\setminus H$ and linear part $L:=L_{P_{\times \fm}}$.} 
 Set $w:=\wt(P)$;
 if $\order L=r$, set $\fv:=\fv(L)$.
We let $a$, $b$ range over $H$ and $\fn$ over $H^\times$.

\begin{definition}  
Call $(P,\fm,\hat a)$   {\bf repulsive-normal} if $\order L=r$, and\index{slot!repulsive-normal}\index{repulsive-normal}
\begin{itemize}
\item[(RN1)] $\fv\prec^\flat 1$; 
\item[(RN2)]  $(P_{\times\fm})_{\geq 1}=Q+R$ where $Q,R\in H\{Y\}$, $Q$ is homogeneous of degree~$1$ and order~$r$,  $L_Q$ splits $\hat a/\fm$-repulsively over $K$, and $R\prec_{\Delta(\fv)} \fv^{w+1} (P_{\times\fm})_1$. 
\end{itemize}
\end{definition}

\noindent
Compare this with ``split-normality'' from Definition~\ref{SN}: clearly repulsive-normal implies split-normal, and hence normal.
If $(P,\fm,\hat a)$ is normal and $L$ splits $\hat a/\fm$-repulsively over $K$, then~$(P,\fm,\hat a)$ is repulsive-normal. 
If~$(P,\fm,\hat a)$ is  repulsive-normal, then so are $(bP,\fm,\hat a)$ for $b\neq 0$ and $(P_{\times\fn},\fm/\fn,\hat a/\fn)$.

\begin{lemma}\label{lem:repulsive-normal comp conj}
Suppose $(P,\fm,\hat a)$ is repulsive-normal and  $\phi$ is active in $H$ such that~${0<\phi\prec 1}$, and $\hat a-a\prec^\flat \fm$
for some $a$. Then the slot $(P^\phi,\fm,\hat a)$ in $H^\phi$ is re\-pul\-sive-nor\-mal.
\end{lemma}
\begin{proof}
First arrange $\fm=1$, and let $Q$, $R$ be as in (RN2) for $\fm=1$. 
Now $(P^\phi,1,\hat a)$ is split-normal by Lemma~\ref{lem:split-normal comp conj}. In fact, 
$P^\phi_{\geq 1} = Q^\phi+R^\phi$, and
the proof
of this lemma shows that $R^\phi \prec_{\Delta(\fw)} \fw^{w+1} P_1^\phi$
where $\fw:=\fv(L_{P^\phi})$.
By Lemma~\ref{lem:hata-repulsive splitting compconj}, $L_{Q^\phi}=L_Q^\phi$ splits $\hat a$-repulsively over $K^\phi$.
So  $(P^\phi,1,\hat a)$ is repulsive-normal.
\end{proof}

\noindent
If $\order L=r$, $\fv\prec^\flat 1$, and $\hat a-a\prec_{\Delta(\fv)} \fm$, then $\hat a - a\prec^\flat \fm$. Thus we obtain from
Lemmas~\ref{lem:good approx to hata} and \ref{lem:repulsive-normal comp conj} the following result:

\begin{cor}\label{cor:repulsive-normal comp conj}  
Suppose  $(P,\fm,\hat a)$ is $Z$-minimal, deep, and repulsive-normal. Let~$\phi$ be active in $H$ with $0<\phi\prec 1$. Then the slot $(P^\phi,\fm,\hat a)$ in $H^\phi$ is repulsive-normal.
\end{cor}

\noindent
Before we turn to the task of obtaining repulsive-normal slots, we deal with the preservation of repulsive-normality under refinements.

\begin{lemma}\label{lemrnlem}
Suppose $(P,\fm,\hat a)$ is repulsive-normal, and let $Q$, $R$ be as in~\textup{(RN2)}.
Let $(P_{+a},\fn,\hat a-a)$ be a steep refinement of $(P,\fm,\hat a)$ where $\fn\prec\fm$ or
$\fn=\fm$. Suppose
$$(P_{+a,\times\fn})_{\geq 1} - Q_{\times\fn/\fm} \prec_{\Delta(\fw)} \fw^{w+1}(P_{+a,\times\fn})_1
\qquad\text{where $\fw:=\fv(L_{P_{+a,\times\fn}})$.}$$
Then  $(P_{+a},\fn,\hat a-a)$ is repulsive-normal.
\end{lemma}
\begin{proof}
By (RN2), $L_Q$ splits $\hat a/\fm$-repulsively over $K$, so $L_Q$ also splits
$(\hat a-a)/\fm$-repulsively over $K$. We have $(\hat a-a)/\fm\prec \fn/\fm \prec 1$ or
$(\hat a-a)/\fm\prec 1=\fn/\fm$, so~$L_Q$ splits $(\hat a-a)/\fn$-repulsively over $K$ by the first part of Corollary~\ref{cor:repulsive hata}, and hence~$L_{Q_{\times\fn/\fm}}=L_Q\cdot(\fn/\fm)$  splits $(\hat a-a)/\fn$-repulsively over $K$ by the second part of that
Corollary~\ref{cor:repulsive hata}. Thus $(P_{+a},\fn,{\hat a-a})$ is repulsive-normal.
\end{proof}

\noindent
The proofs of Lemmas~\ref{splitnormalrefine},~\ref{easymultsplitnormal},~\ref{splnq} give the following repulsive-normal analogues of these lemmas, using also Lemma~\ref{lemrnlem}; for Lemma~\ref{lem:5.21 repulsive-normal} below we adopt the notational conventions about $\fn^q$ ($q\in\Q^>$) stated
before Lemma~\ref{splnq}.

\begin{lemma}\label{lem:5.18 repulsive-normal}
If $(P,\fm,\hat a)$ is repulsive-normal and $(P_{+a},\fm,\hat a-a)$ is a refinement of $(P,\fm,\hat a)$, then $(P_{+a},\fm,\hat a-a)$ is also repulsive-normal.
\end{lemma}

\begin{lemma}\label{lem:5.19 repulsive-normal}
Suppose   $(P,\fm,\hat a)$ is repulsive-normal, $\hat a \prec \fn\prec\fm$, and~${[\fn/\fm]\le\big[\fv]}$. Then the refinement $(P,\fn,\hat a)$ of $(P,\fm,\hat a)$
is repulsive-normal: if $\fm$, $P$, $Q$, $\fv$ are as in \textup{(RN2)}, then \textup{(RN2)} holds with $\fn$, $Q_{\times \fn/\fm}$, $R_{\times \fn/\fm}$, $\fv(L_{P_{\times \fn}})$ in place of~$\fm$,~$Q$,~$R$,~$\fv$. 
\end{lemma}

\begin{lemma}\label{lem:5.21 repulsive-normal}
Suppose $\fm=1$, $(P,1,\hat a)$ is repulsive-normal, $\hat a \prec \fn\prec 1$, and for~$\fv:=\fv(L_P)$ 
we have $[\fn^\dagger]<[\fv]<[\fn]$;
then $(P,\fn^q,\hat a)$ is a repulsive-normal refinement
of~$(P,1,\hat a)$ for all but finitely many~$q\in \Q$ with $0<q<1$. 
\end{lemma}

%\noindent
%We say that $(P^\phi,\fm,\hat a)$ is {\it eventually repulsive-normal}\/ if there is an active $\phi_0\preceq 1$ in~$H$ such that $(P^\phi,\fm,\hat a)$ is repulsive-normal for all active $\phi$ in $H$ with $0<\phi\prec\phi_0$; likewise for $(P^\phi,\fm,\hat a)$ being {\it eventually deep and repulsive-normal}.

\subsection*{Achieving repulsive-normality}
In this subsection we adopt the setting of the subsection {\it Achieving split-normality}\/ 
of Section~\ref{sec:split-normal holes}: {\it $H$ is $\upo$-free and $(P,\fm,\hat a)$ is a minimal hole in $K$ of order $r\geq 1$, $\fm\in H^\times$, and $\hat a\in \hat K\setminus K$, with
$\hat a = \hat b + \hat c \imag$, $\hat b, \hat c\in \hat H$.}\/ 
We let $a$ range over $K$, $b$, $c$ over $H$, and~$\fn$ over $H^\times$. We prove here the following variant of Theorem~\ref{thm:split-normal}:

\begin{theorem}\label{thm:repulsive-normal}   
Suppose the constant field $C$ of $H$ is archimedean and $\deg P>1$. Then one of the following conditions is satisfied:
\begin{enumerate}
\item[$\mathrm{(i)}$] $\hat b\notin H$ and some $Z$-minimal  slot $(Q,\fm,\hat b)$ in $H$ has a special refinement ${(Q_{+b},\fn,\hat b-b)}$ such that $(Q^\phi_{+b},\fn,\hat b-b)$ is eventually deep and repulsive-normal; 
\item[$\mathrm{(ii)}$] $\hat c\notin H$ and some $Z$-minimal  slot $(R,\fm,\hat c)$ in $H$ has a special  refinement ${(R_{+c},\fn,\hat c-c)}$ such that   $(R^\phi_{+c},\fn,\hat c-c)$ is eventually deep and repulsive-normal. 
\end{enumerate}
\end{theorem}

\noindent
To establish this theorem we need to take up the approximation arguments in the proof of Theorem~\ref{thm:split-normal} once again. While in that proof we  treated the cases~${\hat b\in H}$ and~${\hat c\in H}$ separately to obtain stronger results in those cases (Lem\-mas~\ref{lem:hat c in K}, \ref{lem:hat b in K}), 
here  we proceed  differently and first show a repulsive-normal version of  
Proposition~\ref{evsplitnormal} which also applies to those cases. {\it In the rest of
 this subsection we assume that  $C$ is archimedean.}\/

\begin{prop}\label{evrepnormal}
Suppose  the hole $(P,\fm,\hat a)$ in $K$ is special and $v({\hat b-H})\subseteq v({\hat c-H})$ \textup{(}so $\hat b\notin H$\textup{)}. Let~$(Q,\fm,\hat b)$ be a $Z$-minimal deep normal slot in~$H$. Then $(Q,\fm,\hat b)$ has a  repulsive-normal refinement.
\end{prop}
\begin{proof}
As in the proof of Proposition~\ref{evsplitnormal} we first arrange $\fm=1$, and set
$$\Delta\ :=\ \big\{ \delta\in\Gamma:\  \abs{\delta}\in v\big(\hat a-K\big) \big\},$$
a convex subgroup of~$\Gamma$ which is cofinal in $v(\hat a-K) = v({\hat b-H})$,
so $\hat b$ is special over~$H$. 
Lemma~\ref{lem:good approx to hata} applied to $(Q,1,\hat b)$ and $\fv(L_{Q})\prec^\flat 1$ gives that 
$\Gamma^\flat$ is strictly contained in $\Delta$. To show that~$(Q,1,\hat b)$ has a repulsive-normal refinement, we follow the proof of Proposition~\ref{evsplitnormal}, skipping the initial compositional conjugation, and arranging first that $P,Q\asymp 1$. Recall from that proof that $\dot{\hat a}\in \dot{K}^{\operatorname{c}}=\dot{H}^{\operatorname{c}}[\imag]$ and~$\Re \dot{\hat a}=\dot{\hat b}\in \dot{H}^{\operatorname{c}}\setminus \dot{H}$, with $\dot{\hat b} \prec 1$, $\dot{Q}\in \dot{H}\{Y\}$, and so~$\dot{Q}_{+\dot{\hat b}}\in \dot{H}^{\operatorname{c}}\{Y\}$. Let~$A\in \dot{H}^{\operatorname{c}}[\der]$ be the linear part
of~$\dot Q_{+\dot{\hat b}}$. Recall from that proof that $1\le s:=\order Q = \order A \le 2r$ and that 
$A$ splits over $\dot{K}^{\operatorname{c}}$. Then Lemma~\ref{hkspl} gives a {\em real\/} splitting
 $(g_1,\dots,g_s)$ of~$A$ over $\dot K^{\operatorname{c}}$:
$$A\ =\ f (\der-g_1)\cdots (\der-g_s), \qquad 0\ne f\in \dot{H}^{\operatorname{c}},\  g_1,\dots, g_s\in \dot K^{\operatorname{c}}.$$
It follows easily from [ADH, 10.1.8] that the real closed $\d$-valued field $\dot{H}$ is an $H$-field, and so its completion $\dot{H}^{\operatorname{c}}$ is also a real closed $H$-field by [ADH, 10.5.9]. Recall also that $\Delta=v(\dot H^\times)$ is the value group  of $\dot H^{\operatorname{c}}$ and properly contains $\Gamma^\flat$. Thus we can apply Corollary~\ref{cor:repulsive splitting} with~$\dot{H}^{\operatorname{c}}$ in the role of $H$ to get $\fn\in\dot{\mathcal O}$ with~$0\neq \dot\fn\prec^\flat 1$ and~$m$ such that for all $n > m$,
$(h_1,\dots,h_s):=({g_1-n\dot\fn^\dagger},\dots,g_s-n\dot\fn^\dagger)$ is a $\Delta$-repulsive splitting of~$A\dot\fn^n$ over $\dot K^{\operatorname{c}}$, so $\Re h_1,\dots, \Re h_s\ne 0$. For any $n$, $A\dot\fn^n$ is the 
linear part  of  $\dot Q_{+\dot{\hat b},\times\dot\fn^n}\in \dot H^{\operatorname{c}}\{Y\}$,
and  $(h_1,\dots,h_s)$ is also a real splitting
of~$A\dot\fn^n$ over $\dot K^{\operatorname{c}}$:
$$A\dot\fn^n\ =\ \dot{\fn}^n f (\der-h_1)\cdots (\der-h_s). $$
By increasing $m$ we arrange that for all $n > m$  we have
$g_j\not\sim n\dot\fn^\dagger$  ($j=1,\dots,s$), and also $\fv(A\dot\fn^n) \preceq \fv(A)$
provided $\big[\fv(A)\big] < [\dot\fn]$;
for the latter part use Lemma~\ref{lem:nepsilon}.
Below we assume $n > m$. Then $\fv(A\dot\fn^n) \prec 1$: to see this use 
Corollary~\ref{corAfv1}, $\fv(A)\prec 1$, and $g_j \preceq h_j$ $(j=1,\dots,s$).  Note  that
$h_1,\dots, h_s\succeq 1$.
We now apply 
Co\-rol\-lary~\ref{cor:approx LP+f, real, general} to
$\dot H$, $\dot K$, $\dot Q$, $s$, $\dot\fn^n$, $\dot{\hat b}$, $\dot\fn^n f$, $h_1,\dots,h_s$
in place of $H$, $K$, $P$, $r$, $\fm$, $f$, $a$, $b_1,\dots,b_r$, respectively, and any $\gamma\in\Delta$
with $\gamma>v(\dot\fn^n), v(\Re h_1),\dots,v(\Re h_s)$.
This  gives $a,b\in\dot{\mathcal O}$  and $b_1,\dots,b_s\in\dot{\mathcal O}_K$ such that $\dot a, \dot b \ne 0$ in
$\dot H$ and such that for the linear part $\tilde{A}\in \dot{H}[\der]$ of $\dot Q_{+\dot b,\times\dot\fn^n}$ we have
$$\dot b\ - \dot{\hat b} \ \prec\ \dot\fn^n,\qquad \tilde{A}\ \sim\ A\dot{\fn}^n,\qquad \order \tilde{A}\ =\ s, \qquad
\fw\ :=\ \fv(\tilde{A})\ \sim\ \fv(A\dot\fn^n),$$
and such that for $w:=\wt(Q)$ and with $\Delta(\fw)\subseteq \Delta$: 
\begin{align*} \tilde{A}\ =\  \tilde{B} +  \tilde{E},& \quad \tilde{B}\ =\ \dot a(\der-\dot b_1)\cdots(\der-\dot b_s)\in \dot{H}[\der],\quad   \tilde{E}\in \dot{H}[\der],\\
   v(\dot b_1-  h_1),&\dots,v(\dot b_s- h_s)\ >\ \gamma,\quad   \tilde{E}\ \prec_{\Delta(\fw)}\  \fw^{w+1} \tilde{A},  
   \end{align*}
and $(\dot{b}_1,\dots, \dot{b}_s)$ is a real splitting of $\tilde{B}$ over $\dot{K}$.  This real splitting over $\dot{K}$ has a consequence that will be crucial at the end of the proof: by changing $b_1,\dots, b_s$ if necessary, without changing $\dot{b}_1,\dots, \dot{b}_s$ we arrange that $B:=a(\der-b_1)\cdots (\der-b_s)$ lies in $\dot{\mathcal{O}}[\der]\subseteq H[\der]$ and that
$(b_1,\dots, b_s)$ is a real splitting of $B$ over $K$.   (Lemma~\ref{lem:lift real splitting}.)

 Since $\Re\dot{b}_1\sim \Re h_1,\dots, \Re \dot{b}_s\sim \Re h_s$,
 the implication just before Lemma~\ref{lem:repulsive 1} gives that $(\dot{b}_1,\dots, \dot{b}_s)$ is a $\Delta$-repulsive splitting of $\tilde{B}$ 
over $\dot K$.   
Now $\hat b-b\prec\fn^n\prec 1$, so~$(Q_{+b},1,\hat b-b)$ is a refinement of the normal slot
$(Q,1,\hat b)$ in~$H$, hence $(Q_{+b},1,\hat b-b)$ is normal by Proposition~\ref{normalrefine}.
We claim  that the refinement $(Q_{+b},\fn^n,\hat b-b)$ of~$(Q_{+b},1,\hat b-b)$ is also normal.
If $[\fn] \leq \big[\fv(L_{Q_{+b}})\big]$, this claim holds by Corollary~\ref{corcorcor}.
From Lemma~\ref{lem:linear part, new} and~\ref{lem:dotfv} we obtain:
\begin{align*} \order L_{Q_{+b}}\ &=\ \order L_{Q}\ =\ \order L_{Q_{+\hat b}}\ =\ s,\\
 \fv(L_{Q_{+b}})\ \sim\ \fv(L_Q)\ &\sim\ \fv(L_{Q_{+\hat b}}),\qquad v\big(\fv(L_{Q_{+ \hat b}})\big)\ =\ v\big(\fv(A)\big),
 \end{align*} 
so $v\big(\fv(L_{Q_{+b}})\big)=v\big(\fv(A)\big)$.    Moreover,  by Lemma~\ref{lem:dotfv} and the facts about $\tilde{A}$, 
% using $Q_{+b, \times \fn^n}=Q_{\times \fn^n,+b/\times \fn^n}$ and $\fw\sim\fv(A\dot\fn^n)$ {\bf ok}:
$$   v\big(\fv(L_{Q_{+ b,\times\fn^n}})\big)\ =\  v\big(\fv(\tilde{A})  \big)\ =\ v\big(\fv(A\dot\fn^n)\big)\ =\ v(\fw). $$
Suppose  $ \big[\fv(L_{Q_{+b}})\big] < [\fn]$. Then $[\fv(A)] < [\dot{\fn}]$, so $\fv(A\dot\fn^n)\preceq\fv(A)$ using $n>m$. Now the asymptotic relations among the various $\fv(\dots)$ above
give $$\fv(L_{Q_{+b,\times\fn^n}})\ \preceq\ \fv(L_{Q_{+b}}), $$ 
hence $(Q_{+b},\fn^n,\hat b-b)$  is normal by Corollary~\ref{cor:normal for small q, prepZ} applied to $H$ and the normal slot~$(Q_{+b},1,\hat b-b)$ in $H$ in the role of $K$ and $(P,1,\hat a)$, respectively. 
Put~$\fv:=\fv(L_{Q_{+b,\times\fn^n}})$, so $\fv\asymp\fw$. Note that $Q_{+b,\times \fn^n}\in \dot{\mathcal{O}}\{Y\}$, so the image of~$L_{Q_{+b,\times \fn^n}}\in \dot{\mathcal{O}}[\der]$ in 
$\dot{H}[\der]$ is $\tilde{A}$. Thus in $H[\der]$ we have:
$$L_{Q_{+b,\times\fn^n}}\ =\ B +   E \qquad \text{where $E\in \dot{\mathcal O}[\der]$, $\ E\prec_{\Delta(\fv)} \fv^{w+1} L_{Q_{+ b,\times\fn^n}}$.}$$
Now $\dot{b}_1,\dots, \dot{b}_s$ are $\Delta$-repulsive, so $b_1,\dots, b_s$ are $\Delta$-repulsive, hence $$B\ =\ a(\der-b_1)\cdots(\der-b_s)$$ splits $\Delta$-repulsively,  and thus $(\hat b-b)/\fn^n$-repulsively.
Therefore $(Q_{+b},\fn^n,\hat b-b)$ is re\-pul\-sive-normal.
\end{proof}

\noindent
Instead of assuming in the above proposition that $(P,\fm,\hat a)$ is special and $(Q,\fm,\hat b)$ is deep and normal, we can assume, as with  Corollary~\ref{cor:evsplitnormal, 1}, that~$\deg P>1$:

\begin{cor}\label{cor:evrepnormal, 1}
Suppose $\deg P>1$ and $v({\hat b}-H)\subseteq v({\hat c}-H)$. Let~$Q\in Z(H,\hat b)$ have minimal complexity.
Then the $Z$-minimal slot~$(Q,\fm,\hat b)$ in $H$ has a special  refinement~$(Q_{+b},\fn,\hat b-b)$ such that 
$(Q^\phi_{+b},\fn,{\hat b-b})$ is eventually deep and repulsive-normal. 
\end{cor}
\begin{proof}
The beginning of the subsection {\it Achieving split-normality}\/
of Section~\ref{sec:split-normal holes} and~$\deg P >1$ give that $K$ is $r$-linearly newtonian.
Lemmas~\ref{lem:quasilinear refinement} and~\ref{ufm} yield a quasilinear refinement~${(P_{+a},\fn,\hat a-a)}$ of 
our hole $(P,\fm,\hat a)$ in $K$.  
Set $b:=\Re a$. By  Lemma~\ref{lem:same width} we have
$$v\big((\hat a-a)-K\big)\ =\ v({\hat a} - K)\ =\ v\big({\hat b} -H\big)\ =\ {v\big((\hat b-b)-H\big)}.$$
Replacing~$(P,\fm,\hat a)$ and $(Q,\fm,\hat b)$ by
$(P_{+a},\fn,\hat a-a)$ and $(Q_{+b},\fn,\hat b-b)$, respectively, we arrange that~$(P,\fm,\hat a)$ is quasilinear.
Then
by Proposition~\ref{prop:hata special} and~$K$ being $r$-linearly newtonian,  $(P,\fm,\hat a)$ is special;
%Hence if  $(Q,\fm,\hat b)$ is normal, then the rest of the corollary follows from the preceding proposition. To reduce to this case 
hence so is~$(Q,\fm,\hat b)$. Proposition~\ref{varmainthm} gives a refinement~$(Q_{+b},\fn,\hat b-b)$ of $(Q,\fm,\hat b)$ and an active $\phi_0\in H^{>}$ such that~$(Q^{\phi_0}_{+b},\fn,{\hat b-b})$ is deep and normal. 
Refinements of $(P,\fm,\hat a)$ remain quasilinear by Corollary~\ref{cor:ref 2n}.  Since $v(\hat b-H)\subseteq v(\hat c -H)$, Lemma~\ref{lem:same width}(ii) gives 
a refinement~$(P_{+a},\fn,\hat a-a)$ of~$(P,\fm,\hat a)$ with $\Re a=b$.
By Lemma~\ref{speciallemma} the minimal hole~$(P^{\phi_0}_{+a},\fn,\hat a-a)$ in~$K^{\phi_0}$ is special.  Proposition~\ref{evrepnormal} applied to~$(P^{\phi_0}_{+a},\fn,{\hat a-a})$,  $(Q^{\phi_0}_{+b},\fn,\hat b-b)$ in place of~$(P,\fm,\hat a)$, $(Q,\fm,\hat b)$, respectively,   gives us~${b_0\in H}$, $\fn_0\in H^\times$ and a repulsive-normal
refinement $\big(Q^{\phi_0}_{+(b+b_0)},\fn_0, {\hat b - (b+b_0)}\big)$ of~$(Q^{\phi_0}_{+b},\fn,{\hat b-b})$. 
This refinement is steep and hence deep by Corollary~\ref{cor:steep refinement}, since~$(Q^{\phi_0}_{+b},\fn,{\hat b-b})$ is deep.  Thus 
by Corollary~\ref{cor:repulsive-normal comp conj},    $\big(Q_{+(b+b_0)},\fn_0, {\hat b - (b+b_0)}\big)$ is a refinement   
of~$(Q,\fm,\hat b)$ such that that
  $\big(Q^\phi_{+(b+b_0)}, \fn_0,{ \hat b - (b+b_0)}\big)$ is eventually deep and re\-pul\-sive-normal. As a refinement of $(Q, \fm, \hat b)$, it is special. 
\end{proof}

\noindent
In the same way that Corollary~\ref{cor:evsplitnormal, 1}   gave rise to Corollary~\ref{cor:evsplitnormal, 2},  
Corollary~\ref{cor:evrepnormal, 1} gives rise to the following:

\begin{cor} \label{cor:evrepnormal, 2}
If  $\deg P>1$, $v(\hat c-H)\subseteq v(\hat b-H)$, and $R\in Z(H,\hat c)$ has minimal complexity,  then the $Z$-minimal slot~$(R,\fm,\hat c)$ in $H$ has a special refinement~$(R_{+c},\fn,\hat c-c)$   such that $(R^\phi_{+c},\fn,{\hat c-c})$ is eventually deep and repulsive-normal.
\end{cor}

\noindent
By Lemma~\ref{lem:same width} we have  $v(\hat b-H)\subseteq v(\hat c-H)$ or $v(\hat c-H)\subseteq v(\hat b-H)$, hence
the two corollaries above yield Theorem~\ref{thm:repulsive-normal}, completing its proof.  \qed

\subsection*{Strengthening repulsive-normality}
In this subsection we adopt the setting of the subsection {\it Strengthening split-normality}\/ of Section~\ref{sec:split-normal holes}.
Thus~$(P,\fm,\hat a)$ is a slot in $H$ of order $r\geq 1$ and weight $w:=\wt(P)$, and $L:=L_{P_{\times\fm}}$. If $\order L=r$,
we set~$\fv:=\fv(L)$. We let $a$, $b$ range over $H$ and $\fm$, $\fn$ over~$H^\times$.

\begin{definition}\label{def:strongly repulsive-normal}
We say that 
$(P,\fm,\hat a)$ is {\bf almost strongly repulsive-normal} if $\order L=r$,\index{slot!almost strongly repulsive-normal}\index{almost strongly!repulsive-normal}\index{repulsive-normal!almost strongly}  
$\fv\prec^\flat 1$, and there are $Q, R\in H\{Y\}$ such that
\begin{list}{*}{\addtolength\itemindent{-3.5em}\addtolength\leftmargin{0.5em}}
\item[(RN2as)]  $(P_{\times\fm})_{\geq 1}=Q+R$, $Q$ is homogeneous of degree~$1$ and order~$r$,  $L_Q$ has a strong $\hat a/\fm$-repulsive splitting over $K$, and $R\prec_{\Delta(\fv)} \fv^{w+1} (P_{\times\fm})_1$. 
\end{list}
We say that $(P,\fm,\hat a)$ is {\bf strongly repulsive-normal}\index{slot!strongly repulsive-normal}\index{strongly!repulsive-normal} \index{repulsive-normal!strongly} 
 if $\order L =r$, $\fv\prec^\flat 1$, and there are $Q, R\in H\{Y\}$ such that:
\begin{enumerate}
\item[(RN2s)]  
$P_{\times\fm}=Q+R$, $Q$ is homogeneous of degree~$1$ and order~$r$, 
$L_Q$ has a strong $\hat a/\fm$-repulsive splitting over $K$, and $R\prec_{\Delta(\fv)} \fv^{w+1} (P_{\times\fm})_1$.
\end{enumerate}
\end{definition}

\noindent
If $(P,\fm,\hat a)$ is almost strongly repulsive-normal, then $(P,\fm,\hat a)$ is almost strongly split-normal; likewise without ``almost''. 
Thus we can augment our diagram from Sec\-tion~\ref{sec:split-normal holes} as follows, the implications holding for slots of order $\ge 1$
in real closed $H$-fields with small derivation and asymptotic integration: $$\xymatrix@L=6pt{	\parbox{4.5em}{strongly repulsive-normal} \ar@{=>}[r] \ar@{=>}[d]&  \ \parbox{4.5em}{almost strongly repulsive-normal} \ar@{=>}[r] \ar@{=>}[d] &  \ \parbox{4.5em}{repulsive-normal} \ar@{=>}[d] \\ 
\parbox{3.5em}{strongly split-normal} \ar@{=>}[r] \ar@{=>}[d]&  \ \parbox{3.5em}{almost strongly split-normal} \ar@{=>}[r] & \ \parbox{3em}{split-normal} \ar@{=>}[d] \\ 
\parbox{3.5em}{strictly normal} \ar@{=>}[rr] & & \ \parbox{3em}{normal}}  $$
\newline
Adapting the proof of Lemma~\ref{lem:char strong split-norm} gives:

\begin{lemma}\label{lem:char strong rep-norm}
The following are equivalent:
\begin{enumerate}
\item[\textup{(i)}] $(P,\fm,\hat a)$ is strongly repulsive-normal;
\item[\textup{(ii)}] $(P,\fm,\hat a)$ is almost strongly repulsive-normal and strictly normal;
\item[\textup{(iii)}] $(P,\fm,\hat a)$ is almost strongly repulsive-normal and $P(0)\prec_{\Delta(\fv)} \fv^{w+1} (P_1)_{\times\fm}$.
\end{enumerate}
\end{lemma}

\begin{cor}\label{cor:strongly rep splitting => strongly rep-normal} 
If $L$ has a strong $\hat a/\fm$-repulsive splitting  over $K$, then: 
\begin{align*}
\text{$(P,\fm,\hat a)$ is almost strongly repulsive-normal }&\ \Longleftrightarrow\   \text{$(P,\fm,\hat a)$ is  normal,} \\
\text{$(P,\fm,\hat a)$ is strongly repulsive-normal }&\ \Longleftrightarrow\   \text{$(P,\fm,\hat a)$ is strictly normal.}
\end{align*}
\end{cor}

\noindent
If $(P,\fm,\hat a)$ is almost strongly repulsive-normal, then  so are
$(bP,\fm,\hat a)$ for $b\neq 0$ and $(P_{\times\fn},\fm/\fn,\hat a/\fn)$, and likewise with ``strongly''
in place of ``almost strongly''.  
The proof of the next lemma is like that of Lemma~\ref{stronglysplitnormalrefine},
 using Lemmas~\ref{lem:5.18 repulsive-normal} and~\ref{lem:char strong rep-norm} in place of
 Lemmas~\ref{splitnormalrefine} and~\ref{lem:char strong split-norm}, respectively.

\begin{lemma}\label{stronglyrepnormalrefine}
Suppose $(P_{+a},\fm,\hat a-a)$ refines $(P,\fm,\hat a)$.
If $(P,\fm,\hat a)$ is almost strongly repulsive-normal, then so is $(P_{+a},\fm,\hat a-a)$.
If  $(P,\fm,\hat a)$ is   strongly repulsive-normal, $Z$-minimal, and 
$\hat a - a \prec_{\Delta(\fv)} \fv^{r+w+1}\fm$, then $(P_{+a},\fm,\hat a-a)$ is strongly repulsive-normal. 
\end{lemma}

\noindent
Here is the key to achieving almost strong repulsive-normality; its proof is similar to that of Lemma~\ref{stronglysplitnormalrefine, q}:

\begin{lemma}\label{stsprr}
Suppose that $(P,\fm,\hat a)$ is repulsive-normal and $\hat a \prec_{\Delta(\fv)} \fm$. Then for all sufficiently small $q\in\Q^>$,
any $\fn\asymp\fv^q\fm$ yields an almost strongly repulsive-normal refinement $(P,\fn,\hat a)$ of $(P,\fm,\hat a)$.
\end{lemma}
\begin{proof}
First arrange $\fm=1$. Take $Q$, $R$ as in (RN2) for $\fm=1$.  
%and   $q_0\in\Q^>$ with~$\hat a\prec\fv^{q_0}$.
Then Lemma~\ref{lem:achieve strong repulsive splitting}  gives $q_0\in \Q^{>}$ such that $\hat a\prec\fv^{q_0}$ and
for all $q\in\Q$ with~${0<q\leq q_0}$ and~$\fn\asymp\fv^q$,
$L_{Q_{\times\fn}}=L_Q\fn$ has a strong $\hat a/\fn$-repulsive splitting over $K$. 
%Let $(g_1,\dots,g_r)$ be an $\hat a$-repulsive splitting of $L_Q$ over $K$. Then for any $q\in\Q$ with $0<q\leq q_0$ and any $\fn\asymp\fv^q$, $(g_1-\fn^\dagger,\dots,g_r-\fn^\dagger)$ is an $\hat a/\fn$-repulsive splitting of $L_{Q_{\times\fn}}=L_Q\fn$ over $K$, by Corollary~\ref{cor:repulsive hata}. \marginpar{replaced references to  Lemma~\ref{lem:repulsive hata, 1} and   Lemma~\ref{lem:repulsive hata, 2} and its proof} The proof of Lemma~\ref{lem:split strongly multconj} also shows that we can decrease $q_0$ so that for all~$q\in\Q$ with $0<q\leq q_0$ and $\fn\asymp\fv^q$, $(g_1-\fn^\dagger,\dots,g_r-\fn^\dagger)$ is a strong $\hat a/\fn$-repulsive splitting of $L_{Q_{\times\fn}}$ over $K$.
Now  Lem\-ma~\ref{lem:5.19 repulsive-normal}  yields that~$(P,\fn,\hat a)$ is almost strongly repulsive-normal for such~$\fn$.
\end{proof}

\noindent
Using this lemma we now adapt the proof of Corollary~\ref{cor:deep and almost strongly split-normal} to obtain:

\begin{cor}\label{cor:deep and almost strongly repulsive-normal}
Suppose $(P,\fm,\hat a)$ is $Z$-minimal, deep, and repulsive-nor\-mal.
Then $(P,\fm,\hat a)$ has a deep and almost strongly repulsive-nor\-mal refinement.
\end{cor}
\begin{proof} Lemma~\ref{lem:good approx to hata} gives $a$ such that $\hat a - a \prec_{\Delta(\fv)} \fm$. By Corollary~\ref{cor:deep 2, cracks}, the refinement~$(P_{+a},\fm,{\hat a-a})$ of
$(P,\fm,\hat a)$ is deep  with $\fv(L_{P_{+a,\times \fm}})\asymp_{\Delta(\fv)} \fv$, and by Lem\-ma~\ref{lem:5.18 repulsive-normal} it is also
repulsive-normal.
Now apply Lemma~\ref{stsprr} to~${(P_{+a},\fm,\hat a-a)}$ in place of~$(P,\fm,\hat a)$ and again use Corollary~\ref{cor:deep 2, cracks} to preserve being deep.
\end{proof}

\noindent
Next we adapt the proof of Lemma~\ref{lem:strongly split-normal compconj}
to obtain a result about the behavior of (almost) repulsive-normality under compositional conjugation:

\begin{lemma}\label{pqr}
Suppose $\phi$ is active in $H$ with $0<\phi\prec 1$, and there exists $a$ with~$\hat a-a\prec^\flat\fm$. If $(P,\fm,\hat a)$ is almost strongly repulsive-normal, then so is the slot~$(P^\phi,\fm,\hat a)$ in $H^\phi$.  Likewise with ``strongly'' in place of ``almost strongly''.
\end{lemma}
\begin{proof}
We arrange $\fm=1$, 
assume $(P,\fm,\hat a)$ is almost strongly repulsive-normal, and take~$Q$,~$R$ as in~(RN2as).
The proof of Lemma~\ref{lem:split-normal comp conj} shows that with $\fw:=\fv(L_{P^\phi})$
we have $\fw\prec^\flat_\phi 1$ and  $(P^\phi)_{\geq 1} = Q^\phi + R^\phi$ where
$Q^\phi\in H^\phi\{Y\}$ is homogeneous of degree~$1$ and order~$r$, $L_{Q^\phi}$ splits over~$K^\phi$, and
$R^\phi \prec_{\Delta(\fw)} \fw^{w+1} (P^\phi)_1$. By Lemma~\ref{lem:hata-repulsive splitting compconj}, $L_{Q^\phi}=L_Q^\phi$ has even a strong  $\hat a$-repulsive splitting over~$K$.
Hence~$(P^\phi,\fm,\hat a)$ is almost strongly repulsive-normal.
For the rest we use Lem\-ma~\ref{lem:char strong rep-norm} and 
the fact that  if $(P,\fm,\hat a)$ is strictly normal, then so is  $(P^\phi,\fm,\hat a)$.
\end{proof}

%Follow the proof of Lemma~\ref{} with
%Lemmas~\ref{lem:hata-repulsive splitting compconj} and \ref{lem:char strong rep-norm} in place of 
%Lemmas~\ref{lem:split strongly compconj} and \ref{lem:char strong split-norm}, 
%respectively.
%\end{proof}

\noindent
Lemma~\ref{lem:good approx to hata}, the remark preceding Corollary~\ref{cor:repulsive-normal comp conj},  and Lemma~\ref{pqr} yield:

\begin{cor}\label{cor:strongly repulsive-normal compconj}
Suppose $(P,\fm,\hat a)$ is 
$Z$-minimal and deep, and $\phi$ is active in~$H$ with $0<\phi\prec 1$.
If $(P, \fm, \hat a)$ is almost strongly repulsive-normal, then so
is the slot~$(P^\phi,\fm,\hat a)$ in $H^\phi$. Likewise with ``strongly'' in place of ``almost strongly''.
\end{cor}

%\noindent
%We say that {\em $(P^\phi,\fm, \hat a)$ is eventually  strongly repulsive-normal\/} if for some  active~$\phi_0$ in~$H$, $(P^\phi,\fm, \hat a)$ is strongly repulsive-normal for all active $\phi\preceq \phi_0$ in~$H$ with~$\phi>0$; likewise with ``almost strongly'' in place of ``strongly''. 

\noindent
In the case $r=1$, ultimateness yields almost strong repulsive-normality, under suitable assumptions;
more precisely:  

\begin{lemma}\label{lem:rep-norm ultimate}
Suppose $H$ is Liouville closed and of Hardy type, and $\I(K)\subseteq K^\dagger$. 
Assume also that   $(P,\fm,\hat a)$ is normal and special, of order  $r=1$. Then~
$$\text{$(P,\fm,\hat a)$ is ultimate} \quad\Longleftrightarrow\quad\text{$L$ has a strong $\hat a/\fm$-repulsive splitting over $K$,}$$ in which case  $(P,\fm,\hat a)$ is almost strongly repulsive-normal.
\end{lemma}
\begin{proof}
By Lemma~\ref{lem:ultimate normal}, $(P,\fm,\hat a)$ is ultimate iff
$\exc^{\operatorname{u}}(L)\cap v\big((\hat a/\fm)-H\big)\leq 0$, and the latter is
equivalent to $L$ having a strong $\hat a/\fm$-repulsive splitting over $K$, by Corollary~\ref{cor:split rep ultimate}.
For the rest use  Corollary~\ref{cor:strongly rep splitting => strongly rep-normal}.
\end{proof}

\noindent
Liouville closed $H$-fields are $1$-linearly newtonian by Corollary~\ref{cor:Liouville closed => 1-lin newt}, so in view of Lemma~\ref{lem:special dents}  and Corollary~\ref{cor:normal=>quasilinear} we may replace the hypothesis ``$(P,\fm,\hat a)$ is special'' in the previous lemma
by ``$(P,\fm,\hat a)$ is $Z$-minimal or a hole in $H$''. 
This leads to repulsive-normal analogues of Lemma~\ref{5.30real} and Corollary~\ref{cor:5.30real} for $r=1$:

\begin{lemma}  \label{5.30real rep-norm} 
Assume $H$ is Liouville closed   and of Hardy type, and $\I(K)\subseteq K^\dagger$.  
Suppose $(P,\fm,\hat a)$ is $Z$-minimal and quasilinear of order $r=1$. Then there is a refinement $(P_{+a},\fn,\hat a-a)$  of $(P,\fm,\hat a)$ and an active $\phi$ in $H$
with $0<\phi\preceq 1$   such that~$(P^\phi_{+a},\fn,{\hat a-a})$  is  deep,  strictly normal, and ultimate  \textup{(}so $(P^\phi_{+a},\fn,\hat a-a)$  is strongly repulsive-normal by Lem\-mas~\ref{lem:rep-norm ultimate} and~\ref{lem:char strong rep-norm}\textup{)}.
\end{lemma}
\begin{proof}
For any active $\phi$ in $H$ with $0<\phi\preceq 1$ we may replace~$H$,~$(P,\fm,\hat a)$  by~$H^\phi$,~$(P^\phi,\fm,\hat a)$. We may
 also replace~$(P,\fm,\hat a)$ by any of its refinements. 
Since~$H$ is $1$-linearly newtonian,
Corollary~\ref{mainthm, r=1} gives a refinement~$(P_{+a},\fn,\hat a-a)$ of $(P,\fm,\hat a)$ and an active~$\phi$ in $H$ such that~$0<\phi\preceq 1$ and $(P^\phi_{+a},\fn,\hat a-a)$ is   normal.
Replacing~$H$,~$(P,\fm,\hat a)$ by~$H^\phi$,~$(P^\phi_{+a},\fn,{\hat a-a})$,  
we arrange that~$(P,\fm,\hat a)$ itself is normal. 
Then $(P,\fm,\hat a)$ has an ultimate refinement by Proposition~\ref{prop:achieve ultimate},  and 
applying Corollary~\ref{mainthm, r=1} to this refinement and using Lemma~\ref{lem:ultimate refinement}, we  obtain an ultimate refinement  $(P_{+a},\fn,\hat a-a)$ of~$(P,\fm,\hat a)$  and an active $\phi$ in $H$ with $0<\phi\preceq 1$
 such that the $Z$-minimal slot~$(P_{+a}^\phi,\fn,\hat a-a)$ in $H^\phi$ is deep, normal, and ultimate.
 % and hence its linear part splits $(\hat a-a)/\fn$-repulsively over $K^\phi$ by Lem\-ma~\ref{lem:rep-norm ultimate} and the remark following its proof.
 Again replacing~$H$,~$(P,\fm,\hat a)$ by~$H^\phi$,~$(P^\phi_{+a},\fn,{\hat a-a})$,  
we   arrange that~$(P,\fm,\hat a)$  is deep, normal, and ultimate.
Corollary~\ref{cor:achieve strong normality, 1} yields a deep and strictly normal refinement~$(P_{+a},\fm,{\hat a-a})$ of
$(P,\fm,\hat a)$; this refinement is still ultimate by Lemma~\ref{lem:ultimate refinement}.
Hence $(P_{+a},\fm,{\hat a-a})$ is a refinement of $(P,\fm,\hat a)$ as required,  with~$\phi=1$. 
\end{proof}

\noindent
Combining Lemmas~\ref{lem:quasilinear refinement} and~\ref{5.30real rep-norm} with Corollary~\ref{cor:strongly repulsive-normal compconj} yields:

\begin{cor}\label{cor:5.30real rep-norm}
Assume $H$ is Liouville closed, $\upo$-free, and of Hardy type, and $\I(K)\subseteq K^\dagger$.  
Then  every $Z$-minimal slot in $H$ of order $r=1$ has a refinement~$(P,\fm,\hat a)$ such that $(P^\phi,\fm,\hat a)$ is eventually deep, ultimate, and strongly re\-pul\-sive-normal. 
\end{cor}

\noindent
In the next subsection we show how minimal holes of degree $>1$ in $K$ give rise to deep, ultimate,  strongly repulsive-normal,   $Z$-minimal slots in $H$.
% (possibly of order~$>1$) from minimal holes of degree~$>1$ in $K$.

\subsection*{Achieving strong repulsive-normality}
Let $H$ be an $\upo$-free Liouville closed $H$-field with small derivation and constant field $C$, and $(P,\fm,\hat a)$ a minimal hole of order~$r\ge 1$ in~$K:=H[\imag]$. Other conventions are as in the subsection {\it Achieving repulsive-normality.}\/
Our goal is to prove a version of Theorem~\ref{thm:repulsive-normal} with ``repulsive-normal'' improved to 
``strongly repulsive-normal~+~ultimate'': 

{\samepage
\begin{theorem}\label{thm:strongly repulsive-normal} 
Suppose  $C$ is archimedean, $\I(K)\subseteq K^\dagger$,  and $\deg P>1$. Then one of the following conditions is satisfied:
\begin{enumerate}
\item[$\mathrm{(i)}$] $\hat b\notin H$ and some $Z$-minimal slot $(Q,\fm,\hat b)$ in $H$ has a special refinement ${(Q_{+b},\fn,\hat b-b)}$ such that $(Q^\phi_{+b},\fn,\hat b-b)$ is eventually deep, strongly re\-pul\-sive-normal, and ultimate; 
\item[$\mathrm{(ii)}$] $\hat c\notin H$ and some $Z$-minimal slot $(R,\fm,\hat c)$ in $H$ has a special refinement ${(R_{+c},\fn,\hat c-c)}$ such that $(R^\phi_{+c},\fn,\hat c-c)$ is eventually deep, strongly re\-pul\-sive-normal, and ultimate. 
\end{enumerate}
\end{theorem}}

\noindent
The proof of this theorem rests on the following two lemmas, where the standing assumption that $H$ is Liouville closed can be dropped.

\begin{lemma}\label{lem:refine to almost strongly repulsive-normal, Q}
Suppose $\hat b\notin H$ and $(Q,\fm,\hat b)$ is  a $Z$-minimal slot in $H$ with a refinement~${(Q_{+b},\fn,\hat b-b)}$  such that $(Q^\phi_{+b},\fn,\hat b-b)$ is eventually deep and repulsive-normal.
Then $(Q,\fm,\hat b)$ has a refinement~${(Q_{+b},\fn,\hat b-b)}$  such that $(Q^\phi_{+b},\fn,\hat b-b)$ is eventually  deep  and almost strongly repulsive-normal.
\end{lemma}
\begin{proof}
We adapt the proof of Lemma~\ref{lem:refine to almost strongly split-normal, Q}.
% using Corollaries~\ref{cor:deep and almost strongly repulsive-normal} and~\ref{cor:strongly repulsive-normal compconj} 
%in place of Corollary~\ref{cor:deep and almost strongly split-normal} and Lemma~\ref{lem:strongly split-normal compconj}, respectively; for the convenience of the reader, we now give the details.
 Let ${(Q_{+b},\fn,\hat b-b)}$ be a refinement of~$(Q,\fm,\hat b)$  and
let $\phi_0$ be active in $H$ such that $0<\phi_0\preceq 1$
and $(Q^{\phi_0}_{+b},\fn,\hat b-b)$ is   deep and repulsive-normal.
Then  Corollary~\ref{cor:deep and almost strongly repulsive-normal}
yields a refinement $$\big((Q^{\phi_0}_{+b})_{+b_0},\fn_0,(\hat b-b)-b_0\big)$$ of $(Q^{\phi_0}_{+b},\fn,\hat b-b)$ which is deep and
almost strongly repulsive-normal. Hence 
$$\big((Q_{+b})_{+b_0},\fn_0,(\hat b-b)-b_0\big)\ =\ \big( Q_{+(b+b_0)},\fn_0,\hat b - (b+b_0) \big)$$
is  a refinement of $(Q,\fm,\hat b)$, and $\big( Q^\phi_{+(b+b_0)},\fn_0,\hat b - (b+b_0) \big)$ is eventually deep and almost strongly repulsive-normal by Corollary~\ref{cor:strongly repulsive-normal compconj}. 
\end{proof}

\noindent
In the same way we obtain: 

\begin{lemma}\label{lem:refine to almost strongly repulsive-normal, R}
Suppose  $\hat c\notin H$ and $(R,\fm,\hat c)$ is  a $Z$-minimal slot in $H$ with a refinement~${(R_{+c},\fn,\hat c-c)}$  such that $(R^\phi_{+c},\fn,\hat c-c)$ is eventually deep and repulsive-normal.
Then $(R,\fm,\hat c)$  has a refinement~${(R_{+c},\fn,\hat c-c)}$  such that $(R^\phi_{+c},\fn,\hat c-c)$ is eventually deep and almost strongly repulsive-normal.
\end{lemma}

\noindent
Theorem~\ref{thm:repulsive-normal} and the two lemmas above give  Theorem~\ref{thm:repulsive-normal} with ``re\-pul\-sive-normal'' improved to ``almost strongly repulsive-normal''. 
We now upgrade this further to ``strongly repulsive-normal~+~ultimate'' (under an extra assumption).  

\medskip
\noindent
Recall from Lemma~\ref{lem:same width}  that 
$v(\hat b-H)\subseteq v(\hat c-H)$ or $v(\hat c-H)\subseteq v(\hat b-H)$.
Thus the next two lemmas finish the proof of Theorem~\ref{thm:strongly repulsive-normal}.

\begin{lemma}\label{Zdsrnu1}
Suppose   $C$ is archimedean,  $\I(K)\subseteq K^\dagger$, $\deg P>1$, and $$v(\hat b-H)\ \subseteq\ v({\hat c-H}).$$
Let~$Q\in Z(H,\hat b)$ have minimal complexity.
Then the $Z$-minimal slot~$(Q,\fm,\hat b)$ in~$H$ has a special refinement~$(Q_{+b},\fn,{\hat b-b})$ such that $(Q^\phi_{+b},\fn,{\hat b-b})$ is eventually  deep, strongly repulsive-normal, and ultimate.
\end{lemma}
\begin{proof}  Here are two ways of modifying $(Q, \fm, \hat b)$. First, let  $(Q_{+b}, \fn, \hat b -b)$ be a refinement of 
$(Q,\fm, \hat b)$. Lemma~\ref{lem:same width} gives $c\in H$ with $v(\hat a -a)=v(\hat b-b)$ with $a:=b+c\imag$, and so the minimal
hole $(P_{+a}, \fn, \hat a -a)$  in $K$ is a refinement of $(P,\fm,\hat a)$ that relates to $(Q_{+b}, \fn, \hat b -b)$
as $(P, \fm, \hat a)$ relates to $(Q, \fm, \hat b)$. So we can replace~$(P, \fm, \hat a)$ and $(Q, \fm, \hat b)$ by $(P_{+a}, \fn, \hat a -a)$
and $(Q_{+b}, \fn, \hat b -b)$, whenever convenient. Second, let $\phi$ be active in $H$ with $0<\phi\preceq 1$. 
Then we can likewise replace $H$, $K$, $(P, \fm, \hat a)$, $(Q, \fm, \hat b)$ by $H^\phi$, $K^\phi$, $(P^\phi, \fm, \hat a)$, $(Q^\phi, \fm, \hat b)$. 

In this way we first arrange as in the proof of Corollary~\ref{cor:evrepnormal, 1}  that~$(Q,\fm,\hat b)$  is special. Next, we use Proposition~\ref{varmainthm} likewise to arrange  that
$(Q, \fm, \hat b)$ is also normal.  By Propositions~\ref{prop:achieve ultimate} (where the assumption $\I(K)\subseteq K^\dagger$ comes into play) and~\ref{normalrefine}  we arrange that $(Q, \fm, \hat b)$ is
ultimate as well.  The properties ``special'' and ``ultimate''  persist under further refinements and compositional conjugations. 

Now Corollary~\ref{cor:evrepnormal, 1} and Lemma~\ref{lem:refine to almost strongly repulsive-normal, Q} 
give a refinement 
$(Q_{+b},\fn,{\hat b-b})$ of the slot~$(Q,\fm,\hat b)$ in $H$ and an active $\phi_0$ in $H$ with $0<\phi_0\preceq 1$
such that the slot~$(Q^{\phi_0}_{+b},\fn,{\hat b-b})$ in~$H^{\phi_0}$ is  deep and almost strongly repulsive-normal. 
%It is also special, by Lemma~\ref{lem:special refinement}.
Corollary~\ref{cor:achieve strong normality, 1} then yields a deep and strictly normal refinement
$$\big( (Q^{\phi_0}_{+b})_{+b_0},\fn, (\hat b - b)-b_0 \big)$$
of $\big( Q^{\phi_0}_{+b},\fn,{\hat b - b }\big)$.  
This refinement is still almost  strongly re\-pul\-sive-normal by
Lem\-ma~\ref{stronglyrepnormalrefine},   and therefore strongly repulsive-normal by
Lemma~\ref{lem:char strong rep-norm}. 
Co\-rol\-lary~\ref{cor:strongly repulsive-normal compconj} then gives that~$\big( Q_{+(b+b_0)},\fn,{\hat b - (b+b_0) }\big)$
is  a special refinement of our slot~$(Q,\fm,\hat b)$ such that~$\big( Q^\phi_{+(b+b_0)},\fn,{\hat b - (b+b_0) }\big)$ is eventually 
deep and strongly re\-pul\-sive-nor\-mal. 
\end{proof}

\noindent
Likewise: 
%Corollary~\ref{cor:evrepnormal, 2} and Lemma~\ref{lem:refine to almost strongly repulsive-normal, R} in place of Corollary~\ref{cor:evrepnormal, 1} and Lemma~\ref{lem:refine to almost strongly repulsive-normal, Q}, respectively, we prove in a similar way:

\begin{lemma}\label{adersu} 
Suppose $C$ is archimedean,  $\I(K)\subseteq K^\dagger$, $\deg P >1$, and 
$$v(\hat c-H)\ \subseteq\ v(\hat b-H).$$ Let~$R\in Z(H,\hat c)$ have minimal complexity. Then the $Z$-minimal slot~$(R,\fm,\hat c)$ in~$H$ has a special refinement~$(R_{+c},\fn,\hat c-c)$   such that $(R^\phi_{+c},\fn,{\hat c-c})$ is eventually  deep, strongly repulsive-normal, and ultimate. 
\end{lemma} 

\newpage 

\part{Hardy Fields and their Universal Exponential Extensions}\label{part:Hardy fields univ exp ext}

\medskip

\noindent
In this part we turn to Hardy fields.
Section~\ref{sec:germs} contains basic definitions and facts about germs
of one-variable (real- or complex-valued) functions, and in Section~\ref{sec:second-order} we collect
the main facts we need about   linear differential equations.
In Section~\ref{sec:Hardy fields} we introduce Hardy fields and review some extension results  due to
Boshernitzan~\cite{Boshernitzan81, Boshernitzan82, Boshernitzan86} and Rosenlicht~\cite{Ros}. In
Section~\ref{sec:upper lower bds} we discuss upper and lower bounds on the growth of germs in Hardy fields from \cite{Boshernitzan86,Boshernitzan82,Rosenlicht83}, and   Section~\ref{sec:order 2 Hardy fields} contains a first study of second-order linear differential
equations over Hardy fields (to be be completed in Section~\ref{sec:perfect applications}, with our main theorem available).
Section~\ref{sec:upo-free Hardy fields} contains the proof of a significant result about maximal Hardy fields, Theorem~\ref{upo}:
every such Hardy field is $\upo$-free. (See the beginning of that section for a review of this important
property of $H$-asymptotic fields, introduced in~[ADH, 11.7].)
The rest of Section~\ref{sec:upo-free Hardy fields} contains refinements and applications of this fact.
In Section~\ref{sec:bounding} we then prove a general fact about bounding  the derivatives of solutions to linear differential equations,
based on~\cite{Esc, HaLi, Landau}.
In Section~\ref{sec:ueeh} we give an analytic description of the universal exponential extension $\Univ=\Univ_K$, introduced in Part~\ref{part:universal exp ext}, of the algebraic closure~$K$ of a Liouville closed Hardy field extending~$\R$. The elements of~$\Univ$ are exponential sums with coefficients and exponents in $K$.
To extract asymptotic information about the summands in such a sum we use  results of Boshernitzan~\cite{BoshernitzanUniform} about uniform distribution mod~$1$ over Hardy fields. We include   proofs of these results
in Section~\ref{sec:udmod1}, preceded by a development of the required classical facts concerning almost periodic functions in Section~\ref{sec:almost periodic}. (None of the material in Sections~\ref{sec:almost periodic} and~\ref{sec:udmod1} is original, we only aim for an efficient and self-contained exposition.)  

\section{Germs of Continuous Functions}\label{sec:germs}

\noindent
Hardy fields consist of germs of one-variable differentiable real-valued functions. In this section we first consider
the ring $\c$ of germs of {\it continuous}\/ real-valued functions, and its complex counterpart $\c[\imag]$. With an eye towards
applications to Hardy fields, we pay particular attention to extending subfields of $\c$.
 
\subsection*{Germs} As in [ADH, 9.1]  we let $\mathcal{G}$ be the ring
of germs at $+\infty$ of real-valued functions whose domain is
a subset of $\R$
containing an interval~$(a, +\infty)$, $a\in \R$; the domain may vary 
and the ring operations are defined as usual.\index{germ} 
If $g\in \mathcal{G}$ is the germ of a real-valued function 
on a subset of $\R$ containing an interval $(a, +\infty)$, $a\in \R$, then we simplify notation
by letting $g$ also denote this function if the resulting ambiguity is harmless.
With this convention, given a property~$P$ of real numbers
and $g\in \mathcal{G}$ we say that {\em $P\big(g(t)\big)$ holds eventually\/} if~$P\big(g(t)\big)$ holds for all sufficiently large real $t$. 
Thus for $g\in \cal G$ we have $g=0$ iff $g(t)=0$ eventually (and so $g\neq 0$ iff $g(t)\neq 0$ for arbitrarily large $t$). 
Note that the multiplicative group $\mathcal{G}^\times$ of units of $\mathcal{G}$ consists of the $f\in \mathcal{G}$ such that $f(t)\ne 0$, eventually.
We identify each real number~$r$ with the germ at~$+\infty$ of the function~$\R\to \R$ that takes the constant value~$r$. This  makes the field~$\R$ into a subring of $\mathcal{G}$. 
Given $g,h\in\cal G$, we set
\begin{equation}\label{eq:germs partial ordering}
g \leq h \quad :\Longleftrightarrow\quad \text{$g(t)\leq h(t)$, eventually.}
\end{equation}
This defines a partial ordering $\leq$ on $\mathcal G$ which restricts to the usual ordering of $\R$.

Let $g,h\in\mathcal{G}$. Then $g,h\geq 0\Rightarrow g+h,\,g\cdot h,\,g^2\geq 0$, and $g\geq r\in \R^{>}\Rightarrow g\in\mathcal{G}^\times$.
%We say that~$g,h\in\cal G$ are {\it comparable} if~$g\leq h$ or $h\leq g$.
We define $g<h:\Leftrightarrow g\leq h$ and~$g\neq h$. Thus
if $g(t)<h(t)$, eventually, then $g < h$; the converse is not generally valid.

\subsection*{Continuous germs}
We call a germ $g\in \mathcal{G}$ {\em continuous}\/ if it is the germ of
a continuous function $(a,+\infty)\to \R$ for some $a\in \R$, and we let $\mathcal{C}\supseteq \R$ be the subring of $\mathcal{G}$ consisting of the continuous germs $g\in \mathcal{G}$.\index{germ!continuous}\label{p:cont} 
We have $\c^\times=\mathcal{G}^\times\cap\c$;
thus for~$f\in \c^\times$, we have $f(t)\neq 0$, eventually,   hence either $f(t)>0$, eventually, or~$f(t)<0$, eventually, and so $f>0$ or $f<0$.
More generally, if $g,h\in\cal C$ and~$g(t)\neq h(t)$, eventually, then 
 $g(t)<h(t)$, eventually, or  $h(t)<g(t)$, eventually.
 %; in particular,~$g<h$ or $h>g$.
 We let $x$ denote the germ at $+\infty$ of the identity function on $\R$, so $x\in\c^\times$.
 
\subsection*{The ring $\c[\imag]$} In analogy with $\c$ we define its complexification $\c[\imag]$ as the ring of
germs at $+\infty$ of $\C$-valued continuous functions whose domain is a subset of $\R$ containing an interval $(a,+\infty)$, $a\in \R$. It has 
$\c$ as a subring. Identifying each complex number $c$ with the germ at $+\infty$ of the function $\R\to \C$ that takes the constant value $c$ makes $\C$ also a subring of $\c[\imag]$ with $\c[\imag]=\c+\c\imag$, justifying the notation $\c[\imag]$. 
The ``eventual'' terminology for germs $f\in \c$ (like ``$f(t)\ne 0$, eventually'') is extended in the obvious way to germs $f\in \c[\imag]$. Thus for $f\in \c[\imag]$ we have:
$f(t)\ne 0$, eventually, if and only if $f\in \c[\imag]^\times$. 
In particular $\c^\times=\c[\imag]^\times\cap \c$.

\medskip
\noindent
Let $\Phi\colon U\to\C$ be a continuous function where $U\subseteq\C$, and let $f\in\c[\imag]$ be such that~$f(t)\in U$, eventually; then $\Phi(f)$ denotes the germ in $\c[\imag]$ with
$\Phi(f)(t)=\Phi\big(f(t)\big)$, eventually.
For example, taking $U=\C$, $\Phi(z)=\ex^z$, we obtain for $f\in\c[\imag]$ the germ  $\exp f = \ex^f\in\c[\imag]$
with $(\ex^f)(t)=\ex^{f(t)}$, eventually. 
Likewise, for~$f\in \c$ with~$f(t)>0$, eventually, we have the germ $\log f\in\c$. 
For $f\in\c[\imag]$  we have~$\overline{f}\in\c[\imag]$ with  $\overline{f}(t)=\overline{f(t)}$, eventually; the map $f\mapsto\overline{f}$ is an automorphism of the ring $\c[\imag]$
with $\overline{\overline{f}}=f$ and~$f\in\c\Leftrightarrow \overline{f}=f$. 
For $f\in\c[\imag]$ we also have~$\Re f,\Im f,\abs{f}\in\c$ with~$f(t)=(\Re f)(t)+(\Im f)(t)\imag$ and $\abs{f}(t)=\abs{f(t)}$, eventually.

\subsection*{Asymptotic relations on $\c[\imag]$.} 
Although $\c[\imag]$ is not a valued field, it will be convenient to equip $\c[\imag]$ with the asymptotic relations $\preceq$,~$\prec$,~$\sim$ (which are defined on any valued field [ADH, 3.1]) as follows: for $f,g\in \c[\imag]$,
\begin{align*} f\preceq g\quad &:\Longleftrightarrow\quad \text{there exists $c\in \R^{>}$ such that $|f|\le c|g|$,}\\
f\prec g\quad &:\Longleftrightarrow\quad \text{$g\in \c[\imag]^\times$ and $\lim_{t\to \infty} f(t)/g(t)=0$} \\
 &\phantom{:} \Longleftrightarrow\quad \text{$g\in \c[\imag]^\times$ and $\abs{f}\leq c\abs{g}$ for all $c\in\R^>$},\\
f\sim g\quad &:\Longleftrightarrow\quad \text{$g\in \c[\imag]^\times$ and
$\lim_{t\to \infty} f(t)/g(t)=1$}\\ 
\quad&\phantom{:} \Longleftrightarrow\quad f-g\prec g.
\end{align*}
We also use these notations for continuous functions~${[a,+\infty)\to\C}$,  $a\in\R$; for example, for continuous $f\colon [a,+\infty)\to\C$ and
$g\colon [b,+\infty)\to\C$~($a,b\in \R$), $f\preceq g$ means: $\text{(germ of $f$)}\preceq \text{(germ of $g$)}$. 
If $h\in\c[\imag]$ and $1\preceq h$, then~${h\in\c[\imag]^\times}$. Also, for~$f,g\in\c[\imag]$ and $h\in\c[\imag]^\times$ we have
$$f\preceq g\ \Leftrightarrow\ fh\preceq gh, \qquad f\prec g\ \Leftrightarrow\ fh\prec gh, \qquad f \sim g\ \Leftrightarrow\ fh \sim gh.$$
The binary relation $\preceq$ on
$\c[\imag]$ is reflexive and transitive, and $\sim$ is an equivalence relation on $\c[\imag]^\times$.
Moreover, for $f,g,h\in \c[\imag]$ we have
$$f\prec g\ \Rightarrow\ f\preceq g, \qquad f\preceq g \prec h\ \Rightarrow\ f\prec h, \qquad f\prec g \preceq h\ \Rightarrow\ f\prec h.$$ 
Note that $\prec$ is a transitive binary relation
on $\c[\imag]$.  
For $f,g\in \c[\imag]$ we also set
$$f\asymp g:\ \Leftrightarrow\ f\preceq g \ \&\ g\preceq f,\qquad f\succeq g:\ \Leftrightarrow\ g\preceq f,\qquad f\succ g:\ \Leftrightarrow\ g\prec f,
$$
so $\asymp$ is an equivalence relation on $\c[\imag]$, and $f\sim g\Rightarrow f\asymp g$.  
Thus for $f,g,h\in\c[\imag]$, 
$$f\preceq g\ \Rightarrow\  fh\preceq gh,\quad f\preceq h\ \&\  g\preceq h\ \Rightarrow\ f+g \preceq h,
\quad f\preceq 1\ \&\ g\prec 1 \ \Rightarrow\ fg\prec 1,$$
hence\label{p:contpreceq} 
$$\c[\imag]^{\preceq}\ :=\ \big\{f\in\c[\imag]: f\preceq 1\big\}\ = \ \big\{ f\in\c[\imag]: \text{$|f|\leq n$ for some $n$}\big\}$$
is a subalgebra of the $\C$-algebra $\c[\imag]$ and
$$\c[\imag]^{\prec}\	:=\ \big\{f\in\c[\imag]: f\prec 1\big\}\ 
					 =\ \left\{f\in\c[\imag]: \lim_{t\to\infty} f(t)=0 \right\}$$
is an ideal of $\c[\imag]^{\preceq}$. The group of units of $\c[\imag]^{\preceq}$ is
$$\c[\imag]^{\asymp} \ := \ \big\{ f\in\c[\imag]: f\asymp 1\big\} \ = \ \big\{ f\in\c[\imag]: \text{$1/n\leq |f|\leq n$ for some $n\geq 1$}\big\} $$
and has the subgroup 
$$\C^\times\big(1+\c[\imag]^{\prec}\big)\ =\ \left\{f\in\c[\imag]: \lim_{t\to\infty} f(t)\in\C^\times \right\}.$$ 
We set 
$\c^{\preceq}:=\c[\imag]^{\preceq}\cap\c$, and similarly with~$\prec$,~$\asymp$ in place of $\preceq$.

\begin{lemma}\label{lem:sim props}  
Let $f,g,f^*,g^*\in\c[\imag]^\times$ with $f\sim f^*$ and $g\sim g^*$. Then
$1/f\sim 1/f^*$ and~$fg\sim f^*g^*$. Moreover, $f\preceq g\Leftrightarrow f^*\preceq g^*$, and similarly with
$\prec$, $\asymp$, or $\sim$ in place of~$\preceq$.
\end{lemma}

\noindent
This follows easily from the observations above. For later reference we also note:

\begin{lemma}\label{lem:log preceq} 
Let $f,g\in\c^\times$ be such that $1\prec f\preceq g$; then $\log |f| \preceq \log|g|$.
\end{lemma}
\begin{proof}
Clearly $\log|g| \succ 1$.
Take $c\in\R^>$ such that $|f| \leq c|g|$. Then $\log|f|\leq \log c+\log|g|$ where $\log c+\log|g|\sim\log|g|$;
hence $\log |f| \preceq \log|g|$.
\end{proof}

\begin{lemma}\label{lem:fsimg criterion}  
Let $f,g,h\in\c^\times$ be such that  $f-g\prec h$ and
$(f-h)(g-h)=0$. Then~$f\sim g$.
\end{lemma}
\begin{proof}
Take $a\in\R$ and representatives $(a,+\infty)\to\R$ of $f$, $g$, $h$, denoted by the same symbols, such that
for each $t > a$ we have $f(t),g(t),h(t)\neq 0$, and $f(t)=h(t)$ or $g(t)=h(t)$. Let $\varepsilon\in\R$ with $0<\varepsilon \leq 1$ be given, and choose $b\geq a$ such that~$|f(t)-g(t)| \leq \frac{1}{2}\varepsilon |h(t)|$ for all $t > b$.
Set $q:=f/g$ and
let $t > b$; we claim that then $|q(t)-1|\leq \varepsilon$.
This is clear if $g(t)=h(t)$, so suppose otherwise; then~$f(t)=h(t)$, and  
$|1-1/q(t)| \leq \frac{1}{2}\varepsilon\leq \frac{1}{2}$. In particular, $0<q(t)\leq 2$ and so~$|1-q(t)| = |1-1/q(t)|\cdot q(t) \leq \varepsilon$ as claimed. 
\end{proof}

\subsection*{Subfields of $\mathcal{C}$}
Let $H$ be a {\em Hausdorff field}, that is, a subring
of $\mathcal{C}$ that happens to be a field; see \cite{ADH2}.\index{Hausdorff field} 
Then $H$ has the subfield~$H\cap \R$. 
If $f\in H^\times$, then~$f(t)\neq 0$ eventually, hence
either $f(t)<0$ eventually or $f(t)>0$ eventually. The partial ordering of $\cal G$ from \eqref{eq:germs partial ordering} thus restricts to a total
ordering on $H$ making~$H$ an ordered field in the usual sense of that term. 
By \cite[Propositions~3.4 and~3.6]{Boshernitzan81}: 

\begin{prop}\label{b1} 
Let $H^{\operatorname{rc}}$ consist of the
$y\in \mathcal{C}$ with $P(y)=0$ for some $P(Y)$ in $H[Y]^{\ne}$. Then $H^{\operatorname{rc}}$ is the
unique real closed  Hausdorff field that
extends $H$ and is algebraic over~$H$. In particular, $H^{\operatorname{rc}}$ is a real closure of the ordered field $H$.  
\end{prop} 

\noindent
Boshernitzan~\cite{Boshernitzan81} assumes $H\supseteq \R$ for this result, but this is not really needed in the proof, much of
which already occurs in Hausdorff~\cite{Hau09}.  
For the reader's convenience we include a proof of Proposition~\ref{b1}, after some lemmas. Let
$$P(Y)\ =\ P_0Y^n+P_{1}Y^{n-1}+\cdots+P_n\in H[Y] \qquad (P_0,\dots,P_{n}\in H),$$
and take $a\in\R$ such that~$P_0,\dots,P_n$ have representatives in $\c_a$, also
denoted by~$P_0,\dots,P_n$. This yields for $t\geq a$  the polynomial
$$P(t,Y)\  :=\  P_0(t)Y^n+ P_1(t)Y^{n-1}+\cdots+P_n(t)\in \R[Y].$$
For any other choice of $a$ and representatives of $P_0,\dots, P_n$ in 
$\c_a$ this gives for large enough $t$ the same polynomial $P(t,Y)\in \R[Y]$, so the ``eventual'' terminology
makes sense for properties mentioning $P(t,Y)$ with $t$ ranging over $\R$. For example, for~$y\in \c[\imag]$, we have:
$P(y)=0\Leftrightarrow y(t)\in \C \text{ is a zero of $P(t,Y)$, eventually}$. 

\begin{lemma}\label{lem:parametrize zeros}
Suppose $P$ is  irreducible \textup{(}in $H[Y]$\textup{)} of degree $n$, so $n\ge 1$. 
Then there are $y_1,\dots, y_m\in \c$ such that $y_1(t) < \cdots < y_m(t)$, eventually, and the
distinct real zeros
of the polynomial $P(t,Y) \in\R[Y]$ are exactly $y_1(t), \dots, y_m(t)$, eventually.  
Thus $P(y_1)=\cdots=P(y_m)=0$, and if $n$ is odd, then $m\ge 1$.  
\end{lemma}
\begin{proof}
Take $A,B\in H[Y]$ with $1=AP+BP'$. Then $$1\ =\ A(t,Y)P(t,Y)+B(t,Y)P(t,Y)',\quad \text{ eventually}.$$ 
Hence $P(t,Y)\in \R[Y]$ has exactly $n$ distinct complex zeros, eventually.
Now use ``continuity of roots'' as used for example in~\cite[Chapter~II, (2.4)]{Domin}. 
\end{proof}

\begin{lemma}\label{lemPQyz}
Let $P, Q\in H[Y]$ be monic and irreducible with $P\neq Q$, and let~$y,z\in\mathcal C[\imag]$, $P(y)=Q(z)=0$.
Then $y(t)\neq z(t)$, eventually. In particular, if~$y,z\in \mathcal C$, then either $y(t) < z(t)$ eventually, or $y(t)>z(t)$ eventually. 
\end{lemma}
\begin{proof}
Take $A,B\in H[Y]$ such that $1=AP+BQ$. Then 
$$1\ =\ A(t,Y)P(t,Y)+B(t,Y)Q(t,Y),\quad \text{ eventually}.$$ 
Hence $Q(t, y(t))\ne 0$, eventually, and thus $y(t)\ne z(t)$, eventually. 
\end{proof}

\begin{cor}\label{corPQ}
Let $y\in\c$ and $P(y)=0$, $P\in H[Y]^{\neq}$.
Then~${Q(y)=0}$ for some monic irreducible $Q\in H[Y]$. 
\end{cor}
\begin{proof} We have $P=hQ_1^{e_1}\cdots Q_n^{e_n}$ where $h\in H^{\ne}$,~$e_1,\dots,e_n\in\N^{\geq 1}$, and
$Q_1,\dots,Q_n$ in $H[Y]$ are distinct,  and monic irreducible. 
Lemmas~\ref{lem:parametrize zeros} and~\ref{lemPQyz} now yield germs~$y_1,\dots,y_m\in\c$ such that, eventually,  
$y_1(t)<\cdots<y_m(t)$ are
the real zeros of the polynomials $Q_1(t,Y),\dots,Q_n(t,Y)\in \R[Y]$, and thus of $P(t,Y)$, 
and such that for each $i\in \{1,\dots,m\}$ there is a unique $j\in \{1,\dots,n\}$ with $Q_j\big(t, y_i(t)\big)=0$, eventually. 
Continuity arguments and the connectedness of halflines $[a,+\infty)$ 
yields a single~$i$ with~$y_i(t)=y(t)$, eventually,  and thus $Q_j(y)=0$ for some~$j$.
\end{proof}

\begin{proof}[Proof of Proposition~\ref{b1}] Given $y\in H^{\operatorname{rc}}$ we have by Corollary~\ref{corPQ} a monic
irreducible $Q\in H[Y]$ with $Q(y)=0$, so $H[y]\subseteq H^{\operatorname{rc}}$ is a Hausdorff field and algebraic over $H$. Since ``algebraic over'' is  transitive, it follows that $H^{\operatorname{rc}}$ is a Hausdorff field and algebraic over $H$. Such  transitivity also gives
$(H^{\operatorname{rc}})^{\operatorname{rc}}=H^{\operatorname{rc}}$.  
Obviously, any algebraic Hausdorff field extension of $H$ is contained in $H^{\operatorname{rc}}$. So it only
remains to show that the ordered field $H^{\operatorname{rc}}$ is real closed. 
First, if $y\in H^{\operatorname{rc}}$ and $y\ge 0$, then clearly $\sqrt{y}\in \c$ is algebraic over $H^{\operatorname{rc}}$, and
thus in it. Next, let~$P(Y)\in H^{\operatorname{rc}}[Y]$ have odd degree. Then
$P$ has a zero in $H^{\operatorname{rc}}$: this follows from Lemma~\ref{lem:parametrize zeros} by considering an irreducible
factor of $P$ in $H^{\operatorname{rc}}[Y]$  of odd degree. 
\end{proof}

\noindent
We record the following useful consequence of  Corollary~\ref{corPQ} and its proof:

\begin{cor}\label{corb1} Let $P\in H[Y]^{\ne}$ and let $y_1,\dots, y_m$ be the distinct zeros of $P$ in~$H^{\operatorname{rc}}$.
Then $y_1(t),\dots y_m(t)$ are the distinct real zeros of $P(t,Y)$, eventually.
\end{cor}

\noindent
Note that $H[\imag]$ is a subfield of $\c[\imag]$, and by Proposition~\ref{b1} and [ADH, 3.5.4], the subfield~$H^{\operatorname{rc}}[\imag]$ of $\c[\imag]$ is an algebraic closure of the field $H$.
If $f\in\c[\imag]$ is integral over $H$, then so is $\overline{f}$, and hence so are the elements $\Re f=\frac{1}{2}(f+\overline{f})$ and
$\Im f=\frac{1}{2\imag}(f-\overline{f})$ of $\c$ [ADH, 1.3.2]. Thus
$H^{\operatorname{rc}}[\imag]$  consists of the
$y\in \c[\imag]$ with~$P(y)=0$ for some $P(Y)\in H[Y]^{\ne}$. 

\medskip
\noindent
The ordered field $H$ has a convex subring 
$$ \mathcal{O}\ =\ 
\big\{f\in H :\ \text{$|f|\le n$ for some $n$}\big\}\ =\ \c^{\preceq}\cap H,$$
which is a valuation ring of $H$, and we
consider $H$ accordingly as a valued ordered field. The maximal ideal of $\mathcal{O}$ is
$\smallo = \c^{\prec}\cap H$. The residue morphism $\mathcal O\to\res(H)$ restricts to an ordered field
embedding $H\cap\R\to \res(H)$, which is bijective if $\R\subseteq H$.
Restricting the binary relations
$\preceq$,~$\prec$,~$\sim$ from the previous subsection to $H$ gives exactly the asymptotic relations $\preceq$,~$\prec$,~$\sim$ on $H$
that it comes equipped with as a valued field. By [ADH, 3.5.15], 
$$\mathcal{O}+\mathcal{O}\imag\ =\ \big\{f\in H[\imag]: \text{$|f|\le n$ for some $n$}\big\}\ =\ \c[\imag]^{\preceq}\cap H[\imag]$$ is the unique valuation ring of $H[\imag]$ whose intersection with $H$ is $\mathcal{O}$. In this way we consider
$H[\imag]$ as a valued field extension of $H$. The maximal ideal of $\mathcal{O}+\mathcal{O}\imag$ is~$\smallo+\smallo \imag=\c[\imag]^{\prec}\cap H[\imag]$.
The asymptotic relations $\preceq$, $\prec$, $\sim$ on $\c[\imag]$
restricted to~$H[\imag]$ are exactly the asymptotic relations $\preceq$, $\prec$, $\sim$ on $H[\imag]$ that $H[\imag]$ has as a valued field.
Moreover, $f \asymp \abs{f}$ in $\c[\imag]$ for all $f \in H[\imag]$.

\subsection*{Composition} 
Let $g\in \c$, and suppose that $\lim\limits_{t\to +\infty} g(t)=+\infty$; equivalently,  $g\geq 0$ and $g\succ 1$. Then the 
composition operation
$$f\mapsto f\circ g\ :\ \c[\imag] \to \c[\imag], \qquad (f\circ g)(t)\ :=\  f\big(g(t)\big)\ \text{ eventually},$$
is an injective endomorphism of the ring $\c[\imag]$
that is the identity on the subring~$\C$. For $f_1, f_2\in \c[\imag]$ we have:  $f_1\preceq f_2 \Leftrightarrow f_1\circ g \preceq f_2\circ g$, and likewise with $\prec$, $\sim$. 
This endomorphism of $\c[\imag]$ commutes with the automorphism $f\mapsto\overline{f}$ of $\c[\imag]$, and maps each
subfield $K$ of $\c[\imag]$ isomorphically 
onto the subfield $K\circ g=\{f\circ g:f\in K\}$ of~$\c[\imag]$. Note that if the subfield~$K$ of $\c[\imag]$ contains $x$, then $K\circ g$ contains $g$. 
Moreover,   $f\mapsto f\circ g$ restricts to an endomorphism of the subring $\c$ of~$\c[\imag]$ such that if  $f_1, f_2\in \c$ and $f_1\leq f_2$, then $f_1\circ g \leq f_2\circ g$. This endomorphism of $\c$
maps each Hausdorff field $H$ isomorphically (as an ordered field)
onto the Hausdorff field~$H\circ g$.

\medskip\noindent
Occasionally it is convenient to extend the composition operation on~$\c$ to the ring~$\mathcal G$ of all (not necessarily continuous) germs. 
Let~$g\in\mathcal G$ with~$\lim\limits_{t\to+\infty} g(t)=+\infty$. Then for $f\in\mathcal G$ we have the
germ $f\circ g\in\mathcal G$ with  
$$(f\circ g)(t)\ :=\  f\big(g(t)\big)\ \text{ eventually.}$$ 
The map $f\mapsto f\circ g$ is an endomorphism of the $\R$-algebra $\mathcal G$. 
Let $f_1, f_2\in \mathcal{G}$. Then~$f_1\leq f_2\Rightarrow f_1\circ g\leq f_2\circ g$, and likewise with $\preceq$ and $\prec$ instead of $\leq$, where
we extend the binary relations  $\preceq$, $\prec$
from $\c$ to $\mathcal G$ in the natural way: 
\begin{align*} f_1\preceq f_2\quad &:\Longleftrightarrow\quad \text{there exists $c\in \R^{>}$ such that $|f_1(t)|\le c|f_2(t)|$, eventually;}\\
f_1\prec f_2\quad &:\Longleftrightarrow\quad \text{$f_2\in {\mathcal G}^\times$ and $\lim_{t\to \infty} f_1(t)/f_2(t)=0$.}  
\end{align*}

\subsection*{Compositional inversion}
Suppose that $g\in \c$ is eventually strictly increasing such that~$\lim\limits_{t\to +\infty} g(t)=+\infty$.\label{p:ginv}
Then its compositional inverse $g^{\inv}\in \c$ is given by~$g^{\inv}\big(g(t)\big)=t$, eventually,
and  $g^{\inv}$ is also  eventually strictly increasing with~$\lim\limits_{t\to +\infty} g^{\inv}(t)=+\infty$.
Then $f\mapsto f\circ g$  is an automorphism of the ring  $\c[\imag]$, with inverse~$f\mapsto f\circ g^{\inv}$. 
In particular, $g\circ g^{\inv}=g^{\inv}\circ g=x$. Moreover, 
$f\mapsto f\circ g$ restricts to an automorphism of $\c$, and
if  $h\in\c$ is eventually strictly increasing with  $g\leq h$, then~$h^{\operatorname{inv}}\leq g^{\operatorname{inv}}$.

Let now  $f,g\in\c$ with $f,g\geq 0$, $f,g\succ 1$.  
It is not true in general that if~$f$,~$g$ are eventually strictly increasing and $f\sim g$, then $ f^{\operatorname{inv}}\sim g^{\operatorname{inv}}$. (Counterexample:  $f=\log x$, $g=\log 2x$.) Corollary~\ref{cor:Entr} below gives a useful condition on $f$, $g$ under which this implication does hold. 
In addition, let $h\in\c^\times$ be eventually monotone and continuously differentiable with~$h'/h\preceq 1/x$.

\begin{lemma}[Entringer~\cite{Entringer}]\label{lem:Entr}
Suppose ${f\sim g}$.  Then~$h\circ f\sim h\circ g$.
\end{lemma}
\begin{proof}
Replacing $h$ by $-h$ if necessary we arrange that $h\geq 0$, so $h(t)>0$ eventually. 
Set $p:=\min(f,g)\in\c$ and $q:=\max(f,g)\in\c$. Then $0\le p\succ 1$ and $f-g\prec p$.
The Mean Value Theorem gives $\xi\in\mathcal G$ such that 
$p\leq \xi \leq q$ (so $0\le \xi\succ 1$) and
$$h\circ f - h \circ g\  =\  (h'\circ \xi)\cdot (f-g).$$
From  $h'/h\preceq 1/x$ we obtain $h'\circ \xi \preceq ({h\circ \xi})/\xi\preceq (h\circ \xi)/p$, hence
$h\circ f - h \circ g  \prec h\circ \xi$. Set $u:=\max(h\circ p,h\circ q)$. Then 
$0\leq h\circ\xi \leq u$, hence~$h\circ f - h \circ g  \prec u$.
Also~$(u-h\circ f)(u-h\circ g)=0$, so Lemma~\ref{lem:fsimg criterion} yields~$h\circ f  \sim h\circ g$.
\end{proof}

\begin{cor}\label{cor:Entr}
Suppose $f$, $g$ are eventually strictly
increasing with  $f\sim g$  and~$f^{\operatorname{inv}}\sim h$. Then~$g^{\operatorname{inv}}\sim h$.
\end{cor}
\begin{proof}
By the lemma above we have $h\circ f\sim h\circ g$, and from $f^{\operatorname{inv}}\sim h$
we obtain~$x=f^{\operatorname{inv}}\circ f\sim h\circ f$. Therefore
$g^{\operatorname{inv}}\circ g= x\sim h\circ g$ and thus $g^{\operatorname{inv}}\sim h$.
\end{proof}

\begin{cor}\label{cor:Entr, 2}
If $g$, $h$ are eventually strictly increasing, $0\le h\succ 1$, and
 $g\sim h^{\operatorname{inv}}$, then
$g^{\operatorname{inv}}\sim h$.
\end{cor}
\begin{proof}
Take $f= h^{\operatorname{inv}}$ in Corollary~\ref{cor:Entr}.
\end{proof}

\noindent
Sometimes we prefer ``big O'' and ``little o'' notation instead of $\preceq$ and $\prec$: for~$\phi,\xi,\theta\in \c[\imag]$,
$\phi=\xi+O(\theta):\Leftrightarrow \phi-\xi\preceq \theta$ and $\phi=\xi+o(\theta):\Leftrightarrow \phi-\xi\prec \theta$. For use in Section~\ref{sec:perfect applications} we note:
 
\begin{cor}\label{cor:Entr, 3} 
Suppose  $g = x+cx^{-1}+o(x^{-1})$, $c\in\R$, and $g$ is eventually strictly increasing.
Then
$g^{\operatorname{inv}} = x-cx^{-1}+o(x^{-1})$.
\end{cor}
\begin{proof}
We have  $g^{\operatorname{inv}}\sim x$ by Corollary~\ref{cor:Entr, 2} (for $h=x$),
so $g^{\operatorname{inv}}=x(1+\varepsilon)$ where~${\varepsilon\in \c}$, $\varepsilon \prec 1$. 
Now $(1+\varepsilon)^{-1}=1+\delta$ with $\delta\in\c$, $\delta\prec 1$.
Then
$$x\ =\ g\circ g^{\operatorname{inv}}\ =\ x(1+\varepsilon)+cx^{-1}(1+\delta)+o(x^{-1})$$
and thus
$\varepsilon = -cx^{-2}(1+\delta)+o(x^{-2})=-cx^{-2}+o(x^{-2})$. This yields
  $g^{\operatorname{inv}}=x(1+\varepsilon)=x-cx^{-1}+o(x^{-1})$, as claimed.
\end{proof}

\subsection*{Extending ordered fields inside an ambient partially ordered ring}
Let $R$ be a commutative ring with $1\ne 0$, equipped with a translation-invariant
partial ordering $\le$ such that $r^2\ge 0$ for all $r\in R$, and
$rs\ge 0$
for all $r,s\in R$ with $r,s\ge 0$. It follows that for $a,b,r\in R$ we have: 
\begin{enumerate}
\item if $a \le b$ and $r\ge 0$, then $ar\le br$; 
\item if $a$ is a unit and $a>0$, then $a^{-1}=a\cdot (a^{-1})^2>0$; 
\item if $a$,~$b$ are units and $0 < a \le b$, then $0 < b^{-1}\le a^{-1}$. 
\end{enumerate}
Relevant cases: 
$R=\mathcal{G}$ and $R=\mathcal{C}$, with partial ordering given by \eqref{eq:germs partial ordering}. 

\medskip
\noindent
%Call a subset $S$ of $R$ {\em totally ordered\/} if the partial ordering of $R$ induces a total ordering on $S$.
An {\em ordered subring of $R$} is a subring  
of $R$ that is totally ordered by the partial ordering of $R$. An {\em ordered subfield of $R$} is an ordered subring $H$ of $R$
which happens to be a field; then $H$ equipped with the induced ordering
is indeed an ordered field, in the usual sense of that term.
(Thus any Hausdorff field is an ordered subfield of the partially ordered ring $\mathcal{C}$.) 
We identify $\Z$ with its image in $R$ via the unique ring embedding $\Z \to R$, and this makes $\Z$ with its usual ordering into an ordered subring of $R$. 

\begin{lemma}\label{pr0} Assume $D$ is an ordered
subring of $R$ and every nonzero element of $D$ is a unit of $R$.
Then $D$ generates an ordered subfield $\Frac{D}$ of $R$.
\end{lemma}
\begin{proof} It is clear that $D$ generates a subfield $\Frac{D}$ of $R$. For $a\in D$, $a>0$, we have $a^{-1}>0$. It follows that $\Frac{D}$ is totally ordered.
\end{proof}

\noindent
Thus if every $n\ge 1$ is a unit of $R$, then we may identify $\Q$ with its image in $R$ via the unique
ring embedding $\Q\to R$, making $\Q$ into an ordered subfield of $R$.

\begin{lemma}\label{pr2} Suppose $H$ is an ordered subfield of $R$,
all $g\in R$ with $g>H$ are units of $R$, and $H<f\in R$. Then we have an ordered subfield $H(f)$ of $R$. 
\end{lemma}
\begin{proof} 
For $P(Y)\in H[Y]$ of degree $d\ge 1$
with leading coefficient $a>0$ we have~$P(f)=af^d(1+\epsilon)$
with $-1/n < \epsilon < 1/n$ for all $n\ge 1$, in particular,
$P(f)>H$ is a unit of $R$. It remains to appeal to Lemma~\ref{pr0}.
\end{proof}

\begin{lemma}\label{pr1} Let $H$ be a real closed ordered subfield of $R$.
Let $A$ be a nonempty downward closed subset of $H$ such that
$A$ has no largest element and $B:=H\setminus A$ is nonempty and
has no least element. Let $f\in R$ be such that $A<f<B$. Then the subring $H[f]$ of $R$ has the following properties: \begin{enumerate}
\item[$\rm{(i)}$] $H[f]$ is a domain;
\item[$\rm{(ii)}$] $H[f]$ is an ordered subring of $R$;
\item[$\rm{(iii)}$] $H$ is cofinal in $H[f]$;
\item[$\rm{(iv)}$] for all $g\in H[f]\setminus H$ and $a\in H$, if $a<g$, then $a < b<g$ for some $b\in H$, and if $g<a$, then 
$g<b<a$ for some $b\in H$. 
\end{enumerate}
\end{lemma} 
\begin{proof} Let
$P\in H[Y]^{\neq}$; to obtain (i) and (ii) it suffices to show that then
$P(f) < 0$ or $P(f)>0$.  We have $$P(Y)\ =\ c\,Q(Y)\,(Y-a_1)\cdots(Y-a_n)$$ where
$c\in H^{\ne}$,  $Q(Y)$  is  a product of monic quadratic irreducibles in $H[Y]$, and  $a_1,\dots, a_n\in H$. This gives $\delta\in H^{>}$ such that $Q(r)\ge \delta$ for all $r\in R$.  Assume~$c>0$. (The case $c<0$ is handled similarly.) We can arrange that $m\le n$ is such that~$a_i\in A$ for $1\le i \le m$ and $a_j\in B$ for $m < j\le n$.
Take $\epsilon>0$ in $H$ such that
$a_i + \epsilon \le f$ for~$1\le i\le m$ and 
$f\le a_j-\epsilon$ for~$m < j\le n$.
Then $$P(f)\ =\ c\,Q(f)\,(f-a_1)\cdots (f-a_m)(f-a_{m+1})\cdots(f-a_n),$$
and $(f-a_1)\cdots(f-a_m) \ge \epsilon^m$. If
$n-m$ is even, then $(f-a_{m+1})\cdots(f-a_n)\ge \epsilon^{n-m}$, so $P(f)\ge c\delta\epsilon^n >0$. If $n-m$ is odd, then
$(f-a_{m+1})\cdots(f-a_n)\le -\epsilon^{n-m}$,
so~$P(f) \le -c\delta\epsilon^n < 0$. These estimates also yield (iii) and (iv). 
\end{proof}

\begin{lemma}\label{pr3} With $H$,~$A$,~$f$ as in Lemma~\ref{pr1}, suppose all $g\in R$ with $g\ge 1$ are units of $R$. Then we have an ordered subfield $H(f)$ of $R$ such that {\rm(iii)} and {\rm(iv)} of Lemma~\ref{pr1} go through for $H(f)$ in place of $H[f]$.
\end{lemma} 
\begin{proof} 
%It is clear that we have a subfield $K(f)$ of $R$. 
Note that if $g\in R$ and $g\ge \delta\in H^{>}$, then 
$g\delta^{-1}\ge 1$, so $g$ is a unit of $R$ and
$0<g^{-1}\le \delta^{-1}$. For   
$Q\in H[Y]^{\ne}$ with $Q(f)>0$ we can take $\delta\in H^{>}$ such that~$Q(f)\ge \delta$, so $Q(f)\in R^\times$  and $0 <Q(f)^{-1}\le \delta^{-1}$. Thus we have an
ordered subfield
$H(f)$ of $R$ by Lemma~\ref{pr0}, and 
the rest now follows easily.  
\end{proof}

\subsection*{Adjoining pseudolimits and increasing the value group} Let $H$ be a real closed Hausdorff field and  view $H$ as an ordered valued field as
before. Let $(a_{\rho})$ be a strictly increasing divergent pc-sequence in $H$.
Set 
$$A\ :=\ \{a\in H:\ \text{$a< a_{\rho}$ for some $\rho$}\}, \qquad
B\ :=\ \{b\in H:\ \text{$b>a_{\rho}$ for all $\rho$}\},$$ so
$A$ is nonempty and downward closed without a largest element.
Moreover,~$B=H\setminus A$ is nonempty and has no least element, since a
least element of $B$ would be a limit and thus a pseudolimit of
$(a_{\rho})$. Let $f\in \mathcal{C}$ satisfy $A<f<B$.
Then by Lem\-ma~\ref{pr3} for $R=\c$ we have an ordered subfield $H(f)$ of $\mathcal{C}$, and:

\begin{lemma} \label{ps1} $H(f)$ is an immediate valued field extension of $H$ with $a_{\rho} \leadsto f$.
\end{lemma}
\begin{proof} We can assume that $v(a_{\tau}-a_{\sigma})>v(a_{\sigma}-a_{\rho})$ for all indices $\tau>\sigma>\rho$. Set~$d_{\rho}:= a_{s(\rho)}-a_{\rho}$ ($s(\rho):=$~successor of $\rho$).
Then $a_{\rho}+2d_{\rho}\in B$ for all indices~$\rho$; see the discussion
preceding~[ADH,~2.4.2]. It then follows from that lemma that~$a_{\rho} \leadsto f$. Now $(a_{\rho})$ is a divergent pc-sequence in the henselian valued field $H$, so it is of transcendental type over $H$, and thus $H(f)$ is an immediate
extension of~$H$.  
\end{proof}

\begin{lemma}\label{ps2} Let $H$ be a Hausdorff field with divisible value group $\Gamma:=v(H^\times)$. Let $P$ be a nonempty upward closed subset of $\Gamma$, and let $f\in \mathcal{C}$ be such that 
$a < f$ for all~$a\in H^{>}$ with 
$va\in P$, and $f<b$ for all $b\in H^{>}$ with $vb < P$.
Then $f$ generates a Hausdorff field~$H(f)$ with 
$P >vf > Q$, $Q:= \Gamma\setminus P$.  
\end{lemma}
\begin{proof} For any positive $a\in H^{\text{rc}}$ there is
$b\in H^{>}$ with $a\asymp b$ and $a< b$, and also an element
$b\in H^{>}$ with $a\asymp b$ and $a>b$. Thus by Proposition~\ref{b1} we can replace~$H$ by~$H^{\text{rc}}$ and arrange in this way that $H$ is real closed.
 Set 
$$A\ :=\ \{a\in H:\ \text{$a\le 0$  or $va\in P$}\}, \qquad B:= H\setminus A.$$
Then we are in the situation of Lemma~\ref{pr1} for $R=\mathcal{C}$,   so by that lemma and Lemma~\ref{pr3} we have a Hausdorff field
$H(f)$. Clearly then $P> vf >Q$.     
\end{proof}

\subsection*{Non-oscillation} A germ $f\in \c$ is said to {\bf oscillate\/}  if $f(t)=0$ for arbitrarily large~$t$ and $f(t)\ne 0$ for
arbitrarily large $t$. \index{germ!oscillating} Thus for $f,g\in \c$,
$$ f-g \text{  is non-oscillating }\quad \Longleftrightarrow\quad \begin{cases} &\parbox{17em}{ either $f(t)< g(t)$ eventually, or $f=g$, or $f(t)>g(t)$ eventually.}\end{cases}$$
%\begin{align*} f-g \text{  is non-oscillating }\ \Longleftrightarrow\ &\text{ either $f(t)< g(t)$ eventually, or $f=g$}, \\
%&\text{ or $f(t)>g(t)$ eventually.}
%\end{align*} 
In particular, $f\in\c$ does not oscillate iff $f=0$ or $f\in\c^\times$.
If $g\in\c$ and~$g(t)\to +\infty$ as $t\to+\infty$, then $f\in\c$ oscillates iff $f\circ g$ oscillates. 

\begin{lemma}\label{lem:nonosc, 1}
Let $f\in\c$ be such that for every $q\in\Q$ the germ $f-q$ is non-oscillating. 
Then~$\lim\limits_{t\to\infty} f(t)$ exists in $\R\cup\{-\infty,+\infty\}$.
\end{lemma}
\begin{proof} Set $A\ :=\ \{q\in\Q: f(t)> q \text{ eventually}\}$.
If $A=\emptyset$, then $\lim\limits_{t\to\infty} f(t)=-\infty$, whereas if $A=\Q$, then 
$\lim\limits_{t\to\infty} f(t)=+\infty$.
If $A\neq\emptyset,\Q$, then for $\ell:=\sup A\in\R$ we have
%Let $\varepsilon\in\R^>$. Choose
%$a,b\in\Q$ with $\ell-\varepsilon\leq a<\ell<b\leq \ell+\varepsilon$; then 
%$a\in A$, $b\notin A$, so~$\ell-\varepsilon\leq a\leq f\leq b\leq \ell+\varepsilon$. Hence~
$\lim\limits_{t\to\infty} f(t)=\ell$.
\end{proof}

\begin{lemma}\label{lem:nonosc, 2}
Let $H$ be a real closed Hausdorff field and $f\in\c$. Then $f$ lies in a Hausdorff field extension of $H$ iff 
$f-h$ is non-oscillating for all $h\in H$. 
\end{lemma}
\begin{proof}
The forward direction is clear. For the converse, suppose $f-h$ is non-oscillating  for all $h\in H$.
We assume $f\notin H$, so $h<f$ or $h>f$ for all $h\in H$. 
Set~${A:=\{h\in H:h < f\}}$, a downward closed subset of $H$.
If $A=H$, then we are done by Lemma~\ref{pr2} applied to $R=\c$; if $A=\emptyset$ then 
we apply the same lemma to~$R=\c$ and $-f$ in place of $f$.
Suppose $A\neq\emptyset,H$.
If~$A$ has a largest element $a$, then we replace $f$ by~$f-a$ to arrange $0 < f (t)< h(t)$ eventually, for all $h\in H^>$, and
then  Lemma~\ref{pr2} applied to $R=\c$, $f^{-1}$ in place of $f$ yields that $f^{-1}$, and hence also $f$, lies in a Hausdorff field
extension of $H$. The case that $B:=H\setminus A$ has a least element is handled in the same way. 
If $A$ has no largest element and $B$ has no least element, then
we are done by Lemma~\ref{pr3}.
\end{proof}

\section{Linear Differential Equations}\label{sec:second-order}

\noindent
In this section we fix notations and conventions concerning differentiable functions and summarize well-known
results on linear differential equations as needed later, focusing on the case of order $2$.
We also discuss disconjugate linear differential equations, mainly following \cite[Chapter~3]{Coppel}, as well as work by Lyapunov and Perron on ``bounded'' matrix differential equations;  this material is only used in
Section~\ref{sec:lin diff applications} on applications and can be skipped upon first reading.

\subsection*{Differentiable functions} 
Let $r$ range over $\N\cup\{\infty\}$, and let $U$ be a nonempty open subset of $\R$. \label{p:Cr(U)}
Then $\c^r(U)$ denotes the $\R$-algebra of $r$-times
continuously differentiable functions $U\to \R$, with the usual pointwise defined algebra operations. (We use ``$\c$'' instead of ``$C$'' since $C$ will often denote the constant field of a
differential field.) For $r=0$ this is the $\R$-algebra~$\c(U)$ of continuous real-valued functions on~$U$, so 
$$\c(U)\ =\ \c^0(U)\ \supseteq\ \c^1(U)\ \supseteq\ \c^2(U)\ \supseteq\ 
\cdots\ \supseteq\ \c^{\infty}(U).$$ For $r\ge 1$ we have the derivation $f\mapsto f'\colon \c^r(U)\to \c^{r-1}(U)$ (with $\infty-1:=\infty$). This makes $\c^{\infty}(U)$ a differential ring, with
its subalgebra
$\c^{\omega}(U)$ 
of real-analytic functions $U\to \R$ as a differential subring. The algebra operations
on the algebras below are also defined pointwise.  
Note that
$$\c^r(U)^\times\  =\  \big\{ f\in\c^r(U):\text{$f(t)\neq 0$ for all $t\in U$}\big\},$$ 
also for $\omega$ in place of $r$ \cite[(9.2), ex.~4]{Dieudonne}. 

\medskip
\noindent
Let $a$ range over $\R$. \label{p:Car}
% and $r$ range over $\N$. 
Then $\Car$ denotes the $\R$-algebra of functions $[a,+\infty) \to \R$ that extend to a function in $\c^r(U)$ for some open $U\supseteq [a,+\infty)$. Thus $\Caz$ (also denoted by $\c_a$) is the $\R$-algebra of real-valued continuous functions on $[a,+\infty)$, and
$$\Caz\ \supseteq\ \Cao\ \supseteq\ \Cat\ \supseteq\ \cdots\ \supseteq \Cainf.$$  
We have the subalgebra $\Caom$ of $\Cainf$, consisting of the functions~${[a,+\infty)\to \R}$ that extend to a real-analytic function
$U \to \R$ for some open $U\supseteq [a,+\infty)$. 
For ${f\in\Cao}$ and~$g\in\c^1(U)$ extending $f$  
with open $U\subseteq \R$ containing $[a,+\infty)$, the restriction of~$g'$ to $[a,+\infty)\to\R$ depends only on $f$, not on $g$, so
we may define~${f':=g'|_{[a,+\infty)}\in\c_a}$.
For~$r\ge 1$ this gives the derivation $f\mapsto f'\colon \Car\to \Carl$. This makes $\Cainf$ a differential ring with~$\c^{\omega}_a$  as a differential subring.

\medskip
\noindent
For each of the algebras $A$ above we also consider its complexification $A[\imag]$ which consists by definition of the
$\C$-valued functions $f=g+h\imag$ with $g,h\in A$, so
$g=\Re f$ and $h= \Im f$ for such $f$. We consider $A[\imag]$ as a $\C$-algebra with respect to the natural pointwise defined algebra operations. We identify each complex number with
the corresponding constant function to make $\C$ a subfield of
$A[\imag]$ and~$\R$ a subfield of~$A$. (This justifies the notation $A[\imag]$.) 
We have $\Car[\imag]^\times=\c_a[\imag]^\times\cap\Car[\imag]$ and
$(\Car)^\times=\c_a^\times\cap\Car$, and likewise with $r$ replaced by $\omega$.

\medskip
\noindent
For $r\ge 1$ we extend $g\mapsto g'\colon \Car\to \Carl$ to the derivation 
$$g+h\imag\mapsto g'+h'\imag\ :\  \Car[\imag] \to 
\Carl[\imag] \qquad (g,h\in \Car[\imag]),$$
which for $r=\infty$ makes $\Cainf$ a differential subring of $\Cainf[\imag]$.  
We shall use the map
$$f \mapsto f^\dagger:=f'/f\ \colon \ \Cao[\imag]^\times=\big(\Cao[\imag]\big)\!^\times \to \Caz[\imag],$$
with 
$$(fg)^\dagger=f^\dagger+g^\dagger\qquad\text{ for $f,g\in \Cao[\imag]^\times$,}$$
in particular the fact that $f\in \Cao[\imag]^\times$ 
and $f^\dagger\in \Caz[\imag]$ are related by
$$ f(t)\ =\ f(a)\exp\!\left[\int_a^t f^\dagger(s)\,ds\right] \qquad (t\ge a).$$ 
For $g\in \Caz[\imag]$, let $\exp\int g$  denote the function   $t\mapsto \exp\!\left[\int_a^tg(s)\,ds\right]$ in $\Cao[\imag]^\times$.
Then
$$(\exp\textstyle\int g)^\dagger=g\quad\text{ and }\quad \exp\int (g+h)=(\exp\int g)\cdot (\exp\int h)\qquad\text{ for $g,h\in \Caz[\imag]$.}$$ 
Therefore~$f\mapsto f^\dagger\colon \Cao[\imag]^\times\to \Caz[\imag]$ is surjective. 

\begin{notation}
For $b\geq a$ and $f\in \c_a[\imag]$ we set  $f|_b:=f|_{[b,+\infty)}\in \c_a[\imag]$. 
\end{notation}
% is a $\C$-algebra endomorphism of  $\Car[\imag]$, which
%restricts to an $\R$-algebra endomorphism of $\Car$,  
%such that $(y')|_{b}=(y|_{b})'$ if~$r\geq 1$ and~$y\in \Car[\imag]$. \marginpar{added sentence} 

\subsection*{Differentiable germs} 
Let $r\in \N\cup\{\infty\}$ and let $a$ range over $\R$.\index{germ!differentiable}\label{p:Gr}
Let $\Gr$ be the partially ordered subring of 
$\c$ consisting of the germs at $+\infty$ of the functions in~$\bigcup_a \Car$; thus $\Gz=\c$ consists of the germs at~$+\infty$ of the continuous real valued functions on intervals $[a,+\infty)$, $a\in \R$. Note that $\Gr$ with its partial ordering satisfies the conditions on $R$ from  Section~\ref{sec:germs}. Also, every~$g\ge 1$ in~$\Gr$ is a unit of $\Gr$, so
Lemmas~\ref{pr2} and~\ref{pr3} apply to ordered subfields of $\Gr$. We have 
$$\Gz\ \supseteq\ \Go\ \supseteq\ \Gt\ \supseteq\ \cdots\ \supseteq\ \Ginf.$$
Each subring $\Gr$ of $\c$ yields the subring $\Gr[\imag]=\Gr+\Gr\imag$ of~$\Gz[\imag]=\c[\imag]$, with
$$\Gz[\imag]\ \supseteq\ \Go[\imag]\ \supseteq\ \Gt[\imag]\ \supseteq\ \cdots\ \supseteq\ \Ginf[\imag].$$
Suppose $r\geq 1$; then for $f\in\Car[\imag]$ the germ of $f'\in\Carl[\imag]$ only depends on the germ of $f$, and we thus obtain a derivation $g\mapsto g'\colon \Gr[\imag]\to\c^{r-1}[\imag]$
with $\text{(germ of $f$)}'=\text{(germ of $f'$)}$ for $f\in\bigcup_a \Car[\imag]$.
This derivation restricts to a derivation  $\Gr\to\c^{r-1}$. % and for  $r\geq 2$,  $g\mapsto g'\colon \c^r[\imag]\to\c^{r-1}[\imag]$ restricts to $h\mapsto h'\colon \c^{r-1}[\imag]\to\c^{r-2}[\imag]$.
Note that
$\c[\imag]^\times\cap \c^r[\imag] = \c^r[\imag]^\times$, and hence 
$\c^\times\cap \c^r = (\c^r)^\times$. 

\medskip
\noindent
For open $U\subseteq\C$ and $\Phi\colon U\to\C$ of class $\c^r$ (that is, its real and imaginary parts are of class $\c^r$), if $f\in\c^r[\imag]$ and $f(t)\in U$, eventually, then
$\Phi(f)\in\c^r[\imag]$. For example, if $f\in\c^r$, then $\exp f\in\c^r$, and if in addition $f(t)>0$, eventually, then $\log f\in\c^r$.

\medskip
\noindent
We set 
$$\Gi[\imag]\ :=\ \bigcap_{n}\, \Gn[\imag].$$ Thus $\Gi[\imag]$ is naturally a
differential ring with $\C$ as its ring of constants. \label{p:Gi} 
We also have the differential subring
$$\Gi \ :=\ \bigcap_{n}\, \Gn$$ of $\Gi[\imag]$, 
with $\R$ as its ring of constants and $\Gi[\imag]=\Gi+\Gi\imag$.  
Note that~$\Gi[\imag]$ has~$\Ginf[\imag]$ as a differential subring. Similarly, $\Gi$ has $\Ginf$ as a differential subring,
and the differential ring $\Ginf$ has in turn
the differential subring~$\Gom$, whose elements are the
germs at~$+\infty$ of the functions in $\bigcup_a \Caom$.\index{germ!analytic}\index{germ!smooth} 
We have~$\c[\imag]^\times\cap\Calinf[\imag]=(\Calinf[\imag])^\times$ and  
$\c^\times  \cap \Calinf = (\Calinf)^\times$, and likewise with~$\Gom$ in place of $\Calinf$.  If $R$ is a subring of~$\Go$ such that $f'\in R$ for all~$f\in R$, then~$R\subseteq\Gi$ is a differential subring of~$\Gi$.

\subsection*{Basic facts about linear differential equations}
In this subsection we review the main analytic facts about linear differential equations used later.
Let~$a\in\R$, $r\in\N^{\geq 1}$, and $f_1,\dots,f_{r}\in\c_a[\imag]$.
This gives the $\C$-linear map
$$y\mapsto A(y)\ :=\ y^{(r)}+f_1 y^{(r-1)}+\cdots+f_r y\ \colon\ \Car[\imag]\to\c_a[\imag].$$
We now have the classical existence and uniqueness theorem (see, e.g., \cite[(10.6.3)]{Dieudonne} or \cite[\S{}19, I, II]{Walter}):

\begin{prop} \label{prop:existence and uniqueness}
Let $t\in\R^{\geq a}$ be given. Then for any $b\in\c_a[\imag]$ and $c\in\C^r$ there is a unique~$y=y(b,c)\in \Car[\imag]$ such that 
$$A(y)\ =\ b, \qquad \big(y(t),y'(t),\dots,y^{(r-1)}(t)\big)\ =\ c.$$ 
The map 
$c\mapsto y(0,c)\colon\C^r\to \ker A$ is an isomorphism of $\C$-linear spaces, and so in particular, $\dim_{\C} \ker A=r$.  
\end{prop}

\begin{cor} \label{cor:existence and uniqueness}
Let $y\in \ker A$. If for some $t\in\R^{\geq a}$ we have $y^{(j)}(t)=0$ for~$j=0,\dots,r-1$, then $y=0$.
\end{cor}

\noindent
Proposition~\ref{prop:existence and uniqueness} and $\Re A(y)=A(\Re y)$ for  $f_1,\dots, f_r\in \c_a$ and $y\in \Car[\imag]$ give:

\begin{cor} \label{cor:existence and uniqueness, real} Suppose $f_1,\dots,f_r\in\c_a$ and $t\in\R^{\geq a}$. Then for any $b\in\c_a$ and~$c\in\R^r$ we have $y=y(b,c)\in \Car$, and the map 
$c\mapsto y(0,c)\colon\R^r\to \Car\cap\ker A$ is an isomorphism of $\R$-linear spaces.
\end{cor}

\noindent  
Let $b\in\c_a[\imag]$.
Using $y^{(r)}=b-\sum_{i=1}^rf_iy^{(r-i)}$ for  $y\in A^{-1}(b)\subseteq\Car[\imag]$ gives
$$b,f_1,\dots,f_r\in \Can[\imag]\ \Longrightarrow\ A^{-1}(b)\subseteq \c^{n+r}_a[\imag]$$
by induction on $n$.  
Hence  $b,f_1,\dots,f_r\in \c^\infty_a[\imag]\Rightarrow A^{-1}(b)\subseteq \c^\infty_a[\imag]$, in particular,
$f_1,\dots,f_r\in \c^\infty_a[\imag]\Rightarrow \ker A\subseteq \c^\infty_a[\imag]$. 
Also $b,f_1,\dots,f_r\in \Caom[\imag]\Rightarrow A^{-1}(b)\subseteq \Caom[\imag]$ by Lemma~\ref{smo} below (see also~\cite[(10.5.3)]{Dieudonne}), so
$f_1,\dots,f_r\in \Caom[\imag]\Rightarrow \ker A\subseteq \Caom[\imag]$.

\medskip\noindent
Let $y_1,\dots,y_r\in \Car[\imag]$. The {\it Wronskian}\/ $w=\wr(y_1,\dots,y_r)$ of $y_1,\dots,y_r$ is 
$$\wr(y_1,\dots,y_r)\ :=\ \det \begin{pmatrix} y_1 & \cdots & y_r \\
y_1' & \cdots & y_r' \\
\vdots & & \vdots \\
y_1^{(r-1)} & \cdots & y_r^{(r-1)}
\end{pmatrix} \in \Cao[\imag].$$
Hence if $w\neq 0$ (that is, $w(t)\ne 0$ for some $t\ge a$), then $y_1,\dots,y_r$ are $\C$-linearly independent.
The converse does not hold in general, even for $r=2$ and $y_1,y_2\in\c^\infty_a$, see~\cite{Bocher}, 
% \marginpar{result from \cite{Bocher} skipped for now, not used later} 
but we do have:

\begin{lemma}\label{lem:A from yjs}
The following are equivalent:
\begin{enumerate}
\item[\textup{(i)}] $w\in\c_a[\imag]^\times$ \textup{(}that is, $w(t) \neq 0$ for all $t\geq a$\textup{)};
\item[\textup{(ii)}]  $y_1,\dots,y_r$ is a basis of $\ker A$ for some choice of $f_1,\dots,f_r\in \c_a[\imag]$.
\end{enumerate}
\end{lemma}
\begin{proof}
For  (i)~$\Rightarrow$~(ii), assume $w\in\c_a[\imag]^\times$, and use that
$y_1,\dots, y_r\in \ker A$, where $A$ is the $\C$-linear differential operator given by
$$y\mapsto A(y)\ :=\ \wr(y_1,\dots,y_r,y)/w\ \colon\ \Car[\imag] \to\c_a[\imag].$$
For (ii)~$\Rightarrow$~(i), assume (ii) and suppose towards a contradiction that $t\ge a$ is such that $w(t)=0$. This gives $c_1,\dots, c_r\in \C$, not all $0$, such that for $y=\sum_{k=1}^r c_ky_k$ we have $y^{(j)}(t)=0$ for $j=0,\dots,r-1$. 
Hence $y=0$ by Corollary~\ref{cor:existence and uniqueness}.
\end{proof}

\noindent
Let now $y_1,\dots,y_r\in\ker A$. Then by the above
$$
w \neq 0 	\ \Longleftrightarrow\ w\in\Cao[\imag]^\times \ \Longleftrightarrow\ \text{$y_1,\dots,y_r$ are $\C$-linearly independent}.
$$
Moreover,
$w'=-f_1 w$ ({\it Abel's Identity}\/, see \cite[\S{}19, p.~200]{Walter}) and hence
$$w(t)\ =\ w(a)\exp\!\left({-\textstyle\int_a^t f_1(s)\,ds}\right)\quad\text{ for $t\geq a$.}$$
In particular,  $w=w(a)\in\C$ if $f_1=0$. 

\medskip
\noindent
In the next corollary  we let $g_1,\dots,g_r\in\c_a[\imag]$ and consider the $\C$-linear map
$$y\mapsto B(y)\ :=\ y^{(r)}+g_1 y^{(r-1)}+\cdots+g_r y\ \colon\ \Car[\imag]\to\c_a[\imag].$$
 
\begin{cor}\label{cor:uniqueness of coeffs}
$f_1=g_1,\dots,f_{r-1}=g_{r-1}\ \Longleftrightarrow\ 
 A=B \ \Longleftrightarrow \ \ker A=\ker B$.
\end{cor} 
\begin{proof}
Suppose $\ker A=\ker B$.
Let $y_1,\dots,y_r$ be a basis of $\ker A$, and
set $h_j:=f_j-g_j$  ($j=1,\dots,r-1$). Towards a contradiction
suppose $h_j\neq 0$ for some $j$, and take $j$ minimal with this property.
Take a nonempty open interval $I\subseteq\R^{\geq a}$ with~$h_j\in\c(I)[\imag]^\times$.
(Here and below we denote the restrictions of $h_1,\dots,h_{r-1}$ to functions $I\to\C$ by the same symbols.)
Then   $y_1,\dots,y_r$  restrict to
$\C$-linearly independent functions in $\c^r(I)[\imag]$ each satisfying the equation
$$y^{(r-j)}+(h_{j+1}/h_j)y^{(r-j-1)}+\cdots+(h_r/h_j)y=0,$$ 
contradicting \cite[\S{}19, II]{Walter}.
\end{proof}

\noindent
Next some basic properties of Wronskians:

\begin{lemma}\label{lem:Wronskian mult}
Let $u\in\Car[\imag]$. Then $\wr(uy_1,\dots,uy_r) = u^r \wr(y_1,\dots,y_r)$.
In particular, if $y_1\in \Car[\imag]^\times$, $r\geq 2$, and $z_j:=(y_{j+1}/y_1)'\in\Carl[\imag]$ \textup{(}$j=1,\dots,r-1$\textup{)}, then
$\wr(y_1,\dots,y_r) = y_1^r \wr(z_1,\dots,z_{r-1})$.
\end{lemma}
\begin{proof}
For the first identity, use that there are $u_{ij}\in\c_a[\imag]$ ($0\leq i\leq j<r$) with~${u_{0j}=u}$ such that
for all $y\in\Car[\imag]$ we have
$$(uy)^{(j)}\  =\  u_{0j} y^{(j)} + u_{1j} y^{(j-1)} + \cdots + u_{jj} y.$$
The first identity yields the second by taking $u:=y_1^{-1}$.
\end{proof}

\begin{lemma}\label{lem:Wronskian diff} 
Suppose $v:=\wr(y_1,\dots,y_{r-1})$ and $w:=\wr(y_1,\dots,y_r)$ lie in $\c_a^1[\imag]^\times$, with $v:= 1$ if $r=1$. 
Then we have for all $y\in\Car[\imag]$,
$$\big( \!\wr(y_1,\dots,y_{r-1},y)/w \big)'\  =\  (v/w^2)\wr(y_1,\dots,y_r,y).$$
\end{lemma}
\begin{proof} Expand the determinants $\wr(y_1,\dots, y_{r-1},y)$ and $\wr(y_1,\dots, y_r,y)$ according to their last column to get functions $g_1,\dots,g_{r-1},h_1,\dots,h_{r-1}\in\c_a[\imag]$ such that 
\begin{align*}
\big( \!\wr(y_1,\dots,y_{r-1},y)/w \big)' \ &=\  (v/w)y^{(r)}+g_1y^{(r-1)}+\cdots+g_ry, \\
(v/w^2)\wr(y_1,\dots,y_r,y)\ &=\  (v/w)y^{(r)}+h_1y^{(r-1)}+\cdots+h_ry
\end{align*}
for all $y\in\Car[\imag]$. Both left hand sides have  the $\R$-linearly independent $y_1,\dots, y_r$ among their zeros.
Now use Corollary~\ref{cor:uniqueness of coeffs}.
\end{proof}
 
\noindent
Let  $y\in\Car[\imag]$ and $t\geq a$.  The {\bf multiplicity} of $y$ at $t$\index{multiplicity!function}\label{p:multt}
is the largest~$m\leq r$   such that  
$y(t)=y'(t)=\cdots=y^{(m-1)}(t)=0$; notation: $\operatorname{mult}_t(y)$, or $\operatorname{mult}^r_t(y)$ if we need to indicate the dependence on $r$. So for $a\le 0$ we have $\operatorname{mult}^2_0(x^3)=2$, but~$\operatorname{mult}^r_0(x^3)=3$ for $r\geq 3$.)
Thus $t$ is a zero of $y$ (that is, $y(t)=0$) iff~$\operatorname{mult}_t(y)\geq 1$.\index{zero!multiplicity}
If $y\in\ker A$   has a zero of multiplicity~$r$, then~$y=0$ by Corollary~\ref{cor:existence and uniqueness}. Note that~$\operatorname{mult}_t(y) = \min\!\big\{ \!\operatorname{mult}_t(\Re y)
\operatorname{mult}_t(\Im y) \big\}$. For  $z\in\Car[\imag]$ we have
\[%\begin{equation}\label{eq:mult fms}
\operatorname{mult}_t(y+z) \geq \min\! \big\{ \!\operatorname{mult}_t(y),
\operatorname{mult}_t(z) \big\},
\]%\end{equation}
and  using the Product Rule:
$$\operatorname{mult}_t(yz)\ =\ \min\big\{r, \operatorname{mult}_t(y)+\operatorname{mult}_t(z)\big\}.$$
If $r\geq 2$ and $y(t)=0$, then $y'\in\Carl[\imag]$ and $\operatorname{mult}^{r-1}_t(y') = \operatorname{mult}^r_t(y)-1$.
The following is obvious:

\begin{lemma}\label{lem:Ch=>Ma}
Let  $y_1,\dots,y_r\in \Car[\imag]$, $w:=\wr(y_1,\dots,y_r)$, and $t\geq a$.
If~${w(t)=0}$, then  $\operatorname{mult}_t(y)= r$
for some   $\C$-linear combination  $y=c_1y_1+\cdots+c_ry_r$ of~$y_1,\dots,y_r$,
where~$c_1,\dots,c_r\in\C$ are not all zero.
\end{lemma}

\noindent
We also call the sum $$\operatorname{mult}(y)\ :=\ \sum_{t\geq a} \operatorname{mult}_t(y)\in \N\cup\{\infty\}$$
of the multiplicities of all zeros of $y$ the (total) {\bf multiplicity} of $y$, and we denote it by $\operatorname{mult}^r(y)$ if we need to exhibit the dependence on $r$.\index{multiplicity!total}\label{p:multtotal}
Note that~$\operatorname{mult}(y)<\infty$ iff $y$ has finitely many zeros.
% If $y$ or $z$ has only finitely many zeros, then the inequality~\eqref{eq:mult fms} holds with~$\operatorname{mult}$ in place of $\operatorname{mult}_t$, and
 If $z\in\Car[\imag]^\times$, then
$\operatorname{mult}(yz)=\operatorname{mult}(y)$.

\begin{lemma}\label{lem:mult of der}
Suppose $y\in\Car$, $r\geq 2$. Then \textup{(}with $\infty-1:=\infty$\textup{)}:
$$\operatorname{mult}^{r-1}(y')\  \geq\  \operatorname{mult}^r(y) - 1.$$
\end{lemma}
\begin{proof}
Let $m\leq\operatorname{mult}^r(y)$; it is enough to show that then $m-1\leq \operatorname{mult}^{r-1}(y')$. 
Let~$t_1<\cdots<t_n$ be zeros of $y$ such that $\sum_i \operatorname{mult}_{t_i}(y) \geq m$.
For $i=1,\dots,n-1$, Rolle's Theorem yields $s_i\in (t_i,t_{i+1})$ such that $y'(s_i)=0$.
Hence
\begin{align*}
m\  \leq\ \sum_{i=1}^n \operatorname{mult}^{r}_{t_i}(y)\ &=\  n+\sum_{i=1}^n \operatorname{mult}^{r-1}_{t_i}(y')\\ & \leq\
1 + \sum_{i=1}^{n-1} \operatorname{mult}^{r-1}_{s_i}(y') +  \sum_{i=1}^n \operatorname{mult}^{r-1}_{t_i}(y')\  \leq\ 1 + 
\operatorname{mult}^{r-1}(y')
\end{align*}
as required.
\end{proof}

%\medskip
%\noindent
%The Lagrange form of the remainder in the Taylor expansion yields:

%\begin{lemma}\label{lem:local min}
%Suppose $y\in\Car$,  $m=\operatorname{mult}_t(y)$ is even, and~$y^{(m)}(t)>0$. Then there is an open interval $I\subseteq\R$ containing $t$ such that~$y(s)>0$ for $s\in (I\cap\R^{\geq a})\setminus\{t\}$.
%\end{lemma}

\subsection*{Oscillation} 
Let  $y\in\c_a$. 
We say that $y$ {\bf oscillates} if its germ in $\c$ oscillates. So~$y$ does not oscillate iff~$\sgn y(t)$ is  constant, eventually.
If $y$ oscillates, then so does~$cy$ for $c\in\R^\times$.
%Note: $g\in\c$ does not oscillate iff~$g=0$ or $g\in \c^\times$. If $g\in\c$ and~$g(t)\to +\infty$ as $t\to+\infty$, then $y\in\c$ oscillates iff $y\circ g$ oscillates. 
If $y\in\Cao$ oscillates, then so does $y'\in\c_a$, by Rolle's Theorem.

\medskip\noindent
Let now   $r\in\N^{\geq 1}$ and $f_1,\dots,f_r\in\c_a$, and consider the $\R$-linear map
$$y\mapsto A(y)\ :=\ y^{(r)}+f_1 y^{(r-1)}+\cdots+f_r y\ \colon\ \Car \to\c_a.$$
By Corollary~\ref{cor:existence and uniqueness, real}, the $\R$-linear subspace
$\Car\cap\ker A$ of $\Car$ has dimension~$r$. 

\medskip\noindent
Let 
$y\in {\Car}\cap{\ker A}$, $y\ne 0$, and  let $Z:=y^{-1}(0)$ be the set of zeros of~$y$, so $Z\subseteq [a,+\infty)$ is closed in $\R$. By a {\em limit point\/}
of a set $S\subseteq \R$ we mean a point $b\in \R$ such that for every  real $\epsilon>0$ we have $0<|s-b|<\epsilon$ for some $s\in S$. \index{limit point}
  
\begin{lemma}\label{lem:no limit pt}
$Z$ has no limit points. 
\end{lemma}
\begin{proof}
For $j=0,\dots,r$
let $Z_j:=(y^{(j)})^{-1}(0)$ be the set of zeros of $y^{(j)}$, so  $Z=Z_0$. Each $Z_j$ is closed and hence
contains its limit points. If $t_0<t_1$ are in $Z_j$, $0\le j<r$, then $Z_{j+1}\cap (t_0,t_1)\neq\emptyset$,
by Rolle, 
hence each limit point of $Z_j$ is a limit point of $Z_{j+1}$.
Thus if  $t$ is a limit point of~$Z$, then  $t\geq a$ and 
$y(t)=y'(t)=\cdots=y^{(r-1)}(t)=0$, hence $y=0$ by Corollary~\ref{cor:existence and uniqueness}, a contradiction.
\end{proof}

\noindent
By Lemma~\ref{lem:no limit pt}, $Z\cap[a,b]$  is finite for every $b\ge a$. Thus
$$\text{$y$ does not oscillate}
\quad \Longleftrightarrow\quad   \text{$Z$ is finite}\quad \Longleftrightarrow\quad \text{$Z$ is bounded.}$$
If $t_0\in Z$    is not the largest element of $Z$, then $Z\cap (t_0,t_1)=\emptyset$ for some $t_1> t_0$ in~$Z$.
We say that a pair of zeros~$t_0<t_1$ of $y$ is {\bf consecutive} if $Z\cap(t_0,t_1)=\emptyset$.\index{zero!consecutive}

\medskip
\noindent
Next we consider the set~$Z_1:=(y')^{-1}(0)$ of stationary points of $y$.

\begin{lemma}\label{lem:stationary pts}
Suppose $f_r\in\c_a^\times$. Then $Z_1$ has no limit points.
\end{lemma}
\begin{proof}
The proof of Lemma~\ref{lem:no limit pt} shows that if $t$ is a limit point of~$Z_1$, then
$t\geq a$ and $y'(t)=y''(t)=\cdots=y^{(r)}(t)=0$,  and as $y\in \ker A$, this gives  $0=A(y)(t)=f_r(t)y(t)$, so $y(t)=0$,
and thus $y=0$, a contradiction.
\end{proof}

\noindent
Thus if $f_r\in\c_a^\times$, then $Z_1\cap [a,b]$ is finite for all $b\ge a$. 
  
\subsection*{Second-order differential equations} Let $f\in \c_a$, that is, $f\colon [a,\infty) \to \R$ is continuous. We consider the differential equation
\begin{equation}\label{eq:2nd order}\tag{L}
 Y'' + fY\ =\ 0. 
 \end{equation}
The solutions~${y\in \Cat}$ of \eqref{eq:2nd order} form
an $\R$-linear subspace $\Sol(f)$ of $\Cat$. \label{p:Sol} The
solutions~${y\in \Cat[\imag]}$ of \eqref{eq:2nd order} are the $y_1+y_2\imag$ with
$y_1, y_2\in \Sol(f)$ and form
a $\C$-linear subspace $\Sol_{\C}(f)$ of~$\Cat[\imag]$. For any complex numbers
$c$,~$d$ there is a unique solution~$y\in \Cat[\imag]$ of   \eqref{eq:2nd order} with~${y(a)=c}$ and $y'(a)=d$, and the map that assigns to $(c,d)$ in $\C^2$
this unique solution is an isomorphism $\C^2\to \Sol_{\C}(f)$ of 
$\C$-linear spaces; it restricts to an $\R$-linear bijection~${\R^2\to \Sol(f)}$. 
We have $f\in \Can\Rightarrow\Sol(f)\subseteq \c^{n+2}_a$ (hence~$f\in \c^\infty_a\Rightarrow {\Sol(f)\subseteq \c^\infty_a}$)  and $f\in \Caom\Rightarrow \Sol(f)\subseteq \Caom$.  
From~\cite[\S{}27, XI]{Walter}:

\begin{lemma}[Sonin-P\'olya]\label{lem:extremal pts}
Suppose $f\in(\c_a^1)^\times$, $y\in\Sol(f)^{\ne}$, and $t_0<t_1$ are stationary points of $y$.
If $f$ is increasing, then~$|y(t_0)|\geq |y(t_1)|$. If $f$ is decreasing, then~$|y(t_0)|\leq |y(t_1)|$.
 If $f$ is strictly increasing, respectively strictly decreasing, then these inequalities are strict. 
\end{lemma}
\begin{proof}
Put $u:=y^2+ \big((y')^2/f\big)\in\c_a^1$. Then $u'=-f'(y'/f)^2$.
Thus if $f$ is increasing, then  $u$ is decreasing, and
as $u(t_i)=y(t_i)^2$ for $i=0,1$, we get $|y(t_0)|\geq |y(t_1)|$.
The other cases are similar, using also Lemma~\ref{lem:stationary pts} for the strict inequalities.
\end{proof}

%\noindent
%By Lemmas~\ref{lem:stationary pts} and~\ref{lem:extremal pts}, 

\begin{lemma}\label{uniquestat}  Suppose  $f\in(\c_a^1)^\times$, $y\in\Sol(f)$,  and  $t_0<t_1$ are consecutive zeros of
$y$. Then there is exactly one stationary point of $y$ in the interval $(t_0,t_1)$.
\end{lemma}
\begin{proof} If $s_0< s_1$ were stationary points of $y$ in the interval $(t_0, t_1)$, then by Rolle~$y''$ and thus $y$ (in view of $y''=-fy$) would have a zero in the interval $(s_0, s_1)$. 
\end{proof}

\noindent
Let $y_1,y_2\in\Sol(f)$, with 
Wronskian $w=y_1y_2'-y_1'y_2$. Then $w\in\R$, and
$$w\neq 0\ \Longleftrightarrow\ \text{$y_1$,~$y_2$ are $\R$-linearly independent}.$$ 
%\marginpar{Bellman's lemma 2 has a mistake}
By~\cite[Chapter~6, Lemmas~2 and 3]{Bellman} we have:

\begin{lemma}\label{lem:2nd order inhom}
Let $y_1,y_2\in\Sol(f)$ be $\R$-linearly independent and $g\in\c_a$. Then
$$t\mapsto y(t)\ :=\ -y_1(t)\int_a^t \frac{y_2(s)}{w}g(s)\,ds + y_2(t)\int_a^t \frac{y_1(s)}{w} g(s)\,ds\ \colon\  [a,+\infty)\to\R$$
lies in $\Cat$ and satisfies $y''+fy=g$, $y(a)=y'(a)=0$. 
\end{lemma}

\begin{lemma}\label{lem:2nd lin indep sol}
Let $y_1\in\Sol(f)$ with $y_1(t)\neq 0$ for $t\geq a$. Then the function
$$t\mapsto y_2(t):=y_1(t)\int_a^t \frac{1}{y_1(s)^2}\,ds\colon [a,+\infty)\to\R$$
also lies in $\Sol(f)$, and $y_1$, $y_2$ are $\R$-linearly independent.
\end{lemma}

\noindent
From~\cite[Chapter~2, Lemma~1]{Bellman} we also recall:

\begin{lemma}[Gronwall's Lemma] \label{lem:gronwall}
Let   $C\in\R^{\ge}$,~$v,y\in\mathcal C_a$ satisfy~$v(t),y(t)\geq 0$ for all $t\geq a$ and
$$y(t)\ \leq\ C+ \int_a^t v(s)y(s)\,ds\quad\text{for all $t\geq a$.}$$
Then
$$y(t)\ \leq\ C\exp\!\left[ \int_a^t v(s)\,ds\right]\quad\text{for all $t\geq a$.}$$
\end{lemma}

\noindent
Here is a simpler differential version:

\begin{lemma}\label{grondiff} Let $u \in \c_a$ and $y \in \c^1_a$ satisfy $y’(t) \leq u(t) y(t)$ for all $t \geq a$. Then~$y(t) \leq y(a) \exp\!\left(\int_a^t u(s)\,ds\right)$ for all $t \geq a$.
\end{lemma}
\begin{proof} Put $z(t) := y(t) \exp\!\left(-\int_a^t u(s)\, ds\right)$ for $t\geq a$. Then $z\in \c^1_a$, and $z’(t) \leq 0$ for all $t\geq a$, so
$z(t) \leq z(a)=y(a)$ for all $t\geq a$. This yields the desired result. 
\end{proof}

\noindent
{\it In the rest of this subsection we assume that
$a\geq 1$ and that $c\in\R^>$  is such that~$|f(t)| \leq c/t^2$ for all $t\geq a$.}
Under this hypothesis, Lemma~\ref{lem:gronwall} yields the following bound on the growth of the solutions $y\in\Sol(f)$; the proof we give is similar to that of~\cite[Chap\-ter~6, Theorem~5]{Bellman}.

\begin{prop}\label{prop:bound} Let $y\in\Sol(f)$. Then there is $C\in\R^\geq$ such that $|y(t)| \leq Ct^{c+1}$ 
and $|y'(t)|\leq Ct^c$
for all $t\geq a$.
\end{prop}
\begin{proof}
Let $t$ range over $[a,+\infty)$.
Integrating $y''=-fy$ twice between $a$ and $t$, we obtain constants $c_1$, $c_2$ such that
for all $t$,
$$y(t)\ =\ c_1+c_2t -\int_a^t\int_a^{t_1} f(t_2)y(t_2)\,dt_2\,dt_1\ =\ c_1+c_2t-\int_a^t (t-s)f(s)y(s)\, ds$$
and hence, with $C:=|c_1|+|c_2|$,
$$ |y(t)|\ \leq\ Ct + t\int_a^t |f(s)|\cdot|y(s)|\, ds,\ 
\text{ so }\ 
\frac{|y(t)|}{t}\ \leq\ C + \int_a^t s|f(s)|\cdot\frac{|y(s)|}{s}\, ds.$$
Then by the lemma above, 
$$\frac{|y(t)|}{t}\ \leq\ C\exp\!\left[\int_a^t s|f(s)|\,ds\right]\ \leq\
C\exp\!\left[\int_1^t c/s\,ds\right]\ =\ Ct^c$$
and thus
$|y(t)|\leq Ct^{c+1}$. Now  
\begin{align*} y'(t)\ &=\ c_2-\int_a^t f(s)y(s)\,ds, \text{ so}\\
|y'(t)|\ &\leq\ |c_2|+\int_a^t |f(s)y(s)|\,ds\ \leq\ C+Cc\int_1^t s^{c-1}\,ds\\ &=\
C+Cc\left[\frac{t^c}{c}-\frac{1}{c}\right]\ =\ Ct^c. \qedhere
\end{align*}
\end{proof} 

\noindent
Let $y_1, y_2\in \Sol(f)$ be $\R$-linearly independent. Recall that $w=y_1y_2'-y_1'y_2\in\R^\times$. It follows that $y_1$ and~$y_2$ cannot be simultaneously very small:

\begin{lemma}\label{lem:bound} There is a positive constant $d$ such that
$$\max\!\big(|y_1(t)|, |y_2(t)|\big)\ \ge\ dt^{-c} \quad \text{ for all $t\ge a$.}$$
\end{lemma}
\begin{proof} 
Proposition~\ref{prop:bound} yields $C\in \R^>$ such that  $|y_i'(t)|\le Ct^c$ for
$i=1,2$ and all~$t\ge a$. Hence $|w|\le 2\max\!\big(|y_1(t)|, |y_2(t)|\big)Ct^c$ for $t\ge a$, so
\equationqed{ \max\!\big(|y_1(t)|, |y_2(t)|\big)\ \ge\ \frac{|w|}{2C}t^{-c}\qquad(t\ge a).}
\end{proof}

\begin{cor}\label{cor:bound}
Set $y:=y_1+y_2\imag$ and $z:=y^\dagger$. Then for some $D\in\R^>$,  
$$|z(t)|\ \leq\ Dt^{2c}\quad\text{ for all $t\geq a$.}$$
\end{cor}
\begin{proof}
Take $C$ as in the proof of Lemma~\ref{lem:bound}, and $d$ as in that lemma. Then
$$|z(t)|\ =\ \frac{|y_1'(t)+y_2'(t)\imag|}{|y_1(t)+y_2(t)\imag|}\ \leq\  \frac{|y_1'(t)|+|y_2'(t)|}{\max\!\big(|y_1(t)|, |y_2(t)|\big)}\ \leq\ \left(\frac{2C}{d}\right)t^{2c}$$
for $t\geq a$.
\end{proof}

\subsection*{More on oscillation} We continue with the study of \eqref{eq:2nd order}.
Sturm's Separation Theorem says that if $y,z\in\Sol(f)$ are $\R$-linearly independent and   $t_0<t_1$ are consecutive zeros of $z$,  then $(t_0,t_1)$ contains a unique 
zero  of $y$ \cite[\S{}27, VI]{Walter}.  Thus:

\begin{lemma}\label{lem:some vs all oscillate}
Some $y\in\Sol(f)^{\neq}$ oscillates  $\ \Longleftrightarrow\ $ every $y\in\Sol(f)^{\neq}$ oscillates.
\end{lemma}
%\begin{proof}
%Let $y,z\in\Sol(f)^{\neq}$ where $z$ is oscillating. If $y\in\R z$, then $y$ clearly is also oscillating. If $y\notin\R z$, then for each pair $t_0<t_1$ of consecutive zeros of $z$, the interval $(t_0,t_1)$ contains a zero of $y$, by Sturm; hence $y$ has infinitely many zeros, so $y$ is oscillating.
%\end{proof}

\noindent
We say that {\bf $f$ generates oscillations} if some element of $\Sol(f)^{\ne}$  oscillates.\index{germ!generates oscillations}\index{generates oscillations!germ}

\begin{lemma}
Let $b\in\R^{\geq a}$. Then 
$$\text{$f$ generates oscillations}\quad\Longleftrightarrow\quad\text{$f|_{b}\in\c_b$ generates oscillations.}$$
\end{lemma}
\begin{proof}
The forward direction is obvious. For the backward direction, use that every~$y\in\c^2_b$ with $y''+gy=0$ for $g:=f|_{b}$
extends uniquely to a solution of~\eqref{eq:2nd order}.
\end{proof}

\noindent
By this lemma, whether $f$ generates oscillations depends only on its germ in $\c$. So this induces the notion of an element of $\c$ generating oscillations. 
%By the previous lemma we may define: T a germ in~$\c$ {\bf generates oscillations} if it has a representative in $\Caz$ (for some $a\in\R$) which
%generates oscillations. 
Here is another result of Sturm~\cite[loc.~cit.]{Walter} that we use below:

\begin{theorem}[Sturm's Comparison Theorem] \label{thm:Sturm Comp}
Let $g\in\c_a$ with $f(t)\geq g(t)$ for all $t\geq a$. Let $y\in\Sol(f)^{\ne}$ and $z\in\Sol(g)^{\neq}$, and let $t_0<t_1$ be consecutive zeros of~$z$. Then either $(t_0,t_1)$ contains a zero of~$y$,  or on $[t_0, t_1]$ we have
 $f=g$ and~$y=cz$ for some constant $c\in\R^\times$. 
\end{theorem}

\noindent
Here is an immediate consequence:

\begin{cor} \label{cor:gen osc closed upward}
If $g\in\c_a$  generates oscillations and   $f(t)\geq g(t)$, eventually, then~$f$ also generates oscillations.
\end{cor}

\begin{example}
For $k\in\R^\times$ we have the differential equation of the harmonic oscillator,
$$y'' + k^2 y\ =\ 0.$$
A function $y\in \Cat$ is a solution iff for some real constants $c,t_0$ and all $t\ge a$,
$$y(t)\ =\  c\sin k(t-t_0).$$
For $c\ne 0$, any function $y\in \Cat$ as displayed oscillates.  
Thus
if $f(t)\geq \varepsilon$, eventually, for some constant~$\varepsilon>0$, then $f$ generates oscillations. 
\end{example}

\noindent
To \eqref{eq:2nd order} we associate the corresponding {\bf Riccati equation}
\begin{equation}\label{eq:Riccati}\tag{R}
z'+z^2+f\ =\ 0.
\end{equation}
Let $y\in\Sol(f)^{\neq}$ be a non-oscillating solution to \eqref{eq:2nd order}, 
and take $b\geq a$ with~${y(t)\neq 0}$ for $t\geq b$. Then the function
$$t\mapsto z(t)\ :=\ y'(t)/y(t)\ \colon\ [b,+\infty)\to\R$$ 
in $\mathcal C^1_b$ satisfies \eqref{eq:Riccati}.
Conversely, if $z\in \mathcal C^1_b$ ($b\geq a$) is a solution to \eqref{eq:Riccati}, then 
$$t\mapsto y(t)\ :=\ \exp \left(\int_b^t z(s) \,ds\right)\ \colon\ [b,+\infty)\to\R$$
is a non-oscillating solution to \eqref{eq:2nd order} with $y\in(\c^1_b)^\times$ and $z=y^\dagger$.

\medskip
\noindent
Let $g\in \Cao$, $h\in\Caz$ and consider the second-order linear differential equation 
\begin{equation}\label{eq:2nd order, gen}\tag{$\tilde{\operatorname{L}}$}
y''+gy'+hy\ =\ 0.
\end{equation}

\begin{cor}\label{coroscgen}
Set $f:=-\frac{1}{2}g'-\frac{1}{4}g^2+h\in \c_a$. Then the following are equivalent:
\begin{enumerate}
\item[\textup{(i)}] some nonzero solution of \eqref{eq:2nd order, gen}  oscillates;
\item[\textup{(ii)}] all nonzero solutions of \eqref{eq:2nd order, gen} oscillate; 
\item[\textup{(iii)}] $f$ generates oscillations.
\end{enumerate}
\end{cor}
\begin{proof}
Let $G\in (\Cat)^\times$ be given by  
$G(t):=\exp\!\left(-\frac{1}{2}\int_a^t g(s)\,ds\right)$. 
Then $y\in\Cat$ is a solution to \eqref{eq:2nd order} iff $Gy$ is a solution to \eqref{eq:2nd order, gen}; cf. [ADH, 5.1.13].
\end{proof}

\subsection*{More on non-oscillation}  We continue with \eqref{eq:2nd order}. 
Let~$y_1$,~$y_2$ range over elements of $\Sol(f)$, and recall that its  Wronskian~$w=y_1y_2'-y_1'y_2$ lies in $\R$.

\begin{lemma}\label{y1y2}
Suppose $b\geq a$ is such that $y_2(t) \neq 0$ for $t\geq b$. Then for~$q:=y_1/y_2\in\c^2_b$ we have
$q'(t)=-w/y_2(t)^2$  for $t\geq b$,
so $q$ is monotone and $\lim_{t\to\infty} q(t)$ exists in $\R\cup\{-\infty,+\infty\}$.
\end{lemma}

\noindent
This leads to:

\begin{cor}\label{cor:I1 I2}
Suppose  $b\geq a$ and $y_1(t),y_2(t)\neq 0$ for $t\geq b$. For $i=1,2$, set
$$h_i(t)\  :=\  \int_b^t \frac{1}{y_i(s)^2}  \, ds\quad \text{ for $t\ge b$, so $h_i\in \mathcal{C}_b^3$}.$$
Then:  $y_1\prec y_2\ \Longleftrightarrow\ h_1 \succ 1 \succeq h_2$.
\end{cor}
\begin{proof} Suppose $y_1\prec y_2$. Then $y_1$, $y_2$ are $\R$-linearly independent, so $w\ne 0$. Moreover, $q\prec 1$ with $q$ as in in Lemma~\ref{y1y2}, and $q'=-wh_2'$ by that lemma, so $q+wh_2$ is constant, and thus $h_2\preceq 1$. 
Note that $h_1$ is strictly increasing. If $h_1(t)\to r\in \R$ as~$t\to +\infty$,  then
$z:=(r-h_1)y_1\in \Sol(f)$ by Lemma~\ref{y1y2} with $y_1$ and $y_2$ interchanged, and $z\prec y_1$, so $z=0$, hence $h_1=r$, a contradiction. Thus $h_1\succ 1$. 

For the converse, suppose $h_1\succ 1\succeq h_2$. Then $y_1$, $y_2$ are $\R$-linearly independent, so $w\ne 0$. From $h_2\preceq 1$  and $q+wh_2$ being constant we obtain $q\preceq 1$. If $q(t)\to r\ne 0$ as $t\to +\infty$, then
$y_1=qy_2\asymp y_2$, and thus $h_1\asymp h_2$, a contradiction. Hence $q\prec 1$, and thus $y_1\prec y_2$. 
\end{proof}

\noindent
The pair $(y_1,y_2)$ is said to be a {\bf principal system} of solutions of \eqref{eq:2nd order} if\index{solution!principal system}\index{system!principal}\index{principal!system}
\begin{enumerate}
\item $y_1(t),y_2(t)>0$ eventually, and
\item $y_1 \prec y_2$.
\end{enumerate}
Then $y_1$, $y_2$ form a basis of the $\R$-linear space $\Sol(f)$, and $f$ does not generate oscillations, by Lemma~\ref{lem:some vs all oscillate}. Moreover, for  
$y=c_1 y_1+c_2 y_2$ with $c_1,c_2\in\R$, $c_2\neq 0$ we have $y\sim c_2 y_2$.
Here are some facts about this notion:

\begin{lemma}\label{lem:admissible pair, unique}
If $(y_1,y_2)$, $(z_1,z_2)$ are principal systems of solutions of \eqref{eq:2nd order}, then
there are $c_1,d_1,d_2\in\R$ such that
$z_1 = c_1 y_1$, $z_2 = d_1 y_1+d_2 y_2$, and $c_1,d_2>0$.
\end{lemma}

\begin{lemma}\label{lem:admissible pair}
Suppose $f$ does not generate oscillations. Then \eqref{eq:2nd order} has a principal system of solutions.
\end{lemma}
\begin{proof}
It suffices to find a basis $y_1$, $y_2$ of $\Sol(f)$ with $y_1\prec y_2$.  
Suppose $y_1$, $y_2$ is any basis of $\Sol(f)$, and set   $c:=\lim_{t\to\infty} y_1(t)/y_2(t)\in\R\cup\{-\infty,+\infty\}$.
If $c=\pm\infty$, then    interchange $y_1$, $y_2$, otherwise   replace $y_1$ by $y_1-cy_2$. Then $c=0$, so $y_1\prec y_2$.
\end{proof}

\noindent
One calls $y_1$  a {\bf principal} solution of \eqref{eq:2nd order} if $(y_1,y_2)$ is a principal system of solutions of \eqref{eq:2nd order}  for some $y_2$.\index{solution!principal}\index{principal!solution} (See \cite[Theorem~XI.6.4]{Hartman} and \cite{Leighton,LeightonMorse}.) By the previous two lemmas, \eqref{eq:2nd order} has a principal solution iff $f$ does not generate oscillations, and any two principal solutions differ by a multiplicative factor in $\R^>$.
If $y_1\in (\c_a)^\times$ and $y_2$ is given as in Lemma~\ref{lem:2nd lin indep sol}, then~$y_2$  is a non-principal solution of  \eqref{eq:2nd order} and~$y_1\notin \R y_2$.  

\subsection*{Chebyshev systems and Markov systems\astr} 
Let~$r\in\N^{\geq 1}$ and $y_1,\dots,y_r\in\Car$, and  let $V$ be the $\R$-linear subspace of $\Car$ spanned by $y_1,\dots,y_r$. We call~$y_1,\dots,y_r$
 a {\bf Chebyshev system} (on $\R^{\geq a}$) if  for  all~$y=c_1y_1+\cdots+c_ry_r$ with $c_1,\dots,c_r\in\R$ not all zero, we have~$\operatorname{mult}^r(y)<r$.
Note that if~$y_1,\dots,y_r$ is  a Chebyshev system, then~$y_1,\dots,y_r$ are $\R$-linearly independent, and
every  basis of $V$ is a Chebyshev system. 
Chebyshev systems can be used for interpolation:\index{system!Chebyshev}\index{Chebyshev system}

\begin{lemma}\label{lem:interpolation}
Suppose $y_1,\dots,y_r$ are $\R$-linearly independent.
Let~$t_1,\dots,t_n\in\R^{\geq a}$ be pairwise distinct and let~$r_1,\dots,r_n\in\N$ satisfy~$r_1+\cdots+r_n=r$. $($So $n\ge 1)$. Then the following are equivalent:
\begin{enumerate}
\item[\textup{(i)}]  the only   
$y\in V$ with $\operatorname{mult}^r_{t_i}(y)\geq r_i$ for $i=1,\dots,n$ 
is $y=0$;
\item[\textup{(ii)}] for all $b_{ij}\in\R$ \textup{(}$i=1,\dots,n$, $j=1,\dots,r_i$\textup{)}, there
exists $y\in V$ with
$$y^{(j-1)}(t_i)\ =\ b_{ij}\qquad (i=1,\dots,n,\ j=1,\dots,r_i).$$
\end{enumerate}
Moreover, in this case, given any $b_{ij}$ as in  \textup{(ii)}, the element $y\in V$ in  \textup{(ii)}
is unique.
\end{lemma}
\begin{proof}
Each
$y\in V$ equals $c_1y_1+\cdots+c_ry_r$ for a unique~$(c_1,\dots,c_r)\in\R^r$.
Let the~$b_{ij}$ in $\R$ be as in (ii), and set
$$M:=\begin{pmatrix}
\hfill y_1(t_1) 			& \cdots	& \hfill y_r(t_1) \\
\ \ \vdots				&			&\ \ \vdots \\
y_1^{(r_1-1)}(t_1) 	& \cdots	& y_r^{(r_1-1)}(t_1) \\
\hfill y_1(t_2)			& \cdots	& \hfill y_r(t_2) \\
\ \ \vdots				&			&\ \  \vdots \\
y_1^{(r_n-1)}(t_n) 	& \cdots	& y_r^{(r_n-1)}(t_n)
\end{pmatrix}\in\R^{r\times r}, \quad 
b:=\begin{pmatrix} b_{11} \\ \vdots \\ b_{1r_1} \\ b_{21} \\ \vdots \\ b_{nr_n}
\end{pmatrix}\in\R^r.$$
Then given $c=(c_1,\dots,c_r)^{\operatorname{t}}\in\R^r$, the element $y=c_1y_1+\cdots+c_ry_r$ of $V$
satisfies the  inequalities in (i) iff $Mc=0$, and
the displayed equations in~(ii) iff~$Mc=b$.  Thus (i) means injectivity of $M: \R^r\to \R^r$, and (ii)
its surjectivity. 
%$(Hence the claim follows from $\dim \ker M+\dim \operatorname{im} M=r$.
\end{proof}

\noindent
In particular, if   $y_1,\dots,y_r$ is a Chebyshev system, then  for all  $t_i$, $r_i$ ($i=1,\dots,n$) as in the previous lemma and for   all~$b_{ij}\in\R$ \textup{(}$i=1,\dots,n$, $j=1,\dots,r_i$\textup{)}, there
is a unique $y\in V$ with
$y^{(j-1)}(t_i)=b_{ij}$ ($i=1,\dots,n$, $j=1,\dots,r_i$).

\begin{remark}
Suppose $y_1,\dots,y_r$ are $\R$-linearly independent.
If $y_1,\dots,y_r$ is a Chebyshev system, then each $y\in V^{\neq}$ has~$<r$  zeros. Remarkably, the converse of 
this implication also holds; this is due to Aram\u{a}~\cite{Arama} and (in greater generality) Hartman~\cite{Hartman58};   a simple proof, from~\cite{Opial}, is in 
%\marginpar{converse skipped} 
\cite[Chapter~3, Proposition~3]{Coppel}. This links the notion of Chebyshev system considered here
  with the concept of the same name   in approximation theory~\cite[Chapter~3, \S{}4]{Cheney}. (These remarks are not used later.)
\end{remark} 

\noindent
If $y_1,\dots,y_r$ is a Chebyshev system,   then $\wr(y_1,\dots,y_r)\in\c_a^\times$ by Lemma~\ref{lem:Ch=>Ma}. If~$\wr(y_1,\dots,y_j)\in\c_a^\times$   for $j=1,\dots,r$, then $y_1,\dots,y_r$ is  called a {\bf Markov system} (on $\R^{\ge a}$).\index{system!Markov}\index{Markov system}
Thus by Lemma~\ref{lem:Ch=>Ma}, if 
$y_1,\dots,y_j$ is a Chebyshev system for~$j=1,\dots,r$, then $y_1,\dots, y_r$ is a Markov system.   Here is a partial converse:

\begin{lemma}\label{lem:Ma=>Ch}
If $y_1,\dots,y_r$ is a Markov system, then it is  a Chebyshev system.
\end{lemma}
\begin{proof}
The case $r=1$ is trivial, so let~$r\geq 2$ and let
$y_1,\dots,y_r$ be a Markov system; in particular,~$y_1\in\c_a^\times$. Set $z_j:=(y_{j+1}/y_1)'\in\Carl$ for~$j=1,\dots,r-1$. Then~$z_1,\dots,z_{r-1}$ is a Markov system by Lemma~\ref{lem:Wronskian mult}. Assume inductively that it is a Chebyshev system, and let~$y=c_1y_1+\cdots+c_ry_r$, $c_1,\dots,c_r\in\R$ not all zero; we need to show $\operatorname{mult}(y)<r$.  Towards a contradiction, assume~$\operatorname{mult}(y)\geq r$. 
Then~$z:=(y/y_1)'$ satisfies~$\operatorname{mult}(z)\geq  r-1$, by Lemma~\ref{lem:mult of der} and the remarks before it.
Moreover,~$z=c_2z_1+\cdots+c_rz_{r-1}$, and so  
$c_2=\cdots=c_r=0$ and hence~$y=c_1y_1$, and thus $c_1=0$, a contradiction.
\end{proof}

\noindent
If $y_1,\dots,y_r$ is a Markov system and $b\geq a$, then $y_1|_b,\dots,y_r|_b$ is a Markov system on $\R^{\ge b}$, and likewise with ``Chebyshev'' in place of ``Markov''.

\subsection*{Disconjugacy\astr}  
Let $r\in\N^{\geq 1}$ and $f_1,\dots,f_{r}\in\c_a$, and consider the linear differential equation
\begin{equation}\label{eq:disconj}\tag{D}
y^{(r)}+f_1y^{(r-1)}+\cdots+f_ry\ =\ 0
\end{equation}
on $\R^{\geq a}$. Let $\Sol\eqref{eq:disconj}$ be its set of solutions in $\Car$, so $\Sol\eqref{eq:disconj}$ is the kernel of the $\R$-linear map
$$y\mapsto  A(y):=y^{(r)}+f_1y^{(r-1)}+\cdots+f_ry\ \colon\ \Car \to\c_a.$$
Recall that by Corollary~\ref{cor:existence and uniqueness, real} we have $\dim_{\R}\Sol\eqref{eq:disconj}=r$. 
The linear differential equation \eqref{eq:disconj}
is said to be {\bf disconjugate} if $\Sol\eqref{eq:disconj}$ contains a Chebyshev system; that is, every nonzero $y\in\Sol(D)$
has multiplicity~$<r$.\index{disconjugate}\index{linear differential operator!disconjugate}
If \eqref{eq:disconj} is disconjugate, then it has no oscillating solutions.

\begin{example}
The equation $y^{(r)}=0$ is disconjugate, since its solutions in $\Car$ are the polynomial functions  $c_0+ c_1 x+ \cdots + c_{r-1}x^{r-1}$ with $c_0,\dots, c_{r-1}\in \R$.
%, \dots, x^{r-1}$ as aevery one-variable polynomial with real coefficients of degree~$< r$ has less than $r$ zeros
%in $\R$, counted with multiplicities.
\end{example}

\noindent
From Lemma~\ref{lem:interpolation} we obtain:

\begin{cor}[{de la Vall\'ee-Poussin~\cite{dVP}}]
Suppose \eqref{eq:disconj} is disconjugate. Then for all  pairwise distinct~$t_1,\dots,t_n\geq a$, all
  $r_1,\dots,r_n\in\N$ with~$r_1+\cdots+r_n=r$, and  all $b_{ij}\in\R$ \textup{(}$i=1,\dots,n$, $j=1,\dots,r_i$\textup{)}, there
is a unique $y\in \Sol\eqref{eq:disconj}$ such that
$$y^{(j-1)}(t_i)\ =\ b_{ij}\qquad (i=1,\dots,n,\ j=1,\dots,r_i).$$
\end{cor}

 %Combining Lemma~\ref{lem:Ma=>Ch} with the following lemma shows that the linear differential equation
%\eqref{eq:disconj} is disconjugate iff
% $\Sol\eqref{eq:disconj}$ contains a Markov system:
\noindent
Let $b\geq a$ and set $g_j:=f_j|_b\in\c_b$ for $j=1,\dots,r$. This yields the linear
differential equation
\begin{equation}\label{eq:disconj, b}\tag{D$_b$}
y^{(r)}+g_1y^{(r-1)}+\cdots+g_ry\ =\ 0 
\end{equation}
on $\R^{\geq b}$ with the $\R$-linear isomorphism $y\mapsto y|_b\colon \Sol\eqref{eq:disconj}\to \Sol\eqref{eq:disconj, b}$.

\begin{cor}\label{lem:disconj => Ma}
If \eqref{eq:disconj} is disconjugate, then some basis $y_1,\dots, y_r$ of the $\R$-linear space $\Sol\eqref{eq:disconj}$ yields for every
$b>a$  a Markov system $y_1|_{b},\dots, y_r|_{b}$ on $\R^{\ge b}$. 
\end{cor}
\begin{proof}
Let $y_1,\dots,y_r\in\Car$ be solutions of \eqref{eq:disconj} such that
$$y_j(a)\ =\ y_j'(a)\ =\ \cdots\ =\ y_j^{(r-j-1)}(a)\ =\ 0,\ y_j^{(r-j)}(a)\neq 0\qquad\text{for $j=1,\dots,r$.}$$
Then $\wr(y_1,\dots,y_r)(a)\neq 0$, so $y_1,\dots,y_r$ are $\R$-linearly independent.
Suppose~\eqref{eq:disconj} is disconjugate.  Let $j\in\{1,\dots,r\}$, $t\in\R^{> a}$. Then
$\wr(y_1,\dots,y_j)(t)\neq 0$:
otherwise  Lemma~\ref{lem:Ch=>Ma} yields an   $\R$-linear combination $y\neq 0$ of
$y_1,\dots,y_j$
with~${\operatorname{mult}_t(y)\geq j}$, but also $\operatorname{mult}_a(y)\geq r-j$ by choice of $y_1,\dots,y_r$, hence~${\operatorname{mult}(y)\geq r}$, contradicting disconjugacy of \eqref{eq:disconj}. Thus $y_1,\dots, y_r$ has the desired property. 
\end{proof}

\noindent
With $n\geq 1$ understood from the context, let $\der$ denote the $\R$-linear map $$y\mapsto y'\ \colon\ \Can\to\c_a^{n-1},$$ identify $f\in\c_a^{n-1}$ with the $\R$-linear operator~$y\mapsto fy\colon\c_a^{n-1}\to\c_a^{n-1}$, and
for maps~$\alpha\colon\c_a^{n}\to \c_a^{n-1}$,   $\beta\colon \c_a^{n+1}\to\c_a^{n}$, denote $\alpha\circ\beta\colon\c_a^{n+1}\to\c_a^{n-1}$ simply by~$\alpha\beta$.
With these conventions we can state an analytic version of Lemma~\ref{lem:Polya fact}: 

\begin{lemma}
If  $g_j\in  (\c_a^{r-j+1})^\times$ for $j=1,\dots,r$ and we set
\begin{equation}\label{eq:Polya fact}
A\ =\  g_1\cdots g_r (\der g_r^{-1}) \cdots (\der g_2^{-1})( \der g_1^{-1})\ :\  \c_a^r\to \c_a,
\end{equation}
then  $A = (\der-h_r)\cdots(\der-h_1)$ with $h_j:=(g_1\cdots g_j)^\dagger\in\c_a^{r-j}$ for $j=1,\dots,r$.
Conversely, if $h_j\in\c_a^{r-j}$ for $j=1,\dots,r$ and
$A:=(\der-h_r)\cdots(\der-h_1)\colon \c_a^r\to \c_a$ and $h_0:=0$,  then \eqref{eq:Polya fact} holds for $g_j\in  (\c_a^{r-j+1})^\times$ given by
$$g_j(t)\ :=\ \exp \int_a^t \big(h_j(s)-h_{j-1}(s)\big)\,ds\qquad (j=1,\dots,r),$$
and $h_j=(g_1\cdots g_j)^\dagger$ for those $j$.
\end{lemma}

\noindent
We now link the notion of disconjugacy with factorization of the operator $A=\der^r+ f_1\der^{r-1} + \cdots + f_r\colon \c_a^r\to \c_a$
considered earlier in connection with \eqref{eq:disconj}.

\begin{prop}[{Frobenius~\cite{Frobenius}, Libri~\cite{Libri}}]\label{prop:Frobenius}
Suppose $y_1,\dots,y_r\in\Sol\eqref{eq:disconj}$ is a Markov system. Set $w_0:=1$, $w_j:=\wr(y_1,\dots,y_j)\in(\c_a^{r-j+1})^\times$ for $j=1,\dots,r$, and
$$g_1\ :=\ w_1,\qquad g_j\ :=\ w_jw_{j-2}/w_{j-1}^2\quad (j=2,\dots,r).$$
Then $g_j\in (\c_a^{r-j+1})^\times$ for $j=1,\dots,r$ and  $A = g_1\cdots g_r (\der g_r^{-1})  \cdots (\der g_2^{-1})( \der g_1^{-1})$.
\end{prop}
\begin{proof}
It is clear that $g_j\in  (\c_a^{r-j+1})^\times$  and easy to check that
$w_j/w_{j-1}=g_1\cdots g_j$ for $j=1,\dots,r$.
We define for $j=0,\dots,r$ the 
$\R$-linear map
$$y\mapsto A_j(y):=\wr(y_1,\dots,y_j,y)/w_j\ \colon\ \c_a^r\to \c_a^{r-j}.$$
We claim that
$A_j = g_1\cdots g_j \der g_j^{-1} \der \cdots \der g_1^{-1}$.
The case $j=0$ is trivial. Suppose the claim holds for a certain $j<r$. Then
$$ g_1\cdots g_{j+1} \der g_{j+1}^{-1}\der \cdots \der g_2^{-1} \der g_1^{-1} \ =\
g_1\cdots g_{j+1} \der (g_1\cdots g_{j+1})^{-1} A_{j},$$
which sends $y\in \Car$ to 
$$\frac{w_{j+1}}{w_j}\left(\frac{w_j}{w_{j+1}}\frac{\wr(y_1,\dots,y_j,y)}{w_j}\right)'=
\frac{w_{j+1}}{w_j}\left(\frac{\wr(y_1,\dots,y_j,y)}{w_{j+1}}\right)',$$ and this in turn equals
$A_{j+1}(y)=\wr(y_1,\dots,y_j, y_{j+1},y)/w_{j+1}$ by Lemma~\ref{lem:Wronskian diff}.
\end{proof}

\noindent
Here is a converse, with $A$ still the operator $\der^r +f_{1}\der^{r-1}+ \dots +f_r$ figuring in \eqref{eq:disconj}:

\begin{theorem}[{P\'olya \cite{Polya}}]\label{thm:Polya} Suppose 
$$A = g_1\cdots g_r (\der g_r^{-1} ) \cdots (\der g_2^{-1})( \der g_1^{-1})\quad\text{ with $g_j\in  (\c_a^{r-j+1})^\times$  for $j=1,\dots,r$.}$$ 
Then $\Sol\eqref{eq:disconj}$ contains a Markov system $y_1,\dots, y_r$.  
\end{theorem}
%\begin{multline*}
%\text{\eqref{eq:disconj} is disconjugate}\quad\Longleftrightarrow\quad \\
%\text{there are $g_j\in  (\c_a^{r-j+1})^\times$   \textup{(}$j=1,\dots,r$\textup{)} with
 %$A = g_1\cdots g_r \der g_r^{-1} \der \cdots \der g_2^{-1} \der g_1^{-1}$.}
%\end{multline*}
%\end{theorem}
\begin{proof}
%The forward direction follows from Lemma~\ref{lem:disconj => Ma} and Proposition~\ref{prop:Frobenius}.
%For the converse suppose
Let $t_1,\dots, t_r$ range over $\R^{\ge a}$ and define~$y_1,\dots,y_r\in\Car$ by
\begin{align*}
y_1(t_1)\ &:=\ g_1(t_1),\\
y_2(t_1)\ &:=\ g_1(t_1)\int_a^{t_1} g_2(t_2)\,dt_2,\\
&\quad\vdots\\
y_r(t_1)\ &:=\ g_1(t_1)\int_a^{t_1} g_2(t_2) \int_a^{t_2} \cdots
\int_a^{t_{r-1}} g_{r}(t_{r})\, dt_{r}\cdots dt_2. 
\end{align*}
For $j=1,\dots,r$ we have   $A(y_j)=0$, and by an induction using Lem\-ma~\ref{lem:Wronskian mult}, $\wr(y_1,\dots,y_j)=g_1^j g_2^{j-1}\cdots g_j$. So
$y_1,\dots,y_r\in\Sol\eqref{eq:disconj}$ is a Markov system.
% so \eqref{eq:disconj} is disconjugate by Lem\-ma~\ref{lem:Ma=>Ch}.
\end{proof}

\begin{remark}
Suppose $\Sol\eqref{eq:disconj}$ contains a Markov system $y_1,\dots, y_r$. If   $f_1,\dots,f_r\in\c_a^n$, then
$y_1,\dots, y_r \in  \c_a^{n+r}$, so $g_j\in  (\c_a^{n+r-j+1})^\times$ for $j=1,\dots,r$ where $g_1,\dots,g_r$ are as in Proposition~\ref{prop:Frobenius}. 
Likewise, if $f_1,\dots,f_r\in\c_a^\infty$, then those $g_j$ lie in  $(\c_a^{\infty})^\times$, and the same with $\omega$ in place of $\infty$.
\end{remark}

\begin{cor}
With $A=\der^r +f_{1}\der^{r-1}+ \dots +f_r:\c_a^r\to \c_a$,
$\Sol\eqref{eq:disconj}$ contains a Markov system iff 
there are~$g_j\in  (\c_a^{r-j+1})^\times$ \textup{(}$j=1,\dots,r$\textup{)} such that
$$A \ =\  g_1\cdots g_r (\der g_r^{-1})  \cdots (\der g_2^{-1})( \der g_1^{-1}).$$
Moreover,  $\Sol\eqref{eq:disconj, b}$ contains a Markov system for all~$b>a$ iff \eqref{eq:disconj, b} is disconjugate  
for all $b>a$.
\end{cor}
\begin{proof}
The first equivalence follows from Proposition~\ref{prop:Frobenius} and Theorem~\ref{thm:Polya}, and the second equivalence follows
from Lemma~\ref{lem:Ma=>Ch} and Corollary~\ref{lem:disconj => Ma},
\end{proof}

\noindent
We say that~\eqref{eq:disconj} is {\bf eventually disconjugate} if \eqref{eq:disconj, b} is disconjugate for some $b\geq a$.
If \eqref{eq:disconj} is disconjugate, then so is \eqref{eq:disconj, b} for all~${b\geq a}$, and likewise with ``eventually disconjugate'' in place
of ``disconjugate''.\index{disconjugate!eventually}\index{linear differential operator!eventually disconjugate} If~\eqref{eq:disconj} is eventually disconjugate, then no solution of  \eqref{eq:disconj} in $\Car$ oscillates.
If $r=1$, then~\eqref{eq:disconj} is always disconjugate, since its solutions are the functions $t\mapsto c\exp\!\left(-\int_a^t f_1(s)\,ds\right)$ with $c\in\R$.   Returning to the special case   where $r=2$ we have:
%For $r=2$  we have:

\begin{cor}
Suppose \eqref{eq:2nd order} has a  non-oscillating solution $y\neq 0$.
Then   \eqref{eq:2nd order} is eventually disconjugate.
\end{cor}
\begin{proof} Here $f_1=0$, $f_2=f$, and $f$ does not generate oscillations by
 Lemma~\ref{lem:some vs all oscillate}. 
 Let $y_1,y_2\in\Sol(f)$ be non-oscillating and $\R$-linearly independent. Then $\wr(y_1,y_2)\in\R^\times$. 
Take $b\geq a$ such that $y_1|_b\in (\c_b)^\times$. Then~$y_1|_b,y_2|_b$ is a Markov system.
\end{proof}

\begin{remark}
By \cite{Gustafson},
%\marginpar{remark skipped}
there is for each $r>2$ a linear differential equation \eqref{eq:disconj} with only non-oscillating solutions
in $\Car$, but which is not eventually disconjugate. 
(This will not be used later but motivates Corollary~\ref{cor:Trench} below.) 
\end{remark}

\noindent 
Passing to germs instead of functions, we now consider a monic operator
$$A\ =\ \der^r+\phi_1\der^{r-1} + \cdots + \phi_r\in \Calinf[\der]\qquad (\phi_1,\dots, \phi_r\in \Calinf).$$
It gives rise to the $\R$-linear map 
$$y\mapsto A(y)=y^{(r)} + \phi_1y^{(r-1)} + \cdots + \phi_ry\  :\  \Calinf\to \Calinf,$$
whose kernel we denote by $\ker A$.

\begin{lemma}\label{funtogerm} $\dim_{\R}\ker A=r$,  and if $\theta_1,\dots, \theta_r\in \Calinf$ and $\ker A = \ker B$ for~$B=\der^r+ \theta_1\der^{r-1}+\cdots + \theta_r\in \Calinf[\der]$, then $A=B$, that is $\phi_i=\theta_i$ for $i=1,\dots,r$.
\end{lemma}
\begin{proof} Take $a\in \R$ and $f_1,\dots, f_r\in \c_a$ representing $\phi_1,\dots, \phi_r$.  This gives an equation \eqref{eq:disconj}. Let $y_1,\dots, y_r$ be a basis of the $\R$-linear space $\Sol\eqref{eq:disconj}$. Then the germs of $y_1,\dots, y_r$ lie in $\Calinf$, and denoting these germs also by $y_1,\dots, y_r$ one verifies easily that then $y_1,\dots, y_r$ is a basis of
$\ker A$. The second part of the lemma follows in a similar way from Corollary~\ref{cor:uniqueness of coeffs}. 
\end{proof} 

\noindent
We call $A$ as above
{\bf  disconjugate} if for some $a\in \R$ the germs $\phi_1,\dots,\phi_r$ have representatives $f_1,\dots,f_r$ in $\c_a$ such that
the linear differential equation \eqref{eq:disconj} 
on~$\R^{\geq a}$ is disconjugate.\index{disconjugate}\index{linear differential operator!disconjugate}

{\samepage
\begin{lemma}\label{lem:disconj germs}
For $A$ as above, the following are equivalent:
\begin{enumerate}
\item[\textup{(i)}] $A$ is disconjugate;
\item[\textup{(ii)}] $A=g_1\cdots g_r (\der g_r^{-1})  \cdots (\der g_2^{-1})( \der g_1^{-1})$ for some $g_1,\dots,g_r\in(\Calinf)^\times$;
\item[\textup{(iii)}] $A=(\der-h_r)\cdots (\der-h_1)$ for some $h_1,\dots,h_r\in\Calinf$.
\end{enumerate}
Thus if monic $A_1,A_2\in  \Calinf[\der]$  of order~$\geq 1$ are disconjugate, then so is~$A_1A_2$.
\end{lemma}}
\begin{proof}
Assume (i). Then Corollary~\ref{lem:disconj => Ma}
yields $a\in \R$, representatives $f_1,\dots,f_r\in \c_a$ of $\phi_1,\dots, \phi_r$, and a Markov system
 $y_1,\dots,y_r\in\Sol\eqref{eq:disconj}$. 
Let~$g_1,\dots,g_r$ be as in Proposition~\ref{prop:Frobenius}.
Then for $b\geq a$ with $f_1|_b,\dots,f_r|_b\in \c_a^n$ we have $g_j|_{b}\in (\c_b^{n+r-j+1})^\times$ for $j=1,\dots,r$.
So the germs of~$g_1,\dots,g_r$ are in $(\Calinf)^\times$, and denoting these germs also by
$g_1,\dots, g_r$ gives $A=g_1\cdots g_r (\der g_r^{-1}) \cdots (\der g_2^{-1})( \der g_1^{-1})$ by
Proposition~\ref{prop:Frobenius} and Lemma~\ref{funtogerm}. We have now shown (i)~$\Rightarrow$~(ii). For the converse,  we reverse the argument
using Theorem~\ref{thm:Polya}.  The equivalence~(ii)~$\Leftrightarrow$~(iii) is shown just like Lemma~\ref{lem:Polya fact}, using also that $f\mapsto f^\dagger\colon (\Calinf)^{\times} \to \Calinf$ is surjective. 
\end{proof}

\begin{remark}
Lemma~\ref{lem:disconj germs} goes through for monic $A\in \Ginf[\der]$ of order $r$, with $\Ginf$  in place of $\Calinf$ everywhere. Likewise for monic $A\in \Gom[\der]$ of order $r$, with $\Gom$ in place of $\Calinf$ everywhere.
\end{remark}

\noindent
A {\bf principal system} of solutions of
\eqref{eq:disconj} is a tuple $y_1,\dots,y_r$ in $\Sol\eqref{eq:disconj}$ such that
\begin{enumerate}
\item $y_1(t),\dots,y_r(t)>0$ eventually, and
\item $y_1\prec\cdots\prec y_r$ (in $\c$).
\end{enumerate}
Note that then $y_1,\dots,y_r$ are $\R$-linearly independent, and   $z_1,\dots,z_r\in\Car$ is a principal system of
solutions of \eqref{eq:disconj} iff there are $c_{ij}\in\R$ \textup{(}$1\leq j\leq i\leq r$\textup{)} such that
$$z_i\ =\ c_{ii}y_i+c_{i,i-1}y_{i-1}+\cdots+c_{i1}y_1\quad\text{and}\quad c_{ii}>0.$$
%A solution $y\in\Car$ of \eqref{eq:disconj} is {\bf principal} if there is a principal system of solutions~$y_1,\dots,y_r$
%of \eqref{eq:disconj} with $y_1=y$. Thus \eqref{eq:disconj} has at most one principal solution, up to multiplication by
%a positive constant. \index{solution!principal}\index{system!principal}\index{principal!solution}\index{principal!system}
%{above commented out lines checked, but he notion ``principal'' seems superfluous} 
The next result  generalizes Lemma~\ref{lem:admissible pair}. It seems slightly stronger than
a similar result by  Hartman~\cite{Hartman69} and Levin~\cite{Levin}:
 
%\begin{prop}
%If \eqref{eq:disconj} is eventually disconjugate, then it has a principal solution.\marginpar{skipped}
%\end{prop}

\begin{prop} Suppose \eqref{eq:disconj} has no oscillating solutions. Then it has a principal system of solutions.
\end{prop} 
\begin{proof} Let
$y, z\in \Sol\eqref{eq:disconj}$, $y(t), z(t)>0$ eventually. Claim:  $\lim\limits_{t\to +\infty} y(t)/z(t)$ exists in~$[0,+\infty]$. 
Suppose this limit doesn't exist. Then we have $c\in \R^{>}$ such that $$\liminf_{t\to +\infty} y(t)/z(t)\  < c\  < \limsup_{t\to+\infty} y(t)/z(t),$$ so $y(t)/z(t)=c$ for arbitrarily large $t$, but then $y-cz=0$, a contradiction. In particular, for such $y,z$ we have
either $y\prec z$, or $y\sim cz$ for some $c\in \R^{>}$, or~$y\succ z$. If~$y_1,\dots, y_n\in \Sol\eqref{eq:disconj}^{\ne}$ and 
$y_1\prec\cdots \prec y_n$, then $y_1,\dots, y_n$ are $\R$-linearly independent, so $n\le r$, and for any
nonzero $z\in \R y_1+\cdots + \R y_n$ we have~$z\sim cy_j$ for some $j\in \{1,\dots,n\}$ and $c\in \R^\times$. 
Now take such $y_1,\dots, y_n$ with maximal~$n$, so~$n\ge 1$. We claim that then 
$\Sol\eqref{eq:disconj}=\R y_1+\cdots + \R y_n$ (so $n=r$). Let $z\in  \Sol\eqref{eq:disconj}^{\neq}$. We cannot have $z\prec y_1$, nor $y_j\prec z \prec y_{j+1}$ with $1\le j\le n-1$, nor $z\succ y_n$;  hence~$z\sim cy_j$ where $1\le j\le n$ and
$c\in \R^\times$. Then $z-cy_j\prec y_j$. If $z\ne cy_j$, we take~$z-cy_j$ as our new $z$
and obtain likewise 
$z-cy_j\sim dy_i$ with~$1\le i<j$ and $d\in \R^\times$. Continuing this way leads in a finite number of steps to
 $z\in \R y_1+ \cdots + \R y_n$.
 \end{proof}

%\medskip\noindent
%Since we establish a version of this fact in the Hardy field context later, we do not give the proof here; see
%\cite[Chapter~3, Theorem~13]{Coppel}. 

\noindent
The next result is due to Trench~\cite{Trench}. We do not give a proof, since we shall establish in Section~\ref{sec:lin diff applications}  a version of it in the Hardy field context; see also Proposition~\ref{prop:distinguished splitting, Trench}. 

\begin{prop}\label{prop:Trench}
Suppose $\Sol\eqref{eq:disconj}$ contains a Markov system. Then there are 
% \marginpar{skipped}
$g_j\in  (\c_a^{r-j+1})^\times$   \textup{(}$j=1,\dots,r$\textup{)} such that for $A=\der^r+f_1\der^{r-1} + \cdots + f_r:\c_a^r\to \c_a$,
$$A\ =\  g_1\cdots g_r (\der g_r^{-1}) \cdots (\der g_2^{-1})( \der g_1^{-1})\quad\text{and}\quad 
\int_a^\infty \abs{g_j(s)}\,ds = \infty\text{ for $j=2,\dots,r$.}$$
Moreover, such $g_1,\dots,g_r$ are unique up to multiplication by nonzero constants. 
\end{prop}

\noindent
An application of l'H\^{o}pital's Rule shows that for $g_1,\dots,g_r$ as in Proposition~\ref{prop:Trench} the 
%\marginpar{skipped} 
tuple $y_1,\dots,y_r$  in the proof of Theorem~\ref{thm:Polya} is a principal system of solutions of \eqref{eq:disconj}.

\subsection*{Lyapunov exponents\astr}
In this subsection $f$, $g$, $h$ range over $\c[\imag]$. Consider the downward closed subset
$$\Lambda=\Lambda(f):=\big\{ \lambda\in \R: f\ex^{\lambda x}\preceq 1\big\}$$
of $\R$.
If $\lambda<\mu\in \Lambda$, then $f\ex^{\lambda x}\prec 1$. Also 
$$\Lambda(f)\ =\ \Lambda(\bar f)\ =\ \Lambda(\abs f),\quad f\preceq g\  \Rightarrow\  \Lambda(f)\supseteq \Lambda(g).$$

\begin{notation}Set $\R_{\pm\infty}:=\R\cup\{-\infty,+\infty\}$. Then for $S\subseteq\R$ we have $\sup S\in \R_{\pm\infty}$ with
%If $S\neq\emptyset$ is bounded from above,  then we have $\sup S\in\R$.
%If $S=\emptyset$ we define 
$\sup \emptyset:=-\infty$.
%if $S$ is not bounded from above we set $\sup S:=+\infty$.
%Thus in each case $\sup S\in\R_{\pm\infty}$
\end{notation} 

\noindent
The {\bf Lyapunov exponent}\/ \index{Lyapunov!exponent}  of $f$ is
$\lambda(f):=\sup \Lambda(f)\in\R_{\pm\infty}$.
(See \cite[\S{}3.12]{Cesari}.) Note:
$$\lambda(f)=+\infty\quad\Longleftrightarrow\quad \Lambda(f)=\R
\quad\Longleftrightarrow\quad f\prec \ex^{\lambda x}\text{ for all $\lambda\in\R$,}$$
and
$$\lambda(f)\ =\ \lambda(\bar f)\ =\ \lambda(\abs{f}),\quad f\preceq g\ \Rightarrow\ \lambda(f)\ \geq\ \lambda(g), \quad 
f\asymp g\ \Rightarrow\ \lambda(f)\ =\ \lambda(g).$$
If $\lambda=\lambda(f)\in\R$, then for each $\varepsilon\in\R^>$ we have $f\ex^{(\lambda-\varepsilon)x}\prec 1$
and~$f\ex^{(\lambda+\varepsilon) x}\not\preceq  1$.
One also verifies easily that for $f\in \c[\imag]^\times$,
\begin{equation}\label{eq:Lyap limsup}
\lambda(f)\  =\  - \limsup_{t\to+\infty} \frac{\log\abs{f(t)}}{t}.
\end{equation}
If $f=\ex^g$, then 
$\lambda(f)=-\limsup\limits_{t\to+\infty}   \Re g(t)/t$.
Thus $\lambda(c\ex^{\imag \phi})=0$ for~$c\in\C^\times$, $\phi\in\c$. 

%To avoid annoying case distinctions, we assume $\lambda(f),\lambda(g)>-\infty$
%in the next lemma.

{\samepage
\begin{lemma}\label{lem:Lyap exp} Assume $\lambda(f),\lambda(g)>-\infty$. Then:
%\mbox{}
\begin{enumerate}
\item[\textup{(i)}]
$\lambda(f+g)\geq \min\!\big\{ \lambda(f),\lambda(g) \big\}$,  with equality if $\lambda(f)\neq\lambda(g)$;  
\item[\textup{(ii)}]
$\lambda(fg)\geq\lambda(f)+\lambda(g)$;
\item[\textup{(iii)}]
$\lambda(f^m)=m\lambda(f)$ for all $m$.
\end{enumerate}
\end{lemma}}
\begin{proof}
For (i) suppose $\lambda(f)\leq\lambda(g)$. Then 
%$\Lambda(f)\subseteq \Lambda(g)$, 
 for each $\lambda\in \Lambda(f)$ and $\varepsilon\in \R^{>}$ we have~$(f+g)\ex^{(\lambda-\varepsilon) x}\preceq 1$ and so $\lambda-\varepsilon \in \Lambda(f+g)$. This shows $\lambda(f+g)\geq\lambda(f)$,  and~$\lambda(f+g)=\lambda(f)$
if~$\lambda(f)<\lambda(g)$ then follows using $f=(f+g)-g$. Parts (ii) and (iii) follow in a similar way.
\end{proof}

\noindent
By Lemma~\ref{lem:Lyap exp}(ii),
if $f\in\c[\imag]^\times$ and $\lambda(f),\lambda(f^{-1})\in\R$, then $\lambda(f^{-1})\leq-\lambda(f)$.

\begin{example}
If $f=\ex^g$ and $\Re g-cx\prec x$, $c\in\R$, then $\lambda(f)=c$, $\lambda(f^{-1})=-c$.
\end{example}

\noindent Set
$$\c[\imag]^{\flattereq}\ :=\ \big\{f:\lambda(f)>-\infty\big\}\ =\ \big\{f:\text{$f\preceq\ex^{nx}$ for some $n$}\big\}.$$
Then $\c[\imag]^{\flattereq}$ is a subalgebra of the $\C$-algebra $\c[\imag]$  and
$$\c[\imag]^{\flatter}:=\big\{f:\lambda(f)=+\infty\big\}=\big\{f:\text{$f\preceq\ex^{-nx}$ for all $n$}\big\}$$
is an ideal of $\c[\imag]^{\flattereq}$. The group of units of $\c[\imag]^{\flattereq}$ is
$$\c[\imag]^{\comp}:=\big\{ f\in\c[\imag]^\times:\lambda(f),\lambda(f^{-1})\in\R\big\}=\big\{f :
\text{$\ex^{-nx}\preceq f\preceq\ex^{nx}$  for some $n$}
\big\}.$$

%\begin{remark}  
%By Lemma~\ref{lem:Lyap exp}, $f\mapsto \ex^{-\lambda(f)}\colon\c[\imag]^{\flattereq}\to\R^{\geq}$ is a (non-archimedean) pseudovaluation on the \marginpar{remark not checked}
%ring $\c[\imag]^{\flattereq}$ in the sense of Mahler~\cite{Mahler}.
%\end{remark}

\begin{lemma}
Suppose $f\in\c^1[\imag]$. If $\lambda(f')\leq 0$, then $\lambda(f')\leq\lambda(f)$. If~$\lambda(f')>0$, then $c:=\lim\limits_{s\to \infty}f(s)\in\C$ exists and $\lambda(f')\leq\lambda(f-c)$.
\end{lemma}
\begin{proof}
Let $\lambda\in \Lambda(f')$. 
Take $a\in\R$ and a representative of $f$ in $\c^1_a[\imag]$,
also denoted by $f$, as well $C\in\R^>$, such that~$\abs{f'(t)}\leq C\ex^{-\lambda t}$ for $t\geq a$.
If $\lambda<0$, then for $t\geq a$:
\begin{multline*}
\abs{f(t)}-\abs{f(a)}\ \leq\  \abs{f(t)-f(a)}\ =\  \left|\int_a^t f'(s)\,ds\right|\  \leq\   
\int_a^t \abs{f'(s)}\,ds\  \leq \\  C\int_a^t \ex^{-\lambda s}\,ds\ =\  -\frac{C}{\lambda}(\ex^{-\lambda t}-\ex^{-\lambda a}),
\end{multline*}
hence $f\preceq \ex^{-\lambda x}$.
This yields $\lambda(f')\leq\lambda(f)$ if $\lambda(f')\leq 0$.
Suppose $\lambda> 0$.
Then for~$a\leq s\leq t$:
%$s\geq t\geq a$:
$$\abs{f(t)-f(s)}\  =\ \left|\int_s^t f'(u)\,du\right|\ \leq\
\int_s^t \abs{f'(u)}\,du\  \leq\  -\frac{C}{\lambda}(\ex^{-\lambda t}-\ex^{-\lambda s}).$$
Therefore $c:=\lim\limits_{s\to \infty}f(s)$ exists 
and $\abs{c-f(s)}\leq\frac{C}{\lambda}\ex^{-\lambda s}$ for $s\geq a$. Hence~$f-c\preceq \ex^{-\lambda x}$,
so $\lambda\in \Lambda(f-c)$. This yields $\lambda(f')\leq\lambda(f-c)$.
\end{proof}

\noindent
Let $y=(y_1,\dots,y_n)\in\c[\imag]^n$, $n\geq 1$. Put $\lambda(y):=\min\!\big\{\lambda(y_1),\dots,\lambda(y_n)\big\}$.
Then the function $\lambda\colon\c[\imag]^n\to\R_{\pm\infty}$ on the product ring $\c[\imag]^n$
also satisfies (i)--(iii) in Lemma~\ref{lem:Lyap exp} with $f$, $g$ replaced by  $y,z\in\c[\imag]^n$ with $\lambda(y),\lambda(z)>-\infty$. Thus:

\begin{cor}\label{cor:Lyap exp lin indep}
If $m\ge 2$, $y_1,\dots,y_m\in\c[\imag]^n$ are $\C$-linearly dependent, and $\lambda(y_1),\dots, \lambda(y_m)>-\infty$, then
$\lambda(y_i)=\lambda(y_j)$ for some $i\neq j$.
\end{cor}

\noindent
We define $y\preceq g:\Leftrightarrow y_1,\dots,y_n\preceq g\ \big({\Rightarrow \lambda(y)\geq\lambda(g)}\big)$.  
 Note that $\lambda(y) \in\R$ iff~$y\preceq\ex^{mx}$ 
and $y\not\preceq\ex^{-mx}$ for some~$m$. 

Let $\dabs{\,\cdot\,}$ be a norm on the $\C$-linear space~$\C^n$, and accordingly, let $\dabs{y}$ denote
the germ of $t\mapsto \dabs{\big(y_1(t),\dots, y_n(t)\big)}$, so $\dabs{y}\in \c$. 

\begin{cor}\label{cor:Lyap exp ReIm}
%Let $\dabs{\,\cdot\,}$ be a norm on the $\C$-linear space~$\C^n$. Then 
$y\preceq\dabs{y}$, 
$y\preceq g\Leftrightarrow\dabs{y}\preceq g$, and~$\lambda(\dabs{y})=\lambda(y)$.  
\end{cor}
\begin{proof}
Any two norms on $\C^n$ are equivalent, so we may arrange  $\dabs{\,\cdot\,}=\dabs{\,\cdot\,}_1$.
Then $y\preceq g \Rightarrow\dabs{y}= \abs{y_1}+\cdots+\abs{y_n}\preceq g$.
From $\abs{y_j}\leq\dabs{y}$ we get $y_j\preceq\dabs{y}$ for~$j=1,\dots,n$
and thus $y\preceq\dabs{y}$.
Thus~$\dabs{y}\preceq g\Rightarrow y\preceq g$; also $\lambda(\dabs{y})\leq\lambda(y)$.
Finally, 
  Lem\-ma~\ref{lem:Lyap exp}(i) and $\lambda(\abs{f})=\lambda(f)$ yield
 $\lambda(\dabs{y})\geq\lambda(y)$.
\end{proof}

\noindent
In particular, $\lambda(y)=\lambda(\dabs{y})=\lambda(\dabs{y}_2)=\lambda\big((\Re y_1,\dots, \Re y_n, \Im y_1,\dots, \Im y_n)\big)$.

\subsection*{Remarks on matrix differential equations\astr} 
In this subsection  $N$ is an $n\times n$ matrix with entries in $\c[\imag]$, $n\ge 1$.
We consider tuples
$y\in\c^1[\imag]^n$ as  column vectors~$y=(y_1,\dots, y_n)^{\text{t}}$ with entries $y_j$ in $\c^1[\imag]$.
Later in this subsection and in Section~\ref{sec:lin diff applications} we shall tacitly use the following:  
\begin{enumerate}
\item The $\C$-linear space of $y\in\c^1[\imag]^n$ such that $y'=Ny$ has dimension $n$.
\item If all entries of $N$ are in $\Calinf[\imag]$ and $y\in \c^1[\imag]^n$, $y'=Ny$, then $y\in \Calinf[\imag]^n$.
\end{enumerate} 
Classical existence and uniqueness results on matrix linear differential equations give (1), and  induction on the degree of smoothness of $y$ yields (2). 

\medskip
\noindent
Call a matrix $(f_{ij})$ over $\c[\imag]$  {\bf bounded}\/ if~$f_{ij}\preceq 1$ for all~$i$,~$j$.\index{bounded!matrix}
Similarly with~$\c_a[\imag]$ (${a\in\R}$) in place of $\c[\imag]$.
In the proof of Corollary~\ref{cor:bddflattereq} we shall use the following (cf.~\cite[\S{}3.13]{Cesari}, \cite[\S{}A.3.11]{Kamke}):

\begin{lemma}[{Lyapunov~\cite{Lyap}, Perron~\cite{Perron30}}]\label{lem:Lyap}
Suppose $N$ is bounded. If $y\in\c^1[\imag]^n$, $y'=Ny$, and $y\neq 0$, then $\lambda(y)\in\R$.
\end{lemma}
\begin{proof}
We have  $N=A+B\imag$ where $A$, $B$ are 
$n\times n$ matrices over $\c$. Consider the bounded $2n\times 2n$ matrix $M:=\left(\begin{smallmatrix} A & -B \\ B & A\end{smallmatrix}\right)$
over $\c$. For $y=(y_1,\dots,y_n)^{\text{t}}\in\c^1[\imag]^n$, set~$v:=(\Re y_1,\dots,\Re y_n,\Im y_1,\dots,\Im y_n)^{\text{t}}\in (\c^1)^{2n}$; then~$y'=Ny$ iff $v'=Mv$. Now assume $y'=Ny$ and $y\ne 0$. Then by the remark after Corollary~\ref{cor:Lyap exp ReIm} we may
replace $N$, $n$, $y$ by~$M$,~$2n$,~$v$ to arrange that the entries of~$N$ are in $\c$ and~$y\in (\c^1)^n$. 

Let $\lambda,\mu \in\R$ and consider $z:=\ex^{-\lambda x}y$. Then
$z'=(N-\lambda I_n)z$ where $I_n$ is the $n\times n$ identity matrix over $\c$, and thus
$$ \langle z,z\rangle'\  =\ 2\langle z,z'\rangle\  =\  2\big\langle z,\big(N-(\lambda-\textstyle\frac{1}{2}) I_n\big)z\big\rangle- \langle z,z\rangle.$$
The lemma below gives   $\lambda$ such that $\big\langle z,\big(N-(\lambda-\frac{1}{2}) I_n\big)z\big\rangle\le 0$, so $\langle z, z\rangle\in\c^1$   and~$\langle z,z\rangle'\leq 0$, and thus $\langle z,z\rangle   \preceq 1$. % by Lemma~\ref{grondiff}.   
Corollary~\ref{cor:Lyap exp ReIm}   yields $z\preceq 1$, so~$y\preceq\ex^{\lambda x}$.  
Likewise, set $w:=\ex^{\mu x}y$; then $w'=(N+\mu I_n)w$, and apply Lemma~\ref{fbliap} to a representative $F$ of
$-N$ to get $\mu$ with~$\langle w,w\rangle'\geq \langle w,w\rangle$, so~$\langle w,w\rangle\succeq \ex^x$ by Lemma~\ref{grondiff}, hence $w\not\preceq 1$, and thus~$y\not\preceq\ex^{-\mu x}$. So $\lambda(y)\in \R$. 
\end{proof}

\noindent
In the next lemma $F=(f_{ij})$  is an $n\times n$ matrix  over  $\c_a$, $a\in\R$. For $t\in\R^{\ge a}$ this yields the $n\times n$ matrix $F(t):=\big(f_{ij}(t)\big)$ over $\R$.
Let~$I_n$ also be the~$n\times n$ identity matrix over $\R$.

\begin{lemma}\label{fbliap}
Suppose   $F$ is bounded.
Then there exists~$\mu\in\R^>$ such that for all real $\lambda\geq \mu$, $t\ge a$, and   $z\in\R^n$:
$\big\langle z,\big(F(t)-\lambda I_n\big)z\big\rangle \leq 0$.
\end{lemma}
\begin{proof}
Put $G:=\frac{1}{2}(F+F^{\operatorname{t}})$, a symmetric bounded $n\times n$ matrix over $\c_a$ such that~$\big\langle z,F(t)z\big\rangle=\big\langle z,G(t)z\big\rangle$ for $t\ge a$ and $z\in\R^n$, and
replace $F$ by $G$ to arrange that $F$ is symmetric.
Let 
$$P(Y)\ :=\ \det(YI_n-F)\ =\ Y^n+P_1Y^{n-1}+\cdots+P_n\in\c_a[Y]\qquad (P_1,\dots,P_n\in\c_a),$$ 
and for $t\in\R^{\ge a}$ put $$P(t,Y)\ :=\  Y^n+P_1(t)Y^{n-1}+\cdots+P_n(t)\in\R[Y],  $$ so 
for each $\lambda\in\R$,
$P(t,Y+\lambda)$ is the
characteristic polynomial of the symmetric $n\times n$ matrix $F(t)-\lambda I_n$ over~$\R$. Now~$P_1,\dots,P_n\preceq 1$ since $F$ is bounded, 
so~[ADH, 3.5.2]   yields~$\mu\in\R^>$ such that for all $t\in\R^{\ge a}$, all zeros of $P(t,Y)$ in $\R$ are in~$[-\mu,\mu]$.
Let $\lambda\geq \mu$. Then for $t\ge a$, no real zero of $P(t,Y+\lambda)$, and thus no eigenvalue of~$F(t)-\lambda I_n$, is positive. Hence $\big\langle z,(F(t)-\lambda I_n)z\big\rangle   \leq 0$ for all~$z\in\R^n$.
\end{proof}

\noindent
Let $V:=\big\{ y\in \c^1[\imag]^n: y'=Ny \big\}$, an $n$-dimensional $\C$-linear subspace of $\c^1[\imag]$.
Suppose $N$ is bounded. Then $S:=\lambda(V^{\neq}) \subseteq\R$ by Lemma~\ref{lem:Lyap}, and $S$,  called the {\bf Lyapunov spectrum}\/ of~$y'=Ny$, has at most $n$ elements by Corollary~\ref{cor:Lyap exp lin indep}.\index{Lyapunov!spectrum}\index{matrix differential equation!Lyapunov spectrum}\index{spectrum!Lyapunov} According to
 [ADH, 2.3] the surjective map
$$y\mapsto\lambda(y)\colon V\to S_\infty:=S\cup\{\infty\}$$
makes $V$ a valued vector space over $\C$. Thus by [ADH, remark before 2.3.10]:

\begin{cor}[{Lyapunov~\cite{Lyap}}]
If $N$ is bounded, then $V$ has a basis $y_1,\dots,y_n$  such that for all~$c_1,\dots,c_n\in\C$, not all zero, and $y=c_1y_1+\cdots+c_ny_n$,
$$\lambda(y)\ =\ \min\!\big\{\lambda(y_j):c_j\neq 0\big\}.$$
\end{cor}

\noindent
Whether or not $N$ is bounded, a {\bf Lyapunov
fundamental system of solutions} of $y'=Ny$ is a basis  $y_1,\dots,y_n$ of $V$ as in the  corollary above.\index{matrix differential equation!Lyapunov fundamental system of solutions}\index{Lyapunov!fundamental system of solutions}
(In \cite[\S{}3.14]{Cesari}  this  is called  a {\it normal}\/ fundamental system of solutions of $y'=Ny$.)
A {\bf Lyapunov fundamental matrix} for $y'=Ny$ is  an $n\times n$ matrix with entries in $\c^1[\imag]$ whose columns form a Lyapunov fundamental system of solutions of $y'=Ny$.\index{matrix differential equation!Lyapunov fundamental matrix}\index{Lyapunov!fundamental matrix}

\medskip
\noindent
Lemma~\ref{lem:Lyap} also gives:

\begin{cor}\label{cor:Lyap}
Let $f_1,\dots,f_n\in\c[\imag]$ be such that  $f_1,\dots,f_n\preceq 1$. 
Then there are~$\lambda_1,\dots,\lambda_m\in\R$ \textup{(}$1\leq m\leq n$\textup{)} such that
for all $y\in\c^n[\imag]^{\neq}$ such that
$$y^{(n)}+f_1y^{(n-1)}+\cdots+f_ny\ =\ 0$$
we have $\lambda(y,y',\dots,y^{(n-1)})\in\{\lambda_1,\dots,\lambda_m\}$. 
%$y^{(j)}\preceq\ex^{m x}$ for $j=0,\dots,n-1$ but not $y^{(j)} \preceq \ex^{-m x}$ for all $j=0,\dots,n-1$.
\end{cor}

\section{Hardy Fields}\label{sec:Hardy fields}

\noindent
Here we introduce Hardy fields and review some classical extension theorems for Hardy fields. 

%The Product Rule for derivatives yields a fact used later:

%\begin{lemma}\label{lem:product fm asymptotics}
%If $n\ge 1$, $f\in\Caln[\imag]^\times$, $f^{(i)}\prec f$ for $i=1,\dots,n$, and $g\in\Caln[\imag]$,
%$g^{(n)}\in\c[\imag]^\times$, $g^{(i)}\preceq g^{(n)}$ for $i=0,\dots,n$, 
%then $(fg)^{(n)} \sim fg^{(n)}$.
%\end{lemma}
%\begin{proof}
%We have
%\begin{multline*}
%(fg)^{(n)}\Big/\big(fg^{(n)}\big) = 1+ \sum_{i=1}^n {n\choose i} \big(f^{(i)}/f\big) \big( g^{(n-i)}\big/ g^{(n)} \big) \\ \text{where $f^{(i)}/f\prec 1$ and $g^{(n-i)}/ g^{(n)}\preceq 1$ for $i=1,\dots,n$,}
%\end{multline*}
%and so $\lim\limits_{t\to\infty} (fg)^{(n)}(t)\big/(fg^{(n)}) (t) = 1$.
%\end{proof}

\subsection*{Hardy fields}
A {\em Hardy field\/} is a subfield of 
$\Gi$ that is closed under the derivation of $\Gi$; see also~[ADH, 9.1].\index{Hardy field} Let~$H$ be a Hardy field.
Then $H$ is considered as an ordered valued differential field in the obvious way; see Section~\ref{sec:germs} for the ordering and valuation on $H$.
The field of constants of $H$ is $\R\cap H$. Hardy fields are pre-$H$-fields, and $H$-fields if they contain~$\R$;
see [ADH, 9.1.9(i),~(iii)].
As in Section~\ref{sec:germs} we equip the differential subfield $H[\imag]$  of~$\Calinf[\imag]$ with the unique valuation ring 
lying over that of~$H$. Then $H[\imag]$ is a pre-$\d$-valued field of $H$-type with small derivation and constant field $\C\cap H[\imag]$;
if $H\supseteq\R$, then $H[\imag]$ is $\d$-valued   with constant field $\C$.

\medskip
\noindent
We also consider variants:  a {\em $\Ginf$-Hardy field\/} is a Hardy field $H\subseteq \Ginf$, and a {\em $\Gom$-Hardy field\/} (also called an {\em analytic Hardy field\/}) is a Hardy field~${H\subseteq \Gom}$.\index{Cr-Hardy field@$\c^r$-Hardy field}\index{Hardy field!analytic}\index{Hardy field!smooth} 
Most Hardy fields arising in practice are actually $\Gom$-Hardy fields.  Boshernitzan~\cite{Boshernitzan81} (with details worked out in \cite{Gokhman}) first suggested a Hardy field~${H\not\subseteq\Ginf}$, and~\cite[Theorem~1]{Grelowski} shows that each Hardy field with a largest
comparability class extends to a Hardy field $H\not\subseteq\Ginf$. 
%\marginpar{facts mentioned in this paragraph not checked, nor used}
Rolin, Speissegger, Wilkie~\cite{RSW} construct o-minimal expansions $\tilde{\R}$ of the ordered field of real numbers such that $H\subseteq\Ginf$ and $H\not\subseteq\Gom$ for the Hardy field $H$ consisting of the germs of functions 
$\R \to\R$ that are definable in $\tilde{\R}$. Le Gal and Rolin~\cite{LeGalRolin} construct such expansions such that  
$H\not\subseteq\Ginf$ for the corresponding Hardy field $H$.  

\subsection*{Hardian germs} Let $y\in\mathcal G$. 
Following~\cite{Sj}  we call $y$  {\bf hardian} if it lies in a Hardy field (and thus $y\in\Calinf$).
We also say that $y$ is {\bf $\Ginf$-hardian} if $y$ lies in a $\Ginf$-Hardy field,
equivalently, $y\in\Ginf$ and $y$ is hardian; 
likewise with $\Gom$ in place of~$\Ginf$.\index{germ!hardian}\index{hardian}
Let $H$ be a Hardy field. Call $y\in\mathcal G$ {\bf $H$-har\-dian} (or {\bf hardian over $H$}) if $y$ lies in a Hardy field extension
of~$H$. (Thus $y$ is hardian iff~$y$ is $\Q$-hardian.) 
If $H$ is a $\Ginf$-Hardy field and $y\in\Ginf$ is hardian over $H$, then $y$ generates a $\Ginf$-Hardy field extension~$H\langle y\rangle$ of $H$; likewise with $\Gom$ in place of $\Ginf$.

\subsection*{Maximal and perfect Hardy fields}
Let $H$ be a Hardy field. Call~$H$  {\it maximal}\/ if no Hardy field properly contains $H$.\index{Hardy field!maximal}\label{p:E(H)}
Following Boshernitzan~\cite{Boshernitzan82} we denote by $\Ex(H)$ the intersection of all maximal Hardy fields containing $H$;
thus~$\Ex(H)$ is a Hardy field extension of $H$, and a maximal Hardy field contains $H$ iff it contains~$\Ex(H)$, so $\Ex(\Ex(H))=\Ex(H)$. 
If $H^*$ is a Hardy field extension of~$H$, then $\Ex(H)\subseteq\Ex(H^*)$; hence if
$H^*$ is a Hardy field with $H\subseteq H^*\subseteq\Ex(H)$, then~$\Ex(H^*)=\Ex(H)$. 
Note 
that $\Ex(H)$ consists of the $f\in \mathcal G$ that are hardian over each Hardy field  $E\supseteq H$. Hence $\Ex(\Q)$ consists of the germs in $\mathcal G$  that are  hardian over each Hardy field. 
As in \cite{Boshernitzan82}  we also say that~$H$ is {\bf perfect} if~$\Ex(H)=H$.\index{Hardy field!perfect} (This terminology is slightly unfortunate, since Hardy fields,
being of characteristic zero, are perfect as  fields.)
Thus $\Ex(H)$ is the smallest perfect Hardy field extension of $H$. Maximal Hardy fields are perfect.

\subsection*{Differentially maximal Hardy fields} 
Let $H$ be a Hardy field.
We now define differentially-algebraic variants of the above:
call $H$   {\bf differentially maximal}, or {\bf $\d$-maximal}\/ for short, if $H$ has no proper $\d$-algebraic Hardy field extension.\index{Hardy field!d-maximal@$\d$-maximal}\index{Hardy field!differentially maximal} 
  Every maximal Hardy field is $\d$-maximal, so each Hardy field is contained in a $\d$-maximal one; in fact,
by Zorn, each Hardy field~$H$ has a $\d$-maximal Hardy field extension which is $\d$-algebraic over $H$. Let~$\Dx(H)$ be the intersection of all 
$\d$-maximal Hardy fields containing $H$. Then~$\Dx(H)$ is a $\d$-algebraic Hardy field extension of~$H$ with~$\Dx(H)\subseteq\Ex(H)$.
By the next lemma, $\Dx(H)=\Ex(H)$ iff $\Ex(H)$ is $\d$-algebraic over~$H$:\label{p:D(H)}

\begin{lemma}\label{lem:Dx Ex}
$\Dx(H)=\big\{f\in\Ex(H):\text{$f$ is $\d$-algebraic over $H$}\big\}$.
\end{lemma}
\begin{proof}
We only need to show the inclusion ``$\supseteq$''. For this let~$f\in\Ex(H)$ be $\d$-algebraic over $H$, and let $E$ be a $\d$-maximal Hardy field extension of $H$; we need to show $f\in E$. To see this extend $E$
 to a maximal Hardy field $M$; then $f\in M$, hence $f$ generates a Hardy field extension $E\langle f\rangle$ of $E$.
 Since $f$ is $\d$-algebraic over~$H$ and thus over~$E$, this yields $f\in E$ by $\d$-maximality of $E$, as required.
\end{proof}

\noindent
A $\d$-maximal Hardy field contains $H$ iff it contains~$\Dx(H)$, hence $\Dx(\Dx(H))=\Dx(H)$.
If~$H^*$ is a Hardy field extension of~$H$, then $\Dx(H)\subseteq \Dx(H^*)$; hence for each Hardy field~$H^*$  with $H\subseteq H^*\subseteq\Dx(H)$ we have~$\Dx(H^*)=\Dx(H)$.
We say that $H$ is {\bf $\d$-perfect}\index{Hardy field!d-perfect@$\d$-perfect} if $\Dx(H)=H$.  
Thus~$\Dx(H)$ is the smallest $\d$-perfect Hardy field extension of $H$. Every perfect Hardy field is $\d$-perfect, as is every $\d$-maximal Hardy field. The following diagram summarizes the various implications among these properties of Hardy fields:
$$\xymatrix{ \text{maximal} \ar@{=>}[r]\ar@{=>}[d]  &  \text{perfect} \ar@{=>}[d] \\ 
\text{$\d$-maximal} \ar@{=>}[r]  &  \text{$\d$-perfect}  }$$
We call $\Dx(H)$ the {\bf $\d$-perfect hull} of $H$, and $\Ex(H)$ the {\bf perfect hull} of $H$.\index{Hardy field!d-perfect hull@$\d$-perfect hull}\index{Hardy field!perfect hull} It seems that the following question  asked by Boshernitzan~\cite[p.~144]{Boshernitzan82} is still open:

\begin{question}
Is $\Ex(H)$ $\d$-algebraic over $H$, in other words, is $\Dx(H)=\Ex(H)$?
\end{question}

\noindent
Boshernitzan gave support for a positive answer:
Lemma~\ref{lem:bounded Hardy field ext, 1}, Corollary~\ref{cor:Bosh13.10}, and
Theorem~\ref{thm:Bosh 14.4} below. Our
Theorems~\ref{thm:coinitial in E(H)} and~\ref{thm:upo-freeness of the perfect hull} (in combination
with Theorem~\ref{thm:ADH 13.6.1}) can be seen as further support.

\subsection*{Variants of the perfect hull} 
Let  $H$ be a $\c^r$-Hardy field where   $r\in\{\infty,\omega\}$. We say that  $H$ is {\bf $\c^r$-maximal} if no $\c^r$-Hardy field properly contains it. By Zorn, $H$ has a $\c^r$-maximal extension. In analogy with $\Ex(H)$,  define the {\bf $\c^r$-perfect hull~$\Ex^r(H)$} of $H$ to be the
intersection of all $\c^r$-maximal Hardy fields containing~$H$. \label{p:Er(H)} We say that  $H$ is
  {\bf $\c^r$-perfect} if $\Ex^r(H)=H$. The penultimate subsection goes through
  with {\em Hardy field}, {\em maximal}, {\em hardian}, {\em $\Ex(\,\cdot\,)$}, and {\em perfect} replaced by {\em $\c^r$-Hardy field},
  {\em $\c^r$-maximal},  {\em $\c^r$-hardian}, {\em $\Ex^r(\,\cdot\,)$}, and {\em $\c^r$-perfect}, respectively.    (Corollary~\ref{cor:D(H) smooth} shows that no analogue of $\Dx(H)$ is needed
  for the $\c^r$-category.)\index{Cr-Hardy field@$\c^r$-Hardy field!Cr-maximal@$\c^r$-maximal}\index{Cr-Hardy field@$\c^r$-Hardy field!Cr-perfect@$\c^r$-perfect}\index{Cr-Hardy field@$\c^r$-Hardy field!Cr-perfect@$\c^r$-perfect hull}

\subsection*{Some basic extension theorems} 
We summarize some well-known extension results for Hardy fields:

\begin{prop}\label{prop:Hardy field exts} Any Hardy field $H$ has the following Hardy field extensions: \begin{enumerate}
\item[\textup{(i)}] $H(\R)$, the subfield of $\Gi$ generated by $H$ and $\R$;
\item[\textup{(ii)}] $H^{\operatorname{rc}}$, the real closure of $H$ as defined in Proposition~\ref{b1};
\item[\textup{(iii)}] $H(\ex^f)$ for any $f\in H$;
\item[\textup{(iv)}] $H(f)$ for any $f\in \Go$ with $f'\in H$;
\item[\textup{(v)}] $H(\log f)$ for any $f\in H^{>}$.
\end{enumerate}
If $H$ is contained in $\Ginf$, then
so are the Hardy fields in {\rm (i), (ii), (iii), (iv), (v)}; likewise with $\Gom$ instead of $\Ginf$. 
\end{prop}

\noindent
Note that (v) is a special case of (iv), since $(\log f)'=f^\dagger\in H$
for $f\in H^{>}$. Another special case of (iv) is that $H(x)$ is a Hardy field. A consequence of the proposition is that any Hardy field $H$ has a smallest real closed Hardy field extension~$L$ with~$\R\subseteq L$ such that for all $f\in L$ we have $\ex^f\in L$
and $g'=f$ for some~$g\in L$. Note that then~$L$ is a Liouville closed $H$-field as defined in~[ADH, 10.6]. 
% and~$H^*$ is $\d$-algebraic over $H$. \marginpar{added ``$\d$-algebraic'' part} 
Let $H$ be a Hardy field with~$H\supseteq\R$.  As in~\cite{AvdD2} and [ADH, p.~460] we  then  denote the above
$L$ by  $\Li(H)$; so $\Li(H)$ is the smallest   Liouville closed Hardy field containing $H$,
called the {\it Hardy-Liouville closure}\/ of $H$ in \cite{ADH2}. 
We have  $\Li(H)\subseteq \Dx(H)$, hence if $H$ is $\d$-perfect, then $H$ is a Liouville closed $H$-field. Moreover, if $H\subseteq\Ginf$ then $\Li(H)\subseteq\Ginf$, and
similarly with $\Gom$ in place of~$\Ginf$.\index{Hardy field!Hardy-Liouville closure}\label{p:HL}

\medskip
\noindent
The next more general result in Rosenlicht~\cite{Ros} is attributed there to M. Singer: 

\begin{prop}\label{singer} Let $H$ be a Hardy field and $p(Y),q(Y)\in H[Y]$, $y\in \mathcal{C}^1$, such that 
%\marginpar{still needs check}
$y'q(y)=p(y)$ with $q(y)\in\mathcal{C}^\times$. Then~$y$ generates a Hardy field~$H(y)$ over~$H$.
\end{prop} 

\noindent
Note that for $H$, $p$, $q$, $y$ as in the proposition we have~$y\in\Dx(H)$. 

%Thus: 
%\begin{cor}\label{cor:Singer-Rosenlicht}
%Let $H$ be a  $\d$-perfect Hardy field. Then $H$ contains every $y\in\Go$ such that $y'q(y)=p(y)$ with
%$p(Y), q(Y) \in H[Y]$ and $q(y)\in\c^\times$. 
%%Hence $H\supseteq\R$ and~$H$ is a Liouville closed $H$-field.
%\end{cor}

\subsection*{Compositional  conjugation of differentiable germs}  
Let~$\ell\in \Go$, $\ell'(t)>0$ eventually (so $\ell$ is eventually strictly increasing) and $\ell(t)\to+\infty$ as $t\to+\infty$. Then~$\phi:=\ell'\in\c^\times$,  and the compositional inverse 
$\ell^{\inv}\in \Go$ of $\ell$ satisfies 
$$\ell^{\inv}>\R, \qquad (\ell^{\inv})'\ =\ (1/\phi)\circ \ell^{\inv}\in \c.$$
The $\C$-algebra automorphism\label{p:fcirc}
$f\mapsto f^\circ:= f\circ \ell^{\inv}$ of $\mathcal{C}[\imag]$ (with inverse $g\mapsto g\circ\ell$) maps $\Go[\imag]$ onto itself and  
satisfies for $f\in\Go[\imag]$ a useful identity:
$$(f^\circ)'\ =\ (f\circ \ell^{\inv})'\ =\ (f'\circ \ell^{\inv})\cdot (\ell^{\inv})'\ =\ (f'/\ell')\circ \ell^{\inv}
\ =\ (\phi^{-1}f')^\circ.$$
Hence if $n\ge 1$ and $\ell\in\c^n$, then $\ell^{\inv}\in \c^n$ and $f\mapsto f^\circ$ maps~$\Caln[\imag]$ and $\Caln$
onto themselves, for each $n$. Therefore, if $\ell\in\Calinf$, then $\ell^{\inv}\in \Calinf$ and $f\mapsto f^\circ$   maps~$\Calinf[\imag]$ and $\Calinf$ onto themselves;
likewise with $\Ginf$ or $\Gom$ in place of $\Calinf$. 
{\it In the rest of this subsection we assume~$\ell\in\Calinf$.}\/ Denote the differential ring~$\Calinf[\imag]$ by~$R$, and as usual let $R^\phi$ be  $R$  with its derivation~$f\mapsto \der(f)=f'$ replaced by the  derivation~$f\mapsto \derdelta(f)=\phi^{-1}f'$~[ADH, 5.7]. Then  
$f\mapsto f^\circ \colon R^\phi \to R$
is an isomorphism of differential rings by the identity above. 
We extend it to the isomorphism
$$Q\mapsto Q^\circ\ \colon\ R^\phi\{Y\} \to R\{Y\}$$
of differential rings given by $Y^\circ=Y$.  
 Let $y\in R$. Then 
$$Q(y)^\circ \  = \  Q^\circ(y^\circ)\qquad\text{for $Q\in R^\phi\{Y\}$}$$
and thus
$$P(y)^\circ\ =\ P^\phi(y)^\circ\ =\ (P^\phi)^\circ(y^\circ) \qquad\text{for $P\in R\{Y\}$.}$$
This leads to a useful generalization of the   identity for $(f^\circ)'$ above. For this, let~$n\ge 1$ and let~$G^n_k\in\Q\{X\}$ ($1\leq k\leq n$) 
be the differential polynomial introduced in~[ADH, 5.7]; so~$G^n_k$ is homogeneous of degree $n$
and isobaric of weight $n-k$. Viewing the~$G^n_k$ as elements of~$R\{X\}$ and $\derdelta=\phi^{-1}\der$ as an element of
$R[\der]$ we have 
$$\derdelta^n\ =\ G^n_n(\phi^{-1})\der^n+\cdots+G^n_1(\phi^{-1})\der\qquad\text{in the ring $R[\der]$.}$$
Thus
$$\derdelta^2\ =\ \phi^{-2}\der^2-\phi'\phi^{-3}\der, \qquad \derdelta^3\ =\ \phi^{-3}\der^3-3\phi'\phi^{-4}\der^2 + \big(3(\phi')^2-\phi\phi''\big)\phi^{-5} \der , \quad \dots$$
Set $\lambda:=-\phi^\dagger$, and 
let $$R^n_k\ :=\ \operatorname{Ri}(G^n_k)\in\Q\{Z\},\quad\text{ so }  \quad
G^n_k(\phi^{-1})=\phi^{-n}R^n_k(\lambda)
\qquad (0\leq k\leq n).$$ 
Thus
$$\derdelta^n\  =\  \phi^{-n}\big(R^n_n(\lambda)\der^n+\cdots+R^n_1(\lambda)\der\big).$$
For instance,
\begin{align*}
\derdelta^3\	&=\	\phi^{-3}\big(R^3_3(\lambda)\der^3 + R^3_2(\lambda)\der^2 + R^3_1(\lambda)\der\big) \\
			&=\	\phi^{-3}\big(\der^3+3\lambda\der^2+\big(2\lambda^2+\lambda'\big)\der\big).
\end{align*}
We now have:
 
\begin{lemma}\label{lem:fcircn}  
Let $f\in R$ and $n\geq 1$. Then  
$$(f^\circ)^{(n)}\ =\  \left(\phi^{-n}\big(R^n_n(\lambda)f^{(n)}+\cdots+R^n_1(\lambda)f'\big)\right)^\circ.$$
\end{lemma}
\begin{proof}
Let $Q=Y^{(n)}\in R^\phi\{Y\}$, so $Q^\circ=Y^{(n)}\in R\{Y\}$. Then $(f^\circ)^{(n)}=Q^\circ(f^\circ)=Q(f)^\circ= \derdelta^n(f)^\circ$.
Now use the above identity for $\derdelta^n$. 
%the claim follows from the identities following [ADH, 5.7.2] applied to~$R^\phi$, $\phi^{-1}$ in place of~$K$,~$\phi$.
\end{proof}

\noindent
Note also: $(Q_{+f})^\circ=(Q^\circ)_{+f^\circ}$ and $(Q_{\times f})^\circ=(Q^\circ)_{\times f^\circ}$ for $Q\in R^\phi\{Y\}$, $f\in R$.

\subsection*{Compositional  conjugation in Hardy fields} Let now $H$ be a Hardy field, and
let~$\ell\in \Go$ be such that $\ell>\R$ and $\ell'\in H$. Then $\ell\in \Calinf$,  $\phi:=\ell'$
is active in $H$, $\phi>0$, and we have a Hardy field $H(\ell)$. 
%and the compositional inverse 
%$\ell^{\inv}\in \Go$ of $\ell$ satisfies 
%$$\ell^{\inv}>\R, \qquad (\ell^{\inv})'\ =\ (1/\phi)\circ \ell^{\inv}\in H\circ \ell^{\inv}.$$
 The
$\C$-algebra automorphism~$f\mapsto f^\circ:= f\circ \ell^{\inv}$ of~$\mathcal{C}[\imag]$  restricts to an ordered field isomorphism
$$h\mapsto h^\circ\ :\ H \to H^\circ:=H\circ \ell^{\inv}.$$
The identity  
$(f^\circ)'  = (\phi^{-1}f')^\circ$, valid for each $f\in\Go[\imag]$, shows
that  $H^\circ$ is again a Hardy field. Conversely, if $E$ is a subfield of $\Calinf$ with $\phi\in E$ and $E^\circ:= E\circ \ell^{\inv}$ is a Hardy field, then $E$ is a Hardy field. If $H\subseteq \Ginf$ and $\ell\in \Ginf$, then
$H^\circ\subseteq \Ginf$; likewise with $\Gom$ instead of~$\Ginf$. 
If $E$ is a Hardy field extension of $H$, then
 $E^\circ$ is a Hardy field extension of $H^\circ$, and
$E$ is $\d$-algebraic over~$H$ iff $E^\circ$ is $\d$-algebraic over $H^\circ$.
Hence $H$ is maximal iff $H^\circ$ is maximal, and likewise with ``$\d$-maximal'' in place of ``maximal''.
So~$\Ex(H^\circ)=\Ex(H)^\circ$ and~$\Dx(H^\circ)=\Dx(H)^\circ$, and thus $H$ is perfect iff $H^\circ$ is perfect, and likewise with ``$\d$-perfect'' in place of ``perfect''. 
The next lemma is
 \cite[Corollary~6.5]{Boshernitzan81}; see also \cite[Theo\-rem~1.7]{AvdD4}.
 
\begin{lemma}\label{lem:Bosh6.5}
The germ $\ell^{\inv}$ is hardian. Moreover,
if $\ell$ is $\Ginf$-hardian, then $\ell^{\inv}$ is also $\Ginf$-hardian, and
likewise with  $\Gom$  in place of $\Ginf$. 
\end{lemma}
\begin{proof}
By Proposition~\ref{prop:Hardy field exts}(iv) we can arrange that our Hardy field $H$ contains both $\ell$ and~$x$.
Then $\ell^{\inv}=x\circ \ell^{\inv}$ is an element of the Hardy field $H\circ\ell^{\operatorname{inv}}$.
\end{proof}

\noindent
Next we consider the pre-$\d$-valued field $K:=H[\imag]$ of $H$-type, which gives rise to
$$K^\circ:=K\circ \ell^{\inv}=H^\circ[\imag],$$ 
also a pre-$\d$-valued field of $H$-type, and we have the valued field isomorphism
\[h\mapsto h^\circ\ :\ K\to K^\circ.\]
Note: $h\mapsto h^\circ\colon H^{\phi} \to H^\circ$ is an isomorphism 
of pre-$H$-fields, and~${h\mapsto h^\circ\colon K^{\phi} \to K^\circ}$ is an isomorphism of valued differential fields.
Recall that $K$ and $K^\phi$ have the same underlying field.
For~$f,g\in K$ we have
$$f\preceq^\flat_\phi g \text{\ (in $K$)} \quad\Longleftrightarrow\quad 
f\preceq^\flat g \text{\ (in $K^\phi$)} \quad\Longleftrightarrow\quad
f^\circ\preceq^\flat g^\circ\text{\ (in $K^\circ$)},$$
and likewise with $\preceq^\flat_\phi$, $\preceq^\flat$ replaced by $\prec^\flat_\phi$, $\prec^\flat$. 

\begin{lemma}\label{lemlioucomp} From the isomorphisms $H^\phi\cong H^\circ$ and $K^\phi\cong K^\circ$ we obtain:
If $H$ is Liouville closed, then so is $H^\circ$. If $\I(K)\subseteq K^\dagger$, then $\I(K^\circ)\subseteq (K^\circ)^\dagger$.
\end{lemma}

\noindent
So far we focused on pre-composition with  $\ell^{\inv}$. As to pre-composition with $\ell$,
it seems not to be known whether $H\circ\ell\subseteq H$ whenever $H$ is maximal. However, we have the following (cf.~\cite[Lemma~11.6(7)]{Boshernitzan82}):

\begin{lemma}\label{lem:comp with E(Q)}
$\Ex(\Q)\circ\ell\subseteq \Ex(H)$. 
\end{lemma}
 \begin{proof} $\Ex(H^\circ)=\Ex(H)^\circ$ gives
$\Ex(H^\circ) \circ \ell = \Ex(H)$.  Now use~$\Ex(\Q)\subseteq \Ex(H^\circ)$.
 \end{proof}
 
\noindent
Lemma~\ref{lem:comp with E(Q)} gives
$\Ex(\Q)\circ \Ex(\Q)^{>\R}\subseteq\Ex(\Q)$; cf.~\cite[Theo\-rem~6.8]{Boshernitzan81}.
Boshernitzan's conjecture~\cite[\S{}10, Conjecture~3]{Boshernitzan81} that $\Ex(\Q)^{>\R}$ is also closed under compositional inversion seems to be still open. 

\subsection*{Differential algebraicity of compositional inverses\astr}
In the next lemma we let $\ell\in\Calinf$ be hardian with $\ell>\R$.
The argument in the proof of Lemma~\ref{lem:Bosh6.5} shows that  $\ell$ and
$\ell^{\operatorname{inv}}$ are both $\R(x)$-hardian; moreover (cf.~\cite[Lem\-ma~14.10]{Boshernitzan82}):

\begin{lemma}\label{lem:trdegellinv}
We have
\begin{equation}
\operatorname{trdeg}\!\big(\R\langle x, \ell^{\operatorname{inv}}\rangle|\R\big)\  =\ 
\operatorname{trdeg}\!\big(\R\langle x, \ell\rangle|\R\big), \label{eq:trdegellinv}
\end{equation}
hence if $\ell$ is $\d$-algebraic over $\R$, then so is 
$\ell^{\operatorname{inv}}$,  with
$$\operatorname{trdeg}\!\big(\R\langle \ell^{\operatorname{inv}}\rangle|\R\big)\  \leq\
\operatorname{trdeg}\!\big(\R\langle \ell\rangle|\R\big)+1.$$
\end{lemma}
\begin{proof}
Set $H:=\R\langle x, \ell\rangle=\R(x)\langle \ell\rangle$ and $\phi:=\ell'$.
With $\der$ and $\derdelta=\phi^{-1}\der$ denoting the derivations of~$H$ and $H^\phi$,  we have $\phi=1/\derdelta(x)$ and for all~$f\in H$ and $n\geq 1$,
$$\der^n(f)\in \Q\big[\derdelta(f),\derdelta^2(f),\dots,\phi,\derdelta(\phi),\derdelta^2(\phi),\dots\big]$$ 
by [ADH, remarks before~5.7.3].  The differential fields $H$ and $H^\phi$ have the same underlying field, and the former is generated as a field over $\R$ by $x$ and the $\ell^{(n)}$, so  applying the above to $f=\ell$ shows that $H^\phi$ is generated as a differential field~over~$\R$ by $x$ and $\ell$. 
We also have a differential field isomorphism
$h\mapsto h^\circ\colon H^\phi\to H^\circ=H\circ \ell^{\operatorname{inv}}$.
%(See~[ADH, 5.7.1].) 
This yields
$H^\circ=\R\langle\ell^{\operatorname{inv}},x\rangle$ and  
\eqref{eq:trdegellinv}. 
Suppose now that~$\ell$ is $\d$-algebraic over~$\R$; then
by additivity of~$\operatorname{trdeg}$,      
$$\operatorname{trdeg}\!\big(\R\langle x, \ell\rangle|\R\big)\ =\ 
\operatorname{trdeg}\!\big(\R\langle \ell,x\rangle|\R\langle\ell\rangle\big)+
\operatorname{trdeg}\!\big(\R\langle \ell\rangle|\R\big)\ \leq\ 
1+\operatorname{trdeg}\!\big(\R\langle \ell\rangle|\R\big),$$
and so by \eqref{eq:trdegellinv}:
$$\operatorname{trdeg}\!\big(\R\langle \ell^{\operatorname{inv}}\rangle|\R\big)\ \leq\ 
\operatorname{trdeg}\!\big(\R\langle x, \ell^{\operatorname{inv}}\rangle|\R\big)\ \leq\
\operatorname{trdeg}\!\big(\R\langle \ell\rangle|\R\big)+1,$$
hence  $\ell^{\operatorname{inv}}$ is $\d$-algebraic over $\R$.
\end{proof}

\noindent
In Corollary~\ref{cor:trdegellinv} below we prove a uniform version of Lemma~\ref{lem:trdegellinv}.
To prepare for this we prove the next  two lemmas, where $R$ is a differential ring and $x\in R$, $x'=1$. Also,~$\phi\in R^\times$, and we take distinct differential indeterminates $U$, $X$, $Y$ and let 
$G^n_k\in\Q\{U\}\subseteq R^\phi\{U\}$ ($k=1,\dots,n$)   be as in~[ADH, p.~292], so with  $\der$ and~$\derdelta=\phi^{-1}\der$ denoting the derivations
of $R$ and $R^\phi$, we have in $R^\phi[\derdelta]$ for $\der=\phi\derdelta$: 
$$\der^n\  =\  G^n_n(\phi)\cdot\derdelta^n+G^n_{n-1}(\phi)\cdot\derdelta^{n-1}+\cdots+G^n_1(\phi)\cdot\derdelta.$$
Recall that the $G^n_k$ do not depend on $R$, $x$, $\phi$.

\begin{lemma}\label{lem:Hnk}
There are $H^n_k\in\Q\{X'\}\subseteq\Q\{X\}\subseteq R^\phi\{X\}$ \textup{(}$k=1,\dots,n$\textup{)}, independent of $R$,  $x$, 
$\phi$, such that
$G^n_k(\phi)=\phi^{2n-1} H^n_k(x)$.
\end{lemma}
\begin{proof}
By induction on $n\geq 1$.  
For $n=1$ we have $G^1_1=U$, so~$H^1_1:=1$ does the job. Suppose for a certain $n\geq 1$ we have
$H^n_k$ ($k=1,\dots,n$) with the desired properties, and let $k\in\{1,\dots,n+1\}$. 
Now~$G^{n+1}_k = U\cdot\big( \derdelta(G^n_k)+G^n_{k-1} \big)$  by [ADH, (5.7.2)] (with~$G^n_0:=0$), so   using
$\derdelta(\phi)=-\phi^2\derdelta^2(x)$ and setting $H^n_0:=0$, 
\begin{align*}
G^{n+1}_k(\phi)\ 	
					&=\ \phi\cdot \big( (2n-1)	\phi^{2n-2}\derdelta(\phi)H^n_k(x)+ \phi^{2n-1} \derdelta(H^n_k(x)) + \phi^{2n-1} H^n_{k-1}(x) \big)\\
					&=\ \phi^{2n+1}\big( (1-2n)\derdelta^2(x)H^n_k(x) + \derdelta(x)\derdelta(H^n_k(x)) +\derdelta(x)H^n_{k-1}(x) \big).
\end{align*}
Thus we can take \[H^{n+1}_k\ :=\ (1-2n)X''H^n_k + X' (H^n_k)' + X' H^n_{k-1}. \qedhere\] 
\end{proof}

\begin{lemma}\label{pxynq}
Let  $C$ be a subfield of $C_R$ and $P\in C\{X,Y\}\subseteq R\{X,Y\}$. Then
there are~$N\in\N$ and~$Q\in C\{X,Y\}\subseteq R^\phi\{X,Y\}$ such that 
$P(x,Y)^\phi=\phi^N Q(x,Y)$ in $R^\phi\{Y\}$.
Here we can take $N$,~$Q$ independent  of  $x$, $\phi$.
\end{lemma}
\noindent
Note that $C\{X,Y\}$ as a differential subring of $R\{X,Y\}$ is the same as $C\{X,Y\}$ as a differential subring of $R^\phi\{X,Y\}$,
but ``$P\in C\{X,Y\}\subseteq R\{X,Y\}$'' indicates that $P$ is considered as an element of $R\{X,Y\}$ when substituting in $P$, while 
``$Q\in C\{X,Y\}\subseteq R^\phi\{X,Y\}$'' indicates that $Q$ is taken as an element of $R^\phi\{X,Y\}$ when substituting in $Q$.

\begin{proof} For $i=1,2$, let $P_i\in C\{X,Y\}$, $N_i\in\N$, and~$Q_i\in C\{X,Y\}\subseteq R^\phi\{X,Y\}$   
be such that  
$P_i(x,Y)^\phi  = \phi^{N_i}   Q_i(x,Y)$.
Then  
$$(P_1\cdot P_2)(x,Y)^\phi\  =\ \phi^{N_1+N_2}   (Q_1\cdot Q_2)(x,Y).$$
Moreover, $\derdelta(x)=\phi^{-1}$, hence if  $N_1\leq N_2$, then
$$(P_1+ P_2)(x,Y)^{\phi}\ =\  \phi^{N_2}   Q(x,Y)\quad\text{for $Q:= (X')^{N_2-N_1} Q_1 +  Q_2$.}$$
For $P=X$ we have $P(x,Y)^\phi=x=\phi\cdot x\derdelta(x)$, so $N=1$ and $Q=XX'$ works. For~$P=Y$ we can take
$N=0$ and $Q=Y$. 
It is enough to prove the lemma for $P$ such that no monomial in $P$ has any factor $X^{(m)}$ with $m\ge 1$.
 Thus it only remains to do the case $P=Y^{(n)}$ ($n\geq 1$).
With $H^n_k$ as in Lemma~\ref{lem:Hnk} we have
$$(Y^{(n)})^\phi\ =\ G^n_n(\phi)Y^{(n)}+\cdots+G^n_1(\phi)Y'\ =\ \phi^N Q(x,Y)$$
for   $N:=2n-1$ and $Q:=H^n_n(X)Y^{(n)}+\cdots+H^n_1(X)Y'$.
\end{proof}

\noindent
In the next lemma $x$ has its usual meaning as the germ in $\Calinf$ of the identity function on $\R$, we take $R$ as the differential ring
$\Calinf[\imag]$ and $C\{X,Y\}$ as a differential subring of $R\{X,Y\}$ for any subfield $C$ of $\C=C_R$. 

\begin{lemma}\label{lem:pcbullet}
Let $P\in C\{X,Y\}$ where $C$ is a subfield of $\C$. Then there are~${N\in\N}$ and~$P^\bullet\in C\{X,Y\}$  such that  for all $y\in R$ and $\ell\in \Calinf$ with $\ell(t)\to +\infty$ as~${t\to +\infty}$
and $\ell'(t)>0$, eventually, we have for $\phi:= \ell'$:
$$P(x,y)\circ\ell^{\operatorname{inv}}\ =\ \big(\phi\circ\ell^{\operatorname{inv}}\big)^{N} \cdot P^\bullet(\ell^{\operatorname{inv}},y\circ\ell^{\operatorname{inv}})\ \text{ in }R.$$
\end{lemma}
\begin{proof}
Let $\ell\in\Calinf$ be such that  $\ell(t)\to +\infty$ as $t\to +\infty$
and $\ell'(t)>0$, eventually, and set $\phi:= \ell'$. 
For $P_x:= P(x,Y)\in R\{Y\}$ and $y\in R$ we have 
$P(x,y)=P_x(y)=P_x^{\phi}(y)= \phi^NQ(x,y)$, with $N\in \N$ and $Q\in R^\phi\{X,Y\}$ as in Lemma~\ref{pxynq}, so
$$P(x,y)\circ\ell^{\operatorname{inv}}\ =\  \phi^NQ(x,y)\circ \ell^{\operatorname{inv}}\ =\ \big(\phi\circ\ell^{\operatorname{inv}}\big)^N\cdot Q(x,y)\circ\ell^{\operatorname{inv}}.$$
Let $P^\bullet$ be the element of $C\{X,Y\}$ that is mapped to $Q\in R^\phi\{X,Y\}$ under the
ring inclusion $C\{X,Y\}\to R^\phi\{X,Y\}$. The latter is not in general a differential ring morphism, but we have
the differential ring isomorphism
$$y\mapsto y  \circ\ell^{\operatorname{inv}}\ \colon\ R^\phi\to R\circ\ell^{\operatorname{inv}}=R,$$ which gives for $y\in R$ that
\[  Q(x,y)\circ\ell^{\operatorname{inv}}\ =\ P^\bullet(x\circ \ell^{\operatorname{inv}}, y\circ\ell^{\operatorname{inv}})\ =\ P^\bullet(\ell^{\operatorname{inv}},y\circ\ell^{\operatorname{inv}}). \qedhere \]
\end{proof}

\begin{cor}\label{cor:trdegellinv}
For each $P\in\R\{X,Y\}$   there is a $P^\bullet\in\R\{X,Y\}$  such that  for all $\ell\in\Calinf$ with $\ell(t)\to +\infty$ as $t\to +\infty$
and $\ell'(t)>0$, eventually, we have
$$P(x,\ell)=0\ \Longleftrightarrow\ P^\bullet(\ell^{\operatorname{inv}},x)=0.$$
\end{cor}

\noindent
We now indicate how Lemma~\ref{lem:pcbullet} and   Corollary~\ref{cor:trdegellinv} go through for transseries.
Recall from~[ADH, A.7] that there is a unique operation
$$(f,g)\mapsto f\circ g\colon \mathbb T\times\mathbb T^{>\R} \to \mathbb T$$
such that the following conditions hold for all $g\in \mathbb T^{>\R}$:
\begin{enumerate}
\item $x\circ g=g$;
\item $f\mapsto f\circ g\colon\mathbb T\to\mathbb T$ is an $\R$-linear embedding of ordered exponential fields;
\item $f\mapsto f\circ g\colon\mathbb T\to\mathbb T$ is strongly additive.
\end{enumerate}
By \cite[Proposition~6.3]{vdDMM}  the Chain Rule holds: 
$$(f\circ g)'=(f'\circ g)\cdot g'\qquad (f\in\mathbb T,\ g\in \mathbb T^{>\R}).$$
Moreover,   $(f,g)\mapsto f\circ g$ restricts to a binary operation on $\mathbb T^{>\R}$
which makes $\mathbb T^{>\R}$ a group with identity element $x$. For $f\in \mathbb T^{>\R}$ we denote the unique~$g\in \mathbb T^{>\R}$
with~$f\circ g=x$ by $g=f^{\operatorname{inv}}$.  We extend  $\circ$ in a unique way to an operation
$$(f,g)\mapsto f\circ g\colon \mathbb T[\imag]\times\mathbb T^{>\R} \to \mathbb T[\imag]$$
by requiring that for all $g\in \mathbb T^{>\R}$, the operation $f\mapsto f\circ g\colon  \mathbb T[\imag]\to  \mathbb T[\imag]$
is $\C$-linear. It follows that for all $g\in \mathbb T^{>\R}$ the operation $f\mapsto f\circ g\colon  \mathbb T[\imag]\to  \mathbb T[\imag]$ is a field embedding. For $f\in\mathbb T[\imag]$, $g,h\in\mathbb T^{>\R}$ we have $(f\circ g)\circ h=f\circ (g\circ h)$ [ADH, A.7(vi)], so $\mathbb T[\imag]\circ h=\mathbb T[\imag]$.  
For $\ell\in\mathbb T^{>\R}$ and $\phi:=\ell'$ we have a differential field isomorphism
$$y\mapsto y\circ\ell^{\operatorname{inv}}\colon \mathbb T[\imag]^\phi \to \mathbb T[\imag]\circ \ell^{\operatorname{inv}}=\mathbb T[\imag]. $$
Let $P\in C\{X,Y\}$ where $C$ is a subfield of $\C$. Let
$N\in\N$ and $P^\bullet\in C\{X,Y\}$ be as obtained in the proof of
Lemma~\ref{lem:pcbullet}. 
Then that proof gives for all $y\in \mathbb T[\imag]$, $\ell\in \mathbb T^{>\R}$, and  $\phi:= \ell'$:
$$P(x,y)\circ\ell^{\operatorname{inv}}\ =\ \big(\phi\circ\ell^{\operatorname{inv}}\big)^{N} \cdot P^\bullet(\ell^{\operatorname{inv}},y\circ\ell^{\operatorname{inv}})\ \text{ in $\mathbb T[\imag]$.}$$
Hence for $C=\R$ we have $P^\bullet\in\R\{X,Y\}$ and for all
$\ell\in\mathbb T^{>\R}$:
$$P(x,\ell)=0\ \Longleftrightarrow\ P^\bullet(\ell^{\operatorname{inv}},x)=0.$$

\section{Upper and Lower Bounds on the Growth of Hardian Germs\astr}\label{sec:upper lower bds}

\noindent
This section elaborates on \cite{Boshernitzan82,Boshernitzan86,Rosenlicht83}. It is not used for proving our main theorem, but some of it is needed later, in the proofs of Corollary~\ref{cor:17.11},    Proposition~\ref{prop:translog}, and Theorem~\ref{thm:coinitial in E(H)}.
 
\subsection*{Generalizing logarithmic decomposition}
In this subsection $K$ is a differential ring and~${y\in K}$. 
In [ADH, p.~213] we defined the $n$th iterated logarithmic derivative of $y^{\langle n\rangle}$ when $K$ is a differential field.
Generalizing this, set $y^{\langle 0\rangle}:=y$, and recursively, if $y^{\langle n\rangle}\in K$ is defined and a unit in $K$, then
 $y^{\langle n+1\rangle}:=(y^{\langle n\rangle})^\dagger$, while otherwise $y^{\langle n+1\rangle}$ is not defined. (Thus if~$y^{\langle n\rangle}$ is defined, then so are $y^{\langle 0\rangle},\dots,y^{\langle n-1\rangle}$.)
  With~$L_n$ in~$\Z[X_1,\dots,X_n]$ as in~[ADH, p.~213], if $y^{\langle n\rangle}$ is defined, then
$$y^{(n)}\  =\  y^{\langle 0\rangle} \cdot L_n(y^{\langle 1\rangle},\dots, y^{\langle n\rangle}).$$
 If $y^{\langle n\rangle}$ is defined and $\i=(i_0,\dots,i_n)\in\N^{1+n}$, we set
 $$y^{\langle \i \rangle}\ :=\  (y^{\langle 0\rangle})^{i_0} (y^{\langle 1\rangle})^{i_1} \cdots (y^{\langle n\rangle})^{i_n}\in K.$$
Hence if $H$ is a differential subfield of $K$, $P\in H\{Y\}$ has order at most~$n$ and
logarithmic decomposition $P=\sum_{\i} P_{\langle \i\rangle} Y^{\langle \i\rangle}$ 
($\i$ ranging over $\N^{1+n}$, all~$P_{\langle \i\rangle}\in H$, and~$P_{\langle \i\rangle}=0$ for all but finitely many $\i$),
and $y^{\langle n\rangle}$ is defined, then~$P(y)=\sum_{\i}  P_{\langle \i\rangle} y^{\langle \i\rangle}$.
Below we apply these remarks to~$K=\Calinf$, where for $y\in K^\times$ we have $y^\dagger=(\log |y|)'$, hence
$y^{\langle n+1\rangle}=(\log |y^{\langle n\rangle}|)'$ if $y^{\langle n+1\rangle}$ is defined.

\subsection*{Transexponential germs}
For  $f\in\c$ we recursively define the germs~$\exp_n f$ in $\c$ by $\exp_0 f:=f$ and $\exp_{n+1} f:=\exp (\exp_n f)$.
Following \cite{Boshernitzan82} we say that a germ~$y\in\c$ is {\bf transexponential} if $y\geq\exp_n x$ for all~$n$.\index{germ!transexponential}\index{transexponential}
{\it In the rest of this subsection $H$ is a Hardy field.}\/
By Corollary~\ref{cor:val at infty} and Proposition~\ref{prop:Hardy field exts}:

\begin{lemma}\label{lem:bounded Hardy field ext, 1}
If the $H$-hardian germ $y$ is $\d$-algebraic over $H$, then $y\leq\exp_n h$ for some $n$ and some $h\in H(x)$.
\end{lemma}

\noindent
Thus each transexponential hardian germ is $\d$-transcendental (over $\R$).  
%(Apply Lemma~\ref{lem:expn, 2} to $K=\Li(\R\langle y\rangle)$.) 
%{\it In the rest of this subsection:  $y,z\in\Calinf$, with $y$  transexponential and hardian.}\/ 
{\it In the rest of this subsection:  $y\in\Calinf$ is transexponential and hardian, and $z\in\Calinf[\imag]$.}\/
Then~$y^{\langle n\rangle}$ is defined, and~$y^{\langle n\rangle}$ is also transexponential and hardian,  for all  $n$.
Next some variants of results from Section~\ref{sec:it log derivative}. % (Lemmas~\ref{lem:expn, 1} and \ref{lem:expn, 2}, Corollary~\ref{cor:expn}, and Proposition~\ref{prop:val at infty}). 
For this, let $n$ be given and let $f\in\Calinf$, not necessarily hardian, be such that~$f\succ 1$, $f\geq 0$, and~$y\succeq\exp_{n+1}f$.

\begin{lemma}\label{lem:expn variant}
We have   $y^\dagger\succeq\exp_n f$ and $y^{\langle n\rangle}\succeq\exp f$.
\end{lemma}
\begin{proof} Since $y\succeq \exp_2 x$, we have
 $\log y \succeq \exp x$ by Lemma~\ref{lem:log preceq},  and thus 
$y^\dagger=(\log y)'\succeq \log y$.  Since $y\succeq\exp_{n+1}f$, the same lemma gives $\log y \succeq \exp_n f$. Thus~$y^\dagger\succeq\exp_n f$.
Now the second statement follows by an easy induction. 
%from the first in the same way that Lemma~\ref{lem:expn, 1} implied Lemma~\ref{lem:expn, 2}.
\end{proof}

\begin{cor}\label{cor:expn variant}
Let $\i\in\Z^{1+n}$  and suppose $\i>0$ lexicographically. Then $y^{\langle\i\rangle} \succ f$.
\end{cor}
\begin{proof}
Let~$m\in\{0,\dots,n\}$ be minimal such that $i_m\neq 0$; so $i_m\geq 1$.
The remarks after Corollary~\ref{cor:itpsi} then give $y^{\langle\i\rangle} \succ 1$ and
$[v(y^{\langle\i\rangle})]=[v(y^{\langle m\rangle})]$, so we have~$k\in\N$, $k\geq 1$, such that $y^{\langle\i\rangle} \succeq (y^{\langle m\rangle})^{1/k}$.
Then Lemma~\ref{lem:expn variant} gives
$y^{\langle\i\rangle} \succeq (y^{\langle m\rangle})^{1/k}\succeq (\exp f)^{1/k}\succ f$ as required.
\end{proof}

\noindent
In the next proposition and lemma~$P\in H\{Y\}^{\neq}$ has order at most~$n$, and~$\i$,~$\j$,~$\k$ range over~$\N^{1+n}$.
Let  $\j$ be lexicographically maximal such that~$P_{\<\j\>}\neq 0$, and 
choose~$\k$ so that
$P_{\<\k\>}$ has   minimal valuation.
If $P_{\<\k\>}/P_{\<\j\>} \succ x$,  set
$f:=|P_{\<\k\>}/P_{\<\j\>}|$; otherwise set $f:=x$.
 Then~$f\in H(x)$, $f>0$, $f\succ 1$, and $f\succeq P_{\<\i\>}/P_{\<\j\>}$ for all~$\i$.
 
\begin{prop}\label{prop:val at infty variant}
We have $P(y)\sim P_{\<\j\>}y^{\langle \j\rangle}$ and thus 
$$P(y)\in \big(\Calinf\big)^\times, \qquad \sgn P(y)\ =\ \sgn P_{\<\j\>}\neq 0.$$
\end{prop}
\begin{proof}
For $\i<\j$ we have  $y^{\langle \j-\i \rangle}\succ f\succeq P_{\<\i\>}/P_{\<\j\>}$ by Corollary~\ref{cor:expn variant},
therefore~$P_{\<\j\>}y^{\langle \j\rangle} \succ P_{\<\i\>}y^{\langle \i\rangle}$.
Thus $P(y)\sim P_{\<\j\>}y^{\langle \j\rangle}$.
\end{proof}

\begin{lemma}\label{lem:transexp}
Suppose  that $z^{\langle n\rangle}$  is defined and $y^{\langle i\rangle} \sim z^{\langle i\rangle}$ for $i=0,\dots,n$.
Then~$P(y) \sim P(z)$.
\end{lemma}
\begin{proof}
For all $\i$ with $P_{\langle \i\rangle}\neq 0$ we have $P_{\langle \i\rangle}y^{\langle\i\rangle} \sim P_{\langle \i\rangle}z^{\langle\i\rangle}$, by Lem\-ma~\ref{lem:sim props}.
Now use that for   $\i\neq\j$ we have
$P_{\langle \i\rangle} y^{\langle\i\rangle} \prec P_{\langle \j\rangle} y^{\langle\j\rangle}$ by the proof of Proposition~\ref{prop:val at infty variant}.
\end{proof}

\noindent
From here on $n$ is no longer fixed. 

\begin{cor}[{Boshernitzan~\cite[Theorem~12.23]{Boshernitzan82}}] \label{cor:Bosh 12.23}
Suppose   $y\geq \exp_n h$ for all~$h\in H(x)$  and all $n$. Then $y$ is $H$-hardian.
\end{cor}

\noindent
This is an immediate consequence of Proposition~\ref{prop:val at infty variant}. (In \cite{Boshernitzan82}, the proof of this fact  is only indicated.)
From Lemma~\ref{lem:transexp} we also obtain: 

\begin{cor}\label{cor:transexp}
Suppose that $y$ is as in Corollary~\ref{cor:Bosh 12.23} and $z\in\Calinf$, and 
$z^{\langle n\rangle}$ is defined and $y^{\langle n\rangle} \sim z^{\langle n\rangle}$, for all $n$. Then
$z$ is $H$-hardian, and there is a unique ordered differential field isomorphism $H\langle y\rangle \to H\langle z\rangle$ over $H$ which sends~$y$ to $z$.
\end{cor}

\noindent
Lemma~\ref{lem:y H-hardian crit} below contains another criterion for $z$ to be $H$-hardian. This involves a certain binary relation $\sim_\infty$ on germs defined in the next subsection.
Lemma~\ref{lem:transexp} also yields a complex version of Corollary~\ref{cor:transexp}:

\begin{cor}\label{cor:transexp, complex}
Suppose that $y$ is as in Corollary~\ref{cor:Bosh 12.23} and that
$z^{\langle n\rangle}$ is defined and $y^{\langle n\rangle} \sim z^{\langle n\rangle}$, for all $n$. Then
$z$ generates a differential subfield $H\langle z\rangle$ of $\Calinf[\imag]$, 
and there is a unique  differential field isomorphism $H\langle y\rangle \to H\langle z\rangle$ over $H$ which sends~$y$ to $z$.
Moreover, the binary relation $\preceq$ on $\c[\imag]$ restricts to a dominance relation on $H\langle z\rangle$
which makes this an isomorphism of valued differential fields.
\end{cor}

\subsection*{A useful equivalence relation}
We set
$$\Calinf[\imag]^{\preceq}\ :=\ \big\{ f\in \Calinf[\imag]: f^{(n)}\preceq 1\text{ for all $n$} \big\}\ \subseteq\ \c[\imag]^{\preceq},$$
a differential $\C$-subalgebra of  $\Calinf[\imag]$, and
 $$\mathcal I\ :=\ \big\{ f\in \Calinf[\imag]: f^{(n)}\prec 1\text{ for all $n$} \big\}\ \subseteq\ \Calinf[\imag]^{\preceq},$$
a differential ideal of $\Calinf[\imag]^{\preceq}$
(thanks to the Product Rule).
Recall from the remarks preceding   Lemma~\ref{lem:sim props} that~$(\c[\imag]^{\preceq})^\times = \c[\imag]^{\asymp}$. 

\begin{lemma}\label{lem:units of cinf[i]preceq}
The group of units of $\Calinf[\imag]^{\preceq}$ is
$$\Calinf[\imag]^{\asymp}\ :=\ \Calinf[\imag]^{\preceq}\cap\c[\imag]^{\asymp}\ =\ \big\{ f\in \Calinf[\imag]: f \asymp 1,\ f^{(n)}\preceq 1\text{ for all $n$} \big\}.$$
Moreover, $1+\mathcal I$ is a subgroup of~$\Calinf[\imag]^\asymp$.
\end{lemma}
\begin{proof}
It is clear that
$$(\Calinf[\imag]^{\preceq})^\times\ \subseteq\ \Calinf[\imag]^{\preceq}\cap  (\c[\imag]^{\preceq})^\times\ =\  \Calinf[\imag]^{\preceq}\cap\c[\imag]^{\asymp}\ =\ \Calinf[\imag]^{\asymp}.$$
Conversely, suppose $f\in \Calinf[\imag]$ satisfies $f\asymp 1$ and $f^{(n)}\preceq 1$ for all $n$.
For each~$n$ we have $Q_n\in\Q\{X\}$ such that
$(1/f)^{(n)} =  Q_n(f)/f^{n+1}$, hence~$(1/f)^{(n)} \preceq 1$. Thus~$f\in (\Calinf[\imag]^{\preceq})^\times$.
This shows the first statement.  
Clearly~${1+\mathcal I}\subseteq \Calinf[\imag]^{\asymp}$, and  $1+\mathcal I$ is closed under multiplication.
If $\delta\in\mathcal I$, then~$1+\delta$ is a unit of $\Calinf[\imag]^{\preceq}$ and~$(1+\delta)^{-1}=1+\varepsilon$ where $\varepsilon=-\delta(1+\delta)^{-1}\in\mathcal I$. 
%Hence $1+\mathcal I$ is a subgroup of~$\Calinf[\imag]^\asymp$.
\end{proof}

\noindent
For $y,z\in\c[\imag]^\times$ we  define
$$y\sim_\infty z \quad:\Longleftrightarrow\quad y\in z\cdot (1+\mathcal I);$$
hence $y\sim_\infty z\Rightarrow y\sim z$.
Lemma~\ref{lem:units of cinf[i]preceq} yields that $\sim_\infty$ is  an equivalence relation on~$\c[\imag]^\times$, and
for $y_i,z_i\in\c[\imag]^\times$ ($i=1,2$) we   have
$$y_1\sim_\infty y_2\quad \&\quad z_1\sim_\infty z_2 \qquad\Longrightarrow\qquad y_1z_1 \sim_\infty y_2z_2,\quad y_1^{-1}\sim_\infty y_2^{-1}.$$

\begin{lemma}\label{lem:y' siminfty z'}
Let $y,z\in\c^1[\imag]^\times$ with $y\sim_\infty z$ and $z\in z'\,\Calinf[\imag]^{\preceq}$. Then 
$$y', z'\in \c[\imag]^\times, \qquad y'\sim_\infty z'.$$
\end{lemma}
\begin{proof}
Let $\delta\in\mathcal I$ and $f\in\Calinf[\imag]^{\preceq}$  with $y=z(1+\delta)$ and $z=z'f$.
Then $z'\in\c[\imag]^\times$ and
$y'=z'(1+\delta)+z\delta'=z'(1+\delta+f\delta')$ where $\delta+f\delta'\in\mathcal I$, so $y'\sim_\infty z'$.
\end{proof}

\noindent
If  $\ell\in\c^n[\imag]$ and $f\in\c^n$ with $f\geq 0$, $f\succ 1$, then $\ell\circ f\in\c^n[\imag]$.
In fact, for $n\geq 1$ and~$1\le k\le n$ we have
a differential polynomial~$Q^n_k\in\Q\{X'\}\subseteq\Q\{X\}$ of order~$\le n$,  isobaric of weight~$n$, and homogeneous of degree $k$, 
such that for all such $\ell$, $f$, 
$$(\ell\circ f)^{(n)}\  =\  (\ell^{(n)}\circ f)\,Q^n_n(f)+\cdots+(\ell'\circ f)\,Q^n_1(f).$$
For example,
$$Q^1_1=X',\quad Q^2_2=(X')^2,\ Q^2_1=X'',\quad Q^3_3=(X')^3,\ Q^3_2=3X'X'',\ Q^3_1=X'''.$$
%\begin{align*}
%(\ell\circ f)'		&= (\ell'\circ f)f',\\
%(\ell\circ f)''		&= (\ell''\circ f)(f')^2+(\ell'\circ f)f'',\\
% (\ell\circ f)'''	&= (\ell'''\circ f)(f')^3+(\ell''\circ f)\,3f'f''+(\ell'\circ f)\,f'''.
%\end{align*}
The following Lemma is only used in the proof of Theorem~\ref{thm:coinitial in E(H)} below.

\begin{lemma}\label{lem:difference in I} Let $f,g\in\Calinf$ be such that $f,g\geq 0$ and $f,g\succ 1$, and set~$r:={g-f}$.
Suppose  $P(f)\cdot Q(r)\prec 1$ for all $P,Q\in\Q\{Y\}$ with $Q(0)=0$,
and let~$\ell\in\Calinf[\imag]$ be such that $\ell'\in\mathcal I$. Then $\ell\circ g - \ell\circ f\in \mathcal I$.
\end{lemma}
\begin{proof}
Treating real and imaginary parts separately we arrange $\ell\in\Calinf$. Note that $r\prec 1$. 
Taylor expansion [ADH, 4.2] for $P\in \Q\{X\}$ of order~$\leq n$ gives
$$P(g)-P(f)\  =\  \sum_{|\i|\geq 1} \frac{1}{\i!} P^{(\i)}(f) \cdot r^{\i} \qquad (\i\in \N^{1+n}),$$
and thus $P(g)-P(f)\prec 1$ and $rP(g)\prec 1$. 
The Mean Value Theorem yields a germ~$r_n\in\mathcal G$ such that
$$\ell^{(n)}\circ g - \ell^{(n)}\circ f\ =\ \big(\ell^{(n+1)}\circ (f+r_n)\big)\cdot r\quad\text{and}\quad |r_n| \leq |r|.$$
Now $r_0\prec 1$, so  $\ell'\circ (f+r_0)\prec 1$, hence $\ell\circ g-\ell \circ f\prec 1$.
For $1\le k\le n$,
\begin{multline*}
(\ell^{(k)}\circ g)\,Q^n_k(g) - (\ell^{(k)}\circ f)\,Q^n_k(f) = \\
 (\ell^{(k)}\circ f)\, \big( Q^n_k(g)-Q^n_k(f) \big) +  \big(\ell^{(k+1)}\circ(f+r_k)\big)\cdot r    Q^n_k(g),
\end{multline*}
so $(\ell^{(k)}\circ g)\,Q^n_k(g) - (\ell^{(k)}\circ f)\,Q^n_k(f)\prec 1$, and thus $\big( \ell\circ g - \ell\circ f \big)^{(n)}\prec 1$.
\end{proof}

\noindent
We consider next the differential $\R$-subalgebra $$(\Calinf)^{\preceq}\ :=\ \Calinf[\imag]^{\preceq}\cap\Calinf \ \subseteq\ \c^{\preceq}$$ of $\Calinf$. In the rest of this subsection $H$ is a Hardy field and $y,z\in\Calinf$, $y,z\succ 1$. 
Note that $(\Calinf)^{\preceq}\cap H=\mathcal O_H$ and $\mathcal I\cap H=\smallo_H$. This yields:

\begin{lemma}\label{lem:y siminf z}
Suppose $y-z\in (\Calinf)^{\preceq}$ and $z$ is hardian. Then $y\sim_\infty z$.
\end{lemma}
\begin{proof} From $y=z+f$ with $f\in (\Calinf)^{\preceq}$ we obtain $y=z(1+fz^{-1})$. Now $z^{-1}\in \mathcal{I}$, so
$fz^{-1}\in \mathcal{I}$, and thus $y\sim_\infty z$. 
\end{proof}

\noindent
We now formulate a sufficient condition involving $\sim_\infty$ for $y$  to be  $H$-hardian.

\begin{lemma}\label{lem:y H-hardian crit}
Suppose  
$z$ is $H$-hardian  with $z\geq\exp_n h$ for all $h\in H(x)$ and all~$n$, and~$y\sim_\infty z$. Then  $y$ is  $H$-hardian, and  there is a unique ordered differential field isomorphism~$H\langle y\rangle\to H\langle z\rangle$ which is the identity on $H$ and sends  $y$ to $z$.
\end{lemma}
\begin{proof}
By Lemma~\ref{lem:bounded Hardy field ext, 1} we may replace $H$ by the Hardy subfield
 $\Li\!\big(H(\R)\big)$ of~$\Ex(H)$ to arrange that $H\supseteq\R$ is Liouville closed.
By Corollary~\ref{cor:transexp} (with the roles of~$y$,~$z$ reversed) it is enough to show that for each $n$,
$y^{\langle n\rangle}$ is defined, $y^{\langle n\rangle}\succ 1$, and~$y^{\langle n\rangle}\sim_\infty z^{\langle n\rangle}$. This holds by hypothesis for $n=0$. By Lemma~\ref{lem:expn, 1}, $z>H$ gives~$z^\dagger>H$, so~$z=z'f$ with $f \prec 1$ in the Hardy field $H\langle z\rangle$, hence $f^{(n)}\prec 1$ for all $n$. So by Lemma~\ref{lem:y' siminfty z'}, $y^{\langle 1\rangle}=y^\dagger$ is defined, $y^{\langle 1\rangle}\in (\Calinf)^\times$, $y^{\langle 1\rangle} \sim_\infty z^{\langle 1\rangle}$, and thus $y^{\langle 1\rangle}\succ 1$. 
Assume for a certain $n\geq 1$ that $y^{\langle n\rangle}$ is defined,  $y^{\langle n\rangle}\succ 1$, and~$y^{\langle n\rangle}\sim_\infty z^{\langle n\rangle}$. Then~$z^{\langle n\rangle}$ is $H$-hardian and $H<z^{\langle n\rangle}$ by Lemma~\ref{lem:expn, 3}.
Hence  by the case $n=1$ applied to $y^{\langle n\rangle}$, $z^{\langle n\rangle}$ in place of~$y$,~$z$, respectively,  $y^{\langle n+1\rangle}=(y^{\langle n\rangle})^\dagger$ is defined, $y^{\langle n+1\rangle}\succ 1$, and~$y^{\langle n+1\rangle}\sim_\infty z^{\langle n+1\rangle}$.
\end{proof}

\noindent
The next two corollaries are Theorems~13.6 and 13.10, respectively, in  \cite{Boshernitzan82}.

\begin{cor}\label{cor:yhardian}
Suppose $z$ is trans\-ex\-po\-nen\-tial and  hardian, and~$y-z\in (\Calinf)^{\preceq}$.
Then  $y$  is hardian, and there is a unique  isomorphism~$\R\langle y\rangle\to\R\langle z\rangle$ of ordered differential fields that is the identity on $\R$ and sends  $y$ to $z$.
\end{cor}
\begin{proof}
Take $H:=\Li(\R)$. Then $z$ lies in a Hardy field extension of $H$, name\-ly~$\Li\!\big(\R\langle z\rangle\big)$, and $H<z$. So $y\sim_\infty z$ by Lemma~\ref{lem:y siminf z}. Now use Lemma~\ref{lem:y H-hardian crit}. \end{proof}

\begin{cor} \label{cor:Bosh13.10}
If~$z\in\Ex(H)^{>\R}$, then $z \leq \exp_n h$ for some $h\in H(x)$ and some~$n$. \textup{(}Thus if $x\in H$ and
$\exp H \subseteq H$, then $H^{>\R}$ is cofinal in $\Ex(H)^{>\R}$.\textup{)}
\end{cor}
\begin{proof}
Towards a contradiction, suppose $z\in\Ex(H)^{>\R}$ and $z>\exp_n h$ in $\Ex(H)$ for all $h\in H(x)$ and all $n$.   Set~$y:=z+\sin x$. Then 
 $y$ is $H$-hardian   by Lemmas~\ref{lem:y siminf z} and~\ref{lem:y H-hardian crit}, so $y$, $z$ lie in a common Hardy field
 extension of $H$, a contradiction. 
\end{proof}

\noindent
The same proof shows that Corollary~\ref{cor:Bosh13.10} remains true if $H$ is assumed to be a $\Ginf$-Hardy field and
$\Ex(H)$ is replaced by~$\Ex^\infty(H)$; likewise for $\omega$ in place of $\infty$.

\subsection*{Remarks on differential subfields of $\Calinf[\imag]$} 
Let $K$ be a   subfield of $\c[\imag]$. Then the following are equivalent:
\begin{enumerate}
\item The binary relation $\preceq$ on $\c[\imag]$ restricts to a dominance relation on $K$;
\item for all $f,g\in K$: $f\preceq g$ or $g\preceq f$;
\item for all $f\in K$: $f\preceq 1$ or $1\preceq f$.
\end{enumerate}
If $K\subseteq H[\imag]$ where $H$ is a Hausdorff field, then $\preceq$   restricts to a dominance relation on $K$. (See Section~\ref{sec:germs}.)
Moreover, the following are equivalent:
\begin{enumerate}
\item $K=H[\imag]$ for some Hausdorff field $H$; 
\item $\imag\in K$ and $\overline{f}\in K$ for each $f\in K$;
\item $\imag\in K$ and $\Re f, \Im f\in K$ for each $f\in K$.
\end{enumerate}

\noindent
Next a lemma
similar to Lemma~\ref{lem:y H-hardian crit}, but obtained using Corollary~\ref{cor:transexp, complex} instead of Corollary~\ref{cor:transexp}:

\begin{lemma}\label{lem:y H-hardian crit, complex}
Let $H$ be a Hardy field, let  
$z\in\Calinf$ be $H$-hardian  with $z\geq\exp_n h$ for all $h\in H(x)$ and all~$n$, and~$y\in\Calinf[\imag]$ with $y\sim_\infty z$. Then  $y$ generates a differential subfield~$H\langle y\rangle$ of $\Calinf[\imag]$, and  there is a unique   differential field isomorphism~$H\langle y\rangle\to H\langle z\rangle$ which is the identity on $H$ and sends  $y$ to $z$. The binary relation $\preceq$ on $\c[\imag]$ restricts to a dominance relation on $H\langle y\rangle$
which makes this an isomorphism of valued differential fields.
\end{lemma}

\noindent
We use the above at the end of the next subsection to produce a differential subfield of $\Calinf[\imag]$
that is not contained in $H[\imag]$ for any Hardy field $H$.

\subsection*{Boundedness} 
Let $H\subseteq\c$.
We say that $b\in\c$ {\bf bounds $H$} if $h\leq b$ for each~$h\in H$. We call $H$ {\bf bounded}\index{Hardy field!bounded} if some $b\in\c$ bounds $H$,
and we    call $H$ {\bf unbounded}\index{Hardy field!unbounded} if $H$ is not bounded.
If $H_1, H_2\subseteq\c$ and for each $h_2\in H_2$ there is an $h_1\in H_1$ with~$h_2\leq h_1$, then 
any $b\in\c$ bounding $H_1$ also  bounds $H_2$.
Every bounded subset of $\c$ is bounded by a germ in~$\Gom$; this follows from \cite[Lemma~14.3]{Boshernitzan82}:

\begin{lemma} \label{lem:Bosh 14.3} 
%\marginpar{accepted for now on faith}
For every $b\ge 0$ in $\c^\times$ there is a~$\phi\ge 0$ in $(\c^\omega)^\times$ such that $\phi^{(n)}\prec b$ for all $n$.
\end{lemma}

\noindent
Every countable subset of~$\c$ is bounded, by  
du Bois-Reymond~\cite{dBR75}; see also~\cite[Chapter~II]{Ha} or~\cite[Cha\-pi\-tre~V, p.~53, ex.~8]{Bou}.
Thus $H\subseteq\c$   is bounded if it is totally ordered by the partial ordering $\le$ of $\c$ and has countable cofinality. 
If $H$ is a Hausdorff field and $b\in \c$ bounds $H$, then~$b$ also bounds the real closure $H^{\operatorname{rc}}\subseteq\c$ of $H$ [ADH, 5.3.2]. {\it In the rest of this subsection~$H$ is a Hardy field.}\/

\begin{lemma}\label{lem:bounded Hardy field ext, 2}
Let $H^*$ be a  $\d$-algebraic Hardy field extension of $H$ and
suppose~$H$ is bounded. Then $H^*$ is also bounded.
\end{lemma}
\begin{proof} By [ADH, 3.1.11] we have $f\in H(x)^{>}$ such that for all~$g\in H(x)^\times$ there are~$h\in H^\times$ and $q\in \Q$ with $g\asymp hf^q$. Hence $H(x)$ is bounded. 
Replacing $H$, $H^*$ by~$H(x)^{\operatorname{rc}}$,~$\Li\!\big(H^*(\R)\big)$, respectively, we arrange that $H$ is real closed 
with $x\in H$, and
$H^*\supseteq\R$ is Liouville closed.
Let $b\in\c$  bound~$H$. Then any  $b^*\in\c$ such that~$\exp_n b \leq b^*$ for all $n$ 
bounds $H^*$, by Lemma~\ref{lem:bounded Hardy field ext, 1}.
\end{proof}

\begin{lemma}\label{lem:bounded Hardy field ext, 3} Suppose that $H$ is bounded and  $f\in \Calinf$ is hardian over  $H$. Then~$H\<f\>$ is bounded. 
\end{lemma}
\begin{proof}  Lemma~\ref{lem:bounded Hardy field ext, 2} gives that $\Li\!\big(H(\R)\big)$ is bounded; also, $f$ remains hardian over~$\Li\!\big(H(\R)\big)$. Using this we arrange that $H$ is Liouville closed. The case that~$H\<f\>$ has no element $>H$ is trivial, so assume we have $y\in H\langle f \rangle$ with~$y>H$. Then $y$ is $\d$-transcendental over $H$ and the sequence
$y, y^2, y^3,\dots$  is cofinal in $H\<y\>$, by Corollary~\ref{cor:val gp at infty}, so $H\<y\>$ is bounded. Now use that  $f$ is
$\d$-algebraic over $H\<y\>$. 
\end{proof}

\begin{theorem}[{Boshernitzan~\cite[Theorem~14.4]{Boshernitzan82}}]\label{thm:Bosh 14.4}
Suppose $H$ is bounded. Then the perfect hull~$\Ex(H)$ of $H$ is $\d$-algebraic over $H$ and hence  bounded. If
$H\subseteq\Ginf$, then~$\Ex^\infty(H)$ is $\d$-algebraic over $H$; likewise with $\omega$ in place of~$\infty$.
\end{theorem}

\noindent
Using the results above the proof is not difficult. It is omitted in \cite{Boshernitzan82}, but we include it here for the sake of completeness.  First, a lemma:

\begin{lemma}\label{lem:Q bound}
Let  $b\in\c^\times$ bound $H$, let $\phi\ge 0$ in $\Calinf$ satisfy
$\phi^{(n)}\prec b^{-1}$ for all~$n$, and let $r\in \phi\cdot  (\Calinf)^{\preceq}$. Then $Q(r)\prec 1$ for all $Q\in H\{Y\}$ with $Q(0)=0$.
\end{lemma}
 \begin{proof} From  $\phi\in \mathcal I$ we obtain $r\in\mathcal I$, so it is enough that $hr^{(n)}\prec 1$ for all
  $h\in H$ and all~$n$. Now use
the Product Rule and $h\phi^{(n)}\prec hb^{-1}\preceq 1$ for $h\in H^\times$. 
\end{proof}

\begin{proof}[Proof of Theorem~\ref{thm:Bosh 14.4}]
Using Lemma~\ref{lem:bounded Hardy field ext, 2}, replace $H$ by $\Li\!\big(H(\R)\big)$ to arrange that~$H\supseteq\R$ is Liouville closed.
Let $b\in\c$ bound~$H$. Then $b$ also bounds $\Ex(H)$,  by Corollary~\ref{cor:Bosh13.10}. 
 Lemma~\ref{lem:Bosh 14.3} yields~$\phi\ge 0$ in $(\Gom)^\times$ such that~$\phi^{(n)} \prec b^{-1}$ for all~$n$;
 set $r:=\phi\cdot\sin x\in\Gom$. Then~$Q(r) \prec f$ for all~$f\in\Ex(H)^\times$ and $Q\in \Ex(H)\{Z\}$ with~$Q(0)=0$, by Lem\-ma~\ref{lem:Q bound}.

Suppose towards a contradiction that $f\in \Ex(H)$ is   $\d$-transcendental over~$H$, and
set $g:=f+r\in\Calinf$. Then $f$, $g$ are not in a common Hardy field, so $g$ is not hardian over $H$. On the other hand, let $P\in H\{Y\}^{\neq}$.
Then $P(f)\in\Ex(H)^\times$, and by Taylor expansion,
$$P(f+Z)\ =\ P(f)+Q(Z)\quad\text{ where  $Q\in \Ex(H)\{Z\}$ with $Q(0)=0$,}$$ 
so  $P(g)=P(f+r) \sim P(f)$. Hence $g$ is hardian over $H$, a contradiction.

The proof in the case where $H\subseteq\Ginf$ is similar, using the version of Corollary~\ref{cor:Bosh13.10}
for $\Ex^\infty(H)$; similarly for $\omega$ in place of $\infty$.
\end{proof}

\noindent
As to
the existence of   trans\-ex\-po\-nen\-tial hardian germs, we have:

\begin{theorem}\label{thm:Bosh 1.2}
For every $b\in\c$ there is a $\Gom$-hardian germ $y\geq b$. 
%\marginpar{accepted on faith for now} 
\end{theorem}

\noindent
This is Boshernitzan~\cite[Theo\-rem~1.2]{Boshernitzan86}, and leads to~\cite[The\-o\-rem~1.1]{Boshernitzan86}:

\begin{cor}\label{cor:Bosh 1.1} No maximal Hardy field is bounded. 
\end{cor}
\begin{proof}
Suppose $x\in H$, and $b\in\c$ bounds $H$. Take some $b^*\in\c$ such that~$b^*\geq \exp_n b$ for each~$n$.
Now Theorem~\ref{thm:Bosh 1.2} yields a $\Gom$-hardian germ $y\geq b^*$.
By Corollary~\ref{cor:Bosh 12.23}, $y$ is $H$-hardian, so $H\langle y\rangle$ is a proper Hardy field extension of~$H$.
\end{proof}

\noindent
The same proof shows also that no $\Ginf$-maximal Hardy field and no $\Gom$-maximal Hardy field is bounded.
In particular (Boshernitzan \cite[Theo\-rem~1.3]{Boshernitzan86}):

\begin{cor}\label{cor:Bosh 1.3}
Every  maximal Hardy field contains a transexponential  germ. Likewise with ``$\Ginf$-maxi\-mal'' or
``$\Gom$-maxi\-mal'' in place of ``maximal''.
\end{cor}

\begin{remark}
For $\Ginf$-Hardy fields, some of the above is in Sj\"odin's~\cite{Sj}, predating~\cite{Boshernitzan82,Boshernitzan86}: if $H$ is a bounded $\Ginf$-Hardy field, then so is~$\Li\!\big(H(\R)\big)$
\cite[The\-o\-rem~2]{Sj};   no maximal $\Ginf$-Hardy field is bounded \cite[The\-o\-rem~6]{Sj}; and  
$E:=\Ex^\infty(\Q)$ is bounded~\cite[The\-o\-rem~10]{Sj} with
$E\circ E^{>\R}\subseteq E$~\cite[The\-o\-rem~11]{Sj}.
\end{remark}

\noindent
We can now produce  a differential subfield $K$ of $\Gom[\imag]$ 
containing $\imag$ such that  $\preceq$  restricts to a dominance relation on $K$
making $K$ a $\d$-valued field of $H$-type with constant field $\C$, yet $K\not\subseteq H[\imag]$ for every Hardy field $H$: 

Take a transexponential $\Gom$-hardian germ $z$, 
and~$h\in\R(x)$ with $0\neq h\prec 1$. Then 
$\varepsilon:=h\ex^{x\imag}\in\mathcal I$, so $y:=z(1+\varepsilon)\in\Gom[\imag]$ with~$y\sim_\infty z$.
Lemma~\ref{lem:y H-hardian crit, complex} applied with $H=\R$ shows that $y$ generates a differential subfield~$K_0:=\R\langle y\rangle$
of $\Gom[\imag]$, and   $\preceq$ restricts to a dominance relation on $K_0$
making~$K_0$ a $\d$-valued field of $H$-type with constant field $\R$.
Then $K:=K_0[\imag]$ is a differential subfield of $\Gom[\imag]$ with constant field $\C$.
Moreover,~$\preceq$ also restricts to a dominance relation on $K$,   and this dominance
relation makes $K$
a $\d$-valued field of $H$-type~[ADH, 10.5.15]. We cannot have $K\subseteq H[\imag]$ where $H$ is a Hardy field, since
$\Im y=zh\sin x\notin H$.

\subsection*{Lower bounds on $\d$-algebraic hardian germs} 
In this subsection $H$ is a Hardy field. Let $f\in\c$ and $f\succ 1$, $f\ge 0$. Then the germ $\log f\in\c$ also satisfies
$\log f\succ 1$, $\log f\ge 0$. So we may inductively define the germs $\log_n f$ in $\c$ by $\log_0 f:=f$, $\log_{n+1} f:=\log \log_n f$.
Lemma~\ref{lem:bounded Hardy field ext, 1} gives exponential upper bounds on $\d$-algebraic $H$-hardian germs.
The next result leads to logarithmic lower bounds on such germs when $H$ is grounded.

\begin{theorem}[{Rosenlicht~\cite[Theorem~3]{Rosenlicht83}}]\label{thm:Ros83}
Suppose $H$ is grounded,  and let~$E$ be a Hardy field extension of $H$ such that $|\Psi_E\setminus\Psi_H|\leq n$ \textup{(}so $E$ is also 
%\marginpar{accepted on faith for now} 
grounded\textup{)}. Then there are $r,s\in\N$ with $r+s\leq n$ such that
\begin{enumerate}
\item[\textup{(i)}] for any $h\in H^>$ with $h\succ 1$ and $\max\Psi_H=v(h^\dagger)$, 
there exists~$g\in E^>$ such that $g\asymp \log_r h$
and
$\max\Psi_E = v(g^\dagger)$;
\item[\textup{(ii)}] for any $g\in E$ there exists $h\in H$ such that $g<\exp_s h$.
\end{enumerate}
\end{theorem}

\noindent
This theorem is most useful in combination with
the following lemma, which is~\cite[Pro\-po\-sition~5]{Rosenlicht83} (and also \cite[Lemma~2.1]{AvdD3} in the context of pre-$H$-fields).

\begin{lemma}\label{lemro83}
Let  $E$ be a Hardy field extension of $H$ such that $\operatorname{trdeg}(E|H)\leq n$. Then~$|\Psi_E\setminus\Psi_H|\leq n$. 
\end{lemma}

\noindent
From [ADH, 9.1.11] we recall that for $f,g\succ 1$ in a Hardy field we have $f^\dagger\preceq g^\dagger$ iff~$\abs{f}\leq\abs{g}^n$
for some~$n\geq 1$. (See also the discussion before Lemma~\ref{lem:Hardy type}.)
Thus by Lemma~\ref{lemro83} and Theorem~\ref{thm:Ros83}:

\begin{cor}\label{cor:Ros83}
Let $E$ be a Hardy field extension of
$H$ with $\operatorname{trdeg}(E|H)\leq n$, and let
$h\in H^>$ be such that $h\succ 1$ and $\max\Psi_H=v(h^\dagger)$. Then~$E$ is grounded, and
for all $g\in E$ with $g\succ 1$ there is an $m\geq 1$ such that $\log_n h  \preceq g^m$
\textup{(}and hence~$\log_{n+1} h \prec g$ for all $g\in E$ with $g\succ 1$\textup{)}.
\end{cor}

\noindent
Applying Corollary~\ref{cor:Ros83} to  $H=\R(x)$, $h=x$ yields:

\begin{cor}[{Boshernitzan~\cite[Proposition~14.11]{Boshernitzan82}}]\label{cor:Bosh 14.11}
If $y\in\c$ is hardian and $\d$-algebraic over $\R$, then the Hardy field~$E=\R(x)\langle y\rangle$ is grounded, and there is an~$n$ such that $\log_n x \prec g$
for all $g\in E$ with $g\succ 1$.
\end{cor}
 
\section{Second-Order Linear Differential Equations over Hardy Fields}\label{sec:order 2 Hardy fields} 

\noindent
In this section we review Boshernitzan's work~\cite[\S{}16]{Boshernitzan82} on adjoining non-oscillating solutions of second-order linear differential equations to Hardy fields,
deduce some consequences about complex exponentials over Hardy fields used later, and prove a conjecture  from~\cite[\S{}17]{Boshernitzan82}.
{\it Throughout this section~$H$ is a Hardy field.}\/

\subsection*{Oscillation over Hardy fields}
In this subsection we assume $f\in H$ and consider the 
linear differential equation 
\begin{equation}\label{eq:2nd order, 2} \tag{4L}
4Y''+fY\ =\ 0
\end{equation}
over $H$. The factor $4$ is to simplify certain expressions, in conformity with [ADH, 5.2]. 
In [ADH, 5.2]  we defined for any differential field $K$ functions 
$\omega\colon K\to K$ and~$\sigma\colon K^\times\to K$.
We define likewise 
$$\omega\ :\ \Go[\imag]\to \Gz[\imag], \qquad \sigma\ :\ \c^2[\imag]^\times \to \Gz[\imag]$$ 
by
$$\omega(z)\ =\ -2z'-z^2\quad\text{ and }\quad\sigma(y)\ =\ \omega(z)+y^2\text{ for $z:= -y^\dagger$.}$$
Note that $\omega(\Go)\subseteq\Gz$ and $\sigma\big( (\c^2)^\times \big) \subseteq\Gz$, and
$\sigma(y) = \omega(z+y\imag)$ for  $z:= -y^\dagger$. 
To clarify the
role of $\omega$ and $\sigma$ in connection with second-order linear differential equations, suppose
$y\in \Gt$  is a non-oscillating solution to \eqref{eq:2nd order, 2} with $y\neq 0$.
Then~${z:=2y^\dagger\in \Go}$ satisfies
$-2z'-z^2=f$, so $z$ generates a Hardy field~$H(z)$
with $\omega(z)=f$,  by Proposition~\ref{singer}, which in turn
yields a Hardy field $H(z,y)$ with~$2y^\dagger=z$.
Thus $y_1:=y$ lies in a Hardy field extension of $H$. From Lemma~\ref{lem:2nd lin indep sol} and 
Proposition~\ref{prop:Hardy field exts}(iv) we also obtain a solution $y_2$ to \eqref{eq:2nd order, 2}
in a Hardy field
extension of $H\<y_1\>=H(y,z)$ such that~$y_1$,~$y_2$ are $\R$-linearly independent; see also~\cite[Theo\-rem~2, Corollary~2]{Ros}. 
This shows: 
 
\begin{prop}\label{prop:2nd order Hardy field}
If $f/4$ does not generate oscillations, then $\Dx(H)$ contains $\R$-linearly independent solutions $y_1$, $y_2$ to \eqref{eq:2nd order, 2}.  
\end{prop}
%\begin{proof}
%Replacing $H$ by $\Ex(H)$, we  arrange that $H$ is perfect. We may also assume $y\neq 0$, so  $y\in (\Gt)^\times$,  and $z:=2y'/y\in\Go$  satisfies $2z'+z^2+f=0$.   Then $z\in H$ by Corollary~\ref{cor:Singer-Rosenlicht}. Now $y'+yg=0$ where $g:=-z/2\in H$, hence $y\in H$, again by Corollary~\ref{cor:Singer-Rosenlicht}.
%\end{proof}

\noindent
Indeed, if $f/4$ does not generate oscillations, then $\Dx(H)$ contains   solutions $y_1$, $y_2$ to~\eqref{eq:2nd order, 2}
with $y_1,y_2>0$ and $y_1\prec y_2$.  Here $y_1$ is determined up to multiplication by a factor in $\R^>$; we call such $y_1$ a {\bf principal solution} to~\eqref{eq:2nd order, 2}.\index{solution!principal} (Lemmas~\ref{lem:admissible pair, unique},~\ref{lem:admissible pair}.) 
See Section~\ref{sec:uplupo-freeness} for the subsets $\Upg(H)$, $\Upl(H)$ of $H$.

\begin{lemma}\label{lem:omegabar(H)}
Suppose $H$ is $\d$-perfect and $f/4$ does not generate oscillations, and
let $y\in H$ be a principal solution  to~\eqref{eq:2nd order, 2}.
Then $z:=2y^\dagger$ is the unique solution  of the equation $\omega(z)=f$
in~$\Upl(H)$. %\textup{(}In particular, $\omega(H)=\omega\big(\Upl(H)\big)=\bar{\omega}(H)$.\textup{)}
\end{lemma}
\begin{proof} We already know $\omega(z)=f$. The restriction of $\omega$ to $\Upl(H)$ is strictly increasing [ADH, 11.8.20], so it remains to show that $z\in \Upl(H)$. 
Let $h\in H$, $h'=1/y^2$. Then $h\succ 1$ by Corollary~\ref{cor:I1 I2}, hence
$1/y^2\in\Upg(H)$,
so~$z=-(1/y^2)^\dagger\in \Upl(H)$.
\end{proof}

\noindent 
By [ADH, p.~259], with $A=4\der^2+f\in H[\der]$ we have
$$\text{$4y''+fy=0$ for some $y\in H^\times$}\ \Rightarrow\ \text{$A$ splits over $H$} \ \Longleftrightarrow\ f\in\omega(H).$$
To simplify the discussion we now also introduce the subset \label{p:baromega}
$$\bar{\omega}(H):=  \big\{ f\in H: \text{$f/4$ does not generate oscillations} \big\}$$
of $H$. If $E$ is a Hardy field extension of $H$, then $\bar\omega(E)\cap H=\bar\omega(H)$.
By Corollary~\ref{cor:gen osc closed upward}, $\bar{\omega}(H)$ is downward closed, and  $\omega(H)\subseteq \bar{\omega}(H)$
by the discussion following~\eqref{eq:Riccati} in Section~\ref{sec:second-order}.

\begin{cor}  \label{cor:omega(H) downward closed}
If $H$ is $\d$-perfect, then 
\[\omega(H)\	=\ \bar\omega(H) 
		\	=\ \big\{ f\in H :\ \text{$4y''+fy=0$ for some $y\in H^\times$} \big\},\]
and $\omega(H)$ is  downward closed in $H$.
\end{cor}

%Together with Corollary~\ref{cor:gen osc closed upward} and    Proposition~\ref{prop:2nd order Hardy field}, this  yields:

\noindent
If $H$ is $\d$-perfect, then $H^\dagger=H$ by  Proposition~\ref{prop:Hardy field exts}.
%Corollary~\ref{cor:Singer-Rosenlicht}. 
The remarks   after~\eqref{eq:Riccati} show that a part of Corollary~\ref{cor:omega(H) downward closed} holds under
this weaker condition:

\begin{cor} 
If $H^\dagger=H$, then 
$$\omega(H)\ = \ \big\{ f\in H :\ \text{$4y''+fy=0$ for some $y\in H^\times$} \big\}.$$
\end{cor}

\noindent
Lemma~\ref{lem:2nd order inhom} and Proposition~\ref{prop:2nd order Hardy field} also yield:

\begin{cor}\label{cor:2nd order Hardy field, inhom}
If $f\in\bar\omega(H)$, then each~${y\in\Gt}$   such that $4y''+fy\in H$ is in~$\Dx(H)$.
\end{cor}
%\begin{proof}
%We first arrange again that $H$ is perfect. Proposition~\ref{prop:2nd order Hardy field} yields $\R$-linearly independent so\-lu\-tions~$y_1,y_2\in H$ of~\eqref{eq:2nd order, 2}.  Applying   the   variation of constants method to $y_1$, $y_2$  [ADH, 5.5.22] and using Corollary~\ref{cor:Singer-Rosenlicht} yields some~$z\in H$ with $4z''+fz=b$; since $y-z\in \R y_1+\R y_2\subseteq H$, we also have $y\in H$.
%\end{proof}

\noindent
For use in the proof of Corollary~\ref{cor:FLW} we record the following property of $\bar{\omega}(H)$:

\begin{lemma}\label{lem:Upgcapbaromega}
$\Upg(H)\cap\bar{\omega}(H)=\emptyset$.
\end{lemma}
\begin{proof}
We arrange that $H$ is $\d$-perfect. Hence $H\supseteq\R$ is Liouville closed and~$\bar{\omega}(H)=\omega(H)$ by Corollary~\ref{cor:omega(H) downward closed}. From $x^{-1}=x^\dagger\in \Upg(H)$ and $\sigma(x^{-1})=2x^{-2}\asymp (x^{-1})'\prec\ell^\dagger$ for all $\ell\succ 1$ in $H$ we obtain $\Upg(H)\subseteq \sigma\big(\Upg(H)\big){}^\uparrow$,
so $\Upg(H)\cap\omega(H)=\emptyset$ 
by [ADH, remark before 11.8.29].
\end{proof}

\noindent
Next some consequences of Proposition~\ref{prop:2nd order Hardy field} for more general linear differential equations of
order $2$: Let $g,h\in H$, and consider the linear differential equation
\begin{equation}\label{eq:2nd order, 3}\tag{$\tilde{\operatorname{L}}$}
Y''+gY'+hY\ =\ 0
\end{equation}
over $H$. An easy induction on $n$ shows that for a solution $y\in \Gt$ of \eqref{eq:2nd order, 3} we have~$y\in \mathcal{C}^n$ with $y^{(n)}\in Hy+Hy'$ for all $n$, so
$y\in\Calinf$. To reduce \eqref{eq:2nd order, 3}
to an equation~\eqref{eq:2nd order, 2} we take $f:=\omega(g)+4h=-2g'-g^2+4h\in H$, take $a\in \R$, and take a representative of~$g$ in~$\Cao$, also denoted by $g$, and let $G\in (\Gt)^\times$ be 
the germ of 
$$t\mapsto \exp\!\left(-\frac{1}{2}\int_a^t g(s)\,ds\right)\qquad(t\ge a).$$ 
This gives an isomorphism $y\mapsto Gy$ from the $\R$-linear space of 
solutions of \eqref{eq:2nd order, 2} in~$\Gt$ onto the
$\R$-linear space of solutions of \eqref{eq:2nd order, 3} in~$\Gt$,
and $y\in \Gt$ oscillates iff $Gy$ oscillates.
By Proposition~\ref{prop:Hardy field exts}, $G\in\Dx(H)$. 
Using $\frac{f}{4}=-\frac{1}{2}g'-\frac{1}{4}g^2+h$ we now obtain the following
germ version of Corollary~\ref{coroscgen}:
 
\begin{cor}\label{cor:char osc}
The following are equivalent:
\begin{enumerate}
\item[\textup{(i)}] some solution in $\Gt$ of \eqref{eq:2nd order, 3}  oscillates;
\item[\textup{(ii)}] all nonzero solutions in $\Gt$ of \eqref{eq:2nd order, 3} oscillate;
\item[\textup{(iii)}]  $-\frac{1}{2}g'-\frac{1}{4}g^2+h$ generates oscillations.
\end{enumerate}
Moreover, if $-\frac{1}{2}g'-\frac{1}{4}g^2+h$ does not generate oscillations, then all solutions of~\eqref{eq:2nd order, 3} in $\Gt$ belong to $\Dx(H)$. 
\end{cor}

\noindent
Set $A:=\der^2+ g\der + h$, and let $f=\omega(g)+4h$, $G$ be as above. Then $A_{\ltimes G}=\der^2 + \frac{f}{4}$.
Thus by combining Corollary~\ref{cor:2nd order Hardy field, inhom} and Corollary~\ref{cor:char osc} we obtain: 

\begin{cor}\label{cor:char osc, 1}
If \eqref{eq:2nd order, 3} has no oscillating solution in $\Gt$, and ${y\in\Gt}$ is such that~$y''+gy'+hy\in H$, then $y\in\Dx(H)$.  
\end{cor}

\noindent
The next corollary follows from Proposition~\ref{prop:2nd order Hardy field} and [ADH, 5.1.21]: 
 
\begin{cor}\label{cor:char osc, 2}
The following are equivalent, for  $A\in H[\der]$ and $f$ as above:
\begin{enumerate}
\item[\textup{(i)}] $f/4$ does not generate oscillations;
\item[\textup{(ii)}] $A$ splits over some Hardy field extension of $H$;
\item[\textup{(iii)}] $A$ splits over $\Dx(H)$.
\end{enumerate}
\end{cor}

\noindent
For  $A\in H[\der]$ and $f$ as before we have  $A_{\ltimes G}=\der^2 + \frac{f}{4}$ and $G^\dagger=-\frac{1}{2}g\in H$, so: 

\begin{cor}\label{cor:char osc, 3}
$A$ splits over $H[\imag]$ $\Longleftrightarrow$   $\der^2+\frac{f}{4}$ splits over $H[\imag]$.
\end{cor}

\noindent
Proposition~\ref{prop:2nd order Hardy field} and its corollaries \ref{cor:2nd order Hardy field, inhom}--\ref{cor:char osc, 1}
are from \cite[Theorems~16.17, 16.18, 16.19]{Boshernitzan82}, and Corollary~\ref{cor:omega(H) downward closed} is essentially 
\cite[Lemma~17.1]{Boshernitzan82}.

\medskip
\noindent
Proposition~\ref{prop:2nd order Hardy field} applies only when \eqref{eq:2nd order, 2} has a solution in $(\Gt)^\times$.
Such a solution might not exist, but \eqref{eq:2nd order, 2} does have $\R$-linearly independent solutions~${y_1, y_2\in \Gt}$, so~$w:= y_1y'_2-y'_1y_2\in \R^\times$. 
Set $y:= y_1+ y_2\imag$. Then~$4y'' + fy=0$ and~$y\in\Gt[\imag]^\times$, and for $z=2y^\dagger\in \Go[\imag]$   we have $-2z'-z^2=f$. Now
\begin{align*} z\ =\ \frac{2y_1'+2\imag y_2'}{y_1+\imag y_2}\ &=\
\frac{2y_1'y_1+2y_2'y_2- 2\imag(y_1'y_2-y_1y_2')}{y_1^2+y_2^2}\ =\ \frac{2(y_1'y_1+y_2'y_2)+2\imag w}{y_1^2+y_2^2},\\
\text{so }\ \Re z\ &=\ \frac{2(y_1'y_1+y_2'y_2)}{y_1^2+y_2^2}\in \Go, \qquad \Im z\ =\ \frac{2w}{y_1^2+y_2^2}\in \Gt.
\end{align*} 
Thus $\Im z\in (\Gt)^\times$ and $(\Im z)^\dagger=-\Re z$,  and so
$$\sigma(\Im z)\ =\ \omega\big({-(\Im z)^\dagger}+(\Im z)\imag\big)\ =\ 
\omega(z)\ =\ f\qquad\text{in $\Go$.}$$
Replacing~$y_1$ by $-y_1$ changes 
$w$ to $-w$; this way we can arrange ${w>0}$, so $\Im z> 0$.

\medskip
\noindent
Conversely, every $u\in  (\c^2)^\times$ such that $u>0$ and $\sigma(u)=f$ arises in this way. 
To see this, suppose we are given such $u$,  take
$\phi\in \c^3$ with $\phi'=\frac{1}{2}u$, and set
$$y_1\ :=\  \frac{1}{\sqrt{u}}\cos \phi, \qquad y_2\ :=\  \frac{1}{\sqrt{u}}\sin \phi \qquad \text{(elements of $\c^2$)}.$$
Then $\wr(y_1, y_2)=1/2$, and $y_1$, $y_2$ solve  \eqref{eq:2nd order, 2}.
To see the latter, consider  
$$y\ :=\ y_1+y_2\imag\ =\ \frac{1}{\sqrt{u}}\ex^{\phi\imag}\in \c^2[\imag]^\times$$ and note that $z:=2y^\dagger$ satisfies $$\omega(z)\ =\ \omega(-u^\dagger+u\imag)\ =\ \sigma(u)\ =\ f,$$ hence
$4y''+fy=0$. The computation above shows~$\Im z = 1/(y_1^2+y_2^2) = u$.
We have $\phi'> 0$, so either $\phi>\R$ or $\phi-c\prec 1$ for some $c\in\R$, with $\phi>\R$ iff $f/4$ generates oscillations. 
As to uniqueness of the above pair $(y_1,y_2)$, we have:

\begin{lemma}\label{lem:Im z}
Suppose $f\notin\bar\omega(H)$.
Let $\tilde y_1,\tilde y_2\in \c^2$ be  $\R$-linearly independent solutions of \eqref{eq:2nd order, 2} 
with $\wr(\tilde y_1, \tilde y_2)=1/2$.  Set $\tilde y:=\tilde y_1+\tilde y_2\imag$, $\tilde z:=2\tilde y^\dagger$. Then
$$ \Im \tilde z =u\quad\Longleftrightarrow\quad \tilde y=\ex^{\theta\imag}y \text{ for some $\theta\in\R$.}$$
\end{lemma}
\begin{proof}
If $\tilde y=\ex^{\theta\imag} y$ ($\theta\in\R$), then clearly $\tilde z = 2 \tilde y^\dagger = 2 y^\dagger = z$, hence $\Im z = \Im \tilde z$.
For the converse, let $A$ be the invertible $2\times 2$ matrix  with real entries and~$Ay=\tilde y$; here~$y=(y_1, y_2)^{t}$ and
$\tilde y=(\tilde y_1, \tilde y_2)^{t}$, column vectors with entries in $\c^2$.  As in the proof of~[ADH, 4.1.18], 
$\wr(y_1, y_2)=\wr(\tilde y_1, \tilde y_2)$ yields $\det A=1$.

Suppose $\Im \tilde z=u$, so~$y_1^2+y_2^2 = \tilde y_1^2+\tilde y_2^2$.  Choose~$a\in\R$ and representatives for $u$, $y_1$, $y_2$, $\tilde y_1$, $\tilde y_2$ in~$\c_a$, denoted by the same symbols, such that in~$\c_a$ we have $Ay=\tilde y$  and $y_1^2+y_2^2 = \tilde y_1^2+\tilde y_2^2$, and $u(t)\cdot \big(y_1(t)^2+y_2(t)^2\big)=1$ for all $t\ge a$. With $\dabs{\, \cdot\, }$  the usual euclidean norm on $\R^2$, 
we then have $\dabs{Ay(t)}=\dabs{y(t)}=1/\sqrt{u(t)}$ for~$t\geq a$.
Since~$f/4$ generates oscillations, we have
$\phi> \R$, and we conclude that~$\dabs{Av}=1$ for all $v\in\R^2$ with $\dabs{v}=1$. It is well-known that then $A=\left(\begin{smallmatrix} \cos \theta & -\sin\theta \\ \sin\theta & \cos\theta\end{smallmatrix}\right)$ with~$\theta\in\R$  (see, e.g., \cite[Chapter~XV, Exercise~2]{Lang}), so
$\tilde y=\ex^{\theta\imag}y$.
\end{proof}

\noindent
The observations above will be used in the proofs of Theorems~\ref{upo} and~\ref{thm:upo-freeness of the perfect hull} below.
We finish with miscellaneous historical remarks  (not used later):

 \begin{remarks}
The connection between the second-order linear differential equation \eqref{eq:2nd order, 2} and the third-order
non-linear differential equation~$\sigma(y)=f$    
%\marginpar{skipped remarks} 
was first investigated by Kummer~\cite{Kummer} in 1834. Appell~\cite{Appell}
noted that the linear differential equation
$$Y'''+fY'+(f'/2)Y\ =\ 0$$ has $\R$-linearly independent solutions~$y_1^2, y_1y_2, y_2^2\in\Calinf$, though some cases were known earlier~\cite{Clausen,Liouville39}; in particular, $1/u=y_1^2+y_2^2$ is a solution.  See also Lemma~\ref{lem:Appell}.
Hartman~\cite{Hartman61, Hartman73} investigates monotonicity properties of $y_1^2+y_2^2$.
Steen~\cite{Steen} in 1874, and independently Pinney~\cite{Pinney}, remarked that $r:=1/\sqrt{u}=\sqrt{y_1^2+y_2^2}\in\Calinf$ satisfies $4r''+fr=1/r^3$. (See also~\cite{Redheffer}.)
\end{remarks}

\subsection*{Complex exponentials over Hardy fields}
We now use some  of the above  to prove an extension theorem for Hardy fields (cf. \cite[Lem\-ma~11.6(6)]{Boshernitzan82}):

\begin{prop}\label{prop:cos sin infinitesimal, 1}
If $\phi\in H$ and $\phi\preceq 1$, then $\cos \phi,\sin\phi\in\Dx(H)$. 
\end{prop}
\begin{proof}
Replacing $H$ by $\Dx(H)$ we  arrange $\Dx(H)=H$.
Then by Pro\-po\-si\-tion~\ref{prop:Hardy field exts}, $H\supseteq\R$ is a Liouville closed $H$-field, and
by Corollary~\ref{cor:omega(H) downward closed}, $\omega(H)$ is downward closed.
Hence by Lem\-ma~\ref{lem:sc=>tc}, $H$ is trigonometrically closed.
Let now $\phi\in H$ and~${\phi\preceq 1}$.  Then $(\ex^{\phi\imag})^\dagger=\phi'\imag\in K^\dagger$, so
$\cos\phi+\imag\sin\phi=\ex^{\phi\imag}\in K$ using $K\supseteq\C$.  Thus~$\cos\phi,\sin\phi\in H$.
\end{proof}

%\noindent
%First we show:

%\begin{lemma}\label{lemlemsincos}
%Let $\phi\in \Calinf$ with $h:=\phi'\in \c^\times$, and  
%$A:=\der^2 - h^\dagger \der + h^2\in \Calinf[\der]$.
%Then $\big\{f\in \Calinf:\, A(f)=0\big\}\ =\ \R\cos \phi+\R\sin \phi$.
%\end{lemma}
%\begin{proof}
%We have
%$$(\cos \phi)'\ =\  -h\sin \phi,\qquad (\cos \phi)''\ =\ -h'\sin \phi - h^2\cos \phi = h^\dagger (\cos \phi)' - h^2 \cos \phi,$$
%so $A(\cos \phi)=0$. Likewise, $A(\sin \phi)=0$.
%The Wronskian of $\cos\phi$, $\sin\phi$ equals $h$,
%% $\cos\phi(\sin\phi)'-(\cos\phi)'\sin\phi=h$, 
%hence~$\cos\phi$,~$\sin\phi$
%are $\R$-linearly independent.
%\end{proof}

%\begin{proof}[Proof of Proposition~\ref{prop:cos sin infinitesimal, 1}]
%The case $\phi\in \R$ is trivial, so assume $\phi\in H\setminus \R$, $\phi\preceq 1$. Then
%$h:=\phi'\in\I(H)^{\ne}$; we arrange $h > 0$.  
%Now
%$f:=\omega(-h^\dagger)+4h^2=\sigma(2h)$ with~$2h\in H^>\cap \I(H)$ and so $2h\in \Dx(H)^>\setminus\Upg\big(\Dx(H)\big)$
%by [ADH, 11.8.19] and~$\Dx(H)$ being Liouville closed (Proposition~\ref{prop:Hardy field exts}).  So $f\in \bar\omega(H)$ by~[ADH, 11.8.31] and Corollary~\ref{cor:omega(H) downward closed}.
%Then $\cos \phi,\sin \phi\in\Dx(H)$ by Corollary~\ref{cor:char osc} and Lemma~\ref{lemlemsincos}.
%\end{proof}

\begin{cor}\label{cor:cos sin infinitesimal, 1}
Let  $\phi\in H$ and $\phi\preceq 1$. Then $\cos \phi$, $\sin\phi$ generate a $\d$-algebraic Hardy field extension 
$E:=H(\cos\phi,\sin\phi)$  of $H$. If
$H$ is a $\Ginf$-Hardy field, then so is $E$, and likewise with $\Gom$ in place of~$\Ginf$. 
\end{cor}

%\noindent
%Recall: $\cos (\phi + \theta)=\cos(\phi)\cos(\theta)-\sin(\phi)\sin(\theta)$ for $\phi,\theta\in \R$. So for any $a,b\in \R$ there is
%$d\in \R$ such that $a\cos(t) + b\sin(t)=\sqrt{a^2+b^2}\cdot \cos(t+d)$ for all $t\in \R$. 
%For later use we record the following consequence:

\noindent
Recall that for $\phi,\theta\in \R$ we have
\begin{align*} \cos(\phi+\theta)\  &=\ \cos(\phi)\cos(\theta)-\sin(\phi)\sin(\theta), \\  
\cos(\phi-\theta)\  &=\ \cos(\phi)\cos(\theta)+\sin(\phi)\sin(\theta).
\end{align*} 
Recall also the bijection $\arccos\colon [-1,1]\to [0,\pi]$, the inverse of the cosine function on $[0,\pi]$. It follows that for any $a,b\in \R$ we have $d\in \R$ such that $$a\cos(\phi)+b\sin(\phi)\ =\ \sqrt{a^2+b^2}\cdot \cos(\phi+d) \text{ for all $\phi\in \R$:}$$ 
for $a$, $b$ not both $0$ this holds
with $d=\arccos\!\big(a/\sqrt{a^2+b^2}\big)$ when $b\leqslant 0$, and with~$d=-\arccos\!\big(a/\sqrt{a^2+b^2}\big)$
when $b\geqslant 0$. For later use we record some consequences:

\begin{lemma}[Addition of sinusoids]\label{addsin} 
Let $y\in\c$. Then
$$y=a\cos x+b\sin x \text{ for some $a,b\in\R$} \quad \Longleftrightarrow\quad y=c\cos(x+d) \text{ for some $c,d\in\R$.}$$ 
\end{lemma}
%\begin{proof}
%Suppose $y=a\cos x+b\sin x$. Then with $z=a-b\imag\in\mathbb C$ we have
%$y=\Re(z\ex^{x\imag})$; taking $c,d\in\R$ with $z=c\ex^{d\imag}$ we then have
%$y=\Re(c\ex^{(x+d)\imag})=c\cos(x+d)$. Conversely, if $c,d\in\R$
%satisfy $y=c\cos(x+d)$, then $y=\Re(z\ex^{x\imag})$ for $z=c\ex^{d\imag}$, so $y=a\cos x+b\sin x$ for $a=\Re z$, $b=-\Im z$.
%\end{proof}

\begin{cor}\label{cor:sinusoids}
Let $\phi\in\c$. Then  
$\R \cos\phi + \R \sin\phi\, =\,\big\{  c\cos(\phi+d) :\, c,d\in\R \big\}$.
\end{cor}

\begin{cor}\label{arccosH} Suppose $H\supseteq \R$ is real closed and closed under integration, and let $g,h\in H$. Then there is $u\in H$ such that $-\pi \leqslant u \leqslant \pi$ and $g\cos \phi + h\sin \phi = \sqrt{g^2+h^2}\cdot \cos (\phi+ u)$ for all $\phi\in \c$:  if $h < 0$ this holds for
$u=\arccos \big(g/\sqrt{g^2+h^2}\big)$, and if $h>0$ it holds for $u=-\arccos \big(g/\sqrt{g^2+h^2}\big)$.
\end{cor} 
\begin{proof} On the interval $(-1,1)$ the function $\arccos$ is real analytic with derivative~$t\mapsto -1/\sqrt{1-t^2}$. 
Thus $\arccos\big(g/\sqrt{g^2+h^2}\big)\in H$ for $h\ne 0$. \end{proof}

\begin{cor} \label{gphicos} Let $a\in \R$ and let $g,\phi\in \c_a^1$ have germs in $H$ such that $g(t)\ne 0$ eventually, and
$\phi(t)\to +\infty$ as $t\to +\infty$.  Then there is a real $b\geqslant a$ with the property that if
$s_0, s_1\in [b,+\infty)$ with $s_0 < s_1$ are any successive zeros of $y:= g\cos\phi$, then~$y'$ has exactly one zero in the interval $(s_0,s_1)$. 
\end{cor}
\begin{proof} By increasing $a$ we arrange $g(t)\ne 0$ and $\phi'(t)>0$ for all $t\ge a$. Replacing~$g$ by~$-g$ if necessary we
further arrange $g(t)>0$ for all $t\geqslant a$. Let $s_0, s_1\in [a,+\infty)$ with $s_0<s_1$ be successive zeros of $y$. Later we impose a suitable lower bound $b\geqslant a$ on $s_0$.  Then $\phi(s_1)=\phi(s_0)+\pi$, since $s_1$ is the next zero of $\cos \phi$ after $s_0$. Also
\begin{align*} y'\ &=\ g'\cos\phi - g\phi'\sin \phi\  =\  \sqrt{g'^2+(g\phi')^2}\cos(\phi+u),\ \text{ where}\\
u\ &=\ \arccos\big(g'/\sqrt{g'^2+(g\phi')^2}\big),\ \text{ so }0 < u(t)< \pi \text{ for all $t\geqslant a$.}
\end{align*} 
By Rolle, $y'$ has a zero in $(s_0, s_1)$. Let $t\in (s_0, s_1)$ be a zero of $y'$. Then 
$$\phi(s_0) < \phi(t)< \phi(s_0)+ \pi,\quad \phi(t)+u(t)\in\phi(s_0)+ \Z\pi,$$
so $\phi(t)+ u(t)=\phi(s_0)+ \pi$.  Take $b\geqslant a$ in $\R$ so large that $u$ is differentiable on $[b,+\infty)$ and
$\phi'(t)+u'(t)>0$ for all $t\geqslant b$; this is possible because $u\preccurlyeq 1$ is $H$-hardian by Corollary~\ref{arccosH}, and
$\phi(t)+u(t)\to +\infty$ as $t\to +\infty$. Assuming now that $b\leqslant s_0$, we conclude that $t\in (s_0, s_1)$ is uniquely determined by $\phi(t)+u(t)=\phi(s_0)+\pi$. 
\end{proof}

\noindent
The $H$-asymptotic field extension $K:=H[\imag]$ of~$H$ is a differential subring of $\Calinf[\imag]$.  
To handle ultimate dents in $H$ in Section~\ref{sec:ultimate}, we
sometimes assumed $\I(K)\subseteq K^\dagger$, a condition that we consider more closely in 
the next proposition:

\begin{prop}\label{prop:cos sin infinitesimal, 2}
Suppose $H\supseteq\R$ is closed under integration.  
Then the following conditions are equivalent:
\begin{enumerate}
\item[\textup{(i)}] $\I(K)\subseteq K^\dagger$;
\item[\textup{(ii)}] $\ex^{f}\in K$ for all $f\in K$ with $f\prec  1$;
\item[\textup{(iii)}] $\ex^\phi,\cos \phi, \sin\phi\in H$ for all $\phi\in H$ with $\phi\prec 1$.
\end{enumerate}
\end{prop}
\begin{proof}
Assume (i), and let $f\in K$, $f\prec  1$.
Then $f' \in \I(K)$, so we have~$g\in K^\times$ with~$f'=g^\dagger$
and thus $\ex^{f}=cg$ for some $c\in\mathbb C^\times$. Therefore~$\ex^{f}\in K$.
This shows~(i)~$\Rightarrow$~(ii), and (ii)~$\Rightarrow$~(iii) is clear.
Assume (iii), and let $f\in \I(K)$. Then~$f=g+h\imag$, $g, h\in \I(H)$. Taking 
$\phi, \theta\prec 1$ in~$H$
with~$\phi'=g$ and $\theta'=h$,   
$$\exp(\phi + \theta\imag)\ =\ \exp(\phi)\big(\!\cos (\theta) + \sin(\theta)\imag\big)\in H[\imag]\ =\ K$$
has the property that $f=\big(\!\exp(\phi+\theta\imag)\big){}^\dagger\in K^\dagger$. 
This shows (iii)~$\Rightarrow$~(i).
\end{proof}

\noindent
From Propositions~\ref{prop:cos sin infinitesimal, 1} and \ref{prop:cos sin infinitesimal, 2} we obtain:

\begin{cor}\label{cor:cos sin infinitesimal}
If $H$ is $\d$-perfect, then   $\I(K)\subseteq K^\dagger$. 
\end{cor}

\noindent
Next we consider ``polar coordinates'' of nonzero elements of $K$:

\begin{lemma}\label{lem:arg}
Let $f\in\c[\imag]^\times$. Then $\abs{f}\in\c^\times$, and there exists $\phi\in\c$ with~$f=\abs{f}\ex^{\phi\imag}$;
such $\phi$ is unique up to addition of
an element of $2\pi\Z$. If also $f\in\c^r[\imag]^\times$, $r\in\N\cup\{\infty,\omega\}$, then $\abs{f}\in \c^r$
and $\phi\in\c^r$ for such $\phi$.
\end{lemma}
\begin{proof}
The claims about  $\abs{f}$ are clearly true. To show existence of $\phi$ we may replace~$f$ by $f/\abs{f}$
to arrange $\abs{f}=1$.
Take $a\in\R$ and a representative of~$f$ in~$\c_{a}[\imag]$, also denoted by $f$, such that $\abs{f(t)}=1$ for all $t\geq a$.
The proof of \cite[(9.8.1)]{Dieudonne} shows that for $b\in (a,+\infty)$ and $\phi_a\in\R$ with
$f(a)=\ex^{\phi_a\imag}$
there is a unique continuous function
$\phi\colon [a,b]\to\R$ such that $\phi(a)=\phi_a$ and $f(t) = \ex^{\phi(t)\imag}$ for all $t\in [a,b]$, and if also $f|_{[a,b]}$ is of class
$\c^1$, then so is this $\phi$  with $\imag\phi'(t) = f'(t)/f(t)$ for all~$t\in[a,b]$. 
With $b\to +\infty$ this yields the desired result. 
% Applying this fact successively with $a=a_0+n$ and $b=a+1$ we obtain a.
%$\phi\in\c$ as required.
\end{proof}

\begin{lemma}\label{lem:fexphii}
Suppose $H\supseteq\R$ is Liouville closed  and $f\in \c^1[\imag]^\times$. Then~$f^\dagger\in K$ iff 
$\abs{f}\in H^>$ and  $f=\abs{f}\ex^{\phi\imag}$ for some $\phi\in H$.
If in addition $f\in K^\times$, then $f=\abs{f}\ex^{\phi\imag}$ for some $\phi\preceq 1$ in $H$.
\end{lemma}
\begin{proof}
Take $\phi\in\c$ as in Lemma~\ref{lem:arg}. Then $\phi\in\c^1$ and $\Re f^\dagger=\abs{f}^\dagger$, $\Im f^\dagger=\phi'$.
If $f\in K^\times$, then
the remarks preceding Lemma~\ref{lem:W and I(F)} give $\phi'\in \I(H)$, so $\phi\preceq 1$. 
\end{proof}
 
\begin{cor}\label{cor:fexphii}
Suppose $H\supseteq\R$ is Liouville closed with $\I(K)\subseteq K^\dagger$. Let $L$ be a differential subfield of $\Calinf[\imag]$ containing $K$.
Then $L^\dagger \cap K = K^\dagger$.
\end{cor}
\begin{proof}
Let $f\in L^\times$ satisfy $f^\dagger \in K$. Then $f=\abs{f}\ex^{\phi\imag}$ with $\abs{f}\in H^>$ and $\phi\in H$,
by Lemma~\ref{lem:fexphii}. Hence $\ex^{\phi\imag}, \ex^{-\phi\imag}\in L$ and so
$\cos \phi = \frac{1}{2}(\ex^{\phi\imag}+\ex^{-\phi\imag})\in L$. In particular, $\cos \phi$ does not oscillate, so $\phi\preceq 1$ and thus $f=|f|(\cos \phi + \imag \sin \phi)\in K$ by Proposition~\ref{prop:cos sin infinitesimal, 2}.
\end{proof}

%\noindent
%Next we consider ``polar coordinates'' of nonzero elements of $K$:

%\begin{lemma}\label{lem:fexphii}
%Suppose $H\supseteq\R$ is Liouville closed, and
%let $f\in K^\times$. Then there is $\phi\preceq 1$ in $H$ such that $f=\abs{f}\ex^{\phi\imag}$.
%\end{lemma}
%\begin{proof} 
%We have $f=\abs{f}\cdot s$, $s\in K^\times$, $\abs{s}=1$, and so $s^\dagger \in H\imag$ by Lemma~\ref{lem:logder}.
%Take~$\phi\in H$ with $\phi'\imag=s^\dagger$. Then $s=c\ex^{\phi\imag}$ where $c\in\mathbb C^\times$, $\abs{c}=1$.
%Now $c=\ex^{d\imag}$, $d\in\R$;  replacing $\phi$ by $\phi+d$ we   arrange $c=1$, hence
%$f=\abs{f}\ex^{\phi\imag}$.
%By the remarks before Lemma~\ref{lem:W and I(F)} we have $\phi'\imag=s^\dagger\in  \I(H)\imag$, hence $\phi\preceq 1$.
%\end{proof}

\begin{cor}\label{cor:osc => bded}
Let $\phi\in H$, and  suppose $\ex^{\phi\imag}\sim f$ with $f\in E[\imag]^\times$ for some Hardy field extension~$E$ of~$H$.
Then  $\phi\preceq 1$.
\end{cor}
\begin{proof}
We can assume that $E=H$ is Liouville closed and contains $\R$. Towards a contradiction assume $\phi\succ 1$.
Lemma~\ref{lem:fexphii} yields $\theta\preceq 1$ in $H$ such that $f=\abs{f}\ex^{\theta\imag}$. Then
$\ex^{(\phi-\theta)\imag}\sim \abs{f}$ and $\phi-\theta \sim \phi$. Thus replacing $f$, $\phi$ by $\abs{f}$, $\phi-\theta$, respectively,
we  arrange~$f\in H^\times$. Then $\ex^{\phi\imag} = \cos\phi + \imag \sin\phi \sim f$ in $\Calinf[\imag]$ gives $\cos \phi\sim f$, contradicting that
$\cos \phi$ has arbitrarily large zeros.
\end{proof}

\begin{cor}\label{cor:phi preceq 1}
Let $f\in K^\times$, $\phi\in H$, so $y:=f\ex^{\phi\imag}\in\Calinf[\imag]^\times$. Then the following are equivalent:
\begin{enumerate}
\item[\textup{(i)}] $\phi\preceq 1$;
\item[\textup{(ii)}]  $y\in\operatorname{D}(H)[\imag]$;
\item[\textup{(iii)}]  $y\in E[\imag]$ for some Hardy field extension $E$ of $H$;
\item[\textup{(iv)}]  $y\sim g$ for some Hardy field extension $E$ of $H$ and $g\in E[\imag]^\times$.
\end{enumerate}
\end{cor}
\begin{proof} Use Proposition~\ref{prop:cos sin infinitesimal, 2} and Corollaries~\ref{cor:cos sin infinitesimal} and~\ref{cor:osc => bded} to  obtain the chain of implications (i)~$\Rightarrow$~(ii)~$\Rightarrow$~(iii)~$\Rightarrow$~(iv)~$\Rightarrow$~(i). 
\end{proof}

\noindent
Finally, some observations about solutions to linear differential equations  involving trigonometric functions.
 
\begin{lemma}\label{lem:trig surj, complex} 
Let $A\in K[\der]^{\neq}$ and   $\phi\in H$. 
Then~$A(K\ex^{\phi\imag})\subseteq K\ex^{\phi\imag}$.
Moreover, if
$K$ is  $r$-linearly surjective with $r:= \order A$, or $K$ is $1$-linearly surjective and $A$ splits over $K$,
then $A(K\ex^{\phi\imag}) = K\ex^{\phi\imag}$.
\end{lemma} 
\begin{proof}
The  differential operator $B:=A_{\ltimes\!\ex^{\phi\imag}}=\ex^{-\phi\imag}A\ex^{\phi\imag}\in (\Calinf[\imag])[\der]$ of order $r$ has coefficients in $K$. This follows from extending [ADH, 5.8.8] by allowing the element $h$ there (which is $\ex^{\phi\imag}$ here) to be a unit in a differential ring extension of~$K$ instead of a nonzero element in a differential field extension of $K$; the proof of~[ADH, 5.8.8] goes through, {\em mutatis mutandis}, to give this extension. 
Thus if~$y\in K$, then~$A(y\ex^{\phi\imag}) = B(y)\ex^{\phi\imag}$.
Also, if~$A$ splits over $K$, then so does $B$.
Hence if
$K$ is  $r$-linearly surjective, or $K$ is $1$-linearly surjective and $A$ splits over $K$,
then for each~$b\in K$ we obtain~$y\in K$ with $B(y)=b$, and so $A(y\ex^{\phi\imag}) = b \ex^{\phi\imag}$.
\end{proof}
 
\begin{lemma}\label{lem:trig surj}
Let $A\in H[\der]^{\neq}$, and suppose
$K$ is  $r$-linearly surjective with $r:= \order A$, or $K$ is $1$-linearly surjective and $A$ splits over $K$. Let also $h,\phi\in H$. Then
there are $f,g\in H$ such that $A(f\cos\phi + g\sin \phi)=h\cos \phi$.
\end{lemma}
\begin{proof}
Lemma~\ref{lem:trig surj, complex}  gives $y\in K$ such that
$$A(y\ex^{\phi\imag})\ =\   h \ex^{\phi\imag}\ =\ (h \cos \phi)+(h\sin \phi)\imag.$$
Take $f,g\in H$ with $y=f-g\imag$. Then 
$$y\ex^{\phi\imag}\ =\ (f\cos \phi+g\sin \phi) + (-g\cos \phi+f\sin \phi)\imag$$ and hence
$A(f\cos \phi+g\sin \phi)\ =\ h\cos \phi$.
\end{proof}

\begin{lemma}\label{lemtrig}
Let $f,g\in K$, $\phi\in H$,  $\phi\succ 1$, and $f\cos \phi+g\sin\phi\in\C\subseteq \c[\imag]$. Then~$f=g=0$.
\end{lemma}
\begin{proof} 
Take $c\in \C$ such that $f\cos\phi+g\sin \phi=c$. Since $\phi\in H$ and $\phi\succ 1$, there are arbitrarily large $t$ with $\phi(t)\in 2\Z\pi$, so $f(t)=c$, and thus $f=c$. There are also arbitrarily large $t$ with
$\phi(t)\in (2\Z+1)\pi$, and this gives likewise $-f=c$, so~$f=c=0$. 
Hence $g\sin \phi=0$, which easily gives $g=0$.  
%Let $y:=\frac{1}{2}(f-g\imag)$. Then $y\ex^{\phi\imag}+\overline{y}\ex^{-\phi\imag}=f\cos\phi+g\sin\phi\in\C$, hence $z\ex^{\phi\imag}+\overline{z}\ex^{-\phi\imag}=0$ where $z:=y'+\phi'\imag y\in K$. If $y\neq 0$ then $z\neq 0$ since $K^\dagger\cap H\imag \subseteq\I(H)\imag$ by the remarks preceding Lemma~\ref{lem:W and I(F)}, hence $\ex^{2\phi\imag}=-\overline{z}/z\in K^\times$ and thus~$\phi\preceq 1$ by Corollary~\ref{cor:osc => bded}, a contradiction.
\end{proof}

\noindent
Combining Lemmas~\ref{lem:trig surj} and~\ref{lemtrig} gives:

\begin{cor}\label{cor:unique antider}
If $K$ is $1$-linearly surjective, and $h,\phi\in H$, ${\phi\succ 1}$, then there are unique $f,g\in H$ such that $(f\cos\phi+g\sin\phi)'=h\cos\phi$. 
\end{cor}

\subsection*{Behavior of $\sigma$ and $\omega$ under composition} 
In this subsection we fix $\ell\in\c^1$ with~$\ell>\R$ and $\phi:=\ell'\in H$, so $\phi>0$. We use the superscript $\circ$ as in the subsection
on compositional conjugation in Hardy fields of Section~\ref{sec:Hardy fields}.
We refer to~[ADH, 11.8] (or Section~\ref{sec:uplupo-freeness}) for the definition of the sub\-sets~$\Upg(H)$, $\Upl(H)$, and~$\Upd(H)$
of $H$. The  bijection
$$y\mapsto (y/\phi)^\circ\ \colon\  H\to H^\circ$$ restricts to  bijections
$\I(H)\to\I(H^\circ)$ and $\Upg(H)\to\Upg(H^\circ)$, and the
  bijection
$$z\mapsto  \big( (z+\phi^\dagger) /\phi\big){}^\circ\ \colon\ H\to H^\circ$$
restricts to bijections $\Upl(H)\to\Upl(H^\circ)$ and $\Upd(H)\to\Upd(H^\circ)$.
(See the   transformation formulas in [ADH, p.~520].)
Consider   the bijection
$$f\mapsto \Phi(f)\ :=\ \big(\big(f-\omega(-\phi^\dagger)\big)/\phi^2\big){}^\circ\ \colon\ H\to H^\circ.$$ 
Then for $y\in H^\times$, $z\in H$ we have
$$\sigma\big((y/\phi)^\circ\big)\ =\  \Phi\big(\sigma(y)\big),\qquad \omega\big( \big( (z+\phi^\dagger) /\phi\big){}^\circ \big)\  =\  \Phi\big(\omega(z)\big).$$
(See the formulas in [ADH, pp.~518--519].) 
Hence $\Phi$ restricts to  bijections  
$$\sigma(H^\times)\to\sigma\big((H^\circ)^\times\big),\quad  
\sigma\big(\I(H)^{\ne}\big)\to \sigma\big(\I(H^\circ)^{\ne}\big), \quad \sigma\big(\Upg(H)\big)\to \sigma\big(\Upg(H^\circ)\big),$$
and
$$\omega(H)\to\omega(H^\circ),\quad \omega\big(\Upl(H)\big)\to \omega\big(\Upl(H^\circ)\big), \quad 
\omega\big(\Upd(H)\big)\to \omega\big(\Upd(H^\circ)\big).$$

\subsection*{An example of compositional conjugation}    Which ``changes of variable'' preserve the general form of the linear
differential equation \eqref{eq:2nd order, 2}? The next lemma and Corollary~\ref{cor:chvar} below give an answer.

\begin{lemma}\label{chvar} 
Let   $K$ be a differential field, $f\in K$, and   $P(Y):= 4Y''+fY$.
Then for $g\in K^\times$ and $\phi:= g^{-2}$ we have
$$ g^3P_{\times g}^\phi(Y)\ =\ 4Y'' + g^3P(g)Y.$$
\end{lemma} 
\begin{proof} Let $g,\phi\in K^\times$. Then 
\begin{align*} P_{\times g}(Y)\ &=\ 4gY'' + 8g'Y' + (4g'' + fg)Y\ =\ 4gY'' + 8g'Y' + P(g)Y, \quad \text{so}\\
 P_{\times g}^\phi(Y)\ &=\ 4g(\phi^2 Y'' + \phi'Y') + 8g'\phi Y' + P(g)Y\\ 
 &=\ 4g\phi^2 Y'' + (4g\phi' +8g'\phi)Y' + P(g)Y.\end{align*}
Now $4g\phi' + 8g'\phi=0$ is equivalent to $\phi^{\dagger}=-2g^\dagger$, which holds for $\phi=g^{-2}$.
For this~$\phi$ we get $P_{\times g}^\phi(Y)=
g^{-3}\big(4Y'' + g^3P(g)Y\big)$, that is, $g^3P_{\times g}^\phi(Y) = 4Y'' + g^3P(g)Y$. 
\end{proof}

\noindent
Now let  $\ell\in \Go$ be such that $\ell> \R$
and~$\phi:=\ell'\in H$, and let
 $P := 4Y''+fY$ where~$f\in H$. Note that if $y\in \Gt[\imag]$ and   $4y'' + fy=0$,  then $y\in \Calinf[\imag]$. Towards using Lemma~\ref{chvar}, 
suppose $\phi=g^{-2}$, $g\in H^\times$. Using notation from the previous subsection we set $h:= (g^3P(g) )^\circ\in H^\circ$
to obtain the following reduction of solving the differential equation  \eqref{eq:2nd order, 2} to solving a similar equation over $H^\circ$:

\begin{cor}\label{cor:chvar}
Let $y\in\Gt[\imag]$. Then $z:=(y/g)^\circ\in \Gt[\imag]$, and  
$$4y''+fy\ =\ 0\quad \Longleftrightarrow\quad 
4z''+h z\ =\ 0.$$
\end{cor}

\noindent
In particular, $f/4$ generates oscillations iff 
$h/4$
does. In connection with the formulas in the previous subsection, note that
$$g^3P(g)\  =\  g^3(4g''+fg)\ =\ \big(f-\omega(-\phi^\dagger)\big)/\phi^2,$$ so $h=\Phi(f)$.

\begin{lemma}\label{lem:baromega comp conj}
The increasing bijection  $$f\mapsto \Phi(f)\ =\ \big(\big(f-\omega(-\phi^\dagger)\big)/\phi^2\big){}^\circ\ \colon\ H\to H^\circ$$ maps $\bar{\omega}(H)$ onto $\bar{\omega}(H^\circ)$.
\end{lemma}
\begin{proof}
First replace $H$ by its real closure to arrange that $H$ is real closed,   then take~$g\in H^\times$ with $g^{-2}=\phi$, and use
the remarks following Corollary~\ref{cor:chvar}. 
\end{proof}

\noindent
We use the above to prove the Fite-Leighton-Wintner oscillation criterion for self-adjoint second-order linear differential
equations over $H$ \cite{Fite,Leighton50,Wintner49}. (See also \cite[\S{}2]{Hinton} and \cite[p.~45]{Swanson}.)
Let $A\in H[\der]$ be self-adjoint of order~$2$. Then
  $A=f\der^2+f'\der+g$ where~$f,g\in H$, $f\neq 0$, by the example following
Lemma~\ref{lem:Jacobi, 1}.
For $h\in\c$,   let $\int h$ denote a germ in $\c^1$ with~$(\int h)'=h$.

\begin{cor}\label{cor:FLW}
Suppose~$\int f^{-1}>\R$ and $\int g>\R$. Then $A(y)=0$ for some oscillating $y\in\Calinf$.
\end{cor}
\begin{proof}
We arrange that $H\supseteq\R$ is  Liouville closed. Then $f^{-1},g\in\Upg(H)$ by [ADH, 11.8.19].   Note that  $\phi:=f^{-1}$ is active in $H$.
 Put~$B:=4\phi A_{\ltimes \phi^{1/2}}$, so $B=4\der^2+h$ with~$h:=\omega(-\phi^\dagger)+4g\phi$. 
Then~$A(y)=0$ for some oscillating~$y\in\Calinf$ iff~${B(z)=0}$ for some oscillating $z\in\Calinf$ iff $h\notin\bar{\omega}(H)$, by
Corollary~\ref{cor:char osc}.   The latter is equivalent to~$(4g/\phi)^\circ\notin\bar{\omega}(H^\circ)$, by Lemma~\ref{lem:baromega comp conj} applied to $h$ in place of $f$.  Now~$\Upg(H^\circ)\cap\bar{\omega}(H^\circ)=\emptyset$ by Lemma~\ref{lem:Upgcapbaromega}, so it remains to note that~$4g\in\Upg(H)$ yields~$(4g/\phi)^\circ\in\Upg(H^\circ)$, by remarks in the previous subsection.
\end{proof}

\subsection*{More about $\bar\omega(H)$\astr} 
For later use  (in particular, in Section~\ref{sec:perfect applications}) we study here
the downward closed subset $\bar\omega(H)$ of $H$ in more detail. 
Recall that~$\omega(H)\subseteq\bar{\omega}(H)$, with equality 
for $\d$-perfect~$H$. (Corollary~\ref{cor:omega(H) downward closed}.)
In [ADH, 16.3] we introduced the concept of a {\em $\HLO$-cut}\/ in a pre-$H$-field; every
pre-$H$-field has exactly one or exactly two $\HLO$-cuts~[ADH, remark before 16.3.19].
By~[ADH, 16.3.14, 16.3.16]:\index{LambdaOmega-cut@$\HLO$-cut}

\begin{lemma}\label{lem:d-perfect HLO-cut}
Suppose $H$ is $\d$-perfect.  Then $\big(\I(H), \Upl(H), \bar{\omega}(H) \big)$ is a $\HLO$-cut in~$H$, and this is the unique $\HLO$-cut in $H$ iff $H$ is $\upo$-free.
\end{lemma}

\noindent
%If $E$ is a Hardy field extension of $H$, then $\bar{\omega}(H)=\bar{\omega}(E)\cap H$.
Thus in general, %by Lemma~\ref{lem:d-perfect HLO-cut}, 
$$\big(\I\!\big(\!\Dx(H)\big)\cap H, \, \Upl\big(\!\Dx(H)\big)\cap H,\, \bar{\omega}(H) \big)$$   is a $\HLO$-cut  in $H$, and  hence
(see~[ADH, p.~692]):
$$\omega(H)^\downarrow\ \subseteq\ \bar{\omega}(H)\ \subseteq\   H\setminus \sigma\big(\Upg(H)\big){}^\uparrow.$$ 
The classification of $\HLO$-cuts in $H$ from [ADH, 16.3] can  be used to narrow down the  possibilities for
$\bar{\omega}(H)$:

\begin{lemma}\label{lem:baromega(H) for H without as int} 
Let  $\phi\in H^{>}$ be such that $v\phi\notin (\Gamma_H^{\neq})'$. Then 
$$\bar{\omega}(H)\  =\  \omega(-\phi^\dagger) + \phi^2 \smallo_H^\downarrow \quad\text{or}\quad
\bar{\omega}(H)\  =\  \omega(-\phi^\dagger) + \phi^2 \mathcal O_H^\downarrow.
$$
The first alternative holds if $H$ is grounded, and the second alternative holds if 
$v\phi$ is a gap in  $H$ with $\phi\asymp b'$ for some $b\asymp 1$ in $H$.
\end{lemma}
\begin{proof} Either $v\phi=\max \Psi_H$ or $v\phi$ is a gap in $H$, by [ADH, 9.2]. 
The remark before the lemma yields  an $\HLO$-cut $(I,\Upl,\Upo)$ in $H$ where~$\Upo=\overline{\omega}(H)$. Now use the proofs of~[ADH, 16.3.11, 16.3.12, 16.3.13] together with the transformation formulas~[ADH, (16.3.1)] for $\HLO$-cuts.
\end{proof}

\noindent
By~[ADH, 16.3.15] we have:
 
\begin{lemma}\label{lem:baromega(H)}
If $H$ has   asymptotic integration and the set $2\Psi_H$ does  not have a  supremum in $\Gamma_H$, then
$$\bar{\omega}(H)\ =\ \omega\big(\Upl(H)\big){}^\downarrow\ =\ \omega(H)^\downarrow\quad\text{or}\quad \bar{\omega}(H)\ =\ H\setminus \sigma\big(\Upg(H)\big){}^\uparrow.$$
\end{lemma}

\begin{cor}\label{omuplosc}  
Suppose $H$ is $\upo$-free. Then 
$$\bar{\omega}(H)\ =\ \omega\big(\Upl(H)\big){}^\downarrow\ =\ \omega(H)^\downarrow \ = \  H\setminus \sigma\big(\Upg(H)\big){}^\uparrow.$$ 
\end{cor}
\begin{proof} By [ADH, 11.8.30] we have $\omega\big(\Upl(H)\big){}^\downarrow= \omega(H)^\downarrow=  H\setminus \sigma\big(\Upg(H)\big){}^\uparrow$. It follows from [ADH, 9.2.19] that $2\Psi_H$ has no  supremum in $\Gamma_H$. Now use Lemma~\ref{lem:baromega(H)}.
\end{proof}

\noindent
In the next lemma  $L\supseteq\R$ is a Liouville closed $\d$-algebraic Hardy field extension of~$H$ such that $\omega(L)=\bar{\omega}(L)$.
(By Corollary~\ref{cor:omega(H) downward closed}, this holds for $L=\Dx(H)$.)
 Note that then $\bar{\omega}(L)=\omega\big(\Upl(L)\big)$ by [ADH, 11.8.20].
 
\begin{lemma}\label{lem:D(H) upo-free, 4} 
If $H$ is not $\upl$-free or $\bar{\omega}(H)=H\setminus\sigma\big(\Upg(H)\big){}^\uparrow$, then~$L$ is $\upo$-free.
\end{lemma}
\begin{proof}
If $H$ is $\upo$-free or not $\upl$-free, then $L$ is $\upo$-free by Lemmas~\ref{lem:Li(H) upo-free} and~\ref{lem:D(H) upo-free, 3}. Suppose
$H$ is $\upl$-free but not $\upo$-free, and $\bar{\omega}(H)=H\setminus\sigma\big(\Upg(H)\big){}^\uparrow$. So [ADH, 11.8.30] gives $\upo\in H$ with~$\omega\big(\Upl(H)\big) <  \upo  < \sigma\big(\Upg(H)\big)$. Then 
$\upo\in\bar{\omega}(H)\subseteq \bar{\omega}(L)\subseteq\omega\big(\Upl(L)\big)$. Thus $L$ is $\upo$-free by 
Corollary~\ref{cor:D(H) upo-free, 3}.
\end{proof}

\noindent
%As an application of our main theorem, we shall prove
Theorem~\ref{thm:upo-freeness of the perfect hull}, which depends on much of what follows, shows that for~$L=\Dx(H)$ the converse of Lemma~\ref{lem:D(H) upo-free, 4} also holds. 

\subsection*{Proof of a conjecture of Boshernitzan\astr} 
In this last subsection we establish~\cite[Conjecture~17.11]{Boshernitzan82}: Corollary~\ref{cor:17.11}. For this, 
with $\ell_n:=\log_n x$ we set~$\upg_n:=\ell_n^\dagger$, $\upl_n:=-\upg_n^\dagger$, and $\upo_n:=\omega(\upl_n)$, as   in [ADH, 11.5, 11.7], so
$$\upg_n\ =\ \frac{1}{\ell_0\ell_1\cdots\ell_n},\qquad  \upo_n\ =\ \frac{1}{\ell_0^2} + \frac{1}{\ell_0^2\ell_1^2} + \cdots + \frac{1}{\ell_0^2\ell_1^2\cdots\ell_n^2}.$$
 (See also the beginning of Section~\ref{sec:upo-free Hardy fields} below.)
For $c\in\R$,   the germ
$$\frac{\upo_n+c\upg_{n}^2}{4}\ =\ \frac{1}{4}\left(\frac{1}{\ell_0^2}+\frac{1}{(\ell_0\ell_1)^2}+\cdots+\frac{1}{(\ell_0\cdots\ell_{n-1})^2} + \frac{c+1}{(\ell_0\cdots\ell_{n})^2} \right)$$ 
generates oscillations iff $c>0$.
({A.~Kneser~\cite{AKneser}, Riemann-Weber~\cite[p.~63]{Weber}; cf.~\cite{Hille48}}.)
This follows  from the next corollary applied to $f=\upo_n+c\upg_n^2$ and
%Lemma~\ref{lem:baromega(H) for H without as int}
 the grounded Hardy subfield~$H:=\R\langle\ell_n\rangle=\R(\ell_0,\dots,\ell_n)$ of $\Li(\R)$:

\begin{cor}\label{cor:17.7 generalized}
Let $H$ be a grounded Hardy field such that for some $m$ we have~${h\succ \ell_m}$ for all~$h\in H$ with $h\succ 1$. Then for   $f\in H$, the
following are equivalent:
\begin{enumerate}
\item[\textup{(i)}] $f\in\bar{\omega}(H)$; 
\item[\textup{(ii)}]  $f < \upo_n$ for some $n$; 
\item[\textup{(iii)}] there exists $c\in\R^>$ such that for all $n$ we have $f < \upo_n+c\upg_n^2$;
\item[\textup{(iv)}] $f < \upo_n+c\upg_n^2$ for all $n$ and all  $c\in\R^>$.
\end{enumerate}
\end{cor}
\begin{proof}
By [ADH, 10.3.2, 10.5.15] we may replace $H$ by $H(\R)$ to arrange~${H\supseteq\R}$.
By Lemma~\ref{lem:Li(H) upo-free}, $L:=\Li(H)$ is $\upo$-free. With $H_{\upo}$ as in the proof of that lemma,  one verifies easily
that for each $g\in H_\upo$ with~$g\succ 1$ there is an $m$ such that~${g\succ \ell_m}$.
Hence~$(\ell_n)$ is a logarithmic sequence in~$L$, in the sense of 
[ADH, p.~499]. Now the implication~(i)~$\Rightarrow$~(iv) follows from
Corollary~\ref{omuplosc}  and [ADH, 11.8.22], and~(iv)~$\Rightarrow$~(iii) is obvious.  
Since $0<\upg_{n+1}\prec\upg_n$ we obtain for $c\in\R^>$: 
$$\upo_{n+1}+c\upg_{n+1}^2\ =\  \upo_{n}+\upg_{n+1}^2+c\upg_{n+1}^2\ <\ \upo_n+ \upg_n^2\ =\  \sigma(\upg_n).$$
%\sigma(\upg_{n+1})=\upo_{n+1}+\upg_{n+1}^2=\upo_n+2\upg_{n+1}^2<\upo_n+c\upg_n^2$, since~.
In view of [ADH, 11.8.30, 11.8.21] and Corollary~\ref{omuplosc}, this yields~(iii)~$\Rightarrow$~(ii). 
Downward closedness of~$\bar{\omega}(H)$ implies (ii)~$\Rightarrow$~(i).
\end{proof}

\noindent
Using the above equivalence (i)~$\Leftrightarrow$~(ii)  we recover~\cite[Theorem~17.7]{Boshernitzan82}:

\begin{cor}
Suppose $f\in\c$ is hardian and $\d$-algebraic over $\R$. Then
$$ f \text{ generates oscillations}\ \Longleftrightarrow\  f>\upo_n/4 \text{ for all }n.$$
\end{cor}
\begin{proof} By Corollary~\ref{cor:Bosh 14.11} the Hardy field~$H:=\R\langle f\rangle$ satisfies the hypotheses of Corollary~\ref{cor:17.7 generalized}. Also, $f$ generates oscillations iff $4f\notin \bar{\omega}(H)$.
\end{proof}

\noindent
Using the above implication (iii)~$\Rightarrow$~(i) we obtain in the same way: 

\begin{cor}\label{cor:17.11}
Let $f\in\c$  be hardian and $\d$-algebraic over $\R$, and suppose there is a $c\in\R^>$ with~$f<\upo_n+c\upg_n^2$ for all $n$. Then $f/4$ does not generate oscillations.
\end{cor}

\noindent
In the beginning of this subsection we introduced the germs $\ell_n$, and so this may be a good occasion to
observe that the Hardy field $H=\R(\ell_0,\ell_1,\ell_2,\dots)$ they generate over $\R$ is $\upo$-free: since $H$ is ungrounded and $H$ is the union of the grounded Hardy subfields $\R(\ell_0,\dots, \ell_n)$, this follows from
[ADH, 11.7.15]. Thus the Hardy field~$\Li(\R)=\Li(H)$ is $\upo$-free as well.

\section{Maximal Hardy Fields are $\upo$-Free}\label{sec:upo-free Hardy fields} 

\noindent
In this section we discuss the fundamental property of $\upo$-freeness from [ADH] in the context of Hardy fields.
The main result is Theorem~\ref{upo}, from which it follows that every maximal Hardy field is $\upo$-free.
As an application of this theorem, we answer a question from~\cite{Boshernitzan86}.

\subsection*{The property of $\upo$-freeness for Hardy fields}
Let $H\supseteq \R$ be a Liouville closed Hardy field.
Note that then $x\in H$ and $\log f\in H$ for all $f\in H^{>}$.  To work with $\upo$-freeness for $H$ we introduce the ``iterated logarithms''~$\ell_{\rho}$; more precisely, transfinite recursion yields a  sequence
$(\ell_\rho)$ in $H^{>\R}$ indexed by the ordinals~$\rho$ less than some infinite limit ordinal $\kappa$ as follows:  $\ell_0=x$, and
$\ell_{\rho+1}:=\log \ell_\rho$; if $\lambda$ is an infinite limit ordinal
such that all $\ell_\rho$ with $\rho<\lambda$ have already been chosen,
then we pick $\ell_\lambda$ to be any element in $H^{>\R}$ such that $\ell_\lambda\prec \ell_\rho$ for
all $\rho<\lambda$, if there is such an $\ell_\lambda$, while if there is no
such $\ell_\lambda$, we put $\kappa:=\lambda$. From $(\ell_\rho)$ we obtain the sequences~$(\upg_\rho)$ in~$H^{>}$ and~$(\upl_\rho)$ in $H$ as follows:
$$
\upg_\rho\ :=\ \ell_{\rho}^\dagger, \qquad
\upl_\rho\ :=\  -\upg_\rho^\dagger\ =\ -\ell_\rho^\dagger{}^\dagger\ :=\ -(\ell_\rho^\dagger{}^\dagger).$$
% \quad b_\rho :=(1/l_\rho)'{}^\dagger=a_\rho+y_\rho,$$
Then $\upl_{\rho+1}=\upl_{\rho}+\upg_{\rho+1}$ and we have
\begin{align*} \upg_0\ &=\ \ell_0^{-1}, & \upg_1\ &=\ (\ell_0\ell_1)^{-1}, & \upg_2\ &=\ (\ell_0\ell_1\ell_2)^{-1},\\
\upl_0\ &=\ \ell_0^{-1}, & \upl_1\ &=\ \ell_0^{-1} + (\ell_0\ell_1)^{-1}, &
  \upl_2\ &=\ \ell_0^{-1} + (\ell_0\ell_1)^{-1} + (\ell_0\ell_1\ell_2)^{-1},
\end{align*}
and so on. Indeed, $v(\upg_\rho)$ is strictly increasing as a function of $\rho$
and is cofinal in~$\Psi_H=\big\{v(f^\dagger):f\in H,\ 0\neq f\nasymp 1\big\}$; we refer to
[ADH, 11.5, 11.8] for this and some of what follows.
Also, $(\upl_\rho)$ is a strictly increasing pc-sequence which is cofinal in~$\Upl(H)\subseteq H$. We recall here the relevant
descriptions  from~[ADH, 11.8]:
\begin{alignat*}{2}
\Upg(H) &\ =\ \big\{a^\dagger:\ a\in H,\ a\succ 1\big\} &&\ =\ \{b\in H:\ b>\upg_{\rho} \text{ for some $\rho$}\}, \\
\Upl(H)	&\ =\ - \Upg(H)^\dagger &&\ = \ \big\{ {-a^{\dagger\dagger}}:\ a\in H,\ a\succ 1\big\}.
\end{alignat*}
Here, $\Upg(H)\subseteq H^>$ is upward closed  and $\Upl(H)$ is downward closed, since $H$ is Liouville closed. 
The latter also gives that $H$ is
$\upl$-free, that is, 
$(\upl_{\rho})$ has no pseudolimit in $H$. The function $\omega\colon H \to H$ is strictly increasing on~$\Upl(H)$ and setting $\upo_{\rho}:= \omega(\upl_{\rho})$ we obtain a strictly increasing
pc-sequence $(\upo_{\rho})$ which is cofinal in~$\omega\big(\Upl(H)\big)=\omega(H)$:
$$ \upo_0\ =\ \ell_0^{-2},\qquad \upo_1\ =\ \ell_0^{-2} + (\ell_0\ell_1)^{-2},
\qquad \upo_2\ =\ \ell_0^{-2} + (\ell_0\ell_1)^{-2} + (\ell_0\ell_1\ell_2)^{-2},$$
and so on; see  [ADH, 11.7,~11.8] for this and some of what follows. Now $H$ being $\upo$-free is equivalent to $(\upo_\rho)$ having no pseudolimit in $H$.\index{Hardy field!omega-free@$\upo$-free}\index{H-asymptotic field@$H$-asymptotic field!omega-free@$\upo$-free} By 
[ADH, 11.8.30] the pseudolimits of $(\upo_{\rho})$
in $H$ are exactly the $\upo\in H$ such that $\omega(H) < \upo < \sigma\big(\Upg(H)\big)$. Also, $\sigma$ is strictly increasing on $\Upg(H)$. Thus $H$ is not
$\upo$-free if and only if there exists an $\upo\in H$ such that $\omega(H) < \upo < \sigma\big(\Upg(H)\big)$.

\begin{lemma}\label{lem:upo}
Let $\upg\in(\Go)^\times$, $\upg > 0$, and $\upl:=-\upg^\dagger$ with $\upl_\rho < \upl <\upl_\rho+\upg_\rho$ in $\c$, for all $\rho$.
Then $\upg_\rho > \upg > \upg_{\rho}/\ell_{\rho}=(-1/\ell_\rho)'$ in $\c$, for all $\rho$.
\end{lemma}
\begin{proof}
Pick $a\in \R$ (independent of $\rho$) and functions in $\c_a$ whose germs at $+\infty$ are the elements
$\ell_{\rho}$,~$\upg_{\rho}$,~$\upl_{\rho}$ of $H$; denote these functions also by $\ell_{\rho}$,~$\upl_{\rho}$,~$\upg_{\rho}$. From $\ell_{\rho}^\dagger=\upg_{\rho}$ and $\upg_{\rho}^\dagger=-\upl_{\rho}$
 in~$H$ we obtain $c_{\rho}, d_{\rho}>0$ such that for all
 sufficiently large $t\ge a$, 
$$ \ell_{\rho}(t)\ =\ c_{\rho}\exp\!\left[\int_a^t \upg_{\rho}(s)\,ds\right], \quad \upg_{\rho}(t)\ =\ d_{\rho}\exp\!\left[-\int_a^t\upl_{\rho}(s)\,ds\right]. $$
(How large is ``sufficiently large'' depends on $\rho$.) Likewise we pick functions in $\c_a$ whose germ at $+\infty$ are $\upg$, $\upl$, and also denote these functions by $\upg$, $\upl$. From $\upg^\dagger=-\upl$ in $H$ we obtain a real constant $d>0$ such that for all sufficiently large $t\ge a$, 
$$ \upg(t)\ =\ d\exp\!\left[-\int_a^t\upl(s)\,ds\right].$$ 
Also, $\upl_\rho < \upl <\upl_\rho+\upg_\rho$ yields  constants $a_{\rho}, b_{\rho}\in \R$ such that for all $t\ge a$
$$\int_a^t \upl_{\rho}(s)\,ds\ <\ a_{\rho}+ \int_a^t \upl(s)\,ds\ <\ b_{\rho} +\int_a^t \upl_{\rho}(s)\,ds + \int_a^t \upg_{\rho}(s)\,ds,$$
 which by applying $\exp(-*)$ yields that for all sufficiently large $t\ge a$, 
$$\frac{1}{d_{\rho}}\upg_{\rho}(t)\ >\ \frac{1}{\ex^{a_{\rho}}d}\upg(t)\ >\ \frac{c_{\rho}}{\ex^{b_{\rho}}d_{\rho}} \upg_{\rho}(t)/\ell_{\rho}(t).$$
Here the positive constant factors don't matter,
since the valuation of $\upg_{\rho}$ is strictly increasing and that of $\upg_{\rho}/\ell_{\rho}=(-1/\ell_{\rho})'$
is strictly decreasing with $\rho$. Thus 
for all~$\rho$ we have 
$\upg_{\rho} > \upg > \upg_{\rho}/\ell_{\rho}= (-1/\ell_{\rho})'$, in $\c$. 
\end{proof}
 
\noindent
We are now ready to prove a key result:

\begin{theorem}\label{upo} Every Hardy field has a $\d$-algebraic $\upo$-free Hardy field extension.
\end{theorem} 
\begin{proof} 
It is enough to show that every $\d$-maximal Hardy field is $\upo$-free. That reduces to showing that
every non-$\upo$-free Liouville closed Hardy field containing $\R$
has a proper $\d$-algebraic Hardy field extension. So assume $H\supseteq \R$ is a Liouville closed Hardy field and $H$ is not $\upo$-free. We shall construct a proper $\d$-algebraic Hardy field extension
of $H$. We have $\upo\in H$ such that 
$$\omega(H)\ <\  \upo\  <\ \sigma\big(\Upg(H)\big).$$ 
%%In fact, we construct such a $K$ with $\upo\in \sigma(\Upg(K))$. 
%Take $a\in \R$ such that $\upo$ is the germ of a function in 
%$\Cat$, this function also to be denoted by $\upo$.
%%Since $\upo\asymp \ell_0^{-2}=x^{-2}$ in $H$, we can
%%increase $a$, if necessary, and arrange that $0 < \upo(t) < 1$ on 
%%$[a,\infty)$ and $\upo$ is strictly
%%decreasing on $[a,\infty)$. 
With $\upo$ in the role of $f$ in the discussion following Corollary~\ref{cor:char osc, 3}, we have  $\R$-linearly independent
solutions $y_1, y_2\in \Gt$ of the differential equation~${4Y'' + \upo Y=0}$; in fact, $y_1,y_2\in\Gi$.
Then the complex solution $y=y_1+y_2\imag$ is a unit
of~$\Gi[\imag]$, and so we have 
$z:=2y^\dagger\in \Gi[\imag]$.  We shall prove that the elements $\upl:=\Re z$ and~$\upg:= \Im z$ of $\Gi$ generate a Hardy field extension $K=H(\upl, \upg)$ of~$H$
with~$\upo=\sigma(\upg)\in \sigma(K^\times)$. 
We can assume that $w:= y_1y_2'-y_1'y_2\in \R^{>}$, so~$\upg=2w/|y|^2\in(\Gi)^\times$ and $\upg>0$.

We have $\upo_{\rho} \leadsto \upo$, with $\upo-\upo_{\rho}\sim \upg_{\rho+1}^2$ by  [ADH,~11.7.1]. Till further notice we fix $\rho$ and set
$g_{\rho}:=\upg_{\rho}^{-1/2}$, so $2g_{\rho}^\dagger=\upl_{\rho}=-\upg_{\rho}^\dagger$. 
For $h\in H^\times$ we also have~$\omega(2h^\dagger)=-4h''/h$,
hence $P:= 4Y'' + \upo Y\in H\{Y\}$ gives
$$P(g_{\rho})\ =\  
g_{\rho}(\upo-\upo_{\rho})\ \sim\ g_{\rho}\upg_{\rho+1}^2,$$ and so with an eye towards using Lemma~\ref{chvar}:
$$g_{\rho}^3P(g_{\rho})\ \sim\ g_{\rho}^4\upg_{\rho+1}^2\ \sim\ \upg_{\rho+1}^2/\upg_{\rho}^2\ \asymp\ 1/\ell_{\rho+1}^2.$$
Thus with
$g:= g_{\rho}=\upg_{\rho}^{-1/2}$, $\phi=g^{-2}=\upg_{\rho}$ we  have $A_{\rho}\in \R^{>}$ such that
\begin{equation}\label{eq:bound}
g^3P_{\times g}^{\phi}(Y)\ =\ 4Y'' + g^3P(g)Y,\quad 
|g^3P(g)|\ \le\  A_{\rho}/\ell_{\rho+1}^2.
\end{equation}
From  $P(y)=0$ we get $P^{\phi}_{\times g}(y/g)=0$, that is,
$y/g\in \Gi[\imag]^\phi$ is a solution of~$4Y'' + g^3P(g)Y=0$, with $g^3P(g)\in H\subseteq \Gi$. 
Set $\ell:=\ell_{\rho+1}$, so $\ell'=\ell_{\rho}^\dagger=\phi$.
The subsection on compositional conjugation in Section~\ref{sec:Hardy fields}   yields the isomorphism~$h\mapsto h^\circ=h\circ\ell^{\inv}\colon H^{\phi} \to H^\circ$ of $H$-fields, where~$\ell^{\inv}$ is the compositional inverse of~$\ell$. Under this isomorphism the equation $4Y''+g^3P(g)Y=0$ corresponds to the equation
$$4Y'' + f_{\rho}Y\ =\ 0, \qquad f_{\rho}\ :=\ (g^3P(g))^\circ\in H^\circ\ \subseteq\ \Gi.$$ 
By Corollary~\ref{cor:chvar}, the equation $4Y'' + f_{\rho}Y=0$ has the ``real'' solutions
$$y_{j,\rho}\ :=\ (y_j/g)^\circ\in (\Gi)^\circ \ =\ \Gi \qquad(j=1,2),$$  and the ``complex''
solution $$y_\rho\ :=\ y_{1,\rho} + y_{2,\rho}\imag\ =\ (y/g)^\circ,$$ which is a unit of the ring
$\Gi[\imag]$. Set $z_{\rho}:=2 y_{\rho}^\dagger\in \Gi[\imag]$. 
The bound in \eqref{eq:bound} gives 
$|f_{\rho}|\ \le\  A_{\rho}/x^2$,
which by Corollary~\ref{cor:bound} yields  positive constants $B_{\rho}$, $c_\rho$ such that
$|z_{\rho}|\ \le\ B_{\rho}x^{c_\rho}$.  Using 
$(f^\circ)'=(\phi^{-1}f')^\circ$ for $f\in\Calinf[\imag]$ we obtain
$$z_{\rho}\ =\ 2\big((y/g)^\circ\big){}^\dagger \ =\ 2\big(\phi^{-1} (y/g)^\dagger\big)^\circ \ =\ \big((z-2g^\dagger)/\phi\big)^\circ$$
In combination with the  bound on $|z_\rho|$ this yields
\begin{align*} \left|\frac{z-2g^\dagger}{\phi}\right|\ &\le\ B_{\rho}\,\ell_{\rho+1}^{c_\rho}, \quad \text{hence}\\
|z- \upl_{\rho}|\ &\le\ B_{\rho}\,\ell_{\rho+1}^{c_\rho}\,\phi\ =\ 
  B_{\rho}\,\ell_{\rho+1}^{c_\rho}\,\upg_{\rho}, \quad \text{and so} \\ 
z\ &=\ \upl_{\rho}  + R_{\rho}\quad\text{where}\quad |R_{\rho}|\le B_{\rho}\,\ell_{\rho+1}^{c_\rho} \,\upg_{\rho}. 
\end{align*}
We now use this last estimate with $\rho+1$ instead of $\rho$, 
together with 
$$\upl_{\rho+1}\ =\ \upl_{\rho}+\upg_{\rho+1},\quad  \ell_{\rho+1}\upg_{\rho+1}\ =\ \upg_{\rho}.$$
This yields
\begin{align*} z &\, =\, \upl_{\rho} + \upg_{\rho+1} + R_{\rho+1}\\ & \qquad \text{with}\ 
 |R_{\rho+1}|\, \le\, B_{\rho+1}\,\ell_{\rho+2}^{c_{\rho+1}}\, \upg_{\rho+1} \, =\, B_{\rho+1}\big(\ell_{\rho+2}^{c_{\rho+1}}/\ell_{\rho+1}\big)\,\upg_{\rho},\\
\text{so}\quad z &\, =\, \upl_{\rho} + o(\upg_{\rho})\ \text{ that is, }z-\upl_{\rho}\prec \upg_{\rho}, \\
\text{and thus}\quad \upl\, &=\, \Re z\, =\, \upl_{\rho} + o(\upg_{\rho}), \quad \upg\, =\, \Im z\, \prec\, \upg_{\rho}.
\end{align*}
Now varying $\rho$ again, $(\upl_{\rho})$ is a strictly increasing divergent
pc-sequence in $H$ which is cofinal in $\Upl(H)$. By the above,
$\upl= \Re z$ satisfies $\Upl(H) < \upl < \Upd(H)$. This yields an ordered subfield 
$H(\upl)$ of $\Gi$, which by Lemma~\ref{ps1} is an immediate
valued field extension of $H$ with $\upl_{\rho} \leadsto \upl$.
Now $\upl=-\upg^\dagger$ (see discussion before Lemma~\ref{lem:Im z}), so Lemma~\ref{lem:upo} gives~$\upg_{\rho} > \upg > (-1/\ell_{\rho})'$ in $\Gi$, for all $\rho$. In view of Lemma~\ref{ps2} applied to $H(\upl)$,~$\upg$ in the role of~$K$,~$f$ this yields an ordered subfield $H(\upl, \upg)$
of $\Gi$. Moreover, $\upg$ is transcendental over $H(\upl)$ and $\upg$ 
satisfies the second-order differential equation
$2yy''-3(y')^2+y^4-\upo y^2=0$ over $H$ (obtained from the relation 
$\sigma(\upg)=\upo$ by multiplication with $\upg^2$). It follows that $H(\upl,\upg)$ is closed under the derivation of~$\Gi$, and hence~$H(\upl, \upg)=H\<\upg\>$ is a Hardy field that is $\d$-algebraic over $H$. 
\end{proof} 

\noindent
The proof also shows that every $\Ginf$-Hardy field
has an 
$\upo$-free $\d$-algebraic $\Ginf$-Hardy field extension, and the same with 
$\Gom$ instead of $\Gi$. In Section~\ref{sec:perfect applications} below we show that the perfect hull of an  $\upo$-free Hardy field  remains $\upo$-free (Lemma~\ref{lem:upo-freeness of the perfect hull, (ii)=>(iii)}), but that not every perfect Hardy field is $\upo$-free (Example~\ref{ex:counterex}). 

\subsection*{Improving Theorem~\ref{upo}\astr} 
{\it In this subsection   $H\supseteq\R$ is a Liouville closed Hardy field
and $\upo\in H$, $\upg\in (\c^2)^\times$ satisfy
$\omega(H) <  \upo < \sigma\big(\Upg(H)\big)$ and $\sigma(\upg)=\upo$.}\/
Lem\-ma~\ref{lem:Li(H) upo-free} leads to a more explicit version of Theorem~\ref{upo}: 

\begin{cor}\label{cor:upo}
The germ  $\upg$ generates a Hardy field extension~$H\langle\upg\rangle$ of $H$ with a gap $v\upg$,
and so $\operatorname{Li}\!\big(H\langle\upg\rangle\big)$ is an $\upo$-free Hardy field extension of $H$.
\end{cor}
\begin{proof}
Since $\sigma(-\upg)=\sigma(\upg)$, we may arrange $\upg>0$. The discussion before Lem\-ma~\ref{lem:Im z} with $\upo$, $\upg$ in the roles of $f$, $g$, respectively, yields
$\R$-linearly independent solutions $y_1,y_2\in\Calinf$ of the differential equation $4Y''+\upo Y=0$ with
Wronskian~$1/2$ such that~$\upg=1/(y_1^2+y_2^2)$.  
The  proof of Theorem~\ref{upo} shows that~$\upg$ generates a Hardy field extension~$H\langle\upg\rangle=H(\upl, \upg)$ of $H$. Recall that $v(\upg_{\rho})$ is strictly increasing as a function of $\rho$ and cofinal in $\Psi_H$; as $\upg \prec \upg_{\rho}$ for all $\rho$, this gives~$\Psi_H<v\upg$. Also $\upg> (-1/\ell_{\rho})'>0$ for all $\rho$ and
$v(1/\ell_{\rho})'$ is strictly decreasing as a function of $\rho$ and coinitial in $(\Gamma_H^>)'$, and so $v\upg < (\Gamma_H^{>})'$. Then
by [ADH, 13.7.1 and subsequent remark~(2) on p.~626], $v\upg$ is a gap in~$H\langle \upg\rangle$, so $\operatorname{Li}\!\big(H\langle\upg\rangle\big)$ is $\upo$-free by Lemma~\ref{lem:Li(H) upo-free}. 
\end{proof}

\begin{cor} 
Suppose $\upg>0$. Then with~$L:=\operatorname{Li}\!\big(H\langle\upg\rangle\big)$, 
$$\upo\notin\bar{\omega}(H)	    
\ \Longleftrightarrow\ \upg \in\Upg(L),  \qquad
\upo\in\bar{\omega}(H)	    
\ \Longleftrightarrow\ \upg \in\I(L).$$
\end{cor}
\begin{proof} If $\upg\notin \Upg(L)$, then $\upo\in \omega(L)^{\downarrow}$ by [ADH, 11.8.31], hence $\upo\in \bar{\omega}(H)$. If $\upg\in \Upg(L)$, then we can use Corollary~\ref{omuplosc} for $L$ to conclude $\upo\notin \bar{\omega}(H)$. The equivalence on the right now follows from that on the left and [ADH, 11.8.19]. 
\end{proof}

\noindent
We also note that if $\upo/4$ generates oscillations, then we have many choices for $\upg$:

\begin{cor}\label{cor:sigma(upg)=upo}
Suppose~$\upo/4$ generates oscillations. Then there are continuum many $\tilde\upg\in (\Calinf)^\times$
with $\tilde\upg>0$ and $\sigma(\tilde\upg)=\upo$, and no Hardy field extension of~$H$ contains more than one such germ $\tilde\upg$.
\textup{(}In particular, $H$ has continuum many maximal Hardy field extensions.\textup{)}
\end{cor} 
\begin{proof} As before we arrange $\upg>0$ and set $L:=\operatorname{Li}\!\big(H\langle\upg\rangle\big)$.
 Take~$\phi\in L$ with $\phi'= \frac{1}{2}\upg$ and consider the germs
$$y_1 :=  \frac{1}{\sqrt{\upg}}\cos \phi, \quad y_2 := \frac{1}{\sqrt{\upg}}\sin \phi\qquad\text{in $\Calinf$.}$$ 
The remarks preceding Lemma~\ref{lem:Im z} show:  $y_1$, $y_2$ solve the differential equation~$4Y''+\upo Y=0$, their
Wronskian equals $1/2$, and $\phi\succ 1$ (since $\upo/4$ generates oscillations). We now dilate~$y_1$,~$y_2$: let~$r\in\R^>$ be arbitrary and set
$$y_{1r}\ :=\ ry_1,\qquad   y_{2r}\ :=\ r^{-1}y_2.$$
Then $y_{1r}$, $y_{2r}$ still solve the   equation~$4Y''+\upo Y=0$, and their Wronskian is $1/2$.
Put~$\upg_r:=1/(y_{1r}^2+  y_{2r}^2)\in\Calinf$. Then $\sigma(\upg_r)=\upo$.
Let   $r,s\in\R^>$. Then
$$\upg_r=\upg_s \quad\Longleftrightarrow\quad y_{1r}^2+y_{2r}^2 = y_{1s}^2+y_{2s}^2
\quad\Longleftrightarrow\quad (r^2-s^2)\cos^2\phi + \big(\textstyle\frac{1}{r^2}-\frac{1}{s^2}\big)\sin^2\phi=0,$$
and hence $\upg_r=\upg_s$ iff $r=s$.
Next, suppose $M$ is a $\d$-perfect Hardy field extension of $H$ containing both $\upg$ and $\tilde\upg \in (\Calinf)^\times$
with $\tilde\upg>0$ and $\sigma(\tilde\upg)=\upo$.
Corollary~\ref{cor:omega(H) downward closed} gives $\upo\notin\omega(M)$,
hence
$\upg, \tilde\upg\in\Upg(M)$ by [ADH, 11.8.31], and thus $\upg=\tilde\upg$ by [ADH, 11.8.29]. 
\end{proof}

\subsection*{Answering a question of Boshernitzan\astr}
Following~\cite{Boshernitzan86} we say that a germ~$y$ in $\c$ is {\bf translogarithmic} if~$r \leq y \leq \ell_n$ for all $n$ and all $r\in\R$.
Thus for eventually strictly increasing $y\succ 1$ in $\c$,  $y$ is
translogarithmic iff its compositional inverse $y^{\operatorname{inv}}$ is transexponential.
 %If $y\in\c$ is translogarithmic, then so are the germs $cy^r$ ($c,r\in\R^>$), $y\circ \ex^x$, and $y\circ\log x$.
By  Lemma~\ref{lem:Bosh6.5} and Corollary~\ref{cor:Bosh 1.3} there exist $\Gom$-hardian    translogarithmic germs; see also~[ADH, 13.9].\index{germ!translogarithmic}\index{translogarithmic}
 Translogarithmic hardian germs are $\d$-transcendental, by Corollary~\ref{cor:Bosh 14.11}.
In this subsection we use Theorem~\ref{upo} to prove the following analogue of Corollary~\ref{cor:Bosh 1.3} for translogarithmic germs, thus giving a positive answer to Question~4 in~\cite[\S{}7]{Boshernitzan86}: 

\begin{prop}\label{prop:translog}
Every maximal Hardy field contains a translogarithmic germ.
\end{prop}

\noindent
Let $H\supseteq\R$ be a Liouville closed Hardy field; then $H$ has no translogarithmic element iff $(\ell_n)$ is a logarithmic sequence for $H$ in the sense of [ADH, 11.5].
In this case, if~$H$ is also $\upo$-free, then for each  translogarithmic $H$-hardian germ~$y$ the isomorphism type of the ordered  differential field
$H\langle y\rangle$ over $H$ is uniquely determined; more generally, by [ADH, 13.6.7, 13.6.8]:

\begin{lemma}
Let $H$ be an   $\upo$-free $H$-field, with asymptotic couple $(\Gamma,\psi)$, and let $L=H\langle y\rangle$ be a pre-$H$-field extension  of $H$ with $\Gamma^< < vy< 0$. Then for all~$P\in H\{Y\}^{\neq}$ we have
$$v\big(P(y)\big)\ =\ \gamma + \ndeg(P)vy+\nwt(P)\psi_{L}(vy)\qquad\text{where $\gamma=v^{\ev}(P)\in\Gamma$,}$$
and thus
$$\Gamma_L\ =\ \Gamma\oplus \Z vy \oplus \Z \psi_L(vy)\qquad \text{\textup{(}internal direct sum\textup{)}.}$$
Moreover, if $L^*=H\langle y^*\rangle$ is a  pre-$H$-field extension  of $H$ with $\Gamma^< < vy^*< 0$ and~$\sgn y=\sgn y^*$, then
there is a unique pre-$H$-field isomorphism $L\to L^*$ which is the identity on $H$ and sends $y$ to $y^*$.
\end{lemma}
%\begin{proof}
%Since $H$ is real closed, it may be equipped with a monomial group, so the material of [ADH, 13.4] applies, and the lemma now follows from [ADH, 13.4.1, 13.4.2, 13.4.4].
%\end{proof}

\noindent
This lemma   suggests how to obtain Proposition~\ref{prop:translog}: follow the arguments in the proof of [ADH, 13.6.7]. In the rest of this subsection we carry out this plan.  
For this, let $H\supseteq\R$ be a Liouville closed Hardy field and
   $y\in\Calinf$. %We let~$G^n_k\in\Q\{X\}$ ($1\leq k\leq n$) be the differential polynomials introduced in [ADH, 5.7],  viewed as elements of  $H\{X\}$. 

\begin{lemma}\label{lem:translog} 
Suppose $H$ is $\upo$-free and for all $\ell\in H^{>\R}$ we have, in $\c$:
\begin{enumerate}
\item[\textup{(i)}] $1\prec y\prec \ell$;
\item[\textup{(ii)}] $\derdelta^n(y)\preceq 1$ for all $n\ge 1$, where $\derdelta:=\phi^{-1}\der$, $\phi:=\ell'$;
\item[\textup{(iii)}] $y'\in\c^\times$ and $(1/\ell)' \preceq y^\dagger$.
\end{enumerate}
Let~$P\in H\{Y\}^{\neq}$. Then in~$\c$ we have
$$P(y) \sim a\,  y^d \,  (y^\dagger)^w \qquad\text{where
  $a\in H^\times$, $d=\ndeg(P)$, $w=\nwt(P)$.}$$
\textup{(}Hence $y$ is hardian over $H$ and $\d$-transcendental over $H$.\textup{)}
\end{lemma}

\begin{proof}
Since $H$ is real closed, it has a monomial group, so the material of [ADH, 13.3] applies.
Then [ADH, 13.3.3] gives a monic~$D\in\R[Y]^{\neq}$, $b\in H^\times$, $w\in\N$, and an  active element $\phi$ of $H$ with $0<\phi\prec 1$ such that:
$$P^\phi\  =\ b \cdot D\cdot (Y')^w + R,\qquad   R\in H^{\phi}\{Y\},\ R  \prec_{\phi}^\flat  b.$$
Set $d:=\ndeg P$, and note that  by~[ADH, 13.1.9] we have $d=\deg D+w$ and~$w=\nwt P$.
Replace~$P$,~$b$,~$R$ by~$b^{-1}P$,~$1$,~$b^{-1}R$, respectively, to arrange $b=1$. Take~$\ell\in H$ with $\ell'=\phi$, so~$\ell>\R$; we use the superscript $\circ$ as in the subsection on compositional conjugation
of Section~\ref{sec:Hardy fields}; in particular, 
$y^\circ=y\circ\ell^{\operatorname{inv}}$ with~$(y^\circ)'=(\phi^{-1}y')^\circ$, so~$(y^\circ)^\dagger\succeq 1/x^2$
by   hypothesis (iii) of our lemma.
In $H^\circ\{Y\}$ we now have
$$(P^\phi)^\circ\  =\   D\cdot (Y')^w + R^\circ\qquad\text{where}\qquad    R^\circ  \prec^\flat 1.$$
Evaluating at $y^\circ$ we have $D(y^\circ)\big((y^\circ)'\big){}^w \sim (y^\circ)^d \big((y^\circ)^\dagger\big){}^w$ and so
$D(y^\circ)\big((y^\circ)'\big){}^w \succeq x^{-2w} \asymp^\flat 1$.
By (i) we have $(y^\circ)^m\prec x$ for $m\geq 1$, and
by~(ii)  we have
$(y^\circ)^{(n)}\preceq 1$ for $n\geq 1$.
Hence $R^\circ(y^\circ)\preceq  h^\circ$ for some  $h\in H$ with $h^\circ\prec^\flat 1$.
Thus  in $\c$ we have
$$(P^\phi)^\circ(y^\circ) \sim   (y^\circ)^d  \big((y^\circ)^\dagger\big){}^w.$$
Since $P(y)^\circ=(P^\phi)^\circ(y^\circ)$, this yields
$P(y)\sim a\cdot y^d\cdot (y^\dagger)^w$ for $a=\phi^{-w}$.
\end{proof}

\begin{cor}\label{cor:translog}
Suppose $H$ is $\upo$-free and $1\prec y\prec\ell$ for all $\ell\in H^{>\R}$.
Then the following are equivalent:
\begin{enumerate}
\item[\textup{(i)}] $y$ is hardian over $H$;
\item[\textup{(ii)}]  for all $h\in H^{>\R}$ there is an $\ell\in H^{>\R}$ such that $\ell\preceq h$ and
$y$, $\ell$ lie in a common Hardy field;
\item[\textup{(iii)}] for all $h\in H^{>\R}$ there is an $\ell\in H^{>\R}$ such that $\ell\preceq h$ and $y\circ\ell^{\operatorname{inv}}$ is hardian.
\end{enumerate}
\end{cor}
 
\begin{proof}
The implications (i)~$\Rightarrow$~(ii)~$\Rightarrow$~(iii) are obvious. Let~$\ell\in H^{>\R}$ be such that $y^\circ:=y\circ\ell^{\operatorname{inv}}$ lies in a Hardy field~$H_0$; we arrange $x\in H_0$. For $\phi:=\ell'$ we have~$(\phi^{-1}y^\dagger)^\circ  = 
(y^\circ)^\dagger \succ  (1/x)'=-1/x^2$
and thus $y^\dagger \succ  -\phi/\ell^2=(1/\ell)'$.
Also~$y^\circ\prec x$, hence $z:=(y^\circ)'\prec x'=1$ and so~$z^{(n)}\prec 1$ for all $n$.
With $\derdelta:=\phi^{-1}\der$ and $n\geq 1$
we   have 
 $\derdelta^{n}(y)^\circ=z^{(n-1)}$ 
 and thus~$\derdelta^{n}(y)\prec 1$. 
 Moreover, for $h\in H^{>\R}$ with~$\ell\preceq h$ and $\theta:=h'$
 we have~$\theta^{-1}\der=f\derdelta$ where $f:=\phi/\theta\in H$, $f\preceq 1$.  Let $n\ge 1$. Then  
 \begin{multline*} (\theta^{-1}\der)^n=(f\derdelta)^n=G_n^n(f)\derdelta^n + \cdots + G^n_1(f)\derdelta\quad\text{ on $\Calinf$} \\ \text{ where
 $G^n_j\in \Q\{X\}\subseteq H^\phi\{X\}$ for $j=1,\dots,n$.}\end{multline*} 
As $\derdelta$ is small as a derivation on $H$, we have
 $G^n_j(f)\preceq 1$ for~$j=1,\dots,n$, and so~$(\theta^{-1}\der)^n(y) \prec  1$. 
 Thus (iii)~$\Rightarrow$~(i) by Lemma~\ref{lem:translog}.
 \end{proof}

\begin{proof}[Proof of Proposition~\ref{prop:translog}]
Let $H\supseteq\R$ be any $\upo$-free Liouville closed Hardy field not containing any translogarithmic element;
in view of Theorem~\ref{upo} it suffices to show that then some Hardy field extension of $H$ contains a translogarithmic element.  The remark before Proposition~\ref{prop:translog} yields a translogarithmic germ $y$ in a $\Gom$-Hardy field $H_0\supseteq\R$. 
Then for each~$n$, the germs~$y$,~$\ell_n$ are contained in a common Hardy field, namely $\Li(H_0)$.
Hence $y$ generates a proper Hardy field extension of~$H$ by (ii)~$\Rightarrow$~(i) in Corollary~\ref{cor:translog}.
\end{proof}

\noindent
Proposition~\ref{prop:translog}  goes through
when ``maximal'' is replaced by ``$\Ginf$-maximal'' or ``$\Gom$-maximal''. This follows from its proof, using also
remarks after the proof of Theorem~\ref{upo}. 
Here is a conjecture that is much stronger than  Proposition~\ref{prop:translog}; it postulates an analogue of  Corollary~\ref{cor:Bosh 1.1} for infinite ``lower bounds'':

\begin{conjecture}
If $H$ is maximal, then there is no $y\in\c^1$ such that $1 \prec y \prec h$ for all~$h\in H^{>\R}$, and $y'\in\c^\times$.
\end{conjecture}

\noindent 
We observe that in this conjecture we may restrict attention to $\Gom$-hardian germs~$y$:

\begin{lemma}
Suppose  there exists
 $y\in\c^1$ such that $1\prec y\prec h$ for all $h\in H^{>\R}$ and $y'\in\c^\times$. Then there exists such a germ~$y$ which is $\Gom$-hardian.
\end{lemma}
\begin{proof} Take $y$ as in the hypothesis.
Replace $y$ by $-y$ if necessary to arrange~$y>\R$.
Now Theorem~\ref{thm:Bosh 1.2} yields a $\Gom$-hardian germ  $z \geq y^{\operatorname{inv}}$. By 
Lemma~\ref{lem:Bosh6.5}, the germ~$z^{\operatorname{inv}}$ is also  
 $\Gom$-hardian, and~$\R< z^{\operatorname{inv}}\leq y\prec h$ for all $h\in H^{>\R}$.
\end{proof}
 
\subsection*{Generalizing a theorem of Boshernitzan\astr}
{\it In this subsection $H$ is a Hardy field.}\/ Recall from Corollary~\ref{cor:Bosh13.10} that for all $f\in\Ex(H)$ there are~$h\in H(x)$ and~$n$ such that $f\leq \exp_n h$.  In particular, the sequence $(\exp_n x)$ is cofinal in~$\Ex(\Q)$.
By Theorem~\ref{thm:Bosh 14.4} and Corollary~\ref{cor:Ros83}, $(\ell_n)$ is coinitial in~$\Ex(\Q)^{>\R}$; 
see also~\cite[Theorem~13.2]{Boshernitzan82}.
In particular, for the Hardy field $H=\Li(\R)$,
the subset $H^{>\R}$ is coinitial  in  $\Ex(\Q)^{>\R}=\Ex(H)^{>\R}$, equivalently,
$\Gamma_H^<$ is  cofinal   in $\Gamma_{\Ex(H)}^{<}$.  We now generalize this fact, recalling from
the end of Section~\ref{sec:order 2 Hardy fields} that $\Li(\R)$ is $\upo$-free:  

\begin{theorem}\label{thm:coinitial in E(H)}
Suppose $H$ is $\upo$-free. Then $\Gamma_H^<$ is   cofinal   in $\Gamma_{\Ex(H)}^{<}$.
\end{theorem}

\begin{proof}
Replacing $H$ by $\Li\!\big(H(\R)\big)$ and using Theorem~\ref{thm:ADH 13.6.1} we arrange that $H$ is
Liouville closed and $H\supseteq \R$.
Let $y\in\Ex(H)$ and  suppose towards a contradiction that~$\R < y < H^{>\R}$.
Then $f:=y^{\operatorname{inv}}$  is transexponential  and hardian (Lemma~\ref{lem:Bosh6.5}). 
%Take  $b\in\c$ such that $f^n \leq b$ for each~$n$; then $b$ bounds $\R\langle f\rangle$. (See Corollary~\ref{cor:val gp at infty}.)
Lemma~\ref{lem:bounded Hardy field ext, 3} gives a bound $b\in \c$ for $\R\langle f \rangle$. 
Lem\-ma~\ref{lem:Bosh 14.3} gives~$\phi\in (\c^\omega)^\times$ such that~$\phi^{(n)}\prec 1/b$ for all~$n$; 
set $r:=\phi\cdot\sin x\in\c^\omega$. Then by Lem\-ma~\ref{lem:Q bound} (with~$\R\langle f\rangle$ in place of $H$)
we have $Q(r)\prec 1$ for all~$Q\in \R\langle f\rangle\{Y\}$ with~${Q(0)=0}$.
Hence~$g:=f+r$ is eventually strictly increasing with~$g\succ 1$, and~$y=f^{\operatorname{inv}}$ and~$z:=g^{\operatorname{inv}}\in\Calinf$ do not lie in a common Hardy field. 
Thus in order to achieve the desired contradiction   it suffices to show that $z$ 
is $H$-hardian.  For this we use Corollary~\ref{cor:translog}. It is clear that~$f\sim g$, so
%First note that~$g-f=r\in (\Calinf)^{\preceq}$ and so~$f\sim g$ by Lemma~\ref{lem:y siminf z}.
%Moreover~$y\flatter x$, hence~
$y\sim z$ by Corollary~\ref{cor:Entr}, and thus~$1\prec z\prec \ell$ for all~${\ell\in H^{>\R}}$.
Let~${\ell\in H^{>\R}}$ and $\ell \prec x$; we claim that~$z\circ\ell^{\operatorname{inv}}$ is hardian, equivalently, by Lemma~\ref{lem:Bosh6.5}, that~$\ell\circ g=(z\circ\ell^{\operatorname{inv}})^{\operatorname{inv}}$ is hardian.
Now~$\ell\circ f=(y\circ\ell^{\operatorname{inv}})^{\operatorname{inv}}$ is hardian and $\ell\circ f\succ 1$,
and Lemma~\ref{lem:difference in I} gives $\ell\circ f-\ell\circ g\in (\Calinf)^{\preceq}$. 
Hence~${\ell\circ f\sim_{\infty} \ell\circ g}$ by Lemma~\ref{lem:y siminf z}.  For all $n$ we have~${\ell_n\circ \ell=\log_n \ell\in H^{>\R}}$, so~$y\le \ell_n\circ \ell$, hence~${y\circ \ell^{\inv}\le \ell_n}$, which by
compositional inversion gives~${\ell\circ f\ge \exp_n x}$. So $\ell\circ g$ is hardian by Corollary~\ref{cor:yhardian}. 
Thus~$z$ is $H$-hardian by (iii)~$\Rightarrow$~(i) of Corollary~\ref{cor:translog}.
%Applying Lemma~\ref{lem:y H-hardian crit} with $H=\R$ then yields that $\ell\circ g$ is hardian. 
%by Lemma~\ref{lem:y H-hardian crit}. 
\end{proof}

\noindent
If $H\subseteq\Ginf$ is $\upo$-free, then  $\Gamma_H^<$ is also  cofinal   in $\Gamma_{\Ex^\infty(H)}^{<}$, and similarly
with $\omega$ in place of $\infty$. (Same proof as that of the previous theorem.) We also note that if~$\Dx(H)=\Ex(H)$ (e.g., if $H$ is bounded; cf.~Theorem~\ref{thm:Bosh 14.4}), then
Theorem~\ref{thm:coinitial in E(H)} already follows from Theorem~\ref{thm:ADH 13.6.1}.
%Here is a partial converse of the implication in Theorem~\ref{thm:coinitial in E(H)}:

%\begin{cor}
%Suppose  $\Gamma_H^<$ is   cofinal  in $\Gamma_{\Dx(H)}^{<}$. Then $H$ is $\upo$-free, or 
%$H$ is $\upl$-free and $\bar\omega(H)\neq H\setminus\sigma\big(\Upg(H)\big){}^\uparrow$.
%\end{cor}
%\begin{proof}
%If $H$ is not $\upo$-free, then $\Dx(H)$ is also not $\upo$-free, by [ADH, remark after 11.7.19].
%Now the corollary follows from Lemma~\ref{lem:D(H) upo-free, 4} applied to $L=\Dx(H)$.
%\end{proof}

%\noindent
%We do not know whether conversely $\Gamma_H^<$ is   cofinal   in $\Gamma_{\Dx(H)}^{<}$ whenever $H$ is  not $\upo$-free but $\upl$-free with~$\bar\omega(H)\neq H\setminus\sigma\big(\Upg(H)\big){}^\uparrow$.

\section{Bounding Solutions of Linear Differential Equations}\label{sec:bounding}

\noindent
Let   $r\in\N^{\geq 1}$, and with $\i$ ranging over~$\N^{r}$, let  
$$P\  =\ P(Y,Y',\dots,Y^{(r-1)})\ =\ \sum_{\dabs{\i}<r} P_{\i}Y^{\i}\ \in\ \c[\imag]\big[Y,Y',\dots,Y^{(r-1)}\big]$$ 
with $P_{\i}\in\c[\imag]$ for all $\i$ with $\dabs{\i}<r$, and $P_{\i}\ne 0$ for only finitely many such $\i$.
Then~$P$ gives rise to an   evaluation map
$$y  \mapsto P\big(y,y',\dots,y^{(r-1)}\big)\ :\  \mathcal{C}^{r-1}[\imag] \to \c[\imag].$$ 
Let $y\in\c^r[\imag]$ satisfy the differential equation
\begin{equation}\label{eq:Landau}
y^{(r)}\ =\ P\big(y,y',\dots,y^{(r-1)}\big).
\end{equation}
In addition, $\fm$ with $0<\fm \preceq 1$ is a hardian germ, and
$\eta\in\c$ is eventually increasing with~$\eta(t)>0$ eventually, and $n\geq r$.

\begin{prop}\label{prop:EL}
Suppose $P_{\i}\preceq\eta$ for all $\i$, $P(0)\preceq \eta\,\fm^n$, and $y\preceq \fm^n$. 
Then
$$y^{(j)}\ \preceq\ \eta^j\fm^{n-j(1+\varepsilon)}\qquad\text{  for $j=0,\dots,r$ and  all $\varepsilon\in\R^>$,}$$
with $\prec$ in place of $\preceq$ if $y\prec \fm^n$ and $P(0)\prec\eta\,\fm^n$.
\end{prop}

\noindent
The following immediate consequence is used in Section~\ref{sec:ueeh}:

\begin{cor}\label{cor:EL}
Suppose $f_1,\dots,f_r\in\c[\imag]$ and $y\in\c^r[\imag]$ satisfy
$$y^{(r)}+f_1y^{(r-1)}+\cdots+f_ry\ =\ 0,\qquad f_1,\dots,f_r\preceq\eta,\quad y\ \preceq\ \fm^n.$$ 
Then $y^{(j)}\preceq \eta^j\fm^{n-j(1+\varepsilon)}$ for $j=0,\dots,r$ and all $\varepsilon\in\R^>$, with $\prec$ in place of $\preceq$ if~$y\prec \fm^n$. 
\end{cor}

\noindent
We obtain Proposition~\ref{prop:EL} from estimates due to Esclangon and Landau.
To prepare for this we review an argument of Hardy-Littlewood which bounds the derivative~$f'$ of
a twice continuously differentiable function~$f$ in terms of $f$, $f''$. (For another statement in the same spirit see Lemma~\ref{lem:twice diff}.)

\subsection*{Bounding $f'$ in terms of $f$, $f''$}
Let $a\in\R$, let
$\phi,\psi\colon [a,+\infty)\to (0,+\infty)$ be continuous and increasing, and  $f\in\c^2_a[\imag]$.
If $f$ and $f''$ are bounded, then so is $f'$ by the next lemma:

\begin{lemma}[{Hardy-Littlewood \cite{HaLi}}]\label{lem:HL}
Suppose   $\abs{f} \leq \phi$, $\abs{f''} \leq \psi$, and let $\varepsilon\in\R^>$. Then
$\abs{f'(t)} \leq (2+\varepsilon)\sqrt{\phi(t)\psi(t)}$, eventually.
\end{lemma}
\begin{proof}[Proof \textup{(Mordell~\cite{Mordell})}]
First arrange $a=0$ by translating the domain. Let $0<s<t$. Taylor expansion at $t$ yields $\theta=\theta(s,t)\in [0,1]$ such that
$$f(t-s)\ =\  f(t) - sf'(t) + \textstyle\frac{1}{2}s^2 f''(t-\theta s),$$
hence
$$\abs{ f(t-s) - f(t) + sf'(t) }\  \leq\  \textstyle\frac{1}{2}s^2 \psi(t-\theta s) \leq \textstyle\frac{1}{2}s^2 \psi(t).$$
Since $\abs{f(t)}\leq\phi(t)$ and $\abs{f(t-s)}\leq\phi(t-s)\leq\phi(t)$, this yields
$$\abs{f'(t)}\ \leq\   (2/s) \phi(t) +  (s/2) \psi(t).$$
Put $\rho(t):=\sqrt{\phi(t)/\psi(t)}$ for $t>0$.
If $t>2\rho(t)$, then $s:=2\rho(t)$ in the above gives~$\abs{f'(t)} \leq 2\rho(t)$. 
Hence if eventually $t>2\rho(t)$, then we are done.
Suppose otherwise; then~$\rho$ is unbounded,   hence so is $\psi \rho =\sqrt{\phi\psi}$.
Take $b > 0$ such that~$\sqrt{\phi(t)\psi(t)}\geq \abs{f'(0)}/\varepsilon$ for all $t\geq b$.
We claim that~$\abs{f'(t)} \leq (2+\varepsilon)\sqrt{\phi(t)\psi(t)}$ for $t\geq b$.
If  $t>2\rho(t)$, then $\abs{f'(t)}\leq 2\sqrt{\phi(t)\psi(t)}<(2+\varepsilon)\sqrt{\phi(t)\psi(t)}$, so
suppose~$t\leq 2\rho(t)$.
Then
$$\abs{f'(t)-f'(0)}\ =\ \left|\int_0^t f''(s)\,ds\right|\ \leq\ \int_0^t \abs{f''(s)}\,ds\  \leq\  \int_0^t \psi(s)\,ds\  \leq\  t\psi(t)$$
and hence
\[\abs{f'(t)}\ \leq\  \abs{f'(0)} + t\psi(t) \leq \abs{f'(0)} + 2\sqrt{\phi(t)\psi(t)}\ \leq\ (2+\varepsilon)\sqrt{\phi(t)\psi(t)}.\qedhere\]
\end{proof}

\begin{cor}\label{cor:HL}
If   $f\preceq\phi$, $f''\preceq\psi$, then 
$f'\preceq \sqrt{\phi\psi}$, with $f'\prec \sqrt{\phi\psi}$ if also $f\prec\phi$ or $f''\prec\psi$.
\end{cor}

\noindent
In Corollary~\ref{cor:HL} we cannot drop the assumption that $\phi$, $\psi$
are increasing. (Take $a>1$, $f(t)=t\log t$ for $t\geq a$, $\phi=f$, $\psi=f''$.) 
%{\bf next sentence not checked, nor used}
However, Mordell~\cite{Mordell} also shows that if instead of assuming that $\phi$, $\psi$ are increasing, we assume that they are decreasing, then
Lemma~\ref{lem:HL}  holds in a stronger form: $\abs{f} \leq \phi \ \&\ \abs{f''} \leq \psi\, \Rightarrow\, \abs{f'} \leq 2 (\phi\psi)^{1/2}$. The next lemma (not used later) yields a variant of Corollary~\ref{cor:HL} where the germ of~$f$ lies in a complexified Hardy field;
see also~\cite[\S{}7]{HaLi}.

\begin{lemma}
Let  $H$ be a Hardy field, $K=H[\imag]$, and   $g\in K^\times$ such that  $g\prec 1$ or~$g\succ 1$, $g^\dagger\nasymp x^{-1}$. Then~$g'\preceq {\abs{gg''}} ^{1/2}$.
\end{lemma} 
\begin{proof} Arranging that $H$ is real closed and $x\in H$ and using $\abs{h}\asymp h$ for $h\in K$ (see
the remarks before Corollary~\ref{cor:10.5.2 variant}), the
 lemma now follows from parts (1), (2), (4) of \cite[Lemma~5.2]{AvdD3} applied to the asymptotic couple of $K$.
\end{proof}
  
\noindent
We now generalize Corollary~\ref{cor:HL}:

\begin{lemma}[{Hardy-Littlewood  \cite{HaLi}}] \label{lem:HL, n}
Suppose $f\in\c^n_a[\imag]$, $n\geq 1$,
such that $f\preceq\phi$, $f^{(n)}\preceq\psi$. Then for~$j=0,\dots,n$ we have~$f^{(j)}\preceq \phi^{1-j/n}\psi^{j/n}$.
If additionally $f\prec\phi$ or $f^{(n)}\prec\psi$, then $f^{(j)}\prec  \phi^{1-j/n}\psi^{j/n}$  for~$j=1,\dots,n-1$.
\end{lemma}
\begin{proof}
The case $n=1$ is trivial, so let $n\geq 2$. We may also assume $f\neq 0$, and by increasing $a$ we arrange $f(a)\neq 0$.
Let $j$ range over $\{0,\dots,n\}$. Consider  the continuous increasing functions
$$\chi_j\ \colon\ [a, +\infty)\to (0,+\infty),\qquad \chi_j(t)\ :=\  \max_{a\leq s\leq t}\, \abs{f^{(j)}(s)}\,\big/\,\big(\phi(s)^{1-j/n}\psi(s)^{j/n}\big),$$
and set $\chi:=\max\{\chi_0,\dots,\chi_n\}$. Then $\chi(t)\geq\chi_0(t)>0$ for all $t\geq a$. 
We have $$\abs{f^{(j)}} / \big( \phi^{1-j/n}\psi^{j/n} \big)\ \leq\ \chi_j\  \leq\ \chi,$$ therefore
$$f^{(j)}\  \preceq\ \phi^{1-j/n}\psi^{j/n}\chi.$$
By induction on $j=0,\dots,n-2$ we now show
\begin{equation}\label{eq:fj}
f^{(j)}\ \preceq\ \phi^{1-j/n} \psi^{j/n} \chi^{1-1/2^{j}}.
\end{equation}
The case $j=0$ follows from $f\preceq\phi$. Suppose \eqref{eq:fj} holds for a certain   $j<n-2$.
Then Corollary~\ref{cor:HL} with $f^{(j)}$, $\phi^{1-j/n}\psi^{j/n}\chi^{1-1/2^j}$, $\phi^{1-(j+2)/n}\psi^{(j+2)/n}\chi$ 
in the role of   $f$, $\phi$, $\psi$, respectively,
yields:
\begin{align*}
f^{(j+1)}\ 	&\preceq\ \big( \phi^{1-j/n}\psi^{j/n}\chi^{1-1/2^j}\cdot \phi^{1-(j+2)/n}\psi^{(j+2)/n}\chi\big)^{1/2} \\
			&=\  \phi^{1-(j+1)/n}\psi^{(j+1)/n}\chi^{1-1/2^{j+1}}.
\end{align*}
This proves \eqref{eq:fj}. We claim $\chi\preceq 1$. Suppose otherwise; so $\chi(t)\to+\infty$ as~${t\to+\infty}$, since~$\chi$ is increasing, hence
$f^{(n)}\preceq\psi\preceq\psi\chi^{1-1/2^n}\preceq\psi\chi$.  Corollary~\ref{cor:HL} with~$f^{(n-2)}$, $\phi^{2/n}\psi^{1-2/n}\chi^{1-1/2^{n-2}}$, $\psi\chi$ in the role of~$f$,~$\phi$,~$\psi$, respectively, yields 
$$f^{(n-1)}\ \preceq\ \big( \phi^{2/n}\psi^{1-2/n}\chi^{1-1/2^{n-2}}\cdot \psi\chi \big)^{1/2}\ =\ \phi^{1/n}\psi^{1-1/n}\chi^{1-1/2^{n-1}}.$$
So \eqref{eq:fj} then also holds for $j=n-1$, and it clearly holds for $j=n$. But then~$\chi_j\preceq\chi^{1-1/2^n}$ for all~$j$ and so $\chi\preceq\chi^{1-1/2^n}$, contradicting $\chi\succ 1$.

Now suppose $f\prec\phi$. By induction on $j=0,\dots,n-1$ we show  $f^{(j)}\prec  \phi^{1-j/n}\psi^{j/n}$.
The case $j=0$ holds by assumption; suppose it holds for a certain $j\leq n-2$.
Then~$f^{(j+2)}\preceq \phi^{1-(j+2)/n}\psi^{(j+2)/n}$, so Corollary~\ref{cor:HL} with
$$f^{(j)},\quad \phi^{1-j/n}\psi^{j/n},\quad \phi^{1-(j+2)/n}\psi^{(j+2)/n}$$
 in the role of $f$, $\phi$, $\psi$, respectively, yields
$$f^{(j+1)} \prec \big(\phi^{1-j/n}\psi^{j/n}\cdot \phi^{1-(j+2)/n}\psi^{(j+2)/n}\big)^{1/2}=\phi^{1-(j+1)/n}\psi^{(j+1)/n}.$$
If $f^{(n)}\prec \psi$, then likewise~$f^{(n-j)} \prec \phi^{j/n}\psi^{1-j/n}$ for $j=0,\dots,n-1$.
\end{proof}

\begin{cor}\label{cor:HL, bded}
Suppose $f\in\c^n_a[\imag]$  and $f\preceq \phi$, $f^{(n)}\preceq \phi$. Then $f',\dots,f^{(n-1)}\preceq \phi$,
 and if in addition $f\prec \phi$ or $f^{(n)}\prec \phi$, then  $f',\dots,f^{(n-1)} \prec \phi$.
\end{cor}

\subsection*{The theorem of Esclangon-Landau}
In this subsection $n\ge r \ge 1$ and $P$ is as at the beginning of this section, and $y\in\c^r[\imag]$ satisfies \eqref{eq:Landau}. 
Also, $\eta\in\c$ is eventually increasing and positive, so $\eta\succeq 1$. 
The next theorem covers the case~$\fm\asymp 1$ of Proposition~\ref{prop:EL}:

\begin{theorem}[{Landau \cite{Landau}}]\label{thm:Landau}
Suppose $y\preceq 1$ and $P_{\i}\preceq\eta$ for all $\i$. Then $y^{(j)}\preceq \eta^j$ for $j=0,\dots,r$.
Moreover, if   $y\prec 1$, then
$y^{(j)} \prec \eta^j$ for $j=0,\dots,r-1$, and
if in addition $P(0)\prec\eta$, then  also
$y^{(r)} \prec \eta^r$.
\end{theorem}
\begin{proof}
Take $a\in\R$ such that $\eta$ is represented by an increasing continuous function~$\eta\colon [a,+\infty)\to (0,+\infty)$, and
$y$ by a function $y\in \c_a^r[\imag]$. Then
$$t\mapsto \psi(t)\  :=\  \max\left( 1, \max_{a\leq s\leq t} \abs{y^{(r)}(s)} \right)\ \colon\  [a,+\infty)\to [1,+\infty)$$ 
is continuous and increasing with $\abs{y^{(r)}} \leq \psi$.
By Lemma~\ref{lem:HL, n} we have 
$y^{(j)}\preceq\psi^{j/r}$ for $j=0,\dots,r-1$ and thus
$P_{\i}y^{\i} \preceq \eta \psi^{\dabs{\i}/r}\preceq\eta \psi^{1-1/r}$ if $\dabs{\i}<r$.
So~$y^{(r)}=P\big(y,\dots,y^{(r-1)}\big) \preceq  \eta\psi^{1-1/r}$.
Take $C\in\R^{>}$ such that
$$\abs{y^{(r)}(t)}\  \leq\  C\eta(t)\psi(t)^{1-1/r}\qquad\text{for  all  $t\geq a$.}$$
Increasing $C$ we arrange $C\eta(a)\psi(a)^{1-1/r}\geq 1$.
As  $\eta\psi^{1-1/r}$ is   increasing, 
$$\psi(t)\  \leq\  
\max\left( 1, \max_{a\leq s\leq t} C\eta(s)\psi(s)^{1-1/r}  \right)\
 \leq\  C\eta(t)\psi(t)^{1-1/r}\qquad\text{for   $t\geq a$.}$$
Hence~$\abs{y^{(r)}(t)} \leq \psi(t)\leq C^r\eta^r(t)$ for $t\geq a$, so $y^{(r)}\preceq\eta^r$.
By  Lemma~\ref{lem:HL, n} again this yields~$y^{(j)}\preceq \eta^j$ for $j=0,\dots,r$. Assume now that~$y\prec 1$. 
Then by that same lemma,  $y^{(j)}\prec \eta^j$ for $j<r$. We have  $\eta\succeq 1$, so if~$0<\dabs{\i}<r$, then~$y^{\i}\prec\eta^{\dabs{\i}}\preceq\eta^{r-1}$. 
Hence if additionally $P(0)\prec\eta$, then~$y^{(r)}=P\big(y,\dots,y^{(r-1)}\big) \prec \eta^r$.
\end{proof}

\begin{cor}[{Esclangon~\cite{Esc}}]
Suppose $f_1,\dots,f_r,g\in\c[\imag]$ and $y\in\c^r[\imag]$ satisfy 
$$y^{(r)}+f_1y^{(r-1)}+\cdots+f_ry\ =\ g,\qquad f_1,\dots,f_r,g,y\ \preceq\ 1.$$ 
Then $y',\dots,y^{(r)}\preceq 1$. If in addition  $y\prec 1$ and $g\prec 1$, then $y',\dots,y^{(r)}\prec  1$.
\end{cor}

\noindent 
Below  $H$ is a Hardy field and $\fm\in H$, $0<\fm\prec  1$. Recall also that $n\ge r\ge 1$. 
% and suppose $n\geq 1$.
 
\begin{lemma}\label{lem:bd mult conj}
Let    $z\in\c^r[\imag]$. If
$z^{(j)} \preceq \eta^j$ for~$j=0,\dots,r$, then $(z\fm^n)^{(j)}\preceq  \eta^j\fm^{n-j}$ for $j=0,\dots,r$,
and likewise with $\prec$ instead of $\preceq$.
\end{lemma} 
\begin{proof}
Corollary~\ref{cor:kth der of f, preceq} yields $(\fm^{n})^{(m)} \preceq \fm^{n-m}$ for   $m\leq n$.
Thus  if $z^{(j)} \preceq \eta^j$ for~$j=0,\dots,r$, then
$$ z^{(k)}  (\fm^{n})^{(j-k)} \preceq \eta^k \fm^{n-(j-k)}   \preceq \eta^j\fm^{n-j} \qquad (0\leq k\leq j\leq r),$$
so $(z\fm^n)^{(j)}\preceq  \eta^j\fm^{n-j}$ for $j=0,\dots,r$, by the Product Rule. The argument with $\prec$ instead of $\preceq$ is similar.
\end{proof}

\noindent
We now return to the assumptions on $P, y$ in the beginning of this section, so~$y\in\c^r[\imag]$ satisfies \eqref{eq:Landau}.
Suppose also that  $P_{\i}\preceq\eta$ for all $\i$, $P(0)\preceq \eta\,\fm^{n}$, and $y\preceq \fm^n$. Let
 $\varepsilon\in\R^>$ and set for $i=0,\dots,r$,
$$Y_i\ :=\  \sum_{j=0}^i {i\choose j} Y^{(i-j)} (\fm^n)^{(j)}\in H\big[Y,Y',\dots,Y^{(i)}\big]\ \subseteq\ \c[\imag]\big[Y, Y',\dots, Y^{(r)}\big].$$
Then for  $z:=y\,\fm^{-n}\preceq 1$ in $\c^r[\imag]$ the product rule gives
$$Y_i(z,z',\dots,z^{(i)})\ =\ (z\,\fm^n)^{(i)}\ =\ y^{(i)}    \qquad(i=0,\dots,r),$$
so with
$$Q\ :=\ Y^{(r)}-\fm^{-n}\big(Y_r-P(Y_0,\dots,Y_{r-1})\big)\in\c[\imag]\big[Y,Y',\dots,Y^{(r-1)}\big]$$
we have by substitution of $z,\dots, z^{(r)}$ for $Y, Y',\dots, Y^{(r)}$, 
\begin{align*}
z^{(r)}\ &= \ Q\big(z,z',\dots,z^{(r-1)}\big) + \fm^{-n}\big( y^{(r)} - P(y,y',\dots,y^{(r-1)}) \big) \\
&=\  Q\big(z,z',\dots,z^{(r-1)}\big).  
\end{align*}
For $Y_0,\dots,Y_r\in H\{Y\}$ we have
$(Y^{\i})_{\times\fm^n} = Y_0^{i_0}\cdots Y_r^{i_r}$ for $\i=(i_0,\dots,i_r)\in\N^{1+r}$. 
Now $\fm^{-\varepsilon}\in \Li\!\big(H(\R)\big)$. We equip~$\Li\!\big(H(\R)\big)\{Y\}$ with the gaussian extension of the valuation of $\Li\!\big(H(\R)\big)$. Then by  [ADH, 6.1.4],
$$
\fm^{-n} (Y^{\i})_{\times\fm^n}\ \preceq\ \fm^{-\varepsilon}  \qquad\text{for  $\i\in\N^{1+r}\setminus\{ 0 \}$.}$$
Let $\i$ range over $\N^r$ and take $Q_{\i}\in \c[\imag]$
for $\dabs{\i}<r$ such that
$$Q\ =\ \sum_{\dabs{\i}<r}Q_{\i}Y^{\i}, \qquad(Q_{\i}\ne 0 \text{ for only finitely many }\i).$$
Together with $P_{\i}\preceq\eta$ for all $\i$ and $P(0)\preceq \eta\,\fm^{n}$, the remarks above
 yield $Q_{\i}\preceq\eta\,\fm^{-\varepsilon}$ for all~$\i$.
By Theorem~\ref{thm:Landau} applied to $P$, $y$, $\eta$ replaced by~$Q$,~$z$,~$\eta\,\fm^{-\varepsilon}$, respectively,   we now obtain~$z^{(j)}\preceq (\eta\,\fm^{-\varepsilon})^j$   ($j=0,\dots,r$), with~$\prec$ in place of $\preceq$ if~$y\prec\fm^n$ and~$P(0)\prec\eta\,\fm^n$.
Using Lemma~\ref{lem:bd mult conj} with~$\eta\,\fm^{-\varepsilon}$ and $\Li\!\big(H(\R)\big)$ in place of $\eta$ and~$H$ finishes the proof of Proposition~\ref{prop:EL}. \qed

\section{Almost Periodic Functions} \label{sec:almost periodic}

\noindent 
For later use we now discuss  trigonometric polynomials, almost periodic functions, and their mean values; see~\cite{Bohr,Corduneanu} for this material in the case $n=1$.  {\it In this section we assume $n\geq 1$,
and for vectors $r=(r_1,\dots,r_n)$ and $s=(s_1,\dots,s_n)$ in~$\R^n$ we let
$r\cdot s:= r_1s_1+\cdots + r_ns_n\in \R$ be the usual dot product of $r$ and~$s$. We also set 
$rs:=(r_1s_1,\dots, r_ns_n)\in\R^n$, not to be confused with $r\cdot s\in\R$. Moreover, we let~$v,w\colon\R^n\to\mathbb C$ be complex-valued functions on $\R^n$, and let $s$ range over $\R^n$, and $T$ over $\R^{>}$; 
integrals are with respect to the usual Lebesgue measure of $\R^n$.}\/
Set 
$$\dabs{w}\ :=\ \sup_{s}\, \abs{w(s)}\ \in\ [0,+\infty].$$
We shall also have occasion to consider various functions
$\R^n\to\mathbb C$ obtained from~$w$: $\overline{w}$, $\abs{w}$, as well as $w_{+r}$ and  $w_{\times r}$  (for $r\in\R^n$), defined
by
$$\overline{w}(s)\ :=\ \overline{w(s)}, \quad 
\abs{w}(s)\ :=\ \abs{w(s)}, \quad
w_{+r}(s)\ :=\ w(r+s),\quad
w_{\times r}(s)\ :=\ w(rs).$$
We say that $w$ is {\bf $1$-periodic} if
%A trigonometric polynomial $w$ is said to be {\bf periodic} if 
$w_{+k}=w$ for all $k\in\Z^n$. \index{function!1-periodic@$1$-periodic}

\subsection*{Trigonometric polynomials} 
Let $\alpha$ range over $\R^n$.
Call~$w$ a {\bf trigonometric polynomial}\index{function!trigonometric polynomial} if there are~$w_\alpha\in\mathbb C$, with $w_\alpha=0$ for all but finitely many~$\alpha$,
such that for all $s$, 
\begin{equation}\label{eq:trig poly}
w(s)\ =\ \sum_\alpha w_\alpha \ex^{(\alpha\cdot s)\imag}.
\end{equation}
The trigonometric polynomials form a subalgebra of the $\mathbb C$-algebra of
uniformly continuous bounded functions $\R^n\to\mathbb C$. Let $w$ be a trigonometric polynomial. Then~$\overline{w}$
is a trigonometric polynomial, and for $r\in\R^n$, so are  the
functions   $w_{+r}$ and~$w_{\times r}$.
Note that $w$ extends to a complex-analytic function $\C^n\to \C$,  that $\Re w$ and $\Im w$ are real analytic, and that $\partial w/\partial x_j:=(\partial \Re w/\partial x_j) + (\partial \Im w/\partial x_j)\imag$ for  $j=1,\dots,n$ is also a trigonometric polynomial. The functions $s\mapsto \sin(\alpha\cdot s)$ and 
$s\mapsto\cos(\alpha\cdot s)$ on $\R^n$ are real valued trigonometric polynomials. 
 By Corollary~\ref{cor:id thm} below
the coefficients~$w_\alpha$ in~\eqref{eq:trig poly} are uniquely determined by $w$.

%\noindent
%To establish this proposition we use the following well-known fact:

%\begin{lemma}\label{iszero} 
%Let $a_1<\cdots<a_m$ be real numbers, $c_1,\dots,c_m\in\mathbb C$, and $U$  a nonempty open subset of $\R$ such that $\sum_{k=1}^m c_k \ex^{a_kt\imag}=0$ for all $t\in U$. Then $c_1=\cdots=c_m=0$.
%\end{lemma}
%\begin{proof} 
%By induction on $m$. The case $m=0$ being trivial, let $m\ge 1$. Replace~$a_i$ by $a_i-a_1$ ($i=1,\dots,m$) to arrange that $a_1=0$. The function $t\mapsto h(t):=\sum_k c_k \ex^{a_kt\imag}\colon \R\to\mathbb C$ is differentiable and $h= 0$ on $U$, so  $\sum_{k\geq 2} c_ka_k\imag \ex^{a_kt\imag}=h'(t)=0$ for $t\in U$. So $c_2=\cdots=c_m=0$ by inductive hypothesis,   hence also $c_1=0$.
%\end{proof}

%\begin{proof}[Proof of  Proposition~\ref{prop:id thm}] 
%By induction on $n$. The case $n=1$ is Lemma~\ref{iszero}, so let~$n\ge 2$. We arrange $U=(a_1,b_1)\times\cdots\times (a_n,b_n)$ for reals~$a_j<b_j$ ($j=1,\dots,n$). Let $\alpha'$ range over $\R^{n-1}$ and $\alpha_n$ over $\R$. Let $s'\in U':=(a_1,b_1)\times\cdots \times(a_{n-1},b_{n-1})$. Then for all $s_{n}\in (a_n,b_n)$ we have
%$$ \sum_{\alpha_n} \left(\sum_{\alpha'} w_{(\alpha',\alpha_n)} \ex^{(\alpha' \cdot s')\imag}\right)\ex^{\alpha_n s_n\imag}\ =\ 0,$$
%so $\sum_{\alpha'} w_{(\alpha',\alpha_n)}\ex^{(\alpha' \cdot s')\imag}=0$ for all $\alpha_n$ by Lemma~\ref{iszero}. Since $s'\in U'$ is arbitrary, the inductive assumption yields $w_{\alpha}=0$ for all $\alpha$.
%\end{proof}

%A trigonometric polynomial $w$ is said to be {\bf periodic} if 
%$w_{+k}=w$ for all $k\in\Z^n$.

If $w(s)=\ex^{(\alpha\cdot s)\imag}$ for all $s$, then
$w_{+r}=w$ for all $r\in\R^n$ with $\alpha\cdot r\in 2\pi\Z$.
So if~$w$ is a trigonometric polynomial as in \eqref{eq:trig poly} with $w_\alpha=0$ for all $\alpha\notin 2\pi\Z^n$, then~$w$ is $1$-periodic. Next we state a well-known consequence of the Stone-Weierstrass Theorem; 
see \cite[(7.4.2)]{Dieudonne} %or \cite[Theorem~4.25]{BabyRudin} 
for the case $n=1$.

\begin{prop}\label{prop:StoWei}
If $v$ is continuous and $1$-periodic, then
for every~$\varepsilon$ in $\R^>$  there is a $1$-periodic trigonometric polynomial $w$ with $\dabs{v-w}<\varepsilon$.
\end{prop}

\subsection*{Almost periodic functions}  
We call $w$ {\bf almost periodic}\index{function!almost periodic} (in the sense of Bohr) if for every $\varepsilon$ in $\R^>$ there is a trigonometric
polynomial $v$  such that $||v-w|| \le \varepsilon$.  
If~$w$ is almost periodic, then $w$  is uniformly continuous and bounded
(as the uniform limit of a sequence of uniformly continuous bounded functions~$\R^n\to\mathbb C$).
If $w$ is almost periodic, then so are $\overline{w}$, and $w_{+r}$, $w_{\times r}$ for  $r\in\R^n$. 

Note that the $\mathbb C$-algebra of
uniformly continuous bounded functions $\R^n\to \mathbb C$ is a Banach algebra
with respect to $\|\cdot\|$: it is complete with respect to
this norm. The closure of
its subalgebra of trigonometric polynomials with respect to this norm is~${\{w: \text{$w$ is almost perodic}\}}$, which is therefore a Banach subalgebra.
In particular, if $v,w\colon \R^n\to \mathbb C$ are almost periodic, so are
$v+w$ and $vw$. Moreover:

%\begin{lemma}
%Suppose $v$ and $w$ are almost periodic. Then $v+w\colon\R^n\to\mathbb C$ and $vw\colon\R^n\to\mathbb C$  are almost periodic.
%\end{lemma}
%\begin{proof}
%The first claim is obvious. For the second claim, 
%let $\varepsilon\in(0,1]$.
%Take $M\in\R^{\geq}$ with $\dabs{v},\dabs{w}\leq M$, and then take trigonometric
%polynomials $\tilde{v},\tilde{w}\colon\R^n\to\mathbb C$ such that
%$\dabs{\tilde{v}-v},\dabs{\tilde{w}-w} \leq \frac{\varepsilon}{2(M+1)}$.
%Then $\dabs{\tilde{v}}\leq\dabs{\tilde{v}-v}+\dabs{v}\leq 1+M$ and hence
%$$\dabs{vw-\tilde{v}\tilde{w}}\ \leq\ \dabs{w}\cdot \dabs{v-\tilde{v}}+\dabs{\tilde{v}}\cdot\dabs{w-\tilde{w}}\ \leq\
%\frac{M\varepsilon}{2(M+1)}+\frac{(M+1)\varepsilon}{2(M+1)}\ \leq\
%\varepsilon.$$
%It remains to note that $\tilde{v}\tilde{w}$ is a trigonometric polynomial.
%\end{proof}

\begin{cor}\label{apsub} Let $v_1,\dots, v_m\colon \R^n\to \mathbb C$ be almost periodic, let
$X\subseteq {\mathbb C}^m$ be closed, and suppose $F\colon X\to \mathbb C$ is continuous with
$\big(v_1(s),\dots, v_m(s)\big)\in X$ for all $s$. Then the function
$F(v_1,\dots, v_m)\colon \R^n\to \mathbb C$ is almost periodic. 
\end{cor} 
\begin{proof} Since $v_1,\dots, v_m$ are bounded we can arrange that
$X$ is compact. Let $\varepsilon\in \R^{>}$. Then Weierstrass Approximation \cite[(7.4.1)]{Dieudonne}  gives a polynomial $$P(x_1,y_1,\dots, x_m,y_m)\ \in\ \C[x_1,y_1,\dots, x_m,y_m]$$ such that
$|F(z_1,\dots, z_m)-P(z_1,\bar{z}_1,\dots, z_m,\bar{z}_m)|\le \varepsilon$ for all $(z_1,\dots, z_m)\in X$. Hence $\|F(v_1,\dots, v_m)-P(v_1,\bar{v}_1,\dots, v_m\bar{v}_m)\|\le \varepsilon$.
It remains to note that the func\-tion~$P(v_1,\bar{v}_1,\dots, v_m,\bar{v}_m)$ is almost periodic.  
\end{proof} 

\noindent
Call $w$ {\bf normal}\index{function!normal}\index{normal!function} if $w$ is bounded and for every sequence $(r_m)$ in $\R^n$ the se\-quence $(w_{+r_m})$ of
functions $\R^n\to\mathbb C$ has a uniformly converging subsequence.
One verifies easily that if $v$, $w$ are normal, then so are the functions  $v+w$ and $cv$ ($c\in\mathbb C$);
hence by the next lemma, each trigonometric polynomial is normal:

\begin{lemma}
Suppose $w(s)=\ex^{(\alpha\cdot s)\imag}$ for all $s$. Then $w$ is normal.
\end{lemma}
\begin{proof}
Let $(r_m)$ be a sequence in $\R^n$.
%By Bolzano-Weierstrass, after 
Passing to a subsequence of $(r_m)$ we arrange that the sequence
$\big(w(r_m)\big)$ of complex numbers of modulus~$1$ converges. Now use that for all $l$, $m$ and all $s$ we have 
$\abs{ w_{+r_l}(s)-w_{+r_m}(s) } = \abs{ w(r_l)-w(r_m) }$, and thus~$\dabs{w_{+r_l}-w_{+r_m}}\leq \abs{ w(r_l)-w(r_m) }$.
\end{proof}

\begin{lemma}
Let $(w_m)$ be a sequence of normal functions with $\dabs{w_m-w}\to 0$ as~$m\to \infty$. Then $w$ is normal.
\end{lemma}
\begin{proof}  
Let $(r_k)_{k\in \N}$ be a sequence in $\R^n$. Using normality of the $w_m$ we obtain inductively subsequences $(r_{k0}),(r_{k1}),\dots$ of $(r_k)$ such that for all $m$, $\big( (w_m)_{+r_{km}}\big)$ converges uniformly and
$(r_{k,m+1})$  is a subsequence of $(r_{km})$. 
Then for every $m$,
$(r_{m+l,m+l})_{l\geq 0}$ is a subsequence of $(r_{km})$; so $\big( (w_m)_{+r_{kk}}\big)$  converges uniformly. Now let $\varepsilon\in\R^>$ be given. Take $m$ so that~$\dabs{ w_m-w } \le  \varepsilon$, and then take $k_0$ so that~$\dabs{ (w_m)_{+r_{kk}} - (w_m)_{+r_{ll}} } \le \varepsilon$ for all $k,l\geq k_0$.
For such $k$, $l$ we have
\begin{multline*}
\dabs{w_{+r_{kk}} - w_{+r_{ll}} } \leq \\ \dabs{ w_{+r_{kk}} - (w_m)_{+r_{kk}}} +
\dabs{ (w_m)_{+r_{kk}} - (w_m)_{+r_{ll}} } + \dabs{ (w_m)_{+r_{ll}} - w_{+r_{ll}}} \le 3\varepsilon.
\end{multline*}
Thus  $( w_{+r_{kk}})$ converges uniformly.
\end{proof}

\begin{cor}[Bochner]\label{cor:Bochner}
Every almost periodic function $\R^n\to\mathbb C$ is normal.
\end{cor}

\noindent
For $\varepsilon\in\R^>$, we say that $r\in\R^n$ is an {\bf $\varepsilon$-translation vector}\index{translation vector} for $w$ if $\dabs{w_{+r}-w}<\varepsilon$.
 We define  an {\bf $n$-cube}\/ of side length $\ell\in\R^>$  to be
a subset of $\R^n$ of the form~$I=I_1\times\cdots\times I_n$ where each
$I_1,\dots,I_n$ is an open interval of length $\ell$.

\begin{prop}\label{prop:eps-trans}
If $w$ is normal, then for all  $\varepsilon\in\R^>$ there is an $\ell=\ell(w,\varepsilon)\in\R^>$ such that every $n$-cube
of side length $\ell$ contains an $\varepsilon$-translation vector for~$w$.
\end{prop}
\begin{proof}
We assume that $w$ is bounded and show the contrapositive. Let $\varepsilon\in\R^>$ be such that there are $n$-cubes of arbitrarily large sidelength that contain no $\varepsilon$-translation vector for~$w$; to conclude that $w$ is not normal it suffices to have a sequence~$(r_i)_{i\in \N}$ in $\R^n$ such
that  $r_j-r_i$ is not an $\varepsilon$-translation vector for $w$, for all~$i<j$, since then~$\dabs{w_{+r_j}-w_{+r_i}} = \dabs{w_{+(r_j-r_i)}-w}\geq\varepsilon$ for all $i<j$.
Now sup\-pose $r_0,\dots,r_{m}\in\R^n$ are such that 
 $r_j-r_i$ is not an $\varepsilon$-translation vector for $w$, for all $i<j\le m$. 
 Then for $k=1,\dots,n$ we take
 intervals $I_k=(a_k,b_k)$ ($a_k<b_k$ in~$\R$) of equal length~$b_k-a_k>2\max\big\{\abs{r_{0}}_\infty,\dots,\abs{r_{m}}_\infty\big\}$ such that $I:=I_1\times\cdots\times I_n$  does not contain
an $\varepsilon$-translation vector for $w$.  Set $r_{m+1}:=\frac{1}{2}(a_1+b_1,\dots,a_n+b_n)$; then for 
$i\le m$ we have $r_{m+1}-r_i\in I$, hence $r_{m+1}-r_i$ is not an $\varepsilon$-translation vector for~$w$. 
\end{proof}

\noindent
By Corollary~\ref{cor:Bochner}, Proposition~\ref{prop:eps-trans} applies to almost periodic $w$. Bohr~\cite{Bohr} showed conversely 
%\marginpar{not checked} 
that if $w$ is continuous and
satisfies the conclusion of Proposition~\ref{prop:eps-trans}, then $w$ is almost periodic, but we do not use this elegant characterization of almost periodicity below.
We now improve Proposition~\ref{prop:eps-trans} for almost periodic~$w$.
{\it In the rest of this subsection we assume that $w$ is almost periodic.}\/

\begin{lemma}\label{transcube}
Let $\varepsilon\in\R^>$; then there are $\delta,\ell\in\R^>$ such that every $n$-cube of side length $\ell$ contains an $n$-cube 
of side length~$\delta$ consisting entirely of $\varepsilon$-translation vectors for $w$.
\end{lemma}
\begin{proof}
Uniform continuity of $w$ yields $\delta_1\in\R^>$ such that all  $d\in\R^n$ with $\abs{d}_\infty<\delta_1$
are $(\varepsilon/3)$-translation vectors for $w$. Take $\ell_1:=\ell(w,\varepsilon/3)$ as in Proposition~\ref{prop:eps-trans}, and set $\delta:=2\delta_1$, $\ell:=\ell_1+\delta$.
 Let $J=a+(0,\ell)^n$ be a cube of side length $\ell$, where~$a\in\R^n$. Take an $(\varepsilon/3)$-translation vector $r\in a+(\delta_1,\ell_1+\delta_1)^n$  for $w$. The cube $I:=r+(-\delta_1,\delta_1)^n$ of side length $\delta$ is entirely contained in $J$. Let $p\in I$. Then for $d:=p-r$ we have
 $\abs{d}_\infty<\delta_1$, so for all $s$, 
 $$\abs{w(s+p)-w(s)}\ \leq\ \abs{w(s+d+r)-w(s+d)} + \abs{w(s+d)-w(s)}\ <\ \frac{\varepsilon}{3}+\frac{\varepsilon}{3}\ <\ \varepsilon,$$
hence $p$ is an $\varepsilon$-translation vector for $w$.
\end{proof}

\begin{cor}\label{cor:ap lower bd}
Suppose   $w(\R^n)\subseteq\R$, $s_0\in\R^n$, and~$w(s_0)>0$. Then there are~$\delta_1,\ell_1\in\R^>$ such that  every $n$-cube of side length $\ell_1$ contains an $n$-cube~$I$
of side length~$\delta_1$ with $w(s)\geq w(s_0)/3$ for all $s\in I$.
\end{cor}
\begin{proof}
Let  $\delta$, $\ell$ be as in Lemma~\ref{transcube} for $\varepsilon:=w(s_0)/3$. By decreasing $\delta$ we obtain from the uniform continuity of $w$ that all $d\in\R^n$ with $\abs{d}_\infty<\delta/2$ are $\varepsilon$-translation vectors for $w$.
Set  $\delta_1:=\delta$, $\ell_1:=\ell+\delta/2$. Let $J=a-(0,\ell_1)^n$
with $a\in\R^n$ be an $n$-cube of side length~$\ell_1$; we claim that $J$ contains an $n$-cube~$I$ of side length~$\delta$ with $w(s)\geq \varepsilon$ for all $s\in I$.
To prove this claim, consider the $n$-cube~$J_0:=(s_0-a)+(\delta/2,\ell+\delta/2)^n$ of side length $\ell$.
Our choice of $\delta$, $\ell$ gives an $\epsilon$-translation vector $r\in J_0$ for $w$ such that $r+(-\delta/2,\delta/2)^n\subseteq J_0$. Then
$$I\ :=\ (s_0-r)+(-\delta/2,\delta/2)^n\ \subseteq\ s_0-J_0\ =\ a-(\delta/2,\ell+\delta/2)^n\ \subseteq\ J,$$
and for every $s\in I$, setting $d=s-s_0+r$, we have $|d|_{\infty}< \delta/2$, so
$$w(s)\ =\ w(s_0)+\big(w(s_0+d)-w(s_0)\big)-\big(w(s+r)-w(s)\big)\ \geq\ w(s_0)-\varepsilon-\varepsilon\ =\ \varepsilon$$
as required.
\end{proof}

\begin{lemma}\label{lem:limsup vs sup}
Suppose   $w(\R^n)\subseteq\R$. Then with $|s|:= |s|_{\infty}$, 
$$\liminf_{\abs{s}\to \infty} w(s)\ =\ \inf_{s} w(s),\qquad \limsup_{\abs{s}\to \infty} w(s)\ =\ \sup_{s} w(s).$$
\end{lemma}
\begin{proof}
It suffices to prove the second equality: applying it to $-w$ in place of~$w$ gives the first one.
Set $\sigma:= \sup_{s} w(s)$. Let $\varepsilon\in\R^>$, and take $s_0\in\R^n$ with~$w(s_0) > \sigma-\varepsilon$.
By Corollary~\ref{cor:ap lower bd} applied to $s\mapsto v(s):=w(s)+\varepsilon-\sigma$ instead of~$w$
there are~$s$ with arbitrarily large $\abs{s}$ and
$v(s)\geq  v(s_0)/3>0$,   hence $w(s)>\sigma-\varepsilon$.
Thus~$\limsup\limits_{\abs{s}\to \infty} w(s) \geq \sigma$; the reverse inequality holds trivially.
\end{proof}

\subsection*{The mean value}  
{\it In this subsection $v$ and $w$ are bounded and measurable.}\/ 
If 
\begin{equation}\label{eq:mean value}
\lim_{T\to\infty} \frac{1}{T^n} \int_{[0,T]^n} w(s)\,ds
\end{equation}
exists \textup{(}in $\mathbb C$\textup{)}, then we say that $w$ has a mean value, and we
call  the quantity~\eqref{eq:mean value} the {\bf mean value}\index{mean value}\index{function!mean value} of $w$ and denote it by $M(w)$.
One verifies easily that if $v$ and~$w$ have a mean value, then so do the functions $v+w$, $cw$ ($c\in\mathbb C$),
and $\overline{w}$, with
%\marginpar{check $|\int_X v(s)ds|\le \int_X |v(s)|ds$ for bounded measurable $X\subseteq \R^n$} 
$$M(v+w)\ =\ M(v)+M(w),\quad M(cw)\ =\ cM(w),\quad\text{ and }\quad M(\overline{w})\ =\ \overline{M(w)}.$$ 
If  $w$ has a mean value, then $\abs{M(w)}\leq\dabs{w}$. If $w$ and $|w|$ have a mean value, then~$\abs{M(w)}\leq M(\abs{w})$.   
If $w$ has a mean value and $w(\R^n)\subseteq \R$, then $M(w)\in\R$, with $M(w)\geq 0$ if $w\big((\R^{\geq})^n\big)\subseteq\R^{\geq}$.

\begin{lemma}\label{lem:mean value add conj}
Let $d\in\R^n$. Then $w$ has a mean value iff $w_{+d}$ has a mean value,
in which case $M(w)=M(w_{+d})$.
\end{lemma}
\begin{proof}
It suffices to treat the case $d=(d_1,0,\dots,0)$, $d_1\in\R^{>}$.
For~$T > d_1$ we have
\begin{multline*}
\left|\int_{[0,T]^n} w_{+d}(s)\,ds - 
 \int_{[0,T]^n} w(s)\,ds \right|\ =\ \\
  \left| \int_{[T,d_1+T]\times [0,T]^{n-1}} w(s)\,ds - \int_{[0,d_1]\times [0,T]^{n-1}} w(s)\,ds\right|   \ \leq\ 2d_1\dabs{w}T^{n-1},
\end{multline*}
and this yields the claim.
\end{proof}

\begin{cor}
Suppose $w$ has a mean value, and let  $T_0\in\R^>$.
If $w(s)=0$ for all  $s\in (\R^{\geq})^n$ with $\abs{s}\geq T_0$, then $M(w)=0$.
If $w(\R^n)\subseteq\R$ and $w(s)\geq 0$ for all~$s\in (\R^{\geq})^n$ with $\abs{s}\geq T_0$, then $M(w)\geq 0$. \textup{(}As before, $|s|:=|s|_{\infty}$.\textup{)}
\end{cor}

\begin{lemma}\label{lem:mv and sup}
Suppose $w$ has a mean value and $w(\R^n)\subseteq\R$; then
$$\inf_{s} w(s)\ \leq\ \liminf_{\abs{s}\to \infty} w(s)\ \leq\ 
 M(w)\ \leq\ \limsup_{\abs{s}\to \infty} w(s)\ \leq\ \sup_{s} w(s).$$
\end{lemma}
\begin{proof}
The first and last inequalities  are clear.
Towards a contradiction assume~$L:=\limsup_{\abs{s}\to \infty} w(s)<M(w)$, and let $\varepsilon=\frac{1}{2}\big(M(w)-L\big)$. 
Take $T_0\in\R^>$ such that~$w(s)\leq M(w)-\varepsilon$ for all $s$ with $\abs{s}\geq T_0$.
The previous corollary applied to $s\mapsto{ M(w)-\varepsilon-w(s)}$ instead of $w$ implies $M(w)\leq M(w)-\varepsilon$, a contradiction.
This shows the third inequality; the second inequality is proved in a similar way.
\end{proof}

\noindent
Note that if $w$ has a mean value, then so does every $v$ having the same restriction to
$(\R^{\geq})^n$ as $w$, with $M(v)=M(w)$. 
% this allows us to define in the obvious way when
%a measurable function $u\colon U\to\mathbb C$, where $U\subseteq\R^n$ contains~$(\R^{\geq})^n$, has a mean value, and to define $M(u)$.

\begin{lemma}\label{lem:mv sequence}
Let $(v_m)$ be a sequence of bounded measurable functions $\R^n\to \mathbb C$ with a mean value, such that $\lim_{m\to \infty} \dabs{v_m-w}=0$. Then
$w$ has a mean value, and $\lim\limits_{m\to \infty} M(v_m)=M(w)$.
\end{lemma}
\begin{proof}
Let $\varepsilon\in\R^>$ be given, and take $m$ with 
$||v_m-w|| \le \varepsilon$. Since $v:=v_m$ has a mean value, we have $T_0\in\R^>$ such that for all $T_1,T_2\geq T_0$,
$$\left| \frac{1}{T_1^n} \int_{[0,T_1]^n} v(s)\,ds - \frac{1}{T_2^n} \int_{[0,T_2]^n} v(s)\,ds\right|\ \le\ \varepsilon.$$
Then for such $T_1,T_2$ we have
\begin{multline*}
\left| \frac{1}{T_1^n} \int_{[0,T_1]^n} w(s)\,ds - \frac{1}{T_2^n} \int_{[0,T_2]^n} w(s)\,ds\right|\	\leq\   
\frac{1}{T_1^n} \int_{[0,T_1]^n} \big| w(s)-v(s) \big|\,ds +  \\
 \left| \frac{1}{T_1^n} \int_{[0,T_1]^n} v(s)\,ds - \frac{1}{T_2^n} \int_{[0,T_2]^n} v(s)\,ds\right| + 
 \frac{1}{T_2^n} \int_{[0,T_2]^n} \big| w(s)-v(s) \big|\,ds 
\end{multline*}
where each term on the right of $\le$ is~$\le\varepsilon$.
Hence the limit  \eqref{eq:mean value} exists.
To show~$\lim\limits_{m\to \infty} M(v_m)=M(w)$, use
$|M(v_m)-M(w)|=|M(v_m-w)|\le \|v_m-w\|$. 
\end{proof}

\subsection*{The mean value of an almost periodic function}  
{\it In this subsection~$v$ and $w$ are almost periodic. As before, $\alpha$ ranges over $\R^n$.}\/

\begin{lemma}\label{lem:ap mean value}
Suppose $w(s)=\ex^{\imag (\alpha\cdot s)}$ for all $s$. Then $w$ has a mean value, with~$M(w)=1$ if $\alpha=0$ and $M(w)=0$ otherwise.
\end{lemma}
\begin{proof} This is obvious for $\alpha=0$. Assume $\alpha\ne 0$. Then
\begin{align*} \int_{[0,T]^n} \ex^{\imag (\alpha\cdot s)}\,ds\  &=\ T^{|\{j:\alpha_j=0\}|}\cdot \prod_{j, \alpha_j\ne 0} \frac{\ex^{\imag \alpha_j T}-1}{\imag \alpha_j},\quad \text{ so}\\  
\left| \frac{1}{T^n}\int_{[0,T]^n} \ex^{\imag (\alpha\cdot s)}\,ds  \right|\ &\ \leq\ \frac{1}{T^{\abs{\{j:\ \alpha_j\ne 0\}}}}\cdot \prod_{j, \alpha_j\ne 0}\frac{2}{|\alpha_j|},
\end{align*}
and thus $\frac{1}{T^n}\int_{[0,T]^n} \ex^{\imag (\alpha\cdot s)}\,ds \to 0$ as $T\to \infty$.
\end{proof}

\noindent
It follows that every trigonometric polynomial $w$ has a mean value. Using also Lemma~\ref{lem:mv sequence},
every almost periodic function $\R^n\to\mathbb C$ has a mean value.

\begin{lemma}\label{lem:mv periodic}
Suppose $u\colon\R^n\to\mathbb C$ is continuous and $1$-periodic. Then $u$ is almost periodic with mean value
$M(u)=\int_{[0,1]^n} u(s)\,ds$.
\end{lemma}
\begin{proof}
By Proposition~\ref{prop:StoWei}, $u$ is almost periodic. Now use that for $T\in \N^{\ge 1}$,
$$\qquad \qquad \qquad\qquad \int_{[0,T]^n} u(s)\,ds\ =\ T^n \int_{[0,1]^n} u(s)\,ds. \qquad \qquad \qquad \qquad\qquad\ \qedhere $$
\end{proof}

\begin{lemma}\label{lem:M(alpha w)}
Let $r\in (\R^\times)^n$. Then 
the almost periodic function $w_{\times r}$ has the same mean value as $w$.
\end{lemma}
\begin{proof}
Choose a sequence $(w_m)$ of trigonometric polynomials with $\dabs{w_m-w}\to 0$ as~$m\to \infty$.
Then $(w_m)_{\times r}$   is 
a trigonometric polynomial and $\dabs{(w_m)_{\times r}-w_{\times r}}\to 0$ as~$m\to \infty$.  
Lemma~\ref{lem:ap mean value} gives   $M\big((w_m)_{\times r}\big)=M(w_m)$;
now use Lemma~\ref{lem:mv sequence}.
\end{proof}

\begin{prop}[Bohr]\label{prop:Bohr}
Suppose  $w(\R^n)\subseteq\R^{\geq}$. If $M(w)=0$, then $w=0$.
\end{prop}
\begin{proof}
Suppose $s_0\in\R^n$ and $w(s_0)>0$. We claim that then $M(w)>0$.
Take $\delta_1$, $\ell_1$ as in Corollary~\ref{cor:ap lower bd}. Let $k$ range over~$\N^n$. Then
$\int_{\ell_1 k+[0,\ell_1]^n} w(s)\,ds \geq   \delta_1^n w(s_0)/3$
for all $k$, and hence for $m\ge 1$ and  $T=\ell_1 m$:
$$\frac{1}{T^n} \int_{[0,T]^n}w(s)\,ds\ =\ \frac{1}{T^n} \sum_{\abs{k}<m} \int_{\ell_1 k+[0,\ell_1]^n} w(s)\,ds\ \geq\ 
(\delta_1/\ell_1)^n w(s_0)/3.$$
Thus $M(w)\ \geq\ (\delta_1/\ell_1)^n w(s_0)/3\ >\ 0$.
\end{proof}

\noindent
By Proposition~\ref{prop:Bohr}, the map $(v,w)\mapsto \langle v,w\rangle:=M(v\overline{w})$ is  a positive definite hermitian form on the $\mathbb C$-linear space of almost
periodic functions~$\R^n\to\mathbb C$.
Lemma~\ref{lem:mean value add conj} yields $\langle v_{+d},w_{+d}\rangle = \langle v,w\rangle$ for $d\in\R^n$. 
By Lemma~\ref{lem:ap mean value}, the family~$\big(s\mapsto\ex^{(\alpha\cdot s)\imag}\big)_{\alpha}$ of trigonometric polynomials  is 
orthonormal with respect to~$\langle\ \, ,\  \rangle$. In particular, for a trigonometric polynomial $w$ as in \eqref{eq:trig poly} we have~$w_\alpha=\langle w,\ex^{(\alpha\cdot s)\imag}\rangle$, and thus: 

\begin{cor}\label{cor:id thm}
If $w=0$, then $w_\alpha=0$ for all $\alpha$.  
\end{cor}

\begin{cor}\label{cor:1-periodic trig poly} 
Suppose $w$ is a trigonometric polynomial as in \eqref{eq:trig poly}. Then $$w \text{ is $1$-periodic}\ \Longleftrightarrow\ 
w_\alpha=0 \text{ for all }\alpha\notin 2\pi\Z^n.$$
\end{cor}
\begin{proof}
If $w$ is $1$-periodic, then for $k\in\Z^n$ we have
$$w_\alpha\ =\ \langle w,\ex^{(\alpha\cdot s)\imag}\rangle\ =\ \langle w_{+k},(\ex^{(\alpha\cdot s)\imag})_{+k}\rangle\ =\ 
\ex^{-(\alpha\cdot k)\imag}\langle w, \ex^{(\alpha\cdot s)\imag}\rangle\ =\ \ex^{-(\alpha\cdot k)\imag}w_\alpha, $$
which for $w_{\alpha}\ne 0$ gives $\alpha\cdot k\in 2\pi\Z$ for all $k\in \Z^n$, and thus $\alpha\in 2\pi\Z^n$. 
This yields the forward implication, and the backward direction is obvious.
\end{proof}

\noindent
In the next corollary we equip $\R^n$ with the lexicographic ordering.

\begin{cor}
Suppose $w$ is a trigonometric polynomial. Then $w(\R^n)\subseteq\R$ iff there are $c\in\R$ and  $u_\alpha,v_\alpha\in\R$ for $\alpha > 0$,
with $u_\alpha=v_\alpha=0$ for all but finitely many $\alpha>0$, such that for all $s\in \R^n$,
\begin{equation}\label{eq:trig poly real}
w(s)\ =\ c+\sum_{\alpha>0} \big(u_\alpha\cos(\alpha\cdot s)+ v_\alpha\sin(\alpha\cdot s)\big). 
\end{equation}
Moreover, in this case $c$ and the coefficients $u_\alpha$, $v_\alpha$ are unique, and
$w$ is $1$-periodic iff $u_\alpha=v_\alpha=0$ for all $\alpha>0$ with $\alpha\notin 2\pi\Z^n$.
\end{cor}
\begin{proof}
Clearly if $w$ has stated form, then $w(\R^n)\subseteq\R$. For the converse, suppose~$w(\R^n)\subseteq\R$, and $w$ is given as in
\eqref{eq:trig poly}.
Then for $s\in \R^n$, $$\sum_\alpha \overline{w_\alpha}\ex^{-(\alpha\cdot s)\imag}\ =\ \overline{w}(s)\ =\ w(s)\ =\ \sum_\alpha w_\alpha\ex^{(\alpha\cdot s)\imag},$$
and hence $w_0\in\R$ and $\overline{w_\alpha}=w_{-\alpha}$ for $\alpha>0$, by Corollary~\ref{cor:id thm}, so
$$w(s)\ =\ w_0+\sum_{\alpha>0} \big(w_\alpha\ex^{(\alpha\cdot s)\imag} + \overline{w_\alpha}\ex^{-(\alpha\cdot s)\imag}\big)\qquad (s\in\R^n). $$
Put $c:=w_0$ and $u_\alpha=\Re(2w_\alpha)$, $v_\alpha:=\Im(2w_\alpha)$ for $\alpha>0$. 
Then \eqref{eq:trig poly real} holds for~$s\in \R^n$. The rest follows from Corollaries~\ref{cor:id thm} and~\ref{cor:1-periodic trig poly}.
\end{proof}

\section{Uniform Distribution Modulo One}\label{sec:udmod1}

\noindent
In this section we collect some basic facts about uniform distribution modulo $1$ of functions as needed later.
Our main references are \cite[Chapter~1, \S{}9]{KuipersNiederreiter} and \cite{BoshernitzanUniform}.

\subsection*{Natural density}
Below $\R$ is given its usual Le\-bes\-gue measure, {\em measurable\/}   means
{\em Le\-bes\-gue-measurable}, and $\mu$ denotes the Le\-bes\-gue measure on $\R$.
By an ``interval'' we mean here a set $I=[a,b)$ where $a,b\in\R$, $a<b$, so $\mu(I)=b-a$.
{\it In the rest of this subsection $I$ is an interval and $X$, $Y$  are measurable subsets of~$\R$.}\/ We let
$$\rho(I,X)\ :=\ \frac{\mu(I\cap X)}{\mu(I)}\in [0,1]$$
be the {\bf density of $X$ in $I$}.\index{density} So $\rho(I,X)=0$ if $I\cap X=\emptyset$ and $\rho(I,X)=1$ if~${I\subseteq X}$,
and  $\rho(I+d,X+d)=\rho(I,X)$ for~$d\in\R$. 
Clearly $\rho(I,X)\leq\rho(I,Y)$ if $X\subseteq Y$, and 
if $(X_n)$ is a family of pairwise disjoint measurable subsets of $\R$ and
$X=\bigcup_n X_n$, then $\rho(I,X)=\sum_n \rho(I,X_n)$. 

\medskip
\noindent
Let $X\triangle Y:=(X\setminus Y)\cup(Y\setminus X)$ be the symmetric difference of~$X$,~$Y$.
If $\mu(X)<\infty$ and $\mu(Y)<\infty$, then 
$\mu(X)-\mu(Y)\le \mu(X\setminus Y)$ and $\abs{\mu(X)-\mu(Y)}\leq \mu(X\triangle Y)$, so

\begin{lemma}\label{lem:rho symm diff}
$\abs{\rho(I,X)-\rho(I,Y)} \leq \rho(I,X\triangle Y)$.
\end{lemma}

\noindent
Moreover:

\begin{lemma}\label{lemtrd}
Let $d\in \R$; then  $\big|\rho(I,X)-\rho(I+d,X)\big|\leq \abs{d}/\mu(I)$.
\end{lemma}
\begin{proof}
We need to show $\big|\mu(I\cap X)-\mu\big((I+d)\cap X\big)\big| \leq \abs{d}$. Replacing $I$ and $d$ by~$I+d$ and $-d$, if necessary, we
arrange $d\geq 0$. 
Then 
$$-\mu(I)\ =\ -\mu(I+d)\ \leq\ \mu(I\cap X)-\mu\big((I+d)\cap X\big)\ \leq\  \mu(I),$$
hence we are done if $\mu(I)\le d$. Suppose $\mu(I)>d$ and let $I=[a,b)$, $a,b\in\R$; so~$\mu(I)=b-a$.  Then $I\setminus (I+d)=[a,a+d)$ and $(I+d)\setminus I=[b,b+d)$,  hence 
$$-d\ =\ -\mu\big((I+d)\setminus I\big)\ \leq\ \mu(I\cap X)-\mu\big((I+d)\cap X\big)\ \leq\ \mu\big( I\setminus (I+d) \big)\ =\ d$$
as required.
\end{proof}

\noindent
Let $\rho$ range over $[0,1]$ and $T$ over $\R^{>}$. Lemma~\ref{lemtrd} gives:

\begin{cor}
The following conditions on $X$ are equivalent:
\begin{enumerate}
\item[\rm{(i)}] $\lim_{T\to\infty} \rho\big([0,T), X\big)\ =\ \rho$;
\item[\rm{(ii)}] for all $a\in\R$, $\lim_{T\to\infty} \rho\big([a,a+T), X\big)\ =\ \rho$;
\item[\rm{(iii)}] for some $a\in\R$, $\lim_{T\to\infty} \rho\big([a,a+T), X\big)\ =\ \rho$. 
\end{enumerate}
\end{cor}

\noindent
We say that $X$ has  {\bf natural density~$\rho$ at $+\infty$}\/ (short: $X$ has density~$\rho$)\index{density!natural}   if
one of the equivalent conditions in the corollary above holds, and in this case we set~$\rho(X):=\rho$. 
If $X$ has an upper bound in $\R$, then $\rho(X)=0$, whereas if $X$ contains a half-line~$\R^{\geq a}$~($a\in\R$), then $\rho(X)=1$. By
Lemma~\ref{lem:rho symm diff} we have:

\begin{cor}
If $X$ has density $\rho$ and $X\triangle Y$ has density~$0$, then $Y$ has density~$\rho$.
\end{cor}

\noindent
In particular,  the density of $X$ only depends on the germ of~$X$ at $+\infty$, in the following sense:
if  $X\cap\R^{> a}=Y\cap\R^{> a}$ for some~$a\in\R$, then
$X$ has   density~$\rho$ iff~$Y$ has density~$\rho$. 
The collection of measurable subsets of $\R$ that have a density is a boolean algebra of subsets of $\R$,
and $X\mapsto\rho(X)$ is a finitely additive measure on this boolean algebra taking values in $[0,1]$.
If $X$ has a density and $d\in\R$, then~$X+d$ has the same density.

\subsection*{Uniform distribution~$\operatorname{mod}\,1$}
Let $f\colon\R^{\ge a}\to\R$ ($a\in\R$) be measurable.
For $t\in\R$ we let $\{t\}$ be the fractional part of $t$: the element of $[0,1)$ such that $t\in\Z+\{t\}$.
Let  $Y\subseteq [0,1)$ be measurable; then $Y+\Z$ is measurable and hence so is
$$f^{-1}(Y+\Z)\ =\ \big\{t\in\R^{\ge a}: \{f(t)\}\in Y\big\}.$$ 
For $a\le b<c$  we have
$$\mu\big( [b,c) \cap f^{-1}(Y+\Z) \big)\ =\ \int_b^c \chi_Y\big( \{f(t)\} \big) \,dt.$$
Let $\rho\in\R$; then $f^{-1}(Y+\Z)$ has density $\rho$ iff for some $b \ge a$ we have
$$\lim_{T\to\infty} \frac{1}{T} \int_b^{b+T} \chi_Y\big(\{f(t)\}\big)dt\  =\ \rho,$$
and in this case the displayed identity holds for all $b \ge a$. Hence if
$f^{-1}(Y+\Z)$ has density $\rho$ and 
$g\colon\R^{\ge  b}\to\R$ ($b\in\R$) is measurable with the same germ at $+\infty$ as~$f$, then~$g^{-1}(Y+\Z)$ also has density $\rho$.

\begin{definition}
We say that $f$ is {\bf uniformly distributed~$\operatorname{mod}\,1$} (abbreviated: u.d.~$\operatorname{mod}\,1$) if for every interval $I\subseteq [0,1)$ the set
$f^{-1}(I+\Z)$ has density $\mu(I)$.\index{function!uniformly distributed~$\operatorname{mod}\,1$}\index{germ!uniformly distributed~$\operatorname{mod}\,1$}\index{uniformly distributed~$\operatorname{mod}\,1$}
%we have $\rho(f^{-1}\big(I+\Z)\big)=\mu(I)$.
\end{definition}

\noindent
The function $f: \R^{\ge} \to \R$ with $f(t)=t$ for all $t\ge 0$ has $f^{-1}(I+\Z)=I+\N$ for $I$ as above, so $f$ is u.d~mod~$1$.
By the remarks above, if $f$ is u.d.~mod~$1$, then so is any measurable function ${\R^{\ge b}\to\R}$ ($b\in\R$) 
with the same germ at $+\infty$ as $f$. 
%Hence we may define: a germ in~$\mathcal G$ is
% u.d.~$\operatorname{mod}\,1$ if it has a representative $\R^{>b}\to\R$  ($b\in\R$)
%which is u.d.~$\operatorname{mod}\,1$.
If $f$ is  u.d.~$\operatorname{mod}\,1$ and
eventually   increasing or eventually   decreasing,
then~$\abs{f(t)} \to +\infty$ as $t\to+\infty$. 
If $f$ is  u.d.~$\operatorname{mod}\,1$, then so are
the functions $t\mapsto k\cdot f(t)\colon\R^{\ge a}\to\R$ for~$k\in\Z^{\neq}$,  and
$t\mapsto f(d+t)\colon\R^{\ge(a-d)}\to\R$ with $d\in\R$.   

\subsection*{The Weyl Criterion}
In this subsection we fix a measurable function $f\colon\R^{\ge}\to\R$. 
For a bounded measurable function $w\colon [0,1]\to\R$ we consider the relation
\begin{equation}\label{eq:udmod1 w}\tag{W}
\lim_{T\to \infty} \frac{1}{T} \int_0^{T} w\big(\{f(t)\}\big)dt\ =\ \int_0^1 w(s)\,ds.
\end{equation}  
Then $f$ is u.d.~$\operatorname{mod}\,1$ iff \eqref{eq:udmod1 w} holds whenever $w=\chi_I$ is the characteristic function of
some interval $I\subseteq [0,1]$.
It follows that if $f$ is u.d.~$\operatorname{mod}\,1$ and 
$w\colon [0,1]\to \R$ is a step function (that is, an $\R$-linear combination of characteristic functions $\chi_I$ of intervals $I\subseteq [0,1]$), then
\eqref{eq:udmod1 w} holds.

\begin{lemma}\label{lem:udmod1 w}
Let $w\colon [0,1]\to\R$ be bounded and measurable, and suppose that for every $\varepsilon\in\R^>$ there
are bounded measurable functions $w_1,w_2\colon [0,1]\to\R$ such that 
\begin{enumerate}
\item[\rm{(i)}] $w_1\leq w\leq w_2$ on $[0,1)$, 
\item[\rm{(ii)}] $\int_0^1 \big(w_2(s)-w_1(s)\big)\,ds\leq\varepsilon$, and 
\item[\rm{(iii)}] for $i=1,2$,
 \eqref{eq:udmod1 w} holds for $w_i$ instead of $w$.  
\end{enumerate}
Then \eqref{eq:udmod1 w} holds.
\end{lemma}
\begin{proof}
Given $\varepsilon\in\R^>$ and $w_1,w_2\colon [0,1]\to\R$  satisfying (i), (ii), (iii),
we have \begin{align*}
\int_0^1 w(s)\,ds-\varepsilon &\ \leq\ \int_0^1 w_1(s)\,ds \ =\ \lim_{T\to\infty} \frac{1}{T}  \int_0^{T} w_1\big(\{f(t)\}\big)dt \\
&\ \leq\ \liminf_{T\to \infty} \frac{1}{T}  \int_0^{T} w\big(\{f(t)\}\big)dt\ \leq\ 
\limsup_{T\to \infty} \frac{1}{T}  \int_0^{T} w\big(\{f(t)\}\big)dt \\
&\ \leq\ \lim_{T\to \infty} \frac{1}{T} \int_0^T w_2\big(\{f(t)\}\big)dt\ =\ \int_0^1 w_2(s)\,ds \\ &\ \leq \ \int_0^1 w(s)\,ds+\varepsilon. \qedhere
\end{align*}
\end{proof}

\begin{prop} \label{prop:udmod1 test}
 $f$ is u.d.~$\operatorname{mod}\,1$ iff  \eqref{eq:udmod1 w} holds for all con\-ti\-nu\-ous~${w\colon [0,1]\to\R}$. 
\end{prop}
\begin{proof}
If $w\colon [0,1]\to\R$ is continuous, then partitioning $[0,1)$
into intervals as in Riemann integration we obtain for any $\varepsilon\in\R^>$ step functions $w_1,w_2\colon [0,1]\to\R$ such that 
$w_1\leq w\leq w_2$ on $[0,1)$  and $\int_0^1 \big(w_2(s)-w_1(s)\big)\,ds\leq\varepsilon$.  
Moreover, if $I\subseteq [0,1)$ is an interval, then  for any $\varepsilon\in\R^>$ there are continuous functions $w_1,w_2\colon [0,1]\to\R$ with 
$w_1\leq \chi_I \leq w_2$ and $\int_0^1 \big( w_2(s)-w_1(s) \big)ds\leq\varepsilon$. 
The proposition  follows from these facts and Lemma~\ref{lem:udmod1 w}.
\end{proof}

\noindent
It is convenient to extend the notion of mean value to bounded measurable
functions~$g\colon \R^{\ge} \to \mathbb{C}$: for such $g$, if $\lim\limits_{T\to \infty}\frac{1}{T}\int_0^T g(t)\,dt$ exists in $\mathbb{C}$, then we say that~$g$ has mean value 
$M(g):=\lim\limits_{T\to \infty}\frac{1}{T}\int_0^T g(t)\,dt$.
  
\begin{cor} \label{cor:udmod1 test} 
The following conditions on $f$ are equivalent:
\begin{enumerate}
\item[\rm{(i)}] $f$ is u.d.~$\operatorname{mod}\,1$; 
\item[\rm{(ii)}] $\displaystyle\lim_{T\to \infty} \frac{1}{T} \int_0^{T} (w\circ f)(t)\,dt=\int_0^1 w(s)\,ds$ for all continuous $1$-periodic functions $w\colon \R\to\mathbb C$;
\item[\rm{(iii)}] for every continuous $1$-periodic
$w\colon \R\to\mathbb C$, the function  ${w\circ f}\colon\R^{\geq}\to\mathbb C$ 
has mean value $M(w\circ f)=M(w)$.
\end{enumerate}
\end{cor}
\begin{proof}
We first show (i)~$\Leftrightarrow$~(ii).
For the forward direction, apply Proposition~\ref{prop:udmod1 test} to
the real and imaginary parts of $w$, using $w(\{t\})=w(t)$ for $t\in\R$ and $1$-periodic~$w$.
The converse follows from Lemma~\ref{lem:udmod1 w} and the observation that 
if $I\subseteq [0,1)$ is an interval, then  for any $\varepsilon\in\R^>$ we can take 
continuous  functions $w_1,w_2\colon [0,1]\to\R$ with 
$w_1\leq \chi_I \leq w_2$ and $\int_0^1 \big( w_2(s)-w_1(s) \big)ds\leq\varepsilon$ as in the proof of the proposition above, such that
in addition $w_i(0)=w_i(1)$ for $i=1,2$, and then $v_i\colon \R\to\R$ given by $v_i(t)=w_i(\{t\})$ for $t\in\R$ ($i=1,2$)
is continuous and $1$-periodic. 
The equivalence of (ii) and (iii) is immediate from Lemma~\ref{lem:mv periodic}.
\end{proof}

\begin{theorem}[{Weyl~\cite{Weyl}}]\label{thm:Weyl} The function
 $f$ is u.d.~$\operatorname{mod}\,1$ iff for all $n \geq 1$ we have
\begin{equation}\label{eq:Weyl crit}
\lim_{T\to \infty} \frac{1}{T} \int_0^{T} \ex^{2\pi\imag n f(t)}\,dt\ =\ 0.
\end{equation}
\end{theorem}
\begin{proof} 
The forward direction follows from Corollary~\ref{cor:udmod1 test}.
Conversely, sup\-pose that~\eqref{eq:Weyl crit} holds for all $n\geq 1$.
Note that then  for all $k\in \Z^{\ne}$,
\[%\begin{equation}\label{eq:Weyl crit, gen}
\lim_{T\to \infty} \frac{1}{T} \int_0^{T} \ex^{2\pi\imag k f(t)}\,dt\ =\ 0.
\]%\end{equation}
Thus by Corollary~\ref{cor:1-periodic trig poly}, every $1$-periodic trigonometric polynomial $v\colon\R\to\mathbb C$ gives a function $v\circ f$ with mean value  
$M(v\circ f)=M(v)$.
Now let $w\colon \R\to\mathbb C$ be continuous and $1$-periodic.
Proposition~\ref{prop:StoWei} yields a sequence $(v_m)$ of $1$-periodic trigonometric
polynomials $\R\to \mathbb C$ with $\dabs{v_m-w}\to 0$ as $m\to\infty$.
So $M(v_m)\to M(w)$ as~$m\to\infty$, by Lemma~\ref{lem:mv sequence}.
Extend $f$ to a measurable function $\R\to\R$, also denoted by $f$.
Then $\dabs{(v_m\circ f)-(w\circ f)}\to 0$ as $m\to \infty$. 
Hence by Lemma~\ref{lem:mv sequence} again, 
$w\circ f$ has a mean value and
$M(v_m)=M(v_m\circ f)\to M(w\circ f)$ as $m\to \infty$. Therefore~$M(w\circ f)=M(w)$.
Hence  $f$ is u.d.~$\operatorname{mod}\,1$ by Corollary~\ref{cor:udmod1 test}.
%Let $w\colon \R\to\mathbb C$ be continuous and $1$-periodic, and let $\varepsilon\in\R^>$. By a special case of the Stone-Weierstrass  Theorem~\cite[Theorem~4.25]{BabyRudin} there are $c_1,\dots,c_m\in\mathbb C$ and $k_1,\dots,k_m\in\Z$ such that with  $v(t):=c_1\ex^{2\pi\imag k_1 t}+\cdots+c_n\ex^{2\pi\imag k_m t}$ ($t\in\R$) we have  $\sup_{t\in [0,1]} \abs{ v(t) - w(t) } \leq \varepsilon$.
%Then for all $T$ we have
%\begin{align*}
%\left|\frac{1}{T} \int_0^{T} (w\circ f)(t)\,dt - \int_0^1 w(s)ds   \right|\	&\ \leq\ \left|\frac{1}{T} \int_0^{T} \big( (w\circ f) - (v\circ f)\big)(t)\,dt \right|\ +  \\
%&\hskip1.8em\left|\frac{1}{T}\int_0^T (v\circ f)(t)\,dt-\int_0^1 v(s)\,ds\right|\ +  \\
%&\hskip1.8em\left|\int_0^1 \big(v(s)-w(s)\big)\,ds \right|,
%\end{align*}
%where the first and the third summand on the right are~$\leq\varepsilon$ by  choice of $v$, whereas for some $T_0>0$, the second summand is $\leq\varepsilon$ for all $T\geq T_0$, by \eqref{eq:Weyl crit, gen}. Hence  $f$ is u.d.~$\operatorname{mod}\,1$ by Corollary~\ref{cor:udmod1 test}.
\end{proof}

\begin{remark}
Let $w(s)=\ex^{2\pi\imag s}$  ($s\in\R$)  and let $g\colon\R\to\R$ be a continuous function whose restriction to~$\R^{\geq}$ is u.d.~$\operatorname{mod}\,1$.
By  Corollary~\ref{cor:udmod1 test}, $w\circ g$ has a mean value. When is~$w\circ g$  almost periodic? This happens only for very special $g$:
%\marginpar{not checked what follows}
if
$w\circ g$ is almost periodic, then
there are $r\in\R$ and an almost periodic  $h\colon\R\to\R$ such that~$g(t)=rt+h(t)$ for all $t\in\R$, by a theorem of Bohr~\cite{Bohr30}. Moreover,~\cite[Theorem~1]{NS} says that if $h\colon\R\to\R$ is almost periodic, then
for all but countably many  $r\in\R$ the function  $t\mapsto rt+h(t)\colon\R^{\geq}\to\R$ is 
  u.d.~$\operatorname{mod}\,1$. These facts are not used later.
\end{remark}

\subsection*{Uniform distribution~$\operatorname{mod}\,1$ of differentiable functions}
Let $f\colon\R^{\geq a}\to\R$ ($a\in\R$) be continuously differentiable. We
give here sufficient conditions for $f$ to be u.d.~$\operatorname{mod}\,1$
and for $f$ not to be u.d.~$\operatorname{mod}\,1$.
First a lemma in the spirit of Corollary~\ref{cor:HL}:

\begin{lemma}\label{lem:twice diff}
Let $F\colon \R^{>}\to\R$ be twice continuously differentiable such that ${F(t)/t\to 0}$ as $t\to+\infty$. Assume $t\mapsto tF''(t)\colon\R^>\to\R$
is bounded.  Then $F'(t)\to 0$ as~$t\to+\infty$.
\end{lemma}
\begin{proof}
Let $t,\eta>0$.
Taylor's Theorem~\cite[Theorem~19.9]{Bartle} yields $\theta\in [0,1]$ such that 
$$F(t+\eta)-F(t)\ =\ \eta F'(t) + \textstyle\frac{1}{2}\eta^2F''(t+\theta\eta),$$
and thus
$$F'(t)\ =\ \frac{F(t+\eta)-F(t)}{\eta} - \textstyle\frac{1}{2}\eta F''(t+\theta\eta).$$
Take $M\in\R^{>}$ such that $\abs{ t F''(t) } \leq M$ for all $t\in\R^>$.
Let $\varepsilon\in\R^>$, and set $\delta:=\varepsilon/M$. 
Then for all $t>0$, $\eta=\delta t$ yields $\theta=\theta_t\in [0,1]$ with
$$F'(t)\ =\ \left( \frac{F(t+\delta t)}{t+\delta t} \cdot \frac{1+\delta}{\delta} - \frac{F(t)}{t} \cdot \frac{1}{\delta} \right) - \frac{\delta}{2(1+\theta\delta)} (t+\theta\delta t)F''(t+\theta\delta t).$$
The difference in the parentheses tends to zero as $t\to \infty$ whereas
the remaining term is $\le \varepsilon/2$ in absolute value for all $t\in\R^>$.
\end{proof}

\begin{prop}[{Kuipers-Meulenbeld~\cite{KM}}]\label{prop:not udmod1}
Suppose the function  $$t\mapsto f'(t)t\ \colon\ \R^{\ge a}\to\R$$ is boun\-ded. Then $f$ is not 
u.d.~$\operatorname{mod}\,1$.
\end{prop}
\begin{proof}
Replacing $f$ by $t\mapsto f(a+t)\colon\R^{\ge}\to\R$ we arrange $a=0$.
Assume towards a contradiction that  \eqref{eq:Weyl crit} holds for $n=1$, and consider  $F\colon\R^{>}\to\R$ given by
$$F(t)\ :=\ \Re\left(\int_0^t \ex^{2\pi\imag f(s)}\,ds\right)\ =\ \int_0^t \cos\big( 2\pi f(s)\big) \, ds.$$
Then $F$ is twice continuously differentiable with 
$$F'(t)\ =\ \cos\big( 2\pi f(t)\big),\quad F''(t)\ =\ -2\pi f'(t)\sin\big( 2\pi f(t)\big),$$
and 
$F(t)/t\to 0$ as $t\to\infty$ and $t\mapsto tF''(t)\colon\R^>\to\R$
is bounded.
Hence by Lem\-ma~\ref{lem:twice diff} we have $\cos\big( 2\pi f(t)\big)\to 0$ as $t\to\infty$; likewise we show $\sin\big( 2\pi f(t)\big)\to 0$ as $t\to\infty$.
Hence $\ex^{2\pi\imag f(t)}\to 0$ as $t\to\infty$, a contradiction.
\end{proof}

\noindent
In the next proposition we assume $a = 0$ and consider the continuously differentiable function
$t\mapsto g(t):=f(\ex^t)\colon \R\to\R$ (so $f(t)=g(\log t)$ for $t\in\R^>$).

\begin{prop}[{Tsuji~\cite{Tsuji}}]\label{prop:Tsuji}
Suppose  $g$ and  $g'$ are eventually strictly increasing with $g(t)/t\to+\infty$ as $t\to+\infty$. Then $f$ is u.d.~$\operatorname{mod}\,1$.
\end{prop}
\begin{proof}
Let $n\geq 1$; we claim that \eqref{eq:Weyl crit} holds. 
The continuous functions
\begin{align*} t &\mapsto\varphi(t):=2\pi n f(t)\ \colon\ \R^{\ge}\to\R,\\
t&\mapsto \gamma(t):=\varphi'(t)t=2\pi n\, g'(\log t)\ \colon\ \R^>\to\R
\end{align*}  are eventually strictly increasing.
We have $\varphi(t)/\log t\to +\infty$ as $t\to+\infty$. There\-fore~$\gamma(t)\to+\infty$ as $t\to+\infty$:
otherwise $\varphi'(t)\le M/t$ for all $t\ge b$, and some~${b, M>0}$, and
then integration gives $\varphi(t)=O(\log t)$ as $t\to +\infty$, a contradiction.

Take $a_0\in \R^{>}$ such that
$\varphi$ and $\gamma$ are strictly increasing on $\R^{\ge a_0}$ and $\gamma(a_0)>0$.
Set $\rho_0=\varphi(a_0)$, and take
$\eta\colon\R^{\geq\rho_0}\to\R^{\geq a_0}$ so that $(\eta\circ\varphi)(t)=t$ for $t\in\R^{\geq a_0}$.
Then~$\eta'\big(\varphi(t)\big)>0$ and
$\gamma(t)=\eta\big(\varphi(t)\big)/\eta'\big(\varphi(t)\big)$    for $t>a_0$. Hence the function 
$$u\mapsto \eta^\dagger(u)\ :=\ \eta'(u)/\eta(u)\ \colon\ \R^{>\rho_0}\to\R^>$$ is strictly decreasing
with $\lim\limits_{u\to+\infty}\eta^\dagger(u)=0$.
Let now $T > a_0$ and consider
$$I(T)\ :=\ \int_{a_0}^T \sin\varphi(t)\,dt.$$
Set $\rho_T=\varphi(T)$, and let $\tau\in (\rho_0,\rho_T)$.
Substituting $u=\varphi(t)$ gives
$$I(T)\ =\ \int_{\rho_0}^{\rho_T} \eta'(u)\sin u\,du\ =\ \int_{\rho_0}^\tau \eta'(u)\sin u\,du +
\int_\tau^{\rho_T} \eta(u)\eta^\dagger(u)\sin u\,du.$$
Two applications of the Second Mean Value Theorem for Integrals \cite[Theorem~23.7]{Bartle} yield first $\tau_2$ and then~$\tau_1$
such that $\tau\le \tau_1\le \tau_2\le \rho_T$ and
$$\int_\tau^{\rho_T} \eta(u)\eta^\dagger(u)\sin u\,du\ =\ 
\eta^\dagger(\tau) \int_\tau^{\tau_2} \eta(u) \sin u\,du\ =\
\eta^\dagger(\tau) \eta(\tau_2) \int_{\tau_1}^{\tau_2} \sin u\,du,$$
hence
$$\int_\tau^{\rho_T} \eta(u)\eta^\dagger(u)\sin u\,du\ =\ \eta^\dagger(\tau)\,C\qquad\text{where $\abs{C}\leq 2\eta(\rho_T)=2T$.}$$
Let now $\varepsilon\in\R^>$ be given. Take $\tau>\rho_0$ so large that
$\eta^\dagger(\tau)\leq\varepsilon/4$. Then for $T>a_0$ so large that  
$\varphi(T) > \tau$, and
$$\left| \int_{\rho_0}^{\tau} \eta'(u)\sin u\,du \right|\ \leq\ \varepsilon T/2$$
we have $\abs{I(T)} \leq \varepsilon T$. Thus as $T\to \infty$ we have 
$$\frac{1}{T}\int_0^T\sin\varphi(t)\,dt\ =\ \frac{1}{T}\int_0^{a_0}\sin\varphi(t)\,dt + \frac{I(T)}{T}\ \to\  0.$$ 
Likewise, $\frac{1}{T}\int_0^T\cos\varphi(t)\,dt\to 0$
as $T\to\infty$. Thus \eqref{eq:Weyl crit} is satisfied. 
\end{proof}

\begin{theorem}[{Boshernitzan~\cite{BoshernitzanUniform}}]\label{thm:udmod1 Boshernitzan}
Suppose the germ of $f$, also denoted by $f$, lies in a Hardy field $H$. Then:
$f$ is u.d.~$\operatorname{mod}\,1$ $\Longleftrightarrow$ $f\succ \log x$.
\end{theorem}
\begin{proof}
By increasing $H$ we arrange $\log x\in H$.
The claim is obvious if $f\preceq 1$, since then  $f$ is neither $\succ\log x$ nor u.d.~$\operatorname{mod}\,1$.
So suppose $f\succ 1$; then
$f\succ\log x$ iff~$f'\succ 1/x$.
If $f' \preceq 1/x$, then $f$ is not u.d.~$\operatorname{mod}\,1$ by Proposition~\ref{prop:not udmod1}.
Suppose~$f'\succ 1/x$;
to show that $f$  is u.d.~$\operatorname{mod}\,1$ we replace $f$  by $t\mapsto f(a+t)\colon\R^{\ge}\to\R$ and  $H$ by~$H\circ (a+x)$ to
arrange $a=0$. Replacing $f$ by $-f$ if necessary we also arrange $f>0$
in $H$.
Then $f(t)\to+\infty$ as $t\to+\infty$, hence
 the function $t\mapsto g(t):=f(\ex^t)\colon \R\to\R$ is eventually strictly increasing and
its germ, also denoted by $g$, lies in some Hardy field and  satisfies $g\succ x$; thus $g'\succ 1$, hence
$t\mapsto g'(t)\colon\R\to\R$ is also eventually strictly increasing. Thus $f$ is u.d.~$\operatorname{mod}\,1$ by Proposition~\ref{prop:Tsuji}. 
\end{proof}

\noindent 
In particular, if $f$ is u.d.~$\operatorname{mod}\,1$ and its germ lies in a Hardy field, then $\alpha f$ is
u.d.~$\operatorname{mod}\,1$ for every   $\alpha\in\R^\times$.

\subsection*{Uniform distribution~$\operatorname{mod}\,1$  in higher dimensions}
In this subsection $n\geq 1$, $\R^n$ is equipped with its usual Le\-bes\-gue measure $\mu_n$, and {\em measurable\/} for a subset of $\R^n$   means measurable with respect to $\mu_n$. Let $a\in\R$ and consider measurable functions $f_1,\dots, f_n\colon \R^{\ge a} \to \R$, which we combine into a single map $$f\ : =\ (f_1,\dots,f_n)\ \colon\ \R^{\ge a}\to\R^n.$$ 
By a {\bf box} (in $\R^n$) we mean a set $I=I_1\times\cdots\times I_n$ where $I_1,\dots,I_n$ are intervals, so~$\mu_n(I)=\mu(I_1)\cdots\mu(I_n)$ and
$f^{-1}(I+\Z^n)= \bigcap_{j=1}^n f_j^{-1}(I_j+\Z)$ is measurable.

\begin{definition}
We say that $f$ is {\bf uniformly distributed~$\operatorname{mod}\,1$} (abbreviated: u.d.~$\operatorname{mod}\,1$) if for every box $I\subseteq [0,1)^n$ the set $f^{-1}(I+\Z^n)$ has density $\mu_n(I)$.\index{function!uniformly distributed~$\operatorname{mod}\,1$}
\end{definition}

\noindent
For $s=(s_1,\dots,s_n)\in\R^n$, set $\{s\}:=\big(\{s_1\},\dots,\{s_n\}\big)\in [0,1)^n$. With this notation,
$f$ is u.d.~$\operatorname{mod}\,1$ iff for every box $I\subseteq [0,1)^n$ we have
$$\lim_{T\to\infty} \frac{1}{T} \int_a^{a+T} \chi_I\big(\{f(t)\}\big)dt\ =\ \mu_n(I).$$
Let $b\in\R^{\geq a}$, $d\in \R$. Then $f$ is u.d.~$\operatorname{mod}\,1$ iff the restriction of $f$ to $\R^{\ge b}$ is u.d.~$\operatorname{mod}\,1$, and if $f$ is u.d.~$\operatorname{mod}\,1$, then so is $t\mapsto f(d+t)\colon\R^{\geq(a-d)}\to\R^n$.

\medskip
\noindent
{\it In the rest of this subsection we  assume $a=0$.}\/
Proposition~\ref{prop:udmod1 test} and its Corollary~\ref{cor:udmod1 test} generalize to this setting: 

\begin{prop} \label{prop:udmod1 test, n}
The map $f$ is u.d.~$\operatorname{mod}\,1$ if and only if 
$$\displaystyle\lim_{T\to\infty} \frac{1}{T} \int_0^{T} w\big(\{f(t)\}\big)dt\ =\ \int_{[0,1]^n} w(s)\,ds$$ 
for every continuous function $w\colon [0,1]^n\to\R$.
\end{prop}
\begin{proof} 
First we use the proof of Lemma~\ref{lem:udmod1 w} to obtain the  analogue of that lemma for bounded measurable functions $w\colon [0,1]^n\to\R$.
Now, given a continuous function~$w\colon [0,1]^n\to\R$  and  $\varepsilon\in\R^>$,  there are $\R$-linear combinations~${w_1,w_2\colon [0,1]^n\to\R}$ of characteristic functions of pairwise disjoint boxes contained in $[0,1]^n$ such that $w_1\leq w\leq w_2$ on $[0,1)^n$ and $\int_{[0,1]^n} \big({w_2(s)-w_1(s)}\big)\,ds\leq\varepsilon$. This gives one direction.

Next, let $I=I_1\times\cdots\times I_n\subseteq [0,1)^n$ be a box and $\varepsilon\in\R^>$. 
For $j=1,\dots,n$ we have continuous functions $w_{1j},w_{2j}\colon [0,1]\to\R^{\ge}$ such that 
$$0\ \le\ w_{1j}\ \leq\ \chi_{I_j}\ \leq\ w_{2j}\ \le\ 1\quad\text{and}\quad  \int_0^1 \big( w_{2j}(t)-w_{1j}(t) \big)dt\leq \varepsilon/2^n.$$ 
For $s=(s_1,\dots,s_n)\in [0,1]^n$ set $w_{i}(s):=w_{i1}(s_1)\cdots w_{in}(s_n)$.
Then the functions~$w_1,w_2\colon [0,1]^n\to\R$  are continuous with
$$w_1\leq \chi_I \leq w_2\quad\text{and}\quad {\int_{[0,1]^n} \big( w_2(s)-w_1(s) \big)ds}\leq\varepsilon.$$
The proposition follows from these facts just as in the proof of Proposition~\ref{prop:udmod1 test}. 
\end{proof}

\noindent
As Proposition~\ref{prop:udmod1 test} led to Corollary~\ref{cor:udmod1 test}, so does Proposition~\ref{prop:udmod1 test, n} give: 

\begin{cor} \label{cor:udmod1 test, n} 
The following conditions on $f$ are equivalent:
\begin{enumerate}
\item[\rm{(i)}] the map $f$ is u.d.~$\operatorname{mod}\,1$; 
\item[\rm{(ii)}] $\displaystyle\lim_{T\to\infty} \frac{1}{T} \int_0^{T} (w\circ f)(t)\,dt\ =\ \int_{[0,1]^n} w(s)\,ds$ for every continuous $1$-periodic function $w\colon \R^n\to\mathbb C$;
\item[\rm{(iii)}]  for every continuous $1$-periodic
$w\colon \R^n\to\mathbb C$,  the function  ${w\circ f}\colon\R^{\geq}\to\mathbb C$ has mean value  $M(w\circ f)=M(w)$.
\end{enumerate}
\end{cor}

\begin{cor}\label{cor:limsup vs sup, ud mod 1}
Let $w\colon\R^n\to\R^{\geq}$ be $1$-periodic and continuous, and suppose $f$ is  u.d.~$\operatorname{mod}\,1$.   Then
$$\limsup_{t\to \infty} w\big(f(t)\big) = 0\ \Longleftrightarrow\ 
\limsup_{\abs{s}\to\infty} w(s)=0\ \Longleftrightarrow\ w=0.$$
\end{cor}

\begin{proof}
Corollary~\ref{cor:udmod1 test, n} gives $M(w\circ f)=M(w)$, and $\dabs{w}=\limsup_{\abs{s}\to\infty} w(s)$ by
Lemma~\ref{lem:limsup vs sup}.
One verifies easily that $M(w\circ f)\leq \limsup_{t\to\infty} w\big(f(t)\big)$.
The equivalences now follow from these facts and Proposition~\ref{prop:Bohr}.
\end{proof}

\begin{cor}\label{cor:limsup = sup, 1-periodic cont}
Let $w\colon\R^n\to\R$ be $1$-periodic and continuous, and suppose~$f$ is~u.d.~$\operatorname{mod}\,1$. Then
$\limsup\limits_{t\to \infty} w\big(f(t)\big) = \sup_{s} w(s)$ and $\liminf\limits_{t\to \infty} w\big(f(t)\big) = \inf_{s} w(s)$.
\end{cor}
\begin{proof}
Let $a\in \R$, $a < \sup_s w(s)$. Then $w=u+v$ where $u, v\colon \R^n\to \R$
are given by~$u(s):= \min\! \big(a, w(s)\big)$ for all $s$. Then $u$,  and thus  $v$,  is  $1$-periodic and continuous. Now $v\ge 0$, but $v\ne 0$, so $\limsup_{t\to \infty} v\big(f(t)\big)>0$ by 
Corollary~\ref{cor:limsup vs sup, ud mod 1}.  This gives~$\varepsilon>0$ with
$v(f(t))>\varepsilon$ for arbitrarily large $t$. For such $t$
we have~$u\big(f(t)\big)=a$: $u\big(f(t)\big)<a$ would give $u\big(f(t)\big)=w\big(f(t)\big)$, so
$v\big(f(t)\big)=0$. Hence 
$w\big(f(t)\big)= u\big(f(t)\big)+ v\big(f(t)\big) > a+\varepsilon$ for such $t$. The other equality follows likewise.
\end{proof}

\noindent
Weyl's Theorem~\ref{thm:Weyl} also generalizes:

\begin{theorem}\label{nWeyl} 
The map $f$ is u.d.~$\operatorname{mod}\,1$ if and only if for all $k\in (\Z^n)^{\neq}$,
$$\lim_{T\to \infty} \frac{1}{T}\int_0^T \ex^{2\pi\imag (k\cdot f(t))}\,dt\ =\ 0.$$
\end{theorem}
\begin{proof} 
Like that of 
Theorem~\ref{thm:Weyl}, using~\ref{cor:udmod1 test, n} instead of~\ref{cor:udmod1 test}.
% and the following special case of the Stone-Weierstrass theorem: for any continuous $1$-periodic function $w\colon\R^n\to\mathbb C$ and any $\varepsilon\in\R^>$ there are  $c_1,\dots,c_m\in\mathbb C$ and $k_1,\dots,k_m\in\Z^n$, such that with 
% $$v(s)\ :=\ c_1\ex^{2\pi\imag (k_1\cdot s)}+\cdots +c_m\ex^{2\pi\imag (k_m\cdot s)}\qquad (s\in \R^n)$$ we have $\sup_{s\in [0,1]^n}\abs{v(s)-w(s)}\leq\varepsilon$.
\end{proof}

\noindent
Theorems~\ref{thm:Weyl} and~\ref{nWeyl} yield:

\begin{cor}\label{cornWeyl}
The map $f$ is u.d.~$\operatorname{mod}\,1$ if and only if for all $k\in (\Z^n)^{\neq}$, the function $t\mapsto {k\cdot f(t)}\colon\R^{\ge}\to\R$ is u.d.~$\operatorname{mod}\,1$.
\end{cor}

\subsection*{Strengthening uniform distribution}
{\it In this subsection $n\ge 1$, the func\-tions $f_1,\dots, f_n\colon \R^{\ge} \to \R$ are   measurable, $f=(f_1,\dots,f_n)\colon\R^{\geq}\to\R^n$ is the resulting map, and
for $\alpha\in\R^n$  
we set $\alpha f:= (\alpha_1f_1,\dots, \alpha_nf_n)\colon\R^{\geq}\to\R^n$}. 
 
\begin{lemma}\label{lem:ud map fn}
The following conditions on $f$ are equivalent:
\begin{enumerate}
\item[\rm{(i)}] $\alpha f$ is  u.d.~$\operatorname{mod}\,1$ for all $\alpha\in (\R^{\times})^n$;
\item[\rm{(ii)}] $\displaystyle\lim_{T\to\infty} \frac{1}{T}\int_0^T \ex^{2\pi\imag (\beta\cdot f(t))}\,dt\ =\ 0$ for all $\beta\in (\R^n)^{\neq}$;
\item[\rm{(iii)}] for every almost periodic
$w\colon \R^n\to\mathbb C$, the function $w\circ f\colon\R^{\geq}\to\mathbb C$ has mean value $M(w\circ f)=M(w)$.
\end{enumerate}
\end{lemma}
\begin{proof}
Assume (i); let $\beta\in (\R^n)^{\neq}$. For $i=1,\dots,n$ set  $\alpha_i:=1$, $k_i:=0$ if $\beta_i=0$ and 
$\alpha_i:=\beta_i$, $k_i:=1$ if $\beta_i\neq 0$. Then $k=(k_1,\dots,k_n)\in (\Z^n)^{\neq}$, 
$\alpha=(\alpha_1,\dots,\alpha_n)$ is in~$(\R^\times)^n$,  and $\beta\cdot f(t) = k\cdot (\alpha f)(t)$ for all $t\in\R$.
Now (ii) follows from Theorem~\ref{nWeyl} applied to~$\alpha f$ in place of $f$.
%For (ii)~$\Rightarrow$~(iii), assume (ii), and let $w\colon\R^n\to \mathbb{C}$ be almost periodic. Let $\varepsilon\in \R^{>}$ and take a trigonometric polynomial $v\colon \R^n\to \mathbb{C}$ such that $\|v-w\|\le \varepsilon$.  We have
%\begin{align*} \frac{1}{T}\int_0^T (w\circ f)(t)dt -M(w)\
%=\  &\ \frac{1}{T}\int_0^T \big((w\circ f)(t) - (v\circ f)(t)\big)dt\ +\\
%&\ \frac{1}{T}\int_0^T(v\circ f)(t)dt - M(v)\  +\  M(v-w).
%\end{align*}
%Using Lemma~\ref{lem:ap mean value}, argue as at the end of the proof of Theorem~\ref{thm:Weyl} to obtain
%$$\left|\frac{1}{T}\int_0^T (w\circ f)(t)dt -M(w)\right|\ \le\ 3\varepsilon \quad\text{ for all sufficiently large }T.$$
The implication (ii)~$\Rightarrow$~(iii) follows as in the proof of
Theorem~\ref{thm:Weyl}, using the definition of almost periodicity instead of Proposition~\ref{prop:StoWei}. 
Finally, assume (iii), and let $\alpha\in (\R^{\times})^n$; to show that $\alpha f$ is
u.d.~$\operatorname{mod}\,1$ we verify that condition (iii) in Corollary~\ref{cor:udmod1 test, n}
holds for $\alpha f$ in place of~$f$. Thus let $w\colon\R^n\to\mathbb C$ be continuous and $1$-periodic. By (iii) applied to the almost periodic function $$s\mapsto w_{\times\alpha}(s)=w(\alpha s)\ \colon\ \R^n\to\mathbb C$$ in place of $w$,   the function
 $w_{\times\alpha}\circ f=w\circ (\alpha f)\colon\R^{\geq}\to\mathbb C$ has a mean value and~$M(w_{\times\alpha}\circ  f)=M(w_{\times\alpha})$;
 now use that $M(w_{\times\alpha})=M(w)$ by Lemma~\ref{lem:M(alpha w)}.
\end{proof}

\noindent
We say that $f$ is {\bf uniformly distributed} (abbreviated: u.d.) if it satisfies one of the equivalent conditions in Lemma~\ref{lem:ud map fn}.
 This lemma also yields:

\begin{cor}\label{cor:ud map fn}
The map $f$ is u.d.~if and only if for all $\beta\in (\R^n)^{\neq}$, the function~$t\mapsto {\beta\cdot f(t)}\colon\R^{\ge}\to\R$ is  u.d.~$\operatorname{mod}\,1$.
\end{cor}

\noindent
The proof of the next result is like that of Corollary~\ref{cor:limsup vs sup, ud mod 1}, using Lemma~\ref{lem:ud map fn} instead of Corollary~\ref{cor:udmod1 test, n}:

\begin{cor}\label{cor:limsup vs sup}
Suppose $w\colon\R^n\to\R^{\geq}$ is almost periodic and $f$ is  u.d.   Then
$$\limsup_{t\to \infty} w\big(f(t)\big) = 0\ \Longleftrightarrow\ 
\limsup_{\abs{s}\to\infty} w(s)=0\ \Longleftrightarrow\ w=0.$$
\end{cor}
%\begin{proof}
%Lemma~\ref{lem:ud map fn} gives $M(w\circ f)=M(w)$, and $\dabs{w}=\limsup_{\abs{s}\to\infty} w(s)$ by Lemma~\ref{lem:limsup vs sup}. From Lemma~\ref{lem:mv and sup} we obtain $M(w\circ f)\leq \limsup_{t\to\infty} w\big(f(t)\big)$. The equivalences now follow from these facts and Proposition~\ref{prop:Bohr}.
%\end{proof}

\subsection*{Application to Hardy fields} 
{\it In this subsection $f_1,\dots, f_n\colon \R^{\ge}\to \R$ with $n\ge 1$ are continuous, their germs, denoted also by $f_1,\dots, f_n$, lie in a common Hardy field, and $f:=(f_1,\dots,f_n)\colon \R^{\ge}\to \R^n$.}\/  Theorem~\ref{thm:udmod1 Boshernitzan} with Corollary~\ref{cornWeyl} gives:

\begin{cor}[Boshernitzan]\label{cor:Bosh} We have the following equivalence:
$$\text{$f$ is u.d.~$\operatorname{mod}\,1$}\ \Longleftrightarrow\ \text{$k_1f_1+\cdots+k_nf_n \succ \log x$ for all $(k_1,\dots,k_n)\in (\Z^n)^{\neq}$.}$$
\end{cor}

\noindent
Combining  Theorem~\ref{thm:udmod1 Boshernitzan} with Corollary~\ref{cor:ud map fn} yields likewise:

\begin{cor}\label{nbos}
We have the following equivalence:
$$\text{$f$ is u.d.}\ \Longleftrightarrow\ \text{$\alpha_1f_1+\cdots+\alpha_nf_n \succ \log x$ for all $(\alpha_1,\dots,\alpha_n)\in (\R^n)^{\neq}$.}$$
In particular, if $\log x\prec f_1\prec\cdots\prec f_n$, then $f$ is u.d. 
\end{cor}

\noindent
Here is an immediate application of Corollary~\ref{cor:Bosh}:
%From Corollary~\ref{cor:Bosh} we obtain: 

\begin{cor}[Weyl]\label{oweyl}
Let $\lambda_1,\dots,\lambda_n\in\R$. Define $g\colon \R^{\ge} \to \R^n$ by  $g(t)=(\lambda_1 t,\dots,\lambda_n t)$.
Then~$g$ is u.d.~$\operatorname{mod}\,1$
iff $\lambda_1,\dots,\lambda_n$ are $\Q$-linearly independent.
\end{cor}

\noindent
We now get to the result that we actually need in Section~\ref{sec:ueeh}:

\begin{prop}\label{prop:limsup vs sup}
Suppose $w\colon\R^n\to\R^{\geq}$ is almost periodic,  $1\prec f_1\prec\cdots\prec f_n$, and
 $\limsup\limits_{t\to+\infty} w\big(f(t)\big) = 0$. Then $w=0$.
\end{prop}
\begin{proof}
We first arrange $f_1>\R$, replacing $f_1,\dots,f_n$ and $w$ by $-f_1,\dots,-f_n$ and the function $s\mapsto w(-s)\colon\R^n\to\R^{\geq}$, if $f_1<\R$.
Pick $a\ge 0$ such that the restriction of $f_1$ to $\R^{\geq a}$ is strictly increasing, set $b:=f_1(a)$, and let 
$f_1^{\inv}\colon\R^{\geq b}\to\R$ be the compositional inverse of this restriction.
Set $g_j(t):=(f_j \circ f_1^{\inv})(t)$ for $t\ge b$ and~$j=1,\dots,n$
and  consider the map $$g\ =\ (g_1,\dots, g_n)\ =\ f\circ f_1^{\inv}\  :\  \R^{\ge b} \to \R^n.$$
The germs of $g_1,\dots,g_n$, denoted by the same symbols, 
lie in a common Hardy field and satisfy $x=g_1\prec  g_2 \prec  \cdots \prec  g_n$.
Now $f_1^{\inv}$ is strictly increasing and moreover~$f_1^{\inv}(t)\to+\infty$ as $t\to+\infty$, so
$$\limsup_{t\to\infty} w\big(f(t)\big)\ =\ \limsup_{t\to\infty} w\big(f\big(f_1^{\inv}(t)\big)\big)\ =\ 
\limsup_{t\to\infty} w\big(g(t)\big)\ =\ 0.$$
Thus replacing $f_1,\dots,f_n$ by continuous functions $\R^{\geq}\to\R$ with the same
germs as~$g_1,\dots,g_n$, we arrange $x=f_1\prec f_2\prec\cdots\prec f_n$.
Then $f$ is u.d.~by  Corollary~\ref{nbos}. Now use Corollary~\ref{cor:limsup vs sup}.
\end{proof}

\noindent
The next three results are not used later but included for their independent interest. 

\begin{cor}\label{cor:limsup = sup} Assume $w\colon\R^n\to\R$ is almost periodic and  $1\prec f_1\prec\cdots\prec f_n$. Then $\limsup\limits_{t\to \infty} w\big(f(t)\big) = \sup_{s} w(s)$ and $\liminf\limits_{t\to \infty} w\big(f(t)\big) = \inf_{s} w(s)$. 
\end{cor}

\begin{proof}
Let $a\in \R$, $a < \sup_s w(s)$. Then $w=u+v$ where $u, v\colon \R^n\to \R$
are given by~$u(s):= \min\! \big(a, w(s)\big)$ for all $s$. Then $u$,  and thus  $v$,  is   almost periodic by Corollary~\ref{apsub}. 
Now argue as in the proof of Corollary~\ref{cor:limsup = sup, 1-periodic cont}, using
Proposition~\ref{prop:limsup vs sup} instead of  Corollary~\ref{cor:limsup vs sup, ud mod 1}.
\end{proof}

%\begin{proof} 
%Let $a\in \R$, $a < \sup_s w(s)$. Then $w=u+v$ where $u, v\colon \R^n\to \R$ are given by~$u(s):= \min \big(a, w(s)\big)$ for all $s$. Then $u$,  and thus  $v$,  is  almost periodic, by Corollary~\ref{apsub}. Now $v\ge 0$, but $v\ne 0$, so $\limsup_{t\to \infty} v\big(f(t)\big)>0$ by Proposition~\ref{prop:limsup vs sup}. This gives an $\varepsilon>0$ such that $v(f(t))>\varepsilon$ for arbitrarily large $t$. For such $t$ we have~$u\big(f(t)\big)=a$, since $u\big(f(t)\big)<a$ would mean $u\big(f(t)\big)=w\big(f(t)\big)$, so $v\big(f(t)\big)=0$, a contradiction. Hence for such $t$ we have  $w\big(f(t)\big)= u\big(f(t)\big)+ v\big(f(t)\big) > a+\varepsilon$.
%\end{proof} 

\begin{cor}
If $w\colon\R^n\to\C$ is almost periodic and $1\prec f_1\prec\cdots\prec f_n$, then
$$\lim_{t\to\infty} w\big(f(t)\big) \text{ exists in }\C \quad\Longleftrightarrow\quad \text{$w$ is constant.}$$
\end{cor}
\begin{proof} Apply the previous corollary to the real and imaginary part of $w$.
%First suppose $w(\R^n)\subseteq\R$.
%Assume $\lim\limits_{t\to\infty} w\big(f(t)\big)=\ell\in\R$ exists. Then 
%$\limsup\limits_{t\to\infty} w\big(f(t)\big)=\liminf\limits_{t\to\infty} w\big(f(t)\big)=\ell$ and hence
%$\sup_s w(s)=\inf_s w(s)=\ell$ by the previous corollary, so $w=\ell$ is constant. The general case follows from
%this special case applied to $\Re w$, $\Im w$ in place of $w$.
\end{proof}

\noindent
Finally, we use these results to reprove \cite[Theorem~8]{vdHlimsup}. 
Given $\alpha=(\alpha_1,\dots,\alpha_n)\in\mathbb C^n$ we put $\ex^{\alpha}:=(\ex^{\alpha_1},\dots,\ex^{\alpha_n})\in\mathbb C^n$.
Let $m\ge 1$ and set
$$S\ :=\ \big\{(z_1,\dots,z_m)\in\mathbb C^m:\ \abs{z_1}=\cdots=\abs{z_m}=1\big\}.$$

\begin{cor}
Suppose $1\prec f_1\prec\cdots\prec f_n$. Let $\varphi\colon S\to\R$ be continuous and
let $k_1,\dots,k_n\in\N^{\ge 1}$ with $k_1+\cdots+k_n=m$ and
$\lambda_j=(\lambda_{j1},\dots,\lambda_{jk_j})\in\R^{k_j}$ for~$j=1,\dots,n$ be such that $\lambda_{j1},\dots,\lambda_{jk_j}$ are $\Q$-linearly independent.
Then
$$\limsup_{t\to\infty} \varphi\big(\ex^{\imag f_1(t) \lambda_1},\dots,\ex^{\imag f_n(t)\lambda_n}\big)\ =\ 
\max \varphi(S).$$
\end{cor}
\begin{proof}
By Corollary~\ref{apsub}, the function
$$s=(s_1,\dots,s_n) \mapsto w(s):=\varphi\big(\ex^{\imag s_1 \lambda_1},\dots,\ex^{\imag s_n\lambda_n}\big)\ \colon\ \R^n\to\R$$
 is almost periodic. We have
$$w\big(f(t)\big)\ =\ \varphi\big(\ex^{\imag f_1(t) \lambda_1},\dots,\ex^{\imag f_n(t)\lambda_n}\big)\quad\text{for $t\geq 0$,}$$
so by Corollary~\ref{cor:limsup = sup},
$$\limsup_{t\to\infty} \varphi\big(\ex^{\imag f_1(t) \lambda_1},\dots,\ex^{\imag f_n(t)\lambda_n}\big)\ =\ \limsup_{t\to\infty} w\big(f(t)\big)\ =\ \sup_s w(s).$$
For $j=1,\dots,n$ it follows from Corollary~\ref{oweyl}  that the image of the map 
$$t\mapsto \ex^{\imag t\lambda_j}\ :\ \R^{\ge} \to \big\{(z_1,\dots,z_{k_j})\in\mathbb C^{k_j}: \abs{z_1}=\cdots=\abs{z_{k_j}}=1\big\}$$
is dense in its codomain, so the image of the map
$$(s_1,\dots,s_n)\mapsto \big(\ex^{\imag s_1 \lambda_1},\dots,\ex^{\imag s_n\lambda_n}\big)\ \colon\
\R^n\to S$$ 
is dense in $S$. Hence $\sup_s w(s)=\max\varphi(S)$.
\end{proof}

\subsection*{Examples involving real-valued trigonometric polynomials\astr}
The material in this subsection is only used later to justify a remark after Corollary~\ref{lem:gaussian ext dom, general}. 

\begin{exampleNumbered}\label{ex:oweyl}
Let   $a,b\in \R^\times$, and consider the $1$-periodic trigonometric polynomial~$w\colon\R^2\to\R$ given
by $$w(s)\ =\ a\cos(2\pi s_1)+ b\cos(2\pi s_2)\qquad\text{for   $s=(s_1,s_2)\in\R^2$.}$$
Let $\lambda,\mu\in\R$ be $\Q$-linearly independent. Then by Corollaries~\ref{cor:limsup = sup, 1-periodic cont} and~\ref{oweyl}:
$$\limsup_{t\to -\infty}w(\lambda t, \mu t)\ =\ \limsup_{t\to+\infty} w(\lambda t,\mu t)\ =\  \abs{a}+\abs{b}.$$
Note also that $\abs{a}+\abs{b}>w(\lambda t,\mu t)$ for all $t\in\R^{\ne}$.
%if   $\abs{a}+\abs{b}=w(\lambda t,\mu t)$, then~$a\cos(2\pi\lambda t)=\abs{a}$ and $b\cos(2\pi\mu t)=\abs{b}$, so $\lambda t,\mu t\in \Z$,   and thus~$\lambda\in\mu\Q$, a contradiction. 
\end{exampleNumbered}

\noindent
Now let $\alpha\in\R\setminus\Q$ and consider   
$$  v\colon\R\to\R,\qquad v(t):=2-\cos(t)-\cos(\alpha t).$$
Then $v(t)>0$ for all $t\in\R^{\ne}$. Moreover, $\liminf\limits_{t\to+\infty} v(t)=0$, that is,
for each $\varepsilon>0$ there are arbitrarily large $t\in\R$ with $v(t)<\varepsilon$.
With a suitable choice of $\alpha$ we can replace here $\varepsilon$ by any prescribed function $\varepsilon\colon\R\to\R^>$ with $\varepsilon(t)\to 0$ as $t\to+\infty$:

\begin{theorem}[{Basu-Bose-Vijayaraghavan~\cite{BBV}}]\label{thm:BBV}
Let $\phi\colon\R\to\R^>$ %\textup{(}$a\in\R$\textup{)}  be eventually increasing 
be such that $\phi(t)\to+\infty$ as~$t\to+\infty$. Then there exists $\alpha\in\R\setminus\Q$
such that 
$$2-\cos(t)-\cos(\alpha t) < 1/\phi(t)\quad\text{ for arbitrarily large $t\in\R$.}$$
\end{theorem}
\begin{proof}
We first arrange $\phi\geq 1$. We then choose a se\-quence~$(d_n)_{n\geq 1}$ 
of  positive integers  such that with $q_n:=d_1 d_2\cdots d_n$ (so~$q_0=1$):
$$d_n \geq (2\pi +1 )\, \phi(2\pi q_{n-1})\qquad\text{for $n\geq 1$,}$$
and set
$$\alpha\ :=\ \sum_{n=1}^\infty \frac{1}{q_n}.$$
We have $q_{m+1}=q_m d_{m+1} \geq (2\pi+1) q_m \,\phi(2\pi q_m)$,
so if 
$q_{m+n} \geq (2\pi+1)^n q_m \,\phi(2\pi q_m)$, then
$$q_{m+n+1} \geq (2\pi+1) q_{m+n}\, \phi(2\pi q_{m+n}) \geq (2\pi+1) q_{m+n} \geq (2\pi+1)^{n+1} q_m\, \phi(2\pi q_m).$$
Thus by induction on $n$ we obtain
$$q_{m+n} \geq (2\pi+1)^n q_m\, \phi(2\pi q_m)\qquad\text{for $n\geq 1$.}$$
This yields  
\begin{equation}\label{eq:BBV}
\sum_{n=1}^\infty \frac{1}{q_{m+n}} \leq \frac{1}{2\pi q_m\,\phi(2\pi q_m)}\qquad\text{for all $m\geq 1$.}
\end{equation}
Take $p_m\in\N$ ($m\geq 1$) such that
$$\sum_{n=1}^m \frac{1}{q_n} = \frac{p_m}{q_m}.$$
Then    
$$0<\alpha-\frac{p_m}{q_m} = \sum_{n=1}^\infty \frac{1}{q_{m+n}} \leq \frac{1}{2\pi q_m\,\phi(2\pi q_m)}\qquad\text{for all $m\geq 1$.} $$
Suppose $\alpha=p/q$ where $p,q\in\N^{\ge 1}$; then for all $m\geq 1$ we have
$$\frac{pq_m-qp_m}{qq_m} = \alpha-\frac{p_m}{q_m}\leq \frac{1}{2\pi q_m\,\phi(2\pi q_m)}$$
and so
$$1\leq pq_m-qp_m \leq  \frac{q}{2\pi\,\phi(2\pi q_m)},$$
contradicting $\phi(2\pi q_m)\to+\infty$ as $m\to+\infty$. Hence $\alpha\notin\Q$. Next note that
$q_m$ and~$\alpha q_m-q_m\sum_{n\geq 1}\frac{1}{q_{m+n}}$ are integers and so
$$2 - \cos(2\pi q_m) - \cos(2\pi \alpha q_m) = 1 - \cos\left( 2\pi q_m \sum_{n=1}^\infty \frac{1}{q_{m+n}}\right) < \frac{1}{\phi(2\pi q_m)}
$$
using \eqref{eq:BBV}. This yields the theorem.
\end{proof}

\section{Universal Exponential Extensions of Hardy Fields}\label{sec:ueeh}

\noindent
{\it In this section $H\supseteq \R$ is a Liouville closed Hardy field.}\/ Then~$\Calinf[\imag]$ is a differential ring extension
of the $\d$-valued field $K:=H[\imag]$ with the same ring of constants as~$K$, 
namely~$\mathbb C$. Note that for any $f\in \Calinf[\imag]$ we have a $g\in \Calinf[\imag]$
with $g'=f$, and then $u=\ex^g\in \Calinf[\imag]^\times$ satisfies
$u^\dagger =f$.\index{Hardy field!universal exponential extension}\index{universal exponential extension!of a Hardy field} 

\begin{lemma}\label{purelyimag} Suppose $f\in \Calinf[\imag]$ is purely imaginary, that is,
$f\in \imag\Calinf$. Then there is a $u\in \Calinf[\imag]^\times$
such that $u^\dagger=f$ and $|u|=1$. 
\end{lemma}
\begin{proof} Taking $g\in \imag\Calinf$ with $g'=f$, the resulting
$u=\ex^g$ works. 
\end{proof}

\noindent 
We define the subgroup $\ex^{H\imag}$ of $\Calinf[\imag]^\times$ by 
$$\ex^{H\imag}\ :=\ \big\{\!\ex^{h\imag}:\ h\in H\big\}\ =\ \big\{u\in \Calinf[\imag]^\times:\ |u|=1,\ u^\dagger\in H\imag\big\}.$$ 
Then 
$(\ex^{H\imag})^\dagger=H\imag$ by Lemma~\ref{purelyimag}, so $(H^\times\cdot \ex^{H\imag})^\dagger=K$
and thus~$K[\ex^{H\imag}]$ is an exponential extension of $K$ (in the sense of Section~\ref{sec:univ exp ext}) with the same ring of constants $\mathbb C$ as~$K$.

As in the beginning of Section~\ref{sec:ultimate} we fix a complement $\Lambda$ of $K^\dagger$ with $\Lambda\subseteq H\imag$, set~$\Univ:= K\big[\!\ex(\Lambda)\big]$ as usual, and let $\lambda$ range over $\Lambda$. The differential $K$-algebras~$\Univ$ and
$K[\ex^{H\imag}]$ are isomorphic by Corollary~\ref{corcharexp}, but we need something better:

\begin{lemma}\label{exqe} There is an isomorphism $\Univ\to K[\ex^{H\imag}]$ 
of differential $K$-algebras that maps $\ex(\Lambda)$ into $\ex^{H\imag}$. 
\end{lemma} 
\begin{proof} We have a short exact sequence of commutative groups
$$1 \to S \xrightarrow{\ \subseteq\ } \ex^{H\imag} \xrightarrow{\ \ell\ } H\imag \to 0,$$
where $S=\big\{z\in \mathbb{C}^\times:\, |z|=1\big\}$ and
 $\ell(u):=u^\dagger$ for~$u\in \ex^{H\imag}$. Since the subgroup~$S$ of~$\mathbb{C}^\times$ is divisible, this sequence splits: we have a group embedding 
$e\colon H\imag\to \ex^{H\imag}$ such that $e(b)^\dagger=b$ for all $b\in H\imag$. Then
the group embedding $$\ex(\lambda)\mapsto e(\lambda)\ :\ \ex(\Lambda)\to \ex^{H\imag}$$  extends uniquely to a $K$-algebra morphism $\Univ\to K[\ex^{H\imag}]$. Since
$\ex(\lambda)^\dagger=\lambda=e(\lambda)^\dagger$ for all $\lambda$, this is a
differential $K$-algebra morphism, and even an isomorphism by
Lemma~\ref{lem:embed into U} applied to $R=K[\ex^{H\imag}]$.   
\end{proof}

\noindent
Complex conjugation $f+g\imag\mapsto \overline{f+g\imag}=f-g\imag$ ($f,g\in \Calinf$) is an automorphism of the differential ring~$\Calinf[\imag]$ over $H$ and maps
$K[\ex^{H\imag}]$ onto itself, sending each~${u\in \ex^{H\imag}}$ to~$u^{-1}$. Thus any isomorphism
$\iota\colon \Univ\to K[\ex^{H\imag}]$ of differential $K$-algebras with~$\iota\big(\!\ex(\Lambda)\big)\subseteq \ex^{H\imag}$---such $\iota$ exists by Lemma~\ref{exqe}---also satisfies 
$$\iota(\overline{f})\ =\ \overline{\iota(f)}\qquad(f\in \Univ). $$ 
(See Section~\ref{sec:univ exp ext} for the definition of $\overline{f}$ for $f\in\Univ$.
Given such an isomorphism~$\iota$, any differential $K$-algebra isomorphism $\Univ\to K[\ex^{H\imag}]$   mapping $\ex(\Lambda)$ into $\ex^{H\imag}$ equals $\iota\circ\sigma_\chi$ for a unique character $\chi\colon\Lambda\to\C^\times$ with $\abs{\chi(\lambda)}=1$ for all $\lambda$, by Lem\-ma~\ref{autolem, commuting with f*}.) 
Fix such an isomorphism $\iota$ and identify $\Univ$ with its image $K[\ex^{H\imag}]$ via~$\iota$.
We have the asymptotic relations
$\preceq_{\g}$ and~$\prec_{\g}$
on $\Univ$ coming from the gaussian extension $v_{\g}$ of the valuation on $K$. But we also have the asymptotic relations induced on $\Univ=K[\ex^{H\imag}]$
by the relations $\preceq$ and $\prec$ defined on $\mathcal C[\imag]$ in Section~\ref{sec:germs}. It is clear that for $f\in\Univ$:
\begin{align*}
f\preceq_{\g} 1 	&\quad \Longrightarrow \quad f\preceq 1 \quad \Longleftrightarrow \quad \text{for some $n$ we have $\abs{f(t)}\leq n$ eventually,} \\
f\prec_{\g} 1	&\quad \Longrightarrow \quad f\prec 1 \quad \Longleftrightarrow \quad \lim_{t\to+\infty} f(t)=0.
\end{align*}
As a tool for later use we derive a converse of the implication
$f\prec_{\g}1\Rightarrow f\prec 1$: Lemma~\ref{lem:gaussian ext dom} below, where we assume in addition
that~${\I(K)\subseteq K^\dagger}$ and $\Lambda$ is an $\R$-linear subspace
of $K$. 
This requires the material from Section~\ref{sec:udmod1} and some considerations about exponential sums treated in the next subsection.

\subsection*{Exponential sums over Hardy fields}
In this subsection $n\geq 1$.
In the next lemma, $f=(f_1,\dots,f_m)\in H^m$ where $m\ge 1$ and
 $1\prec f_1\prec\cdots\prec f_m$. (In that lemma it doesn't matter which
 functions we use to represent the germs $f_1,\dots, f_m$.)
%Choose~$a\in \R$ so that $f_1,\dots,f_m$ have representatives on $[a,+\infty)$, denoted by the same symbols.
For $r=(r_1,\dots,r_m)\in\R^m$ we set $r\cdot f:= r_1f_1+\cdots + r_mf_m\in H$.
 
\begin{lemma}\label{lem:limsup 1} 
Let $r^1,\dots,r^n\in\R^m$ be distinct and $c_1,\dots,c_n\in\mathbb C^\times$. Then
$$\limsup_{t\to \infty} \left|c_1 \ex^{(r^1\cdot f)(t)\imag}+ \cdots + c_n\ex^{(r^n\cdot f)(t)\imag}\right|\ >\ 0.$$ 
\end{lemma}
\begin{proof} Consider the trigonometric polynomial $w\colon\R^m\to \R^{\ge}$ given by
$$ w(s)\ :=\ \big|c_1 \ex^{(r^1\cdot s)\imag}+\cdots+c_n \ex^{(r^n\cdot s)\imag}\big|^2.$$  
By Corollary~\ref{cor:id thm} we have $w(s)>0$ for some $s\in\R^m$. 
Taking continuous representatives $\R^{\ge} \to \R$ of $f_1,\dots,f_m$, to be denoted also by $f_1,\dots,f_m$, the lemma now follows from Proposition~\ref{prop:limsup vs sup}. 
\end{proof}

\noindent
Next, let  $h_1,\dots,h_n\in H$ be distinct such that
$(\R h_1+\cdots+\R h_n)\cap \I(H)=\{0\}$. Since $H$ is Liouville closed we have
$\phi_1,\dots,\phi_n\in H$  such that $\phi_1'=h_1,\dots, \phi_n'=h_n$.

\begin{lemma}\label{lem:limsup 2}
Let  $c_1,\dots,c_n\in\mathbb C^\times$. Then for $\phi_1,\dots, \phi_n$ as above,
$$\limsup_{t\to \infty} \left| c_1 \ex^{\phi_1(t)\imag}+ \cdots + c_n \ex^{\phi_n(t)\imag}\right|\ >\ 0.$$
\end{lemma}
\begin{proof} The case $n=1$ is trivial, so let $n\ge 2$. Then $\phi_1,\dots, \phi_n$ are not all in~$\R$. 
Set~$V:= \R+\R\phi_1+\cdots+\R\phi_n \subseteq H$, so
$\der V=\R h_1+\cdots+\R h_n$.
We claim that~${V\cap \smallo_H=\{0\}}$.
To see this, let $\phi\in V\cap \smallo_H$; then $\phi'\in \der(V)\cap \I(H)=\{0\}$ and hence $\phi\in\R\cap \smallo_H=\{0\}$,
proving the claim.
Now $H$ is a Hahn space over $\R$ by~[ADH, p.~109], so by [ADH, 2.3.13] we have $f_1,\dots,f_m\in V$  ($1\le m\leq n$) such that~$V = \R + \R f_1 +\cdots + \R f_m$ and $1\prec f_1\prec \cdots \prec f_m$. For $j=1,\dots,n$, $k=1,\dots,m$, 
take $t_j,r_{jk}\in\R$ such that $\phi_j=t_j+\sum_{k=1}^m r_{jk} f_k$
and set $r^j:=(r_{j1},\dots,r_{jm})\in\R^m$.
Since $\phi_{j_1}-\phi_{j_2}\notin\R$ for   $j_1\neq j_2$, we have
$r^{j_1}\neq r^{j_2}$ for   $j_1\neq j_2$.
It remains to apply Lemma~\ref{lem:limsup 1} to $c_1\ex^{t_1\imag},\dots, c_n\ex^{t_n\imag}$ in place of $c_1,\dots, c_n$. 
\end{proof}

\begin{cor}\label{corlimsupresult}
Let $f_1,\dots,f_n\in K$ and set $f:= f_1 \ex^{\phi_1\imag} + \cdots + f_n\ex^{\phi_n\imag}\in\Calinf[\imag]$, and suppose~$f\prec 1$. Then $f_{1},\dots,f_{n}\prec 1$.
\end{cor}
\begin{proof}
We may assume   $0 \neq f_1\preceq \cdots \preceq f_n$.
Towards a contradiction, suppose that~${f_n\succeq 1}$, and
take $m\leq n$ minimal such that $f_m\asymp f_n$. Then
 with~$g_j:=f_j/f_n\in K^\times$ and $g:=g_1 \ex^{\phi_1\imag} + \cdots + g_n\ex^{\phi_n\imag}$ we have $g\prec 1$ and
$g_1,\dots, g_n\preceq 1$, with $g_j\prec 1$ iff~$j<m$.
Replacing $f_1,\dots,f_n$ by $g_m,\dots,g_n$ and $\phi_1,\dots,\phi_n$ by $\phi_m,\dots,\phi_n$
we arrange $f_1\asymp\cdots\asymp f_n\asymp 1$. So
$$f_1=c_1+\varepsilon_1,\,\dots,\, f_n=c_n+\varepsilon_n\quad\text{ with
$c_1,\dots, c_n\in\mathbb C^\times$ and $\varepsilon_1,\dots, \varepsilon_n\in \smallo$.}$$ 
Then $\varepsilon_1\ex^{\phi_1\imag} +\cdots + \varepsilon_n\ex^{\phi_n\imag}\prec 1$, hence
 $$c_1\ex^{\phi_1\imag} +\cdots + c_n\ex^{\phi_n\imag}\ =\ f-\big(\varepsilon_1\ex^{\phi_1\imag} +\cdots + \ex^{\phi_n\imag}\!\big)\ \prec\ 1.$$
Now Lemma~\ref{lem:limsup 2} yields the desired contradiction.
\end{proof}

\noindent
Here is an application of Corollary~\ref{corlimsupresult}: 
 
\begin{lemma}\label{lem:exp sum asymptotics}
Let $f_1,\dots,f_n,g_1,\dots,g_n\in K$ be such that in $\c[\imag]$ we have
$$f_1\ex^{\phi_1\imag}+\cdots+f_n\ex^{\phi_n\imag}\  \sim\  g_1\ex^{\phi_1\imag}+\cdots+g_n\ex^{\phi_n\imag}.$$
Let $j\in\{1,\dots,n\}$ be such that $0\neq f_j\succeq f_k$ for all $k\in\{1,\dots,n\}$. Then $f_j\sim g_j$, and
$f_k-g_k\prec f_j$ for all $k\ne j$. 
\end{lemma}
\begin{proof}
We arrange $j=1$ and $f_1=1$.  Then
$$\ex^{\phi_1\imag}+f_2\ex^{\phi_2\imag}+\cdots+f_n\ex^{\phi_n\imag}\ \sim\ g_1\ex^{\phi_1\imag}+\cdots+g_n\ex^{\phi_n\imag},\qquad
f_2,\dots,f_n\preceq 1.$$
Hence
$$(1-g_1)\ex^{\phi_1\imag}+(f_2-g_2)\ex^{\phi_2\imag}+\cdots+(f_n-g_n)\ex^{\phi_n\imag} \prec \ex^{\phi_1\imag} + f_2\ex^{\phi_2\imag}+\cdots+f_n\ex^{\phi_n\imag} \preceq 1,$$
so $1-g_1\prec 1$  and $f_k-g_k\prec 1$ for all $k\ne j$, by Corollary~\ref{corlimsupresult}.
\end{proof}

\noindent
This leads to a partial generalization of Corollary~\ref{cor:osc => bded}, included for   use in~\cite{ADHld}:

\begin{cor}\label{cor:exp sum asymptotics}
Let $f\in K^\times$, $g_1,\dots,g_n\in K$, and $j\in\{1,\dots,n\}$  such that  in~$\c[\imag]$,
$$f\ex^{\phi_j\imag} \sim g_1\ex^{\phi_1\imag}+\cdots+g_n\ex^{\phi_n\imag}.$$
%  a\sim f\preceq f_1\preceq \cdots \preceq f_n.$$
Then  $f \sim g_j$, and $g_j  \succ g_k$ for all $k\neq j$.
\end{cor}
\begin{proof} Use
Lemma~\ref{lem:exp sum asymptotics} with $f_j:=f$ and $f_k:=0$ for $k\neq j$.
\end{proof}

\noindent
{\em In the rest of this subsection we assume that $\I(K)\subseteq K^\dagger$}. As noted in Sec\-tion~\ref{sec:ultimate} we can then take $\Lambda=\Lambda_H\imag$ where $\Lambda_H$ is an $\R$-linear complement
of $\I(H)$ in $H$. {\em We assume $\Lambda$ has this form,}\/ giving rise to
the valuation $v_{\g}$ on $\Univ=K[\ex^{H\imag}]$ as explained in the beginning of this section.

\begin{lemma}\label{lem:gaussian ext dom} Let $f\in\Univ$ be such that $f\prec 1$. Then $f\prec_{\g} 1$.
\end{lemma}
\begin{proof} We have $f=f_1\ex(h_1\imag)+\cdots + f_n\ex(h_n\imag)$ with
$f_1,\dots, f_n\in K$ and distinct~$h_1,\dots, h_n\in \Lambda_H$, so $(\R h_1+\cdots+\R h_n)\cap \I(H)=\{0\}$. For $h\in \Lambda_H$ we have~$\ex(h\imag)= \ex^{\phi\imag}$ with $\phi\in H$ and $\phi'=h$.
Hence $f=f_1\ex^{\phi_1\imag} + \cdots + f_n\ex^{\phi_n\imag}$ with~$\phi_1,\dots, \phi_n\in H$ such that
$\phi_1'=h_1,\dots, \phi_n'=h_n$. 
Now Corollary~\ref{corlimsupresult}  yields~${f\prec_{\g} 1}$.
\end{proof}

\begin{cor}
\label{cor:gaussian ext dom}
Let $f\in\Univ$ and $\fm\in H^\times$. Then $f\prec \fm$ iff $f\prec_{\g} \fm$.
\end{cor}

\begin{lemma}\label{lem:gaussian ext dom, preceq}
Let $f\in\Univ$ and $\fm\in H^\times$. Then  $f\preceq \fm$ iff $f\preceq_{\g} \fm$.
\end{lemma}
\begin{proof}
Replace $f$, $\fm$ by $f/\fm$, $1$, respectively, to arrange $\fm=1$.
The backward direction was observed earlier in this section.
For the forward direction suppose~${f\preceq 1}$. Then~$f\prec \fn$ for all $\fn\in H^\times$ with $1\prec \fn$, hence
$f\prec_{\g}\fn$ for all $\fn\in H^\times$   with $1\prec_{\g} \fn$, by two applications of Corollary~\ref{cor:gaussian ext dom}, and thus $f\preceq_{\g} 1$. 
%we need to sho $f\prec_{\g}\fn$.Now $f\preceq 1\prec\fn$ and thus~$f\prec\fn$, hence .
\end{proof}

%\noindent
%The next result is included for the sake of completeness and not used later.
 
\begin{cor}\label{lem:gaussian ext dom, general}
Let $f,g\in\Univ$. Then 
\begin{equation}\label{eq:gaussian ext dom, general}
f\preceq g\ \Longrightarrow\ f\preceq_{\g} g,
\end{equation} and likewise with
$(\preceq,\preceq_{\g})$ replaced by $(\asymp,\asymp_{\g})$, $(\prec,\prec_{\g})$, or $(\sim,\sim_{\g})$. In particular, 
$\ex^{\phi\imag}\asymp_{\g} 1$ for all $\phi\in H$. 
%$$f\preceq g\ \Rightarrow\ f\preceq_{\g} g,\qquad f\prec g\ \Rightarrow\ f\prec_{\g} g, \qquad  f \sim g\ \Rightarrow\ f\sim_{\g} g.$$
\end{cor} 
\begin{proof} The case $g=0$ is trivial, so let $g\neq 0$. Then
 $g\asymp_{\g}\fn$ with $\fn\in H^\times$, so $g\preceq_{\g}\fn$ and $\fn\preceq_{\g} g$,
and thus $g\preceq\fn$ by Lemma~\ref{lem:gaussian ext dom, preceq}.
If $f\preceq g$, then $f\preceq\fn$, hence~$f\preceq_{\g}\fn$   by Lemma~\ref{lem:gaussian ext dom, preceq},
so~$f\preceq_{\g} g$.
Likewise, if $f\prec g$, then $f\prec\fn$, so~$f\prec_{\g}\fn$ by Corollary~\ref{cor:gaussian ext dom},
hence~$f\prec_{\g} g$. The rest is now clear.
\end{proof}

\begin{remark}
The converse of \eqref{eq:gaussian ext dom, general} doesn't hold in general, even when we restrict to~$f=1$ and $g\in{\Univ}\cap \c^\times$: let~$\lambda,\mu\in\R$ be $\Q$-linearly independent and set
% to exclude trivialities. To see this let~$a,b\in\R^\times$ and let~$\lambda,\mu\in\R$ be $\Q$-linearly independent, and  
%$$g:=\abs{a}+\abs{b}-a\cos(\lambda x)-b\cos(\mu x)\in\Univ.$$
$$g\ :=\  2- \cos(\lambda x)-\cos(\mu x)\in\Univ;$$
then $1\asymp_{\operatorname{g}} g$, and by Example~\ref{ex:oweyl} we have $g\in\c^\times$ and $1\not\preceq g$.
%This shows that the converse of the displayed implication \eqref{eq:gaussian ext dom, general}   for $({\preceq},{\preceq}_{\operatorname{g}})$ and
%$({\asymp},{\asymp}_{\operatorname{g}})$ does not hold.
Next, take $\phi\in H$ with $\phi>\R$,   choose $\alpha\in\R\setminus\Q$ as in Theorem~\ref{thm:BBV} applied to a representative of the germ $\phi$, and set
$$h\ :=\ \phi\cdot \big(2-\cos(x) - \cos(\alpha x)\big) \in \Univ.$$
Then $h\in\c^\times$ and $h\asymp_{\operatorname{g}}\phi$, so $1  \prec_{\operatorname{g}} h$.
By choice of $\alpha$ we also have~$1\not\prec h$. Hence the converse of \eqref{eq:gaussian ext dom, general} for $({\prec},{\prec}_{\operatorname{g}})$ in place of~$({\preceq},{\preceq}_{\operatorname{g}})$ fails for $f:=1$, $g:=h$.
% for $({\prec},{\prec}_{\operatorname{g}})$, $h$ in place of~$({\preceq},{\preceq}_{\operatorname{g}})$,~$g$, respectively.
\end{remark}

\subsection*{An application to slots in $H$} 
{\it  In   this subsection we assume $\I(K)\subseteq K^\dagger$. 
We take $\Lambda=\Lambda_H\imag$ where~$\Lambda_H$ is an $\R$-linear complement of $\I(H)$ in $H$, and accordingly identify $\Univ$ with $K[\ex^{H\imag}]$ as explained in the beginning of this section.}\/ %and assume $K$ is $\upo$-free if $r\geq 2$. }\/  
Until further notice we let~$(P,1,\hat h)$ be  a slot in~$H$ of 
order $r\geq 1$.}\/ 
We also let $A\in K[\der]$ have order $r$, and we let $\fm$ range over the elements of $H^\times$ such that
$v\fm\in v(\hat h -H)$. We begin with an important consequence of the material in
Section~\ref{sec:bounding}:

\begin{lemma}\label{lem:8.8 refined}
Suppose $(P,1,\hat h)$ is $Z$-minimal, deep, and special,  and $\fv(L_P)\asymp\fv:=\fv(A)$. Let
$y\in\Calr[\imag]$ satisfy $A(y)=0$ and $y\prec \fm$ for all~$\fm$.  
Then $y',\dots,y^{(r)}\prec \fm$ for all $\fm$.
\end{lemma}
\begin{proof} Corollary~\ref{specialvariant} gives an $\fm\preceq \fv$, so it is enough to show $y',\dots,y^{(r)}\prec \fm$ for all $\fm\preceq \fv$. Accordingly we assume $0<\fm\preceq \fv$ below. 
As~$\hat h$ is special over $H$, we have~${2(r+1)v\fm}\in v(\hat h-H)$, so $y\prec\fm^{2(r+1)}$.
Then    Corollary~\ref{cor:EL} with~$n=2(r+1)$,~$\eta=\abs{\fv}^{-1}$, $\varepsilon=1/r$ gives  for $j=0,\dots,r$:
\[y^{(j)}\ \prec\ \fv^{-j} \fm^{n-j(1+\varepsilon)}\ \preceq\ \fm^{n-j(2+\varepsilon)}\ \preceq\ \fm^{n-r(2+\varepsilon)}\  =\ \fm. \qedhere\]
\end{proof}

\noindent
Note that by Proposition~\ref{prop:existence and uniqueness}, if~$\dim_{\C}\ker_{\Univ} A =r$ and~$A(y)=0$, $y\in\Calr[\imag]$, then~$y\in \Univ=K[\ex^{H\imag}]\subseteq\Calinf[\imag]$. Corollary~\ref{cor:gaussian ext dom} is typically used in combination with the ultimate condition. Here is a first easy application:

\begin{lemma}\label{lem3.9 linear}
Suppose $(P,1,\hat h)$ is linear and ultimate, $\dim_{\C}\ker_{\Univ}L_P=r$, and~$y\in\Calr[\imag]$ satisfies $L_P(y)=0$ and~$y\prec 1$. Then~$y\prec \fm$ for all $\fm$.
\end{lemma}
\begin{proof}
We have~$y\in \Univ$, so $y\prec_{\g} 1$ by Lemma~\ref{lem:gaussian ext dom}. If $y=0$ we are done, so assume~${y\ne 0}$. %Set $A:=L_P$. 
 Lemma~\ref{lem:excu properties}(ii)  
%(if $r\ge2$) and Corollary~\ref{cor:excu properties}(ii) (if $r=1$) 
gives~$0 < v_{\g}y\in
v_{\g}(\ker^{\neq}_{\Univ}L_P) = \exc^{\operatorname{u}}(L_P)$, hence~$v_{\g}y > v(\hat h -H)$ by   Lemma~\ref{lem:ultimate deg 1},
so $y\prec_{\g}\fm$ for all $\fm$.  
Now Corollary~\ref{cor:gaussian ext dom} yields the desired conclusion. 
\end{proof}

\begin{cor}\label{cor:3.6 linear}
Suppose that $(P,1,\hat h)$ is $Z$-minimal, deep,  special,  linear, and ultimate,  and that  $\dim_{\C}\ker_{\Univ}L_P=r$.
Let~$f,g\in \Calr[\imag]$ be such that $P(f)=P(g)=0$ and $f,g\prec 1$.
Then~$({f-g})^{(j)}\prec \fm$ for~$j=0,\dots,r$ and all $\fm$.  
\end{cor}
\begin{proof}
Use Lemmas~\ref{lem:8.8 refined} and~\ref{lem3.9 linear} for $A=L_P$ and $y=f-g$.
\end{proof}

\noindent
{\it In the rest of this subsection we assume that
   $(P,1,\hat h)$ is ultimate and normal, $\dim_{\C} \ker_{\Univ}A=r$, and~$L_{P}=A+B$ where}
 $$B\prec_{\Delta(\fv)} \fv^{r+1}A, \qquad \fv:=\fv(A)\prec^\flat 1.$$
Then Lemma~\ref{lem:fv of perturbed op} gives   $\fv(L_P)\sim\fv$,  and by Lem\-ma~\ref{lem:excu properties},
$$v_{\operatorname{g}}(\ker^{\neq}_{\Univ} A)\ =\ \exc^{\operatorname{u}}(A)\ =\ \exc^{\operatorname{u}}(L_P).$$
This yields a variant of Lemma~\ref{lem3.9 linear}:

\begin{prop}\label{lem3.9}  
If $y\in\Calr[\imag]$ and $A(y)=0$, $y\prec 1$,
then~$y\prec \fm$ for all~$\fm$.
\end{prop}
\begin{proof} 
Like that of Lemma~\ref{lem3.9 linear}, using Lemma~\ref{lem:ultimate normal} instead of~\ref{lem:ultimate deg 1}.
\end{proof}
%\noindent
%In the rest of this subsection we refine Proposition~\ref{lem3.9}.

\noindent
The following result will be used in establishing a crucial non-linear version of Corollary~\ref{cor:3.6 linear}, namely Proposition~\ref{prop:notorious 3.6}.

\begin{cor}\label{cor:8.8 refined}
If $(P,1,\hat h)$ is $Z$-minimal, deep, and special, and
$y\in\Calr[\imag]$ is such that $A(y)=0$ and $y\prec 1$,  
then $y,y',\dots,y^{(r)}\prec \fm$ for all $\fm$.
\end{cor}
\begin{proof}
Use first  Proposition~\ref{lem3.9} and then Lemma~\ref{lem:8.8 refined}.
\end{proof}

%\begin{cor}\label{lem3.9 firm} \marginpar{now part of 7.6.7}
%Suppose $(P,1,\hat h)$ is firm.
%Then the only   $y\in\Calr[\imag]$   such that~${A(y)=0}$ and $y\prec 1$ is $y=0$.
%\end{cor}
%\begin{proof}
%Suppose $y\in\Calr[\imag]$ satisfies $A(y)=0$ and $0\neq y\prec 1$. Then $y\in \Univ$, so~$y\prec_{\g} 1$   by Lemma~\ref{lem:gaussian ext dom}. As in the proof of Lemma~\ref{lem3.9 linear} we obtain
%$v_{\g}y\in\exc^{\operatorname{u}}(L_P)$,
%hence~$v_{\g}y\leq 0$ by the remark following 
%Lemma~\ref{lem:firm deg 1}, contradicting $y\prec_{\g}1$.
%\end{proof}

%\begin{cor}\label{cor:8.8 flabby} \marginpar{now part of 7.6.8}
%Suppose $(P,1,\hat h)$ is $Z$-minimal, deep,  special, and flabby. Then there exists $y\in\ker^{\neq}_{\Univ} %A\subseteq\Calinf[\imag]$ such that $y,y',\dots,y^{(r)}\prec  \fm$ for all $\fm$.
%\end{cor}
%\begin{proof}
%%The proof of Proposition~\ref{lem3.9} shows $v_{\operatorname{g}}(\ker^{\neq}_{\Univ} A)=\exc^{\operatorname{u}}
%%(A)=\exc^{\operatorname{u}}(L_P)$. 
%Flabbiness of $(P,1,\hat h)$ gives
%$y\in\ker^{\neq}_{\Univ} A$ such that  $v_{\g}y > v(\hat h -H)$, by Lem\-ma~\ref{lem:firm normal}, so 
%$y\prec_{\operatorname{g}} \fm$ for all $\fm$. Then $y\prec \fm$ for all $\fm$ by Corollary~\ref{cor:gaussian ext dom},
%and thus~${y',\dots,y^{(r)}\prec \fm}$ for all $\fm$ by Lemma~\ref{lem:8.8 refined}. 
%\end{proof}

\noindent
So far we didn't have to name an immediate asymptotic extension of $H$ where $\hat h$ is located, but for the ``complex'' version of the above we need to be more specific. 

As in the beginning of Section~\ref{sec:ultimate}, let $\hat H$ be an immediate asymptotic extension of $H$ and 
$\hat{K}=\hat{H}[\imag]\supseteq \hat{H}$ a corresponding immediate $\d$-valued extension of $K$.  The results in this subsection then go through  if instead of
$(P, 1,\hat h)$ being a slot in~$H$ of order $r\ge 1$  we assume that
$(P, 1,\hat h)$ is a slot in $K$ of order $r\ge 1$ with $\hat h\in \hat K \setminus K$, with $\fm$ now ranging over the elements of~$K^\times$ such that $v\fm\in v(\hat h - K)$. 
%$\hat h-h\preceq\fm$ for some $h\in K$.  

\subsection*{Solution spaces of linear differential operators} 
Recall that $\Lambda\subseteq H\imag$,
$\Univ=K\big[\!\ex(\Lambda)\big]=K[\ex^{H\imag}]$ where $\ex(\Lambda)\subseteq \ex^{H\imag}\subseteq \Calinf[\imag]^\times$.
Hence for each $\lambda$ we have an element $\phi(\lambda)$ of $H$ (unique up to addition of an element of $2\pi\Z$) such that~$\ex(\lambda)=\ex^{\phi(\lambda)\imag}$; we take $\phi(0):=0$. Then $\ex(\lambda)^\dagger=\lambda$ gives $\phi(\lambda)'\imag=\lambda$, and 
% and $\overline{\ex(\lambda)}=\ex(-\lambda)$ for each $\lambda$. 
$$
\phi(\lambda_1+\lambda_2)\ \equiv\ \phi(\lambda_1)+\phi(\lambda_2)\bmod 2\pi\Z \quad\text{for $\lambda_1,\lambda_2\in\Lambda$.}$$
If $\I(K)\subseteq K^\dagger$, then  $\Lambda\cap \I(H)\imag=\{0\}$ (see Lemma~\ref{lem:W and I(F)}), so $\phi(\lambda)\succ 1$ for~$\lambda\neq 0$, hence for $\mu\in\Lambda$: $\lambda=\mu \Leftrightarrow \phi(\lambda)=\phi(\mu)\Leftrightarrow \phi(\lambda)-\phi(\mu) \preceq 1$. 

\begin{lemma}\label{lem:from phi to lambda}
Let $\phi\in H$. Then there exists $\lambda$ such that $\phi-\phi(\lambda)\preceq 1$. If~$\phi\succ 1$, then for any such $\lambda$ we have
$\sgn \phi=\sgn\Im\lambda$.
\end{lemma}
\begin{proof}
From $\ex^{\phi\imag}\in \ex^{H\imag}\subseteq\Univ^\times=K^\times\ex(\Lambda)$ we get $f\in K^\times$ and  $\lambda$ with $\ex^{\phi\imag}=f\ex(\lambda)=f\ex^{\phi(\lambda)\imag}$. Note that $\abs{f}=1$, so $f=\ex^{\theta\imag}$ where $\theta\in H$ with $\theta\preceq 1$, by Lemma~\ref{lem:fexphii}. This yields $\phi-\phi(\lambda)-\theta\in 2\pi\Z$ and so $\phi-\phi(\lambda)\preceq 1$. This proves the first statement.
Now suppose we have any $\lambda$ with $\phi-\phi(\lambda)\preceq 1$. Then $\phi\sim\phi(\lambda)$ if $\phi\succ 1$. So if~$\phi>\R$, then $\phi(\lambda)>\R$ and thus $\Im\lambda=\phi(\lambda)'>0$; likewise,   $\phi<\R$ implies $\Im\lambda<0$.
\end{proof}

\begin{cor}\label{cor:from phi to lambda}
Suppose $\I(K)\subseteq K^\dagger$.
Let $f=f_1 \ex^{\phi_1\imag}+\cdots+f_m\ex^{\phi_m\imag}\in\Univ$ where~$f_1,\dots,f_m\in K$ and~$\phi_1,\dots,\phi_m\in H$ are such that $\phi_j=\phi_k$ or $\phi_j-\phi_k\succ 1$
for~$j,k=1,\dots,m$. Then 
$$f=0 \quad\Longleftrightarrow\quad
\sum_{1\leq k\leq m,\ \phi_k=\phi_j} f_k\  =\  0\ \text{ for $j=1,\dots,m$,}$$
and for $\fm\in H^\times$:
$$f\prec  \fm\quad\Longleftrightarrow\quad
\sum_{1\leq k\leq m,\ \phi_k=\phi_j} f_k\ \prec\ \fm\ \text{for $j=1,\dots,m$,}$$
and likewise with $\preceq$ in place of $\prec$.
\end{cor}
\begin{proof}
We first arrange that $\phi_1,\dots,\phi_m$ are distinct,  and we then need to show: 
${f=0 \Leftrightarrow f_1=\cdots=f_m=0}$, and
$f\prec\fm \Leftrightarrow f_1,\dots,f_m\prec\fm$, and likewise with $\preceq$ in place of $\prec$. To make Corollary~\ref{cor:gaussian ext dom} applicable we also arrange that $\Lambda=\Lambda_H\imag$ with~$\Lambda_H$ an $\R$-linear complement of $\I(H)$ in $H$.  
Lemma~\ref{lem:from phi to lambda} yields~$\lambda_j\in\Lambda$ with~$\phi_j-\phi(\lambda_j)\preceq 1$
for $j=1,\dots,m$; then $\lambda_1,\dots,\lambda_m$ are distinct.
For~${j=1,\dots,m}$, put~$g_j:=f_j\ex^{(\phi_j-\phi(\lambda_j))\imag}\in K$, so
$f_j\ex^{\phi_j\imag}=g_j\ex(\lambda_j)$ and $g_j\asymp\abs{g_j}=\abs{f_j}\asymp f_j$. Now
the claim follows from the $K$-linear independence of $\ex(\lambda_1),\dots,\ex(\lambda_m)$,
Corollary~\ref{cor:gaussian ext dom}, and Lemma~\ref{lem:gaussian ext dom, preceq}.
\end{proof}

\noindent
Let $A\in K[\der]^{\neq}$, $r:=\order A$, and set
$V:=\ker_{\Univ} A$, a $\mathbb C$-linear subspace of~$\Univ$ of dimension at most $r$, with
$\dim_{\mathbb C} V=r$ iff $V=\ker_{\Calinf[\imag]} A$.
We describe in our present setting some consequences of the results 
obtained in Sections~\ref{sec:splitting} and~\ref{sec:eigenvalues and splitting}  about  zeros of linear differential operators in the universal exponential extension.

\begin{lemma}\label{lem:complex basis}
The $\mathbb C$-linear space $V$ has a basis
$$f_1\ex(\lambda_1),\dots,f_d\ex(\lambda_d)\quad\text{ where $f_j\in K^\times$,    
$\lambda_j\in\Lambda$   \textup{(}$j=1,\dots,d$\textup{)}.}$$
For any such basis the set of eigenvalues of $A$ with respect to $\Lambda$ is $\{\lambda_1,\dots,\lambda_d\}$,
and 
$$\operatorname{mult}_{\lambda}(A)\ =\ \abs{ \{ j\in\{1,\dots,d\}: \lambda_j=\lambda \} }\quad\text{ for every $\lambda$.}$$
\end{lemma}

\noindent
This follows from Lemma~\ref{newlembasis} and the considerations preceding it. %, and Corollary~\ref{cor:basis of kerUA}:

Call $\phi_1,\dots,\phi_m\in H$   {\bf apart}\index{apart}\index{germ!apart} 
%{\bf non-interfering} 
if  $\phi_j=0$ or~$\phi_j\succ 1$ for $j=1,\dots,m$, and  
$\phi_j=\phi_k$ or~$\phi_j-\phi_k\succ 1$ for $j,k=1,\dots,m$. (This holds in particular if $\phi_1=\cdots=\phi_m=0$.) If~$\I(K)\subseteq K^\dagger$,
then~$\phi(\lambda_1),\dots,\phi(\lambda_m)$ are apart for any $\lambda_1,\dots,\lambda_m\in\Lambda$.

\begin{cor} \label{cor:complex basis} 
The
 $\mathbb C$-linear space $V$ has a basis
$$f_1\ex^{\phi_1\imag},\dots,f_d\ex^{\phi_d\imag}\quad\text{ where $f_j\in K^\times$,  $\phi_j  \in H$ \textup{(}$j=1,\dots,d$\textup{)}.}$$
If $\I(K)\subseteq K^\dagger$, then
for any such basis the $f_j\ex^{\phi_j\imag}$ with $\phi_j\preceq 1$ form a basis of the $\C$-linear space 
$V\cap K=\ker_K A$, and we can choose the $f_j$, $\phi_j$ 
such that additionally~$\phi_1,\dots \phi_d$ are apart and $(\phi_1,vf_1),\dots,(\phi_d,vf_d)$ are distinct.
\end{cor}
\begin{proof}
The first claim holds by Lemma~\ref{lem:complex basis}. Suppose $\I(K)\subseteq K^\dagger$, and
let a basis of $V$ as in the corollary be given. Then by Lemma~\ref{lem:from phi to lambda} 
we obtain $\lambda_j\in\Lambda$ such that $\phi_j-\phi(\lambda_j)\preceq 1$, and so $f_j\ex^{\phi_j\imag}=g_j\ex(\lambda_j)$
where~$g_j:=f_j\ex^{(\phi_j-\phi(\lambda_j))\imag}\in K^\times$ by Proposition~\ref{prop:cos sin infinitesimal, 2}.
Now~${\lambda_j=0}\Leftrightarrow\phi_j\preceq 1$, by the remarks preceding Lem\-ma~\ref{lem:from phi to lambda}, hence the $f_j\ex^{\phi_j\imag}$ with $\phi_j\preceq 1$ form a basis of~$V\cap K=\ker_K A$.
Moreover,~$g_1\ex^{\phi(\lambda_1)\imag},\dots,g_d\ex^{\phi(\lambda_d)\imag}$ is a basis of $V$ and $\phi(\lambda_1),\dots, \phi(\lambda_d)$ are apart. 

We have~$V=\bigoplus_\lambda V_\lambda$ 
 (internal direct sum of $\C$-linear subspaces) where   $V_\lambda=(\ker_K A_\lambda) \ex(\lambda)$, by the remarks before~\eqref{eq:bd on sum mults}.
For each $\lambda$, the
subspace $\ker_K A_\lambda$ of the $\C$-linear space $K$ is generated by the~$g_j$  with~$\lambda_j=\lambda$. 
  Applying~[ADH, 5.6.6] to each  $A_\lambda$ we  obtain $h_j\in K^\times$  such that
  $h_1\ex^{\phi(\lambda_1)\imag},\dots,h_d\ex^{\phi(\lambda_d)\imag}$ is a basis of~$V$
  where for all $j\ne k$ with $\phi(\lambda_j)=\phi(\lambda_k)$ we have $h_j\nasymp h_k$.
\end{proof}

\noindent 
A {\bf Hahn basis} of $V$ is a basis of $V$ as in Corollary~\ref{cor:complex basis}  such that $\phi_1,\dots,\phi_d$ are apart and 
 $(\phi_1,vf_1),\dots,(\phi_d,vf_d)$ are distinct.\index{basis!Hahn}\index{Hahn basis}\index{linear differential operator!Hahn basis} (It should really be ``Hahn basis with respect to $\phi_1,\dots, \phi_d$'' but in the few cases we use this notion we shall rely on the context as to what tuple $(\phi_1,\dots,\phi_d)\in H^d$ we are dealing with.) If $\I(K)\subseteq K^\dagger$, then for such a Hahn basis
 the $f_j$ with $\phi_j=0$ form a valuation basis of the  subspace~$V\cap K$
  of the valued $\C$-linear space~$K$~[ADH, 2.3]. 

In the next lemma we assume $\I(K)\subseteq K^\dagger$, and we recall that then 
$$d\ \leq\ \sum_\lambda \, \abs{\exc^{\ev}(A_\lambda)}\  \leq\  r$$
by Lemma~\ref{lem:finiteness of excu(A)} and  Proposition~\ref{prop:finiteness of excu(A), real}, and so by Lemma~\ref{lem:finiteness of excu(A)},
$$\sum_\lambda\, \abs{\exc^{\ev}(A_\lambda)}\ =\ d\ \Longrightarrow\ \exc^{\ev}(A_\lambda)\ =\ v(\ker^{\neq}  A_\lambda)\
\text{ for all }\lambda.$$

\begin{lemma}\label{lem:complex basis, excev}
Suppose $\I(K)\subseteq K^\dagger$, and let $ f_1\ex^{\phi_1\imag},\dots,f_d\ex^{\phi_d\imag}$ be a Hahn basis of $V$ 
as in Corollary~\ref{cor:complex basis}. Then for all $\lambda$, 
$$\exc^{\ev}(A_\lambda)\ \supseteq\ v(\ker^{\neq} A_\lambda) \ =\  \big\{ vf_j:\  1\le j\le d,\ \phi_j-\phi(\lambda) \preceq 1 \big\}.$$
and so $\exc^{\operatorname{u}}(A) \supseteq \{ vf_1,\dots,vf_d\}$, with
$\exc^{\operatorname{u}}(A) = \{ vf_1,\dots,vf_d\}$ if  $\sum_\lambda \abs{\exc^{\ev}(A_\lambda)}=d$.
\end{lemma}
\begin{proof}
Take $g_j$, $\lambda_j$ as in the proof of Corollary~\ref{cor:complex basis}. Then 
$g_j\asymp \abs{g_j}=\abs{f_j}\asymp f_j$, and $\lambda_j=\lambda\Leftrightarrow \phi_j-\phi(\lambda)\preceq 1$, for all $\lambda$.  So we can replace $f_j$, $\phi_j$ by $g_j$, $\phi(\lambda_j)$ to arrange $\phi_j=\phi(\lambda_j)$ for $j=1,\dots,d$.
Then for all $\lambda$ the $\C$-linear space $\ker  A_\lambda\subseteq K$ is generated  by the $f_j$ with $\lambda_j=\lambda$, so $$\exc^{\ev}(A_\lambda)\ \supseteq\ v(\ker^{\neq} A_\lambda)\ =\ \{vf_j:\  1\le j\le d,\ \lambda_j=\lambda\}.$$
For the rest use $\exc^{\operatorname{u}}(A)=\bigcup_\lambda \exc^{\ev}(A_\lambda)$ and the remarks preceding the lemma.
\end{proof}

\noindent
Corollaries~\ref{cor:basis of kerUA} and~\ref{spldcr2} yield  conditions on $A$, $K$ that guarantee $\dim_{\C} V=r$:

\begin{lemma}\label{lem:full kernel} 
Suppose $A$ splits over $K$. If~$r\le 1$, or $r=2$, $A\in H[\der]$, or $K$ is
$1$-linearly surjective,  then $$\dim_{\mathbb C} V\ =\ \sum_\lambda \operatorname{mult}_{\lambda}(A)\ =\ r.$$
\end{lemma}

\noindent
Next a complement to Lemma~\ref{lem:trig surj, complex}:

\begin{cor}\label{cor:trig surj, complex} 
Suppose $K$ is $r$-linearly surjective, or $K$ is $1$-linearly surjective and $A$ splits over $K$.
Let $\phi\in H$ be such that $\phi'\imag+K^\dagger$ is not an eigenvalue of $A$. Then $A$ maps $K\ex^{\phi\imag}$ bijectively onto $K\ex^{\phi\imag}$.
%for each $b\in K$ there is a unique $y\in K$ with $A(y\ex^{\phi\imag})=b\ex^{\phi\imag}$.
\end{cor}
\begin{proof}
Let $y\in K$, $A(y\ex^{\phi\imag})=0$. By Lemma~\ref{lem:trig surj, complex} it is enough to show that~${y=0}$.
Suppose towards a contradiction that $y\neq 0$. Then $y\ex^{\phi\imag}\in\Univ^\times=K^\times\ex(\Lambda)$, so~$y\ex^{\phi\imag}=z\ex(\lambda)$, $z\in K^\times$. Then $\phi'\imag-\lambda\in K^\dagger$, so $\lambda$ is not an eigenvalue
of~$A$ with respect to~$\Lambda$. Since $y\in V$, this contradicts Lemma~\ref{lem:complex basis}.
\end{proof}

\noindent
Let $N$ be an $n\times n$ matrix over $K$, $n\ge 1$. We end this subsection with a variant of Lemma~\ref{lem:complex basis}  for the
matrix differential equation $y'=Ny$. Set
$S := \operatorname{sol}_{\Univ}(N)$,
so~$S$ is a $\C$-linear subspace of $\Univ^n$ of dimension $\le n$. \label{p:solU}
%Let $\Omega:= \text{Frac} \Univ$. 

\begin{lemma}\label{matrixversionlem} Suppose $S$ has a basis 
$$\ex^{\phi_1\imag}f_1,\dots, \ex^{\phi_d\imag} f_d\quad\text{ where $\phi_1,\dots, \phi_d\in H$ and $f_1,\dots, f_d\in K^n\subseteq \Univ^n$.}$$ Set $\alpha_j:= \phi_j'\imag+K^\dagger\in K/K^\dagger$ for $j=1,\dots, d$. Then
$$\mult_{\alpha}(N)\ =\ \abs{\{j\in \{1,\dots, d\}:\ \alpha_j=\alpha\}}  \quad \text{ for all }\alpha\in K/K^\dagger.$$
\end{lemma}
\begin{proof} We have $f_j=(f_{1j},\dots, f_{nj})^{\text{t}}\in \Univ^n$ for $j=1,\dots,d$. We first consider the case that $N$ is the companion matrix of a monic $B\in K[\der]$ of order $n$.
Then we have the $\C$-linear isomorphism $z\mapsto (z,z',\dots, z^{(n-1)})^{\text{t}}\colon\ker_{\Univ} B \to S$; its inverse maps the given basis to
 a basis  $\ex^{\phi_1\imag}f_{11},\dots, \ex^{\phi_d\imag}f_{1d}$ of $\ker_{\Univ} B$. 
For $j=1,\dots,d$ we have~$\phi_j-\phi(\lambda_j)\preceq 1$ with $\lambda_j\in \Lambda$, and so this basis has the form
$g_1\ex(\lambda_1),\dots, g_d\ex(\lambda_d)$ with $g_1,\dots, g_d\in K^\times$. Now use Lemmas~\ref{lem:matrix diff equs vs ops} and~\ref{lem:complex basis}, and the fact that~$\alpha_j=\lambda_j+K^\dagger$ for $j=1,\dots,d$.

For the general case, [ADH, 5.5.9] gives the companion matrix $M$ of a monic~$B\in K[\der]$ of order $n$ such that
$y'=Ny$ is equivalent to $y'=My$. This yields $P\in \text{GL}_n(K)$ such that
$f\mapsto Pf\colon S \to \text{sol}_{\Univ}(M)$ is a $\C$-linear isomorphism, and so~$PS= \text{sol}_{\Univ}(M)$.
Since $P\ex^{\phi_j\imag}f_j=\ex^{\phi_j\imag}g_j$ with $g_j\in K^n$ for $j=1,\dots,d$, we obtain a basis~$\ex^{\phi_1\imag}g_1,\dots, \ex^{\phi_d\imag}g_d$ of the $\C$-linear subspace $\text{sol}_{\Univ}(M)$ of $\Univ^n$,
so we are in the special case treated earlier. 
\end{proof}

\subsection*{A relative version of Corollary~\ref{cor:complex basis}\astr} In this subsection $\I(K)\subseteq K^\dagger$. We use an isomorphism as in Lemma~\ref{exqe} to identify $\Univ=K[\ex(\lambda)]$ with $K[\ex^{H\imag}]$. 

Let $F$ be a Liouville closed Hardy field extension of $H$; set~$L:=F[\imag]\subseteq \Calinf[\imag]$. We show here how
various results about $H$, $K$ extend  in a coherent way to $F, L$.  First,  Corollary~\ref{cor:LambdaL, purely imag} yields
a complement~$\Lambda_L$ of the $\Q$-linear subspace~$L^\dagger$ of~$L$
with
$\Lambda\subseteq\Lambda_L\subseteq F\imag$. 
Let $\Univ_L=L\big[\!\ex(\Lambda_L)\big]$ be
the universal exponential extension  
of~$L$ containing $\Univ=K\big[\!\ex(\Lambda)\big]$ as a differential subring  
 described in the remarks following Corollary~\ref{cor:Univ under d-field ext}.
We also have the differential subring~$L[\ex^{F\imag}]$ of $\Calinf[\imag]$ with
$\Univ=K[\ex^{H\imag}]\subseteq L[\ex^{F\imag}]$.  

\begin{lemma}\label{lem:U -> U_L}
There is an isomorphism $\iota\colon\Univ_L\to L[\ex^{F\imag}]$ of differential $L$-algebras with
$\iota\big(\!\ex(\Lambda_L)\big)\subseteq \ex^{F\imag}$ that is the identity on $\Univ$. Thus the diagram below commutes: 
$$\xymatrix@C=5em{\Univ_L \ar@{-->}[r]^\iota_{\cong} & L[\ex^{F\imag}] \\
\Univ \ar@{=}[r] \ar@{-}[u]^{\subseteq} & K[\ex^{H\imag}]\ar@{-}[u]^{\subseteq}}$$
\end{lemma}
\begin{proof}
Lemma~\ref{exqe} yields an isomorphism $\iota_L\colon \Univ_L\to L[\ex^{F\imag}]$ of differential $L$-algebras with
$\iota_L\big(\!\ex(\Lambda_L)\big)\subseteq\ex^{F\imag}$.
By Lemma~\ref{lem:Univ under d-field ext} we have $\iota_L^{-1}\big(K[\ex^{H\imag}]\big)=K[E]$ 
where~$E=\{{u\in\Univ^\times_L:u^\dagger\in K}\}$. From $\Univ^\times_L=L^\times\ex(\Lambda_L)$
we get $E=K^\times\ex(\Lambda)$, so~$K[E]=\Univ$. 
Hence~$\iota_L^{-1}$ restricts to an automorphism of the differential $K$-algebra~$\Univ$.  So this restriction
equals $\sigma_\chi$ where $\chi\in\operatorname{Hom}(\Lambda,\C^\times)$. (Lemma~\ref{autolem}.) 
Extending $\chi$ to $\chi_L\in \operatorname{Hom}(\Lambda_L,\C^\times)$ yields an isomorphism
$$\iota:= \iota_L\circ\sigma_{\chi_L}\  \colon\ \Univ_L\to   L[\ex^{F\imag}]$$ 
 of differential $L$-algebras with the desired property.
\end{proof}

\noindent
Fix an isomorphism $\iota\colon\Univ_L\to L[\ex^{F\imag}]$ as in the previous lemma and identify~$\Univ_L$ with its image via $\iota$; thus $\Univ=K[\ex^{H\imag}]\subseteq L[\ex^{F\imag}]=\Univ_L\subseteq\Calinf[\imag]$.
For each $\mu\in\Lambda_L$ we have an element $\phi(\mu)$ of $F$ (unique up to addition of an element of $2\pi\Z$) such that~$\ex(\mu)=\ex^{\phi(\mu)\imag}$; we take $\phi(0):=0$. The $\phi(\lambda)\in F$ are actually in $H$ and agree with the
$\phi(\lambda)$ defined earlier, up to addition of elements of $2\pi\Z$.
{\it In the rest of this subsection we assume $\I(L)\subseteq L^\dagger$.}\/
So for $\mu_1,\mu_2\in\Lambda_L$: $\mu_1=\mu_2\Leftrightarrow \phi(\mu_1)-\phi(\mu_2)\preceq 1$.  

\begin{lemma}\label{lem:phi(mu)}
Let $\mu\in\Lambda_L$. Then 
$$\phi(\mu)\in H\ \Longleftrightarrow\  
\phi(\mu)\in H+\mathcal O_F  
\ \Longleftrightarrow\  
\mu\in\Lambda.$$
\end{lemma}
\begin{proof}
We have $\mu=\phi(\mu)'\imag$.
So if $\phi(\mu)\in H+\mathcal O_F$, then $\mu\in H\imag+\I(L)\subseteq K+L^\dagger=\Lambda+L^\dagger$, 
and hence $\mu\in \Lambda$.
Conversely, if $\mu\in\Lambda$,
then $\phi(\mu)'\imag\in \Lambda\subseteq H\imag$, so $\phi(\mu)'\in H$, and thus $\phi(\mu)\in H$.
%and since~$H\supseteq\R$ is closed under integration, we get $\phi(\mu)\in H$.
\end{proof}

\noindent
Lemma~\ref{lem:from phi to lambda} with $F$, $L$, $\Lambda_L$ in place of $H$, $K$, $\Lambda$, and Lemma~\ref{lem:phi(mu)} yield:

\begin{cor} \label{cor:phi(mu)} $H+\mathcal O_F=\big\{\phi\in F:\text{$\phi-\phi(\lambda)\preceq 1$ for some $\lambda$}\big\}$.
\end{cor}

\noindent
Let   $A\in K[\der]^{\neq}$, $r:=\order A$, $V:=\ker_{\Univ} A$, $V_L:=\ker_{\Univ_L}A$,  so  
$V  =V_L\cap\Univ$.
Corollary~\ref{cor:complex basis} applied to $F$, $L$ in place of $H$, $K$ then gives  a Hahn basis
 $$f_1\ex^{\phi_1\imag},\ \dots,\ f_d\ex^{\phi_d\imag} \qquad (f_j\in L^\times, \phi_j\in F)$$
of  $V_L$.
Recall from Corollary~\ref{cor:LambdaL, purely imag} that  $\exc^{\ev}(A_\lambda) = \exc^{\ev}_L(A_\lambda)\cap\Gamma$ for all $\lambda$. Applying Lemma~\ref{lem:complex basis, excev}   to  such a Hahn basis of~$V_L$  and $F$, $L$, $\Lambda_L$, $V_L$ in place of~$H$,~$K$,~$\Lambda$,~$V$, and using Corollary~\ref{cor:phi(mu)} we obtain:
$$\exc^{\operatorname{u}}(A)\  \supseteq\ \{vf_j : j=1,\dots,d,\   \phi_j\in H+\mathcal O_F\}\cap\Gamma.$$
Recall   from Corollary~\ref{cor:excev cap GammaOmega, real} that if $A$
is terminal, then $\exc^{\ev}(A_{\lambda}) = \exc^{\ev}_L(A_\lambda)$ for all~$\lambda$, and~$\exc^{\operatorname{u}}(A)  = \exc^{\operatorname{u}}_L(A)$.
We have $d=\dim_{\C} V_L\leq r$, and
   Lemma~\ref{lem:full kernel}  gives  conditions on~$A$,~$F$,~$L$ which guarantee~$d=r$. 
  The next corollary  shows that if $A$ is terminal and $d=r$, 
then the ``frequencies''~$\phi_j$ of the elements of our Hahn basis of~$V_L$ above can be taken in~$H$:

\begin{cor}\label{cor:A ultimate}
Suppose $A$ is terminal and $d=r$. Then $V_L$ has a Hahn basis
 $$f_1\ex^{\phi_1\imag},\ \dots,\ f_r\ex^{\phi_r\imag} \qquad (f_j\in L^\times,\ \phi_j\in H).$$
For any such basis and all $\lambda$ we have
$$\exc^{\ev}(A_{\lambda})\ =\  \exc^{\ev}_L(A_\lambda)\ =\  v(\ker_L^{\neq}A_\lambda)\  =\ 
\big\{vf_j:j=1,\dots,r,\  \phi_j  - \phi(\lambda) \preceq 1 \big\}$$
and
$\exc^{\operatorname{u}}(A) = \exc^{\operatorname{u}}_L(A) = \{vf_1,\dots,vf_r\}\subseteq\Gamma$, and
the eigenvalues of $A$ viewed as element of $L[\der]$ are $\phi_1'\imag+L^\dagger,\dots,\phi_r'\imag+L^\dagger$.
\end{cor}
\begin{proof}
 Lemma~\ref{lem:finiteness of excu(A)}  and Corollary~\ref{cor:excev cap GammaOmega, real}  give
$\exc^{\ev}(A_{\lambda}) = \exc^{\ev}_L(A_\lambda) = v(\ker_L^{\neq}A_\lambda)$ for all $\lambda$,
and $\ker_L^{\neq} A_\mu=\emptyset$ for $\mu\in\Lambda_L\setminus\Lambda$. 
Take any Hahn basis of $V_L$ as described before the corollary.    Lemma~\ref{lem:from phi to lambda} yields $\lambda_j\in\Lambda_L$ with $\phi_j-\phi(\lambda_j)\preceq 1$. We have~$g_j:=f_j\ex^{(\phi_j-\phi(\lambda_j))\imag}\in L^\times$ by Proposition~\ref{prop:cos sin infinitesimal, 2} and
$f_j\ex^{\phi_j\imag}=g_j\ex(\lambda_j)$, so~$g_j\in\ker_L^{\neq} A_{\lambda_j}$.
This yields $\lambda_j\in\Lambda$ for $j=1,\dots,r$. 
Replacing each pair~$f_j$,~$\phi_j$ by~$g_j$, $\phi(\lambda_j)$ we obtain a Hahn basis
of $V_L$ as claimed.
The rest follows from Lemmas~\ref{lem:complex basis, excev} and ~\ref{lem:complex basis}.
\end{proof}

\subsection*{Duality considerations\astr}  As before, $A\in K[\der]^{\ne}$ has order
$r$, and $V:=\ker_{\Univ}A$. 
Recall from Section~\ref{sec:self-adjoint} the bilinear form $[\ ,\, ]_A$ on the $\C$-linear space $\Omega=\Frac(\Univ)$.
As in  the previous subsection we take $\Univ=K[\ex^{H\imag}]$ and fix values $\ex(\lambda)$. 

\begin{cor}\label{cor:conjugate basis}
Suppose $A$ splits over $K$, $\I(K)\subseteq K^\dagger$, and $r\le 1$ or $K$ is $1$-linearly surjective. 
Let $f_j$, $\phi_j$ be  as in Corollary~\ref{cor:complex basis}.
Then the $\C$-linear space~$\ker_{\Calinf[\imag]}A^*$  equals $W:=\ker_{\Univ} A^*$ and has a basis
$$f_1^*\ex^{-\phi_1\imag},\ \dots,\ f_r^*\ex^{-\phi_r\imag}\qquad\text{where $f_k^*\in K^\times$  \textup{(}$k=1,\dots,r$\textup{)}}$$
such that $\big[f_j \ex^{\phi_j\imag}, f_k^* \ex^{-\phi_k\imag}\big]_A=\delta_{jk}$
for $j,k=1,\dots,r$.
\end{cor}
\begin{proof}
By Lemma~\ref{lem:full kernel} we have $\dim_{\C} V = \dim_{\C} W = r$. 
As in the proof of Corollary~\ref{cor:complex basis} we obtain $g_j\in K^\times$, $\lambda_j\in\Lambda$ with~$\phi_j-\phi(\lambda_j)\preceq 1$,
and
$$y_j\ :=\ f_j \ex^{\phi_j\imag} = g_j \ex(\lambda_j)\in \Univ^\times,\quad j=1,\dots,r.$$  The basis $y_1,\dots, y_r$ of $V$ yields by Corollary~\ref{corbasissplit} that~${A=a(\der-a_r)\cdots(\der-a_1)}$ with~$a\in K^\times$ and $(a_1,\dots, a_r)=\operatorname{split}(y_1,\dots, y_r)$.  It is easy to reduce to the case~${a=1}$. 
Then Corollary~\ref{cor:dual split}  provides a basis
$y_1^*,\dots,y_r^*$  of $W$ with~$[y_j, y_k^*]_A=\delta_{jk}$ for all~$j$,~$k$,  $\operatorname{split}(y_r^*,\dots, y_1^*)=(-a_r,\dots, -a_1)$, and
$y_k^*=h_k\ex(-\lambda_k)$, $h_k\in K^\times$, so~$y_k^*=f_k^*\ex^{-\phi_k\imag}$, where~$f_k^*:=h_k\ex^{(\phi_k-\phi(\lambda_k))\imag}\in K^\times$, for $k=1,\dots,r$. 
\end{proof}

\begin{cor}\label{cor:a_r-1} 
Suppose $\dim_{\C}V=r\ge 1$ and $\I(K)\subseteq K^\dagger$, and let $f_j$, $\phi_j$ be  as in Corollary~\ref{cor:complex basis}.
Let $A=\der^r+a_{r-1}\der^{r-1}+\cdots + a_0$ \textup{(}$a_0,\dots, a_{r-1}\in K$\textup{)}. Then
$$\phi_1+\cdots+\phi_r \equiv b\bmod\mathcal O_H\qquad\text{for any $b\in H$ with  $b'=-\Im a_{r-1}$,}$$
and hence   $\phi_1+\cdots+\phi_r\preceq 1  \Longleftrightarrow  a_{r-1}\in K^\dagger$.
In particular, if $A^*=(-1)^r A_{\ltimes a}$ \textup{(}$a\in K^\times$\textup{)} or $a_{r-1}\in H$, 
then $\phi_1+\cdots+\phi_r\preceq 1$.
\end{cor}

\begin{proof}
Take $g_j$, $\lambda_j$ as in the proof of Corollary~\ref{cor:complex basis}.
Then $$\lambda_1+\cdots+\lambda_r\ \equiv\  -a_{r-1}\bmod K^\dagger$$ by Corollary~\ref{cor:sum of evs} and Lemma~\ref{lem:complex basis}. Now $K^\dagger\cap H\imag = \I(H)\imag$ by Lemma~\ref{lem:W and I(F)}  and the remarks preceding it.  Also $\phi(\lambda_j)'\imag=\lambda_j$ for all $j$, and this yields the first claim.
For the rest note that if~$A^*=(-1)^r A_{\ltimes a}$ ($a\in K^\times$), then $a_{r-1}\in K^\dagger$ by the remarks after the proof of Proposition~\ref{prop:Bogner}.
\end{proof}

%\noindent
%In a similar way we obtain from Corollary~\ref{cor:self-dual operator} and \eqref{eq:bd on sum mults}:

\begin{cor}\label{cor:phi_i paired} 
Suppose $A$ is self-dual and $\I(K)\subseteq K^\dagger$. Also assume
$K$ is $1$-linearly surjective and $\dim_{\C}V=r$, or $r\ge 1$ and  $K$ is $(r-1)$-linearly surjective.
Then with the $f_j$, $\phi_j$  as in Corollary~\ref{cor:complex basis} we have
$\phi_1+\cdots+\phi_d\preceq 1$, and there is for each $i\in \{1,\dots,d\}$ a $j\in \{1,\dots, d\}$ with~$\phi_i+\phi_j\preceq 1$.
\end{cor} 
\begin{proof} By Corollary~\ref{cor:self-dual operator} and \eqref{eq:bd on sum mults} we have
$\mult_{\lambda}A=\mult_{-\lambda}A$ for all $\lambda$.  With~$\lambda_1,\dots, \lambda_d$ as in the proof of  Corollary~\ref{cor:complex basis} this gives $\lambda_1+ \cdots + \lambda_d=0$,
by Lemma~\ref{lem:complex basis}, and thus $\phi_1+\cdots+\phi_d\preceq 1$. For $i=1,\dots, d$ we have $\mult_{\lambda_i}A=\mult_{-\lambda_i}A>0$ by Lemma~\ref{lem:complex basis}, so that same lemma gives $j\in \{1,\dots,d\}$ such that $\lambda_i+\lambda_j=0$, hence $\phi_i+\phi_j\preceq 1$. 
\end{proof} 

\begin{cor}\label{cor:phi_i paired, strengthened} 
Let $A$, $K$ be as in Corollary~\ref{cor:phi_i paired}.
Then $V$ has a basis
$$f_1\ex^{\phi_1\imag}, g_1\ex^{-\phi_1\imag},\, \dots,\, f_m\ex^{\phi_m\imag},g_m\ex^{-\phi_m\imag}, \ h_1,\dots,h_n \qquad (2m+n=d)$$
where $f_1,\dots,f_m,g_1,\dots,g_m,h_1,\dots,h_n\in K^\times$, and $\phi_1,\dots,\phi_m\in H^{>\R}$ are apart.
% non-in\-ter\-fe\-ring elements of $H$.
\end{cor} 
\begin{proof} 
By the proof of Corollary~\ref{cor:phi_i paired}, if $\lambda$ is an eigenvalue of $A$, then so is $-\lambda$, with the same multiplicity.
Hence Lemma~\ref{lem:complex basis} yields a basis 
$$f_1\ex(\lambda_1), g_1\ex(-\lambda_1),\, \dots,\, f_m\ex(\lambda_m), g_m\ex(-\lambda_m), \ h_1,\dots,h_n \qquad (2m+n = d)$$
of $V$ where  $f_j,g_j,h_k\in K$ for $j=1,\dots,m$, $k=1,\dots,n$ and $\lambda_j\in\Lambda$ with
$\Im\lambda_j>0$ for $j=1,\dots,m$.
Note $\ex(-\lambda)=\ex(\lambda)^{-1}=\ex^{-\phi(\lambda)\imag}$. Setting $\phi_j:=\phi(\lambda_j)$ for $j=1,\dots,m$ thus
yields a basis of $V$ as claimed.
\end{proof}

\noindent
In Section~\ref{sec:group rings} we defined a ``positive definite hermitian form'' on the
 $K$-linear space~$\Univ=K\big[\!\ex(\Lambda)\big]$,   which via our isomorphism $\iota\colon \Univ\to K[\ex^{H\imag}]$ transfers to a ``positive definite hermitian form'' $\langle\ ,\, \rangle$ on 
the $K$-linear space~$K[\ex^{H\imag}]$. 
Note that~$\langle\ ,\, \rangle$ does not depend
on the initial choice of isomorphism $\iota$ as in Lemma~\ref{exqe} at the beginning of this section,
by the remarks following that lemma and Corollary~\ref{cor:inner prod invariant}. 
Suppose
$$y_1=f_1\ex^{\phi_1\imag},\ \dots,\ y_d=f_d\ex^{\phi_d\imag}$$ is a basis of the $\C$-linear space $V$
as in Corollary~\ref{cor:complex basis} such that for $j,k=1,\dots,d$ we have $\phi_j=\phi_k$ or $\phi_j-\phi_k\succ 1$.
Then by Lemma~\ref{lem:inner prod intrinsic} and Corollary~\ref{cor:osc => bded}, 
$$\langle y_j,y_k\rangle\ =\ 0\ \text{ if $\phi_j\neq\phi_k$,} \qquad \langle y_j,y_k\rangle\ =\ f_j\overline{f_k}\neq 0\ \text{ if $\phi_j=\phi_k$.}$$

\subsection*{The case that $A\in H[\der]$}
{\it In this subsection we assume $A\in H[\der]^{\ne}$ has order~$r$.}\/ Then $V:=\ker_{\Univ}A$ is closed under the complex conjugation  
automorphism of the differential ring $\Calinf[\imag]$.  We have $\Univ_{\operatorname{r}}={\Univ}\cap\Calinf$ and
by Corollary~\ref{cor:U_r} a decomposition of $\Univ_{\operatorname{r}}$ as an internal direct sum of
$H$-linear subspaces:
$$\Univ_{\operatorname{r}}\ =\ H\oplus\bigoplus_{\Im\lambda>0}\big( H\cos \phi(\lambda) \oplus H\sin\phi(\lambda) \big).$$
Set $V_{\operatorname{r}}:=V\cap \Calinf$, an $\R$-linear subspace of $V$ with
$V=V_{\operatorname{r}}\oplus V_{\operatorname{r}}\imag$ (internal direct sum of $\R$-linear subspaces of $V$). Each basis of the $\R$-linear space
$V_{\operatorname{r}}$ is a basis of the $\mathbb C$-linear space $V$; in particular, $\dim_{\mathbb C} V = \dim_{\mathbb R} V_{\operatorname{r}}$.
If $\dim_{\mathbb C} V=r$,  then
$V_{\operatorname{r}} = \ker_{\Calinf} A$. 
If~$\lambda$ is an eigenvalue of~$A$, then so is~$-\lambda$, with $\operatorname{mult}_\lambda(A)=\operatorname{mult}_{-\lambda}(A)$.

\begin{lemma}\label{lem:basis of kerU A}
{\samepage The $\mathbb C$-linear space $V=\ker_{\Univ} A$ has a basis 
$$g_1\ex^{\phi_1\imag},\,g_1\ex^{-\phi_1\imag},\ \dots,\ g_m\ex^{\phi_m\imag},\,g_m\ex^{-\phi_m\imag}, \ h_1,\ \dots,\ h_n \qquad (2m+n\leq r),$$ 
where $g_1,\dots,g_m\in H^{>}$, $\phi_1,\dots,\phi_m\in H$ with $\phi_j-\phi(\lambda_j)\preceq 1$ and 
$\Im\lambda_j>0$ for some $\lambda_j\in \Lambda$ for $j=1,\dots,m$, and $h_1,\dots,h_n\in H^\times$.} For any such basis of $V$,
$$g_1\cos \phi_1,\, g_1\sin \phi_1,  \ \dots, \ 
  g_m\cos \phi_m,\, g_m\sin \phi_m, \  h_1,\ \dots,\  h_n$$
is a basis of the $\R$-linear space $V_{\operatorname{r}}$, and $h_1,\dots,h_n$ is a basis of the $\R$-linear sub\-space~$\ker_H A=V\cap H$ of $H$. 
\end{lemma}
\begin{proof} By Corollary~\ref{cor:complex and real basis} the $\C$-linear space $V$ has a basis
$$f_1\ex(\lambda_1),\ \overline{f_1}\ex(-\lambda_1),\ \dots,\ f_m\ex(\lambda_m),\ \overline{f_m}\ex(-\lambda_m),\ h_1,\ \dots,\  h_n$$
with $f_1,\dots,f_m\in K^\times$, $\lambda_1,\dots,\lambda_m\in \Lambda$ with 
$\Im\lambda_1,\dots, \Im\lambda_m>0$ and $h_1,\dots, h_n$
in~$H^\times$. 
Moreover, for each such basis,
$$\Re\!\big(f_1\ex(\lambda_1)\big),\ \Im\!\big(f_1\ex(\lambda_1)\big),\ \dots,\ \Re\!\big(f_m\ex(\lambda_m)\big),\ \Im\!\big(f_m\ex(\lambda_1)\big),\ h_1,\ \dots,\  h_n$$
is a basis of the $\R$-linear space $V_{\operatorname{r}}$, and $h_1,\dots,h_n$ is a basis of its $\R$-linear subspace~$\ker_H A=V\cap H$.
Set $g_j:=\abs{f_j}=\abs{\overline{f_j}}\in H^{>}$ ($j=1,\dots, m$). 
Lemma~\ref{lem:fexphii} gives~${\phi_j\in H}$ such that~$\phi_j-\phi(\lambda_j)\preceq 1$ and $f_j=g_j\ex^{(\phi_j-\phi(\lambda_j))\imag}$, and thus
$f_j\ex(\lambda_j)=g_j\ex^{\phi_j\imag}$, for~${j=1,\dots,m}$. 
Then
$g_1,\dots, g_m, \phi_1,\dots, \phi_m, h_1,\dots, h_n$ have the desired properties. 
\end{proof}

\begin{cor}\label{cor:basis of kerU A}
Suppose $K$ is $1$-linearly surjective when $r\geq 2$,   $\I(K)\subseteq K^\dagger$, and~$A$ splits over $K$.  Then $V=\ker_{\Calinf[\imag]} A$ and
the  $\mathbb C$-linear space $V$ has a basis 
$$g_1\ex^{\phi_1\imag},\,g_1\ex^{-\phi_1\imag},\ \dots,\ g_m\ex^{\phi_m\imag},\,g_m\ex^{-\phi_m\imag}, \ h_1,\ \dots,\ h_n \qquad (2m+n=r),$$ 
where $g_j,\phi_j\in H^{>}$ with $\phi_j\succ 1$ \textup{(}$j=1,\dots,m$\textup{)} and
$h_k\in H^\times$ \textup{(}$k=1,\dots,n$\textup{)}. For any such basis of $V$, 
the $\R$-linear space $\ker_{\Calinf} A$ has basis
$$g_1\cos \phi_1,\, g_1\sin \phi_1,  \ \dots, \ 
  g_m\cos \phi_m,\, g_m\sin \phi_m, \  h_1,\ \dots,\  h_n,$$
and the $\R$-linear subspace $\ker_H A=H\cap\ker_{\Calinf}A$ of $H$ has basis  $h_1,\dots, h_n$.
\end{cor}
\begin{proof}
By Corollary~\ref{cor:basis of kerUA} we have $\dim_{\C} \ker_{\Univ} A=r$, hence $V=\ker_{\Calinf[\imag]} A$ and~$V_{\operatorname{r}}=\ker_{\Calinf} A$.
Now use Lemmas~\ref{lem:from phi to lambda} and~\ref{lem:basis of kerU A}. 
\end{proof}

\noindent
From Lemma~\ref{lem:basis of kerU A} we obtain likewise, using Lemma~\ref{lem:order 2 eigenvalues} 
and Corollary~\ref{cor:sinusoids}:

\begin{cor}\label{cor:2nd order kernel}
Suppose   $r=2$ and $A$ splits over $K$ but not over~$H$. Then there are $g,\phi\in H^>$
such that   
$$\ker_{\Calinf} A\ =\ \R g\cos\phi + \R g\sin\phi\ =\, 
\big\{  cg\cos(\phi+d) :\  c,d\in\R \big\}.$$
Moreover, if $\I(K)\subseteq K^\dagger$, then we can choose here in addition $\phi\succ 1$. 
\end{cor}

\begin{remark}
Let   $A$, $g$, $\phi$, $r$  be as in Corollary~\ref{cor:2nd order kernel}. If $\phi\succ 1$, then all $y\in \ker^{\ne}_{\Calinf} A$
oscillate. If  $\phi\preceq 1$, then no $y\in\ker_{\Calinf} A$ oscillates.
\end{remark}

\noindent
The following generalizes Corollary~\ref{cor:unique antider}:

\begin{cor}\label{cor:unique antider, generalized}
Suppose $K$ is $1$-linearly surjective and $A$ splits over $K$.
Let $\phi$ be an element of $H$ with $\phi>\R$ such that $\phi'\imag+K^\dagger$ is not an eigenvalue of $A$. Then for
every~$h\in H$ there are unique $f,g\in H$ such that $A(f\cos\phi+g\sin\phi)=h\cos \phi$.
\end{cor}
\begin{proof}
Let $f,g\in H$ and $A(f\cos\phi + g\sin \phi)=0$. By Lemma~\ref{lem:trig surj} it is enough to show $f=g=0$.
Set $y:=\frac{1}{2}(f-g\imag)\in K$, so $y\ex^{\phi\imag}+\overline{y}\ex^{-\phi\imag}=f\cos \phi + g\sin \phi$.
The hypothesis and Corollary~\ref{cor:basis of kerUA} give $V=\ker_{\Calinf[\imag]} A$, so~$f\cos \phi + g\sin \phi\in V$. 
Suppose towards a contradiction that $y\neq 0$.
As in the proof of Corollary~\ref{cor:trig surj, complex} we obtain $y\ex^{\phi\imag}=z\ex(\lambda)$, $z\in K^\times$, where $\lambda$ is not an eigenvalue with respect to $\Lambda$. Also~$\lambda\neq 0$ in view of $K^\dagger\subseteq H+\operatorname{I}(H)\imag$. Hence
$$0\ =\ A(y\ex^{\phi\imag}+\overline{y}\ex^{-\phi\imag})\ =\ A\big(z\ex(\lambda)+\overline{z}\ex(-\lambda)\big)\ =\ A_{\lambda}(z)\ex(\lambda) + A_{-\lambda}(\overline{z})\ex(-\lambda),$$ 
so $A_{\lambda}(z)=0$, contradicting that $\lambda$ is not an eigenvalue of  $A$ with respect to $\Lambda$.
\end{proof}

\noindent
Next a version of Lemma~\ref{lem:gaussian ext dom} for $\Univ_{\operatorname{r}}$.
{\it In the rest of this subsection $\I(K)\subseteq K^\dagger$ and $\Lambda=\Lambda_H\imag$ where $\Lambda_H$ is an
$\R$-linear complement of $\I(H)$ in $H$.}\/

\begin{lemma}\label{lem:gaussian ext dom, real}
Let $\lambda_1,\dots,\lambda_n\in\Lambda$  be distinct,
$\Im \lambda_j>0$ for $j=1,\dots,n$,  and
$$y\ =\ f_1\cos \phi_1 + g_1 \sin  \phi_1 +\cdots+f_n\cos\phi_n + g_n\sin\phi_n+h$$
where $f_1,\dots,f_n,g_1,\dots,g_n,h\in H$ and $\phi_j\in \phi(\lambda_j)+\mathcal O_H$ for $j=1,\dots,n$. Then 
$$y\prec 1\quad\Longrightarrow\quad f_1,\dots,f_n,g_1,\dots,f_n,h\prec 1.$$
\end{lemma}
\begin{proof} Let $j$ range over $\{1,\dots,n\}$. 
Setting $a_j:=\frac{1}{2}(f_j-g_j\imag)\in K$ we have
$$y\ =\ a_1 \ex^{\phi_1\imag} + \overline{a_1} \ex^{-\phi_1\imag}+ \cdots + a_n \ex^{\phi_n\imag} + \bar{a_n} \ex^{-\phi_n\imag} +h,$$
and so with $b_j:= a_j\ex^{\phi_j-\phi(\lambda_j)}\in K$ we have $a_j\asymp b_j$ and 
$$y\ =\ b_1 \ex^{\phi(\lambda_1)\imag} + \overline{b_1} \ex^{-\phi(\lambda_1)\imag}+ \cdots + b_n \ex^{\phi(\lambda_n)\imag} + \bar{b_n} \ex^{-\phi(\lambda_n)\imag} +h,$$
Set $h_j:=\phi(\lambda_j)'\in H$. Then $h_j\imag=\lambda_j$, so  the elements
$$h_1,\, \dots,\, h_n,\, -h_1,\, \dots,\, -h_n,\, 0$$ of~$H$ are distinct, and  
$(\R h_1+\cdots+\R h_n)\cap \I(H)=\{0\}$ in view of $\Lambda\cap \I(H)\imag=\{0\}$. Assuming $y\prec 1$,  Corollary~\ref{corlimsupresult} then yields
$b_1,\bar{b_1},\dots, b_n,\bar{b_n}, h\prec 1$, and thus~$f_1,\dots,f_n,g_1,\dots,f_n,h\prec 1$. 
\end{proof}

\begin{lemma} Recalling that $\Univ_r={\Univ}\cap {\Calinf}$ we have: 
\begin{align*}
H &\ =\ \big\{ y\in \Univ_{\operatorname{r}}: \, 
\text{$y-h$ is non-oscillating for all $h\in H$} \big\} \\
&\ =\ \big\{ y\in \Univ_{\operatorname{r}}:\, \text{$y$ lies in a Hausdorff field extension of $H$} \big\}.
\end{align*}
\end{lemma}
\begin{proof} Let $j$ range over $\{1,\dots,n\}$. 
Suppose $y\in \Univ_{\operatorname{r}}$ and $y-h$ is non-oscillating for all $h\in H$.
Take distinct~$\lambda_1,\dots,\lambda_n\in\Lambda$ with $\Im \lambda_1,\dots, \Im \lambda_n>0$, and take~$f_1,\dots,f_n,g_1,\dots,g_n,h\in H$ 
such that 
$$y\ =\ f_1\cos \phi(\lambda_1) + g_1 \sin \phi(\lambda_1)+\cdots+f_n\cos\phi(\lambda_n) + g_n\sin\phi(\lambda_n)+h.$$
We claim that $y=h$.
To prove this claim, replace $y$ by $y-h$ to arrange $h=0$. Towards a contradiction, assume $y\neq 0$.
Then $f_j\neq 0$ or~$g_j\neq 0$ for some $j$.
Divide~$y$ and $f_1,\dots,f_n,g_1,\dots,g_n$ by a suitable element of $H^\times$ to arrange $f_j,g_j\preceq 1$ for all~$j$ and $f_j\asymp 1$ or $g_j\asymp 1$ for some $j$. Then
$y\preceq 1$ and $y-s$ is non-oscillating for all~$s\in \R$, and so
 Lemma~\ref{lem:nonosc, 1} yields $\ell\in\R$ such that $y-\ell\prec 1$.
 Then Lemma~\ref{lem:gaussian ext dom, real} gives~$f_j,g_j\prec 1$ for all $j$, a contradiction.
This proves the first equality. The second equality follows from Lemma~\ref{lem:nonosc, 2}.
\end{proof}

\noindent
In combination with Corollary~\ref{cor:basis of kerU A} this yields: 

\begin{cor}\label{cor:lindiff nonosc} Recalling that $V_r=\ker_{\Univ}A \cap \Calinf$ we have:
\begin{align*}
\ker_H A  &\ =\ \big\{ y\in V_{\operatorname{r}}: \, 
\text{$y-h$ is non-oscillating for all $h\in H$} \big\} \\
&\ =\ \big\{ y\in V_{\operatorname{r}}:\, \text{$y$ lies in a Hausdorff field extension of $H$} \big\}.
\end{align*}
Hence if $K$ is $1$-linearly surjective in case $r\geq 2$, and $A$ splits over $K$, then every~$y$ in~$\ker_{\Calinf}A$ such that $y-h$ is non-oscillating for all $h\in H$ lies in $H$.
\end{cor}

\subsection*{Connection to Lyapunov exponents\astr}
{\it In this subsection $\I(K)\subseteq K^\dagger$, and we take $\Lambda=\Lambda_H\imag$ where $\Lambda_H$ is an $\R$-linear complement
of $\I(H)$ in $H$. Accordingly, $\Univ=K[\ex^{H\imag}]$. Let also $n\geq 1$.}\/ 
In Section~\ref{sec:second-order} we introduced the Lyapunov exponent~$\lambda(f)\in\R_{\pm\infty}$  
of~$f\in\c[\imag]^n$. For use in Section~\ref{sec:lin diff applications}
we collect here some properties of these exponents $\lambda(f)$ for $f\in\Univ^n\subseteq\c[\imag]^n$.
Recall:~$f,g\in \c[\imag],\ f\preceq g\Rightarrow \lambda(f)\geq\lambda(g)$.

\begin{lemma}\label{lem:Lyup gauss} Let $f,g\in \Univ$. Then
$$f\preceq_{\g}g \Rightarrow \lambda(f)\geq\lambda(g), \qquad f\asymp_{\g}g \Rightarrow \lambda(f)=\lambda(g).$$
\end{lemma}
\begin{proof}
We first treat the special case $g=\fm\in H^\times$. Then
the first statement follows from the remark before the lemma and
 Lem\-ma~\ref{lem:gaussian ext dom, preceq}. Suppose $f\asymp_{\g}\fm$; thanks to the first statement it suffices to show $\Lambda(f)\subseteq\Lambda(\fm)$. Towards a contradiction, suppose $\Lambda(f) \not\subseteq \Lambda(\fm)$. Then we have
 $a\in\R$ with $f\preceq \ex^{-ax}$ and $\fm\not\preceq \ex^{-ax}$, so ~$\ex^{-ax}\prec\fm$ (since~$\fm,\ex^{-ax}\in H$), hence
$f\prec\fm$ and thus $f\prec_{\g}\fm$ by Corollary~\ref{cor:gaussian ext dom}, contradicting~$f\asymp_{\g}\fm$. 
%This shows $\lambda(f)\leq\lambda(\fm)$.

The case $g=0$ being trivial, we now assume $g\ne 0$ for the general case and take~$\fm\in H^\times$ 
 with $g\asymp_{\g}\fm$; then~$\lambda(g)=\lambda(\fm)$ by the special case (with $f=g$), so we may replace~$g$ by $\fm$ to reduce the lemma to the special case.
\end{proof}

\noindent
We turn
$\Univ^n$ into a valued $\C$-linear space with   valuation~$v_{\g}\colon \Univ^n\to\Gamma_\infty$ given by 
$$v_{\g}(f):=\min\!\big\{v_{\g}(f_1),\dots,v_{\g}(f_n)\big\}\quad\text{for~$f=(f_1,\dots,f_n)\in \Univ^n$,}$$
and denote by $\preceq_{\g}$ the associated dominance relation on $\Univ^n$.
In the next four corollaries, $f$, $g$ range over $\Univ^n$.

\begin{cor}\label{cor:Lyup gauss, 1}
$f\preceq_{\g} g\Rightarrow\lambda(f)\geq\lambda(g)$ and $f\asymp_{\g} g\Rightarrow\lambda(f)=\lambda(g)$. 
\end{cor}
\begin{proof}
Suppose $f=(f_1,\dots,f_n)\preceq_{\g} g=(g_1,\dots,g_n)$. Take  $k$ with $v_{\g}g=v_{\g}g_k$.
Then $f_j\preceq_{\g}g_k$ and so $\lambda(f_j)\geq\lambda(g_k)\geq\lambda(g)$, for all $j$, by Lemma~\ref{lem:Lyup gauss}, and thus~$\lambda(f)\geq\lambda(g)$.
\end{proof}

\begin{cor}\label{cor:Lyup gauss, 2}
Let $m\ge 1$, $g_1,\dots,g_m\in\Univ^n$,   and  $g=g_1+\cdots+g_m$
be such that $v_{\g}(g)=\min\!\big\{v_{\g}(g_1),\dots,v_{\g}(g_m)\big\}$.
Then $\lambda(g)=\min\!\big\{\lambda(g_1),\dots,\lambda(g_m)\big\}$.
\end{cor}
\begin{proof}
We may arrange $g_i\preceq_{\g} g_1$ for all $i$, so $v_{\g} g=v_{\g} g_1$.
Then $\lambda(g_1)\leq\lambda(g_i)$ for all~$i$ and $\lambda(g)=\lambda(g_1)$, by Corollary~\ref{cor:Lyup gauss, 1}.
\end{proof}

\noindent
Here is a special case of Corollary~\ref{cor:Lyup gauss, 2}:

\begin{cor}\label{cor:Lyap gauss, 3}
Suppose $m\ge 1$, $f=\ex(h_1\imag)f_1+\cdots+\ex(h_m\imag)f_m$ with
$f_1,\dots,f_m$ in $K^n$ and distinct $h_1,\dots,h_m\in\Lambda_H$.
Then $\lambda(f)=\min\!\big\{\lambda(f_1),\dots,\lambda(f_m)\big\}$.
\end{cor}

\noindent
For the notion of valuation-independence, see [ADH, p.~92].
 
\begin{cor}\label{cor:Lyap val indep}
Suppose $m\ge 1$,  $f=\ex^{\phi_1\imag}f_1+\cdots+\ex^{\phi_m\imag}f_m$,
 $f_1,\dots,f_m\in K^n$  
and $\phi_1,\dots,\phi_m\in H$. Suppose also $\phi_j=\phi_k$ or $\phi_j-\phi_k\succ 1$
for~$j,k=1,\dots,m$, and for $k=1,\dots,m$ the $f_j$ with $1\leq j\leq m$ and $\phi_j=\phi_k$ are 
valuation-independent. Then $v_{\g}(f)=\min\{v(f_1),\dots, v(f_m)\}$, and thus $\lambda(f)=\min\!\big\{\lambda(f_1),\dots,\lambda(f_m)\big\}$.
\end{cor}
\begin{proof}
First arrange  that $l\in\{1,\dots,m\}$ is such that $\phi_1,\dots,\phi_l$ are   distinct and each $\phi_j$ with $l<j\leq m$
equals one of $\phi_1,\dots,\phi_l$.
For $k=1,\dots,l$, take
$\lambda_k\in\Lambda$ with~$\phi_k-\phi(\lambda_k)\preceq 1$ and
put~$g_k:=\sum_{1\leq j\leq m,\ \phi_j=\phi_k} f_j$ and
$h_k:=\ex^{(\phi_k-\phi(\lambda_k))\imag}g_k\in K^n$.  Then $v(g_k)=v(h_k)$, $\ex^{\phi_k\imag}g_k=\ex(\lambda_k)h_k$, and
%$$\lambda(h_j)\ =\ \lambda(g_j)\ =\ {\min\!\big\{\lambda(f_k):1\leq k\leq m,\ \phi_k=\phi_j\}}$$
% by Corollary~\ref{cor:Lyup gauss, 2}.
%Since
$$f\ =\ \ex^{\phi_1\imag}g_1+\cdots+\ex^{\phi_l\imag}g_l\ =\ \ex(\lambda_1)h_1+\cdots+\ex(\lambda_l)h_l$$ 
with distinct $\lambda_1,\dots,\lambda_l$. Hence 
$$v_{\g}(f)\ =\ \min\{v(h_1),\dots, v(h_l)\}\ =\ \min\{v(g_1),\dots, v(g_l)\}.$$ 
Now use $v(g_k)=\min \{v(f_j):\ 1\leq j\leq m,\  \phi_j=\phi_k\}$ for $k=1,\dots,l$.
%    $\lambda(f)=\min\!\big\{\lambda(h_1),\dots,\lambda(h_l)\big\}$ by Corollary~\ref{cor:Lyap gauss, 3}.
 \end{proof}

\noindent
For what  we say below about $\Delta$ and $\Gamma^\flat$, see ~[ADH, 9.1.11]. Set  
$$\Delta\ :=\ \big\{\gamma\in\Gamma:\psi(\gamma)\geq 0\big\}\ =\ \big\{\gamma\in\Gamma:\gamma=O\big(v(\ex^x)\big)\big\},$$ the smallest convex subgroup of $\Gamma=v(K^\times)$ containing $v(\ex^x)\in\Gamma^<$. Then $\Delta$ has the convex subgroup $$\Gamma^\flat\ =\ \big\{\gamma\in\Gamma:\psi(\gamma) > 0\big\}\ =\ \big\{\gamma\in\Gamma:\gamma=o\big(v(\ex^x)\big)\big\},$$ and we have an ordered group isomorphism $r\mapsto v(\ex^{-rx})+\Gamma^\flat\colon\R\to\Delta/\Gamma^\flat$.
Note also that for $f\in K$ we have: $v(f)\in \Gamma^\flat\Leftrightarrow \lambda(f)=0$. 

\begin{lemma}\label{lem:fexlambdax}
Let $f\in\Univ$. Then
$$\lambda(f)\ =\ +\infty\ \Leftrightarrow\ v_{\g}(f)\ >\ \Delta,\qquad
\lambda(f)\ =\ -\infty\ \Leftrightarrow\ v_{\g}(f)\ <\ \Delta,$$
and if~$\lambda(f)\in\R$, then $v_{\g}(f)\in \Delta$ and 
$v_{\g}(f)\equiv v(\ex^{-\lambda(f)x})\bmod\Gamma^\flat$.
\end{lemma}
\begin{proof}
We assume $f\neq 0$, and use Lemma~\ref{lem:Lyup gauss} to replace $f$ by $\fm\in H$ with $f\asymp_{\g}\fm$ so as to arrange $f\in H^\times$. The displayed claims then follow. Suppose~$\lambda(f)\in\R$, and
let $a\in\R^>$. Then 
$f\ex^{(\lambda(f)-\frac{1}{2}a) x}\preceq 1$, so $f\ex^{\lambda(f)x}\prec\ex^{\frac{1}{2}a x}\prec\ex^{ax}$. Also~$f\ex^{\lambda(f) x}\not\preceq\ex^{-a x}$, thus
%$f\ex^{\lambda(f) x}\not\preceq_{\g}\ex^{-a x}$ by Lemma~\ref{lem:gaussian ext dom, preceq} and so
$\ex^{-a x}\prec f\ex^{\lambda(f) x}\prec\ex^{a x}$.
% where we used  Corollary~\ref{cor:gaussian ext dom} to get the second~$\prec_{\g}$,
This holds for all $a\in \R^{>}$, so~${v(f\ex^{\lambda(f) x})\in\Gamma^\flat}$.
\end{proof}

\noindent
 Lemma~\ref{lem:fexlambdax} yields $\lambda(fg)=\lambda(f)+\lambda(g)$ for all  $f,g\in\Univ\cap \c[\imag]^{\flattereq}$. 

\begin{cor}\label{fglambda}
Assume $f,g\in K$, $g\preceq f$, and $\lambda(f)\in\R$. Then~$g'+\lambda(f)g\prec f$.
\end{cor}
\begin{proof} Lemma~\ref{lem:fexlambdax} gives $v(f\ex^{\lambda(f)x})\in \Gamma^\flat$, so we can
replace $f$, $g$ by $f\ex^{\lambda(f) x}$, $g\ex^{\lambda(f) x}$, to arrange $f'\prec f$ and $\lambda(f)=0$.
Now if $f\asymp 1$, then $g\preceq f\asymp 1$ and so~$g'\prec 1\asymp f$, and
if $f\nasymp 1$, then~$g'\preceq f'\prec f$ using [ADH, 9.1.3(iii) and 9.1.4(i)].
\end{proof}

\noindent
From Corollary~\ref{fglambda} we easily obtain: 

\begin{cor}\label{cor:fexlambdax}
Suppose $f\in K^n$  is such that $\lambda(f)\in\R$. 
Then $f'+\lambda(f)f\prec  f$.
\end{cor}

\noindent
Note that $K\cap\c[\imag]^{\flattereq}=\mathcal O_\Delta$ is by Lemma~\ref{lem:fexlambdax} the valuation ring of the coarsening~$v_\Delta$ of the valuation of $K$ by $\Delta$, with maximal ideal $K\cap\c[\imag]^{\flatter}=\smallo_\Delta$, cf.~[ADH, 3.4].  
By Corollary~\ref{cor:Lyap gauss, 3}, the  $\C$-subalgebra $\Univ\cap\c[\imag]^{\flattereq}$ of $\Univ$ satisfies
$${\Univ}\cap\c[\imag]^{\flattereq}  = \bigoplus_{h\in\Lambda_H} \mathcal O_\Delta\ex(h\imag)\quad\text{(internal direct sum of $\mathcal O_\Delta$-submodules of $\Univ\cap\c[\imag]^{\flattereq}$).}$$
We put
$$\Univ^{\flattereq}\  :=\   \bigoplus_{h\in\Lambda_H\cap \mathcal O_H} \mathcal O_\Delta\ex(h\imag),$$
a $\C$-subalgebra of $\Univ\cap\c[\imag]^{\flattereq}$. Then
\begin{align*}
(\Univ^{\flattereq})^\times &= \big\{ g\ex(h\imag): g\in K^\times,\,  h\in\Lambda_H,\ g^\dagger,h\preceq 1\big\} \\
&= \big\{ g\ex(h\imag): g\in K^\times,\, h\in\Lambda_H,\ \lambda(g)\in\R,\, h\preceq 1\big\}.
\end{align*}
In the next lemma $f=g\ex^{\phi\imag}\in\Univ^\times$
where $g\in K^\times$, $\phi\in H$. Then~$\abs{f}=\abs{g}\in H$ and so~$-\lambda(f)= -\lambda(g)=\lim\limits_{t\to\infty} \frac{1}{t}\log\abs{g(t)}$.

\begin{lemma}\label{lem:lambda units}
 $f\in (\Univ^{\flattereq})^\times\ \Leftrightarrow\ g^\dagger,\phi'\preceq 1\ \Leftrightarrow\ f^\dagger \preceq 1$.
If  $f^\dagger \preceq 1$, then
$$-\lambda(f)\ =\  \lim_{t\to\infty} \Re f^\dagger(t)\ =\  \lim_{t\to\infty} \Re g^\dagger(t),\qquad
 \lim_{t\to\infty} \Im f^\dagger(t)\  =\  \lim_{t\to\infty} \phi'(t)$$
 and these limits are in $\R$. 
\end{lemma}
\begin{proof}
Take $h\in\Lambda_H$ with $\phi-\phi(h\imag)\preceq 1$ and put $g_1:=g\ex^{(\phi-\phi(h\imag))\imag}\in K^\times$, so~$f=g_1\ex^{\phi(h\imag)\imag}=g_1\ex(h\imag)$. 
Now $g^\dagger-\abs{g}^\dagger\prec 1$, since
$g\asymp\abs{g}$. Also $\ex(h\imag)=\ex^{\phi(h\imag)\imag}$ gives~$h=\phi(h\imag)'$ by differentiation. Hence
$g_1^\dagger-g^\dagger=(\phi'-h)\imag=(\phi-\phi(h\imag))'\imag \prec 1$, so 
$$\abs{g}^\dagger \preceq 1\ \Leftrightarrow\ g^\dagger\preceq 1\ \Leftrightarrow\ g_1^\dagger\preceq 1,\qquad
\phi'\preceq 1\ \Leftrightarrow\ h\preceq 1.$$ 
This yields the equivalences of the Lemma, using for
$f^\dagger\preceq 1\Rightarrow g^\dagger, \phi'\preceq 1$ that  $f^\dagger=g^\dagger + \phi'\imag$, and
$\Im(g^\dagger)\in \I(H)\imag\subseteq K^{\prec  1}$, the latter a consequence of Lemma~\ref{lem:W and I(F)}  and the remarks preceding it. Now assume $f^\dagger\preceq 1$. Then $g^\dagger, \phi'\preceq 1$, so $vg\in \Delta$, hence by Lemma~\ref{lem:fexlambdax}, $\lambda(f)=\lambda(g)\in \R$ and $v\big(g\ex^{\lambda(g)x}\big)\in \Gamma^\flat$, that is,
$g^\dagger + \lambda(g)\prec 1$, so~$\Re(g^\dagger)+\lambda(f)\prec 1$, and thus
$-\lambda(f)=\ \lim\limits_{t\to\infty} \Re g^\dagger(t)$. Now use $\Re f^\dagger=\Re g^\dagger$ and~$\Im f^\dagger=\Im g^\dagger + \phi'$ and $\Im g^\dagger \prec 1$. 
\end{proof} 

\begin{lemma}\label{lem:Univflattereq}
Let $f\in \Univ^{\flattereq}$. Then $f'\in \Univ^{\flattereq}$ and $\lambda(f)\leq\lambda(f')$. Moreover, if~$\lambda(f)\in\R$, then
 $f'\preceq_{\g} f$, and if $\lambda(f)\in\R^\times$, then $f'\asymp_{\g} f$.
\end{lemma}
\begin{proof}
%By Corollary~\ref{cor:Lyap gauss, 3} 
Suppose first that $f=g\ex(h\imag)$ where $g\in \mathcal O_\Delta^{\neq}$, $h\in\Lambda_H\cap\mathcal O_H$, so $\lambda(f)=\lambda(g)$ and
  $f'=(g'+gh\imag)\ex(h\imag)$.
Then by [ADH, 9.2.24, 9.2.26] we have~$g'\in\mathcal O_\Delta$, with $g'\in\smallo_\Delta$ if $g\in\smallo_\Delta$.
So $f'\in\mathcal O_\Delta\ex(h\imag)$, with $f'\in\smallo_\Delta\ex(h\imag)$ if $g\in\smallo_\Delta$.
This yields~$f'\in \Univ^{\flattereq}$ as well as $\lambda(f')=+\infty$ if~$\lambda(f)=+\infty$, by Lemma~\ref{lem:fexlambdax}.
Now suppose
$\lambda(f)\in\R$. Then~$v(g\ex^{\lambda(g)x})\in\Gamma^\flat$ by Lemma~\ref{lem:fexlambdax}, hence $g^\dagger+\lambda(g)\prec 1$, so $g^\dagger\preceq 1$, and 
thus~$f'={g(g^\dagger+h\imag)\ex(h\imag)}\preceq_{\g}f$, and this yields $\lambda(f')\geq\lambda(f)$ by Lemma~\ref{lem:Lyup gauss}.
If~$\lambda(f)\neq 0$, then~$g^\dagger\sim-\lambda(g)$ and so $g^\dagger+h\imag\sim -\lambda(g)+h\imag\asymp 1$, and thus
$f'\asymp_{\g}f$.

The case $f=0$ is trivial, so we can assume next that $f=f_1+\cdots+f_m$,  $f_j=g_j\ex(h_j\imag)$, $g_j\in\mathcal O_\Delta^{\neq}$, $h_j\in\Lambda_H\cap\mathcal O_H$ for $j=1,\dots,m$, $m\ge 1$, with dis\-tinct~$h_1,\dots,h_m$. 
We arrange $f_1\succeq_{\g}\cdots\succeq_{\g}f_m$, so $f\asymp_{\g} f_1$ and $\lambda(f_1)\leq\cdots\leq\lambda(f_m)$, and
$\lambda(f)=\min\!\big\{\lambda(f_1),\dots,\lambda(f_m)\big\}=\lambda(f_1)$ by Corollary~\ref{cor:Lyap gauss, 3}. 
The special case gives $f_j'\in \Univ^{\flattereq}$ and $\lambda(f_j)\leq\lambda(f'_j)$ for $j=1,\dots,m$, so
$f'\in \Univ^{\flattereq}$ and $\lambda(f)\leq\lambda(f')$. Suppose $\lambda(f)\in\R$. Then $v_{\g}(f_1)\in\Delta$ by Lemma~\ref{lem:fexlambdax}.
If~$\lambda(f_j)=+\infty$, then~$\lambda(f_j')=+\infty$ by the special case, so $v_{\g}(f_j')>\Delta$ and thus~$f_j' \prec_{\g} f_1\asymp_{\g} f$.  
If~$\lambda(f_j)\in\R$, then $f_j'\preceq_{\g}f_j\preceq_{\g} f$, again by the special case. This yields~$f'\preceq_{\g}f$.
  Likewise, if~$\lambda(f)\in\R^\times$, then $f_1'\asymp_{\g} f_1\asymp_{\g}f$ and thus $f'\asymp_{\g} f$.
\end{proof}

\begin{cor}\label{cor:Univflattereq}
If $f\in\Univ^{\flattereq}$ and $\lambda(f)\in\R^\times$, then for all $n$, 
$$f^{(n)}\ \asymp_{\g}\ f, \qquad \lambda(f^{(n)})\ =\ \lambda(f).$$
\end{cor} 

\noindent
For use in the next lemma and then in Section~\ref{sec:lin diff applications}  we also define for $f\in \c^1[\imag]^\times$,   
$$\mu(f)\ :=\  \limsup\limits_{t\to\infty} \Im\big(f'(t)/f(t)\big)\in\R_{\pm\infty}.$$
If $f\in (\Univ^{\flattereq})^\times$, then $\lambda(f),\mu(f)\in\R$ by Lemma~\ref{lem:lambda units}, and $f^\dagger= \Re(f^\dagger) + \Im(f^\dagger)\imag$ then yields $f^\dagger - \big({-\lambda(f)+\mu(f)\imag}\big)\prec 1$.

\medskip\noindent
In the next lemma, suppose
 $f_1,\dots,f_n \in (\Univ^{\flattereq})^\times$ are such that  
$$c_1\ :=\ -\lambda(f_1)+\mu(f_1)\imag,\ \dots,\ c_n\ :=\ -\lambda(f_n)+\mu(f_n)\imag \in\C$$ are distinct.
Also, let $c\in\C$ and suppose
$f:=f_1+\cdots+f_n\in  \c[\imag]^\times$ and~${c-f^\dagger\prec 1}$.

\begin{lemma}\label{smallestrealpart}
Let $i\in\{1,\dots,n\}$  be such that
$f_i\succeq f_k$ for all $k\in\{1,\dots,n\}$.
Then~$c_i=c$  and~$\Re c_k \le \Re c$ for all $k$.
\end{lemma}
\begin{proof}
We let $j$, $k$, $l$ range over $\{1,\dots,n\}$.
Take $g_k\in \mathcal O_\Delta^{\neq}$ and $h_k\in\Lambda_H\cap\mathcal O_H$  such that $f_k=g_k\ex(h_k\imag)$.
Then~$f_k^\dagger=g_k^\dagger+h_k\imag\in \mathcal O$, $c_k-f_k^\dagger\prec 1$,   and
$$f'\ =\ f_1^\dagger g_1 \ex(h_1\imag)+\cdots+ f_n^\dagger g_n \ex(h_n\imag). $$
Suppose $h_j=h_k$ and $g_j\asymp g_k$; then $f_j^\dagger-f_k^\dagger=(g_j/g_k)^\dagger\in\I(K)\subseteq\smallo$ and so $c_j-c_k\prec 1$, hence $j=k$.
We arrange $l\geq i$ so that $h_1,\dots,h_l$ are distinct and the $h_k$ with $k>l$ are in $\{h_1,\dots,h_l\}$.
For $j\leq l$, set
$g_j^* := \sum_{h_k=h_j} g_k$ and~$g_j^{\der} := \sum_{h_k=h_j} f_k^\dagger g_k$, so
$$f\ =\ \sum_{j\leq l} g_j^*\ex(h_j\imag),\qquad f'\ =\ \sum_{j\leq l} g_j^{\der}\ex(h_j\imag).$$
For $j\le  l$ we have a unique $k=k(j)$ with $g_j^*\sim g_k$. Now $g_i\asymp f_i\succeq f_k\asymp g_k$ for all~$k$, so $i=k(i)$, hence
$0\neq g_i^*\sim g_i \succeq g_{k(j)}\asymp g_j^*$  for $j\leq l$. In particular, $f\asymp_{\g} g_i$.
 
Suppose $c\neq 0$. Then $c-f^\dagger\prec 1$ gives $cf\sim f'$. Hence by Lemma~\ref{lem:exp sum asymptotics}
we have~$cg_i^* \sim  g_i^{\der}$ and $\sum_{h_k=h_j}(c-f_k^\dagger)g_k=cg_j^*-g_j^{\der}\prec cg_i^*$
for~$j\neq i$, $j\le l$.  Then~$cg_i\sim g_i^{\der}=\sum_{h_k=h_i} f_k^\dagger g_k$. But if $k\ne i$ and
$h_k=h_i$, then $f_k^\dagger g_k\preceq g_k \prec g_i$, hence~$cg_i\sim f_i^\dagger g_i$, so $c\sim f_i^\dagger$. This proves $c=c_i$. 
Also, if $k\neq i$ and $h_k=h_i$, then  $g_k\prec g_i$, so~$\Re(f_k^\dagger)=\Re(g_k^\dagger)<\Re(g_i^\dagger)=\Re(f_i^\dagger)$ by Corollary~\ref{cor:10.5.2 variant},  hence~$\Re c_k\le \Re c_i=\Re c$.
% which by $\mu(f_k)= \mu(f_j)$ and $c_k\ne c_j$ yields $\Re c_k < \Re c$.
 If $j\leq l$, $j\neq i$ and $h_k=h_j$, then $c\neq c_k$ gives~$g_k\asymp (c-f_k^\dagger)g_k \preceq   cg_j^*-g_j^{\der} \prec cg_i^*\asymp g_i$, and 
as before this yields $\Re c_k \le \Re c$. 
Hence $\Re c_k \le \Re c$ for all $k$.

Next suppose~$c=0$. Then $f'\prec f$ and so $f'\prec_{\g} f \asymp_{\g} g_i$ by Corollary~\ref{lem:gaussian ext dom, general}, hence
$g_j^{\der}\prec g_i$ for $j\leq l$, and this yields
$f_k^\dagger g_k \prec g_i$ for all $k$.
Taking~${k=i}$ now gives~$f_i^\dagger\prec 1$ and so~$c_i=0$, and 
if~$k\neq i$, then $c_k\neq 0$ and  thus $f_k^\dagger\asymp 1$, so~${g_k\asymp f_k^\dagger g_k \prec g_i}$, 
and as in the case $c\neq 0$ this gives $\Re c_k \le \Re c_i =0$.
\end{proof}

\newpage 

\part{Filling Holes in Hardy Fields}\label{part:Hardy fields}

\medskip

\noindent
This part contains in Section~\ref{sec:d-alg extensions} the proof of our main theorem. 
Important tools for this are the normalization and approximation theorems for holes and slots established in Parts~\ref{part:normalization} and~\ref{part:dents in H-fields}.
On the analytic side we need a suitable fixed point theorem proved in Section~\ref{sec:split-normal over Hardy fields}:
Theorem~\ref{thm:fix}. The definition of the operator used there is based on the right-inverses for  linear differential operators
over  Hardy fields constructed in Section~\ref{sec:IHF}.
Section~\ref{sec:smoothness} complements Section~\ref{sec:split-normal over Hardy fields}
 by showing how to recover suitable smoothness for the fixed points obtained this way.

Let $(P,\fm,\hat f)$ be a hole in a Liouville closed Hardy field $H\supseteq \R$ and recall that $\hat f$ lies in an immediate $H$-field extension of $H$ and satisfies
$P(\hat f)=0$, $\hat f\prec\fm$. (This extension is not assumed to be a Hardy field.) Under suitable hypotheses on 
$H$ and $(P,\fm,\hat f)$, our fixed point
  theorem (or rather its ``real'' variant, Corollary~\ref{cor:fix}) produces a  germ $f$ of a one-variable real-valued function such that $P(f)=0$,~$f\prec\fm$; see Section~\ref{secfhhf}.  The  challenge in the proof of our main result is   to show that such an~$f$ generates   a Hardy field extension $H\langle f\rangle$  of $H$ isomorphic to $H\langle\hat f\rangle$ over $H$ (as ordered differential fields).
In particular, we need to demonstrate that this zero $f$ of~$P$ has the same asymptotic properties (relative to $H$) as its formal counterpart~$\hat f$, and
the notion of {\it asymptotic similarity}\/ established in Section~\ref{sec:asymptotic similarity} provides a suitable general framework for doing so.
In order to show that $f$ is indeed asymptotically similar to~$\hat f$ over $H$, we are naturally led to the following task:
{\it given another germ $g$ satisfying $P(g)=0$, $g\prec\fm$,
bound the growth of $h,h',\dots,h^{(r)}$ where  $h:=(f-g)/\fm$ and $r:=\order P$.}\/
Assuming (among other things) that $(P,\fm,\hat f)$ is repulsive-normal in the sense
of Part~\ref{part:dents in H-fields}, this is  accomplished in Section~\ref{sec:weights}, after 
revisiting parts of the material from Sections~\ref{sec:IHF}, \ref{sec:split-normal over Hardy fields}, and~\ref{secfhhf} for certain weighted function spaces. (See Proposition~\ref{prop:notorious 3.6}.)

\section{Inverting Linear Differential Operators over Hardy Fields}\label{sec:IHF}

\noindent
Given a Hardy field $H$ and $A\in H[\der]$ we shall construe $A$ as
a $\C$-linear operator on various spaces of functions. 
%the differential ring $\c^{<\infty}[\imag]$.
%(This ring is defined in the first subsection below.) 
We wish to construct right-inverses to  such operators.
A key assumption here is that $A$ splits over
$H[\imag]$. This reduces the
construction of such inverses mainly to the case of order $1$, and this case is handled in the first two subsections using suitable twisted integration operators.  In the third subsection we put things together and also show how to ``preserve reality'' by taking real parts. In the fourth subsection we introduce damping factors.  Throughout we pay attention to the continuity of various operators with respect to various norms, for use in Section~\ref{sec:split-normal over Hardy fields}.    

\medskip
\noindent
We let $a$ range over $\R$ and $r$ over $\N\cup\{\infty,\omega\}$.  If $r\in \N$, then $r-1$ and $r+1$ have the usual meaning, while for $r\in \{ \infty,\omega\}$ we set $r-1= r+1:=r$. (This convention is just to avoid case distinctions.) We have the usual absolute value on
$\C$ given by~$|a+b\imag|=\sqrt{a^2+b^2}\in \R^{\ge}$ for $a,b\in \R$, so for
$f\in \c_a[\imag]$ we have $|f|\in \c_a$.

\subsection*{Integration and some useful norms}  
For $f\in \c_a[\imag]$ we define $\der_a^{-1}f \in \Cao[\imag]$ by
$$\der_a^{-1}f(t)\ :=\ \int_a^t f(s)\,ds\ :=\ \int_a^t \Re f(s)\, ds+ \imag\int_a^t \Im f(s)\, ds,$$
so $\der_a^{-1}f$  is the unique $g\in \Cao[\imag]$ such that $g'=f$ and $g(a)=0$.  
The integration operator $\der_a^{-1}\colon \c_a[\imag]\to  \Cao[\imag]$
is $\C$-linear and maps  $\Car[\imag]$ into $\Carm[\imag]$. 
For $f\in \c_a[\imag]$ we have
$$%\begin{equation}\label{eq:triangle ineq for int}
\big|\der_a^{-1}f(t)\big|\le \big(\der_a^{-1}|f|\big)(t)\qquad\text{ for all $t\ge a$.}
$$%\end{equation}
Let $f\in \c_a[\imag]$.  Call $f$ {\bf integrable at $\infty$\/}\index{function!integrable at $\infty$} if
$\lim_{t\to \infty} \int_a^t f(s)\,ds$ exists in $\C$. In that case we denote this
limit by $\int_a^{\infty} f(s)\,ds$ and put
$$ \int_{\infty}^a f(s)\,ds\ :=\  -\int_a^\infty f(s)\,ds,$$
and define $\der_{\infty}^{-1}f\in \Cao[\imag]$  by
$$\der_{\infty}^{-1}f(t)\ :=\ \int_{\infty}^t f(s)\,ds\ =\ \int_{\infty}^a f(s)\, ds +  \int_a^t  f(s)\, ds\ =\ \int_{\infty}^a f(s)\, ds +\der_a^{-1}f(t),$$
so  $\der_{\infty}^{-1}f$ is the unique $g\in \Cao[\imag]$  such that $g'=f$ and $\lim_{t\to \infty} g(t)=0$. Note that \label{p:Caint}
\begin{equation}\label{eq:integrable}
\c_a[\imag]^{\inte}\ :=\ \big\{f\in \c_a[\imag]:\ \text{$f$   is integrable at  $\infty$}\big\}
\end{equation}
is a $\C$-linear subspace of $\c_a[\imag]$ and that $\der_{\infty}^{-1}$ defines a $\C$-linear
operator from this subspace into $\Cao[\imag]$ which maps $\Car[\imag]\cap \c_a[\imag]^{\inte}$ into
$\Carm[\imag]$. If $f\in\c_a[\imag]$ and 
$g\in \c_a^{\inte}:=\c_a[\imag]^{\inte}\cap \c_a$ with $\abs{f}\leq g$ as germs in $\c$, then
$f\in \c_a[\imag]^{\inte}$; in particular,
if~$f\in\c_a[\imag]$ and $\abs{f}\in\c_a^{\inte}$, then 
$f\in \c_a[\imag]^{\inte}$. 
Moreover:

\begin{lemma}\label{lem:derinvinf} 
Let $f\in\c_a[\imag]$ and $g\in \c_a^{\inte}$ be such that $\abs{f(t)} \leq g(t)$ for all $t\geq a$. Then
$\abs{\der^{-1}_\infty f(t)} \leq \abs{\der^{-1}_\infty g(t)}$ for all $t\geq a$.
\end{lemma}
\begin{proof}
Let $t\geq a$.
We have $g\geq 0$ on $[a,\infty)$,
hence % $\der_a^{-1} g(t) = \int_a^t g(s)\,ds \leq \int_a^\infty g(s)\,ds$ and so 
$\der^{-1}_\infty g(t) \leq 0$.
Also $\left| \int_t^\infty f(s)\,ds \right| \leq \int_t^\infty \abs{f(s)}\,ds \leq \int_t^\infty g(s)\,ds$.
Thus
$$\abs{\der_\infty^{-1}f(t)}\ =\ \left|\int_t^\infty f(s)\,ds\right| 
\ \leq\ \int_t^\infty g(s)\,ds\ =\ -\der_\infty^{-1} g(t) 
%\int_a^\infty g(s)\,ds - \der_a^{-1}g(t)
\ =\ \abs{\der_\infty^{-1}g(t)} $$
as claimed.
\end{proof}

\noindent
For $f\in \c_a[\imag]$ we set
$$ \|f\|_a\ :=\ \sup_{t\ge a} |f(t)|\ \in\ [0,\infty],$$
so (with $\b$ for ``bounded''): \label{p:Cab}
$$\c_a[\imag]^{\b}\ :=\ \big\{f\in \c_a[\imag]:\ \|f\|_a<\infty\big\}$$ is a $\C$-linear subspace of 
$\c_a[\imag]$, and $f\mapsto \|f\|_a$ is a norm on $\c_a[\imag]^{\b}$ making it a Banach space over $\C$.
It is also convenient to define for $t\ge a$ the seminorm \label{p:absa}
$$\|f\|_{[a,t]}\ :=\ \max_{a\le s\le t} |f(s)|$$ 
on $\c_a[\imag]$. More generally, let $r\in \N$. Then for $f\in \Car[\imag]$ we set
$$ \|f\|_{a;r}\ :=\ \max\big\{\|f\|_a, \dots, \|f^{(r)}\|_a\big\}\ \in\ [0,\infty],$$
so 
$$\Car[\imag]^{\b}:=\big\{f\in \Car[\imag]:\ \|f\|_{a;r}<\infty\big\}$$ 
is a $\C$-linear subspace of 
$\Car[\imag]$, and $f\mapsto \|f\|_{a;r}$ makes $\Car[\imag]^{\b}$ a normed vector space over $\C$.
Note that by Corollary~\ref{cor:HL, bded},
$$\c^r_a[\imag]^{\b}\ =\  \big\{f\in \c^r_a[\imag]:\, \dabs{f}_{a}<\infty \text{ and } \dabs{f^{(r)}}_{a}<\infty\big\},$$
although we do not use this later. 
Note that for $f,g\in \Car[\imag]$ we have $$\|fg\|_{a;r}\ \le\ 2^r\|f\|_{a;r}\|g\|_{a;r},$$
so $\Car[\imag]^{\b}$ is a subalgebra of the $\C$-algebra $\Car[\imag]$. 
If $f\in \Carm[\imag]$, then $f'\in \Car[\imag]$ with~$\|f'\|_{a;r}\le \|f\|_{a;r+1}$. 

\medskip
\noindent
With $\i=(i_0,\dots,i_r)$ ranging over $\N^{1+r}$,  let
$P=\sum_{\i} P_{\i} Y^{\i}$ (all $P_{\i}\in \c_a[\imag]$) be a polynomial in~$\c_a[\imag]\big[Y,Y',\dots,Y^{(r)}\big]$.
For  $f\in \Car[\imag]$ we   set 
$$P(f)\ :=\ \sum_{\i} P_{\i} f^{\i}\in\c_a[\imag]\qquad\text{where $f^{\i}:=f^{i_0}(f')^{i_1}\cdots (f^{(r)})^{i_r}\in \c_a[\imag]$.}$$
We also let
$$\dabs{P}_a\ :=\ \max_{\i} \, \dabs{P_{\i}}_a \in [0,\infty].$$
Then  $\dabs{P}_a<\infty$ iff $P\in\c_a[\imag]^{\b}\big[Y,\dots,Y^{(r)}\big]$,
and $\dabs{\,\cdot\,}_a$ is a norm on the $\C$-linear space $\c_a[\imag]^{\b}\big[Y,\dots,Y^{(r)}\big]$.
In the following     assume $\dabs{P}_a<\infty$.
Then for $j=0,\dots,r$ such that $\partial P/\partial Y^{(j)}\neq 0$ we have
$$\dabs{\partial P/\partial Y^{(j)}}_a\  \leq\  (\deg_{Y^{(j)}}P)\cdot\dabs{P}_a.$$
Moreover:

\begin{lemma}\label{lem:bound on P(f)} 
If $P$ is homogeneous of degree $d\in\N$ and $f\in\Car[\imag]^{\b}$, then
$$\dabs{P(f)}_a\  \leq\  {d+r\choose r} \cdot \dabs{P}_a\cdot\dabs{f}_{a;r}^d.$$
\end{lemma}

\begin{cor}\label{cor:bound on P(f)}
Let $d\leq e$ in $\N$ be such that $P_{\i}=0$ whenever $\abs{\i}<d$ or $\abs{\i}>e$. Then
for~$f\in\Car[\imag]^{\b}$ we have
$$\dabs{P(f)}_a\ \leq\  D \cdot \dabs{P}_a\cdot \big(\dabs{f}_{a;r}^d+\cdots+\dabs{f}_{a;r}^e\big)$$
where $D=D(d,e,r):={e+r+1\choose r+1}-{d+r\choose r+1}\in\N^{\geq 1}$.  
\end{cor}

\noindent
Let $B\colon V \to \Car[\imag]^{\b}$ be a $\C$-linear map from a normed vector space $V$ over $\C$ into~$\Car[\imag]^{\b}$. Then we set 
$$\|B\|_{a;r}\ :=\  \sup\big\{\|B(f)\|_{a;r}:\ f\in V,\ \|f\| \le 1\big\}\ \in\ [0,\infty],$$
the {\bf operator norm of $B$}. 
Hence with the convention $\infty\cdot b:=b\cdot\infty:=\infty$ for~$b\in [0,\infty]$ we have 
$$\|B(f)\|_{a;r}\ \leq\ \|B\|_{a;r}\cdot \|f\|\qquad\text{for $f\in V$.}$$
Note that $B$ is continuous iff $\|B\|_{a;r}<\infty$. 
If the map $D\colon \Car[\imag]^{\b}\to \mathcal C_{a}^{s}[\imag]^{\b}$ ($s\in\N$) is also $\C$-linear, then
$$\| D\circ B \|_{a;s}\ \leq\ \| D \|_{a;s} \cdot \|B\|_{a;r}.$$
For $r=0$ we drop the subscript: $\|B\|_{a}:=\|B\|_{a;0}$.

\begin{lemma}\label{dermphi} Let $r\in \N^{\ge 1}$ and $\phi\in \Carl[\imag]^{\b}$. Then the $\C$-linear operator $$\der-\phi\ :\ \Car[\imag]\to \Carl[\imag], \quad f\mapsto f'-\phi f$$ maps
$\Car[\imag]^{\b}$ into $\Carl[\imag]^{\b}$, and its restriction
$\der-\phi\colon \Car[\imag]^{\b}\to \Carl[\imag]^{\b}$ is continuous with operator norm $\|\der-\phi\|_{a;r-1}\le 1+2^{r-1}\|\phi\|_{a;r-1}$.
\end{lemma}

\noindent
Let $r\in \N$, $a_0\in\R$, and let $a$ range over $[a_0,\infty)$. The 
$\C$-linear map  
$$f\mapsto  f|_{[a,+\infty)}\ \colon\ \mathcal C_{a_0}^r[\imag] \to \mathcal C_{a}^r[\imag]$$
satisfies $\|f|_{[a,+\infty)}\|_{a;r}\leq \|f\|_{a_0;r}$ for $f\in \mathcal C_{a_0}^r[\imag]$, so it maps
$\mathcal C_{a_0}^r[\imag]^{\b}$ into~$\Car[\imag]^{\b}$.
For~${f\in\mathcal C_{a_0}^0[\imag]}$ also denoting its germ at $+\infty$ and its restriction  $f|_{[a,+\infty)}$, we have: 
\begin{align*}
f\preceq 1 &\quad\Longleftrightarrow\quad \|f\|_a < \infty \ \text{for some $a$}  
\quad\Longleftrightarrow\quad  \|f\|_a < \infty \ \text{for all $a$,} \\
f\prec 1 & \quad\Longleftrightarrow\quad \|f\|_a\to 0\ \text{as $a\to\infty$.}
\end{align*}

\subsection*{Twisted integration} 
For  $f\in \c_a[\imag]$ we have the $\C$-linear operator $$g\mapsto fg\ :\ \c_a[\imag]\to \c_a[\imag],$$ which we also denote by $f$. 
We now
fix an element $\phi\in \c_a[\imag]$, and set $\Phi:= \der_a^{-1} \phi$, so
$\Phi\in \Cao[\imag]$, $\Phi(t)=\int_a^t \phi(s)\,ds$ for $t\ge a$, and $\Phi'=\phi$. Thus $\ex^\Phi, \ex^{-\Phi}\in \Cao[\imag]$ with~$(\ex^\Phi)^\dagger = \phi$. Consider the $\C$-linear operator
$$B := \ex^{\Phi}\circ\ \der_a^{-1}\circ \ex^{-\Phi}\ \colon\ \c_a[\imag]\to \Cao[\imag],$$ so
$$ Bf(t)\ =\  \ex^{\Phi(t)}\int_a^t \ex^{-\Phi(s)}f(s)\, ds\quad \text{ for $f\in \c_a[\imag]$.}$$
It is easy to check that $B$ is a right inverse to $\der - \phi\colon\Cao[\imag]\to \c_a[\imag]$ in the sense that~$(\der-\phi)\circ B$ is the identity on~$\c_a[\imag]$. Note that for
$f\in \c_a[\imag]$ we have $Bf(a)=0$, and thus $(Bf)'(a)=f(a)$, using $(Bf)'=f+\phi B(f)$.  
Set $R:= \Re \Phi$ and $S:= \Im \Phi$, so $R, S\in \Cao$, $R'=\Re \phi$, $S'=\Im \phi$, and $R(a)=S(a)=0$. Note also that if~$\phi\in \Car[\imag]$, then $\ex^{\Phi}\in \Carm[\imag]$, so $B$ maps
$\Car[\imag]$ into $\Carm[\imag]$.

\medskip\noindent
Suppose $\epsilon>0$ and $\Re \phi(t)\le -\epsilon$ for all $t\ge a$. Then $-R$ has derivative $-R'(t)\ge \epsilon$ for all $t\ge a$, so $-R$ is strictly increasing with image
$[-R(a),\infty)=[0,\infty)$ and compositional inverse $(-R)^{\inv}\in \Coo$. Making 
the change of variables $-R(s)=u$ for~$s\ge a$, we obtain for $t\ge a$ and $f\in \c_a[\imag]$, and with $s:=(-R)^{\inv}(u)$,
\begin{align*} \int_a^t \ex^{-\Phi(s)}f(s) \,ds\ &=\ \int_{0}^{-R(t)} \ex^{-\Phi(s)}f(s)\frac{1}{-R'(s)}\, du,\ \text{ and thus}\\
 |Bf(t)|\ &\le\  \ex^{R(t)}\cdot \left(\int_{0}^{-R(t)}\ex^u\, du \cdot \|f\|_{[a,t]} \right)\cdot \left\|\frac{1}{\Re \phi}\right\|_{[a,t]} \\
 &=\ \big[1-\ex^{R(t)}\big]\cdot \|f\|_{[a,t]} \cdot \left\|\frac{1}{\Re \phi}\right\|_{[a,t]}\\   &\le\ \|f\|_{[a,t]} \cdot \left\|\frac{1}{\Re \phi}\right\|_{[a,t]}\
  \le\ \|f\|_a\cdot \left\|\frac{1}{\Re \phi}\right\|_a.
  \end{align*}
Thus $B$ maps $\c_a[\imag]^{\b}$ into $\c_a[\imag]^{\b}\cap \Cao[\imag]$ and $B \colon \c_a[\imag]^{\b}\to \c_a[\imag]^{\b}$ 
is continuous with operator norm $\|B\|_{a}\le \big\|\frac{1}{\Re \phi}\big\|_a$.

\medskip\noindent
Next, suppose $\epsilon>0$ and $\Re \phi(t)\ge \epsilon$ for all $t\ge a$. Then $R'(t)\ge \epsilon$ for all $t\ge a$, so~$R(t) \ge \epsilon\cdot(t-a)$ for such $t$. Hence if
$f\in \c_a[\imag]^{\b}$, then $\ex^{-\Phi}f$ is integrable at~$\infty$.
Recall from \eqref{eq:integrable} that~$\c_a[\imag]^{\inte}$ is the $\C$-linear subspace of $\c_a[\imag]$ consisting of the~$g\in \c_a[\imag]$ that are integrable at $\infty$. 
We have
the $\C$-linear maps $$f\mapsto \ex^{-\Phi}f \colon \c_a[\imag]^{\b} \to \c_a[\imag]^{\inte}, \qquad \der_{\infty}^{-1}\colon \c_a[\imag]^{\inte}\to \Cao[\imag],\ \quad f\mapsto \ex^{\Phi}f\colon \Cao[\imag]\to \Cao[\imag]. $$ Composition yields    
the $\C$-linear operator $B\colon \c_a[\imag]^{\b}  \to \Cao[\imag]$,  $$  Bf(t)\ :=\ \ex^{\Phi(t)}\int_\infty^t \ex^{-\Phi(s)}f(s)\, ds \qquad(f\in \c_a[\imag]^{\b}).$$
It is a right inverse to $\der - \phi$ in the sense that
$(\der-\phi)\circ B$ is the identity on~$\c_a[\imag]^{\b}$. 
Note that $R$ is strictly increasing with image
$[0,\infty)$ and compositional inverse~${R^{\inv}\in \Coo}$. 
Making the change of variables $R(s)=u$ for $s\ge a$, we obtain for $t\ge a$ and~${f\in \c_a[\imag]^{\b}}$ with $s:=R^{\inv}(u)$,\begin{align*} \int_\infty^t \ex^{-\Phi(s)}f(s)\, ds\ &=\ -\int_{R(t)}^{\infty} \ex^{-\Phi(s)}f(s)\frac{1}{R'(s)}\, du,\ \text{ and thus}\\
 |Bf(t)|\ &\le\  \ex^{R(t)}\cdot \left( \int_{R(t)}^{\infty}\ex^{-u}\, du\right) \cdot \|f\|_t \cdot \left\|\frac{1}{\Re \phi}\right\|_t \\
  &\le\ \|f\|_t \cdot \left\|\frac{1}{\Re \phi}\right\|_t\ \le\ \|f\|_a \cdot \left\|\frac{1}{\Re \phi}\right\|_a
  \end{align*}
Hence $B$ maps $\c_a[\imag]^{\b}$ into $\c_a[\imag]^{\b}\cap \Cao[\imag]$, and as a $\C$-linear operator $\c_a[\imag]^{\b}\to \c_a[\imag]^{\b}$, $B$ is continuous with operator norm
$\|B\|_a \le \big\|\frac{1}{\Re \phi}\big\|_a$. 
If $\phi\in\Car[\imag]$, then $B$ maps~${\c_a[\imag]^{\b}\cap \Car[\imag]}$ into 
 $\c_a[\imag]^{\b}\cap\Carm[\imag]$.

\medskip\noindent
The case that for some $\epsilon>0$ we have $\Re\phi(t)\le -\epsilon$ for all $t\ge a$ is called the {\em attractive case},\index{function!attractive} and the case that for some
$\epsilon>0$ we have $\Re\phi(t)\ge \epsilon$ for all~$t\ge a$ is called the {\em repulsive case}.\index{function!repulsive}\index{repulsive!function} In both cases the above yields a continuous operator~$B \colon \c_a[\imag]^{\b}\to \c_a[\imag]^{\b}$
with operator norm $\le \big\|\frac{1}{\Re \phi}\big\|_a$
which is right-inverse to the operator $\der-\phi\colon\Cao[\imag]\to \c_a[\imag]$.
We denote this operator $B$ by~$B_{\phi}$ if we need to indicate its dependence on $\phi$. Note also its dependence on $a$. In both the attractive and the repulsive case, $B$ maps
$\c_a[\imag]^{\b}$ into $\c_a[\imag]^{\b}\cap\Cao[\imag]$, and if $\phi\in\Car[\imag]$ then
$B$ maps $\c_a[\imag]^{\b}\cap \Car[\imag]$ into 
 $\c_a[\imag]^{\b}\cap\Carm[\imag]$.  

\medskip
\noindent
Given a Hardy field $H$ and $f\in H[\imag]$ with  $\Re f \succeq 1$ we can choose $a$ and a representative of $f$ in  $\c_a[\imag]$, to be denoted also by $f$, such that $\Re f(t)\ne 0$ for all~$t\ge a$, and then $f\in \c_a[\imag]$ falls either under the attractive case or under the repulsive case. The original germ $f\in H[\imag]$ as well as the function $f\in \c_a[\imag]$ is accordingly said to be attractive, respectively repulsive.  (This agrees with the terminology introduced at the beginning of Section~\ref{sec:repulsive-normal}.)\index{germ!attractive}\index{germ!repulsive}\index{repulsive!germ}\index{repulsive!element}

\subsection*{Twists and right-inverses of linear operators over Hardy fields} 
Let $H$ be a Hardy field, $K:=H[\imag]$, and let $A\in K[\der]$ be a monic operator of
order $r\ge 1$, 
$$A\ =\ \der^r + f_1\der^{r-1}+\cdots + f_r, \qquad f_1,\dots, f_r\in K.$$ 
Take a real number $a_0$ and functions in $\c_{a_0}[\imag]$ that represent
the germs $f_1,\dots, f_r$ and to be denoted also by $f_1,\dots, f_r$. Whenever we increase below
the value of $a_0$, it is understood that we also update the functions $f_1,\dots, f_r$ accordingly, by restriction; the same holds for any function on $[a_0,\infty)$ that gets named. Throughout, $a$ ranges over $[a_0,\infty)$, and $f_1,\dots, f_r$ denote also the restrictions of these functions to $[a,\infty)$, and likewise for any function on $[a_0,\infty)$ that we name. Thus for any $a$ we have the $\C$-linear operator
$$A_a\ \colon\ \Car[\imag]\to\c_a[\imag], \quad y \mapsto y^{(r)} + f_1y^{(r-1)} + \cdots + f_ry.$$
Next, let $\fm\in H^\times$ be given. It gives rise to the twist 
$A_{\ltimes \fm}\in K[\der]$,  
$$A_{\ltimes \fm}\ :=\ \fm^{-1}A\fm\ =\ \der^r + g_1\der^{r-1}+\cdots + g_r, \qquad g_1,\dots, g_r\in K.$$
Now  [ADH, (5.1.1), (5.1.2), (5.1.3)] gives universal expressions for $g_1,\dots, g_r$ in terms of $f_1,\dots, f_r, \fm, \fm^{-1}$; for example, $g_1=f_1+r\fm^\dagger$. 
Suppose the germ $\fm$ is represented by a function 
in $\Cazr[\imag]^\times$, also denoted by $\fm$. Let $\fm^{-1}$ likewise do double duty as the multiplicative inverse of $\fm$ in $\Cazr[\imag]$. The expressions above can be used to show that the germs $g_1,\dots, g_r$ are represented by functions in $\c_{a_0}[\imag]$, to be denoted also by $g_1,\dots, g_r$, such that for all $a$ and all $y\in \Car[\imag]$ we have $$\fm^{-1} A_a(\fm y)\ =\ (A_{\ltimes \fm})_a(y), \text{ where }
  (A_{\ltimes \fm})_a(y)\ :=\  y^{(r)} + g_1y^{(r-1)} + \cdots + g_ry. $$
The operator $A_a\colon \Car[\imag]\to\c_a[\imag]$ is surjective (Proposition~\ref{prop:existence and uniqueness}); we aim to
construct a right-inverse of $A_a$ on the subspace $\c_a[\imag]^{\b}$ of $\c_a[\imag]$. 
For this, we assume given a splitting of $A$ over $K$, 
$$A\ =\ (\der-\phi_1)\cdots (\der-\phi_r), \qquad  \phi_1,\dots, \phi_r\in K.$$
Take functions in $\c_{a_0}[\imag]$, to be denoted also by $\phi_1,\dots, \phi_r$, that represent the germs $\phi_1,\dots, \phi_r$. We increase $a_0$ to arrange
$\phi_1,\dots, \phi_r\in \Cazrl[\imag]$. Note that for~$j=1,\dots,r$
the $\C$-linear map $\der-\phi_j\colon \Cao[\imag]\to \c_a[\imag]$ restricts to a
$\C$-linear map $A_j\colon \Caj[\imag]\to \Cajl[\imag]$, so that we obtain a map
$A_1\circ \cdots \circ A_r\colon \Car[\imag]\to\c_a[\imag]$. It is routine to verify  that for all sufficiently large $a$ we have
$$A_a\ =\ A_1\circ \cdots \circ A_r\ \colon\ \Car[\imag]\to\c_a[\imag].$$
We increase $a_0$ so that
$A_a=A_1\circ \cdots \circ A_r$ for all $a$. Note that $A_1,\dots, A_r$ depend on  $a$, but we prefer not to indicate this dependence notationally.

\medskip
\noindent
Now $\fm\in H^\times$ gives over $K$ the splitting
$$A_{\ltimes \fm}\ =\ (\der - \phi_1+\fm^\dagger)\cdots (\der-\phi_r+\fm^\dagger).$$
Suppose as before that the germ $\fm$ is represented by a function
$\fm\in \Cazr[\imag]^\times$. With the usual notational conventions we 
have $\phi_j-\fm^\dagger\in \Cazrl[\imag]$, giving the $\C$-linear map $\tilde A_j := \der-(\phi_j-\fm^\dagger)\colon \Caj[\imag]\to \Cajl[\imag]$
for $j=1,\dots,r$, which for all sufficiently large $a$ gives, just as for
$A_a$, a factorization
$$(A_{\ltimes \fm})_a\ =\ \tilde A_1\circ \cdots \circ \tilde A_r.$$
To construct a right-inverse of $A_a$ we now assume $\Re\phi_1,\dots, \Re\phi_r\succeq 1$.
Then we increase $a_0$ once more so that for all $t\ge a_0$,  $$\Re \phi_1(t),\dots, \Re \phi_r(t)\ne 0.$$
Recall that for $j=1,\dots,r$ we have
the continuous $\C$-linear operator 
$$B_j\ :=\ B_{\phi_j}\ \colon\ \c_a[\imag]^{\b}\to \c_a[\imag]^{\b}$$ from the previous subsection.   
The subsection on twisted integration now yields:
 
\begin{lemma}\label{cri} The continuous $\C$-linear operator
$$A_a^{-1}\ :=\  B_r\circ \cdots \circ B_1\ \colon\ \c_a[\imag]^{\b}\to \c_a[\imag]^{\b}$$ is a right-inverse of $A_a$: it maps $\c_a[\imag]^{\b}$ into $\c_a[\imag]^{\b}\cap \Car[\imag]$, and $A_a\circ A_a^{-1}$ is the identity on $\c_a[\imag]^{\b}$. For its operator norm we have $\|A_a^{-1}\|_a\ \le\  
\prod_{j=1}^r\big\|\frac{1}{\Re \phi_j}\big\|_a$.
\end{lemma}

\noindent 
Suppose $A$ is real in the sense that $A\in H[\der]$. Then by increasing $a_0$ we arrange that $f_1,\dots, f_{r}\in \c_{a_0}$.  Next, set 
$$\c_a^{\b}\ :=\ \c_a[\imag]^{\b}\cap\c_a \ =\  \big\{f\in \c_a:\, \|f\|_a<\infty\big\},$$
an $\R$-linear subspace of $\c_a$. Then the real part
$$\Re A_a^{-1}\ :\ \c_a^{\b} \to \c_a^{\b},\qquad (\Re A_a^{-1})(f)\ :=\ \Re\!\big(A_a^{-1}(f)\big)$$
is $\R$-linear and maps $\c_a^{\b}$ into $\Car$. Moreover, it is right-inverse to $A_a$ on $\c_a^{\b}$ in the sense that $A_a\circ \Re A_a^{-1}$ is the identity on $\c_a^{\b}$, and for $f\in \c_a^{\b}$,
$$\|(\Re A_a^{-1})(f)\|_a\ \le\ \|A_a^{-1}(f)\|_a.$$

\subsection*{Damping factors} Here $H$, $K$, $A$, $f_1,\dots, f_r$, $\phi_1,\dots, \phi_r$, $a_0$ are as in Lemma~\ref{cri}.
In particular, $r\in \N^{\ge 1}$, $\Re\phi_1,\dots, \Re \phi_r\succeq 1$, and $a$ ranges over $[a_0,\infty)$.
For later use we choose damping factors $u$ 
to make the operator $uA_a^{-1}$ more manageable than~$A_a^{-1}$.
For $j=0,\dots,r$ we set 
\begin{equation}\label{eq:Ajcirc}
A_j^{\circ}\ :=\ A_1\circ \cdots \circ A_j\ \colon\ \Caj[\imag]\to \c_a[\imag],
\end{equation} 
with $A_0^{\circ}$ the identity on $\c_a[\imag]$
and $A_r^\circ=A_a$, and
\begin{equation}\label{eq:Bjcirc}
B_j^\circ\ :=\ B_j\circ \cdots \circ B_1\ \colon\ \c_a[\imag]^{\b}  \to \c_a[\imag]^{\b},
\end{equation}
where $B_0^{\circ}$ is the identity
on $\c_a[\imag]^{\b}$ and $B_r^\circ=A_a^{-1}$. Then $B_j^\circ$ maps $\c_a[\imag]^{\b}$ in\-to~${{\c_a[\imag]^{\b}}\cap { \Caj[\imag]}}$ and $A_j^\circ \circ B_j^\circ$ is the identity on $\c_a[\imag]^{\b}$ by Lemma~\ref{cri}.

\begin{lemma}\label{teq} Let $u\in \Car[\imag]^\times$. Then for $i=0,\dots,r$ and $f\in \c_a[\imag]^{\b}$,
\begin{equation}\label{eq:derivatives of uA^-1}
\big[u\cdot A_a^{-1}(f)\big]^{(i)}\ =\ \sum_{j=r-i}^r u_{i,j} \cdot u\cdot B_j^{\circ}(f) \quad  \text{ in $\Carmi[\imag]$}
\end{equation}
with coefficient functions $u_{i,j}\in \Carmi[\imag]$ given by $u_{i,r-i}=1$, and for $0\le i< r$,
$$ u_{i+1,j}\ =\ \begin{cases}
u_{i,r}'+ u_{i,r}(u^\dagger + \phi_r)	& \text{if $j=r$,} \\
u_{i,j}'+ u_{i,j}(u^\dagger + \phi_j) + u_{i,j+1} & \text{if  $r-i\le j< r$.}\end{cases}$$
\end{lemma} 
\begin{proof} Recall that for $j=1,\dots,r$ and $f\in \c_a[\imag]^{\b}$ we have
$B_j(f)' = f + \phi_jB_j(f)$. 
It is obvious that \eqref{eq:derivatives of uA^-1} holds for $i=0$. Assuming \eqref{eq:derivatives of uA^-1} for a certain $i< r$ we get
$$\big[uA_a^{-1}(f)\big]{}^{(i+1)}\ =\  \sum_{j=r-i}^r u_{i,j}' \cdot uB_j^{\circ}(f) + \sum_{j=r-i}^r u_{i,j} \cdot \big[uB_j^{\circ}(f)\big]',$$ 
 and for $j=r-i,\dots, r$,
$$\big[uB_j^{\circ}(f)\big]'\ =\ u'B_j^{\circ}(f) + u\cdot \big[B_j^{\circ}(f)\big]'\
                    =\ u^\dagger\cdot uB_j^{\circ}(f) + 
                    uB_{j-1}^{\circ}(f) + \phi_j uB_j^{\circ}(f),$$
which gives the desired result. 
\end{proof} 

\noindent
Let $\fv\in \Cazr$ be such that $\fv(t)>0$ for all $t\ge a_0$, $\fv\in H$, $\fv\prec 1$.  Then we have the convex subgroup
 $$\Delta\ :=\ \big\{\gamma\in v(H^\times):\ \gamma=o(v\fv)\big\}$$ 
of $v(H^\times)$. {\em We assume that 
$\phi_1,\dots, \phi_r\preceq_{\Delta} \fv^{-1}$ in the asymptotic field $K$, where~$\phi_j$ and
 $\fv$ also denote their germs.}
For real~$\nu>0$ we have
$\fv^\nu\in (\Cazr)^\times$, so 
$$u\ :=\ \fv^\nu|_{[a,\infty)}\in (\Car)^\times, \qquad \|u\|_a<\infty.$$
In the next proposition $u$ has this meaning, a meaning which accordingly  varies with $a$.    
Recall that $A_a^{-1}$ maps $\c_a[\imag]^{\b}$ into $\c_a[\imag]^{\b}\cap \Car[\imag]$
with $\|A_a^{-1}\|_a<\infty$. 

{\samepage \begin{prop}\label{uban}  Assume $H$ is real closed and $\nu\in \Q$, $\nu > r$. Then:
\begin{enumerate}
\item[\rm(i)] the $\C$-linear operator $u A_a^{-1}\colon \c_a[\imag]^{\b} \to \c_a[\imag]^{\b}$ maps $\c_a[\imag]^{\b}$ into $\Car[\imag]^{\b}$; 
\item[\rm(ii)] $u A_a^{-1}\colon \c_a[\imag]^{\b} \to \Car[\imag]^{\b}$ is continuous;
\item[\rm(iii)] there is a real constant $c\ge 0$ such that $\|u A_a^{-1}\|_{a;r}\le c$ for all $a$;
\item[\rm(iv)] for all $f\in \c_a[\imag]^{\b}$ we have $uA_a^{-1}(f) \preceq \fv^{\nu}\prec 1$; 
\item[\rm(v)]  $\|u A_a^{-1}\|_{a;r}\to 0$ as $a\to \infty$.
\end{enumerate} 
\end{prop}}
\begin{proof} Note that
$\fv^\dagger\preceq_{\Delta} 1$ by [ADH, 9.2.10(iv)]. Denoting the germ of $u$ also by $u$ we have 
$u\in H$ and $u^\dagger=\nu\fv^\dagger \preceq_{\Delta} 1$, in particular, $u^\dagger\preceq \fv^{-1/2}$. Note that the
$u_{i,j}$ from Lemma~\ref{teq}---that is, their germs---lie in $K$. 
Induction on $i$  gives $u_{i,j}\preceq_{\Delta} \fv^{-i}$ for $r-i\le j\le r$. Hence $uu_{i,j}\prec_{\Delta} \fv^{\nu-i}\prec_{\Delta} 1$ for $r-i\le j\le r$. Thus for $i=0,\dots,r$ we have a real constant
$$c_{i,a}\ :=\ \sum_{j=r-i}^r \|u\,u_{i,j}\|_a \cdot \|B_j\|_a\cdots \|B_1\|_a\in [0,\infty)  $$
with $\big\|\big[uA_a^{-1}(f)\big]{}^{(i)}\big\|_a\le c_{i,a}\|f\|_a$ for all $f\in \c_a[\imag]^{\b}$. 
Therefore $uA_a^{-1}$ maps $\c_a[\imag]^{\b}$ into $\Car[\imag]^{\b}$, and the operator $u A_a^{-1}\colon \c_a[\imag]^{\b} \to \Car[\imag]^{\b}$ is continuous with $$\|uA_a^{-1}\|_{a;r}\ \le\ c_a:=\max\{c_{0,a},\dots,c_{r,a}\}.$$
As to (iii), this is because for all $i$,~$j$, $\|u\,u_{ij}\|_a$ is decreasing as a function
of $a$, and $\|B_j\|_a\le \big\|\frac{1}{\Re \phi_j}\big\|_a$ for all $j$. For $f\in \c_a[\imag]^{\b}$ we have $A_a^{-1}(f)\in \c_a[\imag]^{\b}$, so  (iv) holds. As to (v),  $u\,u_{i,j}\prec 1$ gives $\|uu_{ij}\|_a\to 0$ as $a\to \infty$, for all $i$,~$j$. In view of~$\|B_j\|_a\le \big\|\frac{1}{\Re \phi_j}\big\|_a$ for all $j$, this gives
$c_{i,a}\to 0$ as $a\to \infty$ for $i=0,\dots,r$, so~$c_a\to 0$ as $a\to\infty$. 
\end{proof}

%{\bf The next paragraph has been checked but commented out till there is a need for it.} Let $i$ range over $\{0,\dots,r\}$. 
%With the above notations we set
%$u_i(t)=\max\big(|u_{i,r-i}(t)|,\dots, |u_{i,r}(t)|\big)$ for $t\ge a$, so $u_i\in \Caz$ and the germ of
%$u_i$, denoted also by $u_i$, lies in $H$, with $uu_i\prec_{\Delta} \fv^{\nu-i}$.  Let $f\in \Caz[\imag]^{\b}$ and
%suppose for all $i$ there is given a  germ
%$b_i(f)\in H$ such that $B_{r-i}^{\circ}(f),\dots, B_r^{\circ}(f)\preceq b_i(f)$. Then we obtain from Lemma 2.1 that
 %$[uA_a^{-1}(f)]^{(i)}\preceq uu_ib_i(f)$. 

\section{Solving Split-Normal Equations over Hardy Fields}\label{sec:split-normal over Hardy fields}

\noindent
We construct here
solutions of suitable algebraic differential equations over Hardy fields.
These solutions lie in rings $\Car[\imag]^{\b}$ ($r\in \N^{\ge 1}$) and are obtained as
fixed points of certain contractive maps, as is common in solving differential equations. Here we use that~$\Car[\imag]^{\b}$ is a Banach space with respect to the norm $\|\cdot\|_{a;r}$. It will take some effort to
define the right contractions using the operators from Section~\ref{sec:IHF}.  

\medskip\noindent
In this section $H$, $K$, $A$, $f_1,\dots, f_r$, $\phi_1,\dots, \phi_r$, $a_0$ are as in Lemma~\ref{cri}.
In particular, $H$ is a Hardy field, $K=H[\imag]$, and
$$A=(\der-\phi_1)\cdots(\der-\phi_r)\qquad \text{where $r\in\N^{\ge 1}$, $\phi_1,\dots,\phi_r\in K$, $\Re\phi_1,\dots, \Re\phi_r\succeq 1$.}$$
Here $a_0$ is chosen so that we have representatives for $\phi_1,\dots, \phi_r$  in $\Cazrl[\imag]$, denoted also by $\phi_1,\dots,\phi_r$.
We let $a$ range over $[a_0,\infty)$.
In addition we assume that $H$ is real closed, and that
we are given a germ $\fv\in H^{>}$ such that~$\fv\prec 1$ and
$\phi_1,\dots, \phi_r\preceq_{\Delta} \fv^{-1}$ for the convex subgroup
$$\Delta\ :=\ \big\{\gamma\in v(H^\times):\ \gamma=o(v\fv)\big\}$$
of $v(H^\times)$. We increase $a_0$ so
that $\fv$ is represented by a function in $\Cazr$, also denoted by $\fv$, with $\fv(t)>0$ for all $t\ge a_0$.

\subsection*{Constructing fixed points over $H$} Consider a differential equation
\begin{equation}\label{eq:ADE}\tag{$\ast$}
A(y)\ =\ R(y),\qquad y\prec 1,
\end{equation}
where $R\in K\{Y\}$ has order $\le r$, degree $\le d\in \N^{\ge 1}$ and weight $\le w\in \N^{\ge r}$, with
$R\prec_{\Delta}\fv^w$. Now $R=\sum_{\j}R_{\j}Y^{\j}$ with $\j$ ranging here and below over the tuples~$(j_0,\dots, j_r)\in \N^{1+r}$  with $|\j|\le d$ and $\|\j\|\le w$; likewise for $\i$.
For each $\j$ we take a function in $\c_{a_0}[\imag]$
that  represents the germ $R_{\j}\in K$ and let $R_{\j}$ denote
this function as well as its restriction to any $[a,\infty)$.
Thus $R$ is represented on $[a,\infty)$ by a polynomial
$\sum_{\j}R_{\j}Y^{\j}\in \c_a[\imag]\big[Y, \dots, Y^{(r)}\big]$, to be denoted
also by $R$ for simplicity. This yields for each $a$ an evaluation map
$$f\mapsto R(f):=\sum_{\j}R_{\j}f^{\j}\ :\  \Car[\imag]\to \c_a[\imag].$$
As in [ADH, 4.2] we also have for every $\i$ the formal partial derivative
$$ R^{(\i)}\ :=\ \frac{\partial^{|\i|}R}{\partial^{i_0}Y\cdots \partial^{i_r}Y^{(r)}}\ \in\ \c_a[\imag]\big[Y,\dots, Y^{(r)}\big]$$ 
with $R^{(\i)}=\sum_{\j} R_{\j}^{(\i)}Y^{\j}$,  all $R_{\j}^{(\i)}\in \c_a[\imag]$ having their germs in $K$. 

\medskip
\noindent
A {\em solution of \eqref{eq:ADE} on $[a,\infty)$}\/\index{solution!split-normal equation} is a function $f\in \Car[\imag]^{\b}$ such that  $A_a(f)=R(f)$ and $f\prec 1$. One might try to obtain a solution as a fixed point of the operator~$f\mapsto A_a^{-1}\big(R(f)\big)$, but this operator might fail to be contractive
on a useful space of functions. Therefore we twist $A$ and arrange things so that we can use Proposition~\ref{uban}. 
In the rest of this section we fix $\nu\in \Q$ with $\nu > w$ (so $\nu > r$) such that 
$R\prec_{\Delta}\fv^\nu$ and~$\nu\fv^\dagger\not\sim \Re \phi_j$ in $H$
for $j=1,\dots,r$. (Note that such $\nu$ exists.) Then the twist 
$\tilde A:=A_{\ltimes\fv^\nu}=\fv^{-\nu} A\fv^{\nu}\in K[\der]$ splits over $K$ as follows: 
\begin{align*} \tilde A\ &=\ (\der -\phi_1+\nu\fv^\dagger)\cdots (\der-\phi_r+\nu\fv^\dagger), \quad \text{ with }\\ 
\phi_j-\nu\fv^{\dagger}\ &\preceq_{\Delta}\ \fv^{-1}, \quad \Re\phi_j-\nu\fv^\dagger\ \succeq\ 1 \qquad(j=1,\dots,r).
\end{align*}
We also increase $a_0$ so that  
$\Re\phi_j(t)-\nu\fv^\dagger(t)\ne 0$ for all $t\ge a_0$
and such that for all $a$ and $u:=\fv^\nu|_{[a,\infty)}\in (\Car)^\times$ the operator $\tilde{A}_a\colon \Car[\imag]\to \c_a[\imag]$ satisfies $$\tilde{A}_a(y)\ =\ u^{-1}A_a(uy)\qquad(y\in \Car[\imag]).$$
(See the explanations before Lemma~\ref{cri} for definitions of
$A_a$ and $\tilde{A}_a$.) We now increase $a_0$ once more,  fixing it for the rest of the section except in the subsection ``Preserving reality'', so as to obtain as in Lemma~\ref{cri}, with $\tilde{A}$ in the role of $A$, a 
right-inverse $\tilde{A}_a^{-1}\colon \c_a[\imag]^{\b} \to \c_a[\imag]^{\b}$ for such~$\tilde{A}_a$. 
%As in the previous section where we defined $A_a$ (for any $a\ge a_0$), the above factorization of $\tilde A$ gives the operator 
%$\tilde{A}_a\colon \Car[\imag]\to \Caz[\imag]$. It is easy to check that
%for $u:=\fv^\nu|_{[a,\infty)}\in (\Car)^\times$ we have 
%$\tilde{A}_a(f)= u^{-1}A_a(uf)$ for all $f\in \Car[\imag]$. As in the previous section this also gives us the operator
%$\tilde{A}_a^{-1}\colon \Caz[\imag]^{\b} \to \Caz[\imag]^{\b}$ considered in the next
%result.   

\begin{lemma} \label{bdua} We have a continuous operator \textup{(}not necessarily $\C$-linear\textup{)} $$\Xi_a\ :\ \Car[\imag]^{\b}\to \Car[\imag]^{\b},\quad
f\mapsto u\tilde{A}_a^{-1}\big(u^{-1}R(f)\big).$$
It has the property that $\Xi_a(f)\preceq \fv^ \nu\prec 1$ 
and $A_a\big(\Xi_a(f)\big)=R(f)$
for all $f\in \Car[\imag]^{\b}$. 
\end{lemma} 
\begin{proof} We have $\|u^{-1}R_\i\|_a<\infty$ for all $\i$, so
$u^{-1}R(f)=\sum_{\i}u^{-1}R_{\i}f^{\i}\in \c_a[\imag]^{\b}$ for all~$f\in \Car[\imag]^{\b}$, and thus $u\tilde{A}_a^{-1}\big(u^{-1}R(f)\big)\in \Car[\imag]^{\b}$  for such $f$, by Proposition~\ref{uban}(i).
Continuity of $\Xi_a$ follows from Proposition~\ref{uban}(ii) and continuity of
$f\mapsto u^{-1}R(f)\colon \Car[\imag]^{\b} \to \c_a[\imag]^{\b}$.
For $f\in \Car[\imag]^{\b}$ we have  $\Xi_a(f)\preceq \fv^\nu\prec 1$ by Proposition~\ref{uban}(iv),  and {\samepage
$$ u^{-1}A_a\big(\Xi_a(f)\big)\ =\ u^{-1}A_a\big[u\tilde{A}_a^{-1}\big(u^{-1}R(f)\big)\big]\ =\ \tilde{A}_a\big[\tilde{A}_a^{-1}\big(u^{-1}R(f)\big)\big]\ =\ u^{-1}R(f),$$
so $A_a\big(\Xi_a(f)\big)=R(f)$. } 
\end{proof}

\noindent
By Lemma~\ref{bdua},  each $f\in \Car[\imag]^{\b}$ with $\Xi_a(f)=f$  is a solution of \eqref{eq:ADE} on $[a,\infty)$.

\begin{lemma}\label{bdua, bds} 
There is a constant $C_a\in\R^{\geq}$ such that for all $f,g\in \Car[\imag]^{\b}$, 
$$ \|\Xi_a(f+g)-\Xi_a(f) \|_{a;r}\ \le\ C_a\cdot \max\!\big\{1, \|f\|_{a;r}^d\big\}\cdot\big(\|g\|_{a;r} + \cdots + \|g\|_{a;r}^d\big).$$
We can take these $C_a$ such that $C_a\to 0$ as $a\to \infty$, and we do so below. 
\end{lemma}
\begin{proof}
Let $f,g\in \Car[\imag]^{\b}$. We have the Taylor expansion
$$R(f+g)\ =\ \sum_{\i} \frac{1}{\i !}R^{(\i)}(f)g^{\i}\ =\ \sum_{\i}\frac{1}{\i !}\bigg[\sum_{\j}R_{\j}^{(\i)}f^\j\bigg]g^{\i}.$$
Now for all $\i$,~$\j$ we have $R^{(\i)}_{\j}\prec_{\Delta} \fv^\nu$ in $K$, so  
%$u^{-1}R_{\j}^{(\i)}\in \Caz[\imag]^{\b}$ with 
$u^{-1}R_{\j}^{(\i)}\prec 1$.  Hence
$$D_a\ :=\ \sum_{\i,\j} \big\|u^{-1}R^{(\i)}_{\j}\big\|_a\ \in\ [0,\infty)$$
has the property that $D_a\to 0$ as $a\to \infty$, and
$$\big\|u^{-1}\big(R(f+g)-R(f)\big)\big\|_a\ \le\ D_a\cdot\max\!\big\{1, \|f\|_{a;r}^d\big\}\cdot \big(\|g\|_{a;r} + \cdots + \|g\|_{a;r}^d\big).$$
So $h:= u^{-1}\big(R(f+g)-R(f)\big)\in \Caz[\imag]^{\b}$ gives $\Xi_a(f+g)-\Xi_a(f)=u\tilde{A}_a^{-1}(h)$, and 
$$ \|\Xi_a(f+g)-\Xi_a(f) \|_{a;r}= \|u\tilde{A}_a^{-1}(h) \|_{a;r} \le  \|u\tilde{A}_a^{-1} \|_{a;r}\cdot\|h\|_a.$$ 
Thus the lemma holds for $C_a:= \|u\tilde{A}_a^{-1} \|_{a;r}\cdot D_a$. 
\end{proof}

%{\bf In case it is useful later we observe that in the proof of the above lemma we actually have $u^{-1}R_{\j}^{(\i)}\prec_{\Delta} 1$ in $H[\imag]$, for all $\i, \j$; this has been checked}

\noindent 
In the proof of the next theorem we  use the well-known fact that the normed vector space $\Car[\imag]^{\b}$ over $\C$ is actually a Banach space.
Hence if~$S\subseteq \Car[\imag]^{\b}$ is  nonempty and closed and $\Phi\colon S\to S$ is
contractive (that is, there is a real number $\lambda\in [0,1)$ such that $\|\Phi(f)-\Phi(g)\|_{a;r}\leq \lambda \|f-g\|_{a;r}$
for all~$f,g\in S$), then~$\Phi$ has a unique fixed point~$f_0$, and $\Phi^n(f)\to f_0$ as~$n\to\infty$,  for every $f\in S$. (See for example \cite[Chapter~II, \S{}5, IX]{Walter}.)

\begin{theorem} \label{thm:fix} For all sufficiently large $a$ the operator $\Xi_a$ maps the closed ball $$\big\{f\in \Car[\imag]:\ \|f\|_{a;r}\le 1/2\big\}$$
of the Banach space $\Car[\imag]^{\b}$ into itself and has a unique fixed point on this ball. 
\end{theorem}
\begin{proof} We have $\Xi_a(0)=u\tilde{A}_a^{-1}(u^{-1}R_0)$, so $\|\Xi_a(0)\|_{a;r}\le \|u\tilde{A}_a^{-1}\|_{a;r}\|u^{-1}R_0\|_{a}$. Now~$\|u^{-1}R_0\|_a\to 0$ as $a\to \infty$, so by Proposition~\ref{uban}(iii) we can take $a$ 
so large that $\|u\tilde{A}_a^{-1}\|_{a;r}\|u^{-1}R_0\|_{a}\le \frac{1}{4}$. 
For $f$,~$g$ in the closed ball above we have by Lemma~\ref{bdua, bds},
$$ \|\Xi_a(f)-\Xi_a(g) \|_{a;r}\ =\  \|\Xi_a(f+(g-f))-\Xi_a(f) \|_{a;r}\ \le \ C_a\cdot d\|f-g\|_{a;r}.$$
Take $a$ so  large that also $C_a d \le \frac{1}{2}$. Then $\|\Xi_a(f)-\Xi_a(g)\|_{a;r}\le \frac{1}{2}\|f-g\|_{a;r}$. 
Applying this to $g=0$ we see that $\Xi_a$ maps the closed ball above to itself. Thus~$\Xi_a$ has a unique
fixed point on this ball. 
\end{proof}  

\noindent
Note that if $\deg R\leq 0$ (so $R=R_0$), then $\Xi_a(f)=u\tilde{A}_a^{-1}(u^{-1}R_0)$ is independent of~$f\in \Car[\imag]^{\b}$, so for sufficiently large $a$, the fixed point~$f\in \Car[\imag]^{\b}$ of $\Xi_a$ with~${\dabs{f}_{a;r}\leq 1/2}$ is $f=\Xi_a(0)=u\tilde{A}_a^{-1}(u^{-1}R_0)$. 
Here is a variant of Theorem~\ref{thm:fix}:
% produces ``near fixed points'' of $\Xi_a$:

\begin{lemma}\label{lem:fix h} 
Let $h\in\c_{a_0}^r[\imag]$ be such that $\dabs{h}_{a_0;r} \leq 1/8$. Then for sufficiently large $a$ there is a unique
$f\in \Car[\imag]^{\b}$ such that $\|f\|_{a;r}\le 1/2$ and  $\Xi_a(f)=f+h$.
%the closed ball $$\big\{f\in \Car[\imag]:\ \|f\|_{a;r}\le 1/2\big\}$$
%of $\Car[\imag]^{\b}$ contains a unique $f$ such that $\Xi_a(f)=f+h$.
\end{lemma}
\begin{proof}
Consider the operator $\Theta_a=\Xi_a-h\colon \Car[\imag]^{\b} \to \Car[\imag]^{\b}$ given by
$$  \Theta_a(y)\ :=\ \Xi_a(y)-h.$$
Arguing as in the proof of Theorem~\ref{thm:fix} we take $a$ so large that $\dabs{\Xi_a(0)}_{a;r} \leq \frac{1}{8}$.
Then
$\dabs{\Theta_a(0)}_{a;r} \leq \dabs{\Xi_a(0)}_{a;r}+\dabs{h}_{a;r} \leq \frac{1}{4}$.
Also, take $a$ so large that $C_ad \leq \frac{1}{2}$. Then for $f,g\in \Car[\imag]^{\b}$ with $\|f\|_{a;r}, \|g\|_{a;r}\le 1/2$
we have $$\dabs{\Theta_a(f)-\Theta_a(g)}_{a;r}\  =\  \dabs{\Xi_a(f)-\Xi_a(g)}_{a;r}\  \leq\  \textstyle\frac{1}{2}\dabs{f-g}_{a;r}.$$
Now finish as in the proof of Theorem~\ref{thm:fix}.
\end{proof}

\noindent
Next we investigate the difference between solutions of \eqref{eq:ADE} on $[a_0,\infty)$:

\begin{lemma}\label{lem:close} Suppose $f,g\in \Cazr[\imag]^{\b}$ and $A_{a_0}(f)=R(f)$, $A_{a_0}(g)=R(g)$. Then there are positive reals 
$E$,~$\epsilon$ such that for all $a$ 
there exists an $h_a\in \Car[\imag]^{\b}$ with the property that
for $\theta_a:=(f-g)|_{[a,\infty)}$, 
$$A_{a}(h_a)=0, \quad \theta_a-h_a\prec \fv^w, \quad
\|\theta_a-h_a\|_{a;r}\ \le\ E\cdot \|\fv^{\epsilon}\|_{a}\cdot \big(\|\theta_a\|_{a;r}+\cdots + \|\theta_a\|_{a;r}^d \big),$$
and thus $h_a\prec 1$ in case $f-g\prec 1$. 
\end{lemma}
\begin{proof} 
Set $\eta_a:= A_{a}(\theta_a)=R(f)-R(g)$, where $f$ and $g$ stand for their restrictions to $[a,\infty)$. From $R\prec \fv^{\nu}$ we get
$u^{-1}R(f)\in \c_a[\imag]^{\b}$ and $u^{-1}R(g)\in \c_a[\imag]^{\b}$, so
$u^{-1}\eta_a\in \c_a[\imag]^{\b}$. By Proposition~\ref{uban}(i),(iv) we have
$$\xi_a\ :=\ u\tilde{A}_a^{-1}(u^{-1}\eta_a)\in \Car[\imag]^{\b}, \qquad \xi_a\prec \fv^w.$$ 
Now $\tilde{A}_{a}(u^{-1}\xi_a)=u^{-1}\eta_a$, that is,
$A_{a}(\xi_a)=\eta_a$. Note that then $h_a:=\theta_a-\xi_a$ satisfies 
$A_{a}(h_a)=0$.
Now $\xi_a=\theta_a-h_a$ and $\xi_a=\Xi_a(g+\theta_a)-\Xi_a(g)$, hence   by  Lemma~\ref{bdua, bds} and its proof,
\begin{align*} \|\theta_a-h_a\|_{a;r}\ =\ \|\xi_a\|_{a;r}\ &\le\  C_a\cdot \max\!\big\{1,\|g\|_{a;r}^d\big\}\cdot \big(\|\theta\|_{a;r}+\cdots + \|\theta\|_{a;r}^d\big), \text{ with}\\
 C_a\ &:=\ \big\|u\tilde{A}_a^{-1}\big\|_{a;r}\cdot \sum_{\i,\j} \big\|u^{-1}R^{(\i)}_{\j}\big\|_a.
\end{align*}
Take a real $\epsilon>0$ such that $R\prec \fv^{\nu+\epsilon}$. This gives
a real $e>0$ such that $\sum_{\i,\j} \big\|u^{-1}R^{(\i)}_{\j}\big\|_a\le e\|\fv^\epsilon\|_a$ for all $a$. Now use Proposition~\ref{uban}(iii).
\end{proof}

\noindent
The situation we have in mind in the lemma above is that $f$ and $g$ are close at infinity, in the sense that $\|f-g\|_{a;r}\to 0$ as $a\to \infty$.
Then the lemma yields solutions of $A(y)=0$ that are very close to $f-g$ at infinity.  
However,  being very close at infinity as stated in Lemma~\ref{lem:close}, namely $\theta_a-h_a\prec\fv^w$ and the rest, is too weak for later use.
We take up this issue again in Section~\ref{sec:weights} below. 
(In Corollary~\ref{cor:h=0 => f=g} later in the present section we already show: if $f\neq g$ as germs, then $h_a\neq 0$ for sufficiently large $a$.) 

\subsection*{Preserving reality} We now assume in addition that $A$ and $R$ are real, that is,~${A\in H[\der]}$ and $R\in H\{Y\}$.  It is not clear that the fixed points constructed in the proof of Theorem~\ref{thm:fix} are then also real. Therefore we slightly modify this construction using real parts.  
We first apply the discussion following Lemma~\ref{cri} to $\tilde{A}$ as well as to $A$, increasing $a_0$ so  that
for all $a$ the  $\R$-linear real part
$$\Re\tilde{A}_a^{-1}\ \colon\ \c_a^{\b} \to \c_a^{\b}$$ maps~$\c_a^{\b}$ into $\Car$ and is right-inverse to $\tilde{A}_a$ on $(\Caz)^{\b}$, with  
$$\big\|(\Re\tilde{A}_a^{-1})(f)\big\|_a\le \big\|\tilde{A}_a^{-1}(f)\big\|_a\qquad\text{ for all  $f\in \c_a^{\b}$.}$$ 
Next we set
$$(\Car)^{\b}\ :=\ \big\{f\in \Car:\ \|f\|_{a;r}<\infty\big\}\ =\ \Car[\imag]^{\b}\cap \Car,$$
which is a real Banach space with respect to $\|\cdot\|_{a;r}$. 
Finally, this increasing of $a_0$ is done so that the original $R_{\j}\in \c_{a_0}[\imag]$ restrict to updated functions $R_{\j}\in \c_{a_0}$. 
For all $a$, take $u$, $\Xi_a$  as in Lemma~\ref{bdua}. This lemma has the following real analogue as a consequence:

\begin{lemma} \label{realbdua} The operator 
$$\Re\, \Xi_a\ : \ (\Car)^{\b}\to (\Car)^{\b},\quad
f\mapsto \Re\!\big(\Xi_a(f)\big)$$
satisfies $(\Re\, \Xi_a)(f)\preceq\fv^{\nu}$ for $f\in (\Car)^{\b}$, and 
any fixed point  of $\Re\,\Xi_a$ is a solution of \eqref{eq:ADE} on $[a,\infty)$.
\end{lemma} 

\noindent
Below the constants $C_a$ are as in Lemma~\ref{bdua, bds}.

\begin{lemma} \label{realbdua, bds}
For $f,g\in (\Car)^{\b}$,
$$\big\|(\Re\, \Xi_a)(f+g)-(\Re \Xi_a)(f)\big\|_{a;r}\ \le\ C_a\cdot \max\!\big\{1, \|f\|_{a;r}^d\big\}\cdot \big(\|g\|_{a;r} + \cdots + \|g\|_{a;r}^d\big).$$
\end{lemma}

%\noindent
%Any fixed point $f\in (\Car)^{\b}$ of $\Re\,\Xi_a$ is a solution of \eqref{eq:ADE} on $[a,\infty)$: for such $f$ we have $f=(\Re\,\Xi_a)(f)\prec 1$  by Lemma~\ref{realbdua}, and
%$$ A_a(f)\ =\ A_a\big((\Re\,\Xi_a)(f)\big)\ =\ \Re\!\big(A_a(\Xi_a(f))\big)\ =\ \Re\!\big(R(f)\big)\ =\ R(f),$$
%where the third equality uses $A_a\big(\Xi_a(f)\big)=R(f)$ from   Lemma~\ref{bdua}. 

\noindent
The next corollary is derived from Lemma~\ref{realbdua, bds} in the same way as Theorem~\ref{thm:fix} from Lem\-ma~\ref{bdua, bds}:

\begin{cor}\label{cor:fix} For all sufficiently large $a$ the operator $\Re\,\Xi_a$ maps the closed ball $$\big\{f\in \Car:\ \|f\|_{a;r}\le 1/2\big\}$$
of the Banach space $(\Car)^{\b}$ into itself and has a unique fixed point on this ball. 
\end{cor} 

\noindent
Here is the real analogue of Lemma~\ref{lem:fix h}, with a similar proof: 

\begin{cor}\label{cor:fix h}
Let $h\in\c_{a_0}^r$ be such that $\dabs{h}_{a_0;r} \leq 1/8$. Then for sufficiently large $a$  there is a unique $f\in \Car$ such that $ \|f\|_{a;r}\le 1/2$ and $(\Re\Xi_a)(f)=f+h$.
\end{cor} 

\noindent
We also have a real analogue of Lemma~\ref{lem:close}:

\begin{cor} Suppose $f,g\in (\Cazr)^{\b}$ and $A_{a_0}(f)=R(f)$, $A_{a_0}(g)=R(g)$. Then there are positive reals $E$,~$\epsilon$ such that
for all $a$ there exists an $h_a\in (\Car)^{\b}$ with the property that
for $\theta_a:=(f-g)|_{[a,\infty)}$, 
$$A_{a}(h_a)=0, \quad \theta_a-h_a\prec \fv^{w},\quad
\|\theta_a-h_a\|_{a;r}\ \le\ E\cdot \|\fv^{\epsilon}\|_{a}\cdot \big(\|\theta_a\|_{a;r}+\cdots + \|\theta_a\|_{a;r}^d \big).$$
\end{cor}
\begin{proof}
Take $h_a$ to be the real part of an $h_a$ as in Lemma~\ref{lem:close}. 
\end{proof}

\subsection*{Some useful bounds}
To prepare for Section~\ref{sec:weights} we derive in this subsection some bounds from Lem\-mas~\ref{bdua, bds} and~\ref{lem:close}. Throughout we assume $d,r\in \N^{\ge 1}$. 
We begin with an easy inequality:

\begin{lemma}\label{lem:inequ power d}
Let  $(V,\dabs{\, \cdot\, })$ be a normed $\C$-linear space, and $f,g\in V$. Then
$$\dabs{f+g}^d\ 	\leq\  
2^d\cdot \max\!\big\{1,\dabs{f}^d\big\}\cdot\max\!\big\{1,\dabs{g}^d\big\}.$$
\end{lemma}
\begin{proof}
Use that 
$\dabs{f+g}   \leq \dabs{f}+\dabs{g} \leq 2 \max\!\big\{1,\dabs{f}\big\}\cdot\max\!\big\{1,\dabs{g}\big\}$. 
%$$\dabs{f+g}^d\ \leq\  \sum_{i=0}^d {d\choose i} \dabs{f}^{d-i} \dabs{g}^i\   \leq\  2^d\cdot \max\!\big\{1,\dabs{f}^d\big\}\cdot\max\!\big\{1,\dabs{g}^d\big\},$$
%as claimed.
\end{proof}

\noindent
Now let  $u$, $\Xi_a$ be as in Lem\-ma~\ref{bdua}. By that lemma, the operator
$$\Phi_a\ \colon\ \Car[\imag]^{\b}\times\Car[\imag]^{\b} \to \Car[\imag]^{\b}, \quad (f,y)\mapsto \Xi_a(f+y)-\Xi_a(f)$$ 
is continuous. 
Furthermore $\Phi_a(f,0)=0$ for  $f\in \Car[\imag]^{\b}$ and
\begin{equation}\label{eq:Phia difference}
\Phi_a(f,g+y)-\Phi_a(f,g)\ =\ \Phi_a(f+g,y) \qquad\text{for $f,g,y\in \Car[\imag]^{\b}$.}
\end{equation}

\begin{lemma}\label{lem:2.1 summary}
There are $C_a,C_a^+\in\R^{\geq}$ such that for all $f,g,y\in\Car[\imag]^{\b}$,
\begin{equation}\label{eq:2.1, 1}
\dabs{\Phi_a(f,y)}_{a;r}  
 \ \le\ C_a\cdot \max\!\big\{1, \|f\|_{a;r}^d\big\}\cdot \big(\|y\|_{a;r} + \cdots + \|y\|_{a;r}^d\big),
\end{equation}
\vskip-1.5em
\begin{multline}\label{eq:2.1, 2}
\dabs{\Phi_a(f,g+y)-\Phi_a(f,g)}_{a;r}  \ \leq\ \\ C_a^+\cdot\max\!\big\{1,\dabs{f}_{a;r}^d\big\}\cdot\max\!\big\{1, \|g\|_{a;r}^d\big\}\cdot \big(\|y\|_{a;r} + \cdots + \|y\|_{a;r}^d\big).
\end{multline}
We can take these $C_a,C_a^+$ such that $C_a, C_a^+\to 0$ as $a\to\infty$, and do so below.
\end{lemma}

\begin{proof}
The $C_a$ as in Lemma~\ref{bdua, bds} satisfy the requirements on the $C_a$ here.
Now let~$f,g,y\in\Car[\imag]^{\b}$. Then
by \eqref{eq:Phia difference} and \eqref{eq:2.1, 1} we have
$$ \dabs{\Phi_a(f,g+y)-\Phi_a(f,g)}_{a;r}  \ \leq\ C_a \cdot\max\!\big\{1, \|f+g\|_{a;r}^d\big\}\cdot \big(\|y\|_{a;r} + \cdots + \|y\|_{a;r}^d\big).$$
Thus by Lemma~\ref{lem:inequ power d}, 
$C_a^+ := 2^d\cdot C_a$ has the required property.
\end{proof}

\noindent
Next, let $f$, $g$ be as in the hypothesis of Lemma~\ref{lem:close} and take $E$, $\epsilon$, and $h_a$ (for each $a$) as in its conclusion. 
Thus for all $a$ and $\theta_a:=(f-g)|_{[a,\infty)}$, 
$$\|\theta_a-h_a\|_{a;r}\ \le\ E\cdot \|\fv^{\epsilon}\|_{a}\cdot \big(\|\theta_a\|_{a;r}+\cdots + \|\theta_a\|_{a;r}^d \big),$$
and if $f-g\prec 1$, then $h_a\prec 1$. So
$$\|\theta_a-h_a\|_{a;r}\ \le\ E\cdot \|\fv^{\epsilon}\|_{a}\cdot \big(\rho+\cdots+\rho^d\big),\quad \rho:=\dabs{f-g}_{a_0;r}.$$
We let
$$B_a\ :=\ \big\{ y\in\Car[\imag]^{\b}:\ \dabs{y-h_a}_{a;r} \leq 1/2\big\}$$
be the  closed ball  of radius $1/2$ around $h_a$ in $\Car[\imag]^{\b}$.
Using $\fv^\epsilon\prec 1$ we take $a_1\geq a_0$  so that  $\theta_a\in B_a$ for all $a\geq a_1$.
%below we let $a_1$ have this property.  
Then for $a\geq a_1$ we have $$\dabs{h_a}_{a;r}\  \leq\ \dabs{h_a-\theta_a}_{a;r} + \dabs{\theta_a}_{a;r}\ \leq\ \textstyle\frac{1}{2}+\rho,$$
and hence for $y\in B_a$, 
\begin{equation}\label{eq:dabs(y)}
\dabs{y}_{a;r}\ \leq\ \dabs{y-h_a}_{a;r} + \dabs{h_a}_{a;r}\ \leq\ \textstyle\frac{1}{2}+\big(\frac{1}{2}+\rho\big)\ =\ 1+\rho.
\end{equation}
Consider now the continuous operators
$$\Phi_a,\Psi_a\ :\ \Car[\imag]^{\b} \to \Car[\imag]^{\b},\qquad
\Phi_a(y):=\Xi_a(g+y)-\Xi_a(g), \quad \Psi_a(y):=\Phi_a(y)+h_a.$$
In the notation introduced above, $\Phi_a(y)=\Phi_a(g,y)$ for $y\in  \Car[\imag]^{\b}$. With $\xi_a$ as in the proof of Lemma~\ref{lem:close}
we also have~$\Phi_a(\theta_a)=\xi_a$ and  
$\Psi_a(\theta_a)=\xi_a+h_a=\theta_a$. Below we reconstruct the fixed point $\theta_a$ of $\Psi_a$ from $h_a$,  for sufficiently large $a$.   

\begin{lemma}\label{lem:Psin, b}
There exists $a_2\geq a_1$ such that for all $a\geq a_2$ we have
$\Psi_a(B_a)\subseteq B_a$, and $\dabs{{\Psi_a(y)-\Psi_a(z)}}_{a;r} \leq \frac{1}{2} \|y-z\|_{a;r}$ for all $y,z\in B_a$.
\end{lemma}
\begin{proof}  Take $C_a$ as in  Lemma~\ref{lem:2.1 summary}, and let $y\in B_a$. Then  by \eqref{eq:2.1, 1},
\begin{align*}  \|\Phi_a(y) \|_{a;r}\ &\le\ C_a\cdot \max\!\big\{1, \|g\|_{a;r}^d\big\}\cdot \big(\|y\|_{a;r} + \cdots + \|y\|_{a;r}^d\big), \quad \theta_a\in B_a, \text{ so}\\
\|\Psi_a(y)-h_a \|_{a;r}\  &\leq\ C_a M,\quad M:=\max\!\big\{1, \|g\|_{a_0;r}^d\big\}\cdot \big( (1+\rho)+\cdots+(1+\rho)^d\big).
\end{align*} 
Recall that $C_a\to 0$ as $a\to\infty$. Suppose $a\ge a_1$ is so large that $C_a M \le 1/2$. Then~$\Psi_a(B_a)\subseteq B_a$.
With $C_a^+$ as in Lemma~\ref{lem:2.1 summary}, \eqref{eq:2.1, 2} gives for $y,z\in\Car[\imag]^{\b}$,
\begin{multline*}
\dabs{\Phi_a(y)-\Phi_a(z)}_{a;r} \leq \\ C^+_a\cdot\max\!\big\{1, \|g\|_{a;r}^d\big\}\cdot\max\!\big\{1, \|z\|_{a;r}^d\big\}\cdot \big(\|y-z\|_{a;r} + \cdots + \|y-z\|_{a;r}^d\big).
\end{multline*}
{\samepage Hence with $N := \max\!\big\{1, \|g\|_{a_0;r}^d\big\}\cdot (1+\rho)^d\cdot d$ we obtain for $y,z\in B_a$ that
$$\dabs{\Psi_a(y)-\Psi_a(z)}_{a;r}\  \leq\ C^+_a N   \|y-z\|_{a;r}, $$
so
$\dabs{{\Psi_a(y)-\Psi_a(z)}}_{a;r} \leq \frac{1}{2} \|{y-z}\|_{a;r}$ if $C_a^+N\leq 1/2$. }
\end{proof}

\noindent
Below $a_2$ is as in Lemma~\ref{lem:Psin, b}.

\begin{cor}\label{cor:Psin, b}
If $a\geq a_2$, then  $\lim_{n\to\infty} \Psi_a^n(h_a)=\theta_a$ in $\Car[\imag]^{\b}$.
\end{cor}

\begin{proof} Let $a\ge a_2$. Then $\Psi_a$ has a unique fixed point on $B_a$. As $\Psi_a(\theta_a)=\theta_a\in B_a$, this fixed point is~$\theta_a$.
\end{proof}

\begin{cor}\label{cor:h=0 => f=g}
If $f\neq g$ as germs, then $h_a\neq 0$ for  $a\geq a_2$.
\end{cor}
\begin{proof}
Let $a\geq a_2$. Then $\lim_{n\to \infty} \Psi_a^n(h_a)=\theta_a$. If $h_a=0$, then $\Psi_a=\Phi_a$, and hence $\theta_a=0$, since $\Phi_a(0)=0$.
\end{proof}

\section{Smoothness Considerations} \label{sec:smoothness}

\noindent
We assume $r\in \N$ in this section.
We prove here as much smoothness of solutions of algebraic differential equations over Hardy fields as could be hoped for. In particular, the
solutions in $\Car[\imag]^{\b}$ of the equation \eqref{eq:ADE} in Section~\ref{sec:split-normal over Hardy fields} actually have their germs
in $\mathcal{C}^{<\infty}[\imag]$.
To make this precise, consider 
a ``differential'' polynomial $$P\  =\ P(Y,\dots,Y^{(r)})\ \in\ \c^n[\imag]\big[Y,\dots,Y^{(r)}\big].$$ We put {\em differential\/} in quotes since $\c^n[\imag]$  is not naturally 
a differential ring. Nevertheless, $P$ defines an obvious evaluation map
$$f  \mapsto P\big(f,\dots,f^{(r)}\big)\ :\  \Gr[\imag] \to \c[\imag].$$ We also have the ``separant''\index{separant}
$$S_P\ :=\ \frac{\partial P}{\partial Y^{(r)}}\ \in\ \c^n[\imag]\big[Y,\dots,Y^{(r)}\big]$$
of $P$. 
 
\begin{prop}\label{hardysmooth} 
Assume $n\ge 1$. Let $f \in \Gr[\imag]$ be such that  
$$P\big(f,\dots,f^{(r)}\big) = 0\in \c[\imag]\quad\text{ and }\quad S_P\big(f,\dots,f^{(r)}\big)\in\c[\imag]^\times.$$ 
Then~$f \in \c^{n+r}[\imag]$. 
Thus if $P\in  \Gi[\imag]\big[Y,\dots,Y^{(r)}\big]$, then $f\in \Gi[\imag]$.
Moreover, if~$P\in  \Ginf[\imag]\big[Y,\dots,Y^{(r)}\big]$, then~$f\in\Ginf[\imag]$,  and likewise with $\Gom[\imag]$ in place of~$\Ginf[\imag]$.
\end{prop}

\noindent
We deduce this from the lemma below, which has a complex-analytic hypothesis. Let~${U\subseteq \R\times \C^{1+r}}$ be open. Let $t$ range over $\R$ and
$z_0,\dots, z_r$ over $\C$,  and set~$x_j:=\Re z_j$, $y_j:=\Im z_j$ for $j=0,\dots,r$. 
We also set 
$$U(t,z_0,\dots,z_{r-1})\ :=\ \big\{z_r:(t,z_0,\dots,z_{r-1},z_r)\in U\big\},$$ 
an open subset of $\C$, and
we assume $\Phi\colon U \to \C$  and $n\ge 1$ are such that 
$$\Re \Phi,\  \Im \Phi\ :\ U \to \R$$ 
are $\c^n$-functions of
$(t, x_0,y_0,\dots,x_r,y_r)$, and for all~$t,z_0,\dots,z_{r-1}$ the function 
$$z_r\mapsto \Phi(t,z_0,\dots,z_{r-1},z_r)\ \colon\ U(t,z_0,\dots,z_{r-1})\to \C$$ is holomorphic (the complex-analytic hypothesis alluded to).

\begin{lemma}\label{lem:hardysmooth, complex}
{\samepage Let $I\subseteq \R$ be a nonempty open interval and suppose~$f\in\Gr(I)[\imag]$ is such that for all~$t\in I$, 
\begin{itemize}
\item $\big(t, f(t),\dots, f^{(r)}(t)\big)\in U$;
\item $\Phi\big(t, f(t),\dots, f^{(r)}(t)\big)=0$; and
\item $(\partial\Phi/\partial z_r)\big(t, f(t),\dots, f^{(r)}(t)\big)\ne 0$.
\end{itemize}
Then $f\in \c^{n+r}(I)[\imag]$.}
\end{lemma}
\begin{proof}
Set $A:= \Re \Phi,\ B:= \Im \Phi$ and $g:= \Re f,\ h:= \Im f$. Then for all $t\in I$, 
\begin{align*} A\big(t, g(t), h(t), g'(t),h'(t)\dots, g^{(r)}(t), h^{(r)}(t)\big)\ &=\ 0\\
B\big(t, g(t), h(t), g'(t),h'(t)\dots, g^{(r)}(t), h^{(r)}(t)\big)\ &=\ 0.
\end{align*}
Consider the $\c^n$-map $(A,B)\colon U \to \R^2$, with $U$ identified in the usual way with an open subset of $\R^{1+2(1+r)}$. The Cauchy-Riemann equations give
$$\frac{\partial \Phi}{\partial z_r}\ =\ \frac{\partial A}{\partial x_r}+ \imag\frac{\partial B}{\partial x_r}, \qquad \frac{\partial A}{\partial x_r}\ =\ \frac{\partial B}{\partial y_r}, \qquad  \frac{\partial B}{\partial x_r}\ =\ -\frac{\partial A}{\partial y_r}.$$
Thus the Jacobian matrix of the map $(A,B)$ with respect to its last two variables~$x_r$ and $y_r$ has determinant
$$ D\ =\ \left(\frac{\partial A}{\partial x_r}\right)^2 +  \left(\frac{\partial B}{\partial x_r}\right)^2\ =\ \left|\frac{\partial \Phi}{\partial z_r}\right|^2\ :\ U\to \R.$$
Let $t_0\in I$. Then $$D\big(t_0, g(t_0), h(t_0),\dots, g^{(r)}(t_0), h^{(r)}(t_0)\big)\ne 0,$$ so by the Implicit Mapping Theorem \cite[(10.2.2), (10.2.3)]{Dieudonne} we have a connected open neighborhood~$V$ of the point
$$\big(t_0, g(t_0),h(t_0),\dots, g^{(r-1)}(t_0), h^{(r-1)}(t_0)\big)\in \R^{1+2r},$$ 
open intervals $J,K\subseteq\R$ containing $g^{(r)}(t_0)$, $h^{(r)}(t_0)$, respectively, and a $\c^n$-map
$$(G,H)\colon V \to J\times K$$
such that $V\times J\times K\subseteq U$ and the zero set of $\Phi$ on $V\times J\times K$ equals the graph of~$(G,H)$.
Take an open subinterval $I_0$ of $I$ with $t_0\in I_0$ such that for all $t\in I_0$, 
$$ \big(t, g(t), h(t),g'(t), h'(t),\dots, g^{(r-1)}(t), h^{(r-1)}(t), g^{(r)}(t), h^{(r)}(t)\big)\in V\times J\times K.$$
Then the above gives that for all $t\in I_0$ we have
\begin{align*} 
g^{(r)}(t)\ =\ G&\big(t, g(t), h(t),g'(t), h'(t),\dots, g^{(r-1)}(t), h^{(r-1)}(t)\big),\\
h^{(r)}(t)\ =\ H&\big(t, g(t), h(t),g'(t), h'(t),\dots, g^{(r-1)}(t), h^{(r-1)}(t)\big).
\end{align*} 
It follows easily from these two equalities that $g,h$ are of class
$\c^{n+r}$ on $I_0$.
\end{proof}

\noindent
Let $f$ continue to be as in Lemma~\ref{lem:hardysmooth, complex}.  If $\Re\Phi$, $\Im\Phi$ are $\Ginf$, then by taking $n$ arbitrarily high we conclude that $f\in \Ginf(I)[\imag]$. Moreover:

\begin{cor}\label{cor:hardysmooth, complex}
If  $\Re\Phi$, $\Im\Phi$ are  real-analytic, then~$f\in \Gom(I)[\imag]$. 
\end{cor}
\begin{proof}
Same as that of Lemma~\ref{lem:hardysmooth, complex}, with  the reference to~\cite[(10.2.3)]{Dieudonne} replaced  by \cite[(10.2.4)]{Dieudonne} to obtain that $G$, $H$ are real-analytic, and  noting that then the last displayed relations for $t\in I_0$ force $g$, $h$ to be real-analytic on $I_0$ by~\cite[(10.5.3)]{Dieudonne}.
\end{proof}

\begin{lemma}\label{smo} 
Let $I\subseteq\R$ be a nonempty open interval, $n\ge 1$, 
and $$P\  =\ P\big(Y,\dots,Y^{(r)}\big)\ \in\ \Gn(I)[\imag]\big[Y,\dots,Y^{(r)}\big].$$ 
Let $f\in \Gr(I)[\imag]$ be such that 
$$P\big(f,\dots,f^{(r)}\big) = 0\in \c(I)[\imag]\quad\text{ and }\quad
(\partial P/\partial Y^{(r)})\big(f, \dots, f^{(r)}\big) \in \c(I)[\imag]^\times.$$ 
Then $f \in \c^{n+r}(I)[\imag]$. 
Moreover, if $P\in\Ginf(I)[\imag]\big[Y,\dots,Y^{(r)}\big]$, then $f\in \Ginf(I)[\imag]$, and likewise with
$\Gom(I)[\imag]$ in place of $\Ginf(I)[\imag]$.
\end{lemma} 
\begin{proof}
Let $P=\sum_{\i} P_{\i} Y^{\i}$ where all $P_{\i}\in \Gn(I)[\imag]$. Set
$U:=I\times\C^{1+r}$, and consider the map $\Phi\colon U\to\C$ given by 
$$\Phi(t,z_0,\dots,z_r):=\sum_{\i} P_{\i}(t) z^{\i}\quad\text{where $z^{\i}:=z_0^{i_0}\cdots z_r^{i_r}$ for $\i=(i_0,\dots,i_r)\in\N^{1+r}$.}$$
From Lemma~\ref{lem:hardysmooth, complex}  we obtain  $f\in\c^{n+r}(I)[\imag]$. In view of Corollary~\ref{cor:hardysmooth, complex} and the remark preceding it,  and replacing~$n$ by~$\infty$ respectively $\omega$, this argument also gives the second part of the lemma. 
\end{proof}

\noindent
Proposition~\ref{hardysmooth} follows from Lemma~\ref{smo} by taking suitable representatives of the germs involved. 
%\begin{proof}[Proof of Proposition~\ref{hardysmooth}]
%Take $a \in \R$ and a representative of the germ $f$, denoted also by $f$ for convenience, such that $f \in 
%\Car[\imag]$.  Let $\i$ range over the finite set of tuples in $\N^{1+r}$ with~$P_{\i} \ne 0$. Increasing $a$ if necessary and replacing~$f$ by the corresponding restriction we take representatives of the coefficients~$P_{\i}$ of~$P$, denoted also by $P_{\i}$ for convenience, that lie in 
%$\Can[\imag]$, such that  $\sum_{\i} P_{\i}f^{\i} = 0$ on $[a,\infty)$. Increasing $a$ again and taking the corresponding restrictions of the functions~$f$ and~$P_{\i}$ we abuse notation and set
%$$P\ :=\ \sum_{\i}P_{\i}Y^{\i}\in  \Can[\imag]\big[Y,\dots,Y^{(r)}\big],$$
%and arrange that for the corresponding ``separant''
%$ S_P := \frac{\partial P}{\partial Y^{(r)}}\in  \Can[\imag]\big[Y, \dots , Y^{(r)}\big]$ 
%we have $S_P(f,\dots,f^{(r)})\in \Caz[\imag]^\times$.
%Now Lemma~\ref{smo} implies $f\in\c^{n+r}_a[\imag]$; hence the germ $f$ lies in $\c^{n+r}[\imag]$.
%If $P\in\Gi[\imag]\big[Y,\dots,Y^{(r)}\big]$, then we can take $n$ arbitrarily large, and thus 
%$f\in\Gi[\imag]$.
%If $P$ is in $\Ginf[\imag]\big[Y,\dots,Y^{(r)}\big]$ or in $\Gom[\imag]\big[Y,\dots,Y^{(r)}\big]$, then we argue
%in a similar way, with $n$ replaced by $\infty$ respectively $\omega$ (and the notational conventions $\infty+r:=\infty$, $\omega+r:=\omega$).
%\end{proof}
Let now $H$ be a Hardy field and $P\in H[\imag]\{Y\}$ of order~$r$. Then~$P\in \Gi[\imag]\big[Y,\dots,Y^{(r)}\big]$, and so $P(f):=P\big(f,\dots,f^{(r)}\big)\in\c[\imag]$ for $f\in\Gr[\imag]$
 as explained in the beginning of this section.

%\begin{cor}\label{cor:ADE smooth}\marginpar{checked but replaced by more flexible version} 
%Suppose 
%$P_{>1}\prec f_r$ where $P_1 = f_0Y +\cdots +f_rY^{(r)}$, $f_0,\dots,f_r$ in~$H[\imag]$. Let $f \in \mathcal{C}^r[\imag]$ be such that $f,\dots,f^{(r)}\preceq 1$ and $P(f) = 0$. Then $f \in \Gi[\imag]$. If~$H\subseteq
%\Ginf$, then $f\in \Ginf[\imag]$. If~$H\subseteq \Gom$, then~$f\in\Gom[\imag]$.
%\end{cor}
%\begin{proof} 
%We arrange 
%$f_r = 1$, so $P_{>1}\prec  1$.
%Then $S_P=\frac{\partial P}{\partial Y^{(r)}}=1+S$ where $S:=\frac{\partial P_{>1}}{\partial Y^{(r)}}\prec 1$
%and thus $S_P(f,\dots,f^{(r)})=1+S(f,\dots,f^{(r)})$ where $S(f,\dots,f^{(r)})\prec 1$,   so
%$S_P(f,\dots,f^{(r)})\in\c[\imag]^\times$.  Now appeal to Proposition~\ref{hardysmooth}.
%\end{proof}

\medskip\noindent
For notational convenience we set
$$\c^n[\imag]^{\preceq} := \big\{ f\in \c^n[\imag]: f,f',\dots,f^{(n)}\preceq 1\big\},\qquad (\c^n)^{\preceq}:=\c^n[\imag]^{\preceq}\cap\c^n,$$
and likewise with $\prec$ instead of $\preceq$.
Then $\c^n[\imag]^{\preceq}$ is a $\C$-subalgebra of $\c^n[\imag]$ and $(\c^n)^{\preceq}$ is
an $\R$-subalgebra of $\c^n$. Also, $\c^n[\imag]^{\prec}$ is an ideal of $\c^n[\imag]^{\preceq}$,
and likewise with~$\c^n$ instead of $\c^n[\imag]$.
We have $\c^n[\imag]^{\preceq}\supseteq \c^{n+1}[\imag]^{\preceq}$ and 
$(\c^n)^{\preceq}\supseteq (\c^{n+1})^{\preceq}$, and likewise with $\prec$ instead of $\preceq$.
Thus 
in the notation from Section~\ref{sec:upper lower bds}: 
$$\Calinf[\imag]^{\preceq}\ =\ \bigcap_n \c^n[\imag]^{\preceq},\quad
\mathcal I\ =\ \bigcap_n \c^n[\imag]^{\prec},\quad
(\Calinf)^{\preceq}\ =\ \bigcap_n (\c^n)^{\preceq}.$$

\begin{cor}\label{cor:ADE smooth}
Suppose 
$$P=Y^{(r)} + f_1Y^{(r-1)}+ \cdots + f_rY - R\quad\text{ with $f_1,\dots,f_r$ in~$H[\imag]$ and $R_{\ge 1}\prec 1$.}$$  
Let $f \in \mathcal{C}^r[\imag]^{\preceq}$ be such that $P(f) = 0$. Then $f \in \Gi[\imag]$. Moreover, if~$H\subseteq
\Ginf$, then $f\in \Ginf[\imag]$, and if~$H\subseteq \Gom$, then~$f\in\Gom[\imag]$.
\end{cor}
\begin{proof} We have $S_P=\frac{\partial P}{\partial Y^{(r)}}=1-S$ with $S:=\frac{\partial R_{\ge 1}}{\partial Y^{(r)}}\prec 1$
and thus $$S_P(f,\dots,f^{(r)})\ =\ 1-S(f,\dots,f^{(r)}), \qquad S(f,\dots,f^{(r)})\prec 1,$$   so
$S_P(f,\dots,f^{(r)})\in\c[\imag]^\times$.  Now appeal to Proposition~\ref{hardysmooth}.
\end{proof}

\noindent
Thus the germ of any solution on $[a,\infty)$ of the asymptotic equation \eqref{eq:ADE} of Section~\ref{sec:split-normal over Hardy fields} lies in $\Gi[\imag]$, and even in $\Ginf[\imag]$ (respectively $\Gom[\imag]$) if $H$ is in addition a
$\Ginf$-Hardy field (respectively, a $\Gom$-Hardy field).  

\begin{cor}\label{cor:hardynormal} Suppose $(P, 1,\hat a)$ is a normal slot in $H[\imag]$ of order $r\ge 1$, and $f\in \mathcal{C}^r[\imag]^{\preceq}$, $P(f) = 0$. Then $f \in \Gi[\imag]$. If~$H\subseteq
\Ginf$, then $f\in \Ginf[\imag]$. If~$H\subseteq \Gom$, then~$f\in\Gom[\imag]$.
\end{cor}
\begin{proof} Multiplying $P$ by an element of $H[\imag]^\times$ we arrange 
$$P_1\ =\ Y^{(r)}+ f_1Y^{(r-1)}+ \cdots + f_rY,\quad f_1,\dots,f_r\in H[\imag],$$ 
and then the hypothesis of Corollary~\ref{cor:ADE smooth} is satisfied.
\end{proof}

\medskip
\noindent
For the differential subfield $K:=H[\imag]$ of the differential ring~$\Calinf[\imag]$ we have: 

\begin{cor} 
Suppose $f\in\Calinf[\imag]$ is such that $P(f)=0$ and $f$ generates a differential subfield $K\langle f\rangle$ of $\Calinf[\imag]$ over $K$.
If $H$ is a $\Ginf$-Hardy field, then~$K\langle f\rangle\subseteq\Ginf[\imag]$, and likewise with
$\Gom$ in place of $\Ginf$.
\end{cor}
\begin{proof}
Suppose $H$ is a $\Ginf$-Hardy field; it suffices to show $f\in\Ginf[\imag]$.
We may assume that $P$ is a minimal annihilator of $f$ over $K$; then $S_P(f)\neq 0$ in $K\langle f\rangle$
and so $S_P(f)\in\c[\imag]^\times$. Hence the claim follows from Proposition~\ref{hardysmooth}.
\end{proof}

\noindent
With $H$ replacing $K$ in this proof we obtain the  ``real'' version:

\begin{cor} \label{realginfgom} 
Suppose $f\in\Calinf$ is hardian over $H$ and $P(f)=0$ for some~$P\in H\{Y\}^{\ne}$.   
Then $H\subseteq \Ginf\ \Rightarrow\ f\in \Ginf$, and  $H\subseteq \Gom\ \Rightarrow\ f\in \Gom$.
\end{cor}

\noindent
This leads to: 

\begin{cor}\label{cor:Hardy field ext smooth}
Suppose $H$ is a $\Ginf$-Hardy field. Then every $\d$-algebraic Hardy field extension of $H$ is a $\Ginf$-Hardy field; in particular, $\Dx(H)\subseteq\Ginf$. Likewise with~$\Ginf$ replaced by $\Gom$.
\end{cor}

\noindent
In particular, $\Dx(\Q)\subseteq\Gom$~\cite[Theorems~14.3, 14.9]{Boshernitzan82}.

\medskip
\noindent 
Let $H$ be a $\Ginf$-Hardy field.  Then by Corollary~\ref{cor:Hardy field ext smooth},
$H$ is $\d$-maximal iff~$H$ has no proper $\d$-algebraic $\Ginf$-Hardy field extension; thus  every $\Ginf$-maximal Hardy field is $\d$-maximal (so $\Dx(H)\subseteq\Ex^\infty(H)$), and $H$ has a $\d$-maximal $\d$-algebraic $\Ginf$-Hardy field extension.   The same remarks apply with $\omega$ in place of $\infty$.

\subsection*{Existence and uniqueness theorems\astr} 
We finish this section with some existence and uniqueness results for algebraic differential equations.
From this subsection, only Corollary~\ref{cor:uniqueness ADE} is used later (in the proofs of Lemmas~\ref{lem:notorious 3.6 flabby} and~\ref{lem:notorious 3.6 flabby, K}, which are not needed for the proof of our main theorem).
First, let~$U$,~$\Phi$ be as in   Lemma~\ref{lem:hardysmooth, complex} for $n=1$;
the argument in the proof of that lemma combined with
the existence and uniqueness theorem for scalar differential equations \cite[\S{}11]{Walter} yields:

\begin{lemma}\label{lem:hardysmooth, uniqueness}
Let  $(t_0,c_0,\dots,c_r)\in U$ be such that
$$\Phi(t_0,c_0,\dots,c_r) = 0 \quad\text{and}\quad (\partial\Phi/\partial z_r)(t_0,c_0,\dots,c_r) \neq 0.$$
Then for some open interval $I\subseteq\R$ containing $t_0$ there is a unique~$f\in\Gr(I)[\imag]$ such that 
$\big(f(t_0),f'(t_0),\dots, f^{(r)}(t_0)\big)=(c_0,\dots, c_r)$ and
%$f^{(j)}(t_0)=z_j$ for $j=0,\dots,r$ and
for all $t\in I$,
\begin{equation}\label{eq:Phi=0}
\big(t, f(t),\dots, f^{(r)}(t)\big)\in U \quad\text{and}\quad \Phi\big(t, f(t),\dots, f^{(r)}(t)\big)=0.
\end{equation}
\end{lemma}
\begin{proof}
Set $A:=\Re\Phi$, $B:=\Im\Phi$, $a_j:=\Re c_j$, $b_j:=\Im c_j$ ($j=0,\dots,r$). As in the proof of Lemma~\ref{lem:hardysmooth, complex} we identify $U$ with an open subset of
$\R^{1+2(1+r)}$ and consider the $\c^1$-map $(A,B)\colon U\to\R^2$. The Jacobian matrix of the map $(A,B)$ with respect to its last two variables~$x_r$ and $y_r$ has determinant
$$ D\ =\ \left(\frac{\partial A}{\partial x_r}\right)^2 +  \left(\frac{\partial B}{\partial x_r}\right)^2\ =\ \left|\frac{\partial \Phi}{\partial z_r}\right|^2\ :\ U\to \R,$$
with
$D(t_0,a_0,b_0,\dots,a_r,b_r)\neq 0$, hence the Implicit Mapping Theorem \cite[(10.2.2)]{Dieudonne} yields a connected open neighborhood $V$ in $\R^{1+2r}$
of the point
$$u_0:=(t_0,a_0,b_0,\dots,a_{r-1},b_{r-1}),$$ open intervals $J,K\subseteq\R$ containing~$a_r$,~$b_r$,
respectively, such that 
$V\times J\times K\subseteq U$, as well as a $\c^1$-map~$F=(G,H)\colon V\to J\times K$ whose graph is~$\Phi^{-1}(0)\cap (V\times J\times K)$.
Now by  \cite[\S{}11, II]{Walter} we have an open interval $I\subseteq\R$ containing $t_0$ as well as a $\c^r$-map~$u\colon I\to\R^2$  such that~$\big(t_0,u(t_0),u'(t_0),\dots,u^{(r-1)}(t_0)\big)=u_0$
and for all $t\in I$:
$$\big(t,u(t),\dots,u^{(r-1)}(t)\big)\in V\quad\text{and}\quad
u^{(r)}(t)=F\big(t,u(t),\dots,u^{(r-1)}(t)\big).$$
Then the function $f\colon I\to\C$ with $(\Re f,\Im f)=u$ is  an element of~$\Gr(I)[\imag]$  such that
$\big(f(t_0), f'(t_0),\dots, f^{(r)}(t_0)\big)=(c_0,\dots, c_r)$ and \eqref{eq:Phi=0} holds for all $t\in I$.

Let any $f_1\in \Gr(I)[\imag]$ be given with  $\big(f_1(t_0), f_1'(t_0),\dots, f_1^{(r)}(t_0)\big)=(c_0,\dots, c_r)$  such that \eqref{eq:Phi=0} holds for all $t\in I$ with $f_1$ in place of $f$.
The closed subset 
$$S\ :=\ \big\{t\in I:\ \text{$f_1^{(j)}(t)=f^{(j)}(t)$ for $j=0,\dots,r$} \big\}$$ 
of $I$ contains $t_0$; it is enough to show that $S$ is open.
Towards this, let~$t_1\in S$.
The map $u_1:=(\Re   f_1,\Im   f_1)\colon I\to\R^2$ is of class $\c^r$ and  
$$\big(t_1,u_1(t_1),\dots,u^{(r)}_1(t_1)\big)\ =\ \big(t_1,u(t_1),\dots,u^{(r)}(t_1)\big)\in V\times J\times K,$$
which gives an open interval $I_1\subseteq I$ containing $t_1$ such that
$\big(t,u_1(t),\dots,u^{(r)}_1(t)\big)\in V\times J\times K$ for all $t\in I_1$.
Since $\Phi\big(t, f_1(t),\dots, f^{(r)}_1(t)\big)=0$ for $t\in I_1$, this yields
$u_1^{(r)}(t)=F\big(t,u_1(t),\dots,u_1^{(r-1)}(t)\big)$ for $t\in I_1$.
So $u_1=u$ on $I_1$ by the uniqueness part of \cite[\S{}11, III]{Walter}, hence $f_1=f$ on $I_1$, and thus $I_1\subseteq S$.
\end{proof}

\noindent
The second part of the proof gives a bit more: Suppose $I,J\subseteq \R$ are open intervals with $t_0\in I\cap J$ and the functions $f\in\Gr(I)[\imag]$, 
$g\in \Gr(J)[\imag]$ are such that 
$$\big(f(t_0),f'(t_0),\dots, f^{(r)}(t_0)\big)\ =\ (c_0,\dots, c_r)\ =\ \big(g(t_0),g'(t_0),\dots, g^{(r)}(t_0)\big),$$
 \eqref{eq:Phi=0} holds for all $t\in I$, and  \eqref{eq:Phi=0} holds for all $t\in J$ with $g$ instead of $f$. Assume also that
 $(\partial\Phi/\partial z_r)\big(t,f(t),\dots,f^{(r)}(t)\big) \neq 0$ for all $t\in I$. Then 
 $$f(t)\ =\ g(t)\ \text{ for all }t\in I\cap J.$$ 
 Next, let $I\subseteq\R$ be a nonempty open interval and
$$P\  =\ P\big(Y,\dots,Y^{(r)}\big)\ \in\ \Go(I)[\imag]\big[Y,\dots,Y^{(r)}\big].$$ 
Applying Lemma~\ref{lem:hardysmooth, uniqueness} to the map 
$\Phi\colon U:=I\times\C^{1+r}\to\C$  introduced
in the proof of Lemma~\ref{smo}, we obtain:

\begin{lemma}
Let  $t_0\in I$ and $c_0,\dots,c_r\in\C$ be such that
$$\Phi(t_0,c_0,\dots,c_r)=0\quad\text{and}\quad (\partial \Phi/\partial z_r)(t_0,c_0,\dots,c_r)\neq 0.$$
Then there is an open interval $J\subseteq I$ contaning $t_0$
with a unique $f\in \Gr(J)[\imag]$   such that $\big(f(t_0),f'(t_0),\dots, f^{(r)}(t_0)\big)= (c_0,\dots, c_r)$ and
$P\big(f,\dots,f^{(r)}\big) = 0\in \c(J)[\imag]$. 
\end{lemma}

\noindent
This lemma and the remark following the proof of Lemma~\ref{lem:hardysmooth, uniqueness} yield:

\begin{cor}\label{cor:hardysmooth, uniqueness}
Given $t_0\in I$ and $c_0,\dots,c_r\in\C$, there is at most one function~$f\in\Gr(I)[\imag]$ such that 
$\big(f(t_0),f'(t_0),\dots, f^{(r)}(t_0)\big)= (c_0,\dots, c_r)$  as well as
$$P\big(f,\dots,f^{(r)}\big) = 0\in \c(I)[\imag]\quad\text{ and }\quad
(\partial P/\partial Y^{(r)})\big(f, \dots, f^{(r)}\big) \in \c(I)[\imag]^\times.$$ 
\end{cor}

\noindent
Now let $a$ range over $\R$,   $\i$   over $\N^{1+r}$, and
$$P \  =\ P\big(Y,\dots,Y^{(r)}\big)\ =\ \sum_{\i} P_{\i} Y^{\i} \qquad \text{(all $P_{\i}\in \c^1_{a}[\imag]$)}$$
over polynomials in $\c^1_{a}[\imag]\big[Y,\dots,Y^{(r)}\big]$ of degree at most $d\in\N^{\geq 1}$, and set
$P_{\ge 1}:=\sum_{\abs{\i}\ge 1} P_{\i}Y^{\i}=P-P_0$. 
Recall that in Section~\ref{sec:IHF} we defined 
$$P(f)\ :=\ \sum_{\i} P_{\i} f^{\i}\in\c_a[\imag]\qquad (f\in \Car[\imag]).  $$

%\begin{cor} \label{cor:uniqueness ADE} \marginpar{checked but seems not quite enough to justify subsequent remark} 
%There is an $E=E(d,r)\in\N^{\geq 1}$ with the following property:  if 
%$$P_1 = f_0Y +\cdots +f_rY^{(r)}\quad f_0,\dots,f_{r-1}\in\c^1_{a}[\imag],\  f_r\in\c^1_{a}[\imag]^\times,\ \dabs{P_{>1}/f_r}_{a} \leq 1/E,$$
%then for all $t_0\in \R^{\geq a}$ and $c_0,\dots,c_r\in\C$
 %there is at most one~$f\in\Car[\imag]$   such that
%$$P(f)=0,\quad\big(f(t_0),f'(t_0),\dots, f^{(r)}(t_0)\big)= (c_0,\dots, c_r), 
 %\quad\text{ and }\quad \dabs{f}_{a;r}\leq 1.$$
%\end{cor}
%\begin{proof} 
%Set $E:=2d^2D$ where $D=D(1,d,r)$ is as in Corollary~\ref{cor:bound on P(f)}. Let~$P$ be as in the hypothesis and
 %$t_0\in \R^{\geq a}$, $c_0,\dots,c_r\in\C$.
%We can arrange 
%$f_r = 1$, so $\dabs{P_{>1}}_a\leq 1/E$.
%Then $\partial P/\partial Y^{(r)}=1+S$ where~$S:=\partial P_{>1}/\partial Y^{(r)}$
%and therefore $(\partial P/\partial Y^{(r)})(f,\dots,f^{(r)})=1+S(f,\dots,f^{(r)})$.
%We have~$\dabs{S}_a\leq d / E = 1/(2dD)$ and hence by Corollary~\ref{cor:bound on P(f)}:
%$$\dabs{S(f,\dots,f^{(r)})}_a\ \leq\  D \cdot \dabs{S}_a\cdot \big(\dabs{f}_{a;r}^1+\cdots+\dabs{f}_{a;r}^d\big) \leq   1/2.$$
%Thus $(\partial P/\partial Y^{(r)})(f,\dots,f^{(r)})\in\c_a[\imag]^\times$.  Now use Corollary~\ref{cor:hardysmooth, uniqueness}.
%\end{proof}

\begin{cor} \label{cor:uniqueness ADE}
There is an $E=E(d,r)\in\N^{\geq 1}$ with the following property:  if 
$$P = Y^{(r)}+  f_1Y^{(r-1)} +\cdots +f_rY - R,\quad f_1,\dots,f_{r}\in\c^1_{a}[\imag],\   \dabs{R_{\ge 1}}_{a} \leq 1/E,$$
then for any $t_0\in \R^{> a}$ and $c_0,\dots,c_r\in\C$
 there is at most one~$f\in\Car[\imag]$   such that
$$P(f)=0,\quad\big(f(t_0),f'(t_0),\dots, f^{(r)}(t_0)\big)= (c_0,\dots, c_r), 
 \quad\text{ and }\quad \dabs{f}_{a;r}\leq 1.$$
\end{cor}
\begin{proof} 
Set $E:=2d(d+1)D$ with $D=D(0,d,r)$ as in Corollary~\ref{cor:bound on P(f)}. Let~$P$ be as in the hypothesis and
$f\in \Car[\imag]$, $\dabs{f}_{a;r}\le 1$. 
Then $\partial P/\partial Y^{(r)}=1-S$ where~$S:=\partial R_{\ge 1}/\partial Y^{(r)}$
and therefore $(\partial P/\partial Y^{(r)})(f,\dots,f^{(r)})=1-S(f,\dots,f^{(r)})$.
We have~$\dabs{S}_a\leq d / E = 1/\big(2(d+1)D\big)$ and hence by Corollary~\ref{cor:bound on P(f)}:
$$\dabs{S(f,\dots,f^{(r)})}_a\ \leq\  D \cdot \dabs{S}_a\cdot \big(1+\dabs{f}_{a;r}^1+\cdots+\dabs{f}_{a;r}^{d}\big)\ \leq\   1/2.$$
Thus $(\partial P/\partial Y^{(r)})(f,\dots,f^{(r)})\in\c_a[\imag]^\times$.  Now use Corollary~\ref{cor:hardysmooth, uniqueness}.
\end{proof}

\noindent
Thus in the context of Section~\ref{sec:split-normal over Hardy fields}, if $a$ is so large that
the functions $f_1,\dots,f_r$ and the $R_{\j}$ there are~$\c^1$ on $[a,\infty)$ with
$\dabs{R_{\j}}_{a}\leq 1/E(d,r)$, then
for all  $t_0\in \R^{> a}$ and~$c_0,\dots,c_r\in\C$, 
there is at most one $f\in\Car[\imag]$ with $\dabs{f}_{a;r}\leq 1$ such that~$A_a(f)=R(f)$ and $\big(f(t_0),f'(t_0),\dots, f^{(r)}(t_0)\big)= (c_0,\dots, c_r)$.

\subsection*{A theorem of Boshernitzan\astr}
Here we supply a proof of the following result stated in \cite[Theorem~11.8]{Boshernitzan82}, and to be used in Section~\ref{sec:holes perfect}. (The proof in loc.~cit.~is only indicated there.)  Below, $Y$ and $Z$ are distinct indeterminates.  

\begin{theorem}\label{thm:Bosh order 1}
Let $H$ be a Hardy field, $P\in H[Y,Z]^{\neq}$, and suppose $P(y, y')=0$ with  $y\in\c^1$ lying in a Hausdorff field extension of $H$. Then $y\in \Dx(H)$. 
\end{theorem}

\noindent
In particular, if $H$ is a $\d$-perfect Hardy field and $F$ is a  Hardy field properly extending~$H$, then $\operatorname{trdeg}(F|H) \geq 2$.

\medskip
\noindent
For the proof of Theorem~\ref{thm:Bosh order 1} we first observe: 
%note that an argument similar to the proof of Corollary~\ref{realginfgom} shows: \marginpar{should perhaps better be added
%after 6.3.8}

\begin{cor}\label{realginfgom, hausd}
Let $H$ be a Hausdorff field. Then $H\subseteq\c^n\ \Rightarrow\ H^{\operatorname{rc}}\subseteq\c^n$,
and likewise with $<\infty$, $\infty$, and $\omega$ in place of $n$.
\end{cor}
\begin{proof}
If $y\in H^{\operatorname{rc}}$ has minimum polynomial $P\in H[Y]$ over $H$,
then $P(y)=0\in\c$ and $S_P(y)=P'(y)\in\c^\times$. Now use Proposition~\ref{hardysmooth}.
\end{proof}

\begin{lemma}\label{le0} Suppose $f\in \Go$ oscillates.  Then we are in case $\rm{(i)}$ or case $\rm{(ii)}$:  \begin{enumerate}
\item[$\rm{(i)}$] there are arbitrarily large $s$ with $f'(s)=0$ and $f(s) >0$,  
\item[$\rm{(ii)}$] there are arbitrarily large $s$ with $f'(s)=0$ and $f(s) <0$.
\end{enumerate}
In case $\rm{(i)}$ there are also arbitrarily large $s$ with $f'(s)=0$ and $f(s)\le 0$, and in case $\rm{(ii)}$ there are also arbitrarily large $s$ with $f'(s)=0$ and $f(s)\ge 0$.  
\end{lemma}
\begin{proof} Let $f$ be represented by a $\Go$-function on an interval
$(a,+\infty)$, also denoted by $f$. Take $b>a$ such that $f(b)=0$, and then
$c>b$ with $f(c)=0$ such that~${f(t)\ne 0}$ for some $t\in (b,c)$.
Next, take $s\in [b,c]$ such that $|f(s)|=\max_{b\le t\le c} |f(t)|$.
Then~$f(s)\ne 0$ and $f'(s)= 0$. Since $b$ can be taken arbitrarily large, we are in case~(i) or in case (ii) above. (Of course, this is not an exclusive \textit{or}.) For the remaining part of the lemma, use that
in case (i) there are arbitrarily large $s>a$ where $f$ has a local minimum $f(s)\le 0$, and that in case (ii) there are arbitrarily large $s>a$ where $f$ has a local maximum $f(s)\ge 0$.   
\end{proof}

\noindent
In the next two lemmas $H$, $P$, $y$ are as in Theorem~\ref{thm:Bosh order 1}.

\begin{lemma}  \label{lem:order 1 H-hardian}
The germ  $y$ generates a Hardy field extension~$H\langle y\rangle$ of $H$.
If~$H\subseteq\c^\infty$, then $H\langle y\rangle\subseteq\c^\infty$, and likewise with~$\omega$ in place of $\infty$.
\end{lemma}
\begin{proof}
We are done if $y\in H^{\operatorname{rc}}$, since $H^{\operatorname{rc}}$ is a Hardy field
 with~$H\subseteq\c^\infty\Rightarrow H^{\operatorname{rc}}\subseteq\c^\infty$   and  $H\subseteq\c^\omega\Rightarrow H^{\operatorname{rc}}\subseteq\c^\omega$, 
by  Proposition~\ref{prop:Hardy field exts}.

Suppose $y\notin  H^{\operatorname{rc}}$.
We have the Hausdorff field $F:=H(y)\subseteq\c^1$, and its real closure is by Proposition~\ref{b1} the Hausdorff field
$$F^{\operatorname{rc}}\ =\ \big\{ z\in\c:\  \text{$Q(z)=0$ for some $Q\in F[Z]^{\neq}$} \big\}.$$
By Corollary~\ref{realginfgom, hausd} we have $F^{\operatorname{rc}}\subseteq\c^1$.
Set $Q(Z):=P(y,Z)\in F[Z]^{\ne}$. We have~${Q(y')=0}$, so~$y'\in F^{\operatorname{rc}}$, and thus $\der F\subseteq F^{\operatorname{rc}}$. Let now $z\in F^{\operatorname{rc}}$,  and let~$A(Z)$ be the minimum polynomial of $z$  over $F$, say $A=Z^n+A_1Z^{n-1}+\cdots+A_n$ ($A_1,\dots,A_n\in F$, $n\geq 1$). 
With $A^{\der}:=A_1'Z^{n-1}+\cdots+A_n'\in F^{\operatorname{rc}}[Z]$
%and~$A':=S_Q=\frac{\partial Q}{\partial Z}\in F[Z]$,
  we have
$$0\  =\  A(z)'\  =\ A^{\der}(z) + A'(z)\cdot z'$$
with $0\neq  A'(z)\in F^{\operatorname{rc}}$ and so $z'=-A^{\der}(z)/A'(z)\in F^{\operatorname{rc}}$. Hence $F^{\operatorname{rc}}$
is a Hardy field, and so $y$ generates a Hardy field extension~$H\langle y\rangle\subseteq F^{\operatorname{rc}}$ of $H$.
For the rest use Corollary~\ref{cor:Hardy field ext smooth}.
\end{proof}

\noindent
In the proof of the next lemma we encounter an ordered field isomorphism 
$$f\mapsto \tilde f\ :\  E \to \tilde E$$ between
Hausdorff fields $E$ and $\tilde E$. It extends uniquely to an ordered field isomorphism~${E^{\operatorname{rc}}\to\tilde E^{\operatorname{rc}}}$, also denoted by $f\mapsto\tilde f$, 
and to a ring isomorphism
$${Q\mapsto\tilde Q\ :\  E[Y]\to\tilde E[Y]}, \quad \text{ with }\tilde Y=Y.$$  Let $Q\in E[Y]^{\ne}$, and let $y_1< \dots < y_m$ be the distinct zeros of $Q$ in 
$E^{\operatorname{rc}}$.  Then by Corollary~\ref{corb1},  $y_1(t) < \cdots < y_m(t)$ are the distinct real zeros of $Q(t,Y)$, eventually. 
By the isomorphism, $\tilde y_1< \cdots < \tilde y_m$ are the distinct zeros of $\tilde Q$ in 
$\tilde E^{\operatorname{rc}}$, and so~${\tilde y_1(t) < \cdots < \tilde y_m(t)}$ are the distinct real zeros of 
$\tilde Q(t,Y)$, eventually. This has the trivial but useful consequence that, eventually, that is, for all sufficiently large~$t$, 
$$Q(t,Y)=\tilde{Q}(t,Y) \text{ in }\R[Y]\ \Longrightarrow\ y_1(t)=\tilde y_1(t),\dots, y_m(t)=\tilde y_m(t).$$

\begin{lemma}\label{lem:comparable}
Let $u$ be an $H$-hardian germ.  Then $y- u$ is non-oscillating. 
\end{lemma}
\begin{proof}
By Lemma~\ref{lem:order 1 H-hardian}, $y$ is $H$-hardian. Replacing $H$ by $H^{\operatorname{rc}}$
we arrange that~$H$ is real closed.   
Suppose towards a contradiction that $w:=y-u$ oscillates. Then~$w'=y'-u'$ oscillates.  But $y'$ and  $u'$ are $H$-hardian, so $y',u'\notin H$
and for all $h\in H$: $y'>h \Leftrightarrow u'>h$.
This yields an ordered field isomorphism~$H(y')\to H(u')$ over $H$ mapping $y'$ to~$u'$,
which extends uniquely to an ordered field isomorphism
$$f\mapsto \tilde f\ :\   H(y')^{\operatorname{rc}}\to H(u')^{\operatorname{rc}}.$$
Now $P(y,y')=0$ gives $y\in H(y')^{\operatorname{rc}}$, so
$\tilde y\in H(u')^{\operatorname{rc}}\subseteq H\langle u\rangle^{\operatorname{rc}}$. The remarks preceding the lemma applied to $E=H(y')$, $\tilde{E}=H(u')$ and
$Q(Y):=P(Y,y')$ in~$E[Y]$ give that for all sufficiently large $t$ with $y'(t)=u'(t)$ (that is, $w'(t)=0$) we have~$y(t)=\tilde{y}(t)$. 
Now $u, \tilde{y}\in H\langle u\rangle^{\operatorname{rc}}$, so  $\tilde y < u$ or $\tilde y =u$ or $u< \tilde y$. Suppose~$\tilde y <  u$. (The other two cases lead to a contradiction in a similar way.) Then for all sufficiently large $t$ with $w'(t)=0$ we have $y(t)=\tilde{y}(t) <  u(t)$, so $w(t) < 0$, 
contradicting Lemma~\ref{le0} for $f:=w$.  
\end{proof}

\noindent
With these lemmas in place, Theorem~\ref{thm:Bosh order 1} now follows quickly:

\begin{proof}[Proof of Theorem~\ref{thm:Bosh order 1}]
Let $E$ be a $\d$-maximal Hardy field extension of $H$; we show that then $y\in E$. Now
$E$ is real closed by the remarks after Proposition~\ref{prop:Hardy field exts}, and $y-u$ is non-oscillating for all $u\in E$
by Lemma~\ref{lem:comparable}, so $y$ lies in a Hausdorff field extension of $E$ by Lemma~\ref{lem:nonosc, 2},
hence $y$ is $E$-hardian by Lemma~\ref{lem:order 1 H-hardian} with~$E$ in place of $H$, and thus $y\in E$
by $\d$-maximality of~$E$.
\end{proof}

\noindent
As an application of Theorem~\ref{thm:Bosh order 1} we record \cite[Theorem~8.1]{Boshernitzan81}:

\begin{cor}\label{cor:Bosh order 1}
Let $\ell\in\Dx(\Q)$ be such that $\ell>\R$ and $\operatorname{trdeg}\!\big(\R\langle x,\ell \rangle|\R\big) \leq 2$.
Then~$\ell^{\operatorname{inv}}\in\Dx(\Q)$.
\end{cor}
\begin{proof}
By Lemma~\ref{lem:trdegellinv} and the remark preceding it,
 $\ell^{\operatorname{inv}}$ is   $\R(x)$-hardian with
 $$\operatorname{trdeg}\!\big(\R\langle x, \ell^{\operatorname{inv}}\rangle|\R(x)\big) =
 \operatorname{trdeg}\!\big(\R\langle x, \ell^{\operatorname{inv}}\rangle|\R\big)-1=
\operatorname{trdeg}\!\big(\R\langle x,\ell \rangle|\R\big) -1
\leq 1.$$
Now Theorem~\ref{thm:Bosh order 1} with $H:=\R(x)$ and $y:=\ell^{\operatorname{inv}}$ yields
$y\in\Dx(H)=\Dx(\Q)$.
\end{proof}

\section{Application to Filling Holes in Hardy Fields}\label{secfhhf}

\noindent
This section combines the analytic material above with the normalization results of Parts~\ref{part:normalization} and~\ref{part:dents in H-fields}. {\em Throughout
$H$ is a Hardy field with $H\not\subseteq\R$, and $r\in \N^{\ge 1}$}.
Thus~$K:=H[\imag]\subseteq \Calinf[\imag]$ is an $H$-asymptotic extension of $H$.
(Later we impose extra assumptions on~$H$ like being real closed with asymptotic integration.) 
Note that~$v(H^\times)\ne \{0\}$: take $f\in H\setminus \R$; then $f'\ne 0$, and if $f\asymp 1$, then $f'\prec 1$. 

\subsection*{Evaluating differential polynomials at germs} 
Any $Q\in K\{Y\}$ of order $\le r$ can be evaluated at any
germ $y\in \Calr[\imag]$ to give a germ $Q(y)\in \mathcal{C}[\imag]$,
with~$Q(y)\in \mathcal{C}$ for~$Q\in H\{Y\}$ of order $\le r$ and $y\in \Calr$. 
(See the beginning of Section~\ref{sec:smoothness}.)
Here is a variant that we shall need. Let~$\phi\in H^\times$; with~$\der$ denoting the derivation
of $K$, the derivation
of the differential field~$K^\phi$ is then~$\derdelta:= \phi^{-1}\der$. We also let $\derdelta$ denote its ex\-ten\-sion~$f\mapsto \phi^{-1}f'\colon \mathcal{C}^1[\imag] \to \mathcal{C}[\imag]$, which maps $\mathcal{C}^{n+1}[\imag]$ into $\c^n[\imag]$ and
$\c^{n+1}$ into~$\c^n$, for all~$n$. Thus for 
$j\le r$ we have the maps  
$$ \mathcal{C}^r[\imag]\ \xrightarrow{\ \  \derdelta\ \ }\ \mathcal{C}^{r-1}[\imag]\ \xrightarrow{\ \  \derdelta\ \ }\ \cdots\ \xrightarrow{\ \  \derdelta\ \ }\ \mathcal{C}^{r-j+1}[\imag]\ \xrightarrow{\ \  \derdelta\ \ }\  \mathcal{C}^{r-j}[\imag],$$
which by composition yield $\derdelta^j\colon \mathcal{C}^r[\imag]\to \mathcal{C}^{r-j}[\imag]$, mapping $\Calr$ into $\mathcal{C}^{r-j}$. This allows us to define for
$Q\in K^\phi\{Y\}$ of order $\le r$ and $y\in \mathcal{C}^r[\imag]$
the germ $Q(y)\in \mathcal{C}[\imag]$ by
$$Q(y) :=q\big(y, \derdelta(y),\dots, \derdelta^r(y)\big)\quad\text{ where $Q=q\big(Y,\dots, Y^{(r)})\in K^\phi\big[Y,\dots, Y^{(r)}\big]$.}$$ Note that  $H^\phi$ is a  differential subfield of $K^\phi$, and if $Q\in H^\phi\{Y\}$ is of order $\le r$ and~$y\in \Calr$, then $Q(y)\in \mathcal{C}$. 

\begin{lemma}\label{lem:small derivatives of y}
Let $y\in\c^r[\imag]$ and $\fm\in K^\times$. Each of the following conditions implies~$y\in \c^r[\imag]^{\preceq}$:
%$y^{(0)},\dots, y^{(r)}\preceq 1$:
\begin{enumerate}
\item[\textup{(i)}]   $\phi\preceq 1$ and 
$\derdelta^{0}(y),\dots,\derdelta^{r}(y)\preceq 1$; 
\item[\textup{(ii)}] $\fm\preceq 1$ and $y\in \fm\,\c^r[\imag]^{\preceq}$.
%$(y/\fm)^{(0)},\dots,(y/\fm)^{(r)}\preceq 1$. 
\end{enumerate}
Moreover, if $\fm\preceq 1$ and $(y/\fm)^{(0)},\dots,(y/\fm)^{(r)}\prec 1$, then  $y^{(0)},\dots,y^{(r)}\prec 1$.
\end{lemma}
\begin{proof}
For (i), use the smallness of the derivation of $H$ and the transformation formulas in [ADH, 5.7] expressing the iterates of $\der$ in terms of iterates of $\derdelta$. 
For (ii) and the ``moreover'' part, set $y=\fm z$ with
$z=y/\fm$ and use the Product Rule and the smallness of the derivation of $K$.  
\end{proof}

\subsection*{Equations over Hardy fields and over their complexifications}
%{\it In this subsection~$\fm$ ranges over $H^\times$, and $.}\/
Let~$\phi>0$ be active in~$H$.  We recall here from Section~\ref{sec:Hardy fields} how the
 the asymptotic field~$K^\phi=H[\imag]^\phi$ (with derivation $\derdelta=\phi^{-1}\der$) 
 is isomorphic to the asymptotic field $K^\circ:=H^\circ[\imag]$ for a certain Hardy field $H^\circ$: Let
$\ell\in\mathcal C^1$ be such that $\ell'=\phi$; then~$\ell>\R$, $\ell\in \mathcal{C}^{<\infty}$, and $\ell^{\inv}\in \mathcal{C}^{<\infty}$ for the compositional inverse~$\ell^{\inv}$ of $\ell$.
The $\C$-algebra automorphism~$f\mapsto f^\circ:= f\circ \ell^{\inv}$ of $\mathcal{C}[\imag]$ (with inverse $g\mapsto g\circ\ell$) maps~$\Caln[\imag]$ and~$\Caln$
onto themselves, and hence restricts to a $\C$-algebra automorphism of $\Calinf[\imag]$ and~$\Calinf$ mapping
$\Calinf$ onto itself. 
Moreover,
\begin{equation}\label{eq:derc}\tag{$\der, \circ,\derdelta$} (f^\circ)'\ =\  (\phi^\circ)^{-1}(f')^\circ\ =\ \derdelta(f)^\circ \qquad\text{for $f\in \mathcal{C}^1[\imag]$.}
\end{equation}
Thus we have an isomorphism $f\mapsto f^\circ: (\Calinf[\imag])^\phi\to\Calinf[\imag]$ of differential rings, and likewise with $\Calinf$ in place of $\Calinf[\imag]$.
As already pointed out   in Section~\ref{sec:Hardy fields},
$$H^\circ\ :=\ \{h^\circ:h\in H \}\  \subseteq\ \Calinf$$ is a Hardy field,
and $f\mapsto f^\circ$ restricts to  an isomorphism $H^\phi \to H^\circ$ of pre-$H$-fields, and to
an isomorphism $K^\phi\to K^\circ$ of asymptotic fields.
We extend the latter to the isomorphism
$$Q\mapsto Q^\circ\ \colon\ K^\phi\{Y\} \to K^\circ\{Y\}$$
of differential rings given by $Y^\circ=Y$, which restricts to a differential ring isomorphism $H^\phi\{Y\}\to H^\circ\{Y\}$.  Using the identity \eqref{eq:derc} it is routine to check that
for~$Q\in K^\phi\{Y\}$ of order $\le r$ and $y\in \Calr[\imag]$,
$$Q(y)^\circ\  =\  (Q^\circ)(y^\circ). $$
This allows us to translate algebraic differential equations over $K$ 
 into algebraic differential equations over~$K^\circ$: Let $P\in K\{Y\}$ have order $\le r$ and let $y\in \Calr[\imag]$.

\begin{lemma}\label{lem:as equ comp} $P(y)^\circ=P^\phi(y)^\circ=P^{\phi\circ}(y^\circ)$  where $P^{\phi\circ}:=(P^\phi)^\circ\in K^\circ\{Y\}$, hence
$$P(y)=0\ \Longleftrightarrow\ P^{\phi\circ}(y^\circ)=0.$$
\end{lemma}

\noindent
Moreover, $y\prec \fm\ \Longleftrightarrow\ y^\circ \prec \fm^\circ$, for $\fm\in K^\times$, so asymptotic side conditions are automatically taken care of under this ``translation''. Also, if $\phi\preceq 1$ and $y^\circ\in\c^r[\imag]^{\preceq}$, then  $y\in\c^r[\imag]^{\preceq}$, by Lemma~\ref{lem:small derivatives of y}(i) and \eqref{eq:derc}.

\medskip\noindent
{\it In the rest of this section $H\supseteq\R$ is real closed with asymptotic integration.}\/
Then $H$ is an $H$-field, and $K=H[\imag]$ is the algebraic closure of $H$, a $\d$-valued field  with small derivation extending $H$,
constant field $\C$, and value group~$\Gamma:=v(K^\times)=v(H^\times)$.

\subsection*{Slots in Hardy fields and compositional conjugation}
{\it In this sub\-sec\-tion we let~${\phi>0}$ be active in $H$; as in the previous subsection we take $\ell\in\mathcal C^1$ such that~$\ell'=\phi$ and use the superscript $\circ$ accordingly: $f^\circ:= f\circ \ell^{\inv}$ for $f\in \c[\imag]$.}\/

Let  $(P,\fm,\hat a)$ be a slot in $K$ of order $r$, so $\hat a\notin K$ is an element of an immediate asymptotic extension $\hat K$ of~$K$ with~$P\in Z(K,\hat a)$ and $\hat a\prec\fm$. 
We associate to~$(P,\fm,\hat a)$ a slot in $K^\circ$ as follows:
choose an immediate asymptotic extension $\hat K^\circ$ of $K^\circ$ and an
isomorphism $\hat f\mapsto \hat f^\circ\colon\hat K^\phi\to \hat K^\circ$ of asymptotic fields
extending the isomorphism $f\mapsto f^\circ\colon K^\phi\to K^\circ$.
Then $(P^{\phi\circ},\fm^\circ,\hat a^\circ)$ is a slot in $K^\circ$ of the same complexity as~$(P,\fm,\hat a)$. The
equivalence class of   $(P^{\phi\circ},\fm^\circ,\hat a^\circ)$ does not depend on the choice of
$\hat K^\circ$ and the isomorphism $\hat K^\phi\to \hat K^\circ$. 
If $(P,\fm,\hat a)$ is a hole in $K$, then~$(P^{\phi\circ},\fm^\circ,\hat a^\circ)$ is a hole in~$K^\circ$, and likewise with ``minimal hole'' in place of ``hole''. Moreover, by Lemmas~\ref{lem:v(Aphi)}, \ref{lem:normality comp conj}, and \ref{lem:normality comp conj, strong}:

\begin{lemma} \label{lem:Pphicirc, 1}
If $(P,\fm,\hat a)$ is $Z$-minimal, then so is $(P^{\phi\circ},\fm^\circ,\hat a^\circ)$, and likewise with ``quasi-linear'' and ``special'' in place of ``$Z$-minimal''. 
If $(P,\fm,\hat a)$ is steep and~$\phi\preceq 1$, then $(P^{\phi\circ},\fm^\circ,\hat a^\circ)$ is steep,
and likewise with   ``deep'', ``normal'',  and ``strictly normal'' in place of ``steep''. 
%If $H$ is Liouville closed and $(P,\fm,\hat a)$ is
%ultimate, then so is $(P^{\phi\circ},\fm^\circ,\hat a^\circ)$. (But requires $\hat K=\hat{H}[\imag]$.
\end{lemma}

\noindent
Next, let $(P,\fm,\hat a)$ be a slot in $H$ of order $r$, so $\hat a\notin H$ is an element of an immediate asymptotic extension $\hat H$ of~$H$ with~$P\in Z(H,\hat a)$ and $\hat a\prec\fm$. 
We associate to~$(P,\fm,\hat a)$ a slot in $H^\circ$ as follows:
choose an immediate asymptotic extension $\hat H^\circ$ of $H^\circ$ and an
isomorphism $\hat f\mapsto \hat f^\circ\colon\hat H^\phi\to \hat H^\circ$ of asymptotic fields
extending the isomorphism~$f\mapsto f^\circ\colon H^\phi\to H^\circ$.
Then $(P^{\phi\circ},\fm^\circ,\hat a^\circ)$ is a slot in $H^\circ$ of the same complexity as~$(P,\fm,\hat a)$. The
equivalence class of   $(P^{\phi\circ},\fm^\circ,\hat a^\circ)$ does not depend on the choice of
$\hat H^\circ$ and the isomorphism $\hat H^\phi\to \hat H^\circ$. 
If $(P,\fm,\hat a)$ is a hole in $H$, then~$(P^{\phi\circ},\fm^\circ,\hat a^\circ)$ is a hole in~$H^\circ$, and likewise with ``minimal hole'' in place of ``hole''.
Lemma~\ref{lem:Pphicirc, 1} goes through in this setting. Also, recalling Lemma~\ref{lemlioucomp}, if~$H$ is Liouville closed and~$(P,\fm,\hat a)$
is ultimate, then $(P^{\phi\circ},\fm^\circ,\hat a^\circ)$ is ultimate. 

\medskip
\noindent
Moreover, by Lemmas~\ref{lem:split-normal comp conj} and~\ref{lem:strongly split-normal compconj}, and Corollaries~\ref{cor:repulsive-normal comp conj} and~\ref{cor:strongly repulsive-normal compconj}:
% various preservation results in Sections~\ref{sec:split-normal holes} and \ref{sec:repulsive-normal}:

{\samepage
\begin{lemma} \label{lem:Pphicirc, 2}
\mbox{}
\begin{enumerate}
\item[\textup{(i)}]
If $\phi\preceq 1$ and $(P,\fm,\hat a)$ is split-normal, then  $(P^{\phi\circ},\fm^\circ,\hat a^\circ)$ is split-normal;
likewise with ``split-normal'' replaced by   ``\textup{(}almost\textup{)} strongly split-nor\-mal''.
\item[\textup{(ii)}]
If $\phi\prec 1$ and  $(P,\fm,\hat a)$  is $Z$-minimal, deep, and repulsive-normal, then $(P^{\phi\circ},\fm^\circ,\hat a^\circ)$ is repulsive-normal; likewise with ``repulsive-normal'' replaced by
 ``\textup{(}almost\textup{)} strongly re\-pul\-sive-nor\-mal''. 
\end{enumerate}
\end{lemma} 
}

\subsection*{Reformulations}
We reformulate here some results of Sections~\ref{sec:split-normal over Hardy fields} and~\ref{sec:smoothness}   to facilitate their use.
As in Section~\ref{sec:span} we set for $\fv\in K^\times$, $\fv\prec 1$: 
$$\Delta(\fv)\ :=\ \big\{ \gamma\in\Gamma: \gamma=o(v\fv)\big\},$$ 
a proper convex subgroup of $\Gamma$. 
In the next lemma, $P\in  K\{Y\}$ has order $r$ and~$P=Q-R$, where $Q,R\in K\{Y\}$ and $Q$ is homogeneous of degree~$1$ and order~$r$.
We set~$w:=\wt(P)$, so $w\ge r\ge 1$.
%$L:=L_Q\in H[\imag][\der]$, so $L$ has order $r$, and $w:=\wt(P)$, so $w\geq r$.

\begin{lemma}\label{prop:as equ 1}  
Suppose that $L_Q$ splits strongly over $K$, $\fv(L_Q)\prec^\flat 1$, and 
$$R\prec_\Delta \fv(L_Q)^{w+1}Q, \quad \Delta := \Delta\big(\fv(L_Q)\big).$$ Then $P(y) =0$ and $y',\dots, y^{(r)} \preceq 1$
for some $y\prec \fv(L_Q)^w$ in~$\Calinf[\imag]$.  Moreover: \begin{enumerate}
\item[\textup{(i)}] if $P,Q\in H\{Y\}$, then there is such~$y$ in $\Calinf$;
\item[\textup{(ii)}] if $H\subseteq \Ginf$, then for any  $y\in\c^r[\imag]^{\preceq}$ with $P(y)=0$ we have $y\in \Ginf[\imag]$; likewise with $\c^\omega$ in place of $\c^{\infty}$.
\end{enumerate} 
\end{lemma}
\begin{proof}
Set $\fv:= \abs{\fv(L_Q)}\in H^>$, so $\fv\asymp\fv(L_Q)$. Take $f\in K^{\times}$  such that $A:=f^{-1}L_Q$ is monic; then $\fv(A)=\fv(L_Q)\asymp\fv$ and
$f^{-1}R \prec_\Delta f^{-1}\fv^{w+1}Q \asymp \fv^{w}$.
We have~$A=(\der-\phi_1)\cdots(\der-\phi_r)$ where $\phi_j\in K$ and
$\Re\phi_j\succeq\fv^\dagger\succeq 1$ for $j=1,\dots,r$ by the strong splitting assumption. Also $\phi_1,\dots,\phi_r\preceq\fv^{-1}$ by Corollary~\ref{cor:bound on linear factors}. 
The claims now follow from various results in Section~\ref{sec:split-normal over Hardy fields}  applied to the equation~${A(y)=f^{-1}R(y)}$, $y\prec 1$ in the role of \eqref{eq:ADE}, using also Corollary~\ref{cor:ADE smooth}. %Lemma~\ref{bdua}, Theorem~\ref{thm:fix}, and Corollary~\ref{cor:fix}
\end{proof}

%\noindent
%If $H\subseteq\Ginf$, then this lemma goes through with $\Calinf[\imag]$ and $\Calinf$ replaced by 
%$\mathcal{C}^{\infty}[\imag]$ and $\mathcal{C}^{\infty}$, respectively;  this uses the  $\c^\infty$-part of Corollary~\ref{cor:ADE smooth}. 
%Likewise with $\c^\omega$ in place of $\c^\infty$.

\begin{lemma}\label{dentsolver} 
%Let  $(P,\fn,\hat h)$ be a strongly split-normal dent in $H$ of order~$r$.  
Let  $(P,\fn,\hat h)$ be a slot in $H$ of order~$r$ and let $\phi$ be active in~$H$, ${0<\phi\preceq 1}$, such that
$(P^\phi,\fn,\hat h)$ is strongly split-normal. Then for some~$y$ in~$\Calinf$,
% satisfies 
$$P(y)\ =\ 0,\quad y\ \prec\ \fn,\quad y\in \fn\,(\c^r)^{\preceq}.$$
%(y/\fn)^{(j)}\ \preceq\ 1\ 
%\text{ for $j=1,\dots,r$.}$$
%If $\fn\preceq 1$, then $y',\dots, y^{(r)}\preceq 1$ for such $y$.
If $H\subseteq\c^\infty$, then there exists such $y$ in~$\c^\infty$, and likewise with $\c^\omega$ in place of $\c^\infty$.  
\end{lemma}
\begin{proof}
First we consider the case $\phi=1$.
Replace $(P,\fn,\hat h)$ by $(P_{\times\fn},1,\hat h/\fn)$
to arrange $\fn=1$.  
Then $L_{P}$ has order~$r$, $\fv:=\fv(L_{P})\prec^{\flat} 1$, and $P =  Q-R$ where~$Q,R\in H\{Y\}$,
$Q$ is homogeneous of degree $1$ and order~$r$, $L_{Q}\in H[\der]$ splits strongly over~$K$, 
and~$R\prec_\Delta \fv^{w+1} P_1$, where~$\Delta:=\Delta(\fv)$ and $w:= \wt(P)$. 
Now $P_1=Q-R_1$, so~$\fv\sim \fv(L_{Q})$ by Lem\-ma~\ref{lem:fv of perturbed op}(ii), and thus~$\Delta=\Delta\big(\fv(L_{Q})\big)$.
Lemma~\ref{prop:as equ 1} gives $y$ in $\Calinf$ such that~$y \prec \fv^w\prec 1$, 
$P(y)=0$, and
$y^{(j)}\preceq 1$ for $j=1,\dots,r$.
Then $y$ has for $\fn=1$ the properties displayed in the lemma.  
%For the last statement,  use the remark after the proof of Lemma~\ref{prop:as equ 1}.

Now suppose $\phi$ is  arbitrary. Employing $(\phantom{-})^\circ$ as explained earlier in this section, the
slot $(P^{\phi\circ},\fn^\circ,\hat h^\circ)$ in the Hardy field~$H^\circ$ is strongly split-normal, hence by the case $\phi=1$ we have $z\in\Calinf$ with 
$P^{\phi\circ}(z) = 0$, $z\prec\fn^\circ$, and $(z/\fn^\circ)^{(j)}\preceq 1$
for~${j=1,\dots,r}$.
Take $y\in\Calinf$ with $y^\circ=z$. Then
$P(y) = 0$, $y \prec \fn$, and $y\in \fn\,(\c^r)^{\preceq}$ by Lemma~\ref{lem:as equ comp} and a subsequent remark.  
%$\derdelta^{j}(y/\fn) \preceq 1$ for $j=1,\dots,r$ with $\derdelta:=\phi^{-1}\der$. 
%Lem\-ma~\ref{lem:small derivatives of y} yields $(y/\fn)^{(j)}\preceq 1$ for~$j=1,\dots,r$. 
Moreover, if $\phi,z\in\Ginf$, then $y\in\Ginf$, and likewise with $\Gom$ in place of~$\Ginf$.
\end{proof}

\noindent
In the next ``complex'' version, $(P, \fm, \hat a)$ is a slot in $K$ of order $r$ with $\fm\in H^\times$. 

\begin{lemma}\label{lemfillhole} 
Let $\phi$ be active in~$H$, ${0<\phi\preceq 1}$, such that
 the slot $(P^\phi,\fm,\hat a)$ in~$K^\phi$ is  strictly normal, and its linear part splits strongly
over $K^\phi$. Then for some~$y\in \Calinf[\imag]$ we have
$$P(y)\ =\ 0,\quad y\ \prec\ \fm,\quad y\in \fm\,\c^r[\imag]^{\preceq}.$$
%(y/\fm)^{(j)}\ \preceq\ 1\ \text{ for $j=1,\dots,r$.}$$
%If  $\fm\preceq 1$, then $y',\dots, y^{(r)}\preceq 1$ for such $y$.
If $H\subseteq\c^\infty$, then there is such $y$ in $\c^\infty[\imag]$.
If $H\subseteq\c^\omega$, then there is such $y$ in $\c^\omega[\imag]$. 
\end{lemma}
\begin{proof}
Consider first the case $\phi=1$. Replacing~$(P,\fm,\hat a)$  by $(P_{\times\fm},1,\hat a/\fm)$
we arrange $\fm=1$.
Set $L:=L_{P}\in K[\der]$, $Q:=P_1$, and $R:=P-Q$. 
Since~$(P,1,\hat a)$ is strictly normal, we have $\order(L)=r$,  $\fv:=\fv(L)\prec^\flat 1$,
and $R\prec_\Delta \fv^{w+1}Q$  where~$\Delta:=\Delta(\fv)$,  $w := \wt(P)$. As $L$ splits  strongly over~$K$,   Lemma~\ref{prop:as equ 1} gives~$y$ in $\Calinf[\imag]$ such that
$P(y)=0$, $y \prec \fv^w\prec 1$, and $y^{(j)}\preceq 1$
for~$j=1,\dots,r$. For the last part of the lemma, use the last part of  Lemma~\ref{prop:as equ 1}.  
The general case reduces to this special case as in the proof of Lemma~\ref{dentsolver}. 
\end{proof}

\subsection*{Finding germs in holes} {\em In this subsection $\hat H$ is an immediate asymptotic extension of $H$.}\/
This fits into the setting of Section~\ref{sec:split-normal holes} on split-normal slots: 
$K=H[\imag]$ and~$\hat H$ have $H$ as a common asymptotic subfield and~$\hat{K}:=\hat{H}[\imag]$ as a common asymptotic extension, $\hat H$ is an $H$-field, and  $\hat K$ is $\d$-valued.
{\em Assume also that $H$ is $\upo$-free.}\/  Thus $K$ is $\upo$-free by [ADH, 11.7.23].
{\em Let $(P,\fm,\hat a)$ with~$\fm\in H^\times$ and~$\hat a\in 
\hat K\setminus K$ be a minimal hole in $K$ of order~${r\geq 1}$. 
Take $\hat b,\hat c\in\hat H$ so that~$\hat a=\hat b+\hat c\,\imag$.}\/

\begin{prop}\label{propdeg>1}
Suppose  $\deg P>1$. 
Then for some $y\in \Calinf[\imag]$ we have
$$P(y)\ =\ 0,\quad y\ \prec\ \fm,\quad y\in \fm\,\c^r[\imag]^{\preceq}.$$
%(y/\fm)^{(j)}\ \preceq\ 1\ \text{ for $j=1,\dots,r$.}$$
If $\fm\preceq 1$, then $y\prec \fm$ and $y\in \c^r[\imag]^{\preceq}$ for such~$y$.
Moreover, if   $H\subseteq\c^\infty$, then we can take such $y$ in $\c^\infty[\imag]$, and 
if $H\subseteq\c^\omega$, then we can take such~$y$ in $\c^\omega[\imag]$. 
\end{prop}
\begin{proof} Lemma~\ref{lem:achieve strong splitting} gives
a refinement $(P_{+a},\fn,\hat a-a)$  of $(P,\fm,\hat a)$ with $\fn\in H^\times$ and an active $\phi$ in $H$ with $0<\phi\preceq 1$ such that 
the hole $(P^\phi_{+a},\fn,\hat a-a)$ in $K^\phi$ is  strictly normal  and its  linear part splits 
strongly over $K^\phi$. 
Lemma~\ref{lemfillhole} applied to $(P_{+a},\fn,\hat a-a)$ in place of $(P,\fm,\hat a)$
yields $z\in\Calinf[\imag]$ with $P_{+a}(z)=0$, $z\prec\fn$ and~$(z/\fn)^{(j)}\preceq 1$ for $j=1,\dots,r$.
Lem\-ma~\ref{lem:small derivatives of y}(ii) applied to~$z/\fm$,~$\fn/\fm$ in place of~$y$,~$\fm$, respectively, yields~$(z/\fm)^{(j)}\preceq 1$ for $j=0,\dots,r$. Also, $a\prec\fm$ (in $K$), hence~$(a/\fm)^{(j)}\prec 1$ for~$j=0,\dots,r$. Set $y:=a+z$; then $P(y)=0$, $y\prec\fm$, and
$(y/\fm)^{(j)}\preceq 1$ for~$j=1,\dots,r$.
For the rest use Lem\-ma~\ref{lem:small derivatives of y}(ii) and the last statement in Lemma~\ref{lemfillhole}.
\end{proof}

\noindent
Next we treat the linear case:

\begin{prop} \label{prop:deg 1 analytic} 
Suppose $\deg P=1$. Then for some~$y\in \Calinf[\imag]$ we have
$$P(y)\ =\ 0,\quad y\ \prec\ \fm,\quad (y/\fm)'\ \preceq\ 1.$$
If  $\fm\preceq 1$, then $y\prec 1$ and $y'\preceq 1$ for each such~$y$. 
Moreover, if   $H\subseteq\c^\infty$, then we can take such $y$ in $\c^\infty[\imag]$, and 
if $H\subseteq\c^\omega$, then we can take such~$y$ in $\c^\omega[\imag]$.
%If $P\in H\{Y\}$, then we can take such~$y$ in~$\Calinf$. \marginpar{real version needed?}
\end{prop}
\begin{proof} 
We have $r=1$ by Corollary~\ref{cor:minhole deg 1}.
If $\der K=K$ and $\I(K)\subseteq K^\dagger$, then Lemma~\ref{lem:achieve strong splitting, d=1} applies, and we can argue as in the proof of Proposition~\ref{propdeg>1},
using this lemma  instead of Lemma~\ref{lem:achieve strong splitting}.
We reduce the general case to this special case as follows: 
Set $H_1:=\operatorname{D}(H)$; then $H_1$ is an $\upo$-free Hardy field by Theorem~\ref{thm:ADH 13.6.1}, and~$K_1:=H_1[\imag]$ satisfies $\der K_1=K_1$ and $\I(K_1)\subseteq K_1^\dagger$, by Corollary~\ref{cor:cos sin infinitesimal}. Moreover, by Corollary~\ref{cor:Hardy field ext smooth}, if $H\subseteq\c^\infty$, then $H_1\subseteq\c^\infty$, and likewise with~$\c^\omega$ in place
of~$\c^\infty$.
The newtonization $\hat H_1$ of $H_1$ is an immediate asymptotic extension of $H_1$, and~$\hat K_1:=\hat H_1[\imag]$ is   newtonian [ADH, 14.5.7]. 
Corollary~\ref{cor:find zero of P, 2} gives an embedding~$K\langle \hat a\rangle \to \hat K_1$ over $K$; let $\hat a_1$ be the image of $\hat a$ under this embedding.
If  $\hat a_1\in K_1$, then we are done by taking $y:=\hat a_1$, so we may assume $\hat a_1\notin K_1$.   Then $(P,\fm,\hat a_1)$ is a minimal hole in~$K_1$, and the above applies with $H$, $K$, $\hat a$ replaced by $H_1$, $K_1$, $\hat a_1$, respectively.
\end{proof}

\noindent
We can improve on these results in a useful way:

\begin{cor}\label{improap} 
Suppose
$\hat a\sim a\in K$. Then  for some $y\in \Calinf[\imag]$ we have
$$P(y)\ =\ 0,\quad y\ \sim\ a,\quad (y/a)^{(j)}\ \prec\ 1\ 
\text{ for $j=1,\dots,r$.}$$
If   $H\subseteq\c^\infty$, then there is such $y$ in $\c^\infty[\imag]$.
If $H\subseteq\c^\omega$, then there is such~$y$ in $\c^\omega[\imag]$.
\end{cor}
\begin{proof} Take $a_1\in K$ and $\fn\in H^\times$ with 
 $\fn\asymp \hat a-a \sim a_1$, and set~$b:=a+a_1$.
Then~$(P_{+b}, \fn, \hat a-b)$ is a refinement of $(P,\fm,\hat a)$.
Propositions~\ref{propdeg>1} and~\ref{prop:deg 1 analytic} give $z\in \Calinf[\imag]$
with $P(b+z)=0$, $z\prec \fn$ and~$(z/\fn)^{(j)}\preceq 1$ for~$j=1,\dots,r$. 
We have~$(a_1/a)^{(j)}\prec 1$ for~$j=0,\dots,r$, since $K$ has small derivation.
Likewise, $(\fn/a)^{(j)}\prec 1$ for $j=0,\dots,r$, and hence
$(z/a)^{(j)}\prec 1$ for $j=0,\dots,r$, by~$z/a=(z/\fn)\cdot (\fn/a)$ and the Product Rule.
So $y:= b+z$ has the desired property. The rest follows from the ``moreover'' parts of these propositions.
\end{proof}

\begin{remarkNumbered}\label{rem:improap}  Suppose we replace our standing assumption that $H$ is $\upo$-free and~$(P, \fm, \hat a)$ is a minimal hole in $K$ by the assumption that $H$ is $\upl$-free, $\der K=K$, $\I(K)\subseteq K^\dagger$, and $(P, \fm, \hat a)$ is a slot in $K$ of order and degree $1$.
Then Proposition~\ref{prop:deg 1 analytic} and Corollary~\ref{improap} go   through by
the remark following Lemma~\ref{lem:achieve strong splitting, d=1}.
\end{remarkNumbered}

\noindent
Now also drawing upon Theorem~\ref{thm:strongly split-normal}, we arrive at the main result of this section:

\begin{cor}\label{6main}
Suppose $H$ is $1$-linearly newtonian. Then one of the following two conditions is satisfied: 
\begin{enumerate}
\item[\textup{(i)}] $\hat b\notin H$, and there are $Q\in Z(H,\hat b)$ of minimal complexity  
and  $y\in\Calinf$ such that $Q(y)=0$ and $y\prec\fm$; 
\item[\textup{(ii)}] $\hat c\notin H$, and there are $R\in Z(H,\hat c)$  of minimal complexity  
and  $y\in\Calinf$ such that $R(y)=0$ and $y\prec\fm$.
\end{enumerate}
If $H\subseteq\Ginf$, then there is such $y$ in $\Ginf$, and likewise with $\Ginf$ replaced by $\Gom$. 
\end{cor} 
\begin{proof} 
Suppose $\deg P>1$, or $\hat b\notin H$ and  $Z(H,\hat b)$ has an element of order~$1$, 
or~$\hat c\notin H$ and  $Z(H,\hat c)$ has an element of order~$1$. 
Let $\phi$ range over active elements of $H$ with~$0<\phi\preceq 1$.
By the ``moreover'' part of Theorem~\ref{thm:strongly split-normal},
one of the following holds:
\begin{enumerate}
\item  $\hat b\notin H$ and there exist $\phi$ and a $Z$-minimal slot $(Q,\fm,\hat b)$ in $H$ with a refinement~${(Q_{+b},\fn,\hat b-b)}$ such that ${(Q^\phi_{+b},\fn,\hat b-b)}$ is   strongly split-normal; 
\item  $\hat c\notin H$ and there exist $\phi$ and a $Z$-minimal slot $(R,\fm,\hat c)$ in $H$ with a refinement~${(R_{+c},\fn,\hat c-c)}$ such that $(R^\phi_{+c},\fn,\hat c-c)$ is     strongly split-normal. 
\end{enumerate}
Suppose $\hat b\notin H$ and  $\phi, Q, b$  are as in (1); then
 Lemma~\ref{dentsolver} applied to~$(Q_{+b},\fn,\hat b-b)$ in place of $(P,\fn,\hat h)$ 
 yields $z\in\Calinf$ with  $Q_{+b}(z)=0$, $z\prec\fn$; hence~$Q(y)=0$, $y\prec\fm$   for $y:=b+z$, so (i) holds.
 Similarly,~(2) implies~(ii).

Suppose now that $\deg P =1$, that if $\hat b\notin H$, then $Z(H,\hat b)$ has no element of order~$1$, and
that if $\hat c\notin H$, then $Z(H,\hat c)$ has no element of order~$1$.
Since $\deg P =1$, Proposition~\ref{prop:deg 1 analytic} gives $z\in\Calinf[\imag]$ such that $P(z)=0$ and $z\prec\fm$. Recall also that~$P$ has order $1$ by Corollary~\ref{cor:minhole deg 1}.
Consider now the case $\hat b\notin H$. In view of~$P(\hat a)=0$ and
$P(z)=0$ we obtain from Example~\ref{ex:lclm compl conj}  and  Remark~\ref{rem:lclm compl conj} 
a $Q\in H\{Y\}$ of degree~$1$ and  order~$1$ or $2$
such that~$Q(\hat b)=0$ and~$Q(y)=0$ for $y:=\Re z\prec\fm$.
But~$Z(H,\hat b)$ has no element of order~$1$, so $\order Q=2$ and~$Q\in Z(H,\hat b)$ has minimal complexity. Thus (i) holds. 
If $\hat c\notin H$, then the same reasoning shows that~(ii) holds.
\end{proof}

\noindent 
Is $y$ as in~(i) or~(ii) of Corollary~\ref{6main} hardian over $H$?  At this stage we cannot claim this.  
%For example, let~$\hat b$,~$Q$,~$y$ be as in (i); if~$y$ 
%indeed lies in an  immediate Hardy field
%extension of $H$ such that~$(Q,y,\fm)$ is a dent in $H$ equivalent to $(Q,\hat b,\fm)$, then for each~$b\in H$ and~$\fn\in H^\times$ with~$\hat b-b\prec\fn$ we must have
%$y-b\prec\fn$ and thus~$\big((y-b)/\fn\big){}^{(n)}\prec 1$ for each~$n$ (since Hardy fields have small derivation). 
In the next section we introduce weights and their corresponding norms as a more refined tool. This will allow us to
obtain Corollary~\ref{cor:approx y} as a key approximation result for later use. 

\section{Weights}\label{sec:weights}

\noindent
In this section we prove Proposition~\ref{prop:notorious 3.6} to  strengthen
Lemma~\ref{lem:close}. 
This uses the material on repulsive-normal slots from Section~\ref{sec:repulsive-normal}, but
we also need more refined norms for differentiable functions, to which we turn now.

\subsection*{Weighted spaces of differentiable functions}
{\it In this subsection we fix $r\in \N$ and a {\em weight}
function $\w\in\c_a[\imag]^\times$.}\/ 
For $f\in \Car[\imag]$ we set \label{p:absawt}
%We  consider the seminorm $$f\mapsto 
$$\dabs{f}_{a;r}^\w\ :=\ 
\max\big\{ \dabs{\w^{-1}f}_a, \dabs{\w^{-1}f'}_a, \dots, \dabs{\w^{-1}f^{(r)}}_a\big\} \ \in\  [0,+\infty],$$
and $\dabs{f}_{a}^\w:=\dabs{f}_{a;0}^\w$ for $f\in \c_a[\imag]$. 
\noindent
%(For $f\in \Car[\imag]^{\b}$ and $\w\in\R^\times\subseteq\c_a[\imag]^\times$ this is to be distinguished from the short-hand notation
%$\dabs{f}_{a;r}^{\w}=(\dabs{f}_{a;r})^{\w}$.) \marginpar{added this caveat, relevant for the statement of \ref{lem:weighted bd} and other places below}
Then
$$\Car[\imag]^\w\ :=\ \big\{ f\in \Car[\imag]:\, \dabs{f}_{a;r}^\w < +\infty \big\}$$
is a $\C$-linear subspace of 
$$\c_a[\imag]^\w\ :=\ \Caz[\imag]^\w\ =\ 
\w\,\c_a[\imag]^{\b}\ =\ \big\{ f\in \c_a[\imag]:\, f\preceq \w\big\}.$$ 
Below we consider the $\C$-linear space~$\Car[\imag]^\w$ to be equipped with the norm 
$$f\mapsto \dabs{f}_{a;r}^\w.$$  
Recall from Section~\ref{sec:IHF} the convention $b\cdot\infty=\infty\cdot b=\infty$ for $b\in [0,\infty]$. Note that 
\begin{equation}\label{eq:weighted norm, 0}
\dabs{fg}_{a;r}^\w\ \leq\ 2^r \dabs{f}_{a;r}\,\dabs{g}_{a;r}^\w\quad\text{ for $f,g\in \Car[\imag]$,}
\end{equation}
so $\Car[\imag]^\w$ is a $\Car[\imag]^{\b}$-submodule of $\Car[\imag]$. 
Note also that $\dabs{1}_{a;r}^\w=\dabs{\w^{-1}}_{a}$, hence
$$\dabs{f}_{a;r}^\w\ \leq\ 2^r \dabs{f}_{a;r}\,\dabs{\w^{-1}}_{a}\quad\text{ for $f\in \Car[\imag]$}$$
and
$$\w^{-1}\in\c_a[\imag]^{\b} \quad\Longleftrightarrow\quad 1\in \Car[\imag]^{\w} \quad\Longleftrightarrow\quad \Car[\imag]^{\b} \subseteq \Car[\imag]^{\w}.$$
We have
\begin{equation}\label{eq:weighted norm, 1}
\dabs{f}_{a;r}\ \leq\  \dabs{f}^\w_{a;r}\,\dabs{\w}_{a}\quad\text{ for $f\in \Car[\imag]$,}
\end{equation}
and thus
\begin{equation}\label{eq:weighted norm, 2}
\w\in\c_a[\imag]^{\b} \quad\Longleftrightarrow\quad \Car[\imag]^{\b} \subseteq \Car[\imag]^{\w^{-1}} \quad\Longrightarrow\quad \Car[\imag]^{\w}\subseteq \Car[\imag]^{\b}.
\end{equation}
Hence if $\w,\w^{-1}\in \c_a[\imag]^{\b}$, then $\Car[\imag]^{\b} = \Car[\imag]^{\w}$, and the norms
 $\dabs{\,\cdot\,}_{a;r}^\w$ and  $\dabs{\,\cdot\,}_{a;r}$ on this $\C$-linear space are equivalent. 
 (In later use, $\w\in \c_a[\imag]^{\b}$, $\w^{-1}\notin  \c_a[\imag]^{\b}$.)
If~$\w\in \c_a[\imag]^{\b}$, then $\Car[\imag]^{\w}$ is an ideal of the commutative ring~$\Car[\imag]^{\b}$. From \eqref{eq:weighted norm, 0} and~\eqref{eq:weighted norm, 1} we obtain
$$\dabs{fg}_{a;r}^\w\ \leq\ 2^r \,\dabs{\w}_{a} \, \dabs{f}_{a;r}^\w\,\dabs{g}_{a;r}^\w\quad\text{ for $f,g\in \Car[\imag]$.}$$
For $f\in \Carm[\imag]^\w$ we have $\dabs{f}_{a;r}^\w, \dabs{f'}_{a;r}^\w \leq \dabs{f}_{a;r+1}^\w$.
From \eqref{eq:weighted norm, 1} and  \eqref{eq:weighted norm, 2}:

\begin{lemma}\label{lem:w-conv => b-conv}
Suppose $\w\in\c_a[\imag]^{\b}$ $($so $\Car[\imag]^{\w}\subseteq \Car[\imag]^{\b})$ and $f\in \Car[\imag]^{\w}$. If $(f_n)$ is a sequence in~$\Car[\imag]^\w$ and
$f_n\to f$ in $\Car[\imag]^\w$, then also $f_n\to f$ in $\Car[\imag]^{\b}$.
\end{lemma}

\noindent
This is used to show:

\begin{lemma}\label{lem:w-conv}
Suppose $\w\in\c_a[\imag]^{\b}$. Then
the $\C$-linear space $\Car[\imag]^\w$ equipped with the norm
$\dabs{\,\cdot\,}_{a;r}^\w$ is complete.
\end{lemma}
\begin{proof}
We proceed by induction on $r$.
Let $(f_n)$ be a cauchy sequence in the normed space~$\c_a[\imag]^\w$. Then  
the sequence $(\w^{-1}f_n)$ in the Banach space 
$\Caz[\imag]^{\b}$ is cauchy, hence has a limit $g\in \c_a[\imag]^{\b}$, so
with $f:=\w g\in\c_a[\imag]^\w$ we have
$\w^{-1}f_n\to \w^{-1}f$ in~$\c_a[\imag]^{\b}$ and hence $f_n\to f$ in~$\c_a[\imag]^\w$. 
Thus the lemma holds for $r=0$. Suppose the lemma holds for a certain value of
$r$, and let $(f_n)$ be a cauchy sequence in $\Carm[\imag]^\w$.
Then~$(f_n')$ is a cauchy sequence in $\Car[\imag]^\w$ and hence has a limit
$g\in\Car[\imag]^\w$, by inductive hypothesis.
By Lemma~\ref{lem:w-conv => b-conv},    $f_n'\to g$ in
$\c_a[\imag]^{\b}$.
Now~$(f_n)$ is also a cauchy sequence in $\c_a[\imag]^\w$, hence has a limit
$f\in\c_a[\imag]^\w$ (by the case $r=0$), and
by Lemma~\ref{lem:w-conv => b-conv} again, $f_n\to f$ in $\c_a[\imag]^{\b}$.
Thus 
$f$ is differentiable and $f'=g$ by~\cite[(8.6.4)]{Dieudonne}.
This yields~$f_n\to f$ in~$\Carm[\imag]^{\w}$.
\end{proof}

\begin{lemma}\label{lem:Qnk}
Suppose $\w\in \Car[\imag]^{\b}$. If $f\in\Car[\imag]$ and $f^{(k)}\preceq \w^{r-k+1}$ for $k=0,\dots,r$, then $f\w^{-1} \in \Car[\imag]^{\b}$. \textup{(}Thus $\Car[\imag]^{\w^{r+1}}\subseteq \w\,\Car[\imag]^{\b}$.\textup{)}
\end{lemma}
\begin{proof} 
Let $Q^n_k\in\Q\{X\}$ ($0\leq k\leq n$) be as in Lem\-ma~\ref{lem:def of Qnk}. 
Although $\Car[\imag]$ is not closed under taking derivatives, the proof of that lemma and the computation leading to Corollary~\ref{cor:Qnk}  does give for $f\in \Car[\imag]$ and $n\le r$:
$$(f\w^{-1})^{(n)}\ =\ \sum_{k=0}^n Q^n_k(\w) f^{(k)}\w^{k-n-1}.$$
Now use that $Q^n_k(\w)\preceq 1$ for $n\le r$ and $k=0,\dots,n$. 
\end{proof}

\noindent
Next we generalize the inequality \eqref{eq:weighted norm, 0}:

\begin{lemma}\label{lem:weighted norm}
Let $f_1,\dots,f_{m-1},g\in \Car[\imag]$, $m\geq 1$; then
$$\dabs{f_1\cdots f_{m-1}g}_{a;r}^\w\ \leq\ m^r \,\dabs{f_1}_{a;r}\cdots\dabs{f_{m-1}}_{a;r} \,\dabs{g}_{a;r}^\w.$$
\end{lemma}
\begin{proof}
Use the generalized Product Rule [ADH, p.~199] and the well-known identity 
$\sum \frac{n!}{i_1!\cdots i_m!}=m^n$ with the sum over all $(i_1,\dots,i_m)\in\N^m$ with $i_1+\cdots+i_m=n$.
\end{proof}

\noindent
With $\i$ ranging over $\N^{1+r}$,  let
$P=\sum_{\i} P_{\i} Y^{\i}$ (all $P_{\i}\in \c_a[\imag]$) be a polynomial in~$\c_a[\imag]\big[Y,Y',\dots,Y^{(r)}\big]$;
for $f\in \Car[\imag]$ we have  $P(f)=\sum_{\i} P_{\i} f^{\i}\in\c_a[\imag]$.
(See also the beginning of Section~\ref{sec:split-normal over Hardy fields}.)
We set
$$\dabs{P}_a\ :=\ \max_{\i} \, \dabs{P_{\i}}_a \in [0,\infty].$$
{\em In the rest of this subsection we assume $\dabs{P}_a<\infty$}, that is,  $P\in\c_a[\imag]^{\b}\big[Y,\dots,Y^{(r)}\big]$. Hence if~$\w\in \c_a[\imag]^{\b}$, $P(0)\in \c_a[\imag]^\w$, and $f\in\Car[\imag]^\w$, then $P(f)\in \c_a[\imag]^\w$.  %Lemma~\ref{lem:weighted norm} \marginpar{check commented out things} 
% with $0$  in the role of $r$ and $m:=r$ 
Here are weighted versions of Lemma~\ref{lem:bound on P(f)}
and~\ref{cor:bound on P(f)}: 

\begin{lemma}\label{lem:weighted bd}
Suppose $P$ is homogeneous of degree $d\ge 1$, and let $\w\in \c_a[\imag]^{\b}$ and~$f\in\Car[\imag]^\w$. Then
$$\dabs{P(f)}_a^\w\ \leq\ {d+r \choose r} \cdot \dabs{P}_a\cdot  \dabs{f}_{a;r}^{d-1} \cdot  \dabs{f}_{a;r}^\w.$$
%$$\dabs{P(f)}_a^\w\ \leq\ (d+1)^r{d+r \choose r} \cdot \dabs{P}_a\cdot  \dabs{f}_{a;r}^{d-1} \cdot  \dabs{f}_{a;r}^\w.$$
%where $C=C(d,r):=(d+1)^r{d+r \choose r}\in\N^{\geq 1}$.
\end{lemma}
\begin{proof} For $j=0,\dots,r$ we have $\dabs{f^{(j)}}_a\le \dabs{f}_{a;r}$ and $\dabs{f^{(j)}}_a^\w\le \dabs{f}_{a;r}^\w$.
Now $f^{\i}$, where $\i=(i_0,\dots,i_r)\in \N^{r+1}$ and $i_0+\cdots + i_r=d$, is a product of $d$ such factors~$f^{(j)}$, so
Lemma~\ref{lem:weighted norm} with $m:=d$, $r:=0$, gives
$$ \dabs{f^{\i}}_{a}^\tau\  \le\  \dabs{f}_{a;r}^{d-1}\cdot \dabs{f}_{a;r}^{\w}.$$
It remains to note that by  \eqref{eq:weighted norm, 0} we have $\dabs{P_{\i}f^{\i}}_a^\tau\le \dabs{P_{\i}}_a\cdot\dabs{f^{\i}}_a^\tau$. 
\end{proof}

\begin{cor}\label{cor:weighted bd} Let $1\le d \le e$ in $\N$ be such that $P_{\i}=0$ if~$\abs{\i}<d$ or $\abs{\i}>e$.
Then for $f\in\Car[\imag]^\w$ and $\w\in \c_a[\imag]^{\b}$ we have
$$\dabs{P(f)}_a^\w\ \leq\  D \cdot \dabs{P}_a\cdot   \big( \dabs{f}_{a;r}^{d-1}+\cdots+\dabs{f}_{a;r}^{e-1} \big)\cdot \dabs{f}_{a;r}^\w$$
where
$D\ =\ D(d,e,r)\ :=\  {e+r+1 \choose r+1}-{d+r\choose r+1}\in\N^{\geq 1}$.
\end{cor}

%\subsection*{The  Hardy field case} \marginpar{might not need this after all}
%Let $H$ be a Hardy field and $K=H[\imag]$, an asymptotic field with small derivation.  Let $f_1,\dots,f_n\in K$ and $\phi_1,\dots\phi_n\in\Calr$ with $\phi_1',\dots,\phi_n'\in H$, and  let $\fm\in H^\times$ with $\fm\prec 1$. We assume that $f_1,\dots,f_n, \phi_1,\dots,\phi_n,\fm$ have representatives in $\Car[\imag]$, denoted by the same symbols; then $\fm\in(\Car)^{\b}$, and after increasing $a$ if necessary we can also arrange that  $\fm\in(\Caz)^{\times}$.

%\begin{lemma}
%Suppose $f_1,\dots,f_n\prec\fm^{r+1}$ and $\phi_1',\dots,\phi_n'\preceq\fm^{-1}$; then for  $w=\fm$ we have $f_1\ex^{\imag\phi_1}+\cdots+f_n\ex^{\imag\phi_n} \in \Car[\imag]^w$.
%\end{lemma}
%\begin{proof}
%We can assume $n=1$, $f=f_1$, $\phi=\phi_1$. For $j=0,\dots,r$ we have $f \prec \fm^{j+1}$ and hence $(f\ex^{\imag\phi})^{(j)} \prec \fm$ by Lemma~\ref{lem:der of gexphi}, so $\dabs{f\ex^{\imag\phi}}^w_{a;r}<\infty$.
%\end{proof}

\subsection*{Doubly-twisted integration}
In this subsection we adopt the setting in  {\it Twisted integration}\/ of  Section~\ref{sec:IHF}.
Thus $\phi\in\c_a[\imag]$ and $\Phi=\der_a^{-1}\phi$.
Let $\w\in \Cao$ satisfy~${\w(s)>0}$ for $s\geq a$, and set $\widetilde\phi:=\phi-\w^\dagger\in \c_a[\imag]$  and $\widetilde\Phi:=\der_a^{-1}\widetilde\phi$. Thus
$$\widetilde\Phi(t)\ =\ \int_a^t (\phi-\w^\dagger)(s)\, ds\ =\ \Phi(t) - \log \w(t) + \log \w(a)\qquad\text{for $t\geq a$.}$$
Consider the right inverses
$B, \widetilde B\colon  \c_a[\imag] \to \Cao[\imag]$ 
to, respectively, $\der-\phi\colon\Cao[\imag]\to\c_a[\imag]$ and $\der-\widetilde\phi\colon\Cao[\imag]\to\c_a[\imag]$, 
given by
$$B\ :=\ \ex^\Phi \circ\, {\der_a^{-1}} \circ \ex^{-\Phi}, \quad
  \widetilde B\ :=\ \ex^{\widetilde\Phi} \circ\, {\der_a^{-1}} \circ \ex^{-\widetilde\Phi}.$$
For $f\in\c_a[\imag]$ and $t\geq a$ we have
\begin{align*}
\widetilde{B}f(t)\	&=\ \ex^{\widetilde\Phi(t)} \int_a^t \ex^{-\widetilde\Phi(s)}f(s)\,ds \\
			&=\ \w(t)^{-1}\w(a)\ex^{\Phi(t)} \int_a^t  \ex^{-\Phi(s)} \w(s)\w(a)^{-1} f(s)\,ds \\
	&=\ \w(t)^{-1}\ex^{\Phi(t)} \int_a^t \ex^{-\Phi(s)}\w(s) f(s)\,ds\ =\ \w^{-1}(t)\big( B(\w f)\big)(t)
\end{align*}
and so $\widetilde B=\w^{-1}\circ B\circ \w$.
Hence if  $\widetilde\phi$ is attractive, then $B_{\ltimes \w}:=\w^{-1}\circ B\circ \w$ maps~$\c_a[\imag]^{\b}$ into $\c_a[\imag]^{\b}\cap\Cao[\imag]$, and 
the operator $B_{\ltimes \w}\colon \c_a[\imag]^{\b}\to\c_a[\imag]^{\b}$ is continuous with $\dabs{B_{\ltimes \w}}_a \leq \big\| \frac{1}{\Re \widetilde\phi}\big\|_a$; if in addition $\widetilde\phi\in\Car[\imag]$, then $B_{\ltimes \w}$ maps $\c_a[\imag]^{\b}\cap\Car[\imag]$ into
$\c_a[\imag]^{\b}\cap\Carm[\imag]$.
Note that if $\phi\in\Car[\imag]$ and $\w\in\Carm$, then $\widetilde\phi\in\Car[\imag]$.

\medskip
\noindent
Next, suppose $\phi$, $\widetilde\phi$ are both repulsive. Then we have the $\C$-linear op\-er\-a\-tors $B, \widetilde B\colon  \c_a[\imag]^{\b} \to \Cao[\imag]$ given, for $f\in \c_a[\imag]^{\b}$ and $t\geq a$, by
$$Bf(t)\ :=\ \ex^{\Phi(t)} \int_\infty^t \ex^{-\Phi(s)}f(s)\,ds, \qquad 
\widetilde Bf(t)\ :=\ \ex^{\widetilde\Phi(t)} \int_\infty^t \ex^{-\widetilde\Phi(s)}f(s)\,ds.$$
Now assume  $\w\in\c_a[\imag]^{\b}$. Then
we have the $\C$-linear operator $$B_{\ltimes \w}\ :=\ {\w^{-1}\circ B\circ \w}\ \colon\ \c_a[\imag]^{\b} \to \Cao[\imag].$$
A computation as above shows $\widetilde B = B_{\ltimes \w}$; thus $B_{\ltimes \w}$ maps $\c_a[\imag]^{\b}$ in\-to~$\c_a[\imag]^{\b}\cap\Cao[\imag]$, and the operator $B_{\ltimes \w}\colon \c_a[\imag]^{\b}\to\c_a[\imag]^{\b}$ is continuous with $\dabs{B_{\ltimes \w}}_a \leq \big\| \frac{1}{\Re \widetilde\phi}\big\|_a$.
If~$\widetilde{\phi}\in\Car[\imag]$, then~$B_{\ltimes \w}$ maps $\c_a[\imag]^{\b}\cap\Car[\imag]$ into
$\c_a[\imag]^{\b}\cap\Carm[\imag]$. 

\subsection*{More on twists and right-inverses of linear operators over Hardy fields}
In this subsection we adopt the assumptions in force for Lemma~\ref{cri}, which we repeat here. 
Thus $H$ is a Hardy field, $K=H[\imag]$, $r\in\N^{\geq 1}$, and $f_1,\dots,f_r\in K$.
We fix $a_0\in\R$ and functions in $\c_{a_0}[\imag]$ representing the germs $f_1,\dots,f_r$,
denoted by the same symbols. We let~$a$ range over $[a_0,\infty)$, and we denote the restriction of each
$f\in\c_{a_0}[\imag]$ to  $[a,\infty)$ also by $f$.
For each $a$ we then have the $\C$-linear map~$A_a\colon \Car[\imag]\to\c_a[\imag]$ given by
$$  A_a(y)\ =\ y^{(r)}+f_1y^{(r-1)}+\cdots+f_ry .$$ 
We are in addition given a splitting $(\phi_1,\dots,\phi_r)$ of
the linear differential operator~$A=\der^r+f_1\der^{r-1}+\cdots+f_r\in K[\der]$ over $K$
with $\Re\phi_1,\dots,\Re\phi_r\succeq 1$, 
as well as functions in $\c_{a_0}^{r-1}[\imag]$ representing $\phi_1,\dots,\phi_r$, 
denoted by the same symbols and satisfying  
$\Re\phi_1,\dots,\Re\phi_r\in (\c_{a_0})^\times$.
This gives rise to the continuous $\C$-linear operators
$$B_j\ :=\ B_{\phi_j}\ \colon\ \c_a[\imag]^{\b}\to \c_a[\imag]^{\b} \qquad (j=1,\dots,r)$$
and the right-inverse
$$A_a^{-1} \ :=\  B_r\circ \cdots \circ B_1\ \colon\ \c_a[\imag]^{\b}\to \c_a[\imag]^{\b}$$
of $A_a$ with the properties stated in Lemma~\ref{cri}.

\medskip\noindent
Now let $\fm\in H^\times$ with $\fm\prec 1$, and set~$\widetilde{A}:=A_{\ltimes\fm}\in K[\der]$. Let $\w\in (\Cazr)^\times$ be
a representative of $\fm$. Then $\w\in (\Cazr)^{\b}$
and $\widetilde{\phi}_j:=\phi_j-\w^\dagger\in \Cazrl[\imag]$ for $j=1,\dots,r$.
We   have the $\C$-linear maps
$$\widetilde{A}_j\  :=\  \der-\widetilde{\phi}_j\ \colon\ \Caj[\imag]\to \Cajl[\imag] \qquad (j=1,\dots,r)$$
and for sufficiently large $a$ a factorization
$$\widetilde{A}_a\  =\  \widetilde{A}_1\circ\cdots\circ\widetilde{A}_r\ \colon\ \Car[\imag]\to \c_a[\imag].$$
Below we assume this holds for all $a$, as can be arranged by increasing $a_0$. 
We call~$f,g\in\c_a[\imag]$ 
{\bf alike} if  $f$, $g$ are both attractive or
both repulsive. In the same way we define when germs $f,g\in \c[\imag]$
are alike.\index{alike}\index{germ!alike}\index{function!alike} 
Suppose that  $\phi_j$, $\widetilde{\phi}_j$ are alike for~$j=1,\dots,r$.
%(In particular, $\Re\widetilde{\phi}_j\succeq 1$ and $\Re\widetilde{\phi}_j(t)\neq 0$ for each $j=1,\dots,r$.) 
Then we have   continuous $\C$-linear operators
$$\widetilde{B}_j\ :=\ B_{\widetilde{\phi}_j}\ \colon\ \c_a[\imag]^{\b}\to \c_a[\imag]^{\b} \qquad (j=1,\dots,r)$$
and the right-inverse 
$$\widetilde{A}_a^{-1} \ :=\  \widetilde{B}_r\circ \cdots \circ \widetilde{B}_1\ \colon\ \c_a[\imag]^{\b}\to \c_a[\imag]^{\b}$$
of $\widetilde{A}_a$, and the arguments in the previous subsection show that $\widetilde{B}_j=(B_j)_{\ltimes \w}=\w^{-1}\circ B_j\circ \w$ for $j=1,\dots,r$,  and hence
$\widetilde{A}_a^{-1}=\w^{-1}\circ A_a^{-1}\circ \w$. 
For $j=0,\dots,r$ we set, in analogy with~$A_j^\circ$ and $B_j^\circ$ from \eqref{eq:Ajcirc} and \eqref{eq:Bjcirc}, 
$$\widetilde{A}_j^{\circ}\ :=\ \widetilde{A}_1\circ \cdots \circ \widetilde{A}_j\ \colon\ \Caj[\imag]\to \c_a[\imag],\quad
\widetilde{B}_j^\circ\ :=\ \widetilde{B}_j\circ \cdots \circ \widetilde{B}_1\ \colon\ \c_a[\imag]^{\b}  \to\c_a[\imag]^{\b}.$$ 
Then $\widetilde{B}_j$ maps $\c_a[\imag]^{\b}$ into $\c_a[\imag]^{\b}\cap \Caj[\imag]$, $\widetilde{A}_j^\circ \circ \widetilde{B}_j^\circ$ is the identity on $\c_a[\imag]^{\b}$, and~$\widetilde{B}_j^\circ = \w^{-1}\circ B_j^\circ \circ \w$
by the above.

\subsection*{A weighted version of Proposition~\ref{uban}}
We adopt the setting of the subsection {\it Damping factors}\/ of Section~\ref{sec:IHF}, and
make the same assumptions as in the paragraph before Proposition~\ref{uban}.
Thus $H$, $K$, $A$, $f_1,\dots, f_r$, $\phi_1,\dots, \phi_r$, $a_0$ are as in the previous subsection,
$\fv\in \Cazr$ satisfies $\fv(t)>0$ for all $t\ge a_0$, and its germ $\fv$ is in $H$ with $\fv\prec 1$.  
As part of those assumptions we also have $\phi_1,\dots,\phi_r\preceq_\Delta \fv^{-1}$ in the asymptotic field $K$,  
for the convex subgroup 
 $$\Delta\ :=\ \big\{\gamma\in v(H^\times):\ \gamma=o(v\fv)\big\}$$ 
of $v(H^\times)=v(K^\times)$. Also $\nu\in\R^>$ and $u := \fv^\nu|_{[a,\infty)}\in (\Car)^\times$.

\medskip\noindent
To state a weighted version of Proposition~\ref{uban}, let $\fm\in H^\times$, $\fm\prec 1$, and let $\fm$ also denote a representative  in $(\Cazr)^\times$ of the germ $\fm$.  
Set $\w:=\fm|_{[a,\infty)}$, so we have~$\w\in (\Car)^\times\cap (\Car)^{\b}$ and thus $\Car[\imag]^\w\subseteq \Car[\imag]^{\b}$.
(Note that $\w$, like $u$, depends on~$a$, but we do not indicate this dependence notationally.)
With notations as in the previous subsection we assume that for all $a$ we have the factorization $$\widetilde{A}_a\  =\  \widetilde{A}_1\circ\cdots\circ\widetilde{A}_r\ \colon\ \Car[\imag]\to \c_a[\imag],$$ as can be arranged by increasing $a_0$ if necessary.

%Set $\widetilde{\phi}_j:=\phi_j-w^\dagger\in\Carl[\imag]$ for  $j=1,\dots,r$.

\begin{prop} \label{prop:1.5 weighted}
Assume $H$ is real closed, $\nu\in\Q$, $\nu>r$,
and the elements  $\phi_j$,  $\phi_j-\fm^\dagger$ of $\c_{a_0}[\imag]$ are alike for $j=1,\dots,r$.
Then:
\begin{enumerate}
\item[\rm{(i)}] the $\C$-linear operator $u A_a^{-1}\colon \c_a[\imag]^{\b} \to \c_a[\imag]^{\b}$ maps $\c_a[\imag]^{\w}$ into $\Car[\imag]^{\w}$; 
\item[\rm(ii)]  its restriction to a $\C$-linear map $\c_a[\imag]^{\w} \to \Car[\imag]^{\w}$ is continuous; and
%\item[\rm(iii)] there is a real constant $c\ge 0$ such that $\|uw A_a^{-1}\|_{a;r}^w\le c$ for all $a$;
%\item[\rm(iv)] for all $f\in \Caz[\imag]^{w}$ we have $uw A_a^{-1}(f) \preceq \fv^{\nu}\fm^{-1}\prec 1$; 
\item[\rm(iii)] denoting this restriction also by $uA_a^{-1}$, we have $\|uA_a^{-1}\|^\w_{a;r}\to 0$ as $a\to \infty$.
\end{enumerate}
\end{prop}
\begin{proof}
Let $f\in \c_a[\imag]^{\w}$, so $g:=\w^{-1}f\in\c_a[\imag]^{\b}$. Let $i\in\{0,\dots,r\}$; then with $\widetilde{B}^\circ_j$ as in the previous subsection and $u_{i,j}$ as  in Lemma~\ref{teq}, that lemma gives 
$$\w^{-1}\big(uA_a^{-1}(f)\big){}^{(i)}	= \sum_{j=r-i}^r u_{i,j} u\cdot \w^{-1}B^\circ_j(\w g)\  = 
\sum_{j=r-i}^r u_{i,j}u\widetilde{B}^\circ_j(g).$$
The proof of Proposition~\ref{uban} shows $u_{i,j} u\in\c_a[\imag]^{\b}$ with $\dabs{u_{i,j}u}_a\to 0$ as $a\to\infty$.
Set
$$\widetilde{c}_{i,a}\ :=\ \sum_{j=r-i}^r \|u\,u_{i,j}\|_a \cdot \|\widetilde{B}_j\|_a\cdots \|\widetilde{B}_1\|_a\in [0,\infty)  \qquad (i=0,\dots,r).$$
Then
$\big\|\w^{-1}\big[uA_a^{-1}(f)\big]{}^{(i)}\big\|_a\le \widetilde{c}_{i,a}\|g\|_a=\widetilde{c}_{i,a}\|f\|^\w_a$ 
where $\widetilde{c}_{i,a}\to 0$ as $a\to\infty$. This yields (i)--(iii).
\end{proof}

%\subsection*{A weighted variant of Lemma~\ref{lem:2.1 summary}}
\subsection*{Weighted variants of results in Section~\ref{sec:split-normal over Hardy fields}}
In this subsection we adopt the hypotheses in force for Lemma~\ref{bdua}. To summarize those,   
%of Lemma~\ref{lem:2.1 summary}.
$H$, $K$, $A$, $f_1,\dots, f_r$, $\phi_1,\dots, \phi_r$, $a_0$, $\fv$, $\nu$, $u$, $\Delta$ are as in the previous subsection,
$d,r\in \N^{\ge 1}$, $H$ is real closed, $R\in K\{Y\}$ has order~$\leq d$ and weight~$\leq w\in\N^{\geq r}$.
Also $\nu\in \Q$, $\nu>w$,
$R\prec_{\Delta}\fv^\nu$, $\nu\fv^\dagger\not\sim \Re \phi_j$ and
$\Re\phi_j - \nu\fv^\dagger\in (\c_{a_0})^\times$
for $j=1,\dots,r$. Finally, $\tilde A:=A_{\ltimes\fv^\nu}\in K[\der]$ and $\tilde A_a(y)=u^{-1}A_a(uy)$ for $y\in\Car[\imag]$.
%where $A_a$, $\tilde A_a$  are as introduced before  Lemma~\ref{cri}.
Next, let  $\fm$,~$\w$ be as in the previous subsection.  As in Lemma~\ref{lem:2.1 summary} we consider the
continuous operator
$$\Phi_a\colon \Car[\imag]^{\b}\times \Car[\imag]^{\b}\to \Car[\imag]^{\b}$$ given by
$$\Phi_a(f,y)\ :=\ \Xi_a(f+y)-\Xi_a(f)\ =\  u\widetilde{A}_a^{-1}\!\left( u^{-1}\big( R(f+y)-R(f)\big) \right).$$
Here is our weighted version of Lemma~\ref{lem:2.1 summary}:

{\samepage
\begin{lemma}\label{lem:2.1 weighted}
Suppose the elements $\phi_j-\nu\fv^\dagger$,~$\phi_j-\nu\fv^\dagger-\fm^\dagger$ of $\c_{a_0}[\imag]$ are alike, for~$j=1,\dots,r$, and let $f\in\Car[\imag]^{\b}$. Then the operator~$y\mapsto \Phi_a(f,y)$
maps~$\Car[\imag]^\w$ into itself. Moreover, there are~$E_a,E_a^+\in\R^{\geq}$ such that for all $g\in \Car[\imag]^{\b}$ and
$y\in \Car[\imag]^{\w}$, 
$$\dabs{\Phi_a(f,y)}^\w_{a;r}  
 \ \le\ E_a\cdot \max\!\big\{1, \|f\|_{a;r}^d\big\}\cdot  \big( 1+\dabs{y}_{a;r}+\cdots+\dabs{y}_{a;r}^{d-1} \big)\cdot\dabs{y}_{a;r}^\w,$$
\vskip-1.5em
\begin{multline*}
\|\Phi_a(f,g+y)-\Phi_a(f,g) \|_{a;r}^\w\ \le\ \\ E_a^+\cdot \max\!\big\{1, \|f\|_{a;r}^d\big\}\cdot \max\!\big\{1, \|g\|_{a;r}^d\big\}\cdot \big( 1+\dabs{y}_{a;r}+\cdots+\dabs{y}_{a;r}^{d-1} \big)\cdot\dabs{y}_{a;r}^\w.
\end{multline*}
We can take these $E_a$, $E_a^+$ such that $E_a,E_a^+ \to 0$ as $a\to \infty$, and do so below. 
\end{lemma}}
\begin{proof}
Let $y\in \Car[\imag]^{\w}$. 
By Taylor expansion we have
$$R(f+y)-R(f)\ =\ \sum_{\abs{\i}>0} \frac{1}{\i !}R^{(\i)}(f)y^{\i}\ =\ \sum_{\abs{\i}>0} S_{\i}(f)y^{\i}
\quad\text{where $S_{\i}(f):= \frac{1}{\i !}\sum_{\j}R_{\j}^{(\i)}f^\j$,}$$
and $u^{-1}S_{\i}(f)\in\c_a[\imag]^{\b}$. So
$h:=u^{-1}\big(R(f+y)-R(f)\big) \in \c_a[\imag]^\w$, since~$\c_a[\imag]^\w$ is an ideal of  $\c_a[\imag]^{\b}$. 
Applying Proposition~\ref{prop:1.5 weighted}(i) with $\phi_j-\nu\fv^\dagger$ in the role of $\phi_j$ yields~$\Phi_a(f,y)=u\widetilde{A}_a^{-1}(h) \in \Car[\imag]^\w$, establishing the first claim. Next, let $g\in \Car[\imag]^{\b}$. Then $\Phi_a(f,g+y)-\Phi_a(f,g) = \Phi_a(f+g,y)$ by \eqref{eq:Phia difference}. Therefore,  
\begin{align*} \Phi_a(f,g+y)-\Phi_a(f,g)\ =\  u\widetilde{A}_a^{-1}(h),&\quad h\ :=\ u^{-1}\big(R(f+g+y)-R(f+g)\big), \text{ so}\\
\dabs{\Phi_a(f,g+y)-\Phi_a(f,g)}_{a;r}^\w\ &=\ \dabs{ u\widetilde{A}_a^{-1}(h) }_{a;r}^\w\ \leq\ \dabs{u\widetilde{A}_a^{-1}}_{a;r}^\w\cdot
\dabs{h}^\w_a.
\end{align*}
By Corollary~\ref{cor:weighted bd} we have 
$$\dabs{h}_a^\w\ \leq\ D
\cdot\max_{\abs{\i}>0}\, \dabs{u^{-1}S_{\i}(f+g)}_a\cdot \big( 1+\dabs{y}_{a;r}+\cdots+\dabs{y}_{a;r}^{d-1} \big)\cdot\dabs{y}_{a;r}^\w$$
where $D=D(d,r):=\big( \textstyle{d+r+1 \choose r+1}-1\big)$.
Let $D_a$ be as in the proof of Lem\-ma~\ref{bdua, bds}. Then $D_a\to 0$ as $a\to \infty$, and  Lemma~\ref{lem:inequ power d} gives for $\abs{\i}>0$,
\begin{align*} \dabs{u^{-1}S_{\i}(f)}_a\ &\leq\ D_a\cdot\max\!\big\{1,\dabs{f}_{a;r}^d\big\}\\
\dabs{u^{-1}S_{\i}(f+g)}_a\ &\leq\ D_a\cdot\max\!\big\{1,\dabs{f+g}_{a;r}^d\big\}\\
&\le\ 2^dD_a\cdot \max\!\big\{1,\dabs{f}_{a;r}^d\big\}\cdot\max\!\big\{1,\dabs{g}_{a;r}^d\big\}.
\end{align*}
This gives the desired result for 
$E_a:=\dabs{u\widetilde{A}_a^{-1}}_{a;r}^\w\cdot D\cdot D_a$ and $E_a^+:=2^d E_a$, using also Proposition~\ref{prop:1.5 weighted}(iii)
with $\phi_j-\nu\fv^\dagger$ in the role of $\phi_j$.
\end{proof}

\noindent
Lemma~\ref{lem:2.1 weighted} allows us to refine Theorem~\ref{thm:fix} as follows:

\begin{cor}\label{cor:2.3 weighted}
Suppose the elements $\phi_j-\nu\fv^\dagger$,~$\phi_j-\nu\fv^\dagger-\fm^\dagger$ of $\c_{a_0}[\imag]$ are alike, for $j=1,\dots,r$, and $R(0)\preceq \fv^\nu\fm$. Then
for sufficiently large $a$ the operator~$\Xi_a$ maps the closed ball
$B_a := \big\{f\in \Car[\imag]:\, \|f\|_{a;r}\le 1/2\big\}$
of the normed space~$\Car[\imag]^{\b}$ into itself, has a unique fixed point in $B_a$, and this fixed point 
lies in $\Car[\imag]^\w$.
\end{cor}
\begin{proof}
Take $a$ such that $\dabs{\w}_a\leq 1$. Then by \eqref{eq:weighted norm, 1}, $B_a$ contains the closed ball
$$B_a^\w\ :=\ \big\{f\in \Car[\imag]:\, \|f\|^\w_{a;r}\le 1/2\big\}$$
of the normed space $\Car[\imag]^\w$. Let $f,g\in B_a^\w$. Then
$\Xi_a(g)-\Xi_a(f)=\Phi_a(f,g-f)$ lies in $\Car[\imag]^\w$ by Lemma~\ref{lem:2.1 weighted}, and
 with $E_a$ as in that lemma,
\begin{align*}
\dabs{\Xi_a(f)-\Xi_a(g)}_{a;r}^\w\	&=\ \dabs{\Phi_a(f,g-f)}_{a;r}^\w \\
									& \leq\  E_a\cdot
									\max\!\big\{1, \|f\|_{a;r}^d\big\}\cdot  \big( 1+\cdots+\dabs{g-f}_{a;r}^{d-1} \big)\cdot\dabs{g-f}_{a;r}^\w \\
									& \leq\  E_a\cdot d\cdot \dabs{g-f}^\w_{a;r}. 
\end{align*}
Taking $a$ so that moreover $E_ad\leq\frac{1}{2}$ we obtain 
\begin{equation}\label{eq:2.3 weighted}
\dabs{\Xi_a(f)-\Xi_a(g)}_{a;r}^\w\leq\textstyle\frac{1}{2}\dabs{f-g}^\w_{a;r}\quad\text{ for all
$f,g\in B_a^\w$.}
\end{equation} 
Next we consider the case $g=0$. 
Our hypothesis $R(0)\preceq \fv^\nu \fm$ gives $u^{-1}R(0)\in \c_a[\imag]^\w$. Proposition~\ref{prop:1.5 weighted}(i),(ii)
with $\phi_j-\nu\fv^\dagger$ in the role of~$\phi_j$ gives~$\Xi_a(0)\in \Car[\imag]^\w$ and $\dabs{\Xi_a(0)}_{a;r}^\w \leq \dabs{u\widetilde{A}_a^{-1}}_{a;r}^\w\dabs{u^{-1}R(0)}_a^\w$.
Using Proposition~\ref{prop:1.5 weighted}(iii) we now take~$a$ so large that
$\dabs{\Xi_a(0)}_{a;r}^\w \leq \frac{1}{4}$. Then \eqref{eq:2.3 weighted} for $g=0$ gives $\Xi_a(B_a^\w)\subseteq B_a^\w$. 
By Lemma~\ref{lem:w-conv} the normed space~$\Car[\imag]^\w$ is complete,
hence~$\Xi_a$ has a unique fixed point in $B_a^\w$.
\end{proof}

\noindent
Now suppose in addition that $A\in H[\der]$ and $R\in H\{Y\}$.  Set
$$(\Car)^\w\ :=\ \big\{f\in\Car:\dabs{f}^\w_{a;r}<\infty\}\ =\ \Car[\imag]^\w\cap\Car,$$
a real Banach space with respect to $\dabs{\,\cdot\,}_{a;r}^\w$.
Increase $a_0$ as at the beginning of the subsection {\it Preserving reality}\/ of Section~\ref{sec:split-normal over Hardy fields}. Then we have the map 
$$\Re\Phi_a\colon(\Car)^{\b}\times(\Car)^{\b}\to(\Car)^{\b},\qquad (f,y)\mapsto \Re\!\big(\Phi_a(f,y)\big).$$
Suppose the elements $\phi_j-\nu\fv^\dagger$,~$\phi_j-\nu\fv^\dagger-\fm^\dagger$ are alike for  $j=1,\dots,r$, and let~$a$ and $E_a, E_a^{+}$  be as in
Lemma~\ref{lem:2.1 weighted}. Then this lemma yields:

\begin{lemma}\label{lem:2.1 weighted, real}
Let $f,g\in (\Car)^{\b}$ and $y\in(\Car)^\w$. Then $(\Re\Phi_a)(f,y)\in  (\Car)^\w$ and 
$$\dabs{\Re(\Phi_a)(f,y)}^\w_{a;r}  
 \ \le\ E_a\cdot \max\!\big\{1, \|f\|_{a;r}^d\big\}\cdot  \big( 1+\dabs{y}_{a;r}+\cdots+\dabs{y}_{a;r}^{d-1} \big)\cdot\dabs{y}_{a;r}^\w,$$
\begin{multline*}
\|(\Re\Phi_a)(f,g+y)-(\Re\Phi_a)(f,g) \|_{a;r}^\w\ \le\ \\ E_a^{+}\cdot \max\!\big\{1, \|f\|_{a;r}^d\big\}\cdot \max\!\big\{1, \|g\|_{a;r}^d\big\}\cdot \big( 1+\dabs{y}_{a;r}+\cdots+\dabs{y}_{a;r}^{d-1} \big)\cdot\dabs{y}_{a;r}^\w.
\end{multline*}
%We can take these $E_a$ such that $E_a\to 0$ as $a\to \infty$. 
\end{lemma}

\noindent
The same way we derived Corollary~\ref{cor:2.3 weighted} from Lemma~\ref{lem:2.1 weighted}, this leads to: 

\begin{cor}\label{cor:2.3 weighted, real}  If $R(0)\preceq \fv^\nu\fm$, then
for sufficiently large $a$ the operator~$\Re\Xi_a$ maps the closed ball
$B_a := \big\{f\in \Car:\, \|f\|_{a;r}\le 1/2\big\}$
of the normed space $(\Car)^{\b}$ into itself, has a unique fixed point in $B_a$, and this fixed point 
lies in $(\Car)^\w$.
\end{cor}

\subsection*{Revisiting Lemma~\ref{lem:Psin, b}}
Here we adopt the setting of the
previous subsection. As usual, $a$ ranges over~$[a_0,\infty)$.
We continue our investigation of the differences~${f-g}$ between solutions~$f$,~$g$
of the equation \eqref{eq:ADE} on $[a_0,\infty)$ from Section~\ref{sec:split-normal over Hardy fields} which we began in Lemma~\ref{lem:close}, and so we
take $f$, $g$, $E$, $\epsilon$,  $h_a$, $\theta_a$ as in that lemma.
Recall that in the remarks preceding Lemma~\ref{lem:Psin, b} we defined continuous operators~$\Phi_a, \Psi_a\colon \Car[\imag]^{\b}\to \Car[\imag]^{\b}$ by
$$\Phi_a(y)\ :=\ \Phi_a(g,y)\ =\ \Xi_a(g+y)-\Xi_a(g), \quad \Psi_a(y)\ :=\ \Phi_a(y)+h_a\qquad(y\in\Car[\imag]^{\b}).$$
As in those remarks, we set $\rho:=\dabs{f-g}_{a_0;r}$ and
$$B_a\ :=\ \big\{ y\in\Car[\imag]^{\b}:\ \dabs{y-h_a}_{a;r} \leq 1/2\big\},$$
and take $a_1\ge a_0$ so that $\theta_a\in B_a$ for all $a\ge a_1$. Then by \eqref{eq:dabs(y)} we have~$\dabs{y}_{a;r}\leq 1+\rho$ 
for $a\geq a_1$ and $y\in B_a$.
Next, take $a_2\geq a_1$ as in Lemma~\ref{lem:Psin, b}; thus 
for~$a\geq a_2$ and $y,z\in B_a$ we have
$\Psi_a(y) \in B_a$ and $\dabs{{\Psi_a(y)-\Psi_a(z)}}_{a;r} \leq \textstyle\frac{1}{2} \|y-z\|_{a;r}$.
As in the previous subsection,  $\fm\in H^\times$, $\fm\prec 1$, $\fm$ denotes also  a representative  in~$(\Cazr)^\times$ of the germ $\fm$, and $\w:=\fm|_{[a,\infty)}\in (\Car)^\times\cap (\Car)^{\b}$, so $\Car[\imag]^\w\subseteq \Car[\imag]^{\b}$.

\medskip
\noindent
{\it In the rest of this subsection
$\phi_1-\nu\fv^\dagger,\dots,\phi_r-\nu\fv^\dagger\in K$ are $\gamma$-repulsive for $\gamma:=v\fm\in v(H^\times)^>$,
and $h_a\in\Car[\imag]^\w$ for all $a\ge a_2$.}\/
Then Corollary~\ref{cor:repulsive} gives~$a_3\geq a_2$ such that for all~$a\geq a_3$ and~$j=1,\dots,r$, the functions~$\phi_j-u^\dagger, \phi_j-(u\w)^\dagger\in\c_a[\imag]$ are alike and hence  $\Psi_a\big(\Car[\imag]^\w\big)\subseteq  \Car[\imag]^\w$
by Lem\-ma~\ref{lem:2.1 weighted}. 
Thus $\Psi_a^n(h_a) \in \Car[\imag]^\w$ for all  $n$ and $a\geq a_3$.  

For $a\geq a_2$ we have  $\lim_{n\to\infty} \Psi_a^n(h_a)=\theta_a$ in $\Car[\imag]^{\b}$ by Corollary~\ref{cor:Psin, b};
we now aim to strengthen this to ``in  $\Car[\imag]^{\w}$'' (possibly for a larger $a_2$).
Towards this:

%Now put
%$$B_a^w:=\big\{ y\in\Car[\imag]^{\b}: \dabs{y-h_a}^w_{a;r} \leq \widetilde{\rho}\big\}\qquad\text{where $\widetilde{\rho}:=\rho/\dabs{w}_a$.}$$
%Note that $B_a^w\subseteq B_a$ by \eqref{eq:weighted norm, 1}.

\begin{lemma}\label{lem:Psia contract}
There exists $a_4\geq a_3$ such that  
%$\Psi_a(B_a^w)\subseteq B_a^w$, and 
$\dabs{{\Psi_a(y)-\Psi_a(z)}}^\w_{a;r} \leq \frac{1}{2} \|y-z\|^\w_{a;r}$ for all $a\ge a_4$ and $y,z\in B_a\cap \Car[\imag]^\w$.
\end{lemma}
\begin{proof}
For $a\geq a_3$ and $y,z\in\Car[\imag]^\w$, and with $E_a^+$ as in Lemma~\ref{lem:2.1 weighted} we have
\begin{multline*}
\|\Psi_a(y)-\Psi_a(z) \|_{a;r}^\w\ \le\ \\ 
E_a^+\cdot \max\!\big\{1, \|g\|_{a;r}^d\big\}\cdot \max\!\big\{1, \|z\|_{a;r}^d\big\}\cdot \big( 1+\dabs{y-z}_{a;r}+\cdots+\dabs{y-z}_{a;r}^{d-1} \big)\cdot\dabs{y-z}_{a;r}^\w.
\end{multline*}
For each $a\geq a_1$ and $y,z\in B_a$  we then have
$$\max\!\big\{1, \|z\|_{a;r}^d\big\}\cdot \big( 1+\dabs{y-z}_{a;r}+\cdots+\dabs{y-z}_{a;r}^{d-1} \big)\ \leq\ (1+\rho)^d\cdot d,$$
so taking $a_4\geq a_3$ with 
$$E_a^+\max\!\big\{1, \|g\|_{a_0;r}^d\big\} (1+\rho)^dd\le 1/2\quad\text{ for all $a\geq a_4$,}$$ 
we have $\dabs{{\Psi_a(y)-\Psi_a(z)}}^\w_{a;r} \leq \frac{1}{2} \|y-z\|^\w_{a;r}$ for all $a\ge a_4$ and $y,z\in B_a\cap \Car[\imag]^\w$.
\end{proof}

\noindent
Let $a_4$ be as in the previous lemma.

\begin{cor}\label{cor:small theta, 1}
Suppose $a\geq a_4$. Then $\theta_a\in \Car[\imag]^\w$ and $\lim_{n\to \infty} \Psi_a^n(h_a)=\theta_a$ in the normed space $\Car[\imag]^\w$. In particular,~$f-g,(f-g)',\dots, (f-g)^{(r)}\preceq\fm$.
\end{cor}

\begin{proof} 
We have $\Phi_a(h_a)=\Psi_a(h_a)-h_a\in\Car[\imag]^\w$, so  $M:=\dabs{\Phi_a(h_a)}_{a;r}^\w<\infty$. Since~$\Psi_a(B_a)\subseteq B_a$, 
induction on $n$ using Lemma~\ref{lem:Psia contract} gives
$$\dabs{\Psi_a^{n+1}(h_a)-\Psi^n_a(h_a)}^\w_{a;r} \leq M/2^n\qquad\text{ for all $n$.}$$
%Hence for all $m$, $n$ we have
%\begin{align*}
%\dabs{\Psi_a^{m+n}(h_a)-\Psi_a^m(h_a)}^w_{a;r}	&\ \leq\ 
%\sum_{j=0}^{n-1}\, \dabs{\Psi_a^{m+j+1}(h_a)-\Psi_a^{m+j}(h_a)}^w_{a;r} \\ &\ \leq\ 
%\frac{M}{2^m} \sum_{j=0}^{n-1}\frac{1}{2^j}\  \leq \ \frac{M}{2^{m-1}}.
%\end{align*}
Thus $\big( \Psi_a^n(h_a) \big)$ is a cauchy sequence in the normed space~$\Car[\imag]^\w$, and so converges in~$\Car[\imag]^\w$ by Lemma~\ref{lem:w-conv}.
In the normed space $\Car[\imag]^{\b}$ we have $\lim_{n\to\infty}\Psi_a^n(h_a)  =\theta_a$,
by Corollary~\ref{cor:Psin, b}.
Thus $\lim_{n\to \infty}\Psi_a^n(h_a)  =\theta_a$ in $\Car[\imag]^\w$ by Lem\-ma~\ref{lem:w-conv => b-conv}.
\end{proof}

%\noindent
%Combining the above with [{\tt ahardy},~3.1] yields:

%\begin{cor}\label{cor:small theta, 2}
%Suppose $f-g\prec 1$ and
%$h,h',\dots,h^{(r)}\preceq\fm$ for all $h\in\Calinf[\imag]$ with $A(h)=0$ and $h\prec 1$;
%then $(f-g)^{(j)}\preceq\fm$ for $j=0,\dots,r$.
%\end{cor}
%\begin{proof}
%For  sufficiently large $a$ take $h_a\in\Car[\imag]$ with $A_a(h_a)=0$ as in the proof of [{\tt ahardy}, 2.4]. 
%Then $h_a\in\Car[\imag]^{\b}$, and by [{\tt ahardy},~3.1] the germ $h$ of $h_a$ lies in $\Calinf[\imag]$.
%Moreover, $\theta_a-h_a\prec\fv^w\prec 1$ and $f-g\prec 1$ implies $h\prec 1$,  hence $h,h',\dots,h^{(r)}\preceq\fm$,
%so $h_a\in \Car[\imag]^w$. Now the claim follows from Corollary~\ref{cor:small theta, 1}.
%\end{proof}

%\begin{remark}
%Suppose $\deg R\leq 0$; then the previous corollary goes through even without the assumption that
% $\phi_1-\nu\fv^\dagger,\dots,\phi_r-\nu\fv^\dagger\in H[\imag]$ are $\gamma$-repulsive,
% since then the germ  $h$ of $f-g$ satisfies $A(h)=0$, $h\prec 1$.
%\end{remark}

\subsection*{An application to slots in $H$}
Here we adopt the setting of the subsection {\it An application to slots in $H$}\/ in Section~\ref{sec:ueeh}. Thus
$H\supseteq\R$ is a Liouville closed Hardy field,
$K:=H[\imag]$, 
 $\I(K)\subseteq K^\dagger$, and
$(P,1,\hat h)$ is a slot in $H$ of order~$r\geq 1$; we set~$w:= \wt(P)$, $d:= \deg P$. 
{\it Assume also that $K$ is  $1$-linearly surjective if $r\geq 3$.}\/

\begin{prop} \label{prop:notorious 3.6}  
Suppose $(P,1,\hat h)$ is special, ultimate, $Z$-minimal, deep,  and strongly repulsive-normal. 
Let $f,g\in \Calr[\imag]$ and $\fm\in H^\times$ be such that 
$$P(f)\ =\ P(g)\ =\ 0,\qquad f,g\ \prec\ 1, \qquad v\fm\in v(\hat h - H).$$
Then $(f-g)^{(j)}\preceq\fm$ for $j=0,\dots,r$.
\end{prop}
\begin{proof}
We arrange~${\fm\prec 1}$. Let $\fv:=\abs{\fv(L_P)}\in H^>$, so $\fv\prec^\flat 1$, and set $\Delta:=\Delta(\fv)$.
Take $Q,R\in H\{Y\}$ where $Q$ is homogeneous of degree $1$ and order $r$, 
$A:=L_Q\in H[\der]$ has a strong $\hat h$-repulsive splitting over~$K$,
 $P=Q-R$, and $R\prec_\Delta \fv^{w+1}P_1$, so~$\fv(A)\sim \fv(L_P)$ by Lemma~\ref{lem:fv of perturbed op}. Multiplying $P$, $Q$, $R$ by some $b\in H^\times$ we
arrange that $A$ is monic, so $A=\der^r+ f_1\der^{r-1}+\cdots + f_r$ with $f_1,\dots, f_r\in H$ and~$R\prec_\Delta\fv^w$.
Let $(\phi_1,\dots,\phi_r)\in K^r$ be a strong $\hat h$-repulsive splitting of $A$ over $K$, so $\phi_1,\dots,\phi_r$ are $\hat h$-repulsive and
$$A\ =\ (\der-\phi_1)\cdots (\der-\phi_r), \qquad \Re\phi_1,\dots,\Re \phi_r\ \succeq\ \fv^\dagger\ \succeq\ 1. $$
By Corollary~\ref{cor:bound on linear factors} we have $\phi_1,\dots,\phi_r\preceq\fv^{-1}$. Thus we can take $a_0\in \R$ and functions on $[a_0,\infty)$ representing the germs $\phi_1,\dots, \phi_r$, $f_1,\dots, f_r$, $f$, $g$ and the $R_{\j}$ with~$\j\in \N^{1+r}$, $|\j|\le d$, $\|\j\|\le w$ (using the same symbols for the germs mentioned as for their chosen representatives)
so as to be in the situation described in the beginning of Section~\ref{sec:split-normal over Hardy fields}, with $f$ and $g$ solutions on $[a_0,\infty)$ of the differential equation~\eqref{eq:ADE} there. 
As there, we take $\nu\in\Q$ with $\nu > w$ so that~$R \prec_\Delta \fv^\nu$ and~$\nu\fv^\dagger\not\sim \Re\phi_j$
for $j=1,\dots,r$, and then increase $a_0$ to satisfy all assumptions for Lemma~\ref{bdua}. 
Corollary~\ref{specialvariant} gives
$v(\fv^\nu) \in v(\hat h-H)$, so
$\phi_j-\nu\fv^\dagger=\phi_j-(\fv^\nu)^\dagger$ ($j=1,\dots,r$)
is $\hat h$-repulsive by Lemma~\ref{lem:hata-repulsive}(iv), so $\gamma$-repulsive for $\gamma:=v\fm>0$. Now $A$ splits over~$K$, and $K$ is $1$-linearly surjective if $r\ge 3$,  hence 
$\dim_{\C}\ker_{\Univ} A =r$ by Lemma~\ref{lem:full kernel}. Thus 
by Corollary~\ref{cor:8.8 refined} we have $y,y',\dots,y^{(r)}\prec\fm$ for all~$y\in\Calr[\imag]$ with $A(y)=0$, $y\prec 1$. 
In particular, $\fm^{-1}h_a, \fm^{-1}h_a',\dots, \fm^{-1}h_a^{(r)}\prec 1$ for all $a\ge a_0$. Thus the assumptions on $\fm$ and the $h_a$ made just before Lemma~\ref{lem:Psia contract} are satisfied for a suitable choice of $a_2$,
so we can appeal to Corollary~\ref{cor:small theta, 1}.
\end{proof}

\noindent
The assumption that $K$ is $1$-linearly surjective for $r\ge 3$ was only used in the proof above to obtain $\dim_{\C}\ker_{\Univ} A =r$. So if $A$ as in this proof satisfies $\dim_{\C}\ker_{\Univ} A =r$, then we can drop this assumption about $K$, also in the next corollary. 

\begin{cor}\label{cor:notorious 3.6}
Suppose $(P,1,\hat h)$, $f$, $g$, $\fm$  are as in Proposition~\ref{prop:notorious 3.6}. Then 
$$f-g\in \fm\, \c^r[\imag]^{\preceq}.$$
%we have $\big((f-g)/\fm\big){}^{(j)}\preceq 1$ for $j=0,\dots,r$.
\end{cor}
\begin{proof}
If $\fm\succeq 1$,  then  Lemma~\ref{lem:small derivatives of y}(ii) applied with $y=(f-g)/\fm$ and $1/\fm$ in place of~$\fm$ gives what we want.
% $\big((f-g)/\fm\big){}^{(j)}\preceq 1$ for $j=0,\dots,r$.
Now assume $\fm\prec 1$. 
%By [{\tt mN}, 4.17, 4.34], $\hat f$ is special over~$K$. 
Since $\hat h$ is special over $H$, Proposition~\ref{prop:notorious 3.6}
applies to $\fm^{r+1}$ in place of $\fm$, so $(f-g)^{(j)}\preceq\fm^{r+1}$ for $j=0,\dots,r$. Now apply Lemma~\ref{lem:Qnk} to suitable representatives of $f-g$ and~$\fm$.
\end{proof}

\noindent
Later in this section we use Proposition~\ref{prop:notorious 3.6} and its Corollary~\ref{cor:notorious 3.6} to
strengthen some results from Section~\ref{secfhhf}.
In Section~\ref{sec:holes perfect} we give further refinements of that proposition for the case of firm and flabby slots,
but these are not needed for the proof of our main result, given in Section~\ref{sec:d-alg extensions}.

\subsection*{Weighted refinements of results in Section~\ref{secfhhf}}
We now adopt the setting of the subsection {\it Reformulations}\/ of Section~\ref{secfhhf}. Thus $H\supseteq\R$ is a real closed Hardy field
with asymptotic integration, and $K:=H[\imag]\subseteq \Calinf[\imag]$ is its algebraic closure,
with value group $\Gamma:=v(H^\times)=v(K^\times)$.
The next lemma and its corollary refine Lemma~\ref{prop:as equ 1}. Let~$P$,~$Q$,~$R$,~$L_Q$,~$w$ be as introduced before that lemma, set~$\fv:=\abs{\fv(L_Q)}\in H^>$,  and, in case $\fv\prec 1$,  $\Delta:=\Delta(\fv)$.

\begin{lemma}\label{abc}
Let $f\in K^\times$ and $\phi_1,\dots,\phi_r\in K$ be such that
$$ L_Q \ =\  f(\der-\phi_1)\cdots(\der-\phi_r),\qquad \Re\phi_1,\dots,\Re\phi_r\ \succeq\ 1.$$
Assume $\fv\prec 1$ and $R\prec_\Delta \fv^{w+1}Q$. Let $\fm\in H^\times$, 
$\fm\prec 1$, $P(0)\preceq \fv^{w+2}\fm Q$. Suppose that for $j=1,\dots,r$ and all~$\nu\in\Q$ with $w<\nu<w+1$, 
$\phi_j-(\fv^\nu)^\dagger$ and ${\phi_j-(\fv^\nu\fm)^\dagger}$  are alike.
Then $P(y)=0$ and $y,y',\dots,y^{(r)}\preceq\fm$ for some $y\prec\fv^w$ in $\Calinf[\imag]$.  If~$P,Q\in H\{Y\}$,
then there is such $y$ in $\Calinf$.
\end{lemma}
\begin{proof} Note that $\phi_1,\dots, \phi_r\preceq \fv^{-1}$ by Corollary~\ref{cor:bound on linear factors} and that $R\prec_{\Delta} \fv^{w+1}Q$ gives~$f^{-1}R\prec_{\Delta} \fv^w$.
Take~$\nu\in\Q$ such that $w<\nu<w+1$, $f^{-1}R
\prec_{\Delta} \fv^{\nu}$ and~$\nu\fv^\dagger\not\sim\Re\phi_j$
for $j=1,\dots,r$. Set $A\ :=\ f^{-1}L_Q$. From $\nu < w+1$ and 
$$R(0)\ =\ -P(0)\ \prec_\Delta\ \fv^{w+2}\fm Q$$ 
we obtain $f^{-1}R(0) \prec_{\Delta} \fv^{\nu}\fm$.  Thus we can apply successively Corollary~\ref{cor:2.3 weighted}, 
Lemma~\ref{bdua}, and
Corollary~\ref{cor:ADE smooth} to the equation $A(y)=f^{-1}R(y)$, $y\prec 1$ in the role of \eqref{eq:ADE} in Section~\ref{sec:split-normal over Hardy fields} to obtain the first part.  For the real variant, use instead Corollary~\ref{cor:2.3 weighted, real} and
Lemma~\ref{realbdua}.
\end{proof}

\noindent
Lemma~\ref{abc} with $\fm^{r+1}$ for $\fm$ has the following consequence, using Lemma~\ref{lem:Qnk}: 
%we obtain from Lemma~\ref{abc}:

\begin{cor}\label{cor:abc}  Let  $f\in K^\times$ and $\phi_1,\dots,\phi_r\in K$ be such that
$$  L_Q \ =\  f(\der-\phi_1)\cdots(\der-\phi_r),\qquad \Re\phi_1,\dots,\Re\phi_r\ \succeq\ 1.$$
Assume $\fv\prec 1$ and $R\prec_\Delta \fv^{w+1}Q$. Let $\fm\in H^\times$, 
$\fm\prec 1$, $P(0)\preceq \fv^{w+2}\fm^{r+1} Q$. Suppose 
that  for~${j=1,\dots,r}$ and all $\nu\in\Q$ with $w<\nu<w+1$,
$\phi_j-(\fv^\nu)^\dagger$ and~${\phi_j-(\fv^\nu\fm^{r+1})^\dagger}$  are alike.
Then for some $y\prec\fv^w$ in $\Calinf[\imag]$ we have~${P(y)=0}$ and $y\in \fm\, \c^r[\imag]^{\preceq}$.
%$(y/\fm)^{(j)}\preceq 1$ for $j=0,\dots,r$.
 If $P,Q\in H\{Y\}$,
then there is such $y$ in $\Calinf$.
\end{cor} 

\begin{remark}
If $H$ is a $\c^{\infty}$-Hardy field, then Lem\-ma~\ref{abc} and Corollary~\ref{cor:abc} go through with $\Calinf[\imag]$, $\Calinf$ replaced by~$\mathcal{C}^{\infty}[\imag]$, $\mathcal{C}^{\infty}$, respectively. Likewise with $\c^\omega$ in place of $\c^\infty$. (Use Corollary~\ref{cor:ADE smooth}.)  
\end{remark}

\noindent
Next a  variant of Lemma~\ref{dentsolver}. {\it In the rest of this subsection $(P,\fn,\hat h)$ is a deep, strongly repulsive-normal, $Z$-minimal slot in $H$ of order~$r\ge 1$ and weight~$w:=\wt(P)$.  We assume also that $(P,\fn, \hat h)$ is special \textup{(}as will be the case if $H$ is $r$-linearly newtonian,  and $\upo$-free if $r>1$, by Lemma~\ref{lem:special dents}\textup{)}}. 
%From \dots recall $\mathfrak d(\hat h-H) = \{\fn\in H^\times:\ \text{$\hat h - h\asymp \fn$ for some $h\in H$}\}$.

%\begin{notation}
%Let $\hat K$ be an immediate valued field extension of $K$.
%Given $S\subseteq\hat K^\times$ we set
%$$\mathfrak d(S) \ &=  \ \{\fn\in K^\times:\ \text{$s \asymp \fn$ for some $s\in S$} \};$$
%thus $\mathfrak d(S)=v^{-1}(v(S))\cap K^\times$ and $v(\mathfrak d(S))=v(S)$.
%\end{notation}

\begin{lemma}\label{bcd} 
Let $\fm\in H^\times$ be such that $v\fm\in v(\hat h-H)$, $\fm\prec\fn$,  and $P(0)\preceq\fv(L_{P_{\times\fn}})^{w+2}\,(\fm/\fn)^{r+1}\, P_{\times\fn}$.  
Then for some~$y\in \Calinf$,
$$P(y)\ =\ 0,\quad  y\in \fm\,(\c^r)^{\preceq}.$$
%(y/\fm)^{(j)}\ \preceq\ 1\ \text{ for $j=0,\dots,r$.}$$
If $H\subseteq \c^\infty$, then there is such $y$ in $\c^\infty$; likewise with $\c^\omega$ in place of $\c^\infty$.
\end{lemma}
\begin{proof}
Replace $(P,\fn,\hat h)$, $\fm$ by $(P_{\times\fn},1,\hat h/\fn)$, $\fm/\fn$ to arrange $\fn=1$. 
%Take $\ell\in\Calinf$ with $\ell'=\phi$. Using $(\ )^\circ$ as in Section~\ref{sec:Hardy fields} we replace $H^{\phi}$, $P^{\phi}$, $\fm$ by $H^\circ$, $(P^\phi)^\circ$, $\fm^\circ$, respectively, and rename the latter as $H$, $P$, $\phi$ to arrange also $\phi=1$. 
Then $L_{P}$ has order~$r$, $\fv(L_{P})\prec^{\flat} 1$, and $P =  Q-R$ where $Q,R\in H\{Y\}$,
$Q$ is homogeneous of degree~$1$ and order~$r$, $L_{Q}\in H[\der]$ has a  strong $\hat h$-repulsive splitting $(\phi_1,\dots,\phi_r)\in K^r$ over $K=H[\imag]$, 
and $R\prec_{\Delta^*} \fv(L_P)^{w+1} P_1$ with $\Delta^*:=\Delta\big(\fv(L_P)\big)$. 
By Lemma~\ref{lem:fv of perturbed op}(ii) we have $\fv(L_P)\sim \fv(L_Q)\asymp \fv$, so
$\Re \phi_j\succeq \fv^\dagger\succeq 1$ for $j=1,\dots,r$, and  
$\Delta = \Delta^*$.
Moreover, $P(0) \preceq \fv^{w+2}\fm^{r+1}Q$.
 Let $\nu\in\Q$, $\nu>w$, and $j\in\{1,\dots,r\}$. Then $0<v(\fv^\nu) \in v({\hat h-H})$ by Corollary~\ref{specialvariant}, so
 $\phi_j$ is $\gamma$-repulsive for $\gamma=v(\fv^\nu)$, hence~$\phi_j$ and $\phi_j-(\fv^\nu)^\dagger$ are alike by Corollary~\ref{cor:repulsive}.
 Likewise, $0 <v(\fv^\nu\fm^{r+1})\in v({\hat h-H})$ since $\hat h$ is special over $H$, so~$\phi_j$ and $\phi_j -(\fv^\nu \fm^{r+1})^\dagger$ are alike.
Thus $\phi_j-(\fv^\nu)^\dagger$ and $\phi_j-(\fv^\nu\fm^{r+1})^\dagger$ are alike as well. Hence Corollary~\ref{cor:abc} gives $y\prec \fv^w$ in~$\Calinf$ with~$P(y)=0$ and $y\in \fm\, (\c^r)^{\preceq}$. 
%$(y/\fm)^{(j)}\preceq 1$ for $j=0,\dots,r$.
For the rest use the remark following that corollary.  
\end{proof}

\begin{cor}\label{mfhc}
Suppose $\fn=1$,  and let  $\fm\in H^\times$ be such that $v\fm\in v(\hat h-H)$.  Then there are $h\in H$  and $y\in \Calinf$ such that:
$$  \hat h-h\ \preceq\ \fm,\qquad P(y)\ =\ 0,\qquad y\ \prec\ 1,\ y\in (\c^r)^{\preceq},\qquad
%y',\dots,y^{(r)}\ \preceq\ 1,\\
%  \big( (y-h)/\fm\big){}^{(j)}\ &\preceq\ 1\  \text{ for $j=0,\dots,r$.}
y-h\in \fm\,(\c^r)^{\preceq}.$$
 %\end{align*} 
If $H\subseteq \c^\infty$, then we have such $y\in\c^\infty$; likewise with $\c^\omega$ in place of $\c^\infty$. 
\end{cor}
\begin{proof}
Suppose first that $\fm\succeq 1$, and let $h:=0$ and $y$ be as in Lemma~\ref{dentsolver} for~$\phi=\fn=1$. Then $y\prec 1$, 
$y\in (\c^r)^{\preceq}$, so $y\fm\prec 1$, $y/\fm\prec 1$, $y/\fm\in (\c^r)^{\preceq}$ by the Product Rule. 
%$y^{(j)}\preceq 1$ for $j=1,\dots,r$, so  $y/\fm\prec 1$,  $(y/\fm)^{(j)}\preceq 1$ for~$j=1,\dots,r$ by the Product Rule. 
Next assume~$\fm\prec 1$ and
set $\fv:=\abs{\fv(L_{P})}\in H^>$.
By Corollary~\ref{specialvariant} we can take $h\in H$ such that $\hat h-h\prec (\fv\fm)^{(w+3)(r+1)}$, and then by
Lemma~\ref{lem:small P(a)} we have 
$$P_{+h}(0)\ =\ P(h)\ \prec\ (\fv\fm)^{w+3} P\ \preceq\ \fv^{w+3}\fm^{r+1}P_{+h}. $$
By Lemma~\ref{stronglyrepnormalrefine}, $(P_{+h},1,\hat h-h)$ is strongly repulsive-normal, 
and by Corollary~\ref{cor:deep 2, cracks} it is   deep with $\fv(L_{P_{+h}}) \asymp_{\Delta(\fv)} \fv$. Hence Lemma~\ref{bcd} applies to the slot~${(P_{+h},1,\hat h-h)}$   in place of $(P,1,\hat h)$ to yield a $z\in\Calinf$ with
$P_{+h}(z)=0$ and $(z/\fm)^{(j)}\preceq 1$ for~$j=0,\dots,r$. 
  Lemma~\ref{lem:small derivatives of y}    gives $z^{(j)}\prec 1$ for~$j=0,\dots,r$. 
%Expressing the iterates of $\der$ in terms of those of $\derdelta$ as
%before gives $z/\fm, (z/\fm)',\dots, (z/\fm)^{(r)}\preceq 1$, and then $z=(z/\fm)\fm$ gives $z, z',\dots, z^{(r)}\prec 1$. 
Set~$y:=h+z$; then $P(y)=0$, $y^{(j)} \prec 1$  and 
$\big( (y-h)/\fm\big){}^{(j)}\preceq 1$  for $j=0,\dots,r$. 
\end{proof}

\noindent
%{\it In the rest of this subsection we assume that $H$ is Liouville closed and $\I(K)\subseteq K^\dagger$.  Further assume that $\fn=1$ and our dent $(P,1,\hat h)$ in $H$ is ultimate.}\/
We now use the results above to approximate zeros of $P$ in $\Calinf$ by elements of~$H$:

\begin{cor}\label{cor:approx y} Suppose $H$ is Liouville closed, $\I(K)\subseteq K^\dagger$, $\fn=1$, and our slot~$(P,1,\hat h)$ in $H$ is ultimate. Assume also that $K$ is $1$-linearly surjective if~${r\geq 3}$. Let
$y\in\Calinf$ and $h\in H,\ \fm\in H^\times$ be such that~${P(y)=0}$, $y\prec 1$, and $\hat h-h\preceq\fm$.
Then  
$$y-h\in \fm\,(\c^r)^{\preceq}.$$
%$$\left(\frac{y-h}{\fm}\right)^{(j)}\ \preceq\ 1\ \text{ for }j=0,\dots,r.$$
%$$\big( (y-h)/\fm\big){}^{(j)} \ \preceq\ 1\ \text{  for $j=0,\dots,r$.}$$
\end{cor}
\begin{proof}
Corollary~\ref{mfhc} gives  $h_1\in H$, $z\in\Calinf$ with $\hat h-h_1\preceq\fm$, $P(z)=0$, $z\prec 1$,  and $\big( {(z-h_1)/\fm}\big){}^{(j)} \preceq 1$ for $j=0,\dots,r$.  Now 
$$\frac{y-h}{\fm}\ =\ \frac{y-z}{\fm} + \frac{z-h_1}{\fm} + \frac{h_1-h}{\fm}$$
%$$(y-h)/\fm = (y-z)/\fm + (z-h_1)/\fm + (h_1-h)/\fm$$
with 
${\big( (y-z)/\fm \big){}^{(j)} \preceq 1}$ for $j=0,\dots,r$ by Corollary~\ref{cor:notorious 3.6}.
Also $(h_1-h)/\fm\in H$ and $(h_1-h)/\fm\preceq 1$, so
$\big( (h_1-h)/\fm\big){}^{(j)}  \preceq 1$ for all $j\in\N$.
%This yields the claim.
%The dent~$(P^\phi,1,\hat h)$ in $H^\phi$ remains ultimate, $Z$-minimal, deep,  and is by assumption strongly repulsive-normal. Taking $\ell\in\Calinf$ with $\ell'=\phi$ and using $(\ )^\circ$ as before, the Hardy field~${H^\circ\supseteq\R}$ is Liouville closed and $\upo$-free, $K^\circ:=H^\circ[\imag]$  is $1$-linearly surjective if $r\geq 2$, and  $\operatorname{I}(K^\circ)\subseteq (K^\circ)^\dagger$. Thus Corollary~\ref{cor:notorious 3.6} applied to $H^\circ$, $(P^\phi)^\circ$, $y^\circ$, $z^\circ$, $\fm^\circ$ in the role of $H$, $P$, $f$, $g$, $\fm$ yields  $\big( {(y^\circ-z^\circ)/\fm^\circ} \big){}^{(j)} \preceq 1$  for $j=0,\dots,r$, so  $\derdelta^j\big( (y-z)/\fm\big) \preceq 1$  and hence   $\derdelta^j\big( (y-h)/\fm\big) \preceq 1$, $j=0,\dots,r$. 
\end{proof}

\noindent
The above corollary is the only part of this section used  towards establishing our main result, Theorem~\ref{thm:char d-max}. But this use, in proving Theorem~\ref{46}, is essential, and obtaining Corollary~\ref{cor:approx y} required much of the above section.

%\noindent
%The only use of the assumption that $K$ is $1$-linearly surjective for $r\ge 2$ in the proof above is when we appeal to Corollary~\ref{cor:notorious 3.6} to obtain $\big((y-z)/\fm\big){}^{(j)}\preceq 1$ for~$j=0,\dots,r$. This assumption is trivially satisfied for $r=1$; we can also drop this assumption for $r=2$:

%\begin{cor}\label{cor:approx y, r=2}
%Suppose $r=2$, and let $h$, $\fm$, $y$ be as in Corollary~\ref{cor:approx y}. 
%Then   
%$$\big( (y-h)/\fm\big){}^{(j)} \ \preceq\ 1\ \text{  for $j=0,\dots,r$.}$$
%\end{cor}
%\begin{proof} We follow the proof of Corollary~\ref{cor:approx y} above.
%Take $h_1$, $z$ as in the proof of that corollary.  Note that 
%$A\in H[\der]$ as in the proof of Proposition~\ref{prop:notorious 3.6} is of order~$2$ and splits over $K$, and so $\dim_{\C}\ker_{\Univ}A=2$ by Corollary~\ref{spldcr2}. Thus by the remark preceding Corollary~\ref{cor:notorious 3.6} we can still conclude that $\big((y-z)/\fm\big){}^{(j)}\preceq 1$ for~$j=0,1,2$, and then the rest of the proof of Corollary~\ref{cor:approx y} goes through.
%\end{proof}  

\section{Asymptotic Similarity} \label{sec:asymptotic similarity}

\noindent 
Let $H$ be a Hausdorff field and $\hat{H}$ an immediate valued field extension of~$H$. 
Equip~$\hat H$ with the unique field ordering making it an ordered field extension of~$H$ such that $\mathcal O_{\hat H}$  is convex [ADH, 3.5.12]. 
Let~$f\in\c$ and~$\hat{f}\in \hat{H}$  be given. 

\begin{definition} \label{p:simH}
Call $f$  {\bf asymptotically similar\/} to $\hat{f}$ over $H$ (notation: $f\sim_{H}\hat{f}$) if $f\sim \phi$  in $\c$  and~$\phi\sim \hat{f}$ in $\hat{H}$ for some $\phi\in H$. 
(Note that then $f\in\c^\times$ and $\hat f\neq 0$.)\index{germ!asymptotically similar}  
\end{definition}

\noindent
Recall that the binary relations~$\sim$ on $\c^\times$ and~$\sim$ on $\hat H^\times$ are  equivalence relations which 
restrict to the same equivalence relation on $H^\times$. As a consequence, if
$f\sim_H\hat f$, then $f\sim \phi$ in~$\c$ for any~$\phi\in H$ with~$\phi\sim \hat f$ in  $\hat H$, and $\phi\sim\hat f$ in $\hat H$ for any~$\phi\in H$ with $f\sim\phi$ in $\c$. Moreover,
if $f\in H$, then $f \sim_H \hat f \Leftrightarrow f\sim\hat f$ in~$\hat H$, and if $\hat f\in H$, then $f \sim_H \hat f \Leftrightarrow f\sim\hat f$ in $\c$.

\begin{lemma}\label{lem:simH sim}
Let $f_1\in\c$, $f_1\sim f$, let $\hat f_1\in \hat H_1$ for an immediate valued field extension $\hat H_1$ of $H$, and suppose $\hat f\sim\theta$ in $\hat H$ and $\hat f_1\sim \theta$ in $\hat H_1$ for some $\theta\in H$. 
 Then:~$f\sim_H \hat f \Leftrightarrow f_1\sim_H\hat f_1$.
\end{lemma}

\noindent
For $\fn\in H^\times$ we have $f\sim_H \hat{f}\Leftrightarrow \fn f\sim_H \fn \hat{f}$. Moreover, by Lemma~\ref{lem:sim props}:

\begin{lemma}\label{lem:simH}
Let $g\in\c$, $\hat g\in\hat H$, and suppose $f\sim_H\hat f$ and $g\sim_H\hat g$. Then ${1/f\sim_H 1/\hat f}$ and $fg\sim_H\hat f\,\hat g$. Moreover,
$$f\preceq g \text{ in $\c$} \quad\Longleftrightarrow\quad \hat f\preceq \hat g \text{ in $\hat H$,}$$
and likewise with $\prec$, $\asymp$, or $\sim$ in place of $\preceq$.
\end{lemma}

\noindent
Lemma~\ref{lem:simH} readily yields:

{\samepage 
\begin{cor}\label{cor:simHt}
Suppose
$\hat f$ is transcendental over $H$ and  $Q(f)\sim_H Q(\hat f)$ for all~$Q\in H[Y]^{\neq}$. 
Then we have:
\begin{enumerate}
\item[\textup{(i)}]  a subfield $H(f)\supseteq H$ of $\c$ generated by $f$ over~$H$;
\item[\textup{(ii)}] a field isomorphism $\iota\colon H(f)\to H(\hat f)$ over $H$ with~$\iota(f)=\hat f$;  
\item[\textup{(iii)}] with $H(f)$ and $\iota$ as in \textup{(i)} and \textup{(ii)} we have $g\sim_H\iota(g)$ for all $g\in H(f)^{\times}$, hence for all $g_1,g_2\in H(f)$: 
$g_1\preceq g_2$  in~$\c$~$\Leftrightarrow$ $\iota(g_1)\preceq \iota(g_2)$  in~$\hat{H}$.
%$$g_1\preceq g_2 \text{ in $\c$}\quad\Longleftrightarrow\quad \iota(g_1)\preceq\iota(g_2)\text{ in $\hat H$.}$$
\end{enumerate}
Also, $\iota$ in \textup{(ii)} is unique and is an ordered field isomorphism, where  the ordering on~$H(f)$ is its ordering as a Hausdorff field.
\end{cor}}
\begin{proof} To see that $\iota$ is order preserving, use that $\iota$ is a valued field isomorphism by~(iii), 
and apply [ADH, 3.5.12].
\end{proof}

\noindent
Here is the analogue when $\hat f$ algebraic over $H$: 

\begin{cor}\label{cor:simHa}
Suppose  $\hat f$ is algebraic over $H$ with minimum polynomial~$P$ over~$H$ of degree $d$, and
$P(f)=0$, $Q(f)\sim_H Q(\hat f)$ for all~$Q\in H[Y]^{\neq}$ of degree~$<d$.
Then we have:
\begin{enumerate}
\item[\textup{(i)}]  a subfield $H[f]\supseteq H$ of $\c$ generated by $f$ over~$H$;
\item[\textup{(ii)}] a field isomorphism $\iota\colon H[f]\to H[\hat f]$ over $H$ with~$\iota(f)=\hat f$;  
\item[\textup{(iii)}] with $H[f]$ and $\iota$ as in \textup{(i)} and \textup{(ii)} we have $g\sim_H\iota(g)$ for all $g\in H[f]^{\times}$, hence for all $g_1,g_2\in H[f]$: 
$g_1\preceq g_2$  in~$\c$~$\Leftrightarrow$ $\iota(g_1)\preceq \iota(g_2)$  in~$\hat{H}$.
\end{enumerate}
Also, $H[f]$ and $\iota$ in \textup{(i)} and \textup{(ii)} are unique and $\iota$ is an ordered field isomorphism, where  the ordering on $H(f)$ is its ordering as a Hausdorff field.
 \end{cor}

\noindent
If $\hat f\notin H$, then to show that $f-\phi\sim_H \hat f-\phi$ for all $\phi\in H$ it is enough to do this for $\phi$ arbitrarily close to $\hat f$: 

\begin{lemma}\label{lem:simH phi0}
Let $\phi_0\in H$ be such that $f-\phi_0\sim_H\hat f-\phi_0$. Then $f-\phi\sim_H \hat f-\phi$ for all $\phi\in H$ with
$\hat f-\phi_0 \prec \hat f-\phi$.
\end{lemma}
\begin{proof}
Let $\phi\in H$ with $\hat f-\phi_0 \prec \hat f-\phi$. Then $\phi_0-\phi\succ \hat f-\phi_0$, so $\hat f-\phi=(\hat f-\phi_0)+(\phi_0-\phi)\sim\phi_0-\phi$. By Lemma~\ref{lem:simH} we also have $\phi_0-\phi\succ f-\phi_0$, and hence likewise
$f-\phi \sim\phi_0-\phi$. 
\end{proof}

\noindent
We define: $f\approx_H\hat f: \Leftrightarrow f-\phi \sim_H \hat f-\phi$ for all $\phi\in H$. \label{p:approxH}
If $f\approx_H\hat f$, then~$f\sim_H\hat f$ as well as $f,\hat f\notin H$, and 
$\fn f\approx_H \fn\hat f$ for all $\fn\in H^\times$. Hence
$f\approx_H \hat f$ iff $f, \hat f\notin H$ and the isomorphism $\iota\colon H+H f\to H+H\hat f$ of $H$-linear spaces that is the identity on~$H$ and sends $f$ to $\hat f$   satisfies~$g\sim_H \iota(g)$ for all nonzero $g\in H+Hf$.

\medskip
\noindent
Here is an easy consequence of Lemma~\ref{lem:simH phi0}: 

\begin{cor}\label{cor:simH phi0}
Suppose $\hat f\notin H$ and $f-\phi_0 \sim_H \hat f-\phi_0$ for all  $\phi_0\in H$ such that~$\phi_0\sim\hat f$. Then~${f\approx_H \hat f}$.
\end{cor}
\begin{proof}
Take $\phi_0\in H$ with $\phi_0\sim\hat f$, and
let $\phi\in H$ be given. If~${\hat f-\phi\prec\hat f}$, then~${f-\phi \sim_H \hat f-\phi}$ by hypothesis;
otherwise we have~$\hat f-\phi\succeq\hat f\succ \hat f-\phi_0$, and then~${f-\phi} \sim_H \hat f-\phi$ by Lemma~\ref{lem:simH phi0}.
\end{proof}

\noindent
Lemma~\ref{lem:simH sim} yields an analogue  for $\approx_H$:

\begin{lemma}\label{lem:approxH sim}
Let $f_1\in\c$ be such that $f_1-\phi\sim f-\phi$ for all $\phi\in H$, 
and let $\hat f_1$ be an element of an immediate valued field extension of $H$  such that
$v(\hat f-\phi)=v(\hat f_1-\phi)$ for all $\phi\in H$. Then~$f\approx_H \hat f$ iff $f_1\approx_H\hat f_1$.
\end{lemma}

\noindent
Let $g\in\c$ be eventually strictly increasing with $g(t)\to+\infty$ as $t\to+\infty$; we then have the Hausdorff field
$H\circ g=\{h\circ g:h\in H\}$, with ordered valued field isomorphism~$h\mapsto h\circ g\colon H\to H\circ g$. (See Section~\ref{sec:germs}.)
Suppose  
$$\hat h\mapsto \hat h\circ g\ \colon\  \hat H\to\hat H\circ g$$  extends this isomorphism to
a valued field isomorphism, where $\hat H\circ g$ is an immediate valued field extension of the Hausdorff field $H\circ g$.
Then
$$f\sim_H \hat f\quad\Longleftrightarrow\quad f\circ g\sim_{H\circ g} \hat f\circ g,\qquad
f\approx_H \hat f\quad\Longleftrightarrow\quad f\circ g\approx_{H\circ g} \hat f\circ g.$$

\subsection*{The complex version}
We now assume in addition that~$H$ is real closed, with algebraic closure 
$K:=H[\imag]\subseteq\c[\imag]$. We take $\imag$  with $\imag^2=-1$ also as an element of a field
$\hat K:=\hat H[\imag]$ extending both $\hat{H}$ and $K$, and equip $\hat K$ with the unique valuation ring of $\hat K$ lying over
$\mathcal O_{\hat H}$; see the remarks following Lemma~\ref{lem:dotK real closed}. 
Then~$\hat K$ is an immediate valued field extension of~$K$.
Let~$f\in\c[\imag]$ and~$\hat f\in \hat K$ below.  \label{p:simK}

\medskip
\noindent
Call~$f$  {\bf asymptotically similar\/}\index{germ!asymptotically similar} to $\hat{f}$ over~$K$ (notation: $f\sim_{K}\hat{f}$) if for some~${\phi\in K}$ we have $f\sim \phi$  in $\c[\imag]$  and~$\phi\sim \hat{f}$ in~$\hat{K}$. 
Then~$f\in\c[\imag]^\times$ and~$\hat f\neq 0$. As before, 
if~$f\sim_K\hat f$, then $f\sim \phi$ in~$\c[\imag]$ for any~${\phi\in K}$ for which~$\phi\sim \hat f$ in~$\hat K$, and~$\phi\sim\hat f$ in $\hat K$ for any $\phi\in K$ for which~$f\sim\phi$ in~$\c[\imag]$.
Moreover, if~$f\in K$, then~$f \sim_K \hat f$ reduces to $f \sim \hat f$ in $\hat K^\times$. Likewise, if $\hat f\in K$, then
$f \sim_K \hat f$ reduces to $f \sim \hat f$ in~$\c[\imag]^\times$.

\begin{lemma}\label{lem:simH sim K}
Let $f_1\in\c[\imag]$ with $f_1\sim f$. Let  $\hat H_1$ be an immediate valued field extension of $H$, let $\hat K_1:=\hat H_1[\imag]$ be the corresponding immediate valued field extension of $K$ obtained from $\hat H_1$ as $\hat K$ was obtained from $\hat H$. Let  $\hat f_1\in \hat K_1$, and $\theta\in K$ be such that 
$\hat f\sim\theta$ in $\hat K$ and $\hat f_1\sim \theta$ in~$\hat K_1$. Then~$f\sim_K \hat f$ iff $f_1\sim_K\hat f_1$.
\end{lemma}

\noindent
For~$\fn\in K^\times$ we have $f\sim_K \hat{f}\Leftrightarrow \fn f\sim_K \fn \hat{f}$, and~$f\sim_K \hat f\Leftrightarrow\overline{f} \sim_K \overline{\hat f}$ (complex conjugation).
Here is a useful observation relating $\sim_K$ and $\sim_H$:

\begin{lemma}\label{lem:simK Re Im} Suppose $f\sim_K \hat{f}$ and 
$\Re \hat f\succeq \Im\hat f$; then $$\Re f\ \succeq\ \Im f,\qquad \Re f\ \sim_H\ \Re \hat f.$$
\end{lemma}
\begin{proof} Let $\phi\in K$ be such that $f\sim \phi$ in $\c[\imag]$ and
$\phi\sim \hat f$ in $\hat K$. The latter yields~$\Re\phi\succeq \Im\phi$ in $H$ and $\Re\phi\sim \Re\hat f$ in $\hat H$. Using
that $f=(1+\varepsilon)\phi$ with~$\varepsilon\prec 1$ in~$\c[\imag]$
it follows easily that $\Re f\succeq \Im f$ and $\Re f\sim \Re \phi$ in $\c$. 
\end{proof}

\begin{cor}\label{corsimas}
Suppose $f\in\c$ and $\hat f\in\hat H$. Then $f\sim_H \hat f$ iff $f\sim_K \hat f$.
\end{cor}

\noindent
Lemmas~\ref{lem:simH} and~\ref{lem:simH phi0} go through with~$\c[\imag]$,~$K$,~$\hat K$, and~$\sim_K$ in place of $\c$, $H$, $\hat H$,  and $\sim_H$. 
We define:  $f\approx_K\hat f :\Leftrightarrow f-\phi \sim_K \hat f-\phi$ for all $\phi\in K$. 
Now Corollary~\ref{cor:simH phi0} goes through  with $K$, $\sim_K$, $\approx_K$ in place of~$H$,~$\sim_H$,~$\approx_H$. 
%From Lemma~\ref{lem:simH sim K} we obtain a ``complex'' version of Lemma~\ref{lem:approxH sim}:

%\begin{lemma}\label{lem:approxH sim K}
%Let $\hat H_1$ and $\hat K_1$ be as in Lemma~\ref{lem:simH sim K}. Let $f_1\in\c[\imag]$ and $\hat f_1\in \hat K_1$ be such that $f_1-\phi\sim f-\phi$ and 
%$v(\hat f-\phi)=v(\hat f_1-\phi)$ for all $\phi\in K$. Then
%$$f\approx_K \hat f\ \Longleftrightarrow\ f_1\approx_K \hat f_1.$$
%\end{lemma}

\begin{lemma}\label{lemsimimag}
Suppose $f\in\c$, $\hat f\in\hat H$, and $f\sim_H\hat f$.
Then $f+g\imag \sim_K \hat f+g\imag$ for all $g\in H$.
\end{lemma}
\begin{proof}
Let $g\in H$, and take $\phi\in H$ with $f\sim \phi$ in $\c$ and $\phi \sim \hat f$ in $\hat H$. 
Suppose first that $g\prec\phi$. Then $g\imag\prec\phi$, and together with $f-\phi\prec\phi$ this yields~${(f+g\imag)-\phi\prec\phi}$, that is,
$f+g\imag \sim \phi$ in $\c[\imag]$ (cf.~the basic properties of the relation $\prec$ on $\c[\imag]$ stated before Lemma~\ref{lem:sim props}).
Using likewise the analogous properties of $\prec$ on $\hat K$ we obtain~$\phi\sim\hat f+g\imag$ in $\hat K$. If $\phi\prec g$, then~$f\preceq \phi\prec g\imag$ and thus $f+g\imag \sim g\imag$ in~$\c[\imag]$, and
likewise $\hat f+g\imag \sim g\imag$ in $\hat K$.
Finally, suppose $g\asymp\phi$. Take $c\in\R^\times$ and~$\varepsilon\in H$ with $g = c\phi(1+\varepsilon)$ and~$\varepsilon\prec 1$.
We have
$f=\phi(1+\delta)$ where $\delta\in\c$, $\delta\prec 1$, so~$f+g\imag = \phi(1+c\imag)(1+\rho)$ where~$\rho=(1+c\imag)^{-1}(\delta+c\imag\varepsilon) \prec 1$ in $\c[\imag]$, so~$f+g\imag\sim \phi(1+c\imag)$ in $\c[\imag]$.  
Likewise, $\hat f+g\imag\sim \phi(1+c\imag)$ in $\hat K$. 
%In all three cases we obtain  $f+g\imag \sim_K \hat f+g\imag$.
\end{proof}
 
\begin{cor}\label{cor:approxH vs approxK}
Suppose $f\in\c$ and $\hat f\in\hat H$. Then $f\approx_H \hat f$ iff $f\approx_K\hat f$.
\end{cor}
\begin{proof}
If $f\approx_K\hat f$, then for all $\phi\in H$ we have $f-\phi\sim_K \hat f-\phi$, so $f-\phi\sim_H \hat f-\phi$ by Corollary~\ref{corsimas},
hence $f\approx_H \hat f$. Conversely, suppose $f\approx_H \hat f$. Then for all $\phi\in K$
we have $f-\Re \phi \sim_H \hat f-\Re \phi$, so
$f-\phi \sim_K \hat f-\phi$ by Lemma~\ref{lemsimimag}.
\end{proof}

\noindent
Next we exploit that $K$ is algebraically closed: 

\begin{lemma}
$f\approx_K \hat f \ \Longrightarrow\ Q(f) \sim_K Q(\hat f)$ for all $Q\in K[Y]^{\neq}$. 
\end{lemma} 
\begin{proof}
Factor $Q\in K[Y]^{\ne}$ as 
$$Q=a(Y-\phi_1)\cdots (Y-\phi_n),\qquad a\in K^\times,\  \phi_1,\dots,\phi_n\in K$$
and use $f-\phi_j \sim_K \hat f-\phi_j$  ($j=1,\dots,n$) and the complex version of Lemma~\ref{lem:simH}.
\end{proof}

\noindent
This yields a more useful ``complex'' version of Corollary~\ref{cor:simHt}: 

\begin{cor}\label{complexsimHt} 
Suppose $f\approx_K\hat f$. Then $\hat f$ is transcendental over $K$, and:
\begin{enumerate}
\item[\textup{(i)}]  $f$ generates over $K$ a subfield $K(f)$ of $\c[\imag]$;
\item[\textup{(ii)}] we have a field isomorphism $\iota\colon K(f)\to K(\hat f)$ over $K$ with~$\iota(f)=\hat f$;  
\item[\textup{(iii)}]  $g\sim_K\iota(g)$ for all $g\in K(f)^{\times}$, hence  for all $g_1,g_2\in K(f)$:  $$g_1\preceq g_2 \text{  in }\c[\imag]\ \Longleftrightarrow\ \iota(g_1)\preceq \iota(g_2)  \text{ in }\hat{K}.$$
%$$g_1\preceq g_2 \text{ in~$\c[\imag]$} \quad \Longleftrightarrow \quad \iota(g_1)\preceq\iota(g_2) \text{ in $\hat K$.}$$
\end{enumerate}
\textup{(}Thus the restriction of the binary relation $\preceq$ on $\c[\imag]$ to~$K(f)$ is a dominance relation on the field $K(f)$ in the sense of
\textup{[ADH, 3.1.1].)}
\end{cor}

\noindent 
In the next lemma $f=g+h\imag$, $g,h\in \c$, and $\hat f=\hat g+\hat h \imag$, $\hat g, \hat h\in \hat H$. 
Recall from Lemma~\ref{lem:same width} that if $\hat f\notin K$, then  $v(\hat g-H) \subseteq v(\hat h-H)$ or  $v(\hat h-H) \subseteq v(\hat g-H)$.

\begin{lemma}\label{lem:real part approx}
Suppose $f\approx_K \hat f$ and $v(\hat g-H) \subseteq v(\hat h-H)$. Then $g \approx_H  \hat g$.
\end{lemma}
\begin{proof}
Let $\rho\in H$ be such that $\rho\sim\hat g$; by Corollary~\ref{cor:simH phi0} it is enough to show that
then~${g-\rho} \sim_H \hat g-\rho$.
Take $\sigma\in H$ with $\hat g - \rho \succeq \hat h-\sigma$, and set $\phi:=\rho+\sigma\imag\in K$.
Then $$\Re(f-\phi)=g-\rho\quad\text{ and }\quad \Re(\hat f-\phi)=\hat g-\rho\succeq  \hat h-\sigma=\Im(\hat f-\phi),$$
and so by~$f-\phi\sim_H \hat f -\phi$ and Lemma~\ref{lem:simK Re Im} we have $g-\rho\sim_H \hat g-\rho$.
\end{proof}

\begin{cor}
If $f\approx_K \hat f$, then $\Re f \approx_H \Re \hat f$ or $\Im f\approx_H \Im\hat f$.
\end{cor}
 
\noindent
Let $g\in\c$ be eventually strictly increasing with $g(t)\to+\infty$ as $t\to+\infty$; we then have the subfield
$K\circ g=(H\circ g)[\imag]$ of $\c[\imag]$. Suppose the valued field isomorphism $$h\mapsto h\circ g\colon H\to H\circ g$$ is
extended to a valued field isomorphism  
$$\hat h\mapsto \hat h\circ g\ \colon\  \hat H\to\hat H\circ g,$$ 
 where $\hat H\circ g$ is an immediate valued field extension of the Hausdorff field $H\circ g$. 
In the same way we took a common valued field extension $\hat K=\hat H[\imag]$ of $\hat H$ and~$K=H[\imag]$ we now take a common valued field extension $\hat K\circ g=(\hat H\circ g)[\imag]$ of $\hat H\circ g$ and~$K\circ g=(H\circ g)[\imag]$.
This makes $\hat K\circ g$ an immediate extension of $K\circ g$, and we have a %\marginpar{need a cubical commuting diagram to indicate the various inclusions and maps} 
unique valued field isomorphism~$y\mapsto y\circ g\colon \hat K\to \hat K\circ g$ extending the above
map~$\hat h\mapsto \hat h\circ g\colon \hat H\to\hat H\circ g$ and sending $\imag\in \hat{K}$ to $\imag\in \hat K\circ g$.
This map $\hat K\to \hat K\circ g$ also extends $f\mapsto f\circ g\colon K \to K\circ g$ and is the identity on $\C$. 
See the commutative diagram below, where the labeled arrows are valued field isomorphisms and all unlabeled arrows are natural inclusions.

\begin{equation*}
\xymatrix{%
&\hat{H}\circ g\ar@{->}[rr]\ar@{<-}[dl]_(0.6){\hat h\mapsto \hat h\circ g}  \ar@{<-}[dd]|!{[d];[d]}\hole && \hat{K}\circ g  \ar@{<-}[dd] \ar@{<-}[dl]_(0.6){y\mapsto y\circ g} \\
\hat{H}\ar@{->}[rr]\ar@{<-}[dd]&&\hat{K} \ar@{<-}[dd] \\ 
&H\circ g\ar@{->}[rr]|!{[r];[r]}\hole\ar@{<-}[dl]^(0.45){h\mapsto h\circ g} && K\circ g 
\ar@{<-}[dl]^(0.45){f\mapsto f\circ g} \\
H\ar@{->}[rr]&&K
}
\end{equation*}

\noindent
Now
we have
$$f\sim_K \hat f\quad\Longleftrightarrow\quad f\circ g\sim_{K\circ g} \hat f\circ g,\qquad
f\approx_K \hat f\quad\Longleftrightarrow\quad f\circ g\approx_{K\circ g} \hat f\circ g.$$
At various places in the next section we use this for a Hardy field $H$ and active~$\phi>0$ in $H$, with
$g=\ell^{\inv}$, $\ell\in \c^1,\ \ell'=\phi$. In that situation, $H^\circ:=H\circ g$, $\hat H^\circ:=\hat H\circ g$,
and $h^\circ:= h\circ g$, $\hat h^\circ := \hat h\circ g$ for $h\in H$ and
$\hat h\in \hat H$, and likewise with $K$ and $\hat K$ and their elements instead of $H$ and $\hat H$.

\section{Differentially Algebraic Hardy Field Extensions} \label{sec:d-alg extensions}

\noindent
In this section we are finally able to generate under reasonable conditions Hardy field extensions by solutions in $\Calinf$ of algebraic differential equations, culminating in the proof of our main theorem.
We begin with a generality about enlarging differential fields within an ambient differential ring.
Here, a {\em differential subfield\/} of a differential ring  $E$ is a differential subring of $E$ whose underlying ring is a field. 
 
\begin{lemma}\label{dsfe} Let $K$ be a differential field with irreducible $P\in K\{Y\}^{\ne}$ of order~$r\ge 1$, and $E$ a differential ring extension of $K$ with  $y\in E$ such that $P(y)=0$ and $Q(y)\in E^\times$ for all $Q\in K\{Y\}^{\neq}$
of order $<r$. Then $y$ generates over $K$ a differential subfield $K\langle y \rangle\supseteq K$ of $E$. Moreover, $y$ has $P$ as a minimal annihilator over $K$ and $K\langle y \rangle$ equals
$$\left\{\frac{A(y)}{B(y)}:\ A,B\in K\{Y\},\, \order A\leq r,\, \deg_{Y^{(r)}}A <\deg_{Y^{(r)}} P,\, B\ne 0,\, \order B < r\right\}.$$
\end{lemma}
\begin{proof} Let $p\in K[Y_0,\dots, Y_r]$ with distinct indeterminates $Y_0,\dots, Y_r$ be such that $P(Y)=p(Y, Y',\dots, Y^{(r)})$. The $K$-algebra morphism 
$K[Y_0,\dots, Y_r]\to E$ sending $Y_i$ to $y^{(i)}$ for $i=0,\dots,r$
extends to a $K$-algebra morphism $K(Y_0,\dots, Y_{r-1})[Y_r]\to E$  with $p$ in its kernel, and so induces a $K$-algebra morphism 
$$\iota\ :\ K(Y_0,\dots, Y_{r-1})[Y_r]/(p)\to E, \qquad (p)\ :=\ pK(Y_0,\dots, Y_{r-1})[Y_r].$$ 
Now $p$  as an element of $K(Y_0,\dots, Y_{r-1})[Y_r]$ remains irreducible \cite[Chapter~IV, \S{}2]{Lang}. Thus $K(Y_0,\dots, Y_{r-1})[Y_r]/(p)$ is a field, so $\iota$ is injective, and it is routine to check that the image of $\iota$ is $K\langle y \rangle$ as described; see also [ADH, 4.1.6]. 
\end{proof} 

\noindent
In passing we also note  the obvious $\d$-transcendental version of this lemma:

\begin{lemma}\label{dsfe, d-trans}
Let $K$ be a differential field and $E$ be a differential ring extension of $K$ with $y\in E$ such that $Q(y)\in E^\times$
for all $Q\in K\{Y\}^{\neq}$.  Then $y$ generates over~$K$ a differential subfield $K\langle y \rangle$ of $E$. Moreover,  $y$  is $\d$-transcendental over $K$ and 
$$K\<y\>\ =\ \left\{\frac{P(y)}{Q(y)}:\ P,Q\in K\{Y\},\ Q\neq 0\right\}.$$
\end{lemma}

\noindent
We now apply the material above to generate Hardy field extensions.

\subsection*{Application to Hardy fields}  
 {\it In the rest of this section 
$H$ is a real closed Hardy field, $H\supseteq \R$, and~$\hat{H}$ is an immediate $H$-field extension of $H$.}\/
Let $f\in\Calinf$ and~$\hat f\in\hat H$.
Note that if~${Q\in H\{Y\}}$ and $Q(f)\sim_H Q(\hat f)$, then $Q(f)\in\c^\times$.
Hence by Lemma~\ref{dsfe}  with~$E=\Calinf$, $K=H$, we have:

\begin{lemma}\label{lem:dsfeh}
Suppose $\hat{f}$ is $\operatorname{d}$-algebraic over $H$ with minimal annihilator $P$ over~$H$ of order~$r\ge 1$,
and $P(f)=0$ and $Q(f)\sim_H Q(\hat f)$ for all $Q\in H\{Y\}\setminus H$ with ${\order Q < r}$. Then $f\notin H$ and: \begin{enumerate}
\item[\rm{(i)}] $f$ is hardian over $H$; 
\item[\rm{(ii)}] we have a \textup{(}necessarily unique\textup{)}  isomorphism $\iota\colon H\langle f\rangle \to H\langle\hat f\rangle$ of differential fields over $H$ such that $\iota(f)=\hat{f}$.
\end{enumerate}
%Moreover, $H\langle f\rangle$ is a Hardy field.   
\end{lemma}

\noindent
With an extra assumption $\iota$ in Lemma~\ref{lem:dsfeh} is an isomorphism of $H$-fields:

\begin{cor}\label{cor:dsfeh} 
Let $\hat f$, $f$, $P$, $r$, $\iota$ be as in Lem\-ma~\ref{lem:dsfeh}, and suppose also that~${Q(f)\sim_{H} Q(\hat{f})}$ for all~$Q\in H\{Y\}$ with $\order Q = r$ and $\deg_{Y^{(r)}} Q < \deg_{Y^{(r)}} P$.
Then  $g\sim_H \iota(g)$ for all~$g\in H\langle f\rangle^\times$,
hence for
${g_1,g_2\in H\langle f\rangle}$ we have $$g_1\preceq g_2 \text{  in }\c\ \Longleftrightarrow\
\iota(g_1)\preceq \iota(g_2)  \text{  in }\hat{H}.$$
Moreover, $\iota$ is an isomorphism of $H$-fields.
\end{cor}
\begin{proof}
Most of this follows from Lemmas~\ref{lem:simH} and~\ref{lem:dsfeh} and the description of~$H\langle f\rangle$ in Lemma~\ref{dsfe}. For the last statement, use [ADH, 10.5.8]. 
\end{proof}

\noindent
Here is a $\d$-transcendental  version of Lemma~\ref{lem:dsfeh}:

{\samepage
\begin{lemma}\label{dsfeh, d-trans}
Suppose $Q(f) \sim_H Q(\hat f)$ for 
all $Q\in H\{Y\}\setminus H$. Then:
\begin{enumerate}
\item[\textup{(i)}] $f$ is hardian over $H$;
\item[\textup{(ii)}] we have a \textup{(}necessarily unique\textup{)}  isomorphism $\iota\colon H\langle f\rangle
\to H\langle\hat f\rangle$ of differential fields over $H$ with $\iota(f)=\hat f$; and
\item[\textup{(iii)}]  $g\sim_H \iota(g)$ for all~$g\in H\langle f\rangle^\times$,
hence for
 all $g_1,g_2\in H\langle f\rangle$: 
 $$g_1\preceq g_2 \text{  in }\c\ \Longleftrightarrow\
\iota(g_1)\preceq \iota(g_2)  \text{  in }\hat{H}.$$
\end{enumerate}
Moreover, $\iota$ is an isomorphism of $H$-fields.
\end{lemma}}

\noindent
This follows easily from Lemma~\ref{dsfe, d-trans}.

\subsection*{Analogues for $K=H[\imag]$} We have the $\d$-valued extension~$K:= H[\imag]\subseteq \Calinf[\imag]$ of $H$. As before we arrange that $\hat{K}=\hat{H}[\imag]$ is a $\d$-valued extension of $\hat{H}$ as well as an
an immediate extension of $K$. Let~${f\in\Calinf[\imag]}$ and~$\hat{f}\in \hat{K}$.
We now have the obvious ``complex'' analogues of Lemma~\ref{lem:dsfeh} and Corollary~\ref{cor:dsfeh}: % and Lemma~\ref{dsfeh, d-trans}:

\begin{lemma}\label{lem:dsfek}
Suppose $\hat{f}$ is $\d$-algebraic over $K$ with minimal annihilator $P$ over~$K$ of order~$r\ge 1$,
and $P(f)=0$ and $Q(f)\sim_K Q(\hat f)$ for all $Q\in K\{Y\}\setminus K$ with~${\order Q < r}$. Then \begin{enumerate}
\item[\rm{(i)}] $f$ generates over $K$ a differential subfield $K\langle f\rangle$ of $\Calinf[\imag]$;
\item[\rm{(ii)}] we have a \textup{(}necessarily unique\textup{)} isomorphism $\iota\colon K\langle f\rangle \to K\langle\hat f\rangle$ of differential fields over $K$ such that $\iota(f)=\hat{f}$.
\end{enumerate}
\end{lemma}

\begin{cor}\label{cor:dsfek} 
Let $\hat f$, $f$, $P$, $r$, $\iota$ be as in Lem\-ma~\ref{lem:dsfek}, and   suppose also that~$Q(f)\sim_K Q(\hat{f})$ for all $Q\in K\{Y\}$ with $\order Q = r$ and $\deg_{Y^{(r)}} Q < \deg_{Y^{(r)}} P$. 
Then $g\sim_K \iota(g)$ for all $g\in K\langle f\rangle^\times$, so for all $g_1,g_2\in K\langle f\rangle$ we have: 
$$g_1\preceq g_2  \text{ in }\c[\imag]\ \Longleftrightarrow\ \iota(g_1)\preceq \iota(g_2)  \text{ in }\hat{K}.$$
Thus the relation~$\preceq$ on $\c[\imag]$ restricts to a dominance relation on the field $K\langle f\rangle$.
%and equipping $K\<f\>$ with the associated valuation, $\iota$ is a valued field isomorphism.
\end{cor}

%\begin{lemma}\label{dsfek, d-trans}
%Suppose $Q(f) \sim_K Q(\hat f)$ for  all $Q\in K\{Y\}\setminus K$. Then
%\begin{enumerate}
%\item[\textup{(i)}] we have a \textup{(}necessarily unique\textup{)} differential subfield $K\langle f\rangle\supseteq K$ of $\Calinf[\imag]$ that is generated as a differential field over $K$ by $f$;
%\item[\textup{(ii)}] we have a \textup{(}necessarily unique\textup{)} isomorphism $\iota\colon K\langle f\rangle \to H\langle\hat f\rangle$ of differential $K$-algebras with $\iota(f)=\hat f$; and
%\item[\textup{(iii)}] for all $g_1,g_2\in K\langle f\rangle$ we have: $g_1\preceq g_2$ in $\c[\imag]$ $\Leftrightarrow$ $\iota(g_1)\preceq\iota(g_2)$ in $\hat K$.
%\end{enumerate}
%Moreover, $\preceq$ restricts to a dominance relation on $K\< f \>$.
%\end{lemma}

\noindent
From $K$ being algebraically closed we obtain a useful variant of Corollary~\ref{cor:dsfek}: 

\begin{cor}\label{cor:dsfek, order 1}
Suppose $f\approx_K \hat f$, and $P\in K\{Y\}$ is irreducible with 
$$\order P\ =\ \deg_{Y'} P\  =\ 1, \quad P(f)=0\ \text{ in }\ \Calinf[\imag],\quad P(\hat f)=0\ \text{ in }\ \hat K.$$
Then $P$ is a minimal annihilator of $\hat f$ over~$K$, $f$ generates over $K$ a differential subfield~$K\<f\>=K(f)$ of~$\Calinf[\imag]$, and we have an isomorphism $\iota\colon K\<f\>\to K\<\hat f\>$
of differential fields over $K$ such that~${\iota(f)=\hat f}$ and~$g\sim_K\iota(g)$ for all $g\in K\langle f\rangle^\times$. Thus for all $g_1,g_2\in K\langle f\rangle$: 
$g_1\preceq g_2  \text{ in }\c[\imag]\ \Longleftrightarrow\ \iota(g_1)\preceq \iota(g_2)  \text{ in }\hat{K}.$
%Moreover, $\preceq$ restricts to a dominance relation on the field~$K\<f\>=K(f)$, and equipping $K\<f\>$ with the associated valuation, $\iota$ is a valued field isomorphism.
\end{cor}
\begin{proof} By Corollary~\ref{complexsimHt}, $\hat f$ is transcendental over $K$, so $P$ is a minimal annihilator
of $\hat f$ over $K$ by [ADH, 4.1.6]. Now use Lemma~\ref{dsfe} and Corollary~\ref{complexsimHt}. 
\end{proof}

\noindent
This corollary leaves open whether $\Re f$ or $\Im f$ is hardian over $H$. This issue is critical for us and we treat a special case in Proposition~\ref{prop:Re or Im} below.  The example following Corollary~\ref{cor:Bosh 1.3}  shows that that there is a 
differential subfield of $\Calinf[\imag]$ such that the
binary relation~$\preceq$ on $\c[\imag]$ restricts to a dominance relation on it, but which is not contained in $F[\imag]$ for any
Hardy field $F$.
%; cf.~the example following Lemma~\ref{lem:y H-hardian crit, complex}.

\subsection*{Sufficient conditions for asymptotic similarity}
Let~$\hat h$ be an element of our immediate $H$-field extension $\hat H$ of $H$.
% Let $Q\in H\{Y\}\setminus H$ and $Q\notin Z(H,\hat h)$.
%Then there exist $h\in H$ and $\fv\in H^\times$ such that 
%$$h-\hat h\ \prec\ \fv, \qquad \ndeg_{\prec \fv}Q_{+h}\ =\ 0.$$
Note that in the next variant of [ADH, 11.4.3] we use $\operatorname{ddeg}$ instead of $\ndeg$.

\begin{lemma}\label{1142} Let $Q\in H\{Y\}^{\neq}$, $r:=\order Q$, $h\in H$, and
$\fv\in H^\times$ be such that~$\hat h-h\prec \fv$ and $\ddeg_{\prec \fv} Q_{+h}=0$, and
assume $y\in \Calinf$ and $\fm\in H^\times$  satisfies
$$y-h\preceq\fm\prec\fv,\qquad \left(\frac{y-h}{\fm}\right)',\dots,\left(\frac{y-h}{\fm}\right)^{(r)} \preceq 1.$$ Then $Q(y)\sim Q(h)$ in $\Calinf$ and $Q(h)\sim Q(\hat h)$ in $\hat H$; in particular, $Q(y)\sim_H Q(\hat h)$.
\end{lemma} 
\begin{proof} We have $y=h+\fm u$ with $u=\frac{y-h}{\fm}\in \Calinf$ and
$u, u',\dots, u^{(r)}\preceq 1$.
Now 
$$Q_{+h,\times \fm}\ =\ Q(h) + R\qquad\text{ with }R\in H\{Y\},\ R(0)=0,$$ 
which in view of $\operatorname{ddeg} Q_{+h,\times \fm}=0$ gives $R\prec Q(h)$.
Thus 
$$Q(y)\ =\ Q_{+h,\times\fm}(u)\ =\ Q(h) + R(u),\qquad R(u)\ \preceq\ R\ \prec\ Q(h),$$ so $Q(y)\sim Q(h)$ in $\Calinf$. Increasing $|\fm|$ if necessary we arrange $\hat h - h\preceq \fm$, and then a similar computation with $\hat h$ instead of $y$ gives $Q(h)\sim Q(\hat h)$ in $\hat H$.  
\end{proof}

\noindent
{\it In the remainder of this subsection  we assume that $H$ is ungrounded and $H\ne \R$.}\/
 
\begin{cor}  \label{cor:iso hat h, d-alg, weak}
Suppose $\hat h$ is $\d$-algebraic over $H$ with minimal annihilator $P$ over~$H$ of order $r\geq 1$, and let~${y\in\Calinf}$ satisfy $P(y)=0$. Suppose for all~$Q$ in~$H\{Y\}\setminus H$ of order $< r$  there are $h\in H$, $\fm,\fv\in H^\times$, and an active~$\phi>0$ in~$H$ such that $\hat h-h\preceq\fm\prec\fv$,  $\ddeg_{\prec\fv} Q^\phi_{+h}=0$, and
$$ \derdelta^j\!\left(\frac{y-h}{\fm}\right)\preceq 1\qquad\text{ for }j=0,\dots,r-1\text{ and }\derdelta:=\phi^{-1}\der.$$
Then $y\notin H$ and $y$ is hardian over $H$.
\end{cor} 
\begin{proof}
Let $Q\in H\{Y\}\setminus H$ have order $< r$, and take $h$, $\fm$, $\fv$, $\phi$ as in the statement of the corollary. 
By Lemma~\ref{lem:dsfeh} it is enough to show that  then $Q(y) \sim_H Q(\hat h)$. We use $(\phantom{-})^\circ$ as explained at the beginning of Section~\ref{secfhhf}. Thus we have the  Hardy field~$H^\circ$ and the  $H$-field isomorphism $h\mapsto h^\circ\colon H^\phi\to H^\circ$, extended  to an $H$-field isomorphism $\hat f\mapsto \hat f^\circ\colon\hat H^\phi\to\hat H^\circ$, for an immediate $H$-field extension $\hat H^\circ$ of $H^\circ$. 
Set $u:=(y-h)/\fm\in \Calinf$. We have $\ddeg_{\prec\fv^\circ} Q^{\phi\circ}_{+h^\circ}=0$ and  $(u^\circ)^{(j)}\preceq 1$ for  $j=0,\dots,r-1$;
hence $Q^{\phi\circ}(y^\circ)\sim_{H^\circ} Q^{\phi\circ}(\hat h^\circ)$ by Lemma~\ref{1142}. Now $Q^{\phi\circ}(y^\circ)=Q(y)^\circ$ in $\Calinf$ and~$Q^{\phi\circ}(\hat h^\circ)=Q(\hat h)^\circ$ in $\hat H^\circ$, hence $Q(y) \sim_H Q(\hat h)$.
\end{proof}

\noindent
Using Corollary~\ref{cor:dsfeh} instead of Lemma~\ref{lem:dsfeh} we  show likewise:

\begin{cor}\label{cor:iso hat h, d-alg}
Suppose $\hat h$ is $\d$-algebraic over $H$ with minimal annihilator $P$ over~$H$ of order $r\geq 1$, and let~${y\in\Calinf}$ satisfy $P(y)=0$. Suppose for all~$Q$ in~$H\{Y\}\setminus H$ with $\order Q\leq r$ and $\deg_{Y^{(r)}} Q<\deg_{Y^{(r)}} P$ there are $h\in H$, $\fm,\fv\in H^\times$, and an active~$\phi>0$ in $H$ such that $\hat h-h\preceq\fm\prec\fv$,  $\ddeg_{\prec\fv} Q^\phi_{+h}=0$, and
$$ \derdelta^j\!\left(\frac{y-h}{\fm}\right)\preceq 1\qquad\text{ for }j=0,\dots,r \text{ and }  \derdelta:=\phi^{-1}\der.$$
Then $y$ is hardian over $H$ and there is an isomorphism $H\langle y\rangle \to H\langle \hat h\rangle$ of $H$-fields over
$H$ sending~$y$ to $\hat h$.
\end{cor}

\noindent
In the next subsection we use Corollary~\ref{cor:iso hat h, d-alg} to fill in certain kinds of holes in Hardy fields. 
Recall from  [ADH, remark after~11.4.3] that if $\hat h\notin H$ and~$Z(H,\hat h)=\emptyset$, then~$\hat h$ is $\d$-transcendental over $H$.
The next result is a version of  Corollary~\ref{cor:iso hat h, d-alg} for that situation.
(This will not be used until Section~\ref{sec:perfect applications} below.)

\begin{cor}\label{cor:iso hat h, d-trans}
Suppose $\hat h\notin  H$ and  $Z(H,\hat h)=\emptyset$. Let $y\in\Calinf$ be such that for all~$h\in H$, $\fm\in H^\times$ with $\hat h-h\preceq\fm$
and all $n$ there is an active $\phi_0$ in $H$ such that for all active $\phi>0$ in $H$ with $\phi\preceq \phi_0$ we have 
$\derdelta^n\!\big(\frac{y-h}{\fm}\big)\preceq 1$ for~${\derdelta=\phi^{-1}\der}$.
Then~$y$ is hardian over $H$, and there is an
isomorphism ${H\langle y\rangle \to H\langle \hat h\rangle}$ of $H$-fields over $H$
sending~$y$ to $\hat h$.
\end{cor}

\begin{proof}
Let $Q\in H\{Y\}\setminus H$; by
Lemma~\ref{dsfeh, d-trans} it is enough to show that $Q(y)\sim_H Q(\hat h)$.
Since $Q\notin Z(H,\hat h)$, we obtain $h\in H$ and $\fm,\fv\in H^\times$ 
such that $\hat h-h\preceq\fm\prec \fv$ and $\ndeg_{\prec \fv} Q_{+h}=0$.
Let $r:=\order Q$ and choose
an active $\phi>0$ in $H$
such that~$\ddeg_{\prec\fv} Q^\phi_{+h}=0$ and 
$\derdelta^n\!\big(\frac{y-h}{\fm}\big)\preceq 1$ for~${\derdelta=\phi^{-1}\der}$ and $n=0,\dots,r$.
As in the proof of Corollary~\ref{cor:iso hat h, d-alg, weak} this yields $Q(y)\sim_H Q(\hat h)$.
\end{proof}

\subsection*{Generating immediate $\d$-algebraic Hardy field extensions}
{\it In this subsection $H$ is Liouville closed, $(P, \fn, \hat h)$ is a special $Z$-minimal slot 
in $H$ of order~$r\ge 1$,  $K:=H[\imag]\subseteq \Calinf[\imag]$, $\I(K)\subseteq K^\dagger$, and $K$ is $1$-linearly surjective  if~${r\ge 3}$.}\/ We first treat the case where 
 $(P, \fn, \hat h)$ is a hole in $H$ (not just a slot):

\begin{theorem} \label{46}
Assume $(P, \fn, \hat h)$ is a deep, ultimate, and strongly repulsive-normal hole in $H$, and $y\in \Calinf$, $P(y)=0$, $y\prec \fn$.
Then $y$ is hardian over $H$, and there is an
isomorphism $H\langle y\rangle \to H\langle \hat h\rangle$ of $H$-fields over $H$ sending $y$ to $\hat h$.
\end{theorem}
\begin{proof} 
Replacing $(P, \fn, \hat h)$, $y$ by $(P_{\times \fn}, 1, \hat h/\fn)$, $y/\fn$ we arrange $\fn=1$. 
Let $Q$ in~$H\{Y\}\setminus H$, $\order Q\le r$, and
$\deg_{Y^{(r)}} Q< \deg_{Y^{(r)}} P$. Then $Q\notin Z(H,\hat h)$, so we have
$h\in H$  and $\fv\in H^\times$ such that $h-\hat h \prec \fv$ and 
$\ndeg_{\prec \fv} Q_{+h}=0$.  Take any $\fm\in H^\times$ with $\hat h - h\preceq \fm\prec \fv$. 
%We arrange $\fv\preceq \hat h\prec 1$: if $\fv\succ \hat h$, take $\fv^*\in H^\times$ and $h^*\in H$ such that $\hat h-h^*\prec \fv^*\asymp \hat h$, so $\ndeg_{\prec \fv}Q_{+h^*}=0=\ndeg_{\prec \fv^*}Q_{+h^*}$, and replace $\fv$, $h$ by $\fv^*$, $h^*$. 
Take $\fw\in H^\times$ with $\fm\prec \fw\prec\fv$. Then~$\ndeg Q_{+h,\times \fw}=0$, so we have active $\phi$ in $H$, $0<\phi\prec 1$, with $\ddeg Q_{+h, \times \fw}^\phi =0$, and hence $\ddeg_{\prec \fw} Q_{+h}^\phi=0$. Thus renaming $\fw$ as $\fv$ we have arranged $\ddeg_{\prec \fv} Q_{+h}^\phi=0$. 

Set $\derdelta:=\phi^{-1}\der$;
by Corollary~\ref{cor:iso hat h, d-alg} it is enough to show that 
$\derdelta^j\big(\frac{y-h}{\fm}\big)\preceq 1$
for~$j=0,\dots, r$. 
Now using~$(\phantom{-})^\circ$ as before, the hole~$(P^{\phi\circ},1,\hat h^\circ)$ in~$H^\circ$ is special, $Z$-minimal, deep, ultimate, and strongly repulsive-normal,
 by Lemmas~\ref{lem:Pphicirc, 1} and~\ref{lem:Pphicirc, 2}.  It remains to apply
 Corollary~\ref{cor:approx y} to this hole in~$H^\circ$
with $h^\circ$, $\fm^\circ$, $y^\circ$ in place of~$h$,~$\fm$,~$y$.
\end{proof}

\begin{cor}  \label{cor:46}
Let $\phi$ be active in $H$, $0<\phi\preceq 1$, and suppose the slot $(P^\phi,\fn,\hat h)$ in $H^\phi$ is deep, ultimate, and strongly split-normal.
Then $P(y)=0$ and $y\prec\fn$ for some $y\in\Calinf$. If 
$(P^\phi,\fn,\hat h)$ is strongly repulsive-normal, then any such $y$ is hardian over $H$ with $y\notin H$.  
\end{cor}
\begin{proof}
Lemma~\ref{dentsolver}  gives $y\in\Calinf$ with~$P(y)=0$, $y\prec\fn$.
Now suppose $(P^\phi,\fn,\hat h)$ is strongly repulsive-normal, and 
$y\in\Calinf$, $P(y)=0$, $y\prec\fn$.  Using Lemma~\ref{lem:from cracks to holes} we arrange that $(P,\fn,\hat h)$  is a hole in $H$.  The hole~$(P^{\phi\circ},\fn^\circ,\hat h^\circ)$ in $H^\circ$ is
special, $Z$-minimal, deep, ultimate, and
strongly repulsive-normal. Then Theorem~\ref{46} with~$H^\circ$,~$(P^{\phi\circ}, \fn^\circ, \hat h^\circ)$,~$y^\circ$
in place of~$H$,~$(P,\fn,\hat h)$,~$y$ shows
 that $y^\circ$ is hardian over~$H^\circ$ with $y^\circ\notin H^\circ$.
Hence $y$ is hardian over $H$ and $y\notin H$. 
\end{proof}

%\begin{remark} \marginpar{this is commented out: it is just a special case of Corollary~\ref{realginfgom}} 
%Suppose $\phi$ satisfies the hypothesis of Corollary~\ref{cor:46} and $H$ is a $\Ginf$-Hardy field.
%Then there is some $y\in\Ginf$ with $P(y)=0$, $y\prec\fn$. (Same proof as Corollary~\ref{cor:46}, using the $\Ginf$-part of
%Lemma~\ref{dentsolver}.) For each such $y$, we have~$H\<y\>\subseteq\Ginf$.
%Similarly with $\Gom$ in place of~$\Ginf$.
%\end{remark}

%\noindent
%Here is a first significant application of these results:

%\noindent
%In the next subsection we aim to show that the algebraic closures of $\d$-maximal Hardy fields are  
%$1$-linearly newtonian.

\subsection*{Achieving $1$-linear newtonianity}  For the proof of our main theorem we need to show first that for any $\d$-maximal Hardy field $H$ the corresponding $K=H[\imag]$ is $1$-linearly newtonian, the latter being a key hypothesis in Lemma~\ref{53} below. In this subsection we take this vital step: Corollary~\ref{cor:dmax K 1-linnewt}. 

\begin{lemma}\label{prop:dmax 1-newt}
Every $\d$-maximal Hardy field is $1$-newtonian.
\end{lemma}
\begin{proof} Let $H$ be a $\d$-maximal Hardy field. Then $H$ satisfies the conditions at the beginning of the previous subsection, by Corollary~\ref{cor:cos sin infinitesimal} and Theorem~\ref{upo}, for any
special $Z$-minimal slot $(P, \fn, \hat h)$ in $H$ of order $1$. 
By Corollary~\ref{cor:Liouville closed => 1-lin newt}, $H$ is $1$-linearly newtonian. 
Towards a contradiction assume that $H$ is not $1$-newtonian. Then
Lemma~\ref{lem:no hole of order <=r} yields a minimal  hole
$(P,\fn,\hat h)$ in $H$ of order~$r=1$. Using Lemma~\ref{lem:quasilinear refinement} we replace
$(P,\fn,\hat h)$ by a refinement  to  arrange that~$(P,\fn,\hat h)$ is quasilinear.  Then $(P,\fn,\hat h)$  is special,
by Lemma~\ref{lem:special dents}.
Using Corollary~\ref{cor:5.30real rep-norm} we further refine 
$(P,\fn,\hat h)$  to arrange that~$(P^\phi,\fn,\hat h)$ is eventually deep, ultimate, and strongly repulsive-normal.
Now Corollary~\ref{cor:46} gives a proper $\d$-algebraic Hardy field extension of $H$, contradicting $\d$-maximality of $H$.
\end{proof}

\noindent
{\it In the rest of this subsection $H$ has asymptotic integration.}\/  We have the $\d$-valued extension $K:=H[\imag]\subseteq \Calinf[\imag]$ of $H$ and as before we arrange that~$\hat K=\hat H[\imag]$ is a $\d$-valued extension of $\hat H$ as well as an immediate $\d$-valued extension of $K$.

\begin{lemma}\label{lem:deg 1 order 1 approx}
Suppose $H$ is Liouville closed and $\I(K)\subseteq K^\dagger$.
Let $(P,\fn,\hat f)$ be an ultimate linear minimal hole in $K$ of order $r\geq 1$, where $\hat f\in\hat K$, such that~$\dim_{\C}\ker_{\Univ} L_P=r$.
Assume also that $K$ is $\upo$-free if $r\ge 2$. Let $f\in\Calinf[\imag]$ be such that $P(f)=0$, $f\prec\fn$.
Then~$f\approx_K \hat f$.
\end{lemma}
\begin{proof}
Replacing $(P,\fn,\hat f)$, $f$ by $(P_{\times\fn},1,\hat f/\fn)$, $f/\fn$ we arrange $\fn=1$.
Let $\theta\in K^\times$ be such that $\theta\sim\hat f$; we claim that $f\sim\theta$ in $\c[\imag]$ (and so $f\sim_K \hat f$).
Applying Proposition~\ref{prop:deg 1 analytic} and Remark~\ref{rem:improap}  to the linear minimal hole~$(P_{+\theta},\theta,\hat f-\theta)$ in~$K$ gives $g\in\Gi[\imag]$
such that $P_{+\theta}(g)=0$ and $g\prec\theta$. Then~${P(\theta+g)=0}$ and~${\theta+g \prec 1}$, thus 
$L_P(y)=0$ and $y\prec 1$ for~$y:=f-(\theta+g)\in\Calinf[\imag]$.
Hence~$y\prec\theta$ by the  version of Lemma~\ref{lem3.9 linear} for slots in~$K$; see the remark following Corollary~\ref{cor:8.8 refined}.
Therefore $f-\theta=y+g\prec\theta$  and so~$f\sim\theta$, as claimed. 

The refinement $(P_{+\theta},1,\hat f-\theta)$ of $(P,1,\hat f)$ is ultimate thanks to the $K$-version of Lemma~\ref{lem:ultimate refinement}, and $L_{P_{+\theta}}=L_P$, so we can apply the claim to $(P_{+\theta},1,\hat f-\theta)$ instead of~$(P,1,\hat f)$
and $f-\theta$ instead of $f$ to give $f-\theta \sim_K \hat f-\theta$.
Since this holds for all~$\theta\in K$ with $\theta\sim\hat f$, the $K$-version of Corollary~\ref{cor:simH phi0} then yields~$f\approx_K \hat f$.
\end{proof}

\begin{cor}\label{cor:embed K<hatf>}
Let $(P,\fn,\hat f)$ be a linear hole of order $1$ in $K$. \textup{(}We do not assume here that $\hat f\in \hat{K}$.\textup{)} Then there is an embedding $\iota\colon K\langle\hat f\rangle \to \Calinf[\imag]$ of differential $K$-algebras such that $\iota(g)\sim_K g$ for all~$g\in  K\langle\hat f\rangle^\times$. 
\end{cor}
\begin{proof}
 Note that $(P,\fn,\hat f)$ is minimal. We first show how to arrange that $H$ is Liouville closed and $\upo$-free with $\I(K)\subseteq K^\dagger$ and $\hat f\in\hat K$. 
Let $H_1$ be a maximal Hardy field extension of~$H$.
Then $H_1$ is Liouville closed and $\upo$-free, with $\I(K_1)\subseteq K_1^\dagger$ for~$K_1:=H_1[\imag]\subseteq \Calinf[\imag]$. %Moreover, by Corollary~\ref{cor:Hardy field ext smooth}, if $H\subseteq\c^\infty$, then $H_1\subseteq\c^\infty$, and similarly with $\c^\omega$ in place of $\c^\infty$.
Let~$\hat H_1$ be the newtonization of $H_1$; then $\hat K_1:=\hat H_1[\imag]$ is   newtonian [ADH, 14.5.7]. 
Corollary~\ref{cor:find zero of P, 2} gives an embedding $K\langle \hat f\rangle \to \hat K_1$ of valued differential fields over $K$; let $\hat f_1$ be the image  of $\hat f$ under this embedding. If $\hat f_1\in K_1\subseteq\Calinf[\imag]$, then we are done, so assume $\hat f_1\notin K_1$.   Then $(P,\fn,\hat f_1)$ is a  hole in $K_1$, and we replace $H$, $K$, $(P,\fn,\hat f)$ by
 $H_1$, $K_1$, $(P,\fn,\hat f_1)$, and $\hat K$ by $\hat K_1$,  to arrange that $H$ is Liouville closed and $\upo$-free with $\I(K)\subseteq K^\dagger$
 and $\hat f\in \hat K$. 
 
Replacing $(P,\fn,\hat f)$ by a refinement we also arrange that~$(P,\fn,\hat f)$ is ultimate and~$\fn\in H^\times$, by Proposition~\ref{prop:achieve ultimate, K} and Remark~\ref{rem:achieve ultimate, K}. Then
Proposition~\ref{prop:deg 1 analytic} yield an $f\in\Calinf[\imag]$ with~$P(f)=0$, $f\prec \fn$.
Now Lemma~\ref{lem:deg 1 order 1 approx} gives~$f\approx_K\hat f$, and it remains to appeal to Corollary~\ref{cor:dsfek, order 1}.
\end{proof}

\begin{prop}\label{prop:Re or Im}
Suppose $H$ is $\upo$-free and $1$-newtonian.
Let   $(P,\fn,\hat f)$ be a linear hole in $K$ of order $1$ with $\hat f\in \hat K$, 
 and $f\in\Calinf[\imag]$, $P(f)=0$,
and~$f\approx_K\hat f$.
Then~$\Re f$ or $\Im f$ generates a proper $\d$-algebraic Hardy field extension of~$H$.
\end{prop}

\begin{proof}
Let $\hat g:=\Re\hat f$ and $\hat h:=\Im\hat f$.
By Lemma~\ref{lem:same width} we have $v(\hat g-H)\subseteq v(\hat h-H)$ or $v(\hat h-H)\subseteq v(\hat g-H)$.
Below we assume $v(\hat g-H)\subseteq v(\hat h-H)$ (so $\hat g\in \hat H\setminus H$) and show that then $g:=\Re f$ generates a proper $\d$-algebraic Hardy field extension of~$H$. (If $v(\hat h-H)\subseteq v(\hat g-H)$ one shows likewise that $\Im f$ generates a proper $\d$-algebraic Hardy field extension of $H$.)
The hole  $(P,\fn,\hat f)$ in $K$ is minimal, and by arranging~$\fn\in H^\times$ we see that $\hat g$ is $\d$-algebraic over~$H$, by a remark preceding Lemma~\ref{kb}. 
 %Let~$Q$ be a minimal annihilator of~$\hat g$ over $H$.
Every element of $Z(H,\hat g)$ has order~$\geq 2$, by Corollary~\ref{cor:no dent of order <=r} and $1$-newtonianity of~$H$. We arrange that the linear part $A$ of $P$ is monic, so~${A=\der-a}$ with $a\in K$,  $A(\hat f)=-P(0)$ and $A(f)=-P(0)$.
Then Example~\ref{ex:lclm compl conj} and Remark~\ref{rem:lclm compl conj}  applied to
$F=\Calinf$ yields $Q\in H\{Y\}$ with~$1\le \order Q\le 2$ and~$\deg Q=1$ such that~$Q(\hat g)=0$ and $Q(g)=0$. 
Hence~$\order Q=2$ and
$Q$ is a minimal annihilator of~$\hat g$ over $H$.
%an element of  $Z(H,\hat g)$ of minimal complexity.
%In view of~$P(\hat f)=0$ and $P(f)=0$ we obtain from the example preceding Lemma~\ref{lem:lclm compl conj} and the remark following it in combination with [ADH, 4.1.8] that~$Q(g)=0$.

Towards applying Corollary~\ref{cor:iso hat h, d-alg, weak} to $Q$,~$\hat g$,~$g$ in the role of $P$,~$\hat h$,~$y$ there, let $R$ in~$H\{Y\}\setminus H$ have order $< 2$. Then $R\notin Z(H,\hat g)$, so we have $h\in H$ and $\fv\in H^\times$
such that~$\hat g-h\prec\fv$ and $\ndeg_{\prec\fv} R_{+h}=0$. Take any~$\fm\in H^\times$ with~$\hat g-h\preceq\fm\prec\fv$.
By Lemma~\ref{lem:real part approx} we have $g\approx_H\hat g$ and thus~$g-h\preceq\fm$.
After changing $\fv$ as in the proof of Theorem~\ref{46} we obtain an active $\phi$ in $H$, $0<\phi\preceq 1$, such that~$\ddeg_{\prec\fv} R^\phi_{+h}=0$. Set $\derdelta:=\phi^{-1}\der$;
 by Corollary~\ref{cor:iso hat h, d-alg, weak} it is now enough to show that~${\derdelta\big( (g-h)/\fm \big) \preceq 1}$. 
 
 Towards this and using~$(\phantom{-})^\circ$ as before, we have $f^\circ \approx_{K^\circ}\hat f^\circ$,
 and $g^\circ\approx_{H^\circ} \hat g^\circ$
 by the facts about composition in Section~\ref{sec:asymptotic similarity}.
 Moreover, $(g - h)^\circ\preceq \fm^\circ$, and
 $H^\circ$ is $\upo$-free and $1$-newtonian, hence closed under integration by [ADH, 14.2.2].
We now apply Corollary~\ref{cor:dsfek, order 1} with $H^\circ$,~$K^\circ$,~$P^{\phi\circ}$,~$f^\circ$,~$\hat f^\circ$ in the role of
$H$,~$K$,~$P$,~$f$,~$\hat f$ to give $$ \big(f^\circ/\fm^\circ\big)'\  \approx_{K^\circ}\  \big(\hat f^\circ/\fm^\circ\big)',$$   
hence
$(g^\circ/\fm^\circ)' \approx_{H^\circ} (\hat g^\circ/\fm^\circ)'$ by Lemmas~\ref{lem:same width, der} and~\ref{lem:real part approx}. Therefore,
$$\big((g-h)^\circ/\fm^\circ\big)'\ =\ (g^\circ/\fm^\circ)'-(h^\circ/\fm^\circ)' \ \sim_H\ (\hat g^\circ/\fm^\circ)'-(h^\circ/\fm^\circ)' \ =\  \big((\hat g-h)^\circ/\fm^\circ\big)'.$$
Now $(\hat g-h)^\circ/\fm^\circ \preceq 1$, so $ \big((\hat g-h)^\circ/\fm^\circ\big)' \prec 1$, hence $\big((g-h)^\circ/\fm^\circ\big)' \prec 1$ by the last display, and thus ${\derdelta\big( (g-h)/\fm \big) \prec 1}$, which is more than enough.
\end{proof}

\noindent
If $K$ has a linear hole of order $1$, then $K$ has a proper $\d$-algebraic
differential field extension inside  $\Calinf[\imag]$, by Corollary~\ref{cor:embed K<hatf>}. We can now prove a Hardy analogue:

\begin{lemma}\label{linhole1}
Suppose $K$ has a linear hole of order $1$.
Then~$H$ has a proper $\d$-algebraic Hardy field extension.
\end{lemma}
\begin{proof} If $H$ is not $\d$-maximal, then $H$ has indeed a proper $\d$-algebraic Hardy field extension, and if $H$ is $\d$-maximal, then 
$H$ is  Liouville closed,  $\upo$-free, $1$-newtonian, and $\I(K)\subseteq K^\dagger$, 
by Proposition~\ref{prop:Hardy field exts}, Corollary~\ref{cor:cos sin infinitesimal}, Theorem~\ref{upo}, and Lemma~\ref{prop:dmax 1-newt}.  So assume below that $H$ is  Liouville closed,  $\upo$-free, $1$-newtonian, 
and $\I(K)\subseteq K^\dagger$, and that $(P,\fn,\hat f)$ is a linear hole of order $1$  in $K$. 
By Lemma~\ref{lem:hole in hat K} we arrange that $\hat f\in \hat K:= \hat H[\imag]$ where
$\hat H$ is an immediate $\upo$-free newtonian $H$-field extension of $H$. Then $\hat K$ is also newtonian by [ADH, 14.5.7]. 
By Remark~\ref{rem:achieve ultimate, K}  we can replace $(P,\fn,\hat f)$ by a refinement to arrange that $(P,\fn,\hat f)$ is ultimate and~$\fn\in H^\times$. 
Proposition~\ref{prop:deg 1 analytic} now yields $f\in\Gi[\imag]$ with~$P(f)=0$ and~$f\prec \fn$.
Then~$f\approx_K\hat f$ by Lemma~\ref{lem:deg 1 order 1 approx}, and so 
$\Re f$ or $\Im f$ generates a proper $\d$-algebraic Hardy field extension of $H$, by Proposition~\ref{prop:Re or Im}.
\end{proof}

\begin{cor}\label{cor:dmax K 1-linnewt}
If $H$ is $\d$-maximal, then $K$ is $1$-linearly newtonian.
\end{cor}
\begin{proof} Assume $H$ is $\d$-maximal. Then $K$ is $\upo$-free
by Theorem~\ref{upo} and [ADH, 11.7.23].
If $K$ is not $1$-linearly newtonian, 
then it has a linear hole of order $1$, by Lemma~\ref{lem:no hole of order <=r, deg 1}, and so 
$H$ has a proper $\d$-algebraic Hardy field extension, by Lemma~\ref{linhole1}, contradicting $\d$-maximality of~$H$.
\end{proof}

%\noindent
%We can improve on this proposition:

%\begin{cor} 
%Suppose $H$ is $\d$-perfect and $\upl$-free; then $K$ is $1$-linearly newtonian. \marginpar{probably not true}
%\end{cor}
%\begin{proof} 
%By  Corollary~\ref{cor:Singer-Rosenlicht},  $H$ is Liouville closed, and by [ADH, 11.6.8], $K$ is $\upl$-free.  Towards a contradiction assume that $K$ is not $1$-linearly newtonian. As in the proof of Proposition~\ref{prop:dmax K 1-linnewt} we then obtain a hole $(P,1,\hat f)$ in $K$ with $\hat f\in\hat K$ where   $P=Y'-\xi Y-u\in K\{Y\}$  with $\xi\in H\imag$, $u\in K$.   By Lemma~\ref{lem:8.8 order deg 1} there is at most one $f\in\Go[\imag]$ such that~$P(f)=0$ and $f\prec 1$. However, by Proposition~\ref{}, for each $\d$-maximal Hardy field extension $H^*$ of $H$ the $\d$-valued field $K^*:=H^*[\imag]$ is $1$-linearly newtonian and hence $1$-linearly surjective, by [ADH, ];
%\end{proof}

\subsection*{Finishing the story} With one more lemma we will be done.

%In this subsection we assume that $H$ is Liouville closed and $\upo$-free,
%and that~$\hat H$ is a newtonization of $H$; see [ADH, Chapter~14]. Then~$\hat H\ne H$ if~$H$ is not newtonian.  We set $\hat K:= \hat H[\imag]$, so $\hat{K}$ is a newtonian immediate asymptotic extension of $K$ by [ADH, 14.5.7].

\begin{lemma}\label{53} Suppose $H$ is Liouville closed, $\upo$-free, not newtonian, and $K:=H[\imag]$ is $1$-linearly newtonian. Then $H$ has a proper $\d$-algebraic Hardy field extension.
\end{lemma}
\begin{proof} 
By Proposition~\ref{prop:char 1-linearly newt}, $K$ is $1$-linearly surjective and $\I(K)\subseteq K^\dagger$.
Since $H$ is not newtonian, 
neither is $K$, by [ADH, 14.5.6], and so by Lemma~\ref{lem:no hole of order <=r}
we have a minimal hole~$(P,\fm,\hat f)$ in $K$ of order $r\ge 1$, with $\fm\in H^\times$. Then $\deg P>1$ by Corollary~\ref{cor:minhole deg 1}.  
As in the proof of Lemma~\ref{linhole1} we take for $\hat H$ an immediate $\upo$-free newtonian $H$-field extension of $H$ and arrange $\hat f\in\hat K:=\hat{H}[\imag]$.
%The case $r=1$ is taken care of by Theorem~\ref{51}, so assume $r>1$.
Now $\hat f = \hat g + \hat h\imag$ with $\hat g, \hat h \in \hat H$.
By Theorem~\ref{thm:strongly repulsive-normal}, there are two cases:
\begin{enumerate}
\item  $\hat g\notin H$ and some $Z$-minimal slot $(Q,\fm,\hat g)$ in $H$ has a special refinement ${(Q_{+g},\fn,\hat g-g)}$ such that $(Q^\phi_{+g},\fn,\hat g-g)$ is eventually deep, strongly re\-pul\-sive-normal, and ultimate; 
\item  $\hat h\notin H$ and some $Z$-minimal slot $(R,\fm,\hat h)$ in $H$ has a special refinement ${(R_{+h},\fn,\hat h-h)}$ such that $(R^\phi_{+h},\fn,\hat h-h)$ is eventually deep, strongly re\-pul\-sive-normal, and ultimate. 
\end{enumerate}
Suppose $\hat g\notin H$ and $(Q,\fm,\hat g)$ is as in (1).
Then $1\le \order Q \le 2r$
by Lemma~\ref{kb}. {\em Claim}: $Q(y)=0$ for some $y\in\Calinf\setminus H$ that is hardian over $H$. 
To prove this claim, take a special refinement~$(Q_{+g}, \fn, \hat g -g)$ of~$(Q,\fm,\hat g)$ and an active $\phi$ in $H$ 
with $0<\phi\preceq 1$ such that  the slot $(Q^{\phi}_{+g}, \fn, \hat g -g)$ in~$H^\phi$
is deep, strongly repulsive-normal, and ultimate.  
Corollary~\ref{cor:46} applied to $(Q_{+g},\fn,\hat g-g)$ in place of $(P,\fn,\hat h)$ gives a $z\in \Calinf\setminus H$   
that is hardian over $H$ with $Q_{+g}(z)=0$. 
Thus $y:=g+z\in\Calinf$  is as in the Claim.  Case (2) is handled likewise.  
\end{proof}

\noindent
Recall from the introduction that an {\it $H$-closed field}\/ is an $\upo$-free newtonian Liouville closed $H$-field. Recall also that
Hardy fields containing $\R$ are $H$-fields.  The main result of these notes can now be established in a few lines: 

\begin{theorem}\label{thm:char d-max}
A Hardy field  is $\d$-maximal iff it contains $\R$ and is $H$-closed.
\end{theorem}
\begin{proof}
The ``if'' part is a special case of [ADH, 16.0.3].
By Proposition~\ref{prop:Hardy field exts}  and Theorem~\ref{upo}, every $\d$-maximal Hardy field contains~$\R$ and is
Liouville closed and $\upo$-free. 
Suppose~$H$ is $\d$-maximal. Then $K:=H[\imag]$ is $1$-linearly newtonian
by Corollary~\ref{cor:dmax K 1-linnewt}, so $H$ is newtonian by Lemma~\ref{53}.
\end{proof}

\noindent
Theorem~\ref{thm:char d-max} and Corollary~\ref{cor:Hardy field ext smooth} yield Theorem~\ref{thm:B} from the introduction in a refined form:

\begin{cor}\label{thm:extend to H-closed}
Any Hardy field $F$ has a $\d$-algebraic $H$-closed Hardy field extension.
If $F$ is a $\Ginf$-Hardy field, then so is any such extension, and likewise with~$\Gom$ in place of $\Ginf$.
\end{cor}

%Theorem~\ref{thm:B} from the introduction follows from Theorem~\ref{thm:char d-max}, since any Hardy field extends to a $\d$-maximal one. 

\newpage 

\part{Applications}\label{part:applications}

\medskip

\noindent
Here we apply the material in the previous parts.
In Section~\ref{sec:transfer} we show how to transfer first-order logical properties of the differential field $\mathbb T$ of transseries to maximal Hardy fields, proving in particular Theorem~\ref{thm:A} and Corollaries~\ref{cor:elem equiv}--\ref{cor:systems, 2}
as well as the first part of Corollary~\ref{cor:systems, 3} from the introduction.
In Section~\ref{sec:diff closure} we obtain Corollary~\ref{cor:zeros in complexified Hardy field extensions}, elaborate on [ADH, Chapter~16], and relate Newton-Liouville closure to
relative differential closure. In Section~\ref{sec:embeddings into T} we investigate embeddings
of Hardy fields into~$\mathbb T$, and finish the proof of Corollary~\ref{cor:systems, 3}.
There we also determine the universal theory of Hardy fields.
Section~\ref{sec:lin diff applications} contains applications of our main theorem to linear
differential equations over Hardy fields, including proofs of Corollaries~\ref{cor:factorization intro}--\ref{cor:disconjugacy} from the introduction.
The final Corollary~\ref{cor:Boshernitzan intro} from the introduction is established in Section~\ref{sec:perfect applications}, where we
focus on the structure of perfect and $\d$-perfect Hardy fields.

\section{Transfer Theorems}\label{sec:transfer}  

\noindent
From [ADH, 16.3] we recall the notion of a {\it pre-$\HLO$-field~$\boldsymbol  H=(H,\I,\Lambda,\Omega)$}\/: this is a pre-$H$-field $H$ equipped with a $\HLO$-cut $(\I,\Lambda,\Omega)$ of $H$. (See also Section~\ref{sec:upo-free Hardy fields}.)
A {\it $\HLO$-field\/} is a pre-$\HLO$-field $\boldsymbol  H=(H;\dots)$ where $H$ is an $H$-field.\index{LambdaOmega-field@$\HLO$-field}
 If~$\boldsymbol  M=(M;\dots)$ is a pre-$\HLO$-field and $H$ is a pre-$H$-subfield of~$M$, then $H$ has a unique expansion to a  pre-$\HLO$-field $\boldsymbol  H$ such that $\boldsymbol  H\subseteq\boldsymbol  M$. 
By~[ADH, 16.3.19], a pre-$H$-field $H$ has a unique expansion to a pre-$\HLO$-field iff one of the following conditions holds:
\begin{enumerate}
\item $H$ is grounded;
\item there exists $b\asymp 1$ in $H$ such that $v(b')$ is a gap in $H$; 
\item $H$ is $\upo$-free.
\end{enumerate}
In particular, each  $\d$-maximal Hardy field $M$ (being $\upo$-free) has a  unique expansion to a pre-$\HLO$-field~$\boldsymbol  M$,
namely $\boldsymbol  M=\big(M;\I(M),\Upl(M),\omega(M)\big)$,
and then $\boldsymbol  M$ is a $\HLO$-field with constant field $\R$. Below we always view any $\d$-maximal Hardy field
as an $\HLO$-field in this way.

\begin{lemma}\label{lem:canonical HLO}
Let $H$ be a Hardy field. Then $H$ has an expansion to a pre-$\HLO$-field~$\boldsymbol H$ such that $\boldsymbol H\subseteq\boldsymbol M$ for every $\d$-maximal Hardy field $M\supseteq H$.
\end{lemma}
\begin{proof}
Since every $\d$-maximal Hardy field containing $H$ also contains $\operatorname{D}(H)$, it suffices to show this
for $\operatorname{D}(H)$ in place of $H$. So we assume $H$ is $\d$-perfect,
and thus a  Liouville closed $H$-field.
For each $\d$-maximal Hardy field~$M\supseteq H$ we now have
$\I(H)=\I(M)\cap H$ by [ADH, 11.8.2], $\Upl(H)=\Upl(M)\cap H$ by [ADH, 11.8.6], and~$\omega(H)=\bar{\omega}(H)=\bar{\omega}(M)\cap H=\omega(M)\cap H$ by Corollary~\ref{cor:omega(H) downward closed}, as required.
%This yields the claim.
\end{proof}

\noindent
Given a  Hardy field $H$, we call the unique expansion $\boldsymbol H$ of $H$ to a pre-$\HLO$-field  with the property stated in the previous lemma the {\bf canonical $\HLO$-expansion} of~$H$.\index{Hardy field!canonical $\HLO$-expansion}

\begin{cor}\label{cor:canonical HLO}
Let $H$, $H^*$ be Hardy fields, with their canonical $\HLO$-ex\-pan\-sions~$\boldsymbol H$,~$\boldsymbol H^*$, respectively, such that $H\subseteq H^*$. Then $\boldsymbol H\subseteq\boldsymbol H^*$.
\end{cor}
\begin{proof}
Let $M^*$ be any $\d$-maximal Hardy field extension of $H^*$. Then $\boldsymbol H\subseteq\boldsymbol M^*$ as well as $\boldsymbol H^*\subseteq\boldsymbol M^*$, hence $\boldsymbol H\subseteq\boldsymbol H^*$.
\end{proof}

\noindent
{\it In the rest of this section $\mathcal L=\{0,1,{-},{+},{\,\cdot\,},{\der},{\leq},{\preceq}\}$ is the language of ordered valued differential rings}\/ [ADH, p.~678]. We view 
each ordered valued differential field as an $\mathcal L$-structure in the natural way. Given an ordered valued differential field~$H$
and a subset~$A$ of~$H$ we let $\mathcal L_A$ be $\mathcal L$ augmented by names for the elements of~$A$, and expand the $\mathcal L$-structure $H$ 
to an $\mathcal L_A$-structure by interpreting the name of any $a\in A$ as the element $a$ of $H$; cf.~[ADH, B.3].
Let~$H$ be a Hardy field and $\sigma$ be an $\mathcal L_H$-sentence.
We now have our Hardy field analogue of the ``Tarski principle''~[ADH, B.12.14] in real algebraic geometry promised in the introduction:

\begin{theorem}\label{thm:transfer}
The following are equivalent:
\begin{enumerate}
\item[\textup{(i)}] $M\models\sigma$ for some $\d$-maximal Hardy field $M\supseteq H$;
\item[\textup{(ii)}] $M\models\sigma$ for every $\d$-maximal Hardy field $M\supseteq H$;
\item[\textup{(iii)}] $M\models\sigma$ for every maximal Hardy field $M\supseteq H$;
\item[\textup{(iv)}] $M\models\sigma$ for some maximal Hardy field $M\supseteq H$.
\end{enumerate}
\end{theorem}
\begin{proof}
The implications (ii)~$\Rightarrow$~(iii)~$\Rightarrow$~(iv)~$\Rightarrow$~(i) are obvious, since ``maximal~$\Rightarrow$~$\d$-maximal'';
so it remains to show (i)~$\Rightarrow$~(ii).
Let $M$, $M^*$ be $\d$-maximal Hardy field extensions of $H$.
By Lemma~\ref{lem:canonical HLO} and Corollary~\ref{cor:canonical HLO}
expand $M$, $M^*$, $H$ to pre-$\HLO$-fields $\boldsymbol M$, $\boldsymbol M^*$, $\boldsymbol H$, respectively, such that
$\boldsymbol H\subseteq\boldsymbol M$ and $\boldsymbol H\subseteq \boldsymbol M^*$. 
In~[ADH, introduction to Chapter~16] we extended $\mathcal{L}$ to a language $\mathcal L^\iota_{\HLO}$, and explained in~[ADH, 16.5] how each
pre-$\HLO$-field~$\boldsymbol K$ is construed as an  $\mathcal L^\iota_{\HLO}$-structure in such a way that
every extension~${\boldsymbol K\subseteq\boldsymbol L}$  of pre-$\HLO$-fields corresponds to an extension of the associated~$\mathcal L^\iota_{\HLO}$-structures.
By~[ADH, 16.0.1], the  $\mathcal L^\iota_{\HLO}$-theory~$T^{\operatorname{nl},\iota}_{\HLO}$ of $H$-closed
$\HLO$-fields  eliminates quantifiers, and by
Theorem~\ref{thm:char d-max}, the canonical $\HLO$-expansion of each $\d$-maximal Hardy field is a model of $T^{\operatorname{nl},\iota}_{\HLO}$.
Hence~$\boldsymbol M \equiv_H \boldsymbol M^*$~[ADH, B.11.6], so if~$\boldsymbol M  \models \sigma$, then~$\boldsymbol M^*  \models \sigma$.   
\end{proof}

\noindent
Corollaries~\ref{cor:elem equiv} and~\ref{cor:systems, 1} from the introduction are special cases of
Theorem~\ref{thm:transfer}.
By Corollary~\ref{realginfgom},
$\Ginf$-maximal and $\Gom$-maximal Hardy fields are $\d$-maximal,  so the theorem above also yields Corollary~\ref{cor:systems, 2} from the introduction in the following stronger form:

\begin{cor}
If $H\subseteq\Ginf$ and
$M\models\sigma$ for some $\d$-maximal Har\-dy field extension $M$ of $H$, then
$M\models\sigma$ for every $\Ginf$-maximal Hardy field~${M\supseteq H}$. 
Likewise with $\Gom$ in place of $\Ginf$.
\end{cor}
%\noindent
%In particular, if $H$ is $\Ginf$-maximal, then  $H$ is an elementary substructure of each $\d$-maximal Hardy field extension  of~$H$. Similarly with $\Gom$ in place of $\Ginf$.

\subsection*{The structure induced on $\R$}
In the next corollary $H$ is a Hardy field and~$\varphi(x)$ is an $\mathcal L_H$-formula where $x=(x_1,\dots,x_n)$ and
$x_1,\dots, x_n$ are distinct variables. 
Also, $\mathcal L_{\operatorname{OR}}=\{0,1,{-},{+},{\,\cdot\,},{\leq}\}$ is the language of ordered rings, and the ordered field $\R$ of real numbers is interpreted as an $\mathcal L_{\operatorname{OR}}$-structure in the obvious way. 
By Theorem~\ref{thm:char d-max}, $\d$-maximal Hardy fields are $H$-closed fields, so from [ADH, 16.6.7, B.12.13] 
in combination with Theorem~\ref{thm:transfer} we  obtain:

\begin{cor}\label{cor:sa}
There is a quantifier-free
 $\mathcal L_{\operatorname{OR}}$-formula $\varphi_{\operatorname{OR}}(x)$ such that for all
$\d$-maximal Hardy fields $M\supseteq H$ and $c\in\R^n$ we have
$$   M \models \varphi(c) \quad \Longleftrightarrow\quad \R\models \varphi_{\operatorname{OR}}(c).$$
\end{cor}

\noindent
This yields Corollary~\ref{cor:parametric systems} from the Introduction.  
We now justify what we claim about the examples after that corollary. The first of these examples is already covered by~[ADH, 5.1.18, 11.8.25, 11.8.26], so
we only deal  with the second example here:

{ \samepage
\begin{prop}\label{prop:wp}
Let $g_2,g_3\in \R$. Then the following are equivalent:
\begin{enumerate}
\item[\rm(i)] there exists a hardian germ $y\notin \R$  such that $(y')^2=4y^3-g_2y-g_3$;
\item[\rm(ii)] $g_2^3=27g_3^2$ and $g_3 \leq 0$;
\end{enumerate}
\end{prop}}

\noindent
For (i)~$\Rightarrow$~(ii) we take a more general setting, and recycle   arguments used in the proof of [ADH, 10.7.1].
Let $K$ be a valued differential field such that~$\der\mathcal O\subseteq\smallo$ and~$C\subseteq \mathcal O$.
(This holds for any $\d$-valued field with small derivation.) Consider a polynomial $P(Y) = 4Y^3-g_2Y-g_3$ with $g_2,g_3\in C$. Its discriminant is $16\Delta$ where~$\Delta:=g_2^3 - 27g_3^2$.
Take $e_1$, $e_2$, $e_3$ in an algebraic closure of $C$ auch that
$$P(Y)=4(Y-e_1)(Y-e_2)(Y-e_3).$$ 
Then
\begin{equation}\label{eq:e1e2e3}
e_1+e_2+e_3\	=\ 0, \quad
e_1e_2+e_2e_3+e_3e_1\  =\ -\textstyle\frac{1}{4}g_2, \quad
e_1e_2e_3\ =\ \frac{1}{4}g_3,
\end{equation}
and $\Delta\neq 0$ iff $e_1$, $e_2$, $e_3$ are distinct.
In the next two lemmas $y\in K$ and $(y')^2=P(y)$.
Then
$y\preceq 1$: otherwise 
$(y')^2\prec 4y^3\sim P(y)=(y')^2$ by [ADH, 4.4.3], a contradiction. Hence $P(y)=(y')^2\prec 1$. 
Moreover, if $y\asymp 1$, then $\Delta\neq 0$ or $g_3\neq 0$.

\begin{lemma}\label{eee}
Suppose $P'(y)\asymp 1$. Then $y\in \{e_1, e_2, e_3\}$ $($so $y\in C)$. 
\end{lemma}
\begin{proof} The property $\der\mathcal O\subseteq\smallo$ means that the derivation of $K$  is small with trivial induced derivation on its residue field.
By [ADH, 6.2.1, 3.1.9] this property is inherited by any algebraic closure of $K$, and so is the property $C\subseteq \mathcal{O}$
by [ADH, 4.1.2]. Thus by passing to an algebraic closure we arrange  that $K$ is algebraically closed.
%(See also Lemma~\ref{lem:very small der alg ext} below.)
Then~$C$ is also algebraically closed, so $e_1,e_2,e_3\in C$ and thus
$y-e_j\prec 1$ for some~$j\in\{1,2,3\}$,
say $y=e_1+z$ where $z\prec 1$. Since~$P'(y)\asymp 1$ we have~$y-e_2\asymp y-e_3\asymp 1$ and thus
$$z\ \asymp\ 4z(e_1-e_2+z)(e_1-e_3+z)\ =\  P(y)\  =\  (y')^2\  =\  (z')^2.$$
Now if $z\neq 0$, then $(z')^2 \prec z$, again by~[ADH, 4.4.3],
a contradiction. So $y=e_1$. 
\end{proof}

\noindent
In the next lemma  $K$ is in addition equipped with an ordering making $K$ a valued  ordered differential field whose valuation ring is convex. (Any $H$-field with small derivation satisfies the conditions we imposed.) 
Suppose $\Delta=0$. Then $e_1$, $e_2$, $e_3$ lie in the real closure of~$C$, and  after arranging $e_1\geq e_2\geq e_3$,
the first and the  last of the   equations~\eqref{eq:e1e2e3} yield $e_1=e_2\Longleftrightarrow g_3\leq 0$, and $e_2=e_3\Longleftrightarrow g_3\geq 0$. 

\begin{lemma}\label{deltag3}
Suppose $\Delta=0$ and $g_3>0$. Then  $y\in C$.
\end{lemma}
\begin{proof} Passing to the real closure of $K$ with convex valuation extending that of~$K$, cf.~[ADH, 3.5.18], 
we arrange that $K$, and hence $C$, is real closed.
Arranging also~$e_1\geq e_2\geq e_3$, we set $e:=e_2=e_3$. Then $e_1=-2e>0>e$ and $P(Y)=4(Y+2e)(Y-e)^2$.
We have $y\preceq 1$ and $P(y)\prec 1$. Suppose $y\notin C$. Then  $P'(y)\prec 1$ by Lemma~\ref{eee}, so  $y-e\prec 1$. Set $z:=y-e$, so  $0\neq z\prec 1$ and hence
$$12e z^2\ \sim\  4(z+3e)z^2\ =\  P(y)\  =\  (y')^2\  >\ 0, $$
contradicting $e<0$.
\end{proof}

\noindent
{\it Proof of Proposition~\ref{prop:wp}.}\/ 
Suppose $y\notin \R$ is a hardian germ such that $(y')^2=P(y)$ with $P(Y)=4Y^3-g_2Y-g_3$ and $g_2, g_3\in \R$. Then $y\preceq 1$ and $P(y)\prec 1$, but also~$P'(y)\prec 1$ by Lemma~\ref{eee}, hence $\Delta\prec 1$. As $\Delta\in \R$ this gives $\Delta=0$,
so~$g_3\le 0$ by Lemma~\ref{deltag3}. This proves (i)~$\Rightarrow$~(ii).  For the converse, let $K$ be the Hardy field $\R$ in the considerations above
and suppose
$\Delta=0$, so $e_1, e_2, e_3\in \R$. Arrange~$e_1\geq e_2\geq e_3$.
%If $g_3=0$, then  $e_1=e_2=e_3=0$, \marginpar{commented out material has been checked, but duplicates present version}
%$P(Y)=4Y^3$, and the   differential equation 
%$(Y')^2=P(Y)$
%has for each $c\in \R$ the solution 
%$$y\ =\ \frac{1}{(x-c)^2}\in\Ginf(\R^{>c}),$$
%whose germ lies in the Hardy field $\R(x)$. 
%Next, let $g_3 < 0$, so $e_1=e_2 > 0$ and $P(Y)=4(Y-e_1)^2(Y+2e_1)$. Recall 
%that  $\sinh t:=\frac{1}{2}(\ex^t-\ex^{-t})$ for~$t\in\R$.
%Using the change of variables $y=e_1+3e_1/z^2$ we see that~$(Y')^2=P(Y)$
%has solutions
%$$y\ =\ e_1+\frac{3e_1}{\sinh^2\big((x-c)\sqrt{3e_1}\big)}\in\Ginf(\R^{>c}),\qquad(c\in \R)$$
%whose germs lie in the Hardy field $\R(\ex^{x\sqrt{3e_1}})$.
%(See also \cite[\S{}13.15]{EMOT}.) 
%This concludes the proof of Proposition~\ref{prop:wp}. 
%\medskip\noindent
%Let  $g_2, g_3\in \R$ and assume $\Delta=0$ below, so $e_1, e_2,e_3\in \R$. For completeness' sake we indicate all hardian germs $y$ such that $(y')^2=P(y)$. First, the $y\in \R$ such that~$(y')^2=P(y)$ are exactly $e_1$, $e_2$, $e_3$. If there is a hardian $y\notin \R$ with~$(y')^2=P(y)$, then by the proof of
%Proposition~\ref{prop:wp}, either $e_1=e_2=e_3=0$,  or $e_1=e_2>0$.  In Corollaries~\ref{eee=0} and ~\ref{lem:e_1=e_2>0} 
%we deal exhaustively with these two cases.  
If $g_3=0$, then  $e_1=e_2=e_3=0$, and if $g_3<0$, then~$e_1=e_2>0$. In Corollaries~\ref{eee=0} and~\ref{lem:e_1=e_2>0} 
we deal exhaustively with these two cases. In particular, we show there that in each case there is a hardian $y\notin \R$ such
that $(y')^2=P(y)$, thus finishing the proof of  Proposition~\ref{prop:wp}. 

Accordingly we assume below that $g_2, g_3\in \R$ and $\Delta=0$, so $e_1, e_2, e_3\in \R$.  Note that the $y\in \R$ such that~$(y')^2=P(y)$ are exactly $e_1$, $e_2$, $e_3$.

\begin{cor}\label{eee=0}
Suppose $e_1=e_2=e_3=0$ and $y\in\c^1$. Then 
$$\text{$y\in\c^\times$ and $(y')^2=P(y)$}\quad\Longleftrightarrow\quad  \text{$y=\frac{1}{(x-c)^2}$ for some $c\in\R$.}$$
\end{cor}
\begin{proof}
We have $P(Y)=4Y^3$.  The direction $\Leftarrow$ is routine.  
For $\Rightarrow$, suppose $y\in\c^\times$ and $(y')^2=4y^3$.
Then~$y'\in\c^\times$, $y>0$, $z:=y^{-1/2}>0$, and $z'=-\frac{1}{2}y^{-3/2}y'$.
We have~$y'<0$: otherwise~$0<y'=2y^{3/2}$ and so $z'=-1$, hence $z<0$,
a contradiction. Therefore~$y'=-2y^{3/2}$, so $z'=1$  and thus
$y=\frac{1}{z^2}=\frac{1}{(x-c)^2}$ for some $c\in\R$.
\end{proof}

\noindent
Lemma~\ref{lem:e_1=e_2>0} below is an analogue of Lemma~\ref{eee=0} for $e_1=e_2>0$, but
first we make some observations about hyperbolic functions. Recall that for $t\in \R$, 
\begin{align*}  \sinh t\ :&=\ \frac{1}{2}(\ex^t-\ex^{-t}), \qquad \cosh t\ :=\ \frac{1}{2}(\ex^t+\ex^{-t}), \text{ so}\\
\cosh^2t-\sinh^2t\ &=\ 1,\qquad \frac{d}{dt}\sinh t\ =\ \cosh t,\qquad \frac{d}{dt} \cosh t\ =\ \sinh t.
\end{align*}
We also set for $t\in\R$: 
$$
\operatorname{sech} t\ :=\ \frac{1}{\cosh t} \text{ (hyperbolic secant)},\qquad
\tanh t\ :=\ \frac{\sinh t}{\cosh t} \text{ (hyperbolic tangent)}$$
and for $t\neq 0$:
$$\operatorname{csch} t\ :=\ \frac{1}{\sinh t} \text{ (hyperbolic cosecant)},\quad
\coth t\ :=\ \frac{\cosh t}{\sinh t} \text{ (hyperbolic cotangent).}$$
Now $\sinh\colon\R\to \R$ is an increasing bijection, so 
$t\mapsto \operatorname{csch} t\colon \R^{>}\to \R^{>}$ is a decreasing bijection. 
We have $\operatorname{sech}^2t=1-\tanh^2 t$, and for $t\ne 0$,  $\operatorname{csch}^2 t=\coth^2 t-1$.
Moreover, 
$\frac{d}{dt}\operatorname{sech}t=-\tanh t \operatorname{sech} t$ and for~${t\ne 0}$: 
$\frac{d}{dt}\operatorname{csch}t=-\coth t \operatorname{csch} t$.
Hence both $-\operatorname{sech}^2\colon\R\to\R$ and
 $\operatorname{csch}^2\colon\R^\times\to\R$ satisfy the differential equation
$(u')^2 = 4u^2(u+1)$.
We use these facts to prove:

\begin{lemma}\label{wecsch}
Let $w\in \c^1$ and $e\in\R^>$. Then the following are equivalent:
\begin{enumerate}
\item[\rm(i)] $w(t)>0$, eventually, and $(w')^2=ew^2(w+1)$
\item[\rm(ii)]  
 $w=\operatorname{csch}^2\circ (c+\frac{\sqrt{e}}{2}x)$  for some $c\in \R$.
 \end{enumerate}
\end{lemma}
\begin{proof}  Direct computation gives (ii)~$\Leftarrow$~(i). Assume (i). 
Consider the  decreasing bijection $t\mapsto u(t) :=\operatorname{csch}^2(t) \colon \R^>\to\R^>$; let $u^{\operatorname{inv}}: \R^{>}\to \R^{>}$ be its (strictly decreasing) compositional inverse,
so $u^{\operatorname{inv}}\in \c^1(\R^{>})$. Then for $v:=u^{\operatorname{inv}}\circ w\in \c^1$ we have $v(t)>0$, eventually,
$v'=\frac{w'}{u'\circ v}$, and $u\circ v=w$. Thus
$$(v')^2\ =\ \frac{(w')^2}{(u')^2\circ v}\  =\ \frac{e}{4}.$$
Hence $v'=\frac{\sqrt{e}}{2}$, since $v'=-\frac{\sqrt{e}}{2}$ contradicts $v>0$. Now use $w=u\circ v$. 
% we obtain some $c\in\R$ such that~$w=\operatorname{csch}^2(c-\frac{\sqrt{e}}{2}x)$ or~$w=\operatorname{csch}^2(c+\frac{\sqrt{e}}{2}x)$.
 \end{proof}

\noindent
The increasing bijection $t\mapsto \cosh t\colon (0, +\infty)\to (1,+\infty)$ yields the increasing bijection 
$t\mapsto -\operatorname{sech}^2 t\colon  (0, +\infty)\to (-1,0)$. We use this to prove likewise:

\begin{lemma}\label{wesech}
Let $w\in\c^1$ and $e\in\R^>$. Then the following are equivalent: \begin{enumerate}
\item[\rm(i)] $-1 < w(t) < 0$, eventually, and $(w')^2\ =\ ew^2(w+1)$;
\item[\rm(ii)] $ w=-\operatorname{sech}^2\circ (c+\frac{\sqrt{e}}{2}x)$  for some $c\in \R$.
\end{enumerate}
\end{lemma}
%\begin{proof} Direct computation gives (ii)$\Rightarrow$(i). Assume (i).  Consider the increasing bijection $t\mapsto u(t) :=-\operatorname{sech}^2(t) \colon (0,+\infty)\to (-1,0)$; let $u^{\operatorname{inv}}: (-1,0)\to (0,+\infty)$ be its (increasing) compositional inverse,
%so $u^{\operatorname{inv}}\in \c^1\big((-1,0)\big)$. Then for $v:=u^{\operatorname{inv}}\circ w\in \c^1$ we have $v(t)>0$, eventually, {\bf oksofar} 
%$v'=\frac{w'}{u'\circ v}$, and $u\circ v=w$. Thus
%$$(v')^2\ =\ \frac{(w')^2}{(u')^2\circ v}\  =\ \frac{e}{4}.$$
%Hence $v'=\frac{\sqrt{e}}{2}$, since $v'=-\frac{\sqrt{e}}{2}$ contradicts $v>0$. Now use $w=u\circ v$.
%\end{proof} \marginpar{above proof checked but commented out} 

\begin{cor}\label{lem:e_1=e_2>0}
Suppose $e_1=e_2>0$. Then the hardian germs~$y\notin \R$ such that~$(y')^2=P(y)$ are all in $\R(\ex^{x\sqrt{3e_1}})$ and are
given by
$$y\ =\ e_1+3e_1\cdot \operatorname{csch}^2\circ (c+x\sqrt{3e_1}),\quad
y\ =\ e_1-3e_1\cdot\operatorname{sech}^2\circ (c+x\sqrt{3e_1}),$$
where $c\in\R$.
\end{cor}
\begin{proof}
We have $P(Y)=4(Y-e_1)^2(Y+2e_1)$. Let $y\in\c^1$ and
$w:=(y-e_1)/3e_1$. Then $(y')^2=P(y)$ iff $(w')^2=12e_1w^2(w+1)$. There is no hardian
$y< -2e_1$ with~$(y')^2=P(y)$, so we can use Lemmas~\ref{wecsch} and ~\ref{wesech} with $e:=12e_1$. 
\end{proof}

\subsection*{Uniform finiteness}
We now let $H$ be a Hardy field and $\varphi(x,y)$ and $\theta(x)$ be   $\mathcal L_H$-formulas, where $x=(x_1,\dots,x_m)$ and $y=(y_1,\dots,y_n)$.

\begin{lemma}
There is a $B=B(\varphi)\in\N$ such that for all $f\in H^m$: if 
for some $\d$-maximal Hardy field extension $M$ of $H$ there are more than $B$ tuples $g\in M^n$ with $M\models\varphi(f,g)$, then 
for every $\d$-maximal Hardy field extension $M$ of $H$ there are infinitely many $g\in M^n$ with~$M\models\varphi(f,g)$. 
\end{lemma}
\begin{proof}
Fix a $\d$-maximal Hardy field extension $M^*$ of $H$. By \cite[Proposition~6.4]{ADHdim} we have $B=B(\varphi)\in\N$ such that
for all $f\in (M^*)^m$: if $M^*\models\varphi(f,g)$ for more than~$B$ many $g\in (M^*)^n$, then
$M^*\models\varphi(f,g)$ for infinitely many $g\in (M^*)^n$. Now use Theorem~\ref{thm:transfer}. 
\end{proof}

\noindent
In the proof of the next lemma we use that $\mathcal C$ has the cardinality $\mathfrak c=2^{\aleph_0}$ of the continuum, hence
$\abs{H}=\mathfrak c$ if $H\supseteq\R$.

\begin{lemma}
Suppose $H$ is $\d$-maximal and $S:=\big\{f\in H^m:H\models\theta(f)\big\}$ is infinite. Then $\abs{S}=\mathfrak c$.
\end{lemma}
\begin{proof}
Let $d:=\dim(S)$ be the dimension of the definable set $S\subseteq H^m$ as introduced in \cite{ADHdim}.
If $d=0$, then $\abs{S}=\abs{\R}=\mathfrak c$ by remarks following  \cite[Pro\-po\-si\-tion~6.4]{ADHdim}.
Suppose $d>0$, and for $g=(g_1,\dots,g_m)\in H^m$ and $i\in\{1,\dots,m\}$, let~$\pi_i(g):=g_i$. 
Then for some $i\in\{1,\dots,m\}$, the subset $\pi_i(S)$ of $H$ has nonempty interior, by~\cite[Corollary~3.2]{ADHdim}, and hence 
  $\abs{S}=\abs{H}=\mathfrak c$.
\end{proof}

\noindent
The two lemmas above together now yield Corollary~\ref{cor:uniform finiteness} from the introduction. 

\subsection*{Transfer between maximal Hardy fields and transseries}
Let   $\boldsymbol T$ be the unique expansion of  $\mathbb T$ to a pre-$\HLO$-field, so~$\boldsymbol T$ is an $H$-closed 
$\HLO$-field with small derivation and constant field $\R$.

\begin{lemma}\label{lem:unique HLO-expansion, 1}
Let $H$ be a pre-$H$-subfield of $\mathbb T$ with $H\not\subseteq \R$.
Then $H$ has a unique expansion to a pre-$\HLO$-field.
\end{lemma}
\begin{proof} If $H$ is grounded, this follows from [ADH, 16.3.19].
Suppose $H$ is not grounded. Then $H$ has asymptotic integration by the proof of [ADH, 10.6.19] applied to $\Delta:= v(H^\times)$. Starting with an $h_0\succ 1$ in $H$
with $h_0'\asymp 1$ we construct a logarithmic sequence $(h_n)$ in $H$ as in [ADH, 11.5], so $h_n\asymp\ell_n$ for all $n$.
Hence $\Gamma^<$ is cofinal in $\Gamma_{\mathbb T}^<$, so $H$ is $\upo$-free by [ADH, remark before 11.7.20].
Now use [ADH, 16.3.19] again. 
\end{proof}

%\begin{lemma}\label{lem:unique HLO-expansion, 1}
%Let $H$ be a pre-$H$-subfield of $\mathbb T$ properly containing $\R$.
%Then $H$ has a unique expansion to a pre-$\HLO$-field.
%\end{lemma}
%\begin{proof} If $H$ is grounded, this follows from [ADH, 16.3.19].
%Suppose $H$ is not grounded. Then $H$ has asymptotic integration by [ADH, 10.6.19].  Starting with an $h_0\succ 1$ in $H$
%with $h_0'\asymp 1$ we construct a logarithmic sequence $(h_n)$ in $H$ as in [ADH, 11.5], so $h_n\asymp\ell_n$ for all $n$.
%Hence $\Gamma^<$ is cofinal in $\Gamma_{\mathbb T}^<$, so $H$ is $\upo$-free by [ADH, remark before 11.7.20].
%Now use [ADH, 16.3.19] again. 
%\end{proof}

\noindent
In the rest of this subsection $H$ is a Hardy field with canonical $\HLO$-ex\-pan\-sion~$\boldsymbol H$, and $\iota\colon H\to \mathbb T$ is an embedding of ordered differential fields, and thus of pre-$H$-fields.  %Then by  Lemma~\ref{lem:unique HLO-expansion, 1} and a separate treatment of the case $H\subseteq \R$: 

\begin{cor}\label{lem:unique HLO-expansion, 2}
The map $\iota$ is an embedding~$\boldsymbol H\to\boldsymbol T$ of pre-$\HLO$-fields.
\end{cor}
\begin{proof}
If $H\not\subseteq\R$, then this follows from Lemma~\ref{lem:unique HLO-expansion, 1}.
Suppose $H\subseteq\R$. Then $\iota$ is the identity on $H$, so extends to the embedding
$\R(x)\to \mathbb T$ that is the identity on~$\R$ and sends the germ $x$ to $x\in \T$. Now use that $\R(x)\not\subseteq \R$ and Corollary~\ref{cor:canonical HLO}. 
%so $H$ has gap $0$. Then by [ADH, 16.3.13], $H$ has exactly two $\HLO$-cuts $(I,\Upl,\Upo)$,
%one with~$I=\I(H)=\{0\}$ and one with $I=H$. If $M$ is a $H$-field extension of $H$ with small derivation,
%then~$\I(M)\subseteq M^{\prec 1}$,
%thus $\I(M)\cap H=\{0\}$. Hence $\boldsymbol H=(\{0\},\dots)$ since each maximal Hardy field has small derivation,
%and $\I(\T)\cap \iota(H)=\{0\}$  since~$\mathbb T$ also has small derivation; thus $\iota$ is an embedding~$\boldsymbol H\to\boldsymbol T$.
\end{proof}

%\begin{lemma}\label{lem:unique HLO-expansion, 2}
%$\iota$ is an embedding~$\boldsymbol H\to\boldsymbol T$ of pre-$\HLO$-fields.
%\end{lemma}
%\begin{proof}
%To show this we may replace $\iota$ by any extension $H^*\to\mathbb T$ to an embedding of a Hardy field $H^*\supseteq H$ into $\mathbb T$,
%thanks to Corollary~\ref{cor:canonical HLO}. Doing so, using [ADH, 10.5.13] we first arrange that $H$ is an $H$-field, by replacing $\iota$ by its extension to an embedding of the $H$-field hull of $H$
%inside the Hardy field~$H(\R)$ into $\mathbb T$. Next replace~$\iota$ by its extension to an
%embedding~$H(\R)\to\mathbb T$ to arrange that $\R\subseteq H$, by [ADH, remark before 4.6.21] and [ADH, 10.5.16].
%Now if $\R\neq H$, then we are done by Lemma~\ref{lem:unique HLO-expansion, 1}.
%Suppose $\R=H$, so $H$ has gap $0$. Then by [ADH, 16.3.13], $H$ has exactly two $\HLO$-cuts $(I,\Upl,\Upo)$,
%one with $I=\I(H)=\{0\}$ and one with $I=H=\R$. If $M$ is a $H$-field extension of $H$ with small derivation,
%then~$\I(M)\subseteq M^{\prec 1}$,
%thus $\I(M)\cap H=\{0\}$. Now use that each Hardy field as well as~$\mathbb T$ have small derivation.
%\end{proof}

\noindent
%We let $\mathcal L$ as before be the language of valued ordered differential rings.  
Recall from~[ADH, B.4] that for any $\mathcal L_H$-sentence $\sigma$ we obtain an $\mathcal L_{\mathbb T}$-sentence $\iota(\sigma)$ by replacing the name of each~$h\in H$ occurring in $\sigma$ with the name of $\iota(h)$.

\begin{cor}\label{cor:transfer T, 1}
Let  $\sigma$ be an $\mathcal L_H$-sentence. Then \textup{(i)}--\textup{(iv)} in Theorem~\ref{thm:transfer} are also equivalent to:
\begin{enumerate}
\item[\textup{(v)}] $\mathbb T\models\iota(\sigma)$.
\end{enumerate}
\end{cor}
\begin{proof}
Let $M$ be a $\d$-maximal Hardy field extension of $H$; it suffices to show that~$M\models\sigma$ iff $\mathbb T\models\iota(\sigma)$. For this, mimick the proof of (i)~$\Rightarrow$~(ii) in Theorem~\ref{thm:transfer}, using Corollary~\ref{lem:unique HLO-expansion, 2}. 
\end{proof}

\noindent
Corollary~\ref{cor:transfer T, 1} yields the first part of Corollary~\ref{cor:systems, 3} from the introduction, even in a stronger
form.  After an intermezzo on 
differential closure in Section~\ref{sec:diff closure} we prove the second part of that corollary in Section~\ref{sec:embeddings into T}: Corollary~\ref{cor:transfer T, 2}. There we also use:

\begin{lemma}\label{lemhr} $\iota$ extends uniquely to an embedding $H(\R)\to \mathbb{T}$ of pre-$H$-fields. 
\end{lemma}
\begin{proof}  Let $\hat{H}$ be the $H$-field hull of $H$ in $H(\R)$. Then $\iota$ extends uniquely to an $H$-field embedding $\hat{\iota}: \hat{H}\to \T$ by [ADH, 10.5.13]. By [ADH, remark before 4.6.21] and [ADH, 10.5.16] $\hat{\iota}$ extends uniquely to an embedding $H(\R)\to \T$ of $H$-fields. 
\end{proof} 

\noindent
We finish with indicating how Theorem~\ref{thm:A} from the introduction (again, in strengthened form)  follows from \cite{JvdH} and the results above:

\begin{cor}\label{divpcor}
If $P\in H\{Y\}$, $f<g$ in $H$, and $P(f)<0<P(g)$, then each $\d$-maximal Hardy field extension of $H$ contains a
$y$ with $f<y<g$ and $P(y)=0$.
\end{cor}
\begin{proof}
By  \cite{JvdH}, 
%\marginpar{DIVP  and $H$-closed for $\mathbb T_{\text{g}}$ taken on faith for now}
the ordered differential field $\mathbb T_{\text{g}}$ of grid-based transseries is $H$-closed with
small derivation and has the differential intermediate value property  (DIVP).
Hence $\mathbb T$ also has  DIVP, by completeness of $T_H$ (see the introduction). Now use Corollary~\ref{cor:transfer T, 1}.
\end{proof}

\begin{cor}\label{cor:odd degree}
Let $P\in H\{Y\}$ have odd degree. Then there is an $H$-hardian germ $y$   with $P(y)=0$.
\end{cor}
\begin{proof}
This follows from Theorem~\ref{thm:extend to H-closed} and [ADH, 14.5.3]. Alternatively, we can use Corollary~\ref{divpcor}:  Replace $H$ by $\operatorname{Li}\!\big(H(\R)\big)$ to arrange that $H\supseteq\R$ is Liouville closed, and appeal to
 the example following Corollary~\ref{cor:val at infty}.
%we have~$\sgn P(-y)=-\sgn P(y)$ for all sufficiently large $y$ in~$H$.
\end{proof}

\noindent
Note that if $H\subseteq\Ginf$, then in the previous two corollaries we have $H\langle y\rangle\subseteq\Ginf$, by Corollary~\ref{cor:Hardy field ext smooth}; likewise with $\Gom$ in place if $\Ginf$.

\section{Relative Differential Closure}\label{sec:diff closure}

\noindent
Let $K\subseteq L$ be an extension of differential fields, and let $r$ range over $\N$. We say that~$K$ is {\bf $r$-differentially closed} in $L$ 
for every~$P\in K\{Y\}^{\neq}$ of order~$\leq r$, each zero of $P$ in~$L$   lies in~$K$. 
We also say that $K$ is {\bf weakly $r$-differentially closed} in~$L$ if 
 every~$P\in K\{Y\}^{\ne}$ of order~$\leq r$ with a zero in~$L$ has a zero in~$K$.
 We abbreviate ``$r$-differentially closed'' by ``$r$-$\d$-closed.''
 Thus 
$$\text{$K$ is $r$-$\d$-closed in $L$}\quad\Longrightarrow\quad\text{$K$ is  weakly $r$-$\d$-closed in~$L$,}$$ 
and
\begin{align*}
\text{$K$ is $0$-$\d$-closed in $L$}&\quad\Longleftrightarrow\quad 
  \text{$K$ is weakly $0$-$\d$-closed in $L$} \\ &\quad\Longleftrightarrow\quad
  \text{$K$ is algebraically closed in $L$.}
\end{align*}
Hence 
\begin{equation}\label{eq:0-d-closed}
\text{$K$ is weakly $0$-$\d$-closed in $L$}\quad\Longrightarrow\quad\text{$C$ is algebraically closed in $C_L$.}
\end{equation}
Also, if~$K$ is weakly $0$-$\d$-closed in $L$
and $L$ is algebraically closed, then $K$ is algebraically closed, and similarly
with ``real closed'' in place of ``algebraically closed''. In [ADH, 5.8] we defined
$K$ to be {\it  weakly $r$-$\d$-closed}\/ if every $P\in K\{Y\}\setminus K$ of order~$\leq r$ has a zero in $K$.
Thus
$$\text{$K$ is weakly $r$-$\d$-closed}
\ \Longleftrightarrow\  \begin{cases} &\parbox{21em}{$K$ is weakly $r$-$\d$-closed in every differential field extension of $K$.}\end{cases}$$
If
$K$ is weakly   $r$-$\d$-closed in $L$, then $P(K)=P(L)\cap K$ for all $P\in K\{Y\}$ of or\-der~$\leq r$;
in particular, 
\begin{equation}\label{eq:weakly 1-d-closed}
\text{$K$ is weakly $1$-$\d$-closed in $L$}\quad\Longrightarrow\quad \der K=\der L\cap K.
\end{equation}
Also, 
\begin{equation}\label{eq:1-d-closed}
\text{$K$ is $1$-$\d$-closed in $L$}\quad\Longrightarrow\quad C=C_L\text{ and } K^\dagger=L^\dagger\cap K.
\end{equation}
Moreover:

\begin{lemma}\label{lem:weakly r-d-closed}
Suppose $K$ is weakly $r$-$\d$-closed in $L$. If $L$ is $r$-linearly surjective, then so is $K$, and if $L$
is
$(r+1)$-linearly closed, then so is $K$.
\end{lemma}
\begin{proof}
The first statement is clear from the remarks preceding the lemma, and the second statement is shown similarly to  [ADH, 5.8.9].
\end{proof}

\noindent
Sometimes we get more than we bargained for:

\begin{lemma}
Suppose $K$ is not algebraically closed, $C\ne K$, and
$K$ is weakly $r$-$\d$-closed in $L$. Let $Q_1,\dots,Q_m\in K\{Y\}^{\neq}$ of order~$\leq r$ have a common zero
in~$L$,~$m\ge 1$. Then they have a common zero in $K$.
\end{lemma}
\begin{proof}
%We may assume $C\neq K$.
Take a polynomial $\Phi\in K[X_1,\dots,X_m]$ whose only zero in $K^m$ is the origin~$(0,\dots,0)\in K^m$. Then the differential
polynomial~$P:=\Phi(Q_1,\dots,Q_m)\in K\{Y\}$ is nonzero (use [ADH, 4.2.1]) and
has order~$\leq r$. For $y\in L$ we have $$Q_1(y)=\cdots=Q_m(y)=0\quad \Longrightarrow\quad P(y)=0,$$
and for $y\in K$ the converse of this implication also holds. %we have  \[P(y)=0\quad \Longrightarrow\quad Q_1(y)=\cdots=Q_m(y)=0.\qedhere\]
\end{proof}

\noindent
We say that $K$ is {\bf differentially closed}\index{differentially closed} in~$L$ if $K$ is 
$r$-$\d$-closed in~$L$ for each~$r$, and similarly we define when  $K$ is {\bf weakly differentially closed}\index{differentially closed!weakly} in~$L$.
We also use ``$\d$-closed'' to abbreviate ``differentially closed''.
If~$K$, as a differential ring, is an elementary substructure of~$L$, then $K$ is weakly
$\d$-closed in~$L$.
The elements of $L$ that are $\d$-algebraic over $K$ form
the smallest differential subfield of $L$ containing~$K$ which is  $\d$-closed in $L$; we call it the
{\bf differential closure} (``$\d$-closure'' for short) of~$K$ in $L$.\index{differential closure} Thus $K$ is $\d$-closed in $L$ iff   no   $\d$-subfield of $L$
properly containing~$K$ is
 $\d$-algebraic over $K$.
This notion of being differentially closed does not seen prominent in the differential algebra literature, though the definition occurs (as ``differentially algebraic closure'') in~\cite[p.~102]{Kolchin-DA}. Here is a useful fact about it: 

\begin{lemma}\label{lem:dc base field ext}
Let $F$ be a differential field extension of $L$ and $E$ be a   subfield of~$F$ containing $K$ such that $E$ is algebraic over $K$ and $F=L(E)$. 
\[
\xymatrix{& F=L(E) & \\
E \ar@{-}[ur] & & L \ar@{-}[ul] \\
& K \ar@{-}[ul] \ar@{-}[ur] & } 
\]
Then $K$ is $\d$-closed in $L$ iff $E\cap L = K$ and $E$ is $\d$-closed in $F$.
\end{lemma}

\begin{proof}
Suppose $K$ is $\d$-closed in $L$. Then $K$ is algebraically closed in $L$, so $L$ is linearly disjoint from~$E$ over $K$. (See \cite[Chapter~VIII, \S{}4]{Lang}.) In particular~${E\cap L=K}$.
Now let $y\in F$ be $\d$-algebraic over $E$; we  claim that $y\in E$. Note that $y$ is $\d$-algebraic over~$K$. 
Take a field extension $E_0\subseteq E$ of $K$ with $[E_0:K]<\infty$ (so $E_0$ is a $\d$-subfield of~$E$) such that $y\in L(E_0)$;
replacing~$E$,~$F$ by $E_0$, $L(E_0)$, respectively, we arrange that $n:=[E:K]<\infty$.
Let $b_1,\dots,b_n$ be a basis of the $K$-linear space~$E$;
then~$b_1,\dots,b_n$ is also a basis of the $L$-linear space $F$.
Let $\sigma_1,\dots,\sigma_n$ be the distinct field embeddings~$F\to L^{\operatorname{a}}$ over $L$.
Then the vectors
$$\big(\sigma_1(b_1),\dots,\sigma_1(b_n)\big),\dots, \big(\sigma_n(b_1),\dots,\sigma_n(b_n)\big)\in (L^{\operatorname{a}})^n$$
are $L^{\operatorname{a}}$-linearly independent \cite[Chapter~VI, Theorem~4.1]{Lang}. %\marginpar{statement in Lang is probably defective}
Let $a_1,\dots,a_n\in L$ be such that $y=a_1b_1+\cdots+a_nb_n$. Then
$$\sigma_j(y) = a_1\sigma_j(b_1)+\cdots+a_n\sigma_j(b_n)\quad\text{for $j=1,\dots,n$,}$$
hence by Cramer's Rule, 
$$a_1,\dots,a_n\in K\big(\sigma_j(y),\sigma_j(b_i):i,j=1,\dots,n\big).$$
Therefore $a_1,\dots,a_n$ are $\d$-algebraic over $K$, since    $\sigma_j(y)$ and $\sigma_j(b_i)$ for $i,j=1,\dots,n$ are.
Hence $a_1,\dots,a_n\in K$ since $K$ is $\d$-closed in $L$, so $y\in E$ as claimed. This shows the forward implication. The
backward direction is clear.
\end{proof}

\begin{cor}\label{cor:dc base field ext}
If $-1$ is not a square in $L$ and $\imag$ in a differential field extension of $L$ satisfies $\imag^2=-1$,
then:  $K$ is $\d$-closed in $L\ \Leftrightarrow\ K[\imag]$ is $\d$-closed in $L[\imag]$.
\end{cor} 

%\begin{remark}
%Suppose $K$ is real closed and $L=K[\imag]$ ($\imag^2=-1$). Is $L$ is weakly $\d$-closed if and only if every $P\in K\{Y\}$
%of odd degree has a zero in $K$? We do not know.
%\end{remark}

\noindent
In the next lemma {\em extension\/} refers to an extension of valued differential fields. 

\begin{lemma}\label{lem:r-newtonian descends} 
Suppose $K$ is an  $\upl$-free $H$-asymptotic field and is $r$-$\d$-closed in an $r$-newtonian ungrounded $H$-asymptotic extension $L$. Then $K$ is also $r$-newtonian.
\end{lemma}
\begin{proof}
Let $P\in K\{Y\}^{\neq}$ be quasilinear of order~$\leq r$. Then $P$ remains quasilinear when viewed as differential polynomial over $L$, by 
Lemma~\ref{lem:13.7.10}. Hence $P$ has a zero~$y\preceq 1$ in $L$, which lies in $K$ since $K$ is $r$-$\d$-closed in $L$.
\end{proof}

\subsection*{Relative differential closure in $H$-fields}
We now return to the $H$-field setting.  Let $\mathcal L_\der=\{0,1,{-},{+},{\,\cdot\,},{\der}\}$  be the language of differential rings, 
a sublanguage of the language $\mathcal L=\mathcal{L}_\der \cup\{\le, \preceq\}$ of ordered valued differential rings from Section~\ref{sec:transfer}.

 Let now $M$ be an $H$-closed field,  and let
$H$ a pre-$H$-subfield of $M$ whose valuation ring and constant field we denote by $\mathcal{O}$ and $C$.
Construing~$H$ and~$M$ as $\mathcal L$-structures in the usual way,
$H$ is an $\mathcal L$-substructure of $M$.
We also use the sublanguage~$\mathcal L_{\preceq}:=\mathcal L_\der\cup\{ {\preceq} \}$ of~$\mathcal L$, so $\mathcal L_{\preceq}$ is the language of valued differential rings. We
expand the $\mathcal L_\der$-structure~$H[\imag]$ to an   $\mathcal L_{\preceq}$-structure by interpreting $\preceq$ 
as the dominance relation associated to the valuation ring $\mathcal O+\mathcal O\imag$ of $H[\imag]$;
we expand likewise~$M[\imag]$ to an   $\mathcal L_{\preceq}$-structure by interpreting~$\preceq$ 
as the dominance relation associated to the valuation ring $\mathcal O_{M[\imag]}=\mathcal O_M+\mathcal O_M\imag$ of $M[\imag]$. 
Then $H[\imag]$ is an $\mathcal L_{\preceq}$-substructure of~$M[\imag]$.
By~$H\preceq_{\mathcal{L}} M$ we mean that $H$ is an elementary $\mathcal L$-substructure of $M$, and we use expressions like ``$H[\imag]\preceq_{\mathcal{L}_{\preceq}} M[\imag]$''  in the same way; of course, the two
uses of the symbol $\preceq$ in the latter are unrelated. 

By 
Corollary~\ref{cor:dc base field ext}, $H$ is $\d$-closed in~$M$ iff $H[\imag]$ is $\d$-closed in $M[\imag]$.

\begin{lemma}
Suppose $M$ has small derivation. Then 
$$H\ \preceq_{\mathcal L_\der}\ M\ \Longleftrightarrow\ H[\imag]\ \preceq_{\mathcal L_\der}\ M[\imag].$$
Also, if $H\preceq_{\mathcal L_\der} M$, then $H\preceq_{\mathcal L} M$ and  $H[\imag]\preceq_{\mathcal{L}_{\preceq}} M[\imag]$.
%In this case, if in addition~$C=C_M$, then
%$H$ is an elementary $\mathcal L$-substructure of $M$ and
%$H[\imag]$ is an elementary $\mathcal L_{\preceq}$-substructure of $M[\imag]$.
\end{lemma}
\begin{proof}
The  forward direction in the equivalence is obvious. For the converse, let~$H[\imag]\preceq_{\mathcal L_\der}M[\imag]$.
We have  $M\equiv_{\mathcal{L}_{\der}}\mathbb T$ by [ADH, 16.6.3].
Then  [ADH, 10.7.10] yields an $\mathcal L_\der$-formula defining $M$ in $M[\imag]$,  so the same formula defines $M\cap H[\imag]=H$ in~$H[\imag]$, and thus $H\preceq_{\mathcal L_\der} M$. For the ``also'' part, use that the squares of $M$ are the nonnegative elements in its ordering,  that $\mathcal O_M$ is then definable as the convex hull of $C_M$ in $M$ with respect to this ordering, and
if $H\preceq_{\mathcal L_\der} M$, then each $\mathcal L_\der$-formula defining $\mathcal O_M$ in $M$ also defines $\mathcal O=\mathcal O_M\cap H$ in $H$.
 %$H$ is a real closed $H$-field and 
%$\leq$,~$\preceq$ are definable in the $\mathcal L_\der$-reducts of $H$, $M$, and so
 %$H$ is an    elementary $\mathcal L$-substructure of $M$. This also implies that
%$H[\imag]$ is an elementary $\mathcal L_{\preceq}$-substructure of $M[\imag]$.
\end{proof}

\noindent
The next proposition complements~[ADH, 16.0.3, 16.2.5]:  

\begin{prop}\label{prop:16.0.3 converse}
The following are equivalent:
\begin{enumerate}
\item[\textup{(i)}] $H$ is $\d$-closed in $M$;
\item[\textup{(ii)}]  $C=C_M$ and $H\preceq_{\mathcal L} M$;
\item[\textup{(iii)}] $C=C_M$ and $H$ is $H$-closed.
\end{enumerate}
\end{prop}
\begin{proof} Assume (i).  Then $C=C_M$ and $H$ is a Liouville closed $H$-field, by \eqref{eq:0-d-closed}, \eqref{eq:weakly 1-d-closed}, and \eqref{eq:1-d-closed}.
We have $\omega(M)\cap H=\omega(H)$ since $H$ is weakly $1$-$\d$-closed in~$M$, and
$\sigma\big(\Upg(M)\big)\cap H=\sigma\big(\Upg(M)\cap H\big)=\sigma\big(\Upg(H)\big)$ since~$H$ is $2$-$\d$-closed in $M$ and~$\Upg(M)\cap H=\Upg(H)$ by~[ADH, p.~520]. Now~$M$ is
Schwarz closed~[ADH, 14.2.20], so $M=\omega(M)\cup\sigma\big(\Upg(M)\big)$, hence also~$H=\omega(H)\cup\sigma\big(\Upg(H)\big)$, thus 
$H$ is also
Schwarz closed [ADH, 11.8.33]; in particular, $H$ is $\upo$-free.
By Lemma~\ref{lem:r-newtonian descends}, $H$ is newtonian.
This shows (i)~$\Rightarrow$~(iii). The implication  (iii)~$\Rightarrow$~(i) is~[ADH, 16.0.3], and~(iii)~$\Leftrightarrow$~(ii) follows from [ADH, 16.2.5].
\end{proof}

\noindent
Next a consequence of [ADH, 16.2.1], but note first that $H(C_M)$ is an $H$-subfield of $M$
and $\d$-algebraic over $H$, and  recall that each $\upo$-free $H$-field
has a Newton-Liouville closure,   as defined in~[ADH, p.~669].

\[
\xymatrix{
& M \ar@{-}[d] & \\
& H(C_M) & \\
C_M \ar@{-}[ur] & & H \ar@{-}[ul]
}
\]

\begin{cor}\label{cor:nlc, 1} %\marginpar{maybe also with inclusion diagram} 
If $H$ is  $\upo$-free, then the differential closure of~$H$ in $M$ is a New\-ton-Liou\-ville closure of the $\upo$-free $H$-subfield $H(C_M)$ of $M$. 
\end{cor}

\noindent
Let $\boldsymbol M$ be the expansion of $M$ to a $\HLO$-field, and  let   $\boldsymbol H$, $\boldsymbol H(C_M)$ be the   expansions of  $H$, $H(C_M)$, respectively, to   pre-$\HLO$-subfields of $\boldsymbol M$;  then   $\boldsymbol H(C_M)$ is a $\HLO$-field.
By   Pro\-po\-si\-tion~\ref{prop:16.0.3 converse}, the  $\d$-closure 
 $H^{\operatorname{da}}$ of $H$ in~$M$ is $H$-closed and hence has a unique expansion  $\boldsymbol H^{\operatorname{da}}$  
to a $\HLO$-field. 
Then $\boldsymbol H\subseteq \boldsymbol H(C_M)\subseteq \boldsymbol H^{\operatorname{da}}\subseteq\boldsymbol M$.
For the Newton-Liouville closure of a pre-$\HLO$-field, see [ADH, 16.4.8].

\begin{cor}\label{cor:nlc, 2}  
The $\HLO$-field $\boldsymbol H^{\operatorname{da}}$ is a New\-ton-Liouville closure of~$\boldsymbol H(C_M)$.
\end{cor}
\begin{proof}
Let $\boldsymbol H(C_M)^{\operatorname{nl}}$ be a  Newton-Liouville closure of $\boldsymbol H(C_M)$.
Since $\boldsymbol H^{\operatorname{da}}$   is $H$-closed and extends $\boldsymbol H(C_M)$,  there is an embedding $\boldsymbol H(C_M)^{\operatorname{nl}}\to \boldsymbol H^{\operatorname{da}}$ over $\boldsymbol H(C_M)$, and any such embedding is an isomorphism, thanks to  [ADH, 16.0.3].
\end{proof}

\subsection*{Relative differential closure in Hardy fields}
Specializing to Hardy fields, assume below that
$H$ is a Hardy field and set $K:=H[\imag]\subseteq \Calinf[\imag]$, an $H$-asymptotic extension of $H$. 
By definition, $H$ is $\d$-maximal iff~$H$ is $\d$-closed in every Hardy field extension of $H$. The following contains Corollary~\ref{cor:zeros in complexified Hardy field extensions} from the introduction: 

\begin{cor}\label{cor:d-max weakly d-closed}
Suppose $H$ is $\d$-maximal. Then $K$ is weakly $\d$-closed, hence linearly closed  by \textup{[ADH, 5.8.9]}, and linearly surjective.  If~$E$ is a Hardy field extension of~$H$, then
$K$ is $\d$-closed in $E[\imag]$.
\end{cor}
\begin{proof}
By our main Theorem~\ref{thm:char d-max}, $H$ is newtonian, hence  $K$ is  weakly $\d$-closed by [ADH, 14.5.7, 14.5.3],
proving the first statement; the second statement follows from Corollary~\ref{cor:dc base field ext}.
\end{proof}

\noindent
We now strengthen the second part of Corollary~\ref{cor:d-max weakly d-closed}:

\begin{cor}\label{dmaxdom} Suppose $H$ is $\d$-maximal and $L\supseteq K$ is a differential subfield of $\Calinf[\imag]$  such that
$L$ is a $\d$-valued $H$-asymptotic extension of $K$ with respect to some dominance relation on $L$. 
Then $K$ is $\d$-closed in $L$. 
\end{cor}
\begin{proof}  The $\d$-valued field $K$ is $\upo$-free and newtonian by [ADH, 11.7.23, 14.5.7]. Also $L^\dagger\cap K=K^\dagger$ by Corollary~\ref{cor:fexphii}. Now apply Theorem~\ref{noextension}.
\end{proof} 

\noindent
We do not require that the dominance relation on $L$ in Corollary~\ref{dmaxdom} is the restriction to $L$ of the relation $\preceq$ on $\c[\imag]$.  

Recall also that in Section~\ref{sec:Hardy fields} we defined the $\d$-perfect hull $\operatorname{D}(H)$ of $H$   as the intersection of all $\d$-maximal Hardy field extensions of $H$.
By the next result we only need to consider  here  $\d$-algebraic Hardy field extensions of $H$:

\begin{cor}
If  $H$ is   $\d$-closed in some $\d$-maximal Hardy field extension of~$H$, then $H$ is $\d$-maximal.
Hence
$$\operatorname{D}(H)\ =\ \bigcap\big\{ M: \text{$M$ is a $\d$-maximal $\d$-algebraic Hardy field extension of $H$} \big\}.$$
\end{cor}
\begin{proof}
The first part follows from Theorem~\ref{thm:char d-max} and~(i)~$\Rightarrow$~(iii) in Proposition~\ref{prop:16.0.3 converse}. To prove the displayed equality we only  need to show the inclusion~``$\supseteq$''. So let $f$ be an element of every $\d$-maximal $\d$-algebraic Hardy field extension of $H$, and let $M$ be any $\d$-maximal Hardy field extension of $H$; we need to show $f\in M$. Let $E$ be the $\d$-closure of $H$ in $M$.
Then~$E$ is $\d$-algebraic over~$H$, and by the first part, $E$ is $\d$-maximal; thus $f\in E$, hence~$f\in M$ as required.
\end{proof}

\noindent
We can now prove a variant of Lemma~\ref{lem:Dx Ex} for $\Ginf$- and $\Gom$-Hardy fields:

\begin{cor}\label{cor:D(H) smooth}
Suppose $H$ is a $\Ginf$-Hardy field. Then
\begin{align*}
\operatorname{D}(H) &= \bigcap\big\{ M: \text{$M\supseteq H$ $\d$-maximal $\Ginf$-Hardy field} \big\} \\
&= \big\{ f\in \operatorname{E}^\infty(H): \text{$f$ is $\d$-algebraic over $H$}\big\}.
\end{align*}
Likewise with $\omega$ in place of $\infty$.
\end{cor}
\begin{proof}
With both equalities replaced by  ``$\subseteq$'', this follows from the definitions and the remarks
following Corollary~\ref{cor:Hardy field ext smooth}.
Let $f\in \operatorname{E}^\infty(H)$ be $\d$-algebraic over $H$; we claim that $f\in\operatorname{D}(H)$.
To prove this claim, let $E$ be a $\d$-maximal Hardy field extension $E$ of~$H$; it is enough to show that then $f\in E$.
Now $F:=E\cap\Ginf$ is a $\Ginf$-Hardy field extension of~$H$ which is $\d$-closed in $E$, by Corollary~\ref{cor:Hardy field ext smooth}, and hence $\d$-maximal by the previous corollary. Thus we may replace $E$ by $F$ to arrange that $E\subseteq\Ginf$, and then take a $\Ginf$-maximal Hardy field
extension $M$ of $E$. Now~$f\in \operatorname{E}^\infty(H)$ gives $f\in M$, and $E$ being $\d$-maximal and $f$ being $\d$-algebraic over $E$ yields $f\in E$. The proof for~$\omega$ in place of $\infty$ is similar.
\end{proof}

\noindent
Combining Theorem~\ref{thm:Bosh 14.4} with Corollary~\ref{cor:D(H) smooth} yields:

\begin{cor}\label{cor:Bosh Conj 1}
If $H\subseteq\Ginf$ is  bounded, then $\Dx(H)=\Ex(H)=\Ex^\infty(H)$. Likewise with $\omega$ in place of $\infty$.
\end{cor}

%\noindent
%If the question of Boshernitzan stated after Lemma~\ref{lem:Dx Ex} has a positive answer (that is, if  $\Ex(H)$ is always $\d$-algebraic over $H$), then the boundedness hypothesis in Corollary~\ref{cor:Bosh Conj 1} may be dropped.  \marginpar{not so clear anymore what the argument was}

\begin{question}
Do the following implications hold for all $H$?
$$H\subseteq\Ginf\ \Longrightarrow\ \Ex(H)\subseteq\Ex^\infty(H), \qquad H\subseteq\Gom\ \Longrightarrow\ \Ex(H)\subseteq\Ex^\infty(H)\subseteq\Ex^\omega(H).$$
\end{question}

\noindent
Let $\Ex:=\Ex(\Q)$ be the perfect hull of the Hardy field $\Q$. 
Boshernitzan~\cite[(20.1)]{Boshernitzan82} showed that $\Ex\subseteq\Ex^\infty(\Q)\subseteq\Ex^{\omega}(\Q)$.
From Corollary~\ref{cor:Bosh Conj 1} we obtain
$$\Ex\ =\ \Ex^\infty(\Q)\ =\ \Ex^{\omega}(\Q)\ =\ \Dx(\Q),$$
thus establishing \cite[\S{}10, Conjecture~1]{Boshernitzan81}.

Note that   $\Ex$ is $1$-$\d$-closed in all its Hardy field extensions, by Theorem~\ref{thm:Bosh order 1}.
However, $\Ex$ is not $2$-linearly surjective by \cite[Proposition~3.7]{Boshernitzan87}, %{\bf accepted on faith}   
so $\Ex$ is not weakly $2$-$\d$-closed in any  $\d$-maximal Hardy field extension of $\Ex$ (see Lemma~\ref{lem:weakly r-d-closed})
and   $\Ex$ is not $2$-linearly newtonian (see [ADH, 14.2.2]).

\medskip
\noindent
More generally, by Theorem~\ref{thm:Bosh order 1} each $\d$-perfect Hardy field is $1$-$\d$-closed in all its Hardy field extensions.
Together with Lemma~\ref{prop:dmax 1-newt} and~\ref{lem:r-newtonian descends}, this yields a generalization of 
Lemma~\ref{prop:dmax 1-newt}:

\begin{cor}\label{cor:d-perfect => 1-newt}
Every $\d$-perfect Hardy field  is $1$-newtonian.
\end{cor}

\noindent
Let $M$ be a $\d$-maximal Hardy field extension of $H$ and 
$H^{\operatorname{da}}$ the $\d$-closure of $H$ in~$M$, so $H(\R)\subseteq H^{\operatorname{da}}\subseteq M$. 
From Corollary~\ref{cor:nlc, 1} we  obtain a description of $H^{\operatorname{da}}$:

\begin{cor}
If $H$ is $\upo$-free, then $H^{\operatorname{da}}$ is a Newton-Liouville closure of~$H(\R)$.
\end{cor}

\noindent
Next, let $\boldsymbol H(\R)$, $\boldsymbol H^{\operatorname{da}}$, $\boldsymbol M$ be the canonical $\HLO$-expansions of the Hardy fields $H(\R)$, $H^{\operatorname{da}}$, $M$, respectively,
so $\boldsymbol H(\R)\subseteq\boldsymbol H^{\operatorname{da}}\subseteq\boldsymbol M$.     Corollary~\ref{cor:nlc, 2} yields:

\begin{cor}
$\boldsymbol H^{\operatorname{da}}$ is a Newton-Liouville closure of $\boldsymbol H(\R)$.
\end{cor}

\section{Embeddings into Transseries and Maximal Hardy Fields}\label{sec:embeddings into T}

\noindent
We begin with a direct consequence of facts about ``Newton-Liouville closure'' in~[ADH, 14.5, 16.2].  
Let $\boldsymbol H$ be a   $\HLO$-field with underlying $H$-field $H$. 
By [ADH, 14.5.10, 16.4.1, 16.4.8], the constant field of the Newton-Liouville closure of $\boldsymbol H$ is the 
real closure  of $C:= C_{H}$.
Let $\boldsymbol M$ be an $H$-closed $\HLO$-field extension of $\boldsymbol H$,  with underlying $H$-field $M$, and let $\boldsymbol H^{\operatorname{da}}$ be the $\d$-closure of $\boldsymbol H$ in $\boldsymbol M$. 

\begin{prop}\label{prop:embed into H-closed}  Let $\boldsymbol H^*$ be a $\d$-algebraic $\HLO$-field extension of $\boldsymbol H$ such that the constant field of $\boldsymbol H^*$ is algebraic over $C$. Then there is an embedding~$\boldsymbol H^*\to\boldsymbol M$
over $\boldsymbol H$, and the image of any such embedding is contained in~$\boldsymbol H^{\operatorname{da}}$.
\end{prop}
\begin{proof}
The image of any embedding $\boldsymbol H^*\to\boldsymbol M$
over $\boldsymbol H$ is $\d$-algebraic over $H$ and thus contained in $\boldsymbol H^{\operatorname{da}}$. 
For existence, take a Newton-Liouville closure $\boldsymbol M^*$ of $\boldsymbol H^*$. Then $\boldsymbol M^*$ is also a
Newton-Liouville closure of $\boldsymbol H$, by [ADH, 16.0.3],  and thus embeds into $\boldsymbol M$ over $\boldsymbol H$.
\end{proof}

\noindent
Let $\mathcal L$ be the language of ordered valued differential rings, as  in Section~\ref{sec:transfer}.
The second part of Corollary~\ref{cor:systems, 3} in the introduction now follows from the next result: 
 
\begin{cor}\label{cor:transfer T, 2}
Let $H$ be a Hardy field, $\iota\colon H\to\mathbb T$ an ordered differential field embedding, and
$H^*$ a $\d$-maximal $\d$-algebraic Hardy field extension of $H$. Then $\iota$ extends to an ordered valued differential field embedding $H^*\to\mathbb T$, and so for any  $\mathcal L_H$-sentence $\sigma$,  $H^*\models \sigma$ iff $\mathbb T\models \iota(\sigma)$.
\end{cor}
\begin{proof} We have $H(\R)\subseteq H^*$, and so by Lemma~\ref{lemhr} we arrange that $H\supseteq\R$. 
Let~$\boldsymbol H$,~$\boldsymbol H^*$ be the canonical $\HLO$-expansions of $H$, $H^*$, respectively, and let $\boldsymbol T$ be the expansion of~$\mathbb T$ to a $\HLO$-field.  Then $\boldsymbol H\subseteq\boldsymbol H^*$, and by Lemma~\ref{lem:unique HLO-expansion, 2}, $\iota$ is an
embedding $\boldsymbol H\to\boldsymbol T$.  
By Proposition~\ref{prop:embed into H-closed}, $\iota$ extends to an embedding $\boldsymbol H^*\to\boldsymbol T$. 
\end{proof}

\noindent
At the end of Section 5.5 we introduced the Hardy field $H:=\R(\ell_0, \ell_1, \ell_2,\dots)$, and we now mimick this
in $\mathbb T$ by setting $\ell_0:=x$ and $\ell_{n+1}:=\log\ell_n$ in $\mathbb T$. This yields the unique ordered differential
field embedding $H\to \T$ over $\R$ sending $\ell_n\in H$ to $\ell_n \in \T$ for all $n$. Its image is the $H$-subfield
$\R(\ell_0,\ell_1,\dots)$ of $\mathbb T$. 
Since the sequence $(\ell_n)$ in~$\T$ is coinitial in $\T^{>\R}$, each ordered differential subfield of $\mathbb T$
containing $\R(\ell_0,\ell_1,\dots)$ is an $\upo$-free $H$-field, by the remark preceding [ADH, 11.7.20]. 

 From Lemma~\ref{lem:unique HLO-expansion, 1} and Proposition~\ref{prop:embed into H-closed} we obtain:

\begin{cor}\label{cor:embed into H-closed, 1}
If $H\supseteq \R$ is an $\upo$-free $H$-subfield of $\mathbb T$ and $H^*$ is 
a $\d$-algebraic $H$-field extension of $H$ with constant field $\R$, then there exists an $H$-field embedding $H^*\to\mathbb T$ over $H$. 
\end{cor}

\noindent
Corollary~\ref{cor:embed into H-closed, 1} goes through
with $\mathbb T$ replaced by its  $H$-subfield
$$\mathbb T^{\operatorname{da}}\ :=\  \big\{ f\in\mathbb T : \text{$f$ is $\d$-algebraic (over $\Q$)}\big\},$$
a Newton-Liouville closure of  $\R(\ell_0,\ell_1,\dots)$; see  [ADH, 16.6] and Section~\ref{sec:diff closure} above. 
We now apply this observation to o-minimal structures. The {\it Pfaffian closure}\/ of an expansion of the ordered field of real numbers is its smallest expansion  that is closed under taking Rolle leaves of definable $1$-forms of class~$\c^1$. See Speissegger~\cite{Speissegger} for complete definitions, and the proof that 
%\marginpar{these definitions and result  taken on faith for now} 
the Pfaffian closure of an o-minimal expansion of  the ordered field of reals remains o-minimal.  

\begin{cor}\label{cor:embed into H-closed, 2}
The Hardy field $H$ of the Pfaffian closure of the ordered field of real numbers embeds as an $H$-field over $\R$  into $\mathbb T^{\operatorname{da}}$.
\end{cor}
\begin{proof}
Let $f\colon\R\to\R$ be definable in the Pfaffian closure of the ordered field of real numbers.
The proof of \cite[Theorem~3]{LMS} gives  $r\in\N$,  semialgebraic $g\colon\R^{r+2}\to\R$, and  $a\in\R$ 
%\marginpar{existence of such $r,g,a$ taken on faith} 
such that $f|_{(a,\infty)}$ is $\c^{r+1}$ and $f^{(r+1)}(t)=g\big(t,f(t),\dots,f^{(r)}(t)\big)$ for all~$t>a$. Take $P\in \R[Y_1,\dots,Y_{r+3}]^{\neq}$  vanishing identically on the graph of $g$; see~[ADH, B.12.18].
Then $P\big(t,f(t),\dots,f^{(r+1)}(t)\big)=0$ for $t>a$. Hence the germ of $f$ is $\d$-algebraic over $\R$, and so
 $H$ is $\d$-algebraic over $\R$. As $H$ contains the $\upo$-free Hardy field $\R(\ell_0, \ell_1,\dots)$,  we can use the remark following Corollary~\ref{cor:embed into H-closed, 1}.
 % and the remark following it.
\end{proof}

\begin{question}
Let $H$ be the Hardy field of an o-minimal expansion of the ordered field of reals, and 
let $H^*\supseteq H$ be the Hardy field of the Pfaffian closure of this expansion.
Does every embedding $H\to\mathbb T$ extend to an embedding $H^*\to\mathbb T$?
%(By \cite{LMS}, if $\Gamma$ is archimedean, then $H^*$ is {\it levelled}\/ in the sense of \cite{Rosenlicht87}.)
\end{question}

\noindent
We mentioned in the introduction that  an embedding $H\to\mathbb T$ as in Corollaries~\ref{cor:transfer T, 2} and~\ref{cor:embed into H-closed, 2}  can be viewed as an {\it expansion operator}\/ for the Hardy field $H$
and  its inverse as a {\it summation operator.}\/
The corollaries above concern the existence of expansion operators; this relied on the $H$-closedness of $\mathbb T$. Likewise, Theorem~\ref{thm:char d-max} and Proposition~\ref{prop:embed into H-closed}  also give rise to summation operators: 
\begin{cor}\label{cor:embed into H-closed, 3}
Let $H$ be an $\upo$-free $H$-field and let $H^*$ be
a $\d$-algebraic $H$-field extension of $H$ with $C_{H^*}$ algebraic over $C_H$.
Then any $H$-field embedding~$H\to M$ into a $\d$-maximal Hardy field  extends to an $H$-field embedding~$H^*\to M$.
\end{cor}

\noindent
In particular, given any ordered differential subfield $H\supseteq\R(\ell_0,\ell_1,\dots)$ of $\mathbb T$ with $\d$-closure
$H^*$ in $\mathbb T$, 
any $\mathcal{L}$-isomorphism between $H$ and a Hardy field $F$ extends to an $\mathcal{L}$-isomorphism between $H^*$ and a Hardy field extension of $F$. For $H=\R(\ell_0,\ell_1,\dots)\subseteq \mathbb T$ (so $H^*=\mathbb T^{\operatorname{da}}$) we recover the main result of~\cite{vdH:hfsol}:

\begin{cor} \label{cor:transserial}
The $H$-field $\mathbb T^{\operatorname{da}}$ is $\mathcal{L}$-isomorphic to a Hardy field $\supseteq \R(\ell_0, \ell_1,\dots)$. 
\end{cor} 

\noindent
Any Hardy field that is $\mathcal{L}$-isomorphic to~$\mathbb T^{\operatorname{da}}$ is $\d$-maximal, so contains the Hardy field $\operatorname{E}=\operatorname{E}(\Q)=\operatorname{D}(\Q)$; see the remarks following Lemma~\ref{lem:Dx Ex}.
Thus we have an $\mathcal{L}$-embedding $e\colon \operatorname{E}\to \mathbb T^{\operatorname{da}}$, which we can view as an expansion operator for the Hardy field $\operatorname{E}$.
We suspect that  $e(\operatorname{E})$ is independent of the choice of $e$.

\medskip
\noindent
In the remainder of this section we draw some consequences of Corollary~\ref{cor:transserial} for the universal theory of Hardy fields.

\subsection*{The universal theory of Hardy fields}
Recall from Section~\ref{sec:transfer} that $\mathcal L=\{0,1,{-},{+},{\,\cdot\,},\der,{\leq},{\preceq}\}$ is the language of
ordered valued differential rings.  Let
$\mathcal L^\iota$ be $\mathcal L$ augmented by a new unary function symbol~$\iota$.
We view each pre-$H$-field~$H$ as an $\mathcal L^\iota$-structure by interpreting the symbols
from $\mathcal L$ in the natural way and
$\iota$ by the function $\iota\colon H\to H$ given by $\iota(a):=a^{-1}$ for $a\in H^\times$ and $\iota(0):=0$. 

Since every Hardy field extends to a maximal one, each universal 
 $\mathcal L^\iota$-sentence which holds in every maximal Hardy field also  holds in every   Hardy field; likewise with ``$\d$-maximal'', ``perfect'', or ``$\d$-perfect''
 in place of ``maximal''.
We now use  Corollary~\ref{cor:transserial} to show:

\begin{prop}\label{prop:univ theory}
Let $\Sigma$ be the set of universal $\mathcal L^\iota$-sentences true in all   Hardy fields. Then the models of $\Sigma$ are the  pre-$H$-fields with very small derivation.
\end{prop}

\noindent
For this we need  a refinement of [ADH, 14.5.11]:

\begin{lemma}\label{lem:very small der extend}
Let $H$ be a pre-$H$-field with very small derivation. Then $H$ extends to an $H$-closed field with small derivation.
\end{lemma}
\begin{proof}
By Corollary~\ref{cor:dv(K) very small der}, replacing $H$ by its $H$-field hull, we first arrange that $H$ is an $H$-field.
Let   $(\Gamma,\psi)$ be the asymptotic couple of $H$. Then 
$\Psi^{\geq 0}\neq\emptyset$ or~$(\Gamma,\psi)$ has gap~$0$.
Suppose $(\Gamma,\psi)$ has gap $0$. Let $H(y)$ be the $H$-field extension from [ADH, 10.5.11] for $K:=H$, $s:=1$.
Then $y\succ 1$ and $y^\dagger=1/y\prec 1$, so replacing $H(y)$ by $H$
we can arrange that $\Psi^{\geq 0}\neq\emptyset$.
Then
every pre-$H$-field extension of $H$ has small derivation, and so we are done by [ADH, 14.5.11].
\end{proof}

\begin{proof}[Proof of Proposition~\ref{prop:univ theory}]
The natural axioms for pre-$H$-fields with very small derivation formulated in $\mathcal L^\iota$ are universal, so all models of $\Sigma$ are
 pre-$H$-fields with very small derivation. Conversely, given any pre-$H$-field $H$  with very small derivation we show that $H$ is a model of $\Sigma$: use Lemma~\ref{lem:very small der extend}  to extend $H$ to an $H$-closed field
 with small derivation, and note that the $\mathcal L^\iota$-theory of $H$-closed fields with small derivation is complete by [ADH, 16.6.3] and has a Hardy field model by Corollary~\ref{cor:transserial}.
\end{proof}

\noindent
Similar arguments 
allow us to settle a conjecture from \cite{AvdD2}, in slightly strengthened form.
For this, let $\mathcal L_x^\iota$ be $\mathcal L^\iota$ augmented by a constant symbol $x$.
We view each Hardy field containing the germ of the identity function on $\R$ as an
$\mathcal L_x^\iota$-structure by interpreting the symbols from~$\mathcal L^\iota$ as described
at the beginning of this subsection and the symbol $x$ by the germ  of the identity function
on $\R$, which we also denote by $x$ as usual. Each universal $\mathcal L_x^\iota$-sentence which holds in every maximal Hardy field also holds in every  Hardy field containing $x$.

\begin{prop}\label{prop:univ theory x}
Let $\Sigma_x$ be the set of universal $\mathcal L_x^\iota$-sentences true in all  
Hardy fields that contain $x$. Then the models of $\Sigma_x$ are the pre-$H$-fields with
distinguished element $x$  satisfying $x'=1$ and $x\succ 1$.
\end{prop}

\noindent
This follows from [ADH, 14.5.11] and the next lemma just like Proposition~\ref{prop:univ theory} followed from
Lemma~\ref{lem:very small der extend} and [ADH, 16.6.3].

\begin{lemma}\label{lem:complete x}
The  $\mathcal L_x^\iota$-theory of $H$-closed fields  with distinguished element $x$
satisfying $x'=1$ and $x\succ 1$ is complete.
\end{lemma}
\begin{proof}
Let $K_1$, $K_2$ be  models of this theory, and let $x_1, x_2$ be the interpretations of $x$ in $K_1$, $K_2$.
Then [ADH, 10.2.2, 10.5.11] gives an isomorphism $\Q(x_1)\to\Q(x_2)$ of valued ordered differential fields sending $x_1$ to $x_2$.
To show that $K_1\equiv K_2$ as $\mathcal L_x^\iota$-structures we identify $\Q(x_1)$ with $\Q(x_2)$ via this isomorphism. 
View  $\HLO$-fields as $\mathcal L^\iota_{\HLO}$-structures where $\mathcal L^\iota_{\HLO}$ extends $\mathcal L^\iota$ as specified
in~[ADH, Chapter~16]. (See also the proof of Theorem~\ref{thm:transfer}.)
By [ADH, 16.3.19] the $\upo$-free $H$-fields $K_1$, $K_2$ uniquely expand to  $\HLO$-fields $\boldsymbol{K}_1$, $\boldsymbol{K}_2$. The $H$-subfield $\Q(x_1)$ of~$K_1$ is grounded, so expands also uniquely to an $\HLO$-field,
and this $\HLO$-field is a common substructure of both $\boldsymbol{K}_1$ and~$\boldsymbol{K}_2$. 
Hence~$\boldsymbol{K} _1\equiv_{\Q(x_1)} \boldsymbol{K}_2$ by [ADH, 16.0.1, B.11.6]. 
This yields the claim.
\end{proof}

\noindent
From the completeness of the $\mathcal L^\iota$-theory of $H$-closed fields with small derivation
and Lemma~\ref{lem:complete x} in combination with Theorem~\ref{thm:char d-max} we obtain:

\begin{cor}
The set $\Sigma$ of   universal $\mathcal L^\iota$-sentences true in all   Hardy fields is decidable, and so is
the set $\Sigma_x$ of universal $\mathcal L^\iota_x$-sentences true in all   Hardy fields
containing $x$.
\end{cor}

\noindent
We finish with an example of a not-so-obvious property of asymptotic integrals, expressible by universal $\mathcal L^\iota$-sentences, which holds in all Hardy fields. For this, let~$Y=(Y_0,\dots,Y_{n})$ be a
tuple of distinct indeterminates and $P,Q\in\Z[Y]^{\neq}$.

\begin{example}
For all  hardian germs~$\ell_0,\dots,\ell_{n+1},y$, and $\vec y:=(y,y',\dots,y^{(n)})$: 
$$\begin{cases}&
\parbox{32.5em}{if
$\ell_0'=1$, 
$\ell_{j+1}'\ell_j=\ell_j'$ for $j=0,\dots,n$, 
$P(\vec y)=0$, and~$q:=Q(\vec y)\neq 0$,  then 
$(\ell_0\cdots\ell_{n+1}q)'\neq 0$, $(\ell_0\cdots\ell_{n+1}q)^\dagger\not\asymp q$, and $\big(q/(\ell_0\cdots\ell_{n+1}q)^\dagger\big)'\asymp q$.}\end{cases}$$
To see this, let $\ell_0,\dots,\ell_{n+1},y$
be as hypothesized.  Induction on $j$ shows ${\ell_j\in\Li(\R)}$ and $\ell_j\asymp\log_j x$ for $j=0,\dots,n+1$.
Put $H:=\R(x)$ and $E:=H\langle y\rangle$. Then ${\operatorname{trdeg}(E|H) \leq n}$, so  
Theorem~\ref{thm:Ros83} and Lemma~\ref{lemro83}
yield an~$r \in \{0,\dots,n\}$ and~$g\in E^>$ with~$g\asymp\ell_r$ such that~$E$ is grounded with~$\max\Psi_E=v(g^\dagger)$.
Iterating Proposition~\ref{prop:Hardy field exts} and [ADH, 10.2.3 and remark after it], starting  with $g^\dagger$ and
$\log g$ in the role of $s$ and $y$ in [ADH, 10.2.3],  produces
a grounded Hardy field  $F$ with~$E\subseteq F\subseteq\Li(E)$  
and $\max\Psi_F=v(f^\dagger)$ where~$f\in F^\times$, $f\asymp\ell_{n+1}$.
Then
$$\Psi_E\ <\ v(f^\dagger)\ <\ (\Gamma_E^>)',$$ 
so $f^\dagger\nasymp q:=Q(\vec y)$, and thus 
$(f^\dagger/q)^\dagger\asymp (\ell_0\cdots \ell_{n+1}q)^\dagger$. 
Now the conclusion follows from \cite[remarks after Lemma~2.7]{AvdD2} applied to
$q$, $f^\dagger$, $F$ in the role of~$a$,~$b_0$,~$K$.
\end{example}

\section{Linear Differential Equations over Hardy Fields}\label{sec:lin diff applications}

\noindent
In this section we draw some consequences of our main Theorem~\ref{thm:char d-max} for linear differential equations over Hardy
fields. This  also uses results from Section~\ref{sec:ueeh}.
Throughout this section $H$ is a Hardy field and $K:=H[\imag]\subseteq \Calinf[\imag]$.
Recall from Corollary~\ref{cor:d-max weakly d-closed} that if $H$ is $\d$-maximal, then~$K$ is linearly surjective and linearly closed; we use this fact freely below.
Let $A\in K[\der]^{\neq}$ be monic and $r:=\order A$.

\subsection*{Solutions in the complexification of a $\d$-maximal Hardy field}

\begin{theorem}\label{thm:lindiff d-max}
Suppose $H$ is $\d$-maximal. Then $A$ splits over $K$ and the $\C$-linear space  of zeros of $A$ in
$\Calinf[\imag]$ has a basis 
$$f_1\ex^{\phi_1\imag},\ \dots,\ f_r\ex^{\phi_r\imag}\qquad\text{where $f_1,\dots, f_r\in K^\times$, $\phi_1,\dots, \phi_r\in H$.}$$
For any such basis, set $\alpha_j:=\phi_j'\imag+K^\dagger\in K/K^\dagger$ for $j=1,\dots,r$. Then the spectrum of $A$ is~$\{\alpha_1,\dots,\alpha_r\}$,
with 
$$\operatorname{mult}_{\alpha}(A)\ =\ \abs{ \{ j\in\{1,\dots,r\}: \alpha_j=\alpha \} }\quad\text{ for every $\alpha\in K/K^\dagger$,}$$
and for any $a_1,\dots, a_r\in K$ with $A\ =\ (\der-a_r)\cdots (\der-a_1)$ we have
$$\operatorname{mult}_{\alpha}(A)\ =\ \abs{ \{ j\in\{1,\dots,r\}:  a_j+K^\dagger=\alpha \} }\quad\text{ for every $\alpha\in K/K^\dagger$.}$$
%$$A\ =\ (\der-a_r)\cdots (\der-a_1)\qquad\text{where $\alpha_j=a_j+K^\dagger$ for $j=1,\dots,r$.} $$
\end{theorem}

\noindent
(The spectrum of $A$ is as defined in Section~\ref{sec:splitting}, and does not refer to eigenvalues of
the $\C$-linear operator $y\mapsto A(y)$ on~$\Calinf[\imag]$.)

\begin{proof}
The first part follows from Corollaries~\ref{corbasiseigenvalues},~\ref{cor:complex basis}, and 
Lemma~\ref{lem:full kernel}.
For the rest, also use Lemma~\ref{lem:complex basis} and the proof of Corollary~\ref{cor:complex basis}.
% and the remarks before Lemma~\ref{newlembasis}. 
\end{proof}

\begin{remarks}
Suppose $H$ is $\d$-maximal; so $\I(K)\subseteq K^\dagger$ by Corollary~\ref{cor:cos sin infinitesimal}. Hence by 
Corollary~\ref{cor:complex basis} we can choose
the germs~$f_j$,~$\phi_j$ \textup{(}$j=1,\dots,r$\textup{)} in Theorem~\ref{thm:lindiff d-max} 
such that additionally $f_1\ex^{\phi_1\imag},\dots, f_r\ex^{\phi_r\imag}$ is a Hahn basis of $\ker_{\Calinf[\imag]}A$;
%, for~$j,k=1,\dots,r$: $\phi_j=0$ or~$\phi_j\succ 1$, and $\phi_j=\phi_k$ or~${\phi_j-\phi_k\succ 1}$. 
then the~$f_j$ with $\phi_j=0$ form a valuation basis of the valued $\C$-linear space $\ker_K A$. Fix such~$f_j$,~$\phi_j$, and let 
$\langle\ ,\,\rangle$ be the ``positive definite hermitian form''  on the  $K$-linear subspace~$K[\ex^{H\imag}]$
of $\Calinf[\imag]$, as specified in the remarks after Corollary~\ref{cor:phi_i paired, strengthened}, with associated ``norm''~$\dabs{\,\cdot\,}$ on $K[\ex^{H\imag}]$ given by~$\dabs{f}:=\sqrt{\<f,f\>}\in H^{\ge}$. Those remarks give
$$\langle f_j\ex^{\phi_j\imag}, f_k\ex^{\phi_k\imag}\rangle\ =\ \begin{cases} 0 & \text{if $\phi_j\neq\phi_k$,} \\ 
f_j\overline{f_k} & \text{if $\phi_j=\phi_k$,}\end{cases}$$
and so
$\dabs{f_j\ex^{\phi_j\imag}}=\abs{f_j}$. 

Next, let $H_0\supseteq\R$ be a Liouville closed Hardy subfield of $H$, set $K_0:=H_0[\imag]$ and suppose $\I(K_0)\subseteq K_0^\dagger$, $A \in K_0[\der]$, and $A$ splits over $K_0$. Then we can can choose the~$\phi_j$,~$f_j$ in Theorem~\ref{thm:lindiff d-max} 
such that  $f_1\ex^{\phi_1\imag},\dots, f_r\ex^{\phi_r\imag}$ is a Hahn basis of $\ker_{\Calinf[\imag]}A$, 
$\phi_1,\dots,\phi_r\in H_0$, and $vf_1,\dots,vf_r\in v(H_0^\times)$, by Corollaries~\ref{cor:ultimate prod, 2}  and~\ref{cor:A ultimate}.
\end{remarks}

\noindent
For each $\phi\in H$, the $\C$-linear operator $y\mapsto A(y)$ on $\Calinf[\imag]$ maps
the   $\C$-linear subspace $K\ex^{\phi\imag}$ of $\Calinf[\imag]$ into itself (Lemma~\ref{lem:trig surj, complex}); more precisely, by Corollary~\ref{cor:trig surj, complex}:

\begin{cor}\label{cor:lindiff d-max}
Suppose $H$ is $\d$-maximal, and let $\phi\in H$.
Then $A(K\ex^{\phi\imag}) = K\ex^{\phi\imag}$.
Moreover, if $\phi'\imag+K^\dagger\in K/K^\dagger$ is not an eigenvalue of $A$, then
   for each~$b\in K$ there is a unique $y\in K$ with
$A(y\ex^{\phi\imag})=b\ex^{\phi\imag}$.
\end{cor}

\noindent
Can the assumption ``$H$ is $\d$-maximal'' in Theorem~\ref{thm:lindiff d-max} and Corollary~\ref{cor:lindiff d-max} be weakened to ``$H$ is perfect''? The case $H=\Ex:=\Ex(\Q)$ is illuminating:  $\Ex\supseteq\R$ is a Liouville closed $H$-field, so contains the germs $x$ and $\ex^{x^2}$, but Boshernitzan~\cite[Pro\-po\-sition~3.7]{Boshernitzan87} showed that $\Ex$ is not $2$-linearly surjective, as there is no $y\in \Ex$ with~$y''+y=\ex^{x^2}$. 
%\marginpar{result from Boshernitzan accepted on faith for now} 
In fact, the conclusion of Corollary~\ref{cor:lindiff d-max} fails for  $H=\Ex$:

\begin{lemma}
Suppose $H=\Ex$. Then $K$ is not $1$-linearly surjective: there is no~$y\in K$ with $y'-y\imag = \ex^{x^2}$.
\end{lemma}
\begin{proof}
Suppose $y=a+b\imag$ ($a,b\in H$) satisfies $y'-y\imag = \ex^{x^2}$.
Now $$y'-y\imag = (a'+b'\imag) - (-b+a\imag) = (a'+b)+(b'-a)\imag,$$
hence  $a'+b=\ex^{x^2}$ and $b'=a$, so~$b''+b=\ex^{x^2}$, contradicting~\cite[Proposition~3.7]{Boshernitzan87}.
\end{proof}

\noindent
It follows that the conclusion of Theorem~\ref{thm:lindiff d-max} fails for~$H=\Ex$:

\begin{cor}\label{cor:no full basis}
Let $H=\Ex$ and $A=(\der-2x)(\der-\imag)$. Then 
%with~$\Univ:=K[\ex^{H\imag}]$, we have 
$\ker_{K[\ex^{H\imag}]} A = \mathbb C\ex^{x\imag}$.
\end{cor}
\begin{proof}
In Section~\ref{sec:ueeh} we identified the universal exponential extension of $K$
with $K[\ex^{H\imag}]$.
%By Lemma~\ref{newlembasis},   $\ker_{K[\ex^{H\imag}]} A$ has a basis contained in $\Univ^\times=K^\times\ex^{H\imag}$,
We have~$\ex^{x\imag}\in \ker_{K[\ex^{H\imag}]} A$.
Suppose $\dim_{\mathbb C}\ker_{K[\ex^{H\imag}]} A=2$. Then by Corollary~\ref{corbasiseigenvalues}, the  eigenvalues of $A$
are $2x+K^\dagger$ and $\imag+K^\dagger$. Now $2x\in K^\dagger$, which gives $f\in K^\times$ with $A(f)=0$, so  $\ker_{K[\ex^{H\imag}]} A$ has basis
$f,\ex^{x\imag}$. Also $\imag\notin K^\dagger$ by a remark preceding Lemma~\ref{lem:W and I(F)}, so [ADH, 5.1.14(ii)] yields $(\der-\imag)(cf)=\ex^{x^2}$ for some $c\in \mathbb C^\times$, contradicting the lemma above.
\end{proof}

\noindent
A vestige of linear surjectivity is retained by $\d$-perfect Hardy fields:

\begin{cor}\label{cor:d-perfect almost lin surj}
Suppose $H\supseteq\R$ is Liouville closed and $\I(K)\subseteq K^\dagger$, and $A$ splits over $K$.
Then there are $\fm,\fn\in H^\times$ such that for each Hardy field extension~$F$ of~$H$ and $b\in F[\imag]$ with $b\prec\fn$,
there exists $y\in\Dx(F)[\imag]$ that is the unique $y\in \Calinf[\imag]$ with~$A(y)=b$ and~$y\prec\fm$.
\textup{(}So if $H$ is $\d$-perfect, then for such $\fm$, $\fn$ and all $b\in K$ with $b\prec\fn$ there is a unique $y\in K$ with $A(y)=b$
and $y\prec\fm$.\textup{)}
\end{cor}
\begin{proof}
Let $E$ be a $\d$-maximal Hardy field extension of $H$. Theorem~\ref{thm:lindiff d-max} and the remark following it yields a Hahn basis
$$f_1\ex^{\phi_1\imag},\ \dots,\ f_r\ex^{\phi_r\imag}\qquad (f_1,\dots, f_r\in E[\imag]^\times,\ \phi_1,\dots, \phi_r\in E)$$
of $\ker_{\Calinf[\imag]}A$ with $\phi_1,\dots,\phi_r\in H$ and $v f_1,\dots,v f_r\in \Gamma = v(K^\times)$. 
It follows from Corollaries~\ref{cor:ultimate prod, 2} and~\ref{cor:A ultimate} that 
$\exc^{\ev}(A)=\exc^{\ev}_{E[\imag]}(A)=\{vf_j:\  j=1,\dots,r,\ \phi_j=0\}$.
By Corollary~\ref{cor:13.7.10} the quantity $v_A^{\ev}(\gamma)$, for $\gamma\in \Gamma$, does not change when passing from~$K$ to any ungrounded $H$-asymptotic extension of $K$.

Take $\fm,\fn\in H^\times$ with $\fm\prec f_1,\dots,f_r$ and 
$v\fn=v_A^{\ev}(v\fm)$. 
Consider  a Hardy field extension $F$ of $H$. Let $b\in F[\imag]^\times$, $b\prec\fn$, and let $M$ be a $\d$-maximal Hardy field extension of $F$. Then linear newtonianity of $L:=M[\imag]$ and Corollary~\ref{cor:14.2.10, generalized} yields~$y\in L$ with $A(y)=b$, $vy\notin\exc^{\ev}_L(A)=\exc^{\ev}(A)$, and $v^{\ev}_A(vy)=vb$. Then~$v^{\ev}_A(vy)=vb>v\fn=v^{\ev}_A(v\fm)$. Since $v\fm>\exc^{\ev}_L(A)$, this yields $y\prec\fm$ by Lemma~\ref{lem:ADH 14.2.7}.  Suppose~$z\in\Calinf[\imag]$,
$A(z)=b$, and $y\ne z\prec\fm$. Then~$u:=y-z\in\ker_{\Calinf[\imag]}A$ and~$0\ne u\prec\fm$, so $f_j\prec\fm$ for some~$j$ by Corollary~\ref{cor:from phi to lambda} (applied to $E$, $E[\imag]$ in place of $H$, $K$), a contradiction. 
This last argument also takes care of the case~$b=0$: there is no nonzero $u\prec \fm$ in $\Calinf[\imag]$ such that $A(u)=0$. 
\end{proof}

%\marginpar{skipped this paragraph}
\noindent
In \cite{ADHld} we shall prove that if $H$ is $\upo$-free and $\d$-perfect, then $K$ is linearly closed.
(This applies to $H=\Ex$.)
In particular, if  the $\d$-perfect hull $\Dx(H)$ of $H$ is $\upo$-free, then $A$ splits over
the algebraic closure   $\Dx(H)[\imag]$ of $\Dx(H)$.
(In Section~\ref{sec:perfect applications} below we characterize when $\Dx(H)$  is $\upo$-free.)
Now if $A$ splits over $\Dx(H)[\imag]$, then there are~$g,\phi\in\Dx(H)$, $g\neq 0$, such that~$A(g\ex^{\phi\imag})=0$.
The next lemma helps to clarify when for $H\supseteq\R$ we may take here~$g$,~$\phi$ in the Hardy subfield $\operatorname{Li}(H)$ of~$\Dx(H)$.

\begin{lemma}\label{lem:Liouvillian zeros}
Suppose $H\supseteq\R$. The following are equivalent:
\begin{enumerate}
\item[\textup{(i)}]  there exists $y\neq 0$ in a Liouville extension of the differential field $K$ such that~$A(y)=0$;
\item[\textup{(ii)}]  there 
exists $f\in\Li(H)[\imag]$  such that $f'$ is algebraic over $K$ and $A(\ex^{f})=0$;
%are $g,\phi\in\Li(H)$, $g\neq 0$, such that $g^\dagger$, $\phi'$ are algebraic over $H$ and~$A(g\ex^{\phi\imag})=0$;
\item[\textup{(iii)}]  there exists $f\in\Li(H)[\imag]$  such that $A(\ex^{f})=0$.
\end{enumerate}
\end{lemma}
\begin{proof}
Suppose (i) holds. Then~Corollary~\ref{cor:Kolchin}
%\cite[Proposition~1.45]{vdPS} \marginpar{\bf accepted on faith for now} 
gives $y\neq 0$ in a differential field extension~$L$ of $K$   such that $A(y)=0$ and~$g:=y^\dagger$ is algebraic over~$K$. % of degree~$\leq I(r)$. 
We arrange that $L$ contains the algebraic closure $K^{\operatorname{a}}=H^{\operatorname{rc}}[\imag]$ of $K$, where $H^{\operatorname{rc}}\subseteq\Calinf$ is the real closure of the Hardy field $H$.
Thus~$g\in K^{\operatorname{a}}$, and
hence $A=B(\der-g)$ where $B\in K^{\operatorname{a}}[\der]$ by [ADH, 5.1.21].
%and $$\big[H(a,b):H\big] = \big[K(a,b):K\big]=\big[K(z,\overline{z}):K\big]\leq \big[K(z):K\big] \cdot \big[K(\overline{z}):K\big] \leq N,$$
%where for the first equality we used that $H(a,b)$, $K$ are linearly disjoint over $H$. 
Take $f\in\Li(H)[\imag]$ with~$f'=g$ and set~$z:=\ex^f\in\Calinf[\imag]^\times$. Then~$z^\dagger=g$ and thus
$A(z)=0$.
This shows~(i)~$\Rightarrow$~(ii), and
(ii)~$\Rightarrow$~(iii) is trivial. To prove (iii)~$\Rightarrow$~(i),
let $f$ be as in~(iii) and~$y:=\ex^{f}\in \Calinf[\imag]^\times$.
%It is routine to verify that since $\Li(H)$ is a Liouville extension of $H$,
By [ADH, 10.6.6] the differential field
$$L\ :=\operatorname{Li}(H)[\imag]\ \subseteq\  \Calinf[\imag]$$ is a Liouville extension of $K=H[\imag]$. 
Now $L[y]\subseteq \Univ_L:=L\big[\ex^{\Li(H)\imag}\big]  \subseteq   \Calinf[\imag]$ and~$y^\dagger=f'\in L$, so the differential fraction field $L(y)$ of $L$ is a Liouville extension of $L$, and hence of~$K$.
%containing the zero~$y:=g\ex^{\phi\imag}\neq 0$ of~$A$.
\end{proof}

%\begin{remark} 
%Let  $I(r)\in\N$ be as in \marginpar{remark skipped}
%\cite[Proposition~4.18]{vdPS}. Using \cite[Proposition~4.19]{vdPS} one can   show that if one of the equivalent conditions in
%Lemma~\ref{lem:Liouvillian zeros} holds, then
%the germ $f$  in (ii) can be taken so that in addition~$\big[K(f'):K\big]\leq I(r)$.
%\end{remark}

\begin{cor}\label{cor:Liouvillian zeros}
If  $H\supseteq\R$ is real closed and $A(y)=0$ for some $y\neq 0$ in a Liouville extension of $K$,
then $A(\ex^{f})=0$ for some $f \in \Li(H)[\imag]$ with~$f' \in K$. \end{cor}

\noindent
{\em We assume $H\supseteq \R$ in the next three results}. Let  $E$ be a $\d$-maximal Hardy field extension of $H$.
By Theorem~\ref{thm:lindiff d-max},  $A$ splits over~$E[\imag]$. When does $A$ already split over the $\d$-subfield $\Li(H)[\imag]$ of $E[\imag]$? Here is a necessary condition:

\begin{cor}
If $A$ splits over $\Li(H)[\imag]$, then it splits over~$H^{\operatorname{rc}}[\imag]$.
\end{cor}
\begin{proof}
We arrange $H=H^{\operatorname{rc}}$  and proceed by induction on $r$.
The case $r=0$ being trivial, suppose $r\geq 1$ and $A$ splits over~$L:=\Li(H)[\imag]$.
Then [ADH, 5.1.21] yields~$y\neq 0$ in a differential field extension of $L$ with constant field $\C$ and $y^\dagger\in L$
such that~$A(y)=0$. Now~$L\langle y\rangle$ is
 a Liouville extension of $L$ and hence of $K$,
 so Lemma~\ref{lem:Liouvillian zeros} gives $f\in L$ with~$f'\in K$ and $A(\ex^f)=0$.
 Then $A=B(\der-f')$ where~$B\in K[\der]$ by [ADH, 5.1.21], and  $B$ splits over $\Li(H)[\imag]$ by [ADH, 5.1.22].
We can assume inductively that $B$ splits over $K$, and then $A$ does too. 
 \end{proof}

\noindent
With a weaker hypothesis on $A$, we have:

\begin{lemma}\label{lem:alg basis}
Suppose $A(y)=0$ for some~${y\neq 0}$ in a Liouville extension of $K$.
Then there is a monic $B\in K[\der]$ of order $n\geq 1$ such that $A\in K[\der] B$ and
 the $\C$-linear space  of zeros of $B$ in~$\Calinf[\imag]$ has a basis 
$$g_1\ex^{\phi_1\imag},\ \dots,\ g_n\ex^{\phi_n\imag}\qquad\text{where $g_1,\dots, g_n\in \Li(H)^\times$, $\phi_1,\dots, \phi_n\in \Li(H)$}$$
and~$g_1^\dagger,\dots,g_n^\dagger,\phi_1',\dots,\phi_n'\in H^{\operatorname{rc}}$.
Any such $B$ splits over $H^{\operatorname{rc}}[\imag]$.
\end{lemma}

\begin{proof}
Put $L:=\Li(H)[\imag]$ and identify $\Univ_L$ with the differential subring $L\big[\ex^{\Li(H)\imag}\big]$ of~$\Calinf[\imag]$ as explained at the beginning of Section~\ref{sec:ueeh}.
We consider also the differential subfield $K^{\operatorname{a}}:=H^{\operatorname{rc}}[\imag]$ of $L$, and use Lemma~\ref{lem:Univ under d-field ext} to identify
$\Univ:=\Univ_{K^{\operatorname{a}}}$ with a differential subring of $\Univ_L$, so 
  for all $u\in \Calinf[\imag]^\times$ with $u^\dagger\in K^{\operatorname{a}}$ we have~$u\in \Univ^\times$.
%with the differential subring $K^{\operatorname{a}}[E]$ of $\Univ_L$,
%where $E:=\{u\in\Univ_L^\times:u^\dagger\in K^{\operatorname{a}}\}$. (Cf.~Lemma~\ref{lem:Univ under d-field ext}.) 
Corollary~\ref{cor:Liouvillian zeros} yields $f\in L$ such that $f'\in K^{\operatorname{a}}$
and $A(\ex^f)=0$. 
Then $g:=\ex^{\Re f}\in\Li(H)^\times$ and $\phi:=\Im f\in\Li(H)$
with~$\ex^f=g\ex^{\phi\imag}$, so~$g^\dagger=\Re f'$, $\phi'=\Im f'$, and thus $g^\dagger, \phi'\in H^{\operatorname{rc}}$. 
Now $(\ex^f)^\dagger=f'\in K^{\operatorname{a}}$, so $y:=\ex^f\in~\Univ^\times$. 
Let $V$ be the $\C$-linear subspace of $\Univ$ spanned by the $\sigma(y)$ with $\sigma\in \operatorname{Aut}_\der(\Univ|K)$. Then~$y\in V\subseteq\ker_{\Univ} A$ and so $n:=\dim_{\C} V\in\{1,\dots,r\}$.
Corollary~\ref{cor:newlembasis} yields a unique monic $B\in K^{\operatorname{a}}[\der]$ of order $n$ such that
$V=\ker_{\Univ} B$. 
From $\sigma(V)=V$ for all~$\sigma\in\operatorname{Aut}_\der(\Univ|K)$ we get~$B\in K[\der]$ by Corollary~\ref{cor:extend autom}.
 Then $A\in K[\der]B$ by~[ADH, 5.1.15(i), 5.1.11]. To show $V$ has a basis as described in the lemma, let
 $\sigma\in \operatorname{Aut}_\der(\Univ|K)$. Then $\sigma(y)\in \Univ^\times$, so 
 $\sigma(y)^\dagger\in K^{\operatorname{a}}=H^{\operatorname{rc}}+H^{\operatorname{rc}}[\imag]$, hence~$\sigma(y)^\dagger=g_{\sigma}^\dagger+ \phi_{\sigma}'\imag$ with $g_{\sigma}, \phi_{\sigma}\in H^{\operatorname{rc}}$,
 $g_{\sigma}\ne 0$. 
 Also $\big(g_{\sigma}\ex^{\phi_{\sigma}\imag}\big)^\dagger= g_{\sigma}^\dagger+ \phi_{\sigma}'\imag$, and thus
$ \sigma(y)=c_{\sigma}g_{\sigma}\ex^{\phi_{\sigma}\imag}$ with $c_{\sigma}\in \C^\times$. This yields a basis
of $V$ as claimed. The final splitting claim follows from Corollary~\ref{cor:newlembasis}. 
\end{proof}

\noindent
Lemma~\ref{lem:alg basis} yields the following corollary  inspired by \cite[Corollary~3]{Singer76}.

\begin{cor}\label{cor:alg basis}
If  $A$ is irreducible, then the following are equivalent:
\begin{enumerate}
\item[\textup{(i)}] $A(y)=0$ for some~${y\neq 0}$ in a Liouville extension of $K$;
\item[\textup{(ii)}] the $\C$-linear space  of zeros of $A$ in~$\Calinf[\imag]$ has a basis 
$$\hskip2em g_1\ex^{\phi_1\imag},\ \dots,\ g_r\ex^{\phi_r\imag}\qquad\text{where $g_1,\dots, g_r\in \Li(H)^\times$, $\phi_1,\dots, \phi_r\in \Li(H)$}$$
and~$g_1^\dagger,\dots,g_r^\dagger,\phi_1',\dots,\phi_r'$ are algebraic over $H$.
\end{enumerate}
\end{cor}

\noindent
Next we improve the bounds on the derivatives of solutions to linear differential equations from
Corollary~\ref{cor:EL} when the coefficients of the equation are in $K$:
 
\begin{cor}\label{cor:EL sharpened}
Let  $\fm\in H$  with $0<\fm\preceq 1$ and  $y\in\c^r[\imag]$ be such that~$A(y) = 0$ and~$y \preceq \fm^n$. 
Then $y\in\Calinf[\imag]$ and
$$y^{(j)}\ \preceq\ \fm^{n-j}\fv(A)^{-j}\qquad\text{ for $j=0,\dots,n$,}$$ 
with $\prec$ in place of $\preceq$ if~$y\prec \fm^n$.
\end{cor}
\begin{proof}
First arrange that $H$ is $\d$-maximal. Choose a complement
$\Lambda_H$ of the $\R$-linear subspace $\I(H)$ of~$H$, set $\Lambda:=\Lambda_H\imag$, and   identify the universal exponential extension
$\Univ=\Univ_K$ of $K$ with the differential subring $K[\ex^{H\imag}]$ of $\Calinf[\imag]$ as described at the beginning of Section~\ref{sec:ueeh}.
By Lemmas~\ref{lem:complex basis} and~\ref{lem:full kernel} we have $y\in\ker_{\Calinf[\imag]}A=\ker_{\Univ} A$ and
$$y\ =\ f_1\ex^{\phi_1\imag}+\cdots+f_m\ex^{\phi_m\imag},\quad f_1,\dots,f_m\in K,\quad \phi_1,\dots,\phi_m\in H,$$
where $\lambda_1:=\phi_1'\imag,\dots,\lambda_m:=\phi_m'\imag\in\Lambda$ are the  distinct eigenvalues of $A$ with respect to $\Lambda$
and $f_1\ex^{\phi_1\imag},\dots,f_m\ex^{\phi_m\imag}\in\ker_{\Univ}A$.
By Corollary~\ref{cor:gaussian ext dom} and Lemma~\ref{lem:gaussian ext dom, preceq},
we have for $j=1,\dots,m$: $f_j\ex^{\phi_j\imag}\preceq\fm^n$, with
$f_j\ex^{\phi_j\imag}\prec\fm^n$ if $y\prec\fm^n$.
Hence we may arrange that $y=f\ex^{\phi\imag}$ where $f\in K$, $\phi\in H$, and $\lambda:=\phi'\imag\in\Lambda$ is an eigenvalue of~$A$ with respect to $\Lambda$, so $\lambda\preceq\fv^{-1}$ by Corollary~\ref{cor:evs v-small}, where $\fv:=\fv(A)\preceq 1$.
 
Now for each $j\in\N$:  $(\ex^{\phi\imag})^{(j)}=R_j(\lambda)\ex^{\phi\imag}\preceq\fv^{-j}$, using Lem\-ma~\ref{Riccatipower+}
if~${\lambda\succeq 1}$,
and so by the Product Rule: if $g\in\c^j[\imag]^{\preceq}$,
then $(g\ex^{\phi\imag})^{(j)}  \preceq\fv^{-j}$,
and likewise with~$\prec$ in place of $\preceq$. 
If $\fm\asymp 1$, then this observation with $g:=f$ already yields the desired conclusion.
Suppose $\fm\prec 1$. Then with $z:=y\fm^{-n}=f\fm^{-n}\ex^{\phi\imag}$ this same observation with $g:=f\fm^{-n}$  gives for $j=0,\dots,n$: 
$z^{(j)}\preceq\fv^{-j}$, with 
 $z^{(j)}\prec \fv^{-j}$ if~$y\prec\fm^n$. Now $z\in \c^n[\imag]$, so we can  use Lemma~\ref{lem:bd mult conj} for $r=n$ and $\eta=\abs{\fv}^{-1}$.
\end{proof}

\noindent
For $\fm=1$ we obtain from Corollary~\ref{cor:EL sharpened}:

\begin{cor}\label{cor:EL sharpened, m=1}
Let   $y\in\c^r[\imag]$ be such that~$A(y) = 0$ and~$y \preceq 1$. 
Then $y\in\Calinf[\imag]$ and
$y^{(n)}  \preceq \fv(A)^{-n}$ for all $n$,
with $\prec$ in place of $\preceq$ if~$y\prec 1$.
\end{cor}

%\begin{example}
%Suppose $H\supseteq\R$ and $r=2$.  By Lemma~\ref{lem:2nd order, Liouvillian solutions, 1}, if there is some $y\neq 0$ in a Liouville extension of the differential field $K=H[\imag]$ with $A(y)=0$, then $\dim_{C_L}\ker_L A=2$ for some Liouville extension $L$ of $K$. In fact, if $g$, $\phi$ are as in Corollary~\ref{cor:Liouvillian zeros}(iii), then the  differential fraction field $L$ of $K[\ex^{\phi\imag}]$ has this property. \marginpar{check!}
%\end{example}

\noindent
Recall from $\eqref{eq:P_L}$ the {\it concomitant}\/ $P_A\in K\{Y,Z\}$ of $A$. It yields a $\C$-bilinear map
$$(y,z)\mapsto [y,z]_A:=P_A(y,z)\ :\  \Calinf[\imag]\times \Calinf[\imag] \to \Calinf[\imag]$$ used in the next result, which is immediate from Corollaries~\ref{cor:conjugate basis} and~\ref{cor:d-max weakly d-closed}.

\begin{cor}
Suppose $H$ is $\d$-maximal, and let $f_j$, $\phi_j$ be  as in Theorem~\ref{thm:lindiff d-max}.
Then the $\C$-linear space of zeros of the adjoint $A^*$ of $A$ in $\Calinf[\imag]$ has a basis
$$f_1^*\ex^{-\phi_1\imag},\ \dots,\ f_r^*\ex^{-\phi_r\imag}\qquad\text{where $f_j^*\in K^\times$  \textup{(}$j=1,\dots,r$\textup{)}}$$
such that $\big[f_j \ex^{\phi_j\imag}, f_k^* \ex^{-\phi_k\imag}\big]_A=\delta_{jk}$
for $j,k=1,\dots,r$.
\end{cor}
 
\noindent
Recall that $A$ is said to be {\it self-adjoint}\/ if $A^*=A$, and {\it skew-adjoint}\/ if $A^*=-A$. (See Definition~\ref{defsask}.)
Self-adjoint operators play an important role in boundary value problems; see, e.g.,~\cite[Chap\-ter~XIII]{DunfordSchwartz}. 
The next result follows from  Corollaries~\ref{cor:Darboux general},~\ref{cor:phi_i paired} 
and~\ref{cor:d-max weakly d-closed} and applies to such operators:

\begin{cor}\label{corsd}
Suppose $H$ is $\d$-maximal, and $A^*=(-1)^rA_{\ltimes a}$, $a\in K^\times$. Then there are $a_1,\dots,a_r\in K$ such that
$$A\ =\ (\der-a_r)\cdots(\der-a_1)\quad\text{and}\quad a_j+a_{r-j+1}\ =\ a^\dagger\ \text{ for $j=1,\dots,r$.}$$
For  $f_j$, $\phi_j$ as in Theorem~\ref{thm:lindiff d-max} we have: $\phi_1+\cdots+\phi_r\preceq 1$, and
for each $i\in \{1,\dots,r\}$ there is a $j\in \{1,\dots,r\}$ such that~${\phi_i+\phi_j}\preceq 1$.
\end{cor}

\noindent
%Here, the last statement may be seen as an analogue of the fact that the eigenvalues of a real skew-symmetric matrix are purely imaginary,   closed under the operation~$\lambda\mapsto-\lambda$, and    (counted with multiplicity) sum  up to zero.
Operators satisfying the hypothesis of Corollary~\ref{corsd} are {\it self-dual}\/ in the sense of Section~\ref{sec:self-adjoint}.
For sources of such operators in physics, see \cite{BHMW}.
The next result is immediate from Corollary~\ref{cor:self-dual operator} and Theorem~\ref{thm:lindiff d-max} and  gives a sufficient condition for such operators to have nontrivial zeros in complexified Hardy fields:

\begin{cor}\label{cor:nonosc zeros, self-dual op}
If $A$ is self-dual $($which is the case if $A$ is skew-adjoint$)$, $r$ is odd, and $L$ is a $\d$-maximal Hardy field extension of $H$, then there are  $y,z\in L$, not both zero, such that
$A(y+z\imag)=0$. 
\end{cor}

\noindent
 The space of zeros of a self-dual $A$ has a special kind of basis, by Corollary~\ref{cor:phi_i paired, strengthened}: 

\begin{cor}\label{cor:self-dual basis}
Suppose $A$ is self-dual and $H$ is $\d$-maximal. Then   the $\C$-linear space of zeros of $A$ in $\Calinf[\imag]$ has a basis
$$f_1\ex^{\phi_1\imag}, g_1\ex^{-\phi_1\imag},\, \dots,\, f_m\ex^{\phi_m\imag},g_m\ex^{-\phi_m\imag}, \ h_1,\dots,h_n \qquad (2m+n=r)$$
where $f_1,\dots,f_m,g_1,\dots,g_m,h_1,\dots,h_n\in K^\times$, and $\phi_1,\dots,\phi_m\in H^{>\R}$ are apart.
\end{cor}

\subsection*{Bounded operators}
{\it In this subsection  $H$ is $\d$-maximal. \textup{(One can often reduce to this situation by extending a given Hardy field to a $\d$-maximal Hardy field.)} We choose
an $\R$-linear complement~$\Lambda_H$ of~$\I(H)$ in $H$, set
 $\Lambda:=\Lambda_H\imag$,  
and identify~${\Univ:=\Univ_K}$ with $K[\ex^{H\imag}]$ as explained at the beginning of Section~\ref{sec:ueeh}. Also~$A\in\mathcal O[\der]$ \textup{(so~${\fv(A)=1}$)}.}\/ Thus~$\Univ^\times=K^\times\ex^{H\imag}$ and~$V:=\ker_{\Calinf[\imag]}A=\ker_{\Univ}A$. See
Sections~\ref{sec:second-order} and ~\ref{sec:ueeh} for  definitions of Lyapunov exponents and of $\c[\imag]^{\flattereq}$ and~$\Univ^{\flattereq}$. 

\begin {lemma}\label{vuflattereq} $V\subseteq \Univ^{\flattereq}$, and $\lambda(y)=\lambda(y,y',\dots,y^{(r-1)})\in\R$ for all $y\in V^{\neq}$.
\end{lemma}
\begin{proof} Lemma~\ref{lem:ev bded} gives $\Sigma(A)\subseteq [\mathcal{O}]$, and Corollary~\ref{cor:Lyap} yields
$V\subseteq \Univ\cap \c[\imag]^{\flattereq}$. Lemma~\ref{lem:complex basis} gives a basis
$f_1\ex(h_1\imag),\dots, f_r\ex(h_r\imag)$ of the $\C$-linear space $V$ with $f_1,\dots, f_r\in K^\times$ and $h_1,\dots, h_r\in \Lambda_H$, and it says that
then the eigenvalues of $A$ with respect to $\Lambda$ are $h_1\imag,\dots, h_r\imag$. So for $j=1,\dots,r$ we have
$h_j\imag-a\in K^\dagger=H+\I(H)\imag$ with $a\in \mathcal{O}$. Then $h_j-\Im a\in \I(H)\subseteq \mathcal{O}_H$ and
$\Im a\in \mathcal{O}_H$, so~$h_j\in \Lambda_H\cap \mathcal{O}_H$. From $f_j\ex(h_j\imag)\in \c[\imag]^{\flattereq}$
we obtain $f_j\in K\cap\c[\imag]^{\flattereq}=\mathcal{O}_{\Delta}$, so $V\subseteq \Univ^{\flattereq}$. The rest follows from
Corollary~\ref{cor:Lyap} and Lemma~\ref{lem:Univflattereq}.
\end{proof} 

\noindent
For $y\in\c^1[\imag]^\times$, in Section~\ref{sec:ueeh} we also defined
$$\mu(y)\ =\ \limsup_{t\to+\infty} \Im \frac{y'(t)}{y(t)}  \in \R_{\pm\infty}.$$
The zeros of the characteristic polynomial $\chi_A\in\C[Y]$ of~$A$ (defined in Section~\ref{sec:splitting}) contain information about   elements of   $V\cap\Univ^\times$:

\begin{lemma}\label{lem:lin diff applications, A bded}
Let $f\in K^\times$ and $\phi\in H$ be such that $y=f\ex^{\phi\imag}\in V$. Then~$y\in (\Univ^{\flattereq})^\times$, $\lambda:=\lambda(y)\in \R$, $\mu:=\mu(y)\in\R$, $\phi-\mu x\prec x$, and with
$\alpha:=\phi'\imag+K^\dagger$:
$$\chi_A(-\lambda+\mu\imag)=0,\qquad \operatorname{mult}_{\alpha}(A)\  \leq \sum_{c\in\C,\ \Im c=\mu} \operatorname{mult}_c(\chi_A).$$
\end{lemma}
\begin{proof}
Corollary~\ref{cor:ev bded}
gives $y^\dagger\preceq 1$, so $y\in (\Univ^{\flattereq})^\times$, $\lambda,\mu\in\R$ with $y^\dagger-(-\lambda+\mu\imag)\prec 1$ and $\phi'\preceq 1$ by
Lem\-ma~\ref{lem:lambda units} and an observation following Corollary~\ref{cor:Univflattereq}. Then $\phi'-\mu\prec 1$, so  $\phi-\mu x\prec x$.  The rest follows from Corollary~\ref{cor:ev bded} and Lemma~\ref{lem:mult upper bd}.
\end{proof}

\noindent
 A {\bf Lyapunov basis} of $V$ is a basis $y_1,\dots,y_r$ of the $\C$-linear space $V$   such that for all~$c_1,\dots,c_r\in\C$,
not all zero, and $y=c_1y_1+\cdots+c_ry_r$ we have $\lambda(y)={\min\!\big\{\lambda(y_j):c_j\neq 0\big\}}$.\index{Lyapunov!basis}\index{basis!Lyapunov}\index{linear differential operator!Lyapunov basis}
There is a Lyapunov basis of $V$; indeed, by the remarks after Theorem~\ref{thm:lindiff d-max} and  Corollary~\ref{cor:Lyap val indep}:

\begin{cor}\label{cor:lin diff applications, A bded}
The $\C$-linear space $V$ has  a Hahn basis
$$f_1\ex^{\phi_1\imag},\dots,f_r\ex^{\phi_r\imag}\qquad (f_1,\dots,f_r\in K^\times,\ \phi_1,\dots,\phi_r\in H),$$
and every such Hahn basis of $V$ is a Lyapunov basis of $V$.
\end{cor}

\begin{question}
By Perron~\cite[Satz~8]{Perron13b} (see also \cite[Satz~VI]{Spaeth}), $V$ has  a Lyapunov basis~$y_1,\dots,y_r$ such that  for each $\lambda\in\R$,  the number of $j$ with $\lambda(y_j)=\lambda$ is equal to~$\sum_{\mu\in\R} \operatorname{mult}_{-\lambda+\mu\imag}(\chi_A)$.
Can we   choose here $y_1,\dots,y_r$ to be a Hahn basis of~$V$?
\end{question}

\begin{cor}
If $\chi_A$ has no real zeros, then $K^\dagger$ is not an eigenvalue of $A$, and so there is no $y\in K^\times$ such that $A(y)=0$.
\end{cor}
\begin{proof}
Take a Hahn basis of $V$ as in Corollary~\ref{cor:lin diff applications, A bded}. 
Then the eigenvalues of $A$ are $\phi_1'\imag+K^\dagger,\dots,\phi_r'\imag+K^\dagger$, by Theorem~\ref{thm:lindiff d-max}.
Suppose $K^\dagger$ is an eigenvalue of $A$.  Then we have $j$ with $\phi_j'\imag\in K^\dagger=H+\I(H)\imag$, so
$\phi_j'\in I(H)$, hence $\phi_j\preceq 1$, and thus~${\phi_j=0}$. Then  Lemma~\ref{lem:lin diff applications, A bded} yields
 a real zero of $\chi_A$.
For the rest, use that the $f_j$ with $\phi_j=0$ form a basis of~$\ker_K A$ by remarks after Theorem~\ref{thm:lindiff d-max}.
\end{proof}

\begin{example}
The linear differential equation
$$y'''- \left(\imag + \frac{1}{\ex^x}\right) y'' + \left(1-\frac{1}{\log x}\right) y' - \left(\imag + \frac{1}{x^2}\right) y=0$$
has no nonzero  solution in $F[\imag]$ for any Hardy field $F$.
\end{example}

\noindent
We can now prove a strong version of a theorem of Perron~\cite[Satz~5]{Perron13a} in the setting of linear differential equations over complexified Hardy fields.
(A precursor of Perron's theorem for $A\in\C(x)[\der]$ is due to Poincar\'e~\cite{Poincare}.) Perron assumes additionally that the real parts of distinct zeros
of $\chi_A$ are distinct.

\begin{prop}\label{propperron}
Suppose all  \textup{(}complex\textup{)} zeros   of $\chi_A$ are   simple. Let 
$$y_1=f_1\ex^{\phi_1\imag},\ \dots,\ y_r=f_r\ex^{\phi_r\imag}\qquad (f_1,\dots,f_r\in K^\times,\ \phi_1,\dots,\phi_r\in H)$$
be a Hahn basis of $V$. Then the zeros of $\chi_A$
are $$c_1\ :=\ -\lambda(y_1)+\mu(y_1)\imag,\ \dots,\ c_r\ :=\ -\lambda(y_r)+\mu(y_r)\imag,$$
and $\big(y_j^{(n)}/y_j\big)- c_j^n\prec 1$ for $j=1,\dots,r$ and all $n$. 
\end{prop}
\begin{proof}
By Lemma~\ref{lem:lin diff applications, A bded} each $c_j$ is a zero of $\chi_A$, and we claim that there are no other. 
Let $c=-\lambda+\mu\imag$ ($\lambda,\mu\in\R$) be a zero of $\chi_A$. Then Corollary~\ref{cor:simple Riccati zero} and~[ADH, 5.1.21, 5.8.7] yield $A\in K[\der]\big(\der-(p+q\imag)\big)$ with $p,q\in \mathcal{O}_H$, $p+\lambda, q-\mu\prec 1$.
Taking~$f\in H^\times$ and $\phi\in H$ with $f^\dagger=p$ and $\phi'=q$ we have
 $y:=f\ex^{\phi\imag}\in V^{\ne}$  and so~$\lambda(y)=\lambda$ and $\mu(y)=\mu$.
 
 Take $a_1,\dots, a_r\in \C$ such that $y=a_1f_1\ex^{\phi_1\imag}+ \cdots + a_rf_r\ex^{\phi_r\imag}$.
As in the proof of Corollary~\ref{cor:Lyap val indep}  (but with $r$ instead of $m$) we arrange, with $l\in \{1,\dots, r\}$, that~$\phi_1,\dots, \phi_l$ are distinct and each $\phi_j$ with $l<j\le r$ equals one of $\phi_1,\dots, \phi_l$.
For $k=1,\dots,l$ we take (as in  that proof, but with other notation) $h_k\in \Lambda_H$ such that 
$\phi_k-\phi(h_k\imag)\preceq 1$ and put  $g_k:= \sum_{1\le j\le l,\, \phi_j=\phi_k} a_jf_j$ and 
$u_k:= \ex^{(\phi_k-\phi(h_k\imag))\imag}$, and likewise $h\in \Lambda_H$ with $\phi-\phi(h\imag)\preceq 1$, and set $u:=\ex^{(\phi-\phi(h\imag))\imag}$. Then
$$y\ = uf\ex(h\imag)\ =\ \ g_1\ex^{\phi_1\imag} + \cdots + g_l\ex^{\phi_l\imag}\ =\ u_1g_1\ex(h_1\imag)+\cdots + u_l\g_l\ex(h_l\imag).$$
Now $u, u_1,\dots, u_l\in K$, so $h=h_k$ for some $k\in \{1,\dots, l\}$, say $h=h_1$, hence~$uf=u_1g_1$ and so $y=g_1\ex^{\phi_1\imag}$.
Since the $f_j$ with $\phi_j=\phi_1$ are valuation-independent and~$f\asymp g_1$, this yields $j$ with $\phi_j=\phi_1$ and $a_j\ne 0$
such that $f\asymp f_j$. Then~$\lambda=\lambda(y)=\lambda(f)=\lambda(f_j)=\lambda(y_j)$.   
The proof of Lemma~\ref{lem:lin diff applications, A bded} gives
$$\mu\ =\ \mu(y)\ =\ \lim_{t\to \infty}\phi'(t),\qquad
\mu(y_j)\ =\ \lim_{t\to \infty}\phi_j'(t)\ =\ \lim_{t\to \infty}\phi_1'(t).$$ 
But $\phi_1-\phi(h_1\imag)\preceq 1$, 
$\phi-\phi(h\imag)\preceq 1$, and $h=h_1$, so $\phi-\phi_1\preceq 1$, hence $\phi'-\phi_1'\prec 1$, and
thus $\mu=\mu(y_j)$. This yields $c=c_j$. For the last claim, let $j\in \{1,\dots, r\}$. Then for
 $z_j:=y_j^\dagger$ we have $z_j-c_j\prec 1$ (see for example the proof of Lemma~\ref{lem:lin diff applications, A bded}),
 and~$y_j^{(n)}/y_j=R_n(z_j)$.  Now use Lemma~\ref{Riccatipower+} if $c_j\ne 0$. If
  $c_j=0$, then $z_j\prec 1$, so we can use that then $R_n(z_j)\prec 1$ for $n\ge 1$. 
\end{proof}

\begin{cor}\label{cor:Perron 1}
Suppose  the real part of each complex zero of $\chi_A$ is negative.
Then for all~$y\in V$ and all $n$ we have $y^{(n)}\prec 1$.
\end{cor}
\begin{proof}
By Corollary~\ref{cor:lin diff applications, A bded} it is enough to consider the case $y=f\ex^{\phi\imag}\in V$
where~$f\in K^\times$, $\phi\in H$. Then $\lambda:=\lambda(y)=\lambda(f)\in\R^>$ by
  Lemma~\ref{lem:lin diff applications, A bded},  which for~${0<\varepsilon<\lambda}$ gives   $f\prec\ex^{-(\lambda-\varepsilon)x}\prec 1$. 
Now use Corollary~\ref{cor:EL sharpened, m=1}.
\end{proof}

\noindent
We use Corollary ~\ref{cor:Perron 1} to strengthen another theorem of Perron~\cite{Perron20, Perron23} in the Hardy field context:

\begin{cor}\label{cor:Perron 2}
Suppose $a_0:=\chi_A(0)\neq 0$. Let
$b\in K$, $b\preceq 1$. Then there exists~$y\in K$ such that 
$$A(y)\ =\ b,\quad y-(b/a_0)\ \prec\ 1,\quad y^{(n)}\ \prec\ 1\ \text{ for all~${n\geq 1}$.}$$
Moreover, if the real part of each complex zero of~$\chi_A$ is negative, then  all $y\in\Calinf[\imag]$ with~$A(y)=b$ satisfy 
$y-(b/a_0)\prec 1$ and~$y^{(n)}\prec 1$ for all $n\geq 1$.
\end{cor}
\begin{proof}
By Theorem~\ref{thm:char d-max} and [ADH, 14.5.7], $K$ is $r$-linearly newtonian. As~${\der\mathcal O\subseteq\smallo}$,
 for the first part it is  enough 
to find $y\in K$ such that~$A(y)=b$ and~${y-(b/a_0)\prec 1}$.
  Corollary~\ref{cor:ADH 14.2.10 extract} yields such $y$ if $b\asymp 1$, so  suppose~$0\neq b\prec 1$  (since for~${b=0}$ we can take $y=0$). From~$A(1)\asymp 1$ and Proposition~\ref{kerexc} we obtain~$0\notin v(\ker^{\neq}_K A)=\exc^{\ev}(A)$.
Corollary~\ref{cor:14.2.10, generalized} then yields $y\in K^\times$ with $A(y)=b$, $vy\notin\exc^{\ev}(A)$, and $v_A^{\ev}(vy)=vb$.  Now $A(1)\asymp 1$ gives $v_A^{\ev}(0)=0$, so
$y\prec 1$ by  Lemma~\ref{lem:ADH 14.2.7}. The second statement now follows from Corollary~\ref{cor:Perron 1}.
\end{proof}

\begin{example}
Each $y\in\c^2[\imag]$ with $$y''+(2-\imag)(1+x^{-1}\log x)y'+(1-\imag)y\ =\ 2+\ex^{-x^2}$$ satisfies~$y\sim \imag+1$ and $y^{(n)}\prec 1$ for each $n\geq 1$, and there is
such a $y\in F[\imag]$ for some Hardy field $F$.
\end{example}

\noindent
Here is a version of Corollary~\ref{cor:Ri unique, 3} in the Hardy field setting:

\begin{prop}\label{perfAunique}
Let $H_0$ be a $\d$-perfect  Hardy subfield of $H$ such that $A\in K_0[\der]$ for $K_0:=H_0[\imag]$.
Suppose $r:=\order(A)\geq 1$ and $\chi_A$ has distinct zeros $c_1,\dots,c_r\in\C$ with
$\Re c_1 \geq \cdots \geq \Re c_r$. Then there is a unique splitting $(a_1,\dots,a_r)$ of $A$ over $K$
such that $a_1-c_1,\dots,a_r-c_r\prec 1$.
If in addition $\Re c_1 > \cdots > \Re c_r$,
then for this splitting  of $A$ over~$K$  we have~$a_1,\dots,a_r\in\mathcal O_{K_0}$.
\end{prop}
\begin{proof}
The first claim holds by  Corollary~\ref{cor:Ri unique, 3}.  
Suppose $\Re c_1 > \cdots > \Re c_r$; it remains   to show $a_1,\dots,a_r\in\mathcal O_{K_0}$. For this
we proceed by induction on~$r$. The case~$r=1$ being obvious, suppose $r > 1$.
Let $y_1,\dots,y_r$ be a Hahn basis of~$V$  with   $c_j=-\lambda(y_j)+\mu(y_j)\imag$ for~$j=1,\dots,r$
as in Proposition~\ref{propperron}.
Lemma~\ref{lem:fexphii} yields $\theta_j\in H$ with~$y_j=\abs{y_j}\ex^{\theta_j\imag}$
and~$\abs{y_j}\in H^>$, for~$j=1,\dots,r$; these $\theta_j$ might be different from the $\phi_j$ of Proposition~\ref{propperron}.
Let now $F$ be any $\d$-maximal Hardy field extension of~$H_0$; we claim that~$\abs{y_r},\theta_r\in F$.
To see this, use Lemma~\ref{lem:fexphii} and Proposition~\ref{propperron} applied to~$F$  in place of $H$ 
to get
$f\in F^>$, $\theta\in F$ such that~$y:=f\ex^{\theta\imag}$ satisfies~${A(y)=0}$ and $c_r-y^\dagger\prec 1$. 
Then~$y\in V\cap\c[\imag]^\times$. Take~$d_1,\dots,d_r\in\C$ such that~$y=d_1y_1+\cdots+d_ry_r$.
Lemma~\ref{smallestrealpart} applied to the $d_jy_j$ with~${d_j\neq 0}$ in place of~$f_1,\dots,f_n$, and with $c_r$ in the role of $c$, yields $i\in\{1,\dots,r\}$ such that~$d_i\neq 0$, $c_i=c_r$, and $\Re c_j\le \Re c_r$ for all $j$ with $d_j\ne 0$.
Hence $i=r$ is the only one such~$j$ and thus $y=d_ry_r$. 
%Then~$y^\dagger=y_i^\dagger$, so~${c_i-c_r=(c_i-y_j^\dagger)-(c_r-y^\dagger)\prec 1}$, hence $c_i=c_r$ and thus~$i=r$.
This yields~$\abs{y_r}\in \R^> f \subseteq F^>$ and~$\theta_r\in\theta+\R\subseteq F$ by the uniqueness part of Lemma~\ref{lem:arg}, as claimed.  This claim and $\d$-perfectness of~$H_0$
now give~$\abs{y_r},\phi_r\in H_0$, hence~$y_r^\dagger=\abs{y_r}^\dagger+\theta_r'\imag\in K_0$. 
By~[ADH, 5.1.21] we get~$A=B(\der-y_r^\dagger)$ where~$B\in K_0[\der]$ is monic, and by~[ADH, 5.6.3] we have~$B\in\mathcal O_{K_0}[\der]$ with~$\chi_A=\chi_B \cdot (Y-c_r)$, hence
the zeros of $B$ are~$c_1,\dots,c_{r-1}$. Now apply the inductive hypothesis to $B$.
\end{proof}

\begin{example}
The linear differential operator
$$\der^3 - (1-\ex^{-\ex^x})\imag\der^2 - (1+\imag+x^{-2}\log x^2) \der +   (\log \log x)^{-1/2}\in \mathcal O[\der]$$
splits over $\Ex(\Q)[\imag]$. In fact, there is a unique splitting $(a_1,a_2,a_3)$ of this
linear differential operator over $K$ with $a_1 -(1+\imag)\prec 1$, $a_2\prec 1$, and $a_3+1\prec 1$,  and we have~$a_1,a_2,a_3\in\Ex(\Q)[\imag]$.
\end{example}

\begin{question}
Can we drop the assumption $\Re c_1 > \cdots > \Re c_r$  in the last part of  Proposition~\ref{perfAunique}?
\end{question}

\noindent
Next we derive consequences of Theorem~\ref{thm:lindiff d-max}  for matrix differential equations.  

\subsection*{Matrix differential equations} {\em In this subsection $H$ is $\d$-maximal. We take
an $\R$-linear complement~$\Lambda_H$ of~$\I(H)$ in $H$, set
 $\Lambda:=\Lambda_H\imag$,  
and identify $\Univ=\Univ_K$ with $K[\ex^{H\imag}]$ as usual.
Let $N$ be an ${n\times n}$~matrix over~$K$, where~${n\geq 1}$.}
Recall from~[ADH, 5.5]  the definition of {\it fundamental matrix}\/ for  the matrix differential equation $y'=Ny$ over any differential ring extension of $K$.

\begin{cor}\label{cor:shape of fund matrix}
There are $M\in\operatorname{GL}_n(K)$  and $\phi_1,\dots, \phi_n\in H$  with
$\phi_1,\dots, \phi_n$ apart such that for $D :=\operatorname{diag}(\ex^{\phi_1\imag},\dots,\ex^{\phi_n\imag})$,
the $n\times n$ matrix $MD$ over $K[\ex^{H\imag}]$ is a fundamental matrix for $y'=Ny$.
Moreover, for any such $M$ and $\phi_1,\dots\,\phi_n$, 
%$M\in\operatorname{GL}_n(K)$ and~${\phi_1,\dots,\phi_n\in H}$ such that $MD$ is a fundamental matrix for $y'=Ny$ where $D:=\operatorname{diag}(\ex^{\phi_1\imag},\dots,\ex^{\phi_n\imag})$, 
%be given such that $M D$ is a fundamental matrix
%for $y'=Ny$.for any such $M$ and $\phi_1,\dots,\phi_n$,
 setting~$\alpha_j:=\phi_j'\imag+K^\dagger$ for $j=1,\dots,n$,  the spectrum of $y'=Ny$ is~$\{\alpha_1,\dots,\alpha_n\}$, and for all~$\alpha\in K/K^\dagger$,
 $$\operatorname{mult}_{\alpha}(N)\ =\ \abs{ \{ j\in\{1,\dots,n\}: \alpha_j=\alpha \} }.$$
\end{cor}
\begin{proof} The hypothesis of~[ADH, 5.5.14] is satisfied for $R:= \Univ=K[\ex^{H\imag}]$. To see why, let~$L\in K[\der]$ be monic of
order $n$. Then Theorem~\ref{thm:lindiff d-max} and a subsequent remark provide $f_1,\dots, f_n\in K^\times$ and~$\phi_1,\dots, \phi_n\in H$
such that $\phi_1,\dots, \phi_n$ are apart and $f_1\ex^{\phi_1\imag},\dots, f_n\ex^{\phi_n\imag}$ is a basis of $\ker_{\Calinf[\imag]} L=\ker_{\Omega} L$, where
$\Omega:=\Frac \Univ$. Hence $W:=\operatorname{Wr}\!\big(f_1\ex^{\phi_1\imag},\dots,f_n\ex^{\phi_n\imag}\big)\in \operatorname{GL}_n(\Omega)$. Also $\det W\in \ex^{\phi_1\imag + \cdots + \phi_n\imag}K^\times\subseteq \Univ^\times$, hence $W\in \operatorname{GL}_n(\Univ)$. Thus by the remarks preceding [ADH, 5.5.14]: $W$ is a fundamental matrix for $y'=A_L y$ with~$\Univ$ as the ambient differential ring.  
Note that $W=QD$ where $D=\operatorname{diag}(\ex^{\phi_1\imag},\dots,\ex^{\phi_n\imag})$ and $Q\in \operatorname{GL}_n(K)$.

We now follow the proof of [ADH, 5.5.14]: take monic $L\in K[\der]$ of order $n$ such that $y'=Ny$ is equivalent to $y'=A_Ly$, and
take $P\in \operatorname{GL}_n(K)$ such that~$P\operatorname{sol}_R(A_L)=\operatorname{sol}_R(N)$. 
With $W$  the above fundamental matrix for $y'=A_Ly$,  $PW\in  \operatorname{GL}_n\big(R)$ is then a fundamental matrix for $y'=Ny$. So $M:= PQ\in \operatorname{GL}_n(K)$ gives $MD=PW$ as a fundamental matrix for $y'=Ny$.

Let now any $M\in\operatorname{GL}_n(K)$, $\phi_1,\dots,\phi_n\in H$, and $D:=\operatorname{diag}(\ex^{\phi_1\imag},\dots,\ex^{\phi_n\imag})$ be given such that $M D$ is a fundamental matrix
for $y'=Ny$. Let $f_1,\dots,f_n$ be the successive columns of $M$. Then   $\ex^{\phi_1\imag}f_1,\dots, \ex^{\phi_n\imag}f_n$ 
is a basis of the $\C$-linear space~$\operatorname{sol}_{\Univ}(N)$. The ``moreover’’ part now follows from Lemma~\ref{matrixversionlem}.
%Then~$W:=P^{-1}MD\in\operatorname{GL}_n(R)$ is a fundamental matrix for $y'=A_Ly$. Let $(f_1,\dots, f_n)\in K^n$
%be the top row of   $P^{-1}M$. Then $(f_1\ex^{\phi_1\imag},\dots,f_n\ex^{\phi_n\imag})$ is the top row of $W$,
%so by the remark preceding [ADH, 5.5.14] we have 
%$W=\operatorname{Wr}(f_1\ex^{\phi_1\imag},\dots,f_n\ex^{\phi_n\imag})$, $f_1\ex^{\phi_1\imag},\dots,f_n\ex^{\phi_n\imag}\in\ker_R L$, and  $\det W\in R^\times$.
%Hence $f_1\ex^{\phi_1\imag},\dots,f_n\ex^{\phi_n\imag}$ is a basis of the $\C$-linear space $\ker_R L=\ker_{\Calinf[\imag]}L$. 
%The ``moreover'' part now follows from this observation in combination with 
%the fact stated just before Lemma~\ref{lem:der-a onto}, Lemma~\ref{lem:matrix diff equs vs ops}, and
%Theorem~\ref{thm:lindiff d-max}.
\end{proof} 

\noindent
Recall from  Corollary~\ref{cor:Lyap exp ReIm}  that for $f=(f_1,\dots,f_n)\in\c[\imag]^n$,
$$\lambda(f)\ =\ \min\!\big\{\lambda(f_1),\dots, \lambda(f_n)\big\}\ =\ \lambda\big(|f_1|+\cdots + |f_n|\big)\in \R_{\pm \infty}.$$
If $f\in\c^1[\imag]$ and $f\notin \c^1[\imag]^\times$, we set $\mu(f):=-\infty$. With this convention,
$$\mu(f)\ :\ =\max\!\big\{ \mu(f_1),\dots,\mu(f_n) \big\}\quad\text{ for $f=(f_1,\dots,f_n)\in \c^1[\imag]^n$.}$$
We also turn  
$K^n$ into a valued $\C$-linear space with   valuation $v\colon K^n\to\Gamma_\infty$ given by~$v(f):=\min\!\big\{v(f_1),\dots,v(f_n)\big\}$ 
for $f=(f_1,\dots,f_n)\in K^n$. 

\begin{cor}\label{cor:val indep columns}
We can choose $M$, $D$  as in Corollary~\ref{cor:shape of fund matrix}
such that the successive columns $f_1,\dots, f_n$ of $M$ have the property that for $k=1,\dots, n$ the  $f_j$ with $\phi_j=\phi_k$ are valuation-independent. For any such $M, D$, the matrix $MD$ is a Lyapunov fundamental matrix for $y'=Ny$. 
\end{cor}
\begin{proof}
Take $M$, $\phi_1,\dots,\phi_n$ as in Corollary~\ref{cor:shape of fund matrix} such that $\phi_1,\dots,\phi_m$ are distinct, $m\leq n$, and each $\phi_j$ with~$m<j\leq n$ is
equal to some $\phi_k$ with~$1\leq k\leq m$. For~$V := \sol_{\Univ}(N)$ this yields an internal direct sum decomposition
$$V\ =\ \ex^{\phi_1\imag}V_1\oplus\cdots\oplus \ex^{\phi_m\imag}V_m$$ 
into $\C$-linear subspaces of $V$.  
Now~[ADH, remark before 2.3.10] yields for~$k=1,\dots, m$ a valuation basis of $V_k$. Modifying $M$ accordingly, this yields $M, D$ with the desired property. The rest follows from Corollary~\ref{cor:Lyap val indep}.
\end{proof}

\noindent
If the matrix  $N$ is bounded, then the solutions of ${y'=Ny}$ grow only moderately, by Lemma~\ref{lem:Lyap}; their oscillation is also moderate: 

\begin{lemma}\label{lem:mod osc}
Suppose  $N$ is bounded. Let~$M\in\operatorname{GL}_n(K)$ and $\phi_1,\dots,\phi_n\in H$ be such that for $D:=\operatorname{diag}(\ex^{\phi_1\imag},\dots,\ex^{\phi_n\imag})$ the $n\times n$ matrix $MD$ over $K[\ex^{H\imag}]$
 is a fundamental matrix for $y'=Ny$. 
Then $\phi_1,\dots,\phi_n\preceq x$.
\end{lemma} 
\begin{proof}
  Corollary~\ref{cor:shape of fund matrix} yields $\Sigma(N) = \{\phi_1'\imag+K^\dagger,\dots,\phi_n'\imag+K^\dagger\}$.
The differential module over $K$ associated to $N$ (cf.~[ADH, p.~277]) is bounded, by Example~\ref{ex:lattices}(1), hence each 
$\alpha\in\Sigma(N)$   has the form~$a+K^\dagger$ 
with~$a\in \mathcal O$, by Corollary~\ref{cor:bded cyc 4}.  
Together with Lemma~\ref{lem:[O]}
this yields~${\phi_1,\dots,\phi_n\preceq x}$.
\end{proof}

\begin{cor}\label{cor:bddflattereq} Suppose $N$ is bounded. Then $\sol_{\Univ}(N)\subseteq (\Univ^{\flattereq})^n$.
\end{cor} 
\begin{proof} Take $M$ and $D$ as in Lemma~\ref{lem:mod osc}. It suffices to show that the entries of~$MD$ are in $\Univ^{\flattereq}$.  Such an entry equals $g\ex^{\phi\imag}$ where $g$ is an entry of $M$ and~$\phi\in \{\phi_1,\dots, \phi_n\}$. Lemma~\ref{lem:from phi to lambda} gives $h\in \Lambda_H$ such that
$\phi-\phi(h\imag)\preceq 1$. Now $\phi\preceq x$ by Lemma~\ref{lem:mod osc}, so $\phi(h\imag)\preceq x$, and thus
$h=\phi(h\imag)'\preceq 1$. Also $\ex^{\phi\imag}=u\ex(h\imag)$ with~$u=\ex^{(h-\phi(h\imag))\imag}\in K^\times$, so
$g\ex^{\phi\imag}=gu\ex(h\imag)$. 
 Lemma~\ref{lem:Lyap} gives $\lambda(g\ex^{\phi\imag})> -\infty$. 
Now $\lambda(g\ex^{\phi\imag})=\lambda(g)=\lambda(gu)$, so $gu\in \mathcal{O}_\Delta$. Then  $h\in \Lambda_H\cap \mathcal{O}_H$ gives $g\ex^{\phi\imag}\in \Univ^{\flattereq}$. 
\end{proof} 

\noindent
The next result shows how Lemma~\ref{lem:mod osc}  also yields information for un\-bounded~$N$. For $a\in H$ we let ``$N\preceq a$'' stand
for ``$g\preceq a$ for every entry $g$ of $N$''. 

\begin{cor}\label{cor:mod osc}
Suppose $N\preceq \ell'$, $\ell\in H^{>\R}$. Let $M\in\operatorname{GL}_n(K)$, $\phi_1,\dots,\phi_n\in H$, and $D:=\operatorname{diag}(\ex^{\phi_1\imag},\dots,\ex^{\phi_n\imag})$ be such that $MD$ is a fundamental matrix for $y'=Ny$
over $K[\ex^{H\imag}]$. Then ${\phi_1,\dots,\phi_n\preceq \ell}$, and there exists $m\ge 1$ such that
$f\preceq \ex^{m\ell}$ and~$f\not\preceq \ex^{-m\ell}$, for each column $f$ of $M$. 
\end{cor}
\begin{proof}
Put $\phi:=\ell'$ and use the  superscript $\circ$ as in \eqref{eq:derc}, Section~\ref{secfhhf}.
%on compositional conjugation in Hardy fields of Section~\ref{sec:Hardy fields}. 
Then for~$R:=\Calinf[\imag]$ and any fundamental matrix $F\in\operatorname{GL}_n(R)$ for $y'=Ny$, the matrix~$F^\circ\in \operatorname{GL}_n(R)$ is a fundamental matrix for $z'=(\phi^{-1}N)^\circ z$. As $H^\circ$ is $\d$-maximal and
 $(\phi^{-1}N)^\circ$ is bounded,  we can apply Lemmas~\ref{lem:mod osc} and ~\ref{lem:Lyap} (and a remark following
 Corollary~\ref{cor:Lyap exp lin indep}) to $M^\circ$ and $D^\circ$ in the role of $M$ and $D$, and convert this back to
 information about $M$ and $D$ as claimed.
\end{proof}

\noindent 
In the next lemma and its corollary we assume $N$ is bounded and $M$, $D$ are as in Lemma~\ref{lem:mod osc}. Let $\operatorname{st}(N)$ 
(the {\it standard part}\/ of $N$) be the $n\times n$ matrix over $\C$ such that~$N-\operatorname{st}(N)\prec 1$.
For~$f\in K^n$, put $\abs{f}:=\max\!\big\{\abs{f_1},\dots,\abs{f_n}\big\}\in H$, so  $vf=v\abs{f}$.

\begin{lemma}\label{eigenlmh}
Let   $y=\ex^{\phi\imag}f$ where $f=f_k$ is the $k$th column of $M$ and $\phi=\phi_k$, $k\in \{1,\dots,n\}$. Set~$s:=\abs{f}^{-1}f\in\mathcal O^n$. 
Then $\lambda:=\lambda(y)\in \R$, $\mu:=\mu(y)\in\R$, and~$-\lambda+\mu\imag\in\C$ is an eigenvalue of $\operatorname{st}(N)$ with eigenvector $\operatorname{st}(s)\in\C^n$.
\end{lemma}
\begin{proof} Note that $y$ is the $k$th column of $MD$. 
From Lemma~\ref{lem:Lyap} we get $\lambda\in \R$. Let $g$ be a nonzero entry of $f$. Then for the corresponding entry $g\ex^{\phi\imag}$ of $y$ we have
 $(g\ex^{\phi\imag})^\dagger= g^\dagger +  \phi'\imag$, so $\Im\!\big((g\ex^{\phi\imag})^\dagger\big)=\Im(g^\dagger)+\phi'$ with $\Im(g^\dagger)\prec 1$ by a remark preceding Lemma~\ref{lem:W and I(F)}, and $\phi'\preceq 1$ by Lemma~\ref{lem:mod osc}. Hence $\mu(g\ex^{\phi\imag})= \lim_{t\to +\infty} \phi'(t)\in \R$. This gives $\mu=\lim_{t\to +\infty} \phi'(t)\in \R$ and so $\mu-\phi'\prec 1$.  
 
 Next, $y'=Ny$ gives $\phi'\imag f + f'=Nf$,  and then using also Corollary~\ref{cor:fexlambdax},
 $$ Nf\ =\ (-\lambda + \mu\imag)f + (\phi'\imag-\mu\imag)f + \lambda f + f'\ =\  (-\lambda + \mu\imag)f+r,\quad r\in K^n,\ r\prec f.  $$
Dividing by $|f|\in H^\times$ then yields the claim about $-\lambda+\mu\imag$ and $s$. 
\end{proof}

\noindent
The proof of Lemma~\ref{eigenlmh} also gives the next corollary, where 
$I_n$ denotes the $n \times n$ identity matrix over $K$, and~$\operatorname{mult}_c\!\big(\!\operatorname{st}(N)\big):=\dim_\C \ker_{\C^n}\!\big(\!\operatorname{st}(N)-c I_n\big)$ for $c\in\C$:

\begin{cor}
For $k=1,\dots,n$, let $f_k$ be the $k$th column of $M$, so $y_k:=f_k\ex^{\phi_k\imag}$ is the $k$th column of $MD$, and put $c_k:=-\lambda(y_k)+\mu(y_k)\imag$.
If for a certain $k$  the~$f_j$ with $\mu_j=\mu_k$ are valuation-independent, then for this $k$ we have
$$\operatorname{mult}_{c_k}\!\big(\!\operatorname{st}(N)\big)\ \geq\ \abs{\{j:(\lambda_j,\mu_j)=(\lambda_k,\mu_k) \}}.$$ 
\end{cor}

\begin{question}
Suppose $N$ is bounded and $\operatorname{st}(N)$ is the $n\times n$ matrix over $\C$ such that~$N-\operatorname{st}(N)\prec 1$. By   Perron \cite[Satz~13]{Perron29} there is a Lyapunov 
fundamental matrix $F$ for~$y'=Ny$    such that
for each $\lambda\in\R$, the number of columns~$f$ of~$F$ with~$\lambda(f)=\lambda$
equals~$\sum_{\mu\in\R}\operatorname{mult}_{-\lambda+\mu\imag}\!\big(\!\operatorname{st}(N)\big)$.  
 Can one take here~$F$ of the form~$F=MD$ where~$M\in\operatorname{GL}_n(K)$
and $D=\operatorname{diag}(\ex^{\phi_1\imag},\dots,\ex^{\phi_n\imag})$ with~$\phi_1,\dots,\phi_n\in H$?
\end{question}

\noindent
Recall: a column vector $(y_1,\dots,y_n)^{\text{t}}\in\c[\imag]^n$ is said to be  {\it bounded}\/ if $y_1,\dots,y_n\preceq 1$.
%, and {\it unbounded}\/ otherwise.

\begin{lemma}\label{lem:bounded sol}
Suppose $y'=Ny$ where $y\in\c^1[\imag]^n$ is bounded.
Then $$y\ =\ \ex^{\phi_1\imag}z_1+ \cdots+  \ex^{\phi_m}z_m$$ where $m\le n$,
 $\phi_1,\dots, \phi_m\in H$ are distinct and apart,
$z_1,\dots, z_m\in K^n$ are bounded, and $\ex^{\phi_1\imag}z_1,\dots, \ex^{\phi_m\imag}z_m\in \sol_{\Univ}(N).$
% then
%there is such a $y$ with $y=\ex^{\phi\imag}z$ where $\phi\in H$ and $z\in K^n$.
\end{lemma}
\begin{proof}
Let $M\in\operatorname{GL}_n(K)$ and $\phi_1,\dots,\phi_n\in H$ be as in Corollary~\ref{cor:shape of fund matrix}, in particular, $\phi_1,\dots, \phi_n$ are apart.
For $j=1,\dots,n$, let $f_j\in K^n$ be the $j$th column of~$M$. 
Take $c_1,\dots,c_n\in\C$ such that $y=c_1\ex^{\phi_1\imag} f_1+\cdots+c_n\ex^{\phi_n\imag}f_n$.
We arrange that~$\phi_1,\dots, \phi_m$ are distinct, $m\le n$, and each $\phi_j$ with $m<j\le n$ is equal to one of the
$\phi_k$ with $1\le k\le m$. This gives $y=\ex^{\phi_1\imag} z_1+\cdots + \ex^{\phi_m}z_m$ with $z_1,\dots, z_m\in K^n$ and
$\ex^{\phi_1\imag}z_1,\dots, \ex^{\phi_m\imag}z_m\in \sol_{\Univ}(N)$. 
 The $\preceq$-version of
Corollary~\ref{cor:from phi to lambda} with $\fm=1$ then shows that $z_1,\dots, z_m$ are bounded.
% yields $j\in\{1,\dots,m\}$ such that $\sum_{\phi_k=\phi_j} c_kf_k$ is nonzero and bounded.
%Then   $\phi:=\phi_j$ and $z:=\sum_{\phi_k=\phi_j} c_kf_k$ have the required property.
\end{proof}

\noindent
See \cite[Chapter~2]{Bellman} and~\cite[Chapter~II, \S{}3]{Cesari} for classical conditions on a matrix  differential equation to
have only bounded solutions.

\medskip\noindent
Despite Corollary~\ref{cor:shape of fund matrix}, 
the differential fraction field of $K[\ex^{H\imag}]$ is not pv-closed, since it is not even  algebraically closed;
see~[ADH, 5.1.31] and Lemma~\ref{lem:Frac U not alg closed}.
Combining Corollary~\ref{cor:d-max weakly d-closed} and  [ADH, 5.4.2] also yields:

\begin{cor}
For every column $b\in K^n$ the matrix differential equation~$y'=Ny+b$ has a solution in $K^n$. 
\end{cor}

\noindent
We also have a version of Corollary~\ref{cor:nonosc zeros, self-dual op} for matrix differential equations:

\begin{cor}\label{cor:nonosc zeros, self-dual equ}
Suppose $y'=Ny$ is self-dual.
If $\alpha$ is an eigenvalue of $y'=Ny$, then so is  $-\alpha$, with the same multiplicity.
If $n$ is odd, then the matrix differential equation $y'=Ny$  has a 
solution $y\neq 0$ in~$K^n$.
\end{cor}

\noindent
This follows from Corollaries~\ref{cor:self-dual d-equ} and~\ref{cor:shape of fund matrix}.
Note that the hypothesis on $N$ in Corollary~\ref{cor:nonosc zeros, self-dual equ} is satisfied if~$y'=Ny$ is self-adjoint  or hamiltonian. 

If $y'=Ny$ is self-adjoint, and $M\in \operatorname{GL}_n(K)$ and $\phi_1,\dots,\phi_n\in H$ (as in Corollary~\ref{cor:shape of fund matrix}) are such that $MD$ is a fundamental matrix for $y'=Ny$ where~$D:=\operatorname{diag}(\ex^{\phi_1\imag},\dots,\ex^{\phi_n\imag})$, then 
there exists $U\in\operatorname{GL}_n(\C)$ such that
the fundamental matrix $MDU\in \operatorname{GL}_n\big(K[\ex^{H\imag}]\big)$ of $y'=Ny$ is orthogonal as an element of~$\operatorname{GL}_n(\Omega)$, where $\Omega$ is the differential fraction field of $K[\ex^{H^\imag}]$; likewise with ``hamiltonian'' and ``symplectic'' instead of
``self-adjoint'' and ``orthogonal'': Lemmas~\ref{lem:orthogonal fund matrix} and~\ref{lem:symplectic fund matrix}. 

\begin{example}
Any matrix differential equation $y'=Ny$ with 
$$N = \begin{pmatrix} 0 & a & b \\ -a & 0 & c \\ -b & -c & 0 \end{pmatrix} \qquad (a,b,c\in K=H[\imag])$$
has a nonzero solution $y=(y_1, y_2, y_3)^{\text{t}}\in K^3$. 
\end{example}

 \noindent
In the self-dual case we can  improve on
Corollary~\ref{cor:shape of fund matrix}:

\begin{cor}\label{cor:const norm}
Suppose $y'=Ny$ is self-dual. Then there are~$M\in\operatorname{GL}_n(K)$ and
$\phi_1,\dots,\phi_n\in H$ that are apart
such that 
\begin{enumerate}
\item[\textup{(i)}] for each $j\in\{1,\dots,n\}$ there is a $k\in\{1,\dots,n\}$ with $\phi_j=-\phi_k$;
\item[\textup{(ii)}] with $D:=\operatorname{diag}(\ex^{\phi_1\imag},\dots,\ex^{\phi_n\imag})$, the $n\times n$ matrix $MD$ over $K[\ex^{H\imag}]$ is a fundamental matrix for~$y'=Ny$. 
\end{enumerate}
%If $y'=Ny$ is self-adjoint, then for all such $M$ and $\phi_1,\dots,\phi_n$, if  $(f_1,\dots,f_n)^{\operatorname{t}}$ is the $j$th column of $M$ and $(g_1,\dots,g_n)^{\operatorname{t}}$ is the $k$th column of $M$, then  $f_1g_1+\cdots+f_ng_n\in\C$,  with $f_1g_1+\cdots+f_ng_n=0$ if~$\phi_j\neq-\phi_k$.
\end{cor}
\begin{proof}
Corollary~\ref{cor:self-dual equivalence} yields a matrix differential equation $y'=A_Ly$ over $K$, equivalent to $y'=Ny$, where $L\in K[\der]$ is monic
self-dual of order $n$. Then we can use Corollary~\ref{cor:self-dual basis} instead of Theorem~\ref{thm:lindiff d-max} in the proof of Corollary~\ref{cor:shape of fund matrix}.
\end{proof}

\noindent
Let $\Omega$ be the differential fraction field of $K[\ex^{H\imag}]$
and $V:=\operatorname{sol}_\Omega(N)\subseteq K[\ex^{H\imag}]^n$, a $\C$-linear subspace of $\Omega^n$.  
Then $\dim_{\C}V=n$ and $V=\operatorname{sol}_{\Calinf[\imag]}(N)$.
In the   corollary below we assume that $y'=Ny$ is self-adjoint, and we equip $V$ with the symmetric bilinear form $\langle\ ,\,\rangle$ defined after Lemma~\ref{lem:constant norm} (with $\Omega$ instead of $K$).

\begin{cor}\label{cor:const norm, 1}
There are $m\leq n$ and distinct~$\theta_1,\dots,\theta_m$ in $H^{>\R}$ that are apart,
 subspaces~$V_1,\dots,V_m,W$ of the $\C$-linear space $V$ with $W\subseteq K^n$, and for~$j=1,\dots,m$, nonzero subspaces~$V_j^{+},V_j^{-}$ of $K^n$, such that
\begin{align*} V_j\ &=\ V_j^{+}\ex^{\theta_j\imag}\oplus V_j^{-}\ex^{-\theta_j\imag} \qquad \text{\textup{(}internal direct sum of subspaces of $V_j$\textup{)},}\\
 V\ &=\ V_1 \perp \cdots \perp V_m \perp W\qquad\text{\textup{(}orthogonal sum with respect to $\langle\ ,\,\rangle$\textup{)}.}
\end{align*} 
For any such $m$ and $\theta_j$, $V_j$,  $V_j^{+},V_j^{-}$, $W$ we have
$\dim_{\C} V_j^{+} = \dim_{\C} V_j^{-}$ and $\langle\ ,\,\rangle$ restricts to a null form on
$V_j^{+}\ex^{\theta_j\imag}$ and on $V_j^{-}\ex^{-\theta_j\imag}$.
\end{cor}
\begin{proof}
Take $M$ and  $\phi_1,\dots,\phi_n$ as in Corollary~\ref{cor:const norm}, and set
$$D\ :=\ \operatorname{diag}(\ex^{\phi_1\imag},\dots,\ex^{\phi_n\imag}).$$
Let $(f_1,\dots,f_n)^{\operatorname{t}}$ be the $j$th column of $M$ and $(g_1,\dots,g_n)^{\operatorname{t}}$
be the $k$th column of~$M$.  
The $j$th column of $MD$ is~$f=(f_1\ex^{\phi_j\imag},\dots,f_n\ex^{\phi_j\imag})^{\operatorname{t}}$, and
the $k$th column of $MD$  is~$g=(g_1\ex^{\phi_k\imag},\dots,g_n\ex^{\phi_k\imag})^{\operatorname{t}}$,
and $f,g\in\operatorname{sol}_\Omega(N)$ by Corollary~\ref{cor:const norm}(ii).
Thus by Lemma~\ref{lem:constant norm} applied to $\Omega$ in place of $K$ we have
$$\langle f,g\rangle\ =\ (f_1g_1+\cdots+f_ng_n)\ex^{(\phi_j+\phi_k)\imag}\in\C.$$ 
Corollary~\ref{cor:const norm}(i) gives $l\in\{1,\dots,n\}$ with~$\phi_k=-\phi_l$; then $\phi_j+\phi_k=\phi_j-\phi_l$. Hence if~$\phi_j\neq-\phi_k$, then~$\phi_j+\phi_k\succ 1$ by $\phi_j$, $\phi_l$ being apart, so $\ex^{(\phi_j+\phi_k)\imag}\notin K$ by Corollary~\ref{cor:osc => bded} and thus~$\langle f,g\rangle=0$. 
Taking $\theta_1,\dots,\theta_m$ to be the distinct positive elements of $\{\phi_1,\dots,\phi_n\}$,
this yields the existence statement. The rest follows from 
Corollary~\ref{cor:self-dual d-equ}, Lemma~\ref{matrixversionlem}, and again Corollary~\ref{cor:osc => bded}.
%the ``moreover'' part of Corollary~\ref{cor:shape of fund matrix}.
\end{proof}

\begin{remark}
Suppose   $y'=Ny$ is hamiltonian.
Then Corollary~\ref{cor:const norm, 1} remains true with~$\langle\ ,\,\rangle$ replaced by the alternating bilinear form
$\omega$ on $V$ of Lemma~\ref{lem:constant norm, symplectic}.
(Same proof, using Lemma~\ref{lem:constant norm, symplectic} instead of
Lemma~\ref{lem:constant norm}.)
\end{remark}

\noindent
The complex conjugation automorphism of the differential ring $\Calinf[\imag]$ restricts to an automorphism of the differential integral domain
$\Univ=K[\ex^{H\imag}]$, which in turn extends uniquely to an automorphism $g\mapsto \bar{g}$ of the differential field $\Omega$, with $\bar{\bar{g}}=g$ for all $g\in \Omega$. Let 
$\Omega_{\operatorname{r}}$ be the fixed field of this automorphism of $\Omega$. Then $\Omega_{\operatorname{r}}$ is 
a differential subfield of $\Omega$, and $\Omega=\Omega_{\operatorname{r}}[\imag]$. Set $\Univ_{\operatorname{r}} := \Omega_{\operatorname{r}}\cap \Univ$. Then
$$\Omega_{\operatorname{r}}\ =\ \Frac(\Univ_{\operatorname{r}})\ \text{ inside }\Omega,\qquad \Univ_{\operatorname{r}}\  =\ \Univ\cap\ \Calinf\ \text{\ inside }\Calinf[\imag],\quad \Univ\ =\ \Univ_{\operatorname{r}}[\imag]. $$ 
{\em Assume in the rest of this subsection that~$y'=Ny$ is anti-self-adjoint, and equip the $\C$-linear space $V=\operatorname{sol}_{\Omega}(N)$
with the positive definite hermitian form $\langle\ ,\,\rangle$ introduced after Lemma~\ref{lem:constant norm, anti-self-adj}, with
$\Omega$ in the role of $K$}. 
%Then for $H$, $M$, $D$ as in 
%Corollary~\ref{cor:shape of fund matrix} there exists $U\in\operatorname{GL}_n(\C)$ such that
%the fundamental matrix~$MDU\in \operatorname{GL}_n\!\big(K[\ex^{H\imag}]\big)$ of $y'=Ny$ is unitary as an element of $\operatorname{GL}_n(\Omega)$, by Lemma~\ref{lem:unitary fund matrix} applied to $\Omega$, $\Omega_{\operatorname{r}}$ in place of $K$, $H$, respectively.
Then we have the following analogue of Corollary~\ref{cor:const norm, 1}:
% allows us to strengthen this last statement. We now 

\begin{cor}\label{cor:const norm, 2}
There are $m\in \{1,\dots, n\}$, distinct $\theta_1,\dots,\theta_m\in H$ that are apart,  
and nonzero $\C$-linear subspaces~$V_1,\dots,V_m$ of $K^n$ such that $V$ is the following orthogonal sum with respect to $\langle\ ,\,\rangle$:
$$V\ =\ V_1\ex^{\theta_1\imag} \perp \cdots \perp V_m\ex^{\theta_m\imag}.$$
\end{cor}
\begin{proof} Corollary~\ref{cor:shape of fund matrix} gives $M\in \text{GL}_n(K)$ and $\phi_1,\dots,\phi_n\in H$ that are apart such that $MD$ is a fundamental matrix for $y'=Ny$ where $D:=\text{diag}\big(\ex^{\phi_1\imag},\dots, \ex^{\phi_n\imag}\big)$. 
For $j,k=1,\dots,n$, let $f=(f_1\ex^{\phi_j\imag},\dots,f_n\ex^{\phi_j\imag})^{\operatorname{t}}$
be the $j$th column of~$MD$ and~$g=(g_1\ex^{\phi_k\imag},\dots,g_n\ex^{\phi_k\imag})^{\operatorname{t}}$ the $k$th column of $MD$, where 
$f_1,\dots, f_n$ and $g_1,\dots, g_n$ are in~$K$.
Then  by Lemma~\ref{lem:constant norm, anti-self-adj},  
$$\langle f,g\rangle\ =\ (f_1\overline{g_1}+\cdots+f_n\overline{g_n})\ex^{(\phi_j-\phi_k)\imag}\in\C$$ 
and hence $\langle f,g\rangle=0$ if $\phi_j\neq\phi_k$, by Corollary~\ref{cor:osc => bded}.
Taking $\theta_1,\dots,\theta_m$ to be the distinct   elements of $\{\phi_1,\dots,\phi_n\}$,
this yields the desired result. 
\end{proof}

\noindent
 %Suppose $H$ is $\d$-maximal, and take $m$, $\theta_1,\dots, \theta_m, V_1,\dots, V_m$ as in 
% Corollary~\ref{cor:const norm, 2}. Then \cite[XV, Corollary~5.2]{Lang} gives an orthonormal basis 
%$f_1\ex^{\theta_1\imag},\dots, f_n\ex^{\theta_m\imag}$ of $V$ with respect to $\<\cdot,\cdot\>$, where
%$f_1,\dots, f_n\in K^n$; here $\theta_j$ occurs as often as the dimension of $V_j$ for $j=1,\dots,m$, and accordingly we set
%$D=\text{diag}(\ex^{\theta_1\imag},\dots, \ex^{\theta_m\imag})$.
%Let $M$ be the $n\times n$ matrix over $K$  with $k$th column $f_k$ for $k=1,\dots,n$. Then $MD$ has columns
%$f_1\ex^{\theta_1\imag},\dots, f_n\ex^{\theta_m\imag}$, and so is a unitary as an element of
%$ \operatorname{GL}_n(\Omega)$.
Corollary~\ref{cor:const norm, 2} and \cite[Chapter~XV, Corollary~5.2]{Lang} yield~$M\in \text{GL}_n(K)$
 and $\phi_1,\dots, \phi_n\in H$ that are apart, such that $MD$ with~$D:=\text{diag}\big(\ex^{\phi_1\imag},\dots, \ex^{\phi_n\imag}\big)$ is not only a fundamental matrix for $y'=Ny$ but also unitary as an element of 
 $\operatorname{GL}_n(\Omega)$. (Cf.~Lemma~\ref{lem:unitary fund matrix}.)

\begin{example}[Schr\"odinger equation for   quantum systems with $n$ states {\cite[\S{}3.4]{Teschl}}]
This is the matrix differential equation
$y'=-\imag S y$ where the $n\times n$ matrix $S$ over~$K$ (the Hamiltonian of the system) is hermitian, i.e., $S^{\operatorname{t}}=\overline{S}$. 
Then $y'=-\imag S y$ is anti-self-adjoint, so 
we have the positive definite hermitian form 
$$(y,z)\mapsto \langle y,z\rangle = y_1\bar{z_1}+\cdots+y_n\bar{z_n} \quad (y=(y_1,\dots,y_n)^{\operatorname{t}},\  z=(z_1,\dots,z_n)^{\operatorname{t}})$$
on the $\C$-linear space of solutions $W$ of $y'=-\imag S y$ in $\Calinf[\imag]$.   
There are $\phi_1,\dots, \phi_n\in H$ that are apart and $f_1,\dots, f_n\in K^n$ such that
$$f_1\ex^{\phi_1\imag},\dots,f_n\ex^{\phi_n\imag}  \qquad(\text{``wave functions''})$$ 
is an orthonormal basis of $W$ with respect to $\langle\ ,\,\rangle$. Note again the striking fact that~$\langle y, y\rangle$ is a positive real constant, not just an element of $H^{>}$,
for every  $y\in W^{\neq}$. 
\end{example}

\subsection*{Definability} Here we drop the $\d$-maximality assumption from earlier subsections. 
We begin with consequences of our earlier boundedness results for matrix differential equations depending on (constant) parameters. For this, let~$H\supseteq\R$, let $C=(C_1,\dots,C_m)$ be a tuple of distinct indeterminates over~$K$ 
and let $N(C)$ be an $n\times n$ matrix over the polynomial ring~$K[C]$, $n\ge 1$. 
Then for~$c\in\R^m$ we have the matrix differential equation~$y'=N(c)y$   over $K$.
Combining Corollary~\ref{cor:sa} with  Lemmas~\ref{lem:only bounded sols} and~\ref{lem:nonzero unbounded sol}   yields:

\begin{cor}\label{cor:bddsol}
The set of $c\in\R^m$ such that all solutions of $y'=N(c)y$ in  $\c^1[\imag]^n$ are bounded is  semialgebraic, and so is the set of
$c\in \R^m$ such that $y'=N(c)y$ has no nonzero bounded solution in $\c^1[\imag]^n$.
\end{cor}

\noindent
For  matrix differential equations depending analytically on a single parameter, see~\cite[Chap\-ter~VII]{Wasow}. In this connection we record that for $m=1$ it follows from Corollary~\ref{cor:bddsol}:  if $y'=N(c)y$ has for arbitrarily large
$c\in \R$ only bounded solutions in $\c^1[\imag]^n$, then this happens for all sufficiently large $c\in \R$;  if $y'=N(c)y$ has for arbitrarily large
$c\in \R$ no nonzero bounded solution in $\c^1[\imag]^n$, then this happens for all sufficiently large $c\in \R$.

By the next lemma, the property of a matrix differential equation over a complexified Hardy field to have only
bounded solutions is ``uniformly definable'' from the entries in the matrix. 
Here we view the canonical $\HLO$-expansion~$\mathbf H$ of a Hardy field $H$ as a structure
for the language $\mathcal L^\iota_{\HLO}$  from~[ADH, Chapter~16], and we let~${u=(u_{ij})}$, $v=(v_{ij})$ be disjoint multivariables of size~$n\times n$ with $n\ge 1$. 

\begin{lemma}\label{lem:only bounded sols}
There is a quantifier-free $\mathcal L^\iota_{\HLO}$-formula $\beta(u,v)$ such that for every Hardy field $H$ and
$n\times n$ matrix $N$ over $K=H[\imag]$:
$$\mathbf H\models \beta(\Re N,\Im N)\quad\Longleftrightarrow\quad
\text{all solutions of $y'=Ny$  in $\c^1[\imag]^n$ are bounded.}$$
\end{lemma}
\begin{proof}
Let $N$ be an $n\times n$ matrix over $K$, $\phi\in H$, $z\in K^n$.  Then  $\ex^{\phi\imag}z\in K[\ex^{H\imag}]^n$ is a solution of $y'=Ny$ iff
$z'+\phi'\imag z=Nz$. Moreover, for $\d$-maximal $H$ it follows from Corollary~\ref{cor:shape of fund matrix} that
all solutions of $y'=Ny$ in $\c^1[\imag]^n$ are bounded iff all solutions~$\ex^{\phi\imag}z$ with $\phi\in H$ and $z\in K^n$
are bounded. Now use that $\d$-maximal Hardy fields are $H$-closed and that the theory of $H$-closed $H$-fields admits quantifier elimination in the language $\mathcal L^\iota_{\HLO}$.
%Suppose now $H$ is $\d$-maximal, and
%let~$S$ be the set of triples $(M,N,\phi)$ consisting of  an invertible $n\times n$ matrix $M$ over~$K$ with bounded coefficients, an $n\times n$ matrix $N$ over $K$, and~$\phi=(\phi_1,\dots,\phi_n)\in H^n$ such that with~$D:=\operatorname{diag}(\ex^{\phi_1\imag},\dots,\ex^{\phi_n\imag})$,  the $n\times n$ matrix $MD$ over~$K[\ex^{H\imag}]$ 
%is a fundamental matrix for $y'=Ny$.
%Then by the remark at the beginning of this proof, $S$ (viewed as a subset of $H^{2n^2}\times H^{2n^2}\times H^n=H^{n(4n+1)}$ in the natural way)  is definable in~$\mathbf H$ by a quantifier-free $\mathcal L^\iota_{\HLO}$-formula, independent of $H$;
%hence so is the image  of $S$ under the projection $(M,N,\phi)\mapsto N$, by [ADH, 16.0.1] and Theorem~\ref{thm:char d-max}. Together with Corollary~\ref{cor:shape of fund matrix} this yields the lemma.
\end{proof}

\noindent
Using Lemma~\ref{lem:bounded sol}  we obtain in the same way:

\begin{lemma}\label{lem:nonzero unbounded sol}
There is a quantifier-free $\mathcal L^\iota_{\HLO}$-formula $\gamma(u,v)$ such that for every Hardy field $H$ and
$n\times n$ matrix $N$ over $K=H[\imag]$:
$$\mathbf H\models \gamma(\Re N,\Im N)\ \Longleftrightarrow\
\text{some nonzero solution of $y'=Ny$ in $\c^1[\imag]^n$  is bounded.}$$
\end{lemma}
%\begin{proof}
%Similarly to the proof of Lemma~\ref{lem:only bounded sols}, with $S$ replaced by the set of triples~$(N,\phi,z)$ where $N$ is an $n\times n$ matrix over $K$, $\phi\in H$, and $z\in K^n$ is bounded 
%such that $y'=Ny$ for $y:=\ex^{\phi\imag}z$, and using Lemma~\ref{lem:bounded sol} in place of Corollary~\ref{cor:shape of fund matrix}.
%\end{proof}

\begin{example} Let $a,b\in H$, and take
$g,\phi\in\operatorname{Li}\big(H(\R)\big)$ with $g\neq 0$, $g^\dagger=a$, and $\phi'=b$.
Then $\big\{y\in \c^1[\imag]:\, y'=(a+b\imag)y\big\}=\C g\ex^{\phi\imag}$, and $g\ex^{\phi\imag} \asymp g$.  Thus if $H$ is Liouville closed, then by [ADH, 11.8.19]:
\begin{align*}
&\phantom{\Longleftrightarrow}\qquad \text{every $y\in\c^1[\imag]$ with $y'=(a+b\imag)y$ is bounded} \\
& \Longleftrightarrow\quad \text{some $y\in\c^1[\imag]^{\ne}$  with $y'=(a+b\imag)y$  is  bounded}  \\
& \Longleftrightarrow\quad  a\notin\Upg(H)  \\
& \Longleftrightarrow\quad  a\leq 0 \text{ or }a\in\I(H).
\end{align*}
% Hence for $n=1$, $u=u_{11}$ we can take the  quantifier-free $\mathcal L^\iota_{\HLO}$-formula $u\leq 0 \vee \I(u)$
 %for both~$\beta$,~$\gamma$.
\end{example}

\subsection*{The real case}
{\it In this subsection we assume~$A\in H[\der]$.}\/  
Recall that if $H$ is $\d$-maximal, then $K$ is linearly closed, so~[ADH, 5.1.35] yields the following, which includes Corollary~\ref{cor:factorization intro, real} from the introduction:

\begin{cor}
If $H$ is $\d$-maximal, then $A$ is a product of irreducible operators in~$H[\der]$ which are monic of order~$1$ or
monic of order~$2$.
\end{cor}

\noindent
The next result follows from Corollaries~\ref{cor:cos sin infinitesimal} and~\ref{cor:basis of kerU A}, and is a  version of Theorem~\ref{thm:lindiff d-max} in the case of a real operator:

\begin{cor}\label{cor:lindiff, real}
Let $E$ be a $\d$-maximal Hardy field extension of $H$. Then the $\C$-linear space  $V:=\ker_{\Calinf[\imag]}A$ of zeros of $A$ in $\Calinf[\imag]$ has a basis
$$g_1\ex^{\phi_1\imag},\,g_1\ex^{-\phi_1\imag},\ \dots,\ g_m\ex^{\phi_m\imag},\,g_m\ex^{-\phi_m\imag}, \ h_1,\ \dots,\ h_n \qquad (2m+n=r),$$ 
where $g_j,\phi_j\in E^{>}$ with $\phi_j \succ 1$ \textup{(}$j=1,\dots,m$\textup{)} and
$h_k\in E^\times$ \textup{(}$k=1,\dots,n$\textup{)}. For any such basis of $V$, 
the $\R$-linear space $V\cap\Calinf$ of zeros of $A$ in $\Calinf$ has basis
$$g_1\cos \phi_1,\, g_1\sin \phi_1,  \ \dots, \ 
  g_m\cos \phi_m,\, g_m\sin \phi_m, \  h_1,\ \dots,\  h_n,$$
and the $\R$-linear space $V\cap E$ has basis  $h_1,\dots, h_n$.
\end{cor}

\begin{remarks}
Let $E$ be a $\d$-maximal Hardy field extension  of~$H$.
The quantity $n=\dim_{\R} \ker_E A$ in Corollary~\ref{cor:lindiff, real} (and hence also $m=(r-n)/2$) is independent of the choice of $E$, by Theorem~\ref{thm:transfer}.
Likewise, the number of distinct eigenvalues of~$A$ with respect to $E[\imag]$ does not depend on $E$. In more detail, the tuple $(d, \mu_1,\dots, \mu_d)$ where $d$ is the number of distinct eigenvalues of $A$ and $\mu_1\ge \cdots \ge \mu_d\ge 1$ are their multiplicities, with respect to $E[\imag]$, does not depend on $E$.
%(Use Corollary~\ref{corbasiseigenvalues}.) 
\end{remarks}

\noindent
Corollary~\ref{cor:lindiff, real} yields:
%(and Corollary~\ref{cor:odd degree} applied to the Riccati transform of $A$) yields:

\begin{cor}\label{cor:lindiff odd order}
If $r$ is odd, then $A(y)=0$ for some $H$-hardian germ $y\ne 0$.
\end{cor}

\noindent
From Corollary~\ref{cor:unique antider, generalized} we obtain:

\begin{cor}\label{corabcde}
Suppose $H$ is $\d$-maximal, and
let $\phi>\R$ be an element of $H$ such that $\phi'\imag+K^\dagger$ is not an eigenvalue of $A$. Then for every $h\in H$
there are unique~$f,g\in H$ such that $A(f\cos\phi+g\sin\phi)=h\cos\phi$.
\end{cor}

\noindent
Taking $A=\der$, with $K^\dagger$ as the only eigenvalue (Example~\ref{ex:ev order 1}), we recover  the following result due to Shackell~\cite[Theo\-rem~2]{Shackell17}; his proof
%about integrating certain oscillating germs arising from Hardy field elements. 
is based on~\cite{Boshernitzan87}. 
 
\begin{cor} \label{cor:Shackell}
Let $h,\phi\in H$ and $\phi>\R$.
The germ $h\cos \phi\in\Calinf$ has an antiderivative  $f\cos \phi+g\sin \phi\in \Calinf$ 
 with  $f$, $g$ in a Hardy field extension of $H$, and
any Hardy field extension of $H$ contains at most one such pair $(f,g)$.
\end{cor}

\noindent
Besides Corollary~\ref{corabcde} we use here that by a remark preceding Lemma~\ref{lem:W and I(F)} we have $\phi'\imag\notin K^\dagger$ for $\phi\in H$ with $\phi>\R$.

\medskip\noindent
We also record a real version of Corollary~\ref{cor:shape of fund matrix}, which 
follows from Corollary~\ref{cor:lindiff, real} in the same way that Corollary~\ref{cor:shape of fund matrix} followed from
Theorem~\ref{thm:lindiff d-max}.
Let $I_m$ denote the $m\times m$ identity matrix. Recall that $\Univ_{\operatorname{r}}=K[\ex^{H\imag}]\cap \Calinf$.

\begin{cor}\label{cor:shape of fund matrix, real}
Suppose $H$ is $\d$-maximal and $N$ is an $n\times n$ matrix over $H$, $n\ge 1$. Then there are $M\in\operatorname{GL}_n(H)$    as well as $k,l\in\N$ with $2k+l=n$ and
$$D=\begin{pmatrix} \ 
\boxed{D_1}	&			&				&\\
			& \ddots	& 				&\\
			&			& \boxed{D_k}	& \\
			&			&				& \boxed{I_{l}}\ \end{pmatrix}
\quad\text{where $D_j=\left(\begin{smallmatrix}			
\cos \phi_j		& \sin \phi_j \\
{-\sin \phi_j}	& \cos \phi_j
\end{smallmatrix}\right)$, $\phi_j\in H$, $\phi_j>\R$}$$
such that 
the $n\times n$ matrix $MD$ is a fundamental matrix for $y'=Ny$ with respect to $\Univ_{\operatorname{r}}$. In particular, if $n$ is odd, then $y'=Ny$ for some $0\neq y\in H^n$.
\end{cor}

\begin{proof} Let $R$ and $\Omega$ be as in the proof of Corollary~\ref{cor:shape of fund matrix}, and  let
$L\in H[\der]$ be monic of order $n$. Then Corollary~\ref{cor:lindiff, real} yields
$g_1,\dots,g_k,h_1,\dots,h_l$ and $\phi_1,\dots,\phi_k>\R$  in~$H$, where~$2k+l=n$, such that
the $\R$-linear space $\ker_{\Calinf} L$ has basis
$$g_1\cos \phi_1,\, g_1\sin \phi_1,  \ \dots, \ 
  g_k\cos \phi_k,\, g_k\sin \phi_k, \  h_1,\ \dots,\  h_l.$$
  This is also a basis of the $\C$-linear space $\ker_{\Calinf[\imag]} L=\ker_\Omega L$. 
Thus
$$W:=\operatorname{Wr}(g_1\cos \phi_1, g_1\sin \phi_1,   \dots, 
  g_k\cos \phi_k, g_k\sin \phi_k,   h_1, \dots,  h_l)\in\operatorname{GL}_n(\Omega).$$
Note that $\Univ_{\text{r}}=R\cap \Calinf$. It is routine to verify that $W=QD$ where $Q$ is an~$n\times n$~matrix over $H$ and $D$ is the~${n\times n}$ matrix over $\Univ_{\text{r}}$ displayed above.
We have~$\det D=1$, hence $$\det W=\det Q\in H\cap \Omega^\times=H^\times\subseteq \Univ_{\text{r}}^\times$$
and thus $Q\in\operatorname{GL}_n(H)$ and
$W\in\operatorname{GL}_n(\Univ_{\text{r}})$.  So by the remarks before~[ADH, 5.5.14], $W$ is a fundamental matrix
for $y'=A_Ly$ with $\Univ_{\text{r}}$ as the ambient differential ring.

%Consider the differential integral domain $R:=K[\ex^{H\imag}]\cap\Calinf$ extending~$H$, with differential fraction field $\Omega$.
%Let $L\in H[\der]$ be monic of order $n$. Then Corollary~\ref{cor:lindiff, real} yields
%$g_1,\dots,g_k,h_1,\dots,h_l$ and $\phi_1,\dots,\phi_k>\R$  in $H$, where~$2k+l=n$, such that
%the $\R$-linear space $\ker_{\Calinf} L=\ker_\Omega L$ has basis
%$$g_1\cos \phi_1,\, g_1\sin \phi_1,  \ \dots, \ 
 % g_k\cos \phi_k,\, g_k\sin \phi_k, \  h_1,\ \dots,\  h_l.$$
%Thus
%$$W:=\operatorname{Wr}(g_1\cos \phi_1, g_1\sin \phi_1,   \dots, 
 % g_k\cos \phi_k, g_k\sin \phi_k,   h_1, \dots,  h_l)\in\operatorname{GL}_n(\Omega).$$
%One verifies easily that $W=QD$ where $Q$ is an $n\times n$ matrix over $H$ and $D$ is the~${n\times n}$ matrix over $R$ displayed above.
%We have  $\det D=1$, hence $\det W=\det Q\in H\cap \Omega^\times=H^\times\subseteq R^\times$
%and thus $Q\in\operatorname{GL}_n(H)$ and
%$W\in\operatorname{GL}_n(R)$. So by the remarks before~[ADH, 5.5.14], $W$ is a fundamental matrix
%for $y'=A_Ly$ with $R$ as the ambient differential ring.

Now take monic $L\in H[\der]$ such that $y'=Ny$ is equivalent to $y'=A_Ly$, with respect to the differential field $H$.  Then take $P\in\operatorname{GL}_n(H)$
such that $P\operatorname{sol}_{\Univ_{\text{r}}}(A_L)=\operatorname{sol}_{\Univ_{\text{r}}}(N)$, and let $W$ be as above. Then
$PW\in\operatorname{GL}_n(\Univ_{\text{r}})$ is  a fundamental matrix for $y'=Ny$.
With $D$, $Q$ as before such that $W=QD$, and~$M:=PQ\in\operatorname{GL}_n(H)$, we have
$MD=PW$.
\end{proof}

\begin{example} Let $T$ be an $n\times n$ matrix over $H$, $n\ge 1$, and suppose $T$ is skew-symmetric, that is, $T^{\text{t}}=-T$.
Then the purely imaginary matrix $S:=-\imag T$ is hermitian, giving the Schr\"odinger equation $y'= Ty\ (=-\imag Sy)$ as in the example after Corollary~\ref{cor:const norm, 2}.  If $n$ is odd, then this equation has a solution $y\in E^n$ with $y_1^2+\cdots + y_n^2=1$ for some Hardy field extension $E$ of $H$; such~$y$ exhibits no oscillatory behavior and hence is a ``degenerate'' wave.
\end{example}

\noindent
We don't know whether in Corollary~\ref{cor:shape of fund matrix, real} for $n\geq 4$ we can choose~$\phi_1,\dots,\phi_k$ to be apart.
In the next corollary we let  $\langle\ ,\,\rangle: \Omega_{\text{r}}^n\times \Omega_{\text{r}}^n\to \Omega_{\text{r}}$ denote the usual symmetric bilinear form on $\Omega_{\text{r}}^n$, $n\ge 1$,
where $\Omega_{\text{r}}=\Frac(\Univ_{\text{r}})$.

\begin{cor}
Suppose $H$ is $\d$-maximal, $y'=Ny$ is self-adjoint, and~$M$,~$k$,~$l$, and~$D$ are as in Corollary~\ref{cor:shape of fund matrix, real}.
Then $\langle f,f\rangle\in\R^>$ for each column $f$ of $MD$.
Let   $f_1,g_1,\dots,f_k,g_k,h_1,\dots,h_l\in\Univ_{\operatorname{r}}^n$ be the $1$st, $2$nd,\dots, $n$th column of $MD$.
Then for $i=1,\dots,k$, $j=1,\dots,l$ we have $\langle f_i,g_i\rangle = \langle f_i,h_j\rangle = \langle g_i,h_j\rangle = 0$.
\end{cor}
\begin{proof} For any columns $f,g$ of $MD$ we have $\langle f,g\rangle\in\R$, 
by Lemma~\ref{lem:constant norm}. This proves the first claim. 
Let $i\in\{1,\dots,k\}$  and set~$\phi:=\phi_i$.
Then~$f_i=f\cos\phi-g\sin\phi$, $g_i=f\sin\phi+g\cos\phi$ where $f,g\in H^n$.
Hence
$$\langle f_i,g_i \rangle\ =\ (\langle f,f\rangle - \langle g,g\rangle)\cos\phi\sin\phi+\langle f,g\rangle (\cos^2\phi-\sin^2\phi) \in\R.$$
Lemma~\ref{lem:from phi to lambda} gives $\ex^{2\phi\imag}=\theta\ex(\lambda)$, $\ex^{-2\phi\imag} =\theta^{-1}\ex(-\lambda)$ with $\theta\in K^\times$,
$\lambda\in \Lambda^{\ne}$, so the elements $1$, $\ex^{2\phi\imag}$, $\ex^{-2\phi\imag}$ of 
$K[\ex^{H\imag}]$ are $K$-linearly independent. In view of
$$\cos\phi\sin\phi\ =\ \textstyle\frac{1}{4\imag}(\ex^{2\phi\imag}-\ex^{-2\phi\imag}),
\qquad \cos^2\phi-\sin^2\phi\  =\ \frac{1}{2} (\ex^{2\phi\imag}+\ex^{-2\phi\imag}),$$
this yields $\langle f_i,g_i \rangle=0$. 
For $j\in\{1,\dots,l\}$ we have
$$ \langle f_i,h_j\rangle = \langle f,h_j\rangle \cos\phi - \langle g,h_j\rangle\sin\phi  \in\R,
\quad
\langle g_i,h_j\rangle = \langle f,h_j\rangle \sin\phi + \langle g,h_j\rangle\cos\phi  \in\R,$$
and we obtain likewise $\langle f_i,h_j\rangle = \langle g_i,h_j\rangle = 0$.
\end{proof}

\noindent
Corollary~\ref{cor:lindiff, real} holds for $r=1$ with the assumption  ``$E$ is $\d$-maximal'' weakened to ``$E$ is Liouville closed''. The next section has more about the case $r=2$. The next lemma shows that Corollary~\ref{cor:lindiff, real} fails for~$r=3$ 
with the hypothesis  ``$E$ is $\d$-maximal'' replaced by ``$E$ is perfect''.

\begin{lemma}
Suppose $H=\Ex(\Q)$ and $A=(\der-2x)(\der^2+1)$. Then with $\Univ:=K[\ex^{H\imag}]$ we have
$\ker_{\Univ} A = \mathbb C\ex^{-x\imag} \oplus \mathbb C\ex^{x\imag}$.
\end{lemma}

\noindent
The proof is similar to that of Corollary~\ref{cor:no full basis}, using $\der^2+1=(\der-\imag)(\der+\imag)$ in $K[\der]$ and the fact that $y''+y\ne \ex^{x^2}$ for all $y\in K$.

\begin{remark} 
Let $H$ and $A$ be as in the previous lemma. There is an $H$-hardian germ~$y$ with   $y\neq 0$
and $A(y)=0$, but by the lemma, no such $y$ is in $H$. Thus Theorem~1 in~\cite{Ramayyan94} is false.
\end{remark}

\noindent
If $H$ is $\d$-maximal and $A$ has exactly one eigenvalue, then this eigenvalue is~$0$ by 
Corollary~\ref{cor:unique eigenvalue, 2}. 
This situation will be investigated in the next subsection.

\subsection*{Non-oscillation and disconjugacy}
{\it In this subsection we continue to as\-sume that~$A\in H[\der]$.}\/
In light of Corollary~\ref{cor:lindiff, real} one may ask  whether every non-oscillating~$y\in \ker_{\c^{<\infty}}A$  is  
$H$-hardian. The answer is ``no'' for some $A$: Suppose~$y$ in $H$ satisfies~$y''+y=x$. (If $H$ is $\d$-maximal, then $H$ is linearly surjective and such $y$ exists.)  Then~$y\succ 1$, and $y+\sin x\in \ker_{\c^{<\infty}}A$ is non-oscillating, but not $H$-hardian.

Here is a better question: if
$y\in \ker_{\Calinf}A$ and $y-h$ is non-oscillating for all~$h\in H$, does it follow that $y$ is $H$-hardian? The next two results
shows that the answer may depend on $H$. The first is a consequence of Corollary~\ref{cor:lindiff nonosc}. 

\begin{lemma}
Suppose $H$ is $\d$-maximal. Then  every
$y\in\ker_{\Calinf}A$ such that~$y-h$ is non-oscillating for all $h\in H$ lies in $H$.
\end{lemma}

\begin{lemma}
Let $H:=\Ex(\Q)$. Then there is a monic $A\in H[\der]$ of order $5$ and a~$y\in\ker_{\Calinf}A$  such that
 $y-h$ is non-oscillating for all $h\in H$, but $y$ is not hardian. 
\end{lemma}
\begin{proof}
Recall that each $\d$-maximal Hardy field contains an $f$ with~$f''+f=\ex^{x^2}$, by Theorem~\ref{thm:char d-max}.
% (Corollary~\ref{cor:d-max weakly d-closed}.)
Take an $H$-hardian $z\in \Calinf$ with  $z''+z=\ex^{x^2}$. Then~${z-h\succ x^n}$ for all $n$ and all $h\in H$, by \cite[Proposition~3.7 and Theorem~3.9]{Boshernitzan87}. Set
$$B_1\ :=\ \der^3-2x\der^2+\der-2x\in H[\der],\quad  B_2\ :=\ \der^2+4 \in H[\der]. $$
Then $B_1(z)=B_2(\sin 2x)=0$, hence $y:=z+\sin 2x\in\Calinf$ satisfies $A(y)=0$ for some monic $A\in H[\der]$ of order~$5$, by [ADH, 5.1.39].
For all $h\in H$ we have~$y-h=(z-h)+\sin 2x\sim z-h$, so~$y-h$ is non-oscillating.
Towards a contradiction, assume~$y$ is hardian. Then~$y$ is $H$-hardian.
Take an $H\langle y\rangle$-hardian $u\in \Calinf$ such that $u''+u=\ex^{x^2}$. Then~$(u-z)''+(u-z)=0$,  so~$u-z=c\cos(x+d)$, with $c,d\in\R$.
But~$u-y=c\cos(x+d)-\sin 2x$ is not hardian: this is clear if $c=0$, and otherwise
follows from $B_2(u-y)=3c \cos (x+d)$. This is the desired contradiction.
\end{proof}

\noindent
We can also ask: if {\it no}\/ $y\in \ker_{\Calinf}A$ oscillates, does it follow that $\ker_{\Calinf} A$ is contained in some Hardy field extension of $H$? We now extend Corollary~\ref{cor:char osc, 2} to give a positive answer:

\begin{theorem}\label{thm:split over D(H)}
The following are equivalent:
\begin{enumerate}
\item[\textup{(i)}] no $y\in\ker_{\Calinf} A$ oscillates;
\item[\textup{(ii)}] $\ker_{\Calinf} A\subseteq\Dx(H)$;
\item[\textup{(iii)}] $A$ splits over $\Dx(H)$;
\item[\textup{(iv)}] $A$ splits over some Hardy field extension of $H$.
\end{enumerate}
\end{theorem}
\begin{proof}
Corollary~\ref{cor:lindiff, real} gives (i)~$\Rightarrow$~(ii).  
For (ii)~$\Rightarrow$~(iii) use that  $A$ splits over~$\Dx(H)$ whenever $\dim_{\R} \ker_{\Dx(H)}A=r$, by Corollary~\ref{cor:lindiff, real}, and
Corollary~\ref{corbasissplit} with the remark following it.  The implication (iii)~$\Rightarrow$~(iv) is obvious.
Suppose (iv) holds; to show~(i),   arrange that $A$ splits over $H$ and 
$H$ is Liouville closed.
Then $\ker_{\Calinf} A$ is contained in $H$ by Lem\-ma~\ref{lem:basis of kerUA, real, H-field},  so (i) holds.
\end{proof}

{\samepage
\begin{remark} 
The implication (i)~$\Rightarrow$~(ii) in Theorem~\ref{thm:split over D(H)} is
 also claimed in \cite[Theorem~1]{RT97}; but the proof given there
is defective: in the proof of the auxiliary~\cite[Lem\-ma~1]{RT97} it is assumed that
if  $y\in\Calinf$ is non-oscillating and~$A(y)=0$, $y\neq 0$, then~$y^\dagger$ is also non-oscillating; but
$A=\der^3+\der$, $y=2+\sin x$ contradicts this. 
\end{remark}}

%Let $B\in\Calinf[\der]$ be monic; then $B$ is a product of monic operators of order~$1$ in~$\Calinf$ iff $\ker_{\Calinf} B\subseteq\c^\times$. (For a proof see \cite[Chapter~3, \S{}3, Theorem~2 and Lemma~6]{Coppel}.) Hence each of the statements  (i)--(iv) of Theorem~\ref{thm:split over D(H)} is equivalent to ``$A$ is a product of monic operators of order~$1$ in~$\Calinf$''.

\noindent
We say that {\bf $A$ does not generate oscillations} if it satisfies one of the equivalent conditions in the theorem above.
Thus if $r\leq 1$, then $A$ does not generate oscillations, and by Corollary~\ref{cor:char osc},  the operator~$\der^2+g\der+h$ 
($g,h\in H$) generates oscillations iff the germ $-\frac{1}{2}g'-\frac{1}{4}g^2+h$ generates
oscillations in the sense of Section~\ref{sec:second-order}.  
The property of $A$ to not generate oscillations is  uniformly  definable in the canonical $\HLO$-expansion $\boldsymbol H$ of $H$
viewed as  a structure in the language $\mathcal L^\iota_{\HLO}$  from [ADH, Chapter~16]   (see also the proof of Theorem~\ref{thm:transfer}); more precisely:

\begin{cor}
There is a quantifier-free   $\mathcal L^\iota_{\HLO}$-formula $\omega_r(x_1,\dots,x_r)$ such that for every Hardy field $H$ and
all $(h_1,\dots,h_r)\in H^r$:
$$\boldsymbol H\models\omega_r(h_1,\dots,h_r)\quad\Longleftrightarrow\quad \begin{cases}&\parbox{15em}{$\der^r+h_1\der^{r-1}+\cdots+h_r\in H[\der]$
does not generate oscillations.}\end{cases}$$
\end{cor}
\begin{proof}
Note that $A$ does not generate oscillations iff $A$ splits over some $\d$-maximal Hardy field extension of $H$. Now use that  the $\mathcal L^\iota_{\HLO}$-theory of canonical $\HLO$-ex\-pan\-sions of  $\d$-maximal Hardy fields  eliminates quantifiers,
by [ADH, 16.0.1] and
Theo\-rem~\ref{thm:char d-max}.
\end{proof}

\begin{example}
For $\omega_2(x_1,x_2)$ we may take the $\mathcal L^\iota_{\HLO}$-formula $\Upo(-2x_1'-x_1^2+4x_2)$.
Let $\alpha,\beta\in\R$, and let $\upo_n$ be as in Corollary~\ref{cor:17.7 generalized}. Then for $H=\R$ in that corollary, 
\begin{align*}\text{$\der^2+\alpha\der+\beta$ does not generate oscillations} &\ \Longleftrightarrow\ 
\text{$-\alpha^2+4\beta  < \upo_n$ for some $n$} \\
&\ \Longleftrightarrow\ \alpha^2-4\beta\ge 0, 
\end{align*} 
and applying the corollary to $H=\R(x)$ gives
\begin{align*}
\left.\parbox{16em}{$\der^2+\alpha x^{-1}\der+\beta x^{-2}$ does not generate oscillations}\,\right\} &\ \Longleftrightarrow\
\text{$(2\alpha-\alpha^2+4\beta)x^{-2} < \upo_n$ for some $n$} \\
&\ \Longleftrightarrow\ (1-\alpha)^2-4\beta\geq 0,
\end{align*}
in accordance with Corollary~\ref{cor:sa}. 
(By the way, $y''+\alpha x^{-1}y+\beta x^{-2}y=0$ is  
Euler's differential equation of order~$2$, cf.~\cite[\S{}22.3]{Kamke}, \cite[\S{}20, V]{Walter}.)
%\begin{align*} \marginpar{commented out example checked but replaced by more relevant one} 
%\text{$\der^2+\alpha\der+\beta x$ does not generate oscillations} &\ \Longleftrightarrow\ 
%\text{$-\alpha^2+4\beta x < \upo_n$ for some $n$} \\
%&\ \Longleftrightarrow\ \beta\leq 0,
%\end{align*}
%in accordance with Corollary~\ref{cor:sa}.
\end{example}

\noindent
From Corollary~\ref{cor:unique eigenvalue, 1} we obtain:

\begin{cor}\label{cor:dngo mult}
Suppose  $H$ is $\d$-perfect. Then:
$$A \text{ does not generate oscillations}\ \Longleftrightarrow\  \mult_0(A)=r.$$
If $H$ is moreover $\d$-maximal, then:
$$A \text{ does not generate oscillations}\ \Longleftrightarrow\ A \text{ has no eigenvalue different from }0.$$
\end{cor}

\noindent
We say that $B\in H[\der]^{\neq}$ does not generate oscillations if~$bB$ does not generate oscillations, for $b\in H^\times$ such that~$bB$ is monic.\index{generates oscillations!linear differential operator}\index{linear differential operator!generates oscillations} Using~[ADH, 5.1.22] we obtain:

\begin{cor}
Let $B_1,B_2\in H[\der]^{\neq}$. Then $B_1$ and $B_2$ do not generate oscillations iff $B_1B_2$ does not generate oscillations. 
\end{cor}

\noindent
Note also that if $E$ is a Hardy field extension of $H$, then 
$A$ generates oscillations iff $A$ generates oscillations when viewed as element of $E[\der]$. Moreover,
$A$ generates oscillations iff its adjoint $A^*$   does.
%In Corollary~\ref{cor:lclm dngo} below we shall also see that the least common left multiple of
%operators in $H[\der]$ which do not generate oscillations also does not generate oscillations.

In the next corollary $\phi$ ranges over elements of $\Dx(H)^>$ that are active in~$\Dx(H)$.

\begin{cor}\label{cor:split over D(H)}
Suppose $A$ does not generate oscillations.
Then the $\R$-linear space $\ker_{\Calinf}A$ has a basis $y_1,\dots,y_r$ with all $y_j\in\Dx(H)$ and
$y_1\prec\cdots\prec y_r$,
and there is a unique splitting $(a_r,\dots,a_1)$ of $A$ over $\Dx(H)$ 
such that eventually we have~$a_j+\phi^\dagger<a_{j+1}$ for $j=1,\dots,r-1$.
\end{cor}

\begin{proof}
By Theorem~\ref{thm:split over D(H)}, $A$ splits over $\Dx(H)$. Now use Lemma~\ref{lem:basis of kerUA, real, H-field} and Corollary~\ref{cor:distinguished splitting} 
%Proposition~\ref{prop:distinguished splitting} 
applied to the Liouville closed $H$-field $\Dx(H)$ in place of $H$.
\end{proof}

\noindent
 Theorem~\ref{thm:split over D(H)} and Corollary~\ref{cor:split over D(H)} yield Corollary~\ref{cor:unique factorization} from the introduction. The next lemma complements this Corollary~\ref{cor:unique factorization} by taking a look at splittings over the differential ring extension $R:=\Calinf$ of $H$:

\begin{lemma}\label{lem:fact in Calinf}
Suppose $H\supseteq \R$ is Liouville closed, $A$ splits over~$H$, $a_1,\dots,a_r$ lie in $R$, and  $A=(\der-a_r)\cdots(\der-a_1)$ in $R[\der]$. Then $a_1,\dots,a_r\in H$.
\end{lemma}
\begin{proof}
By induction on $r$. The case $r=0$ being trivial, suppose $r\geq 1$. By Lemma~\ref{lem:basis of kerUA, real} we have $\ker_R A\subseteq H$.
Take $y\in R^\times$ with $y^\dagger=a_1$. Then~${A(y)=0}$, hence $y\in H^\times$, so $a_1=y^\dagger\in H$.
Set $B:=(\der-a_r)\cdots(\der-a_2)\in R[\der]$, so~$A=B(\der-a_1)$.
Now $A_{\ltimes y}\in H[\der]$ and
$A_{\ltimes y}=B_{\ltimes y}\der$, so $B_{\ltimes y}\in H[\der]$.
Thus~${B\in H[\der]}$, and~$B$ splits over $H$ by [ADH, 5.1.22]. Hence, inductively, $a_2,\dots,a_r\in H$.
\end{proof}

\begin{cor}\label{lem:inhomog eq in D(H)}
Suppose $A$ does not generate oscillations, and let $b\in H$. Then all $y\in\Calinf$ with $A(y)=b$ are in $\Dx(H)$.  
\end{cor}
\begin{proof}
This follows from Corollary~\ref{cor:split over D(H)} and variation of constants [ADH, 5.5.22], using that $\Dx(H)$ is
closed under integration.
\end{proof}

\noindent
As pro\-mised in the remarks following Corollary~\ref{cor:disconjugacy} from the Introduction we now  strengthen the Trench normal form of
disconjugate linear differential operators. (See  Section~\ref{sec:second-order}, just before Lemma~\ref{lem:disconj germs}, for  ``disconjugate'' in the present context.) Below we use for $h\in \mathcal{C}$ the suggestive notation
$\int^{\infty} h =\infty$ to indicate that for some $a\in \R$ and representative $h\in \mathcal{C}_a$ (and thus for every $a\in \R$ and every representative $h\in \mathcal{C}_a$) we have $\int_a^t h(t)dt\to +\infty$ as $t\to +\infty$.

%{\samepage
%\begin{cor}\label{cor:Trench}
%The linear differential operator $A$ does not generate oscillations iff it is disconjugate.
%Moreover, suppose $A$ is disconjugate. Then there are~$g_1,\dots,g_r\in\Dx(H)$  
%such that 
%\begin{equation}\label{eq:Trench, 1}
%A=g_1\cdots g_r \der g_r^{-1} \der \cdots \der g_2^{-1} \der g_1^{-1}\quad\text{and}\quad g_j\succeq 1\ (j=2,\dots,r),
%\end{equation}
%and if $h_1,\dots,h_r\in (\Calinf)^\times$ satisfy
%\begin{equation}\label{eq:Trench, 2}
%A=h_1\cdots h_r \der h_r^{-1} \der \cdots \der h_2^{-1} \der h_1^{-1}\quad\text{and}\quad \int^\infty \abs{h_j}=\infty\ (j=2,\dots,r),
%\end{equation}
%then $g_j\in \R^> h_j$ for $j=1,\dots,r$.
%\end{cor}}
\begin{cor}\label{cor:Trench} Let $r\ge 1$. Then
$A$ does not generate oscillations iff $A$ is disconjugate.
Suppose $A$ is disconjugate. Then there are~$g_1,\dots,g_r\in \Dx(H)^{>}$   with
\begin{equation}\label{eq:Trench, 1}
A\ =\ g_1\cdots g_r (\der g_r^{-1})  \cdots (\der g_2^{-1})( \der g_1^{-1}), \quad \text{$g_j\in  \Upg\big(\!\Dx(H)\big)$ for $j=2,\dots,r$}.
\end{equation}
Given any such $g_1,\dots, g_r$, if  $h_1,\dots,h_r\in (\Calinf)^\times$  satisfy
\begin{equation}\label{eq:Trench, 2}
A\ =\ h_1\cdots h_r (\der h_r^{-1}) \cdots (\der h_2^{-1})( \der h_1^{-1}),\quad \text{$\int^\infty h_j\ =\ \infty$ for $j=2,\dots,r$},
\end{equation}
then $g_j\in \R^> \cdot h_j$ for $j=1,\dots,r$.
\end{cor}

\begin{proof} 
If $A$ does not generate oscillations, then $A$ is disconjugate by Lemma~\ref{lem:disconj germs} and Theorem~\ref{thm:split over D(H)}. The converse is clear.
Now suppose $A$ is disconjugate. 
Then Proposition~\ref{prop:distinguished splitting, Trench} yields $g_1,\dots,g_r\in\Dx(H)^{>}$
such that \eqref{eq:Trench, 1} holds. Let $h_1,\dots,h_r\in (\Calinf)^\times$ be such that 
\eqref{eq:Trench, 2} holds, and set $a_j:=(h_1\cdots h_j)^\dagger\in\Calinf$ ($j=1,\dots,r$).  Then~$A=(\der-a_r)\cdots(\der-a_1)$, so $a_1,\dots,a_r\in\Dx(H)$
by Lemma~\ref{lem:fact in Calinf}, hence~$h_1,\dots,h_r\in\Dx(H)$ as well. 
Thus $h_1,\dots, h_r>0$ and $h_2,\dots,h_r\in\Upg\big(\Dx(H)\big)$. The uniqueness part of Proposition~\ref{prop:distinguished splitting, Trench}
gives $g_j\in \R^> \cdot h_j$ for~$j=1,\dots,r$.
\end{proof}

\noindent
This yields in particular Corollary~\ref{cor:disconjugacy} from the Introduction.

\begin{example}[{Trench \cite[p.~321]{Trench}}]
Suppose $H=\R$, $r=3$, and
$A=\der^3-\der$. Then $A$ splits as  $(\der-1)\der(\der+1)$ over $H$, so $A$ does not generate oscillations.
%, and we have the splitting $A=of $A$ as in Corollary~\ref{cor:split over D(H)}.
In $\Dx(H)[\der]$,
$$A\ =\ \ex^x \der \ex^{-2x} \der \ex^x \der\ =\ 
\ex^{-x} \der \ex^x \der \ex^x \der \ex^{-x}\ =\ \ex^x \der \ex^{-x}\der \ex^{-x} \der \ex^x,$$
where only the last factorization is as in \eqref{eq:Trench, 1}.
\end{example}

\noindent
Generating oscillations is an invariant of the type of $A$:

\begin{lemma}
Suppose $A$ does not generate oscillations and  $B\in H[\der]$ has the same type as $A$. Then $B$ also does not generate oscillations.
\end{lemma}
\begin{proof}
By [ADH, 5.1.19]:  $r=\order(B)$ and we have $R\in H[\der]$ of order~$<r$ such that~$1\in H[\der]R+R[\der]A$ and $BR\in H[\der]A$.
Now $\ker_{\Calinf} A=\ker_{\Dx(H)} A$, and  [ADH, 5.1.20] gives an isomorphism
$y\mapsto R(y)\colon \ker_{\Dx(H)}A\to\ker_{\Dx(H)} B$ of $\R$-linear spaces. Hence $\ker_{\Calinf} B=\ker_{\Dx(H)} B$, 
so $B$ does not generate oscillations.
\end{proof}

\noindent
Next, let $M$ be a differential module over $H$. Recall from Section~\ref{sec:splitting} the notions of~$M$ splitting, and of
$M$ splitting over a given differential field extension of $H$.

\begin{lemma}
The following are equivalent:
\begin{enumerate}
\item[\textup{(i)}]  $M$ splits over  some Hardy field extension   of $H$;
\item[\textup{(ii)}] $M$ splits over $\Dx(H)$;
\item[\textup{(iii)}] $M$ splits over $\Ex(H)$.
\end{enumerate}
\end{lemma}
\begin{proof}
Let $E$ be a Hardy field extension of $H$
such that  $M$ splits over $E$. To show that  $M$  splits over $\Dx(H)$, 
replace $E$ by $\Dx(E)$ to arrange $E=\Dx(E)$. 
Next,
using~$E\otimes_H M\cong E\otimes_{\Dx(H)} (\Dx(H)\otimes_H M)$ and replacing~$H$,~$M$ by $\Dx(H)$, $\Dx(H)\otimes_H M$, respectively, also arrange~$H=\Dx(H)$.
In particular, $H\not\subseteq\R$, so $M\cong H[\der]/H[\der]B$ for some monic~$B\in H[\der]$, by~[ADH, 5.5.3].
Then $B$ splits over~$E$ by [ADH, 5.9.2] and hence over $H$ (Theorem~\ref{thm:split over D(H)}), so $M$ splits.
This shows (i)~$\Rightarrow$~(ii). 
The implications~(ii)~$\Rightarrow$~(iii)~$\Rightarrow$~(i) are obvious. 
\end{proof}

\noindent
We define:  {\bf $M$ does not generate oscillations} if it satisfies one of the equivalent conditions in the previous lemma.\index{differential module!generates oscillations}\index{generates oscillations!differential module} 
If~$M=H[\der]/H[\der]A$, then $M$ does not generate oscillations iff~$A$ does not generate oscillations.

\begin{cor}\label{cor:dngo base change}
Let $E$ be a Hardy field extension of $H$. Then   $M$  does not generate oscillations iff
its base change~$E\otimes_H M$ to $E$ does not generate oscillations.
\end{cor}
\begin{proof}
Use  
$\Dx(E)\otimes_H M \cong 
\Dx(E)\otimes_E\big(E\otimes_H M)$
as differential modules over $\Dx(E)$.
\end{proof}

\noindent
If $N$ is a differential submodule of $M$, then $M$ does not generate oscillations iff $N$ and $M/N$ do not generate oscillations.
 Hence:

\begin{cor}\label{cor:lclm dngo}
Let $A_1,\dots,A_m\in H[\der]^{\neq}$, $m\geq 1$. Then $A_1,\dots,A_m$ do not generate oscillations iff
$\operatorname{lclm}(A_1,\dots,A_m)\in H[\der]$ does not generate oscillations.
\end{cor}

\noindent
Let now $N$ be an $n\times n$ matrix over $H$, where $n\geq 1$.
We also say that the matrix differential equation~$y'=Ny$ over~$H$ {\bf does not generate oscillations} if 
the differential module over $H$ associated to~$N$ [ADH, p.~277] does  not generate oscillations.\index{matrix differential equation!generates oscillations}\index{generates oscillations!matrix differential equation}
If a matrix differential equation  over $H$    generates oscillations, then 
so does every  equivalent  matrix differential equation over~$H$.
Moreover, given a Hardy field extension $E$ of~$H$, the matrix differential equation $y'=Ny$ over~$H$
does not generate oscillations iff~$y'=Ny$ viewed as matrix differential equation over~$E$ does not generate oscillations (by Corollary~\ref{cor:dngo base change}).

\begin{lemma}\label{lem:dngo}
Suppose $B\in H[\der]$ is monic and $N$ is the companion
matrix of $B$. Then
$y'=Ny$ does not generate oscillations iff $B$ does not generate oscillations.
For each Hardy field extension $E$  of $H(\R)$ we have an isomorphism
$$y\mapsto (y,y',\dots,y^{(n-1)})^{\operatorname{t}} \colon
\ker_{E}(B) \to \operatorname{sol}_{E}(N)$$ of $\R$-linear spaces.
\end{lemma}
\begin{proof}
For the first claim, use that $M_N\cong H[\der]/H[\der]B^*$ by [ADH, 5.5.8], and $B$ does not generate oscillations iff $B^*$ does
not generate oscillations (remark before Corollary~\ref{cor:split over D(H)}). For the second claim, see [ADH, pp.~271--272].
\end{proof}

\begin{cor}
Suppose $H$ is $\d$-perfect. Then $y'=Ny$ does not generate oscillations iff $\operatorname{mult}_0(N)=n$.
\end{cor}
\begin{proof}
Using [ADH, 5.5.9], arrange that $N$ is the companion matrix of the mo\-nic~${B\in H[\der]}$.
Then $\operatorname{mult}_0(B)=\operatorname{mult}_0(N)$
by Lemma~\ref{lem:matrix diff equs vs ops}. Now use   Corollary~\ref{cor:dngo mult} and Lemma~\ref{lem:dngo}.
\end{proof}

{\samepage
\begin{prop}\label{prop:dngo}
The following are equivalent:
\begin{enumerate}
\item[\textup{(i)}] $y'=Ny$ does not generate oscillations;
\item[\textup{(ii)}] $\Dx(H)$ contains a fundamental matrix of solutions for~$y'=Ny$;
\item[\textup{(iii)}] $\Ex(H)$ contains a fundamental matrix of solutions for~$y'=Ny$;
\item[\textup{(iv)}] every maximal Hardy field extension of $H$ contains a fundamental matrix of solutions for~$y'=Ny$;
\item[\textup{(v)}] some Hardy field extension of $H$ contains a fundamental matrix of solutions for~$y'=Ny$.
\end{enumerate}
\end{prop}}
\begin{proof}
Suppose $y'=Ny$ does not generate oscillations.
Then $y'=Ny$ viewed as matrix differential equation over $\Dx(H)$ does not generate oscillations.
Hence to show (ii) we may arrange that $H=\Dx(H)$. Then $H\not\subseteq\R$, so  by [ADH, 5.5.9] we arrange that $N$ is the companion
matrix of the monic~$B\in H[\der]$.
Then $B$ does not generate oscillations, so  $H$ contains a fundamental matrix of solutions for~$y'=Ny$, by Theorem~\ref{thm:split over D(H)} and Lemma~\ref{lem:dngo}.
This proves~(ii).
The implications~(ii)~$\Rightarrow$~(iii)~$\Rightarrow$~(iv)~$\Rightarrow$~(v) are clear.
Suppose (v) holds. To prove (i), we first arrange that $H$ is $\d$-perfect and contains a fundamental matrix 
of solutions for~$y'=Ny$, and as in the proof of
(i)~$\Rightarrow$~(ii) we then arrange that $N$ is the companion matrix of some monic operator in $H[\der]$.
Then $y'=Ny$ does not generate oscillations, by  Theorem~\ref{thm:split over D(H)} and Lemma~\ref{lem:dngo}.
\end{proof}

\noindent
In \cite[Definition~16.14]{Boshernitzan82}, Boshernitzan defines   $y'=Ny$ 
to be {\it $H$-regular}\/ if it satisfies condition~(iii) in the proposition above. In~\cite[Theorem~16.16]{Boshernitzan82} he then notes the following version of Corollary~\ref{lem:inhomog eq in D(H)}, with $\Ex(H)$ in place of $\Dx(H)$:\index{matrix differential equation!$H$-regular}

\begin{cor}
Suppose the matrix differential equation $y'=Ny$ does not generate oscillations. Let   $b\in H^n$ be a column vector.
Then each solution~$y$ in~$(\Calinf)^n$ to the differential equation~$y'=Ny+b$ lies in $\Dx(H)^n$.
\end{cor}
\begin{proof}
By  Proposition~\ref{prop:dngo}, $\Dx(H)$ contains a fundamental matrix of solutions for~$y'=Ny$. Now use
[ADH, 5.5.21] and $\Dx(H)$ being closed under integration.
\end{proof}

\noindent
Here is an application of the material above to the parametrization of curves in euclidean $n$-space, where
for simplicity we only treat the case $n=3$, denoting the usual euclidean norm on $\R^3$ by $|\cdot|$.

\begin{example}[Frenet-Serret formulas]
Let $U\subseteq \R$ be a nonempty open interval  and $\gamma\colon U\to\R^3$ be a $\Ginf$-curve, parametrized by arc length, that  is, $|\gamma'(t)|=1$ for all~${t\in U}$.
Let $T:=\gamma'$   and
$\kappa:=\abs{T'}$ (the {\it curvature}\/ of $\gamma$).
Suppose $\kappa(t)\neq 0$ for each~$t\in U$, set $N:=T'/\abs{T'}$ and $B:=T\times N$ (vector cross product).
Then for~$t\in U$ the vectors~$T(t),N(t),B(t)\in\R^3$ are orthonormal 
% \marginpar{restof example still to check} 
and~$y=(T,N,B)\colon U\to \R^9$ 
is a solution of the matrix differential equation~${y'=Fy}$ in $\mathcal C^\infty(U)$ where
$$F=\begin{pmatrix} 
  			& \kappa I	&   	  \\
-\kappa I	&  			& \tau I  \\
 			& -\tau I	&  
\end{pmatrix}\qquad\text{($I=$ the $3\times 3$ identity matrix),}$$
for some $\Ginf$-function $\tau\colon U\to\R$ (the {\it torsion}\/ of $\gamma$).

Conversely, let $\Ginf$-func\-tions~$\kappa,\tau\colon U\to\R$ such that~$\kappa(t)>0$ for all $t\in U$ be given. Then
there is a   $\Ginf$-curve~$\gamma=(\gamma_1,\gamma_2,\gamma_3)\colon U\to\R^3$, parametrized by arc length, with
curvature~$\kappa$ and torsion $\tau$. In fact, $\gamma$ is unique up to proper euclidean motions in $\R^3$. (See~\cite[Chap\-ter~1]{Spivak} for these facts.)  Fix such $\gamma$ and assume in addition that $U=(c,+\infty)$ with $c\in \R\cup\{-\infty\}$ and that $H$ is a $\d$-maximal $\Ginf$-Hardy field and contains the germs of~$\kappa$,~$\tau$, also  denoted by~$\kappa$,~$\tau$.
Then for some $\alpha\in K/K^\dagger$, the matrix differential equation~$y'=Fy$ over~$K$  
has spectrum~$\{\alpha,-\alpha,0\}$. (Example~\ref{ex:Frenet-Serret}.)
If~$\alpha=0$, then~$y'=Fy$ does not generate oscillations, hence the germs of~$\gamma_1$,~$\gamma_2$,~$\gamma_3$ lie in~$H$ and so~$\gamma$ does not exhibit oscillating behavior. If~$\alpha\neq 0$, then~$\alpha= \phi'\imag + K^\dagger$ where~$\phi\in H$, $\phi>\R$, and then by Corollaries~\ref{cor:Shackell} and~\ref{cor:shape of fund matrix, real}, the 
 germs of~$\gamma_1$,~$\gamma_2$,~$\gamma_3$ lie in $$H\cos\phi+H\sin\phi+H\subseteq\Ginf.$$
 For example, if $\kappa\in\R^>$ and~$\tau\in\R$ are constant, then~$\gamma$ is the helix given by
 $$t\mapsto\big({-a\cos(t/D)},a\sin(t/D),bt/D\big)$$ where~$D=\sqrt{a^2+b^2}$, $\kappa=a/D^2$, $\tau=b/D^2$. \end{example}

\section{Revisiting Second-Order Linear Differential Equations}\label{sec:perfect applications}

\noindent
In this section we analyze the oscillating solutions of second-order linear differential equations over Hardy fields in more detail.
In particular, we prove  Corollary~\ref{cor:Boshernitzan intro} from the introduction. This is connected to the $\upo$-freeness of the perfect hull of a Hardy field, which is characterized in Theorem~\ref{thm:upo-freeness of the perfect hull}. % We finish  with some remarks and conjectures about firm and flabby dents. 
{\it Throughout this section $H$ is a Hardy field and $K:=H[\imag]\subseteq \Calinf[\imag]$.}\/

\subsection*{Parametrizing the solution space}
Let $a,b\in H$. We now continue the study of  the linear differential equation
\begin{equation}\label{eq:2nd order, app}\tag{$\tilde{\operatorname{L}}$}
Y''+aY'+bY\ =\ 0
\end{equation}
over $H$ from Section~\ref{sec:Hardy fields} (with slightly changed notation), and focus on the oscillating case, viewed in the light of our main theorem. (Corollaries~\ref{cor:char osc} and~\ref{cor:char osc, 2} already dealt with the non-oscillating case, which didn't need our main result.)
Most of the following theorem was claimed without proof by Boshernitzan~\cite[Theorem~5.4]{Boshernitzan87}:

\begin{theorem} \label{thm:Bosh}  
Suppose  \eqref{eq:2nd order, app} has an oscillating solution \textup{(}in $\Calinf$\textup{)}. Then there are $H$-hardian germs
$g>0$, $\phi>\R$  such that for all~$y\in\Calinf$,
$$\text{$y$ is a solution of    \eqref{eq:2nd order, app}} \quad\Longleftrightarrow\quad  \text{$y=cg\cos(\phi+d)$ for some $c,d\in\R$.}$$
Any such $H$-hardian germs~$g$,~$\phi$ are $\d$-algebraic over $H$ and lie in a common Hardy field extension of $H$. 
If~$\Dx(H)$ is $\upo$-free, then these properties force $g,\phi\in\Dx(H)$, and determine~$g$   uniquely  up to multiplication by a positive real number and $\phi$  uniquely up to  
addition of a real number.
\end{theorem}

\begin{remarks}
If $H$ is $\upo$-free, then  $\Dx(H)$ is $\upo$-free by Theorem~\ref{thm:ADH 13.6.1}.
Also, if~$H$ is not $\upl$-free or  $\overline{\omega}(H)=H\setminus\sigma\big(\Upg(H)\big){}^\uparrow$,
then $\Dx(H)$ is $\upo$-free, by  Lemma~\ref{lem:D(H) upo-free, 4}.
(See Section~\ref{sec:order 2 Hardy fields} or [ADH, 5.2] 
for the definition of the function~$\sigma$, and recall that $\bar{\omega}(H)$ is the set of all $f\in H$ such that $f/4$ does not generate oscillations. If $H$ is $\upo$-free, then  $\overline{\omega}(H)=H\setminus\sigma\big(\Upg(H)\big){}^\uparrow$, by Corollary~\ref{omuplosc}.) 
%Theorem~\ref{thm:upo-freeness of the perfect hull} below yields the converse of this implication. 
Recall also that $\upl$-freeness includes having asymptotic integration.
In the last sentence of Theorem~\ref{thm:Bosh}  we cannot drop the
hypothesis that $\Dx(H)$ is $\upo$-free; see Remark~\ref{rem:non-uniqueness}.
\end{remarks}

\noindent
 Let $V$ be an $\R$-linear subspace  of $\c$.  A pair  $(g,\phi)$ is said to {\bf parametrize $V$} if  $$g\in \c^\times,\ g>0, \quad \phi\in \c,\ \phi> \R,\qquad
V\ =\  \big\{ cg\cos(\phi+d): c,d\in\R \big\};$$ equivalently, $g\in \c^\times$, $g>0$, $\phi\in \c$, $\phi> \R$, and 
$V= \R g\cos\phi+\R g\sin\phi$,
by Corol\-lary~\ref{cor:sinusoids}.  If $(g,\phi)$ parametrizes~$V$, then so does $(cg,\phi+d)$ for any~$c\in\R^>$,~$d\in\R$.

\begin{example}
The example following Corollary~\ref{cor:gen osc closed upward} shows that for $f\in \R^>$ the pair~$(1,\frac{\sqrt{f}}{2} x)$ parametrizes $\ker_{\Calinf} (4\der^2+f)$. 
\end{example}

\noindent
Suppose~$V=\ker_{\Calinf} (\der^2+a\der+b)$, and let $g\in \c^{\times}$, $g>0$, and~$\phi \in \c$, $\phi>\R$. Then~$(g,\phi)$
parametrizes $V$  iff $g\ex^{\phi\imag}\in\ker_{\Calinf[\imag]}(\der^2+a\der+b)$. 

%\medskip
%\noindent
%(the commented out material that follows has been checked but is no longer needed)
%For later use we record a consequence of Theorem~\ref{thm:lindiff d-max}.
%(Recall from Section~\ref{sec:splitting}  that $\Sigma(A)\subseteq K/K^\dagger$ denotes the set of eigenvalues of  $A\in %K[\der]^{\neq}$.)
 
% \begin{cor}\label{paralpha} If $H$ is maximal and $(g,\phi)\in H^2$ parametrizes $\ker_{\Calinf} (\der^2+a\der +b)$,
% $then $\Sigma(\der^2+a\der + b)=\{-\alpha,\alpha\}$ for $\alpha:= \phi'\imag+K^\dagger$, and $\phi'\imag\notin K^\dagger$.  
 %\end{cor}
 
 \medskip\noindent
For later use we record the next lemma where $V$ is an $\R$-linear subspace of $\c^1$ and~$V':=\{y':y\in V\}$ (an $\R$-linear subspace of $\c$).

\begin{lemma}\label{lem:param V'}
Suppose $H\supseteq\R$ is real closed and closed under integration, and $(g,\phi)\in H\times H$ parametrizes $V$. Set 
$q:=\sqrt{(g')^2+(g\phi')^2}$ and
$u:=\arccos(g'/q)$. Then $q, u\in H$ and $(q, \phi+u )\in H\times H$ parametrizes $V'$.
\end{lemma}
\begin{proof}
Note that $u$  is as in Corollary~\ref{arccosH} with $g'$, $-g\phi'$ in place of~$g$,~$h$.
Let~$y\in V$, so $y=cg\cos(\phi+d)$ where $c,d\in\R$. Then
$$y'\ =\ cg'\cos(\phi+d)-cg\phi'\sin(\phi+d)\ =\ cq\cos(\phi+u+d).$$
Conversely, for $c,d\in \R$ we have $cq\cos(\phi+u+d)=y'$ for $y=cg\cos(\phi+d)\in V$.
\end{proof}

\begin{lemma}\label{parphi} Set $f:= -2a'-a^2+4b$. Let $h$ be an $H$-hardian germ such that~$h>0$ and 
$h^\dagger=-\frac{1}{2}a$. 
Let $g\in \c^\times$, $g>0$ and $\phi\in \c$, $\phi>\R$. Then: 
\begin{enumerate}
\item[\rm(i)] $(g,\phi)$ parametrizes $\ker_{\Calinf}4\der^2+f$ iff $(gh,\phi)$ parametrizes $\ker_{\Calinf}\der^2+a\der +b$.
\end{enumerate}
Assume also that $\phi$ is hardian $($so $\phi' $ is hardian with  $\phi'>0)$. Then: \begin{enumerate}
\item[\rm(ii)]  $(1/\sqrt{\phi'}, \phi)$ parametrizes $\ker_{\Calinf}4\der^2+\sigma(2\phi')$.
\end{enumerate}
\end{lemma} 
\begin{proof}  The arguments leading up to Corollary~\ref{cor:char osc} yield (i). 
As to (ii), the definition of $\sigma$ in [ADH, p. 262] gives
$$ \sigma(2\phi')\ =\ \omega\big({-(2\phi')^\dagger} + 2\phi'\imag\big)\ =\   \omega\big({-\phi'^\dagger}+2\phi'\imag\big)\ =\ 
\omega(2y^\dagger)$$
 where $y:= (1/\sqrt{\phi'})\ex^{\phi\imag}$.  Hence $A(y)=0$ for $A=4\der^2+\sigma(2\phi')$ by the computation in [ADH, p. 258], and thus $\big(1/\sqrt{\phi'}, \phi\big)$ parametrizes $\ker_{\Calinf} A$.
\end{proof} 

\noindent 
Item (i) in Lemma~\ref{parphi} reduces the proof of Theorem~\ref{thm:Bosh} to the case $a=0$, and~(ii) is about
that case. Next we isolate an argument in the proof of [ADH, 14.2.18]:

\begin{lemma}\label{lem:ADH 14.2.18}
Let $E$ be a $2$-newtonian $H$-asymptotic field with asymptotic integration, $e\in E^\times$, $f\in E$, and $\upg$ be active in $E$ such that $e^2=f-\sigma(\upg)$ and $e\succ\upg$.
Then $\sigma(y)=f$ and~$y\sim e$ for some $y\in E^\times$.
\end{lemma}
\begin{proof}
Note that $e$ is active in $E$ since $e\succ\upg$.
By [ADH, 11.7.6] we have
$$\sigma(e)-f\ =\ \sigma(e)-\sigma(\upg)-e^2\  =\  \omega(-e^\dagger) - \omega(-\upg^\dagger) - \upg^2\ \prec\  e^2,$$
and so $\omega(-e^\dagger)-f\sim-e^2$. Eventually $\phi\prec e$, so $(\phi/e)^\dagger\prec e$ by [ADH, 9.2.11].
Hence with $R$, $Q$ as defined before [ADH, 14.2.18], eventually we have $R^\phi\prec e^2$, and thus~$Q^\phi \sim e^2 Y^2 (Y^2-1)$. Now [ADH, 14.2.12] yields $u\in E$ with $u\sim 1$ and~$Q(u)=0$, thus~$\sigma(y)=f$ for $y:=eu\sim e$.
\end{proof}

\noindent
Let  $A=4\der^2+f\in H[\der]$  where~$f/4\in H$ generates oscillations, and set $V:=\ker_{\Calinf} A$. If $H$ is $\upo$-free, then~$f\in \sigma\big(\Upg(H)\big){}^\uparrow$, by Corollary~\ref{omuplosc}. 
Theorem~\ref{thm:Bosh} now follows from Lemmas~\ref{lem:parametrization of ker A}, \ref{lem:asymptotics of phi'}, and \ref{lem:g,phi unique} below, which give more information.

\begin{lemma}\label{lem:parametrization of ker A}
There is a pair of  $H$-hardian germs para\-me\-tri\-zing $V$. For any such pair $(g,\phi)$ we have~$\sigma(2\phi')=f$ and $g^2\phi'\in \R^>$,  so~$g$,~$\phi$ are $\d$-algebraic over $\Q\langle f\rangle$ and lie in a common Hardy field extension of $H$.
If $f\in\Ginf$, then each pair of  $H$-hardian germs  para\-me\-tri\-zing $V$ is in $(\Ginf)^2$;
likewise with $\Gom$ in place of $\Ginf$.
\end{lemma}

\begin{proof}
The first statement follows from Corollaries~\ref{cor:cos sin infinitesimal},~\ref{cor:d-max weakly d-closed}, and~\ref{cor:2nd order kernel}.
Next, let~$(g,\phi)$ be a pair of  $H$-hardian germs parametrizing~$V$. 
Set~$y:=g\ex^{\phi\imag}\in\Calinf[\imag]^\times$; then we have~$A(y)=0$ 
and hence~$\omega(2y^\dagger)=f$ where $y^\dagger=g^\dagger+\phi'\imag\in \Calinf[\imag]$. 
Now for $p,q\in \c^1$ we have $\omega(p+q\imag)=\omega(p)+q^2-2(pq+q')\imag$, so $$\omega(p+q\imag)\in \c\ 
\Leftrightarrow\ pq+q'=0.$$ Therefore~$2g^\dagger=-(2\phi')^\dagger=-(\phi')^\dagger$ and so $g^2\phi'\in \R^>$,
and~$\sigma(2\phi')=f$~[ADH, p.~262]. 
If $f\in\Ginf$, then $y\in\Ginf[\imag]$, so $g^2=\abs{y}^2\in\Ginf$ and hence also $\phi\in\Ginf$ since~$g^2\phi'\in \R^>$;
likewise with $\Gom$ in place of $\Ginf$.
 \end{proof}

\begin{lemma}\label{lem:asymptotics of phi'}
Suppose that $H\supseteq\R$ is real closed with asymptotic integration,  
and that~$f\in\sigma\big(\Upg(H)\big){}^\uparrow$.
Then there is an active $e>0$ in $H$ such that $\phi'\sim e$ for all pairs~$(g,\phi)$ of $H$-hardian germs  parametrizing $V$.
\end{lemma}
\begin{proof}
Choose a  logarithmic sequence $(\ell_\rho)$ for $H$ and set~$\upg_\rho:=\ell_\rho^\dagger$~[ADH, 11.5];
then~$(\upg_\rho)$ is strictly decreasing and coinitial in $\Upg(H)$~[ADH, p.~528].
%We have $f\notin\omega(H)$ since $f$   generates oscillations, hence $f\in \sigma\big(\Upg(H)\big){}^\uparrow$ since~$H$ is $\upo$-free [ADH, 11.8.30].
Take $\rho$ such that~$f>\sigma(\upg_\rho)$. As in the proof of [ADH, 14.2.18], we increase~$\rho$ so that $f-\sigma(\upg_\rho)\succ\upg_\rho^2$, and take $e\in H^>$ with $e^2=f-\sigma(\upg_\rho)$. 
Then~$e\succ\upg_\rho$ and so~$e\in \Upg(H)^\uparrow$.
Let~$(g,\phi)$ be a pair of elements in a Hardy field extension $E$ of~$H$ parametrizing $V$.
We claim that  $\phi'\sim e/2$ (so $e/2$ in place of $e$ has the property desired in the lemma). We arrange that $E$ is  $\d$-maximal. Then~$E$ is Liouville closed and~$\phi>\R$, so
$e,2\phi'\in\Upg(E)$ by~[ADH, 11.8.19].
Now~$E$ is newtonian by Theorem~\ref{thm:char d-max}, so  Lemma~\ref{lem:ADH 14.2.18} yields~$u \sim 1$ in~$E$ such that~$\sigma(eu)=f$. Now the map~$y\mapsto \sigma(y)\colon\Upg(E)\to E$ is strictly increasing
[ADH, 11.8.29], hence $2\phi'=eu$ by Lemma~\ref{lem:parametrization of ker A}, and thus~$\phi'\sim e/2$.
\end{proof}

\begin{lemma}\label{lem:g,phi unique}
Suppose $\Dx(H)$ is $\upo$-free or $f\in\sigma\big(\Upg(H)\big){}^\uparrow$.
Let $H_i$ be a Hardy field extension of $H$ with  $(g_i,\phi_i)\in H_i\times H_i$ parametrizing $V$, for $i=1,2$. 
Then
$${g_1/g_2}\in\R^>, \qquad {\phi_1-\phi_2}\in\R.$$
Thus $g,\phi\in \Dx(H)$ for any pair $(g,\phi)$ of  $H$-hardian germs parametrizing $V$.  
\end{lemma}
\begin{proof}
We arrange that $H_1$, $H_2$ are $\d$-maximal and thus contain $\Dx(H)$. Replacing~$H$ by $\Dx(H)$  we further arrange that $H$ is $\d$-perfect and  $f\in \sigma\big(\Upg(H)\big){}^\uparrow$.
Then $\phi_1'\sim\phi_2'$ by Lemma~\ref{lem:asymptotics of phi'},
and for $i=1,2$ we have $c_i\in\R^>$ with~$\phi_i' = c_i/g_i^2$, by Lemma~\ref{lem:parametrization of ker A}. 
Replacing~$g_i$ by $g_i/\sqrt{c_i}$ we arrange $c_i=1$ ($i=1,2$), so~$g_1\sim g_2$.
Consider now the elements $g_1\cos\phi_1$, $g_1\sin\phi_1$ of~$V$;
take~$a,b,c,d\in\R$ such that
$$g_1\cos \phi_1\  =\  ag_2\cos(\phi_2+b),\qquad g_1\sin \phi_1\  =\  cg_2\cos(\phi_2+d).$$
Then  
\begin{equation}\label{eq:g1g2}
g_1^2\ =\ g_1^2(\cos^2 \phi_1 + \sin^2 \phi_1)\  =\  g_2^2\big( a^2\cos^2(\phi_2+b) + c^2\cos^2(\phi_2+d) \big),
\end{equation}
and hence
$$ a^2\cos^2(\phi_2+b) + c^2\cos^2(\phi_2+d) \ \sim\  1.$$
Thus the $2\pi$-periodic function
$$t\mapsto F(t)\ :=\ a^2\cos^2(t+b)+c^2\cos^2(t+d)\ \colon\ \R\to\R$$ 
satisfies $F(t)\to 1$ as $t\to+\infty$, hence $F(t)=1$ for all $t$, so $g_1=g_2$ by \eqref{eq:g1g2}.
It follows that~$\phi_1'=\phi_2'$, so~$\phi_1-\phi_2\in\R$.

For the final claim, let $(g,\phi)$ be a pair of $H$-hardian germs parametrizing $V$. Let~$M$ be any $\d$-maximal extension of $H$. Then Lemma~\ref{lem:parametrization of ker A} gives a pair $(g_M,\phi_M)\in M^2$ that also parametrizes $V$. By the above, $g/g_M\in \R^{>}$ and $\phi-\phi_M\in \R$, hence~$g, \phi\in M$. Since $M$ is arbitrary, this gives $g,\phi\in \Dx(H)$. 
\end{proof}

\noindent
This finishes the proof of Theorem~\ref{thm:Bosh} (and Corollary~\ref{cor:Boshernitzan intro} from the introduction).   %\qed. 

\begin{cor}\label{cor:perfect Schwarz closed, 1}
Suppose that $H$ is $\d$-perfect. Then $\omega(H)=\bar{\omega}(H)$ is downward closed and $\sigma\big(\Upg(H)\big)$ is upward closed.
\end{cor}
\begin{proof}
By Corollary~\ref{cor:omega(H) downward closed}, $\omega(H)=\bar{\omega}(H)$ is downward closed. 

Let $f\in\sigma\big(\Upg(H)\big){}^\uparrow$. The last part of Lemma~\ref{lem:g,phi unique} gives $g,\phi\in H$ such that $(g,\phi)$ parametrizes $V$. Then $\sigma(2\phi')=f$
by Lemma~\ref{lem:parametrization of ker A}. Now $2\phi'\in \Upg(H)$ by [ADH, 11.8.19], so $f$ lies in
$\sigma\big(\Upg(H)\big)$. Thus $\sigma\big(\Upg(H)\big)$ is upward closed.
%this yields  and~\ref{lem:g,phi unique}, $\sigma\big(\Upg(H)\big)$ is upward closed.
\end{proof}

\noindent
Recall from [ADH, 11.8] that $H\supseteq\R$ is said to be {\it Schwarz closed}\/ if $H$ is Liouville closed and
$H=\omega\big(\Upl(H)\big)\cup\sigma\big(\Upg(H)\big)$.\index{Hardy field!Schwarz closed}  

\begin{cor}  \label{cor:perfect Schwarz closed, 2}
Suppose $H$ is $\d$-perfect. Then the following are equivalent:
\begin{enumerate}
\item[\textup{(i)}] $H$ is Schwarz closed;
\item[\textup{(ii)}] $H$ is $\upo$-free;
\item[\textup{(iii)}] for all $f\in H$ the operator $4\der^2+f\in H[\der]$ splits over $K$;
\item[\textup{(iv)}] for all $a,b\in H$ the operator $\der^2+a\der+b\in H[\der]$ splits over $K$.
\end{enumerate}
\end{cor}
\begin{proof}
The equivalence (iii)~$\Leftrightarrow$~(iv) holds by Corollary~\ref{cor:char osc, 3}, and
the equi\-va\-len\-ces~(i)~$\Leftrightarrow$~(ii)~$\Leftrightarrow$~(iii) follow from [ADH, 11.8.33] and Corollary~\ref{cor:perfect Schwarz closed, 1}.
\end{proof}

\noindent
{\em In the rest of this subsection
$A=\der^2+a\der+b$ \textup{(}$a,b\in H$\textup{)}. We set $V:=\ker_{\Calinf} A$ and
$f:=-2a'-a^2+4b$, and we take $H$-hardian $h>0$ such that $h^\dagger=-\frac{1}{2}a$}. Note the relevance of Lemma~\ref{parphi}(i) in this situation. 

\begin{cor}\label{cor:2nd order, f succ 1/x^2}
Suppose $f>0, f\succ 1/x^2$. Then $f\notin\overline\omega(H)$, and  for some $H$-hardian germ $\phi$ with $\phi'\sim \frac{1}{2}\sqrt{f}$, and~$g:=1/\sqrt{\phi'}$ we have: $(gh,\phi)$ pa\-ra\-me\-tri\-zes~$V$. 
\end{cor}
\begin{proof}
By Theorem~\ref{upo} we arrange that $H\supseteq\R$ is Liouville closed and $\upo$-free. With  notation as at the beginning of Section~\ref{sec:upo-free Hardy fields} we have
$\upo_\rho\sim 1/x^2$ for all~$\rho$; hence~$f/4>\upo_\rho$ for all $\rho$, so $f/4$ generates oscillations by [ADH, 11.8.21] and Corollary~\ref{omuplosc}, and $f\notin \overline\omega(H)$, $f\in \sigma\big(\Upg(H)\big){}^\uparrow$.
Lemma~\ref{lem:parametrization of ker A} gives a pair $(g,\phi)$ parametrizing $\ker_{\Calinf} (4\der^2+f)$ with
$H$-hardian $\phi$ and $g:=1/\sqrt{ \phi'}$.
Now~$\upg:=1/x$ is active in $H$ with $\sigma(\upg)=2\upg^2$ and so $f>\sigma(\upg)$ and  $f-\sigma(\upg)\sim f$. 
%Set $b:=\sqrt{f-\sigma(\upg)}\in H^>$; then $b\sim f^{1/2}\succ\upg$. 
Then $\phi'\sim \frac{1}{2}\sqrt{f}$  by the proof of Lemma~\ref{lem:asymptotics of phi'}, so $\phi$ has the property stated in Corollary~\ref{cor:2nd order, f succ 1/x^2}. 
\end{proof}

\begin{cor}\label{cor:2nd order, phi prec x}
Suppose $f\notin\overline\omega(H)$ and let $(g,\phi)$ be a pair of $H$-hardian germs parametrizing $V$.
Then $\phi\prec x$ iff $f\prec 1$, and the same with $\preceq$ in place of $\prec$.
Also, if $f\sim c\in\R^>$, then $\phi\sim \frac{\sqrt{c}}{2}x$ and $(f<c\Rightarrow
\phi''>0)$, $(f>c\Rightarrow\phi''<0)$.
\end{cor}
\begin{proof}
We arrange $H\supseteq\R$ is $\upo$-free, Liouville closed, and $g,\phi\in H$.
Then $y:=2\phi'\in\Upg(H)$ by [ADH, 11.8.19]. Lemma~\ref{lem:parametrization of ker A} gives $\sigma(y)=f$; also
$\sigma(c)=c^2$ for all~$c\in\R^>$.
As the restriction of $\sigma$ to $\Upg(H)$ is strictly increasing~[ADH, 11.8.29],  this yields the first part. 
Now suppose  $f\sim c\in\R^>$, and take $\lambda\in\R^>$ with $\phi\sim\lambda x$.
Then $y \sim 2\lambda$, and with~$z:=-y^\dagger\prec 1$ we have
$f=\sigma(y)=\omega(z)+y^2\sim 4\lambda^2$.
Hence~$\lambda=\sqrt{c}/2$.
Suppose $f<c$;
then~$f=\sigma(y)<c=\sigma(\sqrt{c})$ yields $y<\sqrt{c}$, so~$\phi'<\lambda$.
With~$g:=\phi-\lambda x$ we have~$g\prec x$, $g'\prec 1$ 
 and   $g'=\phi'-\lambda<0$, so~$g''=\phi''>0$. The case~$f>c$ is similar.
\end{proof}

\noindent
Combining Corollaries~\ref{lem:parametrization of ker A} and \ref{cor:2nd order, phi prec x} yields:

\begin{cor}
Suppose $f\notin\overline\omega(H)$. Then for every~$y\in V^{\neq}$ we have:
\begin{enumerate}
\item[\textup{(i)}] if $f\prec 1$, then $y\not\preccurlyeq h$ \textup{(}so $y$ is unbounded if in addition $a\leq 0$\textup{)};
\item[\textup{(ii)}] if $f\succ 1$, then $y\prec h$; and
\item[\textup{(iii)}] if $f\asymp 1$, then $y\preceq h$.
\end{enumerate}
\end{cor}

\begin{remarks}
See \cite[Chapter~6]{Bellman} for related (though generally weaker) results in a more general setting. For example, if
$g\in\c$  is eventually   increasing with~$g\succ 1$
or~$g\in\c^1$ and $g\sim 1$ with $\int |g'|\preceq 1$,  then every $y\in\c^2$ with $y''+gy=0$ satisfies~$y\preceq 1$;
cf.~\S\S6,~18 in loc.~cit. 
%{\bf these remarks not checked}
\end{remarks}

{\samepage

\begin{cor}\label{cor:2nd order, Liouvillian solutions} 
Suppose $f\notin\overline\omega(H)$, $H\supseteq\R$, and $H$ does not have asymptotic integration or~$H$ is $\upo$-free. Then the following are equivalent: 
\begin{enumerate}
\item[\textup{(i)}] $A(y)=0$ for some $y\neq 0$ in  a Liouville extension of $K$;
\item[\textup{(ii)}]  some pair $(g,\phi)\in \operatorname{Li}(H)^2$ with $g^\dagger$, $\phi'$ algebraic over~$H$ parametrizes $V$;\item[\textup{(iii)}]  some pair in $\operatorname{Li}(H)^2$ parametrizes $V$;
\item[\textup{(iv)}] every pair of $H$-hardian germs pa\-ra\-me\-trizing $V$ lies in $\operatorname{Li}(H)^2$.
\end{enumerate}
\end{cor}}
\begin{proof}
Suppose (i) holds. 
Then Lemma~\ref{lem:Liouvillian zeros} gives $g, \phi\in L:=\operatorname{Li}(H)$, $g\neq 0$, such that
$g^\dagger$, $\phi'$ are algebraic over $H$ and $A(g\ex^{\phi\imag})=0$.
Replacing~$g$ by~$-g$ if necessary we arrange $g>0$.
We have $\phi\succ 1$: otherwise $g\ex^{\phi\imag}\in E[\imag]^\times$ for some Hardy field extension $E$ of $L$, by Corollary~\ref{cor:phi preceq 1}, hence $\Re(g\ex^{\phi\imag})\in V^{\neq}$ does not oscillate, or $\Im(g\ex^{\phi\imag})\in V^{\neq}$ does not oscillate,
a contradiction. Replacing $\phi$ by $-\phi$ if necessary we arrange $\phi>\R$. Then $(g,\phi)$ parametrizes $V$.
This yields~(ii).  The implication~(ii)~$\Rightarrow$~(iii) is trivial.
By the assumptions on $H$, $\operatorname{Li}(H)$, and thus $\Dx(H)$, is $\upo$-free,  so (iii)~$\Rightarrow$~(iv) follows from Theorem~\ref{thm:Bosh}. For~(iv)~$\Rightarrow$~(i), note that the differential fraction field of
$K[\ex^{H\imag}]\subseteq \Calinf[\imag]$ is a Liouville extension of $K$.  
\end{proof}

\begin{cor}\label{cor:parametrization in H}
Suppose $H\supseteq\R$ is Liouville closed and $f\notin\bar{\omega}(H)$. Then the following are equivalent:
\begin{enumerate}
\item[\textup{(i)}]  $g,\phi\in H$ for every pair
$(g,\phi)$ of $H$-hardian germs parametrizing $V$;
\item[\textup{(ii)}] there is a pair of germs in $H$  parametrizing $V$;
\item[\textup{(iii)}]   $f\in\sigma(H^\times)$.
\end{enumerate}
\end{cor}
\begin{proof}
The implications (i)~$\Rightarrow$~(ii)~$\Rightarrow$~(iii) follow from Lemma~\ref{lem:parametrization of ker A} and
the remarks preceding Lemma~\ref{parphi}.
%(and don't need the hypotheses on $H$, $K$).
Suppose $f\in\sigma(H^\times)$.
Since $f\notin\bar\omega(H)$ and~$\omega(H)^\downarrow\subseteq\bar\omega(H)$,
we have $f\notin\omega(H)^\downarrow$, so
  $f\in\sigma\big(\Upg(H)\big)$ by [ADH, 11.8.31].
Also, $4\der^2+f$ splits over $K$ but not over $H$ (cf.~[ADH, pp.~259, 262]) and $f/4$ generates oscillations. Hence   Corollary~\ref{cor:2nd order kernel} and the remark following it yield a pair of germs in $H$ parametrizing $V$. Now (i) follows from Lemma~\ref{lem:g,phi unique}.
\end{proof}

%\begin{cor}\label{cor:D(H) Schwarz closed}
%Suppose $H$ does not have asymptotic integration or $H$ is $\upo$-free. Then $\Dx(H)$ is Schwarz closed.
%\end{cor}
%\begin{proof}
%By the previous corollary it is enough to show that $\Dx(H)$ is $\upo$-free. If $\Dx(H)$ is $\upo$-free, then this follows from  [ADH, 13.6.1]. Otherwise, first replace $H$ by $H(\R)$ and use [ADH, ] to arrange $H\supseteq\R$; then $\operatorname{Li}(H)$ is $\upo$-free by Lemma~\ref{}, hence $\Dx(H)=\Dx(\operatorname{Li}(H))$ is $\upo$-free by the first case.
%\end{proof}

%\noindent
%If $H$ is bounded, then $\Ex(H)=\Dx(H)$ (see the remarks following Lemma~\ref{lem:Dx Ex}).
%Hence if  $H$ is bounded, and   $H$ does not have asymptotic integration or   $H$ is $\upo$-free, then $\Ex(H)$ is Schwarz closed,
%by   Corollary~\ref{cor:D(H) Schwarz closed}. In \dots we show how to remove the hypothesis of boundedness in this statement. \marginpar{to fill in}

\noindent
The case of Theorem~\ref{thm:Bosh} where $a$, $b$ are $\d$-algebraic over $\Q$ is used later. In that case the $\Psi$-set
of the Hardy subfield $H_0:=\Q\langle a,b\rangle$ of $H$ is finite by Lemma~\ref{lemro83}, so $H_0$ has no asymptotic integration. Thus the relevance of the next result:
\medskip

\begin{cor}~\label{naicor}
Suppose $f\notin\bar\omega(H)$ and $H$ has no asymptotic integration. Then there is a pair~$(g,\phi)\in \Dx(H)^2$ parametrizing $V$
such that every pair of $H$-hardian germs parametrizing $V$ equals $(cg, \phi+d)$ for some $c\in \R^{>}$ and $d\in \R$. 
\end{cor}
\begin{proof} The assumption on $H$ gives that $\Dx(H)$ is $\upo$-free. Now use Theorem~\ref{thm:Bosh}. 
\end{proof}

\noindent
{\em In the rest of this subsection $f\notin \bar\omega(H)$, and  $(g,\phi)$ is a pair of $H$-hardian germs parametrizing $V$}.
Then   $\sigma(2\phi')=f$   (cf.~Lemma~\ref{lem:parametrization of ker A}) and  thus $$P(2\phi')\ =\ 0\quad\text{where $P(Y)\ :=\ 2YY''-3(Y')^2+Y^4-fY^2\in H\{Y\}$.}$$ 
Hence Theorem~\ref{thm:Ros83} applied to $E:=H\<\phi\>=H(\phi, \phi', \phi'')$ gives for grounded $H$ elements
%by Corollary~\ref{cor:Ros83}  there are 
$h_0,h_1\in H^{>}$ and $m$, $n$ with  $h_0, h_1\succ 1$ and $m+n\leq 3$, such that
$$\log_{m+1} h_0\ \prec\ \phi\ \preceq\ \exp_n h_1.$$ 
In the next two lemmas we improve on these bounds:

\begin{lemma}\label{boundphi1}
Suppose  $\ell_0\in H^>$, $\ell_0\succ 1$, and $\max\Psi_{H}=v(\upg_0)$ for $\upg_0:=\ell_0^\dagger$. Then
 $f-\omega(-\upg_0^\dagger) \succ \upg_0^2$ and $\phi\succ \log\ell_0$, or
   $f-\omega(-\upg_0^\dagger) \asymp \upg_0^2$ and $\phi\asymp \log\ell_0$.
\end{lemma}
\begin{proof}
By Lemma~\ref{lem:baromega(H) for H without as int}   we have
$f\notin\bar\omega(H) = \omega(-\upg_0^\dagger)+\upg_0^2\smallo_{H}^{\downarrow}$
and hence $f= \omega(-\upg_0^\dagger)+\upg_0^2 u$ where $u\succeq 1$, $u>0$,
so $f-\sigma(\upg_0)=\upg_0^2(u-1)$. Suppose $u\succ 1$; then~$f-\sigma(\upg_0)\sim u\upg_0^2\succ\upg_0^2$, 
and   the proof of Lemma~\ref{lem:asymptotics of phi'} shows that then $\phi'\sim e/2$ where~$e^2=f-\sigma(\upg_0)$, so
$e\succ \upg_0$ and thus $\phi\succ \log\ell_0$. 
Now suppose $u\asymp 1$, and put~$\ell_1:=\log \ell_0$,   $\upg_1:=\ell_1^\dagger$.
Then by    [ADH, 11.7.6], 
$$f-\sigma(\upg_1)\  =\  \omega\big({-\upg_0^\dagger}\big)-\omega\big({-\upg_1^\dagger}\big)+u\upg_0^2-\upg_1^2\ \sim\  u\upg_0^2\ \succ\  \upg_1^2,$$
and arguing as in the proof of Lemma~\ref{lem:asymptotics of phi'} as before gives~$\phi\asymp\log\ell_0$.
\end{proof}

\begin{lemma}\label{boundphi2}
Suppose $f\in \sigma\big(\Upg(H)\big){}^\uparrow$ or $H$ is not $\upl$-free, and  
$u\in H^>$ is such that $u\succ 1$ and  $v(u^\dagger)=\min\Psi_{H}$. Then $\phi\leq u^n$ for some $n\geq 1$.
\end{lemma}
\begin{proof} We have $H$-hardian $\phi\succ 1$,  but this is not enough to get $\theta\in H^\times$ with~$\phi\asymp \theta$. 
That is why we consider first the case that $H\supseteq\R$  is real closed with asymptotic integration, and
$f\in \sigma\big(\Upg(H)\big){}^\uparrow$.
Then Lemma~\ref{lem:asymptotics of phi'} gives $e\in H^>$  such that~$\phi'\sim e$, and as $H$ has asymptotic integration
we obtain $\theta\in H^\times$ with $\phi\asymp \theta$. Hence~$\phi^\dagger\asymp \theta^\dagger \preceq u^\dagger$,
and   thus $\phi\leq u^n$ for some $n\geq 1$, by [ADH, 9.1.11]. 

We now reduce the general case to this special case. Take a $\d$-maximal Hardy field extension $E$ of $H$ with $g,\phi\in E$.
Suppose $H$ is $\upl$-free. Then $f\in \sigma\big(\Upg(H)\big){}^\uparrow$. Also, $H(\R)$ is $\upl$-free with the same value group as $H$, by Proposition~\ref{prop:const field ext}, so~$L:=H(\R)^{\operatorname{rc}}\subseteq E$ has asymptotic integration, with $v(u^\dagger)=\min\Psi_{L}$. Thus $\phi\leq u^n$ for some $n\geq 1$ by the special case applied to $L$ in the role of $H$.
 
For the rest of the proof we  assume $H$ is not $\upl$-free. 
 Then~$H(\R)$ is not $\upl$-free by Lemmas~\ref{notupl1} and~\ref{notupl2}, and  so $L:=H(\R)^{\operatorname{rc}}\subseteq E$ is not $\upl$-free by~[ADH, 11.6.8]. Using [ADH, 10.3.2] we also have
$v(u^\dagger)=\min\Psi_{L}$. Hence replacing~$H$ by $L$
% and noting that~$\Gamma_H^>$ is cofinal in $\Gamma_L^>$,
we arrange that~$H\supseteq\R$ and $H$  is real closed in what follows. 
%If $H$ has asymptotic integration and $f\in \sigma\big(\Upg(H)\big){}^\uparrow$, then we are done.   

Suppose $H$ has no asymptotic integration. As in the proof of Lemma~\ref{lem:Li(H) upo-free} this
yields an $\upo$-free Hardy subfield $L\supseteq H$ of $E$ such that $\Gamma_H^{>}$ is cofinal  in $\Gamma_L^{>}$, so~$v(u^\dagger)=\min \Psi_L$. Moreover, $f\in L\setminus\bar\omega(L)= \sigma\big(\Upg(L)\big){}^\uparrow$ by Corollary~\ref{omuplosc}.
Hence replacing $H$ by~$L^{\operatorname{rc}}$ we have a reduction to the special case. 

Suppose $H$ has asymptotic integration. 
%and~$f\notin \sigma\big(\Upg(H)\big){}^\uparrow$. 
Since $H$ is not $\upl$-free, [ADH, 11.6.1] gives~$s\in H$ creating a gap over $H$.
Take  $y\in E^\times$ with $y^\dagger=s$. Then $vy$ is a gap in $H(y)$ by the remark following [ADH, 11.5.14],
and thus a gap in $L:=H(y)^{\operatorname{rc}}$. Moreover, $\Gamma_H^{>}$ is cofinal in $\Gamma^{>}_{H(y)}$
by [ADH, 10.4.5](i), hence cofinal in $\Gamma_L^{>}$, so $v(u^\dagger)=\min \Psi_L$. Thus replacing
$H$ by $L$ yields a reduction to the ``no asymptotic integration'' case.  
\end{proof}

\begin{cor}\label{cor:upper lower bd for phi}
Let    $a,b\in H:=\R(x)$. Then $g,\phi\in\Dx(\Q)\subseteq\Gom$, and~$\phi \preceq x^n$ for some~$n\geq 1$.
Moreover, $f\succ 1/x^2$ and $\log x\prec\phi$, or $f \asymp 1/x^2$ and~$\log x\asymp\phi$. 
\end{cor}
\begin{proof} Apply Lemma~\ref{boundphi1} with $\ell_0:= x$, and Lemma~\ref{boundphi2} with $u=x$, and note that
$f\prec 1/x^2$ is excluded by the standing assumption $f\notin \bar\omega(H)$.
\end{proof}  

\begin{remark}
Suppose $a=0$. Then $b=f/4$, and $g^2\phi'\in\R^>$ by Lemma~\ref{lem:parametrization of ker A}; hence bounds on~$\phi$ 
give bounds on $g$.
Thus by Corollary~\ref{cor:upper lower bd for phi}, if $b\in H:=\R(x)$, then~$g\succeq x^{-n}$ for some $n\geq 1$,
and either~$f\succ 1/x^2$, $g\prec \sqrt{x}$, or $f\asymp 1/x^2$, $g\asymp\sqrt{x}$.
\end{remark}

\begin{examples} Let $H:=\R(x)$. Then for  $a=0$ and $b=\frac{5}{4}x^{-2}$ the standing assumption 
$f\notin\bar{\omega}(H)$ holds, since $f=5x^{-2}$. The germ $y=x^{1/2}\cos \log x\in\Gom$ solves the corresponding second-order linear differential equation $4Y''+ fY=0$.
Other example: let $H$ contain $x$ and $x^r$ where $r\in\R$, $r>-1$. Then for $a=0$ and~$b:= \frac{1}{4}\big(x^{2r}-r(r+2)x^{-2}\big)\in \Gom$ the standing assumption $f\notin \bar{\omega}(H)$ holds
in view of $f=4b\sim x^{2r}\succ 1/x^2$. Here  $z=x^{-r/2}\cos \left(\frac{x^{r+1}}{2(r+1)}\right)\in\Gom$ 
satisfies~$4z''+fz=0$.
\end{examples}

\noindent
We now set $B:=\der^3+f\der+(f'/2)\in H[\der]$, and observe:

\begin{lemma}\label{lem:B, 1}
$B(1/\phi')=0$.
\end{lemma}
\begin{proof}
We  arrange that $H\supseteq\R$ contains $\phi$ and is Liouville closed, and identify the universal  exponential extension 
$\Univ=\Univ_K$ of~$K=H[\imag]$   with a differential subring of~$\Calinf[\imag]$ as explained
at the beginning of Section~\ref{sec:ueeh}. Then
$$(\phi')^{-1/2}\ex^{\phi\imag}, (\phi')^{-1/2}\ex^{-\phi\imag}\in \ker_{\Univ}4\der^2+f.$$
Thus $B(1/\phi')=0$ by Lemma~\ref{lem:Appell} applied to $\operatorname{Frac}(\Univ)$ in the role of $K$.
\end{proof}

\noindent
For the canonical $\HLO$-expansion of a Hardy field, see Section~\ref{sec:transfer}.

\begin{lemma}\label{lem:B, 2}
Let $\mathbf E$ be a pre-$\HLO$-field extension of the canonical $\HLO$-exp\-an\-sion of $H\langle \phi'\rangle$.
Then $\ker_E B = C_E(1/\phi')$.
\end{lemma}
\begin{proof}
%Thanks to [ADH, remark after 4.1.13] we may replace $\mathbf E$ by any of its pre-$\HLO$-field extensions.
Using [ADH, 16.3.20, remark after 4.1.13] we arrange  $\mathbf E$ to be Schwarz closed.
Then~$f\notin\bar{\omega}(H)=\omega(E)\cap H$,
hence $f\in\sigma(E^\times)$, so $\dim_{C_E} \ker_E B=1$ by Lem\-ma~\ref{lem:kerB}.
\end{proof}

\noindent
We can now complement Corollary~\ref{cor:2nd order, Liouvillian solutions}: 

\begin{cor}\label{cor:phi' quadratic}
Suppose   $\phi'$ is algebraic over $H$. Then   $(\phi')^2\in H$ and $g^\dagger\in H$.
\end{cor}
\begin{proof}
Let $E:=H^{\operatorname{rc}}\subseteq\Calinf$. Then by Corollary~\ref{cor:canonical HLO},
the canonical $\HLO$-expansion  of~$E$ extends that of $H\langle \phi'\rangle$.
Set  $L:=E[\imag]\subseteq \Calinf[\imag]$, so $L$ is an algebraic closure of the differential field~$H$.
Put
$u:=2\phi'\in E$, and
let~$\tau\in\operatorname{Aut}(L|H)$.
Then~$B(\tau(1/u))=0$ by Lemma~\ref{lem:B, 1}. So $\Re\tau(1/u)$ and~$\Im\tau(1/u)$
in $E$ are also zeros of~$B$, hence
Lemma~\ref{lem:B, 2} yields $c\in\C^\times$ with~$\tau(1/u)=c/u$ and thus~$\tau(u)=c^{-1}u$.
Now with $$P(Y)\ :=\  2YY''-3(Y')^2+Y^4-fY^2\in H\{Y\}$$ we have $P(u)=0$, so $P(\tau(u))=0$, hence
$$0\ =\ P(u)-c^2P(\tau(u))\ =\ P(u)-c^2P(c^{-1}u)\ =\ (1-c^{-2})u^4$$
and  thus $c\in\{-1,1\}$, so $\tau(u^2)=u^2$. This proves the first statement.
The second statement follows from the first and $g^2\phi'\in\R^>$ by Lemma~\ref{lem:parametrization of ker A}.
\end{proof}

\subsection*{Distribution of zeros} 
Let   $a,b\in H$ and consider again the differential equation
%~\eqref{eq:2nd order, app}.
\begin{equation}%\label{eq:2nd order, app}
\tag{$\tilde{\operatorname{L}}$}
Y''+aY'+bY\ =\ 0.
\end{equation}
Below we use Theorem~\ref{thm:Bosh}
to show
that for any oscillating solution~$y\in\Calinf$ of \eqref{eq:2nd order, app}
the sequence of successive zeros of $y$ grows very regularly, with growth comparable to that of the sequence of successive relative
maxima of $y$, and also to that of a function whose germ is hardian. 
(For the equation $Y''+fY=0$,   where $f\in\Ex(\Q)$ generates oscillations, this was suggested
after~\cite[\S{}20, Conjecture~4]{Boshernitzan82}.)

To make this precise we first  
define a preordering~$\leq$ on the set  $\R^{\N}$ of sequences of real numbers  by
$$(s_n)\leq (t_n)\quad:\Longleftrightarrow\quad s_n \leq t_n \text{ eventually} 
\quad:\Longleftrightarrow\quad \exists m\,\forall n\geq m\  s_n\leq t_n.$$
(A preordering on a set is a reflexive and transitive  binary relation on that set.)
We say that $(s_n),(t_n)\in\R^{\N}$ are {\it comparable} if $(s_n) \leq (t_n)$ or $(t_n) \leq (s_n)$. 
The induced equivalence relation $\sim_{\ta}$ on $\R^{\N}$ is that of having the same tail: 
$$(s_n)\sim_{\ta}(t_n) \quad:\Longleftrightarrow\quad (s_n) \leq (t_n) \text{ and }(t_n)\leq (s_n)\ \Longleftrightarrow\ s_n=t_n \text{ eventually}. $$
To any germ $f\in \c$ we take a representative in $\c_0$, denoted here also by $f$ for convenience, and associate to this germ the tail of the sequence $\big(f(n)\big)$, noting that this tail is independent of the choice of representative. 

For example, if the germs of $f,g\in\c_0$ are contained in a common Hardy field, then the sequences
$\big(f(n)\big)$, $\big(g(n)\big)$ are comparable. 
%\marginpar{for ``limit point'', see 5.2.10} 
Given an infinite set $S\subseteq\R$ with a lower bound in $\R$ and without a limit point, the {\bf enumeration}\/ of $S$\index{enumeration}
is the strictly increasing sequence~$(s_n)$ with~$S=\{s_0,s_1,\dots\}$ (so $s_n\to+\infty$ as $n\to+\infty$).
 
We take representatives of $a$, $b$ in $\c^1_e$ with $e\in \R$, denoting these by $a$ and $b$ as well, 
and set~$f:=-2a'-a^2+4b\in \c_e$.
Let $y\in \c_e^2$ be oscillating with 
$$y''+ay'+by\ =\ 0,\qquad (\text{so the germ of $f$ does not lie in }\bar{\omega}(H)),$$
and let $(s_n)$ be the enumeration of $y^{-1}(0)$.
(See Lem\-ma~\ref{lem:no limit pt}.)
 Theorem~\ref{thm:Bosh} 
yields $e_0\geq e$, $g\in\c_{e_0}^\times$, and strictly increasing $\phi\in\c_{e_0}$ such that $y|_{e_0}=g\cos \phi$, and~$g$,~$\phi$ lie in a common Hardy field extension of $H$ with $(g,\phi)$ pa\-ra\-me\-tri\-zing $\ker_{\Calinf}(\der^2+a\der+b)$ (where $g$,~$\phi$ also
denote their own germs). 

 \begin{lemma}\label{lem:zeta}
There is a strictly increasing $\zeta\in\c_{n_0}$ \textup{(}$n_0\in\N$\textup{)} 
 such that  $s_n=\zeta(n)$ for all $n\geq n_0$ and the germ of $\zeta$ is hardian with $H$-hardian compositional inverse.
 \end{lemma}
\begin{proof}
Take~$n_0\in \N$ such that $s_n\geq e_0$ for all $n\geq n_0$, and then $k_0\in\frac{1}{2}+\Z$ such that~$\phi(s_n)=(k_0+n)\pi$ for all $n\geq n_0$. Thus $n_0=\big(\phi(s_{n_0})/\pi\big)-k_0$. Let
$\zeta \in\c_{n_0}$ be the compositional
inverse of~$(\phi/\pi)-k_0$ on $[s_{n_0},+\infty)$. Then $\zeta$ has the desired properties: the germ of $\zeta$ is hardian by Lemma~\ref{lem:Bosh6.5}.
\end{proof}

\noindent
If $a,b\in \Ginf$, then we can choose~$\zeta$ in Lemma~\ref{lem:zeta} such that its germ is in~$\Ginf$;
likewise with $\Gom$ in place of $\Ginf$. We do not know whether we can always choose~$\zeta$ in Lemma~\ref{lem:zeta} to have $H$-hardian germ. For  $\phi$ not growing too slowly we can describe the asymptotic behavior of~$\zeta$ in terms of $\phi$:

\begin{cor}\label{asympzeta}
If $\phi\succeq x^{1/n}$ for some $n\ge 1$, then 
in Lemma~\ref{lem:zeta} one can choose $\zeta\sim \phi^{\operatorname{inv}}\circ \pi x$.
\end{cor}
\begin{proof} Let $n_0$, $k_0$, $\zeta$ be as in the proof of Lemma~\ref{lem:zeta}. Then $$\zeta^{\operatorname{inv}}\ \sim\  (\phi/\pi)-k_0\ \sim\  \phi/\pi.$$  Now assume $\phi\succeq x^{1/n}$, $n\ge 1$. Then $\zeta^{\operatorname{inv}}\succeq x^{1/n}$, so $\zeta\preceq x^n$, and thus the condition stated just before Lemma~\ref{lem:Entr} is satisfied for
$h:=\zeta$. We can therefore use Corollary~\ref{cor:Entr, 2} with 
 $\phi/\pi$, $\zeta$ in the role of $g$, $h$ to give $\phi^{\operatorname{inv}}\circ \pi x\sim \zeta$.
\end{proof}

\noindent
Combining Corollaries~\ref{cor:2nd order, phi prec x},~\ref{asympzeta}, and~\ref{cor:Entr, 2} we obtain:

\begin{cor}\label{cor:zeta asymptotics}
If $f\sim c$ \textup{(}$c\in\R^>$\textup{)}, then $s_n\sim \frac{2}{\sqrt{c}}\pi n$ as $n\to\infty$.
\end{cor}

\noindent
Combining Corollary~\ref{cor:upper lower bd for phi} with the proof of Lemma~\ref{lem:zeta}
yields  crude bounds on the growth of $(s_n)$ when $H=\R(x)$:

\begin{cor}
Suppose $a,b\in\R(x)$.
If $f \asymp 1/x^2$, then for some $r\in\R^>$ we have $\ex^{n/r} \leq s_n\leq \ex^{rn}$ eventually.
If $f\succ 1/x^2$, then for some $m\geq 1$ and every~$\epsilon\in\R^>$ we have 
$n^{1/m}\leq s_n\leq \ex^{\epsilon n}$ eventually.
\end{cor}

\noindent
The next lemma is a version of the Sturm Convexity Theorem~\cite[p.~318]{BR} concerning the differences between consecutive zeros of $y$:

\begin{lemma}\label{lem:zeta, differences}
If $f\prec 1$, then 
the sequence $(s_{n+1}-s_n)$ is eventually strictly increasing
with
$s_{n+1}-s_n\to +\infty$ as $n\to\infty$. If
 $f\succ  1$, then $(s_{n+1}-s_n)$ is eventually strictly decreasing with~${s_{n+1}-s_n}\to 0$  as $n\to\infty$.
Now suppose~$f\sim c$~\textup{(}$c\in\R^>$\textup{)}. Then $s_{n+1}-s_n\to 2\pi/\sqrt{c}$ as~$n\to\infty$,
and if $f<c$, then~$(s_{n+1}-s_n)$ is eventually strictly decreasing,
if $f=c$, then~$(s_{n+1}-s_n)$ is eventually constant, and
if $f>c$, then $(s_{n+1}-s_n)$ is eventually strictly increasing.
%Moreover, if~$\phi''>0$, then the sequence $(s_{n+1}-s_n)$ is eventually strictly decreasing; if~$\phi''=0$, then $(s_{n+1}-s_n)$ is eventually constant;  and if $\phi''<0$, then~$(s_{n+1}-s_n)$ is eventually strictly increasing.
\end{lemma}
\begin{proof}
We arrange $\phi\in\c_{e_0}^2$ such that $\phi'(t)>0$ for all $t\ge e_0$. 
Take  $\zeta$ as in the proof of Lemma~\ref{lem:zeta}. Then~$\zeta\in\c_{n_0}^2$ with
$$\zeta'\ =\ \pi \frac{1}{\phi'\circ\zeta},\qquad \zeta''\ =\ -\pi^2\frac{\phi''\circ\zeta}{(\phi'\circ\zeta)^3}.$$
The Mean Value Theorem gives for every $n\geq n_0$ a $t_n\in (n,n+1)$ such that
$$s_{n+1}-s_n\ =\ \zeta(n+1)-\zeta(n)\ =\ \zeta'(t_n).$$
If $f\prec 1$, then Corollary~\ref{cor:2nd order, phi prec x} gives $\phi\prec x$, so $\zeta\succ x$, hence $\zeta'\succ 1$; this proves the first claim of the lemma. The other claims follow likewise using Corollary~\ref{cor:2nd order, phi prec x} and the above remarks on $\zeta'$ and $\zeta''$. 
\end{proof}

\begin{cor}\label{cor:zeta, 3}
Let $h\in\c_0$ and suppose the germ of $h$ is in $\Ex(\Q)$. Then the se\-quen\-ces~$(s_n)$ and $\big(h(n)\big)$ are comparable.
\end{cor}
\begin{proof}
Let $\zeta$ be as in Lemma~\ref{lem:zeta}, and note that the germs of $\zeta$ and $h$ lie in a common Hardy field.
% the compositional inverse $\zeta^{\operatorname{inv}}$ of the germ of $\zeta$
%is $H$-hardian, and by Lemma~\ref{lem:comp with E(Q)}, so is $g\circ \zeta^{\operatorname{inv}}$; hence $g\circ \zeta^{\operatorname{inv}}\leq x$ or~$g\circ \zeta^{\operatorname{inv}}\geq x$ and so $g\leq\zeta$ or $g\geq\zeta$.
\end{proof}

\noindent
Let also $\underline{a}, \underline{b}\in H$, and take representatives of $\underline{a}$, $\underline{b}$ in $\c_{\underline{e}}^1$ ($\underline{e}\in \R$), denoting these by~$\underline{a}$ and~$\underline{b}$ as well. Let $\underline{y}\in \c_{\underline{e}}^2$ be an oscillating solution of the differential equation$$Y''+\underline{a}Y' +\underline{b}Y\ =\ 0,$$
and let $(\underline{s}_n)$ be the enumeration of~$\underline{y}^{-1}(0)$.

\begin{lemma}\label{lem:comp sequ zeros}
The sequences $(s_n)$ and $(\underline{s}_n)$ are comparable.
%Tehere are some $a\in\R$ and some representatives of $y$, $z$ in $\c_a^2$ whose sequences of zeros are comparable.
%whenever $s_1<s_2<\cdots$ are the consecutive zeros of $y|_b$ and
%$t_1<t_2<\cdots$ are the consecutive zeros of $z|_b$, then $(s_n)$,
%, then there is exactly one zero $t\in [t_1,t_2]$ of $z$.
\end{lemma}
\begin{proof}
We arrange that $H$ is maximal and take $\zeta$ as in Lemma~\ref{lem:zeta}. This lemma also  provides a  strictly increasing $\underline{\zeta}\in\c_{\underline{n}_0}$ \textup{(}$\underline{n}_0\in\N$\textup{)} 
 such that  $\underline{s}_n=\underline{\zeta}(n)$ for all~$n\geq \underline{n}_0$ and the germ of $\underline{\zeta}$ is hardian with $H$-hardian compositional inverse. With $\zeta$ and $\underline{\zeta}$ denoting also their germs this gives
 $\zeta^{\operatorname{inv}}\le \underline{\zeta}^{\operatorname{inv}}$ or 
 $\zeta^{\operatorname{inv}}\ge \underline{\zeta}^{\operatorname{inv}}$, hence~$\zeta\ge \underline{\zeta}$ or $\zeta\le \underline{\zeta}$. Thus $(s_n)$ and $(\underline{s}_n)$ are comparable.
\end{proof}

\noindent
Now $y'$  also oscillates, so by Corollary~\ref{gphicos} there is for all sufficiently large $n$ exactly one $t\in (s_n,s_{n+1})$ with $y'(t)=0$.  Also $b\ne 0$ in $H$, since $b=0$ would mean that $z:=y'$  satisfies
$z'+az=0$, so $z$ would be $H$-hardian. 
This leads to the following:   Let $m\geq 1$ and suppose~$y\in\c_e^{m+2}$ (and~$y''+ay'+by=0$ with oscillating~$y$ as before).  Then  the zero sets of $y,y',\dots,y^{(m)}$ are eventually parametrized by hardian germs as follows:
%the germ of $y'$ also satisfies a second-order linear differential equation over $H$. (See the remarks before Lemma~\ref{lem:A^der}.)  Thus Lemma~\ref{lem:zeta} applies to the enumerations of the zero set of~$y'|_{a_1}$ in place of $(s_n)$.

\begin{lemma}\label{lem:param zeros of derivatives} For $i=0,\dots,m$
we have an $n_i\in \N$ and a strictly increasing function~$\zeta_i\in\c_{n_i}$, such that:
\begin{enumerate}
\item[\textup{(i)}]   $\zeta_i(n_i)\ge e$ and $\zeta_i(t)\to +\infty$ as $t\to +\infty$;
\item[\textup{(ii)}] the germ of $\zeta_i$ is hardian with $H$-hardian  compositional inverse; 
\item[\textup{(iii)}] $\big\{\zeta_i(n):\ n\ge n_i\big\}=\big\{t\geq \zeta_i(n_i)\ :\ y^{(i)}(t)=0\big\}$;
\item[\textup{(iv)}]  $\zeta_i^{\operatorname{inv}}-\zeta_{0}^{\operatorname{inv}}\preceq 1$;
\item[\textup{(v)}] if $i<m$, then $\zeta_i(n)<\zeta_{i+1}(n)<\zeta_i(n+1)$ for all $n\ge n_{i+1}$.
\end{enumerate}
\end{lemma} 
\begin{proof}
We arrange that $H$ is maximal. For simplicity we only do the case~${m=1}$; the general case just involves more notation.  
For $\zeta_0$ we take a function $\zeta$ as constructed in the proof of Lemma~\ref{lem:zeta}, and also take
$n_0$ as in that proof, so clauses~(i),~(ii),~(iii) are satisfied for $i=0$.  Set 
$A:=\der^2+a\der +b\in H[\der]$. As $b\ne 0$, we have the monic operator $A^{\der}\in H[\der]$ of order $2$ as defined before Lemma~\ref{lem:A^der}, with~$A^\der(y')=0$. Take  a pair $(g_1,\phi_1)$  of elements of $H$ parametrizing $\ker_{\Calinf} A^{\der}$. Then with $A^\der$, $g_1$, $\phi_1$ instead of $A$, $g$, $\phi$, and taking suitable representatives of the relevant germs, the proof of Lemma~\ref{lem:zeta} provides likewise an $n_1\in \N$, a $k_1\in \frac{1}{2}+\Z$, and a strictly increasing function 
$\zeta_1\in \c_{n_1}$ satisfying clauses (i), (ii), (iii) for $i=1$ and with compositional inverse given by $(\phi_1/\pi)-k_1$. 

Recall that $\Univ=\Univ_K\subseteq\Calinf[\imag]$ is a differential integral domain extending~$K$, and that 
$\ker_{\Calinf[\imag]} B=\ker_{\Univ} B$ for all $B\in K[\der]^{\ne}$, by Theorem~\ref{thm:lindiff d-max}. Therefore~$\ker_{\Calinf[\imag]} A^\der=\big\{y':y\in \ker_{\Calinf[\imag]}A\big\}$ by Lemma~\ref{lem:A^der}, so in view of $A\in H[\der]$,
$$\ker_{\Calinf} A^\der\ =\  \big\{y':y\in \ker_{\Calinf}A\big\}.$$ Now (iv) for $i=1$ follows from
Lemmas~\ref{lem:param V'} and~\ref{lem:g,phi unique}.

%As to (iv), recall that $\Univ=\Univ_K\subseteq\Calinf[\imag]$ is a differential integral domain extending~$K$, 
%and that  for all $B\in K[\der]^{\ne}$ we have $\ker_{\Calinf[\imag]} B=\ker_{\Univ} B$ by Theorem~\ref{thm:lindiff d-max}.  
%Set $\alpha:=\phi'\imag+K^\dagger$;
%then $\Sigma(A)=\{-\alpha,\alpha\}$ by Corollary~\ref{paralpha}, and
%$\Sigma(A^\der)=\Sigma(A)$ by Corollary~\ref{cor:A^der}. Moreover, $\alpha=\phi_1'\imag+ K^\dagger$:
%otherwise~$\alpha=-\phi_1'\imag+ K^\dagger$, so~$(\phi+\phi_1)'\imag\in K^\dagger$, which is false. 
%Thus $(\phi-\phi_1)'\imag\in K^\dagger$, which gives $\phi-\phi_1\preceq 1$. This proves (iv) for $i=1$.  

As to (v), the remark preceding the lemma gives $\ell\in \N$ and $p\in \Z$ such that for all~$n\ge n_1+\ell$ we have:
$n+p\ge n_0$ and $\zeta_1(n)$ is the unique zero of $y'$ in the interval~$\big(\zeta(n+p), \zeta(n+p+1)\big)$.
Set $n_1^*:=n_1+\ell+|p|$, and modify~$\zeta_1$ to~$\zeta_1^*\colon [n_1^*, +\infty)\to \R$ by setting
 $\zeta_1^*(t)=\zeta_1(t-p)$. Then  $\zeta(n)< \zeta_1^*(n) < \zeta(n+1)$ for all $n\ge n_1^*$. 
 The compositional inverse of $\zeta_1^*$ is given by $(\phi_1/\pi)-(k_1-p)$. Thus replacing $\zeta_1$,~$n_1$,~$k_1$ by 
 $\zeta_1^*$,~$n_1^*$,~$k_1-p$, all clauses are satisfied. 
\end{proof}

\noindent
Define $N\colon \R^{\geq e}\to\N$ by
$$N(t)\ :=\  \big| [e,t]\cap y^{-1}(0) \big|\  =\  \min\{n:s_n>t\},$$
so for $n\ge 1$: $N(t)=n\Leftrightarrow s_{n-1} \le t < s_n$. Thus $N(t)\to+\infty$ as $t\to+\infty$; in fact:

\begin{lemma}\label{lem:N(t)} $N \sim \phi/\pi$.
\end{lemma}
\begin{proof}
Take $n_0$, $k$ as in the proof of Lemma~\ref{lem:zeta}, so
$\phi(s_n)=(k+n)\pi$ for $n\geq n_0$.
Let $t\geq e$ be such that $N(t) \geq n_0+1$ ; then
$s_{N(t)-1} \leq t < s_{N(t)}$ and thus
$$N(t)+k-1=\phi(s_{N(t)-1})/\pi  \leq \phi(t)/\pi  < \phi(s_{N(t)})/\pi =N(t)+k.$$
This yields $N \sim \phi/\pi$. 
\end{proof}

\noindent
The quantity $N(t)$ has been studied extensively in connection with second order linear differential equations; see~\cite[Chapter~IX, \S{}5, and the literature quoted on p.~401]{Hartman}.
For example, the lemma below is a consequence of a result due to Wiman~\cite{Wiman} that holds under more general assumptions (see \cite[Chapter~IX, Corollary~5.3]{Hartman}),
but also follows easily using our Hardy field calculus.
Here we assume $a=0$, so $f=4b\in \c_e$.  

\begin{lemma}
Suppose $f(t)>0$ for all $t\geq e$, and $\big(1/\sqrt{f}\big)'\prec 1$. Then 
$$N(t)\ \sim\  \frac{1}{2\pi} \int_e^t  \sqrt{f(s)}\,ds\quad\text{as $t\to +\infty$.}$$
\end{lemma}
\begin{proof} From  $\big(1/\sqrt{f}\big){}'\prec 1$ we get $f^\dagger\prec \sqrt{f}$. Now $f$ is hardian, so $f\preceq 1/x^2$ would give
$f^\dagger\succeq 1/x$, which together with $f\preceq 1/x^2$  contradicts $f^\dagger\prec \sqrt{f}$. Thus~$f\succ 1/x^2$. 
 For the rest of the argument we arrange $H$ is maximal with $(g,\phi)\in H^2$. 
Corollary~\ref{cor:2nd order, f succ 1/x^2} yields a pair $(g_1,\phi_1)\in H^2$ parametrizing $\ker_{\Calinf}(\der^2+b)$ such that~$\phi_1'\sim (1/2)\sqrt{f}$. Then $\phi-\phi_1\in \R$ by Lemma~\ref{lem:g,phi unique}, so $\phi'\sim (1/2)\sqrt{f}$. Let~$\phi_2\in \c_e^1$ be given by $\phi_2(t)=(1/2)\int_e^t \sqrt{f(s)}\, ds$. Then $\phi_2'=(1/2)\sqrt{f}$, so (the germ of) $\phi_2$ lies in $H$ and 
$\sqrt{f}\succ 1/x$, so $\phi_2>\R$. Hence by [ADH, 9.1.4(ii)] we have $\phi\sim \phi_2$. Now
apply Lemma~\ref{lem:N(t)}. 
\end{proof}

\noindent
In view of Lemma~\ref{lem:no limit pt} one may ask to what extent the results in this subsection generalize to higher-order linear
differential equations over Hardy fields.

\subsection*{When is the perfect hull $\upo$-free?}
Here we use the lemmas that made up the proof of Theorem~\ref{thm:Bosh} to characterize $\upo$-freeness of the ($\d$-)~perfect hull of $H$:

\begin{theorem}\label{thm:upo-freeness of the perfect hull}
The following are equivalent: 
\begin{enumerate}
\item[\textup{(i)}] $H$ is not $\upl$-free or  $\bar{\omega}(H)=H\setminus\sigma\big(\Upg(H)\big){}^\uparrow$;
\item[\textup{(ii)}] $\Dx(H)$ is $\upo$-free;
\item[\textup{(iii)}] $\Ex(H)$ is $\upo$-free.
\end{enumerate}
\end{theorem}
 
\noindent
In connection with this theorem recall that by Corollary~\ref{cor:perfect Schwarz closed, 2}, 
a $\d$-perfect Hardy field is Schwarz closed iff it is $\upo$-free, so in (ii), (iii) we could have also written ``Schwarz closed''
instead of ``$\upo$-free''.
The implication (i)~$\Rightarrow$~(ii) was shown already in Lemma~\ref{lem:D(H) upo-free, 4}.
To show the contrapositive of~(iii)~$\Rightarrow$~(i) suppose~$H$ is $\upl$-free and~$\bar{\omega}(H)\neq H\setminus\sigma\big(\Upg(H)\big){}^\uparrow$. Since $\bar{\omega}(H)\subseteq H\setminus\sigma\big(\Upg(H)\big){}^\uparrow$
this yields 
%Corollary~\ref{cor:gen osc closed upward} then yields some 
$\upo\in H$ with~$\bar{\omega}(H)< \upo  < \sigma\big(\Upg(H)\big)$, and so by Lemma~\ref{lem:upo-freeness of the perfect hull} below, $\operatorname{E}(H)$ is not $\upo$-free.  The proof of this lemma relies on Corollary~\ref{cor:perfect Schwarz closed, 1}, but
additionally draws on some results from Sections~\ref{sec:Hardy fields} and~\ref{sec:upo-free Hardy fields}.

\begin{lemma}\label{lem:upo-freeness of the perfect hull}
Suppose $H$ is $\upl$-free, and $\upo\in H$, $\bar{\omega}(H)<\upo<\sigma\big(\Upg(H)\big)$.
Then 
$$\omega(E) < \upo < \sigma\big(\Upg(E)\big)\qquad\text{for $E:=\Ex(H)$,}$$
hence $E$ is not $\upo$-free.
\end{lemma}
\begin{proof}
We may replace~$H$ by any $\upl$-free Hardy subfield $L$ of $E$ containing~$H$ such that $\Gamma^<$ is cofinal in~$\Gamma^<_L$, by [ADH, 11.8.14, 11.8.29]. Using this observation  and Proposition~\ref{prop:const field ext} we replace $H$ by $H(\R)$ to
%first arrange that~$H$ is an $H$-field, by  replacing~$H$ by its $H$-field hull inside~$E$,
%using Lemma~\ref{lem:Gehret}. Next, replacing~$H$ by the Hardy subfield $H(\R)$ of~$E$ and
%using Lemma~\ref{lem:const field ext}, we 
arrange~$H\supseteq\R$.  Next  we
replace~$H$ by $\operatorname{Li}(H)\subseteq E$ to arrange that $H$ is Liouville closed, using Proposition~\ref{prop:Gehret}.  Now~$\omega(E)=\bar{\omega}(E)$ is downward closed and $\bar{\omega}(E)\cap H=\bar{\omega}(H)$,  so $\omega(E) < \upo$. Towards a contradiction, assume~$\upo \in \sigma\big(\Upg(E)\big){}^\uparrow$. 
%With~$\sigma\big(\Upg(E)\big)$ being upward closed  by Corollary~\ref{cor:perfect Schwarz closed, 1}, this yields
Take~$\upg\in \Upg(E)$ with~$\sigma(\upg)=\upo$. 
Corollaries~\ref{cor:upo} and \ref{cor:sigma(upg)=upo} also yield a germ~$\tilde\upg\in (\Calinf)^\times\setminus \{\upg\}$ 
with $\tilde\upg >0$ and~${\sigma(\tilde\upg)=\upo}$, and  a maximal Hardy field extension~$M$ of $H$
containing $\tilde\upg$. 
Since $M$ is $\upo$-free (by Theorem~\ref{upo}) and~$\upo\notin\bar{\omega}(M)$, we have~$\upo\in\sigma\big(\Upg(M)\big){}^{\uparrow}$ by Corollary~\ref{omuplosc} and so~$\tilde \upg\in\Upg(M)$ by~[ADH, 11.8.31].
Since $E\subseteq M$ we have~$\Upg(E)\subseteq\Upg(M)$. Then from~$\sigma(\upg)=\upo=\sigma(\tilde \upg)$ we   obtain $\upg=\tilde \upg$ by [ADH, 11.8.29], a contradiction.
\end{proof}

\begin{remarkNumbered}\label{rem:non-uniqueness}
Suppose $H\supseteq \R$ is Liouville closed and  $\upo\in H$ satisfies 
$$\bar\omega(H)\ <\ \upo\ <\ \sigma\big(\Upg(H)\big).$$    Then 
the uniqueness in Theorem~\ref{thm:Bosh} fails, by Corollaries~\ref{cor:upo} and \ref{cor:sigma(upg)=upo}:  any~$\upg>0$ in any $\d$-maximal
Hardy field extension~$M$ of $H$ with $\sigma(\upg)=\upo$ yields a pair
$(g,\phi)$ parametrizing $V:=\ker_{\Calinf} (4\der^2+\upo)$ where $g:=1/\sqrt{\upg}$ and 
$\phi\in M$, $\phi'= \frac{1}{2}\upg$.
 \end{remarkNumbered}

\noindent
To finish the proof of Theorem~\ref{thm:upo-freeness of the perfect hull}
it remains to show the implication~(ii)~$\Rightarrow$~(iii), which we do in Lemma~\ref{lem:upo-freeness of the perfect hull, (ii)=>(iii)} below. (This implication holds trivially if $H$ is bounded, by Theorem~\ref{thm:Bosh 14.4}.)
We precede this lemma with some observations. 

\medskip
\noindent
If $\phi$ is active in $H$, then the pre-$H$-field $H^\phi$ has small derivation $\derdelta=\phi^{-1}\der$;
so if~$h\in H$, $h \prec 1$, then $\derdelta^n(h) \prec 1$ for all $n$.
The next lemma yields a variant of this when $h$ is multiplied by a germ in~$\Calinf$ with sufficiently small 
derivatives:

\begin{lemma}\label{lem:derdelta y}
Let $y=hz$, $h\in H$, $h\prec 1$ and  $z\in\Calinf$, $z\preceq 1$ and~$z^{(j)}\preceq \ex^{-x}$ for~$j=1,\dots,n$.
Let also $\phi$ be active in $H$ with $0<\phi\preceq 1/x$, and~$\derdelta=\phi^{-1}\der$.
Then~$\derdelta^j(y)\prec 1$ for $j=0,\dots,n$.
\end{lemma}
\begin{proof}
Let $j$, $k$ with $k\leq j$ range over $\{1,\dots,n\}$.
By the Product Rule for the derivation $\derdelta$ and the remark preceding the lemma it is enough to show
that~${\derdelta^j(z)\preceq 1}$. 
Let $R^j_k\in\Q\{X\}$  be as in Lemma~\ref{lem:fcircn}.
Now $\lambda:=-\phi^\dagger\asymp 1/x$, hence $R^j_k(\lambda)\preceq 1$,
and $(\phi^{-j})^\dagger=-j\phi^\dagger\asymp\lambda \prec 1=(\ex^{x})^\dagger$, hence also $\phi^{-j} \prec \ex^{x}$.
This yields
$$\derdelta^j(z)\ =\ \phi^{-j}\big( R^j_j(\lambda)z^{(j)} + \cdots + R^j_1(\lambda)z' \big) \preceq \phi^{-j}\ex^{-x} \prec 1,$$
which is more than enough. 
\end{proof}

\noindent
We have an ample supply of oscillating germs $z$ as in Lem\-ma~\ref{lem:derdelta y}:

\begin{lemma}\label{lem:ex-x sin x}
Let $z:=\ex^{-x}\sin x\in\mathcal C^\omega$; then $\abs{z^{(n)}} \leq 2^n\ex^{-x}$ for all $n$.
\end{lemma}
%\begin{proof}
%For $0\leq j\leq n$  we have  
%$$(\ex^{-x})^{(n-j)} (\sin x)^{(j)}=\begin{cases}
%(-1)^{n-j/2}\ex^{-x}\sin x & \text{if $j$ is even} \\
%(-1)^{n-(j+1)/2}\ex^{-x}\cos x & \text{if $j$ is odd.}
%\end{cases}$$
%Hence by the Product Rule we obtain $\abs{z^{(n)}}  \leq \ex^{-x} \sum_{j=0}^n {n\choose j} = 2^n\ex^{-x}$.
%\end{proof}

\noindent
In the next lemma and its corollary
%, we assume (as in the subsection on Hardy fields of Section~\ref{sec:d-alg extensions}) that 
our Hardy field~$H$ contains $\R$ and is real closed, and
$\hat H$ is an immediate Hardy field extension of $H$. We now have the following perturbation result:
%Let $\hat f\in\hat H\setminus H$. We show how to perturb $\hat f$ to
%an $f\in \Calinf$ with an isomorphism $H\<f\>\cong H\<\hat f\>$:

\begin{lemma}\label{lem:perturb hat f}
Suppose $H$ is ungrounded, $\Psi_H^{>0}:=\Psi_H\cap\Gamma_H^>\neq\emptyset$. Let $\hat f\in \hat H\setminus H$ and~${Z(H,\hat f)=\emptyset}$.
Let $g\in H$, $vg>v(\hat f -H)$ and $z\in\Calinf$, $z\preceq 1$,  $z^{(n)}\preceq\ex^{-x}$ for all~$n\geq 1$.
Then~$f:=\hat f+gz\in\Calinf$ is hardian over $H$, and we have an
iso\-mor\-phism~$H\langle f\rangle \to H\langle \hat f\rangle$ of $H$-fields over $H$
sending~$f$ to~$\hat f$.
\end{lemma}
\begin{proof}
The hypothesis $\Psi_H^{>0}\neq\emptyset$ and [ADH, 9.2.15] yields active $\phi$   in~$H$ with~$\phi^\dagger \asymp \phi$. But also $t^\dagger\asymp t$ for $t:= x^{-1}$ in $H(x)$, so $\phi\asymp t$ by the uniqueness in [ADH, 9.2.15]. Below $\phi$ ranges over the active elements of $H$ such that $0 < \phi\preceq t$, and~$\derdelta:=\phi^{-1}\der$.
Let $h\in H$, $\fm\in H^\times$ be such that $\hat f-h\preceq \fm$; by Corollary~\ref{cor:iso hat h, d-trans}  it is enough to show that then  $\derdelta^n\big(\frac{f-h}{\fm}\big)\preceq 1$ for all $n$.  Now 
$u:=\frac{\hat f-h}{\fm}\in\hat H$, $u\preceq 1$, and the valued differential field~$\hat H^\phi$ has small derivation,
so $\derdelta^n(u)\preceq 1$ for all $n$. Moreover, $g/\fm\in H$, $g/\fm\prec 1$, so
$\derdelta^n\big(\frac{g}{\fm}z\big)\prec 1$ for all $n$, by Lemma~\ref{lem:derdelta y}.
Thus~$\derdelta^n\big(\frac{f-h}{\fm}\big) = \derdelta^n\big(\frac{\hat f-h}{\fm}\big)+
\derdelta^n\big(\frac{g}{\fm}z\big)\preceq 1$, for all $n$.
\end{proof}

\noindent
Using Lemmas~\ref{lem:ex-x sin x} and~\ref{lem:perturb hat f},  and results of [ADH, 11.4] we now obtain:

\begin{cor}\label{cor:d-trans => c-sequ}
Suppose $H$ is ungrounded with $\Psi_H^{>0}\neq\emptyset$.
Let $(f_\rho)$ be  a divergent pc-sequence in $H$ of $\d$-transcendental type over $H$
and with a pseudolimit in~$\Ex(H)$. Then $(f_\rho)$ is a c-sequence.  
\end{cor}
\begin{proof}
Let $f_\rho\leadsto\hat f\in \Ex(H)$. Then by [ADH, 11.4.7, 11.4.13] the Hardy field $H\langle\hat f\rangle$ is an immediate extension of $H$, 
and~$Z(H,\hat f)=\emptyset$
. Suppose $(f_\rho)$ is not a c-sequence. Then we can take $g\in H^\times$ with
$vg>v(H-\hat f)$. By Lemmas~\ref{lem:ex-x sin x} and~\ref{lem:perturb hat f},
the germ $f:=\hat f+g\ex^{-x}\sin x$ generates a Hardy field $H\langle f\rangle$ over $H$; however,
no maximal Hardy field extension of $H$ contains both $\hat f$ and $f$, contradicting $\hat f\in\Ex(H)$.
\end{proof}

\noindent
We can now supply the proof of the still missing implication (ii)~$\Rightarrow$~(iii) in
Theorem~\ref{thm:upo-freeness of the perfect hull}:

\begin{lemma}\label{lem:upo-freeness of the perfect hull, (ii)=>(iii)}
Suppose $H$ is $\upo$-free. Then  $\Ex(H)$ is also $\upo$-free.
\end{lemma}
\begin{proof}
Since~$E:=\Ex(H)$ is Liouville closed  and contains $\R$ we may replace $H$ by the Hardy subfield
$\operatorname{Li}\!\big(H(\R)\big)$ of $E$, which remains $\upo$-free by Theorem~\ref{thm:ADH 13.6.1}, and arrange that $H\supseteq \R$ and $H$ is Liouville closed (so Corollary~\ref{cor:d-trans => c-sequ} applies).
Towards a contradiction, suppose~$\upo\in E$, $\omega(E) < \upo < \sigma\big(\Upg(E)\big)$; then
$\omega(H) < \upo < \sigma\big(\Upg(H)\big)$. Choose a logarithmic sequence $(\ell_\rho)$ for $H$
and define~$\upo_\rho:=\omega(-\ell_\rho^{\dagger\dagger})$. Then $(\upo_\rho)$ is a divergent pc-sequence in $H$
with $\upo_\rho\leadsto\upo$, by [ADH, 11.8.30]. By~[ADH, 13.6.3], $(\upo_\rho)$ is of $\d$-transcendental type over $H$.
Its width is~$\big\{\gamma\in (\Gamma_H)_\infty:\gamma> 2\Psi_H\big\}$ by~[ADH, 11.7.2], which contains $v(1/x^4)=2v\big((1/x)'\big)$, so $(\upo_\rho)$ is not a c-sequence,
contradicting Corollary~\ref{cor:d-trans => c-sequ}.
\end{proof}

\noindent
Next we describe for $j=1,2$ a $\upl$-free Hardy field $H_{(j)}\supseteq\R$ and $\upo_{(j)}\in H_{(j)}$ such that~$\omega\big(\Upl(H_{(j)})\big) <  \upo_{(j)}  < \sigma\big(\Upg(H_{(j)})\big)$  (so $H_{(j)}$ is not $\upo$-free by [ADH, 11.8.30]), and  
\begin{enumerate}
\item  $\upo_{(1)}\in\bar{\omega}(H_{(1)})$;
\item $\upo_{(2)}\notin\bar{\omega}(H_{(2)})$.
\end{enumerate}
It follows that $\bar{\omega}(H_{(1)})=H_{(1)}\setminus\sigma\big(\Upg(H_{(1)})\big){}^\uparrow$ by Lemma~\ref{lem:baromega(H)},
hence   condition (i) in Theorem~\ref{thm:upo-freeness of the perfect hull} is  satisfied for $H=H_{(1)}$, but it is {\it not}\/ satisfied for $H=H_{(2)}$; thus~$\operatorname{E}(H_{(1)})$ is   $\upo$-free, whereas $\operatorname{E}(H_{(2)})$ is not. 

\medskip
\noindent
To construct such $H_{(j)}$ and $\upo_{(j)}\in H_{(j)}$ we start with a hardian translogarithmic germ $\ell_\omega$   (see the remarks before Proposition~\ref{prop:translog}), and
set $$\upg\ :=\ \ell_\omega^\dagger,\quad \upl\ :=\ -\upg^\dagger,\quad \upo_{(1)}\ :=\ \omega(\upl), \quad \upo_{(2)}\ :=\ \sigma(\upg)\ =\ \upo_{(1)} +\upg^2.$$
Using [ADH, Sections 11.5, 11.7] we see that the Hardy field~$E:=\R( \ell_0,\ell_1,\ell_2,\dots)$ is $\upo$-free and  that
the elements
$\upo_{(1)}, \upo_{(2)}\in M:=E\langle\ell_\omega\rangle$  are pseudolimits of the pc-sequence~$(\upo_n)$ in~$E$.  
For $j=1,2$, we consider the  Hardy subfield   
$$H_{(j)}\ :=\  E\langle \upo_{(j)}\rangle$$
of   $M$, an immediate $\upl$-free extension of $E$ by [ADH, 13.6.3, 13.6.4], and therefore~$\omega\big(\Upl(H_{(j)}\big) <  \upo_{(j)}  < \sigma\big(\Upg(H_{(j)})\big)$ by [ADH, 11.8.30]. Moreover 
$$\bar{\omega}\big(H_{(j)}\big)\ =\ \bar{\omega}(M)\cap H_{(j)}\ \text{  for }j=1,2,$$  so (1) holds
since $\upo_1\in\omega(M)\subseteq\bar{\omega}(M)$,  
whereas $\upo_2\in\sigma\big(\Upg(M)\big)\subseteq M\setminus \bar{\omega}(M)$,  hence~(2) holds.

\begin{exampleNumbered}\label{ex:counterex}
Set $H:=\operatorname{E}(H_{(2)})$.  Then the Hardy field $H$ is perfect, so~${H\supseteq\R}$ is a Liouville closed Hardy field with $\I(K)\subseteq K^\dagger$, but $H$ is not $\upo$-free.  
This makes good on a promise made before Lemma~\ref{lem:counterex}.
\end{exampleNumbered} 

\subsection*{Antiderivatives of rational functions as phase functions}
In this sub\-sec\-tion~$H=\R(x)$, so $K=H[\imag]=\C(x)$. If $f\in H\setminus\bar{\omega}(H)$
and $(g,\phi)\in \operatorname{Li}(H)^2$ para\-me\-tri\-zes~$\ker_{\Calinf}(4\der^2+f)$, then $(\phi')^2\in H$ by
Corollaries~\ref{cor:2nd order, Liouvillian solutions}, \ref{naicor}, and~\ref{cor:phi' quadratic}.
In Corollary~\ref{cor:phi' in H, 1} below we give a condition on such~$f$,~$g$,~$\phi$ that ensures $\phi'\in H$, to
be used in Section~\ref{sec:Bessel}.
We precede this with remarks about ramification in quadratic extensions of~$K$.
So let  $L$ be a field extension of $K$ with $[L:K]=2$. 

\begin{lemma}\label{lem:quadratic ext, p}
Up to multiplication by $-1$, there is a unique $y\in L$ such that~$L=K(y)$ 
and $y^2=p(x)$ where   $p\in \C[X]$ is monic and separable.
\end{lemma}
\begin{proof}
By [ADH, 1.3.11],  $A:=\C[x]$ is integrally closed, so
 [ADH, 1.3.12, 1.3.13] yield a $y\in L$ with minimum polynomial $P\in A[Y]$ over $K$ such that~$L=K(y)$.
Take~$a,b\in A$ with $P=Y^2+aY+b$. Replacing $a$, $b$, $y$ by $0$, $b-(a/2)^2$, $y+(a/2)$, respectively, we arrange $a=0$.
Thus $y^2=p(x)$ for $p\in\C[X]$ with $p(x)=-b$, and replacing $y$ by $cy$ for suitable $c\in\C^\times$ we arrange that $p$ is monic. If~$c\in\C$ and~$p\in (X-c)^{2}\C[X]$, then we may also replace $p$, $y$ by $p/(X-c)^{2}$, $y/(x-c)$, respectively. In this way we 
arrange that $p$ is separable. 
Suppose $L=K(z)$ and~$z^2=q(x)$ where~$q\in\C[X]$ is monic and separable.
Take~$r,s\in K$ with~$z=r+sy$. Then~$s\neq 0$, and~$q(x)=z^2=\big(r^2+s^2p(x)\big)+(2rs)y$, hence $r=0$ and so $q(x)=s^2p(x)$.
Since~$p$,~$q$ are monic and separable, this yields $s^2=1$ and thus $z=-y$ or~$z=y$.
\end{proof}

\noindent
In the following $y$, $p$ are as in Lemma~\ref{lem:quadratic ext, p}.
For each $c\in\C$ we have the valuation~$v_c\colon K^\times\to\Z$  that is trivial on $\C$ with ${v_c(x-c)=1}$,
and we also have the valuation~$v_\infty\colon K^\times\to\Z$ that is trivial on~$\C$ with~$v_\infty(x^{-1})=1$~[ADH, 3.1.30]. 
Given~$f\in K^\times$ there are
only finitely many~$c\in\C_\infty:=\C\cup\{\infty\}$ such that~$v_c(f)\neq 0$; moreover, 
$\sum_{c\in\C_\infty} v_c(f)=0$, with $f\in\C^\times$ iff~$v_c(f)=0$ for all~$c\in\C$. 
Let~$c\in\C_\infty$, and  equip $K$ with the valuation ring~$\mathcal O_c$ of $v_c$.
By~[ADH, 3.1.15, 3.1.21], either exactly one or exactly two  valuation rings of $L$ lie over~$\mathcal O_c$.  
The  residue morphism~$\mathcal O_c\to\res(K)$ restricts to an isomorphism $\C\to\res(K)$, and
equipping~$L$ with   a valuation ring lying over $\mathcal O_c$, composition
with the natural inclusion $\res(K)\to\res(L)$ yields an isomorphism $\C\to\res(L)$;
thus the valued field extension~${L\supseteq K}$ is immediate iff~$L$ is unramified over $K$, that is,  $\Gamma_L=\Gamma=\Z$.

\begin{lemma}\label{lem:quadratic ext}
Suppose $c\neq\infty$.
If $p(c)=0$, then only one valuation ring  of $L$ lies over $\mathcal O_c$, and 
equipping $L$ with this valuation ring we have $[\Gamma_L:\Gamma]=2$.
If $p(c)\neq 0$, then there are exactly two valuation rings of $L$ lying over $\mathcal O_c$,
and equipped with any one of these valuation rings, $L$ is  unramified over $K$.
\end{lemma}
\begin{proof}
If $p(c)=0$, then $v_c(p)=1$,
so  by~[ADH, 3.1.28] there is a unique   valuation ring~$\mathcal O_L$ of $L$ lying over $\mathcal O_c$, and equipping $L$ with $\mathcal O_L$  we have $[\Gamma_L:\Gamma]=2$.
Now suppose $p(c)\neq 0$. 
We identify $K$ with its image under the embedding of $K$ into the 
valued  field $K^{\operatorname{c}}:=\C(\!(t)\!)$ of Laurent series over $\C$ which is the identity on~$\C$ and sends $x-c$ to $t$. Take $\alpha\in\C^\times$ with $\alpha^2=p(c)$.
Hensel's Lem\-ma~[ADH, 3.3.5] yields~$z\in K^{\operatorname{c}}$ with $z^2=p(x)$ and $z\sim\alpha$.
Let $\mathcal O_+$, $\mathcal O_-$ be the preimages of the valuation ring of 
 the valued subfield~$K(z)$ of $K^{\operatorname{c}}$ under the field isomorphisms~$L=K(y)\to K(z)$ over~$K$
 with $y\mapsto z$ and $y\mapsto -z$, respectively.
Then~$\mathcal O_+\neq\mathcal O_-$   lie over $\mathcal O_c$, and each turns~$L$ into an immediate extension of~$K$.
\end{proof}

\noindent
In the next lemma we set $d:=\deg p$, so $d\geq 1$. 

\begin{lemma}
If $d$ is odd, then only one valuation ring  of $L$ lies over $\mathcal O_\infty$, and 
equipping $L$ with this valuation ring we have $[\Gamma_L:\Gamma]=2$.
If $d$ is even, then there are exactly two valuation rings of $L$ lying over $\mathcal O_\infty$,
and equipped with any one of these valuation rings, $L$ is unramified over $K$.
\end{lemma}
\begin{proof}
We have $v_\infty(p)=-d$. Hence if $d$ is odd, then we can argue using [ADH, 3.1.28] as in the proof of Lemma~\ref{lem:quadratic ext}. Suppose $d$ is even, so with $e=d/2$ we have~$(y/x^e)^2=p(x)/x^d \sim 1$.
Identify $K$ with its image under the embedding of~$K$ into the 
valued  field $K^{\operatorname{c}}:=\C(\!(t)\!)$ of Laurent series over $\C$ which is the identity on~$\C$ and sends $x^{-1}$ to $t$. Then [ADH, 3.3.5] yields $z\in K^{\operatorname{c}}$
with $z\sim 1$ and $z^2=p(x)/x^d$. 
Let $\mathcal O_+$, $\mathcal O_-$ be the preimages of the valuation ring of 
 the valued subfield~$K(z)$ of~$K^{\operatorname{c}}$ under the field isomorphisms~$L\to K(z)$ over~$K$
 with $y\mapsto x^e z$ and $y\mapsto -x^e z$, respectively.
Then~$\mathcal O_+\neq\mathcal O_-$ are   valuation rings of $L$ lying over~$\mathcal O_\infty$, each of which  turns~$L$ into an immediate extension of $K$.
\end{proof}

\begin{cor}\label{cor:quadratic ext}
There are at least two $c\in\C_\infty$ such that  
some valuation ring of $L$ lying over $\mathcal O_{c}$
makes $L$ ramified over $K$.
\end{cor}

\noindent
Next we let $C$ be any field of characteristic zero and consider the 
$\d$-valued Hahn field~$C(\!(t^{\Q})\!)$   with its strongly additive $C$-linear derivation satisfying~${t'=1}$. 
We let~$q$,~$r$,~$s$ range over $\Q$, and
$z=\sum_q z_qt^q\in C(\!(t^{\Q})\!)^\times$ with all $z_q\in C$.
Put
$$q_0\ :=vz\ =\ \min \operatorname{supp} z\in\Q,$$ so $z\sim z_{q_0}t^{q_0}$.
If $z\notin C(\!(t)\!)$, then we also set~$q_1:=\min\!\big( (\operatorname{supp} z)\setminus\Z\big)\in\Q\setminus\Z$, so~$q_1\geq q_0$.
In Lemmas~\ref{lem:z in C((t)), 1} and~\ref{lem:z in C((t)), 2} below we give sufficient conditions for  $z$ to be in
$C(\!(t)\!)$.
Set
$$w\ :=\ z^2\ =\ \sum_q w_q t^q \quad\text{where}\quad
 w_q\ =\ \sum_{r+s=q} z_rz_s,$$
so $w\sim z_{q_0}^2t^{2q_0}$, and observe:

\begin{lemma}\label{lem:z^2}
If $w\notin C(\!(t)\!)$ \textup{(}and so $z\notin C(\!(t)\!)$\textup{)}, then 
$$\min\!\big( (\operatorname{supp} w)\setminus\Z\big)=q_0+q_1,\quad w_{q_0+q_1}= 2z_{q_0}z_{q_1}.$$
\end{lemma}
%\begin{proof}
%Suppose $r,s\in\operatorname{supp} z$ with $q:=r+s<q_0+q_1$ and $q\notin\Z$.
%We then have $r<q_1$ or~$s<q_1$: otherwise
%$2q_1 \leq r+s<q_0+q_1$ and hence $q_1<q_0$, a contradiction.
%Now if $r<q_1$, then $r\in\Z$ and thus $s=q-r\notin\Z$, so $s\geq q_1$ and thus
%$q=r+s\geq q_0+q_1$, a contradiction; similarly $s<q_1$ leads to a contradiction.
%This yields the lemma.
%\end{proof}

\begin{lemma}\label{lem:z in C((t)), 1}
Suppose $\omega(z)\in t^{-1}C[[t]]$. Then $z\in  C(\!(t)\!)$.
\end{lemma}
\begin{proof}
Put
$u:=z'=\sum_q u_qt^{q}$, 
$u_q=   (q+1)z_{q+1}$.
If $q_0\neq 0$, then $u\sim~q_0z_{q_0} t^{q_0-1}$.
Hence  $q_0\geq -1$: otherwise  $-\omega(z)=2u+w\sim z_{q_0}^2t^{2q_0}$,
contradicting $\omega(z)\preceq t^{-1}\preceq t^{-2}$. Moreover, if $q_0=-1$, then
$(2u+w)-(-2z_{-1}+z_{-1}^2)t^{-2}\prec t^{-2}$ and so $z_{-1}=2$.
Towards a contradiction,  suppose $z\notin   C(\!(t)\!)$.
%and  we claim that $q_1>q_0$: Suppose otherwise; then  $q_0\neq 0$, so $z'\sim q_0z_{q_0} t^{q_0-1}$ where $q_0-1\notin\Z$, and together with $2z'+z^2=-\omega(z)=-f\in C(\!(t)\!)$ this yields the existence of~$r,s\in\operatorname{supp} z$ with $q_0-1=r+s$. Hence $q_0-1\geq 2q_0$ and thus  $q_0=-1$, a contradiction to $q_0=q_1\notin\Z$. This shows $q_1>q_0$.
We have $u\notin C(\!(t)\!)$. Indeed
\begin{equation}\label{eq:u}
\min\!\big( (\operatorname{supp} u)\setminus\Z\big)=q_1-1,\quad u_{q_1-1}= q_1z_{q_1}.
\end{equation}
Also $w=-\omega(z)-2u\notin C(\!(t)\!)$, and by the previous lemma
\begin{equation}\label{eq:v}
\min\!\big( (\operatorname{supp} w)\setminus\Z\big)=q_0+q_1,\quad w_{q_0+q_1}= 2z_{q_0}z_{q_1}.
\end{equation}
From \eqref{eq:u}, \eqref{eq:v}, and $2u+w\in C(\!(t)\!)$ we get $q_1-1=q_0+q_1$ and~$2q_1z_{q_1}=-2z_{q_0}z_{q_1}$,
hence $q_0=-1$, $q_1=-z_{q_0}$.
Thus $q_1=-2<-1=q_0$, a contradiction. 
\end{proof}

\noindent
In [ADH, p. 519] we defined $\omega^\phi:E \to E$ for a differential field $E$ and $\phi\in E^\times$.

\begin{lemma}\label{lem:z in C((t)), 2}
Suppose $\omega^{-1/t^2}(z)\in C[[t]]^\times$. Then $z\in  C(\!(t)\!)$.
\end{lemma}
\begin{proof}
Put $u:=-t^2z'=\sum_q u_qt^{q}$ where $u_q=-(q-1)z_{q-1}$.
If $q_0\neq 0$, then~$u\sim -q_0z_{q_0} t^{q_0+1}$.
We must have $q_0 = 0$: otherwise, if $q_0<1$, then $2q_0<q_0+1$ and so
$-\omega^{-1/t^2}(z)=2u+w \sim z_{q_0}^2 t^{2q_0}$, contradicting $\omega^{-1/t^2}(z)\asymp 1$, whereas if~$q_0\geq 1$
then~$2u+w \preceq t^{q_0+1} \preceq t^2$, again contradicting $\omega^{-1/t^2}(z)\asymp 1$.
Now  suppose $z\notin   C(\!(t)\!)$.
Then 
\begin{equation}\label{eq:u, infty}
\min\!\big( (\operatorname{supp} u)\setminus\Z\big)=q_1+1,\quad u_{q_1+1}= -q_1z_{q_1}
\end{equation}
and by Lemma~\ref{lem:z^2}:
\begin{equation}\label{eq:v, infty}
\min\!\big( (\operatorname{supp} w)\setminus\Z\big)=q_1,\quad w_{q_1}= 2z_{0}z_{q_1}.
\end{equation}
Together with  $2u+w\in C(\!(t)\!)$ this yields a contradiction.
\end{proof}

\noindent
We now apply the above with $C=\C$ to show:
 
\begin{cor}\label{cor:z in C((t))}
Let $z\in L$ be such that $\omega(z)=f\in K$.
If  $v_c(f)\geq -1$ for all~$c\in\C$, or 
$v_c(f)\geq -1$ for all but one $c\in \C$ and 
$v_\infty(f)=0$, then $z\in K$.
\end{cor}
\begin{proof}
Let $c\in\C_\infty$ and let $L$ be equipped with a valuation ring lying over $\mathcal O_c$.
If~$c\in\C$, then
we  have a valued differential field embedding $L\to \C(\!(t^{\Q})\!)$ over~$\C$ with~$x-c\mapsto t$,
and identifying $L$ with its image under this embedding, 
if $v_c(f)\geq -1$, then $f\in t^{-1}\C[[t]]$, hence
$z\in \C(\!(t)\!)$ by Lemma~\ref{lem:z in C((t)), 1}, so $K(z)\subseteq \C(\!(t)\!)$ is unramified over~$K$.
If $c=\infty$, then 
we  have a valued differential field embedding~$L\to \C(\!(t^{\Q})\!)^{-1/t^2}$  over $\C$ with~$x^{-1}\mapsto t$, and again identifying $L$ with its image under this embedding, 
if $v_\infty(f)=0$, then $f\in \C[[t]]^\times$
by Lemma~\ref{lem:z in C((t)), 2}, so~$K(z)$ is unramified over $K$. Now use Corollary~\ref{cor:quadratic ext}. 
\end{proof}

\noindent
In the next two lemmas we fix  $c\in\C_\infty$ and equip $K=\C(x)$ with $v=v_c$. Then the valued differential field
$K$ is $\d$-valued, and for all $z\in K^\times$ with $vz=k\neq 0$ we have~$z^\dagger \sim k(x-c)^{-1}$ if $c\neq\infty$, and
$z^\dagger \sim -kx^{-1}$ if $c=\infty$.
In these two lemmas we let $z\in K^\times$, and set $k:=vz$, $f:=\omega(z)$.

\begin{lemma}\label{lem:half-integer, 1}
Suppose $z\succ 1$.
If~$c=\infty$, then~$f\sim -z^2$. 
If $c\neq\infty$ and  $f\preceq 1$, then $z-2(x-c)^{-1} \preceq 1$.
%$z-2(x-c)^{-1} \preceq x-c$.
\end{lemma}
\begin{proof}
If $c=\infty$, then $x\succ 1$ and so $z^\dagger\sim -kx^{-1}\prec 1\prec z$, hence~$f=\omega(z)=-z(2z^\dagger+z)\sim -z^2$.
Now suppose $c\neq\infty$ and  
$f\preceq 1$. 
Applying the automorphism of the differential field $K$ over $\C$ with $x\mapsto x+c$ we arrange $c=0$.
So~$x\prec 1$ and~$z^\dagger\sim kx^{-1}$.
We   have~$-z(2z^\dagger+z)=\omega(z)=f\preceq 1$, so $2z^\dagger\sim-z$, and thus~$z\sim 2x^{-1}$, that is, $z-2x^{-1}\preceq 1$.
%Put~$u:=2x^{-1}-z$, so $u\prec x^{-1}$;\marginpar{commented out things checked but no longer needed}
%we claim that~$u\preceq x$.  This is clear if $u=0$. Assume~$u\ne 0$. Then, 
%using~[ADH, 5.2.1]  for the second equality below,  
%$$
%u\cdot\big( (u^2x^{4})^\dagger-u \big)\ =\ u\cdot\big(2(u^\dagger+z)+u\big)\ =\ \omega(z)-\omega(2x^{-1})\ =\ 
%\omega(z)\ =\ f.$$
%Also
%$u^2x^4\prec x^2\prec 1$,  hence
%$(u^2x^{4})^\dagger\asymp x^{-1}$, which in view of $x^{-1} \succ u$ and the last display
%yields  $u\cdot x^{-1}\asymp f\preceq 1$ and thus $u\preceq x$ as claimed. 
\end{proof}

\begin{lemma}\label{lem:half-integer, 2}
Suppose
 $c\neq\infty$ and $d\in\C^\times$ is such that $f - d(x-c)^{-2} \preceq 1$. Then~$z\succ 1$, and 
 for some $b\in\C$ with $b(2-b)=d$ we have
 $z-b(x-c)^{-1} \preceq 1$.
\end{lemma}
\begin{proof}
We arrange again $c=0$, so   $\omega(z)=f\sim d x^{-2}$. 
If $z\preceq 1$, then $z'\prec x^{-1}$ and thus~$dx^{-2}\sim \omega(z)=-(2z'+z^2)\prec x^{-1}$,
contradicting~$x\prec 1$. Thus $z\succ~1$, hence~$z^\dagger\sim~kx^{-1}$ with~$k<0$. Together with~$-z(2z^\dagger+z)=\omega(z)\sim dx^{-2}$ this yields~$k=-1$ and $z\sim ~bx^{-1}$ with $b\in \C^\times$, $b(2-b)=d$, so
$z-bx^{-1}\preceq 1$.
\end{proof}

\noindent
Let  $f\in H\setminus\bar{\omega}(H)$, and suppose $(g,\phi)\in \operatorname{Li}(H)^2$ parametrizes $\ker_{\Calinf}(4\der^2+f)$. Here is the promised sufficient condition for $\phi'\in H$:

\begin{cor}\label{cor:phi' in H, 1}
Suppose     $v_c(f)\geq -1$ for all~$c\in\C$, or
$v_c(f)\geq -1$ for all but one  $c\in \C$ and 
$v_\infty(f)=0$. Then~$\phi'\in H$. 
\end{cor}
\begin{proof}
Put $y:=g\ex^{\phi\imag}\in\Calinf[\imag]^\times$. The proof of Lemma~\ref{lem:parametrization of ker A}  gives~${4y''+fy=0}$ and $\omega(z)=f$ for~$z:=2y^\dagger=-\phi'{}^\dagger+2\phi'\imag$ . 
We have the differential field extension~$L:=K[z]=K[\phi']\subseteq\Li(H)[\imag]$ of $K$.
If $\phi'\notin H$, then $[L:K]=2$, and then Corollary~\ref{cor:z in C((t))} gives $z\in K$,
a contradiction. Thus $\phi'\in H$.
\end{proof}

\noindent
In the next section on the Bessel equation the relevant $f$ satisfies an even stronger condition, and this gives more 
information about $\phi$:

\begin{cor}\label{cor:phi' in H, 2}
Suppose $v_c(f)\geq 0$ for all $c\in\C^\times$, $v_\infty(f)=0$, and $d\in \C$, 
 $v_{0}\big( f - dx^{-2} \big) \geq 0$.
Then there are $a,b\in\C$ and distinct $c_1,\dots,c_n\in\C^\times$ such that
$$-\phi'{}^\dagger+2\phi'\imag\ =\ a + bx^{-1} + 2 \sum_{j=1}^n (x-c_j)^{-1} \quad \text{and} \quad b(2-b)=d.$$
\end{cor} 
\begin{proof} Corollary~\ref{cor:phi' in H, 1} and its proof gives $z:=-\phi'{}^\dagger+2\phi'\imag\in K^\times$ and $\omega(z)=f$.  Consider first the case
$d\ne 0$. Then by Lemma~\ref{lem:half-integer, 1} we have  $v_{\infty}(z) \ge 0$ and $$v_c\big(z-2(x-c)^{-1}\big)\ \ge\  0\
\text{ whenever $c\in \C^\times$ and $v_c(z)<0$}.$$
 Lemma~\ref{lem:half-integer, 2} gives $v_0(z-bx^{-1})\ge 0$ with $b\in \C$ such that $b(2-b)=d$. 
Taking~$c_1,\dots, c_n$ as the distinct poles of $z$ in $\C^\times$, this yields the desired result by considering the partial fraction decomposition of $z$ with respect to $\C[x]$.  
Next, suppose $d=0$. Then $v_c(f)\ge 0$ for all $c\in \C_{\infty}$, hence $f\in \C\cap H=\R$. Also $f>0$, since $0\in\bar{\omega}(H)$ and
$f\notin \bar{\omega}(H)$.  The example preceding Lemma~\ref{lem:param V'}, together with Corollary~\ref{naicor},
 gives $\phi=\frac{\sqrt{f}}{2}x+r$ with $r\in \R$, so $z=\sqrt{f}\cdot \imag$, and this gives the desired result with
 $a=\sqrt{f}\cdot \imag$, $b=0$, $n=0$.  
\end{proof}

\section{The Example of the Bessel Equation}\label{sec:Bessel}

\noindent
We are going to use the results from Section~\ref{sec:perfect applications}  to obtain  information about the solutions of the Bessel equation
\begin{equation}\label{eq:Bessel}\tag{$\operatorname{B}_\nu$}
x^2 Y''+xY'+(x^2-\nu^2)Y\ =\ 0
\end{equation}
of order $\nu\in\R$. For solutions in $\c_e^2$ ($e\in \R^{>}$), this is equivalent to the equation \eqref{eq:2nd order, app} in Section~\ref{sec:perfect applications} with~$a=x^{-1}$, $b=1- \nu^2 x^{-2}$, so that $f_{\nu}:=-2a'-a^2+4b$ gives
$$f_{\nu}\  =\ -2(x^{-1})'- (x^{-1})^2+ 4\big(1- \nu^2 x^{-2}\big)\  =\   4+(1-4\nu^2)x^{-2}\ \sim\ 4.$$
Thus $f_{\nu}\notin \overline{\omega}\big(\R(x)\big)$, and we have the isomorphism $y\mapsto x^{1/2}y$ of the $\R$-linear space~$V_{\nu}\subseteq\Calinf$ of solutions of \eqref{eq:Bessel} onto the $\R$-linear space of   solutions in~$\Calinf$ of 
\begin{equation}\label{eq:Bessel selfadj}\tag{$\operatorname{L}_\nu$}
4Y''+f_{\nu}Y\ =\ 0.
\end{equation} 
 The nonzero solutions of~\eqref{eq:Bessel} in~$\c^2(\R^>)$ are known as (real) {\it cylinder functions}; cf.~\cite[\S{}15.22]{Watson}.\index{function!cylinder}

\begin{prop}\label{prop:Bessel}
There is a unique  hardian germ~$\phi_\nu$  such that 
$$\phi_{\nu}-x\ \preceq\ x^{-1} \text{ and }\ V_{\nu}\ =\ \left\{\frac{c}{\sqrt{x\phi_{\nu}'}}\cos(\phi_{\nu}+d):\ c,d\in\R\right\}.$$
This germ $\phi_{\nu}$ lies in $\Dx(\Q)\subseteq\Gom$. $($Recall that $\Dx(\Q)=\Ex(\Q).)$ 
\end{prop}

\noindent
If $\nu^2=\frac{1}{4}$, then $V_{\nu}=\R x^{-1/2}\cos x+ \R x^{-1/2}\sin x$, and Proposition~\ref{prop:Bessel} holds with~$\phi_{\nu}=x$. So suppose $\nu^2\neq\frac{1}{4}$.  Then 
Corollary~\ref{cor:2nd order, f succ 1/x^2} gives a germ $\phi\sim x$ in $\Dx(\Q)$ such that~$(g,\phi)$ parametrizes~$V_{\nu}$,
where~$g:=(x\phi')^{-1/2}$.
%Taking $\upg:=1/x$ we have
%  $$f-\sigma(\upg)=4-(1+4\alpha^2)x^{-2}\sim 4\succ \upg^2$$ and hence $b:=\sqrt{f-\sigma(\upg)}\sim 2$.
%Theorem~\ref{thm:extend to H-closed} and Lemma~\ref{lem:ADH 14.2.18} yield some $\phi$ in a Hardy field containing $\R(x)$ with $2\phi'\sim b$ and $\sigma(2\phi')=f$. Then~$\phi\sim x$, and setting~$g:=1/\sqrt{\phi'}$ we have  
%$$\omega(2g^\dagger+2\phi'\imag)=\omega\big({-(\phi')^\dagger}+2\phi'\imag\big)=\sigma(2\phi')=f,$$
%hence $A(g\ex^{\phi\imag})=0$ for $A:=4\der^2+f$. Thus $$x^{-1/2}g\cos\phi,\quad x^{-1/2}g\sin\phi$$ is a basis of~$V$. 
Using this fact, Proposition~\ref{prop:Bessel} now follows from Corollary~\ref{naicor}, Lemma~\ref{lem:parametrization of ker A},  and the next lemma about any such  pair $(g,\phi)$: 

\begin{lemma}\label{lem:Bessel expansion}
We have $\phi-x-r -   \frac{1}{2}(\nu^2-\frac{1}{4})x^{-1} \preceq x^{-3}$ for some $r\in\R$. \end{lemma}
\begin{proof}
Set  $z:=2\phi'$, so $z=2+\varepsilon$, $\varepsilon\prec 1$. From
 $\sigma(z)=f_{\nu}$ and multiplication by $z^2$,
$$2zz''-3(z')^2+z^2(z^2-f_{\nu})\ =\ 0$$
and thus with $\mu:=4\nu^2-1\in\R^\times$, $u:=-(2zz''-3(z')^2)=3(\varepsilon')^2-2y\varepsilon''$:
$$u\ =\ z^2(z^2-f_{\nu})\sim 4(z^2-f_{\nu})= 4(4\varepsilon + \varepsilon^2+ \mu x^{-2}), \text{ and thus}$$
\begin{equation}\label{eq:Bessel expansion}
 u/4\ \sim\ \varepsilon(4+\varepsilon)+\mu x^{-2}.
\end{equation}
We claim that $u \prec  x^{-2}$.
If $\varepsilon\preceq x^{-2}$, then $\varepsilon' \preceq x^{-3}$, $\varepsilon''\preceq x^{-4}$, and the claim is valid.
If $\varepsilon \succ x^{-2}$, then $\varepsilon^\dagger\preceq (x^{-2})^\dagger=-2x^{-1}$, so
$\varepsilon' \preceq x^{-1}\varepsilon\prec x^{-1}\prec 1$, hence~$\varepsilon''\prec (x^{-1})'=-x^{-2}$, which again yields
$u\prec  x^{-2}$.
The claim and \eqref{eq:Bessel expansion} give~$\varepsilon\sim-\frac{\mu}{4}  x^{-2}$ and hence 
$\delta:=\varepsilon+\frac{\mu}{4}  x^{-2}\prec x^{-2}$. 
Indeed, we have~$\delta \preceq x^{-4}$. To see why, note that~$\varepsilon'\sim \frac{1}{2}\mu x^{-3}$ and
$\varepsilon''\sim -\frac{3}{2}\mu x^{-4}$, so $u\sim 6 \mu x^{-4}$, and  $$\textstyle\frac{3}{2}\mu x^{-4}\ \sim\ g/4\  \sim\  \varepsilon(4+\varepsilon)+\mu x^{-2}\ =\ 4\delta+\varepsilon^2, \qquad \varepsilon^2\ \sim\ \frac{\mu^2}{16}x^{-4}.$$
Now the lemma follows by integration from  
\[\phi'-1+\textstyle\frac{1}{2}(\nu^2-\textstyle\frac{1}{4})x^{-2}\ =\ \frac{1}{2}\varepsilon + \frac{1}{8}\mu x^{-2}\ =\ \delta/2\ \preceq\ x^{-4}. \qedhere\]
\end{proof}

\noindent
With $\phi$ and $r$ as in Lemma~\ref{lem:Bessel expansion}, it is $\phi-r$ that is the germ $\phi_{\nu}$ in 
Proposition~\ref{prop:Bessel}, and till further notice we set $\phi:=\phi_\nu$, $f:=f_{\nu}$, and $V:=V_{\nu}$. Thus $\sigma(2\phi')=f$ and
 $\phi_{\nu}=\phi_{-\nu}$. 
As  mentioned before, we do not know if $\Ex(\Q)^{>\R}$ is closed under compositional  inversion. Nevertheless: 

\begin{lemma}\label{phiinvexq}
$\phi^{\operatorname{inv}}\in\Ex(\Q)$.
\end{lemma}
\begin{proof}
Set $\alpha:=\frac{1}{2}(\nu^2-\frac{1}{4})$. Then $\phi=x+ \alpha x^{-1}+o(x^{-1})$, so
$$\phi^{\operatorname{inv}}\ =\ x-\alpha x^{-1}+o(x^{-1})$$ 
by   Corollary~\ref{cor:Entr, 3}, and $\phi^{\operatorname{inv}}$ is hardian. 
Let~$P\in\R(x)\{Y\}$ be as in the remarks before Lemma~\ref{boundphi1}
with $H=\R(x)$,  so~$P(2\phi')=0$.
Corollary~\ref{cor:trdegellinv} then gives~$\tilde{P}\in\R(x)\{Z\}$ such that for all hardian $y>\R$,
%eventually strictly increasing~$y\in\Calinf$
% with $y\succ 1$:
$$P(2y')\ =\ 0 \quad\Longleftrightarrow\quad \tilde{P}(y^{\operatorname{inv}})\ =\ 0,$$
in particular, $\tilde{P}(\phi^{\operatorname{inv}})=0$. 
Let now $H$ be any maximal Hardy field.   Theorem~\ref{thm:transfer} then yields~$z\in H$ such that~$z=x-\alpha x^{-1}+o(x^{-1})$ and $\tilde{P}(z)=0$, so $y:=z^{\operatorname{inv}}$ is hardian and $P(2y')=0$. 
Then
$\sigma(2y')=f$, so $\big((xy')^{-1/2},y\big)$ parametrizes~$V$ by Lemma~\ref{parphi} and a remark preceding that lemma. Also
~$y=x+ \alpha x^{-1}+o(x^{-1})$ by  Corollary~\ref{cor:Entr, 3}.
Thus~$\phi=y$ by Proposition~\ref{prop:Bessel} and so~$\phi^{\operatorname{inv}}=z\in H$.
\end{proof} 

\noindent
This quickly yields some facts on the distribution of zeros of solutions:
Let~${y\in\c_e^2}$ ($e\in\R^{>}$) be a nonzero solution of \eqref{eq:Bessel} and let $(s_n)$ be the enumeration of its zero set. 
From Corollary~\ref{cor:zeta asymptotics}
%,~\ref{cor:zeta, 3}, 
and Lemma~\ref{lem:zeta, differences} we obtain a well-known result,  see for example  \cite[Chapter~XI, Exercise~3.2(d)]{Hartman}, \cite[\S{}27, XIII]{Walter}:

\begin{cor}\label{corcorasympzeta}
We have $s_n\sim\pi n$ and $s_{n+1}-s_n \to \pi$ as $n\to\infty$. 
\end{cor}

\begin{lemma}\label{zetaeq} 
There is a strictly increasing~$\zeta\in\c_{n_0}$ \textup{(}$n_0\in\N$\textup{)} whose germ is in~$\Ex(\Q)$
such that $s_n=\zeta(n)$ for all $n\geq n_0$. 
\end{lemma}
\begin{proof} Take $e_0\ge e$, a representative of $\phi$ in $\c^1_{e_0}$ denoted also by $\phi$, and $c,d\in \R$, such that
$\phi'(t)>0$ and $y(t)=\big(c/\sqrt{t\phi'(t)}\big)\cdot \cos\big(\phi(t)+d\big)$ for all $t\ge e_0$. So we are in the situation described
before Lemma~\ref{lem:zeta}. 
Next, take $n_0$, $k_0$, $\zeta$ as in the proof of that lemma. Then $\zeta$ is strictly increasing with $s_n=\zeta(n)$ for all $n\geq n_0$, and the germ of $\zeta$, denoted by the same symbol,
satisfies $\zeta=\phi^{\operatorname{inv}}\circ \big(\pi\cdot(x+k_0)\big)$. Now use Lemma~\ref{phiinvexq} and 
$\Ex(\Q)\circ \Ex(\Q)^{>\R}\subseteq\Ex(\Q)$ (see the remark after Lemma~\ref{lem:comp with E(Q)}), to conclude
$\zeta\in \Ex(\Q)$. 
\end{proof} 

\noindent
Lemma~\ref{zetaeq} yields an improvement of Corollary~\ref{cor:zeta, 3} in our (Bessel) case: 

\begin{cor} For any  $h\in\c_0$ with hardian germ the sequences~$(s_n)$ and~$\big(h(n)\big)$ are comparable.
\end{cor}

% (It is unknown whether $\Ex(\Q)^{>\R}$ is closed under compositional inversion,
%see the remarks after the proof of Lemma~\ref{lem:comp with E(Q)}.)
\noindent
Lemma~\ref{lem:zeta, differences} also has the following corollary, the first part of which
was observed by Porter~\cite{Porter} (cf.~al\-so~\cite[\S{}15.8, 15.82]{Watson}).

\begin{cor}
If $\nu^2>\frac{1}{4}$, then the sequence $(s_{n+1}-s_n)$ is eventually strictly decreasing, and if 
$\nu^2<\frac{1}{4}$, then $(s_{n+1}-s_n)$ is eventually strictly increasing.
\end{cor}
 
\noindent
Finally, if $\underline{\nu}\in\R$ and  $\underline{y}\in\c_{\underline{e}}^2$ with $\underline{e}\in \R^{>}$ is a nonzero solution of the Bessel equation of order $\underline{\nu}$, then~$(s_n)$ and the  enumeration $(\underline{s}_n)$ of the zero set of $\underline{y}$ are comparable, by Lemma~\ref{lem:comp sequ zeros}. This is related to classical
results on the ``interlacing of zeros'' of  cylinder functions; cf.~\cite[\S\S{}15.22, 15.24]{Watson}.

\medskip\noindent
In the next lemma, \cite[Chapter~10, Theorem~8]{BR} has
$-\frac{1}{2}(\nu^2-\textstyle\frac{1}{4})x^{-1}$
instead of our~$\frac{1}{2}(\nu^2-\textstyle\frac{1}{4})x^{-1}$.  This sign error originated in an integration on \cite[p.~327]{BR}.

\begin{lemma}\label{cor:BR}
Let $y\in V^{\neq}$. Then there is a pair $(c,d)\in\R^\times\times[0,\pi)$ such that~$y=\frac{c}{\sqrt{x\phi'}}\cos(\phi+d)$, and for any such pair we have
\begin{equation}\label{eq:cd}
y-\frac{c}{\sqrt{x}}\cos\left(x+d+ \textstyle\frac{1}{2}(\nu^2-\textstyle\frac{1}{4})x^{-1}\right)\ \preceq\ x^{-5/2}.
\end{equation}
\end{lemma}
\begin{proof} 
Proposition~\ref{prop:Bessel} yields $(c,d)\in\R\times[0,\pi)$ such that
$y=\frac{c}{\sqrt{x\phi'}}\cos(\phi+d)$. Then~${c\neq 0}$.
From $\phi'-1\preceq x^{-2}$ we get
$\frac{1}{\sqrt{\phi'}}-1\preceq x^{-2}$, and for every $u\in\c$ we have~$\cos({x+u})-\cos(x) \preceq u$.
Using also Lemma~\ref{lem:Bessel expansion} this yields \eqref{eq:cd}.
\end{proof}

\noindent
We complement this with some uniqueness properties:

\begin{lemma}\label{lem:BR} Let  $y\in V^{\ne} $. Then there is a unique $(c,d)\in \R^\times\times [0,\pi)$ such that
\begin{equation}\label{eq:cd+} y-\frac{c}{\sqrt{x}}\cos(x+d)\  \prec\ \frac{1}{\sqrt{x}},
\end{equation}
and this is also the unique $(c,d)\in\R^\times\times[0,\pi)$ such that  $y=\frac{c}{\sqrt{x\phi'}}\cos(\phi+d)$. 
\end{lemma}
\begin{proof} For $(c,d)\in \R^\times\times [0,\pi)$ with $y=\frac{c}{\sqrt{x\phi'}}\cos(\phi+d)$  we have
\eqref{eq:cd}, so $$y-\frac{c}{\sqrt{x}}\cos(x+d)\ \preceq\  x^{-3/2}$$ in view of $\cos({x+u})-\cos(x) \preceq u$ for $u\in \c$.
This gives  \eqref{eq:cd+}. Suppose towards a contradiction that \eqref{eq:cd+} also holds for a pair
$(c^*,d^*)\in \R^\times \times[0,\pi)$ instead of~$(c,d)$, with $(c^*,d^*)\ne (c,d)$. Then $d\ne d^*$, say $d< d^*$, so 
$0 < \theta:= d^*-d < \pi$.  Then~$c\cos(x+d)-c^*\cos(x+d+\theta)\prec  1$,
and hence $c\cos(x)-c^*\cos(x+\theta) \prec  1$, which by a trigonometric identity turns into
$$(c-c^*\cos\theta)\cos(x) + c^*\sin\theta \sin(x)\ =\ \sqrt{(c-c^*\cos\theta)^2+(c^*\sin\theta)^2}\cdot\cos(x+s) \prec 1$$
with $s\in \R$ depending only on $c$, $c^*$, $\theta$; see the remarks preceding Lemma~\ref{addsin}. 
This forces $c^*\sin\theta=0$, but $\sin \theta>0$, so $c^*=0$, contradicting $c^*\in \R^{\times}$.  
\end{proof}

\begin{cor}\label{BR+} For any $(c,d)\in\R^\times\times[0,\pi)$ there is a unique $y\in V^{\ne}$ such that~\eqref{eq:cd+} holds. This
$y$ is given by $y= \frac{c}{\sqrt{x\phi'}}\cos(\phi+d)$. 
\end{cor} 

\begin{remarkNumbered}\label{rem:BR+} Lemmas~\ref{cor:BR}, \ref {lem:BR}, and Corollary~\ref{BR+} remain valid when we replace $\R^\times \times [0,\pi)$ everywhere by $\R^{>}\times [0,2\pi)$. (Use that $\cos(\theta+\pi)=-\cos(\theta)$ for~$\theta\in \R$.)
\end{remarkNumbered}

\noindent
Call a germ in $\c$ {\em  eventually convex\/} if it has a convex representative in $\c_r$ for some~${r\in \R}$;
 likewise with ``concave'' in place of ``convex''.
The two lemmas below comprise a slightly weaker version of \cite[Theorem~2]{Horsley}.
By Lemma~\ref{lem:Bessel expansion} we have $\phi=x+\alpha x^{-1}+O(x^{-3})$ where $\alpha:=\frac{1}{2}(\nu^2-\frac{1}{4})$, so
with~$\phi$ being hardian we obtain
$$\phi'\ =\ 1-\alpha x^{-2}+O(x^{-4}),\qquad
  \phi''\ =\ 2\alpha x^{-3}+O(x^{-5}).$$
Hence $\phi''>0$ if~$\nu^2>\frac{1}{4}$ and $\phi''<0$ if~$\nu^2<\frac{1}{4}$, and thus:

\begin{lemma}
$\phi$ is eventually convex  if $\nu^2>\frac{1}{4}$, and
eventually concave if $\nu^2<\frac{1}{4}$.
\end{lemma}

\noindent
Lemma~\ref{lem:param V'} yields $(q,\theta)\in\Dx(\Q)^2$ parametrizing $V':=\{y':y\in V\}$.

\begin{lemma}\label{lem:ev convex}
We have $\theta-x-r-\textstyle\frac{1}{2}(\nu^2+\frac{3}{4})x^{-1} \preceq x^{-3}$ for some $r\in \R$.
Hence~$\theta$ is eventually convex.
\end{lemma}
\begin{proof}
Set~$g:=(x\phi')^{-1/2}, q:=\sqrt{(g')^2+(g\phi')^2}$, so $g,q\in\Dx(\Q)$.
The proof of Lemma~\ref{lem:param zeros of derivatives}(iv) and Lemmas~\ref{lem:param V'} and~\ref{lem:g,phi unique} give
$\theta=\phi+d+u$ with $d\in\R$ and~$u=\arccos( g'/q)$.  
Now $\phi'{}^\dagger=2\alpha x^{-3}+O(x^{-5})$ and so
$g^\dagger=-\frac{1}{2}x^{-1}+O(x^{-3})$.
Since $(\phi')^2 = 1+O(x^{-2})$, this yields $((g^\dagger)^2+(\phi')^2)^{-1/2} = 1 + O(x^{-2})$ and thus
$$g'/q\ =\ g^\dagger 
\big((g^\dagger)^2+(\phi')^2\big)^{-1/2}\  =\ -\textstyle\frac{1}{2}x^{-1}+O(x^{-3}).$$
 Hence
$$(g'/q)'\ =\ \textstyle\frac{1}{2}x^{-2}+O(x^{-4}),\qquad \big(1-(g'/q)^2\big)^{-1/2}\ =\ 1+O(x^{-2}).$$
 We obtain
$$u'\ =\ -\frac{ (g'/q)' }{\sqrt{1-(g'/q)^2}}\ =\ -\textstyle\frac{1}{2}x^{-2}+O(x^{-4})$$
and thus $u=c+\textstyle\frac{1}{2}x^{-1}+O(x^{-3})$ with $c\in\R$,
and so
$\theta-x-r-(\alpha+\textstyle\frac{1}{2})x^{-1}\preceq x^{-3}$ for $r:=c+d$,
as claimed.
\end{proof}

\subsection*{Asymptotic expansions for $\phi$ and $\phi^{\operatorname{inv}}$}
The arguments in this subsection demonstrate the efficiency of our transfer theorems from Section~\ref{sec:transfer}.
They allow us to produce  hardian   solutions of    algebraic differential equations from transseries  solutions of these equations. 
Such transseries solutions may be
constructed by purely formal computations in $\mathbb T$ (without convergence considerations). 
Our first goal is to improve on the relation $\phi\sim x + \frac{\mu-1}{8}x^{-1}$ from Lemma~\ref{lem:Bessel expansion}, where~$\mu:=4\nu^2$: 

\begin{theorem} \label{thm:asymp exp phi}
The germ $\phi=\phi_\nu$ has an asymptotic expansion
$$\phi\ \sim\ 
x+\frac{\mu-1}{8}x^{-1}+\frac{\mu^2-26\mu+25}{384}x^{-3}+\frac{\mu^3-115\mu^2+1187\mu-1073}{5120} x^{-5}+\cdots $$ 
\end{theorem}

\noindent
Here we use the sign $\sim$ not in the sense of comparing germs, but to indicate an 
{\em asymptotic expansion}: for a sequence $(g_n)$ in $\Calinf[\imag]$ with $g_0\succ g_1\succ g_2 \succ\cdots$ we say that $g\in \Calinf[\imag]$ has
the asymptotic expansion $$g\  \sim\ c_0g_0 + c_1g_1 + c_2g_2+\cdots \qquad(c_0, c_1, c_2,\dots\in\C)$$ if
$g-(c_0g_0+\cdots+ c_ng_n)\prec g_{n}$ for all $n$ (and then the sequence $c_0, c_1, c_2,\dots$ of coefficients is uniquely determined by
$g, g_0, g_1, g_2,\dots$).

\medskip\noindent
In the course of the proof of  Theorem~\ref{thm:asymp exp phi} we also obtain an explicit  formula for the coefficient
of $x^{-2n+1}$ in the asymptotic expansion of the theorem. Towards the proof, set 
$$(\nu,n)\	:=\ \frac{(\mu-1^2)(\mu-3^2)\cdots(\mu-(2n-1)^2)}{n!\,2^{2n}}\qquad  \text{\textup{(}Hankel's symbol\textup{)},}$$
so $(\nu,0)=1$, $(\nu,1)=\frac{\mu -1}{4}$, and $(\nu,n)=(-\nu,n)$. Also, 
if $(\nu,n)=0$ for some $n$, then~$\nu\in\frac{1}{2}+\Z$:
and 
if $\nu=\frac{1}{2}+m$, then~$(\nu,n) =  0$ for  $n\geq m+1$.
Moreover, if~$\nu\notin\frac{1}{2}+\Z$, then in terms of Euler's Gamma function (cf.  \cite[XV, \S{}2, $\Gamma$3, $\Gamma$5]{LangCA}), 
\begin{align*} (\nu,n)\  &=\ \frac{(-1)^n}{\pi n!} \cos(\pi\nu) \,\Gamma(\textstyle\frac{1}{2}+n-\nu)\,\Gamma(\textstyle\frac{1}{2}+n+\nu), \text{ and so}\\
(m,n)\ & =\  \frac{(-1)^{m+n}}{\pi n!} \Gamma(\textstyle\frac{1}{2}+n-m)\,\Gamma(\textstyle\frac{1}{2}+n+m).
\end{align*} 
\big(To prove the first identity, use $\Gamma(z+1)=z\Gamma(z)$ $n$ times to give
\begin{align*} \Gamma(\textstyle\frac{1}{2}+n-\nu)\ &=\ \Gamma(\textstyle\frac{1}{2}-\nu+n)\ =\ 
\Big(\!\prod_{j=0}^{n-1}(\frac{1}{2}-\nu+j)\Big)\cdot \Gamma(\frac{1}{2}-\nu); \text{ likewise}\\
\Gamma(\textstyle\frac{1}{2}+n+\nu)\ &=\  \Gamma(\textstyle\frac{1}{2}+\nu+n)\ =\ 
  \Big(\!\prod_{j=0}^{n-1}(\frac{1}{2}+\nu+j)\Big)\cdot \Gamma(\frac{1}{2}+\nu),
 \end{align*}
 and then use $\Gamma(z)\Gamma(1-z)=\frac{\pi}{\sin(\pi z)}$ to get $\Gamma(\frac{1}{2}-\nu)\Gamma(\frac{1}{2}+\nu)=\frac{\pi}{\cos(\pi \nu)}$.\big)

\medskip
\noindent
Below we consider the $H$-subfield $\R(\!(x^{-1})\!)$ of $\mathbb T$, and set
 $$y\ :=\  \sum_{n=0}^\infty y_n x^{-2n}\in\R(\!(x^{-1})\!)\ \text{ where }y_n\ :=\  
 (2n-1)!! \frac{ (\nu,n) }{2^n}.$$ 
Here $(2n-1)!! :=1\cdot 3\cdot 5\cdots (2n-1)=\frac{(2n)!}{2^n n!}$, so $(-1)!!=1$.
Thus 
\begin{multline*}
y\ =\ 1+ \left(\frac{\mu-1}{8}\right)  x^{-2}+\frac{3!!}{2!}\left( \frac{\mu-1}{8}\right)\left( \frac{\mu-9}{8}\right) x^{-4}+{} \\
\frac{5!!}{3!}\left( \frac{\mu-1}{8}\right)\left( \frac{\mu-9}{8}\right)\left( \frac{\mu-25}{8}\right)x^{-6}
+\cdots.
\end{multline*}
The definition of the $y_n$ yields the recursion
\begin{equation}\label{eq:yn rec}
y_0\ =\ 1\qquad\text{and}\qquad y_{n+1}\ =\ \left(\frac{2n+1}{n+1}\right)\left(\frac{\mu-(2n+1)^2}{8}\right)y_n.
\end{equation}
Using this recursion and $\Gamma(1/2)=\sqrt{\pi}$ for $n=0$,  induction on $n$ yields for $\nu\notin \frac{1}{2}+\Z$: 
$$y_n\  =\   \frac{\Gamma(n+\frac{1}{2})\,\Gamma(\nu+\frac{1}{2}+n)}{n!\,\sqrt{\pi}\,\Gamma(\nu+\frac{1}{2}-n)}.$$
We now verify that $y$ satisfies the linear differential equation
$$Y'''+fY'+(f'/2)Y\ =\ 0.$$
Here
$f=4+(1-\mu)x^{-2}$, so $f'/2=(\mu-1)x^{-3}$.
Thus
$$(f'/2)y\  =\  \sum_n   (\mu-1) y_n   x^{-2n-3}\ =\ (\mu-1)x^{-3}+\sum_{n\geq 1}  (\mu-1) y_n  x^{-2n-3}. $$
We also have
$y'=\sum_{n\geq 1} -2ny_n x^{-2n-1}$
and so
\begin{align*}
fy' \ &=\  \big(4+(1-\mu)x^{-2}\big) \sum_{n\geq 1}  -2ny_n  x^{-2n-1} \\ 
& =\  \sum_{n\geq 1}  2n(\mu-1)y_n x^{-2n-3} -  \sum_{n\geq 1}  8ny_n  x^{-2n-1}   \\
&=\  \sum_{n\geq 1}  2(\mu-1)ny_n x^{-2n-3} - \sum_{m\geq 0} 8(m+1)y_{m+1} x^{-2m-3}\\
& =\  -8y_1 x^{-3} + \sum_{n\geq 1} 
\big(2(\mu-1)n y_n - 8(n+1)y_{n+1}\big)  x^{-2n-3}
\end{align*}
and hence, using $\mu-1=8y_1$:
$$fy'+(f'/2)y\ =\   \sum_{n\geq 1} 
 \big( (2n+1)(\mu-1)y_n-8(n+1)y_{n+1}\big)  x^{-2n-3}.$$
Moreover
$$y''\ =\  \sum_{n\geq 1}  2n (2n+1) y_n   x^{-2n-2},\quad
  y'''\ =\ \sum_{n\geq 1}  -4 n(2n+1)(n+1) y_n   x^{-2n-3}.$$
This yields the claim by \eqref{eq:yn rec}. We now identify the Hardy field $\R(x)$ with an $H$-subfield of
$\mathbb T$ in the obvious way. Then $\R(x)\subseteq \R(\!(x^{-1})\!)$, and the above yields:

\begin{lemma}\label{lem:kerB, Rx}
Let $H$ be an $H$-closed field extending the $H$-field $\R(x)$ and set~$B:=\der^3+f\der+(f'/2)\in \R(x)[\der]$. 
 Then $\dim_{C_H} \ker_H B=1$.
If $H$ extends the $H$-field~$\R\langle x,y\rangle\subseteq \mathbb T$, then   $\ker_{H} B=~C_H y$.
\end{lemma}
\begin{proof}
We have $\sigma(1/x)=2/x^2 \prec 4 \sim f$, so $f\in\sigma\big(\Upg(H)\big){}^\uparrow$,
hence $f\in\sigma(H^\times)\setminus\omega(H)$. Thus $\dim_{C_H} \ker_H B=1$ by Lemma~\ref{lem:kerB}.
Hence if $H$ extends the $H$-field $\R\langle x,y\rangle$, then $\ker_{H} B=C_H y$  since $B(y)=0$ by
the argument preceding the lemma.
\end{proof}

\begin{prop}\label{prop:3rd order sol}
There is a unique hardian germ $\psi=\psi_\nu$ such that 
\begin{equation}\label{eq:psi conditions}
\psi\ \sim\ 1\quad\text{ and }\quad \psi'''+f\psi'+(f'/2)\psi\ =\ 0.
\end{equation}
This $\psi$   satisfies $\psi=1/\phi'\in\Dx(\Q)$ and has the asymptotic expansion
\begin{equation}\label{eq:psi asymp exp}
\psi\ \sim\  1+ \frac{\mu-1}{8}   x^{-2}+ \cdots+
(2n-1)!! \frac{ (\nu,n) }{2^n}x^{-2n}+\cdots.
\end{equation}
Moreover $\psi_{-\nu}=\psi_\nu$, and
if  $\nu=\frac{1}{2}+m$, then 
$$\psi\ =\ 
 1+ \frac{\mu-1}{8}   x^{-2}+\cdots+ (2m-1)!! \frac{ (\nu,m) }{2^m}x^{-2m}.$$
\end{prop}
\begin{proof}
For any $H$-closed field $H\supseteq \R(x)$, consider the statement
$$\begin{cases}
&\parbox{25em}{there is a unique $\psi\in H$ such that \eqref{eq:psi conditions} holds, and this~$\psi$ satisfies $\psi-\displaystyle\sum_{n=0}^{m} y_nx^{-2n}\prec x^{-2m-1}$ for all $m$.}
\end{cases}$$
This
holds for  $H=\mathbb T$
by Lemma~\ref{lem:kerB, Rx}, and hence also for any $\d$-maximal Hardy field $H$ by Corollary~\ref{cor:transfer T, 1} (applied with $\R(x)$ in the role of $H$).
Thus  every $\d$-maximal Hardy field contains a unique germ $\psi$ satisfying~\eqref{eq:psi conditions},
and every such~$\psi$ has an asymptotic expansion~\eqref{eq:psi asymp exp}.

Now let $\psi$ be any hardian germ satisfying~\eqref{eq:psi conditions}.
Take  a $\d$-maximal Hardy field~$H$ containing $\psi$; then $\R\langle x,\phi\rangle\subseteq\Dx(\Q)\subseteq H$.
%Let
%$\Univ$ be the universal  exponential extension of~$K=H[\imag]\subseteq \Calinf[\imag]$, and identify $\Univ$ with a differential subring of $\Calinf[\imag]$ as explained
%at the beginning of Section~\ref{sec:ueeh}. Let $A:=4\der^2+f$ and 
Let $B$ be  as in Lemma~\ref{lem:kerB, Rx}. Then Lemma~\ref{lem:B, 1} gives
% $(\phi')^{-1/2}\ex^{\phi\imag}, (\phi')^{-1/2}\ex^{-\phi\imag}\in\ker_{\Univ}A$,
 $B(1/\phi')=0$. Since~$1/\phi'\sim 1$,
this yields~$\psi=1/\phi'\in\Dx(\Q)$.   
For the rest, use that $\phi_{\nu}=\phi_{-\nu}$ and that for $\nu=\frac{1}{2}+m$ we have $y_n=0$ for~$n\geq m+1$.  
\end{proof}

\begin{cor}\label{cor:3rd order sol}
Let   $\psi=\psi_\nu$   and suppose $\mu\neq 1$.  Then
$$\psi^{(n)}\ \sim\  (-1)^n(n+1)!\left(\frac{\mu-1}{8}\right)x^{-n-2}\quad\text{for $n\geq 1$.}$$
In particular, $\psi$ is eventually strictly increasing if $\nu^2<1/2$ and eventually strictly decreasing if $\nu^2>1/2$.
\end{cor}

\begin{lemma}\label{lem:asympt exp}
Let $H$ be a Hausdorff field extension of $\R(x)$,  $e\colon H\to\mathbb T$ a valued field embedding 
over $\R(x)$, 
and~$h\in H$, $e(h)\in\R(\!(x^{-1})\!)$,
say $$e(h)\ =\ h_{k_0} x^{-k_0}+h_{k_0+1}x^{-k_0-1}+\cdots\qquad\text{where $k_0\in\Z$ and $h_k\in\R$ for $k\geq k_0$.}$$
Then $h$ has an asymptotic expansion 
$$h\ \sim\  h_{k_0} x^{-k_0}+h_{k_0+1}x^{-k_0-1}+\cdots.$$
Moreover, if $e(h)\in\R[x,x^{-1}]$, then $e(h)=h$.
\end{lemma}
\begin{proof}
For  $k\geq k_0$ we set $h_{\leq k}:=h_{k_0} x^{-k_0}+\cdots+h_k x^{-k}\in\R[x, x^{-1}]$. Then    
$$e (h-h_{\leq k} )\ =\  e(h)-h_{\leq k}=h_{k+1}x^{-k-1}+\cdots \prec x^{-k}\ =\ e(x^{-k}),$$
so $h-h_{\leq k}\prec x^{-k}$, giving the  asymptotic expansion. If $e(h)\in \R[x,x^{-1}]$, take
$k\ge k_0$ so large that $e(h)=h_{\leq k}$, and then $e(h)=e(h_{\leq k})$, so $h=h_{\leq k}\in \R[x,x^{-1}]$.
\end{proof}

\noindent
We   now prove Theorem~\ref{thm:asymp exp phi}, using also results and notations from the Appendix to this section.  Corollary~\ref{cor:transfer T, 2} yields an $H$-field embedding $e\colon\Dx(\Q)\to\mathbb T$ over~$\R(x)$. Let $\psi$ be as in Proposition~\ref{prop:3rd order sol}.
Then~$e(\psi)=y$ and so
$$e(\phi')\ =\ e(1/\psi)\ =\ y^{-1}\ =\ z_0+ z_1\frac{x^{-2}}{1!}+z_2\frac{x^{-4}}{2!}+\cdots+z_n\frac{x^{-2n}}{n!}+\cdots$$
where 
$$z_n\ :=\  B_n(-y_1 1!,\dots,-y_n n!)\in\Q[y_1,\dots,y_n]\subseteq\Q[\mu]$$  by Lemma~\ref{lem:inv of 1-y} at the end of this section; here the $B_n$ are as
defined in \eqref{eq:Bell polys}.  
Using Lemma~\ref{lem:asympt exp} and $\phi\sim x$ we obtain the asymptotic expansion
$$\phi\ \sim\  u_0x+u_1\frac{x^{-1}}{1!}+u_2\frac{x^{-3}}{2!}+\cdots+ u_n\frac{x^{-2n+1}}{n!}+\cdots \quad\text{where
$u_n:=  \frac{z_n}{-2n+1}$.}$$
The first few terms of the sequence $(z_n)$ are
\begin{align*}
z_0\	&=\  1, \\
z_1\	&=\  -y_1 =\frac{-(\mu-1)}{8} ,\\
z_2\	&=\  -2y_2+2y_1^2 = \frac{ -3 (\mu-1)(\mu-9)+2(\mu-1)^2}{64}, \\
z_3\	&=\  -6y_3+12y_1y_2-6y_1^3 \\
	&=\  \frac{-15 (\mu-1)(\mu-9)(\mu-25)+18(\mu-1)^2(\mu-9)-6(\mu-1)^3}{512}  
\end{align*}
and so
\begin{align*}
u_0\ &=\  1, \\
u_1\	&=\  \frac{\mu-1}{8} ,\\
u_2\	&=\  \frac{  3 (\mu-1)(\mu-9)-2(\mu-1)^2}{192} = \frac{\mu^2-26\mu+25}{192} \\
u_3\	&=\  \frac{ 15 (\mu-1)(\mu-9)(\mu-25)- 18(\mu-1)^2(\mu-9)+6(\mu-1)^3}{2560}  \\
	&=\  \frac{3(\mu^3-115\mu^2+1187\mu-1073)}{2560}.
\end{align*}
This finishes the proof of Theorem~\ref{thm:asymp exp phi}. \qed

\medskip
\noindent
We  turn to the compositional inverse $\phi^{\operatorname{inv}}$ of $\phi$. Recall: $\phi^{\operatorname{inv}}\in\Dx(\Q)$ by Lemma~\ref{phiinvexq}. To prove the next result we use Corollary~\ref{cor:LIF, x^-1} in the Appendix to this section. 

\begin{cor} \label{cor:phiinv}
We have an asymptotic expansion
$$\phi^{\operatorname{inv}}\ \sim\ x- \frac{\mu-1}{8}x^{-1}- \frac{(\mu-1)(7\mu-31)}{192}  \frac{x^{-3}}{2!}+\cdots$$
\end{cor}
\begin{proof}
Let $e\colon\Dx(\Q)\to\mathbb T$ and 
$(u_n)$ be as above. Set
$$u:=\sum_n u_n\frac{x^{-2n+1}}{n!}=x+ \frac{\mu-1}{8}x^{-1}+\frac{\mu^2-26\mu+25}{384}x^{-3}+\cdots\in\R(\!(x^{-1})\!)\subseteq\mathbb T,$$
so $e(\phi)=u$.
Let~$P\in\R(x)\{Y\}$ be as   in the proof of Lemma~\ref{phiinvexq}, so $P(2u')=e(P(2\phi'))=0$.  
Corollary~\ref{cor:trdegellinv} 
and the remark following it yield a
$\tilde P\in \R(x)\{Z\}$ such that for all hardian $y>\R$,
$$P(2y')\ =\ 0\quad \Longleftrightarrow\quad \tilde P(y^{\operatorname{inv}})\ =\ 0$$
and such that this equivalence also holds for $y\in \mathbb T^{>\R}$ and $y^{\operatorname{inv}}$   the compositional inverse
of $y$ in $\mathbb T$. Hence $\tilde P(e(\phi^{\operatorname{inv}}))=e(\tilde P(\phi^{\operatorname{inv}}))=0$ and $\tilde P(u^{\operatorname{inv}})=0$.
The proof of Lemma~\ref{phiinvexq}  shows that each maximal Hardy field~$H$ contains a unique zero $z$ of~$\tilde P$
such that $z=x-\frac{1}{8}(\mu-1)x^{-1}+o(x^{-1})$.  By Corollary~\ref{cor:transfer T, 1} this remains true with~$\mathbb T$ in place of $H$.
Now~$e(\phi^{\operatorname{inv}})=x-\frac{1}{8}(\mu-1)x^{-1}+o(x^{-1})$, and by the remarks following 
Corollary~\ref{cor:LIF, x^-1} we have
$u^{\operatorname{inv}}=u^{[-1]}= x-\frac{1}{8}(\mu-1)x^{-1}+o(x^{-1})$. Hence~$e(\phi^{\operatorname{inv}})=u^{\operatorname{inv}}$ and thus $\phi^{\operatorname{inv}}$ has an asymptotic expansion as claimed. 
%(See the remarks following Corollary~\ref{cor:McMahon} below
%for a formula giving the terms in this asymptotic expansion.)
\end{proof}

\begin{remark} Corollary~\ref{cor:LIF, x^-1} yields the more detailed asymptotic expansion
$$\phi^{\operatorname{inv}}\ \sim\  x-\sum_{j=1}^\infty g_j\frac{x^{-2j+1}}{j!},\quad\text{where $g_j= \sum_{i=1}^{j}\frac{(2(j-1))!}{(2j-1-i)!}B_{ij}(u_1,\dots,u_{j-i+1})$.}$$
\end{remark} 

\subsection*{Liouvillian phase functions}
The next proposition adds to Corol\-lary~\ref{cor:2nd order, Liouvillian solutions} for the differential  equation~\eqref{eq:Bessel selfadj}. This subsection is not used in the rest of the section. 

\begin{prop}\label{prop:half-integer}
With $\phi=\phi_\nu$, the following are equivalent:
\begin{enumerate}
\item[\textup{(i)}] $\nu\in\frac{1}{2}+\Z$;
\item[\textup{(ii)}] $1/\phi'\in \R[x^{-1}]$;
\item[\textup{(iii)}]  $f\in\sigma\big(\R(x)^>\big)$; recall: $f=4+(1-\mu)x^{-2}$; 
\item[\textup{(iv)}] $\phi\in\Li\!\big(\R(x)\big)$;
\item[\textup{(v)}] $x^2y''+xy'+(x^2-\nu^2)y=0$ for some $y\ne 0$ in a Liouville extension of $\C(x)$;
%\item[\textup{(v)}] $\phi'\in \R(x)$;
\item[\textup{(vi)}] there are $a,b\in\C$ and distinct $c_1,\dots,c_n\in\C^\times$ such that 
$$-\phi'{}^\dagger+2\phi'\imag = a+bx^{-1}+2\sum_{i=1}^n (x-c_i)^{-1}\quad\text{and}\quad b=1+2\nu\text{ or }b=1-2\nu.$$
%\item[\textup{(vii)}] there is some $y\neq 0$ in a Liouville extension of    the differential field~$\C(x)$ with  $x^2y''+xy'+(x^2-\nu^2)y=0$.
\end{enumerate}
\end{prop}
\begin{proof}
The implication (i)~$\Rightarrow$~(ii) follows from
  Proposition~\ref{prop:3rd order sol},  and (ii)~$\Rightarrow$~(iii) from $f=\sigma(2\phi')$.
If $f\in\sigma\big(\R(x)^\times\big)$, then~$\phi\in\Li\!\big(\R(x)\big)$ by  Corollary~\ref{cor:parametrization in H} with~$H:=\Li\!\big(\R(x)\big)$; thus (iii)~$\Rightarrow$~(iv).
For the rest of the proof, recall that by Lemma~\ref{parphi} the pair~$\big(1/\sqrt{\phi'},\phi\big)$ parametrizes~$\ker_{\Calinf}(4\der^2+f)$, so $\big(1/\sqrt{x\phi'},\phi\big)$ parametrizes $V_{\nu}$. Thus (iv)~$\Leftrightarrow$~(v) by Corollary~\ref{cor:2nd order, Liouvillian solutions}. Moreover, if
%$$y\ :=\ (\phi')^{-1/2}\ex^{\phi\imag},\qquad z\ :=\ 2y^\dagger\ =\ -\phi'{}^\dagger+2\phi'\imag\in\C(x).$$
 (iv) holds, then~$\big(1/\sqrt{\phi'},\phi\big)\in\Li\!\big(\R(x)\big)^2$, so (vi) then follows from Corollary~\ref{cor:phi' in H, 2}.
% yields~(iv)~$\Rightarrow$~(v). 

Suppose $a,b,c_1,\dots,c_n$ are as in (vi) and set
$$y\ :=\ (\phi')^{-1/2}\ex^{\phi\imag},\qquad z\ :=\ 2y^\dagger\ =\ -\phi'{}^\dagger+2\phi'\imag\in\C(x),$$
so $4y''+fy=0$, hence $\omega(z)=f$. Then, as germs at $+\infty$, 
%Equip  $\C(x)$ with the valuation~$v_\infty$.
%\colon \C(x)^\times\to\Z$
%which is trivial on $\C$ and $v_\infty(x)=-1$, with valuation ring $\mathcal O$.
$$z\ =\  a+(b+2n)x^{-1} + O(x^{-2}), \text{ so }
z'\ =\ O(x^{-2}),\quad 
z^2\ =\  a^2+2a(b+2n)x^{-1}+O(x^{-2})$$
and hence
$$f\ =\ 4+(1-\mu)x^{-2}\ =\ \omega(z)\ =\ -(2z'+z^2)\ =\  -a^2-2a(b+2n)x^{-1}  +O(x^{-2}),$$
so $b+2n=0$, hence $\nu=-n-\frac{1}{2}$ or $\nu=n+\frac{1}{2}$, and thus $\nu\in\frac{1}{2}+\Z$.
\end{proof}

\begin{remark} 
In the setting of analytic functions, the above equivalence (i)~$\Leftrightarrow$~(v)  goes back to
Liouville~\cite{Liouville41}. For more on this, see \cite[appendix]{Kolchin68}, \cite[\S{}4.2]{Kovacic}, \cite[Chapter~VI]{Ritt}, and \cite[\S{}4.74]{Watson}.)
\end{remark}

\noindent
For the next result, note that $\arctan(g)\in \Li\!\big(\R(x)\big)$ for $g\in \R(x)$. 

\begin{cor}
Suppose $\nu\in\frac{1}{2}+\Z$.
%$\phi\in\Li\!\big(\R(x)\big)$.  
Then there are  distinct   $(a_1,b_1),\dots,(a_m,b_m)$ in~$\R^\times\times\R$ such that
$$\phi\ =\ x+\sum_{i=1}^m \arctan\left(\frac{a_i}{x-b_i}\right).$$
\end{cor}
\begin{proof} Take imaginary parts in the equality of Proposition~\ref{prop:half-integer}(vi), integrate, and appeal to the defining property of $\phi$ in  Proposition~\ref{prop:Bessel} in combination with the fact that for $a,b\in \R$ we have
 $\arctan\big(\frac{a}{x-b}\big)\preceq x^{-1}$. Here we also use that the derivative of $\arctan\big(\frac{a}{x-b}\big)$ is 
 $\frac{-a}{(x-b)^2+a^2}=\Im\big(\frac{1}{x-c}\big)$ for  $a,b\in \R$, $c=b-a\imag$.   
\end{proof}

\noindent
Is $\phi, \phi^{\operatorname{inv}}\in\Li\!\big(\R(x)\big)$ possible? The answer is ``no'' except for $\phi= x$:

\begin{cor}\label{cor:phinv not in Li}
Suppose $\phi\in\Li\!\big(\R(x)\big)$, $\phi\neq x$, and $\nu\geq 0$, so $\nu=\frac{1}{2}+m$ where~$m\geq 1$. Then $\theta:=1/\phi^{\operatorname{inv}}$  satisfies
$$\theta' = -\theta^2(1+y_1\theta^2+\cdots+ y_m\theta^{2m})\quad\text{where $y_i=(2i-1)!!\frac{(\nu,i)}{2^i}$ for $i=1,\dots,m$}$$
and   $\theta\notin  \Li\!\big(\R(x)\big)$.
\end{cor}
\begin{proof}
By the Chain Rule and Proposition~\ref{prop:3rd order sol} we have 
$$\theta'=-\theta^2  (\phi^{\operatorname{inv}})'=-\theta^2 (\psi \circ \phi^{\operatorname{inv}})\quad\text{where
$\psi =  1+ y_1 x^{-2}+\cdots+ y_m x^{-2m}$,}$$
and this yields the first claim. Towards a contradiction, assume $\theta \in  \Li\!\big(\R(x)\big)$.
Then by [ADH, 10.6.6], $\theta$ lies in the Liouville extension $\Li\!\big(\R(x)\big)[\imag]$ of $\C(x)=\R(x)[\imag]$.
%As $\theta\notin \C$, this is impossible by Corollary~\ref{cor:Abel equ, 2}. 
Hence by Corollary~\ref{cor:Srinivasan}, $\theta$ is algebraic over $\C(x)$.
Also $\theta\notin\C$, so Lemma~\ref{lem:Srinivasan} yields $Q\in\C(Y)$ with $Q'=1/P$
where $P:=-Y^2(1+y_1Y^2+\cdots+y_mY^{2m})$. 
Thus 
$$\phi'\ =\ \frac{1}{\psi}\ =\ -\frac{x^{-2}}{P(x^{-1})}\ =\ -x^{-2}Q'(x^{-1})\ =\ Q(x^{-1})'$$ 
and so $\phi\in\R(x)$. This is impossible by the lemma below. 
\end{proof}

\begin{lemma}
$\mu=1$ $\Longleftrightarrow$ $\phi=x$ $\Longleftrightarrow$   $\phi\in\R(x)$.
\end{lemma}
\begin{proof}
The first equivalence is clear from the remarks following Proposition~\ref{prop:Bessel}.
Assume $\phi\in\R(x)$. 
Then by  Proposition~\ref{prop:half-integer},
$$-\phi'{}^\dagger+2\phi'\imag\ =\ a+bx^{-1}+2\sum_{i=1}^n (x-c_i)^{-1}\quad(a,b\in \C, \text{ distinct }c_1,\dots, c_n\in \C^\times),$$
so $2\phi'\imag-a\in \der F\cap \C F^\dagger$ for $F:= \C(x)$, hence $2\phi'\imag=a$ by Corollary~\ref{cor:linindep logs}.  Thus~$\phi\in\R x+\R$.
Since $\phi\sim x$ and~$\phi-x\preceq x^{-1}$, this gives $\phi=x$.
\end{proof}

\begin{question}
Does there exist a $\nu\notin\frac{1}{2}+\Z$ for which
$\phi^{\operatorname{inv}}\in  \Li\!\big(\R(x)\big)$?  
\end{question}

\subsection*{The Bessel functions}
We can now establish some classical facts about distinguished solutions to the Bessel differential 
equation~\eqref{eq:Bessel}: Corollaries~\ref{cor:Bessel 1}, \ref{cor:Bessel 2}, \ref{cor:Nicholson} below.   Our proofs use less complex analysis than those in the literature: we need just one contour integration, for Proposition~\ref{prop:Hankel} below. 
We assume some  basic facts about  Euler's~$\Gamma$-function and recall that $1/\Gamma$
is an entire function with~${-\N}$ as its set of zeros (all simple), so $\Gamma$ is meromorphic on the complex plane
without any zeros and has $-\N$ as its set of poles. Our main reference for these and other properties of $\Gamma$ used below is~\cite{LangCA}. 
Let also~$z\mapsto\log z\colon  \C\setminus\R^{\leq}\to\C$  be the holomorphic extension of the real logarithm function, and
for~$z\in \C\setminus\R^{\leq}$,   set~$z^\nu:=\exp(\nu\log z)$. Let~${\nu\in\C}$ until  further notice, and note that
$(\nu,z) \mapsto z^{\nu}$ is analytic on~$\C\times (\C\setminus \R^{\le})$, and, keeping $\nu$ fixed, has derivative
$z\mapsto \nu z^{\nu-1}$ on $\C\setminus \R^{\le}$. For $z\in \C\setminus \R^{\le}$
(so $z=|z|\ex^{\imag
\theta},\ -\pi < \theta< \pi$), and $\nu, \nu_1, \nu_2\in \C$, $t\in \R^{>}$ we have $$z^{\nu_1+ \nu_2}\ =\ z^{\nu_1}z^{\nu_2},\quad (tz)^\nu\ =\ t^\nu z^\nu,\quad |z^\nu|\ =\ |z|^{\Re \nu}\ex^{-\Im(z)\theta},\quad 
\bar{z^{\nu}}\ =\ \bar{z}^{\bar{\nu}},$$
and for $z_1, z_2\in \C^\times$ with $\Re z_1\ge 0$, $\Re z_2>0$: $\ z_1z_2\in \C\setminus \R^{\le}$, $(z_1z_2)^\nu=z_1^\nu z_2^\nu$. 

\begin{lemma}\label{lem:Bessel conv}
Let $A,B\subseteq\C$  be nonempty and compact. Then 
\begin{align*} \sum_n &\max_{(\nu,z)\in A\times B} \left|\frac{(-1)^n}{n!\,\Gamma(\nu+n+1)}\left(\frac{z}{2}\right)^{2n}\right| < \infty, \text{ so the series}\\
&\sum_n \frac{(-1)^n}{n!\,\Gamma(\nu+n+1)}\left(\frac{z}{2}\right)^{2n}
\end{align*} converges absolutely  and uniformly on $A\times B$.
\end{lemma} 
\begin{proof}
Take $R\in\R^>$ such that $\abs{z}\leq 2R$ for all $z\in B$. Set~$M_n:=\max\limits_{\nu\in A}  \left|\frac{1}{\Gamma(\nu+n+1)}\right|$ and
take $n_0\in \N$   such that $|\nu+n+1|\ge 1$ for all $n\ge n_0$ and $\nu\in A$. Then~$M_{n+1}\le M_n$ for $n\ge n_0$, by 
 the functional equation for
$\Gamma$, so the sequence $(M_n)$ is bounded. Hence $\sum\limits_n \max\limits_{(\nu,z)\in A\times B} \left|\frac{(-1)^n}{n!\,\Gamma(\nu+n+1)}\left(\frac{z}{2}\right)^{2n}\right| \le \sum\limits_n \frac{M_n R^{n}}{n!}<\infty$.
\end{proof}

\noindent
By Lemma~\ref{lem:Bessel conv} and \cite[V, \S{}1, Theorem~1.1]{LangCA} we obtain a holomorphic function
$$z\mapsto J_\nu(z)\ :=\   \sum_n \frac{(-1)^n}{n!\,\Gamma(n+\nu+1)}\left(\frac{z}{2}\right)^{2n+\nu} \colon \C\setminus \R^{\leq} \to\C.$$ 
For example, for $z\in \C\setminus\R^{\leq}$, we have 
$$J_{\frac{1}{2}}(z)\ =\   \sum_{n} \frac{(-1)^n}{n!\,\Gamma(n+\frac{3}{2})}\left(\frac{z}{2}\right)^{2n+\frac{1}{2}}.$$ 
%\quad J_0(z) = \sum_{n\geq 0} \frac{(-1)^n}{(n!)^2}\left(\frac{z}{2}\right)^{2n}, $$
Note also that
\begin{equation}\label{eq:J-m}
J_{-m}(z)\  =\   \sum_{n\geq m} \frac{(-1)^n}{n!\,(n-m)!}\left(\frac{z}{2}\right)^{2n-m}\  =\  (-1)^m J_m(z) 
\end{equation}
and thus
$$J_{-m}(z)\sim \frac{(-1)^m}{m!}\left(\frac{z}{2}\right)^m,\quad J_{m}(z)\sim \frac{1}{m!}\left(\frac{z}{2}\right)^m \qquad\text{as $z\to 0$.}$$
Termwise differentiation shows that  $J_\nu$ satisfies the differential equation  \eqref{eq:Bessel} on ${\C\setminus\R^{\leq}}$. 
The function~$J_\nu$ is known as the {\it Bessel function of the first kind}\/ of order~$\nu$.\index{function!Bessel}\index{Bessel function!first kind}
Note that \eqref{eq:Bessel} doesn't change when replacing~$\nu$ by $-\nu$, so $J_{-\nu}$ is also
a solution of~\eqref{eq:Bessel}. Lem\-ma~\ref{lem:Bessel conv} shows that the function
$$(\nu,z)\mapsto J_{\nu}(z)\ :\ \C\times (\C\setminus\R^{\leq}) \to \C$$ is analytic, and that for fixed $\nu$ the function   $z\mapsto z^{-\nu}J_\nu(z)$  on $\C\setminus \R^{\le}$ extends to an entire function.

\medskip
\noindent
Termwise differentiation gives $(z^\nu J_\nu)'  =  z^\nu J_{\nu-1}$ and 
  $(z^{-\nu} J_\nu)' =  -z^{-\nu} J_{\nu+1}$, so 
\begin{equation}\label{eq:Jnu diff}
 J_{\nu-1}\ =\ \frac{\nu}{z}J_{\nu} + J_{\nu}',  \qquad 
  J_{\nu+1}\ =\ \frac{\nu}{z}J_{\nu} - J_{\nu}',
\end{equation}
 by the Product Rule,
and thus 
\begin{equation}\label{eq:Bessel rec fm}
J_{\nu-1}+J_{\nu+1}\  =\  \frac{2\nu}{z} J_\nu\,\qquad J_{\nu-1}-J_{\nu+1}\  =\ 2J_\nu'.
\end{equation}
Note: if $\nu+1\notin -\N$, then for $z\in \C\setminus \R^{\le}$ and $z\to 0$ we have
\begin{equation}\label{eq:Jnu as t->0+}
J_\nu(z)\sim \frac{1}{\Gamma(\nu+1)}\left( \frac{z}{2} \right)^{\nu},\  \text{ and for }\nu\ne 0:\  
J_\nu'(z) \sim \frac{\nu}{2\Gamma(\nu+1)}\left( \frac{z}{2} \right)^{\nu-1}.
\end{equation}
If $\nu-1\notin\N$, then \eqref{eq:Jnu as t->0+} holds with $-\nu$ in place of $\nu$.
It follows that for $\nu\notin\Z$ the solutions $J_\nu$, $J_{-\nu}$ of \eqref{eq:Bessel} are $\C$-linearly independent.
Set
$$w\ :=\ \operatorname{wr}(J_\nu,J_{-\nu})\ =\ J_\nu J'_{-\nu}-J_\nu' J_{-\nu}.$$ Then $w'(z)=-w(z)/z$ on
$\C\setminus\R^{\leq}$ (cf.~remarks following Lemma~\ref{lem:A from yjs}). This gives~${c\in\C}$ such that
$w(z)=c/z$ on $\C\setminus\R^{\leq}$. 
If $\nu\notin\Z$, then $c\ne 0$, and by \eqref{eq:Jnu as t->0+} and the remark following it we obtain 
$c=-\frac{2\nu}{\Gamma(\nu+1)\Gamma(-\nu+1)}=-\frac{2}{\Gamma(\nu)\Gamma(-\nu+1)}$. Hence
using~${\Gamma(\nu)\Gamma(1-\nu)} = \pi/\sin(\pi\nu)$: 
\begin{equation}\label{eq:wrJnuJ-nu}
\operatorname{wr}(J_\nu,J_{-\nu})(z)\ =\ -\frac{2\sin(\pi\nu)}{\pi z}\quad\text{for $\nu\notin\Z$ and $z\in\C\setminus \R^{\leq}$.}
\end{equation}
Next we express $J_{1/2}$ and $J_{-1/2}$ in terms of $\sin z$ and $\cos z$:
 
\begin{lemma}\label{lem:J1/2}
On  $\C\setminus\R^{\leq}$ we have
$$J_{1/2}(z)\ =\ \sqrt{\frac{2}{\pi z}}\ \sin z, \qquad
  J_{-1/2}(z)\ =\ \sqrt{\frac{2}{\pi z}}\ \cos z.$$
\end{lemma}
\begin{proof}
For $\nu=1/2$, a fundamental system of solutions of \eqref{eq:Bessel} on $\R^>$ is given by~$x^{-1/2}\cos x$,~$x^{-1/2}\sin x$.
This yields $a,b\in\R$ such that on $\R^{>}$,
$$J_{1/2}(t)\  =\  at^{-1/2}\cos t+bt^{-1/2}\sin t.$$
As $t\to 0^+$ we have:
\begin{align*} t^{-1/2}\cos t\ &=\  t^{-1/2}+O(t^{3/2}),\quad t^{-1/2}\sin t\ =\ t^{1/2}+O(t^{5/2}), \text{ and}\\
J_{1/2}(t) &\sim  \displaystyle \frac{1}{\Gamma(3/2)}\left( \frac{t}{2} \right)^{1/2}\ =\ \sqrt{\frac{2}{\pi}}\  t^{1/2}\  \text{ (using  \eqref{eq:Jnu as t->0+}}), 
\end{align*} so $a=0$, $b=\sqrt{\frac{2}{\pi}}$, giving the identity claimed for $J_{1/2}$. 
For $\nu=-1/2$ one can use the left identity \eqref{eq:Jnu diff} for $\nu=1/2$. 
 \end{proof}

\noindent
From Lemma~\ref{lem:J1/2} and \eqref{eq:Jnu diff}  we obtain by induction on $n$:

\begin{cor}
For each $n$ there are $P_n,Q_n\in\Q[Z]$ with $\deg P_n=\deg Q_n=n$, both with positive leading coefficient, such that, with $Q_{-1}:=0$,
we have on $\C\setminus\R^{\leq}$:
\begin{align*}
J_{n+\frac{1}{2}}(z)\	&=\ \sqrt{\frac{2}{\pi z}}\,\big( P_n(z^{-1})\sin z- Q_{n-1}(z^{-1}) \cos z\big), \\
J_{-n-\frac{1}{2}}(z)\	&=\ (-1)^n \sqrt{\frac{2}{\pi z}}\, \big( P_n(z^{-1})\cos z+ Q_{n-1}(z^{-1}) \sin z\big).
\end{align*}
\end{cor}

\noindent
For example,
$$J_{3/2}(z)\ =\ \sqrt{\frac{2}{\pi z}}\left(\frac{\sin z}{z}-\cos z\right).$$
 For $\nu\notin\Z$ we have the solution
$$  Y_\nu\  :=\ \frac{\cos(\pi\nu)J_\nu -J_{-\nu}}{\sin(\pi\nu)}$$
 of \eqref{eq:Bessel} on $\C\setminus\R^{\leq}$.
For fixed $z\in\C\setminus\R^{\leq}$, the entire function $$\nu\mapsto \cos(\pi\nu)J_\nu(z)-J_{-\nu}(z)$$ has a zero at each $\nu\in\Z$, by \eqref{eq:J-m}, so the holomorphic function
$$\nu\mapsto Y_\nu(z)\colon\C\setminus\Z\to\C$$ has a removable singularity at each $\nu\in\Z$, and thus extends to
an entire function whose value at $k\in \Z$ is given by
$$Y_k(z):=\lim_{\nu\in\C\setminus\Z,\ \nu\to k} Y_\nu(z), \quad \text{ for }z\in \C\setminus \R^{\le}.$$
In this way we obtain a two-variable analytic function
$$(\nu,z)\mapsto Y_{\nu}(z)\  :\  \C\times (\C\setminus \R^{\le})\to \C,$$ and thus
for each $\nu\in\C$ a solution $Y_\nu$
 of \eqref{eq:Bessel} on $\C\setminus\R^{\leq}$, called the {\it Bessel function of the second kind}\/ of order~$\nu$.\index{function!Bessel}\index{Bessel function!second kind}
Using \eqref{eq:wrJnuJ-nu} we determine the Wronskian of~$J_\nu$,~$Y_\nu$ (first for $\nu\notin\Z$, and then by continuity for all $\nu$):
\begin{equation}\label{eq:wrjy} \operatorname{wr}(J_\nu,Y_\nu)(z)\  =\  - \frac{\operatorname{wr}(J_\nu, J_{-\nu})(z)}{\sin(\pi\nu)}\ =\ \frac{2}{\pi z}\quad
(z\in\C\setminus\R^{\leq}),
\end{equation}
hence $J_\nu$, $Y_\nu$ are $\C$-linearly independent.
The recurrence formulas \eqref{eq:Bessel rec fm} yield analogous formulas for the Bessel functions of the second kind:
\begin{equation}\label{eq:Bessel rec fm, 2nd}
Y_{\nu-1}+Y_{\nu+1}\  =\  \frac{2\nu}{z} Y_\nu\,\qquad Y_{\nu-1}-Y_{\nu+1}\ =\ 2Y_\nu'.
\end{equation}
Adding and subtracting these identities gives the analogue of \eqref{eq:Jnu diff}
\begin{equation}\label{eq:Bessel rec fm, 3} 
Y_{\nu-1}\ =\  \frac{\nu}{z}Y_{\nu} + Y_{\nu}', \qquad Y_{\nu+1}\ =\ \frac{\nu}{z}Y_{\nu}-Y_{\nu}'.
\end{equation}

\noindent
For $\nu\in\R$ we have $J_\nu(\R^>), Y_\nu(\R^>)\subseteq\R$, and for such $\nu$ 
we let~$J_{\nu}$,~$Y_{\nu}$ denote also the germs (at $+\infty$)  of their restrictions to $\R^>$. We can now state the main result of the rest of this section. It gives rather detailed information about the behavior of $J_{\nu}(t)$ and $Y_{\nu}(t)$ for $\nu\in \R$ and large $t\in \R^{>}$:

\begin{theorem}\label{thm:JnuYnu} Let $\nu\in \R$. Then for the germs $J_{\nu}$ and $Y_{\nu}$ we have
$$J_\nu\  =\  \sqrt{\frac{2}{\pi x \phi_\nu'}}\ \cos\left(\phi_\nu-\frac{\pi\nu}{2}-\frac{\pi}{4}\right),\qquad
Y_\nu\  =\ \sqrt{\frac{2}{\pi x \phi_\nu'}}\ \sin\left(\phi_\nu-\frac{\pi\nu}{2}-\frac{\pi}{4}\right).$$
\end{theorem}

\noindent
The proof will take some effort, especially for $Y_\nu$ with $\nu\in \Z$.

%\begin{lemma}
%Suppose Proposition~\ref{prop:Hankel} holds for  $\nu+1$ in place of $\nu$. Then it also holds for $\nu-1$ in place of $\nu$.
%\end{lemma}
%\begin{proof}
%From  $J_\nu\preceq x^{-1/2}$ and using \eqref{eq:Bessel rec fm} we   obtain:
%\begin{multline*}
%J_{\nu-1}-\sqrt{\frac{2}{\pi x}}\cos\left(x-\frac{\pi(\nu-1)}{2}-\frac{\pi}{4}\right)\  = \\
%\frac{2\nu}{x}J_\nu-J_{\nu+1}+ \sqrt{\frac{2}{\pi x}}\cos\left(x-\frac{\pi(\nu+1)}{2}-\frac{\pi}{4}\right)\   
%\preceq\  x^{-3/2}
%\end{multline*}
%as claimed.
%\end{proof} 

\medskip\noindent
Below   $J$ denotes the analytic function $(\nu,z)\mapsto J_\nu(z)\colon\C\times(\C\setminus\R^{\leq})\to\C$.

\begin{lemma}\label{lem:Yk}
Let $k\in\Z$ and $z\in\C\setminus\R^{\leq}$. Then
$$Y_k(z)\ =\  \frac{1}{\pi}\left(  \left(\frac{\partial J}{\partial\nu}\right)(k,z)+ (-1)^k 
\left(\frac{\partial J}{\partial\nu}\right)(-k,z) 
\right).$$
In particular, $Y_{-k}=(-1)^k Y_k$ and
$Y_0(z) = \displaystyle\frac{2}{\pi}  \left(\frac{\partial J}{\partial\nu}\right)(0,z)$.
\end{lemma}
\begin{proof}
By l'H\^opital's Rule for germs of holomorphic functions at $k$, 
\begin{align*}
&\lim_{\nu\to k}  \frac{\cos(\pi\nu)J(\nu,z)-J(-\nu,z)}{\sin(\pi\nu)}\\ 
=\ &\lim_{\nu\to k} 
 \frac{-\pi\sin(\pi\nu)J(\nu,z)+\cos(\pi\nu)(\partial J/\partial\nu)(\nu,z)+(\partial J/\partial\nu)(-\nu,z)}{\pi\cos(\pi\nu)}\\
 =\ 
 &\frac{1}{\pi} \lim_{\nu\to k} \left(  (\partial J/\partial\nu)(\nu,z) + \frac{(\partial J/\partial\nu)(-\nu,z)}{\cos(\pi\nu)} \right), 
\end{align*}
and this yields the claims.
\end{proof}

\noindent
The following asymptotic relation
is crucial for establishing Theorem~\ref{thm:JnuYnu}.  It is due to Hankel~\cite{Hankel} with earlier special cases provided
by Poisson~\cite{Poisson} ($\nu=0$), Hansen~\cite{Hansen} ($\nu=1$) and Jacobi~\cite{Jacobi} ($\nu\in\Z$).

\begin{prop}[Hankel]\label{prop:Hankel}
Let $\nu\in\R$. Then for the germ $J_{\nu}$
we have: 
$$J_\nu -\displaystyle \sqrt{\frac{2}{\pi x}}\cos\left(x-\frac{\pi\nu}{2}-\frac{\pi}{4}\right)\ \preceq\  x^{-3/2}.$$
\end{prop}

\noindent 
Proposition~\ref{prop:Hankel} with Lemma~\ref{lem:BR} and Remark~\ref{rem:BR+} yield the identity for the germ $J_{\nu}$ in Theorem~\ref{thm:JnuYnu}.  As to $Y_{\nu}$, let us simplify notation by setting $S:=\sqrt{\frac{2}{\pi x \phi_\nu'}}$, $\alpha:=\frac{\pi\nu}{2}$, $\theta:= \phi_{\nu}-\frac{\pi}{4}$. Using the identity for $J_{\nu}$ in Theorem~\ref{thm:JnuYnu},  the numerator in the definition of $Y_{\nu}$ turns into 
$$S\cdot \big[\cos(2\alpha)\cos(\theta-\alpha)-\cos(\theta+\alpha)\big]$$ and the
denominator into $\sin(2\alpha)$. Trigonometric addition formulas yield
$$ \cos(2\alpha)\cos(\theta-\alpha)-\cos(\theta+\alpha)\ =\ \sin(2\alpha)\sin(\theta-\alpha).$$
For $\nu\notin \Z$, we have $\sin(2\alpha)\ne 0$, so this  gives the identity for the germ $Y_{\nu}$ in Theorem~\ref{thm:JnuYnu}. 
The identity for $Y_{\nu}$ with $\nu\in \Z$ will be dealt with after the proof of Proposition~\ref{prop:Hankel}. First a useful reduction step:

\begin{lemma}\label{lem:nu <-> nu+1}
Let $\nu\in\R$, and let $J_{\nu}$,~$J_{\nu+1}$ denote also the germs $($at $+\infty)$ of their restrictions to $\R^>$.
Then the following are equivalent:
\begin{enumerate}
\item[\textup{(i)}]
 $J_{\nu} - \displaystyle \sqrt{\frac{2}{\pi x}}\cos\left(x-\frac{\pi\nu}{2}-\frac{\pi}{4}\right)\ \preceq\ x^{-3/2}$;
\item[\textup{(ii)}]
 $J_{\nu+1} - \displaystyle \sqrt{\frac{2}{\pi x}}\cos\left(x-\frac{\pi(\nu+1)}{2}-\frac{\pi}{4}\right)\ \preceq\ x^{-3/2}$.
\end{enumerate}
\end{lemma}
\begin{proof}
Put $\alpha_\nu:=\frac{\pi\nu}{2}+\frac{\pi}{4}$ and $g_\nu:=\sqrt{\frac{2}{\pi x\phi_{\nu}'}}\in\Dx(\Q)$.
The proof of Lemma~\ref{cor:BR} gives
$\frac{1}{\sqrt{\phi_\nu'}}-1\preceq x^{-2}$, so
$g_\nu-\sqrt{\frac{2}{\pi x}}\preceq x^{-5/2}$, hence $g_\nu\preceq x^{-1/2}$ and  thus $g_{\nu}'\preceq x^{-3/2}$.
Assume now (i). By~\eqref{eq:Jnu diff} we have
\begin{equation}\label{eq:nu, 1}
J_{\nu+1}\ =\  \frac{\nu}{x} J_{\nu} - J_{\nu}'
\end{equation}
and by (i)
\begin{equation}\label{eq:nu, 2}
\frac{1}{x} J_{\nu}\  =\   \sqrt{\frac{2}{\pi}}\ x^{-3/2}\cos(x-\alpha_{\nu}) + O(x^{-5/2})\  =\  O(x^{-3/2}).
\end{equation}
We have $J_{\nu}\in V_{\nu}$, so by
Lemma~\ref{lem:BR} and (ii), 
$$J_{\nu}\ =\  g_{\nu}\cos(\phi_{\nu}-\alpha_{\nu}),$$
and  $\alpha_{\nu+1}=\alpha_{\nu} + \frac{\pi}{2}$ gives $\sin(t-\alpha_{\nu})=\cos(t-\alpha_{\nu+1})$. Thus
\begin{align*}
-J_{\nu}'\  &=\  -g_{\nu}'\cos(\phi_{\nu}-\alpha_{\nu}) + g_{\nu}\phi_{\nu}'\sin(\phi_{\nu}-\alpha_{\nu}) \\
			&=\  g_{\nu}\phi_{\nu}'\cos(\phi_{\nu}-\alpha_{\nu+1})+ O(x^{-3/2}).
\end{align*}
Also $\cos(x+u)-\cos x\preceq u$ for $u\in\c$ and $\phi_{\nu}-x\preceq x^{-1}$, 
so $\cos(\phi_{\nu}-\alpha_{\nu+1})=\cos(x-\alpha_{\nu+1})+O(x^{-1})$.  
Using $\phi'_{\nu}-1\preceq x^{-2}$ and $g_{\nu}-\sqrt{\frac{2}{\pi x}}\preceq x^{-5/2}$
this yields
$$g_{\nu}\phi_{\nu}'\cos(\phi_{\nu}-\alpha_{\nu+1})\  =\  \sqrt{\frac{2}{\pi x}}\cos(x-\alpha_{\nu+1}) + O(x^{-3/2})$$
and so
\begin{equation}\label{eq:nu, 3}
-J_{\nu}'\ =\  \sqrt{\frac{2}{\pi x}}\cos(x-\alpha_\nu) + O(x^{-3/2}).
\end{equation}
Combining \eqref{eq:nu, 1}, \eqref{eq:nu, 2}, \eqref{eq:nu, 3} yields (ii). Likewise one proves
(ii)~$\Rightarrow$~(i), 
using  
$$J_{\nu}\ =\ \frac{\nu+1}{x}J_{\nu+1}  -   J_{\nu+1}'$$
instead of \eqref{eq:nu, 1}.
 \end{proof}

\begin{remark}
Using the identities \eqref{eq:Bessel rec fm, 3} 
instead of \eqref{eq:Jnu diff} shows that Lemma~\ref{lem:nu <-> nu+1} also holds with~$\sin$,~$Y_\nu$,~$Y_{\nu+1}$ in place of $\cos$, $J_\nu$, $J_{\nu+1}$.
\end{remark}

\noindent
Lemma~\ref{lem:nu <-> nu+1} gives
a reduction of Proposition~\ref{prop:Hankel}  to the case~${\nu>-1/2}$. 
(We could also reduce to the case ${\nu>1}$, say, but the choice of $-1/2$ is useful later.)

\begin{lemma}[Poisson representation]\label{lem:Poisson}
Let $\Re\nu>-\frac{1}{2}$ and $z\in\C\setminus\R^{\leq}$. Then
$$J_\nu(z)\ =\  \frac{(\frac{z}{2})^\nu}{\Gamma(\nu+\frac{1}{2})\sqrt{\pi}}\int_{-1}^1 \ex^{tz\imag}(1-t^2)^{\nu-\frac{1}{2}}\,dt.$$
\end{lemma}
\begin{proof} For $p,q\in \C$ with $\Re p,\Re q>0$  and  $B(p,q):=\int_0^1 t^{p-1}(1-t)^{q-1}\,dt$
we have
$$B(p,q)\ =\  \frac{\Gamma(p)\Gamma(q)}{\Gamma(p+q)},$$
see for example \cite[Chapter~2, \S{}1.6]{Olver}. From the definition of $B(p,q)$ we obtain
$$\int_{-1}^1 t^{2n}(1-t^2)^{\nu-\frac{1}{2}}\,dt\ =\  \int_0^1 s^{n-\frac{1}{2}}(1-s)^{\nu-\frac{1}{2}}\,ds\ =\  
B\big(n+\textstyle\frac{1}{2},\nu+\frac{1}{2}\big).$$
In the equalities below we use this for the fourth equality (and one can appeal to a Dominated Convergence Theorem for the second): 
\begin{align*}
\int_{-1}^1 \ex^{tz\imag}(1-t^2)^{\nu-\frac{1}{2}}\,dt\  &=\ \int_{-1}^1\sum_m (1-t^2)^{\nu-\frac{1}{2}}t^m\frac{(\imag z)^m}{m!}\,dt\\
&=\ \sum_m \left(\int_{-1}^1 t^m(1-t^2)^{\nu-\frac{1}{2}}\,dt \right)\frac{(\imag z)^m}{m!} \\
&=\  \sum_n \left(\int_{-1}^1 t^{2n}(1-t^2)^{\nu-\frac{1}{2}}\, dt \right) \frac{(-1)^n z^{2n}}{(2n)!} \\
&=\ \sum_n {B\big(n+\textstyle\frac{1}{2},\nu+\frac{1}{2}\big)}\frac{(-1)^n z^{2n}}{(2n)!}\\
&=\  {\Gamma\big(\nu+\textstyle\frac{1}{2}\big)}\sum_n \frac{(-1)^n2^{2n}\,\Gamma\big(n+\frac{1}{2}\big)}{\Gamma(n+\nu+1)(2n)!}\left(\frac{z}{2}\right)^{2n}\\
&=\   {\Gamma\big(\nu+\textstyle\frac{1}{2}\big)}\sqrt{\pi}\sum_n \frac{(-1)^n}{n!\,\Gamma(n+\nu+1)}\left(\frac{z}{2}\right)^{2n},
\end{align*}
where for the last equality we used $\Gamma\big(n+\frac{1}{2}\big)=\sqrt{\pi}\ (2n)!/(n!2^{2n})$, a consequence of the 
Gauss-Legendre duplication formula for the Gamma function (see~\cite[XV, \S{}2, $\Gamma 8$]{LangCA}). 
\end{proof}

\noindent
We also need the following estimate:
 
\begin{lemma}\label{lem:exp decay}
Let $\lambda,t\in\R$ with $\lambda>-1$ and $t\geq 1$. Then  
$$t^{-(\lambda+1)}\Gamma(\lambda+1)\ =\  \int_0^\infty \ex^{-st}s^\lambda\,ds\  \leq\   
\int_0^1 \ex^{-st}s^\lambda\,ds + \Gamma(\lambda +1)\ex^{1-t},$$
and thus $\displaystyle\int_1^\infty \ex^{-st}s^\lambda\,ds\ \le\ \Gamma(\lambda +1)\ex^{1-t}$. 
\end{lemma}
\begin{proof}
We have
$$\int_0^\infty \ex^{-st}s^\lambda\,ds\ =\  \int_0^\infty \ex^{-u}\left(\frac{u}{t}\right)^\lambda \frac{du}{t}\ =\  
t^{-(\lambda+1)}\Gamma(\lambda+1)$$
and
$$\int_1^\infty \ex^{-st}s^\lambda\,ds\ =\  \ex^{-t}\int_1^\infty \ex^{-t(s-1)}s^\lambda\,ds\ \leq\  \ex^{-t}\int_1^\infty \ex^{-(s-1)}s^\lambda\,ds.$$
Now use
$\displaystyle\int_1^\infty \ex^{-(s-1)}s^\lambda\,ds =\ex\int_1^\infty \ex^{-s}s^\lambda\,ds \leq\ex\int_0^\infty \ex^{-s}s^\lambda\,ds=\ex\Gamma(\lambda+1)$.
\end{proof}

\begin{cor}\label{cor:exp decay}
Let $\lambda\in\C$ with $\Re \lambda >-1$ and $t\in\R^{\ge 1}$. Then
$$\left| \int_0^1 \ex^{-st} s^\lambda\,ds\ -\ t^{-(\lambda+1)}\Gamma(\lambda+1)\right|\  \leq\ 
\Gamma\big(\Re( \lambda)+1\big)\ex^{1-t}.$$
\end{cor}
\begin{proof}
The identities in the beginning of the proof of Lemma~\ref{lem:exp decay} generalize to 
$$\int_0^\infty \ex^{-st} s^\lambda\,ds\ =\  t^{-(\lambda+1)}\Gamma(\lambda+1).$$
Hence
$$\left| \int_0^1 \ex^{-st} s^\lambda\,ds\ -\  t^{-(\lambda+1)}\Gamma(\lambda+1)\right|\  =\ 
\left| \int_1^\infty \ex^{-st} s^\lambda\,ds \right|\  \leq\  
\int_1^\infty \ex^{-st} s^{\Re \lambda}\, ds,$$
and now use Lemma~\ref{lem:exp decay}.
\end{proof}

\noindent
By Lemma~\ref{lem:nu <-> nu+1} the next result is more than enough to give Proposition~\ref{prop:Hankel}. The 
proof is classical and uses Laplace's method, cf.~\cite[Chapter~3, \S{}7]{Olver}.

\begin{lemma}\label{lem:Hankel}
Suppose $\Re\nu>-\frac{1}{2}$, and let $t$ range over $\R^{\ge 1}$. Then
$$J_\nu(t)\  =\  \sqrt{\frac{2}{\pi t}}\ \cos\left(t-\frac{\pi\nu}{2}-\frac{\pi}{4}\right) +O(t^{-\frac{3}{2}}) \quad\text{as $t\to+\infty$.}$$ 
\end{lemma}
\begin{proof} We consider the holomorphic function
$$z\mapsto f_{\nu,t}(z):=\ex^{tz\imag}(1-z^2)^{\nu-\frac{1}{2}} \colon \C\setminus (\R^{\leq -1}\cup\R^{\geq 1})\to\C,$$
and set $I_\nu(t):=\displaystyle\int_{-1}^1 f_{\nu,t}(s)\,ds$. By Lemma~\ref{lem:Poisson} we have $J_{\nu}(t)= \displaystyle\frac{(\frac{t}{2})^\nu}{\Gamma(\nu+\frac{1}{2})\sqrt{\pi}}I_{\nu}(t)$. To determine the asymptotic behavior of $I_{\nu}(t)$
 as $t\to +\infty$ we integrate along
     the contour $\gamma_{R}$ depicted below, where $R$ is a real number $>1$.  

\usetikzlibrary{decorations.markings}

\begin{center}
\begin{tikzpicture}
	\node     (0) at (-3, 5) {};
	\node[above left] at (0) {$-1+R\imag$};
		\node   (1) at (-3, 0) {};
		\node   (2) at (3, 5) {};
			\node[above right] at (2) {$1+R\imag$};

		\node  (3) at (3, 0) {};
		\node   (4) at (-3, 1) {};
			\node[left] at (4) {$-1+\frac{\imag}{R}\ $};
		\node   (5) at (3, 1) {};
		\node[right] at (5) {$\ 1+\frac{\imag}{R}$};

		\node   (6) at (-2, 0) {};
		
				\node[below] at (6) {$\qquad-1+\frac{1}{R}$};

		\node   (7) at (2, 0) {};
		
						\node[below] at (7) {$1-\frac{1}{R}\qquad$};
\node [below] at (-3,0) {$-1$};
\draw[fill] (-3,0) circle [radius=.1em];

\node [below] at (3,0) {$1$};
\draw[fill] (3,0) circle [radius=.1em];

\node [right] at (2,2) {$\gamma_R$};
%		\node  (8) at (4, 1) {};
%		\node  (9) at (-4.75, 1) {};
\draw[->,>=stealth] (-5,0) to (5,0);
\begin{scope}[very thick,decoration={
    markings,
    mark=at position 0.5 with {\arrow{>}}}]
 	 		\draw[postaction={decorate}] (0.center) to (4.center);
		\draw[postaction={decorate}] (5.center) to (2.center);
		\draw[postaction={decorate}] [bend left=45] (4.center) to (6.center);
		\draw[postaction={decorate}] [bend right=-45] (7.center) to (5.center);
		\draw[postaction={decorate}] (6.center) to (7.center);
		\draw[postaction={decorate}] (2.center) to (0.center);
\end{scope}
 \end{tikzpicture}
 \medskip
\end{center}

\noindent 
By Cauchy, $\int_{\gamma_{R}} f_{\nu,t}(z)\,dz=0$ (see \cite[III, \S{}5]{LangCA}), and letting~$R\to+\infty$ we obtain~$I_\nu(t) = I_\nu^-(t) - I_\nu^+(t)$
where 
$$I_\nu^-(t)\ :=\ \imag\int_0^\infty f_{\nu,t}(-1+\imag s) \,ds,\qquad I_\nu^+(t)\ :=\ \imag\int_0^\infty f_{\nu,t}(1+\imag s) \,ds.$$
Now
$$I_\nu^+(t)\ =\ \imag
\int_0^\infty \ex^{t(1+s\imag)\imag}\big(1-(1+s\imag)^2\big)^{\nu-\frac{1}{2}}\,ds\  =\ \imag\ex^{t\imag}\int_0^\infty \ex^{-st} (s^2-2 s \imag)^{\nu-\frac{1}{2}}\,ds.
$$
The complex analytic function
$(\kappa,z)\mapsto \big(1+\frac{\imag}{2}z\big)^{\kappa}-1$ on $\C\times \{z\in \C:\, |z|<2\}$ vanishes on the locus
$z=0$, so we have a complex analytic function $(\kappa,z)\mapsto r_{\kappa}(z)$ on the same region such that
$\big(1+\frac{\imag}{2}z\big)^{\kappa}-1=zr_{\kappa}(z)$ for all $(\kappa,z)$ in this region.
For~$\kappa=\nu-\frac{1}{2}$ this yields a continuous function
 $r\colon [0,1]\to \C$ such that $\big(1+\frac{\imag}{2}s\big)^{\nu-\frac{1}{2}}=1+r(s)s$ for all $s\in [0,1]$,
and $r(0)=(\nu-\frac{1}{2})\frac{\imag}{2}$. In view of an identity stated just before Lemma~\ref{lem:Bessel conv} this yields
for $s\in (0,1]$:
$$(s^2-2s\imag)^{\nu-\frac{1}{2}}\ =\ (-2s\imag)^{\nu-\frac{1}{2}} \big(1+\textstyle\frac{\imag}{2}s\big)^{\nu-\frac{1}{2}}\ =\ (-2s\imag)^{\nu-\frac{1}{2}}+r(s) s (-2s\imag)^{\nu-\frac{1}{2}},$$
so
$$\int_0^1 \ex^{-st} (s^2-2 s \imag)^{\nu-\frac{1}{2}}\,ds\ =\ \int_0^1 \ex^{-st} (-2s\imag)^{\nu-\frac{1}{2}}\,ds+\int_0^1 \ex^{-st} r(s)
s(-2s\imag)^{\nu-\frac{1}{2}}\,ds.$$
By Corollary~\ref{cor:exp decay}, as $t\to +\infty$, 
$$%\label{eq:Hankel, 1}
\int_0^1 \ex^{-st} (-2s\imag)^{\nu-\frac{1}{2}}\,ds\  =\  
  (-2\imag)^{\nu-\frac{1}{2}}
t^{-(\nu+\frac{1}{2})}\Gamma\big(\nu+\textstyle\frac{1}{2}\big)+O(\ex^{-t}).$$
Take $C\in\R^{>}$ such that $\abs{r(s)}\leq C$ for all $s\in [0,1]$, and set  $\lambda:=\Re(\nu)-\frac{1}{2}$, so~$\lambda>-1$.
 Using the identity of Lemma~\ref{lem:exp decay} for the last step,
\begin{align*}
\left|\int_0^1 \ex^{-st} r(s) s(-2s\imag)^{\nu-\frac{1}{2}}\,ds\right|\ &=\ 
\left|\int_0^1 (-2\imag)^{\nu-\frac{1}{2}}r(s)\ex^{-st} s^{\nu+\frac{1}{2}}\,ds\right| \\
&\leq\  2^{\lambda}\ex^{\Im(\nu)\pi/2} C\int_0^1 \ex^{-st}  s^{\lambda+1}\,ds\\   
&\leq\  2^{\lambda}\ex^{\Im(\nu)\pi/2}C\int_0^\infty\ex^{-st}s^{\lambda+1}\,ds\\  &=\ 2^{\lambda}\ex^{\Im(\nu)\pi/2}C\,\Gamma(\lambda+2)t^{-\lambda-2}. 
\end{align*}
Hence, as $t\to +\infty$, 
\begin{equation}\label{eq:I+, 1}
\int_0^1 \ex^{-st} (s^2-2 s \imag)^{\nu-\frac{1}{2}}\,ds\ =\   (-2\imag)^{\nu-\frac{1}{2}}
t^{-\nu-\frac{1}{2}}\Gamma\big(\nu+\textstyle\frac{1}{2}\big)+O(t^{-\lambda -2}).
\end{equation}
Next, take $D\in\R^{>}$ such that $\abs{(1-2\frac{\imag}{s})^{\nu-\frac{1}{2}}}\leq D$ for all $s\in [1,\infty)$.
For such $s$ we have
$\abs{(s^2-2s\imag)^{\nu-\frac{1}{2}}}\leq Ds^{2\lambda}\leq Ds^{2\Re \nu}$ and thus
$$\left| \int_1^\infty \ex^{-st}(s^2-2s\imag)^{\nu-\frac{1}{2}}\,ds\right|\ \leq\  D\int_1^\infty \ex^{-st}s^{2\Re \nu}\,ds,$$
hence by Lemma~\ref{lem:exp decay}:
\begin{equation}\label{eq:I+, 2}
\int_1^\infty  \ex^{-st}(s^2-2s\imag)^{\nu-\frac{1}{2}}\,ds\  =\  O(\ex^{-t}) \quad\text{as $t\to+\infty$.}
\end{equation}
Combining \eqref{eq:I+, 1} and \eqref{eq:I+, 2} yields
$$I_\nu^+(t)\  =\  \imag (-2\imag)^{\nu-\frac{1}{2}} \ex^{t\imag} t^{-\nu-\frac{1}{2}}\Gamma\big(\nu+\textstyle\frac{1}{2}\big) +O(t^{-\lambda-2})\quad\text{as $t\to+\infty$.}$$
In the same way we obtain
$$I_\nu^-(t)\  =\  \imag (2\imag)^{\nu-\frac{1}{2}} \ex^{-t\imag} t^{-\nu-\frac{1}{2}}\Gamma\big(\nu+\textstyle\frac{1}{2}\big) +O(t^{-\lambda-2})\quad\text{as $t\to+\infty$.}$$
 Thus as $t\to +\infty$:
\begin{align*}
I_\nu(t)\	&=\  I_\nu^-(t) - I_\nu^+(t) \\
		&=\   2^{\nu-\frac{1}{2}}\imag\big[    \imag^{\nu-\frac{1}{2}} \ex^{-t\imag} - (-\imag)^{\nu-\frac{1}{2}} \ex^{t\imag}\big]  t^{-\nu-\frac{1}{2}}\Gamma\big(\nu+\textstyle\frac{1}{2}\big) +O(t^{-\lambda-2}).
\end{align*}
Using $\imag^{\nu-\frac{1}{2}}=\ex^{\frac{1}{2}(\nu-\frac{1}{2})\pi\imag}$  and the like we have
\begin{align*} &\imag \big[ \imag^{\nu-\frac{1}{2}} \ex^{-t\imag} - (-\imag)^{\nu-\frac{1}{2}} \ex^{t\imag}\big]\  =\ 
 2\cos\left(t-\frac{\pi\nu}{2}-\frac{\pi}{4}\right), \text{ and thus}\\
J_\nu(t)\  &=\   \frac{(\frac{t}{2})^\nu}{\Gamma(\nu+\frac{1}{2})\sqrt{\pi}}\,I_\nu(t)\ =\ 
\sqrt{\frac{2}{\pi t}}\ \cos\left(t-\frac{\pi\nu}{2}-\frac{\pi}{4}\right)+O(t^{-\frac{3}{2}})
\end{align*} 
as $t\to+\infty$.  \end{proof} 

\noindent
Here is a consequence of Hankel's result, cf. Lommel \cite[p.~67]{Lommel}:

\begin{cor} Let $\nu\in \R$. Then
$J_\nu^2+J_{\nu+1}^2 =~\displaystyle\frac{2}{\pi x}+O(x^{-2})$.
\end{cor}
\begin{proof}
With $\alpha_\nu:=\frac{\pi\nu}{2}+\frac{\pi}{4}$  we have by  Proposition~\ref{prop:Hankel},
$$J_{\nu}\ =\ \sqrt{\frac{2}{\pi x}}\, \cos(x-\alpha_{\nu}) +O(x^{-3/2}),\qquad J_{\nu+1}\ =\ \sqrt{\frac{2}{\pi x}}\, \cos(x-\alpha_{\nu+1}) +O(x^{-3/2}).$$
Now use $\sin(x-\alpha_\nu)=\cos(x-\alpha_{\nu+1})$.
\end{proof}

\noindent
Since Proposition~\ref{prop:Hankel} is now established, so is Theorem~\ref{thm:JnuYnu}, except for the $Y_{\nu}$-identity when $\nu\in \Z$. To treat that case, and also for use in the next subsection, we now prove  a uniform version of Lemma~\ref{lem:Hankel}:

\begin{lemma}\label{lem:Hankel, uniform}
Let $\nu_0\in\C$ and $\Re\nu_0>-\frac{1}{2}$.  Then there are reals $\varepsilon>0$, $t_0\ge 1$, and  $C_0>0$, such that for all $\nu\in \C$ with $\abs{\nu-\nu_0}< \varepsilon$~and all $t\geq t_0$:
$$\Re\nu\ >\ -\frac{1}{2}, \qquad \left| J_\nu(t) -  \sqrt{\frac{2}{\pi t}}\cos\left(t-\frac{\pi\nu}{2}-\frac{\pi}{4}\right)\right|\ \leq\  C_0 t^{-\frac{3}{2}}.$$
\end{lemma} 
\begin{proof}
We follow the proof of Lemma~\ref{lem:Hankel}, where in the beginning we introduced the complex analytic function~$(\nu,z)\mapsto r_{\nu}(z)$ on $\C\times \big\{z\in \C:\, |z|<2\big\}$. Take~$\varepsilon\in \R^{>}$  and $C\in \R^{\ge 1}$ such that $0<\varepsilon<\Re(\nu_0)+\frac{1}{2}$ and  
$\abs{r_{\nu-\frac{1}{2}}(s)}\leq C$ for all $(\nu,s)\in B_0\times [-1,1]$ where $B_0:=\{\nu\in \C:\,\abs{\nu-\nu_0}<\varepsilon\}$. (To handle $I_\nu^+$ we use this for~$s\in [0,1]$, and  
to deal with $I_{\nu}^{-}$ we use $s\in [-1,0]$.) 
Take also $D\in\R^>$ such that
 $\abs{(1-2\frac{\imag}{s})^{\nu-\frac{1}{2}}}\leq D$ for all $\nu\in B_0$ and $s\geq 1$ (and thus also for $\nu\in B_0$ and~$s\le -1$). 
Next, set $\lambda_0:= \Re(\nu_0)-\frac{1}{2}$, and take~$t_0\geq 1$ such that~$\ex^{t-1}\geq t^{\lambda_0+\varepsilon+2}$ for all $t\geq t_0$. 
Below $\nu$ ranges over~$B_0$ and $t$ over~$\R^{\geq t_0}$, and   $\lambda:=\Re(\nu)-\frac{1}{2}$, so~$\lambda>-1$. 
Then, as in the proof of  Lemma~\ref{lem:Hankel}: 
 $$\left|\int_0^1 \ex^{-st} r_{\nu-\frac{1}{2}}(s)s(-2s\imag)^{\nu-\frac{1}{2}}\,ds\right|\ 
 \le\  2^{\lambda}\ex^{\Im(\nu)\pi/2}C\,\Gamma(\lambda+2)t^{-\lambda-2}.$$
Take  
$C_\Gamma\in\R^{>}$ such that  $2^{\lambda}\ex^{\Im(\nu)\pi/2}C\,\Gamma(\lambda+2)\leq C_\Gamma$ 
for all $\nu$. Then 
$$\left|\int_0^1 \ex^{-st} r_{\nu-\frac{1}{2}}(s)s(-2s\imag)^{\nu-\frac{1}{2}}\,ds\right|\ \leq\  C_\Gamma\,t^{-\lambda-2}.$$
By increasing $C_\Gamma$ we arrange that $C_\Gamma\geq C$ and that for all $\nu$,
$$  2^{\lambda}\ex^{\Im(\nu)\pi/2}C\,\Gamma(\lambda+1),\  D\,\Gamma(2\lambda+2)\  \le\  
 C_{\Gamma}.$$
 We have   $\ex^{1-t}\leq t^{-\lambda-2}$,  so by Corollary~\ref{cor:exp decay}:
\[\left|\int_0^1 \ex^{-st} (-2s\imag)^{\nu-\frac{1}{2}}\,ds -
  (-2\imag)^{\nu-\frac{1}{2}}
t^{-(\nu+\frac{1}{2})}\Gamma(\nu+\textstyle\frac{1}{2})\right|\  \leq\  C_\Gamma\,t^{-\lambda-2}.
\]
Combining this with an earlier displayed inequality yields: 
\begin{equation}\label{eq:I+, 1, unif} \left|\int_0^1 \ex^{-st} (s^2-2 s \imag)^{\nu-\frac{1}{2}}\,ds -   (-2\imag)^{\nu-\frac{1}{2}}
t^{-\nu-\frac{1}{2}}\Gamma(\nu+\textstyle\frac{1}{2})\right|\ \leq\  2C_\Gamma\, t^{-\lambda-2}.
\end{equation}
 As in the proof of  Lemma~\ref{lem:Hankel} and using Lemma~\ref{lem:exp decay} we  also have
\begin{align*}
\left|\int_1^\infty \ex^{-st}(s^2-2s\imag)^{\nu-\frac{1}{2}}\,ds\right|\  &\leq\  D\int_1^\infty\ex^{-st}s^{2\Re(\nu)}\,ds\\
&\le\ D\Gamma(2\lambda+2)\ex^{1-t}\ \le\   C_\Gamma\, t^{-\lambda-2}.
\end{align*}
Combining this with \eqref{eq:I+, 1, unif} we obtain:
$$\big|I_{\nu}^+(t)- \imag (-2\imag)^{\nu-\frac{1}{2}} \ex^{t\imag} t^{-\nu-\frac{1}{2}}\Gamma\big(\nu+\textstyle\frac{1}{2}\big)\big|\ \leq\  3C_\Gamma\, t^{-\lambda-2}.$$
In the same way, 
$$\big|I_{\nu}^-(t)- \imag (2\imag)^{\nu-\frac{1}{2}} \ex^{-t\imag} t^{-\nu-\frac{1}{2}}\Gamma\big(\nu+\textstyle\frac{1}{2}\big)\big|\  \leq\  3C_\Gamma\, t^{-\lambda-2}$$
and so as in the proof of Lemma~\ref{lem:Hankel}:  
$$\left|I_\nu(t)-2^{\nu+\frac{1}{2}}t^{-\nu-\frac{1}{2}}\Gamma\big(\nu+\textstyle\frac{1}{2}\big)\cos\left(t-\frac{\pi\nu}{2}-\frac{\pi}{4}\right)\right|\  \leq\  6C_\Gamma\, t^{-\lambda-2}.
$$
Hence
$$\left|J_\nu(t)- \sqrt{\frac{2}{\pi t}}\ \cos\left(t-\frac{\pi\nu}{2}-\frac{\pi}{4}\right)\right|\ \leq\
\frac{6C_\Gamma}{\sqrt{\pi}\,2^{\Re(\nu)}\abs{\Gamma(\nu+\frac{1}{2})}}t^{-\frac{3}{2}}.
$$
Thus $\varepsilon$, $t_0$ as chosen, and a suitable $C_0$ have the required properties.
\end{proof}

\noindent
To finish the proof of Theorem~\ref{thm:JnuYnu}, it suffices by Lemma~\ref{lem:BR} and Remark~\ref{rem:BR+} to show the following:

\begin{lemma} Let $k\in \Z$. Then for the germ $Y_{k}$ we have:
$$Y_k -\displaystyle \sqrt{\frac{2}{\pi x}}\ \sin\left(x-\frac{\pi k}{2}-\frac{\pi}{4}\right)\ \preceq\  x^{-\frac{3}{2}}.$$
\end{lemma}
\begin{proof} By the remark after the proof of Lemma~\ref{lem:nu <-> nu+1} it is enough to treat the case~$k=0$.
Lemma~\ref{lem:Hankel, uniform} with~$\nu_0=0$ yields reals $t_0\geq 1$, $C_0>0$, and $\varepsilon$ with~${0 < \varepsilon < \frac{1}{2}}$ 
such that for all $\nu\in\C$ with~$\abs{\nu}<\varepsilon$ and all $t\geq t_0$:
$$\left| J_\nu(t) -  \sqrt{\frac{2}{\pi t}}\cos\left(t-\frac{\pi\nu}{2}-\frac{\pi}{4}\right)\right|\ \leq\  C_0 t^{-\frac{3}{2}}.$$
Let $t\geq t_0$ be fixed and consider the entire function $d$ given by
\begin{align*} d(\nu)\  &:=\   J(\nu,t) -  \sqrt{\frac{2}{\pi t}}\cos\left(t-\frac{\pi\nu}{2}-\frac{\pi}{4}\right), \text{ so}\\
d'(\nu)\ &=\ \left(\frac{\partial J}{\partial\nu}\right)(\nu,t) - \sqrt{\frac{\pi}{2t}}\sin\left(t-\frac{\pi\nu}{2}-\frac{\pi}{4}\right)
\end{align*}
and hence by Lemma~\ref{lem:Yk}:
$$d'(0)\ =\ \left(\frac{\partial J}{\partial\nu}\right)(0,t) - \sqrt{\frac{\pi}{2t}}\sin\left(t-\frac{\pi}{4}\right)=
\frac{\pi}{2}\left[Y_0(t)- \sqrt{\frac{2}{\pi t}}\sin\left(t-\frac{\pi}{4}\right)\right].
$$
Also $|d'(0)|\le \frac{1}{\varepsilon}\max_{\abs{\nu}=\varepsilon}\abs{d(\nu)}$, a Cauchy inequality, and thus
$$
\left|Y_0(t)- \sqrt{\frac{2}{\pi t}}\sin\left(t-\frac{\pi}{4}\right)\right|\ =\ 
\frac{2}{\pi}\left|d'(0)\right|\ \leq\ \frac{2C_0}{\pi\varepsilon}t^{-\frac{3}{2}},$$
which gives the desired result for $k=0$ and thus for all $k\in \Z$.  
\end{proof}

\noindent
In the rest of this subsection we let $\nu$ range over $\R$ and derive some consequences of Theorem~\ref{thm:JnuYnu}.
Toward showing that the germ $\psi_{\nu}\in \Ex(\Q)$ depends analytically on $\nu$ we introduce the real analytic function $\Psi\colon \R\times \R^{>}\to \R^{>}$ by
$$\Psi(\nu,t)\ :=\ \frac{\pi t}{2}\left[J(\nu,t)^2 + Y(\nu,t)^2\right],$$
and let $\Psi(\nu,-)$ be the function $t\mapsto \Psi(\nu,t)\colon \R^{>}\to \R^{>}$. Then by Theorem~\ref{thm:JnuYnu}:

\begin{cor}\label{jnu2ynu2}
$\psi_{\nu}$ is the germ at $+\infty$ of $\Psi(\nu,-)$. 
\end{cor} 

\noindent
For $\phi_{\nu}$ we consider the real analytic function
%\marginpar{why ``real analytic''}
$$\tilde{\Phi}\ :\ \R\times \R^{>} \to \R, \quad (\nu,t)\mapsto \int_1^t\frac{1}{\Psi(\nu,s)}\,ds.$$
Let $\tilde{\Phi}_{\nu}$ be the germ (at $+\infty$) of $t\mapsto \tilde{\Phi}(\nu,t)$.
Then ${\tilde{\Phi}_{\nu}}'=\frac{1}{\psi_{\nu}}=\phi_{\nu}'$, so
$\phi_{\nu}=\tilde{\Phi}_\nu+ c_{\nu}$ where $c_{\nu}$ is a real constant. To determine this constant we note that by
Proposition~\ref{prop:3rd order sol}  we have
$1-\frac{1}{\psi_{\nu}}\preceq x^{-2}$, which gives the real number
$$\tilde{c}_{\nu}\ :=\ \int_1^\infty\left(1-\frac{1}{\Psi(\nu,s)}\right)\,ds.$$
We also set $\tilde{c}(\nu,t)=\int_1^t\left(1-\frac{1}{\Psi(\nu,s)}\right)\,ds$, so $\tilde{c}(\nu,t)\to \tilde{c}_{\nu}$ as $t\to +\infty$, and
for $t>0$ we have $\tilde{\Phi}(\nu,t)+ \tilde{c}(\nu,t)= t-1$. Taking germs we obtain
$\tilde{\Phi}_{\nu}+\tilde{c}_{\nu}+1-x\prec 1$. Also~$\phi_{\nu}-x\prec 1$, so $\tilde{\Phi}_{\nu}+c_{\nu}-x\prec 1$, and thus
$c_{\nu}=\tilde{c}_{\nu} +1$. This suggests the function
$$\Phi: \R\times \R^{>}\to \R, \quad (\nu,t)\mapsto \tilde{\Phi}(\nu,t) + \tilde{c}_{\nu}+1.$$
The above arguments yield:

\begin{cor} For each $\nu$ the germ of $\Phi(\nu,-)$ is $\phi_{\nu}$.
\end{cor} 

\noindent
Thus as for $\psi_{\nu}$ the germ $\phi_{\nu}$ has a unique real analytic representative on $\R^{>}$, namely~$\Phi(\nu,-)$.  Note that $\Phi$ is real analytic iff the function~$\nu \mapsto \tilde{c}_{\nu}\colon \R \to \R$ is real analytic,  but we don't even know if this last function is continuous. 

\medskip
\noindent
For the next result, cf.~\cite[\S{}13.74]{Watson}:

\begin{cor}\label{cor:Bessel 1}
$J_\nu^2+Y_\nu^2$ is eventually strictly decreasing,    $x(J_\nu^2+Y_\nu^2)$ is eventually strictly increasing 
if $\abs{\nu}<1/2$ and  eventually strictly decreasing if $\abs{\nu}>1/2$.
%% and $(x^2-\nu^2)^{1/2}(J_\nu^2+Y_\nu^2)$ is eventually strictly increasing. (hmm, this may depend on \nu)
\end{cor}
\begin{proof} We have $\psi_{\nu}\sim 1+\frac{\mu-1}{8}x^{-2}$ by Proposition~\ref{prop:3rd order sol}, and 
$\psi_{\nu}$ is hardian, thus ${\psi_{\nu}' \preceq x^{-3}}$, and so $(x^{-1}\psi_{\nu})'  = -x^{-2}\psi_{\nu}+x^{-1}\psi_{\nu}'\sim -x^{-2}$. 
This yields the claims.
\end{proof}

\begin{cor}[{Schafheitlin~\cite[p.~86]{Schafheitlin}}]\label{cor:Bessel 2}
If $\nu > 1/2$, then, as elements of~$\Dx(\Q)$, 
$$\frac{2/\pi}{x}\  <\  J_\nu^2+Y_\nu^2\  <\  \frac{2/\pi}{(x^2-\nu^2)^{1/2}}.$$
\end{cor}

\noindent
Proposition~\ref{prop:3rd order sol} also yields (cf.~\cite[\S{}13.75]{Watson} or \cite[Chapter~9, \S{}9]{Olver}):

\begin{cor} \label{cor:Nicholson}
The germ $J_\nu^2+Y_\nu^2$ has the asymptotic expansion 
$$J_\nu^2+Y_\nu^2\  \sim\   \frac{2}{\pi x} \sum_n  (2n-1)!! \frac{ (\nu,n) }{2^n}  x^{-2n}.$$ 
\end{cor}

\begin{remark}
 Nicholson~\cite{Nicholson11,Nicholson10} (see~\cite[\S\S{}13.73--13.75]{Watson}) established 
 Corollary~\ref{cor:Nicholson},
 but ``the analysis is difficult'' \cite[p.~340]{Olver}. 
 A simpler deduction of the integral representation of $J_\nu^2+Y_\nu^2$ used by Nicholson was given by Wilkins~\cite{Wilkins}, %also employing the third-order
 %linear differential equation from Proposition~\ref{prop:3rd order sol}. 
 see also~\cite[Chap\-ter~9, \S{}7.2]{Olver}.
 %%Our proof of Corollary~\ref{cor:Nicholson} avoids the use of complex analysis. {\bf maybe replace by another characterization?} 
 For more on the history of this result, see \cite[\S{}1]{Kalf}.
\end{remark}

\noindent 
Theorem~\ref{thm:JnuYnu} and the recurrence relations \eqref{eq:Bessel rec fm} and \eqref{eq:Bessel rec fm, 2nd}  yield remarkable identities among the germs $\phi_\nu, \phi_{\nu-1}, \phi_{\nu+1}\in \Dx(\Q)$. For example:

\begin{cor} Recalling that $\psi_\nu=1/\phi_\nu'$, we have
\begin{align*}
-\sqrt{\psi_{\nu-1}}\sin(\phi_{\nu-1}-\phi_\nu) + \sqrt{\psi_{\nu+1}}\sin(\phi_{\nu+1}-\phi_\nu)\ &=\
\frac{2\nu}{x}\sqrt{\psi_\nu}, \\
\sqrt{\psi_{\nu-1}}\cos(\phi_{\nu-1}-\phi_\nu) - \sqrt{\psi_{\nu+1}}\cos(\phi_{\nu+1}-\phi_\nu)\  &=\  0.
\end{align*}
\end{cor}
\begin{proof}
Put $H_\nu:=J_\nu+Y_\nu\imag\in\Calinf[\imag]$. Then 
$$ H_\nu\ =\ \sqrt{\psi_{\nu}}\cdot\sqrt{\frac{2}{\pi x}}\ \ex^{(\phi_\nu- \frac{\pi\nu}{2}-\frac{\pi}{4})\imag}, \qquad 
H_{\nu-1}+H_{\nu+1}\  =\  \frac{2\nu}{x}H_\nu,$$
and dividing both sides of the equality on the right by $\displaystyle\sqrt{\frac{2}{\pi x}}\ \ex^{(\phi_\nu-\frac{\pi\nu}{2}-\frac{\pi}{4})\imag}$ gives
\begin{align*}
&\sqrt{\psi_{\nu-1}}\ex^{(\phi_{\nu-1}-\phi_\nu+\pi/2)\imag} + \sqrt{\psi_{\nu+1}}\ex^{(\phi_{\nu+1}-\phi_\nu-\pi/2)\imag}\\
 &=\ \sqrt{\psi_{\nu-1}}\ \imag\ex^{(\phi_{\nu-1}-\phi_\nu)\imag} - \sqrt{\psi_{\nu+1}}\ \imag\ex^{(\phi_{\nu+1}-\phi_\nu)\imag}\ =\
\frac{2\nu}{x}\sqrt{\psi_\nu}.
\end{align*}
Now take real and imaginary parts in the last identity.
\end{proof}

\subsection*{An asymptotic expansion for the zeros of Bessel functions}
We are going to use Corollary~\ref{cor:phiinv} to strengthen a result of McMahon on
parametrizing the zeros of Bessel functions: Corollary~\ref{cor:McMahon} and the remark following it.
 Lemma~\ref{lem:Fourier} below, due to Fourier~\cite{Fourier} for $\nu=0$ and to Lommel~\cite[p.~69]{Lommel} in general, is only included for completeness; its proof  is based on 
the following useful identity also due to Lommel~\cite{Lommel79}:

\begin{lemma}\label{lem:Lommel}
Let $\alpha,\beta\in\C^\times$, $\mu,\nu\in\C$, and let $y_\mu, y_\nu\colon\C\setminus\R^{\leq}\to\C$ be holomorphic solutions of~\textup{(B$_\mu$)} and \eqref{eq:Bessel}, respectively. Then on 
$\C\setminus\big(\alpha^{-1}\R^{\le}\cup \beta^{-1}\R^{\le}\big)$:
%\begin{multline}\label{eq:albe}
%\frac{d}{dz}\left[z\left(\beta y_\mu(\alpha z)y_\nu'(\beta z) - \alpha y_\nu(\beta z)y_\mu'(\alpha z)  \right)\right]\  = \\
%\left( (\alpha^2-\beta^2)z - \frac{\mu^2-\nu^2}{z}\right)y_\mu(\alpha z)y_\nu(\beta z).
%\end{multline}
\begin{align*}\label{eq:albe}
&\frac{d}{dz}\left[z\left(\beta y_\mu(\alpha z)y_\nu'(\beta z) - \alpha y_\nu(\beta z)y_\mu'(\alpha z)  \right)\right]\ =\\ 
&\left( (\alpha^2-\beta^2)z - \frac{\mu^2-\nu^2}{z}\right)y_\mu(\alpha z)y_\nu(\beta z).
\end{align*}
\end{lemma}
\begin{proof} 
Let $U\subseteq \C$ be open, let $g,\tilde g\colon U\to\C$ be continuous, and let $y,\tilde y\colon U\to\C$ be holomorphic such that $4y''+gy=4\tilde y''+\tilde g\tilde y=0$.
An easy computation gives  
$$\operatorname{wr}(y,\tilde y)'\ =\ \textstyle\frac{1}{4}(g-\tilde g)y\tilde y.$$
Assume for now that $\alpha,\beta\in \R^{>}$ and apply this remark to $U:=\C\setminus\R^{\leq}$ and  
\begin{align*} y(z)\ &:=\ (\alpha z)^{1/2}y_\mu(\alpha z),\qquad g(z)\ :=4\alpha^2+(1-4\mu^2)z^{-2}\\
 \tilde y(z)\ &:=\ (\beta z)^{1/2}y_\nu(\beta z),\qquad \tilde g(z)\ :=\ 4\beta^2+(1-4\nu^2)z^{-2}.
 \end{align*} 
 Then $4y''+gy=0$, since~$z\mapsto z^{1/2}y_\mu(z)\colon \C\setminus\R^{\leq}\to\C$ satisfies ($\operatorname{L}_{\mu}$).
Likewise, $4\tilde y''+\tilde g\tilde y=0$. This yields the claimed identity for $\alpha,\beta\in \R^{>}$ by a straightforward computation using that
$(\alpha z)^{\frac{1}{2}}=\alpha^{\frac{1}{2}}z^{\frac{1}{2}}$ for such $\alpha$ and for $z\in \C\setminus\R^{\le}$. 

Next, for general $\alpha$, $\beta$, $z$ we note that $U_3:=\big\{(\alpha,\beta, z)\in (\C^\times)^3:\, \alpha z, \beta z\notin \R^{\le}\big\}$ is open in $\C^3$.  Moreover, $U_3$ is connected. (Proof sketch:  suppose $(\alpha, \beta, z)\in  U_3$; then so is
 $(\ex^{\theta\imag}\alpha, \ex^{\theta\imag}\beta,\ex^{-\theta\imag} z)$ for $\theta\in \R$, so we can ``connect to'' $\alpha\in \R^{>}$; next keep~$\alpha$,~$z$ fixed, and rotate $\beta$ to a point in $\R^{>}$ while preserving $\beta\notin z^{-1}\R^{\le}$.)  Both sides in the claimed identity define a complex analytic function on $U_3$. Now use analytic continuation as in \cite[(9.4.4)]{Dieudonne}. 
\end{proof}

\noindent
We now define certain improper complex integrals and state some basic facts about them. Let $U$ be an open subset of $\C$ with $0\in \partial U$ 
and let $f\colon U \to \C$ be holomorphic such that for some $\varepsilon\in \R^{>}$ we have
 $f(u)=O\big(|u|^{-1+\varepsilon}\big)$ as $u\to 0$. For $z\in U$ such that
$(0,z]:=\big\{tz:\ t\in (0,1]\big\}\subseteq U$, we set 
$$\int_0^z f(u)\,du\ :=\  \lim_{\delta\downarrow 0}\int_\delta^1zf(tz)\,dt\quad   \text{(the limit exists in $\C$).}$$
Suppose $(0,z]\subseteq U$ for all $z\in U$. Then the function $z\mapsto \int_0^z f(u)\, du$ on $U$ is holomorphic with derivative $f$, and $\lim_{z\to 0} \int_0^z f(u)\,du = 0$; to see this, first show that for any $z_0\in U$, open ball $B\subseteq U$ centered at $z_0$, and $z\in U$, we have $\int_0^z f(u)\,du=\int_0^{z_0} f(u)\,du + \int_{z_0}^z f(u)\,du$, where the last integral is by definition~$\int_0^1(z-z_0)f\big(z_0+t(z-z_0)\big)\,dt$. Thus the integral below makes sense: 
%consequence of 

\begin{cor}\label{lem:Lommel+}
Let $\alpha,\beta,z\in \C^\times$ satisfy $\alpha z,\beta z\notin\R^{\leq}$, and let $\nu\in \R^{\geq-1}$. Then
\begin{equation}\label{eq:alphabeta}
(\alpha^2-\beta^2)\int_0^z uJ_\nu(\alpha u)J_\nu(\beta u)\,du\  =\ z\big(\beta J_\nu(\alpha z)J_\nu'(\beta z)-\alpha J_\nu(\beta z)J_\nu'(\alpha z)\big).
\end{equation}
\end{cor}
\begin{proof}  Fixing $\alpha, \beta\in \C^\times$ and $\nu\in \R^{\geq-1}$, both sides in \eqref{eq:alphabeta} are holomorphic functions of $z$ on the open
subset
$\C\setminus \big(\alpha^{-1}\R^{\le}\cup \beta^{-1}\R^{\le}\big)$ of $\C$, with equal derivatives by  Lemma~\ref{lem:Lommel}. Moreover, both sides  tend to $0$ as $z\to 0$ in $\C\setminus \big(\alpha^{-1}\R^{\le}\cup \beta^{-1}\R^{\le}\big)$, using for the right hand side the first two terms in the power series for $J_\nu$ and $J_\nu'$. 
%This yields \eqref{eq:alphabeta}. 
 \end{proof}

\begin{lemma}\label{lem:Fourier}
Let $\nu\in\R^{\ge -1}$. Then all zeros of $J_\nu$ are contained in~$\R^>$.
\end{lemma}
\begin{proof}
Let $\alpha\in\C\setminus\R^{\leq}$ be a zero of $J_\nu$.
From the power series for~$J_\nu$ we see that 
then
$\bar{\alpha}$ is also a zero of $J_\nu$, and $\alpha\notin\imag\R$.
Putting $\beta=\bar{\alpha}$ and~$z=1$ in~\eqref{eq:alphabeta} yields
$$(\alpha^2-\bar{\alpha}^2)\int_0^1 tJ_\nu(\alpha t)J_\nu(\bar{\alpha} t)\,dt\  =\  \bar{\alpha} J_\nu(\alpha)J_\nu'(\bar{\alpha})-\alpha J_\nu(\bar{\alpha})J_\nu'(\alpha)\ =\ 0.$$
If $\alpha\notin\R$, then this yields
$\displaystyle\int_0^1 tJ_\nu(\alpha t)J_\nu(\bar{\alpha} t)\,dt=0$, but $J_\nu(\alpha t)J_\nu(\bar{\alpha} t)\in \R^{\ge}$ 
for all~$t\in (0,1]$ and $J_\nu(\alpha t)\ne 0$ for some $t\in (0,1]$, a contradiction. 
Thus $\alpha\in\R$.
\end{proof}

\noindent
Taking $\alpha=\beta=1$ in Lemma~\ref{lem:Lommel} yields (for all $\mu,\nu\in \C$): 
\begin{equation}\label{eq:Lommel}
\frac{d}{dz}\left[ z\big(J_\mu'(z)J_\nu(z)-J_\mu(z)J_\nu'(z)\big) \right]\  =\  
(\mu^2-\nu^2) \frac{J_\mu(z)J_\nu(z)}{z}\ \text{ on }\C\setminus\R^{\leq}.
\end{equation}

\noindent
In the next result it is convenient to let $J$ denote the analytic function $(\nu,z)\mapsto J_{\nu}(z)$ on $\C\times (\C\setminus \R^{\le})$, so for $(\nu,z)\in \C\times(\C\setminus \R^{\le})$ we have $J_{\nu}'(z)=\frac{\partial J}{\partial z}(\nu, z)$. In its proof we also use that for any complex analytic functions $A$, $B$ on an open set $U\subseteq \C^2$ and $(\mu,z)$ and $(\nu,z)$ ranging over $U$:
$$\lim_{\mu\to \nu}\frac{\partial}{\partial z}\left(A(\mu,z)\frac{B(\mu,z)-B(\nu,z)}{\mu-\nu}\right)\ =\ \frac{\partial}{\partial z}\left(A(\nu,z)\frac{\partial B}{\partial \nu}(\nu,z)\right),$$
an easy consequence of $\frac{\partial^2B}{\partial \nu\partial z}=\frac{\partial^2B}{\partial z \partial \nu}$.

\begin{cor}\label{cor:Lommel} On $\C\times (\C\setminus \R^{\le})$ we have
%Let $\nu\in\C$, $z\in\C\setminus\R^{\leq}$. Then
$$\frac{d}{dz}\left[ z\left(J_\nu(z)\cdot \frac{\partial ^2J}{\partial\nu\partial z}(\nu,z)- J'_{\nu}(z)\cdot\frac{\partial J}{\partial\nu}(\nu,z) \right) \right]\  =\  
2\nu \frac{J_\nu^2(z)}{z},$$
\end{cor}
\begin{proof}
For $\mu,\nu\in \C$ with $\mu\ne \nu$  we obtain from \eqref{eq:Lommel} that on $\C\setminus \R^{\le}$, 
$$\frac{d}{dz}\left[ z\left(  J_\mu(z) \frac{J_\mu'(z)-J_\nu'(z)}{\mu-\nu} -J_\mu'(z) \frac{J_\mu(z)-J_\nu(z)}{\mu-\nu}\right) \right] \ 
=\ (\mu+\nu) \frac{J_\mu(z)J_{\nu}(z)}{z}.$$
Now let $\mu$ tend to $\nu$ and use the identity preceding the corollary. 
%\cite[V, \S{}1, Theorem~1.2]{LangCA}.
\end{proof}

\noindent
Below $\nu\in\R$, so the set $Z_\nu:=\R^>\cap J_\nu^{-1}(0)$ of positive real zeros of $J_\nu$ is infinite and has no limit point. 
Let $(j_{\nu,n})$ be the enumeration of $Z_{\nu}$. Note that if $t\in Z_{\nu}$, then $J_{\nu+1}(t)=-J_{\nu}'(t)$ by~\eqref{eq:Bessel rec fm}, so $J_{\nu+1}(t)\ne 0$. 

\begin{prop}[{Schl\"afli~\cite{Schlaefli}}]\label{prop:jnun diff} Let $n$  be given. Then the function $$\nu\mapsto j(\nu):=j_{\nu,n}\colon \R^{>-1}\to\R^>$$ is analytic,  and its derivative at $\nu>0$ is given by
$$j'(\nu)\ =\  \frac{2\nu}{j(\nu)\,J_{\nu+1}^2\big(j(\nu)\big)}\int_0^{j(\nu)} J_\nu^2(s) \frac{ds}{s}.$$
In particular, the restriction of $j$ to $\R^>$ is strictly increasing.
\end{prop}

\noindent
In the proof of Proposition~\ref{prop:jnun diff} we use:

\begin{lemma}\label{lem:no small zeros}
Let $\varepsilon\in\R^>$. Then there exists $\delta\in\R^>$ such that $J_\nu(t)\neq 0$ for all~$\nu\geq -1+\varepsilon$ and $t\in (0,\delta]$.
\end{lemma}
\begin{proof}
Take $\delta\in\R^>$ such that $\delta^2<4\log(1+\varepsilon)$. Then for $\nu\geq -1+\varepsilon$ and $0<t\leq\delta$:
$$\left| \frac{\Gamma(\nu+1)J_\nu(t)}{(\frac{1}{2}t)^\nu}-1\right|\  =\  
\left|\sum_{n\geq 1} \frac{(-1)^n(\frac{1}{4}t^2)^n}{n!(\nu+n)\cdots (\nu+1)} \right|\  \leq\  \frac{\exp(\frac{1}{4}\delta^2)-1}{\varepsilon}\ <1\ ,$$
hence $J_\nu(t)\neq 0$.
\end{proof}

\begin{proof}[Proof of Proposition~\ref{prop:jnun diff}]
Let $\nu_0\in\R^{>-1}$. 
For each $m$ we have~$J_{\nu_0}(j_{\nu_0,m})=0$ and $J_{\nu_0}'(j_{\nu_0,m})\neq 0$. 
So IFT (the Implicit Function Theorem~\cite[(10.2.2), (10.2.4)]{Dieudonne}) yields  an interval 
$I=(\nu_0-\varepsilon, \nu_0+\varepsilon)$ with $\varepsilon\in \R^{>}$,  
$-1<\nu_0-\varepsilon$, and for each~$m\leq n$ an analytic function~$j_m\colon I\to\R^>$ with~$J_\nu\big(j_m(\nu)\big)=0$  for $\nu\in I$ and $j_m(\nu_0)=j_{\nu_0,m}$.  Shrinking~$\varepsilon$ if necessary we also arrange to have  $\delta\in\R^>$ such that~$j_{\nu_0,m}>\delta$ and~$|j_m(\nu)-j_{\nu_0,m}|<\delta$ for all $\nu\in I$ and $m\le n$, and
$J_\nu'(t)\neq 0$ for all~$\nu\in I$, $m\le n$ and $t\in\R$ with~${\abs{t-j_{\nu_0,m}}}<\delta$. 
Hence for $\nu\in I$ and $m\le n$ (using IFT at all $\nu\in I$): 
$$Z_\nu\cap {(j_{\nu_0,m}-\delta, j_{\nu_0,m}+\delta)}\ =\ \big\{j_m(\nu)\big\},\qquad j_0(\nu) <j_1(\nu) < \cdots < j_n(\nu).$$
 We claim that for $\nu\in I$ we have  
 $$\big\{t\in Z_{\nu}:\, t\le j_n(\nu)\big\}=\big\{j_0(\nu), j_1(\nu), \dots , j_n(\nu)\big\}.$$ 
Suppose for example that $\nu_1\in I$, $t_1\in Z_{\nu_1}$, $t_1< j_0(\nu_1)$; we shall derive a contradiction. (The assumption 
$\nu_1 \in I$, $t_1\in Z_{\nu}$, $j_m(\nu_1) < t_1 < j_{m+1}(\nu_1)$, $m<n$, leads to a contradiction in the same way.) Using IFT again we see that 
$$U:=\big\{\nu\in I:\, \text{ there is a $t\in Z_\nu$ with $t<j_0(\nu)$}\big\}$$ is open in $I$, and by an easy extra argument
using Lemma~\ref{lem:no small zeros} also closed in~$I$. As~$\nu_1\in U$, this gives $U=I$, so $\nu_0\in U$, a contradiction.   
The claim gives~${j_m(\nu)=j_{\nu,m}}$ for~$\nu\in I$ and  $m\leq n$. Taking $m=n$ it follows that $j$ is analytic.

Next, let $\nu$ range over $\R^{>}$.  Differentiating~$J\big(\nu, j(\nu) \big)=0$ yields
$$\frac{\partial J}{\partial\nu}\big(\nu, j(\nu)\big)+J'_\nu\big(j(\nu)\big)j'(\nu)\ =\ 0.$$
Using the primitive of $s\mapsto \frac{J_{\nu}^2(s)}{s}$ provided by Corollary~\ref{cor:Lommel} gives 
$$\int_0^{j(\nu)} \frac{J_\nu^2(s)}{s}\, ds\  =\  -\frac{j(\nu)}{2\nu}J'_\nu\big(j(\nu)\big) \frac{\partial J}{\partial \nu}\big(\nu, j(\nu)\big).$$
Now combine the two displayed identities with $J_\nu'\big(j(\nu)\big)=-J_{\nu+1}\big(j(\nu)\big)$.
\end{proof}

\noindent
We now bring in  Lemma~\ref{lem:Hankel, uniform} to   bound $j_{\nu,n}$ for $\nu>0$ and sufficiently large $n$:

 \begin{prop}\label{prop:jn0n}
Let $\nu_0\in\R^>$. 
Then there is an $n_0\in\N$ such that for~$n\geq n_0$:
$$(n+\textstyle\frac{1}{2}\nu_0+\frac{1}{2})\pi\  \leq\  j_{\nu_0,n}\  \leq\  (n+\textstyle\frac{1}{2}\nu_0+1)\pi.$$
\end{prop}
\begin{proof}
Compactness and Lemma~\ref{lem:Hankel, uniform} yield $C_0,t_0\in\R^>$ such that for all $\nu$ in the smallest closed interval $I$ containing both $\nu_0$ and $1/2$, and for all~$t\geq t_0$:
$$\left| J_\nu(t) -\displaystyle \sqrt{\frac{2}{\pi t}}\cos\left(t-\frac{\pi\nu}{2}-\frac{\pi}{4}\right)\right|\  \leq\  C_0 t^{-3/2}.$$
We arrange that
$t_0\geq C_0\sqrt{\pi}$. Hence if  $\nu\in I$ and~$j_{\nu,n}\geq t_0$, then
$$\left| \cos\left(j_{\nu,n}-\frac{\pi}{2}\nu-\frac{\pi}{4}\right)\right|\ \leq\   \frac{1}{\sqrt{2}},$$ 
and so we have a unique $k_{\nu,n}\in\Z$ with
$$\textstyle \frac{1}{4}\pi\  \leq\  j_{\nu,n}-(\frac{1}{2}\nu+\frac{1}{4}+k_{\nu,n})\pi\  \leq\  \frac{3}{4}\pi.$$
Let $\nu_1$ be the left endpoint of $I$ (so $\nu_1=1/2$ or $\nu_1=\nu_0$) and  take $n_0\in\N$ such that~$j_{\nu_1,n_0}\geq t_0$. 
Then $j_{\nu,n}\geq   t_0$ for $\nu\in I$ and $n\geq n_0$ by Proposition~\ref{prop:jnun diff}.
Let $n\geq n_0$; we claim that $\nu\mapsto k_{\nu,n}\colon I\to\Z$ is constant.
To see this, note that Proposition~\ref{prop:jnun diff} yields $\delta\in(0,1/4)$ such that
for all $\nu,\tilde\nu\in I$ with $\abs{\nu-\tilde\nu}<\delta$ we have $\abs{j_{\nu,n}-j_{\tilde\nu,n}}<\pi/4$, which in view of
$$\textstyle -\frac{1}{2}\pi\ \leq\  (j_{\nu,n}-j_{\tilde\nu,n})-(\nu-\tilde\nu)\frac{\pi}{2}-(k_{\nu,n}-k_{\tilde\nu,n})\pi\  \leq\ \frac{1}{2}\pi,$$
gives  $k_{\nu,n}=k_{\tilde\nu,n}$. Thus $\nu\mapsto  k_{\nu,n}\colon I\to\Z$ is locally constant, and hence constant.
Let $k_n$ be the common value of $k_{\nu,n}$ for $\nu\in I$.
Now $Z_{1/2}=\{m\pi:m\geq 1\}$ by Lemma~\ref{lem:J1/2}, hence $j_{1/2,n}=(n+1)\pi$ and so 
$$\textstyle \frac{1}{4}\pi\  \leq\ (n+1)\pi-(\frac{1}{2}+k_{n})\pi\  \leq\  \frac{3}{4}\pi.$$
This yields $k_n=n$.
\end{proof}

\begin{cor}\label{cor:McMahon}
Suppose $\nu>0$. There is a strictly increasing~$\zeta\in\c_{n_0}$~\textup{(}$n_0\in\N$\textup{)} whose germ is in~$\Ex(\Q)$
such that $j_{\nu,n}=\zeta(n)$ for all $n\geq n_0$ and which has  for~$s:=\big(x+\textstyle\frac{1}{2}\nu+\frac{3}{4}\big)\pi$
the asymptotic expansion
$$\zeta\  \sim\  s-\left(\frac{\mu-1}{8}\right)s^{-1}-\left(\frac{(\mu-1)(7\mu-31)}{192}\right) \frac{s^{-3}}{2!}+\cdots\ .$$
\end{cor}
\begin{proof}
Take $n_0$ as in Proposition~\ref{prop:jn0n} with $\nu_0=\nu$.
Theorem~\ref{thm:JnuYnu} yields $t_0\in\R^>$ and a representative of  $\phi=\phi_\nu$  in $\c^1_{t_0}$,  also denoted by~$\phi$,
such that for all $t\ge t_0$, 
$$\phi'(t)>0,\qquad J_\nu(t)=\sqrt{\frac{2}{\pi t\phi'(t)}}\ \cos\left(\phi(t)-\frac{\pi}{2}\nu-\frac{\pi}{4}\right).$$
Increasing $n_0$ if necessary we arrange that $j_{\nu,n_0}\geq t_0+\frac{\pi}{2}$.
Then for $n\ge n_0$, $$\phi(j_{\nu,n})-\left(\textstyle\frac{1}{2}\nu+\frac{3}{4}\right)\pi\in \Z\pi.$$ Take $k\in\Z$ with $\phi(j_{\nu,n_0})=\left(\frac{1}{2}\nu+\frac{3}{4}+k\right)\pi$. Then for $n\ge n_0$ we have
$\phi(j_{\nu,n})=\left(n-n_0+\frac{1}{2}\nu+\frac{3}{4}+k\right)\pi$.
By Proposition~\ref{prop:jn0n} we have for $n\ge n_0$,
$$(n+\textstyle\frac{1}{2}\nu+\frac{1}{2})\pi\ \leq\  j_{\nu,n}\ \leq\  (n+\textstyle\frac{1}{2}\nu+1)\pi,$$
and thus for all $n\ge n_0$,
$$\phi\big( (n+\textstyle\frac{1}{2}\nu+\frac{1}{2})\pi\big)\ \leq\ \phi(j_{\nu,n})\ =\ \left(n+\textstyle\frac{1}{2}\nu+\frac{3}{4}+k-n_0\right)\pi\ \leq\ \phi\big( (n+\textstyle\frac{1}{2}\nu+1)\pi\big).$$
Since $\phi-x\preceq x^{-1}$, this yields
$k=n_0$, therefore $\phi(j_{\nu,n})=\left(n+\frac{1}{2}\nu+\frac{3}{4}\right)\pi$ for~${n\geq n_0}$.
Let~$\phi^{\operatorname{inv}}\in\c_{t_1}$ be the compositional inverse of $\phi$, where $t_1:=\phi(t_0)$, and let~$\zeta\in\c_{n_0}$ be given by $\zeta(t):=\phi^{\operatorname{inv}}\left(\big(t+\frac{1}{2}\nu+\frac{3}{4}\big)\pi\right)$
for~$t\geq n_0$. Then $\zeta$ is strictly increasing with $j_{\nu,n}=\zeta(n)$ for~$n\geq n_0$. Taking $\zeta$ and 
$\phi^{\operatorname{inv}}$ as  germs we have $\zeta=\phi^{\operatorname{inv}}\circ s$. Now
$\phi^{\operatorname{inv}}\in \Ex(\Q)$ by Lemma~\ref{phiinvexq}, and $\Ex(\Q)\circ \Ex(\Q)^{>\R}\subseteq \Ex(\Q)$, so $\zeta\in \Ex(\Q)$. 
The claimed asymptotic expansion for $\zeta$ follows from Corollary~\ref{cor:phiinv}. 
\end{proof}

\begin{remark}
The asymptotic expansion for $j_{\nu,n}$ as $n\to\infty$ in Corollary~\ref{cor:McMahon} was obtained by McMahon~\cite{McMahon}.
(For   $\nu=1$, apparently Gauss was   aware of it as early as 1797, cf.~\cite[p.~506]{Watson}.)
What is new here is that we specified a function $\zeta$ with germ in $\Ex(\Q)$ such that $j_{\nu,n}=\zeta(n)$ for all sufficiently large $n$. 
\end{remark}

\noindent
In \cite[p.~247]{Olver}, Olver states: ``No explicit formula is available for the general term'' of the asymptotic expansion
for $j_{\nu,n}$ as $n\to\infty$ in Corollary~\ref{cor:McMahon}. The remark after the proof of Corollary~\ref{cor:phiinv} yields the asymptotic expansion
$$\zeta\ \sim\  s- \sum_{j=1}^\infty \left(\sum_{i=1}^{j}\frac{(2(j-1))!}{(2j-1-i)!}B_{ij}(u_1,\dots,u_{j-i+1})\right)\frac{s^{-2j+1}}{j!},$$
which is perhaps as explicit as possible. 
The values of $u_1$, $u_2$, $u_3$ given before Corollary~\ref{cor:phiinv} yield the first few terms of this
expansion: 
\begin{multline*}
\zeta\ \sim\  s\ -\  \frac{\mu-1}{8}s^{-1}\ -\  \frac{(\mu-1)(7\mu-31)}{192} \frac{s^{-3}}{2!}\ - \\
\frac{(\mu-1)(83\mu^2-982\mu+3779)}{2560} \frac{s^{-5}}{3!}\ -\ \cdots.
\end{multline*}

\subsection*{Appendix: inversion of formal power series}
In this appendix we discuss multiplicative and compositional inversion of power series.
We use~[ADH, 12.5] and its notations. Thus $x,y_1,y_2,\dots,z$ are distinct indeterminates, and
$$R\ :=\ \Q[x,y_1,y_2,\dots],\qquad A\ :=\ R[[z]].$$
We also let $K$ be a field of characteristic zero.
Recall from [ADH, 12.5.1]  the definition of the Bell polynomials $B_{ij}\in \Q[y_1,\dots,y_d]$, where~$i\leq j$ and $d=j-i+1$:
$$B_{ij} :=   \sum_{\substack{\k=(k_1,\dots,k_{d})\in\N^{d} \\ \abs{\k} = i,\ \dabs{\k}=j}} \frac{j!}{k_1!k_2!\cdots k_{d}!}\, \left(\frac{y_1}{1!}\right)^{k_1} \left(\frac{y_2}{2!}\right)^{k_2}\cdots  \left(\frac{y_d}{d!}\right)^{k_d}.$$
(Also $B_{ij}:=0\in \Q[y_1, y_2,\dots]$ for $i>j$.) Let $$y\ :=\ \sum_{n\geq 1} y_n \frac{z^n}{n!}\in zR[[z]].$$ 
By  [ADH, (12.5.2)] we  have in $R[[z]]$: 
$$\frac{y^i}{i!}\ =\  \sum_{j\geq 0} B_{ij} \frac{z^j}{j!}\ =\ \sum_{j\ge i} B_{ij}\frac{z^j}{j!}.$$

\begin{lemma}\label{lem:Bell shift}
Let $i\leq j$ and $d=j-i+1$; then
$$B_{ij}\left(\frac{y_2}{2},\frac{y_3}{3},\dots,\frac{y_{d+1}}{d+1}\right)\  =\  \frac{j!}{(i+j)!} B_{i,i+j}(0,y_2,y_3,\dots,y_{j+1}).$$
\end{lemma}
\begin{proof}
We have 
$$y-y_1z\ =\ z\sum_{n\geq 1} \left(\frac{y_{n+1}}{n+1}\right)\frac{z^n}{n!}\ =\ \sum_{n\geq 2} y_{n}\frac{z^n}{n!},$$ hence
\[\frac{(y-y_1z)^i}{i!}\  =\  \sum_{j\geq 0} B_{ij}\left(\frac{y_2}{2},\frac{y_3}{3},\dots\right) \frac{z^{i+j}}{j!}\  =\  \sum_{k\geq 0} B_{ik}(0,y_2,y_3,\dots)\frac{z^k}{k!}.\qedhere\]
\end{proof}

\noindent
For $j\in\N$ we set 
\begin{equation}\label{eq:Bell polys}
B_j\  :=\  \sum_{i=0}^j i!\,B_{ij}\in\Q[y_1,\dots,y_{j+1}].
\end{equation}
Note that
$$\frac{B_j}{j!} \ =\   \sum_{i=0}^j
\sum_{\substack{\k=(k_1,\dots,k_{d})\in\N^{d} \\ \abs{\k} = i,\ \dabs{\k}=j,\ d=j-i+1}} \frac{i!}{k_1!\cdots k_d!}\, \left(\frac{y_1}{1!}\right)^{k_1}  \left(\frac{y_2}{2!}\right)^{k_2} \cdots  \left(\frac{y_d}{d!}\right)^{k_d}.
$$
%Here
%$$ {i\choose \k} := \frac{i!}{k_1!\cdots k_d!}\quad \text{for $d,i\in\N$, $\k=(k_1,\dots,k_d)\in\N^d$ with $\abs{\k}=i$.}$$
We have $B_0=B_{00}=1$ and 
$B_j = \displaystyle \sum_{i=1}^j i!\,B_{ij}\in \Q[y_1,\dots, y_j]$ for $j\geq 1$. % (since $B_{0j}=0$ for $j\geq 1$).
Using the examples following~[ADH, 12.5.4] we obtain
\begin{align*}
B_1\ &=\ y_1, \\ B_2\ &=\ y_2+2y_1^2, \\ 
B_3\ &=\ y_3+6y_1y_2+6y_1^3, \\ B_4\ &=\ y_4+8y_1y_3+6y_2^2+36y_1^2y_2+24y_1^4, \\ 
B_5\ &=\ y_5+10y_1y_4+20y_2y_3+60y_1^2y_3+90y_1y_2^2+240y_1^3y_2+120y_1^5.
\end{align*}
We have $1-y\in 1+zR[[z]]\subseteq R[[z]]^\times$ with inverse $(1-y)^{-1}=\sum_{i\geq 0} y^i$, so
for~$m\geq 1$:
\begin{equation}\label{eq:1/(1-y)^m}
(1-y)^{-m}\  =\   \sum_{i\geq 0} {i+m-1 \choose m-1}y^i =  \sum_{j\geq 0} \left( \sum_{i=0}^j m^{\bar i}B_{ij} \right) \frac{z^j}{j!}
\end{equation}
where $m^{\bar i}:=m(m+1)\cdots (m+i-1)$ (so $m^{\bar 0}=1$, $1^{\bar i}=i!$). 
In particular,
$$(1-y)^{-1}\  =\  \displaystyle\sum_{j\geq 0} B_j \frac{z^j}{j!}.$$

\begin{cor}\label{lem:inv of 1-y}
Let  $f=\displaystyle\sum_{n\geq 1}f_n\frac{z^n}{n!}\in K[[z]]$ \textup{(}$f_n\in K$\textup{)}. Then 
\begin{align*}
(1-f)^{-1}\	&=\  \sum_{j\geq 0}  B_j(f_1,\dots,f_{j}) \frac{z^j}{j!} \\
			&=\  1+f_1z+(f_2+2f_1^2)\frac{z^2}{2!}+(f_3+6f_1f_2+6f_1^3)\frac{z^3}{3!}+\cdots.
\end{align*}
\end{cor}

\noindent
Next we discuss the compositional inversion of formal power series. From~[ADH, 12.5]  recall that~$zK^\times+z^2K[[z]]$ is a group under formal composition with $z$ as its identity element.
We denote the compositional inverse of any  $f\in zK^\times+z^2K[[z]]$  by $f^{[-1]}$. We equip the field $K(\!(z)\!)$ of Laurent series with the strongly additive and $K$-linear derivation $d/dz$ (so  $z'=1$).

\begin{definition}\label{def:res}
Let $f=\sum_{k\in\Z} f_k z^k\in K(\!(z)\!)$ where $f_k\in K$ for   $k\in\Z$. Then~$\res(f):=f_{-1}\in K$
is the {\bf residue} of $f$. (We also have the residue morphism $f\mapsto f(0)\colon K[[z]]\to K$ of the
valuation ring $K[[z]]$ of $K(\!(z)\!)$.)\index{residue}
\end{definition}

\noindent
%For $f\in K(\!(z)\!)^\times$ as in the definition above we also let $vf=\min\{k\in\Z:f_k\neq 0\}$.
The map $f\mapsto\res(f)\colon K(\!(z)\!)\to K$ is strongly additive and $K$-linear. 

\begin{lemma}\label{lem:resder}
Let $f\in K(\!(z)\!)$. Then  $\res(f')=0$, and if $f\neq 0$, then $\res(f^\dagger) = vf$.
\end{lemma}
\begin{proof}
The first claim is clearly true. For the second, 
let $f=z^kg$ where~$k=vf$, $g\in K[[z]]^\times$. Then
$f^\dagger=kz^{-1}+g^\dagger \in kz^{-1}+K[[z]]$, so $\res(f^\dagger)=k=vf$.
\end{proof}

\begin{cor}\label{cor:Jacobi, 1}
Let $f,g\in  K(\!(z)\!)$. Then $\res(f'g)=-\res(fg')$ and thus, if $g\neq 0$ and $k\in\Z$, then
$\res(f'g^k)=-k\res  (fg^{k-1}g')$.
\end{cor}
\begin{proof}
For the first claim, use the Product Rule and the first part of Lemma~\ref{lem:resder}. The second claim follows from the first.
\end{proof}

\begin{cor}[{Jacobi~\cite{Jacobi1830}}]\label{cor:Jacobi, 2}
Let $f\in zK[[z]]^{\neq}$ and $g\in K(\!(z)\!)$.
Then $$\res\!\big( (g\circ f)f'\big) = vf\,\res(g).$$
\end{cor}
\begin{proof}
By strong additivity and $K$-linearity it is enough to show this for $g=z^k$ ($k\in\Z$).
If~${k\neq -1}$, then
$$\res\!\big( (g\circ f)f'\big) = \res(f^kf')=\res\!\left( \left( f^{k+1}/(k+1)\right)'\right) = 0 = vf\,\res(g),$$
and if $k=-1$, then
$$\res\!\big( (g\circ f)f'\big) = \res(f^\dagger) = vf=vf\,\res(g),$$
using Lemma~\ref{lem:resder}.
\end{proof}

\noindent
We now obtain the Lagrange Inversion Formula, following \cite{Gessel}:

\begin{theorem} \label{thm:LIF}
Let $f\in zK^\times+z^2K[[z]]$, $g\in K[[z]]$. Then
$$g\circ f^{[-1]}\  =\  g(0) + \sum_{n\geq 1}  \frac{1}{n}\res(g'f^{-n})z^n.$$
\end{theorem}
\begin{proof}
Let $h:=g\circ f^{[-1]}=\sum_{n} h_nz^n$ ($h_n\in K$).  Then for $n\geq 1$ we have
\begin{align*} \frac{1}{n} \res (g'f^{-n})\  &=\  \res\!\big( gf^{-n-1}  f' \big)\ =\ \res\!\big((h\circ f)f^{-n-1}f'\big)\\
&=\  \res\left(\big(hz^{-n-1}\circ f\big)\cdot f' \right)\  =\  \res\!\big( hz^{-n-1}\big)\  =\  h_n,
\end{align*}
using Corollary~\ref{cor:Jacobi, 1} for the first equality, and~\ref{cor:Jacobi, 2} for the next to last.
\end{proof}

\noindent 
Taking $g=z^m$  in the above yields:

\begin{cor}\label{cor:LIF}
If $f\in zK^\times+z^2K[[z]]$, then for $m\geq 1$:
$$(f^{[-1]})^m\ =\ \displaystyle\sum_{n\geq m}  \frac{m}{n}\res(z^{m-1}f^{-n})z^n.$$
\end{cor}

\begin{remark}
%For $K=\C$ and convergent $f$, $g$, 
Theorem~\ref{thm:LIF} stems from Lagrange~\cite{Lagrange} and B\"urmann (cf.~\cite{Hindenburg}). The identity in Corollary~\ref{cor:LIF} is from
Jabotinsky~\cite[Theorem~II]{Jabotinsky} and Schur~\cite{Schur}.
%\marginpar{remark skipped}
\end{remark} 

\noindent
We now use Corollary~\ref{cor:LIF} to express the coefficients of $f^{[-1]}$ in terms of those~of~$f$; cf. \cite[\S{}3.8, Theorem~E]{Comtet-Book}. Here $f=z+\sum_{n\geq 2} f_n\frac{z^n}{n!}$ with $f_n\in K$ for $n\ge  2$. Then
$$ g\ :=\ f^{[-1]}\ =\ z+\sum_{n\geq 2} g_n\frac{z^n}{n!} \qquad (g_n\in K  \text{ for }n\ge 2). $$
For $h\in zK[[z]]$, let $\llbracket h\rrbracket$ denote the (upper triangular) iteration matrix of $h$ as in~[ADH, 12.5], so
$\llbracket h\rrbracket_{m,n}=0$ for $m>n$, $\llbracket h\rrbracket_{0,0}= 1$, $\llbracket h\rrbracket_{0,n}=0$ for $n\ge 1$.  

\begin{prop} \label{prop:inverse coeff} For $1\le m\le n$ we have 
$$ \llbracket g\rrbracket_{m,n}\ =\ \sum_{i=0}^{n-m} \frac{(-1)^i(i+n-1)!}{(i+n-m)!(m-1)!} B_{i,i+n-m}(0,f_2,\dots,f_{n-m+1}).$$
\end{prop}
 \begin{proof} By [ADH, remarks before 12.5.5] we have
$g^m / m! = \sum_n  \llbracket g\rrbracket_{m,n} z^n/n!$ and so by Corollary~\ref{cor:LIF},
$\llbracket g\rrbracket_{m,n} = (n-1)!/(m-1)! \res(z^{m-1}f^{-n})$ if $1\le m\le n$.
Set
$$h\ :=\ 1-\frac{f}{z}\ =\ \sum_{n\geq 1} h_n\frac{z^n}{n!} \qquad\text{where $h_n=-f_{n+1}/(n+1)$ for $n\geq 1$.}$$
 Then $\res(z^{m-1}f^{-n})$ is the constant term of $z^mf^{-n}$ and so equals the coefficient of~$z^{n-m}$ in 
 $z^{n-m}z^mf^{-n}= (z/f)^n=\big(\frac{1}{1-h}\big)^n$. Hence by \eqref{eq:1/(1-y)^m} for $n\ge m\ge 1$: 
$$\res(z^{m-1}f^{-n})\  =\  \sum_{i=0}^{n-m} \frac{n^{\bar{i}}}{(n-m)!} B_{i,n-m}(h_1,\dots,h_{n-m-i+1}).$$
Lemma~\ref{lem:Bell shift} gives for $n\ge m$ and $i\le n-m$: 
\begin{align*} B_{i,n-m}\big(h_1,\dots, h_{n-m-i+1}\big)\ &=\ B_{i,n-m}\big(\frac{-f_2}{2},\dots, \frac{-f_{n-m-i+2}}{n-m-i+2}\big)\\
&=\ \frac{(n-m)!}{(i+n-m)!}B_{i,i+n-m}\big(0, -f_2,\dots, -f_{n-m+1}\big)\\
&=\  \frac{(-1)^i(n-m)!}{(i+n-m)!}B_{i,i+n-m}\big(0, f_2,\dots, f_{n-m+1}\big).
\end{align*} 
For $1\le m\le n$ this yields the identity claimed for $\llbracket g\rrbracket_{m,n}$. 
 \end{proof}

\begin{cor}[{Ostrowski~\cite{Ostrowski}}] For $n\ge 2$ we have
 $$g_n = \sum_{\substack{\k=(k_1,\dots,k_{n})\in\N^{n} \\ \abs{\k} = n-1,\ \dabs{\k}=2n-2}} (-1)^{n-k_1-1} \frac{(2n-k_1-2)!}{k_2!k_3!\cdots k_n!}\left(\frac{f_2}{2!}\right)^{k_2} \left(\frac{f_3}{3!}\right)^{k_3}\cdots \left(\frac{f_n}{n!}\right)^{k_n}.$$
\end{cor}
\begin{proof} Let $n\ge 2$. We have $g_n=\llbracket g\rrbracket_{1,n}$, so by
Proposition~\ref{prop:inverse coeff}: 
$$g_n =  \sum_{i=1}^{n-1} (-1)^i B_{i,i+n-1}(0,f_2,f_3,\dots,f_{n})\qquad (n\geq 2).$$
Now use the definition of the Bell polynomials and   reindex.
\end{proof}

\noindent
One now easily determines the first few $g_n$:
\begin{align*}
g_2	&=  -f_2 \\
g_3	&= -f_3+3f_2^2 \\
g_4	&= -f_4 + 10 f_3f_2 - 15f_2^3 \\
g_5 &= -f_5 + 15 f_4f_2 + 10 f_3^2 - 105 f_3 f_2^2 + 105 f_2^4
\end{align*}

\noindent
We now establish analogues of some of these results   for compositional inversion in the field $K(\!(x^{-1})\!)$ of Laurent
series in $x^{-1}$ over $K$. We have
 the usual valuation~$v\colon K(\!(x^{-1})\!)^\times\to\Z$ on  $K(\!(x^{-1})\!)$
with valuation ring~$K[[x^{-1}]]$; so $v(x)=-1$. 
We also have the unique continuous derivation $f\mapsto f'$ on $K(\!(x^{-1})\!)$  such that $a'=0$ for $a\in K$ and~$x'=1$.
This valuation and derivation make $K(\!(x^{-1})\!)$ a $\d$-valued field with small derivation.

\medskip
\noindent
As in $K(\!(z)\!)$, we have a well-behaved notion of composition in~$K(\!(x^{-1})\!)$:
for~$f,g$ in $K(\!(x^{-1})\!)$ with~$f\succ 1$ and  $g=\sum_k g_kx^k$ ($g_k\in K$ for~$k\in\Z$),  
the family~$(g_k f^k)$ is summable in~$K(\!(x^{-1})\!)$, and we denote its sum  by~$g\circ f$. 
For $f\in K(\!(x^{-1})\!)$ with~${f \succ 1}$, the  map $g\mapsto g\circ f$ is a strongly additive and $K$-linear
field embedding, which is bijective if $f \asymp x$.
This can be seen, for example, by relating composition in~$K(\!(x^{-1})\!)$
to composition in  $K(\!(z)\!)$: 
The strongly additive, $K$-linear map
$$\tau \colon K(\!(x^{-1})\!) \to K(\!(z)\!)\quad\text{with $\tau(x^{-k})=z^k$ for all $k\in \Z$},$$
is an isomorphism of valued fields. Let $f\in K(\!(x^{-1})\!)$, $f\succ 1$. Then 
$\tau(1/f)\in zK[[z]]$, and we have a commutative diagram
$$\xymatrix@C+5em{ K(\!(x^{-1})\!) \ar[r]^{g\,\mapsto\, g\circ f} \ar[d]^\tau & K(\!(x^{-1})\!) \ar[d]^\tau \\
K(\!(z)\!) \ar[r]^{h\,\mapsto\, h\circ \tau(1/f)} & K(\!(z)\!)}$$
of strongly additive, $K$-linear maps. Also $\tau(1/f)\in zK^\times +z^2K[[z]]$ if $f\asymp x$. 

\medskip
\noindent
As in Definition~\ref{def:res}, we define:

\begin{definition}
Let $f=\sum_{k\in\Z} f_k x^k\in K(\!(x^{-1})\!)$ where $f_k\in K$ for  $k\in\Z$. Then~$\res(f):=f_{-1}\in K$
is the {\bf residue} of $f$.\index{residue}
\end{definition}

\noindent
The map $f\mapsto\res(f)\colon K(\!(x^{-1})\!)\to K$ is strongly additive and $K$-linear.

\begin{lemma}\label{resres}
Let $f\in K(\!(x^{-1})\!)$. Then  $\res(f')=0$; if $f\neq 0$, then $\res(f^\dagger) = -vf$.
\end{lemma}
\begin{proof} The first claim is clearly true. For the second, 
let $f=x^kg$ where~$k=-vf$, $g\in K(\!(x^{-1})\!)$, $g\asymp 1$. Then
$f^\dagger=kx^{-1}+g^\dagger$ with $g^\dagger\preceq x^{-2}$, so $\res(f^\dagger)=k=-vf$.
\end{proof}

\noindent
Just like Lemma~\ref{lem:resder} led to Corollaries~\ref{cor:Jacobi, 1} and~\ref{cor:Jacobi, 2}, Lemma~\ref{resres} gives:

\begin{cor}\label{cor:Jacobi, x^-1, 1}
If $f,g\in  K(\!(x^{-1})\!)$, then $\res(f'g)=-\res(fg')$ and thus, if also~$g\neq 0$ and $k\in\Z$, then
$\res(f'g^k)=-k\res  (fg^{k-1}g')$.
\end{cor}

\begin{cor}\label{cor:Jacobi, x^-1, 2}
If $f,g\in K(\!(x^{-1})\!)$, $f\succ 1$, then $\res\!\big( (g\circ f)f'\big) = -vf\,\res(g)$.
\end{cor}

\noindent
For $f\in K(\!(x^{-1})\!)$, $f \asymp x$, let $f^{[-1]}$ be the compositional inverse of $f$. One verifies easily
that if $f\in x+x^{-1}K[[x^{-1}]]$, then $f^{[-1]}\in x+x^{-1}K[[x^{-1}]]$. More generally: 

\begin{theorem} \label{thm:LIF, x^-1}
Let $f,g\in   x+x^{-1}K[[x^{-1}]]$. Then
$$g\circ f^{[-1]}\ =\  x - \sum_{n\geq 1}  \frac{1}{n}\res(g'f^{n})x^{-n}.$$
\end{theorem}
\begin{proof}
Let $h:= g\circ f^{[-1]}=\sum_{k} h_kx^k$ ($h_k\in K$).  Then for $k\in\Z^{\neq}$ we have
\begin{align*} \frac{1}{k} \res (g'f^{-k})\ &=\  \res\!\big( gf^{-k-1}  f' \big)\ =\ \res\!\big((h\circ f)f^{-k-1}f'\big)\\
&=\ \res\!\big((hx^{-k-1}\circ f)f'\big)\ =\  \res\!\big( hx^{-k-1}\big)\  =\  h_{k},
\end{align*} 
using Corollary~\ref{cor:Jacobi, x^-1, 1}  for the first equality and ~\ref{cor:Jacobi, x^-1, 2} for the next to last one.
This computation goes through for any $f, g\in K(\!(x^{-1})\!)$ with $f\asymp x$, but under the assumptions of the theorem gives the desired result. 
\end{proof}

%\noindent
%For  use in the subsection above we now study the case 
%$$f=x+\sum_{n\geq 1}f_n \frac{x^{-2n+1}}{n!}\qquad\text{($f_n\in K$ for $n\geq 1$)}$$
%in more detail; for such $f$ we have:

\begin{cor}\label{cor:LIF, x^-1} Let $f=x+\sum_{n\geq 1}f_n \frac{x^{-2n+1}}{n!}$, all $f_n\in K$, and $g=f^{[-1]}$. Then 
$$g\ =\ x-\sum_{j\geq 1} g_j \frac{x^{-2j+1}}{j!}\quad\text{where $g_j= \sum_{i=1}^{j}\frac{(2(j-1))!}{(2j-1-i)!}B_{ij}(f_1,\dots,f_{j-i+1})$.}$$
\end{cor}
\begin{proof}
Put $F:=\sum_{n\geq 1} f_n\frac{z^n}{n!}\in zK[[z]]$, $h:=\sum_{n\geq 1} f_n \frac{x^{-2n}}{n!}$.
Then $f=x({1+h})$, and~$\res(f^n)$ is the  coefficient of $x^{-n-1}$ in $f^n/x^n=(1+h)^n$.
Thus $\res(f^n)=0$ if~$n$ is even. Now suppose $n$ is odd, $n=2j-1$ ($j\geq 1$).
Then
the coefficient of~$x^{-n-1}=x^{-2j}$ in $(1+h)^n$ equals the coefficient of $z^j$ in the power series~$(1+F)^n=\sum_{i=0}^n \frac{n!}{(n-i)!}\frac{F^i}{i!}$, and this coefficient in turn is given by
$$\frac{1}{j!}\sum_{i=1}^j \frac{n!}{(n-i)!} B_{ij}(f_1,\dots,f_{j-i+1}).$$
Now use Theorem~\ref{thm:LIF, x^-1} with $x$ in the role of of $g$ there. 
\end{proof}
 
\noindent
Using the formulas for $B_{ij}$ for small values of $i$, $j$ given on [ADH, p.~554] we readily compute:
\begin{align*}
g_1&=f_1, \\
g_2&=f_2+2f_1^2, \\
g_3&=f_3+12f_1f_2+12f_1^3,\\ 
g_4&=f_4+24f_1f_3+18f_2^2+180f_1^2f_2+120f_1^4\\
g_5&=f_5+40f_1f_4+80f_2f_3+560f_1^2f_3+840f_1f_2^2+3660f_1^3f_2+1830f_1^5
\end{align*}

\begin{remark}
Suppose $K=\R$. Then $\R(\!(x^{-1})\!)$ is a subfield of   $\T$,
and the composition $(g,f)\mapsto g\circ f\colon \T\times\T^{>\R}\to\T$ in $\T$ (see the remarks after Corollary~\ref{cor:trdegellinv}) extends the composition in  $\R(\!(x^{-1})\!)$ defined above.
All $f\in \T^{>\R}$ have a compositional inverse~$f^{\operatorname{inv}}$ in $\T$,
with $f^{\operatorname{inv}}=f^{[-1]}$ if $f\in \R(\!(x^{-1})\!)$. For $f>\R$ in the subfield $\T_{\operatorname{g}}$ of $\T$ consisting of the grid based series we have $f^{\operatorname{inv}}\in \T_{\operatorname{g}}$, and \cite[Section~5.4.2]{JvdH}
has a formula for the coefficients of $f^{\operatorname{inv}}$ in that case. 
%for $f\in\T^{>\R}$
%(strictly speaking, for $f\in\T_{\operatorname{g}}^{>\R}$).
\end{remark}

\section{Holes and Slots in Perfect Hardy Fields}\label{sec:holes perfect}

\noindent
In this section $H\supseteq\R$ is a real closed Hardy field with asymptotic integration. We set
$K:=H[\imag]\subseteq\Calinf[\imag]$,
an algebraically closed  $\d$-valued extension 
of~$H$.
% and let $\fm$, $\fn$ range over~$H^\times$.
Moreover, $\hat H$ is an immediate $H$-field extension of $H$ and $\hat K:=\hat H[\imag]$ is  the corresponding immediate $\d$-valued extension of $K$ as in Section~\ref{sec:d-alg extensions}.
We also fix a $\d$-maximal Hardy field extension~$H_*$ of~$H$. The $H$-field $H_*$ is newtonian, and the $\d$-valued field extension~$K_*:=H_*[\imag]\subseteq \Calinf[\imag]$ of~$K$ is newtonian and linearly closed. 
$$\xymatrix@R=0.25em{K_* & & \hat K  \\
& \ar@{-}[ul]  K  \ar@{-}[ur]  & \\
H_* \ar@{-}[uu] & & \hat H \ar@{-}[uu] \\
& H \ar@{-}[ul]  \ar@{-}[ur]  \ar@{-}[uu] & 
}$$
Recall that if $\I(K)\subseteq K^\dagger$ and $A\in K[\der]^{\neq}$ splits over $K$, then $A$ is terminal. 
In this section we show:

\begin{theorem}\label{thm:d-perfect flabby}
Suppose $H$ is $\d$-perfect and $\upo$-free.
Then every minimal hole in~$K$ of positive order is flabby. Moreover, $H$ has no hole of order~$1$,  every minimal hole in $H$ of order $2$ is flabby,
and if  all $A\in H[\der]^{\neq}$  are terminal, then every minimal hole in $H$ of positive order  is flabby.
\end{theorem}

\noindent
In Corollary~\ref{cor:notorious 3.6 firm, 3, K} below we also show that if $H$ is $\d$-perfect (but not necessarily $\upo$-free), then  every linear minimal hole $(P,\fm,\hat f)$ in $K$ of order $1$ with $\hat f\in\hat K$ is flabby. 
%\marginpar{can we get rid of the assumption $\hat f\in\hat K$?}
(See the discussion after the proof of Lemma~\ref{lem:upo-freeness of the perfect hull, (ii)=>(iii)} for an example
of a $\d$-perfect Hardy field that is not $\upo$-free.) 

\medskip
\noindent
The theorem above originated in an attempt to characterize $\upo$-free $\d$-perfect Hardy fields
among Hardy fields containing $\R$ purely in terms of asymptotic differential algebra. We hope to return to this topic at a later occasion.

%lends credence to the following conjecture, which characterizes $\upo$-free $\d$-perfect Hardy fields among Hardy fields containing $\R$ purely in terms of asymptotic differential algebra.

%\begin{conjecture}
%Suppose $H\supseteq \R$ is Liouville closed and $\upo$-free with $\I(K)\subseteq K^\dagger$. Then
%\begin{align*}
%\text{$H$ is $\d$-perfect}&\quad\Longleftrightarrow\quad\text{every minimal hole in $H$ is flabby} \\
%&\quad\Longleftrightarrow\quad\text{every minimal hole in $K$ is flabby.}
%\end{align*} 
%\end{conjecture}

\medskip
\noindent
With the proof of Corollary~\ref{cor:notorious 3.6 firm, 3, K} we also finish the proof  of Theorem~\ref{thm:d-perfect flabby}.

\subsection*{Asymptotic similarity and equivalence of slots}
Let   $(P,\fm,\hat f)$  be a slot in~$H$ of order $\ge 1$ where $\hat f\in\hat H$. 
If $f\in\Calinf$ is $H$-hardian and
$(P,\fm,f)$ is a slot in~$H$ (we regard this as including the requirement that $f\notin H$ and
 the Hardy field $H\langle f\rangle$ is an immediate extension of $H$), then 
$$f\approx_H \hat f\quad\Longleftrightarrow\quad\text{$(P,\fm,f)$ and $(P,\fm,\hat f)$  are equivalent.}$$
Note that if $f\in\c$, $f\approx_H \hat f$, and $g,h\in H$, $g\neq 0$,  then $fg-h\approx_H\hat fg-h$.
From Corollary~\ref{cor:find zero of P, 2} and newtonianity of $H_*$ we get a useful result
about filling slots in~$H$ by elements of $\d$-maximal Hardy field extensions of~$H$:

\begin{lemma}\label{lem:find zero of P, d-max}
If $H$ is $\upo$-free and
  $(P,\fm,\hat f)$ is $Z$-minimal, then 
  there exists~${f\in H_*}$ 
such that $(P,\fm,f)$ is a hole in $H$ equivalent to $(P,\fm,\hat f)$, in particular, $P(f)=0$, $f\prec\fm$, and   $f\approx_H \hat f$.
\end{lemma}

\noindent
In Lemma~\ref{lem:find zero of P, d-max} we cannot drop the assumption that $H$ is $\upo$-free. 
To see why, suppose $H$ is $\d$-perfect and not $\upo$-free (such $H$ exists by Example~\ref{ex:counterex}), and
take~$\upo\in H$ and $(P,\fm,\upl)$ as in Lemma~\ref{lem:upl-free, not upo-free} for $H$ in the role of $K$ there, so~$P= 2Y'+Y^2+\upo$ and
$(P,\fm,\upl)$ is  a minimal hole  in 
$H$ by Corollary~\ref{cor:upl-free, not upo-free}.
Since~${\upo\notin\omega(H)}$ and $H$ is $1$-$\d$-closed in all its Hardy field extensions,   no $H$-hardian germ~$f$ satisfies~${P(f)=0}$. 
Thus the conclusion of Lemma~\ref{lem:find zero of P, d-max} fails for $\hat f=\upl$.

Corollary~\ref{cor:find zero of P, 3} yields a variant for $P$ of order $1$:

\begin{lemma}\label{lem:find zero of P, d-max, H, d=r=1}
If $H$ is $\upl$-free and $(P,\fm,\hat f)$ is  $Z$-minimal of order~$1$ with a quasi\-linear refinement,
then there exists $f\in H_*$ such that $H\langle f \rangle$ is an immediate extension of $H$ and
$(P,\fm,f)$ is a hole in $H$ equivalent to~$(P,\fm,\hat f)$.
\end{lemma}

\noindent
Here are complex versions of some of the above:
Let  $(P,\fm,\hat f)$  be a slot in~$K$ of order~$\ge 1$ where $\hat f\in\hat K$.
If~$f\in K_*$ and
$(P,\fm,f)$ is a slot in $K$ (so $f\notin K$ and~$K\langle f \rangle\subseteq K_*$ is an immediate extension of $K$),
then 
$$f\approx_K \hat f\quad\Longleftrightarrow\quad\text{$(P,\fm,f)$ and $(P,\fm,\hat f)$  are equivalent.}$$
If $f\in\c[\imag]$, $f\approx_K \hat f$ ,and $g,h\in K$, $g\neq 0$, then~$fg-h\approx_K \hat fg-h$.
Recall that $H$ is $\upo$-free iff $K$ is, by [ADH, 11.7.23].  Again by Corollaries~\ref{cor:find zero of P, 2} and \ref{cor:find zero of P, 3}:
%Here is the complex analogue of Lemma~\ref{lem:find zero of P, d-max}:

\begin{lemma}\label{lem:find zero of P, d-max, K} 
If $H$ is $\upo$-free and
  $(P,\fm,\hat f)$ is $Z$-minimal as a slot in $K$, then   there exists ${f\in K_*}$ 
such that $K\langle f \rangle$ is an immediate extension of $K$ and $(P,\fm,f)$ is a hole in~$K$ equivalent to $(P,\fm,\hat f)$
$($and thus $P(f)=0$, $f\prec\fm$, and  $f\approx_K \hat f)$.
\end{lemma}

\begin{lemma}\label{lem:find zero of P, d-max, H, d=r=1, K}
If $H$ is $\upl$-free and, as a slot over $K$, $(P,\fm,\hat f)$ is  $Z$-minimal of order~$1$ with a quasilinear refinement,
then there exists $f\in K_*$ such that $K\langle f \rangle$ is an immediate extension of $K$ and
$(P,\fm,f)$ is a hole in $K$ equivalent to~$(P,\fm,\hat f)$.
\end{lemma}

 %Next, $(P,\fm,\upl)$ is a $Z$-minimal hole in $K$ by Corollary~\ref{cor:holes in K vs holes in K[i]} and Lemma~\ref{lem:upl-free, not upo-free}. Suppose $\tilde H$ is an immediate Hardy field extension of~$H$  and $P(f)=0$,
 %$f\in \tilde H[\imag]$.
%Then $\upo\in\omega(\tilde H)$ or $\upo\in\sigma(\tilde H^\times)$ by~[ADH, p.~262].
%In the first case we obtain a contradiction to $\upo\notin \omega(H)$ as before, and in the second
%case we get a contradiction from Corollary~\ref{cor:upo}.
%Hence the conclusion of Lemma~\ref{lem:find zero of P, d-max, K}  fails for $\hat f=\upl$. 

\noindent
{\it In the rest of this section $H$ is Liouville closed  and $\I(K)\subseteq K^\dagger$.}\/ 
(These conditions are satisfied if $H$ is  $\d$-perfect.) {\it We take an $\R$-linear complement $\Lambda_H$ of~$\I(H)$
in~$H$, so~$\Lambda:= \Lambda_H\imag$ is a complement of $K^\dagger$ in $K$. Next we take an $\R$-linear complement~$\Lambda_{H_*}$ of~$\I(H_*)$ in $H_*$, so~$\Lambda_{*}:=\Lambda_{H_*}\imag$ is a complement of~$K_*^\dagger$ in $K_*$. Accordingly we identify in the usual way $\Univ:=\Univ_K:= K\big[\!\ex(\Lambda)\big]$ with $K[\ex^{H\imag}]$ and like\-wise~$\Univ_*:= \Univ_{K_*} :=  K_*\big[\!\ex(\Lambda_{*})\big]$ with~$K_*[\ex^{H_*\imag}]$.}   

\subsection*{Zeros of linear differential operators close to the linear part of a slot}
{\it In this subsection $(P,1,\hat h)$ with $\hat h\in\hat H$ is a  normal or linear slot in $H$ of order~$r\geq 1$}. Then $\order(L_P)=r$, so $\dim_{\C} \ker_{\Univ_*}L_P=r$ by Theorem~\ref{thm:lindiff d-max}. Lem\-ma~\ref{lem:excu properties}(ii)
then gives $\exc_{K_*}^{\operatorname{u}}(L_P)=v_{\operatorname{g}}(\ker^{\neq}_{\Univ_*} L_P)$.
If $L_P$ is terminal, then $\exc^{\operatorname{u}}(L_P)=\exc_{K_*}^{\operatorname{u}}(L_P)$, by Corollary~\ref{cor:excev cap GammaOmega}. We use these remarks to deal with firm and flabby cases:

\begin{lemma}\label{lem:8.8 firmLP}
Suppose $(P,1,\hat h)$ is firm and ultimate and $L_P$ is terminal.
Then there is no  $y\in\Calr[\imag]^{\ne}$  such that $L_P(y)=0$ and $y\prec 1$.
\end{lemma}
\begin{proof}
Suppose $y\in\Calr[\imag]^{\ne}$,  $L_P(y)=0$, and $y\prec 1$. Then $y\in \Univ_{*}$, so~$y\prec_{\g} 1$   by Lemma~\ref{lem:gaussian ext dom}. 
The remarks above 
give $ v_{\g}y\in  \exc^{\operatorname{u}}(L_P)$. Then~$v_{\g}y\leq 0$ by Re\-mark~\ref{rem:firm deg 1}, contradicting $y\prec_{\g}1$.
\end{proof}

\begin{lemma}\label{lem:8.8 flabbyLP} 
Suppose $(P,1,\hat h)$ is flabby.
Then there exists   $y\in \Calinf[\imag]^{\ne}$ such that~$L_P(y)=0$  and~$y\prec\fm$ for all $\fm\in H^\times$ with $v\fm\in v(\hat h -H)$.
If in addition~$(P,1,\hat h)$ is  $Z$-minimal, deep,  and special,  then $y',\dots,y^{(r)}\prec  \fm$ for all  such $y$ and~$\fm$. 
\end{lemma}
\begin{proof} 
Flabbiness  of  $(P,1,\hat h)$ and Lemmas~\ref{lem:firm normal} and~\ref{lem:firm deg 1} yield a $\gamma\in\exc^{\operatorname{u}}(L_P)$
with~$\gamma>v(\hat h-H)$.
Then $\gamma\in\exc^{\operatorname{u}}_{K_*}(L_P)$ by 
Corollary~\ref{cor:LambdaL, purely imag}, so a remark above gives~$y\in  \ker^{\neq}_{\Univ_{*}} L_P$ such that~$v_{\g}y=\gamma$. 
Then $y\prec_{\g}\fm$ and thus $y\prec\fm$, for all~$\fm\in H^\times$  with $v\fm\in v(\hat h -H)$.
 For the remainder, use Lemma~\ref{lem:8.8 refined}.
\end{proof}

\noindent
Next we consider a suitable perturbation $A$ of $L_P$: {\em In the rest of this subsection we assume $L_P=A+B$
with~$A,B\in K[\der]$ satisfying}
 $$\order (A)\ =\ r, \qquad  \fv\ :=\ \fv(A)\prec^\flat 1, \qquad B\prec_{\Delta(\fv)} \fv^{r+1}A.$$
Then Lemma~\ref{lem:fv of perturbed op} gives   $\fv(L_P)\sim\fv$.
By Lem\-ma~\ref{lem:excu properties}(ii),(iii),
$$v_{\operatorname{g}}(\ker^{\neq}_{\Univ_*} A)\ =\ \exc_{K_*}^{\operatorname{u}}(A)\ =\ \exc_{K_*}^{\operatorname{u}}(L_P),\qquad
\exc^{\operatorname{u}}(A)\ =\ \exc^{\operatorname{u}}(L_P). $$
If $A$ is terminal, then all five displayed sets are equal by Corollary~\ref{cor:excev cap GammaOmega}.
Recall also from Corollary~\ref{cor:ultimate prod, 2} that if $A$ splits over $K$, then $A$ is terminal, and from Proposition~\ref{prop:finiteness of excu(A), real} that if $\dim_{\C} \ker_{\Univ} A =r$, then $A$ is terminal. 

We can now generalize Proposition~\ref{lem3.9}: 

\begin{prop}\label{lem3.9, generalized}  
Suppose $(P,1,\hat h)$ is ultimate,  $A$ is terminal, and  $y\in\Calr[\imag]$ satisfies $A(y)=0$, $y\prec 1$.
Then~$y\prec \fm$ for all~$\fm\in H^\times$ with~$v\fm\in v(\hat h-H)$.
\end{prop}
\begin{proof}
We have $y\in \Univ_{*}$, so~$y\prec_{\g} 1$  by Lemma~\ref{lem:gaussian ext dom}. 
If $y=0$, then we are done, so suppose $y\neq 0$.
Then $0 < v_{\g}y\in  \exc^{\operatorname{u}}(L_P)$ by remarks before Proposition~\ref{lem3.9, generalized}.  
Hence~$v_{\g}y > v(\hat h -H)$ by  Lemma~\ref{lem:ultimate normal} if~$(P,1,\hat h)$ is normal, and by 
Lemma~\ref{lem:ultimate deg 1} if~$(P,1,\hat h)$ is linear,
so $y\prec_{\g}\fm$ for all $\fm\in H^\times$ with~$v\fm\in v(\hat h-H)$, and thus $y\prec \fm$ for all such $\fm$ by Corollary~\ref{cor:gaussian ext dom}. 
\end{proof}

%\noindent
%The same way Corollary~\ref{cor:8.8 refined} followed from Proposition~\ref{lem3.9}, this yields:

\begin{cor}\label{cor:8.8 refined, generalized}
Suppose $A$ is terminal and $(P,1,\hat h)$ is $Z$-minimal, deep, ultimate, and special. If~$y\in\Calr[\imag]$ satisfies $A(y)=0$ and $y\prec 1$,
then $y,y',\dots,y^{(r)}\prec \fm$ for all $\fm\in H^\times$ with~$v\fm\in v(\hat h-H)$.
\end{cor}
\begin{proof} First use Proposition~\ref{lem3.9, generalized} and then Lemma~\ref{lem:8.8 refined}. 
\end{proof} 

\noindent
Next we turn to firm and flabby cases. 

\begin{lemma}\label{lem:8.8 firm}
Suppose $A$ is terminal and $(P,1,\hat h)$ is firm and ultimate.
Then there is no   $y\in\Calr[\imag]^{\ne}$   such that $A(y)=0$ and $y\prec 1$.
\end{lemma}
\begin{proof}
Suppose $y\in\Calr[\imag]^{\ne}$,  $A(y)=0$, and $y\prec 1$. Then $y\in \Univ_{*}$, so~$y\prec_{\g} 1$   by Lemma~\ref{lem:gaussian ext dom}. 
The remarks before Proposition~\ref{lem3.9, generalized}  
give $ v_{\g}y\in  \exc^{\operatorname{u}}(L_P)$. Hence~$v_{\g}y\leq 0$ by Remark~\ref{rem:firm deg 1}, contradicting $y\prec_{\g}1$.
\end{proof}

\begin{lemma}\label{lem:8.8 flabby} 
Suppose $(P,1,\hat h)$ is flabby.
Then there exists   $y\in \Calinf[\imag]^{\ne}$ such that~$A(y)=0$  and~$y\prec\fm$ for all $\fm\in H^\times$ with $v\fm\in v(\hat h -H)$.
If in addition~$(P,1,\hat h)$ is  $Z$-minimal, deep,  and special,  then $y',\dots,y^{(r)}\prec  \fm$ for all  such $y$ and   $\fm$. 
\end{lemma}
\begin{proof}
Flabbiness    of  $(P,1,\hat h)$ and Lemmas~\ref{lem:firm normal} and~\ref{lem:firm deg 1} yield a $\gamma\in\exc^{\operatorname{u}}(L_P)=\exc^{\operatorname{u}}(A)$
with~$\gamma>v(\hat h-H)$. The rest of the proof is the same as that of Lemma~\ref{lem:8.8 flabbyLP} with $A$ instead of $L_P$. 
%Then $\gamma\in\exc^{\operatorname{u}}_{K_*}(A)$ by 
%Corollary~\ref{cor:LambdaL, purely imag}.
%Since~$K_*$ is $r$-linearly newtonian, Lemma~\ref{lem:excu 1} gives a~$y\in  \ker^{\neq}_{\Univ_{*}} A$ such that~$v_{\g}y=\gamma$. 
%Then $y\prec_{\g}\fm$ and thus $y\prec\fm$, for all $\fm\in H^\times$  with $v\fm\in v(\hat h -H)$.
 %For the remainder, use Lemma~\ref{lem:8.8 refined}.
\end{proof}

\begin{remark}
The material above in this subsection goes through if instead of 
$(P,1,\hat h)$ with $\hat h\in\hat H$ being a normal or linear slot in $H$ of order $r\ge 1$ we assume
$(P,1,\hat h)$ with $\hat h\in\hat K$ is a normal or linear slot in $K$ of order $r\ge 1$, and  ``$\fm\in H^\times$ with~${v\fm\in v(\hat h -H)}$''
is replaced everywhere by ``$\fm\in K^\times$ with $v\fm\in v(\hat h -K)$''. 

To see this, use the $K$-versions of Lemmas~\ref{lem:ultimate normal},~\ref{lem:ultimate deg 1},
\ref{lem:firm normal}, of Lemma~\ref{lem:firm deg 1}  and Remark~\ref{rem:firm deg 1},
and of Lemma~\ref{lem:8.8 refined};
cf.~the discussion at the end of the subsection {\it An application to slots in $H$}\/ of Section~\ref{sec:ueeh}.
\end{remark}

\subsection*{Application to linear slots}
In this subsection  we apply the material in the last subsection to the study of linear slots (in $H$ and in $K$).
Until further notice  $(P,\fm,\hat h)$ with~$\hat h\in\hat H$ is a $Z$-minimal  linear slot in $H$ of order $r\geq 1$. 

\begin{lemma}  \label{lem:P(f)=0, P linear}
There exists $f\in H_*$ such that $P(f)=0$ and $f\prec\fm$.
\end{lemma}
\begin{proof}
We may replace  $(P,\fm,\hat h)$ by a refinement whenever convenient.
Hence by Remark~\ref{rem:achieve isolated} we may arrange that $(P,\fm,\hat h)$ is isolated.
Then $P(0)\neq 0$, and~$\gamma:=v\hat h$ is the unique element of $\Gamma\setminus\exc^{\ev}(L_P)$
such that $v_{L_P}^{\ev}(\gamma)=v(P(0))$, by Lemmas~\ref{lem:from cracks to holes} and~\ref{lem:isolated, d=1}.
Now $H_*$ is linearly newtonian, so Corollary~\ref{cor:14.2.10, generalized} yields $f\in H_*^\times$  
 with~$P(f)=0$, $vf\notin \exc^{\ev}_{H^*}(L_P)$, and $v_{L_P}^{\ev}(vf)=v(P(0))$.
 By Corollary~\ref{cor:13.7.10}, $v_{L_P}^{\ev}(\gamma)$ does not change when passing from $H$ to $H_*$, and $\gamma\notin\exc^{\ev}_{H_*}(L_P)$. Thus~$vf=\gamma$ by Lemma~\ref{lem:ADH 14.2.7}; in particular, $f\prec\fm$.
\end{proof}

\begin{cor}\label{cor:8.8 flabby}
Suppose  $(P,\fm,\hat h)$ is flabby,  and $f\in\c^r$ is such that $P(f)=0$ and $f\prec \fm$.
Then $f\in \Calinf$ and there exists $g\in\Calinf$ such that $P(g)=0$, $g\prec \fm$,     and $0\neq f-g\prec\fn$ for all $\fn\in H^\times$ with
$v\fn\in v(\hat h-H)$.
For any such $g$ we have  
$$f\approx_H\hat h\ \Rightarrow\ g\approx_H\hat h,\quad
  H\subseteq\Ginf\ \Rightarrow\ f,g\in\Ginf, \quad H\subseteq\Gom\ \Rightarrow\ f,g\in\Gom.$$
If~$(P,\fm,\hat h)$ is also deep and special, then $f-g\in \fm(\c^r)^{\prec}$ for any such $g$.
\end{cor}
\begin{proof} Lemma~\ref{smo} gives $f\in \Calinf$. 
Replace $(P,\fm,\hat h)$, $f$, by $(P_{\times\fm},1,\hat h/\fm)$, $f/\fm$ to arrange $\fm=1$.
Lem\-ma~\ref{lem:8.8 flabbyLP} then   yields $y\in \Calinf[\imag]^{\ne}$ such that $L_P(y)=0$  and~$y\prec \fn$
for all $\fn\in H^\times$ with $v\fn\in v(\hat h-H)$.
Replacing $y$ by $\Re y$ or $\Im y$ we arrange $y\in\Calinf$.
Then~$g:=f+y\in\Calinf$ satisfies $f\neq g$, $P(g)=0$, and~$g\prec 1$. 
The rest follows from remarks after Corollary~\ref{cor:existence and uniqueness, real} and from Lemma~\ref{lem:8.8 flabbyLP}.
\end{proof}

\noindent
In the proof of the next corollary we use that if $H$ is $\upo$-free, then
Lemma~\ref{lem:find zero of P, d-max}  yields an~$f\in H_*$ such that~${P(f)=0}$ and~$f\approx_H \hat h$.

\begin{cor}\label{cor:8.8 firm} 
Suppose   $(P,\fm,\hat h)$ is ultimate and $L_P$ is terminal. Then 
$$\text{$(P,\fm,\hat h)$ is firm}\quad\Longleftrightarrow\quad\text{there is a unique $f\in\c^r$ with $P(f)=0$ and $f\prec \fm$.}$$ 
If  $(P,\fm,\hat h)$ is firm, $f\in\c^r$, $P(f)=0$, $f\prec \fm$, then $f\in\Dx(H)$, there is no $g\ne f$ in~$\c^r[\imag]$ with
$P(g)=0$, $g\prec \fm$,  and if in addition $H$ is $\upo$-free, then $f\approx_H\hat h$.
\end{cor}
\begin{proof} We arrange $\fm=1$ as before.
Then Lem\-mas~\ref{lem:8.8 firmLP}  and~\ref{lem:P(f)=0, P linear} yield~``$\Rightarrow$''. 
For~``$\Leftarrow$'' and the rest, use
Corollary~\ref{cor:8.8 flabby} and the remark after its proof, and observe that  our $\d$-maximal Hardy field extension $H_*$ of $H$ was arbitrary.
\end{proof}

\begin{cor}\label{cor:linear flabby}
Suppose $H$ is $\upo$-free and $\d$-perfect, and all $A\in H[\der]\subseteq K[\der]$ of order $r$ are terminal. Then every $Z$-minimal linear slot in $H$ of order~$r$ is flabby.
\end{cor}
\begin{proof}
Given  a firm $Z$-minimal linear slot in $H$ of order $r$ we use Remark~\ref{rem:achieve ultimate} and Lemma~\ref{lem:firm refinement} to refine it to be ultimate.
So we arrive at an ultimate firm $Z$-minimal linear slot  in $H$ of order $r$ with terminal linear part.
This contradicts~$H$ being $\d$-perfect by Corollary~\ref{cor:8.8 firm}.
\end{proof}

\noindent
Next the $K$-versions of Lemma~\ref{lem:P(f)=0, P linear} and its corollaries: Let  
$(P,\fm,\hat f)$ with~$\hat f\in\hat K$ be a $Z$-minimal linear slot in $K$ of order~$r\geq 1$.
Now $K$ is $\upl$-free [ADH, 11.6.8], so we can mimick the proof of Lemma~\ref{lem:P(f)=0, P linear} to obtain:

\begin{lemma}  \label{lem:P(f)=0, P linear, K}
There exists $f\in K_*$ such that $P(f)=0$ and $f\prec\fm$.
\end{lemma}

\noindent
The $K$-version of Lemma~\ref{lem:8.8 flabbyLP} leads to the $K$-version of Corollary~\ref{cor:8.8 flabby}:

\begin{cor}\label{cor:8.8 flabby, K}
Suppose  $(P,\fm,\hat f)$ is flabby,  and $f\in\c^r[\imag]$, $P(f)=0$, and~${f\prec \fm}$.
Then $f\in \Calinf[\imag]$ and  there exists $g\in\Calinf[\imag]$ such that $P(g)=0$, $g\prec \fm$,  and $0\neq f-g\prec\fn$ for all
$\fn\in K^\times$ with $v\fn\in v(\hat f-K)$.
For any such $g$ we have  
$$f\approx_K\hat f\ \Rightarrow\ g\approx_K\hat f,\quad
  H\subseteq\Ginf\ \Rightarrow\ f,g\in\Ginf[\imag], \quad H\subseteq\Gom\ \Rightarrow\ f,g\in\Gom[\imag].$$
If~$(P,\fm,\hat f)$ is also deep and special, then $f-g\in \fm \c^r[\imag]^{\prec}$ for any such $g$.
%we can additionally choose~$g\in \c^r[\imag]^{\prec}$.
\end{cor}

\noindent
If $H$ is $\upo$-free, then
Lemma~\ref{lem:find zero of P, d-max, K}  yields $f\in K_*$ with~${P(f)=0}$ and~$f\approx_K \hat f$.
This remark and $K$-versions of various results yield the $K$-version of
% Lemma~\ref{lem:8.8 firm}, Lemma~\ref{lem:P(f)=0, P linear, K}, and  Corollary~\ref{cor:8.8 flabby, K} instead of
%Lemma~\ref{lem:8.8 firm}, Lemma~\ref{lem:P(f)=0, P linear}, and  Corollary~\ref{cor:8.8 flabby}, respectively,  in
Corollary~\ref{cor:8.8 firm}:

\begin{cor}\label{cor:8.8 firm, K}
Suppose $(P,\fm,\hat f)$ is ultimate and $L_P$ is terminal. Then:
$$\text{$(P,\fm,\hat f)$ is firm} \quad\Longleftrightarrow\quad \text{there is a unique $f\in\c^r[\imag]$ with $P(f)=0$ and $f\prec \fm$.}$$
If $(P,\fm,\hat f)$ is firm, $f\in \c^r[\imag]$, $P(f)=0$, $f\prec \fm$, then 
 $f\in\Dx(H)[\imag]$, and also $f\approx_K \hat f$ in case $H$ is $\upo$-free. 
\end{cor}

\noindent
Using $K$-versions of various results (like Remark~\ref{rem:achieve ultimate, K} instead of Remark~\ref{rem:achieve ultimate}), then yields the $K$-version of Corollary~\ref{cor:linear flabby}:
%Finally, using Remark~\ref{rem:achieve ultimate, K}, the $K$-version of Lemma~\ref{lem:firm refinement}, and Corollary~\ref{cor:8.8 firm, K} instead of
%the remark after   Proposition~\ref{propachieve ultimate}, Lemma~\ref{lem:firm refinement}, and
% Corollary~\ref{cor:8.8 firm}, respectively, similarly to Corollary~\ref{cor:linear flabby}  one shows:

\begin{cor}\label{cor:linear flabby, K}
If $H$ is $\upo$-free and $\d$-perfect, and all $A\in K[\der]$ of order~$r$ are terminal, then every $Z$-minimal linear slot in $K$ of order~$r$ is flabby.
\end{cor}

\noindent
Linear slots in $K$ of order $1$ are $Z$-minimal, and we can say more in this case:

\begin{cor}\label{cor:8.8 firm, K, r=1}
Suppose $r=1$.
If $(P,\fm,\hat f)$ is flabby, then it is ultimate. Moreover, $L_P$ is terminal, so if  $(P,\fm,\hat f)$ is firm and ultimate, then
there is a unique~$f\in\c^r[\imag]$ with $P(f)=0$ and $f\prec \fm$, and for this $f$ we have:  $f\in\Dx(H)[\imag]$, with~$f\approx_K \hat f$ in case $H$ is $\upo$-free or $(P,\fm,\hat f)$ is a hole in $K$.
\end{cor}
\begin{proof} 
Corollary~\ref{cor:flabby, r=1}(i) gives ``flabby $\Rightarrow$ ultimate''.  
Since $L_P$ has order $1$, it is terminal. Now use Corollary~\ref{cor:8.8 firm, K}, and  Lemma~\ref{lem:find zero of P, d-max, H, d=r=1}, noting in connection with that lemma that the linear slot $(P,\fm, \hat f)$ is quasilinear. 
\end{proof}

\noindent
Our next goal is to establish refinements of Proposition~\ref{prop:notorious 3.6} for the case of firm and flabby slots in $H$:
Lemmas~\ref{lem:notorious 3.6 firm}, \ref{lem:notorious 3.6 firm, abs norm}, \ref{lem:notorious 3.6 firm, 1}, and~\ref{lem:notorious 3.6 flabby} below. Towards this goal, we introduce yet another useful concept of normality
for slots.

\subsection*{Absolutely normal slots in $H$}
{\it In this subsection $(P,\fm,\hat h)$ is a slot  in~$H$ of order~$r\geq 1$ with $\hat h\in \hat H$.}\/ Given active   $\phi>0$ in $H$ we take $\ell\in H$ with $\ell'=\phi$, and set $f^\circ:= f\circ \ell^{\inv}$ for
$f\in \c[\imag]$, as usual; see Section~\ref{secfhhf}. Recall from Section~\ref{sec:Hardy fields} that $H^\circ$ is Liouville closed with $K^\circ=H^\circ[\imag]$
and $\I(K^\circ)\subseteq (K^\circ)^\dagger$, and that $H_*^{\circ}$ is a $\d$-maximal Hardy field extension of $H^\circ$, with
$K_*^{\circ}=H_*^{\circ}[\imag]$. 

Since $K_*$ is linearly closed, each linear differential operator~$A\in K[\der]^{\neq}$  splits   over~$K_*$.
If $A\in K[\der]^{\neq}$  splits strongly over~$K_*$,
then by Theorem~\ref{thm:transfer} this remains true when~$K_*$ is replaced by $K_{**}:=H_{**}[\imag]\subseteq \Calinf[\imag]$ for any $\d$-maximal
Hardy field extension $H_{**}$ of $H$. 
We say that $(P, \fm, \hat h)$ is {\bf   absolutely normal} if  it is strictly normal and its linear part splits strongly over $K_*$. 
If~$(P,\fm,\hat h)$ is absolutely normal, then so is~$(P_{\times\fm},1,\hat h/\fm)$.\index{normal!absolutely}\index{absolutely normal}\index{slot!absolutely normal} Moreover:

\begin{lemma}\label{lem:abs norm comp conj}
Suppose $(P,\fm,\hat h)$ is absolutely normal, and  $\phi$ is active in $H$ with~$0<\phi\preceq 1$. Then
the slot $(P^{\phi\circ},\fm^\circ, \hat h^\circ)$ in $H^\circ$ is absolutely normal.
\end{lemma}
\begin{proof}
By Lemma~\ref{lem:normality comp conj, strong}, the slot $(P^\phi,\fm,\hat h)$ in $H^\phi$ is strictly normal, hence so is the slot 
$(P^{\phi\circ},\fm^\circ, \hat h^\circ)$ in $H^\circ$.
By Lemma~\ref{lem:split strongly compconj} the linear part $L_{P_{\times\fm}^\phi}=(L_{P_{\times\fm}})^\phi$ of~$(P^\phi,\fm,\hat h)$ splits strongly over $K_*^\phi=H_*^\phi[\imag]$, hence the linear part of 
$(P^{\phi\circ},\fm^\circ, \hat h^\circ)$ splits strongly over $K_*^{\circ}=H_*^{\circ}[\imag]$. 
\end{proof}

\noindent 
Next we show how to achieve absolute normality:

\begin{prop}\label{prop:achieve abs norm}
Suppose  $(P,\fm,\hat h)$ is $Z$-minimal, deep, and strictly normal, and~$\hat h\prec_{\Delta(\fv)}\fm$ with $\fv:=\fv(L_{P_{\times\fm}})$. 
Then for all sufficiently small~$q\in\Q^>$, $(P,\fn,\hat h)$, for any~$\fn\asymp\fm\abs{\fv}^q$ in $H^\times$, is a deep, absolutely normal refinement of~$(P,\fm,\hat h)$.
\end{prop}
\begin{proof}
Recall that $L_{P_{\times\fm}}$ splits
over the $H$-asymptotic extension~$K_*$ of~$K$. The argument in the proof of Corollary~\ref{cor:split strongly multconj} shows that
for all sufficiently small $q\in\Q^>$, $(P,\fn,\hat h)$, for any~$\fn\asymp \fm\abs{\fv}^q$ in~$H^\times$, is a steep refinement 
of~$(P,\fm,\hat h)$ whose linear part   splits strongly over~$K_*$.
By Corollary~\ref{cor:steep refinement}
any such refinement~$(P,\fn,\hat h)$ of $(P,\fm,\hat h)$ is   deep, and  for all sufficiently small $q\in\Q^>$, any such 
refinement of~$(P,\fm,\hat h)$  also remains strictly normal, by Lemma~\ref{lem:strongly normal refine, 2}
and Remark~\ref{rem:strongly normal refine, 2}. 
\end{proof}

\begin{cor}\label{cor:H*-normal} If $(P, \fm, \hat h)$ is $Z$-minimal, deep, normal, and  special, then $(P, \fm, \hat h)$ has  a   deep, absolutely normal   refinement.
\end{cor}

\noindent
This follows from  Corollary~\ref{cor:achieve strong normality, 1} and Proposition~\ref{prop:achieve abs norm}.
 
\begin{lemma}\label{cor:H*-normal, 1}
Suppose $H$ is $\upo$-free and $(P,\fm,\hat h)$  is $Z$-minimal and special.
Then there are  a refinement 
$(P_{+h},\fn,\hat h-h)$ of $(P,\fm,\hat h)$  and an active $\phi>0$ in~$H$ such that
the slot $(P_{+h^\circ}^{\phi\circ},\fn^\circ,\hat h^\circ-h^\circ)$ in~$H^\circ$ is deep, absolutely normal, and ultimate.
\end{lemma}
\begin{proof}
For any active $\phi>0$ in $H$  we may replace~$H$,~$(P,\fm,\hat h)$  by~$H^\circ$,~$(P^{\phi\circ},\fm^\circ,\hat h^\circ)$, and we may
 also replace~$(P,\fm,\hat h)$ by any of its refinements. 
Since $H$ is $\upo$-free, Proposition~\ref{varmainthm}
yields a refinement $(P_{+h},\fn,\hat h-h)$  of  $(P,\fm,\hat h)$ and an active~$\phi>0$ in $H$
such that the slot $(P_{+h^\circ}^{\phi\circ},\fn^\circ,\hat h^\circ-h^\circ)$ in $H^\circ$ is   normal.
Replacing $H$,~$(P,\fm,\hat h)$  by~$H^\circ$,~$(P_{+h^\circ}^{\phi\circ},\fn^\circ,\hat h^\circ-h^\circ)$
we arrange that $(P,\fm,\hat h)$  is  normal.
 Proposition~\ref{prop:achieve ultimate} now yields an ultimate refinement of $(P,\fm,\hat h)$.
Applying Proposition~\ref{varmainthm}  to this refinement and using Lemma~\ref{lem:ultimate refinement}, we  obtain an ultimate refinement  $(P_{+h},\fn,\hat h-h)$ of~$(P,\fm,\hat h)$  and an active $\phi>0$ in $H$ 
 such that~$(P_{+h^\circ}^{\phi\circ},\fn^\circ,\hat h^\circ-h^\circ)$ is deep, normal, and ultimate.
 Again replacing~$H$,~$(P,\fm,\hat h)$  by~$H^\circ$,~$(P_{+h^\circ}^{\phi\circ},\fn^\circ,\hat h^\circ-h^\circ)$,  
we   arrange that~$(P,\fm,\hat h)$  is deep, normal, and ultimate. Now apply 
Corollary~\ref{cor:H*-normal}   to $(P,\fm,\hat h)$ and use Lemma~\ref{lem:ultimate refinement}.
\end{proof}

\begin{cor}\label{cor:H*-normal, 2}
Suppose $H$ is $\upo$-free and $r$-linearly newtonian, and  $(P,\fm,\hat h)$  is  $Z$-minimal.
Then the conclusion of Lemma~\ref{cor:H*-normal, 1} holds.
\end{cor}
\begin{proof}
As in the beginning of the proof of  Lemma~\ref{cor:H*-normal, 1}, use Theorem~\ref{mainthm}  to arrange that $(P,\fm,\hat h)$  is  normal.
Then~$(P,\fm,\hat h)$ is quasilinear
by Corollary~\ref{cor:normal=>quasilinear} 
and hence special 
by Lemma~\ref{lem:special dents}, so Lemma~\ref{cor:H*-normal, 1} applies to~$(P,\fm,\hat h)$.
\end{proof}

\begin{remarkNumbered}\label{rem:H*-normal, 2}
By Corollary~\ref{degmorethanone} the hypotheses of Corollary~\ref{cor:H*-normal, 2} are satisfied if 
$H$ is $\upo$-free and $(P,\fm,\hat h)$  is a nonlinear minimal hole in $H$. 
\end{remarkNumbered}

\subsection*{Absolutely normal slots in $K$}
{\it Let $(P,\fm,\hat f)$ be a  slot in~$K$ of order $r\geq 1$, with~$\hat f\in \hat K$.}\/
Call $(P,\fm,\hat f)$  {\bf   absolutely normal} if  it is strictly normal and~$L_{P_{\times\fm}}$ splits strongly over $K_*$.\index{normal!absolutely}\index{absolutely normal}\index{slot!absolutely normal}
If~$(P,\fm,\hat f)$ is absolutely normal, then so is~$(P_{\times\fm},1,\hat f/\fm)$.
If $(Q,\fn,\hat h)$ is a slot in $H$ of order $\ge 1$ with $\hat h\in\hat H\subseteq \hat K$ , then it is a slot in~$K$ (Corollary~\ref{cor:holes in K vs holes in K[i]}), and
$(Q,\fn,\hat h)$ is absolutely normal as a slot in $H$ iff it is absolutely normal as a slot in $K$.

\begin{lemma}
Suppose  $(P,\fm,\hat f)$ is absolutely normal.
Then there exists $y$ in $\Gi[\imag]\cap \fm\,\c^r[\imag]^{\preceq}$ such that $P(y)=0$ and $y\prec \fm$.
If $H\subseteq\c^\infty$, then any such $y$ lies in~$\c^\infty[\imag]$, and likewise with $\c^\omega$ in place of $\c^\infty$.  
\end{lemma}
\begin{proof} Use Lemma~\ref{prop:as equ 1} with $Q=(P_{\times\fm})_1$ and $H_*$, $K_*$, $P_{\times\fm}$ in place of $H$, $K$, $P$. For the last part, use Corollary~\ref{cor:ADE smooth} as in the proof of that lemma. 
\end{proof}

\noindent
The $K$-versions of Lemma~\ref{lem:abs norm comp conj},  Proposition~\ref{prop:achieve abs norm}, and Corollary~\ref{cor:H*-normal} (with the same proofs) are as follows:

\begin{lemma}\label{lem:abs norm comp conj,K}
Suppose $(P,\fm,\hat f)$ is absolutely normal, and  $\phi$ is active in $H$ with~$0<\phi\preceq 1$. Then
the slot $(P^{\phi\circ},\fm^\circ, \hat h^\circ)$ in $K^\circ$ is absolutely normal.
\end{lemma}

\begin{prop}\label{prop:achieve abs norm,K}
Suppose  $(P,\fm,\hat f)$ is $Z$-minimal, deep, and strictly normal, and~$\hat f\prec_{\Delta(\fv)}\fm$ with $\fv:=\fv(L_{P_{\times\fm}})$. 
Then for all sufficiently small~$q\in\Q^>$, $(P,\fn,\hat f)$, for any~$\fn\asymp\fm\abs{\fv}^q$ in $K^\times$, is a deep, absolutely normal refinement of~$(P,\fm,\hat f)$.
\end{prop}

\begin{cor}\label{cor:K*-normal} If $(P, \fm, \hat f)$ is $Z$-minimal, deep, normal, and  special, then $(P, \fm, \hat f)$ has  a   deep, absolutely normal  refinement.
\end{cor}

\noindent
The $K$-version of Lemma~\ref{cor:H*-normal, 1}  is as follows (its proof  uses Proposition~\ref{prop:achieve ultimate, K} instead of Proposition~\ref{prop:achieve ultimate}, and the $K$-version of  Lemma~\ref{lem:ultimate refinement}):

\begin{lemma}\label{cor:K*-normal, 1}
Suppose $H$ is $\upo$-free and $(P,\fm,\hat f)$  is $Z$-minimal and special.
Then there are  a refinement 
$(P_{+f},\fn,\hat f-f)$ of $(P,\fm,\hat f)$  and an active $\phi>0$ in~$H$ such that
the slot $(P_{+f^\circ}^{\phi\circ},\fn^\circ,\hat f^\circ-f^\circ)$ in~$K^\circ$ is deep, absolutely normal, and ultimate.
\end{lemma}

\noindent
The $K$-version of Corollary~\ref{cor:H*-normal, 2} now follows in the same way:

\begin{cor}\label{cor:K*-normal, 2}
Suppose $H$ is $\upo$-free, $K$ is $r$-linearly newtonian, and  $(P,\fm,\hat f)$  is  $Z$-minimal.
Then the conclusion of Lemma~\ref{cor:K*-normal, 1} holds.
\end{cor}

%\noindent
%By Corollary~\ref{degmorethanone} the hypotheses of Corollary~\ref{cor:K*-normal, 2} are satisfied if 
%$H$ is $\upo$-free and $(P,\fm,\hat f)$  is a nonlinear minimal hole in $K$. 

\subsection*{Firm slots in $H$}
{\it In this subsection $(P,\fm,\hat h)$ is a slot   in~$H$ of order~$r\geq 1$, with~$\hat h\in \hat H$.}\/
We set~$d:=\deg(P)$, $w:=\wt(P)$, and
begin with a significant strengthening of Proposition~\ref{prop:notorious 3.6} for firm slots in $H$:

\begin{lemma} \label{lem:notorious 3.6 firm} 
Suppose $(P,\fm,\hat h)$ is firm, ultimate, and strongly  split-normal, and  let~$f,g\in \Calr[\imag]$ satisfy~$P(f)=P(g)=0$ and $f,g\prec \fm$. Then $f=g$.
\end{lemma}
\begin{proof}
The proof is similar to that of Proposition~\ref{prop:notorious 3.6}.
We first replace $(P,\fm,\hat h)$, $f$, $g$ by $(P_{\times\fm},1,\hat h/\fm)$, $f/\fm$, $g/\fm$, to arrange $\fm=1$.
We set $\fv:=\abs{\fv(L_P)} \prec^\flat 1$ and   $\Delta:=\Delta(\fv)$,
and take $Q,R\in H\{Y\}$ where $Q$ is homogeneous of degree $1$ and order $r$, 
$A:=L_Q\in H[\der]$ splits strongly  over~$K$,
 $P=Q-R$, and $R\prec_\Delta \fv^{w+1}P_1$, so $\fv(A)\sim \fv(L_P)$. Multiplying $P$, $Q$, $R$ by some $b\in H^\times$ we
arrange that  $A=\der^r+ f_1\der^{r-1}+\cdots + f_r$ with $f_1,\dots, f_r\in H$ and $R\prec_\Delta\fv^w$. We have
\begin{equation}\label{eq:notorious 3.6 firm}
A\ =\ (\der-\phi_1)\cdots (\der-\phi_r), \quad \phi_1,\dots,\phi_r \in K, \quad \Re\phi_1,\dots,\Re \phi_r\ \succeq\ \fv^\dagger\ \succeq\ 1.
\end{equation}
Corollary~\ref{cor:bound on linear factors} yields $\phi_1,\dots,\phi_r\preceq\fv^{-1}$. Take $a_0\in \R$ and functions on $[a_0,\infty)$ representing the germs $\phi_1,\dots, \phi_r$, $f_1,\dots, f_r$, $f$, $g$ and the $R_{\j}$ with $\j\in \N^{1+r}$, $|\j|\le d$, $\|\j\|\le w$ (using the same symbols for the germs mentioned as for their chosen representatives)
so as to be in the situation described in the beginning of Section~\ref{sec:split-normal over Hardy fields}, with $f$ and $g$ solutions on $[a_0,\infty)$ of the differential equation \eqref{eq:ADE} there. 
As there, we take $\nu\in\Q$ with $\nu > w$ so that~$R \prec_\Delta \fv^\nu$ and $\nu\fv^\dagger\not\sim \Re\phi_j$
for~$j=1,\dots,r$, and then increase $a_0$ to satisfy all assumptions for Lemma~\ref{bdua}.  

With $a\ge a_0$ and $h_a\in \Car[\imag]$ as in Lemma~\ref{lem:close} we have $A_a(h_a)=0$ and~${h_a\prec 1}$. 
%Moreover,  $\dim_{\C}\ker_{\Univ} A =r$ by Lemma~\ref{lem:full kernel}. 
Now $A$ splits over $K$, so $A$ is terminal as an element of $K[\der]$, by Corollary~\ref{cor:ultimate prod, 2}. As~$(P,1,\hat h)$ is firm and ultimate, this yields~$h_a=0$ (for all $a\ge a_0$) by Lemma~\ref{lem:8.8 firm} and Corollary~\ref{cor:existence and uniqueness},
and thus $f=g$ by Corollary~\ref{cor:h=0 => f=g}. 
\end{proof}

\noindent
Next we prove variants of Lemma~\ref{lem:notorious 3.6 firm} by modifying the restrictive hypothesis of strong split-normality.

\begin{lemma}\label{lem:notorious 3.6 firm, abs norm}
Suppose  $(P,\fm,\hat h)$ is firm, ultimate, and absolutely normal,
and its linear part~$L_{P_{\times\fm}}\in H[\der]\subseteq K[\der]$ is terminal. 
Then for all $f,g\in \Calr[\imag]$ such that~$P(f)=P(g)=0$ and $f,g\prec \fm$ we have $f=g$.
\end{lemma}
\begin{proof}
Replacing  $(P,\fm,\hat h)$ by $(P_{\times\fm},1,\hat h/\fm)$ we arrange $\fm=1$.
Put $A:=L_P\in H[\der]$ and $R:=P_1-P\in H\{Y\}$, so~$R\prec_\Delta \fv^{w+1}P_1$
 where $\Delta:=\Delta(\fv)$, $\fv:=\fv(A)\prec^\flat 1$.
Multiplying~$A$,~$P$,~$R$ on the left by some $b\in H^\times$ we arrange 
$$A\ =\ \der^r+ f_1\der^{r-1}+\cdots + f_r, \quad f_1,\dots, f_r\in H, \quad R\prec_\Delta\fv^w.$$ Then~\eqref{eq:notorious 3.6 firm} holds with~$\phi_1,\dots,\phi_r \in K_*$ instead of
$\phi_1,\dots,\phi_r \in K$, and $\phi_1,\dots,\phi_r\preceq\fv^{-1}$ by
Corollary~\ref{cor:bound on linear factors}.  Now  argue as at the end of the proof of Lemma~\ref{lem:notorious 3.6 firm} to get~$f=g$, using that $A\in K[\der]$ is terminal by assumption. 
\end{proof}

\noindent
From Corollary~\ref{cor:ultimate prod, 2} and Lemma~\ref{lem:notorious 3.6 firm, abs norm} we obtain:

\begin{cor}\label{cor:notorious 3.6 firm, abs norm}
If  $(P,\fm,\hat h)$ is firm, ultimate, and strictly normal, and its linear part splits strongly over $K$,
then the conclusion of Lemma~\ref{lem:notorious 3.6 firm, abs norm} holds.
\end{cor}

\noindent
{\it In the rest of this subsection we assume that all $A\in H[\der]\subseteq K[\der]$ of order~$r$ are terminal.}\/
Recall that  by  Lemmas~\ref{lem:ultimate refinement} and~\ref{lem:firm refinement},
if $(P,\fm,\hat h)$ is ultimate, then so is any refinement of it, and likewise
with ``firm'' in place of ``ultimate''.

\begin{lemma}\label{lem:notorious 3.6 firm, 1}
Suppose $(P,\fm,\hat h)$ is $Z$-minimal, deep, normal, special, ultimate, firm, and $f,g\in\c^r[\imag]$, $P(f)=P(g)=0$, $f \approx_K \hat h$, $g\approx_K \hat h$.
Then~$f=g$.
\end{lemma}
\begin{proof}
Corollary~\ref{cor:H*-normal}   gives a deep absolutely normal refinement $(P_{+h},\fn,\hat h-h)$ of~$(P,\fm,\hat h)$.
Now apply Lemma~\ref{lem:notorious 3.6 firm, abs norm} to $(P_{+h},\fn,\hat h-h)$, $f-h$, $g-h$ in the role of~$(P,\fm,\hat h)$,~$f$,~$g$.
\end{proof}

\begin{cor}\label{cor:notorious 3.6 firm, 1}
Suppose $H$ is  $\upo$-free and $r$-linearly newtonian, and $(P,\fm,\hat h)$ is  firm and $Z$-minimal. Then 
there is a unique~$f\in\c^r[\imag]$ with $P(f)=0$ and $f\approx_K\hat h$. For this $f$ we have~$f\in \Dx(H)$ and $f\approx_H \hat h$.
\end{cor}
\begin{proof}
Suppose $f,g\in\c^r[\imag]$
satisfy~$P(f)=P(g)=0$, $f \approx_K \hat h$, and~$g\approx_K \hat h$; we claim that $f=g$.
If $\phi>0$ is active in $H$, then by 
the remarks before Lemma~\ref{lem:firm refinement}, Lemma~\ref{lem:Pphicirc, 1}, and the remarks at the end of Section~\ref{sec:asymptotic similarity}, and with the superscript~$\circ$ having the usual meaning, we may  replace $H$, $K$, $\hat H$, $\hat K$, $(P,\fm,\hat h)$, $f$, $g$ by~$H^\circ$,~$K^\circ$,~$\hat H^\circ$,~$\hat K^\circ$,~$(P^{\phi\circ},\fm^\circ,\hat h^\circ)$,~$f^\circ$,~$g^\circ$. Using this observation, Corollary~\ref{cor:H*-normal, 2} and the remarks before Lemma~\ref{lem:notorious 3.6 firm, 1}, we arrange that  $(P,\fm,\hat h)$ is ultimate and absolutely normal. The claim now follows from
Lemma~\ref{lem:notorious 3.6 firm, abs norm}.
From Lemma~\ref{lem:find zero of P, d-max} we obtain~$f\in H_*$ with $P(f)=0$ and $f\approx_H\hat h$;
then $f\approx_K\hat h$ by Corollary~\ref{cor:approxH vs approxK}. Our $\d$-maximal Hardy field extension $H_*$ of $H$ was arbitrary, so the uniqueness statement just proved gives $f\in\Dx(H)$.
\end{proof}

\begin{cor}\label{cor:notorious 3.6 firm, 2} 
Suppose $H$ is $\upo$-free  and $(P,\fm,\hat h)$ is a firm minimal hole in~$H$. Then 
the conclusion of Corollary~\ref{cor:notorious 3.6 firm, 1} holds.
\end{cor}
\begin{proof} If $(P, \fm,\hat h)$ is nonlinear, then $H$ is $r$-linearly newtonian by Corollary~\ref{degmorethanone}, so
the hypotheses of Corollary~\ref{cor:notorious 3.6 firm, 1} are satisfied, and so is its conclusion. 

Suppose $(P, \fm,\hat h)$ is linear. By Remark~\ref{rem:achieve ultimate} we can refine $(P, \fm,\hat h)$ to arrange it to be ultimate. By our standing assumption $L_P$ is terminal, so we can appeal to Corollary~\ref{cor:8.8 firm}.
\end{proof}

%\noindent \marginpar{rest of subsection checked, but obsolete by above version} 
%If $H$ is $\upo$-free  and has a nonlinear minimal hole of order $r$, then $H$ is $r$-linearly newtonian by Corollary~\ref{degmorethanone}.
%Thus by  
%Corollary~\ref{cor:notorious 3.6 firm, 1}:

%\begin{cor}\label{cor:notorious 3.6 firm, 2} 
%If $H$ is $\upo$-free  and $(P,\fm,\hat h)$ is a firm nonlinear minimal hole in~$H$, then 
%the conclusion of Corollary~\ref{cor:notorious 3.6 firm, 1} holds.
%\end{cor}

\subsection*{$Z$-minimal slots in $\d$-perfect Hardy fields}
{\it In this subsection $H$ is $\d$-perfect.}\/
By Corollary~\ref{cor:d-perfect => 1-newt}, $H$ is $1$-newtonian and so
has no 
 quasilinear $Z$-minimal slot  of order $1$, by Corollary~\ref{cor:1-newt => no qlin Zmin dent order 1}.
This allows us to add to the characterization of $\upo$-freeness for $\d$-perfect Hardy fields given in Corollary~\ref{cor:perfect Schwarz closed, 2}:

\begin{cor}\label{cor:upo-free = no order 1 dents}
$$\text{$H$ is $\upo$-free}
\ \Longleftrightarrow \ \text{$H$ has no hole of order $1$}
\ \Longleftrightarrow \  \text{$H$ has no slot of order $1$.}$$ 
\end{cor}
\begin{proof}
The first equivalence follows from Lemma~\ref{lem:no hole of order <=r} and Corollary~\ref{cor:d-perfect => 1-newt}.
For the rest we observe that if $H$ has a slot of order $1$, then it also has a hole of order~$1$:
Given a slot $(P,\fm,\hat h)$ in $H$ of order~$1$, 
take  $Q\in Z(H,\hat h)$ of minimal complexity. Then~$(Q,\fm,\hat h)$ is a $Z$-minimal slot of order~$\leq 1$ in $H$, 
hence is equivalent to a $Z$-minimal hole $(Q,\fm,\hat b)$ in $H$, by
Lemma~\ref{lem:from cracks to holes}, so  $\order Q=1$ by a remark after Lemma~\ref{lem:no hole of order <=r}.
\end{proof}

\noindent
Next, an immediate consequence of Corollaries~\ref{cor:linear flabby} and~\ref{cor:notorious 3.6 firm, 2}:

\begin{cor}\label{cor:d-perfect => min holes flabby, 1}
Let $r\in\N^{\ge 1}$. If $H$ is   $\upo$-free and  all $A\in H[\der]$ of order~$r$ are terminal, then every
minimal hole in $H$ of order $r$ is flabby.
\end{cor}

\noindent
From this we deduce:

\begin{cor}\label{cor:d-perfect => min holes flabby, 2}
Suppose $H$ is $\upo$-free. Then every  minimal hole in $H$ of order~$2$
and  every linear  slot in $H$ of order~$2$ is flabby.
\end{cor}
\begin{proof}
By Corollaries~\ref{cor:ultimate prod, 2} and~\ref{cor:perfect Schwarz closed, 2}, all $A\in H[\der]$ of order~$2$
are terminal. Hence every  minimal hole in $H$ of order~$2$ is flabby, by Corollary~\ref{cor:d-perfect => min holes flabby, 1}.
Every linear slot in $H$ of order~$2$ is $Z$-minimal, by
Corollary~\ref{cor:upo-free = no order 1 dents}, and hence is flabby by
Corollary~\ref{cor:linear flabby}.
\end{proof}

\noindent
In the next subsection we study flabby slots in~$H$ in more detail.

\subsection*{Flabby slots in $H$}
{\it In this subsection $(P,\fm,\hat h)$ is a  slot in~$H$ of order~${r\geq 1}$.}\/
%Set $d:=\deg(P)$, $w:=\wt(P)$. 
Note that if $(P,\fm,\hat f)$ is normal and $f\in \fm\c^r[\imag]^{\preceq}$, $P(f)=0$, then by Corollary~\ref{cor:hardynormal} we have $f\in \Calinf[\imag]$, and $f\in\Ginf[\imag]$ if $H\subseteq\Ginf$,  $f\in \Gom[\imag]$ if $H\subseteq \Gom$.

Next some observations tacitly used in the proof of Lemma~\ref{lem:notorious 3.6 flabby}  below.
For this, suppose $(P,\fm,\hat h)$ is flabby.
Then~$(P_{\times\fm},1,\hat h/\fm)$ is flabby by Lemma~\ref{lem:firm mult conj}, and
if $g\in(\c^r)^\preceq$ and $P_{\times\fm}(g)=0$, $g\prec 1$, 
then $f:=\fm g\in \fm\, (\c^r)^{\preceq}$ satisfies~$P(f)=0$ and~$f\prec\fm$. Likewise, let~${(P_{+h},\fn,\hat h-h)}$ be a refinement of~$(P,\fm,\hat h)$, and suppose~$(P,\fm,\hat h)$ is also linear or normal.
Then the slot~$(P_{+h},\fn,\hat h-h)$ in~$H$ is  flabby by Corollary~\ref{cor:flabby refinement},
and if  $g\in\fn\,(\c^r)^\preceq$ and $P_{+h}(g)=0$, $g\prec \fn$,   then~$f:=h+g\in\fm\,(\c^r)^\preceq$ satisfies
$P(f)=0$ and~$f\prec\fm$.
Finally, let~$\phi$ be active in~$H$, $0<\phi\preceq 1$, and let the superscript~$\circ$ have the usual meaning.
% as explained in  Section~\ref{secfhhf}.
Then the   slot~$(P^{\phi\circ},\fm^\circ,\hat h^\circ)$ in~$H^\circ$ is flabby, and
if $g\in \fm^\circ\,(\c^r)^\preceq$ and $P^{\phi\circ}(g)=0$, $g\prec\fm^\circ$, 
then taking~$f\in\c^r$ with~$f^\circ=g$   we have  $f\in\fm\,(\c^r)^{\preceq}$, $P(f)=0$, and~$f\prec\fm$, using
the remark after Lemma~\ref{lem:as equ comp}.

\begin{lemma}\label{lem:notorious 3.6 flabby}
Suppose $(P,\fm,\hat h)$ is 
$Z$-minimal, deep, normal, special, and flabby.  Then there are~$f\neq g$ in $\fm\, (\c^r)^\preceq$ such that
$P(f)=P(g)=0$, and $f,g\prec \fm$. 
%any such~$f$,~$g$ are in $\Calinf$, and
%if $H\subseteq\Ginf$, then~$f, g\in\Ginf$, and likewise  with $\Gom$ in place of~$\Gin
\end{lemma}

\begin{proof}
Corollary~\ref{cor:H*-normal}   
yields a deep absolutely normal refinement $(P_{+h},\fn,\hat h-h)$ of~$(P,\fm,\hat h)$.
Using the remarks preceding the lemma, we can  replace~$(P,\fm,\hat h)$ by~$(P_{+h},\fn,\hat h-h)$ to arrange that
$(P,\fm,\hat h)$ is absolutely normal, and then
replacing~$(P,\fm,\hat h)$ by $(P_{\times\fm},1,\hat h/\fm)$ we also arrange $\fm=1$. Set $d:=\deg P$ and~$w:=\wt(P)$.
Let~$\fv$,~$\Delta$,~$A$,~$R$ be as in the proof of Lemma~\ref{lem:notorious 3.6 firm, abs norm}.
Multiplying $A$, $P$, $R$ on the left by some $b\in H^\times$ we arrange $A=\der^r+ f_1\der^{r-1}+\cdots + f_r$ with $f_1,\dots, f_r\in H$ and~$R\prec_\Delta\fv^w$. Then~\eqref{eq:notorious 3.6 firm} holds with $\phi_1,\dots,\phi_r\in K_*$ instead of $\phi_1,\dots,\phi_r\in K$, and $\phi_1,\dots,\phi_r\preceq\fv^{-1}$ by
Corollary~\ref{cor:bound on linear factors}. 
 Take $a_0\in \R$ and functions on $[a_0,\infty)$ representing the germs $\phi_1,\dots, \phi_r$, $f_1,\dots, f_r$, and the $R_{\j}$ with $\j\in \N^{1+r}$, $|\j|\le d$, $\|\j\|\le w$ (using the same symbols for the germs mentioned as for their chosen representatives)
so as to be in the situation described in the beginning of Section~\ref{sec:split-normal over Hardy fields}.
Increasing $a_0$ if necessary and choosing $\nu$ as in the proof of Lemma~\ref{lem:notorious 3.6 firm} we arrange that~$f_1,\dots, f_r,$ and the~$R_{\j}$ are in  $\c^{1}_{a_0}$ and $\dabs{R}_{a_0}\leq 1/E$, with~$E=E(d,r)\in\N^{\geq 1}$ as in Corollary~\ref{cor:uniqueness ADE}, and
the hypotheses of Lemma~\ref{bdua} are satisfied.  
Lemma~\ref{lem:8.8 flabbyLP} yields $h\in\Calinf[\imag]$
such that $A(h)=0$, $h\neq 0$, and~$h,h',\dots,h^{(r)}\prec 1$. Replacing~$h$ by $\Re h$ or $\Im h$  we arrange $h\in\Calinf$. 
 Increasing $a_0$ again we arrange that $h$ is represented by a function in $\c^{r}_{a_0}$, denoted by the same symbol,
such that $$A_{a_0}(h)\ =\ 0,\qquad \dabs{h}_{a_0;r}\ \le 1/8,$$ and such that we are in the situation of
Lemma~\ref{realbdua}, with $a$ ranging over $[a_0, +\infty)$. Then
Corollaries~\ref{cor:fix} and~\ref{cor:fix h} 
yield for sufficiently large $a\geq a_0$ functions $f,g\in\Car$ with~$\dabs{f}_{a;r},\dabs{g}_{a;r}\leq 1$ and $(\Re\Xi_a)(f)=f$, $(\Re\Xi_a)(g)=g+h$.   Fix such~$a$,~$f$,~$g$.  Then $A_a(f)=R(f)$ and
$$A_a(g)\ =\ A_a(g+h)\  =\  A_a\big((\Re\Xi_a)(g)\big)\  =\ \Re A_a\big(\Xi_a(g)\big)\ =\ \Re R(g)\ =\  R(g),$$
 and $f\prec 1$ and $g+h\prec 1$ by Lemma~\ref{realbdua}, so $g\prec 1$, hence
 $f$, $g$ are solutions of~\eqref{eq:ADE} on $[a,\infty)$.  
Denoting the germs of~$f$,~$g$ also by~$f$,~$g$ we have  $P(f)=P(g)=0$.  
Moreover, $f\neq g$ as germs: otherwise~$f=g$ in~$\Car$ by the remark after the proof of Corollary~\ref{cor:uniqueness ADE}, and thus 
$h=(\Re\Xi_a)(g)-g=(\Re\Xi_a)(f)-f=0$ in~$\Car$, a contradiction.  
\end{proof}

\begin{cor}\label{cor:notorious 3.6 flabby, 1}
Suppose $H$ is $\upo$-free and $r$-linearly newtonian, and $(P,\fm,\hat h)$ is $Z$-minimal and flabby. Assume also that $(P,\fm, \hat h)$ is linear or normal.
Then the conclusion of Lemma~\ref{lem:notorious 3.6 flabby} holds.
\end{cor}
\begin{proof}
Use Theorem~\ref{mainthm} and the remarks preceding  Lem\-ma~\ref{lem:notorious 3.6 flabby}
to  arrange that $(P,\fm,\hat h)$  is  deep and normal. Then $(P,\fm,\hat h)$ is quasilinear 
by Corollary~\ref{cor:normal=>quasilinear},  and hence special by Lemma~\ref{lem:special dents}.
Now  Lem\-ma~\ref{lem:notorious 3.6 flabby} applies.
\end{proof}

\noindent
Note that the hypotheses of Corollary~\ref{cor:notorious 3.6 flabby, 1} hold if 
$H$ is $\upo$-free and~$(P,\fm,\hat h)$  is a  flabby normal nonlinear minimal hole in~$H$, by
Corollary~\ref{degmorethanone}.  

\medskip
\noindent
Suppose  $H$ is $\upo$-free, all $A\in H[\der]\subseteq K[\der]$ of order $r$ are terminal, and~$(P,\fm,\hat h)$ is a minimal hole in $H$.
Then by Corollary~\ref{cor:notorious 3.6 firm, 2} we have:
$$\text{$(P,\fm,\hat h)$ is firm}\quad\Longrightarrow\quad \text{there is a unique $f\in\c^r$ with $P(f)=0$ and $f\approx_H \hat h$.}$$
Thanks to Corollary~\ref{cor:8.8 flabby}, the converse of this implication also holds if $\deg P=1$,  but
we do not know whether this is still the case when~$\deg P>1$. We now prove a partial generalization of   Corollary~\ref{cor:8.8 firm}: 

\begin{cor}\label{cor:8.8 firm, gen}
Suppose $H$ is $\upo$-free, $(P,\fm,\hat h)$ is an ultimate minimal hole in~$H$ with terminal linear part,
and $(P,\fm,\hat h)$ is linear or absolutely normal.
Then
$$\text{$(P,\fm,\hat h)$ is firm}\quad\Longleftrightarrow\quad \text{there is a unique $f\in\c^r$ with $P(f)=0$ and $f\prec\fm$.}$$
If $(P,\fm,\hat h)$ is firm and $f\in \c^r$, $P(f)=0$, and $f\prec \fm$, then 
$f\in\Dx(H)$ and $f\approx_H\hat h$, and there is no $g\ne f$ in $\c^r[\imag]$ with $P(g)=0$ and $g\prec\fm$.
\end{cor}
\begin{proof}
If  $\deg P=1$, then this follows from  Corollary~\ref{cor:8.8 firm}. Suppose $\deg P>1$.
Lemma~\ref{lem:find zero of P, d-max} yields $f\in H_*$ with $P(f)=0$ and $f\approx_H\hat h$.  This and Lemma~\ref{lem:notorious 3.6 firm, abs norm} yield
the forward direction of the displayed equivalence, as well as the rest in view of $H_*$ being arbitrary.
The converse holds by the remark after Corollary~\ref{cor:notorious 3.6 flabby, 1}. 
\end{proof}

\begin{remark}
Suppose $H$ is $\upo$-free, all $A\in H[\der]\subseteq K[\der]$ of order $r$ are terminal, and~$(P,\fm,\hat h)$ is a  minimal hole in $H$. 
  Then Remarks~\ref{rem:achieve ultimate} and~\ref{rem:H*-normal, 2} give  a refinement 
$(P_{+h},\fn,\hat h-h)$ of~$(P,\fm,\hat h)$  and an active $\phi>0$ in~$H$ such that
the minimal hole $(P_{+h^\circ}^{\phi\circ},\fn^\circ,\hat h^\circ-h^\circ)$ in~$H^\circ$ is ultimate, and linear or absolutely normal.
Therefore the hypotheses of Corollary~\ref{cor:8.8 firm, gen} are satisfied by~$H^\circ$ and
$(P_{+h^\circ}^{\phi\circ},\fn^\circ,\hat h^\circ-h^\circ)$  in place of $H$ and $(P,\fm,\hat h)$.
\end{remark}
  
\noindent
The next two subsections contain analogues of Lemmas~\ref{lem:notorious 3.6 firm, abs norm}, \ref{lem:notorious 3.6 firm, 1}, and~\ref{lem:notorious 3.6 flabby} for slots in~$K$.  

\subsection*{Firm  slots in $K$}
{\it In this subsection $(P,\fm,\hat f)$ is a  slot in~$K$ of order $r\geq 1$, with~$\hat f\in \hat K$.}\/
Here are $K$-versions of Lemma~\ref{lem:notorious 3.6 firm, abs norm} and its Corollary~\ref{cor:notorious 3.6 firm, abs norm} with similar proofs:

\begin{lemma}\label{lem:notorious 3.6 firm, K}
If $(P,\fm,\hat f)$ is firm, ultimate, and absolutely normal, with terminal  linear part, then for any  $f,g\in \c^r[\imag]$ with $P(f)=P(g)=0$, $f,g\prec\fm$ we have $f=g$.
\end{lemma}

\begin{cor}\label{cor:notorious 3.6 firm, K} 
If $(P,\fm,\hat f)$ is firm, ultimate, and strictly normal, and its linear part splits strongly over $K$,
then the conclusion of Lemma~\ref{lem:notorious 3.6 firm, K} holds.
\end{cor}

\noindent
{\it In the rest of this subsection we assume that all $A\in K[\der]$ of order~$r$ are terminal.}\/
(This holds if $r=1$ or $K$ is $\upo$-free with a minimal hole of order $r$ in $K$, because then $K$ is $r$-linearly closed by Corollary~\ref{corminholenewt}.)

\medskip
\noindent
By the $K$-versions of Lemmas~\ref{lem:ultimate refinement} and~\ref{lem:firm refinement},
if $(P,\fm,\hat f)$ is ultimate, then so is each of its refinements, and likewise
with ``firm'' in place of ``ultimate''.
The $K$-version of Corollary~\ref{cor:H*-normal}, and 
Lemma~\ref{lem:notorious 3.6 firm, K} in place of Lemma~\ref{lem:notorious 3.6 firm, abs norm} then yields the $K$-version of Lemma~\ref{lem:notorious 3.6 firm, 1}:

\begin{lemma}\label{lem:notorious 3.6 firm, 1, K}
If $(P,\fm,\hat f)$ is $Z$-minimal, deep, normal, special, ultimate, and firm, then there is at most one~$f\in\c^r[\imag]$
with~$P(f)=0$ and $f \approx_K \hat f$.
\end{lemma}

\noindent
Using the $K$-version of Corollary~\ref{cor:H*-normal, 2},
and  Lemmas~\ref{lem:find zero of P, d-max, K} and~\ref{lem:notorious 3.6 firm, K} instead of Lemmas~\ref{lem:find zero of P, d-max} and~\ref{lem:notorious 3.6 firm, abs norm} we obtain the $K$-version of Corollary~\ref{cor:notorious 3.6 firm, 1}:

\begin{cor}\label{cor:notorious 3.6 firm, 1, K}
If $H$ is  $\upo$-free, $K$ is $r$-linearly newtonian, and $(P,\fm,\hat f)$ is  firm and $Z$-minimal, then 
there is a unique~$f\in\c^r[\imag]$ with $P(f)=0$ and $f\approx_K\hat f$, and this $f$ is in~$\Dx(H)[\imag]$.
\end{cor}

\noindent
Here is a $K$-analogue of Corollary~\ref{cor:notorious 3.6 firm, 2}:

\begin{cor}\label{cor:notorious 3.6 firm, 2, K} 
Suppose $(P,\fm,\hat f)$ is a firm minimal hole in $K$, and~$r=\deg P=1$ or  $H$ is $\upo$-free. 
Then the conclusion of Corollary~\ref{cor:notorious 3.6 firm, 1, K} holds.
\end{cor}  
\begin{proof}
If $(P,\fm,\hat f)$  has complexity~$(1,1,1)$, then by Remark~\ref{rem:achieve ultimate, K} and the $K$-version of Lem\-ma~\ref{lem:firm refinement} we arrange that
$(P,\fm,\hat f)$ is ultimate, so that the desired conclusion follows from Corollary~\ref{cor:8.8 firm, K, r=1}.
If $K$ is $\upo$-free and  $(P,\fm,\hat f)$  has complexity~$>(1,1,1)$,
then $\deg P>1$ by Corollary~\ref{cor:minhole deg 1}, so $K$ is
  $r$-linearly newtonian by Corollary~\ref{degmorethanone}, and the desired conclusion follows from Corollary~\ref{cor:notorious 3.6 firm, 1, K}.
\end{proof}

\begin{cor}\label{cor:notorious 3.6 firm, 3, K} 
Suppose $H$ is $\d$-perfect. 
If $(P,\fm,\hat f)$ is  a minimal hole   in $K$ of complexity $(1,1,1)$, then it is flabby. If
$H$ is $\upo$-free, then every minimal hole in $K$ of positive order is flabby.
\end{cor}
\begin{proof}
The first part is immediate from  Corollary~\ref{cor:notorious 3.6 firm, 2, K}.
For the second part, suppose $H$ is $\upo$-free and we are given a minimal hole in $K$ of positive order.
By Lemma~\ref{lem:hole in hat K} we can  pass to an equivalent hole~$(Q,\fn,\tilde{a})$ in $K$ with $\tilde{a}\in\tilde{H}[\imag]$
for some immediate $H$-field extension $\tilde{H}$ of $H$, so Corollary~\ref{cor:notorious 3.6 firm, 2, K} applies to it.
\end{proof}

\noindent
Theorem~\ref{thm:d-perfect flabby} now follows from Corollaries~\ref{cor:upo-free = no order 1 dents}, \ref{cor:d-perfect => min holes flabby, 1}, \ref{cor:d-perfect => min holes flabby, 2},
and~\ref{cor:notorious 3.6 firm, 3, K}.

\subsection*{Flabby  slots in $K$}
{\it Let $(P,\fm,\hat f)$ be a slot in~$K$ of order $r\geq 1$, $\hat f\in \hat K$.}\/
Note that if $(P,\fm,\hat f)$ is normal and $f\in \fm\c^r[\imag]^{\preceq}$, $P(f)=0$, then by Corollary~\ref{cor:hardynormal} we have~$f\in \Calinf[\imag]$, and $f\in\Ginf[\imag]$ if $H\subseteq\Ginf$,  $f\in \Gom[\imag]$ if $H\subseteq \Gom$.

Suppose $(P,\fm,\hat f)$ is flabby.
The remarks about multiplicative conjugates, refinements, and compositional conjugates 
preceding Lemma~\ref{lem:notorious 3.6 flabby} then go through for the slot $(P,\fm,\hat f)$ in $K$
instead of the slot $(P,\fm,\hat h)$ in $H$ with $\c^r[\imag]$ replacing~$\c^r$ and~$K^\circ$ instead of $H^\circ$; this
uses the $K$-versions of Lemma~\ref{lem:firm mult conj} and Corollary~\ref{cor:flabby refinement}.
It helps in proving a complex version of Lem\-ma~\ref{lem:notorious 3.6 flabby}:

\begin{lemma}\label{lem:notorious 3.6 flabby, K}
Suppose $(P,\fm,\hat f)$ is flabby, special, $Z$-minimal,  deep, and strictly normal. 
Then there are~$f\neq g$ in $\fm\, \c^r[\imag]^\preceq$ such that~$P(f)=P(g)=0$, $f,g\prec \fm$.  
\end{lemma}
\begin{proof}
%Let  $\fv:=\abs{\fv(L_{P_{\times\fm}})}\prec^\flat 1$, $\Delta:=\Delta(\fv)$, $d:=\deg(P)$, $w:=\wt(P)$.
Using Corollary~\ref{cor:K*-normal} and the
remarks preceding the lemma, we arrange
that $\fm=1$ and  $(P,1,\hat f)$ is absolutely normal.
Now argue as in
the proof of Lemma~\ref{lem:notorious 3.6 flabby}, using instead of Lemma~\ref{lem:8.8 flabbyLP} its $K$-version.
We also appeal to
Lemma~\ref{bdua}, Theorem~\ref{thm:fix}, and Lemma~\ref{lem:fix h} instead of to Lemma~\ref{realbdua} and Corollaries~\ref{cor:fix} and~\ref{cor:fix h}. Naturally, we don't need to take real or imaginary parts, and use $\Xi_a$ instead of $\Re\Xi_a$.
\end{proof}

\begin{cor}\label{cor:notorious 3.6 flabby, 1, K}
Suppose $H$ is $\upo$-free, $K$ is $r$-linearly newtonian, and $(P,\fm,\hat f)$ is $Z$-minimal and flabby. Assume also that $(P,\fm, \hat f)$ is linear or normal.
Then the conclusion of Lemma~\ref{lem:notorious 3.6 flabby, K} holds.
\end{cor}
\begin{proof} Like that of Corollary~\ref{cor:notorious 3.6 flabby, 1}, but using Corollary~\ref{cor:achieve strong normality, 2} instead of Theorem~\ref{mainthm}, and  Lemma~\ref{lem:notorious 3.6 flabby, K} instead of Lemma~\ref{lem:notorious 3.6 flabby}. 
\end{proof}

\noindent
In particular,  the conclusion of Lemma~\ref{lem:notorious 3.6 flabby, K} holds 
if $H$ is $\upo$-free and~$(P,\fm,\hat f)$  is a  flabby normal nonlinear minimal hole in~$K$.  

\begin{cor}\label{cor:8.8 firm, gen, K}
Suppose that
$(P,\fm,\hat f)$ is an ultimate minimal hole in~$K$ and, in case the complexity of $(P,\fm,\hat f)$ is $>(1,1,1)$, that
$H$ is $\upo$-free and $(P,\fm,\hat f)$ is absolutely normal. Then
$$\text{$(P,\fm,\hat f)$ is firm}\quad\Longleftrightarrow\quad \text{there is a unique $f\in\c^r[\imag]$ with $P(f)=0$ and $f\prec\fm$.}$$
If $(P, \fm, \hat f)$ is firm, $f\in\c^r[\imag]$, $P(f)=0$, and $f\prec\fm$, then
$f\in\Dx(H)[\imag]$ and $f\approx_K\hat f$.
\end{cor}
\begin{proof}
If $(P,\fm,\hat f)$ has complexity $(1,1,1)$, use Corollaries~\ref{cor:8.8 firm, K} and~\ref{cor:8.8 firm, K, r=1}.
Now suppose
$H$ is $\upo$-free and $(P,\fm,\hat f)$ is absolutely normal of complexity $>(1,1,1)$.
Then $\deg P>1$ by Corollary~\ref{cor:minhole deg 1}, and $L_{P_{\times\fm}}$ is terminal by Corollaries~\ref{corminholenewt}
and~\ref{cor:ultimate prod, 2}.
Thus the forward direction of the displayed equivalence  follows from  
Lemmas~\ref{lem:notorious 3.6 firm, K} and~\ref{lem:find zero of P, d-max, K}, and the backward direction
from the remark after Corollary~\ref{cor:notorious 3.6 flabby, 1, K}. 
The rest follows by applying Lem\-ma~\ref{lem:find zero of P, d-max, K} to all choices of $H_*$. 
\end{proof}

 \begin{remark}
 Suppose $(P,\fm,\hat f)$ is a minimal hole in $K$.
 If $\deg P=1$, then $(P,\fm,\hat f)$ refines to an ultimate hole in $K$ by Remark~\ref{rem:achieve ultimate, K}.
 If $H$ is $\upo$-free and $\deg P>1$,
  then Corollary~\ref{cor:K*-normal, 2}
gives  a refinement 
$(P_{+f},\fn,\hat f-f)$ of~$(P,\fm,\hat f)$  and an active~$\phi>0$ in~$H$ such that
the minimal hole $(P_{+f^\circ}^{\phi\circ},\fn^\circ,\hat f^\circ-f^\circ)$ in~$K^\circ$ is ultimate and absolutely normal. 
\end{remark}

%\bigskip\noindent
%{\bf Notes checked except as indicated in marginal or bold face comments} 

\printindex

\part*{List of Symbols}

\medskip

{\small
\begin{longtable}{m{6.5em} m{30em} m{4em}}

%$A^*$   &  adjoint of the linear differential operator  $A$ \dotfill & \pageref{p:A*} \\[0.4em]

$\mult_a(A)$   & multiplicity of     $A$ at $a\in K$  \dotfill & \pageref{p:multa} \\[0.4em]

$\mult_\alpha(A)$   & multiplicity of    $A$ at~${\alpha\in K/K^\dagger}$  \dotfill & \pageref{p:multalpha} \\[0.4em]

$\Sigma(A)$   &  spectrum of   $A$ \dotfill & \pageref{p:SigmaA} \\[0.4em]

$\exc^{\operatorname{u}}(A)$   &  set of ultimate exceptional values of   $A$ \dotfill & \pageref{p:excu} \\[0.4em]

$\fv(A)$   &  span of  $A$ \dotfill & \pageref{p:span} \\[0.4em]

$\mult_a(M)$   & multiplicity of the   differential module   $M$ at $a\in K$  \dotfill & \pageref{p:multaM} \\[0.4em]

$\mult_\alpha(M)$   & multiplicity of    $M$ at $\alpha\in K/K^\dagger$  \dotfill &  \pageref{p:multalphaM} \\[0.4em]

%$M^*$   &  adjoint of the   differential module $M$ \dotfill & \pageref{p:M*} \\[0.4em]

$\Sigma(M)$   &  spectrum of the   differential module  $M$   \dotfill & \pageref{p:SigmaM} \\[0.4em]

%$S_P$   &  separant of the differential polynomial $P$   \dotfill &  \pageref{p:separant} \\[0.4em]

%$\cc(P)$   &  complexity of  $P$   \dotfill &  \pageref{p:complexity} \\[0.4em]

%$\Ric(P)$   &  Riccati transform of   $P$   \dotfill &  \pageref{p:Ric} \\[0.4em]

$K^\dagger$   &  group of logarithmic derivatives of $K$ \dotfill &  \pageref{p:Kdagger} \\[0.4em]

$\Univ_K$   & universal exponential extension of    $K$  \dotfill &  \pageref{p:UK} \\[0.4em]

$v_{\g}$   & gaussian extension of the valuation of $K$ to $K[G]$  \dotfill &  \pageref{p:vg} \\[0.4em]

$\preceq_{\g}$, $\prec_{\g}$, $\asymp_{\g}$   & dominance relations associated to $v_{\g}$  \dotfill &  \pageref{p:vg} \\[0.4em]

$H^{\trig}$, $H^{\tl}$   &  trigonometric closure of  $H$, trigonometric-Liouville closure of $H$ \dotfill &  \pageref{p:Htrig}, \pageref{p:Htl} \\[0.4em]

 $\Ex(H)$ & perfect hull of $H$  \dotfill &  \pageref{p:E(H)} \\[0.4em]

 $\Dx(H)$ & $\d$-perfect hull of $H$  \dotfill &  \pageref{p:D(H)} \\[0.4em]

 $\Ex^r(H)$ & $\c^r$-perfect hull of $H$  \dotfill &  \pageref{p:Er(H)} \\[0.4em]

$\Li(H)$   &  Hardy-Liouville closure of  $H$  \dotfill &  \pageref{p:HL}  \\[0.4em]

$\bar{\omega}(H)$   &  set of $f\in H$ such that $f/4$ does not generate oscillations  \dotfill &  \pageref{p:baromega} \\[0.4em]

$\sim_H$, $\approx_H$ & asymptotic similarity over $H$ \dotfill &  \pageref{p:simH} \\[0.4em]

$\sim_K$, $\approx_K$ & asymptotic similarity over $K$ \dotfill &  \pageref{p:simK} \\[0.4em]

%$\Upg(H)$, $\Upl(H)$, $\Upd(H)$   & special definable subsets of the   pre-$H$-field $H$  \dotfill &  \pageref{p:special} \\[0.4em]

%$\I(K)$   & special definable $\mathcal O$-submodule of the asymptotic field~$K$  \dotfill &  \pageref{p:I(K)} \\[0.4em]

$(P,\fm,\hat a)$   &  slot in $K$  \dotfill &  \pageref{p:hole}, \pageref{p:slot} \\[0.4em]

$\Delta(\fm)$   & convex subgroup of all $\gamma\in\Gamma$ with $\gamma=o(v\fm)$  \dotfill &  \pageref{p:Delta} \\[0.4em]

$\c^r(U)$ & ring of $\c^r$  functions~$U\to \R$   \dotfill &  \pageref{p:Cr(U)} \\[0.4em]

$\Car$ & ring of    $f|_{\R^{\geq a}}$ where~$f\in\c^r(U)$ for an open $U\supseteq\R^{\geq a}$ \dotfill &  \pageref{p:Car} \\[0.4em]

$\c_a[\imag]^{\inte}$ & ring of  $f\in\c_a[\imag]$ integrable at $\infty$  \dotfill &  \pageref{p:Caint} \\[0.4em]

 ${\|\cdot\|_a}$,   ${\|\cdot\|_{a;r}}$ & supremum norms   \dotfill &  \pageref{p:absa} \\[0.4em]

$\c_a[\imag]^{\b}$, $\c^r_a[\imag]^{\b}$ & rings of functions of bounded norms  \dotfill &  \pageref{p:Cab} \\[0.4em]

 ${\|\cdot\|_a^\w}$,   ${\|\cdot\|_{a;r}^\w}$ & weighted supremum norms   \dotfill &  \pageref{p:absawt} \\[0.4em]

$\c_a[\imag]^{\w}$, $\c^r_a[\imag]^{\w}$ & rings of functions of bounded weighted norm   \dotfill &  \pageref{p:absawt} \\[0.4em]

$\c$, $\c[\imag]$   & rings of continuous germs \dotfill &  \pageref{p:cont} \\[0.4em]

$\c^{\preceq}$, $\c[\imag]^{\preceq}$ & rings of bounded continuous germs \dotfill &  \pageref{p:contpreceq} \\[0.4em]

 $\c^r$, $\c^r[\imag]$ & rings of $\c^r$ germs   \dotfill &  \pageref{p:Gr} \\[0.4em]

 $\Calinf$, $\Calinf[\imag]$ & intersection of
all   $\c^r$  respectively $\c^r[\imag]$ ($r\in\N$)   \dotfill &  \pageref{p:Gi} \\[0.4em]

$g^{\operatorname{inv}}$ & compositional inverse of $g\in\c$ \dotfill &  \pageref{p:ginv} \\[0.4em]

$f^{\circ}$ & the germ $f\circ\ell^{\operatorname{inv}}$ \dotfill &  \pageref{p:fcirc} \\[0.4em]

%$\mult_t(y)$ & multiplicity of $y\in\Car[\imag]$ at $t\geq a$   \dotfill &  \pageref{p:multt} \\[0.4em]

%$\mult(y)$ & total multiplicity of $y\in\Car[\imag]$    \dotfill &  \pageref{p:multtotal} \\[0.4em]

$\Sol(f)$ & solution space of $Y''+fY=0$  \dotfill &  \pageref{p:Sol} \\[0.4em]

$\sol_{\Univ}(N)$ & solution space of $y'=Ny$ \dotfill & \pageref{p:solU}
\end{longtable}}

\newpage

\newlength\templinewidth
\setlength{\templinewidth}{\textwidth}
\addtolength{\templinewidth}{-2.25em}

\patchcmd{\thebibliography}{\list}{\printremarkbeforebib\list}{}{}

\let\oldaddcontentsline\addcontentsline% Store \addcontentsline
\renewcommand{\addcontentsline}[3]{\oldaddcontentsline{toc}{part}{References}}

\def\printremarkbeforebib{\bigskip\hskip1em The citation [ADH] refers to our book \\

\hskip1em\parbox{\templinewidth}{
M. Aschenbrenner, L. van den Dries, J. van der Hoeven,
\textit{Asymptotic Differential Algebra and Model Theory of Transseries,} Annals of Mathematics Studies, vol.~195, Princeton University Press, Princeton, NJ, 2017.
}

\bigskip

} 

\bibliographystyle{amsplain}

\end{document}